\documentclass{amsbook}

\newtheorem{thm}{Theorem}[section]
\newtheorem*{thm*}{Theorem}

\newtheorem{claim}[thm]{Claim}

\newtheorem{cor}[thm]{Corollary}

\newtheorem{lem}[thm]{Lemma}
\newtheorem*{lem*}{Lemma}
\newtheorem{mainthm}{Theorem}
\newtheorem*{mainthm*}{Theorem}
\newtheorem{maincor}[mainthm]{Corollary}

\newtheorem{prop}[thm]{Proposition}

\theoremstyle{definition}

\newtheorem*{case*}{Case}
\newtheorem{conj}[thm]{Conjecture}

\newtheorem{defn}[thm]{Definition}
\newtheorem*{defn*}{Definition}
\newtheorem{exmp}[thm]{Example}
\newtheorem*{exmp*}{Example}

\newtheorem{hyp}[thm]{Hypothesis}

\newtheorem{maindefn}[mainthm]{Definition}

\newtheorem{step}{Step}\renewcommand{\thestep}{}

\theoremstyle{remark}

\newtheorem{case}{Case}\renewcommand{\thecase}{}
\newtheorem{subcase}{Case}
\numberwithin{subcase}{case}
\newtheorem{subsubcase}{Case}
\numberwithin{subsubcase}{subcase}

\newtheorem{rmk}[thm]{Remark}
\newtheorem*{rmk*}{Remark}

\makeatletter
\def\alphenumi{
  \def\theenumi{\alph{enumi}}
  \def\p@enumi{\theenumi}
  \def\labelenumi{(\@alph\c@enumi)}}
\makeatother

\makeatletter
\def\thecase{\@arabic\c@case}
\makeatother

\makeatletter
\def\thestep{\@arabic\c@step}
\makeatother

\makeatletter
\newcommand{\tpitchfork}{%
  \vbox{
    \baselineskip\z@skip
    \lineskip-.52ex
    \lineskiplimit\maxdimen
    \m@th
    \ialign{##\crcr\hidewidth\smash{$-$}\hidewidth\crcr$\pitchfork$\crcr}
  }%
}
\makeatother

\DeclareFontFamily{U}{mathx}{\hyphenchar\font45}
\DeclareFontShape{U}{mathx}{m}{n}{
      <5> <6> <7> <8> <9> <10>
      <10.95> <12> <14.4> <17.28> <20.74> <24.88>
      mathx10
      }{}
\DeclareSymbolFont{mathx}{U}{mathx}{m}{n}
\DeclareFontSubstitution{U}{mathx}{m}{n}
\DeclareMathAccent{\widecheck}{0}{mathx}{"71}
\DeclareMathAccent{\wideparen}{0}{mathx}{"75}

\newcount\hh
\newcount\mm
\mm=\time
\hh=\time
\divide\hh by 60
\divide\mm by 60
\multiply\mm by 60
\mm=-\mm
\advance\mm by \time
\def\hhmm{\number\hh:\ifnum\mm<10{}0\fi\number\mm}

\let\oldmarginpar\marginpar
\renewcommand\marginpar[1]{\-\oldmarginpar[\raggedleft\footnotesize #1]%
{\raggedright\footnotesize #1}}

\newcommand\embed{\hookrightarrow}

\newcommand\barM{{\bar{M}}}

\newcommand\ubarCC{{\underline{\mathbb{C}}}}

\newcommand\ubarRR{{\underline{\mathbb{R}}}}

\newcommand\CC{\mathbb{C}}

\newcommand\KK{\mathbb{K}}

\newcommand\NN{\mathbb{N}}

\newcommand\PP{\mathbb{P}}

\newcommand\RR{\mathbb{R}}
\newcommand\SSS{\mathbb{S}}

\newcommand\UU{\mathbb{U}}

\newcommand\ZZ{\mathbb{Z}}

\newcommand\cF{{\mathcal{F}}}

\newcommand\cH{{\mathcal{H}}}

\newcommand\cK{{\mathcal{K}}}

\newcommand\cM{{\mathcal{M}}}

\newcommand\cP{{\mathcal{P}}}

\newcommand\cS{{\mathcal{S}}}

\newcommand\cU{{\mathcal{U}}}
\newcommand\cV{{\mathcal{V}}}

\newcommand\cX{{\mathcal{X}}}
\newcommand\cY{{\mathcal{Y}}}

\newcommand\fg{{\mathfrak{g}}}

\newcommand\fH{{\mathfrak{H}}}
\newcommand\fk{{\mathfrak{k}}}
\newcommand\fK{{\mathfrak{K}}}

\newcommand\fM{{\mathfrak{M}}}

\newcommand\fN{{\mathfrak{N}}}

\newcommand\fP{{\mathfrak{P}}}

\newcommand\fs{{\mathfrak{s}}}
\newcommand\fS{{\mathfrak{S}}}
\newcommand\ft{{\mathfrak{t}}}

\newcommand\fu{{\mathfrak{u}}}

\newcommand\fV{{\mathfrak{V}}}

\newcommand\sA{{\mathscr{A}}}
\newcommand\sB{{\mathscr{B}}}
\newcommand\sC{{\mathscr{C}}}

\newcommand\sE{{\mathscr{E}}}
\newcommand\sF{{\mathscr{F}}}
\newcommand\sG{{\mathscr{G}}}
\newcommand\sH{{\mathscr{H}}}
\newcommand\sI{{\mathscr{I}}}
\newcommand\sJ{{\mathscr{J}}}
\newcommand\sK{{\mathscr{K}}}
\newcommand\sL{{\mathscr{L}}}
\newcommand\sM{{\mathscr{M}}}
\newcommand\sN{{\mathscr{N}}}
\newcommand\sO{{\mathscr{O}}}
\newcommand\sP{{\mathscr{P}}}
\newcommand\sQ{{\mathscr{Q}}}
\newcommand\sR{{\mathscr{R}}}
\newcommand\sS{{\mathscr{S}}}
\newcommand\sT{{\mathscr{T}}}
\newcommand\sU{{\mathscr{U}}}
\newcommand\sV{{\mathscr{V}}}

\newcommand\sX{{\mathscr{X}}}
\newcommand\sY{{\mathscr{Y}}}
\newcommand\sZ{{\mathscr{Z}}}

\newcommand\beps{{\boldsymbol{\varepsilon}}}

\newcommand\bga{{\boldsymbol{\gamma}}}
\newcommand\bgamma{{\boldsymbol{\gamma}}}
\newcommand\bchi{{\boldsymbol{\chi}}}
\newcommand\biota{{\boldsymbol{\iota}}}
\newcommand\bkappa{{\boldsymbol{\kappa}}}

\newcommand\bmu{{\boldsymbol{\mu}}}

\newcommand\bnu{{\boldsymbol{\nu}}}

\newcommand\bomega{{\boldsymbol{\omega}}}
\newcommand\bOmega{{\boldsymbol{\Omega}}}

\newcommand\bpsi{{\boldsymbol{\psi}}}

\newcommand\bxi{{\boldsymbol{\xi}}}

\newcommand\bB{{\mathbf{B}}}

\newcommand\bC{{\mathbf{C}}}

\newcommand\bD{{\mathbf{D}}}

\newcommand\bE{{\mathbf{E}}}

\newcommand\bF{{\mathbf{F}}}

\newcommand\bg{{\mathbf{g}}}

\newcommand\bh{{\mathbf{h}}}
\newcommand\bH{{\mathbf{H}}}

\newcommand\bJ{{\mathbf{J}}}

\newcommand\bL{{\mathbf{L}}}

\newcommand\bM{{\mathbf{M}}}

\newcommand\bs{{\mathbf{s}}}

\newcommand\bv{{\mathbf{v}}}
\newcommand\bV{{\mathbf{V}}}
\newcommand\bw{{\mathbf{w}}}

\newcommand\bx{{\mathbf{x}}}
\newcommand\bX{{\mathbf{X}}}

\newcommand{\cov}{\nabla}

\newcommand{\rd}{\partial}

\newcommand\eps{\varepsilon}

\newcommand\la{\lambda}
\newcommand\La{\Lambda}
\newcommand\ka{\kappa}
\newcommand\om{\omega}
\newcommand\Om{\Omega}
\newcommand\si{\sigma}

\newcommand\gl{{\mathfrak{g}\mathfrak{l}}}
\newcommand\fsl{{\mathfrak{s}\mathfrak{l}}}

\newcommand\su{{\mathfrak{s}\mathfrak{u}}}

\newcommand\GL{\operatorname{GL}}
\newcommand\Or{\operatorname{O}}

\newcommand\SL{\operatorname{SL}}
\newcommand\SO{\operatorname{SO}}

\newcommand\Spin{\operatorname{Spin}}
\newcommand\SU{\operatorname{SU}}
\newcommand\U{\operatorname{U}}

\newcommand\less{\setminus}

\newcommand\ad{{\operatorname{ad}}}
\newcommand\Ad{{\operatorname{Ad}}}

\DeclareMathOperator{\Aut}{Aut}

\newcommand\ch{\operatorname{ch}}

\newcommand\CCl{\operatorname{{\mathbb{C}\ell}}}
\DeclareMathOperator{\Coh}{Coh}

\newcommand\Coker{\operatorname{Coker}}

\DeclareMathOperator{\Crit}{Crit}

\newcommand\dist{\operatorname{dist}}

\newcommand\End{\operatorname{End}}

\DeclareMathOperator{\Euler}{Euler}
\newcommand\Exp{\operatorname{Exp}}

\newcommand\grad{\operatorname{grad}}

\newcommand\hess{\operatorname{hess}}
\newcommand\Hess{\operatorname{Hess}}

\newcommand\Hol{\operatorname{Hol}}
\newcommand\Hom{\operatorname{Hom}}

\DeclareMathOperator{\Imag}{Im}

\DeclareMathOperator{\ind}{ind}

\newcommand\Ker{\operatorname{Ker}}

\DeclareMathOperator{\Mat}{Mat}

\newcommand\PD{\operatorname{PD}}
\DeclareMathOperator{\Pic}{Pic}
\DeclareMathOperator{\pr}{pr}

\DeclareMathOperator{\pt}{pt}

\newcommand\Ran{\operatorname{Ran}}
\newcommand\rank{\operatorname{rank}}

\newcommand\Real{\operatorname{Re}}
\newcommand\Red{\operatorname{Red}}

\DeclareMathOperator{\Res}{Res}

\DeclareMathOperator{\SFF}{II}

\newcommand\Stab{\operatorname{Stab}}

\newcommand\supp{\operatorname{supp}}

\DeclareMathOperator{\SW}{SW}
\DeclareMathOperator{\Sym}{Sym}

\newcommand\tr{\operatorname{tr}}

\newcommand\vol{\operatorname{vol}}
\newcommand\Vol{\operatorname{Vol}}

\DeclareMathOperator{\Zero}{Zero}

\newcommand\apriori{{\emph{a priori }}}

\newcommand\can{{\mathrm{can}}}

\newcommand\even{{\mathrm{even}}}
\DeclareMathOperator{\expdim}{\mathrm{exp\ dim}}

\newcommand\id{{\mathrm{id}}}

\newcommand\loc{{\mathrm{loc}}}

\newcommand\mutatis{{\emph{mutatis mutandis }}}

\newcommand\odd{{\mathrm{odd}}}

\newcommand\Op{{\mathrm{Op}}}

\newcommand\ps{{\mathrm{ps}}}
\newcommand\red{{\mathrm{red}}}
\newcommand\reg{{\mathrm{reg}}}

\newcommand\spin{\text{spin}}
\newcommand\spinc{\text{$\mathrm{spin}^c$ }}
\newcommand\spinu{\text{$\mathrm{spin}^u$ }}
\newcommand\Spinc{\text{$\mathrm{Spin}^c$}}

\newcommand\symp{\mathrm{symp}}

\newcommand\vir{\mathrm{vir}}

\numberwithin{equation}{section}
\numberwithin{section}{chapter}
\numberwithin{figure}{section}

\usepackage{amscd}
\usepackage{amsfonts}
\usepackage{amssymb}
\usepackage{enumerate}
\usepackage{graphicx}
\usepackage{hyperref}
\usepackage{mathrsfs}
\usepackage[usenames]{color}
\usepackage{paralist}
\usepackage{slashed}
\usepackage{stmaryrd}
\usepackage{svg}
\usepackage{tikz-cd}
\usepackage{url}
\usepackage{verbatim}
\usepackage{xr}
\usepackage[all]{xy}

\usepackage{float}
\usepackage{morefloats}
\usepackage[top=1.3in, bottom=1.1in, left=1.3in, right=1.3in,marginpar=1in]{geometry}

\hypersetup{pdftitle={Virtual Morse--Bott index, moduli spaces of pairs, and applications to the topology of smooth four-manifolds}}
\hypersetup{pdfauthor={Paul M. N. Feehan and Thomas G. Leness}}

\begin{document}

\frontmatter

\title[Virtual Morse--Bott index and moduli spaces of pairs]{Virtual Morse--Bott index, moduli spaces of pairs, and applications to the topology of smooth four-manifolds}

\author[Paul M. N. Feehan]{Paul M. N. Feehan}
\author[Thomas G. Leness]{Thomas G. Leness}

\dedicatory{Paul Feehan dedicates this monograph to his wife, Julia, and children, Odilla and Ois{\'\i}n. \\[2pt] Thomas Leness dedicates this monograph to his wife, Kirsten, children Anne and George, and to the memory of his father, John G. Leness.}






\maketitle

\setcounter{page}{7}

\tableofcontents
\listoffigures

\chapter*{Preface}
\label{chap:Preface}
We previously developed an approach to Bia{\l}ynicki--Birula theory for holomorphic $\mathbb{C}^*$ actions on complex analytic spaces and the concept of virtual Morse--Bott indices for singular critical points of Hamiltonian functions for the induced circle actions (see Feehan \cite{Feehan_analytic_spaces}). For Hamiltonian functions of circle actions on closed, complex K\"ahler manifolds, the virtual Morse--Bott index coincides with the classical Morse--Bott index due to Bott \cite{Bott_1954} and Frankel \cite{Frankel_1959}. A key principle in our approach is that positivity of the virtual Morse--Bott index at a critical point of the Hamiltonian function implies that the critical point cannot be a local minimum even when that critical point is a singular point in the moduli space. In this monograph, we consider our method in the context of the moduli space of non-Abelian monopoles over a closed, complex, K\"ahler surface. We use the Hirzebruch--Riemann--Roch Theorem to compute virtual Morse--Bott indices of all critical strata (Seiberg--Witten moduli subspaces) and we prove that these indices are positive in a setting motivated by the conjecture that all closed, smooth four-manifolds of Seiberg--Witten simple type obey the Bogomolov--Miyaoka--Yau inequality.

\chapter*{Acknowledgments}
\label{chap:Acknowledgments}
We warmly thank Anar Akhmedov, R. Inanc Baykur, Hans Boden, Steven Bradlow, Robert Bryant, Joanna Cirici, Ian Coley, Jean--Pierre Demailly, Simon Donaldson, Tedi Draghici, Mariano Echeverria, Yue Fan, Luis Fernandes, Robert Friedman, Oscar Garc{\'\i}a--Prada, Chris Gerig, Gueo Grantcharov, Daniel Halpern--Leistner, Ian Hambleton, Nigel Hitchin, David Hurtubise, Dominic Joyce, Yi--Jen Lee, Manfred Lehn, Chi Li, Tian--Jun Li, Jason Lotay, Konstantin Mischaikow, Tom Mrowka, Hiraku Nakajima, Karl--Hermann Neeb, Thomas Parker, Duong Phong, Gerd Rudolph, Siddhartha Sahi, Tobias Shin, Ivan Smith, Jian Song, Andras Stipsicz, Dennis Sullivan, Andrei Teleman, Michael Thaddeus, Valentino Tosatti, Charles Weibel, Graeme Wilkin, Scott Wilson, Boyu Zhang, and Zhengyi Zhou for helpful communications. We are especially grateful to Richard Wentworth for his expertise, thoughtful insights, and numerous helpful discussions over many years. We warmly thank the anonymous referees for their comments on our manuscript. We appreciate help with the figure in our monograph provided by John Etnyre, Kyle Hayden, Kristen Hendricks, and Tai Wai Hu. We thank the National Science Foundation for its support of our research via the grants DMS-1510064 and DMS-2104865 (Feehan) and DMS-1510063 and DMS-2104871 (Leness). This monograph is also based in part on work supported by the National Science Foundation under Grant No. 1440140, while one of the authors (Feehan) was in residence at the Mathematical Sciences Research Institute in Berkeley, California, during Fall 2022 as a Research Professor in the program \emph{Analytic and Geometric Aspects of Gauge Theory}.

\bigskip
\bigskip

\leftline{Paul M. N. Feehan}
\leftline{Department of Mathematics}
\leftline{Rutgers, The State University of New Jersey}
\leftline{Piscataway, NJ 08854}
\leftline{United States}
\medskip

\leftline{\texttt{feehan@math.rutgers.edu}}
\leftline{\url{math.rutgers.edu/~feehan}}
\bigskip

\leftline{Thomas G. Leness}
\leftline{Department of Mathematics}
\leftline{Florida International University}
\leftline{Miami, FL 33199}
\leftline{United States}
\medskip

\leftline{\texttt{lenesst@fiu.edu}}
\leftline{\url{fiu.edu/~lenesst}}
\bigskip

\leftline{This version: July 3, 2023}

\mainmatter

\chapter{Introduction}
\label{chap:Introduction}
In Feehan \cite{Feehan_analytic_spaces}, we developed an approach to Morse--Bott theory (more precisely, Bia{\l}ynicki--Birula theory), that applies to singular complex analytic spaces that typically arise in gauge theory --- including moduli spaces of non-Abelian monopoles over closed, complex K\"ahler surfaces (see Feehan and Leness \cite{FL5, FL6, FL7, FL8} and references therein), moduli spaces of projective vortices (see Bradlow \cite{Bradlow_1991} and Garc{\'\i}a--Prada \cite{GarciaPradaDimRedVort}, \cite{BradlowGP}), moduli spaces of stable holomorphic pairs of bundles and sections over closed complex K\"ahler manifolds (see Huybrechts and Lehn \cite{Huybrechts_Lehn_1995jag}), and moduli spaces of Higgs pairs over closed Riemann surfaces (see Hitchin \cite{Hitchin_1987}). When defined over base manifolds that are complex K\"ahler, such moduli spaces are \emph{complex analytic spaces}, equipped with K\"ahler metrics (on the top stratum of smooth points) and Hamiltonian circle actions. Our development of Morse theory for the corresponding Hamiltonian functions extends one pioneered by Hitchin in his celebrated study of the moduli space of Higgs pairs over Riemann surfaces \cite{Hitchin_1987}.

When a manifold is almost Hermitian (as are four-manifolds of Seiberg--Witten simple type) --- where the almost complex structure is not necessarily integrable and the fundamental two-form defined by the almost complex structure and Riemannian metric is not necessarily closed --- one might conjecture that moduli spaces of non-Abelian monopoles, and moduli spaces of projective vortices more generally, can be shown to be almost Hermitian under suitable conditions.

We begin in Section \ref{sec:Frankel_theorem_circle_actions_almost_Hermitian_manifolds} by describing an extension (Theorem \ref{mainthm:Frankel_almost_Hermitian}) of Frankel's Theorem \cite[Section 3]{Frankel_1959} for Hamiltonian circle actions, originally on complex K\"ahler manifolds and extended here to almost Hermitian and almost symplectic manifolds. In Section \ref{sec:Virtual_Morse-Bott_index_Hamiltonian_function_circle_action_complex_analytic_space}, we introduce the concept of a \emph{virtual Morse--Bott index} for the Hamiltonian function of a circle action on a complex analytic space and explain that positivity of the virtual Morse--Bott index of a critical point implies that it cannot be a local minimum. We refer to Feehan \cite{Feehan_analytic_spaces} for precise statements and proofs of those results.

One may ask whether moduli spaces of non-Abelian monopoles over closed, smooth four-manifolds of Seiberg--Witten simple type can be used to constrain the topology of such four-manifolds and this is the primary motivation for our present monograph. An early example of such arguments was provided by Donaldson \cite{DonApplic}, who used the moduli space of anti-self-dual connections to prove that if the intersection form of a closed, simply connected, smooth four-manifold is negative definite, then the intersection form is diagonal. In Section \ref{sec:Morse_theory_existence_anti-self-dual_connections}, we outline how Morse--Bott theory may constrain the geography of closed, smooth four-manifolds that have $b_1=0$, odd $b^+ \geq 3$, and Seiberg--Witten simple type.

In Section \ref{sec:Morse_theory_moduli_space_PU2_monopoles_Kaehler}, we discuss the virtual Morse--Bott index (and nullity and co-index) for critical points of a Hamiltonian function on the moduli space of non-Abelian monopoles over closed, complex K\"ahler surface and introduce the main results of this monograph. Theorem \ref{mainthm:IdentifyCriticalPoints} shows that points in these moduli spaces are \emph{critical points} of the Hamiltonian function (the square of the $L^2$ norm of the spinor section, by analogy with Hitchin's definition in terms of the Higgs field in \cite[Section 7]{Hitchin_1987}) if and only if they represent \emph{Seiberg--Witten monopoles} or \emph{anti-self-dual connections}. Theorem \ref{mainthm:ExistenceOfSpinuForFlow} shows that one can always choose moduli spaces of non-Abelian monopoles to support the strategy outlined in Section \ref{sec:Morse_theory_existence_anti-self-dual_connections}.
Theorem \ref{mainthm:MorseIndexAtReduciblesOnKahler} gives the formula for the virtual Morse--Bott index of the Hamiltonian function for a point represented by a type $1$ Seiberg--Witten monopole, a formula computed by an application of the Atiyah--Singer Index Theorem (more specifically, the Hirzebruch--Riemann--Roch Index Theorem in this monograph).  Corollary \ref{maincor:MorseIndexForFeasibilitySpinuStructure} is an application of Theorem \ref{mainthm:MorseIndexAtReduciblesOnKahler} (and Corollary \ref{maincor:MorseIndexAtReduciblesOnKahlerWithSO3MonopoleCharacteristicClasses}) which shows that the virtual Morse--Bott indices of points representing type $1$ Seiberg--Witten monopoles in the moduli spaces constructed in Theorem \ref{mainthm:ExistenceOfSpinuForFlow} are \emph{positive}\footnote{Hitchin obtained an analogous result \cite[Proposition 7.1, p. 92]{Hitchin_1987} for points in the moduli space of Higgs pairs where the connections are reducible.} and thus cannot be local minima. The Hamiltonian function is zero at points representing anti-self-dual connections and so they are absolute minima. Section \ref{sec:Outline} contains an outline of the remainder of this monograph. A discussion of which of these results hold for four-manifolds that are not complex K\"ahler and how these results might be generalized to smooth four-dimensional manifolds of Seiberg--Witten simple type appears in Remark \ref{rmk:Virtual_Morse-Bott_index_formula_Kaehler_surface_versus_smooth_4-manifold}. In Chapter \ref{chap:Bubbling}, we discuss approaches to the problem that the moduli space of non-Abelian monopoles is non-compact and that its Gieseker or Uhlenbeck compactifications have singularities that are more complicated than those considered in this monograph.

\section[Frankel's Theorem for Hamiltonian circle action on almost symplectic manifold]{Frankel's Theorem for the Hamiltonian function of a circle action on a smooth almost symplectic manifold}
\label{sec:Frankel_theorem_circle_actions_almost_Hermitian_manifolds}
The forthcoming version, Theorem \ref{mainthm:Frankel_almost_Hermitian}, of Frankel's Theorem \cite[Section 3]{Frankel_1959} that we prove in this monograph is more general than that stated in \cite{Frankel_1959} because we allow for circle actions on closed, smooth manifolds $(M,g,J)$ that are only assumed to be \emph{almost Hermitian}\label{page:Almost_Hermitian_manifold}, so the $g$-orthogonal almost complex structure $J$ need not be integrable and the fundamental two-form
\begin{equation}
  \label{eq:Fundamental_two-form}
  \omega = g(\cdot,J\cdot)
\end{equation}
is non-degenerate but not required to be closed, whereas Frankel assumed in \cite[Section 3]{Frankel_1959} that $\omega$ was closed. (Our convention in \eqref{eq:Fundamental_two-form} agrees with that of Kobayashi \cite[Equation (7.6.8), p. 251]{Kobayashi_differential_geometry_complex_vector_bundles} but is opposite to that used elsewhere, for example, Huybrechts \cite[Definition 1.2.13, p. 29]{Huybrechts_2005}.) Frankel notes \cite[p. 1]{Frankel_1959} that the main results of his article hold when $\omega$ is a $g$-harmonic, symplectic form. (If $\omega$ is symplectic, then $d\omega=0$ while if $g$ is \emph{adapted} to $\omega$, it is well-known that $d^{*_g}\omega=0$ --- see Delano\"e \cite{Delanoe_2002} --- and so $\omega$ is harmonic.) Recall that\footnote{If $E$ is a smooth vector bundle over a smooth manifold, we let $\Omega^0(E) = C^\infty(E)$ denote the Fr\'echet space of smooth sections of $E$.} $J \in C^\infty(\End(TM))$ is an \emph{almost complex structure} 
\label{page:Almost_complex_structure} on $M$ if $J^2 = -\id_{TM}$ and $(M,J)$ is thus an \label{page:Almost_complex_manifold} \emph{almost complex manifold}. One says that $J$ is \emph{orthogonal with respect to} or \emph{compatible with} a Riemannian metric $g$ on $M$ if
\begin{equation}
  \label{eq:g_compatible_J}
  g(JX,JY) = g(X,Y),
\end{equation}
for all vector fields $X, Y \in C^\infty(TM)$. Recall that a smooth manifold is called \emph{almost symplectic} \label{page:almost_symplectic_manifold} if it admits a non-degenerate two-form and \emph{symplectic} \label{page:symplectic_manifold} if that two-form is closed (see Libermann and Marle \cite[Definition 12.4]{Libermann_Marle_symplectic_geometry_analytical_mechanics}).

If a smooth manifold $M$ admits a smooth circle action,
\begin{equation}
\label{eq:Circle_action_smooth_manifold}
\rho:S^1\times M\to M,
\end{equation}
we denote the induced circle action on the tangent bundle $TM$ by
\begin{equation}
\label{eq:Circle_action_tangent_bundle}
\rho_*:S^1\times TM\to TM,
\end{equation}
where, using $D_2\rho$ to denote the differential of $\rho$ in directions tangent to $M$,
\begin{equation}
  \label{eq:rho_*_definition}
  \rho_*(e^{i\theta},v) := D_2\rho(e^{i\theta},x)v,
  \quad\text{for all } x \in M, \ v\in T_xM, \text{ and } e^{i\theta}\in S^1.
\end{equation}
A covariant two-tensor field $\varpi \in C^\infty(T^*M\otimes T^*M)$ is called \emph{circle-invariant}
\label{page:circle_invariant} (with respect to \eqref{eq:Circle_action_smooth_manifold}) if it obeys
\begin{equation}
\label{eq:Circle_invariant_covariant_2-tensor}
  \varpi\left(\rho_*(e^{i\theta})v, \rho_*(e^{i\theta})w\right) = \varpi(v,w),
  \quad\text{for all } p \in M, \ v, w \in T_pM, \text{ and } e^{i\theta} \in S^1.
\end{equation}
Similarly, a tensor field $\tau \in C^\infty(TM\otimes T^*M) = C^\infty(\End(TM))$ is called \emph{circle-invariant} (with respect to \eqref{eq:Circle_action_smooth_manifold}) if it obeys
\begin{equation}
\label{eq:Circle_invariant_(1,1)-tensor}
  \tau\left(\rho_*(e^{i\theta})v\right) = \rho_*(e^{i\theta})\tau v,
  \quad\text{for all } p \in M, \ v \in T_pM, \text{ and } e^{i\theta} \in S^1.
\end{equation}
(One could, of course, absorb definitions \eqref{eq:Circle_invariant_covariant_2-tensor} and \eqref{eq:Circle_invariant_(1,1)-tensor} into a general definition of a tensor field invariant under the flow of a vector field, as in Lee \cite[Equation (12.11), p. 323]{Lee_john_smooth_manifolds} and characterize covariant tensor fields that are invariant under the flow of a vector field in terms of vanishing Lie derivatives as in  Lee \cite[Theorem 12.37, p. 324]{Lee_john_smooth_manifolds}.) A smooth two-form $\omega$ or smooth Riemannian metric $g$ on $M$ is called \emph{circle-invariant} (with respect to \eqref{eq:Circle_action_smooth_manifold}) if it obeys \eqref{eq:Circle_invariant_covariant_2-tensor} with $\varpi$ replaced by $\omega$ or $g$, respectively. A circle action is called \emph{Hamiltonian} with respect to a smooth two-form $\omega$ on $M$ if \footnote{By analogy with the usual meaning \cite[Definition 2.1]{Dwivedi_Herman_Jeffrey_van_den_Hurk} of a Hamiltonian vector field and Hamiltonian function.} there exists a smooth function $f:M\to\RR$ such that
\begin{equation}
\label{eq:MomentMap}
  df = \iota_X\omega,
\end{equation}
where $X \in C^\infty(TM)$ is the vector field generated by the circle action, so $X_p = D_1\rho(1,p) \in T_pM$ for all $p\in M$, with $D_1\rho$ denoting the differential of $\rho$ in directions tangent to $S^1$. Adapting Bott \cite[Definition, p. 248]{Bott_1954}, \cite{Bott_1959} and Nicolaescu \cite[Definition 2.41]{Nicolaescu_morse_theory} (see also Feehan \cite[Definition 1.2, p. 3279]{Feehan_lojasiewicz_inequality_all_dimensions}), we make the

\label{page:Morse-Bott_function}
\begin{maindefn}[Morse--Bott function]
\label{maindefn:Morse-Bott_function}  
Let $(M,g)$ be a smooth manifold and $f:M\to\RR$ be a smooth function. We let
\[
  \Crit f := \{p\in M: df(p)=0\}
\]
denote the \emph{critical set} of $f$. The function $f$ is \emph{Morse--Bott at $p$} if there exists a small enough open neighborhood $U\subset M$ of $p$ such that $U\cap \Crit f$ is a smooth submanifold with tangent space
\[
  T_p\Crit f = \Ker\hess f(p),
\]
where $\hess f(p) \in \Hom(T_pM,T_p^*M)$ is the Hessian form defined in Section \ref{sec:Introduction_Hessian_restriction_smooth_function_submanifold_Euclidean_space}. The function $f$ is \emph{Morse--Bott} if it is Morse--Bott at every point of $\Crit f$.
\end{maindefn}  

\begin{mainthm}[Frankel's theorem for circle actions on almost Hermitian manifolds]
\label{mainthm:Frankel_almost_Hermitian}
(Compare Frankel \cite[Section 3]{Frankel_1959}.)
Let $M$ be a smooth manifold endowed with a smooth circle action $\rho$ as in \eqref{eq:Circle_action_smooth_manifold} and a non-degenerate two-form $\om$ that is circle-invariant in the sense of \eqref{eq:Circle_invariant_covariant_2-tensor}. Then
\begin{enumerate}
\item
\label{item:Frankel_almost_Hermitian_FixedPointsAreZerosOfVField}
A point $p\in M$ is a fixed point of the action \eqref{eq:Circle_action_smooth_manifold} if and only if $X_p=0$, where $X \in C^\infty(TM)$ is the vector field generated by the circle action.

\item
\label{item:Frankel_almost_Hermitian_ComponentsOfFixedPointsAreSmoothSubmanifolds}
Each connected component of the fixed-point set of the circle action \eqref{eq:Circle_action_smooth_manifold} is a smooth submanifold of even dimension and even codimension in $M$.
\end{enumerate}
In addition, let $f:M\to\RR$ be a smooth function that is Hamiltonian in the sense of \eqref{eq:MomentMap}. Then
\begin{enumerate}
\setcounter{enumi}{2}
\item
\label{item:Frankel_almost_Hermitian_FixedPointsAreCriticalPoints}
A point $p\in M$ is a critical point of $f$ if and only if $p$ is a fixed point of the circle action \eqref{eq:Circle_action_smooth_manifold}.
\end{enumerate}

\begin{enumerate}
\setcounter{enumi}{3}
\item
  \label{item:Frankel_almost_Hermitian_f_is_MB}
  The function $f$ is Morse--Bott at $p$ in the sense of Definition \ref{maindefn:Morse-Bott_function}.
\end{enumerate}
Furthermore, assume that there is a smooth Riemannian metric on $M$ that is circle-invariant in the sense of \eqref{eq:Circle_invariant_covariant_2-tensor}. Then
\begin{enumerate}
\setcounter{enumi}{4}
\item
\label{item:Frankel_almost_Hermitian_ACisS1Invariant}
There are a smooth almost complex structure $J$ on $TM$ and a smooth Riemannian metric $g$ on $M$ such that  $(\omega,g,J)$ is a compatible triple in the sense that it obeys \eqref{eq:Fundamental_two-form} and \eqref{eq:g_compatible_J}, and $g$ is circle-invariant in the sense of \eqref{eq:Circle_invariant_covariant_2-tensor}, and $J$ is circle-invariant in the sense of \eqref{eq:Circle_invariant_(1,1)-tensor}.

\item
\label{item:Frankel_almost_Hermitian_WeightsAreEigenvalues}
The eigenvalues of the Hessian operator $\Hess_g f(p) \in \End(T_pM)$ defined in Section \ref{sec:Introduction_Hessian_restriction_smooth_function_submanifold_Euclidean_space} are given by the weights of the circle action on $(T_pM,J)$ if the signs of the weights are chosen to be compatible with $J$ in the sense of Definition \ref{defn:Weight_Associated_To_Complex_Structure} and Lemma \ref{lem:S1_Invariant_J_Is_MultiplicationBy_i}.

\item
  \label{item:Frankel_almost_Hermitian_Sylvester}
  The signature $(\lambda_p^+(f),\lambda_p^-(f),\lambda_p^0(f))$ of the Hessian operator $\Hess_g f(p)$ is independent of the Riemannian metric $g$, where $\lambda_p^\pm(f)$ denotes the number of positive (negative) eigenvalues of $\Hess_g f(p)$ and $\lambda_p^0(f)$ denotes the nullity of $\Hess_g f(p)$.
\end{enumerate}
\end{mainthm}

\begin{rmk}[Fixed-point sets of isometric circle actions on Riemannian manifolds are totally geodesic submanifolds]
\label{rmk:Fixed-point_sets_isometric_circle_actions_Riemannian_manifolds_totally_geodesic} 
We recall that if $g$ is a circle-invariant Riemannian metric on $M$, then the vector field $X \in C^\infty(TM)$ generated by the circle action is a Killing field of $g$. A theorem of Kobayashi \cite[p. 63]{Kobayashi_1958} implies that the zero set of a Killing field is a closed, smooth \emph{totally geodesic} submanifold of even codimension and so Item \eqref{item:Frankel_almost_Hermitian_FixedPointsAreZerosOfVField} implies that the fixed-point subset is a closed, smooth \emph{totally geodesic} submanifold of even codimension, strengthening Item  \eqref{item:Frankel_almost_Hermitian_ComponentsOfFixedPointsAreSmoothSubmanifolds}. See Remark \ref{rmk:Simplifications_totally_geodesic_submanifolds} for an explanation of the significance of the totally geodesic property in the context of Morse--Bott theory.
\end{rmk}

We prove Theorem \ref{mainthm:Frankel_almost_Hermitian} in Chapter \ref{chap:Circle_actions_almost_Hermitian_manifolds}. Recall that the \emph{gradient vector field} $\grad_g f \in C^\infty(TM)$ associated to a function $f \in C^\infty(M,\RR)$ is defined by the relation
\begin{equation}
\label{eq:DefineGradient}
  g(\grad_g f,Y) := df(Y), \quad\text{for all } Y \in C^\infty(TM).
\end{equation}  
If $\nabla^g$ is the covariant derivative for the Levi--Civita connection on $TM$ defined by the Riemannian metric $g$, then one can define the Hessian of $f \in C^\infty(M,\RR)$ by (see Petersen \cite[Proposition 2.2.6]{Petersen_2006})
\begin{equation}
\label{eq:DefineHessian}
  \Hess_g f := \nabla^g\grad_g f \in C^\infty(\End(TM)).
\end{equation}
See the forthcoming equation \eqref{eq:Hessian_operator} for an alternative expression for $\Hess_g f$. If $p\in M$ is a critical point of $f$, then Theorem \ref{mainthm:Frankel_almost_Hermitian} implies that subspace $T_p^-M \subset T_pM$ on which the Hessian $\Hess_g f(p) \in \End(T_pM)$ is \emph{negative definite} is equal to the subspace of $T_pM$ on which the circle acts with \emph{negative weight}. Hence, the \emph{Morse--Bott index} \label{page:Morse--Bott_index} of $f$ at a critical point $p$,
\[
  \lambda_p^-(f) := \dim_\RR T_p^-M,
\]
is equal to the dimension of the subspace of $T_pM$ on which the circle acts with \emph{negative weight}; one calls $\lambda_p^+(f)$ and $\lambda_p^0(f)$ the \emph{Morse--Bott co-index} and \emph{nullity}, respectively.
In \cite[Section 3]{Frankel_1959}, Frankel says that the critical set of $f$ in Theorem \ref{mainthm:Frankel_almost_Hermitian} is \emph{non-degenerate at $p$} whereas we use the terminology that $f$ is \emph{Morse--Bott at $p$} following \cite[Definition 1.5]{Feehan_lojasiewicz_inequality_all_dimensions_morse-bott}.

\section[Virtual Morse--Bott index of a critical point for the Hamiltonian function]{Virtual Morse--Bott index of a critical point for the Hamiltonian function of a circle action on a complex analytic space}
\label{sec:Virtual_Morse-Bott_index_Hamiltonian_function_circle_action_complex_analytic_space}
As a well-known application due to Hitchin \cite[Proposition 7.1, p. 92]{Hitchin_1987} illustrates, Theorem \ref{mainthm:Frankel_almost_Hermitian} is remarkably useful, but its hypotheses limit its application to smooth manifolds. One of our goals in this monograph is to indicate how Theorem \ref{mainthm:Frankel_almost_Hermitian} may have useful generalizations to \emph{complex analytic spaces} (see Abhyankar \cite{Abhyankar_local_analytic_geometry}, Fischer \cite{Fischer_complex_analytic_geometry}, Grauert and Remmert \cite{Grauert_Remmert_coherent_analytic_sheaves}, Gunning and Rossi \cite{Gunning_Rossi_analytic_functions_several_complex_variables}, or Narasimhan \cite{Narasimhan_introduction_theory_analytic_spaces}). In this section, we summarize one of the main results of Feehan \cite{Feehan_analytic_spaces} and thus provide motivation for our forthcoming
Theorem \ref{mainthm:MorseIndexAtReduciblesOnKahler} and Corollaries 
\ref{maincor:MorseIndexAtReduciblesOnKahlerWithSO3MonopoleCharacteristicClasses} and \ref{maincor:MorseIndexForFeasibilitySpinuStructure}. 

Suppose that $X$ is a complex, finite-dimensional, K\"ahler manifold with circle action that is compatible with the standard almost complex structure and the induced Riemannian metric. We further assume that the circle action is Hamiltonian with real analytic Hamiltonian function $f:X\to\RR$ such that $df = \iota_\xi\omega$, where $\omega$ is the K\"ahler two-form on $X$ and $\xi$ is the vector field on $X$ generated by the circle action. We let $Y\subset X$ be a circle-invariant, closed, complex analytic subspace, $p \in Y$ be a point, and $F:U\to\CC^r$ be a complex analytic, local defining function for $Y$ on an open neighborhood $U\subset X$ of $p$ in the sense that $Y\cap U = F^{-1}(0)$ and $\sO_{Y\cap U} = (\sO_U/\sI)\restriction Y\cap U$ is the structure sheaf for $Y\cap U$ defined by the ideal $\sI = (f_1,\ldots,f_r) \subset \sO_U$ given by the component functions $f_j$ of $F$. We let the complex linear subspace $\bH_p^2 \subset \CC^r$ be the orthogonal complement of $\Ran dF(p) \subset \CC^r$ and let $\bH_p^1 = \Ker dF(p) \subset T_pX$ denote the (complex linear) 
\emph{Zariski tangent space} to $Y$ at $p$
\label{Zariski_Tangent_Space}. We let $S \subset X$ be the circle-invariant, complex, K\"ahler submanifold given by $F^{-1}(\bH_p^2)$ and observe that $T_pS = \bH_p^1$. If $p$ is a \emph{critical point}
\label{page:Critical_point} of $f:Y\to\RR$ in the sense that $\bH_p^1 \subseteq \Ker df(p)$, then $p$ is a fixed point of the induced circle action on $S$. Because $p\in X$ is also a fixed point of the circle action on $X$, one can show that the holomorphic map $F$ determines a circle action on $\CC^r$ and that the subspace $\bH_p^2$ is circle-invariant with respect to this induced action. We let $\bH_p^{1,\pm} \subset \bH_p^1$ and $\bH_p^{2,\pm} \subset \bH_p^2$ be the subspaces on which the circle acts with positive (respectively, negative) weight and let $\bH_p^{1,0} \subset \bH_p^1$ and $\bH_p^{2,0} \subset \bH_p^2$ be the subspaces on which the circle acts with zero weight. We define the \emph{virtual Morse--Bott nullity, co-index} and \emph{index}, respectively, for $f$ at $p$ by
\begin{subequations}
\label{eq:Virtual_Morse-Bott_signature}  
\begin{align}
  \label{eq:Virtual_Morse-Bott_nullity}
  \lambda_p^0(f) &:= \dim_\RR\bH_p^{1,0} - \dim_\RR\bH_p^{2,0},
  \\
  \label{eq:Virtual_Morse-Bott_coindex}
  \lambda_p^+(f) &:= \dim_\RR\bH_p^{1,+} - \dim_\RR\bH_p^{2,+},
  \\
  \label{eq:Virtual_Morse-Bott_index}
  \lambda_p^-(f) &:= \dim_\RR\bH_p^{1,-} - \dim_\RR\bH_p^{2,-}.
\end{align}
\end{subequations}
One of the main results of Feehan \cite{Feehan_analytic_spaces} asserts that if $\lambda_p^-(f)$ (respectively, $\lambda_p^+(f)$) is positive, then $p$ cannot be a local minimum (respectively, maximum) for $f:Y\to\RR$.

\section{Existence of anti-self-dual connections and the Bogomolov--Miyaoka--Yau inequality}
\label{sec:Morse_theory_existence_anti-self-dual_connections}
Let us now turn to the application that motivates this monograph. For a closed topological four-manifold $X$, we define
\begin{equation}
\label{eq:DefineTopCharNumbersOf4Manifold}
  c_1(X)^2:= 2e(X)+3\sigma(X) \quad\text{and}\quad \chi_h(X):=\frac{1}{4}(e(X)+\sigma(X)),
\end{equation}
where
\begin{equation}
  \label{eq:Euler_characteristic_signature_4-manifold}
  e(X) : =2-2b_1(X)+b_2(X) \quad\text{and}\quad \sigma(X) := b^+(X)-b^-(X)
\end{equation}
are the \emph{Euler characteristic} and \emph{signature} of $X$, respectively. If $Q_X$
\label{page:IntersectionForm}
 denotes the intersection form on $H_2(X;\ZZ)$, then $b^\pm(X)$ are the dimensions of the maximal positive (negative) subspaces of $Q_X$ on $H_2(X;\RR)$.

\begin{defn}[Standard four-manifold]
\label{defn:Standard}
We say that a four-manifold $X$ is \emph{standard} if it is closed, connected, oriented, and smooth with odd $b^+(X) \geq 3$ and $b_1(X) = 0$. 
\end{defn}

We can then formulate the

\begin{conj}[Bogomolov--Miyaoka--Yau (BMY) inequality for four-manifolds with non-zero Seiberg--Witten invariants]
\label{conj:BMY_Seiberg-Witten}
If $X$ is a standard four-manifold of Seiberg--Witten simple type with a non-zero Seiberg--Witten invariant, then
\begin{equation}
  \label{eq:BMY}
  c_1(X)^2 \leq 9\chi_h(X).
\end{equation}
\end{conj}

We refer to \cite{MorganSWNotes, NicolaescuSWNotes, SalamonSWBook, Witten} for the definition of Seiberg--Witten invariants and simple type. If $X$ obeys the hypotheses of Conjecture \ref{conj:BMY_Seiberg-Witten}, then it has an almost complex structure $J$ \cite{KMThom} and in the inequality \eqref{eq:BMY}, which is equivalent to
\[
  c_1(X)^2 \leq 3c_2(X),
\]
the Chern classes are those of the complex vector bundle $(TX,J)$.
Inequality \eqref{eq:BMY} was proved by Miyaoka \cite{Miyaoka_1977} and Yau \cite{YauPNAS, Yau} for compact complex surfaces of general type, which have Kodaira dimension $2$, while a weaker version was proved by Bogomolov \cite{Bogomolov_1978}. Inequality \eqref{eq:BMY} also holds for compact complex surfaces with Kodaira dimension $1$ (elliptic surfaces) and Kodaira dimension $0$ (abelian, Enriques, hyperelliptic, and $K3$ surfaces) \cite{Barth_Hulek_Peters_Van_de_Ven_compact_complex_surfaces}. Conjecture \ref{conj:BMY_Seiberg-Witten} is based on \cite[Problem 4]{Stern_2006} (see also Koll\'ar \cite{Kollar_2008}), though often stated for simply connected, symplectic four-manifolds --- see Gompf and Stipsicz \cite[Remark 10.2.16 (c), p. 401]{GompfStipsicz} or Stern \cite[Problem 2]{Stern_2006}. Taubes \cite{TauSymp, TauSympMore} proved that symplectic four-manifolds have non-zero Seiberg--Witten invariants, generalizing Witten's result for K\"ahler surfaces \cite{Witten}, while Kotschick, Morgan, and Taubes \cite{Kotschick_Morgan_Taubes_1995} and Szab\'o \cite{SzaboNoSymplectic} proved existence of non-symplectic four-manifolds with non-zero Seiberg--Witten invariants. Indeed, Szab\'o constructs a family $\{X_n\}_{n\in\NN}$ of closed, irreducible, oriented, simply connected, smooth, four-dimensional manifolds such that, for each integer $n\geq 1$, one has that $b^+(X_n) = 3$, $b^-(X_n) = 19$, 
and the Seiberg--Witten invariant of $X_n$ is non-zero and divisible by $n$ while the Seiberg--Witten
invariant of $\bar X_n$ is zero. From equations \eqref{eq:DefineTopCharNumbersOf4Manifold}, one finds that $c_1(X)^2 = 0$ and $\chi_h(X) = 2$, so these examples obey the inequality \eqref{eq:BMY}.

We note that the hypotheses of Conjecture \ref{conj:BMY_Seiberg-Witten} include the condition that $X$ be a standard four-manifold, which includes the condition $b_1(X)=0$.  We discuss our use of this condition in the forthcoming Remark \ref{rmk:Onb1Equals0}. The proofs of inequality \eqref{eq:BMY} due to Miyaoka \cite{Miyaoka_1977} and Yau \cite{YauPNAS, Yau} for compact complex surfaces of general type do \emph{not} require that $b_1(X)=0$.

Efforts to prove Conjecture \ref{conj:BMY_Seiberg-Witten}, construct a counterexample, or construct four-manifolds $X$ obeying its hypotheses with $c_1(X)^2$ approaching $9\chi_h(X)$ from below have inspired research by Moishezon and Teicher \cite{Moishezon_Teicher_1987}, Z. Chen \cite{Chen_ZJ_1987,Chen_ZJ_1991}, Akhmedov, Hughes,  Park, Sakalli, and Urzua \cite{Akhmedov_Hughes_Park_2013, Akhmedov_Park_2008jggt, Akhmedov_Park_2010mrl,Akhmedov_Park_Urzua_2010,Akhmedov_Sakalli_2016}, Baldridge, Kirk, and Li \cite{Baldridge_Kirk_2006, Baldridge_Kirk_2007, Baldridge_Li_2005}, Bryan, Donagi, and Stipsicz \cite{Bryan_Donagi_Stipsicz_2001}, Fintushel and Stern \cite{Fintushel_Stern_1994im}, Gompf and Mrowka \cite{GompfNewSymplectic, GompfMrowka}, Hamenst\"adt \cite{Hamenstaedt_2012arxiv}, Park and Stipsicz \cite{Park_Stipsicz_2015, Stipsicz_1998, Stipsicz_2000}, Torres \cite{Torres_2014},
Niepel \cite{Niepel_2005}, and others\footnote{A more complete review of the relevant literature appears in Akhmedov, Hughes and Park \cite{Akhmedov_Park_2008jggt,Akhmedov_Hughes_Park_2013}.}. Conjecture \ref{conj:BMY_Seiberg-Witten} holds for all known examples but inequality \eqref{eq:BMY} can fail for four-manifolds with zero Seiberg--Witten invariants.

For $w\in H^2(X;\ZZ)$ and $4\kappa\in \ZZ$, let $E$ be a rank-two, Hermitian vector bundle over $X$ with $c_1(E)=w$ and Pontrjagin number $p_1(\su(E))=-4\kappa$
\label{page:kaChargeOfBundle}, where $\su(E)\subset\gl(E)$ is the Riemannian vector subbundle of trace-free, skew-Hermitian endomorphisms of $E$. Recall from Donaldson and Kronheimer \cite[Equations (4.2.21) and (4.2.22), p. 137]{DK} that the \emph{expected dimension} of the moduli space $M_\kappa^w(X,g)$ of $g$-anti-self-dual connections on $\su(E)$ is given by
\begin{equation}
\label{eq:Expected_dimension_moduli_space_ASD_connections}
  \expdim M_\kappa^w(X,g) = -2p_1(\su(E)) - 6\chi_h(X),
\end{equation}
noting that equations \eqref{eq:DefineTopCharNumbersOf4Manifold} and \eqref{eq:Euler_characteristic_signature_4-manifold} yield
\[
  \chi_h(X) = \frac{1}{2}\left(1 - b_1(X) + b^+(X)\right).
\]
When $g$ is \emph{generic} in the sense of Donaldson and Kronheimer \cite[Corollary 4.3.19, p. 147]{DK} or Freed and Uhlenbeck \cite[Proposition 3.20]{FU}, then $M_\kappa^w(X,g)$ is a smooth manifold, if non-empty, under the hypotheses of Conjecture \ref{conj:BMY_Seiberg-Witten} on the topology of $X$. (By passing to the smooth blow-up of $X$, we may also assume without loss of generality that $w\pmod{2}$ is \emph{good} in the sense that no integral lift of $w\pmod{2}$ is torsion --- see Definition \ref{defn:Good} --- and so the Uhlenbeck compactification $\bar M_\kappa^w(X,g)$ (Donaldson and Kronheimer\cite[Section 4.4.1, p. 158]{DK}) contains no flat connections.) As $\kappa$ increases relative to $\chi_h(X)$, the expected dimension of $M_\kappa^w(X,g)$ increases and it becomes easier to prove that this moduli space is non-empty. Indeed, existence results due to Taubes \cite[Theorems 1.1 and 1.2]{TauIndef} or Taylor \cite[Theorem 1.1]{Taylor_2002} imply that $M_\kappa^w(X,g)$ is non-empty when $\kappa$ is sufficiently large relative to $\chi_h(X)$.

By contrast, as $\kappa$ becomes smaller (equivalently, as $p_1(\su(E))$ becomes larger), it becomes more difficult to prove that $M_\kappa^w(X,g)$ is non-empty. The \emph{geography} question (see Gompf and Stipsicz \cite{GompfStipsicz}) asks which values of signature and Euler characteristic can be realized by smooth four-manifolds which have a given geometric structure, such as a complex or symplectic structure or an Einstein metric, or admit non-trivial Donaldson or Seiberg--Witten invariants. We consider a bundle $E$ that is constrained by the \emph{fundamental lower bound},
\begin{equation}
\label{eq:p1_lower_bound}
p_1(\su(E)) \geq c_1(X)^2-12\chi_h(X),
\end{equation}
and ask whether existence of a spin${}^c$ structure $\fs$ over $X$ with \emph{non-zero Seiberg--Witten invariant} $\SW_X(\fs)$ (hence the first Chern class of $\fs$ is a Seiberg--Witten \emph{basic class}) implies that $M_\kappa^w(X,g)$ is non-empty?

For now, suppose that $\su(E)$ does admit a $g$-anti-self-dual connection when inequality \eqref{eq:p1_lower_bound} holds and that the metric $g$ on $X$ is generic. Then $M_\kappa^w(X,g)$ with $w=c_1(E)$ and $\kappa=-\frac{1}{4}p_1(\su(E))$ is a non-empty, smooth manifold and so
\[
  \expdim M_\kappa^w(X,g) \geq 0,
\]
noting that in this case the actual and expected dimensions of $M_\kappa^w(X,g)$ coincide. Therefore,
\begin{align*}
  0 &\leq \frac{1}{2}\expdim M_\kappa^w(X,g)
  \\
    &= -p_1(\su(E))-3\chi_h(X) \quad\text{(by \eqref{eq:Expected_dimension_moduli_space_ASD_connections})}
  \\
    &\leq -\left(c_1(X)^2-12\chi_h(X)\right)-3\chi_h(X) \quad\text{(by \eqref{eq:p1_lower_bound})}
  \\
    &= -c_1(X)^2+9\chi_h(X),
\end{align*}
and this yields the Bogomolov--Miyaoka--Yau inequality \eqref{eq:BMY}.

According to Gompf and Stipsicz \cite[Section 1.4.1]{GompfStipsicz}, one has
\[
  c_1(X)^2-12\chi_h(X) = -e(X) = -c_2(X)
\]
and so the fundamental lower bound \eqref{eq:p1_lower_bound} is equivalent to a \emph{fundamental upper bound} on the instanton number for $E$:
\[
  \kappa(E) = -\frac{1}{4}p_1(\su(E)) \leq \frac{1}{4}e(X).
\]
We recall from Donaldson and Kronheimer \cite[Equation (2.1.39), p. 42]{DK} that
\[
  p_1(\su(E)) = c_1(E)^2 - 4c_2(E),
\]
where $c_1(E)$ and $c_2(E)$ are the first and second Chern classes of $E$, respectively.

Of course, if $M_\kappa^w(X,g)$ has expected dimension $2\delta$ and a Donaldson invariant of degree $\delta$ is non-zero, that would also imply that $M_\kappa^w(X,g)$ is non-empty. Witten's formula \cite{FL7} expresses Donaldson invariants in terms of Seiberg--Witten invariants and so this is one way in which Seiberg--Witten invariants yield information about $M_\kappa^w(X,g)$. However, this route to proving that $M_\kappa^w(X,g)$ is non-empty via Witten's formula presumes that the expected dimension of $M_\kappa^w(X,g)$ is non-negative.

For example, suppose that $X$ is a minimal, compact, complex surface of general type. According to Witten \cite[Section 4, Footnote 11, p. 589]{Witten} (see also Morgan \cite[Theorem 7.4.1, p. 122]{MorganSWNotes}), $X$ has just one Seiberg-Witten basic class up to sign, namely the canonical class, $K_X \in H^2(X;\ZZ)$. For such four-manifolds, the formula for the Donaldson series simplifies and it becomes easier to extract expressions for Donaldson invariants in low degrees: see, for example, \cite[Equation (1.2), p. 2]{FL7} for Witten's formula for the Donaldson series in terms of Seiberg--Witten invariants and basic classes first stated in \cite[Equation (2.17), p. 779]{Witten} or see Kronheimer and Mrowka \cite[Theorem 1.7, p. 577]{KMStructure} for the formula for the Donaldson series expressed in terms of Kronheimer--Mrowka basic classes. If $X$ is a standard four-manifold as in Definition \ref{defn:Standard} with one basic class $K$ up to sign and $w\cdot K = \chi_h(K)\pmod {2}$, which is equivalent to $b^+(X) = 3\pmod {4}$ when $w=0$, then Witten's formula yields
\[
  \bD_X^w(h) = 2^{K^2-\chi_h(X)+2}e^{Q_X(h)/2}\SW_X(K)\cosh(K\cdot h), \quad\text{for all } h \in H_2(X;\RR),
\]
for $w \in H^2(X;\ZZ)$: see Kronheimer and Mrowka \cite[Equation (9.7), p. 728]{KMStructure} or Fintushel and Stern \cite[Theorem 5.9, p. 619]{FSStructure}. In particular, the preceding formula implies that
\[
  \bD_X^w(0) = 2^{K^2-\chi_h(X)+2}\SW_X(K),
\]
and $\SW_X(K)$ is non-zero since $K$ is a basic class. Hence, the degree-zero Donaldson invariant $D_X^w(1) = \bD_X^w(0)$ is non-zero and is equal to the count with sign of the number of points in a (necessarily \emph{non-empty}) zero-dimensional, smooth\footnote{For a generic Riemannian metric, $g$.} moduli space $M_\kappa^w(X,g)$ of anti-self-dual connections on an $\SO(3)$ bundle $\su(E)$ over $X$ with $p_1(\su(E)) = 4\kappa$, where $E$ is a complex, rank-two vector bundle $E$ with $c_1(E)=w$. Equation \eqref{eq:Expected_dimension_moduli_space_ASD_connections} thus gives $p_1(\su(E)) = -3\chi_h(X)$. We observe that $p_1(\su(E))$ obeys the constraint \eqref{eq:p1_lower_bound} if and only if $-3\chi_h(X) \geq c_1(X)^2 - 12\chi_h(X)$, which is equivalent to \eqref{eq:BMY}, namely the Bogomolov--Miyaoka--Yau inequality that we hope to prove using gauge theory. This calculation illustrates that the formula for the Donaldson series in terms of Seiberg--Witten invariants does not necessarily provide insight into whether or not the inequality \eqref{eq:BMY} holds.

Our approach to proving that $M_\kappa^w(X,g)$ is non-empty for standard four-manifolds $X$ of Seiberg--Witten simple type under the constraint \eqref{eq:p1_lower_bound} instead relies on the virtual Morse--Bott index \eqref{eq:Virtual_Morse-Bott_index} of a critical point of a Hamiltonian function (see the forthcoming \eqref{eq:Hitchin_function}) on the moduli space of non-Abelian monopoles. We describe some of these ideas and our results in the next section.

\section[Virtual Morse--Bott index and the moduli space of non-Abelian monopoles]{Virtual Morse--Bott index and the moduli space of non-Abelian monopoles over a closed, complex K\"ahler surface}
\label{sec:Morse_theory_moduli_space_PU2_monopoles_Kaehler}
Our strategy to use compactifications of the moduli space $\sM_\ft$ of non-Abelian monopoles, defined in Chapter \ref{chap:Preliminaries},  to prove existence of an anti-self-dual  connection on $\su(E)$ satisfying \eqref{eq:p1_lower_bound} is broadly modeled on that of Hitchin's use of Morse--Bott theory to describe the topology of the moduli space of Higgs pairs over a Riemann surface in \cite[Section 7]{Hitchin_1987}, but there are two essential differences as we explain in this section.

\begin{figure}[H]
\centering
\scalebox{.6}{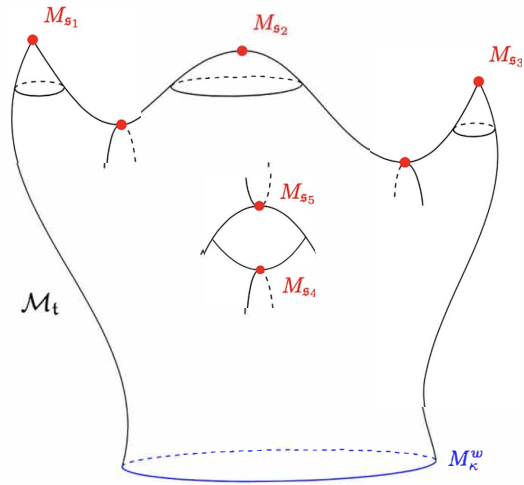}
\caption[\quad Critical sets in the non-Abelian monopole moduli space $\sM_\ft$]{Non-Abelian monopole moduli space $\sM_\ft$ with critical sets of the Hamiltonian function given by Seiberg--Witten moduli subspaces $M_{\fs_i}$ and the moduli subspace $M_\kappa^w(X,g)$ of anti-self-dual connections}
\label{fig:SO3_monopole_cobordism}
\end{figure}

We will introduce an analogue of Hitchin's Hamiltonian function, demonstrate the existence of a \spinu structure $\ft=(\rho,W\otimes E)$ with $p_1(\su(E))$ satisfying the fundamental lower bound \eqref{eq:p1_lower_bound} and whose associated non-Abelian monopole moduli space $\sM_\ft$ is non-empty, and then show that for complex K\"ahler surfaces $X$, all critical points of the Hamiltonian function on $\sM_\ft$ are given by points in $M^w_\ka(X,g) \subset \sM_\ft$ or by points in moduli subspaces $M_\fs \subset \sM_\ft$ (defined by \spinc structures $\fs$) of Seiberg--Witten monopoles with positive virtual Morse--Bott index \eqref{eq:Virtual_Morse-Bott_index}. A schematic illustration of the geometry of the moduli space $\sM_\ft$ is illustrated in Figure \ref{fig:SO3_monopole_cobordism}. Unlike the case of Higgs pairs studied by Hitchin, we must address two fundamental new difficulties in our application that do not arise in \cite{Hitchin_1987}:
\begin{enumerate}
\item\label{item:Obstructions} The strata $M_\fs \subset \sM_\ft$ of moduli subspaces of Seiberg--Witten monopoles are smooth manifolds, but not necessarily smoothly embedded as submanifolds of $\sM_\ft$.
\item\label{item:Bubbling} The moduli space $\sM_\ft$ of non-Abelian monopoles is non-compact due to energy bubbling and the Hamiltonian function \eqref{eq:Hitchin_function} is not proper. While the moduli space $\sM_\ft$ has Uhlenbeck or Gieseker compactifications as described in Chapter \ref{chap:Bubbling}, those compactified moduli spaces have singularities that are more complicated than those which can be analyzed by the methods of this monograph.
\end{enumerate}
Hitchin \cite{Hitchin_1987} assumes that the degree of the Hermitian vector bundle $E$ is odd and that its rank is two and, more generally, it is known (see Fan \cite{Fan_2022} or Gothen \cite{Gothen_1994} and references therein) that if the degree and rank of the Hermitian vector bundle $E$ in the equations \cite[Introduction, p. 60]{Hitchin_1987} for a Higgs pair are \emph{coprime}, then the linearization of those equations at a Higgs pair $(A,\Phi)$ has \emph{vanishing cokernel} and Higgs pairs $(A,\Phi)$ have trivial stabilizer subgroups in $C^\infty(\SU(E))$. However, when the degree and rank of $E$ are not coprime, one may encounter problems similar to those in the preceding Item \eqref{item:Obstructions}, where the linearization of the non-Abelian monopole equations at a Seiberg--Witten pair $(A,\Phi)$ may have a non-vanishing cokernel. (A recent article by Hitchin and Hausel \cite{Hausel_Hitchin_2022} provides generalizations of some of the results in \cite{Hitchin_1987} without an assumption that the degree and rank of $E$ are coprime.) The virtual Morse--Bott index \eqref{eq:Virtual_Morse-Bott_index} that we introduced in Feehan \cite{Feehan_analytic_spaces} may be used to extract useful information regarding the local behavior of the Hamiltonian function near a critical point, where we replace the role of the classical Morse--Bott index \cite{Bott_1954} of a smoothly embedded critical submanifold by that of the virtual Morse--Bott index \eqref{eq:Virtual_Morse-Bott_index}.

Energy bubbling in the preceding Item \eqref{item:Bubbling} is a problem that is characteristic of non-Abelian monopole equations over four-dimensional manifolds \cite{FL1} and does not occur for the moduli space of Higgs monopoles over a Riemann surface. Moreover, for Hitchin \cite{Hitchin_1987}, his analogue of the forthcoming Hamiltonian function \eqref{eq:Hitchin_function} is proper and that is certainly not true for \eqref{eq:Hitchin_function} on $\sM_\ft$ since even each level set may be non-compact due to energy bubbling: for example, the level set $f^{-1}(0)$ corresponds to the moduli subspace $M_\kappa^w(X,g)$ (if non-empty) of anti-self-dual connections whereas in \cite{Hitchin_1987}, the corresponding level set is the (compact) moduli space of flat connections over the Riemann surface. We describe some of the difficulties in our setting caused by this non-compactness and the singularities of the compactifications in Chapter \ref{chap:Bubbling}.

\subsection{An analogue of Hitchin's Hamiltonian function}
We consider the following analogue of Hitchin's Hamiltonian function in \cite[Section 7]{Hitchin_1987},
\label{page:Hitchin_Function}
\begin{equation}
  \label{eq:Hitchin_function}
  f: \sM_\ft \ni [A,\Phi] \mapsto f[A,\Phi] = \frac{1}{2}\|\Phi\|_{L^2(X)}^2 \in \RR.
\end{equation}
This function extends to the compact, smoothly stratified space of ideal non-Abelian monopoles $\sI\!\!\sM_\ft$ containing the \emph{Uhlenbeck compactification} $\bar\sM_\ft$ of $\sM_\ft$ (see Section \ref{sec:PU2Monopoles} or Feehan and Leness \cite{FL1, FL5}). The function $f$ is continuous on $\sI\!\!\sM_\ft$ and smooth but not necessarily Morse--Bott on smooth strata of $\sI\!\!\sM_\ft$. We say that a point $[A,\Phi] \in \sM_\ft$ is a \emph{critical point} of $f$ if the differential $df[A,\Phi]$ is zero on the Zariski tangent space $T_{[A,\Phi]}\sM_\ft$. (See Feehan \cite{Feehan_analytic_spaces} for further discussion of definitions of critical points of functions on singular analytic spaces.) We adapt Hitchin's paradigm in \cite[Section 7]{Hitchin_1987} for Higgs pairs and use \eqref{eq:Hitchin_function} as a Morse function to probe the topology of the moduli space $\sM_\ft$ of non-Abelian monopoles.

\subsection{Feasibility of the non-Abelian monopole cobordism method}
\label{subsec:Feasibility_SO(3)-monopole_cobordism_method}
The first step in our program to prove Conjecture \ref{conj:BMY_Seiberg-Witten} is to construct a  \spinu structure $\ft=(\rho,W\otimes E)$ (see Section \ref{sec:SpincuStr}) with $p_1(\su(E))$ satisfying the fundamental lower bound \eqref{eq:p1_lower_bound} and for which the open subspace of non-split\footnote{We use the terms `split' or `non-split' in this monograph, as given in Definition \ref{defn:Split_trivial_central-stabilizer_spinor_pair}, where in previous publications we used the terms `reducible' or `irreducible' as discussed in Remark \ref{rmk:Old_Definition_of_Reducible_Pairs}.} monopoles $\sM_\ft^* \subset \sM_\ft$ (see Section \ref{sec:PU2Monopoles}) is non-empty. To obtain greater control over the characteristic classes of the \spinu structure, we shall use the \emph{(smooth) blow-up} $\widetilde X:=X\#\overline{\CC\PP}^2$ of $X$,
\label{page:BlowUp}
where $\overline{\CC\PP}^2$ is the smooth four-manifold $\CC\PP^2$ with orientation given by the opposite of the complex orientation, and $\#$ denotes the smooth connected sum (see Gompf and Stipsicz \cite[Definition 2.2.7, p. 43]{GompfStipsicz}).  Let $e\in H^2(\widetilde X;\ZZ)$ be the Poincar\'e dual of the exceptional curve
(see Gompf and Stipsicz \cite[Definition 2.2.7, p. 43]{GompfStipsicz}).  Because $H^2(\widetilde X;\ZZ)\cong H^2(X;\ZZ)\oplus\ZZ e$ (also by \cite[Definition 2.2.7, p. 43]{GompfStipsicz}), we shall treat $H^2(X;\ZZ)$ as a submodule of $H^2(\widetilde X;\ZZ)$.

The passage from $X$ to $\widetilde X$ involves no loss in generality. Indeed, because $c_1(\widetilde X)^2=c_1(X)^2-1$, we can replace the condition \eqref{eq:p1_lower_bound} with the equivalent fundamental lower bound for $E$ over $\widetilde X$,
\begin{equation}
  \label{eq:p1_lower_bound_blow-up}
  p_1(\su(E)) \ge c_1(\widetilde X)^2+1-12\chi_h(\widetilde X).
\end{equation}
As will be explained in the proof of Theorem \ref{mainthm:ExistenceOfSpinuForFlow}, the existence of a Seiberg--Witten basic class on a manifold of simple type implies the existence of an almost complex structure on the manifold as assumed in Item \ref{item:PositiveFormalIndex} of Theorem \ref{mainthm:ExistenceOfSpinuForFlow}.

The characteristic class $c_1(\widetilde X)$, appearing in the expression \eqref{eq:FormalMorseIndexIntroThm} in the forthcoming Theorem \ref{mainthm:ExistenceOfSpinuForFlow}, depends on the choice of an almost complex structure $\tilde J$ on $\widetilde X$. The fact that different choices of an almost complex structure give different cohomology classes $c_1(\widetilde X)$ follows from Gompf and Stipsicz \cite[Theorem 1.4.14, p. 29]{GompfStipsicz}, which states that a class $\ka\in H^2(\widetilde X;\ZZ)$ satisfies $\ka=c_1(T\widetilde X,\tilde J)$ for some almost complex structure $\tilde J$ on $\widetilde X$ if and only if $\ka^2=2e(\widetilde X)+3\si(\widetilde X)$ and $\ka\equiv w_2(\widetilde X)\pmod 2$. A simple example of different choices for $c_1(\widetilde X)$ is given by considering the case where $\ka=c_1(T\widetilde X,\tilde J)$, so that $-\ka=c_1(T\widetilde X,-\tilde J)$ by Milnor and Stasheff \cite[Lemma 14.9, p. 168]{MilnorStasheff}. In Chapter \ref{chap:Feasibility}, we shall prove the

\begin{mainthm}[Feasibility of \spinu structures]
\label{mainthm:ExistenceOfSpinuForFlow}
Let $X$ be a standard four-manifold in the sense of Definition \ref{defn:Standard} with $b^-(X)\ge 2$ and Seiberg--Witten simple type. If $\fs_0$ is a spin${}^c$ structure over $X$ such that $\SW_X(\fs_0)\neq 0$, and $\widetilde X=X\#\overline{\CC\PP}^2$ denotes the smooth blow-up of $X$, and $\tilde g$ is a smooth Riemannian metric on $\widetilde X$, then there exists a Hermitian, rank-two vector bundle $E$ over $\widetilde X$ with associated spin${}^u$ structure $\tilde\ft=(\rho,W\otimes E)$ over $\widetilde X$ such the following hold:
\begin{enumerate}
\item\label{item:ExistenceOfSpinuForFlowNonEmpty}
There is an open, dense subspace of smooth Riemannian metrics in the $C^r$ topology (with $r\ge 3$) such that, if $\tilde g$ belongs to that subspace, then the moduli space $\sM^{*,0}_{\tilde\ft}$ of non-split, non-zero-section non-Abelian monopoles is non-empty.

\item\label{item:ExistenceOfSpinuForFlowp1}
The bundle $E$ over $\widetilde X$ obeys the fundamental lower bound \eqref{eq:p1_lower_bound_blow-up}.

\item\label{item:ExistenceOfACStructure}
There is an almost complex structure $\tilde J$ on $\widetilde X$ with
\[
    c_1(\widetilde X):=c_1(T\widetilde X,\tilde J)=-c_1(\fs_0)-e,
\]
where $e$ is the Poincar\'e dual of the exceptional two-sphere.

\item\label{item:PositiveFormalIndex}
For all spin${}^c$ structures $\tilde\fs$ for which $M_{\tilde\fs}$ continuously embeds in $\sM_{\tilde\ft}$, the \emph{formal Morse--Bott index} is positive:
\begin{equation}
\label{eq:FormalMorseIndexIntroThm}
\lambda^-(\tilde\ft,\tilde\fs)
:=
-2\chi_h(\widetilde X)
-\left( c_1(\tilde\fs)-c_1(\tilde\ft)\right)\cdot c_1(\widetilde X)
-\left( c_1(\tilde\fs)-c_1(\tilde\ft)\right)^2
> 0.
\end{equation}

\item\label{item:Upper_bound_expected_dimension_ASD_moduli_space_blowup}
The expected dimension of the moduli space $M_\kappa^w(\widetilde X,\tilde g)$ of anti-self-dual connections obeys the following inequality:
\begin{equation}
  \label{eq:Upper_bound_expected_dimension_ASD_moduli_space_blowup}
    \frac{1}{2}\expdim M_\kappa^w(\widetilde X,\tilde g) \leq -c_1(X)^2+9\chi_h(X).
  \end{equation}
\end{enumerate}
\end{mainthm}

When the conclusions of Theorem \ref{mainthm:ExistenceOfSpinuForFlow} are obeyed, we call $\tilde\ft$ a \emph{feasible spin${}^u$ structure}. If $X$ is a closed, complex K\"ahler surface, then it is possible to replace the condition in Item \eqref{item:ExistenceOfSpinuForFlowNonEmpty} of Theorem \ref{mainthm:ExistenceOfSpinuForFlow} that the Riemannian metric $\tilde g$ lies in an open, dense subspace of the space of all Riemannian metrics by the condition that the metric be K\"ahler and have a fundamental two-form $\omega$ as in \eqref{eq:Fundamental_two-form} and whose associated cohomology class $[\omega]\in H^{1,1}(X)\cap H^2(X;\RR)$ lies in an open, dense subspace of the K\"ahler cone. We will justify this assertion in the forthcoming Remark \ref{rmk:Proof_of_GenericKahlerMetrics_for_NoZeroPairs}.

\begin{rmk}[On the hypothesis that $b^-(X)\ge 2$]
\label{rmk:OnbMinusEquals1}
The hypothesis that $b^-(X)\ge 2$ in Theorem \ref{mainthm:ExistenceOfSpinuForFlow} is one we hope to remove in future work. We claim that a closed four-manifold $X$ with $b_1(X)=0$ satisfies $b^-(X)\leq 1$ if and only if it satisfies
\begin{equation}
\label{eq:bMinus_ge_2_Region}
 c_1(X)^2\ge 10\chi_h(X)-2
\end{equation}
The inequality \eqref{eq:bMinus_ge_2_Region} can be rewritten as $c_1(X)^2\ge 9\chi_h(X)+(\chi_h(X)-2)$.  Thus, four-manifolds satisfying \eqref{eq:bMinus_ge_2_Region} and $\chi_h(X)>2$ would violate the Bogomolov--Miyaoka--Yau inequality \eqref{eq:BMY}, namely
\[
  c_1(X)^2\le 9\chi_h(X).
\]
Because we are concerned in this monograph with standard four-manifolds that might violate the Bogomolov--Miyaoka--Yau inequality, the assumption $b^-(X)\ge 2$ in Theorem \ref{mainthm:ExistenceOfSpinuForFlow} yields the non-trivial gap defined by \eqref{eq:bMinus_ge_2_Region} in the set of standard four-manifolds to which these results apply.

To prove \eqref{eq:bMinus_ge_2_Region}, recall that $\si(X)=b^+(X)-b^-(X)$ and observe that because $b_1(X)=0$, 
\[
e(X)=2+b^+(X)+b^-(X).
\]
The definitions of $c_1(X)^2$ and $\chi_h(X)$ in \eqref{eq:DefineTopCharNumbersOf4Manifold}
thus give
\[
c_1(X)^2=4+5b^+(X)-b^-(X) \quad\text{and}\quad \chi_h(X)=\frac{1}{2}(1+b^+(X)).
\]
Hence, $X$ satisfies the equality
\begin{equation}
\label{eq:c1Squared_And_bMinus_Inequality}
c_1(X)^2-10\chi_h(X)=-1-b^-(X).
\end{equation}
Combining $b^-(X)\le 1$ with \eqref{eq:c1Squared_And_bMinus_Inequality} yields
\eqref{eq:bMinus_ge_2_Region}.  
\end{rmk}

\begin{rmk}[On the hypothesis that $b_1(X)=0$]
\label{rmk:Onb1Equals0}
The hypothesis that $X$ is a standard manifold in Theorem \ref{mainthm:ExistenceOfSpinuForFlow} includes the condition $b_1(X)=0$. We use that condition --- and the resulting inequality $\chi_h(X)>0$, where $\chi_h(X)$ is defined in \eqref{eq:DefineTopCharNumbersOf4Manifold} --- in the proof of Theorem \ref{mainthm:ExistenceOfSpinuForFlow} to show that the forthcoming constraint \eqref{eq:bGeneralTypeAssumption} does not prevent us from deriving the bound \eqref{eq:Upper_bound_expected_dimension_ASD_moduli_space_blowup} on the expected dimension of the moduli space of anti-self-dual connections. Additional details on the use of the condition $b_1(X)=0$ appear in the forthcoming Remark \ref{rmk:TechnicalRemarkOnb1=0}. It is possible that the condition $b_1(X)=0$ could be removed from Theorem \ref{mainthm:ExistenceOfSpinuForFlow}, although such arguments are likely to be even more delicate.
\end{rmk}

In the proof of the forthcoming Corollary \ref{maincor:MorseIndexAtReduciblesOnKahlerWithSO3MonopoleCharacteristicClasses}, we show that the expression defined in \eqref{eq:FormalMorseIndexIntroThm} is equal to the \emph{virtual Morse--Bott index} \eqref{eq:MorseIndexAtReduciblesOnKahlerSpinNotationType1} of the Hamiltonian function \eqref{eq:Hitchin_function} for \spinu and \spinc structures, $\ft$ and $\fs$, respectively, over a compact, complex, K\"ahler surface $X$ rather than $\tilde\ft$ and $\tilde\fs$ over a smooth four-manifold $\widetilde X$ as in Theorem \ref{mainthm:ExistenceOfSpinuForFlow}. We refer to Section \ref{subsec:SWMonopoles} for an introduction to the moduli space $M_\fs$ of Seiberg--Witten monopoles associated to a \spinc structure $\fs$ and the literature cited there for the definition of the Seiberg--Witten invariant $\SW_X(\fs)$.

\subsection{Critical points of the Hamiltonian function in the moduli space of non-Abelian monopoles over closed complex K\"ahler surfaces}
\label{subsec:CritPointsOfHitchinFunctionOnKahler}
The next step in our program to prove Conjecture \ref{conj:BMY_Seiberg-Witten} is to use Frankel's Theorem \ref{mainthm:Frankel_almost_Hermitian} to identify the critical points of the Hamiltonian function \eqref{eq:Hitchin_function} on $\sM_\ft$. The ultimate goal of our strategy is to produce a zero-section non-Abelian monopole $(A,0)$ on a \spinu structure $\ft$ (and hence a projectively anti-self-dual connection $A$ on $E$) with $p_1(\ft)$ satisfying the inequality \eqref{eq:p1_lower_bound}. Thus, we will often assume that a non-Abelian monopole $(A,\Phi)$ satisfies $\Phi\not\equiv 0$ because the existence of a non-Abelian monopole with $\Phi\equiv 0$ would immediately yield the desired result.

Suppose that $(X,g,J)$ is a complex K\"ahler surface, with K\"ahler form $\omega(\cdot,\cdot)=g(\cdot,J\cdot)$ as in \eqref{eq:Fundamental_two-form}. (In Chapter \ref{chap:Bubbling}, we outline how one might extend the proofs of our main results from K\"ahler surfaces to almost Hermitian manifolds of real dimension four.) A version of the Hitchin--Kobayashi bijection (see the forthcoming Theorems \ref{thm:Lubke_Teleman_6-3-7} and \ref{thm:Lubke_Teleman_6-3-10}) identifies the moduli space $\sM_\ft^0$ of non-zero-section, non-Abelian monopoles (of type $1$) with a moduli space $\fM^0(E,\omega)$ of non-zero-section, \emph{stable holomorphic pairs} in \eqref{eq:Moduli_space_stable_holomorphic_pairs}, so $\sM_\ft^0$ inherits the structure of a complex analytic space and the open subspace of smooth points in $\sM_\ft^0$ is a complex manifold, generalizing results due to Itoh \cite{Itoh_1983, Itoh_1985} for the complex structure on the moduli space $M_\kappa^w(X,g)$ of anti-self-dual  connections via its identification with the moduli space of stable holomorphic bundles; we extend the proofs by Itoh \cite{Itoh_1988} and Kobayashi \cite{Kobayashi_differential_geometry_complex_vector_bundles} of their results for $M_\kappa^w(X,g)$ to prove that the $L^2$ metric $\bg$ and integrable almost complex structure $\bJ$ on the open subspace $\sM_{\ft,\reg}$
\label{page:RegularPointsOfNonAbelianMonopoleModuli} of smooth points of $\sM_\ft$ define a K\"ahler form $\bomega = \bg(\cdot,\bJ\cdot)$ (see Theorem \ref{thm:Kobayashi_7-6-36_pairs}). When $g$ (and other geometric perturbation parameters in the non-Abelian monopole equations) are \emph{generic}, then $\sM_\ft^{*,0} \subset \sM_{\ft,\reg}$ (see Theorem \ref{thm:Transv}) but if $g$ is K\"ahler and thus non-generic, this observation need not hold. 

By generalizing the analysis by Hitchin \cite[Sections 6 and 7]{Hitchin_1987}, we shall see in Section \ref{sec:Marsden-Weinstein_reduction_moduli_space_SO3_monopoles_symplectic_quotient} that the function $f$ in \eqref{eq:Hitchin_function} is a Hamiltonian function for the \emph{circle action} on $\sM_\ft$ in the sense of \eqref{eq:MomentMap}, that is,
\begin{equation}
  \label{eq:Circle_action_on_moduli_space_non-abelian_monopoles_is_Hamiltonian}
  df = \iota_\bxi\bomega \quad\text{on } \sM_{\ft,\reg}^0,
\end{equation}
where $\sM_{\ft,\reg}^0 := \sM_\ft^0\cap \sM_{\ft,\reg}$ and the vector field $\bxi$ on $\sM_{\ft,\reg}^0$ is the generator of the $S^1$ action on $\sM_{\ft,\reg}^0$ given by scalar multiplication on the sections $\Phi$. Because the fundamental two-form $\bomega$ is non-degenerate, $[A,\Phi]\in \sM_{\ft,\reg}^0$ is a critical point if and only if it is a fixed point of the $S^1$ action, and from our previous work on non-Abelian monopoles \cite[Proposition 3.1, p. 86 and Lemma 3.11, p. 93]{FL2a}, we know that this occurs if and only if $(A,\Phi)$ is a \emph{split pair} with $\Phi\not\equiv 0$ (equivalent to a \emph{Seiberg--Witten monopole}) or $\Phi\equiv 0$ (equivalent to a \emph{projectively anti-self-dual  connection}). More broadly, we obtain the

\begin{mainthm}[All critical points of Hitchin's function represent either Seiberg--Witten monopoles or projectively anti-self-dual connections]
\label{mainthm:IdentifyCriticalPoints}
Let $X$ be a closed, complex K\"ahler surface and let $\ft$ be a \spinu structure over $X$ with associated moduli space $\sM_\ft$ of non-Abelian monopoles. If $[A,\Phi]\in\sM_\ft^0$ is a critical point of the Hamiltonian function $f$ in \eqref{eq:Hitchin_function}, then there exists a spin${}^c$ structure $\fs$ over $X$ such that $[A,\Phi]$ belongs to the image of the canonical embedding $M_\fs\hookrightarrow\sM_\ft$.
\end{mainthm}

Theorem \ref{mainthm:IdentifyCriticalPoints} is a corollary of the more general Theorem \ref{thm:Critical_points_Hitchin_Hamiltonian_function_moduli_space_projective_vortices}, which applies to the moduli space of projective vortices on a Hermitian vector bundle $(E,h)$ over a closed, complex K\"ahler manifold $(X,\omega)$. If $[A,\Phi]$ belongs to the open subset $\sM_{\ft,\reg}^0 \subset \sM_\ft$ of smooth, non-zero-section points, then Theorem \ref{mainthm:Frankel_almost_Hermitian} implies that $[A,\Phi]$ is a fixed point of the circle action on $\sM_\ft$ since, as discussed above, $\sM_{\ft,\reg}^0$ is a complex K\"ahler manifold. Thus, $(A,\Phi)$ must be a split pair by the forthcoming Proposition \ref{prop:FixedPointsOfS1ActionOnSpinuQuotientSpace} and hence a Seiberg--Witten monopole by Feehan and Leness \cite[Lemma 3.11, p. 93]{FL2a} (see also Section \ref{subsec:RedPU2Monopole}). More generally, if $[A,\Phi]$ is not a smooth point of $\sM_\ft^0$, then it belongs to a local virtual moduli space $\sM_\ft^\vir$ that is a finite-dimensional, $S^1$-invariant symplectic manifold containing an open neighborhood of $[A,\Phi]$ in $\sM_\ft^0$. (See Sections \ref{sec:Marsden-Weinstein_reduction_moduli_space_SO3_monopoles_symplectic_quotient} and \ref{sec:Complex_Kaehler_structure_moduli_space_projective_vortices_near_singular_points} together with Feehan \cite{Feehan_analytic_spaces} for the construction of such virtual moduli spaces.) The $S^1$ action on $\sM_\ft^\vir$
restricts to the given $S^1$ action on $\sM_\ft^0$. By construction, the tangent space to $\sM_\ft^\vir$ at $[A,\Phi]$ is the Zariski tangent space to $\sM_\ft^0$ at $[A,\Phi]$ and so $[A,\Phi]$ is a critical point of the Hamiltonian function $f$ in \eqref{eq:Hitchin_function} in the sense of Section \ref{sec:Virtual_Morse-Bott_index_Hamiltonian_function_circle_action_complex_analytic_space}. We can therefore apply Theorem \ref{mainthm:Frankel_almost_Hermitian} to conclude that $[A,\Phi]$ is a fixed point of the circle action on $\sM_\ft^\vir$ and hence a fixed point of the circle action on $\sM_\ft^0$. Consequently, $(A,\Phi)$ must again be a split pair by the forthcoming Proposition \ref{prop:FixedPointsOfS1ActionOnSpinuQuotientSpace} and hence a Seiberg--Witten monopole \cite[Lemma 3.11, p. 93]{FL2a}.

\subsection{Virtual Morse--Bott index}
\label{subsec:VMBProperty}
If a Seiberg--Witten fixed point $[A,\Phi]$ is a non-zero-section \emph{smooth} point of $\sM_\ft$ then, by arguments generalizing those of Hitchin \cite[Section 7]{Hitchin_1987}, one can apply Frankel's Theorem \cite{Frankel_1959}
(or Theorem \ref{mainthm:Frankel_almost_Hermitian} in this monograph) to conclude that $f$ is Morse--Bott at $[A,\Phi]$ and compute the Morse--Bott index of $f$ (the dimension of the maximal negative-definite subspace of $\Hess f[A,\Phi]$ on $T_{[A,\Phi]}\sM_\ft$) as the dimension of the negative-weight space $T_{[A,\Phi]}^-\sM_\ft$ for the $S^1$ action on $T_{[A,\Phi]}\sM_\ft$. The dimension of $T_{[A,\Phi]}^-\sM_\ft$ can be computed using the Hirzebruch--Riemann--Roch Index Theorem via the identification of $\sM_\ft^0$ with the moduli space $\fM_{\ps}^0(E,\omega)$ of non-zero-section, polystable holomorphic pairs. 

In \cite[Proposition 7.1, p. 92]{Hitchin_1987}, Hitchin computes the Morse--Bott indices of points that are represented by reducible Higgs pairs with non-vanishing Higgs field; they are always \emph{positive}, thus indicating that they cannot be local minima.

If the Seiberg--Witten fixed point $[A,\Phi]$ is a non-zero-section \emph{singular} point\footnote{Note that Theorem \ref{thm:Transv} does not apply to points represented by split non-Abelian monopoles, that is, Seiberg--Witten monopoles.} of $\sM_\ft$, we use the fact that $M_\fs\subset\sM_\ft$ is a submanifold of a smooth local virtual moduli space $\sM_\ft^\vir \subset \sC_\ft$ implied by the Kuranishi model given by the elliptic deformation complex \eqref{eq:SO3MonopoleDefComplex} for the non-Abelian monopole equations with harmonic spaces $\bH_{A,\Phi}^\bullet$ in \eqref{eq:H_APhi^bullet}. As noted above, the space $\sM_\ft^\vir$ is a symplectic manifold of dimension equal to that of the Zariski tangent space $T_{[A,\Phi]}\sM_\ft$, contains open neighborhoods of $[A,\Phi]$ in $\sM_{\ft,\reg}$ and $M_\fs$ as embedded symplectic submanifolds, and admits an $S^1$ action extending that on $\sM_\ft$ and preserving the symplectic structure. The set of fixed points of the $S^1$ action on $\sM_\ft^\vir$ coincides with the open neighborhood of $[A,\Phi]$ in $M_\fs$ and $f$ is Morse--Bott on $\sM_\ft^\vir$ with critical submanifold given by an open neighborhood of $[A,\Phi]$ in $M_\fs$. (In the simpler setting of Hitchin \cite[Section 7]{Hitchin_1987} --- see also Bradlow and Wilkin \cite[Section 4]{Bradlow_Wilkin_2012}, the critical sets are smooth submanifolds of the moduli space of Higgs pairs and $f$ is Morse--Bott.) Here, $M_\fs$ is a smooth manifold, given by a moduli space of Seiberg--Witten monopoles that may be zero or positive-dimensional (in the latter case $\SW_X(\fs)=0$ when $X$ has Seiberg--Witten simple type), but not necessarily an embedded smooth submanifold of $\sM_\ft$. When $\Phi\not\equiv 0$, an analogue of the virtual Morse--Bott signature \eqref{eq:Virtual_Morse-Bott_signature} of $f$ at $[A,\Phi]$ is given by
\begin{subequations}
\label{eq:Virtual_Morse-Bott_signature_moduli_space_non-abelian_monopoles}  
\begin{align}
  \label{eq:Virtual_Morse-Bott_nullity_moduli_space_non-abelian_monopoles}
  \lambda_{[A,\Phi]}^0(f) &:= \dim_\RR \bH_{A,\Phi}^{0,1} - \dim_\RR \bH_{A,\Phi}^{0,2},
  \\
  \label{eq:Virtual_Morse-Bott_co-index_moduli_space_non-abelian_monopoles}
  \lambda_{[A,\Phi]}^+(f) &:= \dim_\RR \bH_{A,\Phi}^{+,1} - \dim_\RR \bH_{A,\Phi}^{+,2},
  \\
  \label{eq:Virtual_Morse-Bott_index_moduli_space_non-abelian_monopoles}
  \lambda_{[A,\Phi]}^-(f) &:= \dim_\RR \bH_{A,\Phi}^{-,1} - \dim_\RR \bH_{A,\Phi}^{-,2}.
\end{align}
\end{subequations}
In principle, the virtual Morse--Bott index \eqref{eq:Virtual_Morse-Bott_index_moduli_space_non-abelian_monopoles} could be defined as minus the real Euler characteristic of a \emph{negative-weight elliptic subcomplex} (with $\bH_{A,\Phi}^{-,0}=(0)$, since $\bH_{A,\Phi}^0=(0)$ when $\Phi\not\equiv 0$), given by a subcomplex of the elliptic deformation complex for the non-Abelian monopole equations defined by a \emph{negative-weight} subspaces of the $S^1$ action induced on the non-Abelian monopole elliptic deformation complex \eqref{eq:SO3MonopoleDefComplex} by the $S^1$ action in Definition \ref{defn:UnitaryZActionOnAffine}. However, identifying the signs of the weights of the $S^1$ action induced on the elliptic complex \eqref{eq:SO3MonopoleDefComplex} requires an understanding of almost complex structures on the harmonic spaces in this elliptic complex, as discussed in Remark \ref{rmk:UseOfChoiceOfComplexStructure}.

We shall instead compute the indices \eqref{eq:Virtual_Morse-Bott_signature_moduli_space_non-abelian_monopoles} by replacing\footnote{As explained in Remark \ref{rmk:Decoupling_Dirac_equation_and_equivalence_type_I_and_II_solutions}, the restriction to non-Abelian monopoles of type $1$ entails no loss of generality.} the moduli space $\sM_{\ft,1}^0$ of non-zero-section, non-Abelian monopoles of \emph{type $1$} with its image under the Hitchin--Kobayashi correspondence given by the moduli space $\fM_{\ps}^0(E,\omega)$ of non-zero-section, polystable holomorphic pairs in \eqref{eq:Moduli_space_non-zero-section_polystable_holomorphic_pairs}. Because the Hitchin--Kobayashi correspondence in Theorem \ref{thm:HitchinKobayashiCorrespondenceForPairs} is $S^1$-equivariant, we find that $[A,\Phi]\in\sM_{\ft,1}^0$ is a fixed point of the $S^1$ action on $\sM_{\ft,1}^0$ if and only if the corresponding polystable holomorphic pair $[\bar\partial_A,\varphi]$ is a fixed point of the $S^1$ action on $\fM_{\ps}^0(E,\omega)$ induced by the natural $\CC^*$ action $\rho_\CC$ on $\fM_{\ps}^0(E,\omega)$. For such fixed points, the $S^1$ action on $\sM_{\ft,1}^0$ defines an $S^1$ action on the harmonic spaces $\bH_{A,\Phi}^\bullet$ and the $S^1$ action on the quotient space
$\fM(E)$ of holomorphic pairs defines an $S^1$ action on the harmonic spaces $\bH_{\bar\partial_A,\varphi}^\bullet$ for the elliptic complex \eqref{eq:Holomorphic_pair_elliptic_complex} defined by the pair $(\bar\partial_A,\varphi)$. In Chapters \ref{chap:Elliptic_deformation_complex_moduli_space_SO(3)_monopoles_over_almost_Hermitian_four-manifold}
and \ref{chap:VirtualMorseIndexComputation}, we prove that there is an $S^1$-equivariant isomorphism between the harmonic spaces $\bH_{A,\Phi}^\bullet$ and $\bH_{\bar\partial_A,\varphi}^\bullet$ when $X$ has its canonical \spinc structure (see Definition \ref{defn:Canonical_spinc_bundles}) and $\Phi = (\varphi,0) \in \Omega^0(E)\oplus\Omega^{0,2}(E)$ with
$\varphi\not\equiv 0$. We define an analogue of the virtual Morse--Bott signature \eqref{eq:Virtual_Morse-Bott_signature} of $f$, namely
\begin{subequations}
\label{eq:Virtual_Morse-Bott_signature_moduli_space_stable_holomorphic_pairs}  
\begin{align}
  \label{eq:Virtual_Morse-Bott_nullity_moduli_space_stable_holomorphic_pairs}
  \lambda_{[\bar\partial_A,\varphi]}^0(\rho_\CC) := \dim_\CC \bH_{\bar\partial_A,\varphi}^{0,1} - \dim_\CC \bH_{\bar\partial_A,\varphi}^{0,2},
  \\
  \label{eq:Virtual_Morse-Bott_co-index_moduli_space_stable_holomorphic_pairs}
  \lambda_{[\bar\partial_A,\varphi]}^0(\rho_\CC) := \dim_\CC \bH_{\bar\partial_A,\varphi}^{+,1} - \dim_\CC \bH_{\bar\partial_A,\varphi}^{+,2},
  \\
  \label{eq:Virtual_Morse-Bott_index_moduli_space_stable_holomorphic_pairs}
  \lambda_{[\bar\partial_A,\varphi]}^-(\rho_\CC) := \dim_\CC \bH_{\bar\partial_A,\varphi}^{-,1} - \dim_\CC \bH_{\bar\partial_A,\varphi}^{-,2}.
\end{align}
\end{subequations}
We can compute the virtual Morse--Bott index \eqref{eq:Virtual_Morse-Bott_index_moduli_space_stable_holomorphic_pairs} by computing the Euler characteristic of the subcomplex of the elliptic complex \eqref{eq:Holomorphic_pair_elliptic_complex} for the pair $(\bar\partial_A,\varphi)$ with harmonic spaces $\bH_{\bar\partial_A,\varphi}^{-,\bullet}$ on which the $S^1$ action has negative weight, noting that $\bH_{\bar\partial_A,\varphi}^{-,0} = (0)$ and $\bH_{\bar\partial_A,\varphi}^{-,3} = (0)$. The almost complex structure on the vector spaces $\bH_{\bar\partial_A,\varphi}^{-,k}$ (for $k=0,1,2,3$) underlying the elliptic complex \eqref{eq:Holomorphic_pair_elliptic_complex} for the pair $(\bar\partial_A,\varphi)$ is immediately apparent, making it easy to identify the subcomplex on which $S^1$ acts with negative weight, which we do in Section \ref{sec:NormalTangentialSplittingOfDefComplex}, thus facilitating computation of the virtual Morse--Bott index $\lambda_{[A,\Phi]}^-(f)$ using the $S^1$-equivariant, canonical isomorphisms of complex vector spaces,
\[
  \bH_{A,\Phi}^{-,1} \cong \bH_{\bar\partial_A,\varphi}^{-,1}
  \quad\text{and}\quad
  \bH_{A,\Phi}^{-,2} \cong \bH_{\bar\partial_A,\varphi}^{-,2},
\]
noting that the almost complex structures on $\bH_{A,\Phi}^{-,k}$ are defined by those on $\bH_{\bar\partial_A,\varphi}^{-,k}$ for $k=1,2$.

\begin{mainthm}[Virtual Morse--Bott signature of the Hamiltonian function at a non-zero-section, split non-Abelian monopole]
\label{mainthm:MorseIndexAtReduciblesOnKahler}
Let $(\rho_{\can},W_{\can})$ be the canonical \spinc structure of Definition \ref{defn:Canonical_spinc_bundles} over a closed, connected, complex K\"ahler surface $X$, and $E$ be a rank-two Hermitian vector bundle over $X$ that admits a splitting $E=L_1\oplus L_2$ as a direct sum of Hermitian line bundles, and $\ft = (\rho,W_{\can}\otimes E)$ be the corresponding \spinu structure. Assume that $(A,\Phi)$ is a non-zero-section,
type $1$ non-Abelian monopole on $\ft$ that is split in the sense of the forthcoming Definition \ref{defn:Split_trivial_central-stabilizer_spinor_pair} with respect to the decomposition $E=L_1\oplus L_2$, with $\Phi=(\Phi_1,0)$ and $\Phi_1\in\Om^0(W_{\can}^+\otimes L_1)$. Then the virtual Morse--Bott index
\eqref{eq:Virtual_Morse-Bott_index_moduli_space_non-abelian_monopoles} of the Hamiltonian function $f$ in \eqref{eq:Hitchin_function} at the point $[A,\Phi] \in \sM_\ft$ is given by minus twice the Euler characteristic of the negative-weight elliptic complex \eqref{eq:NegativeWeightPairDefComplex} for the holomorphic pair $(\bar\partial_A,\varphi)$, where $\Phi_1=(\varphi,0)\in\Om^0(L_1)\oplus\Om^{0,2}(L_1)$, and equals
\begin{equation}
\label{eq:MorseIndexAtReduciblesOnKahler}
\lambda_{[A,\Phi]}^-(f)
=
-2\chi_h(X)
-\left(  c_1(L_1)-c_1(L_2)\right)\cdot c_1(X)
-\left( c_1(L_1)-c_1(L_2)\right)^2,
\end{equation}
where $\chi_h(X)$ in \eqref{eq:DefineTopCharNumbersOf4Manifold} is the holomorphic Euler characteristic of $X$. The
virtual Morse--Bott co-index of $f$ at $[A,\Phi]$ in \eqref{eq:Virtual_Morse-Bott_co-index_moduli_space_non-abelian_monopoles} is given by
\begin{equation}
\label{eq:MorseCoIndexAtReduciblesOnKahler}
\begin{aligned}
  \lambda_{[A,\Phi]}^+(f)
  &=
  -\left( c_1(L_2)-c_1(L_1)\right)\cdot c_1(X) - \left( c_1(L_2)-c_1(L_1)\right)^2
  \\
  &\qquad +c_1(L_2)\cdot c_1(X) +c_1(L_2)^2,
\end{aligned}
\end{equation}
and equals minus twice the Euler characteristic of the positive-weight elliptic complex \eqref{eq:PositiveWeightPairDefComplex}. The
virtual Morse--Bott nullity of $f$ at $[A,\Phi]$ in \eqref{eq:Virtual_Morse-Bott_nullity_moduli_space_non-abelian_monopoles} is given by
\begin{equation}
\label{eq:MorseNullIndexAtReduciblesOnKahler}
\lambda_{[A,\Phi]}^0(f) = c_1(L_1)^2 + c_1(L_1)\cdot c_1(X),
\end{equation}
and equals minus twice the Euler characteristic of the zero-weight elliptic complex \eqref{eq:ZeroWeightPairDefComplex}.
\end{mainthm}

\begin{rmk}[Invariance of the virtual signature and virtual Morse--Bott index with respect to circle-equivariant isomorphisms of K\"ahler spaces]
\label{rmk:Invariance_virtual_signature}
In the forthcoming Corollary \ref{cor:Invariance_signature_Morse-Bott_function}, we recall the well-known fact that the signature of the Hessian bilinear form of a smooth function at a critical point in a smooth manifold is preserved by diffeomorphisms. A generalization of this result should hold for a suitably-defined virtual signature of the Hessian bilinear form of a Hamiltonian function for a circle action at a critical point in a K\"ahler space, in the sense that the virtual signature (and virtual Morse--Bott index in particular) is preserved by circle-equivariant isomorphisms of K\"ahler spaces. We refer to Moishezon \cite{Moishezon_1975} for a definition of K\"ahler spaces.
\end{rmk}

\begin{rmk}[Hamiltonian circle action on the moduli space of stable, holomorphic pairs]
\label{rmk:Hamiltonian_circle_action_moduli_space_stable_holomorphic_pairs} 
We noted in \eqref{eq:Circle_action_on_moduli_space_non-abelian_monopoles_is_Hamiltonian} that the natural circle action on the moduli space $\sM_\ft^0$ of non-Abelian monopoles $(A,\Phi)$ on a \spinu structure $\ft = (\rho,W\otimes E)$ over a standard Riemannian four-manifold $(X,g)$ induced by complex scalar multiplication on sections $\Phi$ of $W^+\otimes E$ has the Hamiltonian function $f$ defined in \eqref{eq:Hitchin_function}, at least upon restriction to the open subset of smooth points of $\sM_\ft^0$. It is natural to ask whether the circle action induced by the standard $\CC^*$ action on the moduli space $\fM_\ps^0(E,\omega)$ of non-zero-section, polystable holomorphic pairs $(\bar\partial_A,\varphi)$ on a complex vector bundle $E$ over a complex, K\"ahler surface $(X,\omega)$ is also Hamiltonian, at least upon restriction to the subset of smooth points of $\fM_\ps^0(E,\omega)$. If $E$ is equipped with a Hermitian metric $h$, then $\|\varphi\|_{L^2(X)}$ is invariant under the action of the group $C^\infty(\SU(E))$ of determinant-one, \emph{unitary} automorphisms of $E$ but \emph{not} invariant under the action of the group $C^\infty(\SL(E))$ of determinant-one, \emph{complex} automorphisms of $E$. Hence, the obvious analogue of \eqref{eq:Hitchin_function} on $\fM_\ps^0(E,\omega)$ is not well-defined. 

However, the circle action on $\fM_\ps^0(E,\omega)$ can be shown to be Hamiltonian (see Section \ref{subsec:Characterization_critical_points_Hamiltonian_function_circle_action_on_projective_vortices}) by using the circle-equivariant Hitchin--Kobayashi correspondence in Theorem \ref{thm:Lubke_Teleman_6-3-10} between the moduli space $\sM_{\ft,1}^0$ of non-zero-section, type $1$ non-Abelian monopoles and $\fM_{\ps}^0(E,\omega)$. The pull back of the function $f$ in \eqref{eq:Hitchin_function} from $\sM_{\ft,1}^0$ to $\fM_{\ps}^0(E,\omega)$ serves as a Hamiltonian function for the circle action on $\fM_{\ps}^0(E,\omega)$.
\end{rmk}

For reasons explained in the forthcoming Remark \ref{rmk:Virtual_Morse-Bott_index_formula_Kaehler_surface_versus_smooth_4-manifold}, we expect the formula \eqref{eq:MorseIndexAtReduciblesOnKahler} for the virtual Morse--Bott index  in the conclusion of Theorem \ref{mainthm:MorseIndexAtReduciblesOnKahler} to continue to hold for standard four-manifolds.

The Hitchin--Kobayashi correspondence maps $[A,\Phi]\in M_\fs\subset\sM_\ft^0$ to a non-zero-section, polystable holomorphic pair which is split with respect to the decomposition $E=L_1\oplus L_2$, so by expressing the characteristic classes of $L_1$ and $L_2$ in terms of those of the \spinc structure $\fs$ and the \spinu structure $\ft$ we can rewrite
Theorem \ref{mainthm:MorseIndexAtReduciblesOnKahler} as

\begin{maincor}[Virtual Morse--Bott index of the Hamiltonian function at a point represented by a Seiberg--Witten monopole]
\label{maincor:MorseIndexAtReduciblesOnKahlerWithSO3MonopoleCharacteristicClasses}
Let $X$ be a closed, connected, complex K\"ahler surface, $\ft$ be a spin${}^u$ structure over $X$, and $\fs$ be a spin${}^c$ structure over $X$. If $[A,\Phi] \in \sM_\ft$ is a point represented by a non-zero-section, split, type $1$ non-Abelian monopole in the image of the embedding defined in \eqref{eq:DefnOfIotaOnQuotient} of the moduli space $M_\fs$ of Seiberg--Witten monopoles on $\fs$ into the moduli space $\sM_\ft$ of non-Abelian monopoles on $\ft$, then the virtual Morse--Bott index
\eqref{eq:Virtual_Morse-Bott_index_moduli_space_non-abelian_monopoles} of the Hamiltonian function $f$ in \eqref{eq:Hitchin_function} on the moduli space $\sM_\ft$ at $[A,\Phi]$ is given by
\begin{equation}
\label{eq:MorseIndexAtReduciblesOnKahlerSpinNotationType1}
\lambda_{[A,\Phi]}^-(f)
=
-2\chi_h(X)
-\left( c_1(\fs)-c_1(\ft)\right)\cdot c_1(X)
-\left( c_1(\fs)-c_1(\ft)\right)^2,
\end{equation}
where $\chi_h(X)$ is the holomorphic Euler characteristic \eqref{eq:DefineTopCharNumbersOf4Manifold} of $X$, and $c_1(\fs) = c_1(W^+)\in H^2(X;\ZZ)$ is the first Chern class \eqref{eq:DefineChernClassOfSpinc} of the spin${}^c$ structure given by $\fs = (\rho,W)$ with $W=W^+\oplus W^-$, and $c_1(\ft)=c_1(E)+c_1(W^+) \in H^2(X;\ZZ)$ is the first Chern class \eqref{eq:CliffordDeterminantBundle} of the spin${}^u$ structure.
\end{maincor}

If we apply the expression \eqref{eq:MorseIndexAtReduciblesOnKahlerSpinNotationType1} for the virtual Morse--Bott index to the moduli space $\sM_\ft$ described in Theorem \ref{mainthm:ExistenceOfSpinuForFlow}, we obtain the

\begin{maincor}[Positivity of virtual Morse--Bott index of the Hamiltonian function at a point represented by a Seiberg--Witten monopole for a feasible \spinu structure]
\label{maincor:MorseIndexForFeasibilitySpinuStructure}
Let $X$ be a closed, connected, complex K\"ahler surface with $b_1(X)=0$, and $b^-(X)\ge 2$, and $b^+(X)\ge 3$. If $\widetilde X=X\#\overline{\CC\PP}^2$ is the blow-up of $X$ and $\tilde\ft$ is the \spinu structure on $\tilde X$ constructed in Theorem \ref{mainthm:ExistenceOfSpinuForFlow}, then for all non-empty Seiberg--Witten moduli subspaces $M_{\tilde\fs}$ that are continuously embedded in $\sM_{\tilde \ft}$ as non-zero-section, type $1$, split, non-Abelian monopoles, the virtual Morse--Bott index
\eqref{eq:Virtual_Morse-Bott_index_moduli_space_non-abelian_monopoles} of the Hamiltonian function $f$ in \eqref{eq:Hitchin_function} on $\sM_{\tilde \ft}$ is \emph{positive} at all points in $M_{\tilde\fs}$.
\end{maincor}

\begin{rmk}[Sufficiency of type $1$ non-Abelian monopoles]
\label{rmk:DisjointBranchesAndTheFlow}
Although Corollary \ref{maincor:MorseIndexForFeasibilitySpinuStructure} is only stated for type $1$ split, non-Abelian monopoles, a similar result holds for type $2$ split, non-Abelian monopoles. However, as described by Teleman in \cite[Corollary 2.3.3]{TelemanNonabelian} and Dowker in \cite[Theorem 1.4.3]{DowkerThesis}, the moduli space $\sM_\ft$ of  non-Abelian monopoles on a K\"ahler surface contains gauge-equivalence classes of type $1$ and type $2$ pairs as closed subspaces whose intersection is the moduli space of projectively anti-self-dual connections. Thus, if the initial point of the gradient flow of the Hamiltonian function \eqref{eq:Hitchin_function} is in the subspace of points represented by type $1$ pairs, the limit of the flow must be the gauge-equivalence class of a type $1$ pair or a projectively anti-self-dual connection.  Hence, when considering non-Abelian monopoles on a K\"ahler surface, it is not necessary for Corollary  \ref{maincor:MorseIndexForFeasibilitySpinuStructure} to apply to type $2$ non-Abelian monopoles.
\end{rmk}

\begin{rmk}[Extension of the virtual Morse-Bott index formula from complex, K\"ahler surfaces to standard four-manifolds of Seiberg--Witten simple type]
\label{rmk:Virtual_Morse-Bott_index_formula_Kaehler_surface_versus_smooth_4-manifold}
In Theorem \ref{mainthm:MorseIndexAtReduciblesOnKahler} and Corollaries \ref{maincor:MorseIndexAtReduciblesOnKahlerWithSO3MonopoleCharacteristicClasses} and \ref{maincor:MorseIndexForFeasibilitySpinuStructure}, we assume that $X$ is a closed, connected, complex K\"ahler surface. However, we conjecture that Theorem \ref{mainthm:MorseIndexAtReduciblesOnKahler} should extend to standard four-manifolds of Seiberg--Witten simple type. The cohomology calculations used to prove Theorem \ref{mainthm:ExistenceOfSpinuForFlow} do not rely on $X$ being complex or K\"ahler and the corresponding elliptic complexes obtained when closed complex K\"ahler surfaces are replaced by standard four-manifolds are perturbations by compact operators of the elliptic complexes obtained for closed, complex K\"ahler surfaces. Hence, the conclusion for the virtual Morse--Bott index of the Hamiltonian function $f$ in Corollary \ref{maincor:MorseIndexForFeasibilitySpinuStructure} should hold without the hypothesis that $X$ is complex or K\"ahler, provided one can show that the moduli space of non-Abelian monopoles is almost Hermitian.
\end{rmk}

\section{Outline}
\label{sec:Outline}
To motivate our generalization of Frankel's Theorem for almost symplectic manifolds (Theorem \ref{mainthm:Frankel_almost_Hermitian}), we present in Chapter \ref{chap:Hessian_restriction_smooth_function_submanifold_Euclidean_space} a computation of the Hessian of a function restricted to an isometrically embedded submanifold of $\RR^n$, illustrating the difficulties
in any attempt to directly calculate Morse--Bott indices in this monograph without appealing to Frankel's Theorem. In Chapter \ref{chap:Circle_actions_almost_Hermitian_manifolds}, we prove Theorem  \ref{mainthm:Frankel_almost_Hermitian}. In Chapter \ref{chap:Local_indices_vector_fields_manifolds_analytic_spaces}, we describe definitions of the index of a vector field on an analytic space that have appeared in the literature. Chapter \ref{chap:Preliminaries} provides a review of the definitions and basic properties of the moduli spaces of anti-self-dual connections, Seiberg--Witten monopoles, and non-Abelian monopoles over smooth Riemannian four-manifolds. With that background in place, we  prove Theorem \ref{mainthm:ExistenceOfSpinuForFlow} in Chapter \ref{chap:Feasibility}.

We review the relationship between non-Abelian monopoles, projective vortices, and holomorphic pairs over closed, complex K\"ahler surfaces in Chapter \ref{chap:SO(3)monopolesoverKahler surfaces}. The Hitchin--Kobayashi correspondence provides an isomorphism in the sense of real analytic spaces between the moduli space of (type $1$) non-Abelian monopoles and the moduli space of stable, holomorphic pairs over closed, complex K\"ahler surfaces. We begin our exposition of the proof of a more general version of this result in Chapter \ref{chap:Elliptic_complexes} by describing the elliptic complexes and their associated cohomology groups and harmonic spaces that we shall need in this monograph. In Chapter \ref{chap:Elliptic_deformation_complex_moduli_space_SO(3)_monopoles_over_almost_Hermitian_four-manifold}, we compare the cohomology groups of the elliptic complexes for
\begin{inparaenum}[\itshape i\upshape)]
\item projectively Hermitian--Einstein connections and holomorphic structures over complex, K\"ahler manifolds,
\item projective vortices and holomorphic pairs over complex, K\"ahler manifolds, and
\item non-Abelian monopoles and holomorphic pairs over complex, K\"ahler surfaces.
\end{inparaenum}
In order to apply Frankel's Theorem \ref{mainthm:Frankel_almost_Hermitian} to the moduli space of non-Abelian monopoles over a closed, complex K\"ahler surface, we need to show that each non-zero-section point in the latter moduli space has
an $S^1$-invariant open neighborhood that is contained in a local virtual moduli space equipped with an $S^1$-invariant symplectic or complex K\"ahler structure. Thus, in Chapter \ref{chap:Complex_structure_and_Kaehler_metric_moduli_space_HE_connections_complex_Kaehler_manifold} we first review the proof of the existence of complex K\"ahler structures on moduli spaces in the simpler case of the moduli space of projectively Hermitian--Einstein connections over a closed, complex K\"ahler manifold and hence, as a special case, the moduli space of anti-self-dual connections over a closed, complex K\"ahler surface. In Chapter \ref{chap:Complex_structure_and_Kaehler_metric_moduli_space_SO3_monopoles_complex_Kaehler_surface}, we extend those results to the moduli space of projective vortices over a closed, complex K\"ahler manifold and hence, as a special case, the moduli space of non-Abelian monopoles over a closed, complex K\"ahler surface. In Chapter \ref{chap:Analytic_isomorphism_moduli_spaces}, we give an exposition of the proof that the Hitchin--Kobayashi correspondence between projective vortices and holomorphic pairs gives an isomorphism in the sense of real analytic spaces from the moduli space of projective vortices onto the moduli space of polystable holomorphic pairs over a closed, complex K\"ahler manifold and hence, as a special case, an isomorphism in the sense of real analytic spaces from the moduli space of non-Abelian monopoles onto the moduli space of polystable holomorphic pairs over a closed, complex K\"ahler surface.

Chapter \ref{chap:VirtualMorseIndexComputation} contains the proof of Theorem \ref{mainthm:MorseIndexAtReduciblesOnKahler}, which gives the virtual Morse--Bott signature of the Hamiltonian function \eqref{eq:Hitchin_function} at split pairs (corresponding to Seiberg--Witten monopoles over K\"ahler surfaces) in the moduli space of polystable holomorphic pairs (corresponding to non-Abelian monopoles over K\"ahler surfaces). We conclude in Chapter \ref{chap:Bubbling} by discussing potential generalizations of these results to almost Hermitian standard four-manifolds of Seiberg--Witten simple type and some approaches to computing the virtual Morse--Bott index of critical points in the lower levels of the Uhlenbeck compactification \cite[Section 4.5.2, p. 351]{FL1} of the moduli space of non-Abelian monopoles.

We note here a few conventions that apply throughout this monograph unless stated otherwise in specific instances. The base manifolds $X$ will be assumed to be connected, but we will only mention this in statements of results where this property is particularly relevant. We write $C^\infty(V)$ for the Fr\'echet space of smooth sections of a smooth vector bundle $V$ over a smooth manifold. We let $\NN$ denote the set $\{1,2,3,\ldots,\}$ of natural numbers, with zero excluded. 

\chapter{Hessian of the restriction of a smooth function to a submanifold}
\label{chap:Hessian_restriction_smooth_function_submanifold_Euclidean_space}
In the simplest application of the results of this monograph (including those of Hitchin in \cite{Hitchin_1987}), we consider Morse theory for an finite-dimensional, smoothly embedded submanifold $M$ of a smooth Banach manifold $B$ in a setting to which Frankel's Theorem \ref{mainthm:Frankel_almost_Hermitian} applies, where the manifold $B$ is equipped with a circle action $\rho$ and a circle-invariant, non-degenerate two-form $\omega$ that restrict to $M$, and an associated Hamiltonian function $H$ on $B$ that restricts to a Hamiltonian function $h$ on $M$ and serves as a Morse function on $M$. In this situation, Frankel's Theorem \ref{mainthm:Frankel_almost_Hermitian} can be used to compute the signature of the Hessian bilinear form $\hess h(p)$ defined by the Hamiltonian function $h$ at a critical point $p$ (equivalently, a fixed point of the circle action) and if in addition $B$ is equipped with circle-invariant Riemannian metric $g$, then Frankel's Theorem \ref{mainthm:Frankel_almost_Hermitian} can be used to determine the eigenvalues of the Hessian operator $\Hess_gh(p)$ in terms of the weights of the circle action at the fixed point $p$ (equivalently, a critical point of the Hamiltonian function).

Because our finite-dimensional submanifold $M$ can be viewed as the zero locus of a simple, explicit section $F$ of a Banach vector bundle $V$ over $B$ and the Morse function $h$ on $M$ is the restriction of a simple, explicit function $H$ on $B$, it is interesting to ask whether the signature of the Hessian form $\hess_gh(p)$ or the eigenvalues of the Hessian operator $\Hess_gh(p)$ can be computed directly, without appealing to Frankel's Theorem \ref{mainthm:Frankel_almost_Hermitian}. As we shall explain in this chapter, such direct calculations appear to be intractable and these observations nicely illustrate the computational power of Frankel's Theorem \ref{mainthm:Frankel_almost_Hermitian} in indirect calculations of the signature or eigenvalues.

\section{Introduction}
\label{sec:Introduction_Hessian_restriction_smooth_function_submanifold_Euclidean_space}
Following Palais and Terng \cite[Section 9.4]{PalaisTerng}, let $M$ be a smooth manifold, $h:M\to\RR$ be a smooth function, and define its $\RR$-bilinear \emph{Hessian form},
\[
  \hess h:C^\infty(TX)\times C^\infty(TX) \to C^\infty(X),
\]
characterized by
\begin{equation}
  \label{eq:Hessian_bilinear_form}
  \hess h(X,Y) := X(Yh), \quad\text{for all } X, Y \in C^\infty(TM).
\end{equation}
For any point $p \in M$, we see that $\hess h(p)(X,Y) = X_p(Yh)$ depends on $X$ only through its value, $X_p$, at $p$, while if $p$ is a \emph{critical point} of $h$, then
\[
  \hess h(p)(X,Y) - \hess h(p)(Y,X) = X_p(Yh) - Y_p(Xh) = [X,Y]_ph = dh(p)[X,Y] = 0.
\]
Therefore, $\hess h(p)(X,Y) = \hess h(p)(Y,X)$ when $p$ is a critical point of $h$ and, in this case, $\hess h(p)(X,Y)$ depends on the vector fields $X$ and $Y$ only through their values, $X_p$ and $Y_p$, at $p$. Thus, when $p$ is a critical point of $h$, the bilinear form
\[
  \hess h(p) \in \Hom(T_pM\otimes T_pM,\RR)
\]
is \emph{symmetric} and we may compute its values by
\[
  \hess h(p)(X_p,Y_p) = X_p(Yh), \quad\text{for all } X_p, Y_p \in T_pM,
\]
where $Y$ is now any smooth vector field on an open neighborhood $p$ that extends the tangent vector $Y_p$ (see Palais and Terng \cite[Theorem 9.4.3]{PalaisTerng}). If $M$ is equipped with a smooth Riemannian metric $g$, we obtain the $\RR$-linear \emph{Hessian operator},
\[
  \Hess_gh:C^\infty(TX)\to C^\infty(TX),
\]  
characterized by
\begin{equation}
  \label{eq:Hessian_operator}
  \left\langle (\Hess_g h)(X), Y\right\rangle_g := X(Yh), \quad\text{for all } X, Y \in C^\infty(TM).
\end{equation}
If $p$ is a critical point of $h$, then $\Hess_g h(p)$ depends on the vector field $X$ through its value $X_p$ at $p$, and which we emphasize by writing $\Hess_g h(p)X_p \in T_pM$. Moreover, $\Hess_g h(p) \in \End(T_pM)$ is \emph{self-adjoint}. (We reverse the notation of Palais and Terng \cite[p. 197]{PalaisTerng}, writing $\hess h(p) \in \Hom(T_pM\otimes T_pM,\RR)$ for the Hessian bilinear form and $\Hess_g h(p) \in \End(T_pM)$ for the Hessian operator defined by $\hess h(p)$ and a Riemannian metric $g$ on $T_pM$.)

If $\nabla$ is any connection on $TM$, one may also define (see Palais and Terng \cite[p. 198]{PalaisTerng})
\begin{equation}
  \label{eq:Hessian_covariant_derivative_bilinear_form}
  \hess^\nabla h(X,Y) := (\nabla dh)(X,Y) , \quad\text{for all } X, Y \in C^\infty(TM).
\end{equation}
This expression may be usefully rewritten as \cite[p. 198]{PalaisTerng}
\begin{equation}
  \label{eq:Hessian_covariant_derivative_bilinear_form_alternative}
  \hess^\nabla h(X,Y) = X(Yh) - dh(\nabla_XY).
\end{equation}
Consequently, if $p$ is a critical point of $h$, then
\begin{equation}
  \label{eq:Hessian_covariant_derivative_is_basic_bilinear_form}
  \hess^\nabla h(p)(X_p,Y_p) = X_p(Yh) = \hess h(p)(X_p,Y_p) , \quad\text{for all } X_p, Y_p \in T_pM,
\end{equation}
where $Y$ is now any smooth vector field on an open neighborhood $p$ that extends the tangent vector $Y_p$ (see Palais and Terng \cite[Theorem 9.4.4 and Corollary 9.4.5]{PalaisTerng}). If $\nabla$ is a symmetric connection --- that is, a connection with zero torsion, such as the Levi--Civita connection for a Riemannian metric, the expression \eqref{eq:Hessian_covariant_derivative_bilinear_form_alternative} implies that $\hess^\nabla h$ is a covariant, symmetric two-tensor field on $M$ (see Palais and Terng \cite[Theorem 9.4.4]{PalaisTerng}).

\begin{lem}[Expression for the Hessian operator as the Lie derivative with respect to the gradient vector field]
\label{lem:HessianIsBracket}
Let $(M,g)$ be a smooth Riemannian manifold. If $h:M\to\RR$ is a smooth function with gradient vector field $\grad_g h$ on $M$ as in \eqref{eq:DefineGradient} and $p\in M$ is a critical point, so $dh(p)=0$, then for all smooth vector fields $X$ on $M$,
\begin{equation}
  \label{eq:GradientDefnOfHessian}
  \Hess_g h(p)X_p = -[\grad_g h, X]_p.
\end{equation}
\end{lem}

\begin{proof}
Let $\{x^i\}$ be a local coordinate system on an open neighborhood of $p$ with $x^i(p)=0$. Let $X$ and $Y$ be smooth vector fields on $M$ and write them in this coordinate system as
\[
  X=a^i\frac{\partial}{\partial x^i} \quad\text{and}\quad Y=b^j\frac{\partial}{\partial x^j},
\]
where we are using the Einstein summation convention. By definition \eqref{eq:Hessian_bilinear_form} of the Hessian bilinear form, we see that
\[
  \hess h(X,Y)
  =
  a^ib^j\frac{\partial^2h}{\partial x^i\partial x^j}
  + a^i\frac{\partial b^j}{\partial x^i}\frac{\partial h}{\partial x^j}.
\]
At the critical point $p$, we have
\[
  0 = dh(p) = \frac{\partial h}{\partial x^i}(p)dx^i,
\]
and so the preceding expression for $\hess h(X,Y)$ becomes
\begin{equation}
\label{eq:HessianInLocalCoordinates}
\hess h(p)(X_p,Y_p)
=
a^i(p)b^j(p)\frac{\partial^2h}{\partial x^i\partial x^j}(p).
\end{equation}
If the metric $g$ has components $g_{ij}$ in this coordinate system and $(g^{k\ell})$ denotes the inverse of the matrix $(g_{ij})$, then the definition \eqref{eq:DefineGradient} of $\grad_g h$ yields
\[
  \grad_g h = g^{k\ell} \frac{\partial h}{\partial x^k}\frac{\partial}{\partial x^\ell}
\]
in local coordinates. Therefore,
\begin{align*}
  [\grad_g h,X]_p
  &=
    \left[g^{k\ell} \frac{\partial h}{\partial x^k}\frac{\partial}{\partial x^\ell},
    a^i\frac{\partial}{\partial x^i}\right]_p
  \\
  &=
    g^{k\ell}(p) \frac{\partial h}{\partial x^k}(p)\frac{\partial a^i}{\partial x^\ell}(p)\frac{\partial}{\partial x^i}(p)
    -
    a^i(p)\frac{\partial}{\partial x^i}
    \left(g^{k\ell} \frac{\partial h}{\partial x^k}\right)(p) \frac{\partial}{\partial x^\ell}(p)
  \\
  &= -a^i(p) g^{k\ell}(p) \frac{\partial^2 h}{\partial x^i\partial x^k}(p)
    \frac{\partial}{\partial x^\ell}(p),
\end{align*}
where the last equality follows from the fact that $(\partial h/\partial x^k)(p) = 0$ for all $k$. Thus, for any smooth vector field $Y$ on $M$, we obtain
\begin{align*}
\left\langle 
-[\grad_g h,X]_p,Y_p
\right\rangle_g
&=
\left\langle 
a^i(p) g^{k\ell}(p) \frac{\partial^2 h}{\partial x^i\partial x^k}(p)\frac{\partial}{\partial x^\ell}(p),
b^j(p) \frac{\partial}{\partial x^j}(p)
\right\rangle_g
\\
&=
a^i(p) g^{k\ell}(p)\frac{\partial^2 h}{\partial x^i\partial x^k}(p)b^j(p) g_{\ell j}(p)
\\
&=
     a^i(p)b^j(p)\delta_j^k\frac{\partial^2 h}{\partial x^i\partial x^k}(p)
  \\
&= a^i(p)b^j(p)\frac{\partial^2 h}{\partial x^i\partial x^j}(p),
\end{align*}
and consequently,
\begin{equation}
\label{eq:GradientBracketHessian}
\left\langle
-[\grad_g h,X]_p,Y
\right\rangle_p
=
a^i(p)b^j(p)\frac{\partial^2 h}{\partial x^i\partial x^j}(p).
\end{equation}
By comparing \eqref{eq:GradientBracketHessian} and \eqref{eq:HessianInLocalCoordinates}, we see that
\begin{equation}
\label{eq:GradientBracketHessianEquality}
\left\langle
-[\grad_g h,X]_p,Y_p
\right\rangle_g
=
\hess h(p)(X_p,Y_p),
\end{equation}
for all smooth vector fields $X$ and $Y$ on $M$. Consequently, we obtain
\begin{align*}
  \left\langle\Hess h(p)X_p,Y_p\right\rangle
  &= X_p(Yh) \quad\text{(by \eqref{eq:Hessian_operator})}
  \\
  &= \hess h(p)(X_p,Y_p) \quad\text{(by \eqref{eq:Hessian_bilinear_form})}
  \\
  &= \left\langle -[\grad_g h,X]_p,Y_p \right\rangle_g \quad\text{(by \eqref{eq:GradientBracketHessianEquality})},
\end{align*}
for all smooth vector fields $X$ and $Y$ on $M$. We conclude that \eqref{eq:GradientDefnOfHessian} holds and this completes the proof of the lemma.
\end{proof}

\begin{rmk}[Alternative expression for the Hessian operator]
\label{rmk:Alternative_expression_Hessian_operator} 
Given any smooth covariant derivative $\nabla$ on the tangent bundle $TM$, we may also define
\begin{equation}
\label{eq:DefineHessianUsingConnection}
  \Hess_g^\nabla h(X) := \nabla_X\grad_g h \in C^\infty(TM), \quad\text{for all } X \in C^\infty(TM).
\end{equation}
where $\grad_g h$ is as in \eqref{eq:DefineGradient}. Because the torsion of $\nabla$ is a tensor, at a critical point $p$ of $h$ we have the equality
\[
(\nabla_X\grad_g h)_p -(\nabla_{\grad_g h}X)_p-[X,\grad_g h]_p = 0
\]
since $\grad_g h(p) = 0$ and thus
\[
(\nabla_X\grad_g h)_p=-[\grad_g h,X]_p.
\]
Combining the preceding equality with the identity \eqref{eq:GradientDefnOfHessian} and the definition \eqref{eq:DefineHessianUsingConnection} shows that
\[
  \Hess_g^\nabla h(p)X_p = \Hess_gh(p)X_p, \quad\text{for all } X_p \in T_pM,
\]
and so the definitions
\eqref{eq:DefineHessianUsingConnection} and \eqref{eq:Hessian_operator} of Hessian operators agree when $p$ is a critical point.
\end{rmk}

\section[Second fundamental form of a submanifold of Euclidean space]{Second fundamental form of a submanifold of Euclidean space defined by preimage of regular value of a smooth map}
\label{sec:Second_fundamental_form_submanifold_Euclidean_space}
Let $M \subset N$ be a smoothly and isometrically embedded Riemannian submanifold of a smooth Riemannian manifold $N$. Recall that the \emph{second fundamental tensor} $S$ for a submanifold $M\subset N$ is the bilinear map
\[
  S:TM \times T^\perp M \ni (v,\eta) \mapsto \pi_{TM}(\nabla_v^N\eta) \in TM
\]
defined in Jost \cite[Definition 5.1.1]{Jost_riemannian_geometry_geometric_analysis_e7}, where $T^\perp M \subset TN$ is the normal subbundle defined by the isometric embedding $M \hookrightarrow N$, and $\nabla^N$ is the Levi--Civita connection on $TN$, and $\pi_{TM}$ denotes orthogonal projection from $TN = TM\oplus T^\perp M$ onto the subbundle $TM$. The \emph{second fundamental form} is the symmetric bilinear map
\[
  \SFF:TM \times TM  \ni (v,w) \mapsto \pi_{T^\perp M}(\nabla_v^Nw) \in T^\perp M
\]  
by \cite[Definition 5.1.2]{Jost_riemannian_geometry_geometric_analysis_e7} and the proof of \cite[Lemma 5.1.2]{Jost_riemannian_geometry_geometric_analysis_e7} due to Jost, where $\pi_{T^\perp M}$ denotes orthogonal projection from $TN$ onto the normal bundle $T^\perp M$.
(The map $S$ is called the \emph{shape operator} by Palais and Terng \cite[Section 2.1]{PalaisTerng}.) One says that $M$ is a \emph{totally geodesic submanifold} if the second fundamental form vanishes identically along $M$, that is, $\SFF \equiv 0$ (see Jost \cite[Definition 5.2.1 and Theorem 5.2.2, p. 238]{Jost_riemannian_geometry_geometric_analysis_e7}, O'Neill \cite[Chapter 4, Definition 12 and Proposition 13, p. 103]{ONeill_semi-riemannian_geometry_applications_general_relativity}, or Palais and Terng \cite[Section 2.2, p. 33]{PalaisTerng}). In preparation for our proof of the forthcoming Lemma \ref{lem:Hessian_form_restriction_smooth_function_to_submanifold}, we establish the

\begin{lem}[Formula for the second fundamental form of a submanifold of Euclidean space defined by the preimage of a regular value of a smooth map]
\label{lem:Formulae_second_fundamental_form_submanifold_Euclidean_space}
Let $m, n$ be integers with $m\geq 1$ and $n\geq m+1$ and $F:\RR^n\to \RR^{n-m}$ be a smooth map for which the origin $0\in\RR^{n-m}$ is a regular value and $M = F^{-1}(0) \subset \RR^n$ denotes the corresponding embedded smooth submanifold of dimension $m$. Let $p \in M$, fix an orthogonal splitting $\RR^n = T_pM\oplus T_p^\perp M \cong \RR^m\oplus\RR^{n-m}$, and write the coordinates on $\RR^n$ with respect to this splitting as $(x,y)$, with $(x(p),y(p))=(0,0)$. Let $U \subset \RR^m$ be an open neighborhood of the origin and $f:\RR^m\supset U \to \RR^{n-m}$ be the smooth map and $\varphi:\RR^m\supset U \to M$ be the inverse graph coordinate chart, $\varphi(x) = (x,f(x)) \in M$ for all $x\in U$, provided by the Implicit Mapping Theorem with $\varphi(0)=p$. Let $\SFF \in \Hom(T^*M\otimes T^*M, T^\perp M)$ denote the second fundamental form of the embedding $M\hookrightarrow \RR^n$ and induced Riemannian metric $g$ on $M$. Then
\begin{equation}
\label{eq:Formulae_second_fundamental_form_submanifold_Euclidean_space}
  \SFF_g(p) = D^2f(\varphi^{-1}(p)) = - D_2F(p)^{-1}D_1^2F(p) \in \Hom(T_pM\otimes T_pM,T_p^\perp M),
\end{equation}
where $D_1^2F(p)\in \Hom(T_pM\otimes T_pM,T_p^\perp M)$ denotes the second-order partial derivative of $F:\RR^m\times\RR^{n-m}\to \RR^{n-m}$ with respect to $x\in\RR^m$, and $D_2F(p) \in \End(T_p^\perp M)$ denotes the first-order partial derivative of $F$ with respect to $y\in\RR^{n-m}$ that is invertible by hypothesis, and $D^2f(0)$ denotes the second-order derivative of $f:\RR^m\supset U \to \RR^{n-m}$ at $\varphi^{-1}(p) = 0$.
\end{lem}

Before proceeding to the proof, we recall the Chain Rule for second-order derivatives of maps of Banach spaces to help clarify notation guided by Lang \cite[Chapter 1, Sections 1--3]{Lang_introduction_differential_topology}. If $F:\sY\to\sZ$ and $G:\sX\to\sY$ are smooth maps of Banach spaces, then the Chain Rule gives, for all $x\in \sX$,
\begin{equation}
  \label{eq:Chain_rule_second-order_derivatives}
  D^2(F\circ G)(x) = D^2F(G(x))\circ DG(x)^2 + DF(G(x))\circ D^2G(x) \in \Hom(\sX\otimes\sX,\sZ),
\end{equation}
where we identify\footnote{See Ryan \cite{Ryan_introduction_tensor_products_Banach_spaces} for definitions of norms on tensor products of Banach spaces.}
the Banach space of linear operators, $\Hom(\sX\otimes\sX,\sZ)$, with the Banach space of bilinear operators, $\Hom^2(\sX\oplus\sX,\sZ) \cong \Hom(\sX, \Hom(\sX,\sZ))$, in the usual way by
\[
  \widehat T(v_1,v_2) := T(v_1\otimes v_2), \quad\text{for all } v_1, v_2 \in \sX.
\]
When $\widehat T$ is symmetric, its values are determined by its restriction to the diagonal and we write (for example, see Lang \cite[Chapter 1, Sections 1--3]{Lang_introduction_differential_topology} or Whittlesey \cite{Whittlesey_1965})
\[
  \widehat T(v_1,v_2) = T(v_1v_2) = Tv_1v_2 \quad\text{and}\quad \hat T(v,v) = Tv^2,
  \quad\text{for all } v_1, v_2, v \in \sX.
\]
If $A \in \Hom(\sX,\sY)$, we let $A^2 \in \Hom(\sX,\sY\oplus\sY)$ denote the operator defined by
\[
  A^2v := (Av,Av), \quad\text{for all } v \in \sX.
\]  
In particular, the operator $D^2G(x) \in L^2(\sX\oplus\sX,\sZ)$ is symmetric for any $x\in\sX$ and the chain rule \eqref{eq:Chain_rule_second-order_derivatives} can be written more explicitly as
\[
  D^2(F\circ G)(x)v^2 = D^2F(G(x))(DG(x)v)^2 + DF(G(x))D^2G(x)v^2, \quad\text{for all } v \in \sX,
\]
or equivalently, for any $x\in\sX$,
\[
  D^2(F\circ G)(x)(v,v) = D^2F(G(x))(DG(x)v,DG(x)v) + DF(G(x))D^2G(x)(v,v),
  \quad\text{for all } v \in \sX.
\]
We now give the

\begin{proof}[Proof of Lemma \ref{lem:Formulae_second_fundamental_form_submanifold_Euclidean_space}]
We shall write the derivative or partial derivative of a smooth map from one Euclidean space to another by ``$D$'' or ``$D_i$'' in the calculations below, when the domain is given by a product of Euclidean subspaces and $i$ denotes the $i$th factor in the product. We have
\[
  T_pM = T_pF^{-1}(0) = \Ker DF(p),
\]
where $DF(p) = D_1F(p)+D_2F(p) \in \Hom(\RR^m\times \RR^{n-m},\RR^{n-m})$, with
\begin{align*}
  D_1F(p) &= DF(p) \restriction T_pM \in \Hom(T_pM,\RR^{n-m}),
  \\
  D_2F(p) &= DF(p) \restriction T_p^\perp M  \in \Hom(T_p^\perp M,\RR^{n-m}).
\end{align*}
We have $T_pM = \Ker D_1F(p)$, and thus $D_1F(p)=0$, while $D_2F(p)$ is an isomorphism since $0\in \RR^{n-m}$ is a regular value of $F$ by hypothesis. We may assume without loss of generality that $p$ is the origin in $\RR^n = \RR^m\times\RR^{n-m}$. The Implicit Mapping Theorem therefore provides the open neighborhood $U\subset T_pM$ of the origin and the smooth map $f:T_pM \supset U \ni x \mapsto y=f(x) \in T_p^\perp M$ such that $f(0)=0$ and inverse graph coordinate chart $\varphi:U\to M$ with $\varphi(0) = (0,f(0)) = (0,0) = p$ and
\[
  F(x,f(x)) = F(\varphi(x)) = 0 \in \RR^{n-m}, \quad\text{for all } x \in U.
\]
Hence, for $y=f(x)$,
\[
  D_x(F(x,f(x)) =  D_1F(x,y) + D_2 F(x,y)Df(x) = 0, \quad\text{for all } x \in U,
\]
and
\begin{multline*}
  D_x^2(F(x,f(x))  =  D_1^2F(x,y) + 2D_{12}F(x,y)Df(x) + D_2F(x,y)D^2f(x)
  \\
  + D_2^2 F(x,y)Df(x)^2 = 0, \quad\text{for all } x \in U.
\end{multline*}
Evaluating the first-derivative equation at $p$ gives
\[
  D_1F(p) + D_2 F(p)Df(0) = 0 \in \Hom(T_pM,\RR^{n-m}).
\]
Now $D_1F(p)=DF(p)\restriction T_pM=0$ (since $T_pM = \Ker DF(p)$) and thus
\[
  D_2 F(p)Df(0) = 0 \in \Hom(T_pM,\RR^{n-m}).
\]
But $D_2F(p) \in \Hom(T_p^\perp M,\RR^{n-m})$ is an isomorphism and so
\begin{equation}
  \label{eq:Df(0)_is_zero}
  Df(0) = 0 \in \Hom(T_pM,T_p^\perp M).
\end{equation}  
Evaluating the second-derivative equation at $p$ gives
\[
  D_1^2F(p) + 2D_{12}F(p)Df(0) + D_2F(p)D^2f(0) + D_2^2 F(p)Df(0)^2 = 0
  \in \Hom(T_pM\otimes T_pM,\RR^{n-m}).
\]
Using $Df(0) = 0$ by \eqref{eq:Df(0)_is_zero}, the preceding equation simplifies to
\[
  D_1^2F(p) + D_2F(p)D^2f(0) = 0 \in \Hom(T_pM\otimes T_pM,\RR^{n-m}).
\]
Since $D_2F(p) \in  \Hom(T_p^\perp M,\RR^{n-m})$ is an isomorphism, this gives
\begin{equation}
  \label{eq:Explicit_equation_second_derivative_f}
  D^2f(0) = - D_2F(p)^{-1}D_1^2F(p) \in \Hom(T_pM\otimes T_pM,T_p^\perp M).
\end{equation}
According to Salamon and Robbin \cite[Exercise 3.1.11]{Robbin_Salamon_introduction_differential_geometry}, we have
\[
  \SFF_p(v,w) = D^2f(0)vw, \quad\text{for all } v, w \in T_pM,
\]
and thus by \eqref{eq:Explicit_equation_second_derivative_f}, we obtain
\[
  \SFF_p(v,w) = - D_2F(p)^{-1}D_1^2F(p)vw, \quad\text{for all } v, w \in T_pM.
\]
This completes the verification of the identities \eqref{eq:Formulae_second_fundamental_form_submanifold_Euclidean_space}.
\end{proof}

\section{Hessian of the restriction of a smooth function to a submanifold}
\label{sec:Hessian_restriction_smooth_function_submanifold_Euclidean_space}
In this section, we prove the

\begin{lem}[Hessian form of the restriction of a smooth function to a submanifold]
\label{lem:Hessian_form_restriction_smooth_function_to_submanifold}  
Let $M \subset \RR^n$ be a smoothly and isometrically embedded submanifold with induced Riemannian metric $g$. Let $H:\RR^n\to\RR$ be a smooth function and $h:M\to \RR$ denote its restriction to $M$. If $p\in M$ is a critical point of $h$, then
\begin{equation}
\label{eq:Hessian_form_restriction_smooth_function_to_submanifold}   
  \hess_g h(p) = \hess H(p)\circ d\iota(p)^2  + dH(p)\circ \SFF_g(p) \in \Hom(T_p^*M\otimes T_p^*M,\RR), 
\end{equation}
where $\SFF_g \in \Hom(T^*M\otimes T^*M, T^\perp M)$ is the second fundamental form of the embedding $\iota:M\hookrightarrow \RR^n$, and $d\iota:TM\to\RR^n\times \RR^n$ is the differential of that embedding, and $d\iota(p)^2 \in \Hom(T_pM\otimes T_pM,\RR^n)$ is the symmetric operator defined by $d\iota(p)^2(v,w) = (d\iota(p)v, d\iota(p)w)$, for all $v, w \in T_pM$. 
\end{lem}

\begin{proof}
We continue the notation of Lemma \ref{lem:Formulae_second_fundamental_form_submanifold_Euclidean_space}, so $\varphi:\RR^m\supset U \to (x,f(x)) \in M \subset \RR^m\times\RR^{n-m}$ is the graph inverse coordinate chart, and we abbreviate $h(x) = (h\circ\varphi)(x)$ for the function $h\circ\varphi:\RR^m\supset U \to \RR$. We thus compute
\[
  Dh(x) = D_x(H(x,f(x))) = D_1H(x,f(x)) + D_2H(x,f(x))Df(x) \in \Hom(\RR^m,\RR),
  \quad\text{for all } x \in U,
\]
and
\begin{multline*}
  D^2h(x) = D_x^2(H(x,f(x))) = D_1^2H(x,f(x)) + 2D_{12}H(x,f(x))Df(x)
  \\
  + D_2^2H(x,f(x))Df(x)^2 + D_2H(x,f(x))D^2f(x) \in \Hom(\RR^m\otimes\RR^m,\RR),
  \quad\text{for all } x \in U.
\end{multline*}
Evaluating the expression for $Dh(x)$ at the critical point $p=\varphi(0)\in M$ of $h$ gives
\[
  Dh(0) = D_1H(p) + D_2H(p)Df(0) = 0 \in \Hom(T_pM,\RR),
\]
and thus $Dh(0) = D_1H(p) = DH(p) \restriction T_pM = 0 \in \Hom(T_pM,\RR)$, as expected, since $Df(0)=0$ by \eqref{eq:Df(0)_is_zero}. Moreover, evaluating the expression for $D^2h(x)$ at $p=\varphi(0)\in M$ yields
\[
  D^2h(0) = D_1^2H(p) + 2D_{12}H(p)Df(0) + D_2^2H(p)Df(0)^2 + D_2H(p)D^2f(0)
  \in \Hom(T_pM\otimes T_pM,\RR).
\]
Simplifying by using $Df(0)=0$ from \eqref{eq:Df(0)_is_zero} implies that
\[
  D^2h(0) = D_1^2H(p) + D_2H(p)D^2f(0) \in \Hom(T_pM\otimes T_pM,\RR).
\]
Substituting our expression $D^2f(0)$ from \eqref{eq:Explicit_equation_second_derivative_f} leads to the formula
\begin{equation}
  \label{eq:Hessian_restriction_global_function_submanifold_raw}
  D^2h(0) = D_1^2H(p) - D_2H(p)D_2F(p)^{-1}D_1^2F(p) \in \Hom(T_pM\otimes T_pM,\RR).
\end{equation}
Combining the identities \eqref{eq:Hessian_restriction_global_function_submanifold_raw} and \eqref{eq:Formulae_second_fundamental_form_submanifold_Euclidean_space} completes the proof of Lemma \ref{lem:Hessian_form_restriction_smooth_function_to_submanifold}.
\end{proof}

We obtain a Hessian operator $\Hess_g h(p) \in \End(T_pM)$ from its definition \eqref{eq:Hessian_operator} and the expression \eqref{eq:Hessian_restriction_global_function_submanifold_raw} for the Hessian bilinear form $\hess_g h(p) \in T_p^*M\otimes T_p^*M$ using the inverse Riesz isomorphism defined by the Riemannian metric $g$ on $M$, namely $g(p)^{-1}:T_p^*M \to T_pM$, to give
\[
  \Hess_g h(p) = g(p)^{-1}\circ\hess_g h(p) \in T_pM\otimes T_p^*M = \End(T_pM)
\]
and therefore
\begin{equation}
  \label{eq:Hessian_operator_restriction_global_function_submanifold_raw}
   \Hess_g h(p) = g(p)^{-1}D_1^2H(p) - g(p)^{-1}D_2H(p)D_2F(p)^{-1}D_1^2F(p) \in \End(T_pM).
\end{equation}  
If $u\in T_pM$ is a unit eigenvector of $\Hess_g h(p)$ with eigenvalue $\lambda$, then $\hess_g h(p)u^2 = \lambda$ and
\begin{equation}
  \label{eq:Hessian_restriction_global_function_submanifold_raw_eigenvalue}
  \lambda = D_1^2H(p)u^2 - D_2H(p)D_2F(p)^{-1}D_1^2F(p)u^2 \in \RR.
\end{equation}
While the expression \eqref{eq:Hessian_restriction_global_function_submanifold_raw_eigenvalue} appears to give an explicit formula for $\lambda$, it has limited value in practice since it is very difficult to apply it to compute $\lambda$ or even determine its sign. In typical applications in geometric analysis (including within this monograph), the map $F$ of Banach spaces will be a Fredholm map of Banach manifolds or Fredholm section of a Banach vector bundle over a Banach manifold, the derivative $D_2F(p)$ will be an elliptic partial differential operator, and $D_2F(p)^{-1}$ will be the corresponding Green's operator. It is for this reason that Frankel's Theorem \ref{mainthm:Frankel_almost_Hermitian} is so powerful, though $M$ must be equipped with a circle action and a circle-invariant, non-degenerate two-form with Hamiltonian function in order to determine the signature of the Hessian of the Hamiltonian function or, if given a circle-invariant Riemannian metric, the eigenvalues of the Hessian of the Hamiltonian function.

\begin{rmk}[Simplifications for totally geodesic submanifolds]
\label{rmk:Simplifications_totally_geodesic_submanifolds} 
Naturally, Lemma \ref{lem:Hessian_form_restriction_smooth_function_to_submanifold} generalizes to the case where $M$ is an isometrically embedded submanifold of a Riemannian manifold $N$. If $M$ is totally geodesic, the expression \eqref{eq:Hessian_form_restriction_smooth_function_to_submanifold} simplifies to give
\[   
  \hess_g h(p) = \hess H(p)\circ d\iota(p)^2 \in \Hom(T_p^*M\otimes T_p^*M,\RR), 
\]
and the expression \eqref{eq:Hessian_restriction_global_function_submanifold_raw_eigenvalue} for $\lambda$ simplifies to give
\[
  \hess_g h(p)u^2 = \lambda = D_1^2H(p)u^2.
\]
\end{rmk}

\section{Sylvester's law of inertia}
\label{sec:Sylvesters_law_inertia}
As we recalled in Section \ref{sec:Introduction_Hessian_restriction_smooth_function_submanifold_Euclidean_space}, we require a choice of Riemannian metric $g$ on $TM$ in order to define the Hessian operator $\Hess_gh \in \End(TM)$ in terms of the bilinear Hessian form $\hess h \in \Hom(T^*M\otimes T^*M,\RR)$ 
and thus, at a critical point $p$ for the smooth function $h:M\to\RR$, the eigenvalues of the operator $\Hess_g f(p) \in \End(T_pM)$ are well-defined. While those eigenvalues depend on the choice of Riemannian metric $g$, the signature of $\Hess_g f(p)$ and hence that of $\hess_g h(p)$ is independent of $g$. To see this, we recall Sylvester's law of inertia in the forthcoming Theorem \ref{thm:Sylvester_law_inertia}. Let $A \in \Mat(n,\CC)$ be a Hermitian matrix, and let $n_+(A)$, $n_-(A)$, and $n_0(A)$ denote, respectively, the number of positive, negative, and zero eigenvalues of $A$. The triple $(n_+(A),n_-(A),n_0(A))$ is called the \emph{signature} of $A$. We shall only need a version of Sylvester's law for real, symmetric matrices.

\begin{thm}[Sylvester's law of inertia]
\label{thm:Sylvester_law_inertia}  
(See Carrell \cite[Theorem 9.13]{Carrell_groups_matrices_vector_spaces} or Horn and Johnson \cite[Theorem 4.5.8]{Horn_Johnson_matrix_analysis_2013}.)
Let $A$ and $B$ be real symmetric matrices that are congruent\footnote{See Carrell \cite[Definition 9.2]{Carrell_groups_matrices_vector_spaces}.}, so there is a nonsingular matrix $C$ such that $A=C^\intercal BC$. Then $A$ and $B$ have the same signature, that is,
\[
  n_+(A) = n_+(B), \quad n_-(A) = n_-(B), \quad\text{and}\quad n_0(A) = n_0(B).
\]
\end{thm}

Theorem \ref{thm:Sylvester_law_inertia} is most easily applied in the setting of the following

\begin{cor}
\label{cor:Sylvester_law_inertia}
Let $A$ and $P$ be real matrices such that $P$ is positive definite and $A$, $P$, and $AP$ are symmetric. Then $A$ and $AP$ have the same signature.
\end{cor}

\begin{proof}
By hypothesis, we have $AP = (AP)^\intercal = P^\intercal A^\intercal = PA$. Because $P$ is positive definite and symmetric, we may write $P = QDQ^{-1}$, where $D$ is diagonal with positive entries and $Q$ is an orthogonal matrix, so $Q^{-1} = Q^\intercal$ (see Carrell \cite[Corollary 8.28]{Carrell_groups_matrices_vector_spaces}). The square root of $P$ is thus well-defined as $P^{1/2} = QD^{1/2}Q^\intercal$ (see Carrell \cite[Section 8.2.4]{Carrell_groups_matrices_vector_spaces}) and $P^{1/2}$ is clearly symmetric. Since $A$ commutes with $P$, it also commutes with $P^{1/2}$ (see Rudin \cite[p. 326]{Rudin})
and thus $AP = AP^{1/2}P^{1/2} = P^{1/2}A(P^{1/2})^\intercal$. The conclusion now follows from Theorem \ref{thm:Sylvester_law_inertia}.
\end{proof}  

Let $\{x^i\}_{i=1}^m$ be local coordinates on an open neighborhood in $M$ of a critical point $p$ for $h$. From the expression \eqref{eq:Hessian_operator} for $\Hess_gh(p) \in \End(T_pM)$, we have
\begin{align*}
  \frac{\partial^2 h}{\partial x^i\partial x^j}(p)
  &=
  g\left(\Hess_g h(p)\frac{\partial}{\partial x^i}(p), \frac{\partial}{\partial x^j}(p)\right)
  \\
  &=
  g\left((\Hess_g h(p))_i^k\frac{\partial}{\partial x^k}(p),\frac{\partial}{\partial x^j}(p)\right)
  \\
  &=
  (\Hess_g f(p))_i^kg_{kj}(p), \quad\text{for } 1 \leq i,j \leq m. 
\end{align*}
Because the matrix $(g_{kj}(p))$ is positive definite and symmetric and the following matrices are symmetric,
\[
  \left(\frac{\partial^2 h}{\partial x^i\partial x^j}(p)\right) \quad\text{and}\quad
  \left( (\Hess_g f(p))_i^k\right),  
\]
they necessarily have the same signature by Corollary \ref{cor:Sylvester_law_inertia}. In particular, the signature of the bilinear Hessian form $\hess f(p)$ is equal to the signature of the Hessian operator $\Hess_g f(p)$ for any Riemannian metric $g$ on $M$, for a given critical point $p$ for $f$.

Theorem \ref{thm:Sylvester_law_inertia} can also be used to extend the classical observation that the Morse index at a non-degenerate critical point of a Morse function on a smooth manifold is well-defined, independent of the coordinate chart (see Nicolaescu \cite[p. 7]{Nicolaescu_morse_theory}) to the invariance of the signature of a Morse--Bott function under a diffeomorphism. Indeed, one has the well-known

\begin{lem}[Transformation of the Hessian of a smooth function at a critical point under a diffeomorphism]
\label{lem:Del_Pino_21}  
(See del Pino \cite[Lemma 21]{DelPino_gentle_introduction_Morse_theory}.)
Let $f : \RR^n \to \RR$ be a smooth function and $\varphi:\RR^n\to\RR^n$ be a diffeomorphism. If $p$ is a critical point of $f\circ\varphi$, then the following transformation formula holds for the Hessian matrix of $f$:
\[
\Hess (f\circ\varphi)(p) = D\varphi(p)^\intercal \Hess f(\varphi(p)) D\varphi(p).
\]
\end{lem}

\begin{proof}
Given a smooth function $g : \RR^n \to \RR$, its differential is the smooth map $Dg:\RR^n\times\RR^n \to \RR$ and its adjoint with respect to the standard inner product on $\RR^n$ is the gradient map, $\grad g = Dg^\intercal : \RR^n \to \RR^n$. The Hessian matrix of $g$ is thus given by $\Hess g = D\grad g = D(Dg^\intercal)$. Applying this observation and the Chain Rule, we obtain
\[
  D(f\circ\varphi)^\intercal(x) =  D\varphi(x)^\intercal  Df(\varphi(x))^\intercal.
\]
In particular, $p$ is a critical point of $f\circ\varphi$ if and only if $\varphi(p)$ is a critical point of $f$ since $\varphi$ is a diffeomorphism and $D\varphi(p)$ is a diffeomorphism of $\RR^n$. Applying the Chain Rule once more yields
\begin{align*}
  \Hess(f\circ\varphi)(x)
  &= D(D(f\circ\varphi)(x)^\intercal)
  \\
  &=  D\varphi(x)^\intercal  D(Df)^\intercal((\varphi(x)) D\varphi(x)
    + D(D\varphi^\intercal)(x) Df((\varphi(x))^\intercal.
\end{align*}
Since $\varphi(p)$ is a critical point of $f$, the second term vanishes when $x=p$ and the claimed formula holds.
\end{proof}

Theorem \ref{thm:Sylvester_law_inertia} and Lemma \ref{lem:Del_Pino_21} thus yield the

\begin{cor}[Diffeomorphism invariance of the signature of the Hessian bilinear form]
\label{cor:Invariance_signature_Morse-Bott_function}  
Let $M$ be a smooth manifold, $f : M \to \RR$ be a smooth function, and $p$ be a critical point of $f$. Then the signature of the Hessian bilinear form $\hess f(p)$ is well-defined, independent of the local coordinate chart.
\end{cor}

\chapter{Circle actions on almost symplectic manifolds and Frankel's theorem}
\label{chap:Circle_actions_almost_Hermitian_manifolds}
In this chapter, we prove Theorem \ref{mainthm:Frankel_almost_Hermitian}, which extends Frankel's result \cite[Lemma 1]{Frankel_1959} from the category of K\"ahler manifolds to that of almost symplectic manifolds. In Section \ref{sec:Operators_defined_by_vector_fields}, we recall some well-known facts about the action of a vector field on the tangent space at one of its zeros. Section \ref{sec:ReaLAndComplexRepresentations} summarizes results that we require on real and complex representations of $S^1$. In Section \ref{sec:CircleActionsFixedPtsHessians}, we review properties of circle actions on smooth manifolds.  After describing the normal bundle structure of the fixed-point set of such an action, we characterize in Lemma \ref{lem:MultiplicitiesEigenvalues} the circle action on this normal bundle in terms of eigenvalues of a linear map given by Lie bracket with respect to the vector field $X$ generated by the circle action. In Section \ref{sec:CompatTriples}, we prove Lemma \ref{lem:ExistenceOfS1InvarACStructure}, which asserts that the almost complex structure which is compatible with a circle-invariant, non-degenerate two-form $\om$ and circle-invariant Riemannian metric $g$ is also circle-invariant. In Section \ref{sec:HessiansAndMomentMaps}, we prove Lemma \ref{lem:IdentifyEigenvalues}, giving a correspondence between the weights of the circle action and the eigenvalues of the Hessian of the Hamiltonian function, in the sense of \eqref{eq:MomentMap}, and thus complete the proof of Theorem \ref{mainthm:Frankel_almost_Hermitian}.

\section{Operators defined by vector fields}
\label{sec:Operators_defined_by_vector_fields}
If $X$ is a smooth vector field on a smooth manifold, we recall that the Lie derivative of another smooth vector field $Y$ in the direction of $X$ is given by (see Lee \cite[Theorem 9.38]{Lee_john_smooth_manifolds})
\[
  \sL_X Y = [X,Y].
\]
We then have the well-known

\begin{lem}[Action of a vector field on the tangent space at a zero of the vector field]
\label{lem:Action_vector_field_at_a_zero}
(See Bott \cite[pp. 231 and 239--240]{Bott_1967mmj}.)
Let $M$ be a smooth manifold and $X$ be a smooth vector field on $M$. If $X$ is zero at a point $p$, then the Lie derivative $\sL_X:C^\infty(TM)\to C^\infty(TM)$ restricts to an operator $(\sL_X)_p \in \End(T_pM)$. Moreover, if $g$ is a smooth Riemannian metric on $M$, then the following are equivalent:
\begin{enumerate}
\item
  \label{item:Riemannian_metric_under_flow_vector_field}
  $g$ is invariant under the flow of $X$.
\item
  \label{item:Lie_derivative_Riemannian_metric_is_zero}
  $X$ is an infinitesimal isometry of $g$, that is, $\sL_Xg = 0$ on $M$.
\item
  \label{item:Lie_derivative_Riemannian_metric_by_infinitesimal_isometry}
  The following identity holds:
  \begin{equation}
  \label{eq:Lie_derivative_Riemannian_metric_by_infinitesimal_isometry}  
  \langle \sL_X Y, Z \rangle_g + \langle Y, \sL_X Z\rangle_g = X\langle Y,Z\rangle_g,
  \quad\text{for all } Y, Z \in C^\infty(TM).
\end{equation}  
\end{enumerate}
In particular, if $\sL_Xg=0$ then $(\sL_X)_p \in \End(T_pM)$ is skew-adjoint with respect to $g$ on $T_pM$. The eigenvalues of $(\sL_X)_p$ are given by $\pm i\lambda_1,\ldots,\pm i\lambda_n$ (independent of any choice of $g$) and zero with multiplicity $\dim M - 2n$, where $\lambda_j>0$ for $j=1,\ldots,n$.
\end{lem}  

\begin{proof}
Following Bott \cite[p. 231]{Bott_1967mmj}, we observe that if $Y_p \in T_pM$ and $Y$ is any extension of $Y_p$ to a smooth vector field on $M$, then for all smooth functions $f:M\to\RR$ we have
\[
  (\sL_XY)f = [X,Y]f = X(Yf) - Y(Xf),
\]
so that 
\[
  (\sL_XY)_pf = X_p(Yf) - Y_p(Xf) = -Y_p(Xf) \in \RR,
\]
which is independent of the choice of extension of $Y_p$. Thus, $(\sL_XY)_pf$ is uniquely determined by $Y_p$ and $f$ and so $(\sL_X)_p \in \End(T_pM)$ is well-defined by the expression
\[
  (\sL_X)_pY_p := (\sL_XY)_p,
\]
for any extension of $Y_p$ to a smooth vector field $Y$ on $M$.

Items \eqref{item:Riemannian_metric_under_flow_vector_field} and \eqref{item:Lie_derivative_Riemannian_metric_is_zero} are equivalent by Lee \cite[Theorem 12.37]{Lee_john_smooth_manifolds}. According to Lee \cite[Proposition 12.32]{Lee_john_smooth_manifolds}, we have
\[
  X(g(Z,W)) = (\sL_Xg)(Z,W) + g(\sL_XY,W) + g(Y,\sL_XW),
  \quad\text{for all } Y, Z \in C^\infty(TM),
\]
and so the assertion $\sL_Xg = 0$ in \eqref{item:Lie_derivative_Riemannian_metric_is_zero} is equivalent to the identity \eqref{eq:Lie_derivative_Riemannian_metric_by_infinitesimal_isometry} holding in Item \eqref{item:Lie_derivative_Riemannian_metric_by_infinitesimal_isometry}.

Hence, if $\sL_Xg=0$ and $\{e_1,\ldots,e_d\}$ is a local $g$-orthonormal frame for $TM$ on an open neighborhood of $p$, then
\[
  \langle \sL_X e_j, e_k \rangle_g + \langle e_j, \sL_X e_k\rangle_g = 0,
  \quad\text{for } 1 \leq j,k \leq d.
\]
Thus, $(\sL_X)_p \in \End(T_pM)$ is skew-adjoint with respect to $g$.

The final assertion follows from elementary linear algebra. Recall that $A\in\End(\RR^d)$ is a skew-symmetric matrix if and only if there are an orthogonal matrix $Q \in \Or(d)$ and an integer $n\geq 0$ such that $Q^\intercal AQ$ has the block-diagonal form (see Horn and Johnson \cite[Corollary 2.5.11]{Horn_Johnson_matrix_analysis_2013})
\begin{equation}
\label{eq:LieDerivative_Diagonalized_Form}
  \lambda_1\begin{pmatrix}0 & 1 \\ -1 & 0\end{pmatrix}
  \oplus\cdots\oplus
  \lambda_n\begin{pmatrix}0 & 1 \\ -1 & 0\end{pmatrix}
  \oplus
  \underbrace{0\oplus\cdots\oplus 0}_{d-2n}
\end{equation}
for positive real numbers $\lambda_j$, with $j=1,\ldots,n$. Moreover, if $A$ is not identically zero, then its non-zero eigenvalues are $\pm i\lambda_1,\ldots,\pm i\lambda_n$. Independence of the eigenvalues of $(\sL_X)_p$ from a choice of Riemannian metric $g$ or local $g$-orthonormal frame is immediate from the fact that they are defined independently of such choices.
This completes the proof.
\end{proof}

\begin{lem}[$J$-linearity of the Lie derivative at a zero of the vector field]
\label{lem:Action_vector_field_at_a_zero_J_Linear}
Let $M$ be a smooth manifold and $X$ a smooth vector field on $X$. If $X$ is zero at a point $p\in M$ and
there is a smooth almost complex structure $J$ on $M$, then $(\sL_X)_p$ is $J$-linear,
\begin{equation}
  \label{eq:Action_vector_field_at_a_zero_J_Linear}
  J(\sL_X)_p = (\sL_X)_pJ.
\end{equation}
Furthermore, if the almost complex structure $J$ acts by isometries with respect to a smooth Riemannian metric $g$ and if $\sL_X g=0$ as in Lemma \ref{lem:Action_vector_field_at_a_zero}, then $J(\sL_X)_p \in \End(T_pM)$ is self-adjoint with respect to $g$ and is diagonalizable with real eigenvalues $\mu_1,\dots,\mu_n$ (independent of any choice of $g$) each of multiplicity two and zero with multiplicity $\dim M - 2n$, where $|\mu_j|=\la_j$, for $j=1,\ldots,n$, and $\pm i\lambda_1,\dots,\pm i\lambda_n$ are the eigenvalues of $(\sL_X)_p$ appearing in Lemma \ref{lem:Action_vector_field_at_a_zero}.
\end{lem}

\begin{proof}
Observe that for any smooth vector field $Y$ on $M$, the fact that $X_p=0$ yields
\begin{align*}
  J(\sL_X)_pY &= J[X,Y]_p = J(X_pY - Y_pX) = -JY_pX,
  \\
  (\sL_X)_pJY &= [X,JY]_p = X_p(JY) - JY_pX = -JY_pX,
\end{align*}
and this proves \eqref{eq:Action_vector_field_at_a_zero_J_Linear}. (As an aside, we note that the same argument yields
\begin{equation}
\label{eq:JandBracket}
[JX, Y]_p
=
J[X, Y]_p,
\end{equation}
a simple identity that we use later.) If $J$ acts by isometries, then, because $J^2 = -\id$, we obtain
\[
  g(X,JY)=g(JX,J^2Y)=-g(JX,Y),
\]
and thus $J$ is skew-adjoint with respect to $g$. Lemma \ref{lem:Action_vector_field_at_a_zero} implies that $(\sL_X)_p$ is also skew-adjoint with respect to $g$ and so the composition $J(\sL_X)_p$ is symmetric with respect to $g$ and hence, by Horn and Johnson \cite[Corollary 2.5.11]{Horn_Johnson_matrix_analysis_2013}, the operator $J(\sL_X)_p$ is diagonalizable with real eigenvalues. The endomorphisms $J$ and $(\sL_X)_p$ of $T_pM$ commute by \eqref{eq:Action_vector_field_at_a_zero_J_Linear}.  Because they are skew-adjoint, $J$ and $(\sL_X)_p$ are both \emph{normal} in the sense of \cite[Definition 2.5.1]{Horn_Johnson_matrix_analysis_2013}
in that they satisfy $JJ^*=J^*J$ and $(\sL_X)_p(\sL_X)_p^*=(\sL_X)_p^*(\sL_X)_p$.
By \cite[Theorem 2.5.15]{Horn_Johnson_matrix_analysis_2013}, if $\sN$ is a non-empty family of commuting normal matrices, then there is a real, orthogonal matrix $Q$ and a non-negative integer $n$ such that for all $A\in \sN$,
\begin{equation}
\label{eq:QuasiDiagonalForm}
Q^\intercal A Q
=
\begin{pmatrix}
a_1(A) & b_1(A)
\\
-b_1(A) & a_1(A)
\end{pmatrix}
\oplus
\dots
\oplus
\begin{pmatrix}
a_n(A) & b_n(A)
\\
-b_n(A) & a_n(A)
\end{pmatrix}
\oplus
a_{2n+1}(A)
\oplus
\dots
a_d(A)
\end{equation}
where $a_j(A),b_j(A)\in\RR$.  From \eqref{eq:QuasiDiagonalForm}, we see that the eigenvalues of the matrix $A\in\sN$ are given by $a_1(A)\pm i b_1(A),\dots,a_n(A)\pm ib_n(A),a_{2n+1},\dots,a_d(A)$.

Because $(\sL_X)_p$ and $J$ form a non-empty family of commuting
normal endomorphisms of $T_pM$, then
\cite[Theorem 2.5.15]{Horn_Johnson_matrix_analysis_2013} implies that there is an orthonormal basis $\{e_1,\dots,e_d\}$ of $T_pM$ with respect to which $(\sL_X)_p$ and $J$ both have the form appearing on the right-hand-side of \eqref{eq:QuasiDiagonalForm}. The eigenvalues of $(\sL_X)_p$ and $J$ determine the parameters $a_j((\sL_X)_p)$, $b_j((\sL_X)_p)$, $a_j(J)$, and $b_j(J)$. The eigenvalues of $J:T_pM\to T_pM$ are $\pm i$ (with multiplicity $d/2$. The eigenvalues of $(\sL_X)_p$, as identified in Lemma \ref{lem:Action_vector_field_at_a_zero}, are $\pm i\la_1,\dots,\pm i\la_n$ and $0$ with multiplicity $d-2n$. Hence, there is an orthonormal basis $\{e_1,\dots,e_d\}$ of $T_pM$ with respect to which
\begin{equation}
\label{eq:LieDerivativeFamilyDiagonalizedForm}
(\sL_X)_p
=
  \mu_1\begin{pmatrix}0 & 1 \\ -1 & 0\end{pmatrix}
  \oplus\cdots\oplus
  \mu_n\begin{pmatrix}0 & 1 \\ -1 & 0\end{pmatrix}
  \oplus
  \underbrace{0\oplus\cdots\oplus 0}_{d-2n},
\end{equation}
where the $\mu_j \in \RR$ obey $|\mu_j| = \lambda_j$ for $j=1,\ldots,n$ and
\begin{equation}
  \label{eq:JFamilyDiagonalizedForm}
  J
=
\begin{pmatrix}0 & 1 \\ -1 & 0\end{pmatrix}
  \oplus\cdots\oplus
\begin{pmatrix}0 & 1 \\ -1 & 0\end{pmatrix}.
\end{equation}
Equations \eqref{eq:LieDerivativeFamilyDiagonalizedForm} and \eqref{eq:JFamilyDiagonalizedForm}
imply that with respect to the basis $\{e_1,\dots,e_d\}$, 
\[
  J(\sL_X)_p
  =
  -\mu_1\begin{pmatrix}1 & 0 \\ 0 & 1\end{pmatrix}
  \oplus\cdots\oplus
  -\mu_n\begin{pmatrix}1 & 0 \\ 0 & 1\end{pmatrix}
  \oplus
  \underbrace{0\oplus\cdots\oplus 0}_{d-2n}.
\]
This proves the stated relation between the eigenvalues of $(\sL_X)_p$ and $J(\sL_X)_p$.
\end{proof}

\begin{rmk}
\label{rmk:Sign_Of_Parameter}
We now explain why the requirement that $J$ have the form \eqref{eq:JFamilyDiagonalizedForm} means that the
parameters $\mu_j$ appearing in \eqref{eq:LieDerivativeFamilyDiagonalizedForm}
need not be positive, unlike the corresponding parameters $\la_j$ in \eqref{eq:LieDerivative_Diagonalized_Form}. Let $\{e_1,\dots,e_d\}$ again be the orthogonal basis of $T_pM$ with respect to which $(\sL_X)_p$ and $J$ are given by \eqref{eq:LieDerivativeFamilyDiagonalizedForm}  and \eqref{eq:JFamilyDiagonalizedForm} respectively. If $\mu_j<0$, then transposing the basis elements $e_{2j-1}$ and $e_{2j}$ would make $\mu_j$ positive as in \eqref{eq:LieDerivativeFamilyDiagonalizedForm} but doing so would introduce a negative sign in the corresponding factor in the expression \eqref{eq:JFamilyDiagonalizedForm} for $J$. The sign of the parameter $\mu_j$ appearing in \eqref{eq:LieDerivativeFamilyDiagonalizedForm} is determined by whether or not the orientations $\{e_{2j-1},(\sL_X)_pe_{2j-1}\}$ and $\{e_{2j-1},Je_{2j-1}\}$ of the two-plane spanned by $e_{2j-1}$ and $e_{2j}$ agree.

We will revisit this issue in Lemma \ref{lem:S1_Invariant_J_Is_MultiplicationBy_i} where we show how the almost complex structure $J$ determines the sign of the weight of an $S^1$ action at a fixed point.
\end{rmk}

\section{Real and complex representations for linear circle actions}
\label{sec:ReaLAndComplexRepresentations}
In this section, we summarize some relatively standard but less well-known results on real and complex representations for $S^1$.

\begin{defn}[Real and complex representations for groups]
\label{defn:Real_and_complex_representations_groups}
(See Brocker and tom Dieck \cite[Chapter II, Definition 1.1, p. 65]{BrockertomDieck}.)
Let $\KK=\RR$ or $\CC$ and $G$ be a group. A \emph{representation} of $G$ is a group homomorphism $\rho:G\to \GL_\KK(V)$, where $V$ is a $\KK$-vector space and $\GL_\KK(V)$ is the group of $\KK$-linear automorphisms of $V$. The representation is \emph{real} if $\KK=\RR$ and \emph{complex} if $\KK=\CC$.
\end{defn}

\begin{defn}[Trivial representation]
\label{defn:TrivialRepresentation}
The \emph{trivial representation} is given by the group homomorphism $\rho:G\to\GL_\KK(\KK)=\KK^*=\KK\less\{0\}$, defined by $\rho(g)=1$ for all $g\in G$.
\end{defn}

\begin{defn}[Morphism between representations]
\label{defn:Morphism_between_representations}
(See Br\"ocker and tom Dieck \cite[Chapter II, Definition 1.4, p. 67]{BrockertomDieck}.)
Let $\KK=\RR$ or $\CC$ and $V_1, V_2$ be $\KK$-vector spaces. A $\KK$-linear map $T:V_1\to V_2$ is a \emph{morphism between representations} $(\rho_1,V_1)$ and $(\rho_2,V_2)$ if $T\rho_1(g) =\rho_2(g)T$ for all $g\in G$; the two representations are \emph{equivalent} if $T$ is a $\KK$-linear isomorphism of vector spaces.
\end{defn}

Because every complex vector space is also a real vector space and because there is an inclusion 
$r_\RR^\CC:\GL_\CC(W)\to \GL_\RR(W)$, every complex representation $(W,\rho)$ of $G$ defines a real representation which we will denote as $r_\RR^\CC\rho$ (see Br\"ocker and tom Dieck \cite[p. 95]{BrockertomDieck}).

\begin{exmp}[Inequivalence of complex representations of $S^1$ defined by different weights]
\label{exmp:S1ComplexReps}
For each integer $m\in\ZZ$, the map
\[
  \rho_{\CC,m}:S^1 \ni e^{i\theta} \mapsto e^{im\theta} \in \GL_\CC(\CC)=\CC^*
\]
defines a complex representation, $\rho_{\CC,m}$. By the forthcoming Lemma \ref{lem:Isomorphic_S1_Representations}, if $n\in\ZZ$ is such that $m\neq n$, then $\rho_{\CC,m}$ is not equivalent to $\rho_{\CC,n}$.
\end{exmp}

\begin{exmp}[Real and complex representations of $S^1$]
\label{exmp:S1RealReps}
For each $m\in\ZZ$, the map
\[
  \rho_{\RR,m}:S^1 \ni e^{i\theta}
  \mapsto
\begin{pmatrix}
\cos(m\theta) & - \sin(m\theta)
\\
\sin(m\theta) & \cos(m\theta)
\end{pmatrix}
\in
\GL_\RR(\RR^2)
\]
defines a real representation, $\rho_{\RR,m}$.  If we define a linear isomorphism $T:\RR^2\to\CC$ of real vector spaces by
\[
T\begin{pmatrix} a \\ b \end{pmatrix}=a+bi\in\CC,
\]
then $T\rho_{\RR,m}(e^{i\theta}) T^{-1}=\rho_{\CC,m}(e^{i\theta})$ and so $\rho_{\RR^2,m}=r_\RR^\CC\rho_{\CC,m}$.
\end{exmp}

As we shall see in the forthcoming Lemma \ref{lem:Isomorphic_S1_Representations}, the map $\CC\ni v\mapsto \bar v\in\CC$ gives an equivalence between $r_\RR^\CC\rho_{\CC,m}$ and $r_\RR^\CC\rho_{\CC,-m}$ as real representations. Because the map $z\mapsto \bar z$ is not complex linear, however, it does not give an equivalence between $\rho_{\CC,m}$ and $\rho_{\CC,-m}$ as complex representations.

\begin{lem}[Equivalence of representations of $S^1$]
\label{lem:Isomorphic_S1_Representations}
The following hold:  
\begin{enumerate}
\item
\label{item:Isomorphic_S1_Representations_Complex}
The representations $(\CC,\rho_{\CC,m})$ and $(\CC,\rho_{\CC,n})$ in Example \ref{exmp:S1ComplexReps} are equivalent complex representations if and only if $m=n$.
\item
\label{item:Isomorphic_S1_Representations_Real}
The representations $(\RR^2,\rho_{\RR,m})$ and $(\RR^2,\rho_{\RR,n})$ in Example \ref{exmp:S1RealReps} are equivalent real representations if and only if $m=\pm n$.
\end{enumerate}
\end{lem}

\begin{proof}
We first prove Item \eqref{item:Isomorphic_S1_Representations_Complex}.
If $m=n$, then the identity map $\id:\CC\to \CC$ gives the equivalence between $(\CC,\rho_{\CC,m})$ and $(\CC,\rho_{\CC,n})$.  If $(\CC,\rho_{\CC,m})$ and $(\CC,\rho_{\CC,n})$ are equivalent, then there is a complex linear map $T:\CC\to \CC$ with $T\rho_{\CC,m}T^{-1}=\rho_{\CC,n}$.  Because $T$ is linear, there is a constant $\la\in\CC\less\{0\}$ such that $Tz=\la z$ for all $z\in\CC$. Hence, the identity $T\rho_{\CC,m}T^{-1}=\rho_{\CC,n}$ implies that for every $e^{i\theta}\in S^1$, one has $\la e^{im\theta}\la^{-1}=e^{in\theta}$ and thus, $m=n$ as asserted.

We now prove Item \eqref{item:Isomorphic_S1_Representations_Real}. As in the verification of Item \eqref{item:Isomorphic_S1_Representations_Complex}, the identity map gives an equivalence between $(\RR^2,\rho_{\RR,m})$ and itself.  To prove that $(\RR^2,\rho_{\RR,m})$ and $(\RR^2,\rho_{\RR,-m})$ are equivalent, choose
\[
T=
\begin{pmatrix} 1 & 0 \\ 0 & -1 \end{pmatrix},
\quad\text{so that}\quad
T^{-1}=T.
\]
Then
\begin{align*}
T \rho_{\RR,m}(e^{i\theta}) T^{-1}
&=
\begin{pmatrix} 1 & 0 \\ 0 & -1 \end{pmatrix}
\begin{pmatrix}
\cos(m\theta) & - \sin(m\theta)
\\
\sin(m\theta) & \cos(m\theta)
\end{pmatrix}
\begin{pmatrix} 1 & 0 \\ 0 & -1 \end{pmatrix}
\\
&=
\begin{pmatrix}
\cos(m\theta) & - \sin(m\theta)
\\
-\sin(m\theta) & -\cos(m\theta)
\end{pmatrix}
\begin{pmatrix} 1 & 0 \\ 0 & -1 \end{pmatrix}
\\
  &=
\begin{pmatrix}
\cos(m\theta) & \sin(m\theta)
\\
-\sin(m\theta) & \cos(m\theta)
\end{pmatrix}
\\
&=
\rho_{\RR,-m}(e^{i\theta}).
\end{align*}
Hence, for each $m\in\ZZ$, the representations $(\RR^2,\rho_{\RR,m})$ and $(\RR^2,\rho_{\RR,-m})$ are equivalent.

If $(\RR^2,\rho_{\RR,m})$ and $(\RR^2,\rho_{\RR,n})$ are equivalent, then let $T\in\GL_\RR(\RR^2)$ be the linear map with $T\rho_{\RR,m} T^{-1}=\rho_{\RR,n}$.  Because the kernel of the homomorphism $\rho_{\RR,n}:S^1\to\GL_\RR(\RR^2)$
is\label{page:CyclicGroup} $C_n := \{e^{i\theta} \in S^1: e^{in\theta} = 1\}$, the group of $|n|$-th roots of unity, and the kernel of $T\rho_{\RR,m} T^{-1}$ equals the kernel of $\rho_{\RR,m}$, the equality $T\rho_{\RR,m} T^{-1}=\rho_{\RR,n}$ implies that group of $|m|$-th roots of unity is equal to the group of $|n|$-th roots of unity.  Hence, $|m|=|n|$ and so $m=\pm n$ as asserted.  This completes the proof of Item \eqref{item:Isomorphic_S1_Representations_Real} and hence the proof of the lemma.
\end{proof}

\begin{defn}[Subrepresentations and irreducible representations]
\label{defn:SubrepresentationAndIrreducible}
(See Br\"ocker and tom Dieck \cite[Chapter II, Definition 1.8, p. 68]{BrockertomDieck}.)
Let $\KK=\RR$ or $\CC$. If $(V,\rho)$ is a $\KK$-representation and $W\subset V$ is a linear subspace such that $\rho(g)(W)=W$ for all $g\in G$, then $(W,\rho)$ is a \emph{subrepresentation} of $(V,\rho)$.
A representation $(V,\rho)$ is \emph{irreducible} if the only subrepresentations of $(V,\rho)$ are $(0)$ and $V$ itself.
\end{defn}

\begin{defn}[Direct sum of representations]
\label{defn:DirectSumOfRepresentations}
(See Br\"ocker and tom Dieck \cite[Chapter II, paragraph (1.12), p. 69]{BrockertomDieck}.)
The \emph{direct sum} of representations  $(V_1,\rho_1)$ and $(V_2,\rho_2)$ is the representation $(V_1\oplus V_2,\rho_1\oplus \rho_2)$.
\end{defn}

To describe when a representation can be written as a direct sum of irreducible representations, we introduce the 

\begin{defn}[Matrix Lie group]
\label{defn:MatrixLieGroup}
(See Hall \cite[Definition 1.4, p. 4]{Hall_lie_groups_algebras_representations}).
Let $\KK=\RR$ or $\CC$. A \emph{matrix Lie group} is a subgroup $G$ of $\GL_\KK(\KK^n)$ that is a closed subspace.
\end{defn}

We can then state the

\begin{prop}[Representations of compact matrix Lie groups]
\label{prop:Reps_of_Matrix_Lie_Groups_are_semisimple}
(See Br\"ocker and tom Dieck \cite[Chapter II, Propositions 1.9 and 1.14, p. 68 and p. 70]{BrockertomDieck} or Knapp \cite[Corollary 4.7, p. 240]{Knapp_1986} for $\KK=\CC$ and Hall \cite[Proposition 4.28, p. 92]{Hall_lie_groups_algebras_representations} for $\KK=\CC$ or $\RR$.)
Let $\KK=\RR$ or $\CC$. Let $(V,\rho)$ be a $\KK$-representation of a compact matrix Lie group $G$.  Then $(V,\rho)$ is equivalent to a direct sum of irreducible representations of $G$.
\end{prop}

We now describe the irreducible complex representations of the compact matrix Lie group $S^1$.

\begin{prop}[Irreducible complex representations of the circle group]
\label{prop:Irreducible_Complex_Representations_of_S1}
(See \cite[Chapter II, Proposition 8.1, p. 107]{BrockertomDieck}.)
If $(V,\rho)$ is an irreducible complex representation of $S^1$, then $(V,\rho)=(\CC,\rho_{\CC,m})$ for some integer $m\in\ZZ$, where $(\CC,\rho_{\CC,m})$ is the representation defined in Example \ref{exmp:S1ComplexReps}.
\end{prop}

\begin{prop}[Irreducible real representations of the circle group]
\label{prop:Irreducible_Real_Representations_of_S1}
(See Br\"ocker and tom Dieck \cite[Chapter II, Proposition 8.5, p. 109]{BrockertomDieck}.)
If $(V,\rho)$ is a non-trivial, irreducible real representation of $S^1$, then
$(V,\rho)=(\RR^2,\rho_{\RR,m})$ for some integer $m\in\ZZ$, where $(\RR^2,\rho_{\RR,m})$ is the representation defined in Example \ref{exmp:S1RealReps}.
\end{prop}

By combining Propositions \ref{prop:Reps_of_Matrix_Lie_Groups_are_semisimple} and \ref{prop:Irreducible_Complex_Representations_of_S1} or \ref{prop:Irreducible_Real_Representations_of_S1}, we see that any representation of $S^1$ is a direct sum of the representations given in Examples \ref{exmp:S1ComplexReps} or \ref{exmp:S1RealReps}. 

\begin{prop}[Classification of complex representations of $S^1$]
\label{prop:Direct_Sum_Decomposition_of_Complex_S1_Representations}
Every complex representation $(V,\rho)$ of $S^1$ is isomorphic to a direct sum
\[
V
\cong
\bigoplus_{j=1}^r V_j,
\]
where, for $j=1,\dots,r$, there is an integer $m_j$ such that $(V_j,\rho) \cong (\CC,\rho_{\CC,m_j})$. This direct sum decomposition of $(V,\rho)$ is unique up to a reordering of the integers $m_j$.
\end{prop}

\begin{proof}
The existence of such a decomposition follows from  Propositions \ref{prop:Reps_of_Matrix_Lie_Groups_are_semisimple} and \ref{prop:Irreducible_Complex_Representations_of_S1}. The uniqueness of the decomposition follows from Br\"ocker and tom Dieck \cite[Chapter II, Proposition 1.14, p. 70]{BrockertomDieck} or the discussion following Knapp \cite[Corollary 4.16, p. 243]{Knapp_1986}.
\end{proof}

\begin{defn}[Weights]
\label{defn:WeightsOfS1Representation}
(See  Br\"ocker and tom Dieck \cite[Chapter II, Definition 8.2, p. 108]{BrockertomDieck}.)
The \emph{weights} of a complex $S^1$ representation $(V,\rho)$ are the integers $m_j$ appearing in the decomposition of $(V,\rho)$ given in Proposition \ref{prop:Direct_Sum_Decomposition_of_Complex_S1_Representations}.
\end{defn}

\begin{prop}[Classification of real representations of $S^1$]
\label{prop:Direct_Sum_Decomposition_of_Real_S1_Representations}
Every real representation $(V,\rho)$ of $S^1$ is isomorphic to a direct sum
\[
V
\cong
V_0\oplus \bigoplus_{j=1}^r V_j,
\]
where, for $j=1,\dots,r$, there is an integer $m_j$ such that $(V_j,\rho)\cong (\RR^2,\rho_{\RR^2,m_j})$ and $V_0$ is equivalent to the direct sum of $\dim_\RR V- 2r$ copies of the trivial representation in Definition \ref{defn:TrivialRepresentation}. This direct sum decomposition is unique up to the signs of the integers $m_j$ and their ordering.
\end{prop}

\begin{proof}
The existence of the decomposition follows from Propositions \ref{prop:Reps_of_Matrix_Lie_Groups_are_semisimple} and \ref{prop:Irreducible_Real_Representations_of_S1}.
The uniqueness follows from Br\"ocker and tom Dieck \cite[Chapter II, Proposition 6.9, p. 101]{BrockertomDieck} or by applying Proposition \ref{prop:Direct_Sum_Decomposition_of_Complex_S1_Representations} to the complexification, $(V\otimes_\RR\CC,\rho\otimes \id_\CC)$, of $(V,\rho)$.  Note that by  \cite[Proposition 6.1]{BrockertomDieck}, the complexification of $(\RR^2,\rho_{\RR,m})$ is $\rho_{\CC,m}\oplus\rho_{\CC,-m}$, which explains why the uniqueness of the weights of the complexification given by Proposition \ref{prop:Direct_Sum_Decomposition_of_Complex_S1_Representations} only gives the uniqueness up to sign of the weights of the real representation.
\end{proof}

\section{Circle actions and fixed points}
\label{sec:CircleActionsFixedPtsHessians}
Let $M$ be a smooth manifold of dimension $d\ge 2$ and equipped with a smooth circle action \eqref{eq:Circle_action_smooth_manifold},
\[
\rho:S^1\times M\to M.
\]
We recall from \eqref{eq:Circle_action_tangent_bundle} that the induced circle action on the $TM$ bundle,
\[
\rho_*:S^1\times TM\to TM,
\]
is induced by the differential \eqref{eq:rho_*_definition} of $\rho$ in directions tangent to $M$, 
\[
  \rho_*(e^{i\theta})v := D_2\rho(e^{i\theta},p)v, \quad\text{for all } e^{i\theta}\in S^1 \text{ and } v \in T_pM.
\]
We see that $\rho_*$ defines an action of $S^1$ on $TM$ as follows. First, $\rho_*(1)$ is the identity on $TM$ because $\rho(1,\cdot)$ is the identity on $M$. Then,
\begin{align*}
  \rho_*(e^{i\theta})\rho_*(e^{i\theta'})v
  &=
  D_2\rho(e^{i\theta},\rho(e^{i\theta'},p))D_2\rho(e^{i\theta'},p)v \quad\text{(by \eqref{eq:rho_*_definition})}
  \\
  &= D_2\rho(e^{i\theta},\rho(e^{i\theta'},p))v  \quad\text{(by Chain Rule)}
  \\
  &= D_2\rho(e^{i\theta}e^{i\theta'},p)v  \quad\text{(because $\rho$ is an $S^1$ action)}
  \\
  &= \rho_*(e^{i\theta}e^{i\theta'})v, \quad\text{for all } v\in T_pM,
\end{align*}
and thus we have
\begin{equation}
\label{eq:DerivativeOfS1DefinesGroupAction}
\rho_*(e^{i\theta})\rho_*(e^{i\theta'})v
=
\rho_*(e^{i\theta}e^{i\theta'})v, \quad\text{for all } v\in T_pM.
\end{equation}
We now recall a useful special case of a more general result --- see Audin \cite[Corollary 1.2.3]{Audin_torus_actions_symplectic_manifolds}, Bredon \cite[Corollary 6.2.5]{Bredon}, or Duistermaat and Kolk \cite[Theorem 2.2.1 and comments following Definition 2.6.3]{DuistermaatLieGroups} --- on the structure of the subspace of points fixed by the action of a compact Lie group, noting that the fixed-point set of the $S^1$ action is the subspace of points of orbit type $S^1/S^1$ in the terminology of \cite{Audin_torus_actions_symplectic_manifolds, Bredon}.

\begin{thm}[Fixed-point set of a circle action on a smooth manifold is a smooth submanifold]
\label{thm:Fixed-point_circle_action_is_submanifold}
Let $M$ be a smooth manifold endowed with a smooth circle action. Then the set of fixed points of the circle action is an embedded smooth submanifold.
\end{thm}

We now consider the action induced by \eqref{eq:Circle_action_tangent_bundle} on the tangent space of the fixed-point set of $\rho$ on $M$. By Theorem \ref{thm:Fixed-point_circle_action_is_submanifold}, each connected component $F_\alpha \subset M$ of the fixed-point set
\begin{equation}
\label{eq:Fixed-point_set_circle_action_smooth_manifold}
F := \{p\in M:\rho(e^{i\theta},p) = p, \text{for all } e^{i\theta}\in S^1\}
\end{equation}
is a smooth submanifold, possibly of different dimensions for different components. When $p\in F$, the action $\rho_*$ in \eqref{eq:Circle_action_tangent_bundle} gives a circle action on $T_pM$,
\begin{equation}
\label{eq:FixedPointS1Action}
\rho_*:S^1\times T_pM\to T_pM.
\end{equation}
For $p\in F$, we obtain a linear map $\rho_*(e^{i\theta}):T_pM\to T_pM$ and so by \eqref{eq:DerivativeOfS1DefinesGroupAction} the map $\rho_*$ is a real representation of $S^1$ on $T_pM$.
By Proposition \ref{prop:Direct_Sum_Decomposition_of_Real_S1_Representations}, this representation can be written as a direct sum of representations,
\[
T_pM = L_0 \oplus \bigoplus_{j=1}^n L_j,
\]
where $L_0$ is a direct sum of $(d-2n)$ copies of the trivial representation and $L_j$ is equivalent to the representation $(\RR^2,\rho_{\RR,m_j})$ described in Example \ref{exmp:S1RealReps}.  As noted in Item \eqref{item:Isomorphic_S1_Representations_Real} of Lemma \ref{lem:Isomorphic_S1_Representations}, the representations $(\RR^2,\rho_{\RR,m_j})$ and $(\RR^2,\rho_{\RR,-m_j})$ are equivalent, so the representation  $(T_pM,\rho_*)$ does not determine the signs of the integers $m_1,\dots,m_n$.   
The eigenvalues of the Hessian of a Hamiltonian function have no sign ambiguity, so the identification of these eigenvalues with the weights of the $S^1$ action in Item \eqref{item:Frankel_almost_Hermitian_WeightsAreEigenvalues} of Theorem \ref{mainthm:Frankel_almost_Hermitian} requires the weights of the $S^1$ action to be well-defined.
We now show how the choice of an almost complex structure on $T_pM$ determines the signs of the weights by applying  Proposition \ref{prop:Direct_Sum_Decomposition_of_Complex_S1_Representations}.

\begin{lem}[Weights of circle representation at fixed determined by almost complex structure]
\label{lem:S1_Invariant_J_Is_MultiplicationBy_i}
Let $M$ be a smooth manifold and let $p\in M$ be a fixed point of a smooth $S^1$ action $\rho:S^1\times M\to M$. If $J$ is an $S^1$-invariant almost complex structure on $M$ in the sense of \eqref{eq:Circle_invariant_(1,1)-tensor}, that is,
\[
\rho_*(e^{i\theta})J = J\rho_*(e^{i\theta}), \quad\text{for all } e^{i\theta} \in S^1,
\]
and if we define complex multiplication $\CC\times T_pM \to T_pM$ via
\begin{equation}
\label{eq:J_Defines_ComplexMult_On_Lj}
iv:=Jv,
\end{equation}
then the representation $(T_pM,\rho_*)$ admits a decomposition,
\begin{equation}
\label{eq:FixedPointTangentSpaceDecomp}
T_pM =  \bigoplus_{j=0}^n L_j,
\end{equation}
as a direct sum of irreducible, complex representations of $S^1$. Hence, there are integers $m_0,\dots,m_n$ such that
\begin{equation}
\label{eq:WeightOfS1Action}
  \rho_*(e^{i\theta})v = e^{i m_j\theta}v, \quad\text{for all } v\in L_j \text{ and } j = 0,\ldots,n.
\end{equation}
In addition, we have
\begin{equation}
\label{eq:DerivativeOfS1Action}
\left.\frac{d}{d\theta}\rho_*(e^{i\theta})v_j\right|_{\theta=0}
= i m_j v_j.
\end{equation}
\end{lem}

\begin{proof}
Because $J$ satisfies \eqref{eq:Circle_invariant_(1,1)-tensor}, the action $\rho_*$ in \eqref{eq:FixedPointS1Action} on $T_pM$ defines a \emph{complex} representation of $S^1$ on $T_pM$ for each $p\in F$, where $F$ is the fixed-point set of the $S^1$ action. Proposition \ref{prop:Direct_Sum_Decomposition_of_Complex_S1_Representations} implies that the representation
$(T_pM,\rho_*)$ admits a decomposition of the form \eqref{eq:FixedPointTangentSpaceDecomp}
where the subspaces $L_j$ are complex subspaces with respect to the almost complex structure $J$.
The existence of the integers $m_0,\dots,m_n$ satisfying
equation \eqref{eq:WeightOfS1Action} follows from Proposition \ref{prop:Direct_Sum_Decomposition_of_Complex_S1_Representations} and the identification of multiplication by $i$ with $J$.  Equation \eqref{eq:DerivativeOfS1Action} follows by differentiating both sides of Equation \eqref{eq:WeightOfS1Action}.
\end{proof}

\begin{defn}[Weight of a circle action compatible with an almost complex structure]
\label{defn:Weight_Associated_To_Complex_Structure}
Let $L_j\subset T_pM$ be one of the subspaces appearing in the decomposition \eqref{eq:FixedPointTangentSpaceDecomp} and let $J$ be an almost complex structure on $L_j$, so we can write $iv := Jv$ for $v\in L_j$ as in \eqref{eq:DerivativeOfS1Action}.  Then the integer $m_j$ in \eqref{eq:WeightOfS1Action} is the \emph{weight of the $S^1$ action on $L_j$}\label{page:weight_circle_action} with \emph{sign compatible with the almost complex structure $J$}.
\end{defn}

\begin{rmk}[Role of the almost complex structure in determining the signs of the weights]
\label{rmk:UseOfChoiceOfComplexStructure}
Without the requirement of compatibility with an almost complex structure $J$, the representation $(T_pM,\rho_*)$ is only a real representation and so the weights of the $S^1$ action on $T_pM$ are only defined up to sign, as shown by the equivalence of the real representations $\rho_{\RR,m}$ and $\rho_{\RR,-m}$ given in Item \eqref{item:Isomorphic_S1_Representations_Real} of Lemma \ref{lem:Isomorphic_S1_Representations}. However, the equality between the weights of a circle action with the (well-defined) eigenvalues of the Hessian of $f$ appearing in Item \eqref{item:Frankel_almost_Hermitian_WeightsAreEigenvalues} of Theorem \ref{mainthm:Frankel_almost_Hermitian} requires the weights of a circle action to be defined without sign ambiguity.  This equality is proved in the forthcoming Lemma \ref{lem:IdentifyEigenvalues}. In the proof of that lemma, we use the almost complex structure $J$ in Equation \eqref{eq:MomentMapGradientS1VField}, which relates the gradient of $f$ with the vector field generated by the circle action, and in the derivation of Equation \eqref{eq:DerivativeOfS1Action} in the penultimate step of the proof.
\end{rmk}

The following characterization of the weights will be useful when comparing them with the eigenvalues of the Hessian of a Hamiltonian function for the $S^1$ action. The characterization is implicit in Frankel's argument in \cite[p. 3]{Frankel_1959}, but we include the proof in order to show how it generalizes from the case of a complex K\"ahler manifold.

\begin{lem}
\label{lem:MultiplicitiesEigenvalues}
Let $M$ be a smooth manifold and $\rho:S^1\times M\to M$ be a smooth circle action. Let $X$ be the smooth vector field on $M$ induced by $\rho$, thus
\begin{equation}
  \label{eq:VFieldGenByS1Action}
  X_q := D_1\rho(0,q) = \left.\frac{d}{d\theta}\rho(e^{i\theta},q)\right|_{\theta=0} \in T_qM, \quad\text{for all } q \in M,
\end{equation}
where $D_1\rho$ denotes the differential of $\rho$ in the $S^1$ direction. Then the following hold:
\begin{enumerate}
\item
\label{item:MultiplicitiesEigenvaluesFixedPointsAreZeros}
$p\in M$ is a fixed point of the $S^1$ action $\rho$ if and only if $X_p=0$.
\item
  \label{item:MultiplicitiesEigenvalues_IdentifyEigenvalues}
  If $p\in M$ is a fixed point of the $S^1$ action $\rho$, then the linear map $(\sL_X)_p \in \End(T_pM)$ defined in Lemma \ref{lem:Action_vector_field_at_a_zero} satisfies
\begin{equation}
  \label{eq:EigenvectorsOfH_p_endomorphism_tangent_space_fixed_point}
  (\sL_X)_p v=-\left.\frac{d}{d\theta}\rho_*(e^{i\theta})v\right|_{\theta=0},
  \quad\text{for all $v\in T_pM$}.
\end{equation}
\end{enumerate}
\end{lem}

\begin{proof}
The identification of fixed points with zeros of the vector field $X$ in Item \eqref{item:MultiplicitiesEigenvaluesFixedPointsAreZeros} follows from \cite[Lemma 6.1.1]{Bredon}, as we now describe. If $p$ is a fixed point of the $S^1$ action, then $X_p=0$ by the definition of $X$ in \eqref{eq:VFieldGenByS1Action}. Conversely, if $X_p=0$, then the constant curve with value $p$ is an integral curve of the vector field $X$ with initial point $p$. The curve $\RR\ni \theta \mapsto \rho(e^{i\theta},p) \in M$ is also an integral curve of $X$ with initial point $p$. The uniqueness of integral curves for $X$ (see Lee \cite[Theorem 12.9]{Lee_john_smooth_manifolds}) with a given initial point then implies that $p=\rho(e^{i\theta},p)$ for all $\theta\in\RR$ and so $p$ is a fixed point of the action $\rho$. This proves Item \eqref{item:MultiplicitiesEigenvaluesFixedPointsAreZeros}.

We now prove Item \eqref{item:MultiplicitiesEigenvalues_IdentifyEigenvalues}.
If $\tilde v$ is any smooth vector field on $M$ that extends $v\in T_pM$, then by Gallot, Hulin, and Lafontaine \cite[Theorem 1.68]{Gallot_Hulin_Lafontaine_riemannian_geometry}, Spivak \cite[pp. 150--153]{Spivak1}, or Lee \cite[Equation (9.16)]{Lee_john_smooth_manifolds} we obtain
\[
\left.\frac{d}{d\theta}\rho_*(e^{i\theta}) v\right|_{\theta=0}
=
[\tilde v,X]_p = -(\sL_X)_p v,
\]
which yields \eqref{eq:EigenvectorsOfH_p_endomorphism_tangent_space_fixed_point} and proves the lemma.
\end{proof}

\section{Compatible triples}
\label{sec:CompatTriples}
We now describe how the existence of an $S^1$-invariant, non-degenerate two-form $\om$ as assumed in Theorem \ref{mainthm:Frankel_almost_Hermitian} gives the $S^1$-invariant almost complex structure we used to define the weights of the $S^1$ action without sign ambiguity in Lemma \ref{lem:S1_Invariant_J_Is_MultiplicationBy_i}. Recall from \eqref{eq:Circle_invariant_covariant_2-tensor} that a non-degenerate two-form $\om$ is $S^1$-invariant if
\[
\om\left(\rho_*(e^{i\theta})v, \rho_*(e^{i\theta})w\right) = \om(v,w)
  \quad\text{for all } p \in M, v, w \in T_pM, \text{ and } e^{i\theta} \in S^1,
\]
while a Riemannian metric $g$ is $S^1$-invariant (see \eqref{eq:Circle_invariant_covariant_2-tensor}) if
\[
g\left(\rho_*(e^{i\theta})v, \rho_*(e^{i\theta})w\right) = g(v,w)
  \quad\text{for all } p \in M, v, w \in T_pM, \text{ and } e^{i\theta} \in S^1.
\]
According to Bredon \cite[Theorem 6.2.1]{Bredon} (see also Audin \cite[Proof of Theorem 1.2.1]{Audin_torus_actions_symplectic_manifolds}), given a Riemannian metric $g_0$ on $M$, one can define an $S^1$-invariant Riemannian metric $g$ by averaging over $S^1$:
\begin{equation}
\label{eq:Construction_Of_S1_Invariant_Metric}
g(v,w)
=
\int_{S^1} g_0(\rho_*(e^{i\theta})v,\rho_*(e^{i\theta}) w)\ d\theta,
\end{equation}
where $d\theta$ is the Haar measure on $S^1$.

If $\om$ is a non-degenerate two-form, then following McDuff and Salamon \cite[Equation (4.1.2), p. 153]{McDuffSalamonSympTop3}, we say that an almost complex structure $J$ is \emph{$\om$-compatible} if
\begin{equation}
\label{eq:omega_compatible_with_J}
\om(JX,X)>0
\quad\text{and}\quad
\om(JX,JY)=\om(X,Y),\quad\text{for all } X,Y\in TM.
\end{equation}
A \emph{(smooth) compatible triple}\label{page:Compatible_triple} $(\om,g,J)$ on $M$ (see McDuff and Salamon \cite[Section 4.1]{McDuffSalamonSympTop3}) comprises a non-degenerate, smooth two-form $\om$, a smooth Riemannian metric $g$, and a smooth almost complex structure $J$ satisfying
\eqref{eq:g_compatible_J} and \eqref{eq:Fundamental_two-form}, that is,
\begin{equation}
\label{eq:compatible_triple}
\om(X,Y)=g(X, JY)
\quad\text{and}\quad
g(JX,JY)=g(X,Y), \quad\text{for all } X,Y\in C^\infty (TM).
\end{equation}
Regarding \eqref{eq:compatible_triple}, we note the

\begin{rmk}[On the definition of compatible triples]
\label{rmk:CompatibilityConventionExplanation}
Our definition of $\om$-compatible differs from that of McDuff and Salamon \cite[Equation (4.1.2), p. 153]{McDuffSalamonSympTop3} and Cannas da Silva \cite[Definition 12.2, p. 68]{Cannas_da_Silva_lectures_on_symplectic_geometry} where they require $\om(X,JX)>0$.  We follow Kobayashi \cite[Equation (7.6.8), p. 251]{Kobayashi_differential_geometry_complex_vector_bundles} in saying that the metric $g$ and two-form $\om$ are compatible if $\om(X,Y)=g(X,JY)$ in \eqref{eq:Fundamental_two-form}. The convention  $\om(X,Y)=g(X,JY)$ from  \eqref{eq:Fundamental_two-form} implies that $\om(X,JX)=g(X,J^2 X)=-g(X,X)$ while $\om(JX,X)=g(JX,JX)=g(X,X)$.  Thus, we have selected the assumption of $\om(JX,X)>0$ in our definition of $\om$-compatibility to be consistent with the convention \eqref{eq:Fundamental_two-form} from Kobayashi \cite[Equation (7.6.8), p. 251]{Kobayashi_differential_geometry_complex_vector_bundles}.
\end{rmk}

We then have the

\begin{lem}[Existence of circle-invariant compatible triples]
\label{lem:ExistenceOfS1InvarACStructure}
Let $M$ be a smooth manifold with  a non-degenerate two-form $\om$.  
\begin{enumerate}
\item
\label{item:ExistenceOfS1InvarACStructure_CompatTriple}
There is a smooth almost complex structure $J_\om$ and a Riemannian metric $g_\om$ on $M$ such that
$(\om,g_\om,J_\om)$ is a compatible triple.
\item
\label{item:ExistenceOfS1InvarACStructure_S1Invariant}
If $M$ has a smooth $S^1$ action  and $\om$ is $S^1$-invariant, then the almost complex structure $J_\om$ and Riemannian metric $g_\om$ comprising the compatible triple
in Item \eqref{item:ExistenceOfS1InvarACStructure_CompatTriple} are $S^1$-invariant in the sense that \eqref{eq:Circle_invariant_(1,1)-tensor} holds for $J_\om$ and
\eqref{eq:Circle_invariant_covariant_2-tensor} holds for $g_\om$.
\end{enumerate}
\end{lem}

\begin{proof}
Item \eqref{item:ExistenceOfS1InvarACStructure_CompatTriple} follows from McDuff and Salamon \cite[Proposition 2.5.6, p. 67, and Proposition 4.1.1 (iii), p. 153]{McDuffSalamonSympTop3} (see also
Cannas da Silva \cite[Proposition 12.3, p. 68 and Proposition 12.6, p. 70]{Cannas_da_Silva_lectures_on_symplectic_geometry}).  We include details of their construction in order to prove the $S^1$ invariance appearing in Item \eqref{item:ExistenceOfS1InvarACStructure_S1Invariant}. Let $g$ be a smooth Riemannian metric on $M$. From the proof of \cite[Proposition 12.3, p. 68]{Cannas_da_Silva_lectures_on_symplectic_geometry} or by \cite[Step 1 of Proposition 2.5.6, p. 67]{McDuffSalamonSympTop3}, there is a smooth section $A$ of $\End(TM)$ satisfying
\begin{equation}
\label{eq:Omega_and_g_Define_A}
\om(X,Y) = g(X,AY), \quad\text{for all smooth vector fields $X$ and $Y$}.
\end{equation}
We first show that $A$ is skew-symmetric with respect to $g$.  For any smooth vector fields $X$ and $Y$,
\[
g(A^*X,Y)
=
g(X,AY)
=
\om(X,Y)
=
-\om(Y,X)
=
- g(Y,AX)
=
-g(AX,Y).
\]
Thus, $A^*=-A$ and so $A$ is skew-symmetric with respect to $g$.  Because $\om$ is non-degenerate,
\eqref{eq:Omega_and_g_Define_A} implies that $A_x:T_xM\to T_xM$  is injective for all $x\in M$. The $g$-symmetric endomorphism $AA^*$ is positive definite because $g(AA^*X,X)=g(A^*X,A^*X)=g(AX,AX)\ge 0$ for all $X$ and if $g(A_xX_x,A_xX_x) = 0$ at a point $x\in M$, then $A_xX_x = 0$ and thus $X_x = 0$ since $A_x$ is injective. Because $AA^*$ is $g$-symmetric and positive definite, the endomorphism $(AA^*)^{-1/2}$ is well-defined
(see  Bhatia \cite[p. 5]{Bhatia_matrix_analysis}, Rudin \cite[Section 12.24]{Rudin}, or Kato \cite[Equation (5.47)]{Kato}), commutes with $A$ and $A^*$ (see \cite[Section 12.24]{Rudin} or \cite[Equation (5.49)]{Kato}) and $g$-symmetric \cite[Section 12.24, Equation (2)]{Rudin}.

We define a candidate for an almost complex structure $J_\om\in \End(TM)$ by
\begin{equation}
\label{eq:Define_Compatible_J}
J_\om := (AA^*)^{-1/2} A.
\end{equation}
Because $A$ is skew-symmetric and $(AA^*)^{-1/2}$ is symmetric, we observe that $J_\om$ is skew-symmetric with respect to $g$. We see that $J_\om^2=-\id_{TM}$ by computing
\begin{align*}
J_\om^2
&=
-J_\om J_\om^*
\quad\text{(because $J_\om$ is skew-symmetric)}
\\
&=
-\left( (AA^*)^{-1/2} A\right)\left( (AA^*)^{-1/2} A\right)^*
\quad\text{(by definition \eqref{eq:Define_Compatible_J})}
\\
&=
-(AA^*)^{-1/2} A A^* (AA^*)^{-1/2}
\quad\text{(because $(AA^*)^{-1/2}$ is $g$-symmetric)}
\\
&=
-\id_{TM}
\quad\text{(by \cite[Section 12.24]{Rudin} or \cite[Equation (5.49)]{Kato}).}
\end{align*}
We see that $J_\om$ is $g$-orthogonal by computing, for all smooth vector fields $X$ and $Y$,
\[
g(J_\om X,J_\om Y)
=
g(J_\om^*J_\om X,Y)
=
g(-J_\om^2 X,Y)
=
g(X,Y).
\]
We now show that $J_\om$ is $\om$-compatible. First, we observe that
\begin{align*}
\om(J_\om X,J_\om Y)
&=
g(J_\om X,AJ_\om Y)\quad\text{(by \eqref{eq:Omega_and_g_Define_A})}
\\
&=
g(J_\om X,J_\om AY)\quad\text{(because $A$ commutes with $(AA^*)^{-1/2}$)}
\\
&=
g(X,AY)\quad\text{(because $J_\om$ is $g$-orthogonal)}
\\
&=
\om(X,Y)\quad\text{(by \eqref{eq:Omega_and_g_Define_A})},
\end{align*}
and thus,
\begin{equation}
\label{eq:J_and_Om_Compatibility_JInvariant}
\om(J_\om X,J_\om Y)
=
\om(X,Y), \quad\text{for all } X, Y \in C^\infty(TM).
\end{equation}
We continue to verify that $J_\om$ is $\om$-compatible by computing
\begin{align*}
\om(J_\om X,Y)
&=
g(J_\om X,AY)\quad\text{(by \eqref{eq:Omega_and_g_Define_A})}
\\
&=
g(A^* (AA^*)^{-1/2} A X,Y)
\quad\text{(by definition \eqref{eq:Define_Compatible_J})}
\\
&=
g((AA^*)^{1/2}X,Y) \quad\text{(because $A$ commutes with $(AA^*)^{-1/2}$),}
\end{align*}
that is,
\begin{equation}
\label{eq:J_and_gJ_Compatibility_DefiningMetric}
\om(J_\om X,Y)=g((AA^*)^{1/2}X,Y), \quad\text{for all } X, Y \in C^\infty(TM).
\end{equation}
Because $(AA^*)^{1/2}$ is positive definite and symmetric with respect to $g$, then
\eqref{eq:J_and_gJ_Compatibility_DefiningMetric} implies that $\om(J_\om X,X)>0$ for all tangent vectors $X$.
Combining this observation with \eqref{eq:J_and_Om_Compatibility_JInvariant} implies that $J_\om$ is $\om$-compatible in the sense of \eqref{eq:omega_compatible_with_J}. Moreover, the tensor
\begin{equation}
\label{eq:Define_Compatible_gJ}
g_\om(\cdot,\cdot):= g((AA^*)^{1/2}\cdot,\cdot)
\end{equation}
defines a Riemannian metric satisfying (by \eqref{eq:J_and_gJ_Compatibility_DefiningMetric})
\begin{equation}
\label{eq:J_and_gJ_Compatibility_2}
\om(J_\om X,Y)=g_\om(X,Y), \quad\text{for all } X, Y \in C^\infty(TM).
\end{equation}
Thus, $J_\om$ is $\om$-compatible and if we define $g_\om$ by \eqref{eq:Define_Compatible_gJ}, then $\om(J_\om X,Y)=g_\om(X,Y)$ and so the triple $(\om,g_\om,J_\om)$ is compatible in the sense of \eqref{eq:compatible_triple}. This completes the proof of Item \eqref{item:ExistenceOfS1InvarACStructure_CompatTriple}.

We now prove Item \eqref{item:ExistenceOfS1InvarACStructure_S1Invariant}. By \eqref{eq:Construction_Of_S1_Invariant_Metric}, there is an $S^1$-invariant Riemannian metric $g$ on $M$. We see that the section $A$ of $\End(TM)$ defined in \eqref{eq:Omega_and_g_Define_A} by $g$ and $\om$ is $S^1$-invariant by computing
\begin{align*}
g(X,\rho_*(e^{i\theta})A \rho_*(e^{-i\theta}) Y)
  &=
g(\rho_*(e^{-i\theta}) X,A \rho_*(e^{-i\theta}) Y)\quad\text{(by $S^1$ invariance of $g$)}
\\
&=
\om(\rho_*(e^{-i\theta})X, \rho_*(e^{-i\theta})Y)\quad\text{(by \eqref{eq:Omega_and_g_Define_A})}
\\
&=
\om (X,Y)\quad\text{(by $S^1$ invariance of $\om$)}
\\
&=
g(X,AY)\quad\text{(by \eqref{eq:Omega_and_g_Define_A}).}
\end{align*}
Because we have shown
\[
g(X,\rho_*(e^{i\theta})A \rho_*(e^{-i\theta}) Y)
=
g(X,AY)
\]
for all smooth vector fields $X$ and $Y$, the section $A$ satisfies
\begin{equation}
\label{eq:S1_Invariance_of_A}
\rho_*(e^{i\theta})A\rho_*(e^{-i\theta}) = A
\end{equation}
and is thus $S^1$-invariant. Because $A$ is skew-symmetric with respect to $g$, then $A^*=-A$ is also $S^1$-invariant and so
\begin{equation}
\label{eq:S1_Invariance_of_AAt}
\rho_*(e^{i\theta})AA^* \rho_*(e^{-i\theta})=AA^*.
\end{equation}
We then compute
\begin{align*}
\rho_*(e^{i\theta}) J_\om \rho_*(e^{-i\theta})
&=
\rho_*(e^{i\theta}) (AA^*)^{-1/2} A \rho_*(e^{-i\theta})
\quad\text{(by \eqref{eq:Define_Compatible_J})}
\\
  &=
\rho_*(e^{i\theta}) (\rho_*(e^{2i\theta})AA^*\rho_*(e^{-2i\theta}))^{-1/2} \rho_*(e^{-i\theta}) A
\quad\text{(by \eqref{eq:S1_Invariance_of_A} and \eqref{eq:S1_Invariance_of_AAt})}
\\
&=
(AA^*)^{-1/2} A
\quad\text{(by \eqref{eq:DerivativeOfS1DefinesGroupAction})}
\\
&=
J_\om \quad\text{(by \eqref{eq:Define_Compatible_J})}
\end{align*}
and this proves that $J_\om$ is $S^1$-invariant in the sense of \eqref{eq:Circle_invariant_(1,1)-tensor}.

The $S^1$ invariance of the Riemannian metric $g_\om$ 
follows immediately from its definition in \eqref{eq:Define_Compatible_gJ} 
and the $S^1$ invariance of $g$ and $(AA^*)^{-1/2}$.
This completes the proof of Item \eqref{item:ExistenceOfS1InvarACStructure_S1Invariant} and hence the lemma.
\end{proof}

\begin{rmk}[Almost complex structure as alternative to initial choice of non-degenerate two-form]
\label{rmk:Almost_complex_structure_alternative_to_initial_choice_non-degenerate_two-form}
Lemma \ref{lem:ExistenceOfS1InvarACStructure} hypothesizes a non-degenerate smooth two-form $\omega$ and, given any smooth Riemannian metric $g$, constructs a smooth almost complex structure $J$ and smooth Riemannian metric $g_J$ such that $(\omega,g_J,J)$ is a compatible triple. Alternatively, we could begin with a choice of smooth almost complex structure $J$ and, given any smooth Riemannian metric $g$, construct a compatible non-degenerate smooth two-form $\omega_{g,J}$. Indeed, given such a $g$, we may construct a  smooth Riemannian metric $g_J$ such that $J$ is orthogonal with respect to $g_J$ by averaging as in the construction of the Riemannian metric $g_0$ given by McDuff and Salamon in \cite[Proposition 2.5.6, Step 2, p. 69]{McDuffSalamonSympTop3}:
\[
  g_J(X,Y) := \frac{1}{2}\left(g(X,Y) + g(JX,JY)\right), \quad\text{for all } X, Y \in C^\infty(TM).
\]
We may now use \eqref{eq:Fundamental_two-form} to define a non-degenerate smooth two-form on $M$,
\[
  \omega_{g,J}(X,Y) := g_J(X,JY), \quad\text{for all } X, Y \in C^\infty(TM),
\]
and observe that $(\omega_{g,J},g_J,J)$ is a compatible triple in the sense of \eqref{eq:compatible_triple}. Moreover, if $J$ and $g$ are circle-invariant, then $\omega_{g,J}$ is necessarily circle-invariant too.
\end{rmk}  

Thus, given a non-degenerate, smooth two-form $\om$ on a smooth manifold $M$ with a smooth $S^1$ action $\rho$, we can assume that there are an $S^1$-invariant Riemannian metric and an $S^1$-invariant almost complex structure $J$ so that $(M,g,J)$ is an almost Hermitian manifold. 

\section{Hessians of Hamiltonian functions}
\label{sec:HessiansAndMomentMaps}
We now combine Lemmas \ref{lem:MultiplicitiesEigenvalues} and \ref{lem:HessianIsBracket} and identify the eigenvalues of the Hessian of $f$ at a critical point with the weights of the circle action.
The proof of the following lemma is implicit in Frankel \cite[p. 3]{Frankel_1959} and explicit in the proof due to Garc\'{\i}a--Prada, Gothen, and Mu\~{n}oz of \cite[Proposition 3.3]{Garcia-Prada_Gothen_Munoz_2007}, where they consider an $S^1$ action on a complex K\"ahler manifold given by the moduli space of complex rank-three Higgs pairs over a Riemann surface.

\begin{lem}[Critical points of Hamiltonian functions and eigenvalues of their Hessian operators]
\label{lem:IdentifyEigenvalues}
Let $M$ be a smooth manifold with a smooth $S^1$-action $\rho:S^1\times M\to M$ as in \eqref{eq:VFieldGenByS1Action} and smooth non-degenerate two-form $\om$. If $f:M\to \RR$ is a smooth Hamiltonian function for the $S^1$ action as in \eqref{eq:MomentMap}, where $X$ is the smooth vector field on $M$ generated by $\rho$, then the following holds:
\begin{enumerate}
\item
\label{item:IdentifyEigenvaluesCritPointsAreFixedPoints}
$p$ is a critical point of $f$ if and only if it is a fixed point of the $S^1$ action.
\end{enumerate}
If $g$ is a Riemannian metric and $J$ an almost complex structure on $M$ with the property that $(\om,g,J)$ is a compatible triple, then
\begin{equation}
\label{eq:MomentMapGradientS1VField}
\grad_g f = -JX.
\end{equation}
In addition, if the elements $\om$, $g$, and $J$ of the compatible triple $(\om,g,J)$ are $S^1$-invariant, then the following hold:
\begin{enumerate}
\setcounter{enumi}{1}
\item
\label{item:IdentifyEigenvaluesEigenvaluesOfHessianAreWeights}
The eigenvalues of the Hessian of $f$ at $p$ are equal to the weights of the $S^1$ action on $T_pM$ given by Lemma \ref{lem:MultiplicitiesEigenvalues}, with sign compatible with $J$ in the sense of Definition \ref{defn:Weight_Associated_To_Complex_Structure}.
\item
\label{item:IdentifyCritPointsAreFixedPoints_EigenspacesOfHessiansAreL_j}
If $p$ is a critical point of $f$ and $v_j\in L_j$, where $L_j$ is the subspace of $T_pM$ appearing in the decomposition \eqref{eq:FixedPointTangentSpaceDecomp}, then
\begin{equation}
  \label{eq:EigenvectorOfHessian}
  \Hess_g f(p) v_j= m_j v_j, \quad\text{for } j=0,1,\ldots,n,
\end{equation}
where $m_j$ is the weight of the $S^1$ action on $L_j$ as in \eqref{eq:WeightOfS1Action}, with sign compatible with $J$ in the sense of Definition \ref{defn:Weight_Associated_To_Complex_Structure}.
\end{enumerate}
\end{lem}

\begin{proof}
We first prove Item \eqref{item:IdentifyEigenvaluesCritPointsAreFixedPoints}. By the Hamiltonian condition \eqref{eq:MomentMap}, $df_p=(\iota_X\om)_p$. Because $\om$ is non-degenerate, $df_p=(\iota_X\om)_p=0$ if and only if $X_p=0$. Thus $p\in M$ is a critical point, so $df_p=0$, if and only if $X_p=0$. Item \eqref{item:MultiplicitiesEigenvaluesFixedPointsAreZeros} of Lemma \ref{lem:MultiplicitiesEigenvalues} asserts that $X_p = 0$ if and only if $p$ is a fixed point of the $S^1$ action. This completes the proof of Item \eqref{item:IdentifyEigenvaluesCritPointsAreFixedPoints}.

We now prove \eqref{eq:MomentMapGradientS1VField}. Let $(\om,g,J)$ be a compatible triple. Observe that for any smooth vector field $Y$ on $M$, we have
\begin{align*}
  g(\grad_g f,Y)
  &=
    df(Y) \quad\text{(by definition \eqref{eq:DefineGradient} of gradient)}
  \\
  &= (\iota_X\om)(Y) \quad\text{(by Hamiltonian condition \eqref{eq:MomentMap})}
  \\
  &= \om(X,Y)
  \\
  &= g(X,JY) \quad\text{(by compatibility condition \eqref{eq:Fundamental_two-form})}
  \\
  &= g(JX,J^2Y) \quad\text{(by compatibility condition \eqref{eq:g_compatible_J})}
  \\
  &= g(-JX,Y) \quad\text{(by definition of an almost complex structure $J$)},
\end{align*}
and thus, since $Y \in C^\infty(TM)$ was arbitrary, we see that \eqref{eq:MomentMapGradientS1VField} holds. 

Now we assume that the elements of $(\om,g,J)$ are $S^1$-invariant.
Item \eqref{item:IdentifyEigenvaluesEigenvaluesOfHessianAreWeights} follows from Item \eqref{item:IdentifyCritPointsAreFixedPoints_EigenspacesOfHessiansAreL_j}, which we now prove. For $j=0,1,\ldots,n$, let $\tilde v_j$ be any smooth vector field on $M$ that extends $v_j\in L_j\subset T_pM$.  Then
\begin{align*}
  \Hess_g f(p) v_j
  &=
  -[\grad_g f,\tilde v_j]_p \quad\text{(by \eqref{eq:GradientDefnOfHessian})}
  \\
  &
  =
  [JX,\tilde v_j]_p \quad\text{(by \eqref{eq:MomentMapGradientS1VField})}
  \\
  &=
  J[X,\tilde v_j]_p \quad\text{(by \eqref{eq:JandBracket})}
  \\
  &=
  J\left( -\left.\frac{d}{d\theta}\rho_*(e^{i\theta})v_j\right|_{\theta=0}\right)
    \quad\text{(by $\sL_X v_j=[X,\tilde v_j]$) and \eqref{eq:EigenvectorsOfH_p_endomorphism_tangent_space_fixed_point})}
  \\
  &=
  J\left( -im_j v_j\right) \quad\text{(by \eqref{eq:DerivativeOfS1Action})}
  \\
  &=
  m_jv_j \quad\text{(by \eqref{eq:J_Defines_ComplexMult_On_Lj})}.
\end{align*}
This proves equation \eqref{eq:EigenvectorOfHessian} and thus completes the proof of Item \eqref{item:IdentifyCritPointsAreFixedPoints_EigenspacesOfHessiansAreL_j} and the Lemma.
\end{proof}

\begin{cor}
\label{cor:Equality_Of_Hessian_and_CircleAction_Indices}
Let $M$ be a smooth manifold with a smooth $S^1$-action $\rho:S^1\times M\to M$ as in \eqref{eq:VFieldGenByS1Action}.
Assume that $J$ is an $S^1$-invariant almost complex structure on $M$. Let $p\in M$ be a fixed point of the action $\rho$ and let $m_1,\dots,m_n$ be the non-zero weights of the $S^1$ action on $T_pM$ as defined in Lemma \ref{lem:S1_Invariant_J_Is_MultiplicationBy_i} to be compatible with $J$. Let $m_1,\dots,m_r$ be the weights which are positive and $m_{r+1},\dots,m_n$ be the weights which are negative. Define  non-negative integers $m_p^+(\rho,J)$, $m_p^-(\rho,J)$, and $m_p^0(\rho,J)$ as the dimensions of the subspaces of $T_pM$ on which the weights of the action $\rho_*$ are, respectively, positive, negative, and zero:
\begin{equation}
\label{eq:Define Index of S1 Action At Fixed Point}
m_p^+(\rho,J):=2r,\quad
m_p^-(\rho,J):=2(n-r),\quad
m_p^0(\rho,J):=\dim M - 2n.
\end{equation}
If $\om$ is an $S^1$-invariant, non-degenerate two form on $M$, and $f:M\to \RR$ is a smooth Hamiltonian function for the $S^1$ action as in \eqref{eq:MomentMap}, and $(\om,g,J)$ is a compatible triple, then $m_p^-(\rho,J)$ equals the (Morse--Bott) index of $\hess f(p)$, the difference
\[
m_p^+(\rho,J)- m_p^-(\rho,J)
\]
equals the signature of $\hess f(p)$, and $m_p^0(\rho,J)$ equals the nullity of $\hess f(p)$.
\end{cor}

\begin{proof}
The conclusions follow immediately from  equation \eqref{eq:EigenvectorOfHessian}.
\end{proof}

We can now give the

\begin{proof}[Proof of Theorem \ref{mainthm:Frankel_almost_Hermitian}]
Item \eqref{item:Frankel_almost_Hermitian_FixedPointsAreZerosOfVField} of Theorem \ref{mainthm:Frankel_almost_Hermitian} follows  from Item \eqref{item:MultiplicitiesEigenvaluesFixedPointsAreZeros} of Lemma \ref{lem:MultiplicitiesEigenvalues}.

Item \eqref{item:Frankel_almost_Hermitian_FixedPointsAreCriticalPoints} of Theorem \ref{mainthm:Frankel_almost_Hermitian} follows from Item \eqref{item:IdentifyEigenvaluesCritPointsAreFixedPoints} of Lemma \ref{lem:IdentifyEigenvalues}.

Item \eqref{item:Frankel_almost_Hermitian_ACisS1Invariant} of Theorem \ref{mainthm:Frankel_almost_Hermitian} follows from
Lemma \ref{lem:ExistenceOfS1InvarACStructure}.

The equality between the weights of the $S^1$ action at a fixed point $p$ and the eigenvalues of the Hessian of $f$ at $p$ in Item \eqref{item:Frankel_almost_Hermitian_WeightsAreEigenvalues} of Theorem \ref{mainthm:Frankel_almost_Hermitian} follows from Item \eqref{item:IdentifyEigenvaluesEigenvaluesOfHessianAreWeights} of Lemma \ref{lem:IdentifyEigenvalues}.

To prove Item \eqref{item:Frankel_almost_Hermitian_ComponentsOfFixedPointsAreSmoothSubmanifolds} of Theorem \ref{mainthm:Frankel_almost_Hermitian}, we apply Theorem \ref{thm:Fixed-point_circle_action_is_submanifold} to conclude that each connected component of the fixed-point set $F$ is a smooth submanifold of $M$. We prove the remaining assertions in Item \eqref{item:Frankel_almost_Hermitian_ComponentsOfFixedPointsAreSmoothSubmanifolds} after proving Item \eqref{item:Frankel_almost_Hermitian_f_is_MB}.

We now prove Item \eqref{item:Frankel_almost_Hermitian_f_is_MB} of Theorem \ref{mainthm:Frankel_almost_Hermitian}.  By Item \eqref{item:Frankel_almost_Hermitian_FixedPointsAreCriticalPoints}, the subset of critical points $\Crit f \subset M$ is equal to the subset of fixed points $F\subset M$ of the circle action.  The set $F$ is a smooth submanifold by Item \eqref{item:Frankel_almost_Hermitian_ComponentsOfFixedPointsAreSmoothSubmanifolds} of Theorem \ref{mainthm:Frankel_almost_Hermitian}. Let $p\in M$ be a fixed point of the $S^1$ action. By Bredon \cite[p. 306]{Bredon}, the exponential map from $T_pM$ into $M$ defined by an $S^1$-invariant metric is equivariant with respect to the $S^1$ action given by $\rho_*$ on $T_pM$ and the $S^1$ action $\rho$ on $M$.  Hence, this exponential map takes the subspace $L_0$ of $T_pM$, defined in \eqref{eq:FixedPointTangentSpaceDecomp} and on which $\rho_*$ acts as the identity, to the fixed-point submanifold.  Because the derivative of the exponential map at the origin in $T_pM$ is the identity endomorphism of $T_pM$, this derivative identifies the subspace $L_0$  with the tangent space $T_pF$ to the fixed-point submanifold $F$ at the point $p$. The identification of the fixed-point set of the circle action with $\Crit f$ then implies that $T_p\Crit f=L_0$.  The identification of the eigenspaces of the Hessian of $f$ at $p$ with the bundles $L_j$ in the decomposition of $T_pM$ defined in \eqref{eq:FixedPointTangentSpaceDecomp} implies that $L_0$ is the kernel of $\Hess_g f(p)$ and, in particular, that $T_p\Crit f = \Ker\Hess_g f(p)$. Hence, $f$ is Morse--Bott at the point $p$ and this proves Item \eqref{item:Frankel_almost_Hermitian_f_is_MB}.

We now prove the remaining assertions in Item \eqref{item:Frankel_almost_Hermitian_ComponentsOfFixedPointsAreSmoothSubmanifolds}.
The existence of a non-degenerate two-form $\om$ on $M$ implies that $M$ admits an almost complex structure by Lemma \ref{lem:ExistenceOfS1InvarACStructure} and so $M$ is even-dimensional. The subspace $L_0\subset T_pM$ has even codimension by the direct-sum decomposition \eqref{eq:FixedPointTangentSpaceDecomp} of $T_pM$. Hence, each component of the fixed-point set of the circle action is an even-codimension submanifold of an even-dimensional manifold and is thus even-dimensional itself.  This completes the proof of Item \eqref{item:Frankel_almost_Hermitian_ComponentsOfFixedPointsAreSmoothSubmanifolds}.

Finally, Item \eqref{item:Frankel_almost_Hermitian_Sylvester} of Theorem \ref{mainthm:Frankel_almost_Hermitian}  follows from the discussion on the Hessian following Corollary \ref{cor:Sylvester_law_inertia}.
\end{proof}

\chapter{Local indices for vector fields on manifolds and analytic spaces}
\label{chap:Local_indices_vector_fields_manifolds_analytic_spaces}
In this chapter, partly for the sake of disambiguation, we discuss some previous approaches to the question of how to define a local `index' for a vector field when it has a zero at a singular point of an analytic space. Section \ref{sec:Local_indices_vector_fields_manifolds_analytic_spaces_introduction} provides an introduction to the chapter. In Section \ref{sec:Local_index_vector_field_manifold}, we survey definitions of local indices for vector fields on smooth manifolds and in Section \ref{sec:Local_indices_vector_fields_on_singular_varieties}, we outline previous generalizations of those concepts to the case of singular analytic spaces.

\section{Introduction}
\label{sec:Local_indices_vector_fields_manifolds_analytic_spaces_introduction}
In Chapter \ref{chap:Circle_actions_almost_Hermitian_manifolds}, we considered invariants associated with a zero $p$ (not necessarily isolated) of a smooth vector field $X$ on a smooth manifold $M$, namely the eigenvalues of the Lie derivative endomorphism $(\sL_X)_p \in \End(T_pM)$. We observed that these eigenvalues correspond to (integer) weights of a smooth circle action $\rho:S^1\times M \to M$ on an almost Hermitian manifold $(M,J)$ equipped with a circle-invariant almost complex structure $J$, when $X$ is the vector field generated by the circle action (and thus $p$ is a fixed point of that action). Motivated by Morse and Morse--Bott theory, we thus define the \emph{Morse--Bott index of a circle action at a fixed point} by
\begin{equation}
  \label{eq:Morse-Bott_index_circle_action_fixed_point}
  \lambda_p^-(\rho) := \dim_\RR T_p^-M,
\end{equation}
where $T_p^-M$ is the negative-weight subspace of $T_pM$ defined by the circle action. This definition generalizes easily to the case of a real or complex analytic variety $M = F^{-1}(0)$ defined by a circle-equivariant, analytic map $F:\KK^{n+d}\to\KK$, for $\KK=\RR$ or $\CC$. We define the \emph{virtual Morse--Bott index of a circle action at a fixed point} by
\begin{equation}
  \label{eq:Virtual_Morse-Bott_index_circle_action_fixed_point}
  \lambda_p^-(\rho) := \dim_\RR \Ker^- dF(p) - \dim_\RR \Coker^- dF(p),
\end{equation}
where $T_pM = \Ker dF(p)$ is now the Zariski tangent space to $M$ at $p$ and $\Ker^- dF(p) \subset \Ker dF(p)$ and $\Coker^- dF(p) \subset \Coker^- dF(p)$ are the negative-weight subspaces defined by the circle action. Since every analytic space is locally isomorphic to an analytic variety and the preceding definition is independent of the choice of local analytic model space for an open neighborhood of $p$ in $M$, the concept of virtual Morse--Bott index is well-defined on analytic spaces. (See Akhiezer \cite[Sections 1.1 and 1.2]{Akhiezer_lie_group_actions_complex_analysis} for a discussion of group actions on complex analytic spaces and Feehan \cite{Feehan_analytic_spaces} for further development.)

The \emph{Poincar\'e--Hopf index} $\ind_p(X)$ of an isolated zero $p$ of a smooth vector field $X$ on an oriented, smooth manifold $M$ is defined in terms of the degree of an induced map on a sphere of dimension $\dim M-1$ and the \emph{Poincar\'e--Hopf Index Theorem} asserts that if $X$ has finitely many zeros and $M$ is closed, then
\begin{equation}
  \label{eq:Poincare_Hopf_index_theorem}
  \sum_{p\in\Zero(X)}\ind_p(X) = e(M),
\end{equation}
where $e(M)$ is the Euler characteristic of $M$; see Guillemin and Pollack \cite[p. 134]{Guillemin_Pollack}, Milnor \cite[p. 35]{Milnor_topology_from_differentiable_viewpoint}, or Brasselet, Seade, and Suwa \cite[Theorem 1.1.1]{Brasselet_Seade_Suwa_vector_fields_singular_varieties}. The formula \eqref{eq:Poincare_Hopf_index_theorem} is valid for continuous vector fields. Many researchers have explored the question of what is a good notion of index of a vector field if the smooth manifold is replaced by a singular space and that generalize the definition of Poincar\'e--Hopf index, for isolated and non-isolated zeros. These definitions include those of Schwartz (see Section \ref{subsec:Schwartz_index}), G\'omez--Mont, Seade and Verjovsky (see Section \ref{subsec:GSV_index}), and of the \emph{virtual index} introduced Lehmann, Soares and Suwa \cite{Lehmann_Soares_Suwa_1995} for holomorphic vector fields (see Section \ref{subsec:Virtual_index}). Their definitions of ``index'' are very different from our definition of virtual Morse--Bott index, which is only defined when the vector field $X$ is generated by an $S^1$ or $\CC^*$ action, but because there is some overlap of methods and terminology (and thus a potential for confusion), we briefly review these alternative concepts of index in this chapter.

\section{Local indices for vector fields on manifolds}
\label{sec:Local_index_vector_field_manifold}
In this section, we summarize concepts related to the local index of a vector field on a smooth manifold including the Poincar\'e--Hopf index in Section \ref{subsec:Poincare-Hopf_index_Euler_characteristic}, the localization and residue of Chern classes of a vector bundle $E$ given an $\ell$-frame for a $E$ over an open subset in Section \ref{subsec:Localization_residues_Chern_classes}, and Grothendieck residues of holomorphic functions in Section \ref{subsec:Grothendieck_residues}. 

\subsection{Poincar\'e--Hopf index for non-isolated zeros and the Euler characteristic}
\label{subsec:Poincare-Hopf_index_Euler_characteristic}
Following Brasselet, Seade, and Suwa \cite[Section 1.1.2]{Brasselet_Seade_Suwa_vector_fields_singular_varieties}, suppose that $M$ is a closed, topological manifold with a triangulation $K$ and that $S\subset M$ is a compact, connected $K$-subcomplex of $M$. A \emph{cellular tube $T$ around $S$ in $M$} is the union of cells that are duals of simplices in $S$; this notion generalizes the concept of tubular neighborhood of a submanifold $S$. If $S$ is a submanifold without boundary, then $T$ is a bundle on $S$ whose fibers are disks. According to \cite[Theorem 1.1.2]{Brasselet_Seade_Suwa_vector_fields_singular_varieties}, a continuous vector field $X$ without zeros on an open neighborhood of $\partial T$ can be extended to the interior of $T$ with finitely many isolated zeros. The Poincar\'e--Hopf index of $X$ at $S$, denoted $\ind_S(X)$, is defined as the sum of the indices of the extension of $X$ at these points \cite[Equation (1.1.3)]{Brasselet_Seade_Suwa_vector_fields_singular_varieties}. The Poincar\'e--Hopf Index Theorem for a continuous vector field $X$ over a closed, oriented, smooth manifold $M$ now asserts that \cite[Definition 1.1.3]{Brasselet_Seade_Suwa_vector_fields_singular_varieties}
\begin{equation}
  \label{eq:Poincare_Hopf_index_theorem_non-isolated_zeros}
  \sum_{\lambda}\ind_{S_\lambda}(X) = e(M),
\end{equation}
where the sum is over the connected components $\{S_\lambda\}$ of $S = \Zero(X)$.

\subsection{Localization and residues of Chern classes}
\label{subsec:Localization_residues_Chern_classes}
We follow Brasselet, Seade, and Suwa \cite[Sections 1.6.2 and 1.6.3]{Brasselet_Seade_Suwa_vector_fields_singular_varieties} and especially Suwa \cite[Section 2]{Suwa_2003} for the definition of residues of Chern classes of complex vector bundles over smooth manifolds and to which we refer for further details. Let $E$ be a smooth complex vector bundle of rank $r$ over an oriented smooth manifold of dimension $n$. For a connection $\nabla$ on $E$ and $j = 1, \ldots, r$, one denotes by $c^j(\nabla)$ the $j$-th Chern form defined by $\nabla$,
\[
  c^j(\nabla) = \left(\frac{i}{2\pi}\right)^j \sigma_j(F_\nabla),
\]
where $\sigma_j(F_\nabla)$ denotes the $j$-th symmetric form of the curvature matrix $F_\nabla$ of $\nabla$ and is a closed $2j$-form on $M$. Its class
\[
  c^j(\nabla) \in H^{2j}(M,\CC)
\]
is the $j$-th Chern class $c^j(E)$ of $E$. Given two connections $\nabla_0,\nabla_1$ there is a $(2j-1)$-form $c^j(\nabla_0,\nabla_1)$ that is alternating in the two entries and obeys (see Bott \cite{Bott_1972} or Kobayashi and Nomizu \cite[Chapter XII, Section 1, Lemma 5, p. 297]{Kobayashi_Nomizu_v2})
\begin{equation}
  \label{eq:Suwa_2003_1-4}
  c^j(\nabla_0) - c^j(\nabla_1) - dc^j(\nabla_0,\nabla_1) = 0.
\end{equation}
For an open subset $U \subset M$, we denote by $A^q(U)$ the space of complex-valued smooth $q$-forms on $U$. For an open covering $\sU = \{U_\alpha\}$ of $M$, we denote by $A^\bullet(\sU)$ the \v{C}ech--de Rham complex associated to the covering $\sU$ with differential $d$ and by $H^q(A^\bullet(\sU))$ its cohomology (see Bott and Tu \cite[Section 8 and Theorem 8.1]{BT}).

Let $\bs = \{s_1,\ldots,s_\ell\}$ be a smooth $\ell$-frame for $E$ over an open set $U$, thus a set of $\ell$ smooth sections of $E$ that are linearly independent at each point of $U$. A connection $\nabla$ is called $\bs$-\emph{trivial} or \emph{trivial with respect to $\bs$} if $\nabla s_j=0$ for $j=1,\ldots,\ell$. Let $S\subset M$ be a closed subset and $\bs$ be a smooth $\ell$-frame for $E$ over $M\less S$. Let $U_1$ be an open neighborhood of $S$ in $M$ and $U_0 := M \less S$ and write $\sU = \{U_0,U_1\}$. The Chern class $c^j(E) \in H^{2j}(M,\CC)$ is represented by the cocycle \cite[Equation (2.2)]{Suwa_2003}
\begin{equation}
  \label{eq:Suwa_2003_2-2}
  c^j(\nabla_*) := (c^j(\nabla_0),c^j(\nabla_1),c^j(\nabla_0,\nabla_1)) \in A^{2j}(\sU),
\end{equation}
where $\nabla_0$ and $\nabla_1$ are connections for $E$ on $U_0$ and $U_1$, respectively. If one chooses $\nabla_0$ to be $\bs$-trivial, so $c^j(\nabla_0) = 0$ for $j\geq r-\ell+1$ by Suwa \cite[Equation (2.1)]{Suwa_2003}, then $c^j(\nabla_*)$ defines a relative cohomology class,
\[
  c^j(E,\bs) \in H^{2j}(M,M\less S;\CC),
\]
that is independent of the choice of connection $\nabla_1$ or the $\bs$-trivial connection $\nabla_0$. One calls $c^j(E,\bs)$ the \emph{localization of $c^j(E)$ at $S$ with respect to $\bs$}.

Suppose now that $S$ is a compact set admitting a regular neighborhood and let $\{S_\lambda\}$ be the connected components of $S$. Thus, for each $\lambda$, the Alexander dual of the element $c^j(E,\bs)$ (see Hatcher \cite[Proposition 3.46]{Hatcher}) defines a class in $H_{n-2j}(S_\lambda,\CC)$ which one calls the \emph{residue of $\bs$ at $S_\lambda$ with respect to $c^j$}, or simply the \emph{residue of $c^j(E)$}, and denoted by
\begin{equation}
  \label{eq:Residue_section_at_closed_wrt_ci(E)}
  \Res_{c^j}(\bs,E;S_\lambda).
\end{equation} 
For each $\lambda$, one chooses an open neighborhood $U_\lambda$ of $S_\lambda$ in $U_1$, so that the neighborhoods $U_{\lambda}$ are mutually disjoint, and let $R_\lambda$ be a smooth $n$-dimensional manifold with smooth boundary in $U_\lambda$ containing $S_\lambda$ in its interior. We set $R_{0\lambda} := -\partial R_\lambda$. Then the residue \eqref{eq:Residue_section_at_closed_wrt_ci(E)} is represented by an $(n-2i)$-cycle $C$ in $S_\lambda$ such that \cite[Equation (2.3)]{Suwa_2003}
\begin{equation}
  \label{eq:Suwa_2003_2-3}
  \int_C \tau = \int_{R_\lambda} c^j(\nabla_1)\wedge \tau + \int_{R_{0\lambda}} c^j(\nabla_0,\nabla_1)\wedge\tau
\end{equation}  
for every closed $(n-2j)$-form $\tau$ on $U_\lambda$. If $2j=n$, the residue is a number given
by \cite[Equation (2.4)]{Suwa_2003}
\begin{equation}
  \label{eq:Suwa_2003_2-4}
  \Res_{c^j}(\bs,E;S_\lambda)
  = \int_{R_\lambda} c^j(\nabla_1) + \int_{R_{0\lambda}} c^j(\nabla_0,\nabla_1).
\end{equation}
When $\ell=1$, one obtains

\begin{thm}[Localization of the top Chern class]
\label{thm:Brasselet_Seade_Suwa_vector_fields_singular_varieties_1-6-7}  
(See Brasselet, Seade, and Suwa \cite[Theorem 1.6.7]{Brasselet_Seade_Suwa_vector_fields_singular_varieties}.)  
Let $E$ be a smooth complex vector bundle of rank $r$ over a closed, oriented, smooth manifold $M$ of dimension $m$. Let $S\subset M$ be a closed set admitting a regular open neighborhood and $s$ be a smooth section of $E$ such that $s$ has no zeros in $M\less S$. Then for each connected component $S_\lambda$ of $S$, the residue $\Res_{c^r}(s,E;S_\lambda)$ belongs to $H_{m-2r}(S_\lambda,\CC)$ and
\[
  \sum_\lambda (\iota_\lambda)_*\Res_{c^r}(s,E;S_\lambda) = c^r(E)\frown [M],
\]
where $\iota_\lambda:S_\lambda \hookrightarrow M$ denotes the inclusion.
\end{thm}  

Following to Brasselet, Seade, and Suwa \cite[Section 1.6.7(a)]{Brasselet_Seade_Suwa_vector_fields_singular_varieties}), one can define (see Suwa \cite{Suwa_2008})
\[
  \ind_{S_\lambda}(X) := \Res_{c^m}(s,TM;S_\lambda),
\]
if $M$ is a closed, complex manifold of dimension $m$ and $TM$ is the holomorphic tangent bundle and $X$ is a holomorphic vector field, where $\ind_{S_\lambda}(X)$ is the Poincar\'e--Hopf index. If $S_\lambda$ is a single point, this definition coincides with the usual definition of Poincar\'e--Hopf index by Theorem \ref{thm:Brasselet_Seade_Suwa_1-6-14}. By combining Theorems \ref{thm:Brasselet_Seade_Suwa_vector_fields_singular_varieties_1-6-7} and  \ref{thm:Brasselet_Seade_Suwa_1-6-14} with the well-known fact that
\[
  \int_M c^m(TM) = e(M),
\]
one obtains the Poincar\'e--Hopf index formula \eqref{eq:Poincare_Hopf_index_theorem} for holomorphic vector fields with isolated zeros. See Bott \cite[Theorem 1]{Bott_1967mmj}, \cite[Theorem 1]{Bott_1967jdg} and \cite[Theorem 6.1.2]{Brasselet_Seade_Suwa_vector_fields_singular_varieties} (due to Baum and Bott \cite{Baum_Bott_1972} and Bott \cite{Bott_1972}) for generalizations of Theorem \ref{thm:Brasselet_Seade_Suwa_vector_fields_singular_varieties_1-6-7}.

\subsection{Grothendieck residues}
\label{subsec:Grothendieck_residues}
Following Brasselet, Seade, and Suwa \cite[Section 1.6.5]{Brasselet_Seade_Suwa_vector_fields_singular_varieties} or Griffiths and Harris \cite[p. 649]{GriffithsHarris}, let $U\subset\CC^m$ be an open neighborhood of the origin and $f_1,\ldots,f_m \in \sO(U)$ be holomorphic functions such that $f_1^{-1}(0)\cap\cdots\cap f_m^{-1}(0) = \{0\}$. If $\omega$ is a holomorphic $m$-form on $U$, for example
\[
  \omega = g(z)\,dz_1\wedge\cdots\wedge dz_m 
\]  
for $g \in \sO(U)$, one sets
\begin{equation}
  \label{eq:Brasselet_Seade_Suwa_1-6-10_example_1-6-2}
  \Res_0\left[\begin{matrix} \omega \\ f_1,\ldots, f_m\end{matrix}\right]
  :=
  \frac{1}{(2\pi i)^m} \int_\Gamma \frac{\omega}{f_1\cdots f_m},
\end{equation}
where $\Gamma\subset U$ is an $m$-cycle defined by
\[
  \Gamma := \left\{z\in U: |f_1(z)| = \cdots = |f_m(z)| = \eps \right\}
\]
for a constant $\eps \in (0,1]$ and oriented so that $d\theta_1\wedge\cdots\wedge d\theta_m$ is positive, where $\theta_j = \arg f_j$ for $j=1,\ldots,m$.

\begin{thm}[Equivalence of Grothendieck residue and residue of the top Chern class at an isolated zero of a holomorphic section] 
\label{thm:Brasselet_Seade_Suwa_1-6-12}  
(See Brasselet, Seade, and Suwa \cite[Theorem 1.6.12]{Brasselet_Seade_Suwa_vector_fields_singular_varieties} for a statement and Suwa \cite[Theorem 3.1]{Suwa_2000} for a statement and proof.)
Let $E$ be a holomorphic vector bundle of rank $m$ over a complex manifold $M$ of dimension $m$ and $s$ be a holomorphic section of $E$ with an isolated zero at a point $p\in M$. If $\{e_1,\ldots,e_m\}$ is a holomorphic frame for $E$ over an open neighborhood $U\subset M$ of $p$, so that $s = \sum_{j=1}^m f_je_j$, then
\[
  \Res_{c^m}(s,E;p)
  =
  \Res_p\left[\begin{matrix} df_1\wedge\cdots\wedge df_m \\ f_1,\ldots, f_m\end{matrix}\right].
\]
\end{thm}

See also Griffiths and Harris \cite[Section 5.1]{GriffithsHarris} for an introduction to Grothendieck residues.

\begin{thm}[Algebraic expression for the residue of the top Chern class at an isolated zero of a holomorphic section] 
\label{thm:Brasselet_Seade_Suwa_1-6-13}  
(See Brasselet, Seade, and Suwa \cite[Theorem 1.6.13]{Brasselet_Seade_Suwa_vector_fields_singular_varieties} for a statement and Suwa \cite[Corollary 5.3]{Suwa_2005} for a statement and proof.) Continue the hypotheses of Theorem \ref{thm:Brasselet_Seade_Suwa_1-6-12}. Then
\[
  \Res_{c^m}(s,E;p)
  =
  \dim\sO_m/(f_1,\ldots,f_m),
\]
where $\sO_m=\CC\{z_1,\ldots,z_m\}$ is the ring of convergent power series in $m$ complex variables $z_1,\ldots,z_m$ and $(f_1,\ldots,f_m)$ is the ideal of $\sO_m$ generated by $f_1,\ldots,f_m$.
\end{thm}

\begin{thm}[Topological expression for the residue of the top Chern class at an isolated zero of a holomorphic section] 
\label{thm:Brasselet_Seade_Suwa_1-6-14}  
(See Brasselet, Seade, and Suwa \cite[Theorem 1.6.14]{Brasselet_Seade_Suwa_vector_fields_singular_varieties} for a statement and Suwa \cite[Corollary 2.2]{Suwa_2005}, \cite[Corollary 3.4.2]{Suwa_2008} for statements and proofs.)
Continue the hypotheses of Theorem \ref{thm:Brasselet_Seade_Suwa_1-6-12}. Let $\SSS_\eps^{2m-1}(p) \subset U$ denote the sphere with center $p$ and radius $\eps\in(0,1]$ with respect to a choice of Riemannian metric on $U$ and define
\[
  \varphi: \SSS_\eps^{2m-1}(p) \ni z \mapsto \frac{f(z)}{\|f(z)\|} \in \SSS^{2m-1},
\]
where $\SSS^{2m-1} \subset \CC^m$ denotes the unit sphere with center at the origin. Then
\[
  \Res_{c^m}(s,E;p)
  =
  \deg\varphi.
\]
\end{thm}

\section{Local indices for vector fields on singular analytic spaces}
\label{sec:Local_indices_vector_fields_on_singular_varieties}
Brasselet, Seade, and Suwa \cite{Brasselet_Seade_Suwa_vector_fields_singular_varieties} provide a survey of many previous approaches to the problem of how to define the index of a vector field on a singular analytic spaces at one of its zeros and a comprehensive bibliography, with over 170 entries. All of these approaches are quite different to ours and they typically aim to generalize the definition of the Poincar\'e--Hopf index. However, because there is occasional overlap in terminology, if not concepts, we briefly list the approaches here. We omit details and refer instead to \cite{Brasselet_Seade_Suwa_vector_fields_singular_varieties} for definitions, results, and references to the literature. See also Ebeling and Gusein--Zade \cite{Ebeling_Gusein-Zade_2021arxiv} for a more concise introduction.

\subsection{Schwartz or radial index}
\label{subsec:Schwartz_index}
Schwartz considers a special class of vector fields that she called ``radial'', which are obtained
by the process of radial extension. The Schwartz index can be defined for isolated and non-isolated zeros of vector fields on real or complex analytic varieties and versions of the Poincar\'e--Hopf index theorem continue to hold \cite[Chapters 2 and 4]{Brasselet_Seade_Suwa_vector_fields_singular_varieties}.

\subsection{GSV index}
\label{subsec:GSV_index}
The GSV index was introduced by G\'ome--Mont, Seade and Verjovsky, for vector fields on complete intersections. The index has been defined for isolated and non-isolated zeros of vector fields on real or complex analytic varieties and versions of the Poincar\'e--Hopf index theorem continue to hold \cite[Chapters 3 and 4]{Brasselet_Seade_Suwa_vector_fields_singular_varieties}.

\subsection{Virtual index}
\label{subsec:Virtual_index}
The virtual index was introduced Lehmann, Soares and Suwa \cite{Lehmann_Soares_Suwa_1995} for holomorphic vector fields and, despite the name, it is unrelated to our concept of virtual Morse--Bott index. If the variety has only isolated singularities, the virtual index and the GSV index coincide. The virtual index is relatively easy to compute and is defined for vector fields with non-isolated zero set. One can regard the virtual index as being a localization of the top-dimensional Chern class of a \emph{virtual tangent bundle} (the holomorphic tangent bundle minus a holomorphic vector bundle), called the \emph{virtual class}, just as the local Poincar\'e--Hopf index is a localization of the top Chern class of a manifold. Versions of the Poincar\'e--Hopf index theorem continue to hold for this notion of index \cite[Chapter 5]{Brasselet_Seade_Suwa_vector_fields_singular_varieties}. 

\subsection{Homological index and algebraic formulae}
\label{subsec:Homological_index}
For a holomorphic vector field $X$ in $\CC^n$ with an isolated zero at the origin $p$, the local Poincar\'e--Hopf index satisfies
\[
\ind_0(X) = \dim \sO_{\CC^n,0}/(a_1, \ldots , a_n),
\]
where $(a_1, \ldots, a_n)$ is the ideal generated by the components of $X$.

In the real analytic setting, the equivalent statement is given by the formula of Eisenbud, Levine, and Khimshiashvili. According to Eisenbud and Levine \cite[Main Theorem, p. 178]{Eisenbud_1978}, \cite[Theorem 1.2]{Eisenbud_Levine_1977}, the degree of a smooth map $f$ from an open neighborhood in $\RR^n$ to another such neighborhood with $f(0)=0$ is equal to the signature of a certain symmetric bilinear form; see Arnold, Gusein--Zade, and Varchenko \cite[Section 5.12, Theorem 5.12]{Arnold_Gusein-Zade_Varchenko_singularities_differentiable_maps_v1} for another exposition of the proof. Closely related formulae were derived independently by Eisenbud and Levine \cite[Theorem 1.1]{Eisenbud_Levine_1977} and Khimshiashvili \cite[Theorem 5.1]{Khimshiashvili_signature_formulae_topological_invariants}, with certain generalizations due to  Grigor$'$ev and Ivanov \cite{Grigoriev_Ivanov_1980}.

These facts motivated the search for algebraic formulas for indices of vector fields on singular varieties. A major contribution in this direction was given by Arnold for gradient vector fields and also significant contributions by various authors, such as G\'omez--Mont, Gusein--Zade, and Ebeling  \cite[Chapter 7]{Brasselet_Seade_Suwa_vector_fields_singular_varieties}.

\chapter[ASD connections, Seiberg--Witten monopoles, and non-Abelian monopoles]{Moduli spaces of anti-self-dual connections, Seiberg--Witten monopoles, and non-Abelian monopoles}
\label{chap:Preliminaries}
In this chapter, we review some of the notation and definitions from our preceding articles on non-Abelian monopoles \cite{FL2a,FL2b} that we require in this monograph. We begin in Section \ref{sec:Gauge_transformations_stabilizers_unitary_connections} by discussing unitary connections on a complex Hermitian vector bundle, the action of the relevant groups of gauge transformations on these connections, the stabilizers of these connections under this action, and the Lie algebras of these stabilizers. In Section
\ref{sec:SpincuStr}, we recall the definition of a spin${}^u$ structure, its characteristic classes, and the definition of a spin connection.  In Section \ref{sec:SpinuPairsQuotientSpace}, we introduce the quotient space of pairs $\sC_\ft$ for a spin${}^u$ structure on which the non-Abelian monopole equations will be defined, describe the tangent bundle of this quotient space, and construct an embedding of a quotient space of spin${}^c$ pairs into $\sC_\ft$. The definition and basic properties of the moduli space of non-Abelian monopoles from \cite{FL1, FeehanGenericMetric} appear in Section \ref{sec:PU2Monopoles}.  In Section \ref{sec:ASDsingularities}, we describe the stratum of zero-section monopoles, that is, anti-self-dual connections. In Section \ref{sec:Reducibles}, we describe the strata of split, that is, Seiberg--Witten monopoles.  We discuss a circle action on the quotient space of pairs in Section \ref{sec:S1Actions} and the fixed point set of this action in the moduli space of non-Abelian monopoles.  Finally, in Section \ref{sec:IrredMonopolesForNonGenericPerturbations} we prove a criterion for the existence of a non-Abelian monopole which is neither zero-section nor split for a possibly non-generic parameter.

\section{Gauge transformations and stabilizers of unitary connections}
\label{sec:Gauge_transformations_stabilizers_unitary_connections}
Let $(E,h)$ be a smooth Hermitian vector bundle over a smooth manifold $X$ of dimension $d$. Although it plays no role in this section, for consistency with the remainder of this monograph, we shall fix a smooth unitary connection $A_d$ on the smooth Hermitian line bundle $(\det E,h_{\det E})$, where $\det E = \wedge^r E$ with $r$ denoting the complex rank of $E$ and $h_{\det E}$ is the Hermitian metric on the complex line bundle $\det E$ by the Hermitian metric $h$ on $E$. Let $p \in (d/2,\infty)$ be a constant and let $\sA(E,h)$ denote the affine space of $W^{1,p}$ unitary connections $A$ on $E$ that induce $A_d$ on $\det E$ when $E$ has complex rank two or higher (the condition is omitted when $E$ is a Hermitian line bundle).
\label{page:UnitaryConnectionsOverManifoldOfDim_d}
We let
\begin{subequations}
\label{eq:UnitaryAutomorphismBundles}
\begin{align}
  \label{eq:AutomorphismBundles_U(E)}
  \U(E,h)& := \{u\in \End(E): u(x)\in\U(E_x), \text{ for all } x\in X\},
  \\
  \label{eq:AutomorphismBundles_SU(E)}
  \SU(E,h)& := \{u\in \End(E): u(x)\in\SU(E_x), \text{ for all } x\in X\},
\end{align}
\end{subequations}
denote the smooth principal fiber bundles over $X$ with structure groups $\U(r)$ and $\SU(r)$, respectively, when $E$ has complex rank $r$. Gauge transformations are sections of these bundles. When the Hermitian metric $h$ is understood, we shall abbreviate $\U(E):=\U(E,h)$ and $\SU(E):=\SU(E,h)$. 

The Banach Lie group $W^{2,p}(\SU(E))$ of determinant-one, unitary automorphisms of $(E,h)$ acts on (the right of) $\sA(E,h)$
 \begin{equation}
  \label{eq:W2pSUE_action_on_AEh}
  W^{2,p}(\SU(E)) \times \sA(E,h) \ni (u,A) \mapsto u^*A \in \sA(E,h)
\end{equation}
by pullback. By adapting the proofs \mutatis due to Freed and Uhlenbeck of \cite[Appendix A, Proposition A.2, p. 160 and Proposition A.3, p. 161]{FU}, one can show that $W^{2,p}(\SU(E))$ is a Banach Lie group and that its action \eqref{eq:SL(E)_Action_on_(0,1)conn} on $\sA(E,h)$ is smooth.

\begin{rmk}[Pushforward and pullback notation for gauge transformations on connections]
\label{rmk:PushforwardPullbackNotation}
The notation $u^*A$ for the action of a gauge transformation on a connection  in \eqref{eq:GaugeActionOnSpinuPairs} means the \emph{pullback} action which changes the covariant derivative $\nabla_A$ by
\begin{equation}
  \label{eq:Donaldson-Kronheimer_2-1-7_pullback}
  \nabla_{u^*A}=u^{-1}\circ\nabla_A\circ u.
\end{equation}
Observe that this definition implies that $(uv)^*A=v^*u^*A$, justifying the use of the term pullback. The \emph{pushforward} action of a gauge transformation $u$ on a connection $A$ is $u_*A:=(u^{-1})^*A$, so $\nabla_{u_*A}=u\circ\nabla_A\circ u^{-1}$ as in Donaldson and Kronheimer \cite[Equation (2.1.7), p. 34]{DK}.
\end{rmk}

The Sobolev norms implicit in our definition of $\sA(E,h)$ and $W^{2,p}(\SU(E))$ are defined by analogy with those of Aubin \cite[Definition 2.3]{Aubin_1998}, based in turn on those of Adams and Fournier \cite[Section 3.1]{AdamsFournier}. If $(V,h)$ is a smooth Riemannian vector bundle equipped with a compatible smooth covariant derivative $\nabla$ over a smooth Riemannian manifold $(X,g)$ and $s$ is a smooth section of $V$ and $k\geq 0$ is an integer and $p\in [1,\infty)$, then
\begin{equation}
  \label{eq:Aubin_definition_2-3}
  \|s\|_{W_\nabla^{k,p}(X)} := \left(\sum_{l=0}^k \|\nabla^ls\|_{L^p(X)}^p\right)^{1/p},
\end{equation}
where we continue to write $\nabla$ for the covariant derivative on the Riemannian vector bundle $\otimes^lT^*X\otimes V$, when $\l\geq 1$, induced by the Levi-Civita connection for $g$ on $TX$, and thus $T^*X$, and the given covariant derivative on $V$. For $u \in W^{2,p}(\SU(E))$, we consider $\SU(E) \subset \End(E) \cong E\otimes E^*$
and continue to write $\nabla_A$ for the covariant derivative on the Riemannian vector bundle $\otimes^lT^*X\otimes E\otimes E^*$ when $\l\geq 1$ induced by the Levi-Civita connection for $g$ on $TX$, and thus $T^*X$, and the covariant derivative on $E\otimes E^*$ induced by smooth covariant derivative $\nabla_A$ on $E$, yielding
\begin{equation}
  \label{eq:W2p_norm_gauge_transformation}
  \|u\|_{W_A^{2,p}(X)}
  := \left(\|u\|_{L^p(X)}^p + \|\nabla_Au\|_{L^p(X)}^p + \|\nabla_A^2u\|_{L^p(X)}^p\right)^{1/p}.
\end{equation}
Similarly, for $a \in W^{1,p}(T^*X\otimes\su(E))$ and $\nabla_A$ on $E$, we write
\begin{equation}
  \label{eq:W1p_norm_connection}
  \|a\|_{W_A^{1,p}(X)}
  := \left(\|a\|_{L^p(X)}^p + \|\nabla_Aa\|_{L^p(X)}^p\right)^{1/p}.
\end{equation}
We shall use \eqref{eq:Aubin_definition_2-3} as the basis of all our definitions of all Sobolev spaces in this monograph.

\begin{rmk}[On the induced connection on $\End(E)$]
\label{rmk:InducedConn_on_EndE}
Let $A$ be a connection on a vector bundle $E$.  We shall give an explicit expression for the connection induced by $A$ on $\End(E)$ because it often appears in computations in this monograph. Although we shall usually write $A$ for the connection and $d_A$ for the covariant derivative induced by $A$ on an associated bundle, we will occasionally write $A^\ad$ for the connection and $d_{A^\ad}$ for the covariant derivative induced by $A$ on $\End(E)$ when additional clarity is required. By Lawson \cite[Equation (2.4), p. 19]{Lawson}), the covariant derivative $d_{A^{\ad}}$ is defined by
\begin{equation}
  \label{eq:Connection_Induced_On_Endomorphism_Bundle}
  (d_{A^{\ad}} w)s := d_A(w(s)) - w(d_As),
  \quad\text{for all } w\in\Omega^0(\End(E)) \text{ and } s \in \Omega^0(E).
\end{equation}
The expression \eqref{eq:Connection_Induced_On_Endomorphism_Bundle} follows from the identity
\[
d_A(ws) = (d_{A^\ad}w)s + w(d_As),
\]
which can be derived the Leibnitz formula for covariant derivatives of tensor products and the requirement that covariant derivatives commute with contractions.
\end{rmk}

For $A\in\sA(E,h)$, we define the \emph{stabilizer} of $A$ as
\begin{equation}
\label{eq:DefineStabilizerOfConnection_in_SU(E)}
\Stab(A):=\{u\in W^{2,p}(\SU(E)): u^*A=A\}.
\end{equation}

We note the following well-known Lie group properties of stabilizers of unitary connections.

\begin{lem}[Lie group structure of $\Stab(A)$]
\label{lem:LieGroup_Structure_of_Stab(A)}
Let $(E,h)$ be a smooth Hermitian vector bundle over a connected, smooth manifold $X$. If $A\in\sA(E,h)$, then $\Stab(A)$ is a Lie subgroup of $W^{2,p}(\SU(E))$ and (compare the forthcoming \eqref{eq:HE_equation_bHA0})
\[
  \bH_A^0
  :=
  \Ker\left(d_{A^{\ad}}:W^{2,p}(\su(E))\to W^{1,p}(T^*X\otimes \su(E))\right)
\]
is the Lie algebra of $\Stab(A)$.
\end{lem}

\begin{proof}
By Hilgert and Neeb \cite[Lemma 10.1.5 (i), p. 361]{Hilgert_Neeb_structure_geometry_lie_groups}, the stabilizer $\Stab(A)$ is a closed subgroup of $W^{2,p}(\SU(E))$ and thus a Lie group.
The identification of the kernel of $d_A$, and thus of $\bH_A^0$, with the Lie algebra of $\Stab(A)$ appears in the discussion following Donaldson and Kronheimer \cite[Lemma 4.2.8, p. 132]{DK}, the discussion by Freed and Uhlenbeck in \cite[Chapter 4, p. 64]{FU}, and Kronheimer's proof of \cite[Lemma 2.2, p. 64]{Kronheimer_2005}. Because this argument appears frequently in this monograph, we shall give the following additional details (compare the computation in Donaldson and Kronheimer \cite[Equation (2.1.20), p. 37]{DK}).  The Lie algebra of $\Stab(A)$ is generated by vectors tangent to smooth paths $u_t\in W^{2,p}(\SU(E))$ with $u_t^*A=A$.
If $u_t\in W^{2,p}(\SU(E))$ is a smooth path for $t\in(-\eps,\eps)$ with $u_0=\id_E$, and $du_t/dt|_{t=0} = \xi \in W^{2,p}(\su(E))$, and $u_t^*A=A$, then
\footnote{In \cite[Section 2.1.3, p. 280]{FL1}, we chose the action of $W^{2,p}(\SU(E))$ (the Lie group of determinant one, unitary automorphisms of $E$) on $\sA(E,h)$ (the affine space of unitary connections with fixed connection on $\det E$) to be $\nabla_{u_*(A)}s := u\nabla_A(u^{-1}s)$ for all $s\in\Omega^0(E)$, as in the pushforward action in \eqref{eq:PushforwardActionOnUnitaryPairs}. This leads to the expression $d_{A,\Phi}^0\xi = (-\nabla_A\xi,\xi\Phi)$ in \cite[Proposition 2.1 (3), p. 280]{FL1} and in \cite[Equation (2.35), p. 71]{FL2a}. It agrees with the Donaldson--Kronheimer convention in \cite[Equation (2.1.7), p. 34]{DK}, leading to the expression $d_A^0\xi = -\nabla_A\xi$ in \cite[Equation (2.1.20), p. 37]{DK}. By contrast in \eqref{eq:d_APhi^0} we use the Freed--Uhlenbeck convention, $\nabla_{u^*(A)}s := u^{-1}\nabla_A(us)$ for all $s\in\Omega^0(E)$ from \cite[Proof of Proposition A.3]{FU}, leading to the expression $d_A^0\xi = \nabla_A\xi$ in \cite[p. 49]{FU}, and $u(\Phi) := u^{-1}\Phi$ leading to the expression $-\xi\Phi$.}
\begin{equation}
\label{eq:Derivative_Of_GaugeGroupAction}
  0= \left.\frac{d}{dt}u_t^*A\right|_{t=0} = \left.\frac{d}{dt}\left(u_t^{-1}\circ d_A\circ u_t\right)\right|_{t=0}
  =
  [d_A,\xi]
  =
  d_{A^{\ad}}\xi,
\end{equation}
where the final equality follows from \eqref{eq:Connection_Induced_On_Endomorphism_Bundle}. Hence, the Lie algebra of $\Stab(A)$ is identified with the kernel of $d_{A^{\ad}}:W^{2,p}(\su(E))\to W^{1,p}(T^*X\otimes \su(E))$, as asserted.
\end{proof}

Recall that the Lie groups $\SU(r)$ and $\U(r)$ have centers $Z(\SU(r)) = C_r$ and $Z(\U(r)) = S^1$, respectively, where we recall that $C_r = \{\varrho \in \CC: \varrho^r = 1\}$ is the group of $r$-th roots of unity (see Hilgert and Neeb \cite[Example 9.3.13 (b), p. 323]{Hilgert_Neeb_structure_geometry_lie_groups}).
With that in mind, we make the

\begin{defn}[Reducible, split, and central-stabilizer unitary connections]
\label{defn:Reducible_split_trivial-stabilizer_unitary_connection}
Let $A$ be a $W^{1,p}$ unitary connection with $p\in (d/2,\infty)$ on a $W^{2,p}$ Hermitian vector bundle $(E,h)$ of complex rank $r$ over a smooth, connected manifold of dimension $d\geq 2$.
\begin{enumerate}
\item\label{item:central-stabilizer_unitary_connection}
$A$ has \emph{central stabilizer} if the stabilizer $\Stab(A)$ of $A$ is minimal and thus isomorphic to the center $Z(\SU(r))=C_r$ and has \emph{non-central stabilizer} otherwise.

\item\label{item:Split_unitary_connection}  
$A$ is \emph{split} if $A = A_1\oplus A_2$ with respect to a decomposition
\begin{equation}
\label{eq:BasicSplitting}
E = E_1\oplus E_2
\end{equation}
as an orthogonal direct sum of proper, $W^{2,p}$ Hermitian subbundles, where $A_i$ is a $W^{1,p}$ unitary connection on $E_i$, for $i=1,2$, and is \emph{non-split} otherwise. If the pair $(E,A)$ is smooth (respectively, real analytic), then each pair $(E_i,A_i)$ is smooth (respectively, real analytic). 

\item\label{item:Reducible_unitary_connection}
$A$ is \emph{reducible} if the holonomy group of $A$ reduces (in the sense of Kobayashi and Nomizu \cite[Section II.6, p. 81]{Kobayashi_Nomizu_v1}) to a proper subgroup of $\U(r)$ and is \emph{irreducible} otherwise.
\end{enumerate}
\end{defn}

\begin{rmk}[Unitary stabilizers]
\label{rmk:UnitaryStabilizers}
If we define the unitary stabilizer group of a unitary connection $A$ as 
\[
\Stab_{\U(E)}(A):=\{u\in W^{2,p}(\U(E)): u^*A=A\}
\]
then this stabilizer group and that in \eqref{eq:DefineStabilizerOfConnection_in_SU(E)} are related
by
\begin{equation}
\label{eq:Isom_of_StabilizerGroups}
\Stab_{\U(E)}(A)\cong \Stab(A)\times_{C_r}S^1.
\end{equation}
One sees this by observing that if $u\in\Stab_{\U(E)}(A)$, then $\det u$ must be in the stabilizer of
$A^{\det}$, the connection induced by $A$ on the determinant line bundle $\det E$, and thus locally constant.
Because $X$ is connected and $\det u$ is locally constant, $\det u$ is constant so one can choose an $r$-th root $(\det u)^{1/r} \in S^1$ and define the  isomorphism \eqref{eq:Isom_of_StabilizerGroups} by $u\mapsto [u (\det u)^{-1/r}, (\det u)^{1/r}\,\id_E]$.

By \eqref{eq:Isom_of_StabilizerGroups}, we have $\Stab(A)=C_r$ if and only if $\Stab_{\U(E)}(A)=S^1$. The definition of central isotopy in Item \eqref{item:central-stabilizer_unitary_connection} of Definition \ref{defn:Reducible_split_trivial-stabilizer_unitary_connection} is thus equivalent to requiring that
$\Stab_{\U(E)}(A)=S^1$.
\end{rmk}

Recall from Donaldson and Kronheimer \cite[Lemma 4.2.8, p. 132]{DK} or Rudolph and Schmidt \cite{Rudolph_Schmidt_differential_geometry_mathematical_physics_part2} that if $G$ is a compact Lie group and $A$ is a connection on a principal $G$-bundle over a connected base manifold, then the stabilizer $\Stab(A)$ is isomorphic to the centralizer $C_G(\Hol_A)$ of the holonomy group $\Hol_A$ of $A$ in $G$. The structure group $G$ of the principal bundle defining the Hermitian vector bundle $E$ in Definition \ref{defn:Reducible_split_trivial-stabilizer_unitary_connection} is $\U(r)$ or $\SU(r)$. As Donaldson and Kronheimer point out in \cite[Section 4.2, p. 133]{DK}, it is common to abuse terminology and say that connections with central stabilizer are `irreducible' and `reducible' otherwise.  They give the example of a connection on an $\SU(r)$ bundle with holonomy given by $\SO(r)$.  Such a connection would be reducible in the sense of
Definition \ref{defn:Reducible_split_trivial-stabilizer_unitary_connection} \eqref{item:Reducible_unitary_connection}, but would have central stabilizer. Furthermore, the definition of a split connection in Definition \ref{defn:Reducible_split_trivial-stabilizer_unitary_connection}\eqref{item:Split_unitary_connection} requires the holonomy of the connection to be contained in the proper subgroup $\U(r_1)\times \U(r_2)$ of $\U(r)$. L\"ubke and Teleman \cite[Definition 1.1.16, p. 25]{Lubke_Teleman_1995} say that a unitary connection $A$ on $(E,h)$ is `irreducible' if the kernel of $d_A$ on $C^\infty(\fu(E))$ is equal to $i\RR\,\id_E$ and say that $A$ is `reducible' otherwise. Although the concept of split connections and complex vector bundles can be expressed in terms of reducible connections and principal bundles as in Kobayashi and Nomizu \cite[Sections II.7 and II.8]{Kobayashi_Nomizu_v1}, we find it more convenient to work in the setting of complex vector bundles and we use the concept of split connections because that is closely linked to customary usage in the complex geometry literature for holomorphic vector bundles.

\begin{rmk}[Twisted reducibles]
\label{rmk:TwistedReducibles}
The \emph{twisted reducible} unitary connections on a rank-two, smooth Hermitian vector bundle over a smooth manifold of dimension four discussed in Kronheimer and Mrowka \cite[p. 586]{KMStructure}, have
central stabilizer by Definition \ref{defn:Reducible_split_trivial-stabilizer_unitary_connection}
because (as noted in \cite[p. 586]{KMStructure}) their stabilizer is $\{\pm\id_E\}$ but are reducible in the sense of Definition \ref{defn:Reducible_split_trivial-stabilizer_unitary_connection} \eqref{item:Reducible_unitary_connection} as their holonomy is given by $\Or(2)$.
\end{rmk}

We will use the  following algebraic definition to classify the possible stabilizers of a unitary connection.  Let $K(r)$ be the set of ordered pairs of strings $J$  of equal length of positive integers,
\[
J:=((k_1,k_2,\dots,k_\ell),(m_1,m_2,\dots,m_\ell)),\quad\text{for some $\ell$ with $1\le \ell\le r$},
\]
which satisfy $\sum_{i=1}^\ell k_im_i=r$. Each $J\in K(r)$ defines a decomposition,
\[
\CC^r\cong \bigoplus_{i=1}^\ell \left( \CC^{k_i}\otimes \CC^{m_i}\right),
\]
and thus an injective homomorphism of Lie groups
\begin{equation}
\label{eq:Tensor_Embedding}
\rho_J:\prod_{i=1}^\ell \U(k_i) \ni (M_1,\dots,M_\ell)
\mapsto \bigoplus_{i=1}^\ell \left( M_i\otimes\id_{\CC^{m_i}}\right) \in \U(r).
\end{equation}
One defines
\begin{equation}
\label{eq:Definition_SU_J}
\SU_J:= \rho_J\left( \prod_{i=1}^\ell \U(k_i)\right)\cap\SU(r).
\end{equation}
Compare the description of $\Stab(A)$ appearing in Kronheimer \cite[Section 2.1, p. 63]{Kronheimer_2005}.  The following result classifies Howe subgroups of $\SU(r)$.

\begin{lem}[Howe subgroups of $\SU(r)$]
\label{lem:Howe_Subgroups_of_SU(n)}
(See Rudolph and Schmidt \cite[Lemma 8.5.1, p. 665]{Rudolph_Schmidt_differential_geometry_mathematical_physics_part2}.)
Let $H\subset\SU(r)$ be a Howe subgroup of $\SU(r)$, that is, $H$ is the centralizer of some subset of $\SU(r)$.  Then $H$ is conjugate to the subgroup $\SU_J$ defined in \eqref{eq:Definition_SU_J} for some $J\in K(r)$.
\end{lem}

We provide details in addition to those in the proof given by Kronheimer of his result below, noting that Kronheimer's definition of `irreducible' (see \cite[Definition 2.1, p. 62]{Kronheimer_2005}) is equivalent to our definition of central isotopy.

\begin{lem}[Structure of zero-dimensional stabilizers of unitary connections]
\label{lem:Kronheimer_Lemma2.2}
(See Kronheimer \cite[Lemma 2.2, p. 64]{Kronheimer_2005}.)
Let $A$ be a $W^{1,p}$ unitary connection with $p\in (d/2,\infty)$ on a $W^{2,p}$ Hermitian vector bundle $(E,h)$ of complex rank $r\geq 2$ over a smooth, connected manifold of dimension $d\geq 2$. If $\Stab(A)$ is the stabilizer of $A$ in $W^{2,p}(\SU(E))$ and if $\dim\Stab(A)=0$, then $A$ has central stabilizer as in Definition \ref{defn:Reducible_split_trivial-stabilizer_unitary_connection} (\ref{item:central-stabilizer_unitary_connection}).
\end{lem}

\begin{proof}
Because elements of $\Stab(A)$ are parallel sections of $\End(E)$ with respect to the connection induced on $\End(E)$ by $A$, evaluation at a fixed point $x_0\in X$ and a choice of orthonormal frame for $E_{x_0}$ define a monomorphism $\Stab(A)\to \SU(r)$ of Lie groups. We will thus consider $\Stab(A)$ as a subgroup of $\SU(r)$. Because $\Stab(A)$ is the centralizer of the holonomy group of $A$ by Donaldson and Kronheimer \cite[Lemma 4.2.8, p. 132]{DK} or Rudolph and Schmidt \cite[Theorem 6.1.5, p. 467]{Rudolph_Schmidt_differential_geometry_mathematical_physics_part2}, the group $\Stab(A)$ is by definition a \emph{Howe subgroup} as in Rudolph and Schmidt \cite[Definition 8.2.4, p. 640]{Rudolph_Schmidt_differential_geometry_mathematical_physics_part2}. By Lemma \ref{lem:Howe_Subgroups_of_SU(n)}, the group $\Stab(A)$ is conjugate to a group $\SU_J$ for $J\in K(r)$ as defined in \eqref{eq:Definition_SU_J}.

We claim that the dimension of $\SU_J$ will be positive unless $\ell=1$, $k_1=1$, and $m_1=r$. We first consider the case $\ell>1$. Observe that
\begin{equation}
\label{eq:TensorDeterminant}
\det\left( (e^{i\theta}\,\id_{\CC^{k_i}})\otimes\id_{\CC^{m_i}}\right)=e^{i\theta m_ik_i}.
\end{equation}
If $\ell>1$, then
\begin{align*}
&\det \rho_J
\left(
    e^{i\theta m_2k_2}\,\id_{\CC^{k_1}},
    e^{-i\theta m_1k_1}\,\id_{\CC^{k_2}},
    \id_{\CC^{k_3}},\dots,\id_{\CC^{k_\ell}}
\right)
\\
&\quad
=
\det
\left(
    (e^{i\theta m_2k_2}\,\id_{\CC^{k_1}}\otimes\id_{\CC^{m_1}})
\oplus
    (e^{-i\theta m_1k_1}\,\id_{\CC^{k_2}}\otimes\id_{\CC^{m_2}})
\oplus
\oplus_{i=3}^\ell \id_{\CC^{k_im_i}}
\right)
\\
&\quad=
\det\left( (e^{i\theta m_2k_2}\,\id_{\CC^{k_1}})\otimes\id_{\CC^{m_1}}\right)
\det\left( (e^{-i\theta m_1k_1}\,\id_{\CC^{k_2}})\otimes\id_{\CC^{m_2}}\right)
\\
&\quad=
e^{i\theta m_1k_1m_2k_2} e^{-i\theta m_1k_1 m_2k_2}
\quad\text{(by \eqref{eq:TensorDeterminant})}
\\
&\quad=1.
\end{align*}
Hence, if $\ell>1$ the restriction of $\rho_J$ to the circle given by
\[
S_J^1:= \{\left(
    e^{i\theta m_2k_2}\,\id_{\CC^{k_1}},
    e^{-i\theta m_1k_1}\,\id_{\CC^{k_2}},
    \id_{\CC^{k_3}},\dots,\id_{\CC^{k_\ell}}
\right)\}\subset\prod_{i=1}^\ell\U(k_i)
\]
has image contained in $\SU(r)$. Therefore, $\rho_J(S_J^1)$ is a one-dimensional subspace contained in
$\SU_J$. Consequently, if $\ell>1$, then $\SU_J$ has positive dimension.

We now assume that $\ell=1$. If $k_1>1$, then $\dim \SU(k_1)$ is positive and $\rho_J$ embeds $\SU(k_1)$ in $\SU(r)$. Hence, if $\ell=1$ and $k_1>1$, then $\SU_J$ also has positive dimension.  This completes the proof of the claim that the dimension of $\SU_J$ is positive unless $\ell=1$, $k_1=1$, and $m_1=r$.

Finally, if $\ell=1$ and $k_1=1$ (and so $m_1=r$), then because $\U(k_1)\cong S^1$ and 
\[
\det \rho_J(e^{i\theta})=e^{im_1\theta}=e^{ir\theta},
\]
we see that $\rho_J(e^{i\theta})$ lies in $\SU(r)$ if and only if $e^{ir\theta}=1$, that is, $e^{i\theta}$ is an $r$-th root of unity, as asserted.
\end{proof}

We give the following application of Lemma \ref{lem:Howe_Subgroups_of_SU(n)} to recover the well-known classification of stabilizers of $\SU(2)$ connections.

\begin{exmp}[Stabilizers of unitary connections on a complex rank-two Hermitian vector bundle]
\label{exmp:Stab(A)_r=2}
Continue the notation and hypotheses of Lemma \ref{lem:Kronheimer_Lemma2.2}. We now discuss the subgroups $\SU_J$ where $J\in K(2)$ and thus classify the possible stabilizers in $W^{2,p}(\SU(E))$ of unitary connections on a rank-two Hermitian vector bundle. If $r=2$, the set $K(2)$ consists of the following pairs of strings:
\begin{align*}
J(1)&=\left((1,1),(1,1)\right),
\\
J(2) &=\left( (1),(2)\right),
\\
J(3)&=\left( (2),(1)\right).
\end{align*}
For $J=J(1)$, the embedding $\rho_{J(1)}$  given by \eqref{eq:Tensor_Embedding} is
\[
\rho_{J(1)}:\U(1)\times\U(1) \ni \left( e^{i\theta_1},e^{i\theta_2}\right)
\mapsto
\left( e^{i\theta_1}\id_{\CC}\right)\oplus \left( e^{i\theta_2}\id_{\CC}\right)
\in\U(2).
\]
Hence, the intersection $\SU_{J(1)}$ of the image of $\rho_{J(1)}$ with $\SU(2)$ is the standard maximal torus of $\SU(2)$:
\[
  \SU_{J(1)}
  =
\left\{
\begin{pmatrix} e^{i\theta} & 0 \\ 0 & e^{-i\theta}\end{pmatrix}: e^{i\theta}\in S^1
\right\}
\cong S^1.
\]
For $J=J(2)$, the embedding given by \eqref{eq:Tensor_Embedding} is
\[
  \rho_{J(2)}:\U(1) \ni e^{i\theta} \mapsto e^{i\theta}\,\id_{\CC^2} \in \U(2).
\]  
The intersection $\SU_{J(2)}$ of the image of $\rho_{J(2)}$ with $\SU(2)$ is the cyclic group,
\[
\SU_{J(2)} = \{\pm\id_{\CC^2}\} \cong \ZZ/2\ZZ.
\]
For $J=J(3)$, the embedding given by \eqref{eq:Tensor_Embedding} is the identity map $\rho_{J(3)}:\U(2)\to \U(2)$, so the intersection $\SU_{J(3)}$ of the image of $\rho_{J(3)}$ with $\SU(2)$ is
\[
\SU_{J(3)} = \SU(2).
\]
Stabilizers in $W^{2,p}(\SU(E))$ of unitary connections on a rank-two Hermitian vector bundle are given by one of the preceding three subgroups of $\SU(2)$, namely the maximal torus isomorphic to $S^1$, a subgroup isomorphic to the cyclic group $\ZZ/2\ZZ$, and $\SU(2)$ itself. 
\end{exmp}

For the proof of the forthcoming Proposition \ref{prop:FixedPointsOfS1ActionOnSpinuQuotientSpace} and other applications throughout our monograph, we shall require the following refinement of Freed and Uhlenbeck's result \cite[Theorem 3.1, p. 47]{FU}.

\begin{prop}[Characterization of unitary connections with non-discrete stabilizer]
\label{prop:ReducibleUnitaryFromH0}
(Compare Freed and Uhlenbeck \cite[Theorem 3.1, p. 47]{FU} for the case of rank two and Zentner \cite[Lemma 2.2, p. 237]{Zentner_2012} for the case of higher rank.)
Let $(E,h)$ be a smooth Hermitian vector bundle with complex rank $r\geq 2$ over a connected, smooth manifold $X$ and $A$ be a smooth, unitary connection on $E$. If $\xi$ is a $W^{2,p}$ section of $\su(E)$ such that $\xi\not\equiv 0$ and $d_A\xi=0$, then the eigenvalues of $\xi(x)$ are constant with respect to $x\in X$, purely imaginary, and
\begin{equation}
  \label{eq:Splitting_E_constant_eigenvalue_subbundles}
  E = E_1\oplus\cdots\oplus E_s
\end{equation}
as an orthogonal direct sum of proper, smooth Hermitian subbundles of $E$ corresponding to the distinct eigenvalues $i\mu_1,\ldots,i\mu_s$ of $\xi$ for $2\leq s\leq r$, where $\sum_{k=1}^s(\rank_\CC E_k)\mu_k = 0$, and the following hold:
\begin{enumerate}
\item\label{item:ReducibleUnitaryFromH0Splitting_smooth_and_any_rank_differentiable}
The section $\xi$ is smooth and $A$ is split in the sense of Item \eqref{item:Split_unitary_connection} in  Definition \ref{defn:Reducible_split_trivial-stabilizer_unitary_connection} with respect to the decomposition \eqref{eq:Splitting_E_constant_eigenvalue_subbundles}, so $A = A_1\oplus\cdots\oplus A_s$, and each $A_k$ is a 
smooth unitary connection on $E_k$, for $k=1,\ldots,s$.

\item\label{item:ReducibleUnitaryFromH0Splitting_smooth_and_any_rank_analytic}
If $(X,g)$ and $(E,h)$ and $A$ are real analytic, then $\xi$ is real analytic and the subbundles $E_k$ and connections $A_k$ in Item \eqref{item:ReducibleUnitaryFromH0Splitting_smooth_and_any_rank_differentiable} are real analytic, for $k=1,\ldots,s$.
\end{enumerate}
\end{prop}

\begin{rmk}[A non-central element of $\Ker d_A$]
We emphasize that the non-zero element of $\Ker d_A$ is a section of $\su(E)$ rather than one of $\fu(E)$. By writing $\fu(E)=\su(E)\oplus i\ubarRR$, where $\ubarRR$ is the product line bundle $X\times\RR$, we see that every unitary connection $A$ contains non-zero, constant sections of $i\ubarRR$ in $\Ker d_A$.
\label{rmk:NonCentralElementOfStabilizer}
\end{rmk}

\begin{rmk}[Assumption of connectedness in Proposition \ref{prop:ReducibleUnitaryFromH0}]
\label{rmk:ReducibleUnitaryFromH0_connectedness}  
We could remove connectedness from the hypotheses of Proposition \ref{prop:ReducibleUnitaryFromH0} by writing $X$ as a disjoint union of connected components and applying Proposition \ref{prop:ReducibleUnitaryFromH0} to each connected component separately.  However, there could be a different decomposition of $E$ and different eigenvalues of $\xi$ on each connected component.
\end{rmk}

\begin{rmk}[Properties of the connections $A_k$]
\label{rmk:Properties_connections_Ak}
The unitary connections $A_k$ on $E_k$ in Item \eqref{item:ReducibleUnitaryFromH0Splitting_smooth_and_any_rank_differentiable} of Proposition \ref{prop:ReducibleUnitaryFromH0} could be split connections.
\end{rmk}  

\begin{proof}[Proof of Proposition \ref{prop:ReducibleUnitaryFromH0}]
A version of Item \eqref{item:ReducibleUnitaryFromH0Splitting_smooth_and_any_rank_differentiable} appears in the proof of \cite[Theorem 3.1, p. 47]{FU} when $E$ has rank two and without proof in \cite[Lemma 2.2, p. 237]{Zentner_2012} for the case of $E$ with higher rank. The equivalence of Items (a) and (d) in Kobayashi \cite[Proposition 7.4.14, p. 242]{Kobayashi_differential_geometry_complex_vector_bundles} provides a similar result, but is also stated without proof.

Consider Item \eqref{item:ReducibleUnitaryFromH0Splitting_smooth_and_any_rank_differentiable}. Because $\xi$ is a section of $\su(E)$ then, for every point $x\in X$, the linear operator $i\xi(x)\in i\su(E_x) \subset \gl(E_x)$ is Hermitian and thus normal. According to the Spectral Theorem for bounded normal operators on a Hilbert space (see Axler \cite[Theorem 7.24, p. 218]{Axler_linear_algebra_done_right} or Rudin \cite[Theorem 12.23, p. 324]{Rudin} in the finite or infinite-dimensional cases, respectively) and the fact that $i\xi(x)\in\gl(E_x)$ is normal for each $x\in X$, the operator $\xi(x)$ is unitarily equivalent to a diagonal operator $\Lambda(x)$ with eigenvalues $i\lambda_j(x)$ such that $\sum_{j=1}^r\lambda_j = 0$, where $\lambda_j(x)$ is a real-valued function of $x\in X$, for $j=1,\ldots,r$, and $r$ is the complex rank of $E$. Thus, $\xi(x) = u(x)^\dagger\Lambda(x)u(x)$, where $u(x)\in \U(E_x)$ for each $x\in X$.

Each eigenvalue $i\la_j(x)$ is semisimple in the sense of Kato \cite[Section I.5.4, p. 41]{Kato}, for $j=1,\ldots,r$ and so by Kato \cite[Theorem II.5.1 p. 107]{Kato}, the eigenvalue functions $\lambda_j(x)$ are continuous with respect to $x\in X$, for $j=1,\ldots,r$. (Note that $\xi$ is a continuous section of by the Sobolev embedding $W^{2,p}(X) \subset C(X)$ --- see Adams and Fournier \cite[Theorem 4.12, p. 85]{AdamsFournier}.)

We now check that the connection $A$ is split with respect to an orthogonal splitting of $E$ as a direct sum of proper Hermitian subbundles. The equality $d_A^*d_A\xi=0$ and, noting that $A$ is smooth by assumption, regularity theory for the elliptic operator $d_A^*d_A$ with smooth coefficients imply that $\xi$ is smooth. Because $\xi(x)$ is diagonalizable, the eigenvalues $\lambda_j(x)$ are differentiable functions of $x\in X$ by Kato \cite[Theorem II.5.4, p. 111]{Kato}. Let $x_0 \in X$ be a point and suppose that $i\mu_k(x)$, for $k=1,\ldots s$ with $s\leq r$, are the \emph{distinct} eigenvalues of $\xi(x)$, for $x$ in a connected open neighborhood $U_0\subset X$ of the point $x_0$.
According to Kato \cite[Theorem II.5.4, p. 111]{Kato}, the eigenprojections $\pi_k(x)$ onto the $\mu_k(x)$-eigenspaces are differentiable functions of $x\in U_0$. Fix $l\in\{1,\ldots,s\}$ and let $\varphi$ be a smooth section of $E$ such that $\varphi(x)\neq 0$ for all $x\in U_0$. By multiplying $\varphi$ by a positive, smooth function on $U_0$ and relabeling, we may assume that $|\pi_l(x)\varphi(x)| = 1$ for all $x\in U_0$. We may thus differentiate the equalities $\xi\pi_k\varphi = i\mu_k\pi_k\varphi$ on $U_0$ to give
\[
  d_A(\xi\pi_k\varphi) = d_A(i\mu_k\pi_k\varphi) = i(d\mu_k)\pi_k\varphi + i\mu_kd_A\pi_k\varphi
  \quad\text{on } U_0,
  \text{ for } k=1,\ldots,s,
\]
while the identity $d_A\xi=0$ implies that
\[
  d_A(\xi\pi_k\varphi) = (d_A\xi)\pi_k\varphi + \xi d_A\pi_k\varphi =  \xi d_A\pi_k\varphi\quad\text{on } U_0,
  \text{ for } k=1,\ldots,s.
\]
Combining the preceding identities yields
\begin{equation}
\label{eq:CovDerivOfEigensection}
\xi d_A\pi_k\varphi = i(d\mu_k)\pi_k\varphi + i\mu_kd_A\pi_k\varphi\quad\text{on } U_0,
\text{ for } k=1,\ldots,s.
\end{equation}
We also have $\langle\pi_k\varphi,\pi_l\varphi\rangle = \delta_{kl}$ on $U_0$ and thus $d\langle\pi_k\varphi,\pi_l\varphi\rangle = 0$ on $U_0$ for $k=1,\ldots,s$ and so, using the fact that $A$ is unitary,
\[
  \langle d_A\pi_k\varphi,\pi_l\varphi\rangle + \langle \pi_k\varphi,d_A\pi_l\varphi\rangle = 0\quad\text{on } U_0,
  \text{ for } k=1,\ldots,s,
\]
and in particular, $\Real\langle d_A\pi_k\varphi,\pi_l\varphi\rangle = 0$, for $k=1,\ldots,s$. Taking the inner product of \eqref{eq:CovDerivOfEigensection} with $\pi_l\varphi$ and choosing $k=l$ gives
\[
  \langle \xi d_A\pi_l\varphi,\pi_l\varphi\rangle = id\mu_l + i\mu_l\langle d_A\pi_l\varphi,\pi_l\varphi\rangle
  \quad\text{on } U_0.
\]
The imaginary part of the preceding equation, after substituting $\Real \langle d_A\pi_l\varphi,\pi_l\varphi\rangle = 0$, yields
\begin{equation}
\label{eq:ds_j_equation}
  \Imag\langle \xi d_A\pi_l\varphi,\pi_l\varphi\rangle = d\mu_l\quad\text{on } U_0.
\end{equation}
We now observe that, using $\xi^\dagger = -\xi$,
\begin{multline*}
  \Imag\langle \xi d_A\pi_k\varphi,\pi_k\varphi\rangle
  = -\Imag\langle d_A\pi_k\varphi,\xi\pi_k\varphi\rangle
  = -\Imag\langle d_A\pi_k\varphi,i\mu_k\pi_k\varphi\rangle
  \\
  = \Imag\left(i\mu_k\langle d_A\pi_k\varphi,\pi_k\varphi\rangle\right)
  = \mu_k\Real\langle d_A\pi_k\varphi,\pi_k\varphi\rangle
  = 0\quad\text{on } U_0, \text{ for } k=1,\ldots,s.
\end{multline*}
Substituting the resulting identity, $\Imag\langle \xi d_A\pi_k\varphi,\pi_k\varphi\rangle = 0$ on $U_0$, into \eqref{eq:ds_j_equation} with $k=l$ yields $d\mu_l=0$ on $U_0$ and thus $\mu_l$ is constant over $U_0$. Since $l\in\{1,\ldots,s\}$ was arbitrary, the eigenvalues $\mu_k$ are constant over $U_0$, for $k=1,\ldots,s$.

By applying the preceding argument to each point in $X$, we can find an open cover $\{U_\alpha\}$ of $X$ and
for each $\alpha$, real-valued, constant functions $\mu_{\alpha,i}$ on $U_\alpha$ for $i=1,\dots,s_{\alpha}$ such that, for any $x\in X$, the distinct eigenvalues of $\xi(x)$ are $\mu_{\alpha,i}(x)$ for $i=1,\dots,s_{\alpha}$.
By re-indexing if necessary, we can assume that $\mu_{\alpha,i}<\mu_{\alpha,j}$ for each $\alpha$ if $i<j$. If $x\in U_{\alpha}\cap U_{\beta}$, then  $\mu_{\alpha,i}(x)$ and $\mu_{\beta,i}(x)$ are both the $i$-th smallest eigenvalue of $\xi(x)$ and so are equal. Thus, 
\[
\mu_i(x):= \mu_{\alpha,i}(x)\quad\text{if $x\in U_\alpha$}
\]
gives a global, constant function $\mu_i$ on $X$ which is an eigenvalue of $\xi(x)$ for all $x\in X$.

Hence, the associated eigenprojections $\pi_k$ of $E$ are globally defined on $X$, for $k=1,\ldots,s$. Because each $\mu_k$ is constant on $X$, equation \eqref{eq:CovDerivOfEigensection} now implies that for any smooth section $\phi$ of $E$ we have $\xi d_A\pi_k\phi = i\mu_kd_A\pi_k\phi$ over $X$ and thus $d_A\pi_k\phi$ takes values in the differentiable (but not necessarily smooth), proper Hermitian subbundle,
\[
  E_k := \pi_kE \subsetneq E, \quad\text{for } k=1,\ldots,s.
\]
Hence, since $\phi$ is arbitrary, $A$ is split with respect to a pointwise orthogonal splitting of $E$ as a direct sum of differentiable Hermitian subbundles in the sense that
\[
  d_A = \bigoplus_{k=1}^s d_{A_k} \quad\text{on}\quad \bigoplus_{k=1}^s E_k.
\]
According to the cited Spectral Theorem, we may write
\begin{equation}
  \label{eq:Spectral_theorem}
  \xi(x) = \sum_{k=1}^s i\mu_k(x)\pi_k(x), \quad\text{for all } x\in X,
\end{equation}
and since $\xi$ is not identically zero, we must have $\mu_k\neq 0$ for at least once index $k\in\{1,\ldots,s\}$. 
Because $d_A\xi=0$ and the eigenvalues $\mu_k$ are constant, then \eqref{eq:Spectral_theorem} yields
\[
  d_A\xi = \sum_{k=1}^s i\mu_kd_A\pi_k = 0, 
\]
where the preceding equality holds for $\xi \in C^0(\End(E))$ in the sense of distributions. Because $d_A\pi_k\phi$ is a section of $E_k$ by the preceding discussion for any $\phi \in C^\infty(E)$, then $d_A\pi_k = 0$ and
\[
  d_A^*d_A\pi_k = 0,
\]
where both equalities hold for $\pi_k \in C^0(\End(E))$ in the sense of distributions for $k=1,\ldots,s$. Since the Laplace operator $d_A^*d_A$ has smooth coefficients, it follows that $\pi_k \in C^\infty(\End(E))$ for $k=1,\ldots,s$ (see Feehan and Maridakis \cite[Remark 2.2.2, p. 23]{Feehan_Maridakis_Lojasiewicz-Simon_coupled_Yang-Mills} for an explanation). Consequently, the Hermitian subbundles $E_k$ and unitary connections $A_k$ are smooth since $d_{A_k} = \pi_kd_A\pi_k$ for $k=1,\ldots,s$. This verifies Item \eqref{item:ReducibleUnitaryFromH0Splitting_smooth_and_any_rank_differentiable}.

To verify Item \eqref{item:ReducibleUnitaryFromH0Splitting_smooth_and_any_rank_analytic}, we observe that eigenfunctions $\mu_k(x)$ and associated eigenprojections $\pi_k(x)$ are now real analytic functions of $x\in U_0$ by Kato \cite[Section II.1.1, pp. 63--64 and Section II.1.3, pp. 66--67]{Kato}. It follows from the proof of Item \eqref{item:ReducibleUnitaryFromH0Splitting_smooth_and_any_rank_differentiable} that the subbundles $E_k$ and connections $A_k$ are real analytic. This verifies Item \eqref{item:ReducibleUnitaryFromH0Splitting_smooth_and_any_rank_analytic} and completes the proof of Proposition \ref{prop:ReducibleUnitaryFromH0}.
\end{proof}

The following corollary allows us to apply the results of Proposition \ref{prop:ReducibleUnitaryFromH0} to a Sobolev connection that is only gauge equivalent to a smooth connection.

\begin{cor}[Characterization of unitary connections with non-discrete stabilizer that are $W^{2,p}$-gauge-equivalent to a smooth connection]
\label{cor:ReducibleUnitaryFromH0_GaugeEquivalentToSmooth}
Let $(E,h)$ be a smooth Hermitian vector bundle with complex rank $r\geq 2$ over a smooth, connected manifold $X$ and $A$ be a smooth unitary connection on $E$. If $A'=u^*A$, where $u\in W^{2,p}(\SU(E))$, and there is a section $\xi'\in W^{2,p}(\su(E))$ such that $\xi'\not\equiv 0$ and $d_{A'}\xi'=0$, then the eigenvalues of $\xi'(x)$ are constant with respect to $x\in X$, purely imaginary, and $E = E_1'\oplus\cdots\oplus E_s'$ as an orthogonal direct sum of proper, $W^{2,p}$ Hermitian subbundles of $E$ corresponding to the distinct eigenvalues $i\mu_1,\ldots,i\mu_s$ of $\xi'$ for $2\leq s\leq r$, and $A'=A_1'\oplus\cdots\oplus A_s'$, where each $A_k'$ is a $W^{1,p}$ unitary connection on $E_k'$.
\end{cor}

\begin{proof}
A gauge transformation $u\in W^{2,p}(\SU(E))$ acts on a section $\xi$ of $\su(E)$ by conjugation
which we will denote by $\ad(u^{-1})\xi:= u^{-1}\circ\xi\circ u$. The equality $d_{A'}=u^{-1}\circ d_A\circ u$ implies that
\[
  d_{A'}\xi = \ad(u^{-1}) d_A(\ad(u)\xi), \quad\text{for all } \xi \in W^{2,p}(\su(E)).
\]  
If $d_{A'}\xi'=0$, then $d_A(\ad(u)\xi')=\ad(u)d_{A'}\xi'=0$ and so there exists $\xi := \ad(u)\xi' \in \Ker d_A\cap W^{2,p}(\su(E))$ with $\xi\not\equiv 0$.  Proposition \ref{prop:ReducibleUnitaryFromH0} implies that the eigenvalues of $\xi$ are constant and that there are proper, smooth Hermitian subbundles $E_k\subset E$ for $k=1,\dots,s$, an orthogonal decomposition $E=E_1\oplus\cdots\oplus E_s$ corresponding to the distinct eigenvalues $i\mu_1,\ldots,i\mu_s$ of $\xi$, and smooth unitary connections $A_k$ on $E_k$ such that $A=A_1\oplus\cdots\oplus A_s$.  The definition $\xi = \ad(u)\xi'=u\circ\xi'\circ u^{-1}$ implies that the eigenvalues of $\xi'$ are equal to those of $\xi$ and thus are also constant.

Define $E_k':=u^{-1}(E_k)$ for $k=1,\ldots,s$. Because $u$ is a unitary gauge transformation, the decomposition $E=E_1'\oplus\cdots\oplus E_s'$ is an orthogonal direct sum. 
Moreover, if we write $u_k:=u|_{E_k'}$, then $u_k:E_k'\to E_k$ is a unitary isomorphism and
\begin{equation}
\label{eq:U_decomposition}
u = u_1\oplus\cdots\oplus u_s.
\end{equation}
If $v_k\in W^{1,p}(E_k)$, then the definition $\xi=\ad(u)\xi'$ implies that
\[
\xi' u^{-1}v_k
=
u^{-1}\xi v_k
=
i\mu_ku^{-1}v_k, \quad\text{for } k = 1,\ldots,s.
\]
Therefore, $E_k'$ is the subbundle of $E'$ corresponding to the eigenvalue $i\mu_k$ of $\xi'$ for $k=1,\dots,s$.

Define $A_k':=u_k^*A_k$ and observe that $A_k'$ is a $W^{1,p}$ connection on $E_k'$ for $k=1,\dots,s$.
By \eqref{eq:U_decomposition}, we have
\begin{align*}
A' &=u^*A=(u_1\oplus\cdots\oplus u_s)^* (A_1\oplus\cdots\oplus A_s)
  \\
  &= u_1^*A_1\oplus\cdots \oplus u_s^*A_s
  \\
  &= A_1'\oplus\cdots\oplus A_s' \quad\text{on } E_1'\oplus\cdots\oplus E_s'.
\end{align*}
This completes the proof of Corollary \ref{cor:ReducibleUnitaryFromH0_GaugeEquivalentToSmooth}.
\end{proof}

Versions of the following result are well-known, but we include it for ease of reference and precision.

\begin{cor}[Split unitary connections and $\bH_A^0$]
\label{cor:Split_unitary_A_and_Lie_Alg_of_Stab(A)}  
Let $(E,h)$ be a smooth Hermitian vector bundle over a closed, connected, smooth manifold $X$. For $A\in\sA(E,h)$, the vector space $\bH_A^0$, defined in Lemma \ref{lem:LieGroup_Structure_of_Stab(A)} and the forthcoming \eqref{eq:HE_equation_bHA0}, the following hold:
\begin{enumerate}
\item
\label{item:HA0_is_Lie_Alg_of_Stab(A)_Irreducible}
The connection $A$ has central stabilizer as in Definition \ref{defn:Reducible_split_trivial-stabilizer_unitary_connection} if and only if $\bH_A^0=(0)$.

\item
\label{item:A_smooth_and_bHA_non-zero_implies_A_split}
If $A$ is smooth and $\bH_A^0\neq (0)$, then $A$ is split as in Definition \ref{defn:Reducible_split_trivial-stabilizer_unitary_connection}.

\item
\label{item:A_split_implies_bHA_non-zero}
If $A$ is split as in Definition \ref{defn:Reducible_split_trivial-stabilizer_unitary_connection}, then $\bH_A^0\neq (0)$.
\end{enumerate}
\end{cor}

\begin{proof}
Consider  Item \eqref{item:HA0_is_Lie_Alg_of_Stab(A)_Irreducible}. By Lemma \ref{lem:LieGroup_Structure_of_Stab(A)}, the Lie algebra of $\Stab(A)$ is $\bH_A^0$. Observe that $\Stab(A)$ is zero-dimensional if and only if its Lie algebra $\bH_A^0$ is zero-dimensional. If $A$ has central stabilizer, then $\Stab(A)$ and thus $\bH_A^0$ are zero-dimensional. If $\bH_A^0=(0)$, then $\Stab(A)$ is zero-dimensional and $A$ has central stabilizer by Lemma \ref{lem:Kronheimer_Lemma2.2}.

Consider Item \eqref{item:A_smooth_and_bHA_non-zero_implies_A_split}. By assumption, $A$ is smooth and $\bH_A^0\neq (0)$, so there is a non-zero section $\xi\in C^\infty(\su(E))$ in the kernel of $d_A$. The conclusion thus follows from Item \eqref{item:ReducibleUnitaryFromH0Splitting_smooth_and_any_rank_differentiable} in Proposition \ref{prop:ReducibleUnitaryFromH0}.

Consider Item \eqref{item:A_split_implies_bHA_non-zero}. By assumption, $A=A_1\oplus\cdots\oplus A_s$ with respect to a splitting $E = E_1\oplus\cdots\oplus E_s$ as an orthogonal direct sum of proper, smooth Hermitian vector subbundles with $s\geq 2$, where $A_k$ is a smooth unitary connection on $E_k$ for $k=1,\ldots,s$. Define a smooth one-parameter family of gauge transformations $u(\theta) \in W^{2,p}(\U(E))$ by
\begin{equation}
  \label{eq:One-parameter_family_SU(E)_gauge_transformations}
  u(\theta)
  :=
  e^{-i(r-r_1)\theta}\,\id_{E_1}
  \oplus
  e^{ir_1\theta}\,\id_{E_1^\perp}, \quad\text{for all } \theta \in \RR,
\end{equation}
where $r_1$ is the complex rank of $E_1$. Because $E_1^\perp = E_2\oplus\cdots\oplus E_s$ has complex rank $r-r_1$ and
\[
  \det\left(e^{i\theta}\,\id_{E_1}\right) = e^{ir_1\theta}\id_{\det E_1}
  \quad\text{and}\quad
  \det\left(e^{i\theta}\,\id_{E_1^\perp}\right) = e^{i(r-r_1)\theta}\id_{\det E_1^\perp},
\]
we see that
\[
\det u(\theta) = e^{i\theta (-(r-r_1)r_1 + r_1(r-r_1))} = \id_{\det E}, \quad\text{for all } \theta \in \RR.
\]
Hence, the gauge transformations $u(\theta)$ belong to $ W^{2,p}(\SU(E))$ and define a one-parameter family of $\Stab(A)$.  The existence of a smooth one-parameter family in $\Stab(A)$ and the identification of $\bH_A^0$ with the Lie algebra of $\Stab(A)$ imply that $\bH_A^0\neq (0)$. This completes the proof of the corollary.
\end{proof}

\section{Clifford modules and spin connections}
\label{sec:SpincuStr}
Let $V$ be a smooth Hermitian vector bundle over an oriented, smooth, Riemannian manifold $(X,g)$ of even dimension and let $\rho:T^*X\to \End(V)$ be a real-linear map satisfying (see Lawson and Michelsohn \cite[Proposition I.5.10, p. 35]{LM})
\begin{equation}
\label{eq:CliffordMapDefn}
\rho(\alpha)^2 = -g(\alpha,\alpha)\id_{V}
\quad\text{and}\quad
\rho(\alpha)^\dagger = -\rho(\alpha),
\quad\text{for all } \alpha \in C^\infty(T^*X).
\end{equation}
The map $\rho$ uniquely extends to a linear map $\rho:\Lambda^{\bullet}(T^*X)\otimes_\RR\CC\to\End(V)$, and gives $V$ the structure of a \emph{Hermitian Clifford module} for the complex Clifford algebra $\CCl(T^*X)$ (see Lawson and Michelson \cite[Chapter I, Proposition 5.10, p. 35]{LM}).  There is a splitting $V=V^+\oplus V^-$, where $V^\mp$ are the $\pm 1$ eigenspaces of $\rho(\vol)$. A unitary connection $A$ on $V$ is called \emph{spin} if
\begin{equation}
\label{eq:SpinConnection}
[\nabla_A,\rho(\alpha)] =\rho(\nabla\alpha)
\quad\text{on }C^\infty(V),
\end{equation}
for any $\alpha\in C^\infty(T^*X)$, where $\nabla$ is the Levi-Civita connection.

A Hermitian Clifford module $\fs:=(\rho,W)$ over a manifold $X$ of real dimension $2n$ is called a \emph{spin${}^c$ structure}\label{page:Spinc_structure} when $W$ has complex rank equal to $d_{2n}^\CC$, the complex rank of the unique, irreducible complex Clifford module of the complex Clifford algebra of $\RR^{2n}$ (see Lawson and Michelson \cite[Chapter I, Theorem 5.8, p. 33]{LM}).  When the dimension of $X$ is four, a \spinc structure $W$ has complex rank four.  In this case, we define  the associated first Chern class
\begin{equation}
\label{eq:DefineChernClassOfSpinc}
c_1(\fs):=c_1(W^+) \in H^2(X;\ZZ),
\end{equation}
is an integral lift of the second Stiefel-Whitney class of $X$.  If $L$ is a Hermitian line bundle over $X$, we write $\fs\otimes L$ for the \spinc structure $(\rho\otimes\id_L,W\otimes L)$, which we abbreviate as $(\rho,W\otimes L)$.

A Hermitian Clifford module $\ft:=(\rho,V)$ over a manifold of dimension $2n$ is a  \emph{spin${}^u$ structure}\label{page:Spinu_structure} when $V$ has complex rank $d_{2n}^\CC r$ and $r\geq 2$ is an integer. We refer to $r$ as the \emph{rank} of the spin${}^u$ structure.
\label{page:Rank_Of_Spinu}
In particular, a \spinu structure over a manifold of real dimension four is a Hermitian Clifford module of rank $4r$. We let 
$\fu(V) \subset \End(V)$ 
\label{page:u(V)}
denote the subbundle of skew-Hermitian  endomorphisms of $V$,
$\su(V) \subset \End(V)$ 
\label{page:su(V)}
denote the subbundle of skew-Hermitian, traceless endomorphisms of $V$, and $\fg_{\ft}\subset\su(V)$ denote the $\SO(2r-1)$ subbundle given by the span of the sections of the bundle $\su(V)$ which commute with the action of $\CCl(T^*X)$ on $V$. We obtain a splitting,
\begin{equation}
\label{eq:EndSplitting}
\su(V^+)
\cong
\rho(\Lambda^+)\oplus i\rho(\Lambda^+)\otimes_\RR\fg_{\ft}
\oplus \fg_{\ft},
\end{equation}
(see Feehan and Leness \cite[Equation (2.17), p. 67]{FL2a}) and similarly for $\su(V^-)$, where we let $\Lambda^p := \Lambda^p(T^*X)$ denote the vector bundle of $p$-forms over $X$ and, when $X$ has dimension four, recall the splitting $\Lambda^2(T^*X) = \Lambda^+(T^*X)\oplus\Lambda^-(T^*X)$ of self-dual and anti-self-dual two-forms (defined by the orientation and Riemannian metric on $X$) and abbreviate $\Lambda^\pm := \Lambda^\pm(T^*X)$.

When the bundle $V^+$ has complex rank four, there is a complex line bundle,
\begin{equation}
\label{eq:CliffordDeterminantBundle}
{\det}^{\frac{1}{2}}(V^+).
\end{equation}
over $X$  whose tensor-product square is $\det(V^+)$ see Feehan and Leness \cite[Equation (2.9), p. 66]{FL2a}). (If $E$ is a complex vector bundle of $r$ over a manifold, we let $\det E = \wedge^r E$ denote the corresponding determinant line bundle.) A \spinu structure $\ft$ thus defines characteristic classes,
\begin{equation}
\label{eq:SpinUCharacteristics}
c_1(\ft) := \frac{1}{2} c_1(V^+),
\quad
p_1(\ft) := p_1(\fg_{\ft}),
\quad
\text{and}\quad
w_2(\ft) := w_2(\fg_{\ft}).
\end{equation}
Given a \spinc bundle $W$, one has an isomorphism $V\cong W\otimes E$ of Hermitian Clifford modules, where $E$ is a smooth Hermitian vector bundle of complex rank $r$ (see \cite[Lemma 2.3, p. 64]{FL2a} for the case $r=2$); then
\begin{equation}
\label{eq:SpinAssociatedBundles}
\fg_{\ft}
=
\su(E)
\quad\text{and}\quad
{\det}^{\frac{1}{2}}(V^+)
=
\det(W^+)\otimes_\CC\det E.
\end{equation}
Henceforth, we shall assume $V=W\otimes E$ and write $\su(E)$ rather than $\fg_{\ft}$ and $\ft = (\rho,W\otimes E)$ for the \spinu structure defined by a \spinc structure $(\rho,W)$ and Hermitian vector bundle $E$.

\section{Quotient space of pairs for a spin${}^u$ structure}
\label{sec:SpinuPairsQuotientSpace}
Assume now that the base manifold $X$ has real dimension four and let $\ft=(\rho,W\otimes E)$ be a spin${}^u$ structure.  We introduce the action of the group of special unitary gauge transformations on the affine space of spin${}^u$ pairs, define and discuss the stabilizers of spin${}^u$ pairs with respect to this action, and characterize the Lie algebra of such a stabilizer in Section \ref{subsec:SU(E)ActionOnSpinuPairs}.  We define the quotient space $\sC_\ft$ of spin${}^u$ pairs and the subspaces of zero-section and split pairs in Section \ref{subsec:SubspacesOfQuotientSpaceOfSpinuPairs}. In Section \ref{subsection:Tangent_bundle_quotient_space_pairs_for_spinu_structure}, we compute the tangent space of $\sC_\ft$, define an exponential map from this tangent space to $\sC_\ft$, and define a norm on the fibers of this tangent space.
Lastly, in Section \ref{subsec:SpincQuotientSpaceInSpinuSpace}, we define a quotient space $\sC_\fs$ of spin${}^c$ pairs, discuss the embedding of $\sC_\fs$ into $\sC_\ft$, define a normal bundle and tubular neighborhood
of the image of this embedding, and construct a sphere subbundle of this normal bundle.

\subsection{Action of the group of gauge transformations on the affine space of spin${}^u$ pairs}
\label{subsec:SU(E)ActionOnSpinuPairs}
The Banach Lie group $W^{2,p}(\SU(E))$ acts smoothly on the Banach affine space $\sA(E,h)\times W^{1,p}(W^+\otimes E)$ by
\begin{multline}
\label{eq:GaugeActionOnSpinuPairs}
W^{2,p}(\SU(E))\times \left(\sA(E,h)\times W^{1,p}(W^+\otimes E)\right)\ni (u,(A,\Phi))
\\
\mapsto (u^*A,u^{-1}\Phi)
\in \sA(E,h)\times W^{1,p}(W^+\otimes E).
\end{multline}
We denote the quotient space of spin${}^u$ pairs associated to $\ft$ by
\label{page:SpinuPairs}
\begin{equation}
\label{eq:SpinUConfiguration}
\sC_\ft := \left.\left(\sA(E,h)\times W^{1,p}(W^+\otimes E)\right)\right/W^{2,p}(\SU(E)),
\end{equation}
where  $W^{2,p}(\SU(E))$ acts on $\sA(E,h)\times W^{1,p}(W^+\otimes E)$  as in \eqref{eq:GaugeActionOnSpinuPairs}. We write elements of $\sC_\ft$ as $[A,\Phi]$
\label{page:GaugeEquivClassOfSpinuPair}.

\begin{rmk}[Pushforward and pullback notation for gauge transformations on \spinu pairs]
\label{rmk:PushforwardPullbackNotation_spinu_pairs}
The actions of gauge transformations on pairs corresponding to that on connections discussed in Remark \ref{rmk:PushforwardPullbackNotation} is given by
\begin{subequations}
\label{eq:PullPushActionsOnUnitaryPairs}
\begin{align}
\label{eq:PullbackActionOnUnitaryPairs}
u^*(A,\Phi) &:= (u^*A,u^{-1}\Phi),
\\
\label{eq:PushforwardActionOnUnitaryPairs}
u_*(A,\Phi) &:= (u_*A,u\Phi).
\end{align}
\end{subequations}
Because $u^*(A,\Phi)=(u^{-1})_*(A,\Phi)$, the actions \eqref{eq:PullPushActionsOnUnitaryPairs} give the same quotient $\sC_\ft$.
\end{rmk}

We define the \emph{stabilizer group} of $(A,\Phi)\in\sA(E,h)\times W^{1,p}(W^+\otimes E)$ by
\begin{equation}
\label{eq:Define_Stabilizer_of_Spinu_Pair_in_SU(E)}
\Stab(A,\Phi):=\{u\in W^{2,p}(\SU(E)): u^*(A,\Phi)=(A,\Phi)\}.
\end{equation}
We note the following 

\begin{lem}[Lie group structure of $\Stab(A,\Phi)$]
\label{lem:LieGroupStructureOfStab(A,Phi)}
Let $\ft=(\rho,W\otimes E)$ be a spin${}^u$ structure on a closed, oriented, smooth Riemannian manifold $(X,g)$ of real dimension four. For $(A,\Phi)\in \sA(E,h)\times W^{1,p}(W^+\otimes E)$, the stabilizer $\Stab(A,\Phi)$ is a Lie group with Lie algebra given by (compare the forthcoming definition \eqref{eq:H_APhi^bullet})
\[
\bH_{A,\Phi}^0:=
\Ker\left(
d_{A,\Phi}^0: W^{2,p}(\su(E)) \to W^{1,p}(\Lambda^1(\su(E)))\oplus W^{1,p}(W^+\otimes E)
\right)
\]
where $d_{A,\Phi}^0$ is given by
\begin{equation}
  \label{eq:d_APhi^0}
  d_{A,\Phi}^0 \xi = (d_A\xi, -\xi\Phi), \quad\text{for all } \xi \in W^{2,p}(\su(E)).
\end{equation}
\end{lem}

\begin{proof}
By Hilgert and Neeb \cite[Lemma 10.1.5 (i), p. 361]{Hilgert_Neeb_structure_geometry_lie_groups}, the stabilizer $\Stab(A,\Phi)$ is a closed subgroup of $W^{2,p}(\SU(E))$ and thus a Lie group. To identify the Lie algebra of $\Stab(A,\Phi)$, we compute the differential at the identity $\id_E\in C^\infty(\SU(E))$ of the map
\begin{equation}
\label{eq:Unitary_Gauge_Group_Orbit_For_Spinu_Pair}
  W^{2,p}(\SU(E)) \ni u \mapsto u^*(A,\Phi) \in \sA(E,h) \times W^{1,p}(W^+\otimes E),
\end{equation}
where $u^*(A,\Phi)$ is defined in \eqref{eq:GaugeActionOnSpinuPairs}. If $u_t\in W^{2,p}(\SU(E))$ is a path for $t\in(-\eps,\eps)$ with $u_0=\id_E$ and $du_t/dt|_{t=0} = \xi \in W^{2,p}(\su(E))$, then
\begin{equation}
  \label{eq:Differential_of_SU(E)_GaugeAction}
  \left.\frac{d}{dt}u_t^*(A,\Phi)\right|_{t=0} = (d_A\xi, -\xi\Phi).
\end{equation}
Note that we have used the pullback action from \eqref{eq:PullbackActionOnUnitaryPairs} in computing the right-hand-side of \eqref{eq:Differential_of_SU(E)_GaugeAction}, giving the same sign as in \eqref{eq:Derivative_Of_GaugeGroupAction}. Hence, the differential of the map \eqref{eq:Unitary_Gauge_Group_Orbit_For_Spinu_Pair} at $\id_E$ is given by the map $d_{A,\Phi}^0$ in \eqref{eq:d_APhi^0} and hence the kernel of this map is the tangent space at the identity of $\Stab(A,\Phi)$ as asserted.
\end{proof}

\begin{lem}[Non-zero section, rank-two \spinu pairs have trivial stabilizer]
\label{lem:NonZeroSection_Spinu_pairs_Have_Trivial_Stabilizer}
Let $\ft=(\rho,W\otimes E)$ be a spin${}^u$ structure on a closed, connected, oriented, smooth Riemannian manifold $(X,g)$ of real dimension four. If $E$ has rank two and $(A,\Phi)\in \sA(E,h)\times W^{1,p}(W^+\otimes E)$ with $\Phi\not\equiv 0$ and $\Stab(A)\not\cong \SU(2)$, then $\Stab(A,\Phi)=\{\id_E\}$.
\end{lem}

\begin{proof}
Assume, to obtain a contradiction, that $\Stab(A,\Phi)\neq\{\id_E\}$. If $\Phi\not\equiv 0$, then $-\id_E\notin \Stab(A,\Phi)$ because $-\Phi\neq \Phi$. Hence, the facts that $\Stab(A,\Phi)\neq\{\id_E\}$ and $\Stab(A,\Phi)\subset\Stab(A)$ imply that $\Stab(A)\neq\{\pm\id_E\}$, that is, the stabilizer of $A$ is not central.  By Lemma \ref{lem:Kronheimer_Lemma2.2}, we see that $\dim\Stab(A)>0$. Because $\Stab(A)\neq\{\pm\id_E\}$ and $\Stab(A)\neq \SU(2)$, Lemma \ref{lem:Howe_Subgroups_of_SU(n)} implies that $\Stab(A)\cong S^1$ (see also Example \ref{exmp:Stab(A)_r=2}). Proposition \ref{prop:ReducibleUnitaryFromH0} and
the proof of
Corollary \ref{cor:Split_unitary_A_and_Lie_Alg_of_Stab(A)} imply that $A$ is split with respect to an
orthogonal decomposition $E=L_1\oplus L_2$ as a direct sum of Hermitian line bundles and
\[
\Stab(A) = \left\{e^{i\theta}\,\id_{L_1} \oplus e^{-i\theta}\,\id_{L_2}: e^{i\theta}\in S^1\right\}.
\]
Write $\Phi = \phi_1 + \phi_2$, where $\phi_j\in W^{1,p}(L_j)$ for $j=1,2$. Thus, $u(\theta)\Phi = e^{i\theta}\phi_1+e^{-i\theta}\phi_2$ for $u(\theta) = e^{i\theta}\,\id_{L_1} \oplus e^{-i\theta}\,\id_{L_2} \in \Stab(A)$.  Because $\phi_1$ and $\phi_2$ are linearly independent, the equality $u(\theta)\Phi=\Phi$ implies that $e^{i\theta}\phi_1 = \phi_1$ and $e^{-i\theta}\phi_2 = \phi_2$. Thus, $u(\theta)\in\Stab(A,\Phi)$ implies that $u(\theta)=\id_E$, contradicting the assumption that $\Stab(A,\Phi)\neq\{\id_E\}$.
\end{proof}

\subsection{Subspaces of the quotient space of spin${}^u$ pairs}
\label{subsec:SubspacesOfQuotientSpaceOfSpinuPairs}
We now discuss subspaces of the quotient space $\sC_\ft$.

\begin{defn}[Zero-section \spinu pairs]
\label{defn:ZeroSectionSO(3)Pair}
(See Feehan and Leness \cite[Definition 2.2, p. 64]{FL1}.)
We call a pair $(A,\Phi)\in\sA(E,h)\times W^{1,p}(W^+\otimes E)$ a \emph{zero-section} pair if $\Phi\equiv 0$.
\end{defn}

The subspace of $\sC_\ft$ given by gauge-equivalence classes of zero-section pairs is identified with the
\emph{quotient space} $\sB(E,h)$
\label{page:Quotient_space_unitary_connections}
\emph{of fixed-determinant, unitary connections} on $E$. Recall from \eqref{eq:Define_Stabilizer_of_Spinu_Pair_in_SU(E)} that the stabilizer group of $(A,\Phi)\in\sA(E,h)\times W^{1,p}(W^+\otimes E)$ is
\[
  \Stab(A,\Phi):=\{u\in W^{2,p}(\SU(E)): u^*(A,\Phi)=(A,\Phi)\}.
\]
We provide the following partial analogue of Definition \ref{defn:Reducible_split_trivial-stabilizer_unitary_connection} for \spinu pairs.

\begin{defn}[Split, trivial-stabilizer, and central-stabilizer \spinu pairs]
\label{defn:Split_trivial_central-stabilizer_spinor_pair}
(Compare Feehan and Leness \cite[Definition 2.2, p. 64]{FL2a} for the case $r=2$.)
Let $(E,h)$ be a smooth Hermitian vector bundle of complex rank $r\geq 2$ over a smooth, connected Riemannian manifold of real dimension four and let $(A,\Phi)\in\sA(E,h)\times W^{1,p}(W^+\otimes E)$.
\begin{enumerate}
\item\label{item:Trivial-stabilizer_spinor_pair} 
$(A,\Phi)$ has \emph{trivial stabilizer} if the stabilizer group $\Stab(A,\Phi)$ of $(A,\Phi)$ in $W^{2,p}(\SU(E))$ is $\{\id_E\}$ and has \emph{non-trivial stabilizer} otherwise.

\item\label{item:central-stabilizer_spinor_pair}
$(A,\Phi)$ has \emph{central stabilizer} if the stabilizer group $\Stab(A,\Phi)$ of $(A,\Phi)$ in $W^{2,p}(\SU(E))$ is equal to the center $Z(\SU(r)) = C_r$ and has \emph{non-central stabilizer} otherwise.

\item\label{item:Split_spinor_pair}
$(A,\Phi)$ is \emph{split} if the connection $A$ is split as in Definition \ref{defn:Reducible_split_trivial-stabilizer_unitary_connection} \eqref{item:Split_unitary_connection} and $\Phi \in W^{1,p}(W^+\otimes E_1)$. If no such splitting exists, then $(A,\Phi)$ is called \emph{non-split}. We say that $(A,\Phi)$ is \emph{split with respect to the decomposition $E=E_1\oplus E_2$} in \eqref{eq:BasicSplitting} when we wish to specify the bundle splitting.
\end{enumerate}
\end{defn}

\begin{rmk}[On prior definitions of reducible and irreducible pairs]
\label{rmk:Old_Definition_of_Reducible_Pairs}
In our previous articles (for example, Feehan and Leness \cite[Section 2.1.3, p. 284]{FL1} or \cite[Section 2.2, p. 72]{FL2a}), we also defined $\sC_\ft^*$ to be the open subspace of $\sC_\ft$ of points represented by pairs $(A,\Phi)$ where the connection $A$ was non-split. However, we used the term `irreducible' for a non-split pair and `reducible' for a split pair.
\end{rmk}

Zentner \cite{Zentner_2012} has described the subspaces of $\sC_\ft$ arising from different decompositions of the vector bundle $E$ when $E$ has rank $r\geq 2$. Although our focus in this monograph is on \spinu structures $(\rho,W\otimes E)$ where $E$ has $r=2$, we shall explicitly state when this assumption is used. We define subspaces of $\sC_\ft$ by
\begin{subequations}
 \label{eq:SpinuQuotientSpaceSubspaces}
 \begin{align}
   \label{eq:eq:SpinuQuotientSpaceSubspacese_trivial_stabilizer}
   \sC_\ft^{**} &:= \{[A,\Phi]\in\sC_\ft: \Stab(A,\Phi) = \{\id_E\}\},
   \\            
   \label{eq:eq:SpinuQuotientSpaceSubspacese_irreducible}
   \sC_\ft^* &:= \{[A,\Phi]\in\sC_\ft: \text{$(A,\Phi)$ is non-split}\},
   \\
   \label{eq:eq:SpinuQuotientSpaceSubspaces_non-zero-section}
   \sC_\ft^0 &:= \{[A,\Phi]\in\sC_\ft: \Phi\not\equiv 0\},
   \\
   \label{eq:eq:SpinuQuotientSpaceSubspaces_irredicible_non-zero-section}
   \sC_\ft^{*,0} &:= \sC_\ft^*\cap\sC_\ft^0.
\end{align}
\end{subequations}
The quotient space $\sC_\ft^{**}$ in \eqref{eq:eq:SpinuQuotientSpaceSubspacese_trivial_stabilizer} is a real analytic Banach manifold by the analogue for \spinu pairs of the forthcoming Theorem \ref{thm:DK_prop_5-2-9_FU_corollary_page_50_unitary_pairs} for unitary pairs (the proof for \spinu pairs is identical). 
(In Feehan and Leness \cite[Proposition 2.8, p. 287]{FL1}, we proved the weaker result that $\sC_\ft^{*,0}$ is a Banach manifold when $E$ has rank two.) By the analogue of \eqref{eq:sC**(E,h)_subset_sC0(E,h)} in the forthcoming Lemma \ref{lem:sM*0(E,h,omega)_subset_sM**(E,h,omega)} (again, the proof for \spinu pairs is identical), we have the inclusion
\[
  \sC_\ft^{**} \subset \sC_\ft^0.
\]
Recall that when $E$ has rank two and $\Phi\not\equiv 0$, Lemma \ref{lem:NonZeroSection_Spinu_pairs_Have_Trivial_Stabilizer} implies that $\Stab(A,\Phi) = \{\id_E\}$ and so $\sC_\ft^0 \subset \sC_\ft^{**}$, which yields the equality
\begin{equation}
  \label{eq:sC0(ft_equals_sC**(ft)_E_rank2}
  \sC_\ft^0 = \sC_\ft^{**},
\end{equation}
analogous to the equality \eqref{eq:sC0(E,h)_equals_sC**(E,h)_E_rank2} in Lemma \ref{lem:sM*0(E,h,omega)_subset_sM**(E,h,omega)} \eqref{item:Subspace_trivial_stabilizer_equals_non-zero-section_rank-2_unitary_pairs} for unitary pairs.
The map $\sS$ in \eqref{eq:PerturbedSO3MonopoleEquation_map} defines a real analytic Fredholm section of a real analytic vector bundle \eqref{eq:SpinUConfiguration_vector_bundle} over the real analytic Banach manifold $\sC_\ft^{**}$.  

\subsection{Tangent bundle of the quotient space of pairs for a spin${}^u$ structure}
\label{subsection:Tangent_bundle_quotient_space_pairs_for_spinu_structure}
To describe the tangent bundle of $\sC_\ft^{**}$ at the equivalence class of a pair $[A,\Phi]$, we introduce a local slice (see the forthcoming Remark \ref{rmk:Duistermaat_Kolk_2-3-1_local_slice}) for the action \eqref{eq:GaugeActionOnSpinuPairs}, that is, a subspace of the affine space
$\sA(E,h) \times \Omega^0(W^+\otimes E)$ transverse to the orbit given by the image of the map
\[
  W^{2,p}(\SU(E)) \ni u \mapsto u^*(A,\Phi) \in \sA(E,h) \times W^{1,p}(W^+\otimes E),
\]
where $u^*(A,\Phi)$ is defined in \eqref{eq:GaugeActionOnSpinuPairs}. Recall from the proof of Lemma \ref{lem:LieGroupStructureOfStab(A,Phi)} that the operator $d_{A,\Phi}^0$ defined in \eqref{eq:d_APhi^0} is the derivative at the identity of the preceding map. The $L^2$ adjoint of the operator $d_{A,\Phi}^0$  is given by
\[
 d_{A,\Phi}^{0,*}: W^{1,p}\left(T^*X\otimes \su(E) \oplus W^+\otimes E\right) \to L^p(\su(E)),
\]
for which an explicit expression will be given in the forthcoming \eqref{eq:d_APhi^0_star_aphi_identity_and_vanishing}. By Feehan and Maridakis \cite[Theorem 16, p. 18 and Corollary 18, p. 19]{Feehan_Maridakis_Lojasiewicz-Simon_coupled_Yang-Mills}, there is a constant $\zeta = \zeta(A,\Phi) \in (0,1]$ such that the ball,
\begin{multline}
  \label{eq:SliceBall}
  \bB_{A,\Phi}(\zeta)
  :=
  \left\{(a,\phi)\in \Ker d_{A,\Phi}^{0,*}\cap W^{1,p}\left(T^*X\otimes \su(E) \oplus W^+\otimes E\right):\right.
  \\
  \left.\|(a,\phi)\|_{W_A^{1,p}(X)} <\zeta\right\},
\end{multline}
embeds into the quotient $\sC_\ft^{**}$ under the map $(a,\phi)\mapsto [(A,\Phi)+(a,\phi)]$. The tangent bundle of $\sC_\ft^{**}$ is
\begin{equation}
\label{eq:QuotientPairsTangentBundle}
T\sC_\ft^{**}
:=
\left.\left\{(A,\Phi,a,\phi)\in \sA(E,h)\times W^{1,p}(W^+\otimes E)\times\Ker d_{A,\Phi}^{0,*}\right\}\right/W^{2,p}(\SU(E)),
\end{equation}
where the group of gauge transformations acts by
\begin{equation}
\label{eq:GaugeGroupActionOnTangentBundle}
\left( u, (A,\Phi,a,\phi)\right)\mapsto \left( u^*A,u^{-1}\Phi, u^{-1}au,u^{-1}\phi\right).
\end{equation}
The equality
\begin{equation}
\label{eq:GaugeEquivarianceOfTangentBundleMap}
u^*(A,\Phi)+(u^{-1}au,u^{-1}\phi)
=
u^*(A+a,\Phi+\phi)
\end{equation}
implies that the exponential map on the affine space descends to a well-defined map on the quotients,
\begin{equation}
\label{eq:TangentSpaceMap}
\exp:T\sC_\ft^{**} \ni [A,\Phi,a,\phi] \mapsto [(A,\Phi)+(a,\phi)] \in \sC_\ft^{**}.
\end{equation}
(Compare Kondracki and Rogulski \cite[Section 3.3, p. 39]{Kondracki_Rogulski_1986} and our second proof of Lemma \ref{lem:Openness_moduli_subspace_non-split_projectively_Hermitian-Einstein_connections} \eqref{item:Openness_subspace_unitary_connections_minimal_stabilizer}.)
We will use the map \eqref{eq:TangentSpaceMap} in the forthcoming Lemma \ref{lem:TubularNghOfSpincQuotientSpace}. We define a norm on the fibers of $T\sC_\ft^{**}$ by
\begin{equation}
\label{eq;W1pNormOnTangentFibers}
\| [A,\Phi,a,\phi]\|_{W^{1,p}}^p
:=
\|(a,\phi)\|_{W^{1,p}_A}^p
=
\|(a,\phi)\|_{L^p(X)}^p
+
\|\nabla_A (a,\phi)\|_{L^p(X)}^p.
\end{equation}
To see that the preceding norm is independent of the choice of representative of the equivalence class $[A,\Phi,a,\phi]\in T\sC_\ft^{**}$, observe that the equality $\nabla_{u^*A}=u^{-1}\circ\nabla_A\circ u$ implies
that for all $a\in W^{1,p}(T^*X\otimes \su(E))$ and $\phi\in W^{1,p}(W^+\otimes E)$,
\begin{align*}
\nabla_{u^*A} (u^{-1}au)
&=
u^{-1} \left(\nabla_A \left(u\left(u^{-1}au\right)u^{-1}\right)\right)u
=
u^{-1} (\nabla_A a)u,
\\
\nabla_{u^*A}(u^{-1}\phi)
&=
u^{-1}\nabla_A \phi.
\end{align*}
(Note that in the preceding computations, the gauge transformation $u$ acts on sections of $W^+\otimes E$ by $(u,\Phi)\mapsto u\Phi$ and on sections of $T^*X\otimes\su(E)$ by $(u,a)\mapsto \ad(u)a=uau^{-1}$.) Combining the preceding equalities with the fact that $u\in W^{2,p}(\SU(E))$ acts on each fiber of $\su(E)$ and $W^+\otimes E$ by isometries implies that the norm \eqref{eq;W1pNormOnTangentFibers} is well-defined.

\subsection{Embedding quotient spaces of spin${}^c$ pairs in the quotient space of spin${}^u$ pairs}
\label{subsec:SpincQuotientSpaceInSpinuSpace}
We begin by defining the quotient space of spin${}^c$ pairs on a closed, four-dimensional, oriented, smooth Riemannian manifold $(X,g)$. We depart slightly from the notation used for spin${}^c$ pairs in Feehan and Leness \cite{FL2a} in order to match that used for spin${}^u$ pairs in \eqref{eq:SpinUConfiguration}, which in turn was chosen for compatibility with the forthcoming projective vortex equations \eqref{eq:SO(3)_monopole_equations_almost_Hermitian_alpha}.

As in our definition of spin${}^u$ structures, we fix a spin${}^c$ structure $\fs_0=(\rho,W)$ and a unitary connection $A_W$ on $W$ that is a spin in the sense of \eqref{eq:SpinConnection}. Let $(L,h_L)$ be a Hermitian line bundle over $X$ and let $\fs=\fs_0\otimes L:=(\rho\otimes \id_L,W\otimes L)$ be the \spinc structure induced by $L$ and $\fs_0$. In Feehan and Leness \cite[Equation (2.52), p. 76]{FL2a}, we consider the affine space $\sA_{\fs}\times W^{1,p}(W^+\otimes L)$ of \spinc pairs $(B,\Psi)$, where $\sA_{\fs}$ denotes the space of $W^{1,p}$ unitary connections $B$ on $W\otimes L$ that are spin in the sense of \eqref{eq:SpinConnection}.

In this monograph, we instead take
\[
  \sA(L,h_L)\times W^{1,p}(W^+\otimes L)
\]
to be the affine space of spin${}^c$ pairs $(A_L,\Psi)$, where $\sA(L,h_L)$ is the affine space of unitary connections $A_L$ on $L$. While this definition depends on the choice of spin${}^c$ structure $\fs_0$ and connection $A_W$ and thus is not intrinsic to the spin${}^c$ structure $\fs_0\otimes L$, it is more convenient for our current applications.  The map
\begin{equation}
\label{eq:TranslateSWNotation}
\sA(L,h_L)
\ni
A_L
\mapsto
A_W\otimes A_L
\in
\sA_{\fs_0\otimes L}
\end{equation}
identifies the two affine spaces.

Let  $W^{2,p}(X,S^1)$ denote the Banach Lie group of $W^{2,p}$ maps from $X$ to $S^1$. We define the action of $W^{2,p}(X,S^1)$ on the affine space $\sA(L,h_L)\times W^{1,p}(W^+\otimes L)$ by analogy with \eqref{eq:GaugeActionOnSpinuPairs} as
\begin{multline}
\label{eq:SWGaugeGroupAction}
\sA(L,h_L)\times W^{1,p}(W^+\otimes L) \ni
(A_L,\Psi)
\mapsto
s\left(A_L,\Psi\right)
:=
\left(s^*A_L, s^{-1}\Psi\right)
\\
= \left(A_L + (s^{-1}ds), s^{-1}\Psi\right)
\in \sA(L,h_L)\times W^{1,p}(W^+\otimes L).
\end{multline}
The quotient space of spin${}^c$ pairs on $\fs=\fs_0\otimes L$ is
\begin{equation}
\label{eq:SpincConfig}
\sC_\fs := \left.\left(\sA(L,h_L)\times W^{1,p}(W^+\otimes L)\right)\right/W^{2,p}(X,S^1),
\end{equation}
where $W^{2,p}(X,S^1)$ acts by \eqref{eq:SWGaugeGroupAction}, and we recall that $\sC_\fs$ is Hausdorff (see Morgan \cite[Corollary 4.5.4, p. 62]{MorganSWNotes}). We let
\begin{equation}
\label{eq:SpincConfigNonZero}
\sC_\fs^0 := \left\{[A_L,\Psi]\in\sC_\fs: \Psi\not\equiv 0\right\}
\end{equation}
denote the open subspace of $\sC_\fs$ comprising gauge-equivalence classes of non-zero-section \spinc pairs, and which is a Banach manifold (see Morgan \cite[Corollary 4.5.7, p. 64]{MorganSWNotes}) since $W^{2,p}(X,S^1)$ acts freely on the open subspace of non-zero-section \spinc pairs (see Morgan \cite[Lemma 4.5.1, p. 61]{MorganSWNotes}).

We now assume that $\ft=(\rho\otimes\id_E,W\otimes E)$ is a spin${}^u$ structure of rank two and describe an embedding of $\sC_\fs^0$ into $\sC_\ft^0$. If $(A,\Phi)$ is split in the sense of Definition \ref{defn:Split_trivial_central-stabilizer_spinor_pair}, we can write $A=A_1\oplus A_2$, where $A_1$ and $A_2$ are unitary connections on $L_1$ and $L_2$, respectively, and $E = L_1\oplus L_2$ as a direct sum of Hermitian line bundles.
If $h$ is the Hermitian metric on $E$, then we write $h=h_1\oplus h_2$, where $h_i$ is the Hermitian metric on $L_i$ for $i=1,2$. Because $A$ induces the connection $A_d$ on $\det E$ by the definition of $\sA(E,h)$, we have $A_2=A_d\otimes A_1^*$ and $L_2\cong \det E\otimes L_1^*$. If $\fs=(\rho,W\otimes L_1)$, then there is a smooth embedding,
\begin{multline}
\label{eq:DefnOfIota}
\tilde\iota_{\fs,\ft}:
\sA(L_1,h_1)\times W^{1,p}(W^+\otimes L_1)
\ni (A_1,\Psi)
\\
\mapsto
\left(A_1\oplus (A_d\otimes A_1^*),\Psi\oplus 0\right) \in
\sA(E,h)\times W^{1,p}(W^+\otimes E),
\end{multline}
where $A_d$ is the connection on $\det E$ defining $\sA(E,h)$. The embedding \eqref{eq:DefnOfIota} is gauge-equivariant with respect to the homomorphism
\begin{equation}
  \label{eq:GaugeGroupInclusion}
\varrho:W^{1,p}(X,S^1) \ni s
\mapsto
s\,\id_{L_1}\oplus s^{-1}\,\id_{L_2}
\in W^{2,p}(\SU(E)).
\end{equation}
According to \cite[Lemma 3.11, p. 93 and Lemma 3.16, p. 105]{FL2a}, the map \eqref{eq:DefnOfIota} descends to the quotient,
\begin{equation}
\label{eq:DefnOfIotaOnQuotient}
\iota_{\fs,\ft}:\sC_\fs\to \sC_\ft,
\end{equation}
and its restriction to $\sC_\fs^0$ is a smooth embedding $\sC_\fs^0\embed\sC_\ft^0$.

We now introduce a tubular neighborhood of $\iota_{\fs,\ft}(\sC_\fs^0)$ in $\sC_\ft^0$. To define a normal bundle of $\iota_{\fs,\ft}(\sC_\fs^0)$, we begin by reviewing a splitting of the operator $d_{A,\Phi}^{*,0}$ appearing in the definition of $T\sC_\ft^0$ in \eqref{eq:QuotientPairsTangentBundle}.  We write
\[
d_{A,\Phi}^{0,*}: W^{1,p}(F_1)\to L^p(F_0),
\]
where
\[
F_0:= \su(E)
\quad\text{and}\quad
F_1:= T^*X\otimes\su(E) \oplus W^+\otimes E.
\]
The vector bundles $F_0$ and $F_1$ admit orthogonal splittings as 
\begin{equation}
\label{eq:SplittingOfSliceDeformationTerms}
F_j\cong F_j^t\oplus F_j^n,
\end{equation}
where
\[
\begin{array}{lll}
F_0^t:= i\La^0 & \text{and} & F_0^n:=L_2\otimes L_1^*,
\\
F_1^t:= i T^*X\oplus W^+\otimes L_1 & \text{and} & F_1^n:=(T^*X\otimes_\RR L_2\otimes L_1^*) \oplus W^+\otimes L_2.
\end{array}
\]
If $(A,\Phi)=\tilde\iota_{\fs,\ft}(A_1,\Psi)$,  then by \cite[Equations (3.31), (3.32), and (3.33), p. 102]{FL2a}, the operator $d_{A,\Phi}^{0,*}$ admits a splitting with respect to the decomposition \eqref{eq:SplittingOfSliceDeformationTerms},
\[
d_{A,\Phi}^{0,*}=d_{A,\Phi}^{0,t,*}\oplus d_{A,\Phi}^{0,n,*},
\]
where
\[
d^{0,t,*}_{A,\Phi}: W^{1,p}(F_1^t) \to L^p(F_0^t),
\]
is defined in \cite[Equation (3.32), p. 102]{FL2a} and
\[
d_{A,\Phi}^{0,n,*}: W^{1,p}(F_1^n)\to L^p(F_0^n),
\]
is defined in \cite[Equation (3.33), p. 102]{FL2a}.  The direct sum decomposition of $d_{A,\Phi}^{0,*}$ implies an $L^2$-orthogonal decomposition of its kernel,
\begin{equation}
\label{eq:SliceDecompositionAtReducible}
\Ker d_{A,\Phi}^{0,*}
=
\Ker d_{A,\Phi}^{0,t,*} \oplus \Ker d_{A,\Phi}^{0,n,*}.
\end{equation}
As noted following \cite[Equation (3.29), p. 101]{FL2a}, if $(A,\Phi)=\tilde\iota_{\fs,\ft}(A_1,\Psi)$
and if the pair $(B,\Psi)\in\sA_\fs\times W^{1,p}(W^+\otimes L_1)$ is given by $(A_1,\Psi)$ under the correspondence \eqref{eq:TranslateSWNotation}, then the operator $d_{A,\Phi}^{0,t,*}$ equals the operator $d_{B,\Psi}^{0,*}$ given by in \cite[Equation (2.59), p. 77]{FL2a} as the $L^2$-adjoint of the linearization of the $W^{2,p}(X,S^1)$ action in \eqref{eq:SWGaugeGroupAction}. Hence, the tangent bundle of $\sC_\fs^0$ is given by
\begin{multline}
  \label{eq:TangentBundleOfSpincQuotientSpace}
  T\sC_\fs^0
  =
  \left\{(A_1,\Psi,(a_t,\phi_1))\in\sA(L_1,h_1)\times W^{1,p}(W^+\otimes L_1)\times W^{1,p}(F_1^t):\right.
  \\
  \left.\left.d_{A,\Phi}^{0,t,*}(a_t,\phi_1)=0, \quad\text{for } (A,\Phi)=\tilde\iota_{\fs,\ft}(A_1,\Psi)
    \right\}\right/W^{2,p}(X,S^1),
\end{multline}
where $(a_t,\phi_1)\in W^{1,p}(F_1^t)$ and $u \in W^{2,p}(X,S^1)$ acts by
\begin{equation}
\label{eq:GaugeGroupActionOnTangentSpaceOfSpincQuotient}
u\cdot\left(A_1,\Psi,(a_t,\phi_1)\right)
=
\left(u^*A_1,u^{-1}\Phi, (a_t,u^{-1}\phi_1)\right).
\end{equation}
We define a vector bundle over $\iota_{\fs,\ft}(\sC_\fs^0)$ by
\begin{multline}
  \label{eq:NormalBundleOfSpincQuotientInSpinuQuotient}
  \fN(\ft,\fs)
  :=
  \left\{ (A_1,\Psi,(a_n,\phi_2))\in \sA(L_1,h_1)\times W^{1,p}(W^+\otimes L_1)\times W^{1,p}(F_1^n):\right.
  \\
  \left.\left. d^{0,n,*}_{A,\Phi}(a_n,\phi_2)=0, \quad\text{for } (A,\Phi)=\tilde\iota_{\fs,\ft}(A_1,\Psi)
    \right\}\right/W^{2,p}(X,S^1),
\end{multline}
where $(a_n,\phi_2)\in W^{1,p}(F_1^n)$ and $u\in W^{2,p}(X,S^1)$ acts by
\begin{equation}
\label{eq:GaugeGroupActionOnNormalSpaceOfSpincQuotient}
u\cdot\left(A_1,\Psi,(a_n,\phi_2)\right)
=
\left(u^*A_1,u^{-1}\Phi,(u^{-2} a_n,u^{-1}\phi_2)\right).
\end{equation}
The actions of $W^{2,p}(X,S^1)$ in \eqref{eq:GaugeGroupActionOnTangentSpaceOfSpincQuotient} and \eqref{eq:GaugeGroupActionOnNormalSpaceOfSpincQuotient} are given following \cite[Equation (3.40), p. 106]{FL2a}, taking into account the difference in notation described in the forthcoming Remark \ref{rmk:ComparingLineBundlesFromFl2a}.

The following lemma asserts that $\fN(\ft,\fs)$ is a normal bundle for the embedded submanifold $\sC_\fs^0$ in $\sC_\ft^0$ in the sense that its fibers are $L^2$-orthogonal complements of the fibers of $T\sC_\fs^0$ in the restriction of $T\sC_\ft^0$ to $\iota_{\fs,\ft}(\sC_\fs^0)$.

\begin{lem}[Normal bundles for submanifolds of gauge-equivalence classes of split pairs in the configuration space]
\label{lem:DirectSumDecompositionOfTangentBundle}
Let $E$ be a rank-two, smooth Hermitian vector bundle over a closed, oriented, smooth Riemannian four-manifold $(X,g)$.
Let $\fs=(\rho,W\otimes L_1)$ and $\ft=(\rho,W\otimes E)$ be a spin${}^c$ and spin${}^u$ structure on $X$. Assume that $E=L_1\oplus L_2$ as an orthogonal direct sum of Hermitian line bundles. Then the restriction of $T\sC_\ft^0$ to $\iota_{\fs,\ft}(\sC_\fs^0)$ admits a direct-sum decomposition,
\begin{equation}
\label{eq:DirectSumDecompositionOfTangentBundle}
T\sC_\ft^0|_{\iota_{\fs,\ft}(\sC_\fs^0)}
\cong
T\sC_\fs^0 \oplus \fN(\ft,\fs),
\end{equation}
which is orthogonal with respect to the $L^2$ metric on the fibers of $T\sC_\ft^0$, where $T\sC_\fs^0$ is the tangent bundle of $T\sC_\fs^0$ defined in \eqref{eq:TangentBundleOfSpincQuotientSpace} and $\fN(\ft,\fs)$ is as in \eqref{eq:NormalBundleOfSpincQuotientInSpinuQuotient}.
\end{lem}

\begin{proof}
The hypothesis that $E=L_1\oplus L_2$ implies that the embedding $\iota_{\fs,\ft}:\sC_\fs^0\to\sC_\ft^0$ exists.
The restriction of $T\sC_\ft^0$ to $\iota_{\fs,\ft}(\sC_\fs^0)$ admits a reduction of its structure group to $W^{2,p}(X,S^1)$ in the sense that
\begin{multline}
  \label{eq:RestrictionOfSpinuQuotientTangentBundle}
  T\sC_\ft^0|_{\iota_{\fs,\ft}(\sC_\fs^0)}
  =
  \left\{(\tilde\iota_{\fs,\ft}(A_1,\Psi),(a,\phi))
    \in \sA(E,h)\times W^{1,p}(W^+\otimes E)\times W^{1,p}(F_1): \right.
  \\
  \left.\left.
      d^{0,*}_{\tilde\iota_{\fs,\ft}(A_1,\Psi)}(a,\phi)=0\right\}\right/W^{2,p}(X,S^1),
\end{multline}
where $(a,\phi)\in W^{1,p}(F_1)$ and the action of $W^{2,p}(X,S^1)$ is given by the homomorphism \eqref{eq:GaugeGroupInclusion} from $W^{2,p}(X,S^1)$ into $W^{2,p}(\SU(E))$ and the action of $W^{2,p}(\SU(E))$ in \eqref{eq:GaugeGroupActionOnTangentBundle}. The decomposition of the restriction of $T\sC_\ft^0$ then follows from the decomposition of $\Ker d_{A,\Phi}^{0,*}$ in \eqref{eq:SliceDecompositionAtReducible}, the expression for $T\sC_\ft^0$ in \eqref{eq:TangentBundleOfSpincQuotientSpace}, and the definition of $\fN(\ft,\fs)$ in \eqref{eq:NormalBundleOfSpincQuotientInSpinuQuotient}.  That the decomposition is orthogonal follows 
from the orthogonality of the decomposition \eqref{eq:SliceDecompositionAtReducible}.
\end{proof}

While the following construction of a tubular neighborhood is fairly standard (see Lang \cite[Chapter IV, Theorem 5.1, p. 99]{Lang_introduction_differential_topology}), the explicit form of the diffeomorphism $\exp$  will be useful to us.
We first define the following disk and sphere subbundles of $\fN(\fs,\ft)$.
\begin{subequations}
  \label{eq:DefineDiskAndSphereNormalBundle}
\begin{align}
  \label{eq:DefineDiskNormalBundle}
\fN^{\le\eps}(\ft,\fs)
&:=
\left\{[A_1,\Psi,a_n,\phi_2]\in \fN(\ft,\fs):  \| [A_1,\Psi,a_n,\phi_2]\|_{W^{1,p}(X)}\le\eps\right\},
\\
  \label{eq:DefineSphereNormalBundle}
\fN^\eps(\ft,\fs)
&:=
\left\{[A_1,\Psi,a_n,\phi_2]\in \fN(\ft,\fs):  \| [A_1,\Psi,a_n,\phi_2]\|_{W^{1,p}(X)}=\eps\right\},
\end{align}
\end{subequations}
where $\|\cdot\|_{W^{1,p}}$ denotes the fiber norm defined in \eqref{eq;W1pNormOnTangentFibers}.

\begin{lem}[Tubular neighborhoods for submanifolds of gauge-equivalence classes of split pairs in the configuration space]
\label{lem:TubularNghOfSpincQuotientSpace}
Continue the hypotheses of Lemma \ref{lem:DirectSumDecompositionOfTangentBundle}. Then there are an open neighborhood $\sO(\ft,\fs)$ of the zero-section in the bundle $\fN(\ft,\fs)$ in \eqref{eq:NormalBundleOfSpincQuotientInSpinuQuotient} and an open neighborhood $\sU(\ft,\fs)$ of $\iota_{\fs,\ft}(\sC_\fs^0)$ in $\sC_\ft^0$ such that the restriction of the map $\exp$ in \eqref{eq:TangentSpaceMap} to $\sO(\ft,\fs)$ is a diffeomorphism,
\[
\exp: \fN(\ft,\fs) \supset \sO(\ft,\fs)\to \sU(\ft,\fs)\subset\sC_\ft^0.
\]
In addition, for every $[A,\Phi]\in \iota_{\fs,\ft}(\sC_\fs^0)$ there is a constant $\eps=\eps([A,\Phi])\in (0,1]$ such that
\begin{equation}
\label{eq:DiskNormalBundleInclusion}
\fN^{\le\eps}(\ft,\fs)|_{[A,\Phi]} \subset \sO(\ft,\fs)
\end{equation}
where $\fN^{\le\eps}(\ft,\fs)\subset \fN(\ft,\fs)$ is the disk bundle in \eqref{eq:DefineDiskNormalBundle}.
\end{lem}

\begin{proof}
If $(A,\Phi)=\tilde\iota_{\fs,\ft}(A_1,\Psi_1)$, then the tangent space to $\fN(\ft,\fs)$ at a point $[A_1,\Psi,0,0]$ in the zero section of $\fN(\ft,\fs)$ is given by
\begin{align*}
T_{[A_1,\Psi,0,0]}\fN(\ft,\fs)
&\cong
T_{[A_1,\Psi]}\sC_\fs^0 \oplus \fN(\ft,\fs)|_{[A_1,\Psi]}
\\
&\cong
     \Ker d_{A,\Phi}^{0,t,*} \oplus \Ker d_{A,\Phi}^{0,t,*}.
\end{align*}
The tangent space of $\sC_\ft^0$ at $\exp([A_1,\Psi,0,0])=[\tilde\iota_{\fs,\ft}(A_1,\Psi)]=[A,\Phi]$ is
\[
T_{[\tilde\iota_{\fs,\ft}(A_1,\Psi)]}\sC_\ft^0
\cong
\Ker d^{0,*}_{A,\Phi}.
\]
The derivative of $\exp$ at $[A_1,\Psi,0,0]$ is the isomorphism \eqref{eq:SliceDecompositionAtReducible},
\[
\Ker d_{A,\Phi}^{0,t,*} \oplus \Ker d_{A,\Phi}^{0,n,*}
\to
\Ker d^{0,*}_{A,\Phi}.
\]
The existence of the neighborhood $\sO(\ft,\fs)$ on which $\exp$ is a diffeomorphism then follows from the Inverse Mapping Theorem (see Abraham, Marsden, and Ratiu \cite[Theorem 2.5.2, p. 116]{AMR}).

The existence of $\eps = \eps([A,\Phi]) \in (0,1]$ such that the inclusion \eqref{eq:DiskNormalBundleInclusion} holds follows from the openness of $\sO(\ft,\fs)$ and the fact that the topology on the quotient $\sC_\ft$ is defined by
the $W^{1,p}$ topology on the ball \eqref{eq:SliceBall}.
\end{proof}

We now show that given a compact subspace $V\subset\sC_\fs^0$, we can find a constant $\eps = \eps(p,\ft,\fs,V) \in (0,1]$ such that the restriction of the disk bundle $\fN^{\le\eps}(\ft,\fs)$ defined in \eqref{eq:DefineDiskNormalBundle} to $V$ is contained in the neighborhood $\sO(\ft,\fs)$ constructed in Lemma \ref{lem:TubularNghOfSpincQuotientSpace}.  While it is obvious that there are disk bundles contained in $\sO(\ft,\fs)$, we will need to use the specific disk bundle defined in \eqref{eq:DefineDiskNormalBundle} by a value of the fiber norm $\|\cdot\|_{W^{1,p}}$ from \eqref{eq;W1pNormOnTangentFibers} which is independent of the point in $V$. For a subset $V\subset\sC_\fs^0$, we define
\begin{subequations}
  \label{eq:DefineRestrictedDiskAndSphereNormalBundle}
\begin{align}
\label{eq:DefineNormalBundleRestriction}
\fN(\ft,\fs,V)
&:=
\left\{[A_1',\Psi',a_n,\phi_2]\in \fN(\ft,\fs): [A_1',\Psi']\in V\right\},
  \\
  \label{eq:DefineRestrictedDiskNormalBundle}
\fN^{\le\eps}(\ft,\fs,V)
&:=
\fN^{\le\eps}(\ft,\fs)\cap \fN(\ft,\fs,V),
\\
  \label{eq:DefineRestrictedSphereNormalBundle}
\fN^\eps(\ft,\fs,V)
&:=
\fN^\eps(\ft,\fs)\cap \fN(\ft,\fs,V),
\end{align}
\end{subequations}
where $\fN^{\le\eps}(\ft,\fs)$ and $\fN^\eps(\ft,\fs)$ are the disk and sphere bundles defined in
\eqref{eq:DefineDiskNormalBundle} and \eqref{eq:DefineSphereNormalBundle}, respectively. (Note that if $V$ were a closed submanifold of a finite-dimensional, smooth Riemannian manifold $(M,g)$, then the existence of such a constant $\eps$ would follow from the standard construction of a tubular neighborhood --- see, for example, Lee \cite[Theorem 6.24, p. 139]{Lee_john_smooth_manifolds}.)

\begin{lem}[Sphere bundles for submanifolds of gauge-equivalence classes of split pairs in configuration space]
\label{lem:FixedNormsInTubularNgh}
Continue the hypotheses of Lemma \ref{lem:DirectSumDecompositionOfTangentBundle}. If $V\subset\sC_\fs^0$ is a compact subspace, then there is a constant $\eps = \eps(p,\ft,\fs,V) \in (0,1]$ such that
\begin{equation}
\label{eq:FixedNormsInTubularNgh}
\fN^{\le\eps}(\ft,\fs,V)\subset\sO(\ft,\fs),
\end{equation}
where $\fN^{\le\eps}(\ft,\fs,V)$ is the disk bundle over $V$ defined in \eqref{eq:DefineRestrictedDiskNormalBundle}.
\end{lem}

\begin{proof}
If the inclusion \eqref{eq:FixedNormsInTubularNgh} does not hold, then for every\footnote{We adopt the convention that $\NN$ excludes $0$.} $k\in\NN$ the sphere bundle $\fN^{1/k}(\ft,\fs,V)$ is not contained in $\sO(\ft,\fs)$. Thus, for each $k$, there is a point
\[
[A_1(k),\Psi(k),a_n(k),\phi_2(k)]\in \fN^{\le 1/k}(\ft,\fs,V)\less\sO(\ft,\fs)
\]
with $[A_1(k),\Psi(k)]\in V$ and
\[
\| [A_1(k),\Psi(k),a_n(k),\phi_2(k)]\|_{W^{1,p}(X)} \le 1/k,
\]
where $\|\cdot\|_{W^{1,p}(X)}$ denotes the fiber norm defined in \eqref{eq;W1pNormOnTangentFibers}.
Thus,
\begin{equation}
\label{eq:FiberNormLimit}
\lim_{k\to\infty}\| [A_1(k),\Psi(k),a_n(k),\phi_2(k)]\|_{W^{1,p}(X)}=0.
\end{equation}
Since $V$ is a compact subset of $\sC_\fs^0$ and $\{[A_1(k),\Psi(k)]\}_{k\in\NN}\subset V$, there exists a subsequence of
$\{[A_1(k),\Psi(k)]\}_{k\in\NN}$ which converges to a limit $[A_1(\infty),\Psi(\infty)]\in V$. After passing to the subsequence and relabelling, the convergence of $\{[A_1(k),\Psi(k)]\}_{k\in\NN}$ to $[A_1(\infty),\Psi(\infty)]$ and the limit \eqref{eq:FiberNormLimit} imply that the sequence $\{[A_1(k),\Psi(k),a_n(k),\phi_2(k)]\}_{k\in\NN}$ converges to
$[A_1(\infty),\Psi(\infty),0,0]$. However, the complement $\fN(\ft,\fs)\less\sO(\ft,\fs)$ is a closed subspace of $\fN(\ft,\fs)$ containing the sequence $\{[A_1(k),\Psi(k),a_n(k),\phi_2(k)]\}_{k\in\NN}$, so the limit point $[A_1(\infty),\Psi(\infty),0,0]$ must be in $\fN(\ft,\fs)\less\sO(\ft,\fs)$, a contradiction to the inclusion \eqref{eq:DiskNormalBundleInclusion}. Hence, there exists a positive integer $k$ such that for $\eps=1/k$, the inclusion \eqref{eq:FixedNormsInTubularNgh} holds.
\end{proof}

\section{Non-Abelian monopoles}
\label{sec:PU2Monopoles}
We now define the moduli space of non-Abelian monopoles. Fix a spin connection $A_W$ on $W$. We call a pair $(A,\Phi)$ in $\sA(E,h)\times W^{1,p}(W^+\otimes E)$ a solution to the \emph{perturbed non-Abelian monopole equations}\label{SO(3)_monopole} if \cite[Equation (2.15), p. 301]{FL1}, \cite[Equation (2.32), p. 71]{FL2a}
\begin{subequations}
  \label{eq:PerturbedSO3MonopoleEquations}
\begin{align}
  \label{eq:PerturbedSO3MonopoleEquations_curvature}
  (F_A^+)_0 - \tau\rho^{-1}(\Phi\otimes\Phi^*)_{00} &= 0,
  \\
  \label{eq:PerturbedSO3MonopoleEquations_Dirac}
  D_A\Phi + \rho(\vartheta)\Phi &= 0.
\end{align}
\end{subequations}
In \cite{Zentner_2012}, Zentner has studied these equations when $\rank_\CC E>2$.  We will explicitly state when a result only applies in the case $\rank_\CC E=2$.

In \eqref{eq:PerturbedSO3MonopoleEquations}, the operator $D_A=\rho\circ \cov_A:C^\infty(W^+\otimes E)\to C^\infty(W^-\otimes E)$ is the Dirac operator\label{Dirac_operator} associated to the connection $A_W\otimes A$ on $W\otimes E$; the section $\tau$ of $\GL(\Lambda^+)$ is a perturbation close to the identity; the perturbation $\vartheta$ is a complex one-form close to zero; $\Phi^* \in \Hom(W^+\otimes E,\CC)$ is the pointwise Hermitian dual $\langle\cdot,\Phi\rangle_{W^+\otimes E}$ of $\Phi$; and $(\Phi\otimes\Phi^*)_{00}$ is the component of the section $\Phi\otimes\Phi^*$ of $i\fu(W^+\otimes E)$ lying in $\rho(\Lambda^+)\otimes\su(E)$ with respect to the splitting $\fu(W^-\otimes E)=i\underline{\RR}\oplus\su(W^-\otimes E)$ and decomposition \eqref{eq:EndSplitting} of $\su(W^-\otimes E)$, where $\underline{\RR} := X\times\RR$.

For Sobolev exponent $p \in (2,\infty)$, we observe that the perturbed non-Abelian monopole equations \eqref{eq:PerturbedSO3MonopoleEquations} define a smooth (in fact, analytic) map
\begin{multline}
\label{eq:PerturbedSO3MonopoleEquation_map}
  \sS:\sA(E,h)\times W^{1,p}(W^+\otimes E) \ni (A,\Phi)
  \\
  \mapsto
  \left((F_A^+)_0 - \tau\rho^{-1}(\Phi\otimes\Phi^*)_{00}, D_A\Phi + \rho(\vartheta)\Phi \right)
  \in L^p\left(\Lambda^+\otimes\su(E) \oplus W^-\otimes E\right)
\end{multline}
that is equivariant with respect to the action on the domain and codomain induced by the action of the Banach Lie group $W^{2,p}(\SU(E))$ of $W^{2,p}$ determinant-one, unitary automorphisms of $E$. The map $\sS$ in \eqref{eq:PerturbedSO3MonopoleEquation_map} becomes Fredholm upon restriction to a Coulomb-gauge slice, as is clear from Section \ref{sec:Elliptic_deformation_complex_for_SO3_monopole_equations}.

Define the \emph{moduli space of non-Abelian monopoles} by
\begin{equation}
\label{eq:Moduli_space_SO(3)_monopoles}
\sM_\ft := \left\{[A,\Phi] \in \sC_\ft: \sS(A,\Phi) = 0 \right\},
\end{equation}
where $\sS$ is as in \eqref{eq:PerturbedSO3MonopoleEquation_map}. We observe that $\sM_\ft$ may thus be viewed as the zero locus of the induced section $\sS$ from $\sC_\ft$ to the quotient space
\begin{equation}
\label{eq:SpinUConfiguration_vector_bundle}
  \left(\sA(E,h)\times W^{1,p}(W^+\otimes E)\right)\times_{W^{2,p}(\SU(E))}
      \left(L^p\left(\Lambda^+\otimes\su(E) \oplus W^-\otimes E\right)\right).
\end{equation}
Following the definitions of subspaces of $\sC_\ft$ in \eqref{eq:SpinuQuotientSpaceSubspaces}, we set
\begin{subequations}
 \label{eq:PU2MonopoleSubspaces}
 \begin{align}
   \label{eq:PU2MonopoleSubspace_trivial_stabilizer}
   \sM^{**}_{\ft} &:= \sM_{\ft}\cap\sC^{**}_\ft,
   \\
   \label{eq:PU2MonopoleSubspace_irreducible}
   \sM^*_{\ft} &:= \sM_{\ft}\cap\sC^*_\ft,
   \\
   \label{eq:PU2MonopoleSubspace_non-zero-section}
   \sM^0_{\ft} &:= \sM_{\ft}\cap\sC^0_\ft,
   \\
   \label{eq:PU2MonopoleSubspace_irredicible_non-zero-section}
   \sM^{*,0}_\ft &:= \sM_{\ft}\cap\sC_\ft^{*,0},
\end{align}
\end{subequations}
where the subspaces $\sC^{**}_\ft$, $\sC^*_\ft$, $\sC^0_\ft$, and $\sC_\ft^{*,0}$ are defined in \eqref{eq:SpinuQuotientSpaceSubspaces}. We recall the

\begin{prop}[Expected dimension of the moduli space of non-Abelian monopoles]
\label{prop:Expected_dimension_moduli_space_non-Abelian_monopoles}
(See Feehan and Leness \cite[Proposition 2.28, p. 313]{FL1}.)
Let $\ft=(\rho,W\otimes E)$, where $(E,h)$ is a rank-two, smooth Hermitian vector bundle and $(\rho,W)$ is a spin${}^c$ structure over a closed, four-dimensional, oriented, smooth Riemannian manifold $(X,g)$. Then the expected dimension of the moduli space \eqref{eq:Moduli_space_SO(3)_monopoles} of non-Abelian monopoles is
\begin{align}
\label{eq:Transv}
\expdim \sM_{\ft} &\,= d_a(\ft)+2n_a(\ft), \quad\text{where }
\\
\notag
d_a(\ft)
&:= -2p_1(\ft)- \frac{3}{2}(e(X)+\sigma(X)) \text{ and }
n_a(\ft) := \frac{1}{4}(p_1(\ft)+c_1(\ft)^2-\sigma(X)),
\end{align}
and $e(X)$ is the Euler characteristic and $\si(X)$ is the signature of $X$, and is independent of the choice of parameters $(g,\tau,\vartheta)$ in the non-Abelian monopole equations \eqref{eq:PerturbedSO3MonopoleEquations}.
\end{prop}

\begin{thm}[Transversality for the moduli space of non-Abelian monopoles with generic parameters]
\label{thm:Transv}
(See Feehan \cite[Theorem 1.1, p. 910]{FeehanGenericMetric} and Teleman \cite[Theorem 3.19, p. 413]{TelemanGenericMetric}.)
Let $\ft=(\rho,W\otimes E)$, where $(E,h)$ is a rank-two, smooth Hermitian vector bundle and $(\rho,W)$ is a spin${}^c$ structure over a closed, four-dimensional, oriented, smooth Riemannian manifold $(X,g)$. If the parameters $(g,\tau,\vartheta)$ appearing in \eqref{eq:PerturbedSO3MonopoleEquations} are generic in the sense of \cite{FeehanGenericMetric}, then $\sM^{*,0}_{\ft}$ is an embedded smooth submanifold of the Banach manifold $\sC^{*,0}_{\ft}$ that has finite dimension equal to the expected dimension \eqref{eq:Transv}.
\end{thm}

While $\sM_\ft$ is not in general compact, it admits the following Uhlenbeck compactification. For each non-negative integer $\ell$, let $E_\ell$ be a smooth Hermitian vector bundle and $\ft(\ell)=(\rho,W\otimes E_{\ell})$ be a \spinu structure over $X$ characterized by
\begin{subequations}
\label{eq:DefineLowerChargeSpinuStr}
\begin{align}
  \label{eq:DefineLowerChargeSpinuStr_c1}
  c_1(E_\ell) &= c_1(E),
  \\
  \label{eq:DefineLowerChargeSpinuStr_p1}
  p_1(\su(E_\ell)) &= p_1(\su(E))+4\ell,
  \\
  \label{eq:DefineLowerChargeSpinuStr_w2}
  w_2(\su(E_\ell))) &= w_2(\su(E)),
\end{align}
\end{subequations}
where $\ft=(\rho,W\otimes E)$ is a given \spinu structure; when $\ell=0$, then $\ft(0)=\ft$. Note that \eqref{eq:DefineLowerChargeSpinuStr_c1} is equivalent to
\begin{equation}
\label{eq:DefineLowerChargeSpinuStr_c1'}
c_1(\ft(\ell))=c_1(\ft).
\end{equation}
We recall the

\begin{defn}[Uhlenbeck convergence]
\label{defn:UhlenbeckConvergence}
(See Feehan and Leness \cite[Definition 4.19, p. 350]{FL1}.)
For $S_\ell$ the symmetric group on $\ell$ elements,
let $\Sym^\ell(X):=X^\ell/S_\ell$ be the $\ell$-th symmetric product of $X$.
A sequence $\{[A_\alpha,\Phi_\alpha]\}_{\alpha\in\NN} \subset\sC_{\ft}$ is said to \emph{converge in the sense of Uhlenbeck} to a pair $[A_0,\Phi_0,\bx]\in\sC_{\ft(\ell)}\times\Sym^\ell(X)$ if
\begin{itemize}
\item
There is a sequence of $W_\loc^{2,p}$ Hermitian bundle isomorphisms $u_\alpha: E|_{X\less\bx} \to E_\ell|_{X\less\bx}$ such that the sequence $u_\alpha(A_\alpha,\Phi_\alpha)$ converges as $\alpha\to\infty$ to $(A_0,\Phi_0)$ in $W_\loc^{2,p}$ over $X\less\bx$, and
\item
The sequence of measures, $|F_{A_\alpha}|^2\,d\vol_g$ converges as $\alpha\to\infty$ in the weak-star topology to the measure $|F_A|^2\,d\vol_g + 8\pi^2\sum_{x\in \bx} \delta_x$, where $\delta_x$ is the Dirac
measure centered at $x \in X$.
\end{itemize}
Open neighborhoods in the topology defined by this convergence are \emph{Uhlenbeck neighborhoods}.
\end{defn}

We let $\bar\sM_{\ft}$ denote the closure of $\sM_{\ft}$ in the space of \emph{ideal non-Abelian monopoles}\label{ideal_SO(3)_monopoles},
\begin{equation}
\label{eq:idealmonopoles}
\sI\!\!\sM_{\ft} := \bigsqcup_{\ell=0}^\infty\left(\sM_{\ft(\ell)}\times\Sym^\ell(X)\right),
\end{equation}
with respect to the \emph{Uhlenbeck topology} defined by the preceding notion of convergence. We call the intersection of $\bar\sM_{\ft}$ with $\sM_{\ft(\ell)}\times \Sym^\ell(X)$ its \emph{$\ell$-th level}.

\begin{thm}[Uhlenbeck compactness of the moduli space of non-Abelian monopoles]
\label{thm:Compactness}
(See Feehan and Leness \cite[Theorem 1.1, p. 270]{FL1}, Zentner  \cite[Theorem 1.7, p. 235]{Zentner_2012}.)
Let $\ft=(\rho,W\otimes E)$, where $(E,h)$ is a Hermitian vector bundle and $(\rho,W)$ is a spin${}^c$ structure over a closed, four-dimensional, oriented, smooth Riemannian manifold $(X,g)$. Then there is a positive integer $N$, depending at most on the scalar curvature of $(X,g)$, the curvature of the fixed unitary connection $A_d$ on $\det E$, and $p_1(\ft)$, such that the Uhlenbeck closure $\bar\sM_{\ft}$ of $\sM_{\ft}$ in $\sqcup_{\ell=0}^N(\sM_{\ft(\ell)}\times\Sym^\ell(X))$ is a second-countable, compact, Hausdorff space. The space $\bar\sM_{\ft}$ carries a continuous circle action that restricts to the circle action on $\sM_{\ft_\ell}$ for each level.
\end{thm}

We shall need an enhancement of Theorem \ref{thm:Compactness} (more precisely, our sequential compactness result Feehan and Leness \cite[Theorem 4.20, p. 352]{FL1}) that allows for a sequence of non-Abelian monopoles with respect to a sequence of geometric parameters $(g,\tau,\vartheta)$, along the lines of \cite[Theorem 6.1.1, p. 103]{MMR} stated without proof by Morgan, Mrowka, and Ruberman for sequences of anti-self-dual connections and sequences of Riemannian metrics.

\begin{thm}[Uhlenbeck convergence of a sequence of non-Abelian monopoles with respect to a convergent sequence of parameters]
\label{thm:Uhlenbeck_convergence_sequence_non-Abelian_monopoles_parameters}
Continue the hypotheses of Theorem \ref{thm:Compactness}. Let $\{(g_\alpha,\tau_\alpha,\vartheta_\alpha)\}_{\alpha\in\NN}$ be a sequence of parameters defining the non-Abelian monopole equations \eqref{eq:PerturbedSO3MonopoleEquations} that converges in $C^\infty$ to a limit $(g_\infty,\tau_\infty,\vartheta_\infty)$. If $\{(A_\alpha,\Phi_\alpha)\}_{\alpha\in\NN}$ is a sequence of solutions to the corresponding equations \eqref{eq:PerturbedSO3MonopoleEquations} then, after passing to a subsequence and relabelling, this sequence of solutions converges in the sense of Uhlenbeck to a generalized $W^{1,p}$ solution $(A_\infty,\Phi_\infty,\bx)$ to \eqref{eq:PerturbedSO3MonopoleEquations} for the parameters $(g_\infty,\tau_\infty,\vartheta_\infty)$ and a spin${}^u$ structure $\ft(\ell) = (\rho,W\otimes E_\ell)$ for some integer $0 \leq \ell \leq N$, where $E_\ell$ is a smooth Hermitian vector bundle over $X$ related to $E$ by \eqref{eq:DefineLowerChargeSpinuStr}, and $\bx \in \Sym^\ell(X)$.
\end{thm}

\begin{proof}
  Theorem \ref{thm:Compactness} is proved in \cite[Section 4]{FL1}. That proof is quite lengthy, so we shall just outline the main changes here, which are all straightforward modifications of the arguments in \cite{FL1}. The proof of Theorem \ref{thm:Compactness} in \cite{FL1} simplifies when the holonomy perturbations allowed there are replaced by the much simpler geometric perturbations assumed in  \eqref{eq:PerturbedSO3MonopoleEquations}. The \apriori estimates in \cite[Lemma 4.2 p. 336, Lemma 4.3, p. 337, and Lemma 4.4, p. 337]{FL1} hold with the same constants for parameters $(g,\tau,\vartheta)$ in a small $C^r$-open neighborhood (with $r\geq 2$) of the given parameters $(g_\infty,\tau_\infty,\vartheta_\infty)$. The proof of \cite[Proposition 3.18, p. 332]{FL1} (the power $1/2$ on the right-hand side should be omitted) and its corollary as \cite[Lemma 4.9, p. 341]{FL1} (a hypothesis that $\Omega$ is strongly simply connected should be added) extends without change. The proof of Removability of Point Singularities for non-Abelian monopoles, namely \cite[Theorem 4.10, p. 342]{FL1},  is also valid for parameters $(g,\tau,\vartheta)$ in a small $C^r$-open neighborhood (with $r\geq 2$) of $(g_\infty,\tau_\infty,\vartheta_\infty)$. That argument appeals to \apriori estimates in \cite[Proposition 3.12 and Theorem 3.13, p. 328]{FL1} and, again, they hold with the same constants for parameters $(g,\tau,\vartheta)$ in a small $C^r$-open neighborhood of $(g_\infty,\tau_\infty,\vartheta_\infty)$.

 The patching arguments for sequences $\{(A_\alpha,\Phi_\alpha)\}_{\alpha\in\NN}$ relative to fixed parameters $(g,\tau,\vartheta)$ developed in \cite[Lemmas 4.14, 4.15, 4.16, Corollary 4.17, and Proposition 4.18, pp. 348--349]{FL1} carry over without change to allow for sequences $\{(g_\alpha,\tau_\alpha,\vartheta_\alpha)\}_{\alpha\in\NN}$. Finally, the proof of \cite[Theorem 4.20, p. 352]{FL1} carries over without change to this context as well and this completes the proof of Theorem \ref{thm:Uhlenbeck_convergence_sequence_non-Abelian_monopoles_parameters}.
\end{proof}

\section{Stratum of anti-self-dual or zero-section non-Abelian monopoles}
\label{sec:ASDsingularities}
For additional details concerning the construction and properties of the moduli space of anti-self-dual connections, we refer the reader to Donaldson and Kronheimer \cite{DK}, Friedman and Morgan \cite{FrM}, Freed and Uhlenbeck \cite{FU}, Kronheimer and Mrowka \cite{KMStructure}, and Lawson \cite{Lawson}.

As in Section \ref{sec:PU2Monopoles}, we continue to assume that there is a fixed, smooth unitary connection $A_d$ on the smooth Hermitian line bundle $\det E$ defined by the smooth Hermitian vector bundle $(E,h)$ and that $(X,g)$ is a closed, four-dimensional, oriented, smooth Riemannian manifold. While the non-Abelian monopole equations \eqref{eq:PerturbedSO3MonopoleEquations} are defined for a Hermitian vector bundle $E$ of arbitrary rank, we shall restrict to rank two in our discussions of split non-Abelian monopoles unless stated otherwise. From equation \eqref{eq:PerturbedSO3MonopoleEquations}, we see that the stratum of $\sM_{\ft}$ represented by pairs with zero spinor is identified with
\begin{equation}
\label{eq:ASDModuliSpace}
M_\kappa^w(X,g) := \left.\{A\in\sA(E,h): (F_A^+)_0 = 0\}\right/W^{2,p}(\SU(E)),
\end{equation}
the \emph{moduli space of $g$-anti-self-dual connections} on the $\SO(3)$ bundle $\su(E)$, where $\ka:=-\frac{1}{4} p_1(\ft)$ and $w\equiv w_2(\ft)\pmod 2$. If $b^+(X)>0$, then by Donaldson and Kronheimer \cite[Corollary 4.3.18, p. 149]{DK} and Kronheimer and Mrowka \cite[Lemma 2.4 and Corollary 2.5]{KMStructure} for a generic Riemannian metric $g$, the space $M_\kappa^w(X,g)$ is a smooth manifold of the expected dimension \eqref{eq:Expected_dimension_moduli_space_ASD_connections},
\[
  \expdim M_\kappa^w(X,g) = d_a(\ft) = -2p_1(\ft) - \frac{3}{2}(e(X)+\sigma(X)),
\]
where $d_a(\ft)$ was defined in \eqref{eq:Transv}.
As explained in \cite[Section 3.4, p. 96]{FL2a}, it is desirable to choose $w\pmod{2}$ so as to exclude points in $\bar\sM_{\ft}$ with associated flat $\SO(3)$ connections, so we have a \emph{disjoint} union,
\begin{equation}
\label{eq:StratificationCptPU(2)Space}
\bar\sM_{\ft}
\cong
\bar\sM_{\ft}^{*,0} \sqcup \bar M_\kappa^w \sqcup \bar\sM_{\ft}^{\red},
\end{equation}
where $\bar\sM_{\ft}^*\subset\bar\sM_{\ft}$ is the subspace represented by triples whose associated $\SO(3)$ connections are non-split, $\bar\sM_{\ft}^0\subset\bar\sM_{\ft}$ is the subspace represented by triples whose spinors are not identically zero, $\bar\sM_{\ft}^{*,0} = \bar\sM_{\ft}^{*}\cap\bar\sM_{\ft}^{0}$, while $\bar\sM_{\ft}^{\red}\subset\bar\sM_{\ft}$ is the subspace $\bar\sM_{\ft}-\bar\sM_{\ft}^*$ represented by triples whose associated $\SO(3)$ connections are split. We recall the

\begin{defn}[Good cohomology classes]
\label{defn:Good}
(See Feehan and Leness \cite[Definition 3.20, p. 169]{FL2b}.)
A class $v\in H^2(X;\ZZ/2\ZZ)$ is \emph{good} if no integral lift of $v$ is torsion.
\end{defn}

If $w\pmod{2}$ is good, then for \emph{generic} metrics the union \eqref{eq:StratificationCptPU(2)Space} is disjoint, as desired. In practice, rather than constraining $w\pmod{2}$ itself, we use the blow-up trick of Morgan and Mrowka \cite{MorganMrowkaPoly}, replacing $X$ with the smooth blow-up, $X\#\overline{\CC\PP}^2$, and replacing $w$ by $w+\PD[e]$ (where $e\in H_2(X;\ZZ)$ is the exceptional class and $\mathrm{PD}[e]$ denotes its Poincar{\'e} dual), noting that $w+\PD[e]\pmod{2}$ is always good.

\section{Strata of Seiberg--Witten or split non-Abelian monopoles}
\label{sec:Reducibles}
We now describe the relationship between the moduli spaces of Seiberg--Witten invariants and the moduli space of non-Abelian monopoles. For additional details concerning the definition of Seiberg--Witten invariants, we refer the reader to Kronheimer and Mrowka \cite{KMThom, KMBook}, Morgan \cite{MorganSWNotes}, Nicolaescu \cite{NicolaescuSWNotes}, Salamon \cite{SalamonSWBook}, and Witten \cite{Witten}.  

We introduce the Seiberg--Witten equations and define the Seiberg--Witten moduli space in Section \ref{subsec:SWMonopoles}.  Results on compactness, transversality, expected dimension, the existence of zero-section Seiberg--Witten monopoles, and the definitions of the Seiberg--Witten invariant, basic classes, simple type, and a blow-up formula for the Seiberg--Witten invariant appear in Section \ref{subsec:SWInvariants}. Finally, in Section \ref{subsec:RedPU2Monopole}, we characterize subspaces of split non-Abelian monopoles as the image of spaces of Seiberg--Witten monopoles with a particular perturbation under embeddings of the kind discussed in Section \ref{subsec:SpincQuotientSpaceInSpinuSpace}.

\subsection{Seiberg--Witten monopoles}
\label{subsec:SWMonopoles}
As in the definition of spin${}^c$ pairs \eqref{eq:TranslateSWNotation}, we fix a spin${}^c$ structure $\fs_0=(\rho,W)$ and a spin connection $A_W$ on $W$. Let $(L,h_L)$ be a Hermitian line bundle over the closed, oriented, smooth Riemannian four-manifold $(X,g)$. By analogy with Feehan and Leness \cite[Equations (2.55), p. 76, and (2.57), p. 77]{FL2a}, we call a pair $(A_L,\Psi)\in \sA(L,h_L)\times W^{1,p}(W^+\otimes L)$ a solution to the \emph{perturbed Seiberg--Witten monopole equations} on the spin${}^c$ structure $\fs=\fs_0\otimes L$ if
\begin{subequations}
\label{eq:SeibergWitten}
\begin{align}
  \label{eq:SeibergWitten_curvature}
  \tr_{W^+}F^+_{A_W} + 2F^+_{A_L} - \tau\rho^{-1}(\Psi\otimes\Psi^*)_{0} - F^+_{A_{\Lambda}} &=0,
  \\
  \label{eq:SeibergWitten_Dirac}
  D_{A_L}\Psi + \rho(\vartheta)\Psi &=0,
\end{align}
\end{subequations}
where $\tr_{W^+}:\fu(W^+)\to i\ubarRR$ is defined by the trace on $2\times 2$ complex matrices, $(\Psi\otimes\Psi^*)_0$ is the component of the section $\Psi\otimes\Psi^*$ of $i\fu(W^+)$ contained in $i\su(W^+)$, and $D_{A_L}:C^\infty(W^+\otimes L)\to C^\infty(W^-\otimes L)$ is the Dirac operator defined by the spin connection $A_W\otimes A_L$,
\label{page:SpincDiracOperator}
and $A_\Lambda := A_d\otimes A_{\det W^+}$ is a unitary connection on the Hermitian line bundle $\det W^+\otimes \det E$ with first Chern class $\Lambda\in H^2(X;\ZZ)$. The perturbations --- in particular the term $F_{A_\Lambda}^+$ in \eqref{eq:SeibergWitten_curvature} --- are chosen so that if $\ft$ is a spin${}^u$ structure with $c_1(\ft)=\Lambda$ and the embedding $\tilde\iota_{\fs,\ft}$ defined in \eqref{eq:DefnOfIota} exists, then $\tilde\iota_{\fs,\ft}$ maps solutions to equation \eqref{eq:SeibergWitten} to solutions to the non-Abelian monopole equations \eqref{eq:PerturbedSO3MonopoleEquations} that are split in the sense of Definition \ref{defn:Split_trivial_central-stabilizer_spinor_pair} (see Feehan and Leness \cite[Lemma 3.12, p. 95]{FL2a}).

The perturbed Seiberg--Witten equations \cite[Equation (2.55), p. 76, and Equation (2.57), p. 77]{FL2a} are expressed in terms of the spin connection $B := A_W\otimes A_L$ on $(\rho\otimes \id_L,W\otimes L)$ given by the isomorphism of affine spaces \eqref{eq:TranslateSWNotation}. The equalities,
\[
  \tr_{W^+\otimes L} F_B^+
  =
  \tr_{W^+\otimes L} F^+_{A_W\otimes A_L}
  =
 \tr_{W^+} F^+_{A_W}+2F^+_{A_L},
\]
show that the curvature terms appearing in \eqref{eq:SeibergWitten_curvature} are equivalent to those in
\cite[Equation (2.55), p. 76]{FL2a}. The Dirac operator $D_{A_L}$ in \eqref{eq:SeibergWitten_Dirac} is written as $D_B$ in \cite[Equation (2.55), p. 76]{FL2a}. We let
\begin{equation}
\label{eq:Moduli_space_Seiberg-Witten_monopoles}
  M_{\fs}
  :=
  \left.\left\{(A_L,\Psi)\in\sA(L,h_L)\times W^{1,p}(W^+\otimes L): \text{Equation \eqref{eq:SeibergWitten} holds} \right\}
    \right/W^{2,p}(X,S^1)
\end{equation}
denote the \emph{moduli space of Seiberg--Witten monopoles}, where $W^{2,p}(X,S^1)$ acts on the affine space $\sA(L,h_L)\times W^{1,p}(W^+\otimes L)$ as in \eqref{eq:SWGaugeGroupAction}.

\subsection{Seiberg--Witten invariants, basic classes, and simple type}
\label{subsec:SWInvariants}
We recall the following basic results for the moduli space of Seiberg--Witten monopoles.

\begin{prop}[Existence of zero-section Seiberg--Witten monopoles]
\label{prop:SWModuliSpaceZeroSectionCriterion}
(See Morgan \cite[Proposition 6.3.1, p. 91]{MorganSWNotes}.)
Let $\fs$ and $\ft$ be a spin${}^c$ and a rank-two spin${}^u$ structure, respectively, over a closed, four-dimensional,oriented, smooth Riemannian manifold $(X,g)$. If $b^+(X)>0$ and $c_1(\fs)-c_1(\ft)$ is not a torsion class in $H^2(X;\ZZ)$, then there is an open, dense subspace of Riemannian metrics on $X$ in the $C^r$ topology (where $r\geq 3$) such that if $g$ belongs to this subspace, then the moduli space $M_\fs$ does not contain the gauge-equivalence class of a zero-section pair.
\end{prop}

\begin{rmk}[Open and dense subspace of Riemannian metrics]
In Morgan \cite[Proposition 6.3.1, p. 91]{MorganSWNotes}, the conclusion of Proposition \ref{prop:SWModuliSpaceZeroSectionCriterion} is stated to hold for generic Riemannian metrics, meaning Riemannian metrics lying in a Baire set, as described by Freed and Uhlenbeck in \cite[Theorem 3.17, p. 59]{FU}). Applying the argument appearing in Donaldson and Kronheimer \cite[Corollary 4.3.15, p. 148]{DK} gives the conclusion of Proposition \ref{prop:SWModuliSpaceZeroSectionCriterion} for an open and dense subspace of Riemannian metrics.
\end{rmk}

\begin{thm}[Compactness for the moduli space of Seiberg--Witten monopoles]
\label{thm:Compactness_SW}
(See Feehan and Leness \cite[Proposition 2.15, p. 78]{FL2a} and Morgan \cite[Corollary 5.3.7, p. 85]{MorganSWNotes}.)
Let $\fs=(\rho,W)$ be a spin${}^c$ structure over a closed, four-dimensional, oriented, smooth Riemannian manifold $(X,g)$. Then the moduli space $M_\fs$ of Seiberg--Witten monopoles \eqref{eq:Moduli_space_Seiberg-Witten_monopoles} is a compact subset of $\sC_\fs$.
\end{thm}

\begin{prop}[Expected dimension of the moduli space of Seiberg--Witten monopoles]
\label{prop:Expected_dimension_moduli_space_Seiberg-Witten_monopoles}
(See Morgan \cite[Corollary 4.6.2, p. 67]{MorganSWNotes}.)
Continue the hypotheses of Theorem \ref{thm:Compactness_SW}. Then the expected dimension of the moduli space $M_\fs$ of Seiberg--Witten monopoles \eqref{eq:Moduli_space_Seiberg-Witten_monopoles} is
\begin{equation}
\label{eq:DimSW}
d_s(\fs)
=
\dim M_{\fs}
=
\frac{1}{4}(c_1(\fs)^2 -2e(X) -3\sigma(X)).
\end{equation}
where $e(X)$ is the Euler characteristic and $\sigma(X)$ is the signature of $X$, and is independent of the choice of parameters $(g,\tau,\vartheta)$ in the Seiberg--Witten monopole equations \eqref{eq:Moduli_space_Seiberg-Witten_monopoles}.
\end{prop}

\begin{prop}[Transversality for the moduli space of Seiberg--Witten monopoles with generic parameters]
\label{prop:Transv_SW}
(See Feehan and Leness \cite[Proposition 2.16, p. 79]{FL2a}.)
Continue the hypotheses of Theorem \ref{thm:Compactness_SW}. If the parameters $(g,\tau,\vartheta)$ in the Seiberg--Witten monopole equations \eqref{eq:SeibergWitten} are generic in the sense of \cite{FeehanGenericMetric}, then $M_\fs^0:=M_\fs\cap\sC_\fs^0$ is an embedded, orientable, smooth submanifold of the Banach manifold $\sC_\fs^0$ that has finite dimension equal to the expected dimension \eqref{eq:DimSW}.
\end{prop}

Finally, we have the following analogue of Theorem \ref{thm:Uhlenbeck_convergence_sequence_non-Abelian_monopoles_parameters}.

\begin{thm}[Uhlenbeck convergence of a sequence of Seiberg--Witten monopoles with respect to a convergent sequence of parameters]
\label{thm:Uhlenbeck_convergence_sequence_SW_monopoles_parameters}
Continue the hypotheses of Theorem \ref{thm:Compactness_SW}. Let $\{(g_\alpha,\tau_\alpha,\vartheta_\alpha)\}_{\alpha\in\NN}$ be a sequence of parameters defining the Seiberg--Witten monopole equations \eqref{eq:SeibergWitten} that converges in $C^\infty$ to a limit $(g_\infty,\tau_\infty,\vartheta_\infty)$. If $\{(A_{L,\alpha},\Psi_\alpha)\}_{\alpha\in\NN}$ is a sequence of solutions to the corresponding equations \eqref{eq:SeibergWitten} then, after passing to a subsequence and relabelling, this sequence of solutions converges in $W^{1,p}$ modulo $W^{2,p}$ gauge transformations to a $W^{1,p}$ solution $(A_{L,\infty},\Psi_\infty)$
to \eqref{eq:SeibergWitten} for the parameters $(g_\infty,\tau_\infty,\vartheta_\infty)$.
\end{thm}

\begin{proof}
The proof of Theorem \ref{thm:Uhlenbeck_convergence_sequence_SW_monopoles_parameters} follows from that of Theorem \ref{thm:Compactness_SW} by modifications exactly analogous to and simpler than those described in our proof of Theorem \ref{thm:Uhlenbeck_convergence_sequence_non-Abelian_monopoles_parameters}.
\end{proof}

Let $\tilde X=X\#\overline{\CC\PP}^2$ denote the smooth blow-up of $X$ with exceptional class $e\in H_2(\tilde X;\ZZ)$ and denote its Poincar\'e dual by $\PD[e]\in H^2(\tilde X;\ZZ)$ (see Gompf and Stipsicz \cite[Definition 2.2.7, p. 43]{GompfStipsicz}). Let $\fs^\pm=(\tilde\rho,\tilde W)$ denote the \spinc structure on $\tilde X$ with $c_1(\fs^\pm)=c_1(\fs)\pm \PD[e]$ obtained by splicing the \spinc structure $\fs=(\rho,W)$ on $X$ with the \spinc structure on $\overline{\CC\PP}^2$ with first Chern class $\pm \PD[e]$. (See Feehan and Leness \cite[Section 4.5, p. 200]{FL2b} or Salamon \cite[Section 12.4]{SalamonSWBook} for an explanation of the relation between \spinc structures on $X$ and $\tilde X$.) Now
$$
c_1(\fs)\pm \PD[e]-\Lambda \in H^2(\tilde X;\ZZ)
$$
is not a torsion class and so --- for $b^+(X)>0$, generic Riemannian metrics on $X$ and related metrics on the connected sum $\tilde X$ --- the moduli spaces $M_{\fs^\pm}$ contain no zero-section pairs. Thus, for our choice of generic perturbations, the moduli spaces $M_{\fs^\pm}$ are compact, oriented, smooth manifolds, both of dimension $\dim M_{\fs}$.

The \emph{Seiberg--Witten invariant} $\SW_{X}(\fs)$
\label{page:SW_Invariant}
defined by a \spinc structure $\fs$ over $X$ is most easily defined when $b_1(X)=0$ and $b^+(X)>1$ is odd; we refer the reader to Feehan and Leness \cite[Section 4.1]{FL2b} for its definition and references to the development of its properties. One says that $c_1(\fs)$ is a \emph{Seiberg--Witten basic class}\label{page:Basic_class} if $\SW_{X}(\fs) \neq 0$.   We will denote the set of basic classes by
\begin{equation}
\label{eq:SetOfBasicClasses}
B(X):=\{c_1(\fs): \SW_X(\fs)\neq 0\}.
\end{equation}
A four-manifold $X$ has \emph{Seiberg--Witten simple type} \label{page:Seiberg_Witten_simple_type} if all basic classes satisfy
\begin{equation}
  \label{eq:SWSimpleType}
  d_s(\fs)=0,
\end{equation}
where $d_s(\fs)$ is as in \eqref{eq:DimSW}. Versions of the following  have appeared in Fintushel and Stern \cite[Theorem 1.4]{FSTurkish}, Nicolaescu \cite[Theorem 4.6.7]{NicolaescuSWNotes}, and
Fr{\o}yshov \cite[Corollary 14.1.1]{Froyshov_2008}.

\begin{thm}[Blow-up formula for Seiberg--Witten invariants]
\label{thm:SWBlowUp}
Let $X$ be a standard four-manifold and let $\widetilde X=X\#\overline{\CC\PP}^2$. Write $\PD[e]\in H^2(\widetilde X;\ZZ)$ for the Poincar\'e dual of the exceptional sphere. Then $X$ has Seiberg--Witten simple type if and only if $\widetilde X$ does. If $X$ has Seiberg--Witten simple type then
\begin{equation}
\label{eq:BlowUpBasics}
B(\widetilde X)=\{c_1(\fs_0)\pm \PD[e]: c_1(\fs_0)\in B(X)\}.
\end{equation}
If $c_1(\fs_0)\in B(X)$ and $\fs^\pm$ is a spin${}^c$ structure on $\widetilde X$ with
$c_1(\fs^\pm)=c_1(\fs_0)\pm \PD[e]$, then $\SW_{\widetilde X}(\fs^\pm)=\SW_X(\fs)$.
\end{thm}

\subsection{Split non-Abelian monopoles}
\label{subsec:RedPU2Monopole}
Let $\ft$ be a rank-two spin${}^u$ structure on a closed, smooth Riemannian four-manifold. We now describe the equivalence classes of pairs in $\sM_\ft$ which are split in the sense of Definition \ref{defn:Split_trivial_central-stabilizer_spinor_pair}. By \cite[Lemma 3.13, p. 96]{FL2a} the restriction of the map $\iota_{\fs,\ft}$ defined in \eqref{eq:DefnOfIota} to $M_\fs^0$ is a topological embedding $M^0_\fs\embed\sM_\ft$ where $M_{\fs}^0:=M_{\fs}\cap\sC_{\fs}^0$. If $w_2(\ft)\neq 0$ or $b_1(X)=0$, then \cite[Lemma 3.13, p. 96]{FL2a} implies that $\iota_{\fs,\ft}$ is an embedding. The image of $\iota_{\fs,\ft}$ is represented by points in $\sM_\ft$ represented by pairs which are split with respect to the splitting \eqref{eq:BasicSplitting}. Henceforth, we shall not distinguish between $M_\fs$ and its image in $\sM_{\ft}$ under this embedding.

If $\su(E_\ell) \cong i\ubarRR\oplus L$ for some non-negative integer $\ell$ and Hermitian line bundle $L$ with $c_1(L)=c_1(\ft)-c_1(\fs)$, then $p_1(\ft(\ell))=(c_1(\ft)-c_1(\fs))^2$.  Hence, for
\begin{equation}
\label{eq:ReducibleLevel}
\ell(\ft,\fs) := \frac{1}{4}\left((c_1(\ft)-c_1(\fs))^2- p_1(\ft)\right),
\end{equation}
the smooth embedding \eqref{eq:DefnOfIotaOnQuotient} gives a continuous embedding
\begin{equation}
\label{eq:LowerLevelInclusionOFReducibles}
M_{\fs}\times\Sym^\ell(X) \to \sI\!\!\sM_{\ft},
\end{equation}
where $\ell=\ell(\ft,\fs)$.  (Note that the embedding \eqref{eq:DefnOfIotaOnQuotient} is given by \eqref{eq:LowerLevelInclusionOFReducibles} with $\ell(\ft,\fs)=0$). Let $\Spinc(X)$ \label{page:Spinc(X)} denote the set of all \spinc structures $\fs = (\rho,W)$ as in Section \eqref{sec:SpincuStr} on an even-dimensional, oriented, smooth Riemannian manifold $(X,g)$. We record the following cohomological criterion for when the continuous embedding \eqref{eq:LowerLevelInclusionOFReducibles} exists.

\begin{lem}[Existence of embeddings of moduli spaces of Seiberg--Witten monopoles into moduli space of ideal non-Abelian monopoles]
\label{lem:SetOfReducibles}
Let $\ft=(\rho,W\otimes E)$ be a rank-two spin${}^u$ structure over a closed, oriented, smooth, Riemannian four-manifold $(X,g)$. If $\ell$ is a non-negative integer and
\begin{equation}
\label{eq:DefineReduciblesEmbedded}
\Red_\ell(\ft) := \left\{c_1(\fs): \fs \in \Spinc(X) \text{ such that } (c_1(\fs)-c_1(\ft))^2=p_1(\su(E))+4\ell\right\},
\end{equation}
then the continuous embedding $M_\fs\times\Sym^\ell(X)\to \sI\!\!\sM_{\ft}$ in \eqref{eq:LowerLevelInclusionOFReducibles} exists if and only if $c_1(\fs)\in\Red_\ell(\ft)$.
\end{lem}

\begin{proof}
The embedding \eqref{eq:LowerLevelInclusionOFReducibles} exists if and only  the bundles $E_\ell$ and $E':=L_1\oplus \det(E)\otimes L_1^*$ are isomorphic. By Gompf and Stipsicz \cite[Theorem 1.4.20 (a), p. 31]{GompfStipsicz}, $E_\ell\cong E'$ if and only if the first two Chern classes of the bundles $E_\ell$ and $E'$ are equal. We compute
\begin{align*}
c_1(E')
&=
c_1\left( L_1\oplus \det E\otimes L_1^*\right)
=
c_1(L_1)+c_1(\det E)+c_1(L_1^*)
=
c_1(E)
=
c_1(E_\ell),
\end{align*}
where the final equality holds by  \eqref{eq:DefineLowerChargeSpinuStr_c1}. Thus $c_1(E_\ell)=c_1(E')$ always holds. We thus only need to show that $c_2(E_\ell)=c_2(E')$ if and only if $c_1(\fs)\in\Red_\ell(\ft)$  or equivalently,
\begin{equation}
\label{eq:ReducibleLevelEquality}
\left(c_1(\fs)-c_1(\ft)\right)^2=p_1(\su(E))+4\ell.
\end{equation}
We begin by computing,
\begin{align*}
p_1(\su(E))+4\ell
&= p_1(\su(E_\ell)) \quad\text{(by \eqref{eq:DefineLowerChargeSpinuStr_p1})}
\\
&=
c_1(E_\ell)^2 - 4 c_2(E_\ell) \quad\text{(by \cite[Equation (2.1.39), p. 42]{DK})},
\end{align*}
so
\begin{equation}
\label{eq:p1E_ell}
p_1(\su(E))+4\ell=c_1(E_\ell)^2 - 4 c_2(E_\ell)
\end{equation}
Next, observe that
\[
c_2(E')
=
c_2\left( L_1\oplus \det E\otimes L_1^*\right)
=
c_1(L_1)\left( c_1(E)-c_1(L_1)\right)
\]
so
\begin{equation}
\label{eq:c2E'ForReducibles}
c_2(E')
=
c_1(L_1)c_1(E)-c_1(L_1)^2.
\end{equation}
Because $c_1(\fs)=c_1(W^+\otimes L_1)=c_1(W^+)+2c_1(L_1)$ and $c_1(\ft)=c_1(W^+)+c_1(E)$,
\begin{align*}
\left(
c_1(\fs)-c_1(\ft)\right)^2
&=
\left( 2c_1(L_1)-c_1(E)\right)^2
\\
& =4c_1(L_1)^2-4c_1(L_1)c_1(E)+c_1(E)^2
\\
&=-4c_2(E')+c_1(E)^2\quad\text{(by \eqref{eq:c2E'ForReducibles})}
\\
&=c_1(E_\ell)^2-4c_2(E')\quad\text{(by \eqref{eq:DefineLowerChargeSpinuStr_c1})}
\end{align*}
so
\begin{equation}
\label{eq:c1s-c1t_Squared}
\left(c_1(\fs)-c_1(\ft)\right)^2
=
c_1(E_\ell)^2-4c_2(E')
\end{equation}
By comparing \eqref{eq:p1E_ell} and \eqref{eq:c1s-c1t_Squared}, we see that $c_2(E_\ell)=c_2(E')$ if and only if \eqref{eq:ReducibleLevelEquality} holds, completing the proof of the lemma.
\end{proof}

\begin{rmk}[Conventions on splitting of bundles in \spinu structures]
\label{rmk:ComparingLineBundlesFromFl2a}
As we shall often refer to constructions from Feehan and Leness \cite{FL2a, FL2b}, it is worth noting the following difference between the notation used in this monograph and that in \cite{FL2a,FL2b}.  In \cite{FL2a}, the splitting of $W\otimes E$ is written as $W'\oplus W'\otimes L$, whereas the splitting of $W\otimes E$ given by \eqref{eq:BasicSplitting} would be $W\otimes L_1\oplus W\otimes L_2$.  To convert between the two conventions, use $W'=W\otimes L_1$ and $L=L_2\otimes L_1^*$.
\end{rmk}

\section{Circle actions on the affine space of spin${}^u$ pairs}
\label{sec:S1Actions}
The affine space $\sA(E,h)\times W^{1,p}(W^+\otimes E)$ carries a circle action \label{Circle_action_on_SO(3)_pairs} induced by scalar multiplication on the Hermitian vector bundle $E$:

\begin{defn}[Standard $S^1$ action on affine and quotient spaces of spin${}^u$ pairs]
\label{defn:UnitaryZActionOnAffine}
The \emph{standard $S^1$ action on the affine space of spin${}^u$ pairs} is
\begin{multline}
  \label{eq:S1ZAction}
  S^1\times \sA(E,h)\times W^{1,p}(W^+\otimes E)
  \ni (e^{i\theta},A,\Phi)
  \\
  \mapsto (A,e^{i\theta}\Phi) \in \sA(E,h)\times W^{1,p}(W^+\otimes E).
\end{multline}
Because the circle action \eqref{eq:S1ZAction} commutes with that of $W^{2,p}(\SU(E))$ in \eqref{eq:GaugeActionOnSpinuPairs}, the action \eqref{eq:S1ZAction} also defines a circle action on the quotient,
\begin{equation}
\label{eq:S1ZActionOnQuotientSpace}
S^1\times\sC_\ft\to\sC_\ft,
\end{equation}
which we call \emph{standard $S^1$ action on the quotient space $\sC_\ft$}.
\end{defn}

Note that if $r$ is the rank of $\ft$, and $\varrho$ is an $r$-th root of unity, then $\varrho$ acts trivially
on $\sC_{\ft}$. We can see this by observing that for all $e^{i\theta}\in S^1$, the central gauge transformation
$e^{i\theta}\,\id_E\in W^{2,p}(U(E))$ acts trivially on $\sA(E,h)$ and so the action \eqref{eq:S1ZAction} is just the action \eqref{eq:GaugeActionOnSpinuPairs} extended to the elements of $W^{2,p}(\U(E))$ given by $e^{i\theta}\,\id_E$. Because $\varrho\,\id_E$ is a section of $\SU(E)$, it acts trivially on the $W^{2,p}(\SU(E))$-quotient $\sC_\ft$, although non-trivially on the affine space $\sA(E,h)\times W^{1,p}(W^+\otimes E)$.

We now identify the fixed points of the $S^1$ action in Definition \ref{defn:UnitaryZActionOnAffine} on $\sC_\ft$.
Since equation \eqref{eq:PerturbedSO3MonopoleEquations} is invariant under the circle action induced by scalar multiplication on $W^+\otimes E$, the subspaces \eqref{eq:PU2MonopoleSubspaces} of $\sC_{\ft}$ are also invariant under this action.  In \cite[Proposition 3.1, p. 86]{FL2a}, we characterized the fixed points of the action $S^1$ action in Definition \ref{defn:UnitaryZActionOnAffine} on the moduli space $\sM_\ft$ of non-Abelian monopoles.  The following proposition generalizes those results to the quotient space $\sC_\ft$.

\begin{prop}[Fixed points of the $S^1$ action in $\sC_\ft$]
\label{prop:FixedPointsOfS1ActionOnSpinuQuotientSpace}
Let $\ft=(\rho,W\otimes E)$ be a spin${}^u$ structure on a connected, oriented, smooth Riemannian four-manifold. If $[A,\Phi] \in \sC_\ft$ obeys
\begin{enumerate}
\item
\label{item:SO(3)MonopoleFixedPointsZeroSection}
$\Phi\equiv 0$, so $(A,0)$ is a zero-section pair as in Definition \ref{defn:ZeroSectionSO(3)Pair}, or
\item
\label{item:SO(3)MonopoleFixedPointsReducible}
$(A,\Phi)$ is
split as in Definition \ref{defn:Split_trivial_central-stabilizer_spinor_pair},
\end{enumerate}
then $[A,\Phi]$ is a fixed point of the $S^1$ action in Definition \ref{defn:UnitaryZActionOnAffine} on $\sC_\ft$. Conversely, if a fixed point of that $S^1$ action on $\sC_\ft$ is represented by smooth pair $(A,\Phi)$, then it obeys Condition \eqref{item:SO(3)MonopoleFixedPointsZeroSection} or \eqref{item:SO(3)MonopoleFixedPointsReducible}. Lastly, if $E$ has rank two  and the characteristic class $w_2(\ft)$ defined in \eqref{eq:SpinUCharacteristics} is good as in Definition \ref{defn:Good}, then the moduli space $\sM_\ft$ of non-Abelian monopoles \eqref{eq:Moduli_space_SO(3)_monopoles} contains no points represented by a pair that is both split and zero-section.
\end{prop}

\begin{proof}
If $(A,\Phi)$ satisfies Condition \eqref{item:SO(3)MonopoleFixedPointsZeroSection}, so $\Phi\equiv 0$, then $[A,0]$ is a fixed point of the $S^1$ action on $\sC_\ft$ in Definition \ref{defn:UnitaryZActionOnAffine} since $S^1\subset \Stab(A)$.

If $(A,\Phi)$ satisfies Condition \eqref{item:SO(3)MonopoleFixedPointsReducible}, so $(A,\Phi)$ is a split pair with respect to an orthogonal decomposition $E = E_1\oplus E_1^\perp$ as in Definition \ref{defn:Split_trivial_central-stabilizer_spinor_pair}, then because the one-parameter family of gauge transformations $u(\theta) \in W^{2,p}(\SU(E))$ defined by \eqref{eq:One-parameter_family_SU(E)_gauge_transformations} is in $\Stab(A)$ and thus acts trivially on $A$ while $u(\theta)^{-1}\Phi = e^{i(r-r_1)\theta}\Phi$ since $\Phi \in W^{1,p}(W^+\otimes E_1)$, the pair $(A,\Phi)$ satisfies $(A,e^{i(r-r_1)\theta}\Phi) = u(\theta)^*(A,\Phi)$ and thus
\[
  (A,e^{i\theta}\Phi) = u(\theta/(r-r_1))^*(A,\Phi), \quad\text{for all } \theta \in \RR.
\]
Since $u(\theta/(r-r_1)) \in W^{2,p}(\SU(E))$, we see that $[A,\Phi] = [A,e^{i\theta}\Phi]$. Hence,
if $(A,\Phi)$ satisfies Condition \eqref{item:SO(3)MonopoleFixedPointsZeroSection} or  \eqref{item:SO(3)MonopoleFixedPointsReducible}, then $[A,\Phi]$ is a fixed point of the $S^1$ action on $\sC_\ft$ in Definition \ref{defn:UnitaryZActionOnAffine}.

Conversely, suppose that $(A,\Phi)$ is a smooth pair with $\Phi$ not identically zero and $[A,\Phi]\in \sC_\ft$ is a fixed point of the $S^1$ action on $\sC_\ft$ in Definition \ref{defn:UnitaryZActionOnAffine}. We claim that $(A,\Phi)$ satisfies Condition \eqref{item:SO(3)MonopoleFixedPointsReducible}. Because $[A,\Phi]\in \sC_\ft$ is a fixed point of the $S^1$ action, the smooth path
\[
  \RR\ni\theta\mapsto (A,e^{i\theta}\Phi) \in \sA^\infty(E,h)\times C^\infty(W^+\otimes E)
\]
is in the $C^\infty(\SU(E))$-orbit of $(A,\Phi)$, where $\sA^\infty(E,h) \subset \sA(E,h)$ denotes the affine subspace of smooth unitary connections on $E$ that induce $A_d$ on $\det E$.  Therefore, the tangent vector $(0,i\Phi)$ to this curve at $\theta=0$ is in the range of the linearization of the $C^\infty(\SU(E))$ action \eqref{eq:GaugeActionOnSpinuPairs},
\[
C^\infty(\su(E)) \ni \xi \mapsto (d_A\xi,-\xi\Phi) \in \Omega^1(\su(E)) \times C^\infty(W^+\otimes E),
\]
given in \eqref{eq:d_APhi^0}. Hence, there exists $\xi\in C^\infty(\su(E))$ such that $(-d_A\xi,\xi\Phi)=(0,i\Phi)$ and so
$d_A\xi=0$ and $\xi\Phi=i\Phi$.  The identity $\xi\Phi=i\Phi$ and the assumption $\Phi\not\equiv 0$ imply that $\xi\not\equiv 0$. By Item \eqref{item:ReducibleUnitaryFromH0Splitting_smooth_and_any_rank_differentiable} in Proposition \ref{prop:ReducibleUnitaryFromH0}, the connection $A$ is split, so $A=A_1\oplus\cdots\oplus A_s$ with respect to a decomposition $E=E_1\oplus\cdots\oplus E_s$ as an orthogonal direct sum of proper, smooth Hermitian subbundles for some $s\geq 2$.
Observe that
\[
  i\Phi = \xi\Phi = \sum_{k=1}^s \xi\pi_k\Phi = \sum_{k=1}^s i\mu_k\pi_k\Phi,
\]
where $\mu_1,\ldots,\mu_s$ are the constant, distinct eigenvalues of $\xi$ and thus, noting that $\Phi = \sum_{k=1}^s \pi_k\Phi$,
\[
  \pi_k\Phi = \mu_k\pi_k\Phi, \quad\text{for } k = 1,\ldots,s.
\]
We must have $\mu_k = 1$ for all $k$ such that $\pi_k\Phi$ is non-zero and because the $\mu_k$ are distinct and $\Phi\not\equiv 0$, we must have $\pi_l\Phi\not\equiv 0$ for exactly one index $l$. We may assume without loss of generality, after re-indexing if necessary, that $l=1$ and $\mu_1 = 1$. Hence, we must have $\mu_k=0$ and $\pi_k\Phi \equiv 0$ for $k=2,\ldots,s$, but the fact that the $\mu_k$ are distinct forces $s=2$. Therefore, $\Phi\in C^\infty(W^+\otimes E_1)$ and so $(A,\Phi)$ satisfies Condition \eqref{item:SO(3)MonopoleFixedPointsReducible}, as claimed.

By \cite[Lemma 3.2, p. 87]{FL2a} or Morgan and Mrowka \cite[p. 226]{MorganMrowkaPoly}, if the characteristic class $w_2(\ft)$ defined in \eqref{eq:SpinUCharacteristics} is good as in Definition \ref{defn:Good}, then $\sM_\ft$ contains no points represented by a pair which is both split and a zero-section pair.
\end{proof}

When the vector bundle $E$ has rank two, the action \eqref{eq:S1ZAction} is defined by the homomorphism from $S^1$ into the group of unitary gauge transformations of $E$,
\begin{equation}
\label{eq:DefineUnitaryS1ZActionsAtReducible}
\rho_Z(e^{i\theta}):=e^{i\theta}\,\id_{L_1}\oplus e^{i\theta}\,\id_{L_2},
\quad\text{for all } e^{i\theta} \in S^1.
\end{equation}
If the bundle $E$ admits an orthogonal splitting as a direct sum of Hermitian line bundles, $E=L_1\oplus L_2$, then there is an additional $S^1$ action defined by complex scalar multiplication on $L_2$.  We define two more homomorphisms from $S^1$ into the group of unitary gauge transformations of $E$ by
\begin{subequations}
\label{eq:DefineUnitaryS1ActionsAtReducible}
\begin{align}
\label{eq:DefineUnitaryS1sActionsAtReducible}
\rho_{\SU}(e^{i\theta})&:=e^{i\theta}\,\id_{L_1}\oplus e^{-i\theta}\,\id_{L_2},
\\
\label{eq:DefineUnitaryS12ActionsAtReducible}
\rho_2(e^{i\theta})&:=\id_{L_1}\oplus e^{i\theta}\,\id_{L_2}, \quad\text{for all } e^{i\theta} \in S^1.
\end{align}
\end{subequations}
Compare the $\GL(E)$-valued homomorphisms to be defined in \eqref{eq:DefineS1ActionsAtReducible}.

\begin{rmk}[Choice of line bundles in splitting]
\label{rmk:ChoiceOfLineBundlesInDefinitionOfSplitting}
We note that the definition of the homomorphisms \eqref{eq:DefineUnitaryS1ActionsAtReducible} depend not only on the existence of a splitting $E=L_1\oplus L_2$ but also on the choice of line bundles giving the splitting.  In the interest of legibility, we do not include the choice of line bundles $L_1$ and $L_2$ in the notation for the homomorphisms \eqref{eq:DefineUnitaryS1ActionsAtReducible}. 
\end{rmk}

Observe that the homomorphisms \eqref{eq:DefineUnitaryS1ZActionsAtReducible} and \eqref{eq:DefineUnitaryS1ActionsAtReducible} (when the latter exist) are related by
\begin{equation}
\label{eq:UnitaryS1ActionsRelation}
\rho_2(e^{2i\theta})=\rho_{\SU}(e^{-i\theta})\rho_Z(e^{i\theta}),
\quad\text{for all } e^{i\theta} \in S^1.
\end{equation}
The homomorphism \eqref{eq:DefineUnitaryS1ZActionsAtReducible} and, when $E$ admits a splitting $E=L_1\oplus L_2$, the homomorphism \eqref{eq:DefineUnitaryS12ActionsAtReducible}  define actions of $S^1$
on $\sA(E,h)\times W^{1,p}(W^+\otimes E)$ by
\begin{subequations}
\label{eq:S1ActionsOnSpinuPreQuotient}
\begin{align}
\label{eq:S1ZActionsOnSpinuPreQuotient}
\left( e^{i\theta},(A,\Phi)\right)
&\mapsto
( \rho_Z(e^{i\theta})^*A,\rho_Z(e^{-i\theta})\Phi),
\\
\label{eq:S12ActionsOnSpinuPreQuotient}
\left( e^{i\theta},(A,\Phi)\right)
&\mapsto
(\rho_2(e^{i\theta})^*A,\rho_2(e^{-i\theta})\Phi)
\end{align}
\end{subequations}
and on $\sC_\ft$ by
\begin{subequations}
\label{eq:S1ActionsOnSpinuQuotient}
\begin{align}
\label{eq:S1ZActionsOnSpinuQuotient}
\left( e^{i\theta},[A,\Phi]\right)
&\mapsto
[ \rho_Z(e^{i\theta})^*A,\rho_Z(e^{-i\theta})\Phi],
\\
\label{eq:S12ActionsOnSpinuQuotient}
\left( e^{i\theta},[A,\Phi]\right)
&\mapsto
[ \rho_2(e^{i\theta})^*A,\rho_2(e^{-i\theta})\Phi].
\end{align}
\end{subequations}
Because $\rho_{\SU}$ takes values in $C^\infty(\SU(E))$, an action defined by $\rho_{\SU}$ as in \eqref{eq:S1ActionsOnSpinuQuotient} would be trivial. The action \eqref{eq:S1ZActionsOnSpinuQuotient} is the action \eqref{eq:S1ZActionOnQuotientSpace}.  The equality
\eqref{eq:UnitaryS1ActionsRelation} implies that the actions \eqref{eq:S1ZActionsOnSpinuQuotient} and \eqref{eq:S12ActionsOnSpinuQuotient}  are equal up to multiplicity but each has its advantages.  The action \eqref{eq:S1ZActionsOnSpinuQuotient} is defined even when the bundle $E$ does not admit a splitting into line bundles
while the gauge transformations in the image of $\rho_2$ are in the stabilizer of pairs in the image of $\tilde\iota_{\fs,\ft}$ as we now describe.

\begin{lem}[Circle equivariance of the non-Abelian monopole map]
\label{lem:S1ActionsForReduciblePairGerm}
Let $\fs=(\rho,W\otimes L_1)$ and $\ft=(\rho,W\otimes E)$ be a spin${}^c$ and rank-two spin${}^u$ structure, respectively, over a closed, oriented, smooth, Riemannian four-manifold $(X,g)$. If $E$ admits a splitting $E=L_1\oplus L_2$ as an orthogonal direct sum of Hermitian line bundles, then the following hold
\begin{enumerate}
\item If $(A,\Phi)\in\sA(E,h)\times W^{1,p}(W^+\otimes E)$ is in the image of the map $\tilde\iota_{\fs,\ft}$ defined in \eqref{eq:DefnOfIota}, then $(A,\Phi)$ is a fixed point of the action \eqref{eq:S12ActionsOnSpinuPreQuotient}.
\item The non-Abelian monopole map \eqref{eq:PerturbedSO3MonopoleEquation_map},
\[
\sS:\sA(E,h)\times W^{1,p}(W^+\otimes E) \to
L^p(\La^+\otimes\su(E))\oplus L^p(W^-\otimes E),
\]
is $S^1$-equivariant with respect to the action \eqref{eq:S12ActionsOnSpinuPreQuotient} on the domain of $\sS$ and the action
\begin{multline}
\label{eq:S1L2ActionOnnonAbelianMonopoleMapCodomain}
S^1\times L^p(\La^+\otimes\su(E))\oplus L^p(W^-\otimes E)
\ni
\left(e^{i\theta},(\om,\Psi)\right)
\\
\mapsto
\left( \rho_2(e^{-i\theta})\om \rho_2(e^{i\theta}),\rho_2(e^{-i\theta})\Psi\right),
\end{multline}
on the codomain of $\sS$.
\end{enumerate}
\end{lem}

\begin{proof}
Because $(A,\Phi)$ is in the image of $\tilde\iota_{\fs,\ft}$, the definition of $\tilde\iota_{\fs,\ft}$ in \eqref{eq:DefnOfIota} implies that we can write $A=A_1\oplus A_2$ and $\Phi=\varphi_1\oplus 0$, where $A_j$ is a unitary connection on $L_j$, for $j=1,2$, and $\varphi_1\in W^{1,p}(W^+\otimes L_1)$.  Then $\rho_2(e^{i\theta})$ is in the stabilizer of $A$ and $\rho_2(e^{i\theta})\varphi_1=\varphi_1$, which proves the first assertion.

The $S^1$-equivariance of the non-Abelian monopole map $\sS$ in the second assertion follows immediately from the equivariance of $\sS$ with respect to gauge transformations.
\end{proof}

The homomorphism \eqref{eq:DefineUnitaryS1ZActionsAtReducible} and, when $E$ admits a splitting $E=L_1\oplus L_2$, the homomorphism \eqref{eq:DefineUnitaryS12ActionsAtReducible} also define $S^1$ actions on  $T\sC_\ft^0$ and its subbundles $T\sC_\fs^0$ and $\fN(\ft,\fs)$:
\begin{subequations}
\label{eq:S1ActionOnTangentBundle}
\begin{align}
\label{eq:S1ZActionOnTangentBundle}
\left( e^{i\theta}, [A,\Phi,a,\phi]\right)
&\mapsto
[ \rho_Z(e^{i\theta})^*A,\rho_Z(e^{-i\theta})\Phi,\rho_Z(-e^{i\theta})a \rho_Z(e^{i\theta}),\rho_Z(e^{-i\theta})\phi],
\\
\label{eq:S12ActionOnTangentBundle}
\left( e^{i\theta}, [A,\Phi,a,\phi]\right)
&\mapsto
[ \rho_2(e^{i\theta})^*A,\rho_2(e^{-i\theta})\Phi,\rho_2(e^{-i\theta})a \rho_2(e^{i\theta}),\rho_2(e^{-i\theta})\phi],
  \\
  \notag
  &\qquad\text{for all}\ e^{i\theta}\in S^1\ \text{and }  [A,\Phi,a,\phi]\in T\sC_\ft^0.
\end{align}
\end{subequations}
The map $\exp:T\sC_\ft^0\to\sC_\ft$ defined in \eqref{eq:TangentSpaceMap} is $S^1$-equivariant with respect to the $S^1$ action  on $T\sC_\ft^0$ given in \eqref{eq:S1ZActionOnTangentBundle} and the action given by scalar multiplication on the spinor component in \eqref{eq:S1ZActionOnQuotientSpace}) on $\sC_\ft^0$. As is true for the actions \eqref{eq:S1ActionsOnSpinuQuotient}, the actions \eqref{eq:S1ZActionOnTangentBundle} and \eqref{eq:S12ActionOnTangentBundle} are equal up to multiplicity.

\begin{lem}[Circle equivariance of the exponential map]
\label{lem:S1EquivarianceOfg}
Let $\ft=(\rho,W\otimes E)$ be a rank-two spin${}^u$ structure over a closed, oriented, smooth, Riemannian four-manifold $(X,g)$. Then the map $\exp:T\sC_\ft^0\to \sC_\ft$ defined in \eqref{eq:TangentSpaceMap} is equivariant with respect to the $S^1$ action \eqref{eq:S1ZActionOnTangentBundle} on $T\sC_\ft^0$ and the $S^1$ action \eqref{eq:S1ZActionsOnSpinuQuotient} on $\sC_\ft$ (both implied by $\rho_Z$). If $E$ admits a splitting $E=L_1\oplus L_2$ as an orthogonal direct sum of Hermitian line bundles, then the map $\exp$ is also equivariant with respect to the $S^1$ action \eqref{eq:S12ActionOnTangentBundle}  on $T\sC_\ft^0$ and the $S^1$ action \eqref{eq:S12ActionsOnSpinuQuotient} on $\sC_\ft$ (both implied by $\rho_2$).
\end{lem}

\begin{proof}
For any unitary gauge transformation $u\in W^{2,p}(\U(E))$ and any point $[A,\Phi,a,\phi]\in T\sC_\ft^0$,
the computation
\begin{align*}
\exp([ u^*A,u^{-1}\Phi,u^{-1}a u,u^{-1}\phi])
&=
[u^*A + u^{-1}au,u^{-1}\Phi+u^{-1}\phi]
  \\
  &=
[u^*(A+a),u^{-1}(\Phi+\phi)]
\end{align*}
shows that $\exp$ has the asserted equivariance.
\end{proof}

We now describe the $S^1$ action on the restriction of $T\sC_\ft^0$ to $\iota_{\fs,\ft}(\sC_\fs^0)$.

\begin{lem}[Circle action on the restriction of the tangent bundle for configuration space to a submanifold of gauge-equivalence classes of split pairs]
\label{lem:FixedPointsOfS1ActionOnTangentSpaceToQuotientSpace}
Let $\fs=(\rho,W\otimes L_1)$ and $\ft=(\rho,W\otimes E)$ be a spin${}^c$ and rank-two spin${}^u$ structure, respectively, over a closed, oriented, smooth, Riemannian four-manifold $(X,g)$ with $c_1(\fs)\in\Red_0(\ft)$ as defined in \eqref{eq:DefineReduciblesEmbedded}. If $E$ admits a splitting $E=L_1\oplus L_2$  as an orthogonal direct sum of Hermitian line bundles, then the following hold:
\begin{enumerate}
\item
\label{item:FixedPointsOfS1ActionOnTangentSpaceToQuotientSpace_ClosedUnderAction}
The subspace $T\sC_\ft^0|_{\iota_{\fs,\ft}(\sC_\fs^0)}$ of $T\sC_\ft^0$ is closed under the $S^1$ actions \eqref{eq:S1ZActionOnTangentBundle} and \eqref{eq:S12ActionOnTangentBundle}.
\item
\label{item:FixedPointsOfS1ActionOnTangentSpaceToQuotientSpace_IdentifyFixedPoints}
The fixed-point set of the $S^1$ actions \eqref{eq:S1ActionOnTangentBundle} on $T\sC_\ft^0|_{\iota_{\fs,\ft}(\sC_\fs^0)}$ is the subspace given by $T\sC_\fs^0$ under the isomorphism \eqref{eq:DirectSumDecompositionOfTangentBundle}.
\item
\label{item:FixedPointsOfS1ActionOnTangentSpaceToQuotientSpace_ScalarMultOnNormalFiber}
The $S^1$ action \eqref{eq:S12ActionOnTangentBundle} on $\fN(\ft,\fs)$ is given by scalar multiplication on the fiber. 
\end{enumerate}
\end{lem}

\begin{proof}
If $\pi:T\sC_\ft^0\to\sC_\ft^0$ is the obvious projection map, then $\pi$ is equivariant with respect to the actions \eqref{eq:S1ActionOnTangentBundle} on $T\sC_\ft^0$ and the actions \eqref{eq:UnitaryS1ActionsRelation} on $\sC_\ft^0$ by inspection of the definitions. Because $\iota_{\fs,\ft}(\sC_\fs^0)$ is in the fixed-point set of the action \eqref{eq:S1ZActionOnQuotientSpace}, the equivariance of $\pi$ implies that the set
\[
\pi^{-1}\left( \iota_{\fs,\ft}(\sC_\fs^0)\right)
=
T\sC_\ft^0|_{\iota_{\fs,\ft}(\sC_\fs^0)}
\]
is closed under the $S^1$ actions \eqref{eq:S1ActionOnTangentBundle}. This establishes Item \eqref{item:FixedPointsOfS1ActionOnTangentSpaceToQuotientSpace_ClosedUnderAction}.

Because the actions \eqref{eq:S1ZActionOnTangentBundle} and \eqref{eq:S12ActionOnTangentBundle} are equal up to multiplicity, they have the same fixed-point sets in $T\sC_\ft^0|_{\iota_{\fs,\ft}(\sC_\fs^0)}$.
Gauge transformations in the image of the homomorphism $\rho_2$ of \eqref{eq:DefineUnitaryS12ActionsAtReducible} are in the stabilizer of pairs in the image of the embedding $\tilde\iota_{\fs,\ft}$ of \eqref{eq:DefnOfIota}. If
$(A,\Phi)=\tilde\iota_{\fs,\ft}(A_1,\Psi)$ and we use the decomposition \eqref{eq:SliceDecompositionAtReducible} to write  $(a,\phi)\in \Ker d_{A,\Phi}^{0,*}$ as
\[
(a,\phi)=\left((a_t,\phi_t),(a_n,\phi_2)\right)
\in
\Ker d_{A,\Phi}^{0,t,*} \oplus \Ker d_{A,\Phi}^{0,n,*},
\]
then the action of $\rho_2$ in \eqref{eq:S12ActionOnTangentBundle} on the fiber
$\pi^{-1}([\tilde\iota_{\fs,\ft}(A_1,\Psi)])$ is given by
\begin{equation}
\label{eq:S12ActionOnFiberOfFixedPoint}
\left( e^{i\theta}, (a,\phi)\right)
\mapsto
\left( \rho_2(e^{-i\theta}) a \rho_2(e^{i\theta}),\rho_2(e^{-i\theta}) \phi\right)
=
\left( (a_t,\phi_1), (\rho_2(e^{-i\theta})a_n,\rho_2(e^{-i\theta})\phi_2)\right).
\end{equation}
The fixed-point set of this action is $\Ker d_{A,\Phi}^{0,t,*}$, which is the fiber of $T\sC_\fs^0$ over
$[\tilde\iota_{\fs,\ft}(A_1,\Psi)]$.  Hence, $T\sC_\fs^0$ is the fixed-point set of the $S^1$ actions, verifying Item \eqref{item:FixedPointsOfS1ActionOnTangentSpaceToQuotientSpace_IdentifyFixedPoints}. The fiber of the vector bundle $\fN(\ft,\fs)$ over $\pi^{-1}([\tilde\iota_{\fs,\ft}(A_1,\Psi)])$ is $\Ker d_{A,\Phi}^{0,n,*}$ and the description of the $S^1$ action in \eqref{eq:S12ActionOnFiberOfFixedPoint} yields Item \eqref{item:FixedPointsOfS1ActionOnTangentSpaceToQuotientSpace_ScalarMultOnNormalFiber}.
\end{proof}

We use the description of the $S^1$ actions in Lemma \ref{lem:FixedPointsOfS1ActionOnTangentSpaceToQuotientSpace} to prove the following

\begin{lem}[Circle actions on the image of the sphere bundle in the configuration space of pairs]
\label{lem:S1ActionsOnImageOfSphereBundle}
Let $\fs=(\rho,W\otimes L_1)$ and $\ft=(\rho,W\otimes E)$ be a spin${}^c$ and rank-two spin${}^u$ structure, respectively, over a closed, oriented, smooth, Riemannian four-manifold $(X,g)$  with $c_1(\fs)\in\Red_0(\ft)$ as defined in \eqref{eq:DefineReduciblesEmbedded}. Assume that $E$ admits a splitting $E=L_1\oplus L_2$  as an orthogonal direct sum of Hermitian line bundles. For a compact subset $V\subset\sC_\fs^0$, let $\fN^\eps(\ft,\fs,V)\subset\fN(\ft,\fs)$ be the sphere bundle \eqref{eq:DefineRestrictedSphereNormalBundle} over $V$,
where $\eps\in (0,1]$ is the constant given in Lemma \ref{lem:FixedNormsInTubularNgh} such that $\fN^\eps(\ft,\fs,V)\subset\sO(\ft,\fs)$.  Then the $S^1$ actions in \eqref{eq:S1ActionsOnSpinuQuotient} are free on $\exp(\fN^\eps(\ft,\fs,V))$.
\end{lem}

\begin{proof}
The conclusion follows from the equivariance of $\exp$ given in Lemma \ref{lem:S1EquivarianceOfg} and the description of the $S^1$ action on $\fN(\ft,\fs)$ in Item \eqref{item:FixedPointsOfS1ActionOnTangentSpaceToQuotientSpace_ScalarMultOnNormalFiber} of Lemma \ref{lem:FixedPointsOfS1ActionOnTangentSpaceToQuotientSpace}.
\end{proof}

\section[Existence of non-split, non-zero-section, non-Abelian monopoles]{Existence of non-split, non-zero-section  non-Abelian monopoles for non-generic perturbations}
\label{sec:IrredMonopolesForNonGenericPerturbations}
The parameter space of perturbations for the non-Abelian monopole equations \eqref{eq:PerturbedSO3MonopoleEquations} is
\begin{multline}
\label{eq:PerturbationSpace}
\fP:=
\left\{(f,\tau,\vartheta)\in C^r(\GL(T^*X))\times C^r(\GL(\La^+))\times C^r(T^*X
  \otimes_\RR\CC): \right.
\\
\left.\|\tau-\id_{\La^+}\|_{C^0} < 1/64\right\},
\end{multline}
where $r\ge 2$ is an integer. The perturbations $\tau$ and $\vartheta$ in \eqref{eq:PerturbationSpace} appear explicitly in the non-Abelian monopole equations, \eqref{eq:PerturbedSO3MonopoleEquations}, while the automorphism $f\in C^r(\GL(T^*X))$ in \eqref{eq:PerturbationSpace} defines the metric and Clifford multiplication map appearing in \eqref{eq:PerturbedSO3MonopoleEquations}, by acting on the metric and on the Clifford multiplication map by pullback  (see Freed and Uhlenbeck \cite[Theorem 3.4]{FU} and Feehan \cite[Section 2.1, p. 916]{FeehanGenericMetric}).

\begin{defn}[Moduli spaces of Seiberg--Witten and non-Abelian monopoles for a given perturbation parameter]
\label{defn:ModuliSpacesWithPerturbationSpecified}
If $\fs$ is a spin${}^c$ structure, $\ft$ is a spin${}^u$ structure on $X$, and $p\in\fP$ is a perturbation as in \eqref{eq:PerturbationSpace}, then we write
\begin{subequations}
\label{eq:ModuliSpacesWithPerturbationSpecified}
\begin{align}
M_\fs(p)&:=\{[A_1,\Psi]\in\sC_\fs: \text{$(A_1,\Psi)$ satisfies \eqref{eq:SeibergWitten} for the perturbation $p$}\},
\\
\sM_\ft(p)&:=\{[A,\Phi]\in\sC_\ft: \text{$(A,\Phi)$ satisfies \eqref{eq:PerturbedSO3MonopoleEquations} for the perturbation $p$}\},
\end{align}
\end{subequations}
for the Seiberg--Witten and non-Abelian monopole moduli spaces defined by the perturbation $p\in \fP$.
\end{defn}

\begin{defn}[Generic perturbations]
\label{defn:GenericPerturbations}
We call a perturbation $p\in\fP$ \emph{generic} for a spin${}^u$ structure $\ft$ if the map  defining the non-Abelian monopole equations \eqref{eq:PerturbedSO3MonopoleEquations} with that perturbation $p\in\fP$ vanishes transversely on $\sC_\ft^{*,0}$. We call a perturbation $p\in\fP$ \emph{generic} for a spin${}^c$ structure $\fs$ if the map  defining the Seiberg--Witten equations \eqref{eq:SeibergWitten} with that perturbation $p\in\fP$ vanishes transversely on $\sC_\fs^0$.
\end{defn}

In the forthcoming proposition, we give a criterion for the existence of a non-split non-Abelian monopole for any perturbation.  While the proof of Proposition \ref{prop:NonEmptyModuliSpace} is straightforward for a generic perturbation, we will need to work with non-generic perturbations such as perturbations $p$ where the metric is K\"ahler.

\begin{prop}[Existence of non-split non-Abelian monopoles for non-generic perturbations]
\label{prop:NonEmptyModuliSpace}
Let $\fs$ and $\ft$ be a spin${}^c$ and rank-two spin${}^u$ structure on a standard four-manifold $X$ (as in Definition \ref{defn:Standard}) with $c_1(\fs)\in\Red_0(\ft)$ as defined in \eqref{eq:DefineReduciblesEmbedded}. Assume that there are no zero-section pairs in $M_\fs(p_\infty)$, where $p_\infty\in\fP$, the space of perturbations defined in \eqref{eq:PerturbationSpace}. If $\SW_X(\fs)\neq 0$ and
\[
  \expdim \sM_\ft > 0,
\]
then $\sM_\ft^{*,0}(p_\infty) := \sM_\ft(p_\infty)\cap\sC_\ft^{*,0}$ is non-empty.
\end{prop}

We first prove Proposition \ref{prop:NonEmptyModuliSpace} in the special case that $p_\infty$ is a generic parameter.  When the parameter $p$ is generic for $\ft$ and for $\fs$ as in Definition \ref{defn:GenericPerturbations} and $M_\fs(p)\subset\sC_\fs^0$, then $\iota_{\fs,\ft}(M_\fs(p))$ is a smooth submanifold of $\sC_\ft^0$ with a tubular neighborhood.  A sphere bundle in this tubular neighborhood defines a link of $\iota_{\fs,\ft}(M_\fs(p))$ in $\sC_\ft^0$. The link $\bL_{\ft,\fs}(p)$ of $M_\fs(p)$ in $\sM_\ft(p)$ is defined in Feehan and Leness \cite[Equation (3.47)]{FL2a} as the $S^1$-quotient (by the action \eqref{eq:S1ZActionOnQuotientSpace}) of the intersection of $\sM_\ft(p)$ with the link of $\iota_{\fs,\ft}(M_\fs(p))$ in $\sC_\ft^0$.  In addition, the inclusion $\bL_{\ft,\fs}(p)\subset\sM_\ft^{*,0}(p)/S^1$ follows from showing that $\bL_{\ft,\fs}(p)$ contains no fixed points of the $S^1$ action and applying Proposition \ref{prop:FixedPointsOfS1ActionOnSpinuQuotientSpace}.
We recall the definition  \cite[Equation (3.47), p. 112]{FL2a} of the link $\bL_{\ft,\fs}(p)$ of $M_\fs(p)$ in $\sM_\ft(p)$ as
\begin{equation}
\label{eq:IdentifyLinkWithSphereBundle}
\bL_{\ft,\fs}(p):= \left(\exp(\fN^\eps(\ft,\fs,M_\fs(p)))\cap \sM_\ft(p)\right)/S^1,
\end{equation}
where $\fN^\eps(\ft,\fs,M_\fs(p))$ is the sphere bundle defined in \eqref{eq:DefineRestrictedSphereNormalBundle}. We then have the

\begin{lem}[Existence of non-split, non-zero-section non-Abelian monopoles for generic perturbations]
\label{lem:NonEmptyModuliSpaceForGenericParameter}
Continue the assumptions and notation of Proposition \ref{prop:NonEmptyModuliSpace} and in addition assume that the perturbation $p_\infty\in\fP$ is generic for $\ft$ and for $\fs$ in the sense of Definition \ref{defn:GenericPerturbations}. Then the link $\bL_{\ft,\fs}(p_\infty)$ of $M_\fs(p_\infty)$ in $\sM_{\ft}(p_\infty)/S^1$ is non-empty and so $\sM_\ft^{*,0}(p_\infty)$ is non-empty.
\end{lem}

\begin{proof}
Let $p_\infty\in\fP$ be generic and let $\bL_{\ft,\fs}(p_\infty)\subset\sM_\ft^{*,0}(p_\infty)$ be the link of $M_\fs(p_\infty)$ in $\sM_\ft(p_\infty)/S^1$ defined in \eqref{eq:IdentifyLinkWithSphereBundle}.
Because $X$ is standard in the sense of Definition \ref{defn:Standard}, then $H^1(X;\ZZ)$ is trivial and so for all $\alpha,\alpha'\in H^1(X;\ZZ)$, the condition $\alpha\smile\alpha'=0$ in Feehan and Leness \cite[Theorem 4.13, p. 194]{FL2b} is trivially satisfied. We can thus apply the link-pairing formula provided by  \cite[Theorem 4.13]{FL2b} to $\bL_{\ft,\fs}(p_\infty)$. Because $\expdim\sM_{\ft}>0$ and $\SW_{X}(\fs)\neq 0$, then \cite[Theorem 4.13, p. 194]{FL2b} implies that the fundamental class of $\bL_{\ft,\fs}(p_\infty)$ defines a non-zero element of $H_\bullet(\sC_{\ft}^{*,0}/S^1;\RR)$ and thus $\bL_{\ft,\fs}(p_\infty)$ is non-empty.  Because
$\bL_{\ft,\fs}(p_\infty)\subset\sM_\ft^{*,0}(p_\infty)$ is non-empty, then $\sM_\ft^{*,0}(p_\infty)$ is also non-empty.
\end{proof}

Next, we have the

\begin{lem}
\label{lem:NoZeroSectionNeighborhood}
Let $\fs$ be a spin${}^c$ structure on a standard four-manifold $X$. If $M_\fs(p_\infty)$ contains no zero-section pairs, then there is an open neighborhood $U_\infty$ of $p_\infty$ in $\fP$ such that for all $p\in U_\infty$, the moduli space $M_\fs(p)$ contains no zero-section pairs.
\end{lem}

\begin{proof}
If the conclusion were false, then there would be a sequence of perturbations $\{p_k\}_{k\in\NN}\subset\fP$ converging to $p_\infty$ and, for each $k\in\NN$, a zero-section pair $[A_1(k),0]$ with $[A_1(k),0]\in M_\fs(p_i)$. By Theorem \ref{thm:Uhlenbeck_convergence_sequence_SW_monopoles_parameters}, the sequence $\{[A_1(k),0]\}_{k\in\NN}$ would have a subsequence converging to $[A_1(\infty),0]\in M_\fs(p_\infty)$, a contradiction to the hypothesis that $M_\fs(p_\infty)$ contains no zero-section pairs.
\end{proof}

We are now ready to complete the 

\begin{proof}[Proof of Proposition \ref{prop:NonEmptyModuliSpace}]
Let $U_\infty$ be the open neighborhood of $p_\infty$ in $\fP$ given by Lemma \ref{lem:NoZeroSectionNeighborhood} with the property that $M_\fs(p)\subset\sC_\fs^0$ for $p\in U_\infty$.  The set of generic parameters is the complement of a first category set by Feehan \cite[Theorem 1.3]{FeehanGenericMetric} and is thus dense in $\fP$.  Hence, we can find a sequence of generic parameters $\{p_k\}_{k\in\NN}\subset U_\infty$ converging to $p_\infty$. By Theorem \ref{thm:Uhlenbeck_convergence_sequence_SW_monopoles_parameters}, the space $V := \cup_{k\in\NN}M_\fs(p_k)\subset\sC_\fs^0$ is compact. Let $\fN^\eps(\ft,\fs,V)$ be the sphere bundle defined in 
\eqref{eq:DefineRestrictedSphereNormalBundle}. By Lemma \ref{lem:NonEmptyModuliSpaceForGenericParameter}, the link $\bL_{\ft,\fs}(p_k)$ defined in \eqref{eq:IdentifyLinkWithSphereBundle} is non-empty for all $k\in\NN$. Hence, there is a sequence of points $\{[A(k),\Phi(k)]\}_{k\in\NN}$ such that the image of $[A(k),\Phi(k)]$ in the quotient $\sC_\ft^0/S^1$ by the $S^1$ action \eqref{eq:S1ZActionOnQuotientSpace} lies in $\bL_{\ft,\fs}(p_k)$.  By \eqref{eq:IdentifyLinkWithSphereBundle} we can then write
\begin{equation}
\label{eq:LinkPointsInNormalBundleImage}
[A(k),\Phi(k)]=\exp([A_1(k),\Psi(k),a_n(k),\phi_2(k)]), \quad\text{for all } k \in \NN,
\end{equation}
where $[A_1(k),\Psi(k)]\in M_\fs(p_k)$ and $[A_1(k),\Psi(k),a_n(k),\phi_2(k)]\in \fN^\eps(\ft,\fs,M_\fs(p_k))$.
By Theorem \ref{thm:Uhlenbeck_convergence_sequence_non-Abelian_monopoles_parameters}, the sequence $\{[A(k),\Phi(k)]\}_{k\in\NN}$ has a subsequence which converges to a point in $\sI\!\!\sM_{\ft}(p_\infty)$.
After relabelling and passing to a subsequence, we can assume that $\{[A(k),\Phi(k)]\}_{k\in\NN}$ converges to a point in $\sI\!\!\sM_{\ft}(p_\infty)$. By Theorem \ref{thm:Uhlenbeck_convergence_sequence_SW_monopoles_parameters},
after relabelling and passing to a subsequence. the sequence $\{[A_1(k),\Psi(k)]\}_{k\in\NN}$ converges to $[A_1(\infty),\Psi(\infty)]\in M_\fs(p_\infty)$. Because the sequence $\{[A_1(k),\Psi(k)]\}_{k\in\NN}$ converges to $[A_1(\infty),\Psi(\infty)]$, the $L^p$ norms of the curvatures of $A_1(k)$ are uniformly bounded. If we write $(A(k)',\Phi(k)')=\tilde\iota_{\fs,\ft}(A_1(k),\Psi(k))$, then continuity of $\tilde\iota_{\fs,\ft}$ implies that the sequence $\{(A(k)',\Phi(k)')\}_{k\in\NN}$ also converges and so the $L^p$ norms of the curvatures of $A(k)'$ are uniformly bounded. The difference between the $L^p$ norms of the curvatures of $A(k)$ and $A(k)'$ is bounded by the $W^{1,p}_{A(k)}$-norm of $(a_{21}(k),\phi_2(k))$, which equals $\eps$ by the definition of $\fN^\eps(\ft,\fs,M_\fs(p_k))$ in
\eqref{eq:DefineRestrictedSphereNormalBundle}. Hence, the $L^p$ norms of the curvatures of $A(k)$ are bounded uniformly with respect to $k$, and therefore the  sequence $\{[A(k),\Phi(k)]\}_{k\in\NN}$ does not exhibit Uhlenbeck bubbling (our definition of $\sC_\ft$ required that $p \in (2,\infty)$) and its limit is in $\sM_\ft(p_\infty)\subset\sC_\ft$. Because $\{[A(k),\Phi(k)]\}_{k\in\NN}$ is contained in the closed subspace $\exp(\fN^\eps(\ft,\fs,V))$, the limit of the sequence must be in
\[
\exp(\fN^\eps(\ft,\fs,V))\cap\sM_\ft.
\]
By Lemma \ref{lem:S1ActionsOnImageOfSphereBundle}, the subspace $\exp(\fN^\eps(\ft,\fs,V))$ does not contain any fixed points of the $S^1$ action \eqref{eq:S12ActionsOnSpinuQuotient}. Proposition \ref{prop:FixedPointsOfS1ActionOnSpinuQuotientSpace} thus implies that
\[
\exp(\fN^\eps(\ft,\fs,V))\cap\sM_\ft(p_\infty)
\subset
\sM_\ft(p_\infty)\cap\sC_\ft^{*,0}.
\]
Therefore, the limit point $[A(\infty),\Phi(\infty)]$ lies in $\sM_\ft(p_\infty)\cap\sC_\ft^{*,0}$, and this completes the proof of Proposition \ref{prop:NonEmptyModuliSpace}.
\end{proof}

\chapter{Feasibility}
\label{chap:Feasibility}
In this chapter, we prove Theorem \ref{mainthm:ExistenceOfSpinuForFlow} and Corollary \ref{maincor:MorseIndexForFeasibilitySpinuStructure}. In Section \ref{sec:Intermediate_calculations_feasibility}, we construct the \spinu structure used in these proofs, identifying the  properties of the characteristic classes of the \spinu structure used in these proofs in Lemma \ref{lem:SpinuStrRequirements} and then showing that \spinu structures with such characteristic classes exist on the smooth blow-up of a standard four-manifold.  We apply these results to give the proofs of Theorem \ref{mainthm:ExistenceOfSpinuForFlow} and Corollary \ref{maincor:MorseIndexForFeasibilitySpinuStructure} in Section \ref{sec:Proofs_main_results_feasibility}.

\section{Intermediate calculations for feasibility}
\label{sec:Intermediate_calculations_feasibility}
In this section, we assemble the results on the existence of \spinu structures needed for the proof of Theorem \ref{mainthm:ExistenceOfSpinuForFlow} and Corollary \ref{maincor:MorseIndexForFeasibilitySpinuStructure}.  To motivate the computations in this chapter, we begin with Lemma \ref{lem:SpinuStrRequirements} by identifying cohomological conditions on the characteristic classes $p_1(\tilde\ft)$ and $c_1(\tilde\ft)$ which ensure that the moduli space $\sM_{\tilde\ft}$ of non-Abelian monopoles defined by a \spinu structure $\tilde\ft$ satisfies the conclusions \eqref{item:ExistenceOfSpinuForFlowNonEmpty} and \eqref{item:ExistenceOfSpinuForFlowp1} of Theorem \ref{mainthm:ExistenceOfSpinuForFlow}.  Lemma \ref{lem:ExistenceOfSpinu} gives conditions on cohomology classes which ensure they are the characteristic classes $p_1(\tilde\ft)$ and $c_1(\tilde\ft)$ of a \spinu structure $\tilde\ft$. In Lemmas \ref{lem:SemiAbundanceCondition}, \ref{lem:Existence}, and \ref{lem:SpinuCharSatisfyingCriteria}, we construct cohomology classes satisfying the conditions of Lemmas \ref{lem:SpinuStrRequirements} and \ref{lem:ExistenceOfSpinu}, thus constructing a \spinu structure whose associated moduli space of non-Abelian monopoles satisfies the conclusions of Theorem \ref{mainthm:ExistenceOfSpinuForFlow}. We apply our formula \eqref{eq:MorseIndexAtReduciblesOnKahler} for the virtual Morse--Bott index given by Theorem \ref{mainthm:MorseIndexAtReduciblesOnKahler} (which we prove in Chapter \ref{chap:VirtualMorseIndexComputation}) to prove Corollary \ref{maincor:MorseIndexForFeasibilitySpinuStructure} for this \spinu structure.

\begin{lem}
\label{lem:SpinuStrRequirements}
Let $X$ be a standard four-manifold as in Definition \ref{defn:Standard}  and let $\widetilde X=X\#\overline{\CC\PP}^2$ denote the smooth blow-up of $X$ at a point. Assume that there are an integer $b$, a rank-two spin${}^u$ structure $\tilde\ft$, and a spin${}^c$ structure $\tilde\fs$ over $\widetilde X$ that obey the following conditions:
\begin{subequations}
\label{eq:SpinuStrRequirements_c1Squared_equality_p1_bound_LaDotK_Requirement}
\begin{align}
\label{eq:SpinuStrRequirements_c1Squared_equality}
c_1(\tilde\fs)^2&=c_1(\tilde X)^2,
\\
\label{eq:SpinuStrRequirements_p1_bound}
p_1(\tilde\ft)&= c_1(\widetilde X)^2-12\chi_h(\widetilde X)+b,
\\
\label{eq:SpinuStrRequirements_LaDotK_Requirement}
c_1(\tilde \ft)\cdot c_1(\tilde\fs)&=\frac{1}{2}c_1(\tilde \ft)^2+6\chi_h(\widetilde X)-\frac{1}{2}b,
\end{align}
\end{subequations}
where $c_1(\widetilde X)^2$ and $\chi_h(\widetilde X)$ are as in \eqref{eq:DefineTopCharNumbersOf4Manifold}, $c_1(\tilde\fs)$ is as in \eqref{eq:DefineChernClassOfSpinc}, and $p_1(\tilde\ft)$ and $c_1(\tilde\ft)$ are as in \eqref{eq:SpinUCharacteristics}. Then the following hold:
\begin{enumerate}
\item\label{item:SpinuStrRequirements_SWModuliIncluded}
$\sM_{\tilde\ft}$ contains the image of the moduli space $M_{\tilde \fs}$ of Seiberg--Witten monopoles \eqref{eq:Moduli_space_Seiberg-Witten_monopoles} under the map $\iota_{\tilde\fs,\tilde\ft}$ in \eqref{eq:DefnOfIotaOnQuotient}.

\item\label{item:SpinuStrRequirements_Basic_Lower_Bound}
If $E$ is a complex rank-two vector bundle over $\widetilde X$ with $\su(E)\cong \fg_{\tilde\ft}$, where $\fg_{\tilde\ft}$ is as in \eqref{eq:SpinAssociatedBundles}, then
\[
p_1(\su(E))=c_1(\widetilde X)^2-12\chi_h(\widetilde X)+b,
\]
and so, if $b\ge 1$, it obeys the fundamental lower bound \eqref{eq:p1_lower_bound_blow-up}.

\item\label{item:SpinuStrRequirements_ASDModuliIncluded}
If $\tilde g$ is a Riemannian metric on $\widetilde X$ and $M_\kappa^w(\widetilde X,\tilde g)$  is the moduli space \eqref{eq:ASDModuliSpace} of anti-self-dual $\SO(3)$ connections identified with the stratum of $\sM_{\tilde\ft}$ represented by zero-section pairs, where $w\equiv w_2(\tilde\ft)\pmod 2$ and $\ka=-(1/4)p_1(\tilde\ft)$, then the expected dimension of $M_\kappa^w(\widetilde X,\tilde g)$ obeys
\begin{equation}
\label{eq:ASDDimInTermsOfBMY}
\expdim M_\kappa^w(\widetilde X,\tilde g) = 2\left(-c_1(\widetilde X)^2+9\chi_h(\widetilde X)\right)-2b
\end{equation}
and so, if $b\geq 1$, then
\begin{equation}
  \label{eq:ASDDimInTermsOfBMY_lower_bound}
  \expdim M_\kappa^w(\widetilde X,\tilde g) \leq 2\left(-c_1(X)^2+9\chi_h(X)\right).
\end{equation}

\item\label{item:SpinuStrRequirementsEmptyZeroSection}
If $M_{\tilde \fs}$ does not contain the gauge-equivalence class of a zero-section pair and, in addition to the hypotheses
\eqref{eq:SpinuStrRequirements_c1Squared_equality_p1_bound_LaDotK_Requirement}, the spin${}^u$ structure $\tilde\ft$ and spin${}^c$ structure $\tilde\fs$ satisfy
\begin{subequations}
\label{eq:SpinuStrRequirements_Good_NonZeroSW_LaSquared_requirement}
\begin{gather}
\label{eq:SpinuStrRequirementsNonZeroSW}
\SW_{\widetilde X}(\tilde\fs) \neq 0,
  \\
  \label{eq:SpinuStrRequirements_LaSquared_requirement}
c_1(\tilde\ft)^2 > 4\left( c_1(\widetilde X)^2-8\chi_h(\widetilde X)\right) +3b,
\end{gather}
\end{subequations}
then the moduli space $\sM^{*,0}_{\tilde\ft}$ in \eqref{eq:PerturbedSO3MonopoleEquations} of non-split, non-zero-section
non-Abelian monopoles is not empty.
\end{enumerate}
\end{lem}

\begin{rmk}[On the choice of the integer $b$ in Lemma \ref{lem:SpinuStrRequirements}]
\label{rmk:ProblemWithb}
As discussed in Sections \ref{sec:Morse_theory_existence_anti-self-dual_connections} and \ref{sec:Morse_theory_moduli_space_PU2_monopoles_Kaehler}, the goal of the program that motivates this monograph is to prove that for a \spinu structure $\tilde\ft$ satisfying \eqref{eq:SpinuStrRequirements_c1Squared_equality_p1_bound_LaDotK_Requirement}
and \eqref{eq:SpinuStrRequirements_Good_NonZeroSW_LaSquared_requirement}, the moduli space  $M_\ka^w(\widetilde X,\tilde g)$ of projectively anti-self-dual connections is non-empty, thus implying that $\dim  M_\ka^w(\widetilde X,\tilde g)\ge 0$. For $b\in\ZZ$ and $p_1(\tilde\ft)$ as in \eqref{eq:SpinuStrRequirements_p1_bound}, if
\[
  \expdim M_\ka^w(\widetilde X,\tilde g) \geq 0,
\]
then from \eqref{eq:ASDDimInTermsOfBMY} one can  derive the bound $c_1(\widetilde X)^2\le 9\chi_h(X)-b$ or, using $c_1(\widetilde X)^2=c_1(X)^2-1$,
\begin{equation}
\label{eq:BMYwithb}
c_1(X)^2\le 9\chi_h(X)-(b-1).
\end{equation}
The four-manifolds constructed in the literature discussed in Section \ref{sec:Morse_theory_existence_anti-self-dual_connections}  will violate the inequality \eqref{eq:BMYwithb} if $b$ is sufficiently large. However, the only upper bound on $b$ provided by this monograph is $b< 10\chi_h(X)$ from the forthcoming \eqref{eq:bGeneralTypeAssumption} in Lemma \ref{lem:IndexInAlgebraicCase}, but this upper bound is too weak to prevent \eqref{eq:BMYwithb} from contradicting known examples.

However, preliminary computations suggest that stronger upper bounds on the possible values of $b$ arise from a constraint that virtual Morse--Bott indices are also positive for ideal Seiberg--Witten points in lower levels of the Uhlenbeck compactification of $\sM_\ft$.
\end{rmk}

\begin{proof}
We observe that, by Lemma \ref{lem:SetOfReducibles}, to prove Item \eqref{item:SpinuStrRequirements_SWModuliIncluded} we only need to show that $c_1(\tilde\fs)\in\Red_0(\tilde\ft)$, as defined in \eqref{eq:DefineReduciblesEmbedded}, and thus we claim that
\begin{equation}
\label{eq:CohomologicalLevelZeroCondition}
  (c_1(\tilde\fs)-c_1(\tilde\ft))^2 = p_1(\tilde\ft).
\end{equation}
To prove claim \eqref{eq:CohomologicalLevelZeroCondition}, we compute:
\begin{align*}
(c_1(\tilde\fs)-c_1(\tilde\ft))^2
&=
c_1(\tilde\fs)^2 - 2c_1(\tilde\fs)\cdot c_1(\tilde\ft) + c_1(\tilde\ft)^2
\\
&=
c_1(\tilde X)^2 -c_1(\tilde\ft)^2-12\chi_h(\widetilde X)+b + c_1(\tilde\ft)^2
\quad\text{(by \eqref{eq:SpinuStrRequirements_c1Squared_equality} and \eqref{eq:SpinuStrRequirements_LaDotK_Requirement})}
\\
&=   
c_1(\tilde X)^2-12\chi_h(\widetilde X)+b
\\
&=
p_1(\tilde \ft)
\quad\text{(by \eqref{eq:SpinuStrRequirements_p1_bound}).}
\end{align*}
This proves claim \eqref{eq:CohomologicalLevelZeroCondition} and thus Item \eqref{item:SpinuStrRequirements_SWModuliIncluded}.

Item \eqref{item:SpinuStrRequirements_Basic_Lower_Bound} follows from the definition \eqref{eq:SpinUCharacteristics} of $p_1(\tilde \ft)=p_1(\fg_{\tilde \ft})$ and from the expression \eqref{eq:SpinuStrRequirements_p1_bound} for $p_1(\tilde \ft)$.

Consider Item \eqref{item:SpinuStrRequirements_ASDModuliIncluded}. We compute:
\begin{align*}
\expdim M_\ka^w(\widetilde X,\tilde g)
  &=
    -2p_1(\tilde \ft)-6\chi_h(\widetilde X)
    \quad\text{(by substituting \eqref{eq:SpinUCharacteristics} and \eqref{eq:SpinAssociatedBundles} into \eqref{eq:Expected_dimension_moduli_space_ASD_connections})}
\\
&=2\left( -c_1(\widetilde X)^2 +9\chi_h(\widetilde X)\right)-2b
\quad\text{(by \eqref{eq:SpinuStrRequirements_p1_bound})},
\end{align*}
and this proves the equality \eqref{eq:ASDDimInTermsOfBMY}. Because $c_1(\tilde X)^2=c_1(X)^2-1$ and $\chi_h(\widetilde X)=\chi_h(X)$, we observe that \eqref{eq:ASDDimInTermsOfBMY} yields
\begin{align*}
  \expdim M_\ka^w(\widetilde X,\tilde g) &= 2\left( -c_1(X)^2+1 +9\chi_h(X)\right)-2b
  \\
  &= 2\left( -c_1(X)^2 +9\chi_h(X)\right)-2(b-1),
\end{align*}
and so the inequality \eqref{eq:ASDDimInTermsOfBMY_lower_bound} follows if $b\geq 1$. This proves
Item \eqref{item:SpinuStrRequirements_ASDModuliIncluded}.

Finally, we consider Item \eqref{item:SpinuStrRequirementsEmptyZeroSection}. The definitions of $c_1(X)^2$ and $\chi_h(X)$ in \eqref{eq:DefineTopCharNumbersOf4Manifold} yield the equalities
\begin{equation}
\label{eq:ReverseDefinitionOfTopCharNumbersOf4Manifolds}
e(X)=12\chi_h(X)-c_1(X)^2\quad\text{and}\quad \si(X)=c_1(X)^2-8\chi_h(X),
\end{equation}
where $e(X)$ is the Euler characteristic of $X$ and $\si(X)$ the signature of $X$. We obtain
\begin{align*}
2\expdim\sM_{\tilde\ft}
&=
2d_a(\tilde\ft)+4n_a(\tilde\ft) \quad\text{(by \eqref{eq:Transv})}
\\
&=
-3p_1(\tilde\ft) +c_1(\tilde\ft)^2 -3e(X) - 4\si(X)
\quad\text{(by \eqref{eq:Transv})}
\\
&=-3 p_1(\tilde\ft) +c_1(\tilde\ft)^2 -  c_1(\widetilde X)^2 - 4 \chi_h(\widetilde X)
\quad\text{(by \eqref{eq:ReverseDefinitionOfTopCharNumbersOf4Manifolds})}
\\
&=c_1(\tilde\ft)^2-4c_1(\widetilde X)^2+32\chi_h(\widetilde X)-3b
\quad\text{(by \eqref{eq:SpinuStrRequirements_p1_bound}).}
\end{align*}
Hence, if $c_1(\tilde\ft)^2$ obeys \eqref{eq:SpinuStrRequirements_LaSquared_requirement}, then $\expdim\sM_{\tilde\ft}>0$.
If $M_{\tilde\fs}$ contains no zero-section points, then \ref{eq:SpinuStrRequirementsNonZeroSW} and Proposition \ref{prop:NonEmptyModuliSpace} ensure that $\sM_\ft^{*,0}$ is non-empty, which proves the conclusion in Item
\eqref{item:SpinuStrRequirementsEmptyZeroSection}. This completes the proof of Lemma \ref{lem:SpinuStrRequirements}.
\end{proof}

If $X$ has Seiberg--Witten simple type, then in Lemma \ref{lem:SpinuStrRequirements}, condition \eqref{eq:SpinuStrRequirementsNonZeroSW} implies condition \eqref{eq:SpinuStrRequirements_c1Squared_equality}. To establish the existence of \spinu and \spinc structures whose characteristic classes satisfy conditions
\eqref{eq:SpinuStrRequirements_p1_bound}, \eqref{eq:SpinuStrRequirements_LaDotK_Requirement}, and \eqref{eq:SpinuStrRequirements_LaSquared_requirement} in Lemma \ref{lem:SpinuStrRequirements}, we will use the following criterion for the existence of \spinu structure with given characteristic classes.

\begin{lem}[Existence of \spinu structures]
\label{lem:ExistenceOfSpinu}
(See Feehan and Leness \cite[Lemma 2.1, p. 7]{FL8}.)
Let $Y$ be a standard four-manifold as in Definition \ref{defn:Standard}. Given $\wp\in H^4(Y;\ZZ)$ and $\La\in H^2(Y;\ZZ)$ and $\mathfrak{w}\in H^2(Y;\ZZ/2\ZZ)$, there exists a rank-two spin${}^u$ structure $\ft$ on $Y$ with $p_1(\ft)=\wp$ and $c_1(\ft)=\La$ and $w_2(\ft)=\mathfrak{w}$ if and only if the following hold:
\begin{enumerate}
\item
\label{item:ExistenceOfSpinu1}
There is a class $w\in H^2(Y;\ZZ)$ with $\mathfrak{w}=w\pmod 2$; and
\item
\label{item:ExistenceOfSpinu2}
$\La\equiv \mathfrak{w}+w_2(X)\pmod 2$; and
\item
\label{item:ExistenceOfSpinu3}
$\wp\equiv \mathfrak{w}^2\pmod 4$.
\end{enumerate}
\end{lem}

Next, we adapt the methods and results of Feehan and Leness \cite[Appendix A]{FL2a} to establish the existence of cohomology classes on $X$ which we will use to construct characteristic classes satisfying the conditions of Lemmas \ref{lem:ExistenceOfSpinu} and \ref{lem:SpinuStrRequirements}. Recall from Gompf and Stipsicz, \cite[Definition 1.2.8 (c), p. 10]{GompfStipsicz} that a class $\ka\in H^2(X;\ZZ)$ is \emph{characteristic} \label{page:characteristic} if $Q_X(\alpha,\alpha)\equiv Q_X(\alpha,\ka)\pmod 2$ for all $\alpha\in H^2(X:\ZZ)$.

\begin{lem}[Existence of hyperbolic summands]
\label{lem:SemiAbundanceCondition}
Let $X$ be a standard four-manifold as in Definition \ref{defn:Standard} with $b^-(X)\ge 2$. Let $\ka\in H^2(X;\ZZ)$ be characteristic with $\ka^2=c_1(X)^2$. Then there is a hyperbolic summand of $H^2(X;\ZZ)$ which is $Q_X$-orthogonal to $\ka$, that is, there are classes $f_1,f_2\in H^2(X;\ZZ)$ obeying the following properties:
\begin{equation}
\label{eq:HyperbolicSummandConditions}
f_1\cdot f_2 = 1 \quad\text{and}\quad f_i\cdot \kappa = 0 = f_i\cdot f_i,
\quad \text{for } i=1,2.
\end{equation}
If $\ka$ is not a torsion class, then the cohomology classes $\ka$, $f_1$, $f_2$ are linearly independent over $\ZZ$ in the sense that if $r_1,r_2,r_3\in\ZZ$ satisfy $r_1f_1+r_2f_2+r_3\ka=0$, then $r_1=r_2=cr_3=0$.
\end{lem}

\begin{proof}
The linear independence of $\ka$, $f_1$, and $f_2$ follows from \eqref{eq:HyperbolicSummandConditions}. To see this, observe that if $r_1f_1+r_2f_2+r_3\ka=0$, then by \eqref{eq:HyperbolicSummandConditions} we obtain
\begin{align*}
0&= f_1\cdot\left( r_1f_1+r_2f_2+r_3\ka=0\right)= r_2,
\\
0&= f_2\cdot\left( r_1f_1+r_2f_2+r_3\ka=0\right)= r_1.
\end{align*}
Because $r_1=r_2=0$, we have $r_3\ka=0$. Because $\ka$ is not torsion, the equality $r_3\ka=0$ implies that $r_3=0$ and this completes the proof of the linear independence of $\ka$, $f_1$, and $f_2$.

We now prove \eqref{eq:HyperbolicSummandConditions}. We first treat the case when $Q_X$ is even and then the case when $Q_X$ is odd.

\begin{case}[$Q_X$ even]\label{case:Q_XEven}
From Gompf and Stipsicz \cite[Theorem 1.2.21, p. 14]{GompfStipsicz}, we have $Q_X=mE_8\oplus (2n+1)H$ when $Q_X$ is even, where $E_8$ is the rank-eight, positive-definite form defined following \cite[Lemma 1.2.20, p. 13]{GompfStipsicz} and $H$ is a hyperbolic summand, so
\[
  H \cong \begin{pmatrix} 0 & 1 \\ 1 & 0\end{pmatrix}.
\]
Note that while $n$ is a non-negative integer, $m$ can be negative in which case $mE_8$ means $|m|$ summands of the negative-definite form $-E_8$. We claim that $n\ge 1$. Because $\si(H)=0$, then $\si(X)=8m$ as $E8$ has real rank eight.  Next, we show that $n\ge 1$.
\begin{subcase}
If $\si(X)\le 0$, then $m\le 0$ and $2n+1=b^+(X)$ and so the hypothesis that $b^+(X)\ge 3$ implies that $n\ge 1$.
\end{subcase}
\begin{subcase}
If $\si(X)> 0$, then $m=\si(X)/8$ by Gompf and Stipsicz \cite[Theorem 1.2.21, p. 14]{GompfStipsicz} and $2n+1=b^-(X)$.  Because $b^-(X)\ge 2$ by hypothesis, we must again have $n\ge 1$, as claimed.
\end{subcase}

Because $n\ge 1$, $Q_X$ contains a summand $R=H\oplus H$. Let $h_1,h_2$ and $h_1',h_2'$ be bases for the subspace $V_R\subset H^2(X;\ZZ)$ on which $R$ is non-zero. Thus, $Q_X(h_1,h_2)=Q_X(h_1',h_2')=1$ and all other pairings of these four elements vanish. As in the proof of Feehan and Leness \cite[Lemma A.2, p. 126]{FL2a}, we shall write $\ka_R$ for the component of $\ka$ in $V_R$.  Because $\ka_R$ is characteristic for $(V_R,Q_X|_R)$, then $Q_X(\ka_R,\ka_R)=2d$ for some $d\in\ZZ$. Define $v := h_1+dh_2$, so $Q_X(v,v)=Q_X(\ka_R,\ka_R)$. By a result of Wall \cite[Theorem 1]{WallUnimodQuadForms}, there is an isometry $A$ of $(V_R,Q_X|_R)$ that maps $v$ to $\ka_R$. Let $V_{R'}\subset V_R$ be the hyperbolic summand spanned by $h_1'$ and $h_2'$. Because $V_{R'}$ is $Q_X$-orthogonal  to $v$ and $A$ is an isometry, $A(V_{R'})$ is a hyperbolic summand that is $Q_X$-orthogonal to $A(v)=\ka_R$. Because $A(V_{R'})\subset V_R$ and $A(V_{R'})$ is $Q_X$-orthogonal to $\ka_R$, which is the component of $\ka$ in $V_R$, then $A(V_{R'})$ is $Q_X$-orthogonal to $\ka$. This completes the treatment of Case \eqref{case:Q_XEven}.
\end{case}

\begin{case}[$Q_X$ odd]
\label{case:Q_XOdd}
By Feehan and Leness \cite[Lemmas A.4, p. 126]{FL2a}, when $Q_X$ is odd the $Q_X$-orthogonal complement of $\ka$ contains a hyperbolic summand if $b^+(X)\ge 5$ and $b^-(X)\ge 3$.  By the assumption that $\ka^2=c_1(X)^2$, the same result holds if $Q_X$ is odd with $b^+(X)=3$ and $b^-(X)\ge 2$ by \cite[Lemmas A.6, p. 129]{FL2a}. Thus, to complete the proof of Lemma \ref{lem:SemiAbundanceCondition} in the case when $Q_X$ is odd we only need to consider the case $b^+(X)\ge 5$ and $b^-(X)=2$.  In the latter case, we can write
\begin{equation}
\label{eq:OddIntersectioNFormDecomposition}
H^2(X;\ZZ) \cong  \ZZ g_1 \oplus \ZZ g_2 \oplus \bigoplus_{i=1}^{b^+(X)} \ZZ e_i,
\end{equation}
where $Q_X(e_i,e_i)=1$ and $Q_X(g_i,g_i)=-1$ and all other pairings of $e_i$ and $g_j$ vanish for $i=1,\ldots,b^+(X)$ and $j=1,2$. We will next prove the

\begin{claim}
\label{claim:ExistenceOfPrimitive}
Let $\ka\in H^2(X;\ZZ)$ be the characteristic element in the statement of the lemma. There is a primitive $\ka_p\in H^2(X;\ZZ)$ with $\ka=d\ka_p$ and $d$ an odd, positive integer and a primitive class $\tilde\ka\in H^2(X;\ZZ)$ of the form
\begin{equation}
\label{eq:PrimitiveEll}
\tilde\ka = a_1e_1+a_2e_2+a_3e_3+a_4e_4+ce_5+\left(\sum_{j=6}^{b^+(X)} e_j\right) + cg_1+cg_2,
\end{equation}
such that $\tilde\ka^2=\ka_p^2$, where $g_i,e_j\in H^2(X;\ZZ)$ are as in \eqref{eq:OddIntersectioNFormDecomposition} and $a_j$ and $c$ are odd integers.
\end{claim}

Assuming Claim \ref{claim:ExistenceOfPrimitive}, we can proceed to complete the proof of Lemma \ref{lem:SemiAbundanceCondition} in Case \ref{case:Q_XOdd}. Using $r^2\equiv r\pmod 2$ for any $r\in\ZZ$, an elementary computation shows that an element $\tilde \ka$ of the form \eqref{eq:PrimitiveEll} will be characteristic if the coefficients $a_j$ and $c$ are odd. Thus, a class $\tilde\ka$ satisfying the conclusions of Claim \ref{claim:ExistenceOfPrimitive} is characteristic. If we define $h_i\in H^2(X;\ZZ)$ by $h_i=e_5+g_i$ for $i=1,2$, then $Q_X(h_i,\tilde\ka)=0$ and the span of $\{h_1,h_2\}$ gives a hyperbolic summand $Q_X$-orthogonal to $\tilde\ka$. Because $\tilde\ka$ is primitive, the previously cited result of Wall \cite[Theorem 1]{WallUnimodQuadForms} implies that there is an automorphism of $(H^2(X;\ZZ),Q_X)$ mapping $\tilde\ka$ to $\ka_p$ and thus the hyperbolic summand spanned by $h_1$ and $h_2$ to $\ka_p^\perp$ and hence to $\ka^\perp$, as required. This completes the proof of Lemma \ref{lem:SemiAbundanceCondition} in Case \ref{case:Q_XOdd}.

It remains to complete the

\begin{proof}[Proof of Claim \ref{claim:ExistenceOfPrimitive}]
Recall from Gompf and Stipsicz \cite[Definition 1.2.8 (c), p. 10]{GompfStipsicz} that $v\in H^2(X;\ZZ)$ is \emph{primitive} if $v=dv'$ for $d\in\ZZ$ and $v'\in H^2(X;\ZZ)$ implies that $d=\pm 1$. As noted in \cite[Definition 1.2.8 (c), p. 8]{GompfStipsicz}, we can write $\ka=d\ka_p$ where $\ka_p$ is primitive and we note that by replacing $\ka_p$ with $-\ka_p$ if necessary, we can assume that $d$ is a positive integer. Because $Q_X$ is odd by the hypothesis of Case \ref{case:Q_XOdd}, there is a class $v\in H^2(X;\ZZ)$ with $v^2$ odd.  As $\ka$ is characteristic by hypothesis of the lemma, we have $v^2\equiv Q_X(v,\ka)\pmod 2$ so $Q_X(v,\ka)$ is odd. By $\ka=d\ka_p$, we have $Q_X(v,\ka)=d Q_X(v,\ka_p)$ so $Q_X(v,\ka)$ odd implies that $d$ is odd.

\begin{subcase}[Proof of Claim \ref{claim:ExistenceOfPrimitive} when $b^+(X)=5$ and $b^-(X)=2$]
\label{subcase:ClaimProofWhenb+5b-2}
With $b^+(X)=5$ and $b^-(X)=2$, the hypothesis that $\ka^2=c_1(X)^2$ and the definition of $c_1(X)^2$ in \eqref{eq:DefineTopCharNumbersOf4Manifold} imply that
\[
\ka^2=2e(X)+3\si(X)=2(2+b^+(X)+b^-(X))+3(b^+(X)-b^-(X))=27.
\]
As noted above, we can write $\ka=d\ka_p$, where $d\in\ZZ$ and $\ka_p$ is primitive. The equality $27=\ka^2=d^2\ka_p^2$ implies that $d=\pm 1$ or $d=\pm 3$.  We consider the cases $d=\pm 1$ and $d=\pm 3$ separately.

\begin{subsubcase}[$d=\pm 3$]
If $d=\pm 3$, then $\ka_p^2=3$.  In this case, we define $\tilde\ka\in H^2(X;\ZZ)$ of the form \eqref{eq:PrimitiveEll} by
\[
\tilde\ka = e_1+e_2+e_3+e_4+e_5+g_1+g_2.
\]
Then $\tilde\ka^2=3=\ka_p^2$ and $\tilde\ka$ is primitive since the coefficient of $e_1$ is one. Hence, this $\tilde\ka$ satisfies Claim \ref{claim:ExistenceOfPrimitive}.
\end{subsubcase}

\begin{subsubcase}[$d=\pm 1$]
If $d=\pm 1$, so that $\ka_p^2=27$, we define $\tilde\ka\in H^2(X;\ZZ)$ of the form \eqref{eq:PrimitiveEll} by
\[
\tilde\ka = 5e_1+e_2+e_3+e_4+e_5+g_1+g_2.
\]
This $\tilde\ka$ is primitive because the coefficient of $e_2$ is one and this $\tilde\ka$ satisfies $\tilde\ka^2=27=\ka_p^2$. Hence, this $\tilde\ka$ satisfies Claim \ref{claim:ExistenceOfPrimitive}.
\end{subsubcase}
This completes the proof of Claim \ref{claim:ExistenceOfPrimitive} for Case \eqref{subcase:ClaimProofWhenb+5b-2}. 
\end{subcase}

\begin{subcase}[Proof of Claim \ref{claim:ExistenceOfPrimitive} when $b^+(X)>5$ and $b^-(X)=2$]
\label{subcase:ClaimProofWhenb+>5b-2}
As above, we can write $\ka=d \ka_p$ where $d$ is an odd, positive integer and $\ka_p$ is primitive. Because $d$ is odd, we have $d=2p+1$ for some $p\in\ZZ$. Because $\ka^2=(4p^2+4p+1)\ka_p^2$ and $p^2\equiv p\pmod 2$, we have
\begin{equation}
\label{eq:PrimitiveSquaredMod8}
\ka^2\equiv \ka_p^2\pmod 8.
\end{equation}
Because we are considering the case with $b^-(X)=2$, we have
\begin{align*}
\ka^2
&=c_1(X)^2
\\
&=2e(X)+3\si(X)
\\
&=4+2b^+(X)+2b^-(X)+3b^+(X)-3b^-(X)
\\
&=5b^+(X)+2.
\end{align*}
Equation \eqref{eq:PrimitiveSquaredMod8} then implies that there exists a $t\in\ZZ$ such that
\begin{equation}
\label{eq:PrimitiveSquaredValue}
\ka_p^2 = 5b^+(X)+2+8t.
\end{equation}
If $\tilde\ka\in H^2(X;\ZZ)$ has the form \eqref{eq:PrimitiveEll}, then $\tilde\ka$ will be primitive since the coefficient of $e_6$ is one. Thus $\tilde\ka$ will satisfy Claim \ref{claim:ExistenceOfPrimitive} if $\tilde\ka^2=\ka_p^2$, that is, by the expression \eqref{eq:PrimitiveSquaredValue} for $\ka_p^2$ and the expression for $\tilde\ka$ in \eqref{eq:PrimitiveEll}, if we can find odd integers $a_1,\dots,a_4,c$ such that
\[
5b^+(X)+2+8t
=
a_1^2+a_2^2+a_3^2+a_4^2 - c^2 +(b^+(X)-5),
\]
or equivalently,
\begin{equation}
\label{eq:ell2=ka_p2}
a_1^2+a_2^2+a_3^2+a_4^2=4b^+(X)+c^2+7+8t.
\end{equation}
Because $c$ is odd, we can write $c=2s+1$ for $s\in\ZZ$ and so the right-hand side of \eqref{eq:ell2=ka_p2} becomes
\begin{align*}
4b^+(X)+(2s+1)^2+7+8t
&=
4b^+(X)+4s^2+4s+8+8t
\\
&\equiv
4b^+(X)\pmod 8
\end{align*}
If we define $A$ by the right-hand-side of \eqref{eq:ell2=ka_p2},
\[
A:=4b^+(X)+c^2+7+8t,
\]
then because $b^+(X)$ is odd by hypothesis, the equivalence $4b^+(X)+(2s+1)^2+7+8t\equiv 4b^+(X)\pmod 8$ implies that
\[
A\equiv 4\pmod 8.
\]
By taking $s$ sufficiently large, the right-hand side of \eqref{eq:ell2=ka_p2} and thus $A$ will become positive. Because $A$ is a positive integer with $A\equiv 4\pmod 8$, then there are odd integers $a_1,\dots,a_4$ such that $A=a_1^2+a_2^2+a_3^2+a_4^2$ by Feehan and Leness \cite[Lemma A.5, p. 128]{FL2a}. The definition of $A$ and \eqref{eq:ell2=ka_p2} then imply that such a choice of $a_1,\dots,a_4,c$ defines a cohomology class $\tilde\ka$ satisfying $\tilde\ka^2=\ka_p^2$, and thus $\tilde\ka$ satisfies the conclusions of Claim \ref{claim:ExistenceOfPrimitive} for Case \eqref{subcase:ClaimProofWhenb+>5b-2}.
\end{subcase}

This completes the proof of Claim \ref{claim:ExistenceOfPrimitive}.
\end{proof}

This completes the treatment of Case \eqref{case:Q_XOdd}.
\end{case}

Because Claim \ref{claim:ExistenceOfPrimitive} implies Lemma \ref{lem:SemiAbundanceCondition} in Case \eqref{case:Q_XOdd}, this completes the proof of Lemma \ref{lem:SemiAbundanceCondition}.
\end{proof}

In the following series of lemmas, we let $\widetilde X=X\#\overline{\CC\PP}^2$ denote the smooth blow-up of $X$ at a point and let $e\in H^2(\widetilde X;\ZZ)$ be the Poincar\'e dual of the exceptional sphere. From Gompf--Stipsicz \cite[Definition 2.2.7, p. 43]{GompfStipsicz}, we have $e^2=-1$ and $H^2(\widetilde X;\ZZ)\cong H^2(X;\ZZ)\oplus H^2(\overline{\CC\PP}^2;\ZZ)$. We will use the preceding isomorphism and $Q_{\widetilde X}$-orthogonal decomposition to allow us to view $H^2(X;\ZZ)$ as a submodule of $H^2(\widetilde X;\ZZ)$.

\begin{lem}[Existence of \spinu structures]
\label{lem:Existence}
Let $X$ be a standard four-manifold with $b^-(X)\ge 2$. Let $\fs_0$ be a spin${}^c$ structure on $X$ with $c_1(\fs_0)^2=c_1(X)^2$ and let $f_1,f_2\in H^2(X;\ZZ)$ be cohomology classes satisfying the conclusions of Lemma \ref{lem:SemiAbundanceCondition} with $\ka=c_1(\fs_0)$ in \eqref{eq:HyperbolicSummandConditions}. Let $\widetilde X$ be the smooth blow-up of $X$ at a point and $e\in H^2(\widetilde X;\ZZ)$ be the Poincar\'e dual of the exceptional sphere. If $x,b\in\ZZ$ satisfy $x\equiv b\pmod 2$, then
\begin{equation}
\label{eq:Definey}
y_\pm := \frac{1}{2}x^2\mp x+\frac{1}{2}b-6\chi_h(\widetilde X)
\end{equation}
is an integer and there are rank-two spin${}^u$ structures $\ft_b^+$ and $\ft_b^-$ on $\widetilde X$ satisfying
\begin{equation}
\label{eq:CharClassesOfSpinutb}
p_1(\ft_b^\pm )=\left(c_1(\widetilde X)^2 -12\chi_h(\widetilde X) +b\right)\PD[\pt]
\quad\text{and}\quad
c_1(\ft_b^\pm )= xe+y_\pm f_1+f_2,
\end{equation}
where $\pt\in X$ is a point and $\PD[\pt]\in H^4(X;\ZZ)$ is its Poincar\'e dual. If we further assume that $x\equiv 0\pmod 2$, then the characteristic class $w_2(\ft_b^\pm)$ defined in \eqref{eq:SpinUCharacteristics} is good in the sense of Definition \ref{defn:Good}.
\end{lem}

\begin{rmk}[On the utility of goodness]
\label{rmk:GoodnessOfw2}
We shall not use the goodness of $w_2(\ft_b^\pm)$ in this monograph but include it for possible future use.
\end{rmk}

\begin{proof}
Because $x\equiv b\pmod 2$, we can write $x=2k+b$ for some $k\in\ZZ$.  Then
\begin{align*}
  y_\pm
  &
    =\frac{1}{2}(4k^2+4kb+b^2)\mp (2k+b)+\frac{1}{2}b-6\chi_h(\widetilde X)
    \quad\text{(by \eqref{eq:Definey})}
  \\
  &= 2k^2+2kb\mp 2k+\frac{1}{2}\left( b^2 + b\right)\mp b -6\chi_h(\widetilde X).
\end{align*}
Because $b^2+ b\equiv 0\pmod 2$ for any integer $b$, the last expression is an integer and so $y_\pm$ in \eqref{eq:Definey} is an integer.

For $\fs_0\in\Spin^c(X)$ as in the hypotheses with 
\begin{equation}
\label{eq:c_1SquaredCondition}
c_1(\fs_0)^2=c_1(X)^2, 
\end{equation}
we will show that Lemma \ref{lem:ExistenceOfSpinu} applies to the classes
$\La_\pm :=xe+y_\pm f_1+f_2\in H^2(\widetilde X;\ZZ)$, $w := \La_\pm-(c_1(\fs_0)+e)\in H^2(\widetilde X;\ZZ)$, and 
\[
\wp_\pm:=\left(c_1(\widetilde X)^2 -12\chi_h(\widetilde X) +b\right)\PD[\pt].
\]
By Gompf and Stipsicz, \cite[Proposition 2.4.16, p. 56]{GompfStipsicz}, $c_1(\fs_0)$ is characteristic for $X$ and an integer lift of $w_2(X)$ and so $c_1(\fs_0)+e$ is characteristic for $\widetilde X$ and is an integer lift of $w_2(\widetilde X)$. (One can check that $c_1(\fs_0)+e$ is characteristic by writing $y\in H^2(\widetilde X;\ZZ)$ as $y_0+ke$ where $y\in H^2(X;\ZZ)$ and $k\in \ZZ$.  Because $c_1(\fs_0)$ is characteristic, $y_0^2\equiv c_1(\fs_0)\cdot y_0\pmod 2$.  Then, $y^2=y_0^2-k^2$  while $(c_1(\fs_0)+e)\cdot (y_0+ke)=c_1(\fs_0)\cdot y_0 -k\equiv y_0^2-k^2\pmod 2$, implying that $c_1(\fs_0)+e$ is characteristic.) By the preceding definition of $w$, we have $\La_\pm-w=c_1(\fs_0)+e$ and so $\La_\pm \equiv w+w_2(X)\pmod 2$ as required in Item \eqref{item:ExistenceOfSpinu2} in Lemma \ref{lem:ExistenceOfSpinu}. Because $e$ is $Q_{\widetilde X}$-orthogonal to $H^2(X;\ZZ)$ and $f_i\cdot c_1(\fs_0)=0$ for $i=1,2$ by \eqref{eq:HyperbolicSummandConditions}, we have
\begin{equation}
\label{eq:IntersectionPairingsOnBlowUp}
e^2=-1,\quad
e\cdot f_i = f_i\cdot c_1(\fs_0) = f_i\cdot f_i=0 \quad\text{for } i=1,2,\quad
f_1\cdot f_2=1,\quad
e\cdot c_1(\fs_0)=0.
\end{equation}
Thus,
\begin{align*}
w^2 &=
    \left( (x-1)e+y_\pm f_1+f_2 -c_1(\fs_0)\right)^2
    \quad\text{(by definition of $w$ and $\La_\pm$)}
\\
&=
    -x^2+2x -1+2y_\pm+c_1(X)^2
     \quad\text{(by \eqref{eq:IntersectionPairingsOnBlowUp} and \eqref{eq:c_1SquaredCondition})}
\\
&=
    -x^2+2x +2y_\pm +c_1(\widetilde X)^2
    \quad\text{(by $c_1(\widetilde X)^2=c_1(X)-1$)}
\\
&
\equiv
    b-12\chi_h(\widetilde X) +c_1(\widetilde X)^2\pmod 4
    \quad\text{(by \eqref{eq:Definey}).}
\end{align*}
Hence, the cohomology classes $\La_\pm$, $w$, and $\wp_\pm$ defined above satisfy Item \eqref{item:ExistenceOfSpinu3} in Lemma \ref{lem:ExistenceOfSpinu}. Lemma  \ref{lem:ExistenceOfSpinu} then implies that there are a \spin${}^u$ structures $\ft_b^\pm$ with characteristic classes satisfying \eqref{eq:CharClassesOfSpinutb}.

We now show that if $x$ is even, then $w_2(\ft_b^\pm)$ is good in the sense of Definition \ref{eq:StratificationCptPU(2)Space}.  By Lemma \ref{lem:ExistenceOfSpinu}, $w_2(\ft_b^\pm)\equiv w\pmod 2$.  Any integer lift $w'$ of $w_2(\ft_b^\pm)$ can then be written as $w'=w+2u$ where $u\in H^2(\widetilde X;\ZZ)$. If we assume $x$ is even, then for $e^*$ the homology class of the exceptional curve,
\begin{align*}
\langle w',e^*\rangle
& \equiv
\langle w,e^*\rangle\pmod 2
\\
&\equiv 
\langle \La_\pm -(c_1(\fs_0)+e),e^*\rangle \pmod 2
\\
&\equiv
\langle -e+y_\pm f_1+f_2-c_1(\fs_0),e^*\rangle\pmod 2\quad\text{(by the assumption that $x\equiv 0\pmod 2$)}
\\
&\equiv 1\pmod 2
\end{align*}
Thus, if $x$ is even, $\langle w',e^*\rangle$ is odd so $w'$ cannot be a torsion class, implying that $w_2(\ft_b^\pm)$ is good. This completes the proof of Lemma \ref{lem:Existence}.
\end{proof}

Next, we show that the \spinu structures constructed in Lemma \ref{lem:Existence} obey the hypotheses of Lemma \ref{lem:SpinuStrRequirements}.

\begin{lem}[Properties of characteristic classes of \spinu structures]
\label{lem:SpinuCharSatisfyingCriteria}
Continue the notation and hypotheses of Lemma \ref{lem:Existence}. If $\fs^\pm$ are  spin${}^c$ structures on $\widetilde X$ satisfying $c_1(\fs^\pm)=c_1(\fs_0)\pm e$
\footnote{The existence of such spin${}^c$ structures follows from Gompf and Stipsicz \cite[Proposition 2.4.16, p. 56]{GompfStipsicz}}, then the following hold:
\begin{enumerate}
\item
If $x\in\ZZ$ obeys
\begin{equation}
\label{eq:xRequirementPlus}
x< -2\left( c_1(\widetilde X)^2 -5\chi_h(\widetilde X)\right) +b,
\end{equation}
then the characteristic classes $p_1(\ft_b^+)$ and $c_1(\ft_b^+)$ of the spin${}^u$ structure $\ft_b^+$ constructed in Lemma \ref{lem:Existence} and $c_1(\fs^+)$ obey the hypotheses  \eqref{eq:SpinuStrRequirements_c1Squared_equality_p1_bound_LaDotK_Requirement} and \eqref{eq:SpinuStrRequirements_LaSquared_requirement} of Lemma \ref{lem:SpinuStrRequirements}.

\item
If $x$ obeys
\begin{equation}
\label{eq:xRequirementMinus}
x> 2\left( c_1(\widetilde X)^2 -5\chi_h(\widetilde X)\right) -b,
\end{equation}
then the characteristic classes $p_1(\ft_b^-)$ and $c_1(\ft_b^-)$ of the spin${}^u$ structure $\ft_b^-$ constructed in Lemma \ref{lem:Existence} and $c_1(\fs^-)$ obey the hypotheses \eqref{eq:SpinuStrRequirements_c1Squared_equality_p1_bound_LaDotK_Requirement} and \eqref{eq:SpinuStrRequirements_LaSquared_requirement} of Lemma \ref{lem:SpinuStrRequirements}.

\item
\label{eq:xRequirement_SWnonTorsion}
The class $c_1(\fs^\pm)-c_1(\ft_b^\pm)$ is not torsion.
\end{enumerate}
\end{lem}

\begin{proof}
By Gompf and Stipsicz \cite[Proposition 2.4.16, p. 56]{GompfStipsicz}, if $\tilde\ka\in H^2(\widetilde X;\ZZ)$ is characteristic, then there is a spin${}^c$ structure $\tilde\fs$ on $\widetilde X$ with $c_1(\tilde \fs)=\tilde\ka$.  Because $c_1(\fs_0)$ is characteristic on $X$, we see that $c_1(\fs_0)\pm e$ is characteristic on $\widetilde X$ (as discussed in the proof of Lemma \ref{lem:Existence}). Hence there are spin${}^c$ structures $\fs^\pm$ on $\widetilde X$ with $c_1(\fs^\pm)=c_1(\fs_0)\pm e$.

Because $c_1(\fs_0)^2=c_1(X)^2$ by hypothesis, we compute that $c_1(\fs^\pm)^2=\left( c_1(\fs_0)\pm e\right)^2=c_1(\fs_0)^2-1=c_1(\widetilde X)^2$ by Gompf and Stipsicz \cite[Definition 2.2.7, p. 43]{GompfStipsicz}.
Hence, the classes $c_1(\fs^\pm)$ obey condition \eqref{eq:SpinuStrRequirements_c1Squared_equality} in Lemma \ref{lem:SpinuStrRequirements} while the classes $p_1(\ft_b^\pm)$ obey condition \eqref{eq:SpinuStrRequirements_p1_bound} by their expressions in \eqref{eq:CharClassesOfSpinutb}.

We proceed to show that $c_1(\ft_b^\pm)$ and $c_1(\fs^\pm)$ satisfy conditions \eqref{eq:SpinuStrRequirements_LaDotK_Requirement} and \eqref{eq:SpinuStrRequirements_LaSquared_requirement} in the hypotheses of Lemma \ref{lem:SpinuStrRequirements}. We note that the Poincar\'e dual $e\in H^2(\widetilde X;\ZZ)$ of the exceptional sphere satisfies $e^2=-1$ (see Gompf and Stipsicz \cite[Definition 2.2.7, p. 43]{GompfStipsicz}).  As in Lemma \ref{lem:Existence}, we let $f_1,f_2\in H^2(X;\ZZ)$ be the cohomology classes satisfying the conclusions of Lemma \ref{lem:SemiAbundanceCondition} with $\ka=c_1(\fs_0)$. Because $e$ is orthogonal to $H^2(X;\ZZ)\subset H^2(\widetilde X;\ZZ)$ we have $f_i\cdot e=0$ for $i=1,2$. The expression for $c_1(\ft_b^\pm )$  in \eqref{eq:CharClassesOfSpinutb} then yields
\begin{align*}
c_1(\ft_b^\pm)^2&=\left( xe+y_\pm f_1+f_2\right)^2= -x+2y_\pm
\end{align*}
so
\begin{equation}
\label{eq:CharacteristicPairings_c1t}
c_1(\ft_b^\pm)^2=-x^2+2y_\pm
\end{equation}
Because $c_1(\fs_0)\cdot e=0$, $c_1(\fs_0)\cdot f_i=0$,  and $f_i\cdot e=0$ for $i=1,2$ we can compute
\[
c_1(\ft_b^\pm)\cdot c_1(\fs^+)
=
\left( xe+y_\pm f_1+f_2\right)\cdot
\left( c_1(\fs_0)+ e\right)
=-x
\]
so
\begin{equation}
\label{eq:CharacteristicPairings_c1Dotc_s+}
c_1(\ft_b^\pm)\cdot c_1(\fs^+)=-x.
\end{equation}
A similar computation gives
\begin{equation}
\label{eq:CharacteristicPairings_c1Dotc_s-}
c_1(\ft_b^\pm)\cdot c_1(\fs^-)=x.
\end{equation}
To verify that $c_1(\ft_b^\pm)$ obeys condition \eqref{eq:SpinuStrRequirements_LaSquared_requirement} in Lemma \ref{lem:SpinuStrRequirements}, we compute
\begin{align*}
&c_1(\ft_b^\pm)^2
-
\left( 4c_1(\widetilde X)^2 -32\chi_h(\widetilde X) +3b\right)
\\
&\quad =-x^2+2y_\pm-
\left( 4c_1(\widetilde X)^2 -32\chi_h(\widetilde X) +3b\right)
\quad\text{(by \eqref{eq:CharacteristicPairings_c1t})}
\\
&\quad =
\mp 2x+b-12\chi_h(\widetilde X)-
     \left( 4c_1(\widetilde X)^2 -32\chi_h(\widetilde X) +3b\right)
     \\&
\qquad\text{(by definition \eqref{eq:Definey} of $y_\pm$)}
\\
&\quad =
\mp 2x -\left( 4c_1(\widetilde X)^2 -20\chi_h(\widetilde X) +2b\right).
\end{align*}
The preceding quantity is positive by our hypotheses \eqref{eq:xRequirementPlus} and \eqref{eq:xRequirementMinus}. Thus $c_1(\ft_b^\pm)$ satisfies  condition \eqref{eq:SpinuStrRequirements_LaSquared_requirement} in Lemma \ref{lem:SpinuStrRequirements}.

To prove that $c_1(\ft_b^\pm)$ and $c_1(\fs^\pm)$ satisfy Item  \eqref{eq:SpinuStrRequirements_LaDotK_Requirement} in Lemma \ref{lem:SpinuStrRequirements}, we compute
\begin{align*}
\frac{1}{2}c_1(\ft_b^\pm)^2+6\chi_h(\widetilde X) -\frac{1}{2}b
&=
-\frac{1}{2}x^2+y_\pm +6\chi_h(\widetilde X) -\frac{1}{2}b
\quad\text{(by \eqref{eq:CharacteristicPairings_c1t})}
\\
&=
\mp x
\quad\text{(by the definition of $y_\pm$ in \eqref{eq:Definey})}
\\
&=
c_1(\ft_b^\pm)\cdot c_1(\fs^\pm)
\quad\text{(by \eqref{eq:CharacteristicPairings_c1t} and \eqref{eq:CharacteristicPairings_c1t})}.
\end{align*}
Thus, $c_1(\ft_b^\pm)$ and $c_1(\fs^\pm)$ satisfy Item  \eqref{eq:SpinuStrRequirements_LaDotK_Requirement} in Lemma \ref{lem:SpinuStrRequirements}.

Consider Item \eqref{eq:xRequirement_SWnonTorsion}. By the definition of $c_1(\ft_b^\pm)$ in \eqref{eq:CharClassesOfSpinutb}, we have
\begin{align*}
c_1(\fs^\pm)-c_1(\ft_b^\pm)
&=
c_1(\fs_0)\pm e -xe-y_\pm f_1-f_2
\\
&=
c_1(\fs_0)+(-x\pm 1)e-y_\pm f_1-f_2.
\end{align*}
If there is an integer $n\in\ZZ$ such that $n(c_1(\fs^\pm)-c_1(\ft_b^\pm))=0$ then, by substituting the preceding identity, we obtain
\[
0=
n c_1(\fs_0)+n(-x\pm 1)e- n y_\pm f_1-nf_2.
\]
Taking the product of both sides of the preceding identity with $f_1$ gives
\[
0
=
f_1\cdot\left( n c_1(\fs_0)+n(-x\pm 1)e- n y_\pm f_1-nf_2\right)
=-n.
\]
Hence, the identity $n(c_1(\fs^\pm)-c_1(\ft_b^\pm))=0$ implies $n=0$ and so $c_1(\fs^\pm)-c_1(\ft_b^\pm)$ is not torsion. This completes the proof of Item \eqref{eq:xRequirement_SWnonTorsion} and thus of the lemma.
\end{proof}

The expression $\la^-(\tilde\ft,\tilde\fs)$ for the formal Morse--Bott index in \eqref{eq:FormalMorseIndexIntroThm} contains the first Chern class $c_1(\widetilde X)$ which depends on a choice of an almost complex structure on the four-manifold.  We begin our discussion of the formal Morse--Bott index by choosing almost complex structures on $X$ and on its smooth blow-up $\widetilde X$.

\begin{lem}
\label{lem:ACStructuresWithGivenc1}
Let $X$ be a closed, oriented, smooth Riemannian four-manifold and let $\fs_0$ be a spin${}^c$ structure on $X$ with $c_1(\fs_0)^2=c_1(X)^2$. Then there are almost complex structures $J$ on $X$ and $\tilde J$ on $\widetilde X=X\#\overline{\CC\PP}^2$ such that
\begin{subequations}
\label{eq:ACStructuresForIndex}
\begin{align}
\label{eq:ACStructureForIndexOnX}
c_1(X)&:=c_1(T_J^{1,0}X)=-c_1(\fs_0),
\\
\label{eq:ACStructureForIndexOnBlowUp}
c_1(\widetilde X)&:=c_1(T\widetilde X,\tilde J) =-c_1(\fs_0)-e,
\end{align}
\end{subequations}
where $e\in H^2(\widetilde X;\ZZ)$ is the Poincar\'e dual of the exceptional two-sphere, $T_J^{1,0}X$ is the subbundle of $TX\otimes_\RR\CC$ defined by $J$ (see Wells, \cite[Chapter I, Section 3, pp. 31-32]{Wells3}), and similarly for $(T\widetilde X,\tilde J)$.
\end{lem}

\begin{proof}
If $Y$ is a closed four-manifold and $\ka\in H^2(Y;\ZZ)$, then by Gompf and Stipsicz \cite[Theorem 1.4.15, p. 29]{GompfStipsicz}, there is an almost complex structure $J$ on $Y$ with $c_1(T_J^{1,0}Y)=\ka$ if and only if $\ka^2=c_1(Y)^2$ and $\ka\equiv w_2(Y)\pmod 2$. The spin${}^c$ structure $\fs_0$ on $X$ satisfies $c_1(\fs_0)^2=c_1(X)^2$ by hypothesis and $c_1(\fs_0)\equiv w_2(X)\pmod 2$ by Gompf and Stipsicz \cite[Proposition 2.4.16, p. 56]{GompfStipsicz}. The equalities $(-c_1(\fs_0))^2=c_1(\fs_0)^2$ and $-c_1(\fs_0)\equiv c_1(\fs)\pmod 2$ then imply that $(-c_1(\fs_0))^2=c_1(X)^2$ and $-c_1(\fs_0)\equiv w_2(X)\pmod 2$. Hence, there is an almost complex structure $J$ on $X$ with $c_1(T_J^{1,0}X)=-c_1(\fs_0)$.  Similarly, $-c_1(\fs_0)-e\equiv w_2(\widetilde X)\pmod 2$ and $(-c_1(\fs_0)-e)^2=c_1(\fs_0)^2-1=c_1(\widetilde X)^2$, so there is an almost complex structure $\tilde J$ on $\widetilde X$ satisfying $c_1(T\widetilde X,\tilde J)=-c_1(\fs_0)-e$.  This completes the proof of \eqref{eq:ACStructuresForIndex}.
\end{proof}

In the forthcoming lemma, we will give criteria which imply that the virtual Morse--Bott index $\la^-_{[A,\Phi]}(f)$ in
\eqref{eq:Virtual_Morse-Bott_index_moduli_space_non-abelian_monopoles} of the Hamiltonian function $f$ in \eqref{eq:Hitchin_function} is positive at any point $[A,\Phi]\in \iota_{\tilde\fs,\tilde\ft}(M_{\tilde\fs})\subset\sM_{\tilde\ft}$
by proving that the formal expression for this index, given by the function $\la^-(\tilde\ft,\tilde\fs)$ in \eqref{eq:FormalMorseIndexIntroThm} is positive. The forthcoming Lemma \ref{lem:IndexInAlgebraicCase} applies 
not just at points in $M_{\tilde\fs}$ where $\SW_X(\tilde\fs)\neq 0$ but to points in the following expansion of the set $B(X)$ of basic classes defined in \eqref{eq:SetOfBasicClasses},
\begin{equation}
\label{eq:NonEmptyReducibleMonopoles}
\widehat B(X)
:=
\{c_1(\fs): \fs\in\Spinc(X) \text{ and } M_{\fs}\neq \varnothing\}.
\end{equation}
The inclusion $B(X)\subset\widehat B(X)$ may be strict if there are solutions to the Seiberg--Witten equations for a \spinc structure with zero Seiberg--Witten invariant.

\begin{lem}[Positivity of an expression for the formal Morse--Bott index]
\label{lem:IndexInAlgebraicCase}
Continue the hypotheses of Lemma \ref{lem:Existence}.
Let $c_1(X)$ and $c_1(\widetilde X)$ be the Chern classes defined in \eqref{eq:ACStructuresForIndex}.
Let $f_1,f_2\in H^2(X;\ZZ)\subset H^2(\widetilde X:\ZZ)$ be the cohomology classes satisfying equation \eqref{eq:HyperbolicSummandConditions} in Lemma \ref{lem:SemiAbundanceCondition} with
$\ka=c_1(X)$. Then there are integers $b$ satisfying
\begin{equation}
\label{eq:bGeneralTypeAssumption}
b< 10\chi_h(\widetilde X),
\end{equation}
and for every such $b$, there are integers $x$ satisfying
\begin{subequations}
\label{eq:IndexInAlgebraicCase_xAssumptions}
\begin{align}
\label{eq:GeneralTypeAssumptionOnx}
&
x>
	\max_{c_1(\tilde\fs)\in \widehat B(\widetilde X)}
	\left( c_1(\tilde\fs) +c_1(\widetilde X)\right)\cdot c_1(\widetilde X)
\quad\text{where $\widehat B(\widetilde X)$ is defined in \eqref{eq:NonEmptyReducibleMonopoles}},
\\
\label{eq:IndexInAlgebraicCase_xbRequirementMinus}
&\text{$x$ and $b$ satisfy \eqref{eq:xRequirementMinus}},
\\
  \label{eq:IndexInAlgebraicCase_x_NonNeg}
    &x \geq 0,
\\
\label{eq:IndexInAlgebraicCase_EqualParity}
& x\equiv b\pmod 2.
\end{align}
\end{subequations}
Moreover, for any ordered pair of integers $(b,x)$ satisfying \eqref{eq:bGeneralTypeAssumption} and \eqref{eq:IndexInAlgebraicCase_xAssumptions}, the following holds. Let $\ft_b^-$ be the rank-two spin${}^u$ structure
constructed in Lemma \ref{lem:Existence} such that $p_1(\ft_b^-)$ and $c_1(\ft_b^-)$ satisfy \eqref{eq:CharClassesOfSpinutb} for the cohomology classes $f_1$ and $f_2$ and the ordered pair of integers $(b,x)$. Then for all $\tilde\fs\in\Spinc(\widetilde X)$ with $c_1(\tilde \fs)\in \widehat B(\widetilde X)$ and $\iota_{\tilde\fs,\ft_b^-}(M_{\tilde \fs})\subset\sM_{\ft^-_b}$, we have $\la^-(\ft_b^-,\tilde\fs)>0$, where $\la^-(\ft_b^-,\tilde\fs)$ is defined in \eqref{eq:FormalMorseIndexIntroThm}.
\end{lem}

\begin{rmk}[On the role of the hypothesis that $b_1(X)=0$]
\label{rmk:TechnicalRemarkOnb1=0}
As discussed in Remark \ref{rmk:Onb1Equals0}, the assumption that $b_1(X)=0$ (implicit in the hypothesis that $X$ is standard) implies that $\chi_h(\widetilde X)>0$ and that allows us to find a positive integer $b$ satisfying \eqref{eq:bGeneralTypeAssumption}.  The positivity of $b$ gives the desired upper bound on the expected dimension of the moduli space of anti-self-dual connections stated in \eqref{eq:ASDDimInTermsOfBMY_lower_bound}.
\end{rmk}

\begin{proof}[Proof of Lemma \ref{lem:IndexInAlgebraicCase}]
That there is an integer $b$ satisfying \eqref{eq:bGeneralTypeAssumption} is obvious and so we now fix such a $b$ and verify the existence of $x$ satisfying \eqref{eq:IndexInAlgebraicCase_xAssumptions}. Because the set $\widehat B(\widetilde X)$ is finite by
Morgan \cite[Theorem 5.2.4, p. 79]{MorganSWNotes}, the quantity on the right-hand-side of \eqref{eq:GeneralTypeAssumptionOnx} is finite and so we may chose an integer $x$ to satisfy \eqref{eq:GeneralTypeAssumptionOnx}. Then the inequalities \eqref{eq:GeneralTypeAssumptionOnx} and \eqref{eq:xRequirementMinus} (for the fixed integer $b$) are both lower bounds for $x$, we may also choose $x$ to be non-negative and so satisfy \eqref{eq:IndexInAlgebraicCase_x_NonNeg} and the parity constraint of  \eqref{eq:IndexInAlgebraicCase_EqualParity}. Hence, there are integers $b$ and $x$ satisfying
\eqref{eq:bGeneralTypeAssumption} and \eqref{eq:IndexInAlgebraicCase_xAssumptions}.

We next prove that $\la^-(\ft_b^-,\tilde\fs)>0$, as asserted in the final conclusion of Lemma \ref{lem:IndexInAlgebraicCase}.
Let $\fs^-$ be the \spinc structure on $\widetilde X$ defined in Lemma \ref{lem:SpinuCharSatisfyingCriteria} with $c_1(\fs^-)=c_1(\fs_0)-e$. By \eqref{eq:ACStructureForIndexOnX}, we have
\begin{equation}
\label{eq:SpincMinus_and_c1X}
c_1(\fs^-)=-c_1(X)-e.
\end{equation}
The equalities \eqref{eq:ACStructureForIndexOnX} and \eqref{eq:ACStructureForIndexOnBlowUp}
imply that $c_1(\widetilde X)=-c_1(\fs_0)-e=c_1(X)-e$, which we record as
\begin{equation}
\label{eq:BlowUpC1}
c_1(\widetilde X)=c_1(X)-e.
\end{equation}
We compute that
\begin{align*}
c_1(\fs^-)\cdot c_1(\widetilde X)
&=
                                    \left(-c_1(X)-e\right)\left(c_1(X)-e\right)
                                    \quad\text{(by \eqref{eq:BlowUpC1} and \eqref{eq:SpincMinus_and_c1X})}
\\
&=-c_1(X)^2-1
     \quad\text{(by $e^2=-1$)},
\end{align*}
and so, using $c_1(\widetilde X)^2=c_1(X)^2-1$,
\begin{equation}
\label{eq:IndexCompCohomEqualitySpinCanDotc1}
c_1(\fs^-)\cdot c_1(\widetilde X) = -c_1(\widetilde X)^2-2.
\end{equation}
By hypothesis, the cohomology classes $f_1,f_2\in H^2(X;\ZZ)$ satisfy equation \eqref{eq:HyperbolicSummandConditions} with $\kappa = c_1(\fs_0)$, that is,
\[
f_1\cdot f_2=1 \quad\text{and}\quad c_1(\fs_0)\cdot f_i=0=f_i^2, \quad\text{for } i=1,2.
\]
Because $c_1(\fs_0)\in H^2(X;\ZZ)$ and $e$ is $Q_{\widetilde X}$-orthogonal to $H^2(X;\ZZ)$, we have $c_1(\fs_0)\cdot e=0$. By combining the preceding equalities with the expression for $c_1(\ft_b^-)$ in \eqref{eq:CharClassesOfSpinutb} and that for $c_1(\widetilde X)$ in \eqref{eq:ACStructureForIndexOnBlowUp} , we obtain
\[
c_1(\ft_b^-)\cdot c_1(\widetilde X)=(xe+y_-f_1+f_2)\cdot (-c_1(\fs_0)-e)=x,
\]
and so
\begin{equation}
\label{eq:IndexCompCohomEqualityc1tDotc_1X}
c_1(\ft_b^-)\cdot c_1(\widetilde X)=x.
\end{equation}
We compute
\[
c_1(\fs^-)^2=\left( -c_1(X)-e\right)^2=c_1(X)^2 -1=c_1(\widetilde X)^2.
\]
The resulting equality $c_1(\fs^-)^2=c_1(\widetilde X)^2$ and Lemma \ref{lem:SpinuCharSatisfyingCriteria} imply that the \spinc structure $\fs^-$ and \spinu structure $\ft_b^-$ constructed in Lemma \ref{lem:Existence} satisfy the conditions
\eqref{eq:SpinuStrRequirements_c1Squared_equality_p1_bound_LaDotK_Requirement} in Lemma \ref{lem:SpinuStrRequirements}. Therefore, Lemma \ref{lem:SpinuStrRequirements} implies that $\iota_{\fs^-,\ft_b^-}(M_{\fs^-})\subset \sM_{\ft_b^-}$. The preceding inclusion and \eqref{eq:LowerLevelInclusionOFReducibles} imply that
\begin{equation}
\label{eq:Levelofs-}
\ell(\ft_b^-,\fs^-)=0.
\end{equation}  
The formula for $\ell(\ft_b^-,\fs^-)$ in \eqref{eq:ReducibleLevel} then implies that
\begin{equation}
\label{eq:IndexCompCohomEqualityReducibleLevel}
(c_1(\fs^-)-c_1(\ft_b^-))^2=p_1(\ft_b^-).
\end{equation}
Combining \eqref{eq:IndexCompCohomEqualityReducibleLevel} with the expression for $p_1(\ft_b^-)$ in \eqref{eq:CharClassesOfSpinutb} yields
\begin{equation}
\label{eq:IndexCompCohomEqualityReducibleLevel2}
(c_1(\fs^-)-c_1(\ft_b^-))^2
=
c_1(\widetilde X)^2 -12\chi_h(\widetilde X) +b.
\end{equation}
Substituting the equalities indicated below into the expression \eqref{eq:FormalMorseIndexIntroThm}
yields,
\begin{align*}
\la^-(\ft_b^-,\fs^-)
&=
-2\chi_h(\widetilde X) -\left( c_1(\fs^-)-c_1(\ft_b^-)\right)\cdot c_1(\widetilde X)
-\left( c_1(\fs^-)-c_1(\ft_b^-)\right)^2
  \\
  &\qquad\text{(by \eqref{eq:FormalMorseIndexIntroThm})}
\\
&=
-2\chi_h(\widetilde X) -\left( c_1(\fs^-)-c_1(\ft_b^-)\right)\cdot c_1(\widetilde X)
-\left(c_1(\widetilde X)^2 -12\chi_h(\widetilde X) +b\right)
  \\
  &\qquad\text{(by \eqref{eq:IndexCompCohomEqualityReducibleLevel2})}
\\
&=
-2\chi_h(\widetilde X)+c_1(\widetilde X)^2+2+x
     -c_1(\widetilde X)^2 +12\chi_h(\widetilde X)-b,
  \\
  &\qquad\text{(by \eqref{eq:IndexCompCohomEqualitySpinCanDotc1} and \eqref{eq:IndexCompCohomEqualityc1tDotc_1X})},
\end{align*}
and thus
\begin{equation}
\label{eq:KXIndex}
\la^-(\ft_b^-,\fs^-)
=
x +10\chi_h(\widetilde X)+2-b.
\end{equation}
Our hypothesis \eqref{eq:bGeneralTypeAssumption} on $b$ then ensures that $\la^-(\ft_b^-,\fs^-)>x$. By inequality \eqref{eq:IndexInAlgebraicCase_x_NonNeg}, we have $x$ is non-negative and so $\la^-(\ft_b^-,\fs^-)>0$.

We next prove that $\la^-(\ft_b^-,\tilde\fs)>0$ for any $\tilde\fs\in\Spinc(\widetilde X)$ for which there is an embedding $\iota_{\tilde\fs,\ft_b^-}(M_{\tilde \fs})\subset\sM_{\ft^-_b}$ as in \eqref{eq:DefnOfIotaOnQuotient}. By \eqref{eq:LowerLevelInclusionOFReducibles}, the level of the set $M_{\tilde\fs}$ of points represented by split monopoles in $\bar\sM_{\ft_b^-}$ is zero, that is, $\ell(\ft^-_b,\tilde\fs)=0$.  The expression for  $\ell(\ft^-_b,\tilde\fs)$ in \eqref{eq:ReducibleLevel} and \eqref{eq:Levelofs-} then imply that
\begin{equation}
\label{eq:SameP1FromSameLevel}
\left(c_1(\tilde \fs)-c_1(\ft^-_b)\right)^2
=
\left(c_1(\fs^-)-c_1(\ft^-_b)\right)^2.
\end{equation}
By comparing the expressions for $\la^-(\ft_b^-,\tilde\fs)$ and $\la^-(\ft_b^-,\fs^-)$ in \eqref{eq:FormalMorseIndexIntroThm}, we see that
\begin{multline*}
\la^-(\ft_b^-,\tilde\fs)-\la^-(\ft_b^-,\fs^-)
=
-\left( c_1(\tilde\fs)-c_1(\ft_b^-)\right)\cdot c_1(X)
-\left( c_1(\tilde\fs)-c_1(\ft_b^-)\right)^2
\\
-
\left(-\left( c_1(\fs^-)-c_1(\ft_b^-)\right)\cdot c_1(X)
-\left( c_1(\fs^-)-c_1(\ft_b^-)\right)^2\right).
\end{multline*}
Thus,
\begin{align*}
		\la^-(\ft_b^-,\tilde\fs)
		&=
		\la^-(\ft_b^-,\fs^-)
		+
           \left(c_1(\fs^-)- c_1(\tilde\fs)\right)\cdot c_1(\widetilde X)
  \\
                &\qquad
        -
        \left( \left(c_1(\tilde \fs)-c_1(\ft^-_b)\right)^2 -\left(c_1(\fs^-)-c_1(\ft^-_b)\right)^2\right)
        \\
        &=
        \la^-(\ft_b^-,\fs^-)
		+
           \left(c_1(\fs^-)- c_1(\tilde\fs)\right)\cdot c_1(\widetilde X)
  \\
                &
    \qquad\text{(by \eqref{eq:SameP1FromSameLevel})}
		\\
		&=
        \left( x +10\chi_h(\widetilde X)+2-b\right)
        +
        \left( -c_1(\widetilde X)-c_1(\tilde \fs)\right)\cdot c_1(\widetilde X) -2
  \\
                &\qquad\text{(by \eqref{eq:KXIndex} and \eqref{eq:IndexCompCohomEqualitySpinCanDotc1})},
\end{align*}
and so we obtain
\begin{equation}
\label{eq:LowerBoundForIndexAtTopStratumSW}
\la^-(\ft_b^-,\tilde\fs)=
\left(
x-  \left( c_1(\widetilde X)+c_1(\tilde \fs)\right)\cdot c_1(\widetilde X)
\right)
+
\left( 10\chi_h(\widetilde X) - b\right).
\end{equation}
The preceding expression is positive for $c_1(\tilde\fs)\in\widehat B(\widetilde X)$ by the hypotheses \eqref{eq:bGeneralTypeAssumption} and \eqref{eq:GeneralTypeAssumptionOnx}. This completes the proof of Lemma \ref{lem:IndexInAlgebraicCase}.
\end{proof}

\begin{rmk}[Virtual Morse--Bott index at Seiberg-Witten points in lower levels of the Uhlenbeck compactification]
\label{rmk:VirtualMBIndexAtLowerStrataReducibles}
We now discuss how the requirements on $t_b^-$ in \eqref{eq:bGeneralTypeAssumption} and \eqref{eq:IndexInAlgebraicCase_xAssumptions} must change to make the formal Morse--Bott index \eqref{eq:FormalMorseIndexIntroThm}  positive not just for \spinc structures $\tilde\fs$ where $M_{\tilde\fs}$ admits a continuous embedding in $\sM_\ft$ but for \spinc structures $\tilde\fs_\ell$ where $M_{\tilde\fs_\ell}$ admits a continuous embedding in a lower level of the Uhlenbeck compactification, $\bar\sM_\ft$.

Let $\ft_b^-$ be the rank-two \spinu structure constructed in Lemma \ref{lem:Existence} and appearing again in Lemma \ref{lem:IndexInAlgebraicCase} and let $\fs^-$ be a \spinc structure satisfying \eqref{eq:SpincMinus_and_c1X}.
For a non-negative integer $\ell$, let $\ft_b^-(\ell)$ be the \spinu structure on $\widetilde X$ satisfying, as in \eqref{eq:DefineLowerChargeSpinuStr},
\[
c_1(\ft_b^-(\ell))=c_1(\ft_b^-)
\quad\text{and}\quad
p_1(\ft_b^-(\ell))=p_1(\ft_b^-)+4\ell.
\]
Let $\tilde\fs_\ell$ be another \spinc structure on $\widetilde X$ with $c_1(\tilde\fs_\ell)\in \widehat B(\widetilde X)$ (as defined in \eqref{eq:NonEmptyReducibleMonopoles}) and assume that there is a continuous embedding,
\begin{equation}
\label{eq:LowerLevelInclusionOFReduciblesRemark}
M_{\tilde\fs_\ell}\times\Sym^\ell(\widetilde X) \hookrightarrow
\sM_{\ft_b^-(\ell)}\times\Sym^\ell(\widetilde X)
\subset
\sI\!\!\sM_{\ft_b^-},
\end{equation}
as in \eqref{eq:LowerLevelInclusionOFReducibles}. By Lemma \ref{lem:SetOfReducibles} and \eqref{eq:IndexCompCohomEqualityReducibleLevel}, we see that
\[
\left( c_1(\tilde\fs_\ell)-c_1(\ft_b^-(\ell))\right)^2
=
p_1(\ft_b^-)+4\ell
=
\left( c_1(\fs^-)-c_1(\ft_b^-)\right)^2+4\ell.
\]
(Recall from \eqref{eq:ReducibleLevel} that the preceding equality defines the level as a function $\ell=\ell(\ft_b^-,\tilde\fs_\ell)$.)

By modifying the calculation of $\la^-(\ft_b^-,\tilde\fs)$ given at the end of Lemma \ref{lem:IndexInAlgebraicCase} and using the equality $c_1(\ft_b^-)=c_1(\ft_b^-(\ell))$ from \eqref{eq:DefineLowerChargeSpinuStr_c1'}, we would obtain the following generalization of \eqref{eq:LowerBoundForIndexAtTopStratumSW}:
\begin{equation}
\label{eq:LowerBoundForIndexAtLowerLevelSW}
\la^-(\ft_b^-(\ell),\tilde\fs_\ell)=
\left(
x-  \left( c_1(\widetilde X)+c_1(\tilde \fs_\ell)\right)\cdot c_1(\widetilde X)
\right)
+
\left( 10\chi_h(\widetilde X) - b\right)
-4\ell.
\end{equation}
The assumption \eqref{eq:GeneralTypeAssumptionOnx} implies that the term
\[
x-  \left( c_1(\widetilde X)+c_1(\tilde \fs_\ell)\right)\cdot c_1(\widetilde X)
\]
on the right-hand-side of \eqref{eq:LowerBoundForIndexAtLowerLevelSW} remains positive. However, the assumption \eqref{eq:bGeneralTypeAssumption} is no longer strong enough to ensure that
the term
\[
\left( 10\chi_h(\widetilde X) - b\right)
-4\ell
\]
on the right-hand-side of \eqref{eq:LowerBoundForIndexAtLowerLevelSW} remains positive. The integer $\ell=\ell(\ft_b^-,\tilde\fs_\ell)$ depends on $\fs_\ell$ and so, to ensure that the expression \eqref{eq:LowerBoundForIndexAtLowerLevelSW} remains positive for all $\tilde \fs_\ell$ with $c_1(\tilde\fs_\ell)\in\widehat B(\widetilde X)$, one would need to replace the assumption \eqref{eq:bGeneralTypeAssumption} with the inequality
\[
b<  10\chi_h(\widetilde X)-4\ell_{\max},
\]
where
\[
  \ell_{\max}
  := \max\left\{\ell(\ft_b^-,\tilde\fs_\ell): c_1(\tilde\fs_\ell)\in\widehat B(\widetilde X), \text{ for some } \ell \geq 0\right\}.
\]
Hence, the requirement that $\la^-(\ft_b^-(\ell),\tilde\fs_\ell)$ be positive for all integers $\ell\geq 0$ with \spinc structures $\tilde\fs_\ell$ obeying \eqref{eq:LowerLevelInclusionOFReduciblesRemark} is more demanding than a requirement that $\la^-(\ft_b^-,\tilde\fs_\ell)$ be positive when $\tilde\fs_\ell$ obeys \eqref{eq:LowerLevelInclusionOFReduciblesRemark} with $\ell=0$.
\end{rmk}

\section{Proofs of main results on feasibility}
\label{sec:Proofs_main_results_feasibility}
We can now complete the

\begin{proof}[Proof of Theorem \ref{mainthm:ExistenceOfSpinuForFlow}]
Let $\fs_0$ be the spin${}^c$  structure on $X$ in the hypotheses of Theorem \ref{mainthm:ExistenceOfSpinuForFlow}. Because $X$ is of Seiberg--Witten simple type as in \eqref{eq:SWSimpleType}, Theorem \ref{thm:SWBlowUp} implies that $\widetilde X$ is also Seiberg--Witten simple type.  Write $\fs^\pm$ for the spin${}^c$ structures on $\widetilde X$ with $c_1(\fs^\pm)=c_1(\fs_0)\pm e$. Theorem \ref{thm:SWBlowUp} implies that
\[
\SW_{\widetilde X}(\fs^\pm)
=
\SW_X(\fs_0)
=
\SW_{\widetilde X}(\tilde\fs).
\]
Thus, $c_1(\fs^-)\in B(\widetilde X)$, the set of basic classes defined in \eqref{eq:SetOfBasicClasses}. Because $\widetilde X$ has Seiberg--Witten simple type and $c_1(\fs^-)\in B(\widetilde X)$, then $c_1(\fs^-)$ satisfies conditions \eqref{eq:SpinuStrRequirements_c1Squared_equality} and \eqref{eq:SpinuStrRequirementsNonZeroSW} in Lemma \ref{lem:SpinuStrRequirements}.

The hypothesis in Theorem \ref{mainthm:ExistenceOfSpinuForFlow} that $X$ is standard implies that $b_1(X)=0$ and $b^+(X)\ge 3$. Therefore, $\chi_h(X)=(1+b^+(X))/2$ by \eqref{eq:DefineTopCharNumbersOf4Manifold} and so $\chi_h(X)\geq 2$ since $b^+(X)\ge 3$. Because $\chi_h(\widetilde X) = \chi_h(X)$, we see that $10\chi_h(\widetilde X)=10\chi_h(X) \geq 20$. Thus, we can chose a \emph{positive} integer $b$ satisfying \eqref{eq:bGeneralTypeAssumption}. For this choice of $b$, let $x$ be an integer satisfying \eqref{eq:IndexInAlgebraicCase_xAssumptions}. Let $\ft_b^-$ be the spin${}^u$ structure constructed in Lemma \ref{lem:Existence} using $c_1(\fs_0)$ and the above choice of $b$ and $x$.

We proceed to verify that $\fs^-$ and $\ft_b^-$ satisfy the conditions \eqref{eq:SpinuStrRequirements_c1Squared_equality_p1_bound_LaDotK_Requirement}
and \eqref{eq:SpinuStrRequirements_Good_NonZeroSW_LaSquared_requirement} required of $\tilde\fs$ and $\tilde\ft$ in Lemma \ref{lem:SpinuStrRequirements}. As noted above, $\fs^-$ satisfies conditions \eqref{eq:SpinuStrRequirements_c1Squared_equality} and \eqref{eq:SpinuStrRequirementsNonZeroSW}. Because $b$ and $x$ satisfy  \eqref{eq:IndexInAlgebraicCase_xbRequirementMinus} in Lemma \ref{lem:IndexInAlgebraicCase}, Lemma \ref{lem:SpinuCharSatisfyingCriteria} implies that $\fs^-$ and $\ft_b^-$ satisfy conditions \eqref{eq:SpinuStrRequirements_p1_bound}, \eqref{eq:SpinuStrRequirements_LaDotK_Requirement}, and \eqref{eq:SpinuStrRequirements_LaSquared_requirement}. Thus, $\fs^-$ and $\ft_b^-$ satisfy the conditions \eqref{eq:SpinuStrRequirements_c1Squared_equality_p1_bound_LaDotK_Requirement}
and \eqref{eq:SpinuStrRequirements_Good_NonZeroSW_LaSquared_requirement} in Lemma \ref{lem:SpinuStrRequirements}.

By Item \eqref{eq:xRequirement_SWnonTorsion} of Lemma \ref{lem:SpinuCharSatisfyingCriteria}, 
$c_1(\fs^-)-c_1(\ft_b^-)$ is not torsion. Proposition \ref{prop:SWModuliSpaceZeroSectionCriterion} then implies that for a smooth Riemannian metric $\tilde g$ in an open, dense subspace of the space of Riemannian metrics in the $C^r$ topology (with $r\ge 3$) on $\widetilde X$, the moduli space $M_{\fs^-}$ does not contain the gauge-equivalence class of a zero-section pair.

Next, we explain how Lemma \ref{lem:SpinuStrRequirements} yields the conclusions of Theorem \ref{mainthm:ExistenceOfSpinuForFlow}. Because $M_{\fs^-}$ does not contain a gauge-equivalence class of a zero-section pair for a generic metric, Item \eqref{item:SpinuStrRequirementsEmptyZeroSection} in Lemma \ref{lem:SpinuStrRequirements} implies that Item \eqref{item:ExistenceOfSpinuForFlowNonEmpty} of Theorem \ref{mainthm:ExistenceOfSpinuForFlow} holds. Item \eqref{item:SpinuStrRequirements_Basic_Lower_Bound} in Lemma \ref{lem:SpinuStrRequirements} and the fact that $b$ is a positive integer yield Item \eqref{item:ExistenceOfSpinuForFlowp1} of Theorem \ref{mainthm:ExistenceOfSpinuForFlow}. The almost complex structure $\tilde J$ in equation \eqref{eq:ACStructureForIndexOnBlowUp} satisfies the requirement of Item \eqref{item:ExistenceOfACStructure} in Theorem \ref{mainthm:ExistenceOfSpinuForFlow}.
Item \eqref{item:Upper_bound_expected_dimension_ASD_moduli_space_blowup} in Theorem \ref{mainthm:ExistenceOfSpinuForFlow} follows from the fact that $b\geq 1$ and from the resulting inequality \eqref{eq:ASDDimInTermsOfBMY_lower_bound} in Item \eqref{item:SpinuStrRequirements_ASDModuliIncluded} in Lemma \ref{lem:SpinuStrRequirements}.

Because $b$ and $x$ satisfy the conditions \eqref{eq:bGeneralTypeAssumption} and \eqref{eq:IndexInAlgebraicCase_xAssumptions}, Lemma \ref{lem:IndexInAlgebraicCase} implies that $\la^-(\ft_b^-,\tilde\fs')>0$ for all  $\tilde\fs'\in\Spinc(\widetilde X)$ with $c_1(\tilde \fs')\in \widehat B(\widetilde X)$ and $\iota_{\tilde\fs',\ft_b^-}(M_{\tilde \fs})\subset\sM_{\ft^-_b}$.  Thus, $\ft_b^-$ satisfies Item \eqref{item:PositiveFormalIndex} of Theorem \ref{mainthm:ExistenceOfSpinuForFlow}.  This completes the proof of Theorem \ref{mainthm:ExistenceOfSpinuForFlow}.
\end{proof}

\begin{rmk}[Non-emptiness of $\sM_{\tilde\ft}^{*,0}$ for generic K\"ahler metrics]
\label{rmk:Proof_of_GenericKahlerMetrics_for_NoZeroPairs}
As noted following the statement of Theorem \ref{mainthm:ExistenceOfSpinuForFlow}, when $X$ is a closed, complex K\"ahler surface, it is possible to show that $\sM_{\tilde\ft}^{*,0}$ is non-empty when the metric $\tilde g$ is K\"ahler with fundamental form lying in an open, dense subspace of the K\"ahler cone. The requirement of Item \eqref{item:ExistenceOfSpinuForFlowNonEmpty} of Theorem \ref{mainthm:ExistenceOfSpinuForFlow} that the metric $g$ lie in an open, dense subspace of the space of Riemannian metrics is only used in the proof of Theorem \ref{mainthm:ExistenceOfSpinuForFlow} to use Proposition \ref{prop:SWModuliSpaceZeroSectionCriterion} to prove that $M_{\fs^-}$ does not contain the equivalence class of a zero-section pair. If $X$ is a closed K\"ahler surface,
we can use the forthcoming Lemma \ref{lem:ExistenceOfZeroSectionSWonKahler} instead of Proposition \ref{prop:SWModuliSpaceZeroSectionCriterion} to prove that $M_{\fs^-}$ does not contain the gauge-equivalence class of a zero-section pair for any K\"ahler metric whose fundamental form lies in an open dense subspace of the K\"ahler cone.
\end{rmk}

Next, we give the

\begin{proof}[Proof of Corollary \ref{maincor:MorseIndexForFeasibilitySpinuStructure}]
Let $\tilde\ft$ be the \spinu structure on $\widetilde X$ provided by Theorem \ref{mainthm:ExistenceOfSpinuForFlow}. If $X$ is a closed, complex K\"ahler surface, then $\widetilde X$ is as well and so for all gauge-equivalence classes of type $1$ non-Abelian monopoles $[A,\Phi]$ in $\iota_{\tilde\fs,\tilde\ft}(M_{\tilde\fs})\subset\sM_{\tilde\ft}$, equation \eqref{eq:MorseIndexAtReduciblesOnKahlerSpinNotationType1} in Corollary \ref{maincor:MorseIndexAtReduciblesOnKahlerWithSO3MonopoleCharacteristicClasses} and the definition \eqref{eq:FormalMorseIndexIntroThm} imply that $\la_{[A,\Phi]}^-(f)=\la^-(\tilde\ft,\tilde\fs)$. The conclusion of Corollary \ref{maincor:MorseIndexForFeasibilitySpinuStructure} now follows immediately from Item \eqref{item:PositiveFormalIndex} of Theorem \ref{mainthm:ExistenceOfSpinuForFlow} and Corollary \ref{maincor:MorseIndexAtReduciblesOnKahlerWithSO3MonopoleCharacteristicClasses}.
\end{proof}

\chapter[Monopoles, projective vortices, and  holomorphic pairs over K\"ahler surfaces]{Non-Abelian monopoles, projective vortices, and  holomorphic pairs over K\"ahler surfaces}
\label{chap:SO(3)monopolesoverKahler surfaces}
In \cite{Witten}, Witten proved that solutions to the unperturbed version of the Seiberg--Witten equations \eqref{eq:SeibergWitten} over a complex K\"ahler surface are always of type $1$ or type $2$ and thus can be viewed as solutions to projective vortex equations in the sense of Bradlow \cite{Bradlow_1990, Bradlow_1991}. In this chapter, we give an exposition of the proof of the corresponding facts for the unperturbed version of the non-Abelian monopole equations \eqref{eq:PerturbedSO3MonopoleEquations} over a complex K\"ahler surface, so $\tau=\id_{\Lambda^+}$ and $\vartheta=0$ and $g$ is not assumed to be generic. For this purpose, we shall require elements from the theory of complex manifolds and so, for ease of reference throughout this monograph, we collect the facts we need in Section \ref{sec:Almost_Hermitian_manifolds}. In Section \ref{sec:Spincu_structures_Dirac_operators_over_almost_Hermitian_manifolds}, we review the definition of the canonical \spinc structure on an almost Hermitian four-manifold and describe the corresponding Dirac operator.  We proceed in Section \ref{sec:SO3_monopole_equations_over_almost_Hermitian_four-manifolds} to apply this background to develop the form of the non-Abelian monopole equations over \emph{almost K\"ahler} (that is, symplectic) manifolds\label{page:Almost_Kaehler_manifold} and, with a minor modification, to almost Hermitian manifolds more generally. In Section \ref{sec:Witten_dichotomy_non-Abelian_monopoles}, we give an exposition of the proof that solutions to the unperturbed non-Abelian monopole equations over complex K\"ahler surfaces are always type $1$ or $2$ and may thus be viewed as solutions to a projective vortex equation. We also comment on how those solutions may be distinguished based on properties of Hermitian line bundles $L$ and or vector bundles $E$ that twist the canonical \spinc structure $(\rho_\can,W_\can)$ over an almost Hermitian manifold (see the forthcoming Definition \ref{defn:Canonical_spinc_bundles}).

\section{Almost Hermitian manifolds}
\label{sec:Almost_Hermitian_manifolds}
In order to make our conventions and notation clear, in this section we review background material on almost complex manifolds, almost Hermitian manifolds, almost K\"ahler (symplectic) manifolds, complex manifolds, and K\"ahler manifolds from Cannas da Silva \cite{Cannas_da_Silva_lectures_on_symplectic_geometry}, Griffiths and Harris \cite{GriffithsHarris}, H\"ormander \cite{Hormander_introduction_complex_analysis_several_variables}, Huybrechts \cite{Huybrechts_2005}, Kobayashi and Nomizu \cite{Kobayashi_Nomizu_v2}, and Wells \cite{Wells3}. 

\subsection{Hodge star, Lefshetz, and dual Lefshetz operators}
\label{subsec:Hodge_star_Lefshetz_dual_Lefshetz_operators}
We shall follow the exposition and notation of Huybrechts \cite{Huybrechts_2005} for the sake of consistency. Let $(X,g)$ be a smooth Riemannian manifold of dimension $d$. The \emph{Hodge star operator}
\label{page:Hodge_star}
\[
  \star:\Omega^\bullet(X) \to \Omega^\bullet(X)
\]
is defined in the usual way (see Huybrechts \cite[Section 1.2, p. 32]{Huybrechts_2005}) by the relation
\[
  \alpha\wedge\star\beta = g(\alpha,\beta)\,d\vol = \langle\alpha,\beta\rangle_{\Lambda^\bullet(X)}\,d\vol,
  \quad\text{for all } \alpha, \beta \in \Omega^\bullet(X),
\]
where we write $\Omega^\bullet(X) = \oplus_{k=0}^d \Omega^k(X)$ and $\Lambda^\bullet(X) = \oplus_{k=0}^d \Lambda^k(X)$ (pointwise orthogonal direct sums), with $\Omega^k(X) := C^\infty(\Lambda^k(X))$ and $\Lambda^k(X) := \wedge^k(T^*X)$ and $\Omega^\bullet(X) := C^\infty(\Lambda^\bullet(X))$ and $\Lambda^\bullet(X) := \wedge^\bullet(T^*X)$, write $g$ for the Riemannian metric $g$ on $TX$ and all associated tensor product and dual vector bundles, and let $d\vol \in \Omega^d(X)$ denote the volume form defined by $g$.

Suppose that $(X,J)$ is an almost complex manifold of even real dimension $d = 2n$. As usual (see Huybrechts \cite[Section 1.3, pp. 42--43]{Huybrechts_2005} or Wells \cite[Chapter I, Section 3, pp. 31--32]{Wells3}), we let $T^{1,0}X$ and $T^{0,1}X$ denote the subbundles of $T_\CC X := TX\otimes_\RR\CC$ defined by the $\pm i$ eigenspaces, respectively, of the $\CC$-linear extension $J \in \End(T_\CC X)$ of $J\in \End(TX)$, let $T^*X^{1,0}$ and $T^*X^{0,1}$ denote the dual bundles of $T^{1,0}X$ and $T^{0,1}X$, respectively, write
\[
  \wedge^{p,q}(T^*X) := \wedge^p(T^*X^{1,0})\otimes \wedge^q(T^*X^{0,1}),
\]
which is customarily abbreviated by $\Lambda^{p,q}(X)$, and denote
\[
  \Omega^{p,q}(X):=C^\infty(\Lambda^{p,q}(X)).
\]
(Recall from Wells \cite[Chapter I, Section 3, p. 32, lines 6--8]{Wells3} that there is a $\CC$-linear isomorphism $(TX,J) \cong T^{1,0}X$, where $(TX,J)$ is the complex vector bundle constructed from $TX$ by defining scalar multiplication by $i$ on the fibers of $TX$ via $J$.) We usually write $\Lambda^0(X,\CC)$, or $\Lambda^0(X)$ when complex coefficients are understood, instead of $\Lambda^{0,0}(X)$; similarly, we usually write $\Omega^0(X,\CC)$, or $\Omega^0(X)$ when complex coefficients are understood, instead of $\Omega^{0,0}(X)$.

Suppose further that $(X,g,J)$ is a smooth almost Hermitian manifold, so now $J\in \End(TX)$ is a $g$-orthogonal smooth almost complex structure, and $\omega = g(\cdot,J\cdot)$ is the fundamental two-form as in \eqref{eq:Fundamental_two-form}. The \emph{Lefschetz operator} \label{page:Lefschetz_operator} is given by (see Huybrechts \cite[Remark  1.2.19, p. 31]{Huybrechts_2005})
\[
  L:\Omega^\bullet(X) \ni \alpha \mapsto \omega\wedge\alpha \in \Omega^\bullet(X)
\]
and has a complex linear extension (see Huybrechts \cite[Definition 1.2.18, p. 31]{Huybrechts_2005})
\[
  L:\Omega^\bullet(X,\CC) \ni \alpha \mapsto \omega\wedge\alpha \in \Omega^\bullet(X,\CC),
\]
where $\Omega^\bullet(X,\CC) = \Omega^\bullet(X)\otimes\CC$ and $\Lambda^\bullet(X,\CC) = \Lambda^\bullet(X)\otimes\CC$. The \emph{dual Lefshetz operator}
\[
  \Lambda:\Omega^\bullet(X) \to \Omega^\bullet(X)
\]
is defined (see Huybrechts \cite[Definition 1.2.21, p. 33]{Huybrechts_2005}) by the relation
\[
  g(\Lambda\alpha,\beta) = g(\alpha,L\beta),
  \quad\text{for all } \alpha, \beta \in \Omega^\bullet(X).
\]
The dual Lefshetz operator also has a complex linear extension
\[
  \Lambda:\Omega^\bullet(X,\CC) \to \Omega^\bullet(X,\CC)
\]
and obeys, with respect to the Hermitian metric $h$ on $TX\otimes\CC$ and all associated tensor product and dual vector bundles,
\[
  h(\Lambda\alpha,\beta) = h(\alpha,L\beta),
  \quad\text{for all } \alpha, \beta \in \Omega^\bullet(X,\CC).
\]
The Hermitian metric $h$ on $TX\otimes\CC$ extends $g$ on $TX$ as in Huybrechts \cite[Section 1.2, p. 30]{Huybrechts_2005}. The Hodge star operator has a complex linear extension (see Huybrechts \cite[Paragraph prior to Lemma 1.2.24, p. 33]{Huybrechts_2005})
\[
  \star:\Omega^\bullet(X,\CC) \to \Omega^\bullet(X,\CC),
\]
and obeys
\begin{multline}
\label{eq:Huybrechts_page_33_inner_product_complex-valued_forms}
\alpha\wedge\star\bar{\beta}
= h(\alpha,\beta)\,d\vol
= \langle\alpha,\beta\rangle_{\Lambda^\bullet(X,\CC)}\,d\vol,
\\
\text{for all } \alpha,\beta \in\Omega^\bullet(X,\CC)\quad\text{(see \cite[Section 1.2, p. 33]{Huybrechts_2005})}.
\end{multline}
The properties of $\star$ for forms over real manifolds are implied by Huybrechts \cite[Proposition 1.2.20, p. 32]{Huybrechts_2005} on the level of cotangent spaces, including the following:
\begin{gather}
\label{eq:Huybrechts_proposition_1-2-20_i}
\star 1 = d\,\vol \quad\text{(by \cite[Proposition 1.2.20 (i), p. 32]{Huybrechts_2005})},  
\\  
\label{eq:Huybrechts_proposition_1-2-20_ii}
\langle\alpha, \star\beta\rangle_{\Lambda^k(X)} = (-1)^{k(d-k)}\langle\star\alpha, \beta\rangle_{\Lambda^k(X)}, \quad\text{for all } \alpha\in\Omega^k(X) \text{ and } \beta\in\Omega^{d-k}(X)
\\
\notag\qquad\text{(by \cite[Proposition 1.2.20 (ii), p. 32]{Huybrechts_2005})},
\\
\label{eq:Huybrechts_proposition_1-2-20_iii}
\star^2 = \star\circ\star = (-1)^{k(d-k)} \quad\text{on } \Omega^k(X) \quad\text{(by \cite[Proposition 1.2.20 (iii), p. 32]{Huybrechts_2005})},
\\
\label{eq:Huybrechts_proposition_1-2-20_Bonus}
\star d\vol = 1
\quad\text{(by \eqref{eq:Huybrechts_proposition_1-2-20_i} and \eqref{eq:Huybrechts_proposition_1-2-20_iii}) }
\\
\label{eq:Huybrechts_page_33_and_lemma_1-2-24}
\star:\Omega^{p,q}(X) \ni \varphi \mapsto \star\varphi \in  \Omega^{n-q,n-p}(X)
\\
\notag\qquad\text{is $\CC$-linear (see \cite[Lemma 1.2.24 (ii), p. 33]{Huybrechts_2005})}.
\end{gather}
The operators $L$ and $\Lambda$ enjoy the following properties:
\begin{gather}
\label{eq:Huybrechts_remark_1-2-19-ii}  
L:\Omega^{p,q}(X) \ni \alpha \mapsto \omega\wedge\alpha \in \Omega^{p+1,q+1}(X) \quad\text{(by \cite[Remark 1.2.19 (ii), p. 31]{Huybrechts_2005})},
\\
\label{eq:Huybrechts_definition_1-2-21_and_lemmas_1-2-23_and_24}
\Lambda = L^* = \star^{-1}\circ L\circ\star:\Omega^{p,q}(X) \to \Omega^{p-1,q-1}(X)
\\
\notag\qquad\text{(by \cite[Definition 1.2.21 and Lemmas 1.2.23 and 1.2.24, p. 33]{Huybrechts_2005})},
\end{gather}
We recall (see Huybrechts \cite[Proposition 1.2.24 (i). p. 33]{Huybrechts_2005} or Wells \cite[Chapter V, Proposition 2.2, p. 165]{Wells3}) that for any closed, smooth almost Hermitian manifold $(X,g,J)$ and pairs of non-negative integers $(p,q)\neq (r,s)$, one has the following
pointwise orthogonal property:
\begin{equation}
  \label{eq:Almost_hermitian_manifold_pq_forms_pointwise_orthogonal_rs_forms_unless_pq_equals_rs}
  \langle\kappa,\mu\rangle_{\Lambda^{\bullet,\bullet}(X)} = 0 \quad\text{for all } \kappa \in \Omega^{p,q}(X) \text{ and } \mu \in \Omega^{r,s}(X).
\end{equation}
Suppose that we are also given a Hermitian vector bundle $(E,\langle\cdot,\cdot\rangle_E)$ over $X$. The Hermitian metric $E\cong E^*$ is interpreted as a $\CC$-antilinear isomorphism, $s \mapsto \langle\cdot,s\rangle_E$, and $E^*=\Hom(E,\CC)$ is the complex dual vector bundle (following the notation implied by Huybrechts in \cite[Lemma 1.2.6, p. 26]{Huybrechts_2005}). Then
\begin{gather}
\label{eq:Huybrechts_definition_4-1-11_and_example_4-1-2}
\langle\varphi\otimes s, \psi\otimes t\rangle_{\Lambda^{p,q}(X)\otimes E}
= \langle\varphi, \psi\rangle_{\Lambda^{p,q}(X)} \langle s,t\rangle_E,
\\
\notag\qquad\text{for all } \varphi,\psi \in \Omega^{p,q}(X) \text{ and } s, t \in \Omega^0(E)
\\
\notag\qquad\text{(see \cite[Example 4.1.2, p. 166, and Definition 4.1.11, p. 169]{Huybrechts_2005})}.
\end{gather}
We will often employ the abbreviations
\begin{equation}
\label{eq:VectorBundleValued(pq)forms}
\Lambda^k(E):=\Lambda^k(X)\otimes E \quad\text{and}\quad \Lambda^{p,q}(E):=\Lambda^{p,q}(X)\otimes E,
\end{equation}
and write $\Omega^k(E) = \Omega^0(\Lambda^k(E))$ and $\Omega^{p,q}(E) = \Omega^0(\Lambda^{p,q}(E))$. We also note the following identities
\begin{gather}
L:\Omega^{p,q}(E) \ni \alpha \mapsto \omega\wedge\alpha \in \Omega^{p+1,q+1}(E),
\\
\Lambda = L^* = \star^{-1}\circ L\circ\star:\Omega^{p,q}(E) \to \Omega^{p-1,q-1}(E),
\end{gather}
that are implied by properties of $L$ and $\Lambda$ for the case $E = \underline{X} = X\times\CC$.

\subsection{Nijenhuis tensor and components of the exterior derivative}
\label{subsec:Nijenhuis_tensor_components_exterior_derivative}
Recall that the \emph{Nijenhuis tensor} $N$ is defined by (see Salamon \cite[Equation (3.5), p. 62]{SalamonSWBook})
\begin{equation}
\label{eq:Nijenhuis_tensor}
  N(v,w) = [v,w] + J[Jv,w] + J[v, Jw] - [Jv, Jw], \quad\text{for all } v, w \in C^\infty(TX),
\end{equation}
and that (see Salamon \cite[Exercise 3.9, p. 62]{SalamonSWBook}) this expression defines a tensor, that is, the value of $N(v,w)$ at $x\in X$ depends only on the values of $v_x$ and $w_x$. Moreover (see Salamon \cite[Exercise 3.10, p. 62]{SalamonSWBook}), the tensor $N$ is skew-symmetric, so
\[
  N \in C^\infty(\wedge^2(T^*X)\otimes TX),
\]  
and $J$-anti-linear in each variable, so
\[
  N(Jv,w) = -JN(v,w), \quad\text{for all } v, w \in C^\infty(TX).
\]
We may also write
\[
  N \in C^\infty(\Hom(T^*X, \wedge^2(T^*X)\otimes TX)),
\]
where we use the canonical identification $TX = T^{**}X$ to give the isomorphism 
\[
  C^\infty(\wedge^2(T^*X)\otimes TX)) \cong C^\infty(\Hom(T^*X, \wedge^2(T^*X)\otimes TX)),
\]
If $\alpha \in \Omega^{1,0}(X)$, then (see Salamon \cite[Exercise 3.11, p. 62]{SalamonSWBook})
\begin{equation}
\label{eq:Nijenhuis_tensor_(0,2)_component_exterior_derivative_(1,0)_forms}
 (d\alpha)^{0,2} = \frac{1}{4}\alpha\circ N.
\end{equation}
If $\beta \in \Omega^{0,1}(X)$, then $\bar\beta \in \Omega^{1,0}(X)$ (see Wells \cite[Chapter I, Section 3, p. 32]{Wells3}). (Note that $J\beta(v) = \beta(Jv) = -i\beta(v)$ for all $v \in C^\infty(T^{0,1}X)$ and thus $J\beta = -i\beta$ if $Jv=-iv$ for $v\in C^\infty(T^{0,1}X)$.) Hence
\[
  (d\bar\beta)^{0,2} = \frac{1}{4}\bar\beta\circ N.
\]
From Wells \cite[Chapter I, Proposition 3.6, p. 33]{Wells3}, we have $(d\bar\beta)^{0,2} = \overline{(d\beta)^{2,0}}$ and thus
\begin{equation}
\label{eq:Nijenhuis_tensor_(2,0)_component_exterior_derivative_(0,1)_forms}
 (d\beta)^{2,0} = \frac{1}{4}\beta\circ N.
\end{equation}
See also Bryant \cite[Section 2]{Bryant_2006}, Donaldson \cite[Section 4]{DonYangMillsInvar}, McDuff and Salamon \cite[Appendix C.7]{McDuffSalamonSympTop3}, and Nicolaescu \cite[Section 1.4]{NicolaescuSWNotes}.

More generally, for forms of higher degree, one may proceed as follows. Following Salamon \cite[Exercise 3.12, p. 63]{SalamonSWBook}, for any integer $r\geq 0$ one defines
\[
  \iota(N):\Omega^r(X,\CC) \to \Omega^{r+1}(X,\CC)
\]
for $\alpha \in \Omega^1(X)$ by $\iota(N)\alpha = \alpha\circ N$ and in general by
\[
  \iota(N)(\sigma\wedge\tau) := (\iota(N)\sigma)\wedge\tau + (-1)^{\deg\sigma}\sigma\wedge\iota(N)\tau
\]
for $\sigma\in\Omega^k(X,\CC)$ and $\tau\in\Omega^l(X,\CC)$. In particular,
\[
  d\tau - \partial\tau - \bar\partial\tau = \frac{1}{4}\iota(N)\tau \in \Omega^{p+2,q-1}(X)\oplus \Omega^{p-1,q+2}(X)
\]
for all $\tau \in \Omega^{p,q}(X)$. Following Cirici and Wilson \cite[Section 2]{Cirici_Wilson_2020_harmonic}, we denote
\[
  \mu\tau = d^{2,-1}\tau \quad\text{and}\quad \bar\mu\tau = d^{-1,2}\tau,
\]
so that
\[
  \mu\tau =  \frac{1}{4}\pi_{p+2,q-1}(\iota(N)\tau) \quad\text{and}\quad \bar\mu\tau = \frac{1}{4}\pi_{p-1,q-2}(\iota(N)\tau).
\]
In particular, for $\alpha \in \Omega^{1,0}(X)$ and $\beta \in \Omega^{0,1}(X)$, we have
\begin{align*}
  \mu\alpha = 0 \quad&\text{and}\quad \bar\mu\alpha = \frac{1}{4}\alpha\circ N \in \Omega^{0,2}(X),
  \\
  \mu\beta = \frac{1}{4}\beta\circ N \in \Omega^{2,0}(X) \quad&\text{and}\quad \bar\mu\beta = 0,
\end{align*}
in agreement with Salamon \cite[Exercise 3.11, p. 62]{SalamonSWBook}.

Given a connection $A$ on $E$, the corresponding exterior covariant derivative operator (see, for example, Donaldson and Kronheimer \cite[Equation (2.1.12), p. 35]{DK})
\[
  d_A:\Omega^{p,q}(E) \to \Omega^{p+q+1}(E)
\]
thus splits into four components.
\begin{subequations}
  \label{eq:d_A_components_almost_complex_manifold}
  \begin{align}
    \label{eq:del_A_old}
    \partial_A = d_A^{1,0} &= \pi_{p+1,q}d_A:\Omega^{p,q}(E) \to \Omega^{p+1,q}(E),
    \\
    \label{eq:mu_A}
    \mu_A = d_A^{2,-1} &= \pi_{p+2,q-1}d_A:\Omega^{p,q}(E) \to \Omega^{p+2,q-1}(E),
    \\
    \label{eq:mu_bar_A}
    \bar\mu_A = d_A^{-1,2} &= \pi_{p-1,q+2}d_A:\Omega^{p,q}(E) \to \Omega^{p-1,q+2}(E),
    \\
    \label{eq:del_bar_A_old}
    \bar\partial_A = d_A^{0,1} &= \pi_{p,q+1}d_A:\Omega^{p,q}(E) \to \Omega^{p,q+1}(E),
\end{align}  
\end{subequations}
so that (see Donaldson \cite[p. 34]{DonYangMillsInvar})
\begin{equation}
  \label{eq:d_A_sum_components_almost_complex_manifold}
  d_A = \partial_A + \mu_A + \bar\mu_A + \bar\partial_A \quad\text{on } \Omega^{p,q}(E).
\end{equation}
In particular, when $N\equiv 0$, one has
\begin{equation}
  \label{eq:d_A_sum_components_almost_complex_manifold_integrable}
  d_A = \partial_A + \bar\partial_A \quad\text{on } \Omega^{p,q}(E),
\end{equation}
as usual over complex manifolds.

\section{Canonical spin${}^c$ and spin${}^u$ structures and Dirac operators over almost Hermitian manifolds}
\label{sec:Spincu_structures_Dirac_operators_over_almost_Hermitian_manifolds}
We review the definition of the canonical \spinc structure on an almost Hermitian four-manifold in Section \ref{subsec:Spinc_and_spinu_structures_over_almost_Hermitian_manifolds} and describe the Dirac operator for this \spinc structure in Section \ref{subsec:Dirac_operator_almost_Hermitian}.

\subsection{Canonical spin${}^c$ and spin${}^u$ structures over almost Hermitian manifolds}
\label{subsec:Spinc_and_spinu_structures_over_almost_Hermitian_manifolds}
We review Clifford multiplication over almost Hermitian four-manifolds $(X,g,J)$, first for the canonical \spinc structure and then for any \spinc structure. We closely follow Salamon \cite{SalamonSWBook}, though we also rely on Donaldson \cite{DonSW}, Kotschick \cite{KotschickSW}, Morgan \cite{MorganSWNotes}, Nicolaescu \cite{NicolaescuSWNotes}, Taubes \cite[Section 1]{TauSymp}, and Witten \cite{Witten}. We begin with a discussion of certain identities for the interior product of complex differential forms that we shall need. For a smooth almost Hermitian manifold $(X,g,J)$ of dimension $2n$, smooth form $\alpha \in \Omega^{0,k}(X)$, and integers $\ell\geq k\geq 0$, we recall that left multiplication by $\alpha$ and interior multiplication,
\begin{align}
  L_\alpha:\Omega^{0,\ell-k}(X) \ni \varpi
  &\mapsto \alpha\wedge\varpi \in \Omega^{0,\ell}(X),
  \\
  \iota(\bar\alpha):\Omega^{0,\ell}(X) \ni \chi
  &\mapsto \iota(\bar\alpha)\chi \in \Omega^{0,\ell-k}(X),
\end{align}
are related by (see Salamon \cite[Lemma 3.5, p. 60]{SalamonSWBook}) 
\[
  \langle L_\alpha\varpi, \chi\rangle_{\Lambda^{0,\ell}(X)}
  = \langle \alpha\wedge\varpi, \chi\rangle_{\Lambda^{0,\ell}(X)}
  = 2^k(-1)^{k(k-1)/2}\langle\varpi, \iota(\bar\alpha)\chi\rangle_{\Lambda^{0,\ell-k}(X)},
\]
where $\iota(\bar\alpha)$, for $\alpha = \sum \alpha_I e_I''$, is defined for $\ell\geq 1$ by
\begin{equation}
  \label{eq:Salamon_3-3}
  \iota(\bar\alpha)\chi := \sum \bar\alpha_I\iota(e_I)\chi, \quad e_I'' = e_{i_1}\wedge\cdots\wedge e_{i_\ell}'',
\end{equation}
where the sum runs over all multi-indices $I = \{i_1,\ldots,i_\ell\}$ with $i_1<\cdots<i_\ell$ and
\[
  \iota(e_I)\chi := \iota(e_{i_1})\cdots\iota(e_{i_\ell})\chi,
\]
and $\iota(v):\Omega^\ell(X,\CC)\to\Omega^{\ell-1}(X,\CC)$ denotes interior product,
\begin{equation}
\label{eq:Interior_product}
\iota(v)\chi(v_2,\ldots v_\ell) := \chi(v,v_2,\ldots v_\ell),
\quad\text{for all } v,v_2\ldots,v_\ell \in C^\infty(TX),
\end{equation}  
and $\{e_1,Je_1,\ldots,e_n,J_ne_n\}$ is a local orthonormal frame for $TX$ and $\{e_1^*,Je_1^*,\ldots,e_n^*,J_ne_n^*\}$ is the corresponding local orthonormal frame for $T^*X$ defined by (see Salamon \cite[Section 3.1, pp. 58--59]{SalamonSWBook}) for any $v \in C^\infty(TX)$ by 
\begin{multline*}
  v^* := g(\cdot,v) \in \Omega^1(X,\RR)
  \\
  v' := 2\pi_{1,0}v^* = \langle \cdot,v\rangle_{T_\CC X} \in \Omega^{1,0}(X),
  \\
  \quad\text{and}\quad
  v'' := 2\pi_{0,1}v^* = \langle v,\cdot\rangle_{T_\CC X} = \overline{v'} \in \Omega^{0,1}(X),
\end{multline*}
and, for any $\theta \in \Omega^1(X,\RR)$,
\[
  \pi_{1,0}\theta := \frac{1}{2}(\theta + iJ\theta) \quad\text{and}\quad \pi_{0,1}\theta := \frac{1}{2}(\theta - iJ\theta).
\]
giving $J\pi_{1,0}\theta = -i\theta$ and $J\pi_{0,1}\theta = i\theta$. (Note that Salamon \cite[Section 3.1, p. 55]{SalamonSWBook} adopts the \emph{opposite convention} for the Hermitian inner product on a complex Hilbert space $\fH$, taking $\langle x,y \rangle_\fH$ to be complex linear in $y$ and complex anti-linear in $x$.)

Recall that the Hermitian inner product on $\Lambda_\CC^\bullet(X) := \wedge^\bullet(T_\CC^*X)$ is defined by the relation (see Huybrechts \cite[Section 1.2, p. 33]{Huybrechts_2005})
\begin{equation}
  \label{eq:Huybrechts_Hermitian_inner_product_complex_forms_p_33}
  \langle \kappa, \mu \rangle_{\Lambda_\CC^\bullet(X)}\cdot\vol := \kappa \wedge \star \bar\mu,
  \quad\text{for all } \kappa, \mu \in \Omega^\bullet(X,\CC).
\end{equation}
Therefore, given $\alpha \in \Omega^{k,l}(X)$, we have (using $\star 1 = \vol$ by \cite[Proposition 1.2.20 (i)]{Huybrechts_2005})
\[
  \langle \alpha\wedge\varpi, \chi\rangle_{\Lambda_\CC^\bullet(X)}\cdot\vol
  =
  \langle \alpha\wedge\varpi, \chi\rangle_{\Lambda_\CC^\bullet(X)}\cdot(\star 1)
  =
  \alpha\wedge\varpi\wedge\star\bar\chi,
\]
which gives (using $\star^2 1 = 1$ by Huybrechts \cite[Proposition 1.2.20 (iii), p. 32]{Huybrechts_2005})
\begin{multline}
  \label{eq:Huybrechts_Hermitian_inner_product_complex_forms_p_33_Hodge_star}
  \langle\alpha\wedge\varpi, \chi\rangle_{\Lambda^{p,q}(X)}
  =
  \star(\alpha\wedge\varpi \wedge \star\bar\chi),
  \\
  \text{for all } \alpha \in \Omega^{k,l}(X), \varpi \in \Omega^{p,q}(X),
  \text{ and } \chi \in \Omega^{p+k,q+l}(X).
\end{multline}
If $\Lambda_\alpha:\Omega^{p+k,q+l}(X)\to\Omega^{p,q}(X)$ denotes the adjoint of left exterior multiplication by $\alpha \in \Omega^{k,l}(X)$,
\[
  L_\alpha:\Omega^{p,q}(X) \ni \varpi \mapsto \alpha\wedge\varpi \in \Omega^{p+k,q+l}(X),
\]
with respect to the Hermitian inner product on $\Lambda_\CC^\bullet(X)$, that is,
\begin{equation}
  \label{eq:Left_multiplication_of_pq_forms_by_kl_form_Hermitian_adjoint}
  \langle L_\alpha\varpi, \chi\rangle_{\Lambda^{p+k,q+l}(X)}
  =
  \langle\varpi, \Lambda_\alpha\chi\rangle_{\Lambda^{p,q}(X)},
  \quad\text{for all } \varpi \in \Omega^{p,q}(X), \chi \in \Omega^{p+k,q+l}(X),
\end{equation}
then, noting that $\deg\bar\chi=(q+l,p+k)$ and $\deg\star\bar\chi=(n-p-k,n-q-l)$ and thus $\deg(\alpha \wedge \star\bar\chi)=(n-p,n-q)$ and $\star^2=(-1)^{(2n-p-q)(p+q)}$ on $\Omega^{n-p,n-q}(X)$ by Huybrechts \cite[Proposition 1.2.20 (iii), p. 32]{Huybrechts_2005},
\begin{align*}
  \langle\varpi, \Lambda_\alpha\chi\rangle_{\Lambda^{p,q}(X)}
  &= \langle \alpha\wedge\varpi, \chi\rangle_{\Lambda^{p+k,q+l}(X)} \quad\text{(by \eqref{eq:Left_multiplication_of_pq_forms_by_kl_form_Hermitian_adjoint})}
  \\
  &= \star(\alpha\wedge\varpi \wedge \star\bar\chi) \quad\text{(by \eqref{eq:Huybrechts_Hermitian_inner_product_complex_forms_p_33_Hodge_star})}
  \\
  &= (-1)^{(k+l)(p+q)}\star(\varpi \wedge \alpha \wedge \star\bar\chi)
    \quad\text{(as $\deg\alpha = (k,l)$ and $\deg\varpi=(p,q)$)}
  \\
  &= (-1)^{(k+l)(p+q)}(-1)^{(2n-p-q)(p+q)}\star(\varpi \wedge \star^2(\alpha \wedge \star\bar\chi))
  \\
  &\qquad\text{(as $\deg(\alpha \wedge \star\bar\chi)=(n-p,n-q)$ and $\star^2=(-1)^{(2n-p-q)(p+q)}$)}
  \\
  &= (-1)^{(k+l)(p+q)}(-1)^{(p+q)^2}\left\langle\varpi, \overline{\star(\alpha \wedge \star\bar\chi)} \right\rangle_{\Lambda^{p,q}(X)}
   \quad\text{(by \eqref{eq:Huybrechts_Hermitian_inner_product_complex_forms_p_33_Hodge_star})}    
  \\
  &= (-1)^{(k+l)(p+q)}(-1)^{(p+q)^2}\left\langle\varpi, \star(\bar\alpha \wedge \star\chi) \right\rangle_{\Lambda^{p,q}(X)}
  \\
  &\qquad \text{(by $\CC$-linearity \cite[Section 1.2, p. 33]{Huybrechts_2005} of Hodge $\star:\Lambda_\CC^\bullet(X)\to \Lambda_\CC^\bullet(X)$)},
\end{align*}
and therefore\footnote{Compare the identity given by Warner \cite[Equation (2.11.2), p. 61 and Exercise 2.14, p. 80]{Warner} in the case of a Riemannian manifold and $\deg\alpha=1$.}
\begin{equation}
  \label{eq:Left_multiplication_of_p+kq+l_forms_by_kl_form_Hermitian_adjoint_formula}
  \Lambda_\alpha\chi
  =
  (-1)^{(k+l)(p+q)}(-1)^{(p+q)^2} \star(\bar\alpha \wedge \star\chi),
  \quad\text{for all } \chi \in \Omega^{p+k,q+l}(X).
\end{equation}
Note that if $k+l\equiv 1 \pmod{2}$, then the identity \eqref{eq:Left_multiplication_of_p+kq+l_forms_by_kl_form_Hermitian_adjoint_formula} simplifies to give
\[
  \Lambda_\alpha\chi = \star(\bar\alpha\wedge*\chi), \quad\text{for all } \chi \in \Omega^{p+k,q+l}(X).
\]
We recall the important

\begin{defn}[Canonical \spinc structure]
\label{defn:Canonical_spinc_bundles}
(See Donaldson \cite[Equation (15)]{DonSW}, Kotschick \cite[Fact 2.1]{KotschickSW}, Salamon \cite[Lemma 4.52, p. 141]{SalamonSWBook}, and Taubes \cite[Section 1]{TauSymp}.)
Let $(X,g,J)$ be an almost Hermitian four-manifold. Then the \emph{canonical \spinc structure} $\fs_{\can}=(\rho_\can,W_\can)$ on $X$ is defined by the bundle
\begin{equation}
  \label{eq:Canonical_spinc_bundles}
  W_\can^+ := \Lambda^{0,0}(X) \oplus \Lambda^{0,2}(X), \quad W_\can^- := \Lambda^{0,1}(X),
  \quad W_\can = W_\can^+\oplus W_\can^-,
\end{equation}
and the map $\rho_\can:TX\to \Hom(W_\can)$ defined for all $Y\in C^\infty(TX)$ and $\phi \in \Omega^{0,\bullet}(X)$ by
\begin{equation}
  \label{eq:Canonical_Clifford_multiplication}
    \rho_\can(Y)\phi = \frac{1}{\sqrt{2}}Y''\wedge\phi - \sqrt{2}\iota(Y)\phi,
\end{equation}
where $\iota(Y):\Omega^{0,q}(X) \to \Omega^{0,q-1}(X)$ denotes interior product for $q\geq 1$ and
\[
  Y'' = 2\pi_{0,1}Y^* = \langle Y,\cdot\rangle_{T_\CC X} \in \Omega^{0,1}(X).
\]
\end{defn}
The \emph{canonical line bundle} of an almost complex manifold $(X,J)$ of real dimension $2n$ is\footnote{If $X$ is a complex manifold, then the holomorphic tangent bundle $\sT_X$ (see Huybrechts \cite[Definition 2.2.14, p. 71]{Huybrechts_2005}) is naturally isomorphic as a complex vector bundle to $T^{1,0}X$ by \cite[Proposition 2.6.4 (ii), p. 104]{Huybrechts_2005}.} (compare Donaldson \cite[Section 4, p. 56]{DonSW}, Fine \cite[Definition 1.6]{Fine_2012}, Griffiths and Harris \cite[Chapter 1, Section 2, p. 146]{GriffithsHarris}, or Wells \cite[Chapter VI, Section 1, p. 218]{Wells3})
\begin{equation}
  \label{eq:DefineCanonicalLineBundle}
  K_X := \Lambda^{n,0}(X) = \wedge^n(T^{1,0}X),
\end{equation}
and the \emph{anti-canonical line bundle} is its dual line bundle (see Fine \cite[Definition 1.6]{Fine_2012} or
Gauduchon \cite[Section 2, Remark 5, p. 275]{Gauduchon_1997}),
\begin{equation}
  \label{eq:DefineAntiCanonicalLineBundle}
  K_X^* := \Hom_\CC(K_X,\CC).
\end{equation}
Then by Kotschick \cite[Fact 2.1]{KotschickSW}, the first Chern class of $\fs_{\can}$ is
\begin{equation}
\label{eq:c1_of_CanonicalSpinc}
c_1(\fs_{\can}) = -c_1(K_X) = c_1(X).
\end{equation}
According to Salamon \cite[Lemma 3.4, p. 141]{SalamonSWBook},
\[
  2\langle\iota(Y)\phi, \chi\rangle_{\Lambda^{0,q-1}(X)} = \langle \phi, Y''\wedge\chi\rangle_{\Lambda^{0,q}(X)},
  \quad\text{for all } \chi \in \Omega^{0,q-1}(X) \text{ and } \phi \in \Omega^{0,q}(X).
\]
But then, as in the derivation of \eqref{eq:Left_multiplication_of_p+kq+l_forms_by_kl_form_Hermitian_adjoint_formula},
\begin{align*}
  2\langle\chi, \iota(Y)\phi\rangle_{\Lambda^{0,q-1}(X)}
  &= \langle Y''\wedge\chi, \phi\rangle_{\Lambda^{0,q}(X)}
  \\
  &= \star (Y''\wedge\chi \wedge \star\bar\phi)
    \quad\text{(by \eqref{eq:Huybrechts_Hermitian_inner_product_complex_forms_p_33_Hodge_star})}
  \\
  &= (-1)^{q-1}\star (\chi \wedge Y''\wedge \star\bar\phi) \quad\text{(as $\deg\chi = q-1$)}
  \\
  &= (-1)^{q-1}(-1)^{(q-1)(2n-q+1)} \star (\chi \wedge \star (\star(Y''\wedge \star\bar\phi)))
  \\
  &= (-1)^{q-1}(-1)^{(q-1)^2} \star (\chi \wedge \star (\star(Y''\wedge \star\bar\phi)))
  \\
  &= (-1)^{2(q-1)} \star (\chi \wedge \star (\star(Y''\wedge \star\bar\phi)))
  \\
  &= \star (\chi \wedge \star (\star(Y''\wedge \star\bar\phi)))  
  \\
  &= \left\langle \chi, \overline{\star(Y''\wedge \star\bar\phi)} \right\rangle_{\Lambda^{0,q-1}(X)}
    \quad\text{(by \eqref{eq:Huybrechts_Hermitian_inner_product_complex_forms_p_33_Hodge_star})}
  \\
  &= \left\langle \chi, \star(\overline{Y''}\wedge \star\phi) \right\rangle_{\Lambda^{0,q-1}(X)},
  \\
  &\qquad \text{(by $\CC$-linearity \cite[Section 1.2, p. 33]{Huybrechts_2005} of Hodge $\star:\Lambda_\CC^\bullet(X)\to \Lambda_\CC^\bullet(X)$)},    
\end{align*}
where to obtain the fourth equality, we used the fact that $\deg\bar\phi = \deg\phi = q$, so $\deg\star\bar\phi = 2n-q$, and $\deg(Y''\wedge \star\bar\phi) = 2n-q+1$, and thus $\star^2=(-1)^{(q-1)(2n-q+1)}$ on $\Omega^{2n-q+1}(X)$ by Warner \cite[Equation (6.1.1)]{Warner}. Therefore,
\begin{equation}
  \label{eq:Salamon_lemma_3-4_formula}
  2\iota(Y)\phi = \star(Y'\wedge \star\phi), \quad\text{for all } Y\in C^\infty(TX)
  \text{ and } \phi \in \Omega^{0,q}(X).
\end{equation}
Hence, the Clifford multiplication \eqref{eq:Canonical_Clifford_multiplication} may be rewritten as
\begin{equation}
  \label{eq:Canonical_Clifford_multiplication_Hodge_star}
  \rho_\can(Y)\phi = \frac{1}{\sqrt{2}}\left(Y''\wedge\phi - \star(Y'\wedge \star\phi)\right),
  \quad\text{for all } Y\in C^\infty(TX) \text{ and } \phi \in \Omega^{0,q}(X).
\end{equation}
Note that $Y^* = \pi_{1,0}Y^* + \pi_{0,1}Y^* = \frac{1}{2}(Y'+Y'')$. For $a \in \Omega^1(X,\RR)$, we apply the formula \eqref{eq:Canonical_Clifford_multiplication_Hodge_star} with $Y \in C^\infty(TX)$ uniquely defined by $a = g(\cdot,Y) = Y^*$ and $a = a'+a'' = \frac{1}{2}(Y'+Y'') \in \Omega^{1,0}(X)\oplus\Omega^{0,1}(X)$ to give
\begin{equation}
  \label{eq:Canonical_Clifford_multiplication_Hodge_star_one_forms}
  \rho_\can(a)\phi = \sqrt{2}\left(a''\wedge\phi - \star(a'\wedge \star\phi)\right),
  \quad\text{for all } a \in \Omega^1(X,\RR) \text{ and } \phi \in \Omega^{0,q}(X).
\end{equation}  
Choosing $\phi = (\sigma,\tau) \in \Omega^0(W^+) = \Omega^{0,0}(X)\oplus\Omega^{0,2}(X)$ and $\nu \in \Omega^0(W^-) = \Omega^{0,1}(X)$ and specializing \eqref{eq:Canonical_Clifford_multiplication} to the case of the positive and negative spinors in \eqref{eq:Canonical_spinc_bundles} gives the identities\footnote{Morgan \cite[p. 109]{MorganSWNotes} gives the expression \eqref{eq:Canonical_Clifford_multiplication_positive_spinors} but omits the factor $\sqrt{2}$ in the first term. Kotschick \cite[Fact 2.1]{KotschickSW} gives expressions for Clifford multiplication that agree with \eqref{eq:Canonical_Clifford_multiplication_positive_spinors} and \eqref{eq:Canonical_Clifford_multiplication_negative_spinors}.}
\begin{subequations}
  \begin{align}
    \label{eq:Canonical_Clifford_multiplication_positive_spinors}
    \rho_\can(a)(\sigma,\tau) &= \sqrt{2}\left(a''\wedge\sigma - \star(a'\wedge\star\tau)\right)
                                \in \Omega^{0,1}(X),
  \\
    \label{eq:Canonical_Clifford_multiplication_negative_spinors}
    \rho_\can(a)\nu &= \sqrt{2}\left(- \star(a'\wedge\star\nu) + a''\wedge\nu\right)
                      \in \Omega^{0,0}(X)\oplus\Omega^{0,2}(X).
  \end{align}  
\end{subequations}
To deduce \eqref{eq:Canonical_Clifford_multiplication_positive_spinors} from \eqref{eq:Canonical_Clifford_multiplication_Hodge_star_one_forms}, observe that for $\sigma \in \Omega^{0,0}(X)$ we have  $\star\sigma \in \Omega^{2,2}(X)$ and $a'\wedge\star\sigma \in \Omega^{3,2}(X)$, that is, $a'\wedge\star\sigma = 0$, and thus by \eqref{eq:Canonical_Clifford_multiplication_Hodge_star_one_forms},
\[
  \rho_\can(a)\sigma = \sqrt{2}a''\wedge\sigma \in \Omega^{0,1}(X).
\]
For $\tau \in \Omega^{0,2}(X)$ we have $a''\wedge\tau \in \Omega^{0,3}(X)$, and so $a''\wedge\tau=0$, while $\star\tau \in \Omega^{0,2}(X)$ and $a' \wedge \star\tau \in \Omega^{1,2}(X)$ and $\star(a' \wedge \star\tau)  \in \Omega^{0,1}(X)$, and thus by \eqref{eq:Canonical_Clifford_multiplication_Hodge_star_one_forms},
\[
  \rho_\can(a)\tau = - \sqrt{2}\star(a' \wedge \star\tau)  \in \Omega^{0,1}(X),
\]
which verifies \eqref{eq:Canonical_Clifford_multiplication_positive_spinors}. To deduce \eqref{eq:Canonical_Clifford_multiplication_negative_spinors} from \eqref{eq:Canonical_Clifford_multiplication}, observe that for $\nu \in \Omega^{0,1}(X)$ we have $\star\nu \in \Omega^{1,2}(X)$ by Huybrechts \cite[Lemma 1.2.24, p. 33]{Huybrechts_2005} and $a'\wedge\star\nu \in \Omega^{2,2}(X)$ and $\star(a'\wedge\star\nu) \in \Omega^{0,0}(X)$, so that
\[
  \rho_\can(a)\nu = \sqrt{2}\left(a''\wedge\nu - \star(a'\wedge\star\nu)\right)
  \in \Omega^{0,0}(X)\oplus\Omega^{0,2}(X),
\]
as claimed.

We refer to Salamon \cite[Lemma 4.52, p. 141]{SalamonSWBook} for a verification that \eqref{eq:Canonical_Clifford_multiplication} satisfies the axioms of a \spinc structure \cite[Definition 4.32, p. 125]{SalamonSWBook}, namely that
\begin{equation}
  \label{eq:Salamon_4-18}
  \rho(a)^\dagger = -\rho(a) \in \End(W) \quad\text{and}\quad \rho(a)^\dagger\rho(a) = g(a,a)\id_W,
\end{equation}
for all $a \in \Omega^1(X,\RR)$. If $E$ is a smooth, complex vector bundle over $(X,g,J)$, we obtain a \spinu structure  over $(X,g,J)$ corresponding to a Clifford multiplication bundle map
\[
  \rho:T^*X \to \Hom(W^+, W^-)
\]
via the bundle map
\[
  \rho:T^*X \to \Hom(W^+\otimes E, W^-\otimes E)
\]
defined on elementary tensors by
\begin{multline*}
  \rho(a)(\phi\otimes s) := (\rho(a)\phi)\otimes s \in \Omega^0(W^-\otimes E),
  \\
  \text{for all } a \in \Omega^1(X,\RR), \phi \in \Omega^0(W^+), \text{ and } s \in \Omega^0(E).
\end{multline*}
See Feehan and Leness \cite[Definition 2.2, p. 64, and Lemma 2.3, p. 64]{FL2a} for a more invariant definition of \spinu structure.

\subsection{Dirac operators over complex K\"ahler, almost K\"ahler, complex Hermitian, and almost Hermitian manifolds}
\label{subsec:Dirac_operator_almost_Hermitian}
We now describe the Dirac operator on the canonical \spinc structure $(\rho_\can,W_\can)$ of Definition \ref{defn:Canonical_spinc_bundles} on a closed, smooth almost Hermitian four-manifold $(X,g,J)$. From Donaldson \cite[Equation (15) and p. 56]{DonSW} and Gauduchon \cite[Equations (3.3.4), (3.7.1), and (3.7.2)]{Gauduchon_1997}, the Dirac operator $D:C^\infty(W_\can)\to C^\infty(W_\can)$ defined by the canonical \spinc structure and \emph{Chern connection} \cite[p. 273 and Equation (3.6.1)]{Gauduchon_1997} (see Remark \ref{rmk:Chern_connection}) on the Hermitian line bundle $L = K_X^*$ (the anti-canonical line bundle in \eqref{eq:DefineAntiCanonicalLineBundle} and dual of the canonical line bundle $K_X$ in \eqref{eq:DefineCanonicalLineBundle}), is related to the normalized Dolbeault operator $\sqrt{2}(\partial + \partial^*)$ acting on the spinor bundle $W_\can = W_\can^+\oplus W_\can^-$, where as in
\eqref{eq:Canonical_spinc_bundles},
\[ 
  W_\can^+ = \Lambda^{0,0}(X)\oplus \Lambda^{0,2}(X) \quad\text{and}\quad W_\can^- = \Lambda^{0,1}(X),
\]
by the identity
\begin{equation}
\label{eq:Gauduchon_3-7-2}
D\Phi = \sqrt{2}(\bar\partial + \bar\partial^*)\Phi + \frac{1}{4}\lambda\cdot(\Phi^+ - \Phi^-),
\quad\text{for all } \Phi \in C^\infty(W_\can),
\end{equation}
where $\lambda \in \Omega^1(X,\RR)$ is the \emph{Lee form} defined by $(g,J)$, so \cite[Equation (1.1.4)]{Gauduchon_1997}
\begin{equation}
\label{eq:Gauduchon_1-1-4}
  \lambda = \Lambda(d\omega),
\end{equation}
and $\omega(\cdot,\cdot) = g(J\cdot,\cdot)$ from Gauduchon \cite[p. 259]{Gauduchon_1997}, and ``$\cdot$'' denotes Clifford multiplication \cite[Equation (3.1.4) and p. 276]{Gauduchon_1997}. Note that $d\omega = \lambda\wedge\omega$ by \cite[p. 259]{Gauduchon_1997} when $X$ has real dimension four. When $(X,g,J)$ is almost K\"ahler, so $d\omega=0$, then the identity \eqref{eq:Gauduchon_3-7-2} reduces to the identity $D\Phi = \sqrt{2}(\bar\partial + \bar\partial^*)\Phi$ stated by Donaldson in \cite[Equation (15)]{DonSW}.

When $A$ is a unitary connection on a Hermitian vector bundle $E$ over $X$ and the Levi-Civita connection $\nabla$ on $TX$ is replaced by the connection $\nabla_A$ on $TX\otimes E$, then Gauduchon's proof of \eqref{eq:Gauduchon_3-7-2} now \mutatis yields
\begin{equation}
\label{eq:Gauduchon_3-7-2_auxiliary_Hermitian_bundle_E}
  D_A\Phi = \sqrt{2}(\bar\partial_A + \bar\partial_A^*)\Phi + \frac{1}{4}\lambda\cdot(\Phi^+ - \Phi^-),
\end{equation}
for all $\Phi \in C^\infty(W_\can\otimes E)$.

\begin{rmk}[Chern connection]
\label{rmk:Chern_connection}  
If $(E,h)$ is a \emph{Hermitian} vector bundle with a holomorphic structure $\bar\partial_E$ over a complex manifold $X$, we may choose $A$ to be the \emph{Chern connection} on $E$ --- the \emph{unique} unitary connection on $E$ defined by the holomorphic structure $\bar\partial_E$ and Hermitian metric $h$ (see Kobayashi \cite[Proposition 1.4.9, p. 11]{Kobayashi_differential_geometry_complex_vector_bundles}).
\end{rmk}

\section{Non-Abelian monopole equations over almost Hermitian four-manifolds}
\label{sec:SO3_monopole_equations_over_almost_Hermitian_four-manifolds}
In this section, we develop the structure of an unperturbed version of the non-Abelian monopole equations \eqref{eq:PerturbedSO3MonopoleEquations} when specialized from Riemannian four-manifolds to almost K\"ahler (that is, symplectic) manifolds and, with a minor modification, to almost Hermitian manifolds more generally, following Dowker \cite[Chapter 1]{DowkerThesis} and L\"ubke and Teleman \cite[Section 6.3]{Lubke_Teleman_2006}. Let $(E,h)$ be a Hermitian vector bundle over a closed, oriented, smooth Riemannian manifold $(X,g)$. Note that we are allowing the rank of $E$ to be greater than two.  By analogy with \eqref{eq:PerturbedSO3MonopoleEquations}, we say that a pair $(A,\Phi) \in \sA(E,h)\times\Omega^0(W^+\otimes E)$ is a solution to the \emph{unperturbed non-Abelian monopole equations} if
\begin{subequations}
\label{eq:SO(3)_monopole_equations}
\begin{align}
  \label{eq:SO(3)_monopole_equations_curvature}
  (F_A^+)_0 - \rho^{-1}(\Phi\otimes\Phi^*)_{00} &= 0,
  \\
  \label{eq:SO(3)_monopole_equations_Dirac}
  D_A\Phi &=0,
\end{align}  
\end{subequations}
where, abbreviating $\Lambda^\pm(X) = \Lambda^\pm(T^*X)$ and $\Lambda^2(X) = \Lambda^2(T^*X)$ and noting that $\Lambda^2(X) = \Lambda^+(X) \oplus \Lambda^-(X)$, we recall from \eqref{eq:EndSplitting} that
\begin{equation}
  \label{eq:rho_bundle_isomorphism_Lambda+_with_suW+}
  \rho:\Lambda^+(X) \cong \su(W^+)
\end{equation}
is an isomorphism of Riemannian vector bundles. For any complex vector space $\cV$, we let
\begin{equation}
  \label{eq:Complex_dual_space}
  \cV^* = \Hom(\cV,\CC)
\end{equation}
denote its complex dual space. In equation \eqref{eq:SO(3)_monopole_equations_curvature}, we write
\begin{equation}
  \label{eq:Phi_star}
  \Phi^* = \langle\cdot,\Phi\rangle_{V^+} \in \Omega^0(V^{+,*})
\end{equation}
for $V := W\otimes E$ and $V^\pm := W^\pm\otimes E$, so $V=V^+\otimes V^-$, and Hermitian metric $\langle\cdot,\cdot\rangle_{V^+}$ on $V^+$ given by the tensor product of the Hermitian metrics $\langle\cdot,\cdot\rangle_{W^+}$ on $W^+$ and $\langle\cdot,\cdot\rangle_E$ on $E$, and
\[
  \Phi\otimes\Phi^* \in \Omega^0(i\fu(V^+)) \subset \Omega(\gl(V^+)),
\]
where $\gl(V^+) = V^{+,*}\otimes V^+$. Recall from \cite[Equation (2.13), p. 66]{FL2a} that we have a fiberwise orthogonal decomposition,
\[
  i\su(V^+) = i\su(W^+)\oplus \su(W^+)\otimes\su(E)\oplus i\su(E),
\]
and $i\fu(V^+) = i\underline{\RR}\oplus i\su(V^+)$, where $\underline{\RR} := \RR\,\id_{V^+} \subset \gl(V^+)$. We define the section
\[
  (\Phi\otimes\Phi^*)_{00} \in \Omega^0(\su(W^+)\otimes_\RR\su(E))
\]
to be the image of $\Phi\otimes\Phi^*$ under the orthogonal projection from the fibers of $i\fu(V^+)$ onto the fibers of the subbundle $\su(W^+)\otimes\su(E) \subset i\fu(V^+)$. Because $A$ is a unitary connection on $E$, then $F_A \in \Omega^2(\fu(E))$ and $(F_A^+)_0 \in \Omega^0(\Lambda^+(X)\otimes\su(E))$ denotes the image of  $F_A^+ \in \Omega^0(\Lambda^+(X)\otimes\fu(E))$ induced by the orthogonal projection from the fibers of $\fu(E)$ onto the fibers of the subbundle $\su(E) \subset \fu(E)$. (See Section \ref{sec:PU2Monopoles} for a discussion of the perturbed non-Abelian monopole equations \eqref{eq:PerturbedSO3MonopoleEquations}.)

\begin{rmk}[Typographical errors in literature on the non-Abelian monopole equations]
\label{rmk:SO(3)_monopole_equations_typographical_errors}  
There are typographical errors in some expositions in the literature for the system of equations \eqref{eq:SO(3)_monopole_equations}. In Bradlow and Garc{\'\i}a--Prada \cite[Equations (4.1), (4.3), and (4.4), pp. 575--576]{BradlowGP}, the term $(\Phi\otimes\Phi^*)_{00}$ is incorrectly multiplied by $i=\sqrt{-1}$ and this error is repeated in Dowker \cite[Section 1.1]{DowkerThesis}. 
\end{rmk}

\begin{rmk}[Perturbations of the non-Abelian monopole equations via generic geometric parameters]
\label{rmk:Perturbations_SO3_monopole_equations_generic_geometric_parameters}
While \eqref{eq:SO(3)_monopole_equations} are the unperturbed non-Abelian monopole equations, we recall from Theorem \ref{thm:Transv} that in order to obtain a regular moduli subspace $\sM_\ft^{*,0} \subset \sM_\ft$, one needs to use the perturbed non-Abelian monopole equations \eqref{eq:PerturbedSO3MonopoleEquations} and allow a choice of generic Riemannian metric $g$, perturbation $\tau\in\GL(\Lambda^+(X))$ (near the identity) of the term $\rho^{-1}(\Phi\otimes\Phi^*)_{00}$ in \eqref{eq:SO(3)_monopole_equations_curvature}, and perturbation $\vartheta \in \Omega^1(X,\CC)$ (near zero) of the Dirac operator $D_A$ in \eqref{eq:SO(3)_monopole_equations_Dirac}. In this monograph we primarily restrict our attention to the case where $X$ is a complex, K\"ahler surface and drop the assumption that the metric $g$ be generic, choose $\tau$ to be the identity endomorphism of $\Lambda^+(X)$, and set $\vartheta=0$.
\end{rmk}

Suppose now that $(X,g,J)$ is a closed, smooth, almost Hermitian four-manifold. Our forthcoming Lemma \ref{lem:SO3_monopole_equations_almost_Kaehler_manifold} and the description of the Dirac operator in Section \ref{subsec:Dirac_operator_almost_Hermitian} establish the special form of the non-Abelian monopole equations over $(X,g,J)$, by analogy with the special form of the Seiberg--Witten monopole equations over $(X,g,J)$ (see Donaldson \cite[p. 59]{DonSW}, Kotschick \cite[Section 2.1]{KotschickSW}, Morgan \cite[Equations (7.1) and (7.2), p. 112]{MorganSWNotes}, Nicolaescu \cite[Equation (3.2.9)]{NicolaescuSWNotes}, Salamon \cite[Proposition 12.2, p. 367]{SalamonSWBook}, or Witten \cite[Equation (4.1)]{Witten}). In preparation for our proof of Lemma \ref{lem:SO3_monopole_equations_almost_Kaehler_manifold}, we recall that the spin bundles associated to the \emph{canonical} \spinc structure $\fs_\can = (\rho_\can, W_\can)$ on $X$ are (see equations \eqref{eq:Canonical_spinc_bundles} for $W_\can$ and \eqref{eq:Canonical_Clifford_multiplication} for $\rho_\can$)
\[
  W_\can^+ := \Lambda^{0,0}(X) \oplus \Lambda^{0,2}(X), \quad W_\can^- := \Lambda^{0,1}(X),
  \quad\text{and}\quad W_\can = W_\can^+\oplus W_\can^-.
\]
The direct sum decomposition of $W_\can^+$ is orthogonal with respect to the Hermitian metric on $W_\can$, via the standard Hermitian metric on $\Lambda^\bullet(X,\CC)$ (see Section \ref{subsec:Hodge_star_Lefshetz_dual_Lefshetz_operators}). Consequently,
\begin{multline}
  \label{eq:Canonical_spinu_bundles}
  W_\can^+\otimes E = E \oplus \Lambda^{0,2}(E), \quad W_\can^-\otimes E = \Lambda^{0,1}(E),
  \\
  \text{and}\quad W_\can\otimes E = (W_\can^+\otimes E) \oplus (W_\can^-\otimes E),
\end{multline}
where we abbreviate $\Lambda^{p,q}(E) := \Lambda^{p,q}(X)\otimes E$ for all non-negative integers $p,q$, and note that the preceding direct sum decompositions are orthogonal with respect to the induced (tensor product) Hermitian metric on $W_\can\otimes E$ (see Section \ref{subsec:Hodge_star_Lefshetz_dual_Lefshetz_operators}). (Recall that we write $\Lambda^{0,0}(X) =\Lambda^0(X,\CC) = X\times\CC$ and thus $\Lambda^{0,0}(E) = \Lambda^0(E) = E$.) 

\begin{rmk}[On the restriction to the canonical \spinc structure]
\label{rmk:Restriction_to_canonical_spinc_structure}
There is no loss of generality in restricting to the canonical \spinc structure $(\rho_\can,W_\can)$ since we shall be considering \spinu structures $(\rho_\can,W_\can,E)$ and thus varying $(\rho_\can,W_\can)$ by tensoring it with Hermitian line bundles $L$.
\end{rmk}  

Thus, using \eqref{eq:Canonical_spinu_bundles} we can write $\Phi \in \Omega^0(W^+\otimes E) = \Omega^0(E)\oplus \Omega^{0,2}(E)$ as $\Phi=(\varphi, \psi)$, where $\varphi\in\Omega^0(E)$ and $\psi\in\Omega^{0,2}(E)$. When $(X,g,J)$ is almost K\"ahler, then the Dirac operator can be expressed in the form (see Donaldson \cite[p. 59]{DonSW} and Kotschick \cite[Section 2.1]{KotschickSW})
\begin{equation}
  \label{eq:Coupled_Dirac_operator_almost_Kaehler}
  D_A\Phi = \bar{\partial}_A\varphi + \bar{\partial}_A^*\psi \in \Omega^{0,1}(E)\oplus \Omega^{1,0}(E) = \Omega^{1,1}(E) = \Omega^0(W^-\otimes E).
\end{equation}
See Section \ref{subsec:Dirac_operator_almost_Hermitian}, based on \cite{Gauduchon_1997}, for the modification of this expression for the Dirac operator when $(X,g,J)$ is only almost Hermitian.

By substituting $\Phi = (\varphi,\psi) = \varphi+\psi$ into the expression \eqref{eq:Phi_star} for $\Phi^*$, we obtain
\[
  (\varphi+\psi)^* = \langle \cdot, \varphi+\psi \rangle_{E\oplus \Lambda^{0,2}(E)}
  = \langle \cdot, \varphi \rangle_E + \langle \cdot, \psi \rangle_{\Lambda^{0,2}(E)} = \varphi^* + \psi^* \in E^* \oplus (\Lambda^{0,2}(E))^*.
\]
To help understand the term $\psi^* \in (\Lambda^{0,2}(E))^*$, we observe more generally that
\[
  (\Lambda^{p,q}(E))^* = (\Lambda^{p,q}(X))^*\otimes E^*.
\]
We first note that the Riesz map
\[
  \Lambda^{p,q}(X) \ni \beta \mapsto \langle\cdot,\beta\rangle_{\Lambda^{p,q}(X)} \in (\Lambda^{p,q}(X))^*
\]
is a \emph{complex anti-linear}, isometric isomorphism of Hermitian vector bundles. Recall from Huybrechts \cite[Proposition 1.2.8, p. 27]{Huybrechts_2005} that complex conjugation defines a complex anti-linear isomorphism,
\[
  \Lambda^{q,p}(X) \ni \bar\gamma \mapsto \gamma \in \Lambda^{p,q}(X),
\]
of complex vector bundles. Hence, the composition of the Riesz map and complex conjugation,
\[
  \Lambda^{q,p}(X) \ni \bar\gamma \mapsto \langle\cdot,\gamma\rangle_{\Lambda^{p,q}(X)} \in (\Lambda^{p,q}(X))^*,
\]
is a \emph{complex-linear}, isometric isomorphism of Hermitian vector bundles. In particular, we have a complex-linear, isometric isomorphism of Hermitian vector bundles,
\begin{equation}
  \label{eq:Complex_linear_isomorphism_Lambda_qp_Edual_and_Lambda_pq_E_alldual}
  \Lambda^{q,p}(E^*) = \Lambda^{q,p}(X)\otimes E^* \cong (\Lambda^{p,q}(X))^* \otimes E^*
  = (\Lambda^{p,q}(X)\otimes E)^* = (\Lambda^{p,q}(E))^*,
\end{equation}
induced by the following map on elementary tensors
\begin{equation}
  \label{eq:Complex_linear_isomorphism_Lambda_qp_Edual_and_Lambda_pq_E_alldual_elementary_tensors}
  \bar\gamma\otimes s^* \mapsto \langle\cdot,\gamma\rangle_{\Lambda^{p,q}(X)}\otimes \langle\cdot, s\rangle_E = \langle\cdot,\gamma\otimes s\rangle_{\Lambda^{p,q}(E)},
  \quad\text{for all } \gamma \in \Omega^{p,q}(X) \text{ and } s \in \Omega^0(E),
\end{equation}
where
\[
  s^* := \langle\cdot, s\rangle_E \in \Omega^0(E^*).
\]  
Hence, if $\chi = \sum_{k=1}^r \chi_k\otimes e_k \in \Omega^{p,q}(U,E)$ with respect to a local orthonormal $\RR$-frame $\{e_k\}$ over an open subset $U\subset X$
for the Hermitian vector bundle $E$ of complex rank $r$, where $\chi_k \in \Omega^{p,q}(U,\CC)$, then using \eqref{eq:Complex_linear_isomorphism_Lambda_qp_Edual_and_Lambda_pq_E_alldual_elementary_tensors} the isomorphism \eqref{eq:Complex_linear_isomorphism_Lambda_qp_Edual_and_Lambda_pq_E_alldual} is given by
\begin{equation}
\label{eq:Complex_conjugate_sections_E_and_Lambda02E}  
\sum_{k=1}^r \bar\chi_k\otimes e_k^* \mapsto \sum_{k=1}^r \left\langle\cdot,\chi_k\otimes e_k\right\rangle_{\Lambda^{p,q}(E)}
= \langle\cdot, \chi\rangle_{\Lambda^{p,q}(E)} \in \Omega^{q,p}(U,E^*).
\end{equation}
We may thus interpret
\[
  \chi^* = \langle\cdot,\chi\rangle_{\Lambda^{p,q}(E)} \in \Omega^0((\Lambda^{p,q}(E))^*)
  \cong \Omega^0(\Lambda^{q,p}(E^*)) = \Omega^{q,p}(E^*)
\]
and we write
\begin{equation}
  \label{eq:Hermitian_duals_sections_LambdapqE}
  \chi^* = \langle\cdot,\chi\rangle_{\Lambda^{p,q}(E)} \in \Omega^{q,p}(E^*), \quad\text{for all } \chi \in \Omega^{p,q}(E).
\end{equation}
Now substituting $\Phi = \varphi+\psi$ into the expression $\Phi\otimes\Phi^*$ in equation \eqref{eq:SO(3)_monopole_equations_curvature}, noting that $\varphi\in\Omega^0(E)$ and $\psi\in\Omega^{0,2}(E)$, we obtain
\[
  \Phi\otimes\Phi^* = (\varphi+\psi)\otimes(\varphi+\psi)^* = \varphi\otimes\varphi^* + \varphi\otimes\psi^* + \psi\otimes\varphi^* + \psi\otimes\psi^*,
\]
where, for elementary tensors
\begin{equation}
  \label{eq:varphi_psi_elementary_tensors}
  \varphi = \alpha\otimes s \quad\text{and}\quad \psi = \beta\otimes s,
\end{equation}
with $\alpha\in\Omega^{0,0}(X) = C^\infty(X,\CC)$ and $\beta\in\Omega^{0,2}(X)$ and $s \in \Omega^0(E)$, we have
\begin{equation}
  \label{eq:Hermitian_duals_sections_E_and_Lambda02E}
  \varphi^* = \bar\alpha\otimes s^* \in \Omega^0(E^*) \quad\text{and}\quad \psi^* = \bar\beta\otimes s^* \in \Omega^{2,0}(E^*).
\end{equation}
Using $(\cdot)_0$ to denote the orthogonal projection onto the trace-free part with respect to $\gl(E)$, we have 
\begin{align*}
  \varphi\otimes\varphi^* \in \Omega^0(i\fu(E)) &\quad\text{and}\quad i(\varphi\otimes\varphi^*)_0 \in \Omega^0(\su(E)),
  \\
  \varphi\otimes\psi^* \in \Omega^{2,0}(\gl(E)) &\quad\text{and}\quad (\varphi\otimes\psi^*)_0 \in \Omega^{2,0}(\fsl(E)),
  \\
  \psi\otimes\varphi^* \in \Omega^{0,2}(\gl(E)) &\quad\text{and}\quad (\psi\otimes\varphi^*)_0 \in \Omega^{0,2}(\fsl(E)),
  \\
  \psi\otimes\psi^* \in \Omega^{2,2}(i\fu(E)) &\quad\text{and}\quad i\star(\psi\otimes\psi^*)_0 \in \Omega^0(\su(E)).
\end{align*}
Taking the trace-free part of $\Phi\otimes\Phi^*$ with respect to $\gl(E)$ thus yields
\begin{align*}
  (\Phi\otimes\Phi^*)_0
  &= (\varphi\otimes\varphi^*)_0 \in \Omega^0(i\su(E))
  \\
  &\quad + (\varphi\otimes\psi^*)_0 \in \Omega^{2,0}(\fsl(E))
  \\
  &\quad + (\psi\otimes\varphi^*)_0 \in \Omega^{0,2}(\fsl(E))
  \\
  &\quad + (\psi\otimes\psi^*)_0 \in \Omega^{2,2}(i\su(E)).
\end{align*}
We analyze the terms appearing above by first writing them under the assumption that $\varphi$ and $\psi$ are elementary tensors of the form \eqref{eq:varphi_psi_elementary_tensors}:
\begin{subequations}
\label{eq:TensorEqualitiesOnHermitian}
\begin{align}
  \label{eq:TensorEqualitiesOnHermitian_varphi_varphi*}
  \varphi\otimes\varphi^* &= \alpha\bar\alpha\, s\otimes s^* \in \Omega^0(i\fu(E)),
  \\
  \label{eq:TensorEqualitiesOnHermitian_varphi_psi*}
  \varphi\otimes\psi^* &= \alpha\bar\beta\otimes s\otimes s^* \in \Omega^{2,0}(\gl(E)),
  \\
  \label{eq:TensorEqualitiesOnHermitian_psi_varphi*}
  \psi\otimes\varphi^* &= \beta\bar\alpha\otimes s\otimes s^* \in \Omega^{0,2}(\gl(E)),
  \\
  \label{eq:TensorEqualitiesOnHermitian_psi_psi*}
  \psi\otimes\psi^* &= \beta\wedge\bar\beta\otimes s\otimes s^* \in \Omega^{2,2}(i\fu(E)).
\end{align}
\end{subequations}
The equality \eqref{eq:TensorEqualitiesOnHermitian_varphi_varphi*} simplifies to give
\begin{equation}
\label{eq:TensorEqualityOmega(0,0)Term}
\varphi\otimes\varphi^*\otimes\om=|\alpha|^2\omega \otimes s\otimes s^*.
\end{equation}
Applying the dual Lefshetz operator $\Lambda$ in \eqref{eq:Huybrechts_definition_1-2-21_and_lemmas_1-2-23_and_24} to the identity \eqref{eq:TensorEqualitiesOnHermitian_psi_psi*} yields
\begin{equation}
\label{eq:TensorEqualityOmega(0,2Term1}
\La(\psi\otimes\psi^*)=\La(\beta\wedge\bar\beta)\otimes s \otimes s^*,
\end{equation}
To further analyze the expression $\La(\beta\wedge\bar\beta)$ in \eqref{eq:TensorEqualityOmega(0,2Term1}, we first note that self-duality of $\bar\beta$ and the expression \eqref{eq:Huybrechts_page_33_inner_product_complex-valued_forms} for the inner product on complex differential forms imply that
\[
\beta\wedge\bar\beta=\beta\wedge\star\bar\beta=|\beta|^2 d\vol,
\]
and so
\begin{equation}
\label{eq:(0,2)FormWedgeAndVolume}
\beta\wedge\bar\beta=|\beta|^2 d\vol.
\end{equation}
Next, we compute that
\begin{align*}
\La ( d\vol )
&=
(\star^{-1}\circ L \circ\star)\star 1
\quad\text{(by \eqref{eq:Huybrechts_proposition_1-2-20_i} and \eqref{eq:Huybrechts_definition_1-2-21_and_lemmas_1-2-23_and_24})}
\\
&=
\star^{-1}L\, 1
\quad\text{(by \eqref{eq:Huybrechts_proposition_1-2-20_iii})}
\\
&=\om
\quad\text{(by self-duality of $\om$ and expression \eqref{eq:Huybrechts_remark_1-2-19-ii} for $L$).}
\end{align*}
Substituting the preceding equality into \eqref{eq:(0,2)FormWedgeAndVolume} yields
\[
\La(\beta\wedge\bar\beta)
=
|\beta|^2\om.
\]
Combining the latter identity with \eqref{eq:TensorEqualityOmega(0,2Term1} yields
\begin{equation}
\label{eq:TensorEqualityOmega(0,2Term2}
\La(\psi\otimes\psi^*)
=
|\beta|^2 \om.
\end{equation}
To take the trace-free part of $(\Phi\otimes\Phi^*)_0$ with respect to $\gl(W^+)$ and compute $\rho^{-1}(\Phi\otimes\Phi^*)_{00}$, we shall use Salamon \cite[Lemma 4.62, p. 150 and Remark 4.63, p. 151 and proof of Proposition 12.2, p. 367]{SalamonSWBook} (see also Donaldson \cite[Section 4]{DonSW}, Kotschick \cite[Section 2.1 and derivation of Equations (8) and (9)]{KotschickSW}, Morgan \cite[Section 7.1]{MorganSWNotes}, Nicolaescu \cite[Equations (3.2.6) and (3.2.7) and derivation of Equation (3.2.9)]{NicolaescuSWNotes}, and Witten \cite[Section 4]{Witten}). We begin by recalling the following two lemmas from Salamon \cite{SalamonSWBook}.

\begin{lem}[Clifford multiplication by self-dual two-forms for the canonical \spinc structure]
\label{lem:Salamon_4-61}
(See Salamon \cite[Lemma 4.61, p. 149, together with Lemma 4.60, p. 148, and Remark 3.7, p. 61]{SalamonSWBook}.)
Let $(X,g,J)$ be a four-dimensional, smooth, almost Hermitian manifold with fundamental two-form $\omega = g(\cdot,J\cdot)$ as in \eqref{eq:Fundamental_two-form} and $(\rho_\can,W_\can)$ be the canonical spin${}^c$ structure over $X$ (see equations \eqref{eq:Canonical_spinc_bundles} for $W_\can$ and \eqref{eq:Canonical_Clifford_multiplication} for $\rho_\can$). If $\phi \in \Omega^0(W_\can^+)$ is given by $\phi = (\alpha,\beta) \in \Omega^0(X,\CC)\oplus \Omega^{0,2}(X)$, then\footnote{In this monograph, we use the convention that a Hermitian inner product $\langle\cdot,\cdot\rangle$ is complex anti-linear in the second variable, whereas Salamon uses the convention that it is complex anti-linear in the first variable.}
\begin{equation}
  \label{eq:Salamon_lemma_4-61}
  \rho_\can(\eta)(\alpha,\beta)
  =
  2\left(\eta_0\alpha + \langle \beta,\eta_2 \rangle, \alpha\eta_2 - \eta_0\beta\right)
  \in \Omega^0(X,\CC) \oplus \Omega^{2,0}(X),
  \quad\text{for all } \eta \in \Omega^+(X,i\RR),
\end{equation}
where $\eta = \eta_{2,0} + i\eta_0\omega + \eta_{0,2}$, with $\eta_2 := \eta_{0,2}\in\Omega^{0,2}(X)$, and $\eta_{2,0} = -\overline{\eta_{0,2}} \in \Omega^{2,0}(X)$, and $\eta_0 \in \Omega^0(X,\RR)$.
\end{lem}

The complex linearity of the Clifford multiplication map gives the following

\begin{cor}[Clifford multiplication by  $(0,2)$ and $(2,0)$ forms for the canonical \spinc structure]
\label{cor:CanonicalSpinStructureCliffordMultBy2Forms}
Let $(X,g,J)$ be a four-dimensional, smooth, almost Hermitian manifold with fundamental two-form $\omega = g(\cdot,J\cdot)$ as in \eqref{eq:Fundamental_two-form} and $(\rho_\can,W_\can)$ be the canonical spin${}^c$ structure over $X$ (see equations \eqref{eq:Canonical_spinc_bundles} for $W_\can$ and \eqref{eq:Canonical_Clifford_multiplication} for $\rho_\can$). If $\phi = (\alpha,\beta) \in \Omega^0(X,\CC)\oplus \Omega^{0,2}(X)$, $\eta_{2,0}\in\Om^{2,0}(X)$, and $\eta_{0,2}\in\Om^{0,2}(X)$,
then
\begin{subequations}
  \label{eq:CanonicalSpinStructureCliffordMultByComplex2Forms}
\begin{align}
  \label{eq:CanonicalSpinStructureCliffordMultBy(0,2)Forms}
  \rho_\can(\eta_{0,2})(\alpha,\beta)
  &=
  \left(0,2\alpha\eta_{0,2}\right),
  \\
  \label{eq:CanonicalSpinStructureCliffordMultBy(2,0)Forms}
  \rho_\can(\eta_{2,0})(\alpha,\beta)
  &=
  \left(-2\langle\beta,\bar\eta_{2,0}\rangle,0\right).
\end{align}
\end{subequations}
Equivalently, recalling that $W_\can^+ = \Lambda^0(X)\oplus\Lambda^{0,2}(X)$,
\begin{subequations}
  \label{eq:CanonicalSpinStructureCliffordMultByComplex2FormsCompressedVersion}
\begin{align}
  \label{eq:CanonicalSpinStructureCliffordMultBy(0,2)FormsCompressedVersion}
  \rho_\can(\eta_{0,2})
  &= \begin{pmatrix} 0 & 0 \\ 2\eta_{0,2}\wedge(\cdot) & 0\end{pmatrix} \in \End(W_\can^+),
  \\
  \label{eq:CanonicalSpinStructureCliffordMultBy(2,0)FormsCompressedVersion}
  \rho_\can(\eta_{2,0})
  &= \begin{pmatrix} 0 & -2\langle\cdot,\bar\eta_{2,0}\rangle \\ 0 & 0\end{pmatrix} \in \End(W_\can^+).
\end{align}
\end{subequations}
\end{cor}

\begin{proof}
If $\eta_{0,2}\in\Om^{0,2}(X)$, then
\begin{equation}
\label{eq:Write02FormAsSumOfRealForms}
\eta_{0,2}
=
\frac{1}{2}\left( \eta_{0,2}-\bar\eta_{0,2}\right)
+
\frac{1}{2i}\left( i\eta_{0,2}+i\bar\eta_{0,2}\right).
\end{equation}
Because  $\eta_{0,2}-\bar\eta_{0,2}\in\Om^+(X,i\RR)$, then equation \eqref{eq:Salamon_lemma_4-61} implies that
\begin{equation}
\label{eq:CliffordMultByOneRealComponentOf2Form}
\rho_\can\left( \eta_{0,2}-\bar\eta_{0,2}\right)(\alpha,\beta)
=
2\left( \langle \beta,\eta_{0,2}\rangle,\alpha\eta_{0,2}\right).
\end{equation}
Similarly, because $i\eta_{0,2}+i\bar\eta_{0,2}\in\Om^+(X,i\RR)$, then equation \eqref{eq:Salamon_lemma_4-61}
and the complex anti-linearity of the Hermitian structure in the second argument imply that
\begin{equation}
\label{eq:CliffordMultByOtherRealComponentOf2Form}
\rho_\can\left( i\eta_{0,2}+i\bar\eta_{0,2}\right)(\alpha,\beta)
=
2\left( -i\langle \beta,\eta_{0,2}\rangle,i\alpha\eta_{0,2}\right).
\end{equation}
As noted following \eqref{eq:CliffordMapDefn}, the extension of $\rho_\can$ from real to complex forms is complex linear. Thus,
\begin{align*}
\rho_\can(\eta_{0,2})(\alpha,\beta)
&=
\rho_\can\left( \frac{1}{2}\left( \eta_{0,2}-\bar\eta_{0,2}\right)
+
\frac{1}{2i}\left( i\eta_{0,2}+i\bar\eta_{0,2}\right) \right)(\alpha,\beta)
\quad\text{(by \eqref{eq:Write02FormAsSumOfRealForms})}
\\
&=
\frac{1}{2}\rho_\can\left( \eta_{0,2}-\bar\eta_{0,2}\right)(\alpha,\beta)
+
\frac{1}{2i}\rho_\can\left( i\eta_{0,2}+i\bar\eta_{0,2}\right)(\alpha,\beta)
  \\
  &\qquad\text{(by complex linearity of $\rho_\can$)}
  \\
  &=
\left( \langle \beta,\eta_{0,2}\rangle,\alpha\eta_{0,2}\right)
+\frac{1}{i}\left( -i\langle \beta,\eta_{0,2}\rangle,i\alpha\eta_{0,2}\right)
\quad\text{(by \eqref{eq:CliffordMultByOneRealComponentOf2Form} and \eqref{eq:CliffordMultByOtherRealComponentOf2Form})}
\\
&=
\left( 0, 2\eta_{0,2}\alpha\right).
\end{align*}
This proves \eqref{eq:CanonicalSpinStructureCliffordMultBy(0,2)Forms}.  Equation \eqref{eq:CanonicalSpinStructureCliffordMultBy(0,2)FormsCompressedVersion} follows immediately from 
\eqref{eq:CanonicalSpinStructureCliffordMultBy(0,2)Forms}.
To prove \eqref{eq:CanonicalSpinStructureCliffordMultBy(2,0)Forms}, we observe that for
$\eta_{2,0}\in\Om^{2,0}(X)$,
\begin{equation}
\label{eq:Write20FormAsSumOfRealForms}
\eta_{2,0}
=
\frac{1}{2}\left( -\bar\eta_{2,0}+\eta_{2,0}\right)
+
\frac{1}{2i}\left( i\bar\eta_{2,0}+i\eta_{2,0}\right).
\end{equation}
Because $-\bar\eta_{2,0}+\eta_{2,0}\in\Om^+(X,i\RR)$, then equation \eqref{eq:Salamon_lemma_4-61} implies that
\begin{equation}
\label{eq:CliffordMultByOneRealComponentOf20Form}
\rho_\can\left(-\bar\eta_{2,0}+\eta_{2,0} \right)(\alpha,\beta)
=
2\left( -\langle \beta,\bar\eta_{2,0}\rangle,-\alpha\bar\eta_{2,0}\right).
\end{equation}
Similarly, because $i\bar\eta_{2,0}+i\eta_{2,0} \in\Om^+(X,i\RR)$, then equation \eqref{eq:Salamon_lemma_4-61}
and the complex anti-linearity of the Hermitian structure in the second argument imply that
\begin{equation}
\label{eq:CliffordMultByOtherRealComponentOf20Form}
\rho_\can\left(i\bar\eta_{2,0}+i\eta_{2,0} \right)(\alpha,\beta)
=
2\left( -i\langle \beta,\bar\eta_{2,0}\rangle,i\alpha\bar\eta_{2,0}\right).
\end{equation}
Then,
\begin{align*}
\rho_\can(\eta_{2,0})(\alpha,\beta)
&=
\rho_\can\left( \frac{1}{2}\left( -\bar\eta_{2,0}+\eta_{2,0}\right)
+
\frac{1}{2i}\left( i\bar\eta_{2,0}+i\eta_{2,0}\right)\right)(\alpha,\beta)
\quad\text{by \eqref{eq:Write20FormAsSumOfRealForms})}
\\
&=
\frac{1}{2}\rho_\can\left(  -\bar\eta_{2,0}+\eta_{2,0}\right)(\alpha,\beta)
+
\frac{1}{2i}\rho_\can\left( i\bar\eta_{2,0}+i\eta_{2,0}\right)(\alpha,\beta)
  \\
  &\qquad
\quad\text{(by complex linearity of $\rho_\can$)}
\\
&=
\left( -\langle \beta,\bar\eta_{2,0}\rangle,-\alpha\bar\eta_{2,0}\right)
+\frac{1}{i}\left( -i\langle \beta,\bar\eta_{2,0}\rangle,i\alpha\bar\eta_{2,0}\right)
\quad\text{(by \eqref{eq:CliffordMultByOneRealComponentOf20Form} and \eqref{eq:CliffordMultByOtherRealComponentOf20Form})}
\\
&=
\left(-2\langle \beta,\bar\eta_{2,0}\rangle, 0\right).
\end{align*}
This proves \eqref{eq:CanonicalSpinStructureCliffordMultBy(2,0)Forms}. Equation \eqref{eq:CanonicalSpinStructureCliffordMultBy(2,0)FormsCompressedVersion} follows immediately from
\eqref{eq:CanonicalSpinStructureCliffordMultBy(2,0)Forms}. This completes the proof of the lemma.
\end{proof}

\begin{lem}[Quadratic term in the Seiberg--Witten monopole curvature equation over a K\"ahler surface ]
\label{lem:Salamon_4-62_and_remark_4-63}
(See Salamon \cite[Lemma 4.62, p. 150, and Remark 4.63, p. 151 and proof of Proposition 12.2, p. 367]{SalamonSWBook}.)
Continue the hypotheses of Lemma \ref{lem:Salamon_4-61}. Then $(\phi\otimes\phi^*)_0 \in \Omega^0(i\su(W^+))$ and
\begin{equation}
  \label{eq:Salamon_lemma_4-62}
  \rho_\can^{-1}(\phi\otimes\phi^*)_0
  =
  \frac{i}{4}\left(|\alpha|^2 - |\beta|^2\right)\omega + \frac{1}{2}\left(\beta\bar\alpha - \alpha\bar\beta\right)
  \in C^\infty(X,i\RR)\omega \oplus \Omega^{0,2}(X) \oplus \Omega^{2,0}(X),
\end{equation}
where $|\cdot|$ indicates the pointwise norm on $\Omega^0(X,\CC)$ or $\Omega^{0,2}(X)$ and $(\cdot)_0$ denotes the trace-free part of a section of $\gl(W^+)$.
\end{lem}

When $(X,g,J)$ is an almost K\"ahler four-manifold, the non-Abelian monopole equations \eqref{eq:SO(3)_monopole_equations} take the form described in the forthcoming Lemma \ref{lem:SO3_monopole_equations_almost_Kaehler_manifold} --- compare Bradlow and Garc\'ia--Prada \cite[Equation (5.3), p. 577]{BradlowGP},  Dowker \cite[p. 10]{DowkerThesis}, Labastida and Mari\~no \cite[Equation (3.1)]{Labastida_Marino_1995nam4m}, and Okonek and Teleman \cite[Proposition 2.6, p. 900]{OTVortex}. Bradlow and Garc\'ia--Prada and Dowker use slightly different conventions in their presentation of the equations, Labastida and Mari\~no use quite different notation, and the presentation by Okonek and Teleman is brief. Hence, to avoid ambiguity and make our conventions clear, we include our derivation of the non-Abelian monopole equations over an almost K\"ahler four-manifold, obtaining exact agreement with the equations described by Okonek and Teleman in \cite[Proposition 2.6, p. 900]{OTVortex}.

\begin{lem}[Non-Abelian monopole equations over an almost K\"ahler manifold]
\label{lem:SO3_monopole_equations_almost_Kaehler_manifold}
(See Okonek and Teleman \cite[Proposition 2.6, p. 900]{OTVortex}.)  
Let $(X,g,J)$ be a smooth, almost K\"ahler manifold of real dimension four, let $(\rho_\can,W_\can)$ be the canonical spin${}^c$ structure over $X$ (see equations \eqref{eq:Canonical_spinc_bundles} for $W_\can$ and \eqref{eq:Canonical_Clifford_multiplication} for $\rho_\can$), and let $(E,h)$ be a smooth, Hermitian vector bundle over $X$. If $A$ is a smooth unitary connection on $E$ and $\Phi$ is a smooth section of $W_\can^+\otimes E$, then $(A,\Phi)$ is solution to the unperturbed non-Abelian monopole equations \eqref{eq:SO(3)_monopole_equations} with $\Phi=(\varphi,\psi) \in \Omega^0(E)\oplus\Omega^{0,2}(E)$ if and only if $(A,\varphi,\psi)$ obeys
\begin{subequations}
\label{eq:SO(3)_monopole_equations_Kaehler}    
\begin{align}
  \label{eq:SO(3)_monopole_equations_(1,1)_curvature}
  (\Lambda F_A)_0 &= \frac{i}{2}\left( (\varphi\otimes\varphi^*)_0 - \star(\psi\otimes\psi^*)_0\right),
  \\
  \label{eq:SO(3)_monopole_equations_(0,2)_curvature}
  (F_A^{0,2})_0 &= \frac{1}{2}(\psi\otimes\varphi^*)_0,
  \\
  \label{eq:SO(3)_monopole_equations_(2,0)_curvature}
  (F_A^{2,0})_0 &= -\frac{1}{2}(\varphi\otimes\psi^*)_0,
  \\
  \label{eq:SO(3)_monopole_equations_Dirac_almost_Kaehler}
  \bar{\partial}_A\varphi + \bar{\partial}_A^*\psi &= 0.
\end{align}
\end{subequations}
\end{lem}

Before proving Lemma \ref{lem:SO3_monopole_equations_almost_Kaehler_manifold}, we make some remarks. To interpret the left-hand sides of equations \eqref{eq:SO(3)_monopole_equations_(1,1)_curvature}, \eqref{eq:SO(3)_monopole_equations_(0,2)_curvature}, and \eqref{eq:SO(3)_monopole_equations_(2,0)_curvature}, we recall that
\[
  F_A^{0,2} \in \Omega^{0,2}(\fu(E)) \quad\text{and}\quad (F_A^{0,2})_0 \in \Omega^{0,2}(\su(E)),
\]
with $F_A^{2,0} = -(F_A^{0,2})^\dagger$, while
\[
  \Lambda F_A \in \Omega^0(\fu(E)) \quad\text{and}\quad (\Lambda F_A)_0 \in \Omega^0(\su(E)),
\]
as desired.

\begin{rmk}[Equations \eqref{eq:SO(3)_monopole_equations_Kaehler} defined on Hermitian vector bundle of arbitrary complex rank over an almost Hermitian manifold of any real even dimension]
\label{rmk:SO(3)_monopole_equations_Kaehler_any_dimension}
It is useful to note that, despite their origin as the non-Abelian monopole equations \eqref{eq:SO(3)_monopole_equations} defined by a Hermitian vector bundle $E$ with complex rank two and \spinc structure $(\rho,W)$ over a smooth Riemannian manifold $(X,g)$ of real dimension four, the equations \eqref{eq:SO(3)_monopole_equations_Kaehler} are well-defined for a Hermitian vector bundle $E$ of arbitrary complex rank over an almost Hermitian manifold $(X,g,J)$ of any real even dimension. 
\end{rmk}

\begin{rmk}[Pre-holomorphic pairs]
\label{rmk:Pre-holomorphic_pairs}  
We call equations \eqref{eq:SO(3)_monopole_equations_(0,2)_curvature} and \eqref{eq:SO(3)_monopole_equations_Dirac_almost_Kaehler} the    
\emph{pre-holomorphic pair equations} and call their solutions $(A,\Phi)$ \emph{pre-holomorphic pairs}.
\end{rmk}  

\begin{rmk}[Fixed unitary connection on $\det E$]
\label{rmk:Fixed_unitary_connection_detE}
In writing \eqref{eq:SO(3)_monopole_equations_Kaehler}, we recall (see Section \ref{sec:SpinuPairsQuotientSpace}) that the unitary connection on the Hermitian line bundle $\det E$ induced by the unitary connection $A$ on $E$ is fixed once and for all and thus has no dynamical role in this system.
\end{rmk}  

\begin{rmk}[Alternative form of the $(1,1)$ component of the curvature equation \eqref{eq:SO(3)_monopole_equations_(1,1)_curvature}]
\label{rmk:Alternative_form_1-1_component_SO3_monopole_equations}
As we recall in the proof of Lemma \ref{lem:SO3_monopole_equations_almost_Kaehler_manifold}, the volume form on $(X,g,J)$ is given by $\vol = \frac{1}{2}\omega^2 = \frac{1}{2}\omega\wedge\omega = \star 1$ (see the forthcoming \eqref{eq:Pseudovolume_almost_Hermitian_manifold}), where $\omega=\omega(\cdot,J\cdot)$ is the fundamental two-form \eqref{eq:Fundamental_two-form} defined by $(g,J)$, so
\[
  \frac{1}{2}\Lambda^2(\psi\otimes\psi^*)_0 = \frac{1}{2}\langle(\psi\otimes\psi^*)_0,\omega^2\rangle_{\Lambda^{2,2}(X)} = \langle(\psi\otimes\psi^*)_0,\vol\rangle _{\Lambda^{2,2}(X)}= \star(\psi\otimes\psi^*)_0 \in \Omega^0(i\su(E)).
\]
Equation \eqref{eq:SO(3)_monopole_equations_(1,1)_curvature} can be equivalently written as
\begin{equation}
  \label{eq:SO(3)_monopole_equations_(1,1)_curvature_omega_forms}
  (F_A^\omega)_0 = \frac{i}{4}\left( \omega\otimes(\varphi\otimes\varphi^*)_0 - \Lambda(\psi\otimes\psi^*)_0\right),
\end{equation}
where $F_A^\omega$ is the component of $F_A \in \Omega^2(\fu(E))$ that belongs to $\Omega^0(\fu(E)\otimes \omega)$ and $(F_A^\omega)_0$ is the component of $F_A \in \Omega^2(\fu(E))$ that belongs to $\Omega^0(\su(E))\otimes \omega$ or, equivalently, $F_A^\omega$ is the component of $F_A \in \Omega^+(\fu(E))$ that belongs to $\Omega^0(\fu(E))\otimes \omega$ and $(F_A^\omega)_0$ is the component of $F_A \in \Omega^+(\fu(E))$ that belongs to $\Omega^0(\su(E))\otimes \omega$. Dowker \cite[p. 10]{DowkerThesis} writes $\Lambda^2(\psi\otimes\psi^*)_0$ in \eqref{eq:SO(3)_monopole_equations_(1,1)_curvature} rather than $\star(\psi\otimes\psi^*)_0$ as in Okonek and Teleman \cite[Proposition 2.6, p. 900]{OTVortex}, so there is a minor typographical error in his presentation of the non-Abelian monopole equations over a complex, K\"ahler surface.
\end{rmk}

\begin{rmk}[Reality condition for sections of the complex endomorphism bundle and redundancy of equation \eqref{eq:SO(3)_monopole_equations_(2,0)_curvature}]
\label{rmk:Reality_condition_sections_complex_endomorphism_bundle_and_20_curvature_equation}
The condition that the following element is skew-Hermitian, 
\begin{equation}
  \label{eq:Decompose_a_in_Omega1suE_into_10_and_01_components}
  a = \frac{1}{2}(a'+a'') \in \Omega^1(\su(E))=\Omega^{1,0}(\fsl(E))\oplus \Omega^{0,1}(\fsl(E)),
\end{equation}
is equivalent to the condition that (see Itoh \cite[p. 850]{Itoh_1985} or Kobayashi \cite[Equation (7.6.11), p. 251]{Kobayashi_differential_geometry_complex_vector_bundles})
\begin{equation}
\label{eq:Kobayashi_7-6-11}
  a' = -(a'')^\dagger \in \Omega^{1,0}(\fsl(E)),
\end{equation}
where $a^\dagger = \bar{a}^\intercal$ is defined by taking the complex, conjugate transpose of complex, matrix-valued representatives with respect to any local frame for $E$ that is orthonormal with respect to the Hermitian metric on $E$. Using \eqref{eq:Kobayashi_7-6-11}, we see that equation \eqref{eq:SO(3)_monopole_equations_(2,0)_curvature} is the complex conjugate transpose of \eqref{eq:SO(3)_monopole_equations_(0,2)_curvature} and so can be omitted without loss of generality. 
\end{rmk}

\begin{rmk}[Non-Abelian monopole equations over almost Hermitian four-manifolds]
\label{rmk:SO3_monopole_equations_over_almost_Hermitian_four-manifolds}
Note that we only need to assume that $(X,g,J)$ is \emph{almost K\"ahler} for this decomposition of the non-Abelian monopole equations to hold, rather than K\"ahler as assumed by Dowker \cite{DowkerThesis} and Okonek and Teleman \cite{OTVortex}. When $(X,g,J)$ is only \emph{almost Hermitian}, equations \eqref{eq:SO(3)_monopole_equations_(1,1)_curvature} and \eqref{eq:SO(3)_monopole_equations_(0,2)_curvature} are the same as in the K\"ahler case, but as we describe in Section \ref{subsec:Dirac_operator_almost_Hermitian}, the Dirac equation \eqref{eq:SO(3)_monopole_equations_Dirac_almost_Kaehler} is modified by the addition of a zeroth-order term.
\end{rmk}

We now give the

\begin{proof}[Proof of Lemma \ref{lem:SO3_monopole_equations_almost_Kaehler_manifold}]
For an elementary tensor $\Phi = \phi\otimes s$, where $\phi \in \Omega^0(W_\can^+)$ is given by $\phi = (\alpha,\beta) \in \Omega^0(X,\CC)\oplus \Omega^{0,2}(X)$ and $s \in \Omega^0(E)$, then
by \eqref{eq:Complex_linear_isomorphism_Lambda_qp_Edual_and_Lambda_pq_E_alldual_elementary_tensors}
\[
  \Phi^* = \langle\cdot,\Phi\rangle_{W_\can^+\otimes E} = \langle\cdot,\phi\otimes s\rangle_{W_\can^+\otimes E} = \langle\cdot,\phi\rangle_{W_\can^+}\langle\cdot,s\rangle_E = \phi^*\otimes s^*,
\]
and thus
\[
  \Phi\otimes\Phi^* = \phi\otimes\phi^*\otimes s\otimes s^*.
\]
Therefore, by \eqref{eq:Salamon_lemma_4-62} we obtain
\begin{align*}
  \rho_\can^{-1}(\Phi\otimes\Phi^*)_0
  &= \rho_\can^{-1}(\phi\otimes\phi^*)_0\otimes s\otimes s^*
  \\
  &= \frac{i}{4}\left(|\alpha|^2 - |\beta|^2\right)\omega\otimes s\otimes s^* + \frac{1}{2}\left(\beta\bar\alpha - \alpha\bar\beta\right)\otimes s\otimes s^*,
\end{align*}
where $(\cdot)_0$ denotes the trace-free part of a section of $\gl(W_\can^+)$. We denote $\varphi = \alpha\otimes s$ and $\psi = \beta\otimes s$ as in \eqref{eq:varphi_psi_elementary_tensors} and observe that $\varphi^* = \bar\alpha\otimes s^*$ and $\psi^* = \bar\beta\otimes s^*$ by \eqref{eq:Hermitian_duals_sections_E_and_Lambda02E}, to obtain
\begin{align*}
  \rho_\can^{-1}(\Phi\otimes\Phi^*)_0
  &= \frac{i}{4}\left(|\alpha|^2 - |\beta|^2\right)\omega\otimes s\otimes s^* + \frac{1}{2}\beta\bar\alpha\otimes s\otimes s^* - \frac{1}{2}\alpha\bar\beta\otimes s\otimes s^*
  \\
  &= \frac{i}{4}\left(|\alpha|^2\omega - |\beta|^2\omega\right)\otimes s\otimes s^* + \frac{1}{2}\psi\otimes\varphi^* - \frac{1}{2}\varphi\otimes \psi^*
  \\
  &= \frac{i}{4}\left(\varphi\otimes\varphi^*\otimes \omega - \Lambda(\psi\otimes\psi^*)\right)
    + \frac{1}{2}\psi\otimes\varphi^*
    - \frac{1}{2}\varphi\otimes \psi^*,
\end{align*}
where the last equality follows from \eqref{eq:TensorEqualityOmega(0,0)Term} and
\eqref{eq:TensorEqualityOmega(0,2Term2}. Therefore, taking the trace-free part of sections of $\gl(E)$ yields, when $\Phi$ is an elementary tensor,
\begin{multline}
\label{eq:QuadraticTermAC}
  \rho_\can^{-1}(\Phi\otimes\Phi^*)_{00}
  = \frac{i}{4}\left( \omega\otimes(\varphi\otimes\varphi^*)_0 - \Lambda(\psi\otimes\psi^*)_0\right)
  -\frac{1}{2}(\varphi\otimes\psi^*)_0 + \frac{1}{2}(\psi\otimes\varphi^*)_0
  \\
  \in \Omega^{1,1}(\su(E)) \oplus \Omega^{2,0}(\fsl(E)) \oplus \Omega^{0,2}(\fsl(E)).
\end{multline}
We now show that the identity \eqref{eq:QuadraticTermAC} continues to hold when $\Phi\in\Omega^0(W_\can^+\otimes E)$ is not assumed to be an elementary tensor, that is, not necessarily of the form $\phi\otimes s$, for $\phi\in\Omega^0(W_\can^+)$ and $s\in\Omega^0(E)$. Since $W^+_\can$ is expressed as a direct sum of complex line bundles, $\Lambda^0(X)\oplus\Lambda^{0,2}(X)$, we may always write $\Phi = \varphi + \psi$ with $\varphi = \alpha\otimes s = s$ and $\psi = \beta\otimes t$, for $\alpha=1$ and $\beta\in\Omega^{0,2}(\CC)$ and $s, t\in\Omega^0(E)$. Now $\varphi\in\Om^0(E) \subset \Omega^0(W_\can^+)$ and $\psi\in\Om^{0,2}(E) \subset \Omega^0(W_\can^+)$, and we observe that
\[
\Phi\otimes\Phi^*
=
\varphi\otimes\varphi^*
+
\left(\varphi\otimes\psi^*+\psi\otimes\varphi^*\right)
+
\psi\otimes\psi^* \in \Omega^0(i\fu(W_\can^+\otimes E)).
\]
By the linearity of the isometric isomorphism $\rho_{\can}^{-1}:\su(W_\can^+) \to \Lambda^0(X)\oplus\Lambda^{0,2}(X)$ and the orthogonal projection $(\cdot)_{00}:i\fu(W_\can^+\otimes E) \to \su(W_\can^+)\otimes\su(E)$ of Riemannian vector bundles, we obtain
\begin{equation}
\label{eq:QuadraticDecomposition}
\rho_{\can}^{-1}(\Phi\otimes\Phi^*)_{00}
=
\rho_{\can}^{-1}(\varphi\otimes\varphi^*)_{00}
+
\rho_{\can}^{-1}\left(\varphi\otimes\psi^*+\psi\otimes\varphi^*\right)_{00}
+
\rho_{\can}^{-1}(\psi\otimes\psi^*)_{00}.
\end{equation}
Since $\varphi$ and $\psi$ are elementary tensors, we can apply \eqref{eq:QuadraticTermAC} with $\Phi=\varphi$ and $\Phi = \psi$ to give
\begin{subequations}
\label{eq:QuadraticDecomp2}
\begin{align}
\label{eq:QuadraticDecomp2_varphi}
\rho_{\can}^{-1}(\varphi\otimes\varphi^*)_{00}
&=
\frac{i}{4}\omega\otimes(\varphi\otimes\varphi^*)_0,
\\
\label{eq:QuadraticDecomp2_psi}
\rho_{\can}^{-1}(\psi\otimes\psi^*)_{00}
&=
-\frac{i}{4}\La(\psi\otimes\psi^*)_0.
\end{align}
\end{subequations}
It remains to compute the terms
\[
\rho_{\can}^{-1}\left(\varphi\otimes\psi^*+\psi\otimes\varphi^*\right)_{00}
\]
in \eqref{eq:QuadraticDecomposition}.
For $\alpha\in\Om^{0,0}(\CC)$ and $\beta\in\Om^{0,2}(\CC)$, the endomorphisms $\alpha\otimes\beta^*$ and $\beta\otimes\alpha^*$ of $W_\can^+$ are trace-free because $\alpha$ and $\beta$ are orthogonal as elements of $W_\can^+$. Hence, the $\fu(W_\can^+)$-trace component of
\[
\varphi\otimes\psi^*+\psi\otimes\varphi^*
=
(\alpha\otimes\beta^*)\otimes (s\otimes t^*)
+
(\beta\otimes\alpha^*)\otimes (t\otimes s^*)
 \in \fu(W_\can^+)\otimes\fu(E)
\]
is zero. This implies that
\[
(\varphi\otimes\psi^*+\psi\otimes\varphi^*)_{00}
=
(\varphi\otimes\psi^*+\psi\otimes\varphi^*)_0,
\]
where $(\cdot)_0$ in the right-hand side of the preceding identity denotes the orthogonal projection $\fu(E) \to \su(E)$. Explicitly, we have
\begin{equation}
\label{eq:CanonSpinStructureQuadraticForm}
  (\varphi\otimes\psi^*+\psi\otimes\varphi^*)_{00}
  =
  \beta^*\otimes(s\otimes t^*)_0 + \beta\otimes (t\otimes s^*)_0.
\end{equation}
Note that in the preceding identity we view $\beta$ and $\beta^*$ as linear maps,
\begin{align*}
  \beta &:=\beta\wedge(\cdot) \in \Hom(\Lambda^0(X),\Lambda^{0,2}(X)) \subset \End(W_\can^+),
  \\
  \beta^* &:=\langle \cdot,\beta\rangle_{\Lambda^{0,2}(X)}
            \in \Hom(\Lambda^{0,2}(X),\Lambda^0(X)) \subset \End(W_\can^+).
\end{align*}
With this understanding, equation \eqref{eq:CanonicalSpinStructureCliffordMultByComplex2FormsCompressedVersion} then gives
\begin{equation}
  \label{eq:CanonicalSpinStructureCliffordMultByEndEValuedForms}
\rho_{\can}\left( \beta\otimes (t\otimes s^*)_0 -\bar\beta\otimes (s\otimes t^*)_0\right)
=
2\beta\otimes (t\otimes s^*)_0
+
2 \beta^*\otimes (s\otimes t^*)_0.
\end{equation}
Consequently,
\begin{align*}
\rho_{\can}^{-1}\left(\varphi\otimes\psi^*+\psi\otimes\varphi^*\right)_{00}
&=
\rho_{\can}^{-1}\left(\beta^*\otimes(s\otimes t^*)_0 + \beta\otimes (t\otimes s^*)_0\right)
\quad\text{(by \eqref{eq:CanonSpinStructureQuadraticForm})}
\\
&=
\frac{1}{2}\beta\otimes (t\otimes s^*)_0 -\frac{1}{2}\bar\beta\otimes (s\otimes t^*)_0
\quad\text{(by \eqref{eq:CanonicalSpinStructureCliffordMultByEndEValuedForms})}
\\
&=
\frac{1}{2}(\psi\otimes\varphi^*)_0
-\frac{1}{2}(\varphi\otimes\psi^*)_0,
\end{align*}
where, following \eqref{eq:Hermitian_duals_sections_E_and_Lambda02E}, in the last equality we used $\psi = \beta\otimes t$ and $\varphi = s$ and $\varphi^* = \langle\cdot,\varphi\rangle_E = \langle\cdot,s\rangle_E = s^*$ and
\[
  \psi^* = \langle\cdot,\psi\rangle_{\Lambda^{0,2}(E)}
  = \langle\cdot,\beta\otimes t\rangle_{\Lambda^{0,2}(E)}
  = \bar\beta\otimes\langle\cdot,t\rangle_E
  = \bar\beta\otimes t^*.
\]
Thus, we obtain
\begin{equation}
\label{eq:QuadraticDecomp3}
\rho_{\can}^{-1}\left(\varphi\otimes\psi^*+\psi\otimes\varphi^*\right)_{00}
=
\frac{1}{2}(\psi\otimes\varphi^*)_0
-\frac{1}{2}(\varphi\otimes\psi^*)_0
\in
\Omega^{0,2}(\fsl(E)) \oplus \Omega^{2,0}(\fsl(E)),
\end{equation}
By combining the identities \eqref{eq:QuadraticDecomposition}, \eqref{eq:QuadraticDecomp2}, and \eqref{eq:QuadraticDecomp3}, we see that \eqref{eq:QuadraticTermAC} continues to hold when the assumption that $\Phi$ is an elementary tensor is relaxed. Hence, by choosing $(\rho,W) = (\rho_\can,W_\can)$ and substituting the expression \eqref{eq:QuadraticDecomposition} for $\rho_\can^{-1}(\Phi\otimes\Phi^*)_{00}$ in equation \eqref{eq:SO(3)_monopole_equations_curvature} and applying the decomposition of $\Lambda^2(E)$ implied by the forthcoming decomposition \eqref{eq:Donaldson_Kronheimer_lemma_2-1-57} of $\Lambda^2(X,\CC)$, we obtain the equivalent system of equations \eqref{eq:SO(3)_monopole_equations_(1,1)_curvature_omega_forms} in Remark \ref{rmk:Alternative_form_1-1_component_SO3_monopole_equations} together with \eqref{eq:SO(3)_monopole_equations_(0,2)_curvature} and \eqref{eq:SO(3)_monopole_equations_(2,0)_curvature}. Equation \eqref{eq:SO(3)_monopole_equations_Dirac_almost_Kaehler} follows by substituting the expression \eqref{eq:Coupled_Dirac_operator_almost_Kaehler} for $D_A\Phi$ in \eqref{eq:SO(3)_monopole_equations_Dirac}.

It remains to show that equation \eqref{eq:SO(3)_monopole_equations_(1,1)_curvature} is equivalent to equation \eqref{eq:SO(3)_monopole_equations_(1,1)_curvature_omega_forms}. Because $\star 1 = d\vol$ and $1 = \star^{-1}(d\vol)$, we have
\[
  \Lambda\omega = \star^{-1}(L\star\omega) = \star^{-1}(\omega\wedge\star\omega) = \star^{-1}(\omega\wedge\omega) = 2\star^{-1}(d\vol) = 2\star^{-1}(\star 1) = 2,
\]
and therefore
\begin{equation}
  \label{eq:Lambda_SO(3)_monopole_equations_coupled_spinor_quadratic}
  \Lambda\left(\frac{i}{4}\left( \omega\otimes(\varphi\otimes\varphi^*)_0 - \Lambda(\psi\otimes\psi^*)_0\right)\right)
  = \frac{i}{4}\left(2(\varphi\otimes\varphi^*)_0 - \Lambda^2(\psi\otimes\psi^*)_0\right).
\end{equation}
Finally, note that if $\gamma \in \Omega^{2,2}(X)$, then
\begin{multline*}
  \Lambda^2\gamma = \star^{-1}(L\star(\Lambda\gamma)) = \star^{-1}(L\star(\star^{-1}(L\star\gamma))) = \star^{-1}(\omega^2(\star\gamma))
  \\
  = 2(\star^{-1}(d\,\vol))\star\gamma = 2(\star^{-1}\star 1)\star\gamma = 2\star\gamma.
\end{multline*}
Consequently, we obtain
\begin{equation}
  \label{eq:Lambda_SO(3)_monopole_equations_coupled_spinor_quadratic_simplified}
  \frac{i}{4}\left(2(\varphi\otimes\varphi^*)_0 - \Lambda^2(\psi\otimes\psi^*)_0\right)
  = \frac{i}{2}\left((\varphi\otimes\varphi^*)_0 - \star(\psi\otimes\psi^*)_0\right).
\end{equation}
Thus, by applying $\Lambda$ to equation \eqref{eq:SO(3)_monopole_equations_(1,1)_curvature_omega_forms} and substituting the identities \eqref{eq:Lambda_SO(3)_monopole_equations_coupled_spinor_quadratic} and \eqref{eq:Lambda_SO(3)_monopole_equations_coupled_spinor_quadratic_simplified}, we obtain equation \eqref{eq:SO(3)_monopole_equations_(1,1)_curvature}. This completes the proof of Lemma \ref{lem:SO3_monopole_equations_almost_Kaehler_manifold} and the verification of equation \eqref{eq:SO(3)_monopole_equations_(1,1)_curvature_omega_forms} in Remark \ref{rmk:Alternative_form_1-1_component_SO3_monopole_equations}.
\end{proof}

\section{Witten's dichotomy for solutions to the non-Abelian monopole equations}
\label{sec:Witten_dichotomy_non-Abelian_monopoles}
In Section \ref{subsec:Witten_dichotomy_pre-holomorphic_pairs}, we prove an analogue of \emph{Witten's dichotomy} \label{page:Wittens_Dichotomy} for solutions to the Seiberg--Witten monopole equations (see Donaldson \cite[Section 4]{DonSW}, Kotschick \cite[Section 2.1 and derivation of Equations (8) and (9)]{KotschickSW}, Morgan \cite[Lemma 7.2.1, p. 113]{MorganSWNotes}, Nicolaescu \cite[Proposition 3.2.4]{NicolaescuSWNotes}, Salamon \cite[Proposition 12.3, p. 368 and Corollary 12.4, p. 369]{SalamonSWBook}, and Witten \cite[Section 4]{Witten}). In Section \ref{subsec:Distinguishing_type_1_and_2_solutions_by_degree}, we comment on the extent to which type $1$ and $2$ solutions can be distinguished by properties of the Hermitian vector bundle $E$.

\subsection{Witten's dichotomy for pre-holomorphic pairs on Hermitian vector bundles over complex Hermitian manifolds}
\label{subsec:Witten_dichotomy_pre-holomorphic_pairs}
The version of Witten's dichotomy that we prove as Lemma \ref{lem:Okonek_Teleman_1995_3-1} generalizes analogous results in \cite{BradlowGP, DowkerThesis, Labastida_Marino_1995nam4m,OTVortex, OTQuaternion} by allowing $(E,h)$ to be a Hermitian vector bundle of arbitrary complex rank over an almost Hermitian manifold $(X,g,J)$ of any real even dimension and by omitting the assumption that $(A,\Phi)$ solves \eqref{eq:SO(3)_monopole_equations_(1,1)_curvature}.

\begin{lem}[Witten's dichotomy for pre-holomorphic pairs on Hermitian vector bundles over complex Hermitian manifolds]
\label{lem:Okonek_Teleman_1995_3-1}
(Compare Bradlow and Garc\'ia--Prada \cite[Proposition 5.1, p. 577 and Proposition 5.2, p. 579]{BradlowGP}, Dowker \cite[p. 10]{DowkerThesis}, Labastida and Mari\~no \cite{Labastida_Marino_1995nam4m}, and Okonek and Teleman \cite[Lemma 3.1]{OTVortex}, \cite[Proposition 4.1]{OTQuaternion}.)
Let $(X,g,J)$ be a closed, connected, complex Hermitian manifold, $(\rho_\can,W_\can)$ be the canonical spin${}^c$ structure over $X$ (see equations \eqref{eq:Canonical_spinc_bundles} for $W_\can$ and \eqref{eq:Canonical_Clifford_multiplication} for $\rho_\can$), and $(E,h)$ be a smooth Hermitian vector bundle over $X$. If $(A,\Phi)$ is a smooth solution to the pre-holomorphic pair equations \eqref{eq:SO(3)_monopole_equations_(0,2)_curvature} and \eqref{eq:SO(3)_monopole_equations_Dirac_almost_Kaehler} with $\Phi=(\varphi,\psi) \in \Omega^0(E)\oplus\Omega^{0,2}(E)$ and the induced unitary connection $A_{\det}$ on the Hermitian line bundle $\det E$ also obeys
\begin{equation}
  \label{eq:FdetA02_is_zero}
  F_{A_{\det}}^{0,2} = 0 \in \Omega^{0,2}(\det E),
\end{equation}
then one of $\varphi$ or $\psi$ must be identically zero.
\end{lem}

\begin{proof}
The argument is broadly similar to the proof of the Witten's dichotomy for Seiberg--Witten monopoles (see Donaldson \cite[Section 4]{DonSW}, Morgan \cite[Lemma 7.2.1, p. 113]{MorganSWNotes}, or Nicolaescu \cite[Proposition 3.2.4]{NicolaescuSWNotes} (Witten's own argument \cite[Section 4]{Witten} is more indirect). The argument presented by Dowker in \cite[p. 10]{DowkerThesis} for non-Abelian monopoles is incomplete and the arguments presented by Bradlow and Garc\'ia--Prada, Labastida and Mari\~no, and Okonek and Teleman are indirect and address special cases of the lemma, so we include a direct proof here of a more general result. Applying the operator $\bar\partial_A$ to equation \eqref{eq:SO(3)_monopole_equations_Dirac_almost_Kaehler} yields
\begin{align*}
  0
  &=
    \bar\partial_A^2\varphi+\bar\partial_A\bar\partial_A^*\psi
  \\
  &=
    F^{0,2}_A\varphi + \bar\partial_A\bar\partial_A^*\psi \in \Omega^{0,2}(E),
\end{align*}
where the second equality follows from the forthcoming decomposition \eqref{eq:Decomposition_FA_bitype} of the curvature $F_A$ over a complex manifold and that fact that $J$ is a integrable complex structure, so $d = \partial+\bar\partial$ on $\Omega^1(X,\CC)$ and thus $d_A = \partial_A + \bar\partial_A$ on $\Omega^1(E)$
(in contrast with the more general expression \eqref{eq:d_A_sum_components_almost_complex_manifold} for $d_A$ on an almost Hermitian manifold). By hypothesis, $A$ induces a unitary connection $A_{\det}$ on $\det E$ that obeys \eqref{eq:FdetA02_is_zero} and so we have $(F_A^{0,2})_0=F_A^{0,2}$.  Thus, by combining the preceding equality with
\eqref{eq:SO(3)_monopole_equations_(0,2)_curvature} we obtain
\[
0=
\frac{1}{2}(\psi\otimes\varphi^*)_0\varphi +\bar\partial_A\bar\partial_A^*\psi  \in \Omega^{0,2}(E).
\]
Because $(\psi\otimes\varphi^*)_0=\psi\otimes\varphi^*-\frac{1}{r}\langle\psi,\varphi\rangle_E\,\id_E$,
where $r=\rank_\CC E$, we can rewrite the preceding identity as
\[
0=
\frac{1}{2}\psi |\varphi|_E^2 -\frac{1}{2r}\langle \psi,\varphi\rangle_E\varphi
+\bar\partial_A\bar\partial_A^*\psi \in \Omega^{0,2}(E).
\]
Taking the pointwise $\Lambda^{0,2}(E)$-inner product of the preceding equation with $\psi\in\Omega^{0,2}(E)$ yields
\[
  0
  =
    \frac{1}{2}\langle \psi, \psi \rangle_{\Lambda^{0,2}(E)}|\varphi|_E^2
    - \frac{1}{2r} \left\langle \psi, \langle \psi,\varphi\rangle_E\varphi \right\rangle_{\Lambda^{0,2}(E)}
    + \langle\bar\partial_A\bar\partial_A^*\psi,\psi\rangle_{\Lambda^{0,2}(E)}.
\]
Substituting the following identity into the preceding equation,
\begin{equation}
  \label{eq:Pointwise_norm_quadratic_identity}
  \left\langle \psi, \langle \psi,\varphi\rangle_E\varphi \right\rangle_{\Lambda^{0,2}(E)}
  =
  |\langle \psi,\varphi\rangle_E|_{\Lambda^{0,2}(X)}^2,
  \quad\text{for all } \varphi \in \Omega^0(E) \text{ and } \psi \in \Omega^{0,2}(E),
\end{equation}
gives
\begin{equation}
  \label{eq:Witten_dichotomy_pre-Cauchy-Schwartz}
  \frac{1}{2}|\psi|_{\Lambda^{0,2}(E)}^2|\varphi|_E^2
    - \frac{1}{2r} |\langle \psi,\varphi\rangle_E|_{\Lambda^{0,2}(X)}^2
    + \langle\bar\partial_A\bar\partial_A^*\psi,\psi\rangle_{\Lambda^{0,2}(E)} = 0.
\end{equation}
To understand the origin of the equality \eqref{eq:Pointwise_norm_quadratic_identity}, observe that $\Lambda^{0,2}(X)$ is a complex line bundle, so we may assume without loss of generality that $\psi$ is an elementary tensor, $\psi = \beta\otimes s$ for $\beta \in \Omega^{0,2}(X)$ and $s \in \Omega^0(E)$. Thus
\begin{align*}
  \left\langle \psi, \langle \psi,\varphi\rangle_E\varphi \right\rangle_{\Lambda^{0,2}(E)}
  &=
    \left\langle \beta\otimes s, \langle \beta\otimes s,\varphi\rangle_E\varphi \right\rangle_{\Lambda^{0,2}(E)}
  \\
  &=
    \left\langle \beta\otimes s, \langle s,\varphi\rangle_E \beta\otimes \varphi \right\rangle_{\Lambda^{0,2}(E)}
  \\
  &= \langle\beta,\beta\rangle_{\Lambda^{0,2}(X)}
    \langle s,\varphi\rangle_E \overline{\langle s,\varphi\rangle}_E
    =
    |\beta|_{\Lambda^{0,2}(X)}^2 |\langle s,\varphi\rangle_E|^2,
\end{align*}
while
\begin{align*}
  |\langle \psi,\varphi\rangle_E|_{\Lambda^{0,2}(X)}^2
  &= \langle \langle \psi,\varphi\rangle_E, \langle \psi,\varphi\rangle_E \rangle_{\Lambda^{0,2}(X)}
  \\
  &= \langle \langle \beta\otimes s,\varphi\rangle_E, \langle \beta\otimes s,\varphi\rangle_E \rangle_{\Lambda^{0,2}(X)}
  \\
  &= \langle\beta\langle s,\varphi\rangle_E,
    \beta \langle s,\varphi\rangle_E\rangle_{\Lambda^{0,2}(X)}
  \\
  &= \langle\beta,\beta\rangle_{\Lambda^{0,2}(X)}
    \langle s,\varphi\rangle_E \overline{\langle s,\varphi\rangle}_E
    =
    |\beta|_{\Lambda^{0,2}(X)}^2 |\langle s,\varphi\rangle_E|^2.
\end{align*}
Therefore, by combining the preceding equalities, we obtain \eqref{eq:Pointwise_norm_quadratic_identity}, as claimed. By applying the pointwise Cauchy--Schwartz inequality to $\langle s,\varphi\rangle_E$, we therefore see that
\begin{align*}
  |\langle \psi,\varphi\rangle_E|_{\Lambda^{0,2}(X)}^2
  &=
    |\beta|_{\Lambda^{0,2}(X)}^2 |\langle s,\varphi\rangle_E|^2
  \\
  &\leq
  |\beta|_{\Lambda^{0,2}(X)}^2 |s|_E^2 |\varphi|_E^2
  \\
  &=
  |\beta\otimes s|_{\Lambda^{0,2}(E)}^2 |\varphi|_E^2
  \\
  &=
  |\psi|_{\Lambda^{0,2}(E)}^2 |\varphi|_E^2,
\end{align*}
that is,
\[
  |\langle \psi,\varphi\rangle_E|_{\Lambda^{0,2}(X)}^2
  \leq
  |\psi|_{\Lambda^{0,2}(E)}^2 |\varphi|_E^2. 
\]
Substituting the preceding inequality into \eqref{eq:Witten_dichotomy_pre-Cauchy-Schwartz} gives
\[
  \frac{r-1}{2r}|\psi|_{\Lambda^{0,2}(E)}^2|\varphi|_E^2 + \langle\bar\partial_A\bar\partial_A^*\psi,\psi\rangle_{\Lambda^{0,2}(E)} \leq 0.
\]
Integrating the preceding inequality over $X$ and integrating by parts yields
\[
\frac{r-1}{2r}\int_X|\psi|_{\Lambda^{0,2}(E)}^2|\varphi|_E^2\,d\vol_g
+
\|\bar\partial^*_A\psi\|_{L^2(X)}^2
\leq
0.
\]
Thus, $\bar\partial^*_A\psi=0$ on $X$ and if $\varphi\not\equiv 0$, we must have $\psi\equiv 0$ on an open subset of $X$ and hence $\psi \equiv 0$ on $X$ by the unique continuation property for solutions to second-order elliptic equations defined by differential operators with smooth coefficients and scalar principal symbol (see Aronszajn \cite{Aronszajn}).
\end{proof}

\begin{rmk}[Projective vortices and type $1$ non-Abelian monopoles]
\label{rmk:Projective_vortices_type_1_monopoles} 
Continue the hypotheses of Lemma \ref{lem:Okonek_Teleman_1995_3-1}. If $\psi=0$, then the equations \eqref{eq:SO(3)_monopole_equations_Kaehler} reduce to
\begin{subequations}
\label{eq:SO(3)_monopole_equations_almost_Hermitian_alpha}    
\begin{align}
  \label{eq:SO(3)_monopole_equations_(1,1)_curvature_alpha}
  (\Lambda F_A)_0 &= \frac{i}{2}(\varphi\otimes\varphi^*)_0,
  \\
  \label{eq:SO(3)_monopole_equations_(0,2)_curvature_zero_alpha}
  (F_A^{0,2})_0 &= 0,
  \\
  \label{eq:SO(3)_monopole_equations_Dirac_almost_Hermitian_alpha}
  \bar{\partial}_A\varphi &= 0.
\end{align}
\end{subequations}
Equations \eqref{eq:SO(3)_monopole_equations_(0,2)_curvature_zero_alpha} and \eqref{eq:FdetA02_is_zero} assert that $A$ defines an integrable holomorphic structure $\bar\partial_A$ on $E$ since $\bar\partial_A\circ\bar\partial_A = 0$ by the forthcoming \eqref{eq:dbarE_squared_is_zero} and \eqref{eq:Holomorphic_curvature}. The section $\varphi$ is then holomorphic via \eqref{eq:SO(3)_monopole_equations_Dirac_almost_Hermitian_alpha} and \eqref{eq:SO(3)_monopole_equations_(1,1)_curvature_alpha} is the ($1,1)$-component of the) \emph{projective vortex equation}. The system of equations \eqref{eq:SO(3)_monopole_equations_almost_Hermitian_alpha} is defined on a Hermitian vector bundle of arbitrary complex rank over an almost Hermitian manifold of any real even dimension. 
\label{page:Type1_and_Type2_NonAbelianMonopole}
When $X$ has complex dimension two, the pairs $(\bar\partial_A,\varphi)$ are called \emph{type $1$} solutions to the unperturbed non-Abelian monopole equations \eqref{eq:SO(3)_monopole_equations_Kaehler} and when $X$ has arbitrary complex dimension, they are called \emph{projective vortices}.
\end{rmk}

\begin{rmk}[Holomorphic pairs]
\label{rmk:Holomorphic_pair}  
Continue the hypotheses of Lemma \ref{lem:Okonek_Teleman_1995_3-1}. Recall that because the unitary connection $A_{\det}$ on the Hermitian line bundle $\det E$ induced by $A$ on $E$ obeys $F_{A_{\det}}^{0,2} = 0$ by hypothesis \eqref{eq:FdetA02_is_zero}, the equation \eqref{eq:SO(3)_monopole_equations_(0,2)_curvature_zero_alpha} is strengthened to $F_A^{0,2}=0$
and thus $\bar\partial_A\circ\bar\partial_A = 0$, just as in Remark \ref{rmk:Projective_vortices_type_1_monopoles}. A solution $(A,\varphi)$ to the resulting system, \eqref{eq:SO(3)_monopole_equations_(0,2)_curvature_zero_alpha} with \eqref{eq:FdetA02_is_zero} and \eqref{eq:SO(3)_monopole_equations_Dirac_almost_Hermitian_alpha}, namely
\[
  F_A^{0,2}=0 \quad\text{and}\quad\bar{\partial}_A\varphi = 0,
\]
is called a \emph{holomorphic pair}. This system of equations is defined on a Hermitian vector bundle of arbitrary complex rank over an almost Hermitian manifold of any real even dimension.
\end{rmk}

\begin{rmk}[Projective vortices and type $2$ monopoles]
\label{rmk:Projective_vortices_type_2_monopoles} 
In the case $\varphi = 0$, then the equations \eqref{eq:SO(3)_monopole_equations_Kaehler} reduce to
\begin{subequations}
  \label{eq:SO(3)_monopole_equations_almost_Hermitian_beta}
  \begin{align}
  \label{eq:SO(3)_monopole_equations_(1,1)_curvature_beta}  
  (i\Lambda F_A)_0 &= \frac{i}{2}\star(\psi\otimes\psi^*)_0,
  \\
  \label{eq:SO(3)_monopole_equations_(0,2)_curvature_zero_beta}  
  (F_A^{0,2})_0 &= 0,
    \\
  \label{eq:SO(3)_monopole_equations_Dirac_almost_Hermitian_beta}
  \bar{\partial}_A^*\psi &= 0.
\end{align}
\end{subequations}
In this case, Equation \eqref{eq:SO(3)_monopole_equations_(0,2)_curvature_zero_alpha} asserts that $A$ defines a holomorphic structure $\bar\partial_A$ on $E$. The unitary connection $A$ on $E$ induces a unitary connection $A'$ on $E^*\otimes K_X$ and Equation \eqref{eq:SO(3)_monopole_equations_Dirac_almost_Hermitian_beta} implies that $\bar\psi$ (equivalently, $\psi^*$) is a holomorphic section of the Hermitian vector bundle $F := E^*\otimes K_X$ with respect to\footnote{We shall abuse notation and simply write this as $\bar\partial_A$ if the meaning is unambiguous.} $\bar\partial_{A'}$, and $(A',\bar\psi)$ satisfies the projective vortex equation \eqref{eq:SO(3)_monopole_equations_(1,1)_curvature_beta}; the pairs $(\bar\partial_A,\psi)$ are called \emph{type 2} solutions to the unperturbed non-Abelian monopole equations \eqref{eq:SO(3)_monopole_equations_Kaehler}. See Remark \ref{rmk:Decoupling_Dirac_equation_and_equivalence_type_I_and_II_solutions} for a more detailed explanation of how a type $2$ solution may be regarded as a type $1$ solution for the Hermitian vector bundle $F$ in place of $E$ and vice versa.
\end{rmk}

\begin{rmk}[Decoupling of the Dirac equation and equivalence of type $1$ and $2$ solutions]
\label{rmk:Decoupling_Dirac_equation_and_equivalence_type_I_and_II_solutions}  
For suitable perturbations (see Biquard \cite[Lemma 2.1]{Biquard_1998} and Witten \cite[Equation (4.11)]{Witten}) for the Seiberg--Witten monopole equations for the canonical \spinc structure over complex Hermitian surfaces, the single Dirac-type equation $\bar\partial_A\varphi + \bar\partial_A^*\psi = 0$ decouples to give a pair of equations,
\[
  \bar\partial_A\varphi = 0, \quad \bar\partial_A^*\psi = 0.
\]
In the setting of the non-Abelian monopole equations over almost Hermitian four-manifolds, the second equation is equivalent to $\partial_A\psi = 0$ and thus $\overline{\partial_A\psi} = \bar\partial_A\bar\psi = 0$. But using the isomorphism of Hermitian vector bundles $\bar E \cong E^*$ yields $\bar\psi \in \Omega^{2,0}(\bar E) = \Omega^{2,0}(E^*) = \Omega^0(\Lambda^{2,0}(X)\otimes E^*)$, where $K_X $ is the canonical line bundle \eqref{eq:DefineCanonicalLineBundle}, and so $\bar\psi \in \Omega^0(K_X\otimes E^*)$. If we now define the Hermitian vector bundle $F := K_X\otimes E^*$, then we obtain an isomorphism of Hermitian vector bundles $E^* \cong K_X^*\otimes F$ so that, if $\varphi\in\Omega^0(E)$, then $\bar\varphi\in\Omega^0(E^*) = \Omega^0(K_X^*\otimes F) = \Omega^0(\Lambda^{0,2}(X)\otimes F) = \Omega^{0,2}(F)$, where $K_X^* = \Lambda^{0,2}(X)$ is the anti-canonical line bundle over $X$ defined in \eqref{eq:DefineAntiCanonicalLineBundle}.

In other words, if $(\bar\partial_A,(\varphi,\psi))$ (where $\bar\partial_A$ denotes the partial connection on the complex vector bundle $E$ induced by the unitary connection $A$ on  the Hermitian vector bundle $E$) is a solution to the non-Abelian monopole equations defined by $E$ over an almost Hermitian manifold $(X,g,J)$ if and only if $(\bar\partial_A,(\bar\psi,\bar\varphi))$  (where $\bar\partial_A$ now denotes the partial connection on the complex vector bundle $E^*\otimes K_X$ induced by $A$) is a solution to the non-Abelian monopole equations defined by $F = K_X\otimes E^*$ over $(X,g,J)$.
\end{rmk}

\subsection{Distinguishing type 1 and 2 monopoles}
\label{subsec:Distinguishing_type_1_and_2_solutions_by_degree}
The statement of Witten's dichotomy (for example, see Morgan \cite[Lemma 7.2.1, p. 113]{MorganSWNotes}) for solutions to the (unperturbed) Seiberg--Witten monopole equations \eqref{eq:SeibergWitten} (with $\tau$ equal to the identity endomorphism of $\Lambda^+$ and $\vartheta=0$) corresponding to a \spinc structure $(\rho,W) = (\rho_\can\otimes\id_L, W_\can\otimes L)$ over a closed, complex K\"ahler surface $(X,\omega)$ indicates that the degree of the Hermitian line bundle $L$ (see the forthcoming definition  \eqref{eq:DegreeIntegral}) determines whether the solutions are type $1$ or type $2$; see also Morgan 
\cite[Corollaries 7.2.2, p. 114 and 7.2.3, p. 115]{MorganSWNotes}. In the case of solutions of the (unperturbed) non-Abelian monopole equations \eqref{eq:SO(3)_monopole_equations}, the question of how to distinguish the solutions becomes more subtle.

\subsubsection{Distinguishing type 1 and 2 Seiberg--Witten monopoles}
\label{subsubsec:Distinguishing_type_1_and_2_solutions_SW}
The relationship between the degree of $L$ and the type of the solutions of the Seiberg--Witten equations follows from the analogue of equation \eqref{eq:SO(3)_monopole_equations_(1,1)_curvature} for the Seiberg--Witten equations. 
We consider the unperturbed version of the Seiberg--Witten equations \eqref{eq:SeibergWitten}
(see also Morgan \cite[p. 55]{MorganSWNotes}) on a spin${}^c$ structure $\fs=(\rho,W)$ over a complex  K\"ahler surface $(X,\omega)$.  Let $(\rho_{\can},W_{\can})$ be the canonical spin${}^c$ structure in Definition \ref{defn:Canonical_spinc_bundles}. By Kronheimer and Mrowka \cite[Proposition 1.1.1, p. 3]{KMBook}, there is a Hermitian line bundle $L$ on $X$ such that $(\rho,W)=(\rho_{\can}\otimes\id_L,W_\can\otimes L)$.  Fix the spin connection $A_{\can}$ on $W_{\can}^+=\Lambda^0(X)\oplus \Lambda^{0,2}(X)$ given by the direct sum of the product connection on $\La^0(X)= X\times \CC$ and the connection $A_{K^*}$ induced on $\Lambda^{0,2}(X)= K_X^*$ (see \eqref{eq:DefineCanonicalLineBundle} and \eqref{eq:DefineAntiCanonicalLineBundle}) by the Levi--Civita connection on $TX$, and $A_L$ is a unitary connection on $L$.  The unperturbed versions of the Seiberg--Witten equations \eqref{eq:SeibergWitten} for the pair $(A_L,\Psi)$ are given by
\begin{subequations}
\label{eq:UnperturbedSeibergWitten}
\begin{align}
  \label{eq:UnperturbedSeibergWitten_curvature}
  \tr_{W^+_{\can}} F^+_{A_\can}+2F_{A_L}^+ - \rho^{-1}(\Psi\otimes\Psi^*)_{0}  &=0,
  \\
  \label{eq:UnperturbedSeibergWitten_Dirac}
  D_{A_L}\Psi  &=0.
\end{align}
\end{subequations}
Note that the term $F_{A_\Lambda}^+$ appearing in \eqref{eq:SeibergWitten_curvature} does not appear in
\eqref{eq:UnperturbedSeibergWitten_curvature} because we are discussing the \emph{unperturbed} equations. Write $\Psi=(\alpha,\beta)$, where $\alpha\in\Omega^0(L)$ and $\beta\in\Omega^{0,2}(L)$. The expression \eqref{eq:Coupled_Dirac_operator_almost_Kaehler} for the Dirac operator over a complex K\"ahler manifold shows that \eqref{eq:UnperturbedSeibergWitten_Dirac} is equivalent to
\begin{equation}
\label{eq:SWKahlerDirac1}
  \bar\partial_{A_L}\alpha + \bar\partial_{A_L}^*\beta = 0.
\end{equation}
We now rewrite \eqref{eq:UnperturbedSeibergWitten_curvature}, decomposing it into $(1,1)$, $(0,2)$, and $(2,0)$ components. 
Because $A_{\can}$ is the direct sum of the product connection on $\Lambda^0(X)$, whose curvature vanishes, and the connection $A_{K^*}$ on $\Lambda^{0,2}(X)$, we see that
\begin{equation}
\label{eq:Curvature_of_Canonical_SpinConnection}
\tr_{W_{\can}^+} F_{A_\can}^+ = F_{K^*} = -F_K,
\end{equation}
where $F_{K^*}$ is the curvature of $A_{K^*}$ and $F_K$ is the curvature of the dual connection $A_K$ on the canonical line bundle $K_X$ in \eqref{eq:DefineCanonicalLineBundle}. Thus,
\begin{equation}
\label{eq:TraceOfSpinConn}
\tr_{W^+_{\can}} F^+_{A_\can}+2F_{A_L}^+
=
2 F_{A_L}^+-F_K^+.
\end{equation}
Because the connection $A_{K^*}$ is compatible with the holomorphic structure on $K_X$ induced by the complex structure on $X$, we have $F_K^{0,2}=0$ and $F_K^{2,0}=0$ and so
\begin{equation}
\label{eq:ComponentsOfSpinConn}
2 F_{A_L}^+-F_K^+
=
2F_{A_L}^{1,1}-F_K^{1,1}
+
2F_{A_L}^{2,0} + 2F_{A_L}^{0,2}.
\end{equation}
According to \eqref{eq:Salamon_lemma_4-62}, we have
\begin{equation}
\label{eq:ComponentsOfQuadraticTerm}
\rho^{-1}(\Psi\otimes\Psi^*)_0
=
 \frac{i}{4}\left(|\alpha|_L^2 - |\beta|_L^2\right)\omega + \frac{1}{2}\left(\beta\bar\alpha - \alpha\bar\beta\right).
\end{equation}
By combining \eqref{eq:TraceOfSpinConn}, \eqref{eq:ComponentsOfSpinConn}, and \eqref{eq:ComponentsOfQuadraticTerm}, we see that
 \eqref{eq:UnperturbedSeibergWitten_curvature} is equivalent to
\begin{align*}
2F_{A_L}^{1,1}-F_K^{1,1}&=\frac{i}{4}\left(|\alpha|_L^2 - |\beta|_L^2\right)\omega,
\\
  2F_{A_L}^{0,2} &= \frac{1}{2}\beta\bar\alpha,
\\
 2F_{A_L}^{2,0} &= -\frac{1}{2}\alpha\bar\beta.
\end{align*}
Because $A_L$ is a unitary connection on a Hermitian line bundle, $F_{A_L}^{2,0}=-\bar F_{A_L}^{0,2}$ (see Kobayashi \cite[Equation (1.4.8), p. 11]{Kobayashi_differential_geometry_complex_vector_bundles}). Hence, the third of the preceding equations is redundant and may be omitted. By combining the first two of the preceding equations with \eqref{eq:SWKahlerDirac1}, we see that the unperturbed Seiberg--Witten equations are equivalent to
\begin{subequations}
\label{eq:Seiberg_Witten_equations_complex_Kaehler_surface}
\begin{align}
  \label{eq:SWEquivalent(1,1)Component}
  2F_{A_L}^{1,1}-F_K^{1,1} &= \frac{i}{4}(|\alpha|_L^2 -|\beta|_L^2)\om,
  \\
  \label{eq:SWEquivalent(0,2)Component}
    2F_{A_L}^{0,2} &= \frac{1}{2}\beta\bar\alpha,
  \\
  \label{eq:SWEquivalentDiracComponent}
  \bar\partial_{A_L}\alpha + \bar\partial_{A_L}^*\beta &= 0.
\end{align}
\end{subequations}
This completes our decomposition of the Seiberg--Witten equations on a complex K\"ahler surface.

\begin{rmk}[Other versions of the decomposition of the Seiberg--Witten equations on complex K\"ahler surfaces]
\label{rmk:Decompositions_SW_equations_Kaehler_surfaces}
Decompositions of the Seiberg--Witten equations on K\"ahler surfaces similar to \eqref{eq:Seiberg_Witten_equations_complex_Kaehler_surface} appear in many sources, including Kotschick \cite[Equations (10) and (11), p. 202]{KotschickSW}, Morgan \cite[p. 110 and Equations (7.1) and (7.2), p. 112]{MorganSWNotes}, Nicolaescu \cite[Equations (3.2.8) and (3.2.9), p. 231]{NicolaescuSWNotes}, Salamon \cite[Equation (12.1), p. 367]{SalamonSWBook}, and of course Witten \cite[Equation (4.1), p. 784]{Witten}, who gave the first derivation. The variety of conventions and notation in these sources makes comparison difficult. For example, our equations \eqref{eq:Seiberg_Witten_equations_complex_Kaehler_surface} differ from those in Morgan
\cite[Equations (7.1) and (7.2), p. 112]{MorganSWNotes} only in so far as that the connection $A_L$ appearing in
\eqref{eq:Seiberg_Witten_equations_complex_Kaehler_surface} is a unitary connection on the Hermitian line bundle $L$, whereas Morgan uses a unitary connection $A$ on the Hermitian line bundle $\det(W_{\can}^+\otimes L)\cong K_X^*\otimes L^{\otimes 2}$.  Morgan writes $\sL_0$ for the line bundle that we denote by $L$ --- see \cite[Section 7.1, p. 110]{MorganSWNotes}.
Hence, the curvatures of $A_L$ and $A$ are related by $F_A=2F_{A_L}-F_K$ (compare \eqref{eq:TraceOfSpinConn}).

The connection $B$ in Kotschick \cite[Equations (10) and (11), p. 202]{KotschickSW} is the same as the connection $A_L$ appearing in \eqref{eq:Seiberg_Witten_equations_complex_Kaehler_surface}, while the connection $A_0$
in \cite[Equations (10) and (11), p. 202]{KotschickSW} is a connection on $\La^{0,2}(X)$ so $F_{A_0}=F_{K^*}=-F_K$, as we have calculated. With this understanding, one can derive equations
\eqref{eq:SWEquivalent(1,1)Component} and \eqref{eq:SWEquivalent(0,2)Component} by multiplying both sides of \cite[Equations (10) and (11), p. 202]{KotschickSW} by two.

Similar translations are also possible for the expressions in Nicolaescu \cite[Equations (3.2.8) and (3.2.9), p. 231]{NicolaescuSWNotes} and Salamon \cite[Equation (12.1), p. 367]{SalamonSWBook}, where Salamon uses $F_B$ to denote the curvature of a \emph{virtual connection} in the sense of \cite[Remark 6.5, p. 191]{SalamonSWBook}.
\end{rmk}

Keeping in mind that equations \eqref{eq:Seiberg_Witten_equations_complex_Kaehler_surface} are valid over almost Hermitian manifolds of any dimension and focusing for now just on the pair of equations \eqref{eq:SWEquivalent(0,2)Component} and \eqref{eq:SWEquivalentDiracComponent}, we see that it is possible to generalize the statement of Witten's Dichotomy in Morgan \cite[Lemma 7.2.1, p. 113]{MorganSWNotes} to one along the lines of Lemma \ref{lem:Okonek_Teleman_1995_3-1}:

\begin{lem}[Witten's dichotomy for pre-holomorphic pairs on Hermitian line bundles over complex Hermitian manifolds]
\label{lem:Okonek_Teleman_1995_3-1_SW}
Let $(X,g,J)$ be a closed, connected, complex Hermitian manifold, $(\rho_\can,W_\can)$ be the canonical spin${}^c$ structure over $X$ (see equations \eqref{eq:Canonical_spinc_bundles} for $W_\can$ and \eqref{eq:Canonical_Clifford_multiplication} for $\rho_\can$), and $(L,h)$ be a smooth Hermitian line bundle over $X$. If $(A_L,\Psi)$ is a smooth solution to the pre-holomorphic pair equations \eqref{eq:SWEquivalent(0,2)Component} and \eqref{eq:SWEquivalentDiracComponent} with $\Psi=(\alpha,\beta) \in \Omega^0(L)\oplus\Omega^{0,2}(L)$, then one of $\alpha$ or $\beta$ must be identically zero.
\end{lem}

\begin{proof}
We proceed as in the proof of Lemma \ref{lem:Okonek_Teleman_1995_3-1}, although the argument is simpler. Applying the operator $\bar\partial_{A_L}$ to equation \eqref{eq:SWEquivalentDiracComponent} yields
\begin{align*}
  0 &= \bar\partial_{A_L}^2\alpha+\bar\partial_{A_L}\bar\partial_{A_L}^*\beta
  \\
  &= F^{0,2}_{A_L}\alpha + \bar\partial_{A_L}\bar\partial_{A_L}^*\beta \in \Omega^{0,2}(L),
\end{align*}
where the second equality follows from the forthcoming decomposition \eqref{eq:Decomposition_FA_bitype} of the curvature $F_{A_L}$ over a complex manifold and that fact that $J$ is a integrable complex structure, so $d = \partial+\bar\partial$ on $\Omega^1(X,\CC)$ and thus $d_{A_L} = \partial_{A_L} + \bar\partial_{A_L}$ on $\Omega^1(L)$. By combining the preceding equality with \eqref{eq:SWEquivalent(0,2)Component} we obtain
\[
  0 = \frac{1}{4}(\beta\bar\alpha)\alpha +\bar\partial_{A_L}\bar\partial_{A_L}^*\beta  \in \Omega^{0,2}(L).
\]
Noting that $\bar\alpha\alpha = \langle\alpha,\alpha\rangle_L = |\varphi|_L^2$ and taking the pointwise $\Lambda^{0,2}(L)$-inner product of the preceding equation with $\beta\in\Omega^{0,2}(L)$ yields
\[
  0 = \frac{1}{4}|\beta|_{\Lambda^{0,2}(L)}|\varphi|_L^2 + \langle\bar\partial_{A_L}\bar\partial_{A_L}^*\beta,\beta\rangle_{\Lambda^{0,2}(L)}.
\]
Integrating the preceding inequality over $X$ and integrating by parts yields
\[
\frac{1}{4}\int_X|\beta|_{\Lambda^{0,2}(L)}^2|\alpha|_L^2\,d\vol_g
+
\|\bar\partial^*_{A_L}\beta\|_{L^2(X)}^2
=
0.
\]
Thus, $\bar\partial^*_{A_L}\beta=0$ on $X$ and if $\alpha\not\equiv 0$, then we must have $\beta\equiv 0$ on an open subset of $X$ and hence $\beta\equiv 0$ on $X$ by the unique continuation property for solutions to second-order elliptic equations defined by differential operators with smooth coefficients and scalar principal symbol (see Aronszajn \cite{Aronszajn}). A similar argument shows that if $\beta\not\equiv 0$, then we must have $\alpha\equiv 0$.
\end{proof}

To distinguish when $(\rho_{\can}\otimes\id_L,W_{\can}\otimes L)$ admits Seiberg--Witten monopoles of type $1$ or $2$, we recall the definition (see Kobayashi \cite[Equation (3.1.17), p. 52]{Kobayashi_differential_geometry_complex_vector_bundles}) of the degree of a complex vector bundle $E$ over a closed, complex K\"ahler manifold $(X,\omega)$ with complex dimension $n$ by
\begin{equation}
  \label{eq:Degree}
  \deg E := \int_X c_1(E)\wedge\omega^{n-1} = \langle c_1(E)\smile\omega^{n-1}, [X]\rangle.
\end{equation}
where $c_1(E) \in H^2(X;\RR)$ is the first Chern class of $E$. Recall that the first Chern class $c_1(E) = c_1(\det E) \in H^2(X,\ZZ)$ of a complex vector bundle $E$ over a smooth manifold $X$ is represented by the real two-form (see Donaldson and Kronheimer \cite[Section 2.1.4, p. 39]{DK} or Kobayashi \cite[Equation (2.2.14), p. 37]{Kobayashi_differential_geometry_complex_vector_bundles})
\begin{equation}
  \label{eq:Chern-Weil_formula}
  \frac{i}{2\pi}\tr_E F_A \in \Omega^2(X,\RR),
  \quad\text{with } c_1(E) = \left[\frac{i}{2\pi}\tr_E F_A\right] \in H^2(X,\RR).
\end{equation}
By applying the Chern--Weil formula \eqref{eq:Chern-Weil_formula}, the definition \eqref{eq:Degree} of the degree of $E$ yields
\begin{equation}
  \label{eq:DegreeIntegral}
  \deg E = \frac{i}{2\pi } \int_X (\tr_E F_A)\wedge\omega^{n-1}.
\end{equation}
The following result is well-known in algebraic geometry,
see Kobayashi \cite[Theorem 3.1.24, p. 56]{Kobayashi_differential_geometry_complex_vector_bundles}; we record it here for convenience.

\begin{lem}[Holomorphic line bundles with negative degree over K\"ahler manifolds have no nontrivial holomorphic sections]
\label{lem:NegDegLineBundlesNoSections}
Let $\sL$ be a holomorphic line bundle over a closed, smooth K\"ahler manifold $(X,\om)$. If $\deg\sL<0$, then $H^0(X;\sL)=(0)$, that is, $\sL$ admits no non-vanishing holomorphic sections.
\end{lem}

\begin{proof}
The zero-locus, $s^{-1}(0)$, of  a non-zero, holomorphic section $s$ of $\sL$ is a divisor,
\[
D=\sum_i n_i V_i,
\]
where $V_i\subset X$ is a complex analytic hypersurface and each integer $n_i\in\ZZ$ records the order to which the section $s$ vanishes along $V_i$ (see Griffiths and Harris \cite[Chapter 1, Section 1, p. 11]{GriffithsHarris}). Because $s$ is holomorphic, $D$ is an \emph{effective} divisor (see \cite[Chapter 1, Section 1, p. 130]{GriffithsHarris}) in the sense that each $n_i$ is strictly positive. Let $n$ denote the complex dimension of $X$. From the discussion in 
\cite[Chapter 0, Section 2, p. 33 and Chapter 0, Section 4, pp. 60-61]{GriffithsHarris} or Voisin \cite[Lemma 11.14, p. 270]{Voisin_hodge_theory_complex_algebraic_geometry_I}, we see that because the singular set of each complex analytic hypersurface $V_i$ has complex codimension at least two in $X$, integration over the locus $V_i^{\reg}$ of smooth points of $V_i$ defines a functional
\[
H^{2n-2}(X;\RR) \ni [\Theta]\mapsto \int_{V_i^{\reg}} \Theta \in\RR.
\]
Because $X$ is a closed, oriented manifold of real dimension $2n$, Poincar\'e duality (see Bott and Tu \cite[Equation (5.4) p. 44]{BT}) gives a cohomology class $\eta_D\in H^2(X;\RR)$ representing the above functional in the sense that
\begin{equation}
\label{eq:PoincareDual_of_Divisor}
\int_X \eta_D\wedge \Theta
=
\sum_i n_i \int_{V_i^{\reg}}\Theta.
\end{equation}
By the identification of $\sL$ with the line bundle defined by $D$ on \cite[Chapter 1, Section 1, p. 136]{GriffithsHarris} and the resulting equality $c_1(\sL)=\eta_D$ in \cite[Chapter 1, Section 1, p. 141]{GriffithsHarris} or Voisin \cite[Theorem 11.33, p. 288]{Voisin_hodge_theory_complex_algebraic_geometry_I}, we can rewrite \eqref{eq:PoincareDual_of_Divisor} as
\begin{equation}
\label{eq:c1_is_PD_of_Divisor}
\int_X c_1(\sL)\wedge \omega^{n-1}
=
\sum_i n_i \int_{V_i^{\reg}} \omega^{n-1}.
\end{equation}
By Wirtinger's Theorem (see \cite[Chapter 0, Section 2, p. 31]{GriffithsHarris}),
\begin{equation}
\label{eq:WirtingerThm}
\Vol(V_i)=\frac{1}{(n-1)!}\int_{V_i^{\reg}}\omega^{n-1}.
\end{equation}
Therefore, we have
\begin{align*}
\deg\sL
&=\int_X c_1(\sL)\wedge\omega^{n-1} \quad\text{(by \eqref{eq:Degree})}
\\
&= \sum_i n_i \int_{V_i^{\reg}} \omega^{n-1} \quad\text{(by \eqref{eq:c1_is_PD_of_Divisor})}
\\
&=
\sum_i n_i(n-1)!\Vol(V_i) \quad\text{(by \eqref{eq:WirtingerThm})}.
\end{align*}
The resulting equality,
\[
\deg\sL
=
\sum_i n_i(n-1)!\Vol(V_i),
\]
and the fact that the zero locus of the holomorphic section $s$ is an effective divisor, so that $n_i> 0$ for all $i$, imply that $\deg\sL\ge 0$.  Thus, if $\sL$ admits a non-vanishing holomorphic section, then $\deg\sL\ge 0$. The contrapositive of the preceding statement is the conclusion of the lemma.
\end{proof}

\begin{rmk}[Non-existence of non-trivial holomorphic sections on negative line bundles]
\label{rmk:NegativityOfALineBundle}
We now describe a different method from that of Lemma \ref{lem:NegDegLineBundlesNoSections} for proving that a holomorphic line bundle has no non-trivial holomorphic sections. Recall from Griffiths and Harris \cite[Definition and Proposition, p. 148]{GriffithsHarris}, Kobayashi \cite[Section 3.3]{Kobayashi_differential_geometry_complex_vector_bundles}, and Wells \cite[Chapter VI, Definition 2.1, p. 233]{Wells3} that a complex line bundle $L$ over a closed, complex K\"ahler manifold $M$ of dimension $n$ is \emph{positive} if and only if $c_1(L) \in H^2(M;\RR)$ is represented by a positive $(1,1)$ form, in which case one writes $c_1(L) > 0$, and $L$ is \emph{negative} if $L^*$ is positive, expressed as $c_1(L) < 0$. We write $c_1(L)=0$ if $L$ is neither positive nor negative. A closed, real $(1,1)$ form
\[
  \Omega = \frac{i}{2\pi} \sum_{\alpha,\beta=1}^n\Omega_{\alpha\bar\beta}\,dz^\alpha\wedge\bar z^\beta
\]
is \emph{positive} if at each point $x\in X$, the Hermitian matrix $(\Omega_{\alpha\bar\beta}(x))$ is positive definite \cite[Section 3.3]{Kobayashi_differential_geometry_complex_vector_bundles}.

If $(L,\partial_L)$ is a \emph{holomorphic} line bundle with $c_1(L) < 0$, then by the vanishing theorems of Kodaira and Nakano (see Kobayashi \cite[Theorems 3.3.1 and 3.3.2, pp. 64--65]{Kobayashi_differential_geometry_complex_vector_bundles} or Wells \cite[Chapter II, Example 2.12, p. 48 and Chapter VI, Theorem 2.4 (b), p. 226]{Wells3}) we have 
\[
  H^q(X;\bOmega^p(L)) = 0, \quad\text{for } p+q<n,
\]
where $\bOmega^p(L) = \sO(X,\wedge^p(T^*X)\otimes_\CC L) = \bOmega^p\otimes L$, the sheaf of (global) holomorphic $p$-forms on $X$ with coefficients in $L$, and  $\bOmega^0(L) = \sO(X,L)$. Dolbeault's Theorem (see Wells \cite[Chapter II, Theorem 3.17, p. 61]{Wells3} when $L=X\times\CC$ and \cite[Chapter II, Theorem 3.20, p. 63]{Wells3} when $L$ is replaced by an arbitrary Hermitian vector bundle $E$) asserts that
\[
 H^q(X;\bOmega_X^p(L)) \cong H^{p,q}(X;L),
\]
where $H^{p,q}(X;L)$ is Dolbeault cohomology (see Wells \cite[Chapter II, Theorem 3.20, p. 63, Chapter IV, Example 2.7, p. 117 and Example 5.7, p. 151, and Chapter V, Theorem 2.7, p. 170]{Wells3}. In particular, if $c_1(L)<0$ then
\[
  H^{0,0}(X;L) \cong H^0(X;\sO(L)) = 0
\]
and so $L$ has no non-zero global holomorphic sections.
\end{rmk}

Arguing as in the proofs of \cite[Lemma 7.2.1, p. 113 and Corollaries 7.2.2, p.114 and 7.2.3, p. 115]{MorganSWNotes} and Salamon \cite[Corollary 12.4, p. 369]{SalamonSWBook}, one can distinguish between the two types of solutions that arise in Lemma \ref{lem:Okonek_Teleman_1995_3-1_SW}
by combining equation \eqref{eq:SWEquivalent(1,1)Component}, the definition of the degree of a line bundle \eqref{eq:Degree}, and Lemma \ref{lem:NegDegLineBundlesNoSections}
as follows.

\begin{lem}[Degree constraint on the type of a Seiberg--Witten monopole over a complex K\"ahler surface]
\label{lem:DegreeConstraints_On_Type_of_SW_Solution}
Let $L$ be a Hermitian line bundle over a complex K\"ahler surface $(X,\om)$ and $(A_L,\Psi)$ be a solution to the Seiberg--Witten equations \eqref{eq:Seiberg_Witten_equations_complex_Kaehler_surface} for the spin${}^c$ structure $(\rho_{\can}\otimes\id_{L},W_\can\otimes L)$ where $(\rho_{\can},W_\can)$ is the canonical spin${}^c$ structure as in Definition \ref{defn:Canonical_spinc_bundles}. If $\Psi=(\alpha,\beta)$ with $\alpha\in\Om^0(L)$ and $\beta\in\Om^{0,2}(L)$, then the following hold:
\begin{enumerate}
\item
\label{item:DegreeConstraints_On_Type_of_SW_Solution_Type1}
The pair $(A_L,\Psi)$ is type $1$ with $\alpha\not\equiv 0$ and $\beta\equiv 0$ if and only if $0\le 2\deg L<\deg K_X$.
\item
\label{item:DegreeConstraints_On_Type_of_SW_Solution_Type2}
The pair $(A_L,\Psi)$ is type $2$ with $\alpha\equiv 0$ and $\beta\not\equiv 0$ if and only if $\deg K_X< 2\deg L\le 2\deg K_X$.
\item
\label{item:DegreeConstraints_On_Type_of_SW_Solution_ZeroSection}
The section $\Psi$ vanishes identically if and only if $2\deg L=\deg K_X$.
\end{enumerate}
\end{lem}

\begin{proof}
By Lemma \ref{lem:Okonek_Teleman_1995_3-1_SW}, we know that $\alpha\equiv 0$ or $\beta\equiv 0$ and thus exactly one of the three cases in Items \eqref{item:DegreeConstraints_On_Type_of_SW_Solution_Type1},
\eqref{item:DegreeConstraints_On_Type_of_SW_Solution_Type2}, and \eqref{item:DegreeConstraints_On_Type_of_SW_Solution_ZeroSection} holds.  Moreover, in all three cases, equation \eqref{eq:SWEquivalent(0,2)Component} implies that $(L,\bar\partial_{A_L})$ is a holomorphic line bundle.

By combining the expression \eqref{eq:DegreeIntegral} with $E$ replaced by $L$ for the degree of the line bundle $L$ with equation \eqref{eq:SWEquivalent(1,1)Component}, we obtain
\begin{equation}
\label{eq:DegreeConstraintOnSectionNorms}
  2\deg L-\deg K_X
  =
  \frac{1}{8\pi}\int_X ( |\beta|^2-|\alpha|^2) \vol_\omega,
\end{equation}
recalling that $\vol_\omega = \omega^2/2$ (see the forthcoming \eqref{eq:Pseudovolume_almost_Hermitian_manifold}).

Consider Item \eqref{item:DegreeConstraints_On_Type_of_SW_Solution_Type1}. Assume that $(A_L,\Psi)$ is type $1$ with $\alpha\not\equiv 0$ and $\beta\equiv 0$. The expression \eqref{eq:DegreeConstraintOnSectionNorms} implies that $2\deg L-\deg K_X<0$ and so $2\deg L<\deg K$. Because $(A_L,\Psi)$ is type $1$, equation \eqref{eq:SWEquivalentDiracComponent} implies that $\alpha$ is a holomorphic section of the holomorphic line bundle $(L,\bar\partial_{A_L})$. Because $\alpha\not\equiv 0$ by assumption, Lemma \ref{lem:NegDegLineBundlesNoSections} ensures that $\deg L\geq 0$. Thus, we obtain $0\leq 2\deg L<\deg K$, which verifies the ``only if'' assertion in Item \eqref{item:DegreeConstraints_On_Type_of_SW_Solution_Type1}.

Consider Item \eqref{item:DegreeConstraints_On_Type_of_SW_Solution_Type2}. Assume that
$(A_L,\Psi)$ is type $2$ with $\varphi\equiv 0$ and $\beta\not\equiv 0$.  The expression \eqref{eq:DegreeConstraintOnSectionNorms} implies that $2\deg L-\deg K_X>0$ and so $2\deg L>\deg K_X$.
Because $X$ has complex dimension two, equation \eqref{eq:SWEquivalentDiracComponent}
and the definition of Dolbeault cohomology (see Kobayashi \cite[Equation (3.2.45), p. 63]{Kobayashi_differential_geometry_complex_vector_bundles}) imply that $\beta$ defines a non-trivial element of the Dolbeault cohomology $H^{0,2}(X;L)$ if $L$ has the holomorphic structure given by $\bar\partial_{A_L}$. The isomorphisms
\begin{align*}
  H^{0,2}(X;L)
  &\cong
    H^{2,0}(X;L^*)^*\quad\text{(Serre Duality)}
  \\
  &\cong
    H^0(X;K_X\otimes L^*)^* \quad\text{(Dolbeault isomorphism)}
\end{align*}
imply that $H^0(X;K_X\otimes L^*)$ is non-trivial. (For statements of Serre duality, see Kobayashi \cite[Theorem 3.2.50, p. 64]{Kobayashi_differential_geometry_complex_vector_bundles}, Nicolaescu \cite[Equation (3.1.8), p. 204]{NicolaescuSWNotes},  Voisin \cite[Theorem 5.32, p. 135]{Voisin_hodge_theory_complex_algebraic_geometry_I}, or Wells \cite[Chapter V, Theorem 2.7, p. 170]{Wells3}. For a statement of the Dolbeault isomorphism, see Kobayashi \cite[Equation (3.2.46), p. 63]{Kobayashi_differential_geometry_complex_vector_bundles}.)
By Lemma \ref{lem:NegDegLineBundlesNoSections}, the non-triviality of $H^0(X;K_X\otimes L^*)$  yields the inequality $\deg K_X\otimes L^*\geq 0$ and so $\deg K_X\geq \deg L$. Thus, we obtain $\deg K_X < 2\deg L \leq 2\deg K$, which verifies the ``only if'' assertion in Item \eqref{item:DegreeConstraints_On_Type_of_SW_Solution_Type2}.

Consider Item \eqref{item:DegreeConstraints_On_Type_of_SW_Solution_ZeroSection}. If $\Psi\equiv 0$, then the equality $2\deg L=\deg K_X$ follows immediately from equation \eqref{eq:DegreeConstraintOnSectionNorms} and this verifies the ``only if'' assertion in Item \eqref{item:DegreeConstraints_On_Type_of_SW_Solution_ZeroSection}.

We verify the ``if'' assertions in Items \eqref{item:DegreeConstraints_On_Type_of_SW_Solution_Type1},
\eqref{item:DegreeConstraints_On_Type_of_SW_Solution_Type2}, and \eqref{item:DegreeConstraints_On_Type_of_SW_Solution_ZeroSection} as follows.
By Lemma \ref{lem:Okonek_Teleman_1995_3-1_SW} one and only one of the three cases,
\begin{inparaenum}[\itshape i\upshape)]
\item
$(A_L,\Psi)$ is type $1$ with $\alpha\not\equiv 0$, or
\item
 $(A_L,\Psi)$ is type $2$ with $\beta\not\equiv 0$, or
\item
 $\Psi\equiv 0$
\end{inparaenum} can hold.
Because only one of the three conditions 
\begin{inparaenum}[\itshape i\upshape)]
\item
$0\le 2\deg L<\deg K_X$, or 
\item
$\deg K_X<2\deg L\le 2\deg K_X$, or 
\item
$2\deg L=\deg K_X$
\end{inparaenum}
can hold, we see that the case where $0\le 2\deg L<\deg K_X$ must correspond to the case where $(A_L,\Psi)$ is type $1$ with $\varphi\not\equiv 0$. The remaining two ``if'' assertions are proved in the same way.
\end{proof}

\subsubsection{Distinguishing type 1 and 2 quaternionic monopoles}
\label{subsec:Distinguishing_type_1_and_2_solutions_quaternionic}
A dichotomy also arises for solutions to the quaternionic monopole equations due to Okonek and Teleman \cite[Lemma 3.1]{OTVortex}. They consider pairs $(A,\Phi)$, where $A$ is a unitary connection on $E$ and $\Phi$ is a section of $W_{\can}^+\otimes E$, satisfying
\begin{equation}
\label{eq:U2MonopoleEquations}
\begin{aligned}
\rho(F^+_A)&=(\Phi\otimes\Phi^*)_0,
\\
D_A\Phi&=0,
\end{aligned}
\end{equation}
where $(\Phi\otimes\Phi^*)_0$ denotes the projection of $\Phi\otimes\Phi^*$ onto $\Om^0(\su(W^+)\otimes\fu(E))$.
Because equation \eqref{eq:U2MonopoleEquations} implies an equality between the $\gl(E)$-traces of $F_A^+$ and $(\Phi\otimes\Phi^*)_0$, an argument similar to that in Morgan \cite[Section 7.2]{MorganSWNotes} shows that the degree of $E$ determines the type of the pair.

\subsubsection{Distinguishing type 1 and 2 non-Abelian monopoles}
\label{subsec:Distinguishing_type_1_and_2_solutions_nonabelian}
The argument in the proof of Okonek and Teleman \cite[Lemma 3.1]{OTVortex} does not apply the non-Abelian monopole equations \eqref{eq:SO(3)_monopole_equations_Kaehler} because, unlike equation \eqref{eq:U2MonopoleEquations}, the non-Abelian monopole equations do not constrain $\tr_E F_A^+$.  Rather, the central component of the curvature of $A$ is constrained by the requirement that $A$ induces a fixed unitary connection $A_d$ on the Hermitian line bundle $\det E$. Because the non-Abelian monopole equations do not give a relationship between  $\tr_E F_A^+$ and the section $\Phi$, one cannot adapt the argument in \cite[Section 3]{OTVortex} to show that the degree of $E$ determines whether non-Abelian monopoles for the \spinu structure $\ft  = (\rho_{\can}\otimes \id_E, W_{\can}\otimes E)$ will be type $1$ or type $2$. Instead, we see in Section \ref{sec:HE_connections} that the connection $A_d$, whose curvature is $F_{A_d} = \frac{1}{2}\tr_E F_A$, must be Hermitian--Einstein (see Definition \ref{defn:HE_connection}), meaning that $F_{A_d}^{0,2}=0$ and $\La F_{A_d}=-i2\lambda$ where $\la$ is a constant, determined by the degree of $E$ in the forthcoming \eqref{eq:Einstein_factor}, and $E$ has complex rank two.

As described by Teleman in \cite[Corollary 2.3.3]{TelemanNonabelian} or by Dowker in \cite[Theorem 1.4.3]{DowkerThesis}, the moduli space of non-Abelian monopoles over a complex K\"ahler surface contains gauge-equivalence classes of type $1$ and type $2$ pairs as closed subspaces, whose intersection is the moduli space of projectively anti-self-dual connections (projectively Hermitian--Einstein connections, more generally). Both types of solutions may appear in the moduli space and in \cite[Theorem 7.1]{OTQuaternion}, Okonek and Teleman observe that for the preceding \spinu structure $\ft$ with $c_1(E)=K_X$ (the canonical class for $X$), the involution described in Remark \ref{rmk:Decoupling_Dirac_equation_and_equivalence_type_I_and_II_solutions}
will interchange the two types of pairs.

However, there are analogues of the Kodaira Vanishing Theorem for holomorphic vector bundles $E$ over complex K\"ahler manifolds $M$ of dimension $n$. Following Kobayashi \cite[Equation (2.1.1), p. 30]{Kobayashi_differential_geometry_complex_vector_bundles}, let $P(E)$ be the fiber bundle over $M$ with fibers $P(E_x) = (E_x\less\{0\})/\CC^*$ over $x\in M$ and let $L(E)$ be the tautological line bundle over $P(E)$. One calls $E$ \emph{negative} if $c_1(L(E)) <0$ (as discussed in Remark \ref{rmk:NegativityOfALineBundle})
and \emph{positive} if $E^*$ is negative  \cite[Section 3.5, p. 78]{Kobayashi_differential_geometry_complex_vector_bundles}. According to Kobayashi \cite[Theorem 3.5.9 or Corollary 3.5.10, pp. 78--79]{Kobayashi_differential_geometry_complex_vector_bundles} with $k=n$ and $E$ negative of rank $r\leq n$, one has
\[
  H^q(X;\bOmega^p(E)) = 0, \quad\text{for } p+q \leq n-r.
\]
In particular, choosing $p=q=0$ and $n=2$ and $r=2$, we see that
\[
  H^0(X;\sO(E))=0.
\]
Hence, the criterion that $E$ be negative can be used to rule out non-trivial type $1$ solutions $(A,\Phi)$ to \eqref{eq:SO(3)_monopole_equations_Kaehler} with $\Phi=(\varphi,0)$ and $\varphi$ a holomorphic section of $(E,\bar\partial_A)$, just as in Section \ref{subsubsec:Distinguishing_type_1_and_2_solutions_SW}.

Similarly, the criterion that $K_X\otimes E^*$ be negative can be used to rule out non-trivial type $2$ solutions $(A,\Phi)$ to \eqref{eq:SO(3)_monopole_equations_Kaehler} with $\Phi=(0,\psi)$ and $\bar\psi$ a holomorphic section of the complex vector bundle $K_X\otimes E^*$ with holomorphic structure induced by $\bar\partial_A$, just as in Section \ref{subsubsec:Distinguishing_type_1_and_2_solutions_SW}.

\subsection{Zero-section Seiberg--Witten monopoles on complex K\"ahler surfaces}
We note the following adaptation of Proposition \ref{prop:SWModuliSpaceZeroSectionCriterion} from smooth Riemannian metrics to complex, K\"ahler metrics.  Recall from Huybrechts \cite[Definition 3.2.14, p. 130]{Huybrechts_2005} that the \emph{K\"ahler cone} on a compact, complex manifold is the subset
\begin{equation}
\label{eq:KahlerCone}
\sK\subset H^{1,1}(X)\cap H^2(X;\RR),
\end{equation}
of cohomology classes $[\omega]$, where $\omega$ is \emph{positive} in the sense that, with respect to local holomorphic coordinates $z=(z_1,\ldots,z_n)$ on coordinate domains $U\subset X$, one has $\omega = \frac{i}{2}\sum_{i,j}h_{ij}(z)dz_i\wedge d\bar z_j$ and $h_{ij}(z)$ is a positive definite, Hermitian matrix for any $z\in U$ (see Huybrechts \cite[Lemma 3.1.7, p. 130]{Huybrechts_2005}). If $[\omega]\in\sK$ is defined by a Hermitian metric $h$ on $X$, then $h$ is K\"ahler.

\begin{lem}[Non-existence of zero-section Seiberg--Witten monopoles over complex K\"ahler surfaces]
\label{lem:ExistenceOfZeroSectionSWonKahler}
Let $\fs$ be a \spinc structure and $\ft=(\rho,W\otimes E)$ a rank-two \spinu structure on a closed, complex K\"ahler surface. Consider the moduli space $\sM_\ft$ of non-Abelian monopoles defined by the K\"ahler metric $g$
and the perturbations $\tau$ and $\vartheta$ in \eqref{eq:PerturbedSO3MonopoleEquations} given by $\tau=\id_{\La^+}$ and $\vartheta=0$. Assume that $c_1(\fs)\in\Red_0(\ft)$, where $\Red_0(\ft)$ is defined in \eqref{eq:DefineReduciblesEmbedded}, so there is a continuous embedding $\iota_{\fs,\ft}:M_\fs\hookrightarrow \sM_\ft$. 
Assume that the unitary connection $A_\La$ on $\det W^+\otimes \det E$ appearing in the perturbed Seiberg--Witten equations \eqref{eq:SeibergWitten} satisfies $F_{A_\Lambda}^{0,2}=0$. If $c_1(\fs)-c_1(\ft)$ is not torsion, then there is an open, dense subspace $\sU$ of the K\"ahler cone \eqref{eq:KahlerCone} such that if the metric on $X$ is K\"ahler with fundamental two-form $[\om]\in\sU$, then $\iota_{\fs,\ft}(M_\fs)\subset\sM_\ft$ does not contain the gauge-equivalence class of a zero-section pair.
\end{lem}

\begin{proof}
As in Section \ref{subsubsec:Distinguishing_type_1_and_2_solutions_SW}, we write the spin${}^c$ structure $\fs$ as $\fs=\fs_{\can}\otimes L$ where $\fs_\can=(\rho_{\can},W_{\can})$ is the canonical spin${}^c$ structure defined in Definition \ref{defn:Canonical_spinc_bundles}.  We fix the connection $A_{\can}$ on $W_{\can}$ as defined in Section \ref{subsubsec:Distinguishing_type_1_and_2_solutions_SW}. If $[A,0]\in\sM_\ft$ is in the image of $\iota_{\fs,\ft}$, then by the discussion in Section \ref{subsec:RedPU2Monopole}, we have $(A,0)=\tilde\iota_{\fs,\ft}(A_L,0)$ where
$(A_L,0)$ is a solution to the perturbed Seiberg--Witten equations \eqref{eq:SeibergWitten}. Let $\La_\om$ be the dual Lefschetz operator \eqref{eq:Huybrechts_definition_1-2-21_and_lemmas_1-2-23_and_24} associated to the fundamental two-form $\om$ of a K\"ahler metric.   By writing the $(1,1)$ and $(0,2)$ components of the curvature equation  \eqref{eq:SeibergWitten_curvature} separately, we see that the curvature of $A_L$ satisfies
\begin{subequations}
\begin{align}
\label{eq:SWCurvature(1,1)Component}
\La_\om\left( \tr_{W^+_{\can}} F_{A_{\can}}+2 F_{A_L}- F_{A_\Lambda}\right)&=0,
\\
\label{eq:SWCurvature(0,2)Component}
\tr_{W^+_{\can}} F_{A_{\can}}^{0,2} +2F_{A_L}^{0,2} - F_{A_\Lambda}^{0,2}&=0,
\end{align}
\end{subequations}
where $A_{\can}$ is the spin connection on $(\rho_\can,W_{\can})$ introduced in Section \ref{subsubsec:Distinguishing_type_1_and_2_solutions_SW}. 
The $(0,2)$ component of \eqref{eq:Curvature_of_Canonical_SpinConnection} implies that
\[
\tr_{W^+_{\can}} F_{A_{\can}}^{0,2}
=
-F_K^{0,2}
=0.
\]
The equality $\tr_{W^+_{\can}} F_{A_{\can}}^{0,2}=0$, equation \eqref{eq:SWCurvature(0,2)Component}, and the hypothesis $F_{A_\Lambda}^{0,2}=0$ imply that $F_{A_L}^{0,2}=0$.  Hence, because $A_L$ and $A_\Lambda$ are unitary connections, we see that
\begin{equation}
\label{eq:PerturbedCurvature}
\frac{1}{2\pi i}\left( \tr_{W^+_{\can}} F_{A_{\can}}+2 F_{A_L} - F_{A_\Lambda}\right)
\end{equation}
is a closed, real $(1,1)$-form representing the cohomology class $c_1(\fs)-c_1(\ft)$. By Huybrechts \cite[Corollary 3.1.8, p. 117]{Huybrechts_2005}, the set of K\"ahler forms is an open convex cone in the real vector space
\begin{equation}
\label{eq:Closed(1,1)Forms}
\{\om\in \Om^{1,1}(X): d\om=0\}\cap \Om^2(X;\RR).
\end{equation}
Because $c_1(\fs)-c_1(\ft)$ is not a torsion class, the form \eqref{eq:PerturbedCurvature} is non-zero. Equation \eqref{eq:SWCurvature(1,1)Component} is equivalent to $\om$ being orthogonal to the non-zero, real $(1,1)$-form \eqref{eq:PerturbedCurvature} and so equation \eqref{eq:SWCurvature(1,1)Component} defines a hyperplane in the vector space of closed, real $(1,1)$-forms \eqref{eq:Closed(1,1)Forms} and thus in $\sK$, the K\"ahler cone \eqref{eq:KahlerCone}.  Let $\sU\subset\sK$ be the complement of this hyperplane. The existence of a solution $(A_L,0)$ to the Seiberg--Witten equations \eqref{eq:SeibergWitten} implies that the K\"ahler form $\om$ lies on this hyperplane and thus $[\om]\in\sK\less\sU$. Hence, if the cohomology class $[\om]$ of the K\"ahler form lies in $\sU$, there are no zero-section solutions to the perturbed Seiberg--Witten equations \eqref{eq:SeibergWitten} as asserted.
\end{proof}

\chapter[Elliptic complexes]{Elliptic complexes}
\label{chap:Elliptic_complexes}
In this chapter, we describe the elliptic complexes and associated cohomology groups and harmonic spaces that we shall need in this monograph. We begin in Section \ref{sec:Kuranishi_models} by reviewing Kuranishi models defined by analytic Fredholm maps of Banach manifolds and then proving the key Lemma \ref{lem:EquivariantKuranishiLemma} that establishes the existence of an equivariant Kuranishi model and builds on an earlier result \cite[Lemma 4.7]{FU} due to Freed and Uhlenbeck. Section \ref{sec:Elliptic_complexes} provides an introduction to abstract elliptic complexes that underpins much of our work in this monograph. Section \ref{sec:Elliptic_deformation_complex_for_SO3_monopole_equations} reviews the elliptic deformation complex for the non-Abelian monopole equations over a closed, four-dimensional, oriented, smooth Riemannian manifold and is the first of a group of sections in this chapter that discusses elliptic complexes for systems of nonlinear equations arising in gauge theory relevant to the main results of this monograph. In Section \ref{sec:Elliptic_deformation_complex_for_projective_vortex_equations}, we describe the elliptic complex for the projective vortex equations over an almost Hermitian manifold. Section \ref{sec:Elliptic_deformation_complex_for_pre-holomorphic_pair_equations_complex_surface} introduces the elliptic deformation complex for the pre-holomorphic pair equations over a complex K\"ahler surface --- equations that arise from the non-Abelian monopole equations over a four-dimensional Riemannian manifold. In Section \ref{sec:Elliptic_deformation_complex_holomorphic_pair_equations}, we review the elliptic complex for the holomorphic pair equations over a complex manifold. Section \ref{sec:Fredholm_maps_elliptic_deformation_complexes} provides an introduction to abstract complexes defined by Fredholm maps and sections of sheaves. In Section \ref{sec:Fredholm_complexes_Euler_characteristics}, we define the Euler characteristic of these Fredholm complexes and prove a lemma that identifies the homologies of the elliptic complexes of the pre-holomorphic pair equations and holomorphic pair equations. Lastly, in Section \ref{sec:Witten_dichotomy_for_linearization_SO3_monopole_equations}, we discuss Witten's dichotomy for solutions to the linearized non-Abelian monopole equations and the associated vanishing of cohomology groups.

\section{Kuranishi models defined by analytic Fredholm maps of Banach manifolds}
\label{sec:Kuranishi_models}
Given a smooth map $F:\sX\to\sY$ of smooth, real or complex Banach manifolds and a smoothly embedded submanifold $\sZ \subset \sY$, the \emph{Kuranishi method} \cite{Kuranishi} aims to describe the local structure of $F^{-1}(\sZ)$ near a point $x_0$. A standard procedure (see Guillemin and Pollack \cite[p. 28]{Guillemin_Pollack}) allows us to reduce this question to one where $\sZ$ is a single point, say $y_0 \in \sY$. Furthermore, since we are interested in the local structure of $F^{-1}(y_0)$, we may further reduce this question to one where $\sX$ and $\sY$ are replaced by real or complex Banach spaces, $\cX$ and $\cY$, by choosing local coordinate charts $\phi:\sX \subset \sU \to \cX$ and  $\psi:\sY \subset \sV \to \cY$ and considering the composed map $\psi\circ F \circ \phi^{-1}:\cX \supset \phi(\sU) \to \cY$. Thus, it suffices to consider the case of a smooth map $F:\cX \supset \cU \to \cY$ and the preimage $F^{-1}(y_0)$ in an open neighborhood $\cU$ of a point $x_0$. We may assume without loss of generality that $x_0$ is the origin in $\cX$ and $y_0$ is the origin in $\cY$.

If $dF(0):\cX\to\cY$ is surjective then, after possibly shrinking $\sU$, the preimage $F^{-1}(0)$ is a smoothly embedded submanifold of $\cU$ (see Abraham, Marsden, and Ratiu \cite[Theorem 3.5.12, p. 202]{AMR}, Lang \cite[Proposition 2.5, p. 29]{Lang_introduction_differential_topology}, or Margalef Roig and Outerelo Dom{\'\i}nguez \cite[Proposition 7.1.14, p. 319]{Margalef-Roig_Outerelo-Dominguez_differential_topology}).
With our applications in mind, we focus on the case where $dF(0):\cX\to\cY$ is not necessarily surjective but $\cY_0 := \Ran dF(0)$ has finite codimension in $\cY$ and is thus a closed subspace of $\cY$ with closed complement $C$, yielding a topological direct sum $\cY = \cY_0\oplus C$. By construction, $F$ is transverse to $C$ at the point $0 \in \cX$ and thus, again after possibly shrinking $\cU$, the preimage $F^{-1}(C)$ is a smoothly embedded submanifold of $\cU$ (see the references cited above) with tangent space
\[
  T_0F^{-1}(C) = dF(0)^{-1}(T_0C) =  dF(0)^{-1}(C) = dF(0)^{-1}(0) = \Ker dF(0),
\]
equal to the Zariski tangent space to $F^{-1}(0)$ at $0$. Denote $K := \Ker dF(0)$, a closed subspace of $\cX$ and assume, again with our applications in mind, that $K$ is complemented in $\sX$ by a closed subspace $\cX_0$, and therefore $\cX = \cX_0 \oplus K$. The Implicit Mapping Theorem (see, for example, Abraham, Marsden, and Ratiu \cite[Theorem 2.5.7, p. 121]{AMR}), provides a smooth embedding $\gamma:K\cap\cU \to \cX$ such that $\gamma(K\cap\cU) = F^{-1}(C)$ and thus we obtain a smooth map $\chi:K\cap\cU \to C$ defined by
\[
  \chi(x) = (F\circ\gamma)(x), \quad\text{for all } x\in K\cap\cU, 
\]  
with the resulting property that
\[
  F^{-1}(0) = \chi^{-1}(0) \subset \cU.
\]
If we also assume that $F$ is analytic, the Implicit Mapping Theorem ensures that $\gamma$ and hence $\chi$ are analytic too (see, for example, Feehan and Maridakis \cite[Theorem F.1 and Remark F.2, p. 127]{Feehan_Maridakis_Lojasiewicz-Simon_coupled_Yang-Mills}). Finally, if we further assume that $K$ is finite-dimensional, then $F^{-1}(0)$ has the structure of a finite-dimensional analytic variety over a field $\KK=\RR$ or $\CC$.

The preceding paragraphs provide an intrinsic approach to constructing the Kuranishi model for an open neighborhood of a singular point in the zero locus $F^{-1}(0)$ that does not rely on choices of local coordinate charts. Although a brief outline of a proof of the following lemma that does rely on choices of local coordinate charts appears in Freed and Uhlenbeck \cite[Lemma 4.7]{FU} and Donaldson and Kronheimer \cite[Proposition 4.2.19, p. 136]{DK}, we provide a detailed proof here in order to specify the exact appearance of equivariance in their argument.

\begin{lem}[Equivariant Kuranishi model]
\label{lem:EquivariantKuranishiLemma}
(See Freed and Uhlenbeck \cite[Lemma 4.7]{FU}.)
Let $G$ be a Lie group, $E_1$ and $E_2$ be Banach spaces over $\KK=\RR$ or $\CC$, and $\rho_i:G\to \GL(E_i)$ be a homomorphism of groups for $i=1,2$. Let $k$ be a positive integer or $\infty$ or\footnote{As customary, we use $C^\omega$ to indicate a $\KK$-analytic map.} $\omega$ and $S:U\to E_2$ be a $G$-equivariant, $C^k$  map from an open neighborhood $U$ of the origin in $E_1$ such that $S(0)=0$ and the derivative $L = DS(0) \in \Hom(E_1,E_2)$ is Fredholm.  Assume that the splittings\footnote{The splittings in the sense of Banach spaces necessarily exist since $\Ker L$ has finite dimension over $\KK$ and $\Ran L$ has finite codimension over $\KK$ --- see Rudin \cite[Lemma 4.12, p. 106]{Rudin}.} of Banach spaces
\begin{subequations}
\label{eq:FredholmSplitting}
\begin{align}
\label{eq:FredholmSplittingDomain}
  E_1&= \Ker L\oplus F_1,
  \\
  \label{eq:FredholmSplittingCoDomain}
  E_2&=\Ran L\oplus C_2,
\end{align}
\end{subequations}
are such that the subspaces $F_1$ and $C_2$ are closed under the $G$ actions defined by $\rho_1$ and $\rho_2$, respectively. Then, after possibly shrinking $U$ and taking its $G$-orbit, there are
\begin{enumerate}
\item
  \label{item:EquivariantKuranishiExistenceOfEmbedding}
  A $G$-equivariant, $C^k$ embedding $\bga:U \to E_1$ onto an open neighborhood of the origin in $E_1$, and
\item
\label{item:EquivariantKuranishiExistenceOfObstruction}
A $G$-equivariant, $C^k$ map $\bchi:U\to C_2$
\end{enumerate}
that obey the following additional properties:
\begin{enumerate}
\setcounter{enumi}{2}
\item
\label{item:EquivariantKuranishiOriginPreserving}
$\bga(0)=0$;
\item
\label{item:EquivariantKuranishiDecompositionOfMap}
$S\circ \bga^{-1}=L+\bchi$;
\item
\label{item:EquivariantKuranishiObstructionMapValues}
$\bchi(0)=0$ and $D\bchi(0)=0$; and
\item
  \label{item:EquivariantKuranishiEmbeddingHomeomorphismOfZeroSet}
$\bgamma:U\to E_1$ restricts to a homeomorphism
\[
\bga: S^{-1}(0) \to \bchi^{-1}(0)\cap\Ker L.
\]
\end{enumerate}
\end{lem}

One calls $\bga$ and $\bchi$ in Lemma \ref{lem:EquivariantKuranishiLemma} the \emph{embedding map}
\label{page:Embedding_obstruction_Kuranishi} and \emph{obstruction map}, respectively.

\begin{rmk}[On the hypotheses of Lemma \ref{lem:EquivariantKuranishiLemma}]
\label{rmk:KuranishiModelEquivariantSplitting}
Because the derivative $L$ in Lemma \ref{lem:EquivariantKuranishiLemma} is Fredholm, there is always a splitting of $E_1$ and $E_2$ in the sense of Banach spaces as in \eqref{eq:FredholmSplitting} since $\Ker L$ and $\Coker L$ are finite-dimensional (see Rudin \cite[Lemma 4.21, p. 106]{Rudin}). The subspaces $\Ker L \subset E_1$ and $\Ran L \subset E_2$ will always be closed under the $G$ actions but it is not immediately apparent that the complements, $F_1$ and $C_2$, can be chosen to be closed under the $G$ actions. If $E_1$ and $E_2$ are Hilbert spaces and the representations $\rho_1$ and $\rho_2$ are orthogonal or unitary,
as can always be assumed if $G$ is compact by replacing the inner products with their averages with respect over $G$ (using integration with respect to the Haar measure), as done in Br\"ocker and tom Dieck \cite[Chapter II, Theorem 1.7, p. 68]{BrockertomDieck} or Knapp \cite[Proposition 4.6, p. 240]{Knapp_1986} in the case of finite-dimensional representations (the same averaging argument extends without change to infinite-dimensional representations) and we choose $F_1 := (\Ker L)^\perp$ and $C_2 := (\Ran L)^\perp$ (orthogonal complements with respect to the inner products in $H_1$ and $H_2$, respectively), then these orthogonal complements are also closed under the $G$ actions.
\end{rmk}

\begin{proof}[Proof of Lemma \ref{lem:EquivariantKuranishiLemma}]
Because $S$ is $G$-equivariant, the derivative $L=DS(0)\in\Hom(E_1,E_2)$ is also $G$-equivariant and so $\Ker L$ and $\Ran L$ are closed under the $G$ actions defined by $\rho_1$ and $\rho_2$, respectively. Let $\pi_K:E_1\to \Ker L$ and $\pi_C:E_2\to C_2$ be the projections
associated with the splittings \eqref{eq:FredholmSplitting}. The hypothesis that $F_1$ and $C_2$ are closed under the $G$ actions defined by $\rho_1$ and $\rho_2$, respectively, implies that for all $g\in G$, the operator $\rho_1(g)\in\End(E_1)$ is diagonal with respect to the splitting \eqref{eq:FredholmSplittingDomain} and the operator $\rho_2(g)\in\End(E_2)$ is diagonal with respect to the splitting \eqref{eq:FredholmSplittingCoDomain}. The projections $\pi_K$ and $\pi_C$ are thus $G$-equivariant.

By the Open Mapping Theorem, the restriction of $L$ to $F_1$ has a unique bounded inverse, $L^{-1} \in \Hom(\Ran L, F_1)$.  Taking the inverse of both sides of the equality $\rho_2(g^{-1})\circ L\circ\rho_1(g)=L$ gives $\rho_1(g^{-1})\circ L^{-1}\circ \rho_2(g)=L^{-1}$, for all $g\in G$, implying that $L^{-1}$ is $G$-equivariant. Extend $L^{-1}$ from $\Ran L$ to $E_2$ using the composition $L^{-1}\circ (\id_{E_2}-\pi_C)$. Because $\pi_C$ and hence $\id_{E_2}-\pi_C$ are $G$-equivariant, the operator $L^{-1}\circ (\id_{E_2}-\pi_C)$ is $G$-equivariant. We claim that, after possibly shrinking $U$, the $C^k$ map
\begin{equation}
  \label{eq:DefineDiffeomorphism}
  \bgamma = \pi_K+L^{-1}\circ(\id_{E_2}-\pi_C)\circ S: U \to E_1
\end{equation}
satisfies the properties claimed for $\bgamma$ in Item \eqref{item:EquivariantKuranishiExistenceOfEmbedding}. By its definition in terms of the $G$-equivariant maps $\pi_K$ and $L^{-1}\circ (\id_{E_2}-\pi_C)\circ S$, the map $\bgamma$ is $G$-equivariant. Observe that $\bga(0)=0$, so $\bga$ also satisfies Item \eqref{item:EquivariantKuranishiOriginPreserving}. The derivative of $\bga$ at
origin in $E_1$ is
\[
D\bgamma(0) = \pi_K + L^{-1}\circ (\id_{E_2}-\pi_C)\circ L = \id_F \in \End(E_1),
\]
because $L^{-1}\circ (\id_{E_2}-\pi_C)\circ L = \id_F-\pi_K$. By the Inverse Mapping Theorem for $C^k$ maps of Banach spaces (see, for example, Feehan and Maridakis \cite[Theorem F.1 and Remark F.2, p. 127]{Feehan_Maridakis_Lojasiewicz-Simon_coupled_Yang-Mills}), after possibly shrinking the open neighborhood $U$ of the origin in $E_1$, the map $\bgamma$ satisfies Item \eqref{item:EquivariantKuranishiExistenceOfEmbedding}.

If $G$ acts by isometries, then every ball around the origin is $G$-invariant and we can assume that the neighborhood $U$ of the origin is such a ball.  More generally, we can replace $U$ by its $G$-orbit,
\[
G\cdot U:=\{\rho_1(g)u: g\in G,\, u\in U\}.
\]
We show that $\bgamma$ remains injective on $G\cdot U$ as follows.  If $\bgamma(\rho_1(g_1)u_1)=\bgamma(\rho_1(g_2)u_2)$ for $g_1,g_2\in G$ and $u_1,u_2\in U$, then the equivariance of $\bgamma$ implies that $\rho_2(g_2^{-1}g_1)\bgamma(u_1)=\bgamma(u_2)$. Thus, $\rho_2(g_2^{-1}g_1)\bgamma(u_1)\in\bgamma(U)$ and so $\bgamma^{-1}(\rho_2(g_2^{-1}g_1)\bgamma(u_1))$ is well-defined, since we know that $\gamma:U \to \gamma(U)$ is a bijection. Because $\bgamma^{-1}$ is $G$-equivariant, the equality $\rho_2(g_2^{-1}g_1)\bgamma(u_1)=\bgamma(u_2)$
implies that
\[
u_2= \bgamma^{-1}\left(\rho_2(g_2^{-1}g_1)\bgamma(u_1) \right)
=
\rho_1(g_2^{-1}g_1)\bgamma^{-1}(\bgamma(u_1))
=
\rho_1(g_2^{-1}g_1)u_1.
\]
The resulting equality, $u_2=\rho_1(g_2^{-1}g_1)u_1$, implies that $\rho_1(g_2)u_2=\rho_1(g_1)u_1$ and we have shown that $\bgamma$ is injective on $G\cdot U$.  The fact that the derivative of $\bgamma$ is an isomorphism on the larger set $G\cdot U$ follows by differentiating the equality $\rho_2(g^{-1})\bgamma(\rho_1(g) u)=\bgamma(u)$ with respect to $u$ to give
\[
D\rho_2(g^{-1})\circ D\bgamma(\rho_1(g)u)\circ D\rho_1(g)
=
D\bgamma(u),
\]
and by observing that the differentials $D\rho_2(g^{-1})$ and $D\rho_1(g)$ are isomorphisms.

We claim that the $C^k$ map
\begin{equation}
  \label{eq:KuranishiLemmaDefineObstructionMap}
  \bchi =\pi_C\circ S\circ \bgamma^{-1}:U\to C_2
\end{equation}
satisfies the properties in Item \eqref{item:EquivariantKuranishiExistenceOfObstruction}. By the argument given earlier to prove that $L^{-1}$ is $G$-equivariant, we see that $\bgamma^{-1}$ is  $G$-equivariant as well.
The definition of $\bchi$ as the composition of the $G$-equivariant, $C^k$ maps $\pi_C$ and $S$ and $\bgamma^{-1}$ implies that $\bchi$ is a $G$-equivariant, $C^k$ map and thus satisfies Item \eqref{item:EquivariantKuranishiExistenceOfObstruction}. Because $S(0)=0$, the definition \eqref{eq:DefineDiffeomorphism} of $\bgamma$ implies that $\bga(0)=0$ and so $\bga$ satisfies Item \eqref{item:EquivariantKuranishiOriginPreserving}.

We prove that $\bgamma$ and $\bchi$ satisfy Item \eqref{item:EquivariantKuranishiDecompositionOfMap} as follows. Because $L\circ\pi_K=0$, the definition \eqref{eq:DefineDiffeomorphism} of $\bgamma$ yields
\[
L\circ \bgamma=(\id_{E_2}-\pi_C)\circ S,
\]
and thus $L=(\id_{E_2}-\pi_C)\circ S\circ \bgamma^{-1}$.  Combining this equality with the definition \eqref{eq:KuranishiLemmaDefineObstructionMap} of $\bchi$
gives
\begin{equation}
\label{eq:KuranishiLemmaDecompositionOfMap}
S\circ \bgamma^{-1}=L+\bchi,
\end{equation}
so $\bchi$ satisfies Item \eqref{item:EquivariantKuranishiDecompositionOfMap}. The definition \eqref{eq:KuranishiLemmaDefineObstructionMap} of $\bchi$ and the equalities $S(\bga^{-1}(0))=S(0)=0$ imply that $\bchi(0)=0$. Because $\pi_C\circ L=0$ and
$(D\bga^{-1})(0) = D\bga(0)^{-1} = \id_{E_1}$,
the definition \eqref{eq:KuranishiLemmaDefineObstructionMap} of $\bchi$ implies that $D\bchi(0)=0$. Hence, $\bchi$ satisfies Item  \eqref{item:EquivariantKuranishiObstructionMapValues}.

Finally, we verify the assertion in Item \eqref{item:EquivariantKuranishiEmbeddingHomeomorphismOfZeroSet}. Equation \eqref{eq:KuranishiLemmaDecompositionOfMap} implies that $\bgamma$ defines a homeomorphism between $S^{-1}(0)$ and $(L+\bchi)^{-1}(0)$. Because the images of $L$ and $\bchi$ are contained in the factors $\Ran L$ and $C_2$ of the Banach-space splitting \eqref{eq:FredholmSplittingCoDomain} of $E_2$, the definition of $\bchi$ implies that
\[
  (L+\bchi)^{-1}(0) = L^{-1}(0)\cap \bchi^{-1}(0) = \Ker L\cap\bchi^{-1}(0).
\]
Thus, $\bgamma$ defines a homeomorphism between $S^{-1}(0)$ and $\Ker L\cap\bchi^{-1}(0)$, as asserted in Item \eqref{item:EquivariantKuranishiEmbeddingHomeomorphismOfZeroSet}.
\end{proof}

\section{Elliptic complexes}
\label{sec:Elliptic_complexes}
Since we shall be discussing many different elliptic complexes and how they are related, it is convenient to recall the general theory of elliptic complexes from Atiyah and Bott \cite{Atiyah_Bott_1967, Atiyah_Bott_1968}, Gilkey \cite[Section 1.5]{Gilkey2}, H\"ormander \cite[Chapter 14]{Hormander_v3}, Kumano-go \cite[Section 9.1]{Kumanogo}, and Wells \cite[Section 4.5]{Wells3}. We partly adapt the exposition due to Wells here.

Let $\{E_k\}_{k=0}^n$ (with $n\geq 1$) be a sequence of smooth (real or complex) vector bundles defined over a closed, smooth manifold $X$ and $\{d_k\}_{k=0}^{n-1}$ of (real or complex) linear differential operators of some fixed order \cite[Chapter IV, Section 2, p. 114]{Wells3} (but always first-order in our applications), to form a sequence \cite[Chapter IV, Equation (5.1), p. 144]{Wells3}:
\begin{equation}
  \label{eq:Elliptic_complex}
  0 \xrightarrow{} C^\infty(E_0)  \xrightarrow{d_0}  C^\infty(E_1)  \xrightarrow{d_1} \cdots \xrightarrow{d_{n-2}}   C^\infty(E_{n-1}) \xrightarrow{d_{n-1}}  C^\infty(E_n)  \xrightarrow{} 0
\end{equation}
This sequence of operators and spaces of smooth sections of vector bundles \eqref{eq:Elliptic_complex} is called a \emph{(cochain) complex} if $d_k\circ d_{k-1} = 0$ for $k = 1, \ldots, n-1$. We define $E_k := 0$ for all integers $k<0$ or $k>n$ and $d_k := 0$ for all integers $k<0$ or $k>n-1$ and define the \emph{cohomology groups} (real or complex vector spaces) associated to the complex \eqref{eq:Elliptic_complex} by \cite[Chapter IV, Equation (5.3), p. 145]{Wells3}
\begin{multline}
  \label{eq:Cohomology_groups_elliptic_complex}
  H^k(E) := \left.\Ker\left(d_k:C^\infty(E_k)\to C^\infty(E_{k+1})\right)\right/
  \Ran\left(d_{k-1}:C^\infty(E_{k-1})\to C^\infty(E_k)\right),
  \\
  \text{for } k = 0, \ldots,n.
\end{multline}
Let $X\times\{0\} \subset T^*X$ denote the zero section and let $\pi^*E_k$ denote the pullback of $E_k$ to 
$T^*X\less (X\times\{0\})$ via the projection map $\pi:T^*X\less(X\times\{0\}) \to X$ for $k=0,\ldots,n$. To each $d_k$, one can associate its symbol $\sigma(d_k) \in \Hom(\pi^*E_k, \pi^*E_{k+1})$ (see \cite[Chapter IV, Equation (2.1), p. 115]{Wells3}) and thus obtain the associated symbol sequence \cite[Chapter IV, Equation (5.2), p. 144]{Wells3}:
\begin{equation}
  \label{eq:Symbol_complex}
  0 \xrightarrow{} \pi^*E_0  \xrightarrow{\sigma(d_0)}  \pi^*E_1  \xrightarrow{\sigma(d_1)} \cdots \xrightarrow{\sigma(d_{n-2})} \pi^*E_{n-1} \xrightarrow{\sigma(d_{n-1})} \pi^*E_n  \xrightarrow{} 0
\end{equation}
The complex \eqref{eq:Elliptic_complex} is called \emph{elliptic} \label{page:Elliptic_complex} if the associated symbol sequence \eqref{eq:Symbol_complex} is \emph{exact} at each non-zero $\xi\in T^*X$ (see Wells \cite[Chapter IV, Definition 5.1, p. 144]{Wells3}). The real (respectively, complex) vector bundles $E_k$ may be equipped with smooth, Riemannian (respectively, Hermitian) metrics, with respect to which one may define the $L^2$ or \emph{formal adjoint differential operators} \cite[Chapter IV, Proposition 2.8, p. 117]{Wells3} \label{page:Adjoint_Differential_Operator}
\[
  d_k^*:C^\infty(E_{k+1}) \to C^\infty(E_k), \quad\text{for } k = 0, \ldots,n-1,
\]
and \emph{Laplace operators} \cite[Chapter IV, Section 5, p. 145]{Wells3}
\[
  \Delta_k := d_k^*d_k + d_{k-1}d_{k-1}^*:C^\infty(E_k) \to C^\infty(E_k), \quad\text{for } k = 0, \ldots,n.
\]
The partial differential operators $\Delta_k$ are elliptic and self-adjoint \cite[Chapter IV, Section 5, p. 145]{Wells3}, with finite-dimensional kernels \cite[Chapter IV, Theorem 4.12, p. 141]{Wells3},
\[
  \dim\Ker\Delta_k < \infty, \quad\text{for } k = 0, \ldots,n.
\]  
If the complex $(d_\bullet,C^\infty(E_\bullet))$ in \eqref{eq:Elliptic_complex} is elliptic, then \cite[Chapter IV, Theorem 5.2 and Proposition 5.3, p. 147]{Wells3}, there are canonical isomorphisms 
\[
  H^k(E) \cong \bH^k(E), \quad\text{for } k = 0, \ldots,n,
\]
between the cohomology groups \eqref{eq:Cohomology_groups_elliptic_complex} and the harmonic spaces (real or complex vector spaces)
\begin{equation}
  \label{eq:Harmonic_spaces_elliptic_complex}
  \bH^k(E) := \Ker\left(d_k+d_{k-1}^*:C^\infty(E_k)\to C^\infty(E_{k+1}\oplus E_{k-1})\right), \quad\text{for } k = 0, \ldots,n,
\end{equation}
and orthogonal (Hodge) decompositions
\begin{equation}
  \label{eq:Hodge_decomposition_all_k-forms}
  C^\infty(E_k) = \bH^k(E) \oplus \Ran d_{k-1} \oplus \Ran d_{k+1}^*, \quad\text{for } k = 0, \ldots,n.
\end{equation}
Note that by Rudin \cite[Theorem 4.12 and Section 4.6]{Rudin} (with orthogonal complements replacing annihilators), we have
\begin{equation}
  \label{eq:Rudin_theorem_4-12}
  \Ker d_{k-1}^* = (\Ran d_{k-1})^\perp \quad\text{and}\quad \Ker d_{k+1} = (\Ran d_{k+1}^*)^\perp.
\end{equation}
Moreover, the harmonic spaces $\bH^k(E)$ may be alternatively expressed as \cite[Chapter IV, Proposition 5.3, p. 147]{Wells3}
\[
  \bH^k(E) = \Ker\Delta_k, \quad\text{for } k = 0, \ldots,n.
\]  
When the operators $d_k$ are first-order, it is standard (see, for example, Feehan \cite{Feehan_yang_mills_gradient_flow_v4}) that the complex \eqref{eq:Elliptic_complex} and preceding orthogonal decompositions extend from the category of $C^\infty$ sections of $E_k$ to $W^{n-k,p}$ sections of $E_k$, for $k=0,\ldots,n$ and $p \in (1,\infty)$. (Similar remarks apply to sequences of differential operators $d_k$ of order $m\geq 1$.) The \emph{index} (or \emph{Euler characteristic}) of an elliptic complex \eqref{eq:Elliptic_complex} is defined by (see Gilkey \cite[Equation (1.5.17)]{Gilkey2})
\begin{equation}
  \label{eq:Index_elliptic_complex}
  \Euler_\KK(d_\bullet) := \sum_{k\in\ZZ} (-1)^k\dim_\KK H^k(E),
\end{equation}
where $\KK = \RR$ or $\CC$. We recall that (see Gilkey \cite[Equation (1.5.18)]{Gilkey2})
\begin{equation}
  \label{eq:Index_rolled-up_elliptic_complex}
  \Euler_\KK(d_\bullet) = \Euler_\KK(d_\even),
\end{equation}
where
\[
  d_\even : C^\infty(E_\even) \to C^\infty(E_\odd),
\]
is defined by
\[
  d_\even := \bigoplus_{k\text{ even}} (d_k + d_{k-1}^*),
  \quad
  E_\even := \bigoplus_{k\text{ even}} E_k,
  \quad\text{and}\quad
  E_\odd := \bigoplus_{k\text{ odd}} E_k.
\]
According to Gilkey \cite[Lemma 1.5.1]{Gilkey2}, the complex $(d_\bullet,C^\infty(E_\bullet))$ is elliptic if and only if the Laplacians $\Delta_k$ are elliptic for all $k$ or the operator $d_\even$ is elliptic.

Furthermore, by Gilkey \cite[Lemma 1.4.5]{Gilkey2}, the index of an elliptic partial differential operator\footnote{Or pseudodifferential operator, more generally.} depends only the homotopy type of the leading symbol of $P$ within the class of elliptic pseudodifferential operators of the same order. Consequently, by the equality \eqref{eq:Index_rolled-up_elliptic_complex}, the index of the elliptic complex \eqref{eq:Elliptic_complex} depends only the homotopy type of on the leading symbols $\sigma(d_k)$ of the partial differential operators $d_k$ for $k\in\ZZ$ within the class of pseudodifferential operators of the same order.

We shall often encounter decompositions of a complex \eqref{eq:Elliptic_complex} into a direct sum of subcomplexes of the form (compare Weibel \cite[Section 1.2]{Weibel_introduction_homological_algebra})
\begin{equation}
  \label{eq:Elliptic_complex_pm}
  0 \xrightarrow{} C^\infty(E_0^\pm)  \xrightarrow{d_0^\pm}  C^\infty(E_1^\pm)  \xrightarrow{d_1^\pm} \cdots \xrightarrow{d_{n-2}^\pm}   C^\infty(E_{n-1}^\pm) \xrightarrow{d_{n-1}^\pm}  C^\infty(E_n^\pm)  \xrightarrow{} 0
\end{equation}
where the vector bundles and differentials are related by
\[
  E_k = E_k^+\oplus E_k^- \quad\text{and}\quad d_k = d_k^+\oplus d_k^-, \quad\text{for all } k \in \ZZ.
\]  
Their cohomology groups clearly satisfy (see Rotman \cite[Exercise 10.30]{RotmanAlgebra_2002})
\begin{equation}
\label{eq:EllipticCohomologyDecomposition}
  H^k(E) = H^k(E^+)\oplus H^k(E^-), \quad\text{for all } k \in \ZZ.
\end{equation}
Because the leading symbols obey
\[
  \sigma(d_k) = \sigma(d_k^+) \oplus \sigma(d_k^-),   
\]
where $\sigma(d_k^+) \in \Hom(\pi^*E_k^\pm, \pi^*E_{k+1}^\pm)$, for all $k\in\ZZ$, we see that the complex $(d_\bullet,C^\infty(E_\bullet))$ in \eqref{eq:Elliptic_complex} is elliptic if and only if each of the complexes $(d_\bullet^\pm,C^\infty(E_\bullet^\pm))$ in \eqref{eq:Elliptic_complex_pm} are elliptic.

We shall frequently also have to relate cohomology groups defined by two different complexes. Thus, suppose we are given another complex in the sense of \eqref{eq:Elliptic_complex}
\[
  0 \xrightarrow{} C^\infty(E_0')  \xrightarrow{d_0'}  C^\infty(E_1')  \xrightarrow{d_1'} \cdots \xrightarrow{d_{n-2}'}   C^\infty(E_{n-1}') \xrightarrow{d_{n-1}'}  C^\infty(E_n')  \xrightarrow{} 0
\]
and homomorphisms $f_k:E_k\to E_k'$ of vector bundles inducing a \emph{morphism} or \emph{cochain map} between the complexes $(C^\infty(E_\bullet),d_\bullet)$ and $(C^\infty(E_\bullet'),d_\bullet')$, so that the induced homomorphisms $f_k:C^\infty(E_k)\to C^\infty(E_k')$ obey
\[
  d_k' \circ f_k = f_{k+1} \circ d_k, \quad\text{for all } k \in \ZZ,
\]
by the usual definitions of morphism or chain map in Weibel \cite[p. 3 or pp. 16--17]{Weibel_introduction_homological_algebra}, and so the following diagrams commute:
\[
  \begin{tikzcd}
  C^\infty(E_{k-1}) \arrow[r, "d_{k-1}"] \arrow[d, "f_{k-1}"]
  &C^\infty(E_k) \arrow[r, "d_k"] \arrow[d, "f_k"]
  &C^\infty(E_{k+1}) \arrow[d, "f_{k+1}"]
  \\
  C^\infty(E_{k-1}') \arrow[r, "d_{k-1}'"]
  &C^\infty(E_k')  \arrow[r, "d_k'"]
  &C^\infty(E_{k+1}')
  \end{tikzcd}
\]
A cochain map $f = \{f_k\}_{k\in\ZZ}$ induces homomorphisms $f_{k,*}:H^k(E)\to H^k(E')$ of the cohomology groups \cite[Exercise 1.1.2]{Weibel_introduction_homological_algebra}. The morphism of complexes is called a \emph{quasi-isomorphism} if the homomorphisms $f_{k,*}:H^k(E)\to H^k(E')$ are isomorphisms for all $k\in\ZZ$ \cite[Definition 1.1.2]{Weibel_introduction_homological_algebra}. If $f$ and $g$ are cochain homotopic, then they induce the same homomorphisms, $f_{k,*} = g_{k,*}:H^k(E)\to H^k(E')$, of cohomology groups.

Two cochain maps $f,g$ are \emph{cochain homotopic}, $f\sim g$, if there exist homomorphisms $s_k:E_k\to E_{k-1}'$ of vector bundles so that the induced homomorphisms $s_k:C^\infty(E_k)\to C^\infty(E_{k-1}')$ obey
\[
  f_k - g_k = d_{k-1}' \circ s_k + s_{k+1} \circ d_k, \quad\text{for all } k \in \ZZ,
\]
and the collection of maps $s =  \{s_k\}_{k\in\ZZ}$ is called a \emph{cochain homotopy} from $f$ to $g$; a cochain map $f$ is a \emph{cochain homotopy equivalence} if there is a cochain map $g$ defined by homomorphisms $g_k:E_k'\to E_k$ of vector bundles such that $g\circ f$ and $f\circ g$ are cochain homotopic to the identity maps of $(C^\infty(E_\bullet),d_\bullet)$ and $(C^\infty(E_\bullet'),d_\bullet')$, respectively \cite[Definition 1.4.4]{Weibel_introduction_homological_algebra}. The induced morphisms, $f_*$ and $g_*$, of complexes are then inverses of one another and so are quasi-isomorphisms.

\section{Elliptic deformation complex for the non-Abelian monopole equations over a Riemannian four-manifold}
\label{sec:Elliptic_deformation_complex_for_SO3_monopole_equations}
In this section, we review the elliptic deformation complex for the non-Abelian monopole equations over a closed, four-dimensional, oriented, smooth Riemannian manifold $(X,g)$. Let $E$ be a smooth,
Hermitian vector bundle over $X$ and recall from Feehan and Leness  \cite[Equation (2.37), p. 312]{FL1} or \cite[Equation (2.47), p. 75]{FL2a} that the elliptic deformation complex for the moduli space of non-Abelian monopoles modulo the group $C^\infty(\SU(E))$ of determinant one, smooth, unitary gauge transformations at a $C^\infty$ pair $(A,\Phi)$ obeying the unperturbed non-Abelian monopole equations \eqref{eq:SO(3)_monopole_equations} is given by
\begin{equation}
\label{eq:SO3MonopoleDefComplex}
\begin{CD}
0
@>>>
\Omega^0(\su(E))
@> d_{A,\Phi}^0>>
\begin{matrix}
\Omega^1(\su(E))
\\
\oplus
\\
\Omega^0(W^+\oplus E)
\end{matrix}
@> d_{A,\Phi}^1 >>
\begin{matrix}
\Omega^+(\su(E))
\\
\oplus
\\
\Omega^0(W^-\otimes E)
\end{matrix}
@>>>
0
\end{CD}
\end{equation}
where $d_{A,\Phi}^0$ in \eqref{eq:d_APhi^0} is the differential at the identity $\id_E\in C^\infty(\SU(E))$ of the map
\[
  C^\infty(\SU(E)) \ni u \mapsto u^*(A,\Phi) \in \sA(E,h) \times \Omega^0(W^+\otimes E),
\]
where $u^*(A,\Phi)$ is defined in \eqref{eq:GaugeActionOnSpinuPairs}.  Thus,
\[
  d_{A,\Phi}^0 \xi= (d_A\xi, -\xi\Phi), \quad\text{for all } \xi \in \Omega^0(\su(E)).
\]
The differential $d_{A,\Phi}^1$ is the differential at a pair $(A,\Phi)$ of the map
\[
  \sA(E,h) \times \Omega^0(W^+\otimes E) \ni (A,\Phi) \mapsto \sS(A,\Phi) \in \Omega^+(\su(E))\oplus \Omega^0(W^-\otimes E)
\]
defined by the unperturbed non-Abelian monopole equations \eqref{eq:SO(3)_monopole_equations} (compare \cite[Equation (2.36), p. 312]{FL1} or \cite[Equation (2.47), p. 75]{FL2a} for the perturbed non-Abelian monopole equations), so that
\begin{equation}
\label{eq:d1OfSO3MonopoleComplex}
d_{A,\Phi}^1(a,\phi)
:=
\begin{pmatrix}
d^+_Aa -(\Phi\otimes\phi^*+\phi\otimes\Phi^*)_{00}
\\
D_A\phi+\rho(a)\Phi
\end{pmatrix},
\quad\text{for all } (a,\phi) \in \Omega^1(\su(E))\oplus \Omega^0(W^+\otimes E).
\end{equation}
Because the non-Abelian monopole equations \eqref{eq:SO(3)_monopole_equations} are invariant under the action of $C^\infty(\SU(E))$, if $(A,\Phi)$ is a non-Abelian monopole, then
we have
\begin{equation}
  \label{eq:d_squared_zero_SO3_monopole_complex}
  d_{A,\Phi}^1 \circ d_{A,\Phi}^0 = 0 \quad\text{on } \Omega^0(\su(E)),
\end{equation}
and so the sequence \eqref{eq:SO3MonopoleDefComplex} is a complex. Moreover, the complex \eqref{eq:SO3MonopoleDefComplex} is elliptic, since the operators $d_A^+ + d_A^*$ and $D_A$ are elliptic.

The cohomology groups of the complex \eqref{eq:SO3MonopoleDefComplex} are defined by the real vector spaces
\[
  H_{A,\Phi}^k := \left.\Ker d_{A,\Phi}^k\right/\Ran d_{A,\Phi}^{k-1},
  \quad\text{for } k=0,1,2,
\]
while their harmonic representative spaces are defined by the real vector spaces
\begin{subequations}
\label{eq:H_APhi^bullet}  
\begin{align}
  \label{eq:H_APhi^0}
  \bH_{A,\Phi}^0 &:= \Ker\left(d_{A,\Phi}^0:\Omega^0(\su(E)) \to \Omega^1(\su(E))\oplus \Omega^0(W^+\otimes E)\right),
  \\
  \label{eq:H_APhi^1}
  \bH_{A,\Phi}^1 &:= \Ker\left(d_{A,\Phi}^1 + d_{A,\Phi}^{0,*}:\Omega^1(\su(E))\oplus \Omega^0(W^+\otimes E)\right.
  \\
  \notag
  &\qquad \left. \to \Omega^+(\su(E))\oplus \Omega^0(W^-\otimes E)\oplus \Omega^0(\su(E))\right),
  \\
  \label{eq:H_APhi^2}
  \bH_{A,\Phi}^2 &:=\Ker \left(d_{A,\Phi}^{1,*}:\Omega^+(\su(E))\oplus \Omega^0(W^-\otimes E) \to \Omega^1(\su(E))\oplus \Omega^0(W^+\otimes E)\right).
\end{align}
\end{subequations}
We note the following regularity result for solutions to the non-Abelian monopole equations.

\begin{lem}[Regularity for $W^{1,p}$ non-Abelian monopole equations]
\label{lem:Regularity_For_Non-Abelian_Monopoles}
Let $\ft=(W\otimes E,\rho)$ be a spin${}^u$ structure over a closed, connected, four-dimensional, oriented, smooth Riemannian manifold $(X,g)$, where $(E,h)$ is a smooth, Hermitian vector bundle.  If $(A,\Phi) \in \sA(E,h)\times W^{1,p}(W^+\otimes E)$ is a solution to the (perturbed) non-Abelian monopole equations \eqref{eq:PerturbedSO3MonopoleEquations} and $p\in(2,\infty)$, then there is gauge transformation $u\in W^{2,p}(\SU(E))$ such that $(u^*A,u^{-1}\Phi)$ is a smooth solution to the non-Abelian monopole equations.
\end{lem}

\begin{proof}
We recall that if $(A,\Phi)$ is a $W^{2,2}$ solution to the non-Abelian monopole equations, then there is a $W^{2,3}$ unitary gauge transformation $u$ such that $u^*(A,\Phi)$ is smooth by Feehan and Leness \cite[Proposition 3.3, p. 320 and Proposition 3.7, p. 323]{FL1}. The argument for a $W^{1,p}$ non-Abelian monopole and existence of a $W^{2,p}$ determinant one, unitary gauge transformation $u$ is very similar. 
\end{proof}

We will use the following lemma in Chapter \ref{chap:VirtualMorseIndexComputation} to relate the index of elliptic complexes with the virtual Morse--Bott index of a Hamiltonian function. Lemma \ref{lem:H0_Of_NonAbelianMonopoleComplex_Vanishes} is a relatively rare example of a result whose proof relies on the complex vector bundle $E$ having rank two.

\begin{lem}[Condition for vanishing of $\bH_{A,\Phi}^0$]
\label{lem:H0_Of_NonAbelianMonopoleComplex_Vanishes}
Let $\ft=(W\otimes E,\rho)$ be a spin${}^u$ structure over a closed, four-dimensional, oriented, smooth Riemannian manifold $(X,g)$, where $E$ is a smooth, rank-two Hermitian vector bundle.  If $(A,\Phi) \in \sA(E,h)\times W^{1,p}(W^+\otimes E)$ is a solution to the (perturbed) non-Abelian monopole equations \eqref{eq:PerturbedSO3MonopoleEquations} with $\Phi\not\equiv 0$, then
\[
\bH_{A,\Phi}^0=(0),
\]
where the harmonic space $\bH_{A,\Phi}^0$ is defined in \eqref{eq:H_APhi^0}.
\end{lem}

\begin{proof}
By Lemma \ref{lem:Regularity_For_Non-Abelian_Monopoles}, there is a gauge transformation $u \in W^{2,p}(\SU(E))$ such that $(A',\Phi') := (u^*A,u^{-1}\Phi)$ is a smooth solution to the non-Abelian monopole equations \eqref{eq:PerturbedSO3MonopoleEquations} and $\Phi'\not\equiv 0$ since $\Phi' = u^{-1}\Phi$. By its definition in \eqref{eq:H_APhi^0}, we have
\[
  \bH_{A,\Phi}^0 = \Ker d_{A,\Phi}^0 \cap W^{2,p}(\su(E)).
\]
By writing $\xi' = u^*\xi = u^{-1}\xi u$ for $\xi \in W^{2,p}(\su(E))$, we see that the definition of $d_{A,\Phi}^0$ in \eqref{eq:d_APhi^0} and the pullback action of $u$ on $\sA(E,h)\times W^{1,p}(W^+\otimes E)$ in \eqref{eq:Donaldson-Kronheimer_2-1-7_pullback} and \eqref{eq:PullbackActionOnUnitaryPairs} gives
\begin{multline*}
  d_{A',\Phi'}^0\xi'
  = (d_{A'}\xi', -\xi'\Phi')
  = (d_{u^*(A)}u^*\xi, -(u^*\xi)u^{-1}\Phi)
  \\
  = (u^{-1}(d_A)(uu^{-1}\xi u), -(u^{-1}\xi u)u^{-1}\Phi
  = (u^{-1}(d_A)(\xi u), -u^{-1}\xi\Phi)
  = (u^*(d_A\xi),-u^{-1}\xi\Phi)
  \\
  = u^*(d_{A,\Phi}^0\xi)
  \in W^{1,p}(T^*X\otimes\su(E)) \oplus W^{1,p}(W^+\otimes E).
\end{multline*}
Thus, as expected, $d_{A',\Phi'}^0\xi' = 0 \iff d_{A,\Phi}^0\xi = 0$ and $\bH_{A',\Phi'}^0 \cong \bH_{A,\Phi}^0$, so we may assume without loss of generality that $(A,\Phi)$ is a smooth pair.

We now argue by contradiction and suppose that $\bH_{A,\Phi}^0 \neq (0)$, so there exists $\xi\in W^{2,p}(\su(E))$ which is not identically zero such that $d_{A,\Phi}^0\xi = (d_A\xi, -\xi\Phi) = 0$. By Proposition \ref{prop:ReducibleUnitaryFromH0}, the existence of $\xi \in \Ker d_A \cap W^{2,p}(\su(E))$ with $\xi\not\equiv 0$ implies that $E = L_1\oplus L_2$ as an orthogonal direct sum of smooth Hermitian line bundles corresponding to the constant, distinct eigenvalues $i\mu_1,i\mu_2$ of $\xi$. Since $\mu_1$ or $\mu_2$ must be non-zero and $\tr\xi = i\mu_1 + i\mu_2 = 0$ (because $\xi$ is a section of $\su(E)$), we must have $\mu_1 = \mu$ and $\mu_2 = -\mu$ for a non-zero, real constant $\mu$. From \eqref{eq:Spectral_theorem}, we obtain
\[
  \xi = i\mu\pi_1 - i\mu\pi_2,
\]
where the $\pi_k\in C^\infty(\gl(E))$ are orthogonal projections from $E$ to $L_k$ for $k=1,2$. Because $\det\xi = -\mu^2 \neq 0$, we have $\xi\in C^\infty(\GL(E))$ since $\xi$ is pointwise invertible. Therefore, $\xi\Phi \equiv 0$ implies that $\Phi\equiv 0$, contradicting our assumption that $\Phi\not\equiv 0$. Hence, if $\Phi\not\equiv 0$, then $\Ker d_{A,\Phi}^0 = (0)$ and this completes the proof.
\end{proof}

We record the following identification of the Euler characteristic of the elliptic deformation complex for non-Abelian monopoles \eqref{eq:SO3MonopoleDefComplex} with the expected dimension of the moduli space $\sM_\ft$.

\begin{prop}[Expected dimension of the moduli space of non-Abelian monopoles in terms of Euler characteristic  of the elliptic deformation complex]
\label{prop:IdentifyExpDimOfnonAbelianMonopoleModuliSpaceWithEulerCharacteristic}
(See Feehan and Leness \cite[Section 2.6, p. 310]{FL1}.) Let $\ft$ be a \spinu structure over a closed, oriented, smooth, Riemannian four-manifold $(X,g)$. If $(A,\Phi)$ is a non-Abelian monopole on $\ft$, then the negative of Euler characteristic of the elliptic deformation complex for non-Abelian monopoles \eqref{eq:SO3MonopoleDefComplex} is equal to the expected dimension of $\sM_\ft$,
\begin{equation}
\label{eq:IdentifyExpDimOfnonAbelianMonopoleModuliSpaceWithEulerCharacteristic}
\exp\,\dim\sM_\ft
=
-\sum_{i=0}^2 (-1)^k\dim_\RR \bH_{A,\Phi}^k,
\end{equation}
where $\bH_{A,\Phi}^k$ is defined in \eqref{eq:H_APhi^bullet} for $k=0,1,2$.
\end{prop}

This completes our discussion of the elliptic deformation complex \eqref{eq:SO3MonopoleDefComplex} and the harmonic space representatives for its cohomology groups.

\section{Elliptic complex for the  projective vortex equations over an almost Hermitian manifold}
\label{sec:Elliptic_deformation_complex_for_projective_vortex_equations}
Let $E$ be a smooth Hermitian vector bundle over an almost Hermitian manifold $(X,g,J)$ with fundamental two-form $\omega = g(\cdot,J\cdot)$ as in \eqref{eq:Fundamental_two-form} and fixed smooth, unitary connection on the complex line bundle $\det E$. As noted in Section \ref{sec:SO3_monopole_equations_over_almost_Hermitian_four-manifolds}, the projective vortex equations \eqref{eq:SO(3)_monopole_equations_almost_Hermitian_alpha} are defined for a bundle $E$ of any complex rank over an almost Hermitian manifold $X$ of any even real dimension.  When $X$ has complex dimension $n\geq 3$, we assume in addition that $(X,J)$ is a complex manifold\footnote{L\"ubke and Teleman \cite[Section 6.2]{Lubke_Teleman_2006} require $(X,g,J)$ to be a complex Gauduchon manifold.}, so the higher-order differentials in the forthcoming sequence \eqref{eq:Projective_vortex_elliptic_deformation_complex} form a complex.

The projective vortex equations \eqref{eq:SO(3)_monopole_equations_almost_Hermitian_alpha} are invariant under the action of the group $C^\infty(\SU(E))$ of determinant one, unitary gauge transformations (but not invariant under the larger group of determinant-one, not necessarily unitary gauge transformations
$C^\infty(\SL(E))$ defined following the forthcoming \eqref{eq:AutomorphismBundles}) and so the sequence of partial differential operators at a $C^\infty$ pair $(A,\varphi) \in \sA(E,h)\times\Omega^0(E)$ for the moduli space of projective vortices modulo $C^\infty(\SU(E))$ is given by
\begin{multline}
\label{eq:Projective_vortex_elliptic_deformation_complex}
\begin{CD}
0
@>>>
\Omega^0(\su(E))
@> d_{A,\varphi}^0>>
\begin{matrix}
\Omega^1(\su(E))
\\
\oplus
\\
\Omega^0(E)
\end{matrix}
@> d_{A,\varphi}^1 >>
\begin{matrix}
\Omega^0(\su(E))
\\
\oplus
\\
\Omega^{0,2}(\fsl(E))
\\
\oplus
\\
\Omega^{0,1}(E)
\end{matrix}
@> d_{A,\varphi}^2 >>
\begin{matrix}
\Omega^{0,3}(\fsl(E))
\\
\oplus
\\
\Omega^{0,2}(E)
\end{matrix}
\cdots
\end{CD}
\\
\begin{CD}
\cdots
@>d_{A,\varphi}^{n-1} >>
\begin{matrix}
\Omega^{0,n}(\fsl(E))
\\
\oplus
\\
\Omega^{0,n-1}(E)
\end{matrix}
@> d_{A,\varphi}^n >>
\begin{matrix}
0
\\
\oplus
\\
\Omega^{0,n}(E)
\end{matrix}
@>>>
0
\end{CD}
\end{multline}
where the definition of the differential $d_{A,\varphi}^0$ is implied by that of $d_{A,\Phi}^0$ in \eqref{eq:d_APhi^0} as the derivative at the identity $\id_E\in C^\infty(\SU(E))$ of the map
\[
  C^\infty(\SU(E)) \ni u \mapsto u(A,\varphi) \in \sA(E,h) \times \Omega^0(E),
\]
so that
\begin{equation}
\label{eq:d0_projective_vortex_elliptic_deformation_complex}
  d_{A,\varphi}^0\xi = (d_A\xi,-\xi\varphi) \in \Omega^1(\su(E))\oplus\Omega^0(E),
  \quad\text{for all } \xi \in \Omega^0(\su(E)).
\end{equation}
For convenience, we note that if $u_t\in W^{2,p}(\SU(E))$ is a path for $t\in(-\eps,\eps)$ with $u_0=\id_E$ and $du_t/dt|_{t=0} = \xi \in W^{2,p}(\su(E))$, then the argument giving \eqref{eq:Differential_of_SU(E)_GaugeAction} implies that
\begin{equation}
\label{eq:d0_projective_vortex_is_diff_of_gauge_action}
  d_{A,\varphi}^0\xi =
  \left. \frac{d}{dt} u_t^*(A,\varphi)\right|_{t=0},
  \quad\text{for all } \xi \in \Omega^0(\su(E)).
\end{equation}
Also, $d_{A,\varphi}^1$ is the differential at a pair $(A,\varphi)$ of the map $\sS$ in the forthcoming \eqref{eq:Projective_vortex_map}, 
\[
  \sS:\sA(E,h) \times \Omega^0(E) \to \Om^0(\su(E))\oplus \Om^{0,2}(\fsl(E))\oplus \Om^{0,1}(E),
\]
defined by the projective vortex equations \eqref{eq:SO(3)_monopole_equations_almost_Hermitian_alpha}, so that (see the forthcoming expression \eqref{eq:Itoh_1985_2-18_SO3_monopole_complex_Kaehler} for $d_{A,\Phi}^1$ and its derivation in the proof of the forthcoming Proposition \ref{prop:Itoh_1985_proposition_2-4_SO3_monopole_complex_Kaehler}),
\begin{multline}
\label{eq:d1_projective_vortex_elliptic_deformation_complex}
d_{A,\varphi}^1(a,\sigma)
:=
\begin{pmatrix}
\Lambda(\bar\partial_Aa' + \partial_Aa'') - i(\varphi\otimes\sigma^*+\sigma\otimes\varphi^*)_0
\\
\bar\partial_Aa''
\\
\bar\partial_A\sigma+a''\varphi
\end{pmatrix},
\\
\text{for all } (a,\sigma) \in \Omega^1(\su(E))\oplus \Omega^0(E),
\end{multline}
where we write $a=\frac{1}{2}(a'+a'')$, with $a''\in\Omega^{0,1}(\fsl(E))$ and $a'=-(a'')^\dagger \in\Omega^{1,0}(\fsl(E))$. Because the projective vortex equations \eqref{eq:SO(3)_monopole_equations_almost_Hermitian_alpha} are invariant under the action of $C^\infty(\SU(E))$, we have
\begin{equation}
  \label{eq:d_squared_zero_projective_vortex_elliptic_deformation_complex}
  d_{A,\varphi}^1 \circ d_{A,\varphi}^0 = 0 \quad\text{on } \Omega^0(\su(E)),
\end{equation}
when $\sS(A,\varphi) = (0,0,0)$, that is, when $(A,\varphi)$ is a projective vortex.

The higher-order differentials $d_{A,\varphi}^k$ for $k=2,\ldots,n$, in the sequence \eqref{eq:Projective_vortex_elliptic_deformation_complex} are defined by analogy of those of Kim \cite{Kim_1987} and Kobayashi \cite[Section 7.2, p. 226]{Kobayashi_differential_geometry_complex_vector_bundles} for the Hermitian--Einstein equations over a complex, K\"ahler manifold of complex dimension $n$ and, furthermore, coincide with the forthcoming definition \eqref{eq:dkStablePair} of the differentials $\bar\partial_{A,\varphi}^k$ in the forthcoming complex \eqref{eq:Holomorphic_pair_elliptic_complex} for $k=3,\ldots,n$:
\begin{align}
  \label{eq:d2_projective_vortex}
  d_{A,\varphi}^2 &:= \bar\partial_{A,\varphi}^2\circ \pi_0^\perp,
  \\
  \label{eq:dk_projective_vortex}
  d_{A,\varphi}^k &:= \bar\partial_{A,\varphi}^k, \quad\text{for } k = 3, \ldots, n,
\end{align}
where $\pi_0^\perp := \id - \pi_0$ and $\pi_0$ is the $L^2$ orthogonal projection
\begin{equation}
  \label{eq:pi_0}
  \pi_0: \Omega^0(\su(E)) \oplus \Omega^{0,2}(\fsl(E)) \oplus \Omega^{0,1}(E) \to \Omega^0(\su(E))
\end{equation}
and ``$\id$'' is the identity operator on the vector space $\Omega^0(\su(E)) \oplus \Omega^{0,2}(\fsl(E)) \oplus \Omega^{0,1}(E)$. Our calculations in Section \ref{sec:Elliptic_deformation_complex_holomorphic_pair_equations} for the higher-order differentials immediately yield the identities
\begin{align*}
  d_{A,\varphi}^2 \circ d_{A,\varphi}^1 &= 0 \quad\text{on } \Omega^1(\su(E)) \oplus \Omega^0(E),
  \\
  d_{A,\varphi}^3 \circ d_{A,\varphi}^2 &= 0 \quad\text{on } \Omega^0(\su(E)) \oplus \Omega^{0,2}(\fsl(E)) \oplus \Omega^{0,1}(E),
  \\
  d_{A,\varphi}^k \circ d_{A,\varphi}^{k-1} &= 0 \quad\text{on } \Omega^{0,k-1}(\fsl(E)) \oplus \Omega^{0,k-2}(E),
  \quad\text{for } k = 4, \ldots, n,
\end{align*}
and thus \eqref{eq:Projective_vortex_elliptic_deformation_complex} is indeed a complex. The fact that the complex \eqref{eq:Projective_vortex_elliptic_deformation_complex} is also \emph{elliptic} follows from Kobayashi \cite[Lemma 7.2.20, p. 227]{Kobayashi_differential_geometry_complex_vector_bundles}, which is proved by showing that the \eqref{eq:Projective_vortex_elliptic_deformation_complex} is exact at the symbol level and, in particular exact at the vector spaces
\[
  B^1
  :=
  \begin{matrix}
    \Omega^1(\su(E))
    \\
    \oplus
    \\
    \Omega^0(E)
  \end{matrix},
  \quad B^2
  :=
  \begin{matrix}
    \Omega^0(\su(E))
    \\
    \oplus
    \\
    \Omega^{0,2}(\fsl(E))
    \\
    \oplus
    \\
    \Omega^{0,1}(E)
  \end{matrix},
  \quad\text{and}\quad B^3
  :=
  \begin{matrix}
    \Omega^{0,3}(\fsl(E))
    \\
    \oplus
    \\
    \Omega^{0,2}(E)
  \end{matrix},
\]
since symbol exactness is clear at the remaining vector spaces $B^0$ and $B^4,\ldots,B^n$ in the complex.

The cohomology groups of the complex \eqref{eq:Projective_vortex_elliptic_deformation_complex} are defined by the real vector spaces
\[
  H_{A,\varphi}^k := \left.\Ker d_{A,\varphi}^k\right/\Ran d_{A,\varphi}^{k-1},
  \quad\text{for } k=0,1,\ldots,n,
\]
while their harmonic representative spaces are defined by the real vector spaces
\begin{subequations}
\label{eq:H_Avarphi^bullet}  
\begin{align}
  \label{eq:H_Avarphi^0}
  \bH_{A,\varphi}^0 &:= \Ker d_{A,\varphi}^0\cap \Omega^0(\su(E)),
  \\
  \label{eq:H_Avarphi^1}
  \bH_{A,\varphi}^1 &:= \Ker\left(d_{A,\varphi}^1 + d_{A,\varphi}^{0,*}\right)
                      \cap\left(\Omega^1(\su(E))\oplus \Omega^0(E)\right),
  \\
  \label{eq:H_Avarphi^2}
  \bH_{A,\varphi}^2 &:=\Ker \left(d_{A,\varphi}^2 + d_{A,\varphi}^{1,*}\right)
                      \cap\left(\Omega^0(\su(E))\oplus \Omega^{0,2}(\fsl(E))\oplus \Omega^{0,1}(E)\right),
  \\
  \label{eq:H_Avarphi^k}
  \bH_{A,\varphi}^k &:=\Ker \left(d_{A,\varphi}^k + d_{A,\varphi}^{k-1,*}\right)
                      \cap\left( \Omega^{0,k}(\fsl(E))\oplus \Omega^{0,k-1}(E) \right),
                      \quad\text{for } k = 3, \ldots,n+1,
\end{align}
\end{subequations}
with the convention that $\Omega^{0,k}(\fsl(E)) = (0)$ and $\Omega^{0,k}(E) = (0)$ when $k<0$ or $k>n$.

We note the following regularity result for solutions to the projective vortex equations.

\begin{lem}[Regularity for $W^{1,p}$ projective vortices]
\label{lem:Regularity_For_ProjVortices}
Let $(E,h)$ be a smooth, Hermitian vector bundle over a closed, connected, smooth almost Hermitian, real $2n$-dimensional manifold $(X,g,J)$.  If $(A,\varphi)\in\sA(E,h)\times W^{1,p}(E)$ is a solution to the projective vortex equations \eqref{eq:SO(3)_monopole_equations_almost_Hermitian_alpha} and $p\in(n,\infty)$, then there is gauge transformation $u\in W^{2,p}(\SU(E))$ such that $(u^*A,u^{-1}\varphi)$ is a smooth solution to the projective vortex equations \eqref{eq:SO(3)_monopole_equations_almost_Hermitian_alpha}.
\end{lem}

\begin{proof}
The proof of Lemma \ref{lem:Regularity_For_Non-Abelian_Monopoles} applies \mutatis to give the corresponding conclusion for solutions to the projective vortex equations.
\end{proof}

We record the following analogue of Lemma \ref{lem:H0_Of_NonAbelianMonopoleComplex_Vanishes} for projective vortices.

\begin{lem}[Condition for vanishing of $\bH_{A,\varphi}^0$]
\label{lem:H0_Of_ProjectiveComplex_Vanishes}
Let $(E,h)$ be a smooth, rank-two Hermitian vector bundle over a closed, almost Hermitian, real $2n$-dimensional manifold $(X,g,J)$.  If $(A,\varphi)\in\sA(E,h)\times W^{1,p}(E)$ is a solution to the projective vortex equations \eqref{eq:SO(3)_monopole_equations_almost_Hermitian_alpha} on $E$ with $p\in(n,\infty)$ and $\varphi\not\equiv 0$, then
\[
\bH_{A,\varphi}^0=(0),
\]
where $\bH_{A,\varphi}^0$ is the harmonic space defined in \eqref{eq:H_Avarphi^0}.
\end{lem}

\begin{proof}
The proof is identical to that of Lemma \ref{lem:H0_Of_NonAbelianMonopoleComplex_Vanishes}, except that we apply the regularity result given by Lemma \ref{lem:Regularity_For_ProjVortices} for projective vortices instead of the regularity result given by Lemma \ref{lem:Regularity_For_Non-Abelian_Monopoles} for non-Abelian monopoles used in the proof of Lemma \ref{lem:H0_Of_NonAbelianMonopoleComplex_Vanishes}.
\end{proof}

\section[Elliptic deformation complex for pre-holomorphic pairs over K\"ahler surfaces]{Elliptic deformation complex for the  pre-holomorphic pair equations over a complex K\"ahler surface}
\label{sec:Elliptic_deformation_complex_for_pre-holomorphic_pair_equations_complex_surface}
In discussing the elliptic deformation complex for the moduli space of \emph{pre-holomorphic pairs} on a smooth Hermitian vector bundle $E$ over a complex Hermitian surface, we will need to define complex gauge transformations. For this purpose, let
\begin{subequations}
\label{eq:AutomorphismBundles}
\begin{align}
  \label{eq:AutomorphismBundles_GL(E)}
  \GL(E)& := \{u\in \End(E): u(x)\in\GL(E_x), \text{ for all } x\in X\},
  \\
  \label{eq:AutomorphismBundles_SL(E)}
  \SL(E)& := \{u\in \End(E): u(x)\in\SL(E_x), \text{ for all } x\in X\},
\end{align}
\end{subequations}
denote the smooth principal fiber bundles over $X$ with structure groups $\GL(r,\CC)$ and $\SL(r,\CC)$, respectively, when $E$ has complex rank $r$. Complex gauge transformations are sections of these bundles.

The elliptic deformation complex for the moduli space of \emph{pre-holomorphic pairs} on a smooth Hermitian vector bundle $E$ over a complex Hermitian surface $(X,g,J)$ modulo the action of the group $C^\infty(\SL(E))$ of smooth complex, determinant-one gauge transformations of $E$ at a $C^\infty$ pair $(\bar\partial_A,(\varphi,\psi))$ obeying the pre-holomorphic pair equations \eqref{eq:SO(3)_monopole_equations_(0,2)_curvature}
and \eqref{eq:SO(3)_monopole_equations_Dirac_almost_Kaehler} is given by
\begin{equation}
\label{eq:Pre-holomorphic_pair_elliptic_complex}
\begin{CD}
0
@>>>
\Omega^0(\fsl(E))
@> \bar\partial_{A,(\varphi,\psi)}^0 >>
\begin{matrix}
\Omega^{0,1}(\fsl(E))
\\
\oplus
\\
\Omega^0(E)\oplus \Omega^{0,2}(E)
\end{matrix}
@> \bar\partial_{A,(\varphi,\psi)}^1 >>
\begin{matrix}
\Omega^{0,2}(\fsl(E))
\\
\oplus
\\
\Omega^{0,1}(E)
\end{matrix}
@>>>
0
\end{CD}
\end{equation}
where $\varphi\in\Omega^0(E)$ and $\psi\in\Omega^{0,2}(E)$. Here, $\bar\partial_{A,(\varphi,\psi)}^0$ is the differential at the identity $\id_E\in\SL(E)$ of the map
\[
  C^\infty(\SL(E)) \ni u \mapsto u(\bar\partial_A,(\varphi,\psi)) \in \sA^{0,1}(E) \times \Omega^0(E)\oplus\Omega^{0,2}(E),
\]
for a given pair $(\bar\partial_A,(\varphi,\psi)$, so that
\begin{equation}
\label{eq:d0StableComplex}
\bar\partial_{A,(\varphi,\psi)}^0\zeta
:=
\begin{pmatrix}
  \bar\partial_A\zeta
  \\
-\zeta(\varphi,\psi)
\end{pmatrix},
\quad\text{for all } \zeta \in \Omega^0(\fsl(E)).
\end{equation}
We consider the affine space of $(0,1)$-connections (see the forthcoming \eqref{eq:Affine_space_01-connections})
\[
  \sA^{0,1}(E) = \bar\partial_{A_0} + \Omega^{0,1}(\fsl(E)),
\]
where $\bar\partial_{A_0}:\Omega^0(E)\to\Omega^{0,1}(E)$ is a fixed, smooth $(0,1)$-connection as in
the forthcoming \eqref{eq:Donaldson_Kronheimer_2-1-48}. Moreover, $\bar\partial_{A,(\varphi,\psi)}^1$ is the differential at the pair $(\bar\partial_A,(\varphi,\psi))$ of the map
\begin{equation}
\label{eq:PreHolomorphicPairMap}
  \sA^{0,1}(E) \times \Omega^0(E)\oplus\Omega^{0,2}(E) \ni (\bar\partial_A,(\varphi,\psi)) \mapsto \fS(\bar\partial_A,(\varphi,\psi)) \in \Omega^{0,2}(\fsl(E)) \oplus \Omega^{0,1}(E)
\end{equation}
defined by the pre-holomorphic pair equations \eqref{eq:SO(3)_monopole_equations_(0,2)_curvature}
and \eqref{eq:SO(3)_monopole_equations_Dirac_almost_Kaehler}. By substituting $a = (a'+a'')/2$ as in \eqref{eq:Decompose_a_in_Omega1suE_into_10_and_01_components} and using $DF_A$ to denote the derivative of the curvature map, $\sA(E,h) \ni A \mapsto F_A \in \Omega^2(\su(E))$, at a unitary connection $A$ (the Chern connection defined by the $(0,1)$-connection $\bar\partial_A$ and Hermitian metric $h$ on $E$ --- see Remark \ref{rmk:Chern_connection}), we obtain
\[
  (DF_A)a = d_Aa \quad\text{and}\quad (DF_A^{0,2})a = \left((DF_A)a\right)^{0,2} = (d_Aa)^{0,2} = \frac{1}{2}\bar\partial_Aa''.
\]
By applying an overall rescaling of equation \eqref{eq:SO(3)_monopole_equations_(0,2)_curvature} by a factor of $1/2$ (before or after computing the linearization), we thus see that
\begin{multline}
\label{eq:d1StableComplex}
\bar\partial_{A,(\varphi,\psi)}^1(a'',\sigma,\tau)
:=
\begin{pmatrix}
\bar\partial_Aa'' - (\tau\otimes\varphi^* + \psi\otimes\sigma^*)_0
\\
\bar\partial_A\sigma+\bar\partial_A^*\tau +a''\varphi-\star(a'\wedge \star\psi)
\end{pmatrix},
\\
\text{for all } (a'',\sigma,\tau) \in \Omega^{0,1}(\fsl(E))\oplus \Omega^0(E)\oplus \Omega^{0,2}(E),
\end{multline}
where $a' := -(a'')^\dagger \in \Omega^{1,0}(\fsl(E))$ by \eqref{eq:Kobayashi_7-6-11}. Because the pre-holomorphic pair equations \eqref{eq:SO(3)_monopole_equations_(0,2)_curvature}
and \eqref{eq:SO(3)_monopole_equations_Dirac_almost_Kaehler} are
invariant under the action of $\SL(E)$, we have
\begin{equation}
  \label{eq:d_squared_zero_pre-holomorphic_pair_complex}
  \bar\partial_{A,(\varphi,\psi)}^1 \circ \bar\partial_{A,(\varphi,\psi)}^0 = 0 \quad\text{on } \Omega^0(\fsl(E)),
\end{equation}
when $(\bar\partial_A,(\varphi,\psi))$ is a pre-holomorphic pair in the sense that $(\bar\partial_A,(\varphi,\psi))$ is in the zero locus of the map \eqref{eq:PreHolomorphicPairMap}.

Hence, the sequence \eqref{eq:Pre-holomorphic_pair_elliptic_complex} is a complex and is elliptic since the operators $\bar\partial_A+\bar\partial_A^*$ are elliptic (see Wells \cite[Chapter IV, Example 5.7, p. 151]{Wells3}).

The cohomology groups of the complex \eqref{eq:Pre-holomorphic_pair_elliptic_complex} are defined by the real vector spaces
\[
  H_{\bar\partial_{A,(\varphi,\psi)}}^k := \left.\Ker \bar\partial_{A,(\varphi,\psi)}^k\right/\Ran \bar\partial_{A,(\varphi,\psi)}^{k-1},
  \quad\text{for } k=0,1,2,
\]
while their harmonic representative spaces are defined by the real vector spaces
\begin{subequations}
  \label{eq:H_dbar_APhi^0bullet}
\begin{align}
  \label{eq:H_dbar_APhi^00}
  \bH_{\bar\partial_{A,(\varphi,\psi)}}^0 &:= \Ker\left(\bar\partial_{A,(\varphi,\psi)}^0:\Omega^0(\fsl(E)) \to \Omega^{0,1}(\fsl(E))\oplus \Omega^0(E)\oplus \Omega^{0,2}(E)\right),
  \\
  \label{eq:H_dbar_APhi^01}
  \bH_{\bar\partial_{A,(\varphi,\psi)}}^1 &:= \Ker\left(\bar\partial_{A,(\varphi,\psi)}^1 + \bar\partial_{A,(\varphi,\psi)}^{0,*}:\Omega^{0,1}(\fsl(E))\oplus \Omega^0(E)\oplus \Omega^{0,2}(E)\right.
  \\
  \notag
  &\qquad \left. \to \Omega^{0,2}(\fsl(E))\oplus \Omega^{0,1}(E) \oplus \Omega^0(\fsl(E))\right),
  \\
  \label{eq:H_dbar_APhi^02}
  \bH_{\bar\partial_{A,(\varphi,\psi)}}^2 &:=\Ker \left(\bar\partial_{A,(\varphi,\psi)}^{1,*}:\Omega^{0,2}(\fsl(E))\oplus \Omega^{0,1}(E) \to \Omega^{0,1}(\fsl(E))\oplus \Omega^0(E)\oplus \Omega^{0,2}(E) \right).
\end{align}
\end{subequations}
From \eqref{eq:d0StableComplex} and \eqref{eq:H_dbar_APhi^00}, we see that
\begin{equation}
  \label{eq:H_dbar_APhi^00_explicit}
  \bH_{\bar\partial_{A,(\varphi,\psi)}}^0 = \left\{\zeta \in \Omega^0(\fsl(E)): \bar\partial_A\zeta - \zeta(\varphi,\psi) = 0 \in \Omega^{0,1}(\fsl(E))\oplus \Omega^0(E)\oplus \Omega^{0,2}(E)\right\}.
\end{equation}
To obtain a more explicit expression for $\bH_{\bar\partial_{A,(\varphi,\psi)}}^1$ in \eqref{eq:H_dbar_APhi^01}, observe that \eqref{eq:d0StableComplex} yields
\begin{align*}
  \left(\bar\partial_{A,(\varphi,\psi)}^{0,*}(a'',\sigma,\tau), \zeta\right)_{L^2(X)}
  &=
  \left((a'',\sigma,\tau), \bar\partial_{A,(\varphi,\psi)}^0\zeta\right)_{L^2(X)}
  \\
  &=
  \left((a'',\sigma,\tau), \bar\partial_A\zeta - \zeta(\varphi,\psi)\right)_{L^2(X)}
  \\
  &=
  (a'',\bar\partial_A\zeta)_{L^2(X)} - (\sigma, \zeta\varphi)_{L^2(X)} - (\tau, \zeta\psi)_{L^2(X)}
  \\
  &=
    (\bar\partial_A^*a'',\zeta)_{L^2(X)} - (\sigma, R_\varphi\zeta)_{L^2(X)} - (\tau, R_\psi\zeta)_{L^2(X)}
    \quad\text{(by \eqref{eq:Right_multiplication_of_section_slE_by_sections_E_or_02E})}
  \\
  &=
  (\bar\partial_A^*a'' - R_\varphi^*\sigma - R_\psi^*\tau, \zeta)_{L^2(X)}
  \quad\text{for all } \zeta \in \Omega^0(\fsl(E)).
\end{align*}
Hence, the preceding identity and expression \eqref{eq:d1StableComplex} for $\bar\partial_{A,(\varphi,\psi)}^1$ imply that $\bH_{\bar\partial_{A,(\varphi,\psi)}}^1$ is given by the vector space of all
\[
  (a'',\sigma,\tau) \in \Omega^{0,1}(\fsl(E))\oplus \Omega^0(E)\oplus \Omega^{0,2}(E)
\]
such that (again using $a' := -(a'')^\dagger \in \Omega^{1,0}(\fsl(E))$ by \eqref{eq:Kobayashi_7-6-11})
\begin{subequations}
  \label{eq:H_dbar_APhi^01_explicit}
  \begin{align}
     \label{eq:H_dbar_APhi^01_explicit_0}
    \bar\partial_A^*a'' - R_\varphi^*\sigma - R_\psi^*\tau
    &= 0 \in \Omega^0(\fsl(E)),
    \\
     \label{eq:H_dbar_APhi^01_explicit_02}
    \bar\partial_Aa'' - (\tau\otimes\varphi^*+\psi\otimes\sigma^*)_0
    &= 0 \in \Omega^{0,2}(\fsl(E)),
    \\
     \label{eq:H_dbar_APhi^01_explicit_01}
    \bar\partial_A\sigma+\bar\partial_A^*\tau +a''\varphi-\star(a'\wedge \star\psi)
    &= 0 \in \Omega^{0,1}(E).
  \end{align}
\end{subequations}
Lastly, to obtain an explicit expression for $\bH_{\bar\partial_{A,(\varphi,\psi)}}^2$ in \eqref{eq:H_dbar_APhi^02}, observe that
\begin{align*}
  \left(\bar\partial_{A,(\varphi,\psi)}^{1,*}(v'',\nu), (a'',\sigma,\tau)\right)_{L^2(X)}
  &=
  \left((v'',\nu), \bar\partial_{A,(\varphi,\psi)}^1(a'',\sigma,\tau)\right)_{L^2(X)},
  \\
  &=
  \left(v'',\bar\partial_Aa'' - (\tau\otimes\varphi^*+\psi\otimes\sigma^*)_0\right)_{L^2(X)}
  \\
  &\quad + \left(\nu, \bar\partial_A\sigma+\bar\partial_A^*\tau +a''\varphi-\star(a'\wedge \star\psi)\right)_{L^2(X)} \quad\text{(by \eqref{eq:d1StableComplex})}
  \\
  &= (\bar\partial_A^*v'',a'')_{L^2(X)} - \left(v'', (\tau\otimes\varphi^*+\psi\otimes\sigma^*)_0\right)_{L^2(X)}
  \\
  &\quad + (\bar\partial_A^*\nu,\sigma)_{L^2(X)} + (\bar\partial_A\nu,\tau)_{L^2(X)} + \left(\nu, a''\varphi-\star(a'\wedge \star\psi)\right)_{L^2(X)}
  \\
  &\qquad\text{for all } (a'',\sigma,\tau) \in \Omega^{0,1}(\fsl(E))\oplus \Omega^0(E)\oplus \Omega^{0,2}(E).
\end{align*}
Hence, the expression \eqref{eq:d1StableComplex} for $\bar\partial_{A,(\varphi,\psi)}^1$ implies that $\bH_{\bar\partial_{A,(\varphi,\psi)}}^2$ is given by the vector space of all
\[
  (v'',\nu) \in \Omega^{0,2}(\fsl(E))\oplus\Omega^{0,1}(E)
\]
such that, for all $(a'',\sigma,\tau) \in \Omega^{0,1}(\fsl(E))\oplus \Omega^0(E)\oplus \Omega^{0,2}(E)$, with $a' := -(a'')^\dagger$ as above,
\begin{equation}
  \label{eq:H_dbar_APhi^01_part_explicit}
  \begin{aligned}
  (\bar\partial_A^*v'',a'')_{L^2(X)} + \left(\nu, a''\varphi-\star(a'\wedge \star\psi)\right)_{L^2(X)}
  &= 0, 
  \\
  (\bar\partial_A^*\nu,\sigma)_{L^2(X)}  - \left(v'', (\psi\otimes\sigma^*)_0\right)_{L^2(X)} 
  &= 0,
  \\
  (\bar\partial_A\nu,\tau)_{L^2(X)} - \left(v'', (\tau\otimes\varphi^*)_0\right)_{L^2(X)}
  &= 0.
\end{aligned}
\end{equation}
According to the forthcoming equations \eqref{eq:_Inner_product_v''_and_tau_otimes_bar_varphi_equals_v''varphi_and_tau} and \eqref{eq:L2_inner_product_nu_and_a''_varphi_and_star_(a'_wedge_star_psi)}, we have the identities
\begin{align*}
  (\nu, a''\varphi)_{L^2(X)}
  &= \left((\nu\otimes\varphi^*)_0, a''\right)_{L^2(X)},
  \\
  \left(\nu, -\star(a'\wedge \star\psi)\right)_{L^2(X)}
  &= \left(\star((\star\psi^*)\wedge\nu)_0), a'\right)_{L^2(X)},
  \\
  \left(v'', (\psi\otimes\sigma^*)_0\right)_{L^2(X)}
  &= (v''\sigma, \psi)_{L^2(X)},
  \\
  \left(v'', (\tau\otimes\varphi^*)_0\right)_{L^2(X)}
  &= (v''\varphi, \tau)_{L^2(X)}.
\end{align*}
When $(A,(\varphi,\psi))$ is a type $1$ solution to the unperturbed non-Abelian monopole equations \eqref{eq:SO(3)_monopole_equations_Kaehler} over a complex K\"ahler surface, so $\psi=0$ and $\bar\partial_A\circ\bar\partial_A = F_A^{0,2}=0$, the expression \eqref{eq:d1StableComplex} for $\bar\partial_{A,(\varphi,\psi)}^1$ simplifies to give
\begin{multline}
\label{eq:d1StableComplex_type1}
\bar\partial_{A,(\varphi,0)}^1(a'',\sigma,\tau)
:=
\begin{pmatrix}
\bar\partial_Aa'' - (\tau\otimes\varphi^*)_0
\\
\bar\partial_A\sigma+\bar\partial_A^*\tau +a''\varphi
\end{pmatrix},
\\
\text{for all } (a'',\sigma,\tau) \in \Omega^{0,1}(\fsl(E))\oplus \Omega^0(E)\oplus \Omega^{0,2}(E).
\end{multline}
Moreover, the defining equations \eqref{eq:H_dbar_APhi^01_part_explicit} for $(v'',\nu) \in \bH_{\bar\partial_{A,(\varphi,0)}}^2$ simplify to the assertion that, for all $(a'',\sigma,\tau) \in \Omega^{0,1}(\fsl(E))\oplus \Omega^0(E)\oplus \Omega^{0,2}(E)$, the pair $(v'',\nu)$ obeys
\begin{align*}
  \left(\bar\partial_A^*v'' + (\nu\otimes\varphi^*)_0, a''\right)_{L^2(X)}
  &= 0, 
  \\
  (\bar\partial_A^*\nu,\sigma)_{L^2(X)} &= 0,
  \\
  (\bar\partial_A\nu - v''\varphi, \tau)_{L^2(X)}
  &= 0.
\end{align*}
Therefore, $(v'',\nu) \in \bH_{\bar\partial_{A,(\varphi,0)}}^2$ if and only if
\begin{subequations}
  \label{eq:H_dbar_APhi^01_part_explicit_type1}
  \begin{align}
  \label{eq:H_dbar_APhi^01_part_explicit_type1_01}  
  \bar\partial_A^*v'' + (\nu\otimes\varphi^*)_0 &= 0 \in \Omega^{0,1}(\fsl(E)), 
    \\
  \label{eq:H_dbar_APhi^01_part_explicit_type1_00_and_02}
  \bar\partial_A\nu + \bar\partial_A^*\nu  - v''\varphi &= 0 \in \Omega^0(E)\oplus \Omega^{0,2}(E).
\end{align}
\end{subequations}
We shall use the preceding observations as check on consistency later in this chapter.

\section{Elliptic complex for the holomorphic pair equations over a complex manifold}
\label{sec:Elliptic_deformation_complex_holomorphic_pair_equations}
Given a smooth complex vector bundle $E$ over a complex manifold $X$, we let the $\CC$-linear operator
\begin{equation}
  \label{eq:Donaldson_Kronheimer_2-1-48}
  \bar\partial_E:\Omega^0(E) \to \Omega^{0,1}(E)
\end{equation}
denote a \emph{partial connection} in the terminology of Donaldson and Kronheimer \cite[Equation (2.1.48), p. 44]{DK} or $(0,1)$-\emph{connection} in the terminology of Friedman and Morgan \cite[Section 4.1.2, Definition 1.5, p. 283]{FrM} and Itoh \cite[Definition 3.1]{Itoh_1985} or \emph{semiconnection} in the terminology of L\"ubke and Teleman \cite[Section 4.3]{Lubke_Teleman_1995}, so it obeys the Leibnitz rule
\begin{equation}
  \label{eq:Donaldson_Kronheimer_2-1-45_i}
  \bar\partial_E(fs) = (\bar\partial f)s + f\bar\partial_Es,
  \quad\text{for all } f \in \Omega^0(X,\CC) \text{ and } s \in \Omega^0(E).
\end{equation}
More generally, the identity \eqref{eq:Donaldson_Kronheimer_2-1-45_i} extends to give (see Huybrechts \cite[Lemma 1.3.6, p. 44, Definition 3.A.1, p. 145, and Section 4.1, p. 170]{Huybrechts_2005})
\begin{equation}
  \label{eq:Donaldson_Kronheimer_2-1-45_i_pq-forms}
  \bar\partial_E(ws) = (\bar\partial_E w)s + (-1)^{p+q} w\bar\partial_Es,
  \quad\text{for all } w \in \Omega^{p,q}(\gl(E)) \text{ and } s \in \Omega^0(E).
\end{equation}
We recall \cite[Equation (2.1.50), p. 44]{DK}, \cite[Section 4.1.2, p. 284]{FrM} that the \emph{$(0,2)$-curvature} of $\bar\partial_E$ is defined by
\begin{equation}
  \label{eq:Donaldson_Kronheimer_2-1-50}
  F_{\bar\partial_E} := \bar\partial_E\circ\bar\partial_E \in \Omega^{0,2}(\gl(E)).
\end{equation}
The $(0,1)$-connection $\bar\partial_E$ is \emph{integrable} in the sense of \cite[p. 45]{DK}, \cite[Section 4.1.2, p. 284]{FrM}, and thus defines a \emph{holomorphic structure} on $E$, if and only if $F_{\bar\partial_E} \equiv 0$, that is,
\begin{equation}
  \label{eq:dbarE_squared_is_zero}
  \bar\partial_E\circ\bar\partial_E = 0.
\end{equation}
Recall that if $A$ is a smooth connection on $E$, then (see Kobayashi \cite[Equation (1.3.2), p. 8]{Kobayashi_differential_geometry_complex_vector_bundles})
\begin{multline}
\label{eq:Decomposition_FA_bitype}  
  F_A = \partial_A\circ\partial_A + (\partial_A\circ\bar\partial_A + \bar\partial_A\circ\partial_A) + \bar\partial_A\circ\bar\partial_A
  \\
  = F_A^{2,0} + F_A^{1,1} + F_A^{0,2}
  \in \Omega^{2,0}(\gl(E))\oplus \Omega^{1,1}(\gl(E))\oplus\Omega^{0,2}(\gl(E)).
\end{multline}
Note that by \eqref{eq:Donaldson_Kronheimer_2-1-50} and \eqref{eq:Decomposition_FA_bitype} we have
\begin{equation}
  \label{eq:Holomorphic_curvature}
  F_{\bar\partial_A} = \bar\partial_A\circ\bar\partial_A = F_A^{0,2}.
\end{equation}
We recall from Kobayashi \cite[Proposition 1.3.5, p. 8]{Kobayashi_differential_geometry_complex_vector_bundles} that given a \emph{holomorphic vector bundle} $(E,\bar\partial_E)$, there is a connection $A$ on $E$ (not necessarily unique) such that $\bar\partial_A = \bar\partial_E$. Conversely, if $A$ is a smooth connection on $E$ such that $F_A^{0,2} = 0$, then there is a unique holomorphic structure on $E$ such that $\bar\partial_E = \bar\partial_A$.

Let $(\bar\partial_E,\varphi)$ with $\varphi\in\Omega^0(E)$ be a \emph{holomorphic pair} on a Hermitian vector bundle $E$ over a complex manifold $X$, so
\begin{subequations}
  \label{eq:Holomorphic_pair}
  \begin{align}
  \label{eq:F_dbarE_is_zero}  
  F_{\bar\partial_E} &= 0,
    \\
  \label{eq:dbarE_varphi_is_zero}  
  \bar\partial_E\varphi &= 0.
\end{align}
\end{subequations}
To preserve consistency with solutions to the unperturbed non-Abelian monopole equations \eqref{eq:SO(3)_monopole_equations_Kaehler} and the corresponding pre-holomorphic pair equations \eqref{eq:SO(3)_monopole_equations_(0,2)_curvature}
and \eqref{eq:SO(3)_monopole_equations_Dirac_almost_Kaehler} over complex K\"ahler surfaces, we shall assume that the complex line bundle $\det E$ has a \emph{fixed holomorphic structure} $\bar\partial_{\det E}$, which we usually denote by $\bar\partial_{E_d}$. We define the \emph{affine space of (smooth) $(0,1)$-connections} by 
\begin{equation}
  \label{eq:Affine_space_01-connections}
  \sA^{0,1}(E) := \bar\partial_{E_0} + \Omega^{0,1}(\fsl(E)),
\end{equation}
where $\bar\partial_{E_0}:\Omega^0(E)\to\Omega^{0,1}(E)$ is a fixed, smooth $(0,1)$-connection as in
\eqref{eq:Donaldson_Kronheimer_2-1-48} and which induces $\bar\partial_{E_d}$ on $\det E$. We shall also write the affine space of $(0,1)$-connections of class $W^{1,p}$ as
\[
  \sA^{0,1}(E) = \bar\partial_{E_0} + W^{1,p}\left(\Lambda^{0,1}(\fsl(E))\right),
\]
and rely on the context to eliminate ambiguity in its meaning.

The equations \eqref{eq:Holomorphic_pair} are preserved by the action of $\SL(E)$ and so the corresponding sequence when $X$ has complex dimension $n$ is given by
\begin{multline}
\label{eq:Holomorphic_pair_elliptic_complex}
\begin{CD}
0
@>>>
\Omega^0(\fsl(E))
@> \bar\partial_{E,\varphi}^0 >>
\begin{matrix}
\Omega^{0,1}(\fsl(E))
\\
\oplus
\\
\Omega^0(E)
\end{matrix}
@> \bar\partial_{E,\varphi}^1 >>
\begin{matrix}
\Omega^{0,2}(\fsl(E))
\\
\oplus
\\
\Omega^{0,1}(E)
\end{matrix}
@> \bar\partial_{E,\varphi}^2 >>
\cdots
\end{CD}
\\
\begin{CD}
\cdots
@> \bar\partial_{E,\varphi}^{n-1} >>
\begin{matrix}
\Omega^{0,n}(\fsl(E))
\\
\oplus
\\
\Omega^{0,n-1}(E)
\end{matrix}
@> \bar\partial_{E,\varphi}^n >>
\begin{matrix}
0
\\
\oplus
\\
\Omega^{0,n}(E)
\end{matrix}
@>>>
0
\end{CD}
\end{multline}
where
\begin{subequations}
  \label{eq:dStablePair}
\begin{align}
  \label{eq:d0StablePair}
\bar\partial_{E,\varphi}^0\zeta
  &:=
    \begin{pmatrix}
      \bar\partial_E\zeta
      \\
      -\zeta\varphi
    \end{pmatrix},
  \quad
  \text{for all } \zeta \in \Omega^0(\fsl(E)),
\\    
\label{eq:d1StablePair}
\bar\partial_{E,\varphi}^1(\alpha,\sigma)
&:=
\begin{pmatrix}
\bar\partial_E\alpha
\\
\bar\partial_E\sigma + \alpha\varphi
\end{pmatrix},
\quad
\text{for all } (\alpha,\sigma) \in \Omega^{0,1}(\fsl(E))\oplus \Omega^0(E),
\\
\label{eq:dkStablePair}
\bar\partial_{E,\varphi}^k(w,\nu)
&:=
\begin{pmatrix}
  \bar\partial_Ew
  \\
  \bar\partial_E\nu  + (-1)^{k-1}w\varphi
\end{pmatrix},  
\quad
  \text{for all } (w,\nu) \in \Omega^{0,k}(\fsl(E))\oplus \Omega^{0,k-1}(E),
\end{align}
\end{subequations}
where $k=2,\ldots,n+1$ in \eqref{eq:dkStablePair}, with the convention that $\Om^{0,k}(E) = (0)$ when $k>n$. As usual, the differential $\bar\partial_{E,\varphi}^0$ is defined by the linearization of the action of $C^\infty(\SL(E))$ on the affine space of pairs $(\bar\partial_E,\varphi)$, while the differential $\bar\partial_{E,\varphi}^1$ is defined by the linearization of the holomorphic pair equations \eqref{eq:Holomorphic_pair} and we again have
\[
  \bar\partial_{E,\varphi}^1\circ\bar\partial_{E,\varphi}^0 = 0 \quad\text{on } \Omega^0(\fsl(E)).
\]  
The definition of $\bar\partial_{E,\varphi}^k$ in \eqref{eq:dkStablePair} for $k\geq 2$ is motivated by the requirement that
\[
  \bar\partial_{E,\varphi}^k\circ\bar\partial_{E,\varphi}^{k-1} = 0
  \quad\text{on } \Omega^{0,k-1}(\fsl(E)) \oplus \Omega^{k-2}(E), \quad\text{for } k=2,\ldots,n.
\] 
Indeed, for any $k=1,\ldots,n$ and $(w,\nu)\in \Omega^{0,k}(\fsl(E)) \oplus \Omega^{k-1}(E)$ we have
\begin{align*}
  \left(\bar\partial_{E,\varphi}^{k+1}\circ \bar\partial_{E,\varphi}^k\right)(w,\nu)
  &=
    \bar\partial_{E,\varphi}^{k+1}\left(\bar\partial_Ew, \bar\partial_E\nu + (-1)^{k-1}w\varphi\right)
    \quad\text{(by \eqref{eq:dkStablePair})}
  \\
  &= \left((\bar\partial_E\circ\bar\partial_E)w, \bar\partial_E(\bar\partial_E\nu + (-1)^{k-1}w\varphi) + (-1)^k(\bar\partial_Ew)\varphi\right)
    \quad\text{(by \eqref{eq:dkStablePair})}
  \\
  &= \left(0, (\bar\partial_E\circ \bar\partial_E)\nu +(-1)^{k-1} \bar\partial_E(w\varphi)
    +(-1)^k(\bar\partial_Ew)\varphi \right)
  \\
  &= \left(0, (-1)^{k-1}(\bar\partial_Ew)\varphi +(-1)^{2k-1} w\wedge\bar\partial_E\varphi +(-1)^k(\bar\partial_Ew)\varphi \right)
  \\
  &= (0,0),
\end{align*}
where the final equalities follow from the
Leibnitz Rule \eqref{eq:Donaldson_Kronheimer_2-1-45_i_pq-forms} for $\bar\partial_E$ and the facts that
\[
  (\bar\partial_E\circ \bar\partial_E)\nu = 0
  \quad\text{and}\quad
  (\bar\partial_E\circ \bar\partial_E)w = 0
  \quad\text{(by \eqref{eq:dbarE_squared_is_zero})}
  \quad\text{and}\quad
  \bar\partial_E\varphi = 0 \quad\text{(by \eqref{eq:dbarE_varphi_is_zero})}.
\]
These observations ensure that the sequence \eqref{eq:Holomorphic_pair_elliptic_complex} is a complex and it is elliptic since the operators $\bar\partial_E + \bar\partial_E^*$ are elliptic (see Wells \cite[Chapter IV, Example 5.7, p. 151]{Wells3}).

The cohomology groups of the elliptic complex \eqref{eq:Holomorphic_pair_elliptic_complex} are defined by the complex vector spaces
\begin{equation}
  \label{eq:H_dbar_Avarphi^0bullet_cohomology}
  H_{\bar\partial_E,\varphi}^k := \left.\Ker \bar\partial_{E,\varphi}^k\right/\Ran \bar\partial_{E,\varphi}^{k-1},
  \quad\text{for } k=0,1\ldots,n+1,
\end{equation}
while their harmonic representative spaces are defined by the complex vector spaces
\begin{multline}
  \label{eq:H_dbar_Avarphi^0bullet}
  \bH_{\bar\partial_E,\varphi}^k := \Ker\left(\bar\partial_{E,\varphi}^k + \bar\partial_{E,\varphi}^{k-1,*}:\Omega^{0,k}(\fsl(E))\oplus \Omega^{0,k-1}(E)\right.
  \\
  \left. \qquad\to \Omega^{0,k+1}(\fsl(E))\oplus \Omega^{0,k}(E) \oplus \Omega^{0,k-1}(\fsl(E))\oplus \Omega^{0,k-2}(E)\right),
  \\
  \text{for } k=0,1,\ldots,n+1,
\end{multline}
with the convention that $\Omega^{0,k}(\fsl(E)) = (0)$ and $\Omega^{0,k}(E) = (0)$ if $k<0$ or $k>n$.

We now wish to compare the elliptic complexes \eqref{eq:Pre-holomorphic_pair_elliptic_complex} and \eqref{eq:Holomorphic_pair_elliptic_complex} when $n=2$. Observe that, for any $\tau \in \Omega^{0,2}(E)$ and all $(w,\nu) \in \Omega^{0,2}(\fsl(E))\oplus \Omega^{0,1}(E)$, we have
\begin{align*}
  (\bar\partial_{E,\varphi}^{2,*}\tau, (w,\nu))_{L^2(X)}
  &= (\tau, \bar\partial_{E,\varphi}^2(w,\nu))_{L^2(X)}
  \\
  &= (\tau, \bar\partial_E\nu  - w\varphi)_{L^2(X)} \quad\text{(by \eqref{eq:dkStablePair} with $k=2$)}
  \\
  &= (\bar\partial_E^*\tau, \nu)_{L^2(X)}  - (\tau, w\varphi)_{L^2(X)}
  \\
  &= (\bar\partial_E^*\tau, \nu)_{L^2(X)}  - ((\tau\otimes\varphi^*)_0, w)_{L^2(X)},  
\end{align*}
where the last inequality follows from \eqref{eq:_Inner_product_v''_and_tau_otimes_bar_varphi_equals_v''varphi_and_tau}, and thus
\[
  \bar\partial_{E,\varphi}^{2,*}\tau
  =
  \bar\partial_E^*\tau  - (\tau\otimes\varphi^*)_0 \in \Omega^{0,2}(\fsl(E))\oplus \Omega^{0,1}(E),
  \quad\text{for all } \tau \in \Omega^{0,2}(E).
\]  
We may thus partially roll up the complex \eqref{eq:Holomorphic_pair_elliptic_complex} to give the equivalent elliptic complex,
\begin{equation}
\label{eq:Holomorphic_pair_elliptic_complex_rolled_up}
\begin{CD}
0
@>>>
\Omega^0(\fsl(E))
@> \bar\partial_{E,\varphi}^0 >>
\begin{matrix}
\Omega^{0,1}(\fsl(E))
\\
\oplus
\\
\Omega^0(E)\oplus\Omega^{0,2}(E)
\end{matrix}
@> \bar\partial_{E,\varphi}^1 + \bar\partial_{E,\varphi}^{2,*}>>
\begin{matrix}
\Omega^{0,2}(\fsl(E))
\\
\oplus
\\
\Omega^{0,1}(E)
\end{matrix}
@>>>
0
\end{CD}
\end{equation}
where
\begin{multline}
\label{eq:d1StableComplex_type1_rolled_up}
(\bar\partial_{E,\varphi}^1+\bar\partial_{E,\varphi}^{2,*})(a'',\sigma,\tau)
=
\begin{pmatrix}
\bar\partial_Ea'' - (\tau\otimes\varphi^*)_0
\\
\bar\partial_E\sigma + \bar\partial_E^*\tau + a''\varphi
\end{pmatrix},
\\
\text{for all } (a'',\sigma,\tau) \in \Omega^{0,1}(\fsl(E))\oplus \Omega^0(E)\oplus \Omega^{0,2}(E).
\end{multline}
The elliptic complex \eqref{eq:Pre-holomorphic_pair_elliptic_complex} with $\psi=0$ coincides with the elliptic complex \eqref{eq:Holomorphic_pair_elliptic_complex_rolled_up}, as one can see by comparing the definition \eqref{eq:d1StableComplex} of $\bar\partial_{E,(\varphi,0)}^1$ with the expression for $\bar\partial_{E,\varphi}^1+\bar\partial_{E,\varphi}^{2,*}$ in \eqref{eq:d1StableComplex_type1_rolled_up}. Their cohomology groups and Euler characteristics are compared in the forthcoming Corollary \ref{cor:Comparison_elliptic_complexes_cohomology_groups_pre-holomorphic_and_holomorphic_pairs}.

\section{Fredholm maps and sections of sheaves}
\label{sec:Fredholm_maps_elliptic_deformation_complexes}
For versions of the elliptic complex \eqref{eq:Elliptic_complex} in the applications that we discuss in this monograph, the differential $d_1$ is the linearization of a map $F$ at a point, as discussed in Section \ref{sec:Kuranishi_models}; the term $C^\infty(E_0)$ corresponds to the Lie algebra of a Lie group of automorphisms (that leaves $F^{-1}(y_0)$ invariant) and the differential $d_0$ is the linearization of automorphism group at that point; and the term $C^\infty(E_2)$ is the codomain of $d_1$. If $n=2$, then $\Coker d_1$ is finite-dimensional and we call \eqref{eq:Elliptic_complex} an \emph{elliptic deformation complex} for a map $F$ whose linearization $d_1$ is Fredholm upon restriction to $(\Ran d_0)^\perp = \Ker d_0^*$. However, if $n\geq 3$, then $\Coker d_1$ need not be finite-dimensional and, in order to obtain a Fredholm map, we need to effectively replace $C^\infty(E_2)$ by $\Ker d_2$. In this section, we shall describe the details of this process for the holomorphic pair map.

The following map, though complex analytic with respect to $W^{1,p}$ and $L^p$ completions of the domain and codomain, defined by the holomorphic pair equations \eqref{eq:Holomorphic_pair} is not Fredholm,
\begin{equation}
  \label{eq:Holomorphic_pair_map}
  \fS:\sA^{0,1}(E)\times \Omega^0(E)
  \ni (\bar\partial_E,\varphi) \mapsto (F_{\bar\partial_E},\bar\partial_E\varphi) \in
  \Omega^{0,2}(\fsl(E)) \oplus \Omega^{0,1}(E),
\end{equation}
where $F_{\bar\partial_E} = \bar\partial_E\circ\bar\partial_E$ as in
\eqref{eq:Donaldson_Kronheimer_2-1-50},
since neither the kernel nor cokernel are finite-dimensional at pairs $(\bar\partial_E,\varphi)$ as we shall explain below.

To obtain a Fredholm map, we fix a pair $(\bar\partial_E,\varphi)$ and first restrict to a Coulomb-gauge slice (see the forthcoming Theorem \ref{thm:Existence_of_complex_gauge_transformation_to_Coulomb_01_pairs}) for the action of $C^\infty(\SL(E))$,
\begin{multline}
  \label{eq:Holomorphic_pair_map_Coulomb_gauge_slice}
  \fS:(\bar\partial_E,\varphi) +
  \Ker\bar\partial_{E,\varphi}^{0,*}\cap \left( \Omega^{0,1}(\fsl(E)) \oplus \Omega^0(E)\right)
  \ni (\bar\partial_E+\alpha,\varphi+\sigma)
  \\
  \mapsto \left(F_{\bar\partial_E+\alpha},(\bar\partial_E+\alpha)(\varphi+\sigma)\right) \in
  \Omega^{0,2}(\fsl(E)) \oplus \Omega^{0,1}(E),
\end{multline}
defined by the kernel of the operator
\[
  \bar\partial_{E,\varphi}^{0,*}: \Omega^{0,1}(\fsl(E)) \oplus \Omega^{0,1}(E) \to \Omega^0(\fsl(E)),
\]
namely, the $L^2$ adjoint of the operator $\bar\partial_{E,\varphi}^0$ in \eqref{eq:d0StablePair}.

Observe that the linearization of the map $\fS$ in \eqref{eq:Holomorphic_pair_map} is the differential $\bar\partial_{E,\varphi}^1$ given by \eqref{eq:d1StablePair} and part of the elliptic deformation complex \eqref{eq:Holomorphic_pair_elliptic_complex}, namely
\[
  \bar\partial_{E,\varphi}^1: \Omega^{0,1}(\fsl(E))\oplus \Omega^0(E)
  \ni (\alpha,\sigma) \mapsto (\bar\partial_E\alpha,\bar\partial_E\sigma + \alpha\sigma) \in
  \Omega^{0,2}(\fsl(E))\oplus \Omega^{0,1}(E).
\]  
The Hodge decomposition for elliptic complexes \eqref{eq:Hodge_decomposition_all_k-forms} yields
\begin{align*}
  \Omega^{0,1}(\fsl(E))\oplus \Omega^0(E)
  &=
    \bH_{\bar\partial_E,\varphi}^1
    \oplus \Ran \bar\partial_{E,\varphi}^0
    \oplus \Ran \bar\partial_{E,\varphi}^{2,*},
  \\
  \Ker \bar\partial_{E,\varphi}^{0,*}
  &= (\Ran \bar\partial_{E,\varphi}^0)^\perp
    = \bH_{\bar\partial_E,\varphi}^1 \oplus \Ran \bar\partial_{E,\varphi}^{2,*}.
\end{align*}
In particular, the restriction of the preceding operator to the Coulomb-gauge slice, 
\begin{multline}
  \label{eq:dbar_Avarphi_1_Coulomb-gauge_slice}
  \bar\partial_{E,\varphi}^1:
  \Ker \bar\partial_{E,\varphi}^{0,*} \cap\left(\Omega^{0,1}(\fsl(E))\oplus \Omega^0(E)\right)
  \ni (\alpha,\sigma)
  \\
  \mapsto (\bar\partial_E\alpha,\bar\partial_E\sigma + \alpha\sigma) \in
  \Omega^{0,2}(\fsl(E))\oplus \Omega^{0,1}(E),
\end{multline}
has finite-dimensional kernel $\Ker(\bar\partial_{E,\varphi}^1 + \bar\partial_{E,\varphi}^{0,*}) = \bH_{\bar\partial_E,\varphi}^1$.

However, the operator \eqref{eq:dbar_Avarphi_1_Coulomb-gauge_slice} has infinite-dimensional cokernel (even when $X$ has complex dimension two) since the Hodge decomposition for elliptic complexes \eqref{eq:Hodge_decomposition_all_k-forms} gives
\begin{align*}
  \Omega^{0,2}(\fsl(E))\oplus \Omega^1(E)
  &=
  \bH_{\bar\partial_E,\varphi}^2 \oplus \Ran \bar\partial_{E,\varphi}^1
    \oplus \Ran \bar\partial_{E,\varphi}^{2,*},
  \\
  \Ker \bar\partial_{E,\varphi}^2
  &= \bH_{\bar\partial_E,\varphi}^2 \oplus \Ran \bar\partial_{E,\varphi}^1,
\end{align*}
where from \eqref{eq:dkStablePair} we have
$\bar\partial_{E,\varphi}^2 = (\bar\partial_E, \bar\partial_E + \cdot\varphi)$ and
\[
  \Ran \bar\partial_{E,\varphi}^{2,*}
  =
  \Ran\left(\bar\partial_{E,\varphi}^{2,*}:\Omega^{0,3}(\fsl(E))\oplus \Omega^{0,2}(E) \to
\Omega^{0,2}(\fsl(E))\oplus \Omega^{0,1}(E)\right).
\]
One approach to non-Fredholm map problems of this kind is to adapt a construction described by Feehan in \cite[Section 4.4, second proof of Theorem 8]{Feehan_lojasiewicz_inequality_ground_state} for the elliptic deformation complex characterizing the local structure of the moduli space of Yang--Mills connections over a closed, smooth Riemannian manifold of dimension three or higher, extending an equivalent local version of those ideas due to Morgan, Mrowka, and Ruberman \cite[Theorem 12.1.1 and its proof]{MMR} for the moduli space of flat connections over a closed, smooth Riemannian manifold of dimension three. Friedman and Morgan \cite[pp. 455--456]{FrMSWInvariants} describe an approach similar to that of \cite{MMR} for a related non-Fredholm map problem arising in their discussion of the Seiberg--Witten monopole equations over a complex, K\"ahler surface. For another approach, see Section \ref{subsec:Kuranishi_model_holomorphic_structure}.

We recall from Friedman and Morgan \cite[Section 4.1.2, p. 284]{FrM} that $F_{\bar\partial_E}$ obeys the \emph{Bianchi identity},
\begin{equation}
  \label{eq:Holomorphic_curvature_Bianchi_identity}
  \bar\partial_EF_{\bar\partial_E} = 0 \in \Omega^{0,3}(\fsl(E)).
\end{equation}
Consequently, by the definition of $\bar\partial_{E,\varphi}^2$ in \eqref{eq:dkStablePair} (with $k=2$) we see that for any $(0,1)$-pair (not necessarily holomorphic)
\[
  \bar\partial_{E,\varphi}^2(F_{\bar\partial_E},\bar\partial_E\varphi)
  =
  (\bar\partial_EF_{\bar\partial_E},\bar\partial_E\bar\partial_E\varphi - F_{\bar\partial_E}\varphi)
  =
  (0,0),
\]
since $\bar\partial_E\bar\partial_E\varphi = F_{\bar\partial_E}\varphi$ by \eqref{eq:Decomposition_FA_bitype}, and thus we obtain the  \emph{Bianchi identity for $(0,1)$-pairs},
\begin{equation}
  \label{eq:Holomorphic_pair_Bianchi_identity}
  \bar\partial_{E,\varphi}^2(F_{\bar\partial_E},\bar\partial_E\varphi) = (0,0)
  \in \Omega^{0,3}(\fsl(E))  \oplus \Omega^{0,2}(E).
\end{equation}
The linear partial differential operator
\[
  \bar\partial_{E,\varphi}^2:\Omega^{0,2}(\fsl(E))  \oplus \Omega^{0,1}(E)
  \to \Omega^{0,3}(\fsl(E))  \oplus \Omega^{0,2}(E)
\]
may be viewed as a homomorphism of sheaves of (germs of local) sections of smooth vector bundles over $X$ and thus $\Ker\bar\partial_{\bullet,\bullet}^2$ may be viewed as a sheaf over $\sA^{0,1}(E) \times \Omega^0(E)$ (see Grauert and Remmert \cite[Section A.2.2]{Grauert_Remmert_coherent_analytic_sheaves}, Kuranishi \cite{Kuranishi_1964}, and Maurin \cite[Theorem 1, p. 268]{Maurin_riemann_legacy}).

Adapting \cite[Section 4.4, second proof of Theorem 8]{Feehan_lojasiewicz_inequality_ground_state}, we thus reinterpret the holomorphic pair map \eqref{eq:Holomorphic_pair_map} as a $C^\infty(\SL(E))$-equivariant section of a sheaf $\Ker\bar\partial_{\bullet,\bullet}^2$ (of germs of local sections of an infinite-rank vector bundle) over the affine space $\sA^{0,1}(E)\times \Omega^0(E)$,
\begin{multline}
  \label{eq:Holomorphic_pair_section_sheaf}
  \fS:\sA^{0,1}(E)\times \Omega^0(E) \ni
  (\bar\partial_E,\varphi) \mapsto
  \left((\bar\partial_E,\varphi), (F_{\bar\partial_E},\bar\partial_E\varphi)\right)
  \\
  \in \Ker\left(\bar\partial_{E,\varphi}^2:\Omega^{0,2}(\fsl(E)) \oplus \Omega^{0,1}(E)
  \to \Omega^{0,3}(\fsl(E)) \oplus \Omega^{0,2}(E) \right),
\end{multline}
where the fiber of the sheaf $\Ker\bar\partial_{\bullet,\bullet}^2$ over $(\bar\partial_E,\varphi) \in \sA^{0,1}(E) \times \Omega^0(E)$ is given by
\[
  \Ker\left(\bar\partial_{E,\varphi}^2:\Omega^{0,2}(\fsl(E)) \oplus \Omega^{0,1}(E)
  \to \Omega^{0,3}(\fsl(E)) \oplus \Omega^{0,2}(E) \right).
\]
We claim that the section $\fS$ in \eqref{eq:Holomorphic_pair_section_sheaf} becomes Fredholm upon restriction to a Coulomb-gauge slice for the action of $C^\infty(\SL(E))$ on its domain:
\begin{multline}
  \label{eq:Holomorphic_pair_section_sheaf_Coulomb-gauge_slice}
  \fS:(\bar\partial_E,\varphi) + \Ker\bar\partial_{E,\varphi}^{0,*}\cap\left(\Omega^{0,1}(\fsl(E)) \oplus \Omega^0(E)\right) \ni (\bar\partial_E+\alpha,\varphi+\sigma)
  \\
  \mapsto
  \left((\bar\partial_E+\alpha,\varphi+\sigma),
    \left(F_{\bar\partial_E+\alpha},(\bar\partial_E+\alpha)(\varphi+\sigma)\right)\right)
  \\
  \in \Ker\left(\bar\partial_{E,\varphi}^2:\Omega^{0,2}(\fsl(E)) \oplus \Omega^{0,1}(E)
  \to \Omega^{0,3}(\fsl(E)) \oplus \Omega^{0,2}(E) \right),
\end{multline}
where the sheaf $\Ker\bar\partial_{\bullet,\bullet}^2$ is now defined by restriction to the Coulomb-gauge slice,
\begin{equation}
  \label{eq:Ker_dbarAvarphi^2_sheaf_Coulomb-gauge_slice}
  \Ker\bar\partial_{\bullet,\bullet}^2
  := \left\{\Ker \bar\partial_{E+\alpha,\varphi+\sigma}^2: (\alpha,\sigma) \in \Ker\bar\partial_{E,\varphi}^{0,*}\cap\left(\Omega^{0,1}(\fsl(E)) \oplus \Omega^0(E)\right) \right\}.
\end{equation}
For $p\in(n,\infty)$, we thus consider the equation
\begin{multline}
  \label{eq:Holomorphic_pair_Fredholm}
  (F_{\bar\partial_E},\bar\partial_E\varphi) = (0,0) \in
  \Ker\left(\bar\partial_{E,\varphi}^2:L^p\left(\Lambda^{0,k}(\fsl(E)) \oplus\Lambda^{0,1}(E)\right)\right.
  \\
    \to \left. W^{-1,p}\left(\Lambda^{0,3}(\fsl(E)) \oplus \Lambda^{0,2}(E)\right) \right)
\end{multline}
for a pair
\[
  (\bar\partial_E,\varphi) \in (\bar\partial_{E_0},\varphi_0) + W^{1,p}\left(\Lambda^{0,1}(\fsl(E) \oplus E\right).
\]
We have seen that the kernel of the differential of the section $\fS$ in \eqref{eq:Holomorphic_pair_section_sheaf_Coulomb-gauge_slice} is finite-dimensional, so it suffices to consider its cokernel or the cokernel of the differential of the simpler section $\fS$ in \eqref{eq:Holomorphic_pair_section_sheaf} (that is, prior to restriction to the Coulomb-gauge slice). For that purpose, we shall appeal to the following analogue of \cite[Lemma 1.2]{FrMSWInvariants} due to Friedman and Morgan that allows us to work with a map of fixed Banach or Fr\'echet spaces rather than a section of a sheaf.

\begin{lem}[Local equivalence of holomorphic pair map with composition of holomorphic pair map and orthogonal projection]
\label{lem:Friedman_Morgan_1999_1-2}
Let $E$ be a smooth, complex vector bundle over a complex manifold $X$ and fixed holomorphic structure on the complex line bundle $\det E$. Let $(\bar\partial_{E_0},\varphi_0)$ be a smooth, holomorphic pair on $E$ that induces the given holomorphic structure on $\det E$. If $p\in (n,\infty)$ is a constant, then there is an open neighborhood
\[
  \UU_{\bar\partial_{E_0},\varphi_0} \subset W^{1,p}\left(\Lambda^{0,1}(\fsl(E)) \oplus E\right)
\]
of the origin such that
\[
  (\bar\partial_E,\varphi) \in (\bar\partial_{E_0},\varphi_0) + \UU_{\bar\partial_{E_0},\varphi_0}
\]
is a solution to \eqref{eq:Holomorphic_pair_Fredholm} if and only if it is a solution to
\begin{equation}
  \label{eq:Holomorphic_pair_projection}
  \Pi_{\bar\partial_{E_0},\varphi_0}^{0,2}(F_{\bar\partial_E},\bar\partial_E\varphi) = (0,0) \in
  \Ker\bar\partial_{E_0,\varphi_0}^2
  \cap L^p\left(\Lambda^{0,2}(\fsl(E)) \oplus \Lambda^{0,1}(E)\right),
\end{equation}
where
\[
  \Pi_{\bar\partial_{E_0},\varphi_0}^{0,2}:L^p\left(\Lambda^{0,2}(\fsl(E)) \oplus \Lambda^{0,1}(E)\right)
  \circlearrowleft
\]
is $L^2$-orthogonal projection onto the closed subspace
\[
  \Ker\bar\partial_{E_0,\varphi_0}^2
  \cap L^p\left(\Lambda^{0,2}(\fsl(E)) \oplus \Lambda^{0,1}(E)\right).
\]
\end{lem}

\begin{rmk}[Local trivializations of sheaves of sections of a vector bundle defined by the kernel of a differential operator]
\label{rmk:Local_trivializations_sheaves}
We note that in Morgan, Mrowka, and Ruberman \cite[Theorem 12.1.1 and its proof]{MMR} and Friedman and Morgan \cite[pp. 455--456]{FrMSWInvariants}, the authors are essentially considering
local trivializations of vector bundles through versions of Lemma \ref{lem:Friedman_Morgan_1999_1-2}.
\end{rmk}

\begin{proof}[Proof of Lemma \ref{lem:Friedman_Morgan_1999_1-2}]
We adapt Friedman and Morgan's proof of \cite[Lemma 1.2]{FrMSWInvariants}. If $(\bar\partial_E,\varphi)$ is a solution to \eqref{eq:Holomorphic_pair_Fredholm}, then it is clearly also a solution to \eqref{eq:Holomorphic_pair_projection}.

Conversely, suppose that $(\bar\partial_E,\varphi)$ is a solution to \eqref{eq:Holomorphic_pair_projection}, so that (see \eqref{eq:Rudin_theorem_4-12}) 
\[
  (F_{\bar\partial_E},\bar\partial_E\varphi)
  \in \left(\Ker\bar\partial_{E_0,\varphi_0}^2\right)^\perp =  \Ran\bar\partial_{E_0,\varphi_0}^{2,*}.
\]
We may thus write
\begin{equation}
  \label{eq:F_dbarA_and_dbarA_varphi_equals_dbarAvarphi2*_gamma_and_psi}
  (F_{\bar\partial_E},\bar\partial_E\varphi)
  = \bar\partial_{E_0,\varphi_0}^{2,*}(\gamma,\psi),
\end{equation}
for some
\[
  (\gamma,\psi) \in W^{1,p}\left(\Lambda^{0,3}(\fsl(E)) \oplus \Lambda^{0,2}(E)\right).
\]
Applying the Bianchi identity for $(0,1)$-pairs \eqref{eq:Holomorphic_pair_Bianchi_identity} to the preceding identity gives
\[
  \bar\partial_{E,\varphi}^2\bar\partial_{E_0,\varphi_0}^{2,*}(\gamma,\psi)
  =
  \bar\partial_{E,\varphi}^2(F_{\bar\partial_E},\bar\partial_E\varphi)
  =
  (0,0) \in L^p\left(\Lambda^{0,3}(\fsl(E)) \oplus \Lambda^{0,2}(E)\right).
\]
Therefore, we obtain
\[
  \Pi_{\bar\partial_{E_0},\varphi_0}^{0,3}\bar\partial_{E,\varphi}^2\bar\partial_{E_0,\varphi_0}^{2,*}(\gamma,\psi)
  =
  (0,0) \in \left(\Ker\bar\partial_{E_0,\varphi_0}^{2,*}\right)^\perp \cap
  W^{-1,p}\left(\Lambda^{0,3}(\fsl(E) \oplus \Lambda^{0,2}(E)\right),
\]
where
\[
  \Pi_{\bar\partial_{E_0},\varphi_0}^{0,3}: W^{-1,p}\left(\Lambda^{0,3}(\fsl(E)) \oplus \Lambda^{0,2}(E)\right) \circlearrowleft,
\]
is $L^2$-orthogonal projection onto the closed subspace $\Ran\bar\partial_{E_0,\varphi_0}^2 = (\Ker\bar\partial_{E_0,\varphi_0}^{2,*})^\perp$. We may assume without loss of generality that
\[
  (\gamma,\psi) \in \left(\Ker\bar\partial_{E_0,\varphi}^{2,*}\right)^\perp.
\]
We write
\[
  (\bar\partial_E,\varphi) - (\bar\partial_{E_0},\varphi_0)
  = (\alpha,\sigma) \in W^{1,p}\left(\Lambda^{0,1}(\fsl(E)) \oplus E\right)
\]
and now claim that the operator
\begin{multline}
  \label{eq:Pi_dbarEvarphi^2_dbarE0varphi0^2*}
  \Pi_{\bar\partial_{E_0},\varphi_0}^{0,3}\bar\partial_{E,\varphi}^2\bar\partial_{E_0,\varphi_0}^{2,*}:
  \left(\Ker\bar\partial_{E_0,\varphi}^{2,*}\right)^\perp
  \cap W^{1,p}\left(\Lambda^{0,3}(\fsl(E)) \oplus \Lambda^{0,2}(E)\right)
  \\
  \to
  \left(\Ker\bar\partial_{E_0,\varphi}^{2,*}\right)^\perp
  \cap W^{-1,p}\left(\Lambda^{0,3}(\fsl(E)) \oplus \Lambda^{0,2}(E)\right)
\end{multline}
is \emph{injective} for $(\alpha,\sigma)$ in a small enough neighborhood of the origin. We shall apply an argument that is dual to our proof of surjectivity of a certain perturbed Laplace operator (analogous to the $L^2$ adjoint $\bar\partial_{E_0,\varphi_0}^2\bar\partial_{E_0+\alpha}^*$) in Feehan and Maridakis \cite[Lemma 2.4.1, p. 30]{Feehan_Maridakis_Lojasiewicz-Simon_coupled_Yang-Mills}, but omit some details that are included in \cite{Feehan_Maridakis_Lojasiewicz-Simon_coupled_Yang-Mills}. (See also Morgan, Mrowka, and Ruberman \cite[Theorem 12.1.1 and its proof]{MMR} for related arguments.)

When $(\alpha,\sigma) = (0,0)$, so $\bar\partial_{E,\varphi}^2 = \bar\partial_{E_0+\alpha,\varphi_0+\sigma}^2 = \bar\partial_{E_0,\varphi_0}^2$, the resulting operator 
$\Pi_{\bar\partial_{E_0},\varphi_0}^{0,3}\bar\partial_{E_0,\varphi_0}^2\bar\partial_{E_0,\varphi_0}^{2,*}$ in \eqref{eq:Pi_dbarEvarphi^2_dbarE0varphi0^2*} is bijective and bounded, thus an isomorphism of Banach spaces by the Open Mapping Theorem (see Rudin \cite[Theorem 2.11]{Rudin}), with left inverse
\begin{multline*}
  G_{\bar\partial_{E_0},\varphi_0}:
  \left(\Ker\bar\partial_{E_0,\varphi}^{2,*}\right)^\perp
  \cap W^{-1,p}\left(\Lambda^{0,3}(\fsl(E)) \oplus \Lambda^{0,2}(E)\right)
  \\
  \to
  \left(\Ker\bar\partial_{E_0,\varphi}^{2,*}\right)^\perp
  \cap W^{1,p}\left(\Lambda^{0,3}(\fsl(E)) \oplus \Lambda^{0,2}(E)\right).
\end{multline*}
Hence, the following operator is invertible,
\[ G_{\bar\partial_{E_0},\varphi_0}\Pi_{\bar\partial_{E_0},\varphi_0}^{0,3}\bar\partial_{E_0,\varphi_0}^2\bar\partial_{E_0,\varphi_0}^{2,*}:
  \left(\Ker\bar\partial_{E_0,\varphi}^{2,*}\right)^\perp
  \cap W^{1,p}\left(\Lambda^{0,3}(\fsl(E)) \oplus \Lambda^{0,2}(E)\right)
  \circlearrowleft
\]
and by openness of the general linear group on a Banach space (see Kadison and Ringrose \cite[Proposition 3.1.6]{KadisonRingrose1}), the operator
\[ G_{\bar\partial_{E_0},\varphi_0}\Pi_{\bar\partial_{E_0},\varphi_0}^{0,3}\bar\partial_{E,\varphi}^2\bar\partial_{E_0,\varphi_0}^{2,*}:
  \left(\Ker\bar\partial_{E_0,\varphi}^{2,*}\right)^\perp
  \cap W^{1,p}\left(\Lambda^{0,3}(\fsl(E)) \oplus \Lambda^{0,2}(E)\right)
  \circlearrowleft
\]
is also invertible for $(\alpha,\sigma)$ in a small enough open neighborhood of the origin,
\[
  \UU_{\bar\partial_{E_0},\varphi_0} \subset W^{1,p}\left(\Lambda^{0,1}(\fsl(E)) \oplus E\right).
\]
In particular, the operator $\Pi_{\bar\partial_{E_0},\varphi_0}^{0,3}\bar\partial_{E,\varphi}^2\bar\partial_{E_0,\varphi_0}^{2,*}$ is injective for $(\alpha,\sigma) \in \UU_{\bar\partial_{E_0},\varphi_0}$, as claimed. Thus $(\gamma,\psi) = (0,0)$ and so $(F_{\bar\partial_E}, \bar\partial_E\varphi) = (0,0)$ by \eqref{eq:F_dbarA_and_dbarA_varphi_equals_dbarAvarphi2*_gamma_and_psi}. Therefore, $(\bar\partial_E,\varphi)$ is a solution to \eqref{eq:Holomorphic_pair_Fredholm}.
\end{proof}

\section{Fredholm complexes and their Euler characteristics}
\label{sec:Fredholm_complexes_Euler_characteristics}
It will be useful to relate the index of a complex with one that is partially or fully rolled up and, for this purpose, we generalize the discussion in Gilkey \cite[Section 1.5]{Gilkey2}.

\begin{lem}[Euler characteristics of Fredholm complexes]
\label{lem:Euler_characteristic_Fredholm_complex}
Let $\KK=\RR$ or $\CC$, let $n\geq 1$ be an integer, let $\fV_i$, for $i=0,\ldots,n+1$, be a finite sequence of $\KK$-Hilbert spaces, and let $d_i:\fV_i\to \fV_{i+1}$ be a sequence of bounded $\KK$-linear operators such that $d_{i+1}\circ d_i = 0$ and $d_i$ has finite-dimensional kernel, closed range in $\fV_{i+1}$, and finite codimension in $\Ker d_{i+1}$ for $i=0,\ldots,n$. Then the Euler characteristic
\begin{equation}
  \label{eq:Euler_characteristic_Fredholm_complex}
  \chi_\KK(\fV) := \sum_{i=0}^{n+1} (-1)^i \dim_\KK H^i(\fV)
\end{equation}
defined by the cohomology groups
\begin{equation}
  \label{eq:Fredholm_complex_cohomology_groups}
  H^i(\fV) := \Ker d_i / \Ran d_{i-1} \quad\text{for } i = 0,\ldots,n+1,
\end{equation}
(with $d_{-1}=0$ and $d_{n+1}=0$) of the complex\footnote{This is an example of a Fredholm complex in the terminology of Segal \cite{Segal_1970}.}
\begin{equation}
  \label{eq:Fredholm_complex}
  0 \xrightarrow{} \fV_0  \xrightarrow{d_0} \fV_1  \xrightarrow{d_1} \cdots \xrightarrow{d_{n-3}} \fV_{n-2} \xrightarrow{d_{n-2}} \fV_{n-1} \xrightarrow{d_{n-1}}  \fV_n \xrightarrow{d_n} \fV_{n+1}  \xrightarrow{} 0
\end{equation}
is equal to the Euler characteristic of each one of the following complexes:
\begin{subequations}
\begin{gather}
  \label{eq:Fredholm_complex_partially_rolled_up}
  0 \xrightarrow{} \fV_0  \xrightarrow{d_0} \fV_1  \xrightarrow{d_1} \cdots \xrightarrow{d_{n-3}} \fV_{n-2} \xrightarrow{d_{n-2}} \fV_{n-1}\oplus\fV_{n+1} \xrightarrow{d_{n-1} + d_n^*} \fV_n  \xrightarrow{} 0
  \\
  \label{eq:Fredholm_complex_partially_rolled-up_modified}
  0 \xrightarrow{} \fV_0  \xrightarrow{d_0} \fV_1  \xrightarrow{d_1} \cdots \xrightarrow{d_{n-3}} \fV_{n-2} \xrightarrow{d_{n-2}} \fV_{n-1} \oplus \Ker d_n^* \xrightarrow{d_{n-1}} \Ker d_n  \xrightarrow{} 0
\end{gather}
\end{subequations}
\end{lem}  

\begin{proof}
By analogy with \eqref{eq:Rudin_theorem_4-12}, we recall that a specialization of Rudin \cite[Theorem 4.12]{Rudin} to the case of bounded operators on Hilbert spaces yields, for $i=0,\ldots,n$,
\begin{equation}
  \label{eq:Rudin_theorem_4-12_Fredholm_complex}
  \Ker d_i^* = (\Ran d_i)^\perp \quad\text{and}\quad \Ker d_i = (\Ran d_i^*)^\perp.
\end{equation}
By applying \eqref{eq:Rudin_theorem_4-12_Fredholm_complex}, we see that for $i = 0,\ldots,n+1$,
\[
  \fV_i = \Ran d_{i-1} \oplus (\Ran d_{i-1})^\perp = \Ran d_{i-1} \oplus \Ker d_{i-1}^*.
\]
We define the harmonic spaces for the complex \eqref{eq:Fredholm_complex} By
\begin{equation}
  \label{eq:Fredholm_complex_harmonic_spaces}
  \bH^i(\fV) := \Ker (d_i + \Ker d_{i-1}^*), \quad\text{for } i = 0,\ldots,n+1,
\end{equation}
and observe that
\[
  \Ker (d_i + \Ker d_{i-1}^*)  = \Ker d_i \cap \Ker d_{i-1}^*
  = \Ker d_i \cap (\Ran d_{i-1})^\perp \cong \Ker d_i/\Ran d_{i-1},
\]
that is, by \eqref{eq:Fredholm_complex_cohomology_groups} and \eqref{eq:Fredholm_complex_harmonic_spaces},
\begin{equation}
  \label{eq:Cohomology_groups_isomorphic_to_harmonic_spaces}
  \bH^i(\fV) \cong H^i(\fV), \quad i = 0,\ldots,n+1.
\end{equation}
We now make the following claim that we shall prove shortly:
\begin{subequations}
  \label{claim:Isomorphisms_last_three_terms_Fredholm_complex}
  \begin{align}
    \label{claim:Isomorphism_last_and_lastminustwo_terms_Fredholm_complex}
    \Ker(d_{n-1} + d_n^*)/\Ran d_{n-2} &\cong H^{n-1}(\fV)\oplus H^{n+1}(\fV),
    \\
    \label{claim:Isomorphism_lastminusone_term_Fredholm_complex}
    \Coker(d_{n-1} + d_n^*) &\cong H^n(\fV).
  \end{align}    
\end{subequations}
Assuming \eqref{claim:Isomorphisms_last_three_terms_Fredholm_complex}, we see that the Euler characteristic of the complex \eqref{eq:Fredholm_complex_partially_rolled_up} is given by
\begin{align*}
  {}&(-1)^{n-1}\dim_\KK\Ker (d_{n-1}+d_n^*)/\Ran d_{n-2} + (-1)^n\dim_\KK\Coker (d_{n-1}+d_n^*)
  \\
  &\quad + \sum_{i=0}^{n-2} (-1)^i \dim_\KK H^i(\fV)
  \\
    &= (-1)^{n-1}\dim_\KK\left(\bH^{n-1}(\fV)\oplus \bH^{n+1}(\fV)\right)
  \\
    &\qquad + (-1)^n\dim_\KK \bH^n(\fV) + \sum_{i=0}^{n-2} (-1)^i \dim_\KK \bH^i(\fV)
      \quad\text{(by \eqref{eq:Cohomology_groups_isomorphic_to_harmonic_spaces} and  \eqref{claim:Isomorphisms_last_three_terms_Fredholm_complex})}
  \\
  &= \chi_\KK(\fV) \quad\text{(by \eqref{eq:Euler_characteristic_Fredholm_complex})},
\end{align*}
so the complexes \eqref{eq:Fredholm_complex} and \eqref{eq:Fredholm_complex_partially_rolled_up} have the same Euler characteristics.

We now verify the claim \eqref{claim:Isomorphisms_last_three_terms_Fredholm_complex}. By \eqref{eq:Rudin_theorem_4-12_Fredholm_complex} and noting that $\Ran d_{n-2} \subset  \Ker d_{n-1}$, we have
\begin{align*}
  \Ker(d_{n-1} + d_n^*)/\Ran d_{n-2}
  &\cong
  \Ker(d_{n-1} + d_n^*)\cap(\Ran d_{n-2})^\perp
  \\
  &=
  \Ker(d_{n-1} + d_n^*)\cap\Ker d_{n-2}^*.
\end{align*}
We observe that the subspaces,
\[
  \Ran d_{n-1} \quad\text{and}\quad \Ran d_n^*  \subset \fV_n
\]
are orthogonal since
\[
  d_n\circ d_{n-1} = 0 \quad\text{on } \fV_{n-1},
\]
and thus
\[
  (d_{n-1}v_{n-2}, d_n^*v_n)_{\fV_{n-1}}  = ((d_n\circ d_{n-1})v_{n-2}, v_n)_{\fV_n} = 0,
\]
for all $v_{n-2} \in \fV_{n-2}$ and $v_n \in \fV_n$. Therefore, 
\begin{align*}
  \Ker(d_{n-1} + d_n^*)/\Ran d_{n-2}
  &\cong \left(\Ker d_{n-1} \cap \Ker d_{n-2}^* \right) \oplus \Ker d_n^*
  \\
  &= \Ker\left(d_{n-1} + d_{n-2}^* \right) \oplus \Ker d_n^*,
\end{align*}
where we used the fact that $d_n^*$ and $d_{n-2}^*$  have domains and codomains given by
\[
  d_{n-2}^* : \fV_{n-1} \to \fV_{n-2}
  \quad\text{and}\quad
  d_n^* : \fV_{n+1} \to \fV_n.
\]
Consequently, we have
\[
  \Ker(d_{n-1} + d_n^*)/\Ran d_{n-2}
  \cong 
  \bH^{n-1}(\fV)\oplus \bH^{n+1}(\fV),
\]
verifying claim \eqref{claim:Isomorphism_last_and_lastminustwo_terms_Fredholm_complex}. To obtain the second equality in \eqref{claim:Isomorphisms_last_three_terms_Fredholm_complex}, we observe that
\begin{align*}
  \Coker(d_{n-1} + d_n^*)
  &\cong \Ker (d_{n-1} + d_n^*)^*
  \\
  &= \Ker (d_n + d_{n-1}^*)
  \\
  &= \bH^n(\fV) \quad\text{(by \eqref{eq:Fredholm_complex_harmonic_spaces})},
\end{align*}
verifying claim \eqref{claim:Isomorphism_lastminusone_term_Fredholm_complex}. 

We now compare the Euler characteristics of \eqref{eq:Fredholm_complex} and \eqref{eq:Fredholm_complex_modified}. Note that
\[
  \fV_{n+1} = \Ran d_n \oplus (\Ran d_n)^\perp = \Ran d_n \oplus \Ker d_n^*
\]
while
\[
  \fV_n = \Ker d_n \oplus (\Ker d_n)^\perp = \Ker d_n \oplus \Ran d_n^*
\]
and so the complex \eqref{eq:Fredholm_complex} is equal to
\begin{equation}
  \label{eq:Fredholm_complex_decomposed}
  0 \xrightarrow{} \fV_0  \xrightarrow{d_0} \fV_1  \xrightarrow{d_1} \cdots \xrightarrow{d_{n-2}} \fV_{n-1} \xrightarrow{d_{n-1}} \Ker d_n \oplus (\Ker d_n)^\perp \xrightarrow{d_n} \Ran d_n \oplus (\Ran d_n)^\perp  \xrightarrow{} 0
\end{equation}
Because $d_n:(\Ker d_n)^\perp \to \Ran d_n$ is an isomorphism of Hilbert spaces by the Open Mapping Theorem \cite[Theorem 2.11]{Rudin}, the complex \eqref{eq:Fredholm_complex_decomposed} has the same cohomology groups as that of the modified complex
\[
  0 \xrightarrow{} \fV_0  \xrightarrow{d_0} \fV_1  \xrightarrow{d_1} \cdots \xrightarrow{d_{n-2}} \fV_{n-1} \xrightarrow{d_{n-1}} \Ker d_n \xrightarrow{\times 0} (\Ran d_n)^\perp  \xrightarrow{} 0
\]
By \eqref{eq:Rudin_theorem_4-12_Fredholm_complex}, the preceding complex is equal to the complex
\begin{equation}
  \label{eq:Fredholm_complex_modified}
  0 \xrightarrow{} \fV_0  \xrightarrow{d_0} \fV_1  \xrightarrow{d_1} \cdots \xrightarrow{d_{n-2}} \fV_{n-1} \xrightarrow{d_{n-1}} \Ker d_n \xrightarrow{\times 0} \Ker d_n^*  \xrightarrow{} 0
\end{equation}
Therefore, the complexes \eqref{eq:Fredholm_complex} and \eqref{eq:Fredholm_complex_modified} have the same Euler characteristics.

We can now compare the Euler characteristics of \eqref{eq:Fredholm_complex} and \eqref{eq:Fredholm_complex_partially_rolled-up_modified}. The Euler characteristic of the complex \eqref{eq:Fredholm_complex_decomposed} is equal to that of the partially rolled-up complex
\[
  0 \xrightarrow{} \fV_0  \xrightarrow{d_0} \fV_1  \xrightarrow{d_1} \cdots \xrightarrow{d_{n-2}}
  \fV_{n-1} \oplus \Ran d_n \oplus (\Ran d_n)^\perp
  \xrightarrow{d_{n-1} + d_n^*} \Ker d_n \oplus (\Ker d_n)^\perp  \xrightarrow{} 0
\]
By \eqref{eq:Rudin_theorem_4-12_Fredholm_complex}, the preceding complex is equal to the complex
\[
  0 \xrightarrow{} \fV_0  \xrightarrow{d_0} \fV_1  \xrightarrow{d_1} \cdots \xrightarrow{d_{n-2}}
  \fV_{n-1} \oplus (\Ker d_n^*)^\perp \oplus \Ker d_n^*
  \xrightarrow{d_{n-1} + d_n^*} \Ker d_n \oplus \Ran d_n^*  \xrightarrow{} 0
\]
The Euler characteristic of the preceding complex is equal to that of the modified, partially rolled-up complex
\[
  0 \xrightarrow{} \fV_0  \xrightarrow{d_0} \fV_1  \xrightarrow{d_1} \cdots \xrightarrow{d_{n-2}}
  \fV_{n-1} \oplus \Ker d_n^*
  \xrightarrow{d_{n-1}} \Ker d_n  \xrightarrow{} 0
\]
and which is the complex \eqref{eq:Fredholm_complex_partially_rolled-up_modified}. The preceding complex is a partially rolled up version of the complex \eqref{eq:Fredholm_complex_modified} and thus has the same Euler characteristic as that of \eqref{eq:Fredholm_complex_modified} and hence also the same Euler characteristic as that of \eqref{eq:Fredholm_complex}, as claimed.
\end{proof}

Lemma \ref{lem:Euler_characteristic_Fredholm_complex} yields the useful

\begin{cor}[Comparison of elliptic complexes and cohomology groups for pre-holomorphic and holomorphic pairs]
\label{cor:Comparison_elliptic_complexes_cohomology_groups_pre-holomorphic_and_holomorphic_pairs}
Let $(E,h)$ be a smooth Hermitian vector bundle over a closed complex K\"ahler surface $(X,g,J)$ and $(\rho_\can,W_\can)$ be the canonical spin${}^c$ structure over $X$ (see \eqref{eq:Canonical_spinc_bundles} for $W_\can$ and \eqref{eq:Canonical_Clifford_multiplication} for $\rho_\can$). If $(A,(\varphi,\psi))$ is a pre-holomorphic pair in the sense of Section \ref{sec:Elliptic_deformation_complex_for_pre-holomorphic_pair_equations_complex_surface} and $\psi = 0$, then $(\bar\partial_A,\varphi)$ is a holomorphic pair in the sense of \eqref{eq:Holomorphic_pair} and the following equalities of complex vector spaces hold,
\begin{equation}
\label{eq:Cohomologies_complexes_pre-holomorphic_pair_type1_and_holomorphic_pair_n_is_2}
  \bH_{\bar\partial_{A,(\varphi,0)}}^0 = \bH_{\bar\partial_A,\varphi}^0,
  \quad
  \bH_{\bar\partial_{A,(\varphi,0)}}^1 = \bH_{\bar\partial_A,\varphi}^1\oplus \bH_{\bar\partial_A,\varphi}^3,
  \quad
  \bH_{\bar\partial_{A,(\varphi,0)}}^2 = \bH_{\bar\partial_A,\varphi}^2,
\end{equation}
where the spaces $\bH_{\bar\partial_{A,(\varphi,0)}}^\bullet$ are as in \eqref{eq:H_dbar_APhi^0bullet} and the spaces $\bH_{\bar\partial_A,\varphi}^\bullet$ are as in \eqref{eq:H_dbar_Avarphi^0bullet}. Moreover, the Euler characteristics over $\CC$ of the elliptic complexes \eqref{eq:Pre-holomorphic_pair_elliptic_complex} with $\psi\equiv 0$ and \eqref{eq:Holomorphic_pair_elliptic_complex} when $n=2$ are equal.
\end{cor}

\begin{proof}
The fact that $(\bar\partial_A,\varphi)$ is a holomorphic pair is immediate from the definitions in Section  \ref{sec:Elliptic_deformation_complex_for_pre-holomorphic_pair_equations_complex_surface} and equation \eqref{eq:Holomorphic_pair}. By definition of the differentials $\bar\partial_{A,(\varphi,\psi)}^\bullet$ and $\bar\partial_{A,\varphi}^\bullet$ in the complexes \eqref{eq:Pre-holomorphic_pair_elliptic_complex} and \eqref{eq:Holomorphic_pair_elliptic_complex}, respectively, we see that the complex \eqref{eq:Pre-holomorphic_pair_elliptic_complex} with $\psi\equiv 0$ is equal to
\begin{equation}
\label{eq:Pre-holomorphic_pair_elliptic_complex_type1}
\begin{CD}
0
@>>>
\Omega^0(\fsl(E))
@> \bar\partial_{A,\varphi}^0 >>
\begin{matrix}
\Omega^{0,1}(\fsl(E))
\\
\oplus
\\
\Omega^0(E)\oplus \Omega^{0,2}(E)
\end{matrix}
@> \bar\partial_{A,\varphi}^1 + \bar\partial_{A,\varphi}^{2,*}>>
\begin{matrix}
\Omega^{0,2}(\fsl(E))
\\
\oplus
\\
\Omega^{0,1}(E)
\end{matrix}
@>>>
0
\end{CD}
\end{equation}
The complex \eqref{eq:Pre-holomorphic_pair_elliptic_complex_type1} is the partially rolled-up version of the complex \eqref{eq:Holomorphic_pair_elliptic_complex} when $n=2$:
\begin{equation}
\label{eq:Holomorphic_pair_elliptic_complex_n_is_2}
\begin{CD}
0
@>>>
\Omega^0(\fsl(E))
@> \bar\partial_{A,\varphi}^0 >>
\begin{matrix}
\Omega^{0,1}(\fsl(E))
\\
\oplus
\\
\Omega^0(E)
\end{matrix}
@> \bar\partial_{A,\varphi}^1 >>
\begin{matrix}
\Omega^{0,2}(\fsl(E))
\\
\oplus
\\
\Omega^{0,1}(E)
\end{matrix}
@> \bar\partial_{A,\varphi}^2 >>
\begin{matrix}
0
\\
\oplus
\\
\Omega^{0,2}(E)
\end{matrix}
@>>>
0
\end{CD}
\end{equation}
Specializing Lemma \ref{lem:Euler_characteristic_Fredholm_complex} to the case $n=2$ gives the following Fredholm complexes, 
\begin{gather*}
  0 \xrightarrow{} \fV_0  \xrightarrow{d_0} \fV_1  \xrightarrow{d_1} \fV_2 \xrightarrow{d_2} \fV_3 \xrightarrow{} 0
  \\
  0 \xrightarrow{} \fV_0  \xrightarrow{d_0} \fV_1\oplus\fV_3 \xrightarrow{d_1+d_2^*} \fV_2  \xrightarrow{} 0
  \\
  0 \xrightarrow{} \fV_0  \xrightarrow{d_0} \fV_1\oplus \Ker d_2^* \xrightarrow{d_1} \Ker d_2  \xrightarrow{} 0
\end{gather*}
all of which have the same Euler characteristic In particular, by the claim \eqref{claim:Isomorphisms_last_three_terms_Fredholm_complex} and definitions \eqref{eq:Fredholm_complex_cohomology_groups} and \eqref{eq:Fredholm_complex_harmonic_spaces} in Lemma \ref{lem:Euler_characteristic_Fredholm_complex} and its proof, the harmonic representatives of the cohomology groups of the complexes \eqref{eq:Pre-holomorphic_pair_elliptic_complex_type1} and \eqref{eq:Holomorphic_pair_elliptic_complex_n_is_2} are related by \eqref{eq:Cohomologies_complexes_pre-holomorphic_pair_type1_and_holomorphic_pair_n_is_2}, as claimed, and by Lemma \ref{lem:Euler_characteristic_Fredholm_complex}, their Euler characteristics are equal over $\CC$.
\end{proof}

\section{Witten's dichotomy for solutions to the linearized non-Abelian monopole equations and vanishing of cohomology groups}
\label{sec:Witten_dichotomy_for_linearization_SO3_monopole_equations}
We observe that
\[
  \Ker\bar\partial_{A,\varphi}^0 \cap \Omega^0(\fsl(E)) = \Ker\bar\partial_{A,(\varphi,0)}^0 \cap \Omega^0(\fsl(E))
\]
and so the degree zero cohomology groups of \eqref{eq:d1StableComplex} and \eqref{eq:Holomorphic_pair_elliptic_complex_rolled_up} are equal. To give insight into how the degree one cohomology groups of \eqref{eq:d1StableComplex} and \eqref{eq:Holomorphic_pair_elliptic_complex_rolled_up} compare, recall that if $\psi\equiv 0$ and $\varphi\not\equiv 0$, then all non-Abelian monopoles defining points in an open neighborhood of $[A,(\varphi,0)]$ in $\sM_\ft$ are also type $1$ since the property that $\varphi$ is not identically zero holds on an open neighborhood of $[A,(\varphi,0)]$ and the Witten's dichotomy holds for non-Abelian monopoles by Lemma \ref{lem:Okonek_Teleman_1995_3-1}. Hence,
one would expect elements of the kernel of $\bar\partial_{A,(\varphi,0)}^1$ to have the form $(a'',\sigma,0)$ and we justify this expectation in the following

\begin{lem}[Witten's dichotomy for solutions to the linearized pre-holomorphic pair equations on Hermitian vector bundles over complex Hermitian manifolds]
\label{lem:Witten_dichotomy_linearized_type1_SO3_monopole_equation}
Continue the hypotheses of Lemma \ref{lem:Okonek_Teleman_1995_3-1}. If $\psi=0\in \Omega^{0,2}(E)$ and $(a'',\sigma,\tau)\in\Ker\partial_{A,(\varphi,0)}^1$, where $\bar\partial_{A,(\varphi,0)}^1$ is the differential in \eqref{eq:d1StableComplex}, then $\tau\equiv 0$.
\end{lem}

\begin{proof}
The argument is similar to the proof of Lemma \ref{lem:Okonek_Teleman_1995_3-1}, to which we refer for additional details. Because $F_A^{0,2}=0$ by \eqref{eq:FdetA02_is_zero} and \eqref{eq:SO(3)_monopole_equations_(0,2)_curvature} and because $\bar\partial_A\varphi=0$ by \eqref{eq:SO(3)_monopole_equations_Dirac_almost_Kaehler}, applying $\bar\partial_A$ to the second component of $\bar\partial_{A,(\varphi,0)}^1(a'',\sigma,\tau)=0$ given by \eqref{eq:d1StableComplex} yields
\[
0=\bar\partial_A^2\sigma+\bar\partial_A\bar\partial_A^*\tau +\bar\partial_A(a''\varphi)
=
\bar\partial_A\bar\partial_A^*\tau+(\bar\partial_Aa'')\varphi.
\]
We substitute the first component
of  $\bar\partial_{A,(\varphi,0)}^1(a'',\sigma,\tau)=0$ given by \eqref{eq:d1StableComplex} into the preceding identity to give
\[
0
=\bar\partial_A\bar\partial_A^*\tau + (\tau\otimes\varphi^*)_0\varphi.
\]
Using $(\tau\otimes\varphi^*)_0=\tau\otimes\varphi^*-\frac{1}{2}\langle\tau,\varphi\rangle_E\,\id_E$,
we can rewrite the preceding identity as
\[
0
=\bar\partial_A\bar\partial_A^*\tau + \tau|\varphi|_E^2 -\frac{1}{2}\langle \tau,\varphi\rangle_E\,\varphi.
\]
Taking the pointwise $\Lambda^{0,2}(E)$-inner product of the preceding equation with $\tau$ and applying the Cauchy--Schwartz inequality yields
\begin{align*}
0 &=\langle \bar\partial_A\bar\partial_A^*\tau,\tau\rangle_{\Lambda^{0,2}(E)}
+
|\tau|_{\Lambda^{0,2}(E)}^2|\varphi|_E^2 - \frac{1}{2}|\langle \tau,\varphi\rangle_E|_{\Lambda^{0,2}(X)}^2
\\
&\ge \langle \bar\partial_A\bar\partial_A^*\tau,\tau\rangle_{\Lambda^{0,2}(E)}
+
\frac{1}{2}|\tau|_{\Lambda^{0,2}(E)}^2|\varphi|_E^2.
\end{align*}
Integrating this inequality over $X$ and integrating by parts yields
\[
0\ge\|\bar\partial^*_A\tau\|_{L^2(X)}^2 +\frac{1}{2}\int_X|\tau|_{\Lambda^{0,2}(E)}^2|\varphi|_E^2\,d\vol.
\]
Thus, $\bar\partial^*_A\tau=0$ on $X$ and because $\varphi\not\equiv 0$ by hypothesis, we must have $\tau\equiv 0$ on an open subset of $X$ and hence we again have $\tau \equiv 0$ on $X$ by the unique continuation property.
\end{proof}

Lemma \ref{lem:Witten_dichotomy_linearized_type1_SO3_monopole_equation} yields the forthcoming vanishing result for a cohomology group of the complex \eqref{eq:Holomorphic_pair_elliptic_complex} when $X$ has dimension two.

\begin{cor}[Vanishing of third-order cohomology group for holomorphic pair elliptic complex over a complex K\"ahler surface]
\label{cor:Vanishing_third-order_cohomology_group_holomorphic_pair_elliptic_complex_Kaehler_surface}
Continue the hypotheses of Lemma \ref{lem:Witten_dichotomy_linearized_type1_SO3_monopole_equation}. Then
\begin{equation}
\label{eq:Cohomology_holomorphic_pair_n_is_2_vanishing}
  \bH_{\bar\partial_A,\varphi}^3 = 0
\end{equation}
and the relationship \eqref{eq:Cohomologies_complexes_pre-holomorphic_pair_type1_and_holomorphic_pair_n_is_2} between the harmonic representatives of the cohomology groups for the elliptic complexes \eqref{eq:Pre-holomorphic_pair_elliptic_complex} with $\psi\equiv 0$ and \eqref{eq:Holomorphic_pair_elliptic_complex} when $X$ has dimension two simplifies to
\begin{equation}
\label{eq:Cohomologies_complexes_pre-holomorphic_pair_type1_and_holomorphic_pair_n_is_2_simplified}
  \bH_{\bar\partial_{A,(\varphi,0)}}^\bullet = \bH_{\bar\partial_A,\varphi}^\bullet,
\end{equation}
where the spaces $\bH_{\bar\partial_{A,(\varphi,0)}}^\bullet$ are as in \eqref{eq:H_dbar_APhi^0bullet} and the spaces $\bH_{\bar\partial_A,\varphi}^\bullet$ are as in \eqref{eq:H_dbar_Avarphi^0bullet}.
\end{cor}

\begin{proof}
By definition \eqref{eq:H_dbar_Avarphi^0bullet}, we have
\begin{align*}
  \bH_{\bar\partial_A,\varphi}^1
  &= \Ker\left(\bar\partial_{A,\varphi}^1 + \bar\partial_{A,\varphi}^{0,*}:\Omega^{0,1}(\fsl(E))\oplus \Omega^0(E)
       \to \Omega^{0,2}(\fsl(E))\oplus \Omega^{0,1}(E) \oplus \Omega^0(\fsl(E))\right),
  \\
  \bH_{\bar\partial_A,\varphi}^3
  &= \Ker\left(\bar\partial_{A,\varphi}^{2,*}:\Omega^{0,2}(E) \to \Omega^{0,2}(\fsl(E)) \oplus \Omega^{0,1}(E)\right),
\end{align*}
while by definition \eqref{eq:H_dbar_APhi^01}, we have
\[
  \bH_{\bar\partial_{A,(\varphi,\psi)}}^1 \subset \Ker\left(\bar\partial_{A,(\varphi,\psi)}^1:\Omega^{0,1}(\fsl(E))\oplus \Omega^0(E)\oplus \Omega^{0,2}(E) \to \Omega^{0,2}(\fsl(E))\oplus \Omega^{0,1}(E)\right).
\]
From the proof of Corollary \ref{cor:Comparison_elliptic_complexes_cohomology_groups_pre-holomorphic_and_holomorphic_pairs}, we have
\[
  \bar\partial_{A,(\varphi,0)}^1
  =
  \bar\partial_{A,\varphi}^1 + \bar\partial_{A,\varphi}^{2,*},
\]
where the domains and codomains of $\bar\partial_{A,\varphi}^1$ and $\bar\partial_{A,\varphi}^{2,*}$ are given by
\begin{align*}
  \bar\partial_{A,\varphi}^1&:\Omega^{0,1}(\fsl(E))\oplus \Omega^0(E) \to \Omega^{0,2}(\fsl(E))\oplus \Omega^{0,1}(E),
  \\
  \bar\partial_{A,\varphi}^{2,*}&:\Omega^{0,2}(E) \to \Omega^{0,2}(\fsl(E))\oplus \Omega^{0,1}(E)
\end{align*}
Lemma \ref{lem:Witten_dichotomy_linearized_type1_SO3_monopole_equation} implies that $(\alpha,\sigma,\tau) \in \bH_{\bar\partial_{A,(\varphi,0)}}^1$ if and only if $\tau\equiv 0 \in \bar\partial_{A,\varphi}^{2,*}$ and because
\[
  \bH_{\bar\partial_{A,(\varphi,0)}}^1
  =
  \bH_{\bar\partial_A,\varphi}^1 \oplus \bH_{\bar\partial_A,\varphi}^3
\]
by the identity \eqref{eq:Cohomologies_complexes_pre-holomorphic_pair_type1_and_holomorphic_pair_n_is_2} from Corollary \ref{cor:Comparison_elliptic_complexes_cohomology_groups_pre-holomorphic_and_holomorphic_pairs}, so $\bH_{\bar\partial_A,\varphi}^3 = 0$, as claimed in \eqref{eq:Cohomology_holomorphic_pair_n_is_2_vanishing} and so the identities \eqref{eq:Cohomologies_complexes_pre-holomorphic_pair_type1_and_holomorphic_pair_n_is_2} simplify to those in \eqref{eq:Cohomologies_complexes_pre-holomorphic_pair_type1_and_holomorphic_pair_n_is_2_simplified}.
\end{proof}

\chapter[Comparison of elliptic complexes]{Comparison of elliptic complexes}
\label{chap:Elliptic_deformation_complex_moduli_space_SO(3)_monopoles_over_almost_Hermitian_four-manifold}
Our goal in this chapter is to compare the cohomology groups (more specifically, harmonic spaces) for the elliptic complexes arising in this monograph, primarily between those for the non-Abelian monopole and holomorphic pair equations over complex K\"ahler surfaces but, with later applications in mind, we do the comparisons in more generality. In Sections \ref{sec:Isomorphisms_between_first-order_cohomology_groups}, \ref{sec:Isomorphisms_between_zeroth-order_cohomology_groups}, and \ref{sec:Isomorphisms_between_second-order_cohomology_groups} we compare the first, zeroth, and second-order harmonic spaces, respectively, defined by the elliptic deformation complexes for the non-Abelian monopole and pre-holomorphic pair equations over complex K\"ahler surfaces, making use of results described in Section \ref{sec:Spincu_structures_Dirac_operators_over_almost_Hermitian_manifolds} for spin${}^c$ and spin${}^u$ structures over almost Hermitian manifolds in order to help compare the second-order harmonic spaces. In Section \ref{sec:Elliptic_deformation_complex_moduli_space_ASD_connections_over_complex_Kaehler_surface}, we specialize the preceding results to give a comparison between the elliptic complexes for the anti-self-dual and holomorphic curvature equations over complex K\"ahler surfaces. Section \ref{sec:Elliptic_complex_moduli_space_HE_connections_over_complex_Kaehler_manifold} contains a generalization of the preceding comparison to one between the elliptic complexes for the projectively Hermitian--Einstein and holomorphic curvature equations over complex K\"ahler manifolds of arbitrary dimension. Lastly, in Section \ref{sec:Elliptic_complex_moduli_space_projective vortices_over_complex_Kaehler_manifold} we compare the elliptic complexes for the projective vortex and holomorphic pair equations over complex K\"ahler manifolds.

\section[Isomorphisms between first-order cohomology groups]{Comparison of elliptic complexes for non-Abelian monopole and pre-holomorphic pair equations over complex K\"ahler surfaces: First-order cohomology groups}
\label{sec:Isomorphisms_between_first-order_cohomology_groups}
Our goal in this section and in Sections \ref{sec:Isomorphisms_between_zeroth-order_cohomology_groups} and \ref{sec:Isomorphisms_between_second-order_cohomology_groups} is to generalize Friedman and Morgan \cite[Section 4.3.3, Proposition 3.7, p. 324]{FrM} for complex K\"ahler surfaces --- which compares the elliptic deformation complex \eqref{eq:Elliptic_deformation_complex_ASD_equation} for the anti-self-dual equation and the elliptic deformation complex for the holomorphic curvature equation \eqref{eq:Elliptic_deformation_complex_holomorphic_curvature_equation} --- to give a comparison of the elliptic deformation complex \eqref{eq:SO3MonopoleDefComplex} for a solution $(A,\Phi)$ to the non-Abelian monopole equation \eqref{eq:SO(3)_monopole_equations} and the elliptic deformation complex \eqref{eq:Pre-holomorphic_pair_elliptic_complex} for the corresponding solution $(\bar\partial_A,(\varphi,\psi)) = \pi_{0,1}(A,\Phi)$ to the pre-holomorphic pair equations \eqref{eq:SO(3)_monopole_equations_(0,2)_curvature}
and \eqref{eq:SO(3)_monopole_equations_Dirac_almost_Kaehler} defined by $\pi_{0,1}A = \bar\partial_A$ and $\Phi=(\varphi,\psi)$.

We shall later specialize this comparison to the case of a type $1$ solution $(A,(\varphi,0))$ to the non-Abelian monopole equations, equivalent to a solution $(\bar\partial_A,\varphi)$ to the projective vortex equations \eqref{eq:SO(3)_monopole_equations_almost_Hermitian_alpha}, and the corresponding solution $(\bar\partial_A,\varphi)$ to the holomorphic pair equations \eqref{eq:Holomorphic_pair}. While our calculations would become simpler if we made this assumption initially, the more general comparison that we describe here will be important for generalizations of the methods and results of this monograph from the category of closed, complex K\"ahler surfaces to the categories of closed, symplectic four-manifolds or even standard four-manifolds of Seiberg--Witten simple type. 

\begin{prop}[Canonical real linear isomorphism between first-order cohomology groups for elliptic complexes over closed K\"ahler surfaces]
\label{prop:Itoh_1985_proposition_2-4_SO3_monopole_complex_Kaehler}
Let $E$ be a smooth Hermitian vector bundle over a closed complex K\"ahler surface $(X,g,J)$ and $(\rho_\can,W_\can)$ denote the canonical spin${}^c$ structure over $X$ (see \eqref{eq:Canonical_spinc_bundles} for $W_\can$ and \eqref{eq:Canonical_Clifford_multiplication} for $\rho_\can$). If
\[
  (A,\Phi) \in \sA(E,h) \times \Omega^0(W^+\otimes E)
\]
is a smooth solution to the unperturbed non-Abelian monopole equations \eqref{eq:SO(3)_monopole_equations_Kaehler}, then
\[
  (\bar\partial_A,\varphi,\psi) \in \sA^{0,1}(E) \times \Omega^0(E)\oplus\Omega^{0,2}(E)
\]
is the solution to the pre-holomorphic pair equations \eqref{eq:SO(3)_monopole_equations_(0,2)_curvature} and \eqref{eq:SO(3)_monopole_equations_Dirac_almost_Kaehler} obtained by using \eqref{eq:Decompose_a_in_Omega1suE_into_10_and_01_components} and \eqref{eq:Canonical_spinc_bundles}, respectively, to write
\begin{align}
  \label{eq:del_bar_A}
  \bar\partial_A &= \frac{1}{2}\pi_{p,q+1}\circ d_A: \Omega^{p,q}(\fsl(E)) \to \Omega^{p,q+1}(\fsl(E)),
  \\
  \label{eq:Phi_coupled_spinor_is_pair_sections_E_oplus_02E}
  \Phi &= (\varphi,\psi) \in \Omega^0(E)\oplus\Omega^{0,2}(E).
\end{align}
The induced real linear map of infinite-dimensional real vector spaces
\begin{equation}
  \label{eq:Real_linear_map_tangent_spaces_SO(3)_pairs_to_preholomorphic_pairs}
  \Omega^1(\su(E))\oplus \Omega^0(W^+\otimes E) \ni (a,\phi)
  \mapsto
  (a'',\sigma,\tau) \in \Omega^{0,1}(\fsl(E)) \oplus \Omega^0(E)\oplus\Omega^{0,2}(E)
\end{equation}
implied by \eqref{eq:Decompose_a_in_Omega1suE_into_10_and_01_components} and \eqref{eq:Canonical_spinc_bundles} determine an isomorphism of finite-dimensional real vector spaces,
\begin{equation}
\label{eq:Itoh_1985_proposition_2-4_SO3_monopole_complex_Kaehler_isomorphism_widehat}
\bH_{A,\Phi}^1 \cong \widehat\bH_{\bar\partial_{A,(\varphi,\psi)}}^1,
\end{equation}
where $\bH_{A,\Phi}^1$ is defined by \eqref{eq:H_APhi^1} and $\widehat\bH_{\bar\partial_{A,(\varphi,\psi)}}^1$ is defined to be the finite-dimensional vector space of all solutions $(a'',\sigma,\tau) \in \Omega^{0,1}(\fsl(E))\oplus \Omega^0(E)\oplus\Omega^{0,2}(E)$ to equations \eqref{eq:H_dbar_APhi^01_explicit_02}, \eqref{eq:H_dbar_APhi^01_explicit_01}, and
\begin{equation}
  \label{eq:widehatH_dbar_APhi^01_explicit_0}
    \bar\partial_A^*a'' - R_\varphi^*\sigma + R_\psi^*\tau = 0 \in \Omega^0(\fsl(E)).
\end{equation}
Finally, if $(A,\Phi)$ is a type $1$ solution in the sense of Remark \ref{rmk:Projective_vortices_type_1_monopoles}, so $\psi \equiv 0$, then the isomorphism \eqref{eq:Itoh_1985_proposition_2-4_SO3_monopole_complex_Kaehler_isomorphism_widehat} simplifies to
\begin{equation}
\label{eq:Itoh_1985_proposition_2-4_SO3_monopole_complex_Kaehler_isomorphism}
\bH_{A,\Phi}^1 \cong \bH_{\bar\partial_{A,(\varphi,0)}}^1,
\end{equation}
where $\bH_{\bar\partial_{A,(\varphi,0)}}^1$ is defined by \eqref{eq:H_dbar_APhi^01_explicit} (or more invariantly, by \eqref{eq:H_dbar_APhi^01}) with $\psi \equiv 0$.
\end{prop}

\begin{rmk}[Isomorphism between first-order harmonic spaces for the non-Abelian monopole and holomorphic pair complexes]
\label{rmk:Difference_widehatH_dbar_APhi^01_and_H_dbar_APhi^01}
We note that the systems of defining equations for $\widehat\bH_{\bar\partial_{A,(\varphi,\psi)}}^1$ and $\bH_{\bar\partial_{A,(\varphi,\psi)}}^1$ both include equations \eqref{eq:H_dbar_APhi^01_explicit_02} and \eqref{eq:H_dbar_APhi^01_explicit_01}: the vector spaces differ when $\psi\not\equiv 0$ because the term $R_\psi^*\tau$ appears in \eqref{eq:widehatH_dbar_APhi^01_explicit_0} with a \emph{positive} sign but appears in \eqref{eq:H_dbar_APhi^01_explicit_0} with a \emph{negative} sign. By applying the identity $\bH_{\bar\partial_{A,(\varphi,0)}}^1 = \bH_{\bar\partial_A,\varphi}^1$ from \eqref{eq:Cohomologies_complexes_pre-holomorphic_pair_type1_and_holomorphic_pair_n_is_2_simplified} in Corollary \ref{cor:Vanishing_third-order_cohomology_group_holomorphic_pair_elliptic_complex_Kaehler_surface} when $\Phi=(\varphi,0)$, the conclusion \eqref{eq:Itoh_1985_proposition_2-4_SO3_monopole_complex_Kaehler_isomorphism} is strengthened to give
\begin{equation}
\label{eq:Itoh_1985_proposition_2-4_SO3_monopole_complex_Kaehler_isomorphism_simplified}
\bH_{A,\Phi}^1 \cong \bH_{\bar\partial_A,\varphi}^1,
\end{equation}
where $\bH_{\bar\partial_A,\varphi}^1$ is as in \eqref{eq:H_dbar_Avarphi^0bullet}.
\end{rmk}

The proof of Proposition \ref{prop:Itoh_1985_proposition_2-4_SO3_monopole_complex_Kaehler} is unfortunately rather long and technical. We begin by establishing a few preparatory lemmas. We first record the following generalization and extension of Wells \cite[Chapter I, Proposition 3.6, p. 33]{Wells3} from the case of $\partial$ and $\bar\partial$ on $\Omega^{\bullet,\bullet}(X)$.

\begin{lem}[Commuting complex conjugation and pointwise Hermitian adjoint with $\bar\partial_A$ and $\bar\partial_A^*$]
\label{lem:PointwiseHermitianAdjoint}
Let $E$ be a smooth Hermitian vector bundle $E$ over a closed, smooth almost complex manifold $X$. If $A$ is a smooth unitary connection on $E$ and $p,q$ are non-negative integers, then
\begin{subequations}
  \label{eq:Commute_dagger_and_del}
  \begin{align}
    \label{eq:Commute_dagger_and_del_E}
    \partial_A(\bar\psi) &= \overline{\left(\bar\partial_A\varphi\right)},
                              \quad\text{for all } \psi \in \Omega^{p,q}(E),
    \\
    \label{eq:Commute_dagger_and_del_glE}
    \partial_A(\eta^\dagger) &= \left(\bar\partial_A\eta\right)^\dagger,
                               \quad\text{for all } \eta \in \Omega^{p,q}(\gl(E)),
  \end{align}
\end{subequations}
and
\begin{subequations}
  \label{eq:AdjointCommute_dagger_and_del}
  \begin{align}
    \label{eq:AdjointCommute_dagger_and_del_E}
    \partial_A^*(\bar\psi) &= \overline{(\bar\partial_A^*\psi)},
                            \quad\text{for all } \psi \in \Omega^{p,q+1}(E),
    \\
    \label{eq:AdjointCommute_dagger_and_del_glE}
    \partial_A^*(\eta^\dagger) &= \left(\bar\partial_A^*\eta\right)^\dagger,
                            \quad\text{for all } \eta \in \Omega^{p,q+1}(\gl(E)).
  \end{align}    
\end{subequations}
\end{lem}

\begin{proof}
Recall from Wells \cite[Chapter I, Section 3, pp. 32--33]{Wells3} that complex conjugation defines a conjugate-linear isomorphism of complex vector bundles,
\begin{equation}
\label{eq:Complex_conjugation_forms}  
  Q:\Lambda^{p,q}(X) \ni \theta \mapsto \bar\theta \in \Lambda^{q,p}(X),
\end{equation}
for all non-negative integers $p, q$ and that by Wells \cite[Chapter I, Proposition 3.6, p. 33]{Wells3}, one has the identity,
\begin{equation}
\label{eq:Wells_proposition_1-3-6}  
  \partial\bar\theta = \overline{\bar\partial\theta}, \quad\text{for all } \theta \in \Omega^{p,q}(X).
\end{equation}
We first prove the identities in \eqref{eq:Commute_dagger_and_del} and since the proof of \eqref{eq:Commute_dagger_and_del_E} is very similar to that of \eqref{eq:Commute_dagger_and_del_glE}, we focus on \eqref{eq:Commute_dagger_and_del_glE}. With respect to a local trivialization of $E|_U \cong U \times \CC^r$ (where $r = \rank_\CC E$) over an open subset $U\subset X$, we have
\begin{align*}
  \bar\partial_A\psi &= \bar\partial\psi + \alpha\wedge\psi, \quad\text{for all } \psi\in\Omega^{p,q}(E),
  \\
  \bar\partial_A\eta &= \bar\partial\eta + [\alpha,\eta], \quad\text{for all } \eta\in\Omega^{p,q}(\gl(E)),
\end{align*}
for $\alpha \in \Omega^{0,1}(U,\gl(r))$ while, noting that $\alpha^\dagger \in \Omega^{1,0}(U,\gl(r))$,
\begin{align*}
  \partial_A\psi &= \partial\psi - \alpha^\dagger\wedge\psi, \quad\text{for all } \psi\in\Omega^{p,q}(E),
  \\
  \partial_A\eta &= \partial\eta - [\alpha^\dagger,\eta], \quad\text{for all } \eta\in\Omega^{p,q}(\gl(E)). 
\end{align*}
Here, we use the fact that $d_A = \partial_A+\bar\partial_A$ and $d_A\xi = d\xi + [a,\xi]$ for all $\xi \in \Omega^0(\fu(E))$, with $a = \frac{1}{2}(\alpha - \alpha^\dagger) \in \Omega^1(U,\fu(r))$. Over $U$, we thus have
\begin{align*}
  \partial_A\eta
  &= \partial\eta - [\alpha^\dagger,\eta]
    = \partial\eta + [\alpha,\eta^\dagger]^\dagger
  \\
  &= \overline{\bar\partial\bar\eta} + [\alpha,\eta^\dagger]^\dagger
    \quad\text{(by \eqref{eq:Wells_proposition_1-3-6})}
  \\
  &= \left(\overline{\bar\partial\bar\eta^\intercal}\right)^\intercal + [\alpha,\eta^\dagger]^\dagger
  \\
  &= \left(\bar\partial\eta^\dagger\right)^\dagger + [\alpha,\eta^\dagger]^\dagger
  \\
  &= \left(\bar\partial\eta^\dagger + [\alpha,\eta^\dagger]\right)^\dagger
    = \left(\bar\partial_A\eta^\dagger\right)^\dagger,
\end{align*}
and therefore, over $X$, we obtain \eqref{eq:Commute_dagger_and_del_glE} and an almost identical argument gives \eqref{eq:Commute_dagger_and_del_E}. 

Second, we prove the identities in \eqref{eq:AdjointCommute_dagger_and_del} and again focus on \eqref{eq:AdjointCommute_dagger_and_del_glE}. If $\eta \in \Omega^{p,q+1}(\gl(E))$, then $\eta^\dagger \in \Omega^{p+1,q}(\gl(E))$ and $\partial_A^*(\eta^\dagger) \in \Omega^{p,q}(\gl(E))$ and, for all $\zeta \in \Omega^{p,q}(\gl(E))$, we have
\begin{align*}
(\zeta,\partial_A^*( \eta^\dagger))_{L^2(X)}
{}&=
(\partial_A\zeta,\eta^\dagger)_{L^2(X)}
\\
{}&=
(\eta,(\partial_A\zeta)^\dagger)_{L^2(X)} \quad\text{(by the pointwise identity \eqref{eq:Inner_product_M_dagger_N_equals_inner_product_N_dagger_M})}
\\
{}&=
(\eta,\bar\partial_A(\zeta^\dagger))_{L^2(X)}
\quad\text{(by \eqref{eq:Commute_dagger_and_del_glE})}
\\
{}&=
(\bar\partial_A^*\eta,\zeta^\dagger)_{L^2(X)}
\\
{}&=
(\zeta,(\bar\partial_A^*\eta)^\dagger)_{L^2(X)}  \quad\text{(by the pointwise identity \eqref{eq:Inner_product_M_dagger_N_equals_inner_product_N_dagger_M})}.
\end{align*}
Hence, $\partial_A^*( \eta^\dagger) = (\bar\partial_A^*\eta)^\dagger$ for all $\eta  \in \Omega^{p,q+1}(\gl(E))$, since the preceding equalities hold for all $\zeta \in \Omega^{p,q}(\gl(E))$. This proves \eqref{eq:AdjointCommute_dagger_and_del_glE} and an almost identical argument gives \eqref{eq:AdjointCommute_dagger_and_del_E}. 
\end{proof}

In setting of Proposition \ref{prop:Itoh_1985_proposition_2-4_SO3_monopole_complex_Kaehler}, suppose that $(a,\phi) \in \bH_{A,\Phi}^1$ and so by \eqref{eq:H_APhi^1} we have
\[
  d_{A,\Phi}^1(a,\phi) = 0 \quad\text{and}\quad d_{A,\Phi}^{0,*}(a,\phi) = 0.
\]
We shall separately develop more explicit versions of the preceding two equations and then combine the results to prove Proposition \ref{prop:Itoh_1985_proposition_2-4_SO3_monopole_complex_Kaehler}. We begin with the

\begin{lem}[Explicit version of equation $d_{A,\Phi}^1(a,\phi)=0$]
\label{lem:Itoh_1985_2-18_SO3_monopole_complex_Kaehler}
Let $(A,\Phi)$ be a smooth solution to the unperturbed non-Abelian monopole equations \eqref{eq:SO(3)_monopole_equations_Kaehler}, where $\Phi = (\varphi,\psi) \in \Omega^0(E)\oplus\Omega^{0,2}(E)$.
Then
\[
  (a,\phi) \in \Omega^1(\su(E))\oplus\Omega^0(W_\can^+\otimes E)
\]  
obeys
\[
  d_{A,\Phi}^1(a,\phi) = 0,
\]
for $d_{A,\Phi}^1$ as in \eqref{eq:d1OfSO3MonopoleComplex} if and only if
\[
  (a',a'',\sigma,\tau)
  \in
  \Omega^{1,0}(\fsl(E))\oplus \Omega^{0,1}(\fsl(E))\oplus \Omega^0(E)\oplus\Omega^{0,2}(E),
\]  
obeys
\begin{subequations}
\label{eq:Itoh_1985_2-18_SO3_monopole_complex_Kaehler}  
\begin{align}
  \label{eq:Itoh_1985_2-18_SO3_monopole_complex_Kaehler_11_component_curvature_equation}
  \Lambda(\bar\partial_Aa' + \partial_Aa'') - i\left(\sigma\otimes\varphi^* + \varphi\otimes\sigma^*\right)_0 + i\star\left(\tau\otimes\psi^* + \psi\otimes\tau^*\right)_0 &= 0 \in \Omega^0(\su(E)),
  \\
  \label{eq:Itoh_1985_2-18_SO3_monopole_complex_Kaehler_02_component_curvature_equation}
  \bar\partial_Aa'' - \left(\tau\otimes\varphi^* + \psi\otimes\sigma^*\right)_0 &= 0 \in \Omega^{0,2}(\fsl(E)),
  \\
  \label{eq:Itoh_1985_2-18_SO3_monopole_complex_Kaehler_20_component_curvature_equation}
  \partial_Aa' + \left(\varphi\otimes\tau^* + \sigma\otimes\psi^*\right)_0 &= 0 \in \Omega^{2,0}(\fsl(E)),
  \\
  \label{eq:Itoh_1985_2-18_SO3_monopole_complex_Kaehler_Dirac_operator}  
  \bar\partial_A\sigma + \bar\partial_A^*\tau + a''\varphi - \star (a'\wedge\star\psi) &= 0 \in \Omega^{0,1}(E),
\end{align}
\end{subequations}
where $(a,\phi)$ and $(a'',\si,\tau)$ are related by \eqref{eq:Real_linear_map_tangent_spaces_SO(3)_pairs_to_preholomorphic_pairs}, with $a = \frac{1}{2}(a'+a'')$ as in\footnote{This convention is not used in standard references such as Donaldson and Kronheimer \cite{DK}, Griffiths and Harris \cite{GriffithsHarris}, Huybrechts \cite{Huybrechts_2005}, or Kobayashi \cite{Kobayashi_differential_geometry_complex_vector_bundles}.} \eqref{eq:Decompose_a_in_Omega1suE_into_10_and_01_components} and $a' = -(a'')^\dagger$ as in \eqref{eq:Kobayashi_7-6-11}.
\end{lem}

\begin{proof}
The derivation of the expression for the derivative $d_{A,\Phi}^1$ in \eqref{eq:d1OfSO3MonopoleComplex} implied by \eqref{eq:SO(3)_monopole_equations_Kaehler} at a pair
$(A,\Phi)$, where $\Phi = (\varphi,\psi) \in \Omega^0(E)\oplus\Omega^{0,2}(E)$, yields \eqref{eq:Itoh_1985_2-18_SO3_monopole_complex_Kaehler}.
In deriving
\eqref{eq:Itoh_1985_2-18_SO3_monopole_complex_Kaehler}, we used the forthcoming expression \eqref{eq:dbar_A_star} in Section \ref{subsubsec:Donaldson_Kronheimer_implied_L2_adjoint_dbar_A_and_delA} to compute the derivative of the operator $\bar\partial_A^* = -\star\circ\partial_A\circ\star$ (which arises in equation \eqref{eq:Itoh_1985_2-18_SO3_monopole_complex_Kaehler_Dirac_operator}) in the direction $a \in \Omega^1(\su(E))$, and we recall from \eqref{eq:Hermitian_duals_sections_LambdapqE} that 
\[
  \varphi^* := \langle\cdot,\varphi\rangle_E \in \Omega^0(E^*) \quad\text{and}\quad
  \psi^* := \langle\cdot,\psi\rangle_{E} \in \Omega^{2,0}(E^*),
\]
and similarly for $\sigma^* \in \Omega^0(E^*)$ and $\tau^* \in \Omega^{2,0}(E^*)$. Equation \eqref{eq:Itoh_1985_2-18_SO3_monopole_complex_Kaehler_Dirac_operator} follows by taking the derivative of the Dirac equation \eqref{eq:SO(3)_monopole_equations_Dirac_almost_Kaehler} and using the forthcoming explicit expression \eqref{eq:dbar_A_star} for $\bar\partial_A^*$. 
\end{proof}

Equations \eqref{eq:Itoh_1985_2-18_SO3_monopole_complex_Kaehler_11_component_curvature_equation} and \eqref{eq:Itoh_1985_2-18_SO3_monopole_complex_Kaehler_02_component_curvature_equation} hold when $(X,g,J)$ is almost Hermitian; Equation \eqref{eq:Itoh_1985_2-18_SO3_monopole_complex_Kaehler_Dirac_operator} holds when $(X,g,J)$ is almost K\"ahler but must be modified when it is almost Hermitian. Note that $a' \in \Omega^{1,0}(\fsl(E))$ and $\star\psi \in \Omega^{0,2}(E)$, so $a'\wedge\star\psi \in \Omega^{1,2}(E)$ and thus $\star(a'\wedge\star\psi) \in \Omega^{0,1}(E)$ in \eqref{eq:Itoh_1985_2-18_SO3_monopole_complex_Kaehler_Dirac_operator}, as desired. By using identity $a' = -(a'')^\dagger$ from \eqref{eq:Kobayashi_7-6-11} and the complex conjugation identities in  Lemma \ref{lem:PointwiseHermitianAdjoint}, we see that the equation \eqref{eq:Itoh_1985_2-18_SO3_monopole_complex_Kaehler_20_component_curvature_equation} is the complex conjugate transpose of \eqref{eq:Itoh_1985_2-18_SO3_monopole_complex_Kaehler_02_component_curvature_equation} and so can be omitted from the system \eqref{eq:Itoh_1985_2-18_SO3_monopole_complex_Kaehler} without loss of generality.

\begin{rmk}[Simplification of Equation \eqref{eq:Itoh_1985_2-18_SO3_monopole_complex_Kaehler} when Witten's dichotomy holds]
\label{rmk:Simplification_SO3_monopole_equations_complex_Kaehler_manifolds_Witten_dichotomy}  
When $\psi \equiv 0$ or $\varphi \equiv 0$, the term $(\psi\otimes\varphi^*)_0$ also vanishes and so the system \eqref{eq:Itoh_1985_2-18_SO3_monopole_complex_Kaehler} simplifies. For example, when $\psi=0$, then
\begin{align*}
  \Lambda(\bar\partial_Aa' + \partial_Aa'') - i\left(\sigma\otimes\varphi^* + \varphi\otimes\sigma^*\right)_0&= 0 \in \Omega^0(\su(E)),
  \\
  \bar\partial_Aa'' &= 0 \in \Omega^{0,2}(\fsl(E)),
  \\
  \bar\partial_A\sigma + a''\varphi &= 0 \in \Omega^{0,1}(E),
\end{align*}
where we applied Lemma \ref{lem:Witten_dichotomy_linearized_type1_SO3_monopole_equation} to conclude that we also have $\tau = 0$.
Similarly, when $\varphi=0$, then\footnote{Note that these equations are not complex linear in $(a'',\tau)$.}
\begin{align*}
  \Lambda(\bar\partial_Aa' + \partial_Aa'') + i\star\left(\tau\otimes\psi^* + \psi\otimes\tau^*\right)_0 &= 0 \in \Omega^0(\su(E)),
  \\
  \bar\partial_Aa'' &= 0 \in \Omega^{0,2}(\fsl(E)),
  \\
  \bar\partial_A^*\tau - \star (a'\wedge\star\psi) &= 0 \in \Omega^{0,1}(E),
\end{align*}
where we applied Lemma \ref{lem:Witten_dichotomy_linearized_type1_SO3_monopole_equation} to conclude that we also have $\sigma = 0$.

We may use the replacement $a' = -(a'')^\dagger$ in the last equation. We recall from Lemma \ref{lem:Okonek_Teleman_1995_3-1} that if $(X,g,J)$ is complex K\"ahler, then the analogue of Witten's dichotomy for non-Abelian monopoles, namely $\psi \equiv 0$ or $\varphi \equiv 0$, does hold. 
\end{rmk}

\begin{rmk}[Complex linearity of equation \eqref{eq:Itoh_1985_2-18_SO3_monopole_complex_Kaehler}]
\label{rmk:Complex_linearity_SO3_monopole_equations_almost_complex_Hermitian_four-manifolds}
While we shall see that the real equation \eqref{eq:Itoh_1985_2-18_SO3_monopole_complex_Kaehler_11_component_curvature_equation} will combine with the real Coulomb gauge condition \eqref{eq:SO3_monopole_Coulomb_gauge_slice_condition} to yield the dramatically simpler and complex linear equation \eqref{eq:dbar_A_slice_condition_SO3_monopoles_complex_Kaehler}, a simplification that continues to hold even when $(X,g,J)$ is only almost K\"ahler, neither equations \eqref{eq:Itoh_1985_2-18_SO3_monopole_complex_Kaehler_02_component_curvature_equation} nor \eqref{eq:Itoh_1985_2-18_SO3_monopole_complex_Kaehler_Dirac_operator} will be complex linear in $(a'',\tau,\sigma)$ when $\psi\not\equiv 0$, due to the presence of $\sigma^*$ in \eqref{eq:Itoh_1985_2-18_SO3_monopole_complex_Kaehler_02_component_curvature_equation} and $a'=-(a'')^\dagger$ in \eqref{eq:Itoh_1985_2-18_SO3_monopole_complex_Kaehler_Dirac_operator}.
\end{rmk}  

We digress to discuss our convention regarding adjoints of operators between real and complex inner product spaces. Let $\fH$ be an inner product space over $\RR$ and $\fK$ be an inner product space over $\CC$. If $S \in \Hom_\RR(\fK,\fH)$ and $T \in \Hom_\RR(\fH,\fK)$, then we define $S^* \in \Hom_\RR(\fH,\fK)$ and $T^* \in \Hom_\RR(\fK,\fH)$ by the identities
\begin{equation}
  \label{eq:Adjoint_operator_real_to_complex_inner_product_space}
  \begin{aligned}
    \Real\langle S^*h, k\rangle_\fK &= \langle h,Sk\rangle_\fH, 
    \\
    \langle T^*k, h\rangle_\fH &= \Real\langle k,Th\rangle_\fK, \quad\text{for all } h \in \fH, k \in \fK.
\end{aligned}
\end{equation}
(Note that if $\fK_\RR$ is the real vector space underlying $\fK$, then according to the convention of Kobayashi \cite[Equation (7.6.5), p. 251]{Kobayashi_differential_geometry_complex_vector_bundles}, its inner product is given by 
\begin{equation}
  \label{eq:Kobayashi_7-6-5_abstract_inner_product_space}
  \langle k_1, k_2 \rangle_{\fK_\RR} = \langle k_1, k_2 \rangle_\fK + \langle k_2, k_1 \rangle_\fK
  = 2\Real\langle k_1, k_2 \rangle_\fK,
\end{equation}
but the choice of whether or not to include the factor $2$ is just a convention.) When $\fH$ and $\fK$ are both real or both complex inner product spaces, then we define adjoints in the usual way and this completes our digression.

\begin{lem}[A second version of the Coulomb gauge slice condition]
\label{lem:SO3_monopole_Coulomb_gauge_slice_condition}
Let $(\rho,W\otimes E)$ be a spin${}^u$ structure over a closed, four-dimensional, oriented, smooth Riemannian manifold $(X,g)$. If $(A,\Phi) \in \sA(E,h) \times \Omega^0(W^+\otimes E)$ is a smooth pair, then
\begin{equation}
\label{eq:d_APhi^0_star_Expression}
d_{A,\Phi}^{0,*}(a,\phi) = d_A^*a - \sR_\Phi^*\phi,
\quad\text{for all } (a,\phi) \in \Omega^1(\su(E))\oplus\Omega^0(W_\can^+\otimes E),
\end{equation}
so the Coulomb gauge slice condition $d_{A,\Phi}^{0,*}(a,\phi) = 0$ is equivalent to
\begin{equation}
  \label{eq:d_APhi^0_star_aphi_identity_and_vanishing}
  d_{A,\Phi}^{0,*}(a,\phi) = d_A^*a - \sR_\Phi^*\phi = 0,
\end{equation}
where the real linear operator
\begin{equation}
  \label{eq:Phi_star_Omega0V+_to_Omega0suE}
  \sR_\Phi^*:\Omega^0(W_\can^+\otimes E) \to \Omega^0(\su(E))
\end{equation}
is the pointwise adjoint (with respect to the Riemannian metric on $\su(E)$ and the real part of the Hermitian metric on $W_\can^+\otimes E$ as suggested by \eqref{eq:Adjoint_operator_real_to_complex_inner_product_space}) of the real linear, right composition operator,
\begin{equation}
  \label{eq:Phi_Omega0suE_to_Omega0V+}
  \sR_\Phi:\Omega^0(\su(E)) \ni \xi \mapsto \xi\Phi \in \Omega^0(W_\can^+\otimes E).
\end{equation}
\end{lem}  

\begin{proof}
The Coulomb gauge slice condition, implied by the expression \eqref{eq:d_APhi^0} for the differential,
\[
d_{A,\Phi}^0\xi = (d_A\xi, -\xi\Phi), \quad\text{for all } \xi \in \Omega^0(\su(E)),
\]
based at the smooth pair $(A,\Phi) \in \sA(E,h)\times\Omega^0(W_\can^+\otimes E)$ in the direction $(a,\phi) \in \Omega^1(\su(E))\oplus \Omega^0(W_\can^+\otimes E)$, is given by
\begin{equation}
  \label{eq:SO3_monopole_Coulomb_gauge_slice_condition}
  d_{A,\Phi}^{0,*}(a,\phi) = 0. 
\end{equation}
The unmodified $L^2$ inner product on the codomain $\Omega^1(\su(E)) \oplus \Omega^0(W_\can^+\otimes E)$ of $d_{A,\Phi}^0$ is given by
\[
  (d_{A,\Phi}^0\xi,(a,\phi))_{L^2(X)} = ((d_A\xi, -\xi\Phi),(a,\phi))_{L^2(X)}
  := \underbrace{(d_A\xi, a)_{L^2(X)}}_{\textrm{real}} - \underbrace{(\xi\Phi,\phi)_{L^2(X)}}_{\textrm{complex}}. 
\]
while the unmodified $L^2$ inner product 
\[
  \underbrace{(\xi, d_{A,\Phi}^{0,*}(a,\phi))_{L^2(X)}}_{\textrm{real}}
\]
on the domain $\Omega^1(\su(E))$ of $d_{A,\Phi}^0$ is real. The $L^2$ adjoint operator $d_{A,\Phi}^{0,*}$ in \eqref{eq:SO3_monopole_Coulomb_gauge_slice_condition} acting on $(a,\phi) \in \Omega^1(\su(E)) \oplus \Omega^0(W_\can^+\otimes E)$ is therefore defined by the relation (compare with \eqref{eq:Adjoint_operator_real_to_complex_inner_product_space})
  \begin{equation}
  \label{eq:Definition_d0APhi_star}
  (\xi, d_{A,\Phi}^{0,*}(a,\phi))_{L^2(X)}
  =
  (d_A\xi, a)_{L^2(X)} - \Real(\xi\Phi,\phi)_{L^2(X)}, \quad\text{for all } \xi \in \Omega^0(\su(E)).
\end{equation}
Equation \eqref{eq:SO3_monopole_Coulomb_gauge_slice_condition} is thus equivalent to the assertion that for all $\xi\in\Omega^0(\su(E))$,
\begin{align*}
  0 &= (\xi,d_{A,\Phi}^{0,*}(a,\phi))_{L^2(X)}
  \\
  &= (d_A\xi,a)_{L^2(X)} - \Real(\xi\Phi,\phi)_{L^2(X)} \quad\text{(by \eqref{eq:Definition_d0APhi_star})}
  \\
  &= (d_A\xi,a)_{L^2(X)} - \Real(\sR_\Phi\xi,\phi)_{L^2(X)} \quad\text{(by \eqref{eq:Phi_Omega0suE_to_Omega0V+})}
  \\
  &= (\xi, d_A^*a)_{L^2(X)} - (\xi,\sR_\Phi^*\phi)_{L^2(X)} \quad \text{(by \eqref{eq:Adjoint_operator_real_to_complex_inner_product_space})}
  \\
  &= (\xi, d_A^*a - \sR_\Phi^*\phi)_{L^2(X)},
\end{align*}
and this yields \eqref{eq:d_APhi^0_star_Expression} and \eqref{eq:d_APhi^0_star_aphi_identity_and_vanishing}.
\end{proof}

We shall need a more explicit version of equation $d_A^*a - \sR_\Phi^*\phi = 0$ in \eqref{eq:d_APhi^0_star_aphi_identity_and_vanishing}. For this purpose, we define the complex linear operator
\begin{equation}
  \label{eq:varphi_psi_star_Omega0E_oplus_Omega02E_to_Omega0slE}
  \Omega^0(E)\oplus \Omega^{0,2}(E) \ni (\sigma,\tau) \mapsto R_{(\varphi,\psi)}^*(\sigma,\tau) := R_\varphi^*\sigma + R_\psi^*\tau \in \Omega^0(\fsl(E))
\end{equation}
to be the pointwise adjoint (with respect to the Hermitian metrics on $E$ and $\gl(E)$) of the complex linear, right composition operator,
\begin{equation}
  \label{eq:varphi_psi_Omega0slE_to_Omega0E_oplus_Omega02E}
  \Omega^0(\fsl(E)) \ni \zeta \mapsto R_{(\varphi,\psi)}\zeta = \zeta(\varphi,\psi) = (\zeta\varphi,\zeta\psi) = (R_\varphi\zeta, R_\psi\zeta) \in \Omega^0(E)\oplus \Omega^{0,2}(E),
\end{equation}
where $R_\varphi$ and $R_\psi$ are the complex linear, right composition operators,
\begin{subequations}
  \label{eq:Right_multiplication_of_section_slE_by_sections_E_or_02E}  
  \begin{align}
  \label{eq:Right_multiplication_of_section_slE_by_section_E}  
  R_\varphi:\Omega^0(\fsl(E)) \ni \zeta &\mapsto \zeta\varphi \in \Omega^0(E),
    \\
  \label{eq:Right_multiplication_of_section_slE_by_section_02E}  
  R_\psi:\Omega^0(\fsl(E)) \ni \zeta &\mapsto \zeta\psi \in \Omega^{0,2}(E),
\end{align}
\end{subequations}
with pointwise adjoints (with respect to the Hermitian metrics on $E$ and $\gl(E)$) given by the complex linear operators,
\begin{align*}
  R_\varphi^*:\Omega^0(E) \to \Omega^0(\fsl(E)),
  \\
  R_\psi^*:\Omega^{0,2}(E) \to \Omega^0(\fsl(E)).
\end{align*}
We note that there is a pointwise-orthogonal, direct-sum decomposition given by
\begin{equation}
  \label{eq:slE_equals_suE_oplus_isuE}
  \fsl(E) \ni \zeta \mapsto \xi_1 + i\xi_2 \in \su(E)\oplus i\su(E)
\end{equation}
defined by $\xi_1=\pi_{\su(E)}\zeta$ and $i\xi_2 = \pi_{\su(E)}^\perp\zeta$, where $\xi_1,\xi_2 \in \su(E)$, using the real-linear, pointwise-orthogonal projections
\begin{equation}
  \label{eq:Pointwise_orthogonal_projections_slE_onto_suE_and_isuE}
  \begin{aligned}
    \pi_{\su(E)}: \fsl(E) \ni \zeta &\mapsto \frac{1}{2}(\zeta-\zeta^\dagger) \in \su(E),
    \\
    \pi_{\su(E)}^\perp = \id_{\fsl(E)} - \pi_{\su(E)} = \pi_{i\su(E)}: \fsl(E) \ni \zeta
    &\mapsto \frac{1}{2}(\zeta+\zeta^\dagger) \in i\su(E).
  \end{aligned}
\end{equation}  
We then obtain
\begin{equation}
\label{eq:su2_Projection_of_AdjointOf_sl2_MultOperator_Is_su2_MultOperator}
  \sR_\Phi^* = \pi_{\su(E)}R_{(\varphi,\psi)}^*:   \Omega^0(E)\oplus \Omega^{0,2}(E) \to \Omega^0(\su(E)).
\end{equation}
Indeed, by definition \eqref{eq:Phi_Omega0suE_to_Omega0V+} of $\sR_\Phi$ and definition \eqref{eq:varphi_psi_Omega0slE_to_Omega0E_oplus_Omega02E} of $R_{(\varphi,\psi)}$ and for all $\pi_{\su(E)}\zeta \in \Omega^0(\su(E))$ defined by $\zeta \in \Omega^0(\fsl(E))$ and writing $\Phi = (\varphi,\psi)$, we have
\[
  \sR_\Phi\pi_{\su(E)}\zeta = (\pi_{\su(E)}\zeta)\Phi = (\pi_{\su(E)}\zeta)(\varphi,\psi) = R_{(\varphi,\psi)}\pi_{\su(E)}\zeta 
\]
and thus, using $\pi_{\su(E)}^*=\pi_{\su(E)}$ and $\pi_{\su(E)}\sR_\Phi^*=\sR_\Phi^*$, we have
$\sR_\Phi^*=(R_{(\varphi,\psi)}\pi_{\su(E)})^* = \pi_{\su(E)}^*R_{(\varphi,\psi)}^* = \pi_{\su(E)}R_{(\varphi,\psi)}^*$, as claimed.

By definition of $\pi_{\su(E)}$ and because $\Phi = (\varphi,\psi)$, we see that
\begin{align*}
  \sR_\Phi^*(\sigma,\tau)
  &= \pi_{\su(E)}R_{(\varphi,\psi)}^*
    \quad\text{(by \eqref{eq:su2_Projection_of_AdjointOf_sl2_MultOperator_Is_su2_MultOperator})}
  \\
  &= \frac{1}{2}\left(R_{(\varphi,\psi)}^*(\sigma,\tau)
    - \left(R_{(\varphi,\psi)}^*(\sigma,\tau)\right)^\dagger\right).
\end{align*}
Moreover, for all $\zeta \in \Omega^0(\fsl(E))$ we obtain
\begin{multline*}
  \langle R_{(\varphi,\psi)}^*(\sigma,\tau), \zeta \rangle_{\fsl(E)}
  = \langle (\sigma,\tau), R_{(\varphi,\psi)}\zeta \rangle_{E\oplus\Lambda^{0,2}(E)}
  \\
  = \langle \sigma, R_\varphi\zeta \rangle_{E} + \langle \tau, R_\psi\zeta \rangle_{\Lambda^{0,2}(E)}
  = \langle R_\varphi^*\sigma, \zeta \rangle_{\fsl(E)} + \langle R_\psi^*\tau, \zeta \rangle_{\fsl(E)},
\end{multline*}  
and so
\[
  R_{(\varphi,\psi)}^*(\sigma,\tau) = R_\varphi^*\sigma + R_\psi^*\tau.
\]
We can therefore compute
\begin{align*}
  \sR_\Phi^*(\sigma,\tau)
  &= \frac{1}{2}\left(R_{(\varphi,\psi)}^*(\sigma,\tau)
    - \left(R_{(\varphi,\psi)}^*(\sigma,\tau)\right)^\dagger\right)
  \\
  &= \frac{1}{2}\left(R_\varphi^*\sigma + R_\psi^*\tau
    - \left(R_\varphi^*\sigma + R_\psi^*\tau\right)^\dagger\right),
\end{align*}
and thus
\begin{equation}
  \label{eq:sR_Phi^*_sigma_tau}
  \sR_\Phi^*(\sigma,\tau)
  =
  \frac{1}{2}\left(R_\varphi^*\sigma - (R_\varphi^*\sigma)^\dagger
    + R_\psi^*\tau - (R_\psi^*\tau)^\dagger\right).
\end{equation}
We are now ready to prove the

\begin{lem}[A third version of the Coulomb gauge slice condition]
\label{lem:Real_Coulomb_slice_SO3_monopole_complex_Kaehler_del_A_and_del-bar_A_adjoints}
Continue the hypotheses of Lemma \ref{lem:SO3_monopole_Coulomb_gauge_slice_condition}. Then the Coulomb gauge slice condition \eqref{eq:d_APhi^0_star_aphi_identity_and_vanishing} is equivalent to 
\begin{equation}
\label{eq:Real_Coulomb_slice_SO3_monopole_complex_Kaehler_del_A_and_del-bar_A_adjoints}
\partial_A^*a' + \bar\partial_A^*a''
- \left(R_\varphi^*\sigma - (R_\varphi^*\sigma)^\dagger
  + R_\psi^*\tau - (R_\psi^*\tau)^\dagger\right) = 0,
\end{equation}
where $\Phi = (\varphi,\psi) \in \Omega^0(E)\oplus\Omega^{0,2}(E)$ and $a=\frac{1}{2}(a'+a'')$ from \eqref{eq:Decompose_a_in_Omega1suE_into_10_and_01_components} and $\phi = (\sigma,\tau) \in \Omega^1(\su(E))\oplus\Omega^0(E)\oplus\Omega^{0,2}(E)$.
\end{lem}

\begin{proof}
By substituting the identity $d_A = \partial_A+\bar\partial_A$ on $\Omega^0(\fsl(E))$ and the expression \eqref{eq:sR_Phi^*_sigma_tau} for $R_{(\varphi,\psi)}^*$ and $a=\frac{1}{2}(a'+a'')$ from \eqref{eq:Decompose_a_in_Omega1suE_into_10_and_01_components} and removing overall scaling factor of $1/2$, we will show that the Coulomb gauge slice condition \eqref{eq:d_APhi^0_star_aphi_identity_and_vanishing} becomes \eqref{eq:Real_Coulomb_slice_SO3_monopole_complex_Kaehler_del_A_and_del-bar_A_adjoints}. We first claim that 
\begin{equation}
\label{eq:d_star_equals_partial_stars}
d_A^*a =\frac{1}{2}\left( \partial_A^*a' + \bar\partial_A^*a''\right), \quad\text{for all } a=\frac{1}{2}(a'+a'')\in\Om^1(\su(E)).
\end{equation}
Indeed, writing $\xi = \pi_{\su(E)}\zeta = \frac{1}{2}(\zeta-\zeta^\dagger)$ by \eqref{eq:Pointwise_orthogonal_projections_slE_onto_suE_and_isuE},
\begin{align*}
  (d_A^*a,\xi)_{L^2(X)}
  &= (a,d_A\xi)_{L^2(X)}
  \\
  &= \frac{1}{2}(a',\partial_A\xi)_{L^2(X)} + \frac{1}{2}(a'',\bar\partial_A\xi)_{L^2(X)}
  \\
  &= \frac{1}{4}(a',\partial_A\zeta)_{L^2(X)} - \frac{1}{4}(a',\partial_A(\zeta^\dagger))_{L^2(X)}
    + \frac{1}{4}(a'',\bar\partial_A\zeta)_{L^2(X)} - \frac{1}{4}(a'',\bar\partial_A(\zeta^\dagger))_{L^2(X)}
  \\
  &= \frac{1}{4}(a',\partial_A\zeta)_{L^2(X)} - \frac{1}{4}(\bar\partial_A\zeta,(a')^\dagger)_{L^2(X)}
    + \frac{1}{4}(a'',\bar\partial_A\zeta)_{L^2(X)} - \frac{1}{4}(\partial_A\zeta, (a'')^\dagger)_{L^2(X)}
  \\
  &= \frac{1}{4}(a',\partial_A\zeta)_{L^2(X)} + \frac{1}{4}(\bar\partial_A\zeta,a'')_{L^2(X)}
    + \frac{1}{4}(a'',\bar\partial_A\zeta)_{L^2(X)} + \frac{1}{4}(\partial_A\zeta, a')_{L^2(X)}
  \\
  &=\frac{1}{2} \Real(a',\partial_A\zeta)_{L^2(X)} + \frac{1}{2}\Real(a'',\bar\partial_A\zeta)_{L^2(X)}
  \\
  &= \frac{1}{2}\Real(\partial_A^*a',\zeta)_{L^2(X)} + \frac{1}{2}\Real(\bar\partial_A^*a'',\zeta)_{L^2(X)}
  \\
  &= \frac{1}{2}\Real((\partial_A^*a'+\bar\partial_A^*a''),\zeta)_{L^2(X)}.
\end{align*}
But $(\partial_A^*a'+\bar\partial_A^*a'')^\dagger = -\bar\partial_A^*a''-\partial_A^*a'$
by \eqref{eq:Kobayashi_7-6-11} and \eqref{eq:AdjointCommute_dagger_and_del_glE}
and thus $\partial_A^*a'+\bar\partial_A^*a'' \in \Omega^0(\su(E))$ and
\begin{align*}
  (d_A^*a,\xi)_{L^2(X)}
  &= \frac{1}{2}\Real((\partial_A^*a'+\bar\partial_A^*a''),\zeta)_{L^2(X)}
  \\
  &= \frac{1}{4}\left(\Real((\partial_A^*a'+\bar\partial_A^*a''),\zeta)_{L^2(X)}
    + \Real((\partial_A^*a'+\bar\partial_A^*a''),\zeta)_{L^2(X)}\right)
  \\
  &= \frac{1}{4}\left(\Real((\partial_A^*a'+\bar\partial_A^*a''),\zeta)_{L^2(X)}
    + \Real(\zeta^\dagger, (\partial_A^*a'+\bar\partial_A^*a'')^\dagger)_{L^2(X)}\right)
  \\
  &\qquad\text{(by forthcoming \eqref{eq:Inner_product_M_dagger_N_equals_inner_product_N_dagger_M})}
  \\
  &= \frac{1}{4}\left(\Real((\partial_A^*a'+\bar\partial_A^*a''),\zeta)_{L^2(X)}
    - \Real(\zeta^\dagger, \partial_A^*a'+\bar\partial_A^*a'')_{L^2(X)}\right)
  \\
  &\qquad\text{(by \eqref{eq:Kobayashi_7-6-11} and \eqref{eq:AdjointCommute_dagger_and_del_glE})}.
\end{align*}
Therefore,
\begin{align*}
  (d_A^*a,\xi)_{L^2(X)}
  &= \frac{1}{4}\Real((\partial_A^*a'+\bar\partial_A^*a''),\zeta-\zeta^\dagger)_{L^2(X)}
  \\
  &= \frac{1}{2}((\partial_A^*a'+\bar\partial_A^*a''),\xi)_{L^2(X)}, \quad\text{for all } \xi \in \Omega^0(\su(E)).
\end{align*}
This completes the proof of \eqref{eq:d_star_equals_partial_stars}. Substituting the equalities \eqref{eq:d_star_equals_partial_stars} and \eqref{eq:sR_Phi^*_sigma_tau} into \eqref{eq:d_APhi^0_star_aphi_identity_and_vanishing} and removing the common factor of $1/2$ 
yields \eqref{eq:Real_Coulomb_slice_SO3_monopole_complex_Kaehler_del_A_and_del-bar_A_adjoints}, completing the proof of Lemma \ref{lem:Real_Coulomb_slice_SO3_monopole_complex_Kaehler_del_A_and_del-bar_A_adjoints}.
\end{proof}

We make the

\begin{claim}
\label{claim:R_varphi_star_and_R_psi_star_tracefree}  
\begin{subequations}
\label{eq:R_varphi_star_and_R_psi_star_tracefree}  
\begin{align}
\label{eq:R_varphi_star_sigma_is_sigma_tensor_varphi_star_tracefree}  
  R_\varphi^*\sigma &= (\sigma\otimes\varphi^*)_0 \in \Omega^0(\fsl(E)),
  \\
  \label{eq:R_psi_star_tau_is_tau_tensor_psi_star_tracefree}
  R_\psi^*\tau &= \star\left(\tau\wedge\psi^*\right)_0 \in \Omega^0(\fsl(E)),
  \\
  \label{eq:R_varphi_star_sigma_dagger_is_varphi_tensor_sigma_star_tracefree}  
  (R_\varphi^*\sigma)^\dagger &= (\varphi\otimes\sigma^*)_0  \in \Omega^0(\fsl(E)),
  \\
  \label{eq:R_psi_star_tau_dagger_is_psi_tensor_tau_star_tracefree}
  (R_\psi^*\tau)^\dagger &= \star(\psi\wedge\tau^*)_0  \in \Omega^0(\fsl(E)).
\end{align}
\end{subequations}
\end{claim}

(Note that $\tau\otimes\psi^* = \tau\wedge\psi^*$ and $\psi\otimes\tau^* = \psi\wedge\tau^*$ since wedge product on two-forms is commutative.) Given Claim \ref{claim:R_varphi_star_and_R_psi_star_tracefree}, we can prove the

\begin{lem}[A fourth and final version of the Coulomb gauge slice condition]
\label{lem:Real_Coulomb_slice_SO3_monopole_complex_Kaehler_del_A_and_del-bar_A_adjoints_tensor_products}
Continue the hypotheses of Lemma \ref{lem:Real_Coulomb_slice_SO3_monopole_complex_Kaehler_del_A_and_del-bar_A_adjoints}. Then the Coulomb gauge slice condition \eqref{eq:Real_Coulomb_slice_SO3_monopole_complex_Kaehler_del_A_and_del-bar_A_adjoints}
is equivalent to  
\begin{equation}
\label{eq:Real_Coulomb_slice_SO3_monopole_complex_Kaehler_del_A_and_del-bar_A_adjoints_tensor_products}
  \partial_A^*a' + \bar\partial_A^*a'' - \left(\sigma\otimes\varphi^* - \varphi\otimes\sigma^*\right)_0 - \star\left(\tau\otimes\psi^* - \psi\otimes\tau^*\right)_0 = 0.
\end{equation}
\end{lem}

\begin{proof}[Proof of Lemma \ref{lem:Real_Coulomb_slice_SO3_monopole_complex_Kaehler_del_A_and_del-bar_A_adjoints_tensor_products}]
Substituting the identities \eqref{eq:R_varphi_star_and_R_psi_star_tracefree} into the Coulomb gauge slice condition \eqref{eq:Real_Coulomb_slice_SO3_monopole_complex_Kaehler_del_A_and_del-bar_A_adjoints} yields \eqref{eq:Real_Coulomb_slice_SO3_monopole_complex_Kaehler_del_A_and_del-bar_A_adjoints_tensor_products}.
\end{proof}  

\begin{proof}[Proof of Claim \ref{claim:R_varphi_star_and_R_psi_star_tracefree}]
We now prove the claimed equations \eqref{eq:R_varphi_star_and_R_psi_star_tracefree}, beginning with the claim \eqref{eq:R_varphi_star_sigma_is_sigma_tensor_varphi_star_tracefree}. For all $\zeta\in\Omega^0(\fsl(E))$ and recalling that $\gl(E)=\fsl(E)\oplus \CC\,\id_E$ is a pointwise orthogonal direct sum, we have
\[
  \langle (\sigma\otimes\varphi^*)_0,\zeta\rangle_{\fsl(E)} = \langle \sigma\otimes\varphi^*,\zeta\rangle_{\gl(E)}.
\]
Denote $r = \rank_\CC E$. For a local orthonormal frame $\{e_i\}$ for $E$ over an open subset $U\subset X$ and corresponding local orthonormal frame $\{e_i^*\}$ for $E^*$ over $U$ given by $e_i^* = \langle\cdot,e_i\rangle_E$, we write
\begin{multline*}
  \zeta = \sum_{i,j=1}^r \zeta_{ij} e_i\otimes e_j^* \in \Omega^0(U,\fsl(E)),
  \\
  \sigma = \sum_{i=1}^r\sigma_i\otimes e_i
  \quad\text{and}\quad \varphi = \sum_{i=1}^r\varphi_i\otimes e_i \in \Omega^0(U,E),
  \\
  \varphi^* = \langle\cdot,\varphi\rangle_E = \sum_{i=1}^r\bar\varphi_i\otimes e_i^* \in \Omega^0(U,E^*),
\end{multline*}
appealing to \eqref{eq:Complex_conjugate_sections_E_and_Lambda02E} and \eqref{eq:Hermitian_duals_sections_LambdapqE} for the expression for $\varphi^*$. We observe that
\begin{align*}
  \langle\sigma\otimes\varphi^*, \zeta\rangle_{\gl(E)}
  &= \sum_{i,j,k,l} \langle\sigma_k\bar\varphi_l\, e_k\otimes e_l^*, \zeta_{ij}\,e_i\otimes e_j^*\rangle_{\gl(E)}
  \\
  &= \sum_{i,j,k,l} \langle\sigma_k\bar\varphi_l\, e_k\otimes e_l^*, \zeta_{ij}\,e_i\otimes e_j^*\rangle_{E\otimes E^*}
  \\
  &= \sum_{i,j,k,l} \langle\sigma_k e_k\otimes e_l^*, \zeta_{ij}\varphi_l\,e_i\otimes e_j^*\rangle_{E\otimes E^*}
  \\
  &= \sum_{i,j,k,l} \langle\sigma_k e_k, \zeta_{ij}\varphi_l\,e_i\rangle_E \langle e_l^*,e_j^* \rangle_{E^*}
  \\
  &= \sum_{i,j,k,l} \langle\sigma_k e_k, \zeta_{ij}\varphi_l\,e_i\rangle_E \langle e_j,e_l \rangle_E
\quad\text{(by \eqref{eq:Inner_product_dual_Hilbert_space})}
  \\
  &= \sum_{i,j,k,l} \langle\sigma_k e_k, \zeta_{ij}\varphi_l\,e_i\rangle_E \overline{\langle e_l,e_j \rangle}_E  
  \\
  &= \sum_{i,j,k,l} \langle\sigma_k e_k, \zeta_{ij}\varphi_l\,e_i \langle e_l,e_j \rangle_E \rangle_E
  \\
  &= \sum_{i,j,k,l} \langle\sigma_k e_k, \zeta_{ij}\,e_i \otimes e_j^*(\varphi_le_l) \rangle_E
  = \langle\sigma, \zeta\varphi\rangle_E.
\end{align*}
Consequently, we obtain
\begin{equation}
  \label{eq:Inner_product_sigma_otimes_bar_varphi_star_and_zeta_equals_sigma_and_zeta_varphi}
  \langle\sigma\otimes\varphi^*, \zeta\rangle_{\gl(E)}
  = \langle\sigma, \zeta\varphi\rangle_E
  = \langle\sigma, R_\varphi\zeta\rangle_E
  = \langle R_\varphi^*\sigma, \zeta\rangle_E,
\end{equation}
and so the claim \eqref{eq:R_varphi_star_sigma_is_sigma_tensor_varphi_star_tracefree} follows.

Similarly, by abbreviating $\Lambda^{0,2}(E) = \Lambda^{0,2}(X)\otimes E$ and 
$\Lambda^{2,2}(\fsl(E)) = \Lambda^{2,2}(X)\otimes \fsl(E)$ as usual, we see that for all $\zeta \in \Omega^0(\fsl(E))$ we have
\begin{align*}
 \langle R_\psi^*\tau,\zeta\rangle_{\fsl(E)} 
  &= \langle \tau,R_\psi\zeta\rangle_{\Lambda^{0,2}(E)}
  \\
  &= \langle \tau, \zeta\psi\rangle_{\Lambda^{0,2}(E)} 
  \\
  &= \langle\tau\wedge\psi^*,\star\zeta\rangle_{\Lambda^{2,2}(\gl(E))} \quad\text{(by the forthcoming \eqref{eq:Inner_product_tau_and_zeta_psi_equals_tau_otimes_bar_psi_star_and_star_zeta})}
  \\
  &= \langle(\tau\wedge\psi^*)_0,\star\zeta\rangle_{\Lambda^{2,2}(\fsl(E))}
  \\
  &= \langle\star(\tau\wedge\psi^*)_0,\zeta\rangle_{\fsl(E)},
\end{align*}
and so the claim \eqref{eq:R_psi_star_tau_is_tau_tensor_psi_star_tracefree} follows, modulo the forthcoming \eqref{eq:Inner_product_tau_and_zeta_psi_equals_tau_otimes_bar_psi_star_and_star_zeta}. In the preceding derivation, we have used the fact that, for elementary tensors $\tau = \alpha\otimes s$ and $\psi = \beta\otimes t$ for $\alpha,\beta \in \Omega^{0,2}(X)$ and $s,t \in \Omega^0(E)$,
\begin{align*}
  \langle \tau, \zeta\psi\rangle_{\Lambda^{0,2}(E)}\vol
  &= \langle \alpha\otimes s, \zeta(\beta\otimes t)\rangle_{\Lambda^{0,2}(E)}\vol
  \\
  &= \langle s, \zeta t \rangle_E\, \langle \alpha, \beta \rangle_{\Lambda^{0,2}(X)}\vol
  \\
  &= \langle s, \zeta t \rangle_E\, \alpha\wedge\star\bar\beta \quad\text{(by Huybrechts \cite[Section 1.2, p. 33 and Definition 3.2.1, p. 125]{Huybrechts_2005})}
  \\
  &= \langle s, \zeta t \rangle_E\, \alpha\wedge\bar\beta \quad\text{(since $\Omega^{2,0}(X)\subset\Omega^+(X,\CC)$)}
  \\
  &= \langle s, \zeta t \rangle_E\, \langle \alpha\wedge\bar\beta, \star 1 \rangle_{\Lambda^{2,2}(X)}\vol
  \\
  &= \langle s\otimes t^*, \zeta \rangle_{\gl(E)}\, \langle \alpha\wedge\bar\beta, \star 1 \rangle_{\Lambda^{2,2}(X)}\vol
    \quad\text{(by \eqref{eq:Inner_product_sigma_otimes_bar_varphi_star_and_zeta_equals_sigma_and_zeta_varphi})}
  \\
  &= \langle (\alpha \otimes s)\wedge (\bar\beta\otimes t^*), \star\zeta \rangle_{\Lambda^{2,2}(\gl(E))}\vol
  \\
  &= \langle (\alpha \wedge \bar\beta)\otimes (s\otimes t^*), \star\zeta \rangle_{\Lambda^{2,2}(\gl(E))}\vol, 
\end{align*}
where\footnote{We use the convention in Kobayashi \cite[Example 1.3.2, p. 28]{Kobayashi_Nomizu_v1} and Warner \cite[Equations (2.10.1) and (2.10.6)]{Warner} that gives $\delta\wedge\gamma = \frac{1}{2}(\delta\otimes\gamma - \gamma\otimes\delta)$, for any complex vector space $V$ and $\delta,\gamma\in V^*$.} to obtain the final equality above, we use the following identities:
\begin{multline*}
  (\alpha\otimes s) \wedge (\bar\beta\otimes t^*)
  =
  \frac{1}{2}\left(\alpha\otimes s\otimes \bar\beta\otimes t^* - \bar\beta\otimes t^*\otimes \alpha\otimes s\right)
  \\
  =
  \frac{1}{2}\left(\alpha\otimes \bar\beta - \bar\beta\otimes \alpha\right)\otimes s\otimes t^*
  =
  (\alpha\wedge\bar\beta)\otimes s\otimes t^*.
\end{multline*}
Hence,
\begin{equation}
  \label{eq:Inner_product_tau_and_zeta_psi_equals_tau_otimes_bar_psi_star_and_star_zeta_elemtary_tensors}
  \langle \alpha\otimes s, \zeta(\beta\otimes t)\rangle_{\Lambda^{0,2}(E)}
  =
  \langle (\alpha \wedge \bar\beta)\otimes (s\otimes t^*), \star\zeta \rangle_{\Lambda^{2,2}(\gl(E))}.
\end{equation}
Writing $\tau = \sum_{k=1}^r\tau_k\otimes e_k \in \Omega^{0,2}(U,E)$ and $\psi = \sum_{l=1}^r\psi_l\otimes e_l \in \Omega^{0,2}(U,E)$ for $\tau_k, \psi_k \in \Omega^{0,2}(U)$ and appealing to \eqref{eq:Complex_conjugate_sections_E_and_Lambda02E} and \eqref{eq:Hermitian_duals_sections_LambdapqE} to write
\[
  \psi^* = \sum_{l=1}^r\bar\psi_l\otimes e_l^* \in \Omega^{2,0}(U,E^*),
\]
we thus obtain
\begin{align*}
  \langle \tau\wedge\psi^*, \star \zeta \rangle_{\Lambda^{2,2}(\gl(E))}
  &= \sum_{k,l=1}^r \langle \tau_k\wedge\bar\psi_l\otimes e_k\otimes e_l^*, \star \zeta \rangle_{\Lambda^{2,2}(\gl(E))}
  \\
  &= \sum_{k,l=1}^r \langle \tau_k \otimes e_k, \zeta(\psi_l\otimes e_l) \rangle_{\Lambda^{0,2}(E)}
    \quad\text{(by \eqref{eq:Inner_product_tau_and_zeta_psi_equals_tau_otimes_bar_psi_star_and_star_zeta_elemtary_tensors})}
  \\
  &= \langle \tau, \zeta\psi \rangle_{\Lambda^{0,2}(E)},
\end{align*}
and therefore
\begin{equation}
  \label{eq:Inner_product_tau_and_zeta_psi_equals_tau_otimes_bar_psi_star_and_star_zeta}
  \langle R_\psi^*\tau,\zeta\rangle_{\fsl(E)}
  =
  \langle \tau, \zeta\psi\rangle_{\Lambda^{0,2}(E)}
  =
  \langle \tau\wedge \psi^*, \star \zeta \rangle_{\Lambda^{2,2}(\gl(E))},
\end{equation}
as claimed. We shall also use the identity \eqref{eq:Inner_product_tau_and_zeta_psi_equals_tau_otimes_bar_psi_star_and_star_zeta} in the form
\begin{equation}
  \label{eq:Inner_product_tau_and_zeta_psi_equals_star_tau_otimes_bar_psi_star_and_zeta}
  \langle R_\psi^*\tau,\zeta\rangle_{\fsl(E)}
  =
  \langle \tau, \zeta\psi\rangle_{\Lambda^{0,2}(E)}
  =
  \langle \star(\tau\wedge\psi^*), \zeta \rangle_{\Lambda^0(\gl(E))},
\end{equation}
obtained using the fact that $\langle\varpi,\star\zeta\rangle_{\Lambda^{2,2}(\gl(E))} = \langle\star\varpi,\zeta\rangle_{\Lambda^0(\gl(E))}$ by \cite[Proposition 1.2.20 (ii), p. 32]{Huybrechts_2005}.

Next, we prove the claim \eqref{eq:R_varphi_star_sigma_dagger_is_varphi_tensor_sigma_star_tracefree}. If $\zeta = \sum_{i,j=1}^r\zeta_{ij}\, e_i\otimes e_j^* \in \Omega^0(U,\fsl(E))$, then by \eqref{eq:Hermitian_dual_elementary_tensor} and \eqref{eq:Adjoint_is_conjugate_linear},
\[
  \zeta^\dagger = \sum_{i,j=1}^r (\zeta_{ij}\, e_i\otimes e_j^*)^\dagger = \sum_{i,j=1}^r \bar\zeta_{ij}\, e_j\otimes e_i^*
   = \sum_{i,j=1}^r \bar\zeta_{ji}\, e_i\otimes e_j^*.
\]
Therefore, for all $\zeta \in \Omega^0(\fsl(E))$ we have
\begin{align*}
  \langle (R_\varphi^*\sigma)^\dagger, \zeta\rangle_{\fsl(E)}
  &=
    \langle \zeta^\dagger, R_\varphi^*\sigma\rangle_{\fsl(E)}
    \quad\text{(by \eqref{eq:Inner_product_M_dagger_N_equals_inner_product_N_dagger_M})}
  \\
  &= \langle R_\varphi(\zeta^\dagger), \sigma\rangle_E
    \quad\text{(by \eqref{eq:varphi_psi_star_Omega0E_oplus_Omega02E_to_Omega0slE})}
  \\
  &= \langle \zeta^\dagger\varphi, \sigma\rangle_E
    \quad\text{(by \eqref{eq:varphi_psi_Omega0slE_to_Omega0E_oplus_Omega02E})},
\end{align*}
while
\begin{multline*}
  \langle (\varphi\otimes\sigma^*)_0, \zeta\rangle_{\fsl(E)}
  =
  \langle \varphi\otimes\sigma^*, \zeta\rangle_{\fsl(E)}
  \\
  =
  \sum_{i,j}\langle \varphi\otimes\sigma^*, \zeta_{ij} e_i\otimes e_j^*\rangle_{\fsl(E)}
  =
  \sum_{i,j} \bar\zeta_{ij}\langle \varphi\otimes\sigma^*, e_i\otimes e_j^*\rangle_{\fsl(E)}
  \\
=
  \sum_{i,j}\bar\zeta_{ij}\langle \varphi, e_i\otimes e_j^*(\sigma)\rangle_E \quad\text{(by \eqref{eq:Inner_product_sigma_otimes_bar_varphi_star_and_zeta_equals_sigma_and_zeta_varphi})}
  \\
  =
  \sum_{i,j}\langle \varphi, \zeta_{ij}e_i\otimes e_j^*(\sigma)\rangle_E
  =
  \langle \varphi, \zeta\sigma \rangle_E
  =
  \langle \zeta^\dagger\varphi, \sigma \rangle_E.
\end{multline*}
By combining the preceding two chains of equalities, we verify the claim \eqref{eq:R_varphi_star_sigma_dagger_is_varphi_tensor_sigma_star_tracefree}.

By modifying the derivation of \eqref{eq:R_varphi_star_sigma_dagger_is_varphi_tensor_sigma_star_tracefree} \mutatis we can also prove the claim \eqref{eq:R_psi_star_tau_dagger_is_psi_tensor_tau_star_tracefree}. Indeed, for all $\zeta \in \Omega^0(\fsl(E))$ we have
\begin{align*}
  \langle (R_\psi^*\tau)^\dagger, \zeta\rangle_{\fsl(E)}
  &=
  \langle \zeta^\dagger, R_\psi^*\tau\rangle_{\fsl(E)} \quad\text{(by \eqref{eq:Inner_product_M_dagger_N_equals_inner_product_N_dagger_M})}
  \\
  &=   \langle R_\psi\zeta^\dagger, \tau\rangle_{\Lambda^{0,2}(E)}
    \quad\text{(by \eqref{eq:varphi_psi_star_Omega0E_oplus_Omega02E_to_Omega0slE})}
  \\
  &=   \langle \zeta^\dagger\psi, \tau\rangle_{\Lambda^{0,2}(E)}
    \quad\text{(by \eqref{eq:varphi_psi_Omega0slE_to_Omega0E_oplus_Omega02E})},
\end{align*}
while 
\begin{align*}
  \langle \star(\psi\wedge\tau^*)_0, \zeta\rangle_{\fsl(E)}
  &=
  \langle \star(\psi\wedge\tau^*), \zeta\rangle_{\fsl(E)}
  \\
  &= \langle \psi\wedge\tau^*, \star\zeta\rangle_{\Lambda^{2,2}(\fsl(E))}
  \\
  &= \langle \psi, \zeta\tau \rangle_{\Lambda^{0,2}(E)} \quad\text{(by \eqref{eq:Inner_product_tau_and_zeta_psi_equals_tau_otimes_bar_psi_star_and_star_zeta} with $\psi,\tau$ interchanged)}
  \\
  &= \langle \zeta^\dagger\psi, \tau \rangle_{\Lambda^{0,2}(E)}.
\end{align*}
By combining the preceding chains of equalities, we see that we have verified the claim \eqref{eq:R_psi_star_tau_dagger_is_psi_tensor_tau_star_tracefree}. This completes the 
proof of Claim \ref{claim:R_varphi_star_and_R_psi_star_tracefree}
\end{proof}

Although Lemma \ref{lem:Real_Coulomb_slice_SO3_monopole_complex_Kaehler_del_A_and_del-bar_A_adjoints_tensor_products} is the final version of the Coulomb slice condition we require for the proof of Proposition \ref{prop:Itoh_1985_proposition_2-4_SO3_monopole_complex_Kaehler}, we record the following result for subsequent comparisons of slice conditions.

\begin{lem}[Comparison of operators defining unitary and complex Coulomb-gauge slice conditions]
\label{lem:RealPartOfPairs_rd_0_star_is_d_APhi_0_*}
Continue the hypotheses of Lemma \ref{lem:SO3_monopole_Coulomb_gauge_slice_condition}.
In addition, assume that $(A,\Phi)$ is type $1$ in the sense that $\Phi=(\varphi,0)\in\Om^0(E)\oplus \Om^{0,2}(E)$.  Then the operators $d_{A,\Phi}^{0,*}$ (namely, the $L^2$ adjoint of the operator $d_{A,\Phi}^0$ defined in \eqref{eq:d_APhi^0}) and $\bar\partial_{A,\varphi}^{0,*}$ (namely, the $L^2$ adjoint of the operator\footnote{See Section \ref{sec:Elliptic_complexes}.} $\bar\partial_{A,\varphi}^0$ defined in \eqref{eq:d0StablePair}) are related by
\begin{equation}
\label{eq:RealPartOfPairs_rd_0_star_is_d_APhi_0_*}
d_{A,\Phi}^{0,*}(a,\phi)
=
\pi_{\su(E)} \bar\partial_{A,\varphi}^{0,*}(a'',\sigma)
\end{equation}
where $(a,\phi)\in\Om^1(\su(E))\oplus\Om^0(W_{\can}^+\otimes E)$ and $(a'',\sigma)\in\Om^{0,1}(\fsl(E))\oplus \Om^0(E)$ are related by $a=\frac{1}{2}(a'+a'')$ as in \eqref{eq:Decompose_a_in_Omega1suE_into_10_and_01_components} and $a' = -(a'')^\dagger$ as in \eqref{eq:Kobayashi_7-6-11} and $\phi=(\si,0)\in\Om^0(E)\oplus\Om^{0,2}(E)$.
\end{lem}

\begin{proof}
We compute
\begin{align*}
\pi_{\su(E)} \bar\rd_A^* a''
&=
\frac{1}{2}\left( \bar\rd_A^* a''-(\bar\rd_A^* a'')^\dagger\right)
\\
&=
\frac{1}{2}\left( \bar\rd_A^* a''+(\rd_A^* a')^\dagger\right)
\quad\text{(by \eqref{eq:Commute_dagger_and_del_glE} and $a'=-(a'')^\dagger$)}
\\
&=
d_A^*a \quad\text{by \eqref{eq:d_star_equals_partial_stars})}.
\end{align*}
Thus, for $a=(1/2)(a'' -(a'')^\dagger)$, we have
\begin{equation}
\label{eq:su(E)_component_of_bar_partial_adjoint_is_d_adjoint}
\pi_{\su(E)} \bar\rd_A^* a''
=
d_A^*a.
\end{equation}
From the definition \eqref{eq:d0StablePair} of $\bar\rd_{A,\varphi}^0$ and the definition \eqref{eq:Right_multiplication_of_section_slE_by_section_E} of $R_\varphi$, we have
\[
\bar\rd_{A,\varphi}^{0,*}(a'',\sigma)
=
\bar\rd_A^*a'' - R_\varphi^*\si.
\]
Therefore,
\begin{align*}
\pi_{\su(E)}\bar\rd_{A,\varphi}^{0,*}(a'',\sigma)
&=
\pi_{\su(E)}\left( \bar\rd_A^*a'' - R_\varphi^*\si\right)
\\
&=\pi_{\su(E)}\bar\rd_A^*a'' - \pi_{\su(E)}R_\varphi^*
\\
&=
d_A^*a - \sR_\Phi^*(\si,0)
\quad
\text{by \eqref{eq:su(E)_component_of_bar_partial_adjoint_is_d_adjoint} and \eqref{eq:su2_Projection_of_AdjointOf_sl2_MultOperator_Is_su2_MultOperator})}
\\
&= d_{A,\Phi}^{0,*}(a,\phi) \quad\text{by \eqref{eq:d_APhi^0_star_Expression})},
\end{align*}
which completes the proof of the lemma.
\end{proof} 

\begin{rmk}[Motivation for the convention $a=(a'+a'')/2$]
\label{rmk:Factor_one_half_motivation_1-forms}
The identity \eqref{eq:RealPartOfPairs_rd_0_star_is_d_APhi_0_*} provides the motivation for our convention of writing $a=\frac{1}{2}(a'+a'')\in\Om^1(\su(E))$, where $a'' \in \Omega^{0,1}(\fsl(E))$ and $a' = -(a'')^\dagger \in \Omega^{1,0}(\fsl(E))$.
\end{rmk}

Because $(X,g,J)$ is assumed to be at least almost K\"ahler in Proposition \ref{prop:Itoh_1985_proposition_2-4_SO3_monopole_complex_Kaehler}, the following \emph{K\"ahler Identities} hold (see Huybrechts \cite[Proposition 3.1.12, p. 120]{Huybrechts_2005} and Wells \cite[Chapter V, Corollary 4.10, p. 193]{Wells3} for the complex K\"ahler case and Cirici and Wilson \cite[Proposition 3.1 (4)]{Cirici_Wilson_2020_harmonic} for the almost K\"ahler case),
\begin{subequations}
  \label{eq:Kaehler_identity_commutator_Lambda_del-bar_and_Lambda_del}
  \begin{align}
    \label{eq:Kaehler_identity_commutator_Lambda_del-bar}
    [\Lambda,\bar\partial] &= -i\partial^* \quad\text{and}
    \\
    \label{eq:Kaehler_identity_commutator_Lambda_del}
    [\Lambda,\partial] &= i\bar\partial^* \quad\text{on } \Omega^\bullet(X,\CC).
  \end{align}
\end{subequations}
According to Demailly \cite[Section 7.1]{Demailly_complex_analytic_differential_geometry} (see also Demailly \cite{Demailly_1986}), Donaldson \cite[Proof of Proposition 3]{DonASD} and Donaldson and Kronheimer \cite[Equation (6.1.8), p. 212]{DK}, and Kobayashi \cite[Section 3.2, p. 62]{Kobayashi_differential_geometry_complex_vector_bundles} the K\"ahler Identities \eqref{eq:Kaehler_identity_commutator_Lambda_del-bar_and_Lambda_del} extend to bundle-valued forms to give
\begin{subequations}
\label{eq:Kaehler_identity_commutator_Lambda_del-bar_A_and_Lambda_del_A}   
\begin{align}
\label{eq:Kaehler_identity_commutator_Lambda_del-bar_A} 
  [\Lambda, \bar\partial_A] &= -i\partial_A^* \quad\text{and}
  \\
  \label{eq:Kaehler_identity_commutator_Lambda_del_A} 
[\Lambda, \partial_A] &= i\bar\partial_A^* \quad\text{on } \Omega^\bullet(\fsl(E)).
\end{align}
\end{subequations}
One can prove these identities first for $\Omega^\bullet(E)$ and then $\Omega^\bullet(\gl(E))$ or $\Omega^\bullet(\fsl(E))$. We can now conclude the

\begin{proof}[Proof of Proposition \ref{prop:Itoh_1985_proposition_2-4_SO3_monopole_complex_Kaehler}]
Suppose that $(a,\phi) \in \bH_{A,\Phi}^1$ and so by \eqref{eq:H_APhi^1} we have
\[
  d_{A,\Phi}^1(a,\phi) = 0 \quad\text{and}\quad d_{A,\Phi}^{0,*}(a,\phi) = 0.
\]
Note that equation \eqref{eq:Itoh_1985_2-18_SO3_monopole_complex_Kaehler_02_component_curvature_equation} is equal to \eqref{eq:H_dbar_APhi^01_explicit_02} and equation \eqref{eq:Itoh_1985_2-18_SO3_monopole_complex_Kaehler_Dirac_operator} is equal to \eqref{eq:H_dbar_APhi^01_explicit_01} in the definition \eqref{eq:H_dbar_APhi^01_explicit} of $\bH_{\bar\partial_{A,(\varphi,\psi)}}^1$. Hence, we need to show that $(a'',\sigma,\tau)$ obeys the remaining defining equation \eqref{eq:H_dbar_APhi^01_explicit_0} for $\bH_{\bar\partial_{A,(\varphi,\psi)}}^1$ in order to conclude that $(a'',\sigma,\tau) \in \bH_{\bar\partial_{A,(\varphi,\psi)}}^1$.

By applying the K\"ahler identities \eqref{eq:Kaehler_identity_commutator_Lambda_del-bar_A_and_Lambda_del_A} to give
\[
  \Lambda(\bar\partial_Aa' + \partial_Aa'') = - i\partial_A^*a' + i\bar\partial_A^*a'',
\]
we see that \eqref{eq:Itoh_1985_2-18_SO3_monopole_complex_Kaehler_11_component_curvature_equation} yields
\[
  - i\partial_A^*a' + i\bar\partial_A^*a'' - i\left(\sigma\otimes\varphi^* + \varphi\otimes\sigma^*\right)_0 + i\star\left(\tau\otimes\psi^* + \psi\otimes\tau^*\right)_0 = 0,
\]
that is,
\[
  \partial_A^*a' - \bar\partial_A^*a'' + \left(\sigma\otimes\varphi^* + \varphi\otimes\sigma^*\right)_0 - \star\left(\tau\otimes\psi^* + \psi\otimes\tau^*\right)_0 = 0.
\]
We subtract the preceding identity from \eqref{eq:Real_Coulomb_slice_SO3_monopole_complex_Kaehler_del_A_and_del-bar_A_adjoints_tensor_products}, namely
\[
  \partial_A^*a' + \bar\partial_A^*a'' - \left(\sigma\otimes\varphi^* - \varphi\otimes\sigma^*\right)_0 - \star\left(\tau\otimes\psi^* - \psi\otimes\tau^*\right)_0 = 0
\]
to give
\[
  2\bar\partial_A^*a'' - 2(\sigma\otimes\varphi^*)_0 + 2\star(\psi\otimes\tau^*)_0
    = 0
\]
and therefore
\begin{equation}
\label{eq:dbar_A_slice_condition_SO3_monopoles_complex_Kaehler}
\bar\partial_A^*a'' -(\sigma\otimes\varphi^*)_0 + \star(\psi\otimes\tau^*)_0 = 0.
\end{equation}
But equation \eqref{eq:dbar_A_slice_condition_SO3_monopoles_complex_Kaehler} for $(a'',\sigma,\tau)$ is 
equal to remaining defining equation for $\widehat\bH_{\partial_{A,(\varphi,\psi)}}^1$, 
namely \eqref{eq:widehatH_dbar_APhi^01_explicit_0},
since \eqref{eq:R_varphi_star_sigma_is_sigma_tensor_varphi_star_tracefree} and \eqref{eq:R_psi_star_tau_is_tau_tensor_psi_star_tracefree} yield
\[
  R_\varphi^*\sigma = (\sigma\otimes\varphi^*)_0
  \quad\text{and}\quad
  R_\psi^*\tau = \star\left(\tau\wedge\psi^*\right)_0.
\]
Finally, when $\psi \equiv 0$, then equation \eqref{eq:dbar_A_slice_condition_SO3_monopoles_complex_Kaehler} for $(a'',\sigma,\tau)$ coincides exactly with the remaining defining equation for $\bH_{\partial_{A,(\varphi,\psi)}}^1$, namely \eqref{eq:H_dbar_APhi^01_explicit_0}. This completes the proof of Proposition \ref{prop:Itoh_1985_proposition_2-4_SO3_monopole_complex_Kaehler}.
\end{proof}

\subsection{Auxiliary results in linear algebra required for proof of isomorphism between first-order cohomology groups}
\label{subsec:Auxiliary_results_linear_algebra_first-order_cohomology}
In this subsection, we gather auxiliary results in linear algebra required by the proof of Proposition \ref{prop:Itoh_1985_proposition_2-4_SO3_monopole_complex_Kaehler}.

\subsubsection{Certain properties of complex vector spaces and complex Hilbert spaces}
\label{subsubsec:Properties_complex_vector_spaces_and_complex_Hilbert_spaces}
We to recall some properties of complex vector spaces and complex Hilbert spaces. We may always view a complex vector space $\fK$ as the complexification $\sK^\CC = \sK\otimes_\RR\CC$ of some real subspace $\sK\subset\fK$, called a \emph{real form} \cite[Definition 4.1]{Conrad_complexification} for $\fK$, and if $\{e_i\}$ is an $\RR$-basis for $\sK$, then $\{e_i\}$ is a $\CC$-basis for $\fK$ (see Conrad \cite[Theorem 3.2]{Conrad_complexification}). Given a real form $\sK$ for $\fK$, there is a unique complex conjugation map, $\sQ\in\End_\RR(\fK)$ \cite[Definition 4.4 and Theorem 4.11]{Conrad_complexification}, which we denote by $\sQ(k)=\bar{k}$ for all $k\in\fK$. Moreover, we may always view a complex Hilbert space $\fH$ as the complexification $\sH^\CC = \sH\otimes_\RR\CC$ of some real Hilbert subspace $\sH\subset\fH$ (a real form for $\fH$) and if $x,y \in \sH$, then $\langle x,y\rangle_\fH = \langle x,y\rangle_\sH$ (see Conway \cite[Exercise 1.1.7]{Conway_course_functional_analysis}). Moreover, if $h=x+iy, k=u+iv \in \fH$ for $x,y,u,v\in\sH$, then
\[
  \langle h,k\rangle_\fH = \langle x+iy, u+iv\rangle_\fH = \langle x,u\rangle_\sH + \langle y, v\rangle_\sH + i\left(\langle y, u\rangle_\sH - \langle x,v\rangle_\sH\right)
\]
and
\[
  \overline{\langle h,k\rangle}_\fH = \langle x,u\rangle_\sH + \langle y, v\rangle_\sH + i\left(\langle x,v\rangle_\sH - \langle y, u\rangle_\sH\right)
  = \langle x-iy, u-iv\rangle_\fH,
\]
and thus, using $\bar h=x-iy, \bar k=u-iv$,
\begin{equation}
  \label{eq:Complex_conjugation_and_inner_product}
  \langle k,h\rangle_\fH = \overline{\langle h,k\rangle}_\fH = \langle \bar h,\bar k\rangle_\fH, \quad\text{for all } h, k \in \fH.
\end{equation}
If $\{e_i\}$ is an orthonormal $\RR$-basis for $\sH$, then $\{e_i\}$ is an orthonormal $\CC$-basis for $\fH$ by \cite[Theorem 3.2]{Conrad_complexification}, \cite[Exercise 1.1.7]{Conway_course_functional_analysis}.

For any complex Hilbert space $\fH$ with dual space $\fH^*=\Hom(\fH,\CC)$, the composition of complex conjugation $\sQ \in \End_\RR(\fH)$ and the Riesz (complex anti-linear, isometric) isomorphism $\sR \in \End_\RR(\fH)$ \cite[Theorem 12.5]{Rudin} given by $\sR:\fH \ni h\mapsto h^* := \langle\cdot,h\rangle_\fH \in \fH^*$ yields a \emph{complex linear} isometry $\sR\circ\sQ:\fH \cong \fH^*$ given by $h \mapsto \bar{h}^* := \langle\cdot,\bar{h}\rangle_\fH$. The inner product, 
\begin{equation}
\label{eq:Inner_product_dual_Hilbert_space}  
  \langle h^*,k^*\rangle_{\fH^*} = \langle k,h\rangle_\fH, \quad\text{for all } h, k \in \fH,
\end{equation}
is implied by the polarization identity \cite[p. 86]{RudinRealComplex} (for an inner product that is complex conjugate linear in the second variable) and the fact that $\sR$ is an isometry, so $\|h^*\|_{\fH^*}=\|h\|_\fH$:
\begin{align*}
  \langle h^*,k^*\rangle_{\fH^*}
  &= \frac{1}{4}\left(\|h^*+k^*\|_{\fH^*}^2 - \|h^*-k^*\|_{\fH^*}^2 + i\|h^*+ik^*\|_{\fH^*}^2 - i\|h^*-ik^*\|_{\fH^*}^2\right)
  \\
  &= \frac{1}{4}\left(\|h^*+k^*\|_{\fH^*}^2 - \|k^*-h^*\|_{\fH^*}^2 + i\|k^*-ih^*\|_{\fH^*}^2 - i\|k^*+ih^*\|_{\fH^*}^2\right)
  \\
  &= \frac{1}{4}\left(\|(k+h)^*\|_{\fH^*}^2 - \|(k-h)^*\|_{\fH^*}^2 + i\|(k+ih)^*\|_{\fH^*}^2 - i\|(k-ih)^*\|_{\fH^*}^2\right)
  \\
  &= \frac{1}{4}\left(\|k+h\|_\fH^2 - \|k-h\|_\fH^2 + i\|k+ih\|_\fH^2 - i\|k-ih\|_\fH^2\right)
   = \langle k,h\rangle_\fH.                                 
\end{align*}
If we replace $\sR$ with the composition $\sR\circ\sQ$, then \eqref{eq:Complex_conjugation_and_inner_product} and \eqref{eq:Inner_product_dual_Hilbert_space} yield
\begin{equation}
\label{eq:Inner_product_dual_Hilbert_space_complex_conjugation}  
  \langle \bar h^*,\bar k^*\rangle_{\fH^*} = \langle h,k\rangle_\fH, \quad\text{for all } h, k \in \fH.
\end{equation}
This concludes our discussion of properties of complex vector spaces and complex Hilbert spaces relevant to the proof of Proposition \ref{prop:Itoh_1985_proposition_2-4_SO3_monopole_complex_Kaehler}.

\subsubsection{Certain properties of adjoints of operators on complex Hilbert spaces}
\label{subsubsec:Certain_properties_adjoints_operators_on_complex_Hilbert_spaces}
We recall properties of adjoints of operators on complex Hilbert spaces relevant to the proof of Proposition \ref{prop:Itoh_1985_proposition_2-4_SO3_monopole_complex_Kaehler}. For complex Hilbert spaces $\fH, \fK$, the complex vector space $\Hom(\fH,\fK) = \fK\otimes\fH^*$ is a Hilbert space with inner product induced from those of $\fH$ and $\fK$ (see Kadison and Ringrose \cite[Section 2.6, p. 125]{KadisonRingrose1}) induced from its definition on elementary tensors:
\begin{equation}
\label{eq:Inner_product_elementary_tensors_K_otimes_H_star}
\langle k_1\otimes h_1^*, k_2\otimes h_2^*\rangle_{\fH^*\otimes\fK} := \langle k_1, k_2\rangle_\fK \langle h_2, h_1\rangle_{\fH^*},
\end{equation}
where we use \eqref{eq:Inner_product_dual_Hilbert_space} to write
\[
  \langle h_1^*, h_2^*\rangle_{\fH^*} = \langle h_2, h_1\rangle_\fH.
\]
As an aside, note that if $M = \sum_{i,j=1}^{r,s} m_{ij} k_i\otimes h_j^*$ and $N = \sum_{i,j=1}^{r,s} n_{ij} k_i\otimes h_j^*$ and $\{k_i\otimes h_j^*\}_{i,j=1}^{r,s}$ is an orthonormal set with respect to the inner product \eqref{eq:Inner_product_elementary_tensors_K_otimes_H_star}, then
\[
  \langle M, N \rangle_{\Hom(\fH,\fK)}
  =
  \sum_{i,j=1}^{r,s} \sum_{a,b=1}^{r,s} \langle m_{ij}k_i\otimes h_j^*, n_{ab}k_a\otimes h_b^* \rangle_{\Hom(\fH,\fK)}
  =
  \sum_{i,j=1}^{r,s} m_{ij}\bar n_{ji},
\]
that is
\begin{equation}
  \label{eq:Frobenius_inner_product_Hom_H_to_K}
  \langle M, N \rangle_{\Hom(\fH,\fK)} = \tr(MN^\dagger),
\end{equation}  
which is the usual expression for the \emph{Frobenius inner product} on $\Hom(\fK,\fH)$ (see \cite[p. 241]{Montero_Gonzalez_Florez_Garcia_Suarez_2002}.) Therefore,
\begin{equation}
  \label{eq:Frobenius_norm_Hom_H_to_K}
  \|M\|_{\Hom(\fH,\fK)}^2 = \langle M, M \rangle_{\Hom(\fH,\fK)} = \tr(MM^\dagger)
\end{equation}
defines the corresponding \emph{Frobenius norm} on $\Hom(\fH,\fK)$. Note that
\begin{equation}
  \label{eq:Operator_norm_Hom_H_to_K}
  \|M\|_\Op := \sup_{h\in\fH\less\{0\}}\frac{\|Mh\|_\fK}{\|h\|_\fH} 
\end{equation}
defines the \emph{operator norm} on $\Hom(\fH,\fK)$, which is well-defined regardless of whether $M$ has finite rank. Observe that if $\fK=\fH$ and $\{e_i\}_{i=1}^n$ is an orthonormal basis for $\fH$ and $M = e_i\otimes e_j^*$, then
\[
  \|M\|_\Op = \sup_{h\in\fH\less\{0\}}\frac{\|e_i\otimes e_j^*(h)\|_\fH}{\|h\|_\fH}
  \geq \frac{\|e_i\otimes e_j^*(e_j)\|_\fH}{\|e_j\|_\fH} = 1,
\]
and if $h_j = e_j^*(h)$, then
\[
  \|M\|_\Op = \sup_{h\in\fH\less\{0\}}\frac{\|e_i\otimes e_j^*(h)\|_\fH}{\|h\|_\fH}
  = \frac{|h_j|}{\|h\|_\fH} \leq 1,
\]
so $\|e_i\otimes e_j^*\|_\Op = 1$. For $M = e_i\otimes e_j^*$, we also have
\begin{align*}
  \|M\|_{\End(\fH)}^2
  &= \tr(MM^\dagger) \quad\text{(by \eqref{eq:Frobenius_norm_Hom_H_to_K})}
  \\
  &= \tr((e_i\otimes e_j^*)(e_i\otimes e_j^*)^\dagger)
  \\
  &=  \tr((e_i\otimes e_j^*)(e_j\otimes e_i^*))
    \quad\text{(by forthcoming \eqref{eq:Hermitian_dual_elementary_tensor})}
  \\
  &=  \tr(e_i\otimes e_i^*) = 1.
\end{align*}
Observe that
\[
  \|M\|_{\gl(\fH)} = \left\|\sum_{i,j=1}^n m_{ij}e_i\otimes e_j^*\right\|_{\gl(\fH)}
  = \left(\sum_{i,j=1}^n m_{ij}\bar m_{ij}\right)^{1/2}, \quad\text{for all } M \in \gl(\fH),
\]
since $\{e_i\otimes e_j^*\}_{i,j=1}^n$ is an orthonormal basis for $\gl(\fH)$, and so we conclude that
\begin{equation}
  \label{eq:Operator_norm_End_H}
  \|M\|_\Op = \|M\|_{\End(\fH)}, \quad\text{for all } M \in \End(\fH). 
\end{equation}
Using
\[
  \Hom(\fH,\fK)^* = (\fK\otimes\fH^*)^* = \fK^*\otimes\fH^{**} = \fH\otimes\fK^* = \Hom(\fK,\fH),
\]
we claim that
\begin{equation}
  \label{eq:Inner_product_M_dagger_N_equals_inner_product_N_dagger_M}
  \langle M_1,M_2\rangle_{\Hom(\fH,\fK)} = \langle M_2^\dagger,M_1^\dagger\rangle_{\Hom(\fH,\fK)},
  \quad\text{for all } M_1, M_2 \in \Hom(\fH,\fK).
\end{equation}
Indeed, any elementary tensor in $\Hom(\fH,\fK)$ obeys
\begin{equation}
  \label{eq:Hermitian_dual_elementary_tensor}
  (k\otimes h^*)^\dagger = h\otimes k^*, \quad\text{for all } h \in \fH, k \in \fK,
\end{equation}
since the Hermitian adjoint $N^\dagger \in \Hom(\fK,\fH)$ of an operator $N \in \Hom(\fH,\fK)$ is defined by the relation
\begin{equation}
  \label{eq:Hermitian_adjoint}
  \langle y,Nx \rangle_\fK = \langle N^\dagger y,x \rangle_\fH, \quad\text{for all } x \in \fH, y \in \fK,
\end{equation}
and thus $M = k\otimes h^*$ obeys
\begin{multline*}
  \langle M^\dagger y, x \rangle_\fH
  = \langle h\otimes k^*(y), x \rangle_\fH
  = \langle h, x \rangle_\fH\langle y, k \rangle_\fK
  = \langle y, \overline{\langle h, x \rangle}_\fH k \rangle_\fK
  \\
  = \langle y, \langle x, h \rangle_\fH k \rangle_\fK
  = \langle y, h^*(x) k \rangle_\fK
  = \langle y, k\otimes h^*(x) \rangle_\fK
  = \langle y, Mx \rangle_\fK.
\end{multline*} 
Therefore, if $M_1 = k_1\otimes h_1^*$ and $M_2 = k_2\otimes h_2^*$ are operators in $\Hom(\fH,\fK)$, then their adjoints $M_1^\dagger = h_1\otimes k_1^*$ and $M_2^\dagger = h_2\otimes k_2^*$ in $\Hom(\fK,\fH)$ obey 
\begin{multline*}
  \langle M_2^\dagger,M_1^\dagger\rangle_{\Hom(\fK,\fH)}
  = \langle  h_2\otimes k_2^*,h_1\otimes k_1^*\rangle_{\fH\otimes\fK^*}
  = \langle  h_2,h_1\rangle_\fH \langle  k_2^*,k_1^*\rangle_{\fK^*}
  \\
  = \langle  h_1^*,h_2^*\rangle_{\fH^*} \langle  k_1,k_2\rangle_\fK
  = \langle  k_1\otimes h_1^*,k_2\otimes h_2^*\rangle_{\fK\otimes\fH^*}
  = \langle  M_1,M_2\rangle_{\Hom(\fH,\fK)}.
\end{multline*}
Since \eqref{eq:Inner_product_M_dagger_N_equals_inner_product_N_dagger_M} holds when $M_1, M_2$ are elementary tensors in $\fk\otimes\fH^*$, then it holds for all elements of $\fk\otimes\fH^*$.
This follows from the equality
\[
(N_1+N_2)^\dagger = N_1^\dagger + N_2^\dagger
\]
which we prove by computing
\begin{align*}
  \langle y,(N_1+N_2)x \rangle_\fK
  &=
  \langle y,N_1x \rangle_\fK + \langle y,N_1x \rangle_\fK
  \\
  &= \langle N_1^\dagger y,x \rangle_\fH+\langle N_2^\dagger y,x \rangle_\fH
  \\
  &=
  \langle (N_1^\dagger +N_2^\dagger)y,x \rangle_\fH.
\end{align*}
If $M \in \gl(\fH) = \fH\otimes\fH^*$, then $U = M-M^\dagger \in \fu(\fH) \subset \gl(\fH) = \fu(\fH)\oplus i\fu(H)$.  Observe that
\[
  \pi_{\fu(\fH)}:\gl(\fH) \ni M \mapsto \frac{1}{2}(M-M^\dagger) \in \fu(\fH)
\]
is an orthogonal projection, since $\pi_{\fu(\fH)}^2 = \pi_{\fu(\fH)}$ and
\begin{multline*}
  \langle \pi_{\fu(\fH)}M,U \rangle_{\fu(\fH)}
  = \frac{1}{2}\langle M-M^\dagger,U \rangle_{\fu(\fH)}
  = \frac{1}{2}\langle M-M^\dagger,U \rangle_{\gl(\fH)}
  = \frac{1}{2}\langle M,U \rangle_{\gl(\fH)} - \frac{1}{2}\langle M^\dagger,U \rangle_{\gl(\fH)}
  \\
  = \frac{1}{2}\langle M,U \rangle_{\gl(\fH)} - \frac{1}{2}\langle U^\dagger, M \rangle_{\gl(\fH)}
  = \frac{1}{2}\langle M,U \rangle_{\gl(\fH)} + \frac{1}{2}\langle U, M \rangle_{\gl(\fH)}
  \\
  = \frac{1}{2}\langle M,U \rangle_{\gl(\fH)} + \frac{1}{2}\overline{\langle M, U \rangle}_{\gl(\fH)}
  = \Real \langle M,U \rangle_{\gl(\fH)}
  = \Real \langle M,\pi_{\fu(\fH)}U \rangle_{\gl(\fH)}.
\end{multline*}
Hence, $\pi_{\fu(\fH)}^* = \pi_{\fu(\fH)}$, as claimed. In general, we have $\|\pi_{\fu(\fH)}\| = 1$
and
\begin{multline*}
  \|M\|_{\gl(\fH)}^2 = \|\pi_{\fu(\fH)}M + \pi_{\fu(\fH)}^\perp M\|_{\gl(\fH)}^2
  = \|\pi_{\fu(\fH)}M\|_{\gl(\fH)}^2 + \|\pi_{\fu(\fH)}^\perp M\|_{\gl(\fH)}^2
  \\
  \geq \|\pi_{\fu(\fH)}M\|_{\gl(\fH)}^2 = \|\pi_{\fu(\fH)}M\|_{\fu(\fH)}^2,
\end{multline*}
so that
\[
  \|\pi_{\fu(\fH)}M\|_{\fu(\fH)} \leq \|M\|_{\gl(\fH)}, \quad\text{for all } M \in \gl(\fH).
\]
Lastly, note that by \eqref{eq:Hermitian_adjoint} the map $\Hom(\fH,\fK) \ni M \mapsto M^* \in \Hom(\fK,\fH)$ is conjugate linear, so
\begin{equation}
  \label{eq:Adjoint_is_conjugate_linear}
  (cM)^\dagger = \bar c M^\dagger, \quad\text{for all } c\in \CC, M \in \Hom(\fH,\fK).
\end{equation}
This concludes our digression on properties of adjoints of operators on complex Hilbert spaces.

\section[Isomorphisms between zeroth-order cohomology groups]{Comparison of elliptic complexes for non-Abelian monopole and pre-holomorphic pair equations over almost Hermitian four-manifolds: Isomorphisms between zeroth-order cohomology groups}
\label{sec:Isomorphisms_between_zeroth-order_cohomology_groups}
In this section, we prove the forthcoming

\begin{prop}[Canonical real linear isomorphism between zeroth-order cohomology groups for elliptic complexes over 
almost Hermitian four-manifolds]  
\label{prop:Donaldson_1985_proposition_3_and_Itoh_1985_proposition_4-1_SO3_monopole_almost_Hermitian}
Let $(X,g,J)$ be a closed, smooth, almost Hermitian manifold of real dimension four, $E$ be a smooth Hermitian vector bundle $E$ over $X$, and $(\rho_\can,W_\can)$ denote the canonical spin${}^c$ structure over $(X,g,J)$ (see \eqref{eq:Canonical_spinc_bundles} for $W_\can$ and \eqref{eq:Canonical_Clifford_multiplication} for $\rho_\can$).
If $(A,\Phi)$ is a solution to the non-Abelian monopole equations \eqref{eq:SO(3)_monopole_equations_Kaehler} with $\Phi=(\varphi,\psi)\in \Om^0(E)\oplus \Om^{0,2}(E)$, then the monomorphism of infinite-dimensional vector spaces
\begin{equation}
\label{eq:E0Inclusion}
\Om^0(\su(E)) \ni \xi \mapsto \xi \in \Om^0(\fsl(E)
\end{equation}
induces a monomorphism of finite-dimensional complex vector spaces,
\begin{equation}
\label{eq:InclusionOfH0Spaces}
\bH_{A,\Phi}^0\subset \widehat\bH_{\bar\partial_{A,(\varphi,\psi)}}^0
:=
\pi_{\su(E)}\left(\bH_{\bar\partial_{A,(\varphi,\psi)}}^0\right),
\end{equation}
where $\bH_{\bar\partial_{A,(\varphi,\psi)}}^0$ is as in \eqref{eq:H_dbar_APhi^00}, and $\bH_{A,\Phi}^0$ is as in \eqref{eq:H_APhi^0}, and
\[
\pi_{\su(E)}:\Om^0(\fsl(E)) \ni \zeta \mapsto \frac{1}{2}\left(\zeta-\zeta^\dagger\right) \in \Om^0(\su(E)),
\]
is the projection defined in \eqref{eq:Pointwise_orthogonal_projections_slE_onto_suE_and_isuE}. In addition, if $(X,g,J)$ is almost K\"ahler and $(A,\Phi)$ is a type $1$ solution as in Remark \ref{rmk:Projective_vortices_type_1_monopoles}, so $\psi \equiv 0$, then the monomorphism \eqref{eq:InclusionOfH0Spaces} is an isomorphism and there is an isomorphism
\begin{equation}
\label{eq:H0Inclusion}
\bH_{A,\Phi}^0\otimes_\RR\CC
\cong
\bH_{\bar\partial_{A,(\varphi,0)}}^0
\end{equation}
defined by the complex linear extension of \eqref{eq:E0Inclusion}.
\end{prop}

We have the following analogue of Remark \ref{rmk:Difference_widehatH_dbar_APhi^01_and_H_dbar_APhi^01}.

\begin{rmk}[Isomorphism between zero-order harmonic spaces for the non-Abelian monopole and holomorphic pair complexes]
\label{rmk:Isomorphism_nonAbelianH0_preholomH0}
If $(X,g,J)$ is almost K\"ahler and $(A,\Phi)=(A,(\varphi,0))$ is a type $1$ solution as in Remark \ref{rmk:Projective_vortices_type_1_monopoles}, then
combining the isomorphism \eqref{eq:H0Inclusion} with the equality
\[
\bH_{\bar\rd_A,\varphi}^0 = \bH_{\bar\rd_A,(\varphi,0)}^0
\]
from equation \eqref{eq:Cohomologies_complexes_pre-holomorphic_pair_type1_and_holomorphic_pair_n_is_2} of Corollary \ref{cor:Comparison_elliptic_complexes_cohomology_groups_pre-holomorphic_and_holomorphic_pairs} gives an isomorphism of complex vector spaces,
\begin{equation}
\label{eq:Isomorphism_nonAbelianH0_preholomH0}
\bH_{A,\Phi}^0\otimes_\RR\CC
\cong
\bH_{\bar\rd_A,\varphi}^0,
\end{equation}
and this is the isomorphism that we use in practice.
\end{rmk}

\begin{proof}
The definitions of the harmonic space $\bH_{A,\Phi}^0$ in \eqref{eq:H_APhi^0}, of $d_{A,\Phi}^0$ in \eqref{eq:d_APhi^0}, and of $\sR_\Phi$ in \eqref{eq:Phi_Omega0suE_to_Omega0V+} imply that
\begin{equation}
\label{eq:H0dIntersection}
\bH_{A,\Phi}^0
=
\Ker d_A \cap \Ker\sR_\Phi\subset \Om^0(\su(E)).
\end{equation}
Similarly, the expression for $\bH_{\bar\rd_A,(\varphi,\psi)}^0$ in \eqref{eq:H_dbar_APhi^00_explicit} and the definition of $R_{(\varphi,\psi)}$ in \eqref{eq:varphi_psi_Omega0slE_to_Omega0E_oplus_Omega02E} imply that
\begin{equation}
\label{eq:H0rdIntersection}
\bH_{\bar\rd_A,(\varphi,\psi)}^0
=
\Ker\bar\rd_A \cap \Ker R_{(\varphi,\psi)}\subset\Om^0(\fsl(E)).
\end{equation}
Assume that $\xi\in\Ker d_A \cap \Ker\sR_\Phi$. Recall from \eqref{eq:d_A_sum_components_almost_complex_manifold_integrable} that\footnote{Observe that we do not need to assume $N_J\equiv 0$ to obtain $d_A\xi=0 \implies\bar\partial\xi=0$ since we could have equally have applied the more general relation \eqref{eq:d_A_sum_components_almost_complex_manifold} that is valid on  any almost complex manifold.}
\[
  d_A\xi = \partial_A\xi + \bar\partial_A\xi \in \Om^{1,0}(\fsl(E))\oplus\Om^{1,0}(\fsl(E)),
\]
so $\rd_A\xi\in\Om^{1,0}(\fsl(E))$ and $\bar\rd_A\xi\in\Om^{0,1}(\fsl(E))$ are linearly independent and so the equality $d_A\xi=0$ implies that $\bar\rd_A\xi=0$ and $\xi\in\Ker\bar\partial_A$. The definitions of $\sR_\Phi$ in \eqref{eq:Phi_Omega0suE_to_Omega0V+} and $R_{(\varphi,\psi)}$ in \eqref{eq:varphi_psi_Omega0slE_to_Omega0E_oplus_Omega02E} imply that
the monomorphism \eqref{eq:E0Inclusion} maps $\Ker\sR_\Phi$ into $\Ker R_{(\varphi,\psi)}$. Consequently, $\xi \in \Ker\bar\rd_A \cap \Ker R_{(\varphi,\psi)}$ and the expressions \eqref{eq:H0dIntersection} for $\bH_{A,\Phi}^0$ and \eqref{eq:H0rdIntersection} for $\bH_{\bar\rd_A,(\varphi,\psi)}^0$ imply that the monomorphism \eqref{eq:E0Inclusion} gives a monomorphism,
\[
\bH_{A,\Phi}^0 \to \bH_{\bar\rd_A,(\varphi,\psi)}^0.
\]
Thus  $\xi\in \bH_{A,\Phi}^0$ implies  that $\xi\in \bH_{\bar\rd_A,(\varphi,\psi)}^0$ and so
\[
\xi=\pi_{\su(E)}\xi \in \widehat \bH_{\bar\rd_A,(\varphi,\psi)}^0,
\]
proving that the monomorphism \eqref{eq:E0Inclusion} induces a monomorphism \eqref{eq:InclusionOfH0Spaces}. We now turn to the question of whether the monomorphism \eqref{eq:InclusionOfH0Spaces} is surjective.

\begin{claim}
\label{claim:Kernel_of_H0_Laplacian_Non_Abelian_Monopoles}
If $\zeta\in \bH_{\bar\rd_A,(\varphi,\psi)}^0$ and we assume in addition that $(X,g,J)$ is almost K\"ahler and $(A,\Phi)$ is a type $1$ solution, then
\begin{equation}
\label{eq:Kernel_of_H0_Laplacian_Non_Abelian_Monopoles}
d_{A,\Phi}^{0,*}d_{A,\Phi}\left(\zeta-\zeta^\dagger\right)=0,
\end{equation}
and so $\xi = \pi_{\su(E)}\zeta = (1/2)(\zeta-\zeta^\dagger)\in \bH_{A,\Phi}^0$.
\end{claim}

\begin{proof}[Proof of Claim \ref{claim:Kernel_of_H0_Laplacian_Non_Abelian_Monopoles}]
Because $\bar\rd_A\zeta=0$ by \eqref{eq:H0rdIntersection}, equation \eqref{eq:Commute_dagger_and_del_glE} implies that
\[
\rd_A(\zeta^\dagger)=(\bar\rd_A\zeta)^\dagger=0.
\]
Thus, using
\[
\La F_A= \La (\bar\rd_A\rd_A + \rd_A\bar\rd_A),
\]
the equalities $\bar\rd_A\zeta=0$ and $\rd_A\zeta^\dagger=0$ give us
\begin{equation}
\label{eq:La_F_A_on_zeta_and_zetaDagger}
\La F_A \zeta
=
\La \bar\rd_A\rd_A \zeta,
\quad\text{and}\quad
\La F_A \zeta^\dagger
=
\La \rd_A\bar\rd_A\zeta^\dagger.
\end{equation}
The preceding equalities imply that
\begin{align*}
\La\left(\bar\rd_A\rd_A - \rd_A\bar\rd_A\right)\left(\zeta-\zeta^\dagger\right)
&=
\La\left( \bar\rd_A\rd_A \zeta + \rd_A\bar\rd_A\zeta^\dagger\right)
\quad\text{(since $\bar\rd_A\zeta=0$ and $\rd_A\zeta^\dagger=0$)}
\\
&=
\La F_A \left(\zeta+\zeta^\dagger\right)
\quad\text{(by \eqref{eq:La_F_A_on_zeta_and_zetaDagger}),}
\end{align*}
that is (note the change from $\zeta-\zeta^\dagger$ to $\zeta+\zeta^\dagger$),
\begin{equation}
\label{eq:Difference_of_(1,1)_On_su(2)_zeta}
\La\left(\bar\rd_A\rd_A - \rd_A\bar\rd_A\right)\left(\zeta-\zeta^\dagger\right)
=
\La F_A \left(\zeta+\zeta^\dagger\right).
\end{equation}
We note that sections in the images of
\[
d_A:\Om^0(\su(E))\to \Om^1(\su(E))\oplus \Om^0(W_\can^+\otimes E)
\]
and
\[
\sR_{(\varphi,\psi)}:\Om^0(\su(E))\to \Om^1(\su(E))\oplus \Om^0(W_\can^+\otimes E)
\]
are pointwise orthogonal because the image of $d_A$ is contained in the  $\Om^1(\su(E))$ while that of $\sR_{(\varphi,\psi)}$ is contained in $\Om^0(W_\can^+\otimes E)$. This pointwise orthogonality implies $L^2$ orthogonality and thus
\begin{equation}
\label{eq:Vanishing_CurvatureMultiplication_CrossTerms_in_d0_Laplacian}
d_A^*\sR_{(\varphi,\psi)}=0
\quad\text{and}\quad
\sR_{(\varphi,\psi)}^* d_A=0.
\end{equation}
Then, we compute
\begin{align*}
d_{A,\Phi}^{0,*}d_{A,\Phi}^0\left(\zeta-\zeta^\dagger\right)
&=
\left( d_A^*-\sR_{(\varphi,\psi)}^*\right)\left( d_A-\sR_{(\varphi,\psi)}\right)\left(\zeta-\zeta^\dagger\right)
\quad\text{(by \eqref{eq:d_APhi^0_star_Expression} and \eqref{eq:d_APhi^0})}
\\
&=
     \left( d_A^*d_A +\sR_{(\varphi,\psi)}^*\sR_{(\varphi,\psi)}\right)\left(\zeta-\zeta^\dagger\right)
\quad\text{(by \eqref{eq:Vanishing_CurvatureMultiplication_CrossTerms_in_d0_Laplacian})}
\\
&=
\left(
    \left( \rd_A^*+\bar\rd_A^*\right)
    \left(\rd_A +\bar\rd_A\right)
 +\sR_{(\varphi,\psi)}^*\sR_{(\varphi,\psi)}\right)\left(\zeta-\zeta^\dagger\right)
\\
 &=
\left(
    \left( i[\La,\bar\rd_A] -i[\La,\rd_A]\right)
    \left(\rd_A +\bar\rd_A\right)
 +\sR_{(\varphi,\psi)}^*\sR_{(\varphi,\psi)}\right)\left(\zeta-\zeta^\dagger\right)
 \quad\text{(by \eqref{eq:Kaehler_identity_commutator_Lambda_del-bar_A_and_Lambda_del_A})}
\\
 &=
\left(
    \left( i\La\bar\rd_A -i\La\rd_A\right)
    \left(\rd_A +\bar\rd_A\right)
 +\sR_{(\varphi,\psi)}^*\sR_{(\varphi,\psi)}\right)\left(\zeta-\zeta^\dagger\right)
  \\
  &\qquad\text{(since $\La$ vanishes on $\Om^1(\fsl(E))$)}
\\
&=
     \left( i\La\bar\rd_A\rd_A -i\La \rd_A\bar\rd_A  +\sR_{(\varphi,\psi)}^*\sR_{(\varphi,\psi)}\right)\left(\zeta-\zeta^\dagger\right)
  \\
  &\qquad\text{(since $\La F_A^{0,2}=0$ and $\La F_A^{2,0}=0$)}
\\
&=
i\La F_A \left(\zeta+\zeta^\dagger\right)
+
\sR_{(\varphi,\psi)}^*\sR_{(\varphi,\psi)}\left(\zeta-\zeta^\dagger\right)
\quad\text{(by \eqref{eq:Difference_of_(1,1)_On_su(2)_zeta})}
\end{align*}
giving us
\begin{equation}
\label{eq:LaplacianComparisonForH0}
d_{A,\Phi}^{0,*}d_{A,\Phi}^0\left(\zeta-\zeta^\dagger\right)
=
i\La F_A \left(\zeta+\zeta^\dagger\right)
+
\sR_{(\varphi,\psi)}^*\sR_{(\varphi,\psi)}\left(\zeta-\zeta^\dagger\right).
\end{equation}
(Note that the above use of \eqref{eq:Kaehler_identity_commutator_Lambda_del-bar_A_and_Lambda_del_A} is the only place in this proof where the assumption that $(X,g,J)$ is almost K\"ahler is used.) The connection $A$ on $E$ defines a connection on $\fsl(E)$ by the commutation relation $d_A\zeta=[d_A,\zeta]$, for all $\zeta \in \Omega^0(\fsl(E))$. Hence in \eqref{eq:LaplacianComparisonForH0}, the curvature $F_A$ of the connection $A$ acts on sections of $\fsl(E)$ by $\zeta\mapsto [F_A,\zeta]$, as described in the comments in Donaldson and Kronheimer following \cite[Lemma 6.1.7, p. 212]{DK}.
We can thus rewrite the curvature term in \eqref{eq:LaplacianComparisonForH0} as
\begin{align*}
i\La F_A(\zeta+\zeta^\dagger)
&=
\left[i\La F_A,\zeta+\zeta^\dagger\right]
\\
&=
\left[i(\La F_A)_0,\zeta+\zeta^\dagger\right] \quad\text{(trace component acts trivially)}
\\  
&=
       -\frac{1}{2}\left[(\varphi\otimes\varphi^*)_0-\star(\psi\otimes\psi^*)_0,\zeta+\zeta^\dagger\right]
       \quad\text{(by \eqref{eq:SO(3)_monopole_equations_(1,1)_curvature}),}
\end{align*}
and thus
\begin{equation}
\label{eq:H0_Laplacian_Curvature_Term}
\begin{aligned}
i\La F_A(\zeta+\zeta^\dagger)
&=
-\frac{1}{2} \left(
(\varphi\otimes\varphi^*)_0\zeta+(\varphi\otimes\varphi^*)_0\zeta^\dagger
-
\zeta (\varphi\otimes\varphi^*)_0 -\zeta^\dagger (\varphi\otimes\varphi^*)_0\right.
\\
&\quad
\left.
-\star(\psi\otimes\psi^*)_0\zeta-\star(\psi\otimes\psi^*)_0\zeta^\dagger
+\zeta \star(\psi\otimes\psi^*)_0 + \zeta^\dagger \star(\psi\otimes\psi^*)_0
\right).
\end{aligned}
\end{equation}
For $\zeta\in\Om^0(\fsl(E))$ and $\Phi\in\Om^0(W_{\can}^+\otimes E)$, we observe that
by \eqref{eq:Phi_star}, we have
\[
(\zeta\Phi)^*
=
\langle \cdot, \zeta\Phi\rangle_{W_{\can}^+\otimes E}
=
\langle \zeta^\dagger(\cdot),\Phi\rangle_{W_{\can}^+\otimes E}
=
\Phi^*\zeta^\dagger.
\]
Thus, for $\zeta\in\Om^0(\fsl(E))$, and $\varphi\in\Om^0(E)$, and $\psi\in\Om^{0,2}(E)$,
\begin{align*}
\zeta (\varphi\otimes\varphi^*)_0 &= ((\zeta \varphi)\otimes\varphi^*),
\\
(\varphi\otimes\varphi^*)_0\zeta^\dagger &= (\varphi\otimes (\zeta \varphi)^*)_0,
\\
\zeta \star(\psi\otimes\psi^*)_0 &= \star((\zeta\psi)\otimes\psi)_0,
\\
\star (\psi\otimes\psi^*)_0\zeta^\dagger &= \star(\psi\otimes (\zeta \psi)^*)_0.
\end{align*}
Applying the preceding equalities and the assumption that $\zeta\in\Ker R_\varphi\cap \Ker R_\psi$ to \eqref{eq:H0_Laplacian_Curvature_Term} yields
\begin{equation}
\label{eq:Laplacian_on_E0_CurvatureTerm}
i\La F_A(\zeta+\zeta^\dagger)
=
-\frac{1}{2} \left(
(\varphi\otimes\varphi^*)_0\zeta
-\zeta^\dagger (\varphi\otimes\varphi^*)_0
-\star(\psi\otimes\psi^*)_0\zeta
+ \zeta^\dagger \star(\psi\otimes\psi^*)_0
\right).
\end{equation}
We now rewrite the other term appearing in the right-hand side of \eqref{eq:LaplacianComparisonForH0}, namely
\[
\sR_{(\varphi,\psi)}^*\sR_{(\varphi,\psi)}\left(\zeta-\zeta^\dagger\right).
\]
To do so, we first observe that from the definitions of the operators $\sR_\Phi$
in \eqref{eq:Phi_Omega0suE_to_Omega0V+} and $R_{(\varphi,\psi)}$ in \eqref{eq:varphi_psi_Omega0slE_to_Omega0E_oplus_Omega02E},
\[
\sR_{(\varphi,\psi)}(\zeta-\zeta^\dagger)
=
(R_\varphi+R_\psi)(\zeta-\zeta^\dagger),
\]
where $R_\varphi$ and $R_\psi$ are the components of $R_{\varphi,\psi}$ defined in
\eqref{eq:varphi_psi_Omega0slE_to_Omega0E_oplus_Omega02E}.  Next, because the subspaces $E$ and $\Lambda^{0,2}(E)$ of $W_\can^+\otimes E$ are fiberwise orthogonal, we have the pointwise equality
\[
  \langle\zeta_1\varphi,\zeta_2\psi\rangle_{\Lambda^{0,\bullet}(E)} = 0,
  \quad\text{for all } \zeta_1,\zeta_2\in\Om^0(\fsl(E)), \varphi\in\Om^0(E), \text{ and } \psi\in\Om^{0,2}(E).
\]
Thus, $\zeta\psi\in\Om^{0,2}(E)$ is pointwise orthogonal to sections in the image of $R_\varphi$ and $\zeta\varphi\in \Om^0(E)$ is pointwise orthogonal to sections in the image of $R_\psi$. This pointwise orthogonality implies $L^2$ orthogonality and thus
\begin{equation}
\label{eq:VanishingCrossTermsInMultiplicationLaplacian}
R_\psi^*(\zeta \varphi)=0
\quad\text{and}\quad
R_\varphi^*(\zeta\psi)=0,
\quad\text{for all } \zeta\in\Om^0(\fsl(E)), \varphi\in\Om^0(E), \text{ and } \psi\in\Om^{0,2}(E).
\end{equation}
From \eqref{eq:su2_Projection_of_AdjointOf_sl2_MultOperator_Is_su2_MultOperator}, we have
\[
\sR_{(\varphi,\psi)}^*
=
\pi_{\su(E)}R_{(\varphi,\psi)}^*
=
\pi_{\su(E)}\left(R_\varphi^*+R_\psi^*\right),
\]
where $\pi_{\su(E))}$ is the projection operator defined in \eqref{eq:Pointwise_orthogonal_projections_slE_onto_suE_and_isuE}.
Thus,
\begin{align*}
\sR_{(\varphi,\psi)}^*\sR_{(\varphi,\psi)}\left(\zeta-\zeta^\dagger\right)
&=
\pi_{\su(E)}\left(R_\varphi^*+R_\psi^*\right)(R_\varphi+R_\psi)(\zeta-\zeta^\dagger)
\\
&=
-\pi_{\su(E)}\left(R_\varphi^*+R_\psi^*\right)(\zeta^\dagger \varphi+\zeta^\dagger\psi)
\quad\text{(since $\zeta\in\Ker R_\varphi\cap\Ker R_\psi$)}
\\
&=
-\pi_{\su(E)}\left( R_\varphi^*(\zeta^\dagger\varphi) + R_\psi^*(\zeta^\dagger\psi)\right)
\quad\text{(by \eqref{eq:VanishingCrossTermsInMultiplicationLaplacian})}
\\
&=
-\pi_{\su(E)}\left(
((\zeta^\dagger\varphi)\otimes \varphi^*)_0
+ \star((\zeta^\dagger\psi)\wedge\psi^*)_0\right)
\quad\text{(by \eqref{eq:R_varphi_star_sigma_is_sigma_tensor_varphi_star_tracefree} and \eqref{eq:R_psi_star_tau_is_tau_tensor_psi_star_tracefree})}
\\
&=
-\pi_{\su(E)}\left( \zeta^\dagger(\varphi\otimes \varphi^*)_0 +
\zeta^\dagger\star(\psi\wedge\psi^*)_0\right),
\end{align*}
and therefore,
\begin{multline}
\label{eq:Laplacian_on_E0_MultiplicationTerm}
\sR_{(\varphi,\psi)}^*\sR_{(\varphi,\psi)}\left(\zeta-\zeta^\dagger\right)
\\
=
-\frac{1}{2}
\left(
\zeta^\dagger(\varphi\otimes \varphi^*)_0
-
(\varphi\otimes \varphi^*)_0\zeta
+
\zeta^\dagger\star(\psi\wedge\psi^*)_0
-
\star(\psi\wedge\psi^*)_0\zeta
\right).
\end{multline}
Combining \eqref{eq:LaplacianComparisonForH0},
\eqref{eq:Laplacian_on_E0_CurvatureTerm}, and
\eqref{eq:Laplacian_on_E0_MultiplicationTerm}
yields
\begin{align*}
d_{A,\Phi}^{0,*}d_{A,\Phi}^0\left(\zeta-\zeta^\dagger\right)
&=
-\frac{1}{2} \left(
(\varphi\otimes\varphi^*)_0\zeta
-\zeta^\dagger (\varphi\otimes\varphi^*)_0
-\star(\psi\otimes\psi^*)_0\zeta
 +\zeta^\dagger \star(\psi\otimes\psi^*)_0
\right)
\\
&\qquad
-\frac{1}{2}
\left(
\zeta^\dagger(\varphi\otimes \varphi^*)_0
-
(\varphi\otimes \varphi^*)_0\zeta
+
\zeta^\dagger\star(\psi\wedge\psi^*)_0
-
\star(\psi\wedge\psi^*)_0\zeta
\right)
\\
&=
\star(\psi\otimes\psi^*)_0\zeta
-\zeta^\dagger \star(\psi\otimes\psi^*)_0.
\end{align*}
Thus, if $(A,\Phi)$ is type $1$ so $\psi=0$ and $\zeta\in \bH_{\bar\partial_{A,(\varphi,0)}}^0$, then
$\zeta-\zeta^\dagger$ satisfies \eqref{eq:Kernel_of_H0_Laplacian_Non_Abelian_Monopoles} and so
\[
0
=
\left( \zeta-\zeta^\dagger,d_{A,\Phi}^{0,*}d_{A,\Phi}^0\left(\zeta-\zeta^\dagger\right) \right)_{L^2}
=
\left(d_{A,\Phi}^0\left(\zeta-\zeta^\dagger\right),d_{A,\Phi}^0\left(\zeta-\zeta^\dagger\right) \right)_{L^2}
=
\| d_{A,\Phi}^0\left(\zeta-\zeta^\dagger\right)\|_{L^2}^2,
\]
and thus
\begin{equation}
\label{eq:su2ProjectionInH0}
\zeta-\zeta^\dagger\in \bH_{A,\Phi}^0.
\end{equation}
The definition of $\pi_{\su(E))}$ in \eqref{eq:Pointwise_orthogonal_projections_slE_onto_suE_and_isuE} and \eqref{eq:su2ProjectionInH0} imply that
\[
\pi_{\su(E)}\left(\bH_{\bar\partial_{A,(\varphi,\psi)}}^0\right)
\subset 
\bH_{A,\Phi}^0,
\]
and so the inclusion \eqref{eq:InclusionOfH0Spaces} is an equality. This completes the proof of Claim \ref{claim:Kernel_of_H0_Laplacian_Non_Abelian_Monopoles}.
\end{proof}

Because we have already shown that \eqref{eq:InclusionOfH0Spaces} is injective, Claim \ref{claim:Kernel_of_H0_Laplacian_Non_Abelian_Monopoles} implies that 
\eqref{eq:InclusionOfH0Spaces} is an isomorphism with the additional hypotheses that $(X,g,J)$ is almost K\"ahler and $(A,\Phi)$ is type $1$.

Lastly, we prove \eqref{eq:H0Inclusion}. If $\zeta\in\bH_{\bar\partial_{A,(\varphi,0)}}^0$ and we write $\zeta=\xi_1+i\xi_2$, where $\xi_j\in\Om^0(\su(E))$ for $j=1,2$, then \eqref{eq:su2ProjectionInH0} implies that
\[
\xi_1=\frac{1}{2}(\zeta-\zeta^\dagger)\in \bH_{A,\Phi}^0.
\]
Because $\bH_{\bar\partial_{A,(\varphi,0)}}^0$ is a complex vector space, the fact that $\zeta\in\bH_{\bar\partial_{A,(\varphi,0)}}^0$ implies $i\zeta\in\bH_{\bar\partial_{A,(\varphi,0)}}^0$ and so \eqref{eq:su2ProjectionInH0} yields
\[
\xi_2
=
-\frac{1}{2}\left( i\zeta - (i\zeta)^\dagger\right)
\in \bH_{A,\Phi}^0.
\]
Thus,
\[
\zeta=\xi_1+i\xi_2
\in \bH_{A,\Phi}^0\otimes_\RR\CC,
\]
completing the proof of \eqref{eq:H0Inclusion} and thus of Proposition \ref{prop:Donaldson_1985_proposition_3_and_Itoh_1985_proposition_4-1_SO3_monopole_almost_Hermitian}.
\end{proof}

\section[Isomorphisms between second-order cohomology groups]{Comparison of elliptic complexes for non-Abelian monopole and pre-holomorphic pair equations over complex K\"ahler surfaces: Second-order cohomology groups}
\label{sec:Isomorphisms_between_second-order_cohomology_groups}
We can now proceed to analyze the second-order cohomology group for the non-Abelian monopole elliptic deformation complex \eqref{eq:SO3MonopoleDefComplex} over a closed, complex K\"ahler surface. 

\begin{prop}[Canonical real linear isomorphism between second-order cohomology groups for elliptic complexes over closed K\"ahler surfaces]
\label{prop:Itoh_1985_proposition_2-3_SO3_monopole_complex_Kaehler}
Let $E$ be a smooth, Hermitian vector bundle $E$ over a closed, complex K\"ahler surface $(X,g,J)$ and $(\rho_\can,W_\can)$ denote the canonical spin${}^c$ structure over $X$ (see \eqref{eq:Canonical_spinc_bundles} for $W_\can$ and \eqref{eq:Canonical_Clifford_multiplication} for $\rho_\can$). If $(A,\Phi)$ is a smooth solution to the unperturbed non-Abelian monopole equations \eqref{eq:SO(3)_monopole_equations_Kaehler} with $\Phi=(\varphi,\psi)$, then the induced real linear map of infinite-dimensional real vector spaces
\begin{equation}
  \label{eq:Real_linear_map_obstruction_spaces_SO(3)_pairs_to_preholomorphic_pairs}
  \Omega^+(\su(E))\oplus \Omega^0(W^-\otimes E) \ni (v,\nu)
  \mapsto
  (v'',\nu,\xi) \in \Omega^{0,2}(\fsl(E)) \oplus \Omega^{0,1}(E)\oplus\Omega^0(\su(E))
\end{equation}
implied by \eqref{eq:Canonical_spinc_bundles} and the forthcoming decomposition \eqref{eq:Donaldson_Kronheimer_lemma_2-1-57_sum_3_terms} determine a canonical isomorphism of finite-dimensional real vector spaces,
\begin{equation}
\label{eq:Itoh_1985_proposition_2-3_SO3_monopole_complex_Kaehler_isomorphism}  
\bH_{A,\Phi}^2 \cong \left\{ (v'',\nu,\xi) \in \Omega^{0,2}(\fsl(E))\oplus \Omega^{0,1}(E) \oplus \Omega^0(\su(E)) \text{ obeys } \eqref{eq:4-1_and_2_and_3}
\right\},
\end{equation}
where $\bH_{A,\Phi}^2$ is as in \eqref{eq:H_APhi^2} (see the forthcoming \eqref{eq:4-1_and_2_and_3} for a more explicit expression). If $\psi=0$, so $(A,\Phi)$ is a type $1$ solution, then there is a canonical isomorphism of finite-dimensional real vector spaces,
\begin{equation}
\label{eq:Itoh_1985_proposition_2-3_SO3_monopole_complex_Kaehler_isomorphism_type1}
 \bH_{A,\Phi}^2 \cong \bH_{\bar\partial_{A,(\varphi,0)}}^2 \oplus \bH_{A,\Phi}^0,
\end{equation}
where $\bH_{\bar\partial_{A,(\varphi,0)}}^2$ is defined in \eqref{eq:H_dbar_APhi^02} (with a more explicit expression in \eqref{eq:H_dbar_APhi^01_part_explicit_type1}) and $\bH_{A,\Phi}^0$ is defined in \eqref{eq:H_APhi^0} (with a more explicit expression implied by \eqref{eq:d_APhi^0} or \eqref{eq:H0dIntersection}).
\end{prop}

We have the following analogue of Remark \ref{rmk:Difference_widehatH_dbar_APhi^01_and_H_dbar_APhi^01}.

\begin{rmk}[Relation between the second-order harmonic spaces for the non-Abelian monopole and holomorphic pair complexes]
\label{rmk:Relation_H_APhi^2_and_H_dbar_APhi^2}
When $X$ has complex dimension two, the identity $\bH_{\bar\rd_A,(\varphi,0)}^2=\bH_{\bar\rd_A,\varphi}^2$ from \eqref{eq:Cohomologies_complexes_pre-holomorphic_pair_type1_and_holomorphic_pair_n_is_2_simplified} in Corollary \ref{cor:Vanishing_third-order_cohomology_group_holomorphic_pair_elliptic_complex_Kaehler_surface} and the isomorphism \eqref{eq:Itoh_1985_proposition_2-3_SO3_monopole_complex_Kaehler_isomorphism_type1} yields
\begin{equation}
  \label{eq:H2_NonAbelianMonopole_HolomorphicPair}
 \bH_{A,\Phi}^2 \cong \bH_{\bar\partial_A,\varphi}^2 \oplus \bH_{A,\Phi}^0,
\end{equation}
which is the isomorphism that we use in practice.
\end{rmk}

\begin{proof}
We continue the notation in the proof of Proposition \ref{prop:Itoh_1985_proposition_2-4_SO3_monopole_complex_Kaehler}. It is convenient to define a differential,
\begin{equation}
  \label{eq:hatd_APhi^1}
  \hat d_{A,\Phi}^1: \Omega^1(\su(E)) \oplus \Omega^0(W^+\otimes E)
  \to
  \Omega^0(\su(E)) \oplus \Omega^{2,0}(\fsl(E)) \oplus \Omega^{0,2}(\fsl(E)) \oplus \Omega^{0,1}(E),
\end{equation}
by the components of the left-hand side of the system \eqref{eq:Itoh_1985_2-18_SO3_monopole_complex_Kaehler}, so that (again writing $\Phi=(\varphi,\psi)$ and $\phi=(\sigma,\tau)$) while $a=\frac{1}{2}(a'+a'')$ with $(a')^\dagger = -a''$ by \eqref{eq:Decompose_a_in_Omega1suE_into_10_and_01_components} and \eqref{eq:Kobayashi_7-6-11})
\begin{equation}
\label{eq:Itoh_1985_2-6_SO3_monopole_complex_Kaehler}  
\begin{aligned}
  \hat d_{A,\Phi}^1(a,\phi)
  &=
  \begin{pmatrix}
      \Lambda(\bar\partial_Aa' + \partial_Aa'') - i\left(\sigma\otimes\varphi^* + \varphi\otimes\sigma^*\right)_0 + i\star\left(\tau\otimes\psi^* + \psi\otimes\tau^*\right)_0
      \\
      \partial_Aa' + \left(\varphi\otimes\tau^* + \sigma\otimes\psi^*\right)_0
      \\
      \bar\partial_Aa'' - \left(\tau\otimes\varphi^* + \psi\otimes\sigma^*\right)_0
      \\
      \bar\partial_A\sigma + \bar\partial_A^*\tau + a''\varphi - \star (a'\wedge\star\psi)
  \end{pmatrix}      
  \\
  &\in \Omega^0(\su(E)) \oplus \Omega^{2,0}(\fsl(E)) \oplus \Omega^{0,2}(\fsl(E)) \oplus \Omega^{0,1}(E).
\end{aligned}
\end{equation}
Before proceeding further, we note the

\begin{rmk}[Equivalent forms of the first component of $d_{A,\Phi}^1(a,\phi)$]
We remark that the restriction of the Lefshetz operator \cite[Definition 1.2.18, p. 31]{Huybrechts_2005},
\begin{equation}
  \label{eq:Huybrechts_definition_1-2-18}
  L = L_\omega:\Omega^0(\su(E)) \ni \xi \mapsto \xi\otimes\omega \in \Omega^{1,1}(\su(E)),
\end{equation}
is a real linear injection \cite[Proposition 1.2.30 (iv), p. 36]{Huybrechts_2005} and so the expression for $\pi_0\hat d_{A,\Phi}^1(a,\phi)$ implied by  \eqref{eq:Itoh_1985_2-6_SO3_monopole_complex_Kaehler} is equivalent to 
\begin{multline}
\label{eq:Itoh_1985_2-6_SO3_monopole_complex_Kaehler_0-component}  
  L\pi_\omega d_{A,\Phi}^1(a,\phi)
  =
  \langle\bar\partial_Aa' + \partial_Aa'',\omega\rangle_{\Lambda^{1,1}(X)}\otimes\omega
  \\
  - \frac{i}{2}\left(\sigma\otimes\varphi^* + \varphi\otimes\sigma^*\right)_0\otimes\omega + \frac{i}{2}\star\left(\tau\otimes\psi^* + \psi\otimes\tau^*\right)_0\otimes\omega
  \\
  \in \Omega^0(\su(E))\otimes\omega \subset \Omega^{1,1}(\fsl(E)),
\end{multline}
where $\pi_0$ denotes the projection from the codomain of $\hat d_{A,\Phi}^1$ onto $\Omega^0(\su(E))$ and $\pi_\omega$ denotes the projection from the codomain of $d_{A,\Phi}^1$ onto $\Omega^0(\su(E))\otimes\omega$. See Remark \ref{rmk:Alternative_form_1-1_component_SO3_monopole_equations} and the proof of Lemma \ref{lem:SO3_monopole_equations_almost_Kaehler_manifold}.
\end{rmk}

Recall from \cite[Lemma 2.1.57, p. 47]{DK} that
\begin{subequations}
\label{eq:Donaldson_Kronheimer_lemma_2-1-57}  
\begin{align}
  \label{eq:Donaldson_Kronheimer_lemma_2-1-57_self_dual_two_forms}  
  \Omega^+(X,\CC) &= \Omega^{2,0}\oplus \Omega^0(X,\CC)\omega \oplus \Omega^{0,2}(X),
  \\
  \label{eq:Donaldson_Kronheimer_lemma_2-1-57_anti_self_dual_two_forms}
  \Omega^-(X,\CC) &= \Omega_0^{1,1}(X),
\end{align}
\end{subequations}
where $\Omega_0^{1,1}(X) \subset \Omega^{1,1}(X)$ is the subspace of forms that are pointwise orthogonal to $\omega$. We shall write $v\in\Omega^+(\su(E))$ with respect to its $L^2$ orthogonal decomposition \eqref{eq:Donaldson_Kronheimer_lemma_2-1-57} as
\begin{equation}
  \label{eq:Donaldson_Kronheimer_lemma_2-1-57_sum_3_terms}
  v = \frac{1}{2}(v' + v'') + \frac{1}{2}\xi\otimes\omega
  \in \Omega^{2,0}(\fsl(E)) \oplus \Omega^{0,2}(\fsl(E))\oplus \Omega^0(\su(E))\otimes\omega,
\end{equation}
where $v''\in\Omega^{0,2}(\fsl(E))$ and $v'$ is defined by the reality condition (see Itoh \cite[p. 850]{Itoh_1985})
\begin{equation}
  \label{eq:Itoh_1985_p_850_reality_condition_su(E)_2-forms}
  v' = -(v'')^\dagger \in\Omega^{2,0}(\fsl(E)),
\end{equation}
and $\xi \in \Omega^0(\su(E))$. We claim that, for each
\[
  v \in\Omega^+(\su(E)) \quad\text{and}\quad \nu \in \Omega^0(W_\can^-\otimes E) = \Omega^{0,1}(E),
\]
the expression for the $L^2$ adjoint
\[
  \hat d_{A,\Phi}^{1,*}(v,\nu) \in \Omega^{1,0}(\fsl(E)) \oplus \Omega^{0,1}(\fsl(E))
  \oplus \Omega^0(E) \oplus \Omega^{0,2}(E)
\]
defined by
\begin{multline}
  \label{eq:Itoh_1985_2-6_SO3_monopole_complex_Kaehler_L2_adjoint} 
  \left(\hat d_{A,\Phi}^{1,*}(v,\nu), (a,\phi)\right)_{L^2(X)}
  =
  \left((v,\nu), \hat d_{A,\Phi}^1(a,\phi)\right)_{L^2(X)},
  \\
  \text{for all } (a,\phi) \in \Omega^1(\su(E)) \oplus \Omega^0(W_\can^+\otimes E)
  = \Omega^1(\su(E)) \oplus \Omega^0(E) \oplus \Omega^{0,2}(E),
\end{multline}
is given by the forthcoming system \eqref{eq:Itoh_1985_2-7_and_8_SO3_monopole_complex_Kaehler}. The Hilbert space inner product on the \emph{left}-hand side of \eqref{eq:Itoh_1985_2-6_SO3_monopole_complex_Kaehler_L2_adjoint} is defined by\footnote{The factor of $1/2$ in \eqref{eq:L2_inner_product_domain_d1APhi} arises by writing $v = \frac{1}{2}(v'+v'')+\frac{1}{2}\xi\otimes\omega$ as in \eqref{eq:Donaldson_Kronheimer_lemma_2-1-57_sum_3_terms} and identifying $v$ with the triple $\frac{1}{2}(v',v'',\xi)$.}
\begin{multline}
\label{eq:L2_inner_product_domain_d1APhi}  
  \left(\hat d_{A,\Phi}^{1,*}(v,\nu), (a,\phi)\right)_{L^2(X)}
  \\
  :=
  \frac{1}{2}\underbrace{\left(\hat d_{A,\Phi}^{1,*}v', (a,\phi)\right)_{L^2(X)}
  + \frac{1}{2}\left(\hat d_{A,\Phi}^{1,*}v'', (a,\phi)\right)_{L^2(X)}}_{\textrm{real}}
  \\
  + \underbrace{\left(\hat d_{A,\Phi}^{1,*}\xi, (a,\phi)\right)_{L^2(X)}}_{\textrm{real}}
  + \Real\underbrace{\left(\hat d_{A,\Phi}^{1,*}\nu, (a,\phi)\right)_{L^2(X)}}_{\textrm{complex}},
\end{multline}
where the Hermitian inner products containing $v'$ and $v''$ will be shown to be complex conjugates of one another (and so the indicated pair is real) and the inner product containing $\xi$ will be shown to be real. Similarly, the Hilbert space inner product on the \emph{right}-hand side of \eqref{eq:Itoh_1985_2-6_SO3_monopole_complex_Kaehler_L2_adjoint} is defined by
\begin{multline}
\label{eq:L2_inner_product_codomain_d1APhi}  
  \left((v,\nu), \hat d_{A,\Phi}^1(a,\phi)\right)_{L^2(X)}
  :=
  \underbrace{\frac{1}{2}\left(v', \hat d_{A,\Phi}^1(a,\phi)\right)_{L^2(X)}
    + \frac{1}{2}\left(v'', \hat d_{A,\Phi}^1(a,\phi)\right)_{L^2(X)}}_{\textrm{real}}
  \\
  + \frac{1}{2}\underbrace{\left(\xi, \hat d_{A,\Phi}^1(a,\phi)\right)_{L^2(X)}}_{\textrm{real}}
  + \Real\underbrace{\left(\nu, \hat d_{A,\Phi}^1(a,\phi)\right)_{L^2(X)}}_{\textrm{complex}}.
\end{multline}
Because $\hat d_{A,\Phi}^1$ in \eqref{eq:Itoh_1985_2-6_SO3_monopole_complex_Kaehler} is not diagonal respect to the real and complex inner products on subspaces of its domain or codomain, we need to take the real parts of complex inner products appearing on the right and left-hand sides of \eqref{eq:Itoh_1985_2-6_SO3_monopole_complex_Kaehler_L2_adjoint} in order to give a meaningful definition of the $L^2$ adjoint $\hat d_{A,\Phi}^{1,*}$ of $\hat d_{A,\Phi}^1$; see our discussion around \eqref{eq:Adjoint_operator_real_to_complex_inner_product_space}. We now make the 

\begin{claim}[Variational identities defining $\hat d_{A,\Phi}^{1,*}$]
\label{claim:Itoh_1985_2-7_and_8_SO3_monopole_complex_Kaehler}
Continue the hypotheses and notation in the proof of Proposition \ref{prop:Itoh_1985_proposition_2-3_SO3_monopole_complex_Kaehler}. Then for each $(v,\nu) \in \Omega^+(\su(E))\oplus\Omega^0(W_\can^-\otimes E)$, the following identities hold for all $(a,\phi) \in \Omega^1(\su(E))\oplus\Omega^0(W_\can^+\otimes E)$:
\begin{subequations}
  \label{eq:Itoh_1985_2-7_and_8_SO3_monopole_complex_Kaehler}
  \begin{align}
    \label{eq:Itoh_1985_2-7_SO3_monopole_complex_Kaehler_20_02_slE_forms}
    \frac{1}{2}\left(\hat d_{A,\Phi}^{1,*}(v'+v''), (a,\phi)\right)_{L^2(X)}
    &= \Real(\bar\partial_A^*v'',a'')_{L^2(X)} - \Real(v''\varphi, \tau)_{L^2(X)}
    \\
    \notag
    &\quad - \Real(v''\sigma, \psi)_{L^2(X)},
    \\
    \label{eq:Itoh_1985_2-8_SO3_monopole_complex_Kaehler_11omega_slE_forms}
    \frac{1}{2}\left(\hat d_{A,\Phi}^{1,*}\xi, (a,\phi)\right)_{L^2(X)}
    &= -\Imag(\partial_A\xi, a')_{L^2(X)} - \Imag(\xi\varphi, \sigma)_{L^2(X)}
    \\
    \notag
    &\quad + \Imag(\xi\psi, \tau)_{L^2(X)},
    \\
    \label{eq:Itoh_1985_2-8_negative_spinor_SO3_monopole_complex_Kaehler_L2}
    \left(\hat d_{A,\Phi}^{1,*}\nu, (a,\phi)\right)_{L^2(X)}
    &= (\bar\partial_A^*\nu,\sigma)_{L^2(X)} + (\bar\partial_A\nu, \tau)_{L^2(X)} + (\nu, a''\varphi)_{L^2(X)}
    \\
    \notag
    &\qquad - (\nu, \star(a'\wedge\star\psi))_{L^2(X)}.                     
\end{align}
\end{subequations}
\end{claim}

\begin{proof}[Proof of Claim \ref{claim:Itoh_1985_2-7_and_8_SO3_monopole_complex_Kaehler}]
To verify \eqref{eq:Itoh_1985_2-7_SO3_monopole_complex_Kaehler_20_02_slE_forms}, we begin by noting that
\begin{align*}
  (\hat d_{A,\Phi}^{1,*}v'', (a,\phi))_{L^2(X)}
  &= (v'', \hat d_{A,\Phi}^1(a,\phi))_{L^2(X)} \quad\text{(by
  \eqref{eq:Itoh_1985_2-6_SO3_monopole_complex_Kaehler_L2_adjoint})}
  \\
  &= \left(v'', \bar\partial_Aa'' - \left(\tau\otimes\varphi^* + \psi\otimes\sigma^*\right)_0\right)_{L^2(X)}
    \quad\text{(by \eqref{eq:Itoh_1985_2-6_SO3_monopole_complex_Kaehler} and \eqref{eq:Almost_hermitian_manifold_pq_forms_pointwise_orthogonal_rs_forms_unless_pq_equals_rs})}
  \\
  &= (\bar\partial_A^*v'', a'')_{L^2(X)} - (v'', \tau\otimes\varphi^*)_{L^2(X)} - (v'', \psi\otimes\sigma^*)_{L^2(X)}.
\end{align*}
For all $v''\in\Omega^{0,2}(\fsl(E))$, and $\tau\in\Omega^{0,2}(E)$, and $\varphi\in\Omega^0(E)$, a modification of the derivation of \eqref{eq:Inner_product_tau_and_zeta_psi_equals_tau_otimes_bar_psi_star_and_star_zeta} yields
\begin{equation}
\label{eq:_Inner_product_v''_and_tau_otimes_bar_varphi_equals_v''varphi_and_tau}  
  \langle v'', \tau\otimes\varphi^* \rangle_{\Lambda^{0,2}(\gl(E))}
  =
  \langle v''\varphi, \tau\rangle_{\Lambda^{0,2}(E)}.
\end{equation}
Indeed, to prove \eqref{eq:_Inner_product_v''_and_tau_otimes_bar_varphi_equals_v''varphi_and_tau}, we compute that
\[
  \langle v'', (\tau\otimes\varphi^*)_0 \rangle_{\Lambda^{0,2}(\fsl(E))}
  =
  \langle v'', (\tau\otimes\varphi^*)_0 \rangle_{\Lambda^{0,2}(\gl(E))}
  =
  \langle v'', \tau\otimes\varphi^* \rangle_{\Lambda^{0,2}(\gl(E))}
\]
and for elementary tensors $v'' = \beta\otimes\zeta$ and $\tau = \alpha\otimes s$ for $\alpha,\beta\in\Omega^{0,2}(X)$ and $s \in \Omega^0(E)$ and $\zeta \in \Omega^0(\fsl(E))$, and recalling that $\varphi^* = \langle\cdot,\varphi\rangle_E \in \Omega^0(E^*)$ by \eqref{eq:Hermitian_duals_sections_LambdapqE}, we see
\begin{align*}
  \langle v'', \tau\otimes\varphi^* \rangle_{\Lambda^{0,2}(\gl(E))}
  &=
  \langle \beta\otimes\zeta, \alpha\otimes s\otimes\varphi^* \rangle_{\Lambda^{0,2}(\gl(E))}
  \\
  &= \langle \beta, \alpha\rangle_{\Lambda^{0,2}(X)}
    \langle \zeta, s\otimes\varphi^* \rangle_{\gl(E)}
  \\
  &= \langle \beta, \alpha\rangle_{\Lambda^{0,2}(X)}
    \langle \zeta\varphi, s \rangle_E \quad\text{(by \eqref{eq:Inner_product_sigma_otimes_bar_varphi_star_and_zeta_equals_sigma_and_zeta_varphi})}
  \\
  &= \langle \beta\otimes\zeta\varphi, \alpha\otimes s\rangle_{\Lambda^{0,2}(E)}
  \\
  &= \langle v''\varphi, \tau\rangle_{\Lambda^{0,2}(E)},    
\end{align*}
as claimed. (In order to apply \eqref{eq:Inner_product_sigma_otimes_bar_varphi_star_and_zeta_equals_sigma_and_zeta_varphi} above, we may write 
\[
  \langle \zeta, s\otimes\varphi^* \rangle_{\gl(E)} = \overline{\langle \zeta, s\otimes\varphi^* \rangle}_{\gl(E)} = \overline{\langle s, \zeta\varphi \rangle}_E = \langle \zeta\varphi, s \rangle_E,
\]
as indicated.) By applying \eqref{eq:_Inner_product_v''_and_tau_otimes_bar_varphi_equals_v''varphi_and_tau} to the right-hand side of the previous defining identity for $\hat d_{A,\Phi}^{1,*}v''$, we obtain
\[
  (\hat d_{A,\Phi}^{1,*}v'', (a,\phi))_{L^2(X)}
  =
  (\bar\partial_A^*v'', a'')_{L^2(X)} - (v''\varphi, \tau)_{L^2(X)} - (v''\sigma, \psi)_{L^2(X)},
\]
for all $(a,\phi)$. Next, we observe that
\begin{align*}
  (\hat d_{A,\Phi}^{1,*}v', (a,\phi))_{L^2(X)}
  &= (v', \hat d_{A,\Phi}^1(a,\phi))_{L^2(X)}
    \quad\text{(by \eqref{eq:Itoh_1985_2-6_SO3_monopole_complex_Kaehler_L2_adjoint})}
  \\
  &= \left(v', \partial_Aa' + \left(\varphi\otimes\tau^* + \sigma\otimes\psi^*\right)_0\right)_{L^2(X)}
    \quad\text{(by \eqref{eq:Itoh_1985_2-6_SO3_monopole_complex_Kaehler} and \eqref{eq:Almost_hermitian_manifold_pq_forms_pointwise_orthogonal_rs_forms_unless_pq_equals_rs})}
  \\
  &= (\partial_A^*v', a')_{L^2(X)} + (v', \varphi\otimes\tau^*)_{L^2(X)} + (v', \sigma\otimes\psi^*)_{L^2(X)}.
\end{align*}
For all $v'\in\Omega^{2,0}(\fsl(E))$, $\tau\in\Omega^{0,2}(E)$, and $\varphi\in\Omega^0(E)$, we can show that
\begin{equation}
\label{eq:_Inner_product_v'_and_varphi_otimes_bar_tau_equals_tau_and_v''varphi}  
  \langle v', \varphi\otimes\tau^* \rangle_{\Lambda^{2,0}(\gl(E))}
  =
  -\langle v''\varphi, \tau\rangle_{\Lambda^{0,2}(E)}.
\end{equation}
Indeed, to prove \eqref{eq:_Inner_product_v'_and_varphi_otimes_bar_tau_equals_tau_and_v''varphi}, we observe that
\[
  \langle v', (\varphi\otimes\tau^*)_0 \rangle_{\Lambda^{2,0}(\fsl(E))}
  =
  \langle v', (\varphi\otimes\tau^*)_0 \rangle_{\Lambda^{2,0}(\gl(E))}
  =
  \langle v', \varphi\otimes\tau^* \rangle_{\Lambda^{2,0}(\gl(E))}
\]
and
\begin{align*}
  \langle v', \varphi\otimes\tau^* \rangle_{\Lambda^{2,0}(\gl(E))}
  &=
    \langle (\varphi\otimes\tau^*)^\dagger, (v')^\dagger \rangle_{\Lambda^{0,2}(\gl(E))}
    \quad\text{(by \eqref{eq:Inner_product_M_dagger_N_equals_inner_product_N_dagger_M})}
  \\
  &=
    -\langle \tau\otimes\varphi^*, v'' \rangle_{\Lambda^{0,2}(E)}
    \quad\text{(by \eqref{eq:Itoh_1985_p_850_reality_condition_su(E)_2-forms})}
  \\
  &=
  -\langle \tau, v''\varphi \rangle_{\Lambda^{0,2}(E)} \quad\text{(by \eqref{eq:_Inner_product_v''_and_tau_otimes_bar_varphi_equals_v''varphi_and_tau})},
\end{align*}
as claimed. By applying \eqref{eq:_Inner_product_v'_and_varphi_otimes_bar_tau_equals_tau_and_v''varphi} to the right-hand side of the previous defining identity for $\hat d_{A,\Phi}^{1,*}v'$, we obtain
\[
  (\hat d_{A,\Phi}^{1,*}v', (a,\phi))_{L^2(X)}
  =
  (\partial_A^*v', a')_{L^2(X)} - (\tau, v''\varphi)_{L^2(X)} - (\psi, v''\sigma)_{L^2(X)},
\]
for all $(a,\phi)$. By combining the defining identities for $\hat d_{A,\Phi}^{1,*}v'$ and $\hat d_{A,\Phi}^{1,*}v''$, we see that
\begin{multline*}
  (\hat d_{A,\Phi}^{1,*}(v'+v''), (a,\phi))_{L^2(X)}
  =
  (\bar\partial_A^*v'', a'')_{L^2(X)} + (\partial_A^*v', a')_{L^2(X)}
  \\
  - (v''\varphi, \tau)_{L^2(X)} - (\tau, v''\varphi)_{L^2(X)}
  \\
  - (v''\sigma, \psi)_{L^2(X)} - (\psi, v''\sigma)_{L^2(X)},
\end{multline*}
for all $(a,\phi)\in\Omega^1(\su(E))\oplus\Omega^0(W_\can^+\otimes E)$. Therefore,
\begin{multline*}
   (\hat d_{A,\Phi}^{1,*}(v'+v''), (a,\phi))_{L^2(X)}
  =
  (\bar\partial_A^*v'', a'')_{L^2(X)} + (\partial_A^*v', a')_{L^2(X)}
  \\
  - 2\Real(v''\varphi, \tau)_{L^2(X)}  - 2\Real(\psi, v''\sigma)_{L^2(X)},
\end{multline*}
for all $(a,\phi)\in\Omega^1(\su(E))\oplus\Omega^0(W_\can^+\otimes E)$. But
\begin{align*}
  (\partial_A^*v', a')_{L^2(X)}
  &= ((a')^\dagger, (\partial_A^*v')^\dagger)_{L^2(X)}
    \quad\text{(by \eqref{eq:Inner_product_M_dagger_N_equals_inner_product_N_dagger_M})}
  \\
  &= -(a'', \bar\partial_A^*((v')^\dagger))_{L^2(X)}
    \quad\text{(by \eqref{eq:AdjointCommute_dagger_and_del_glE})}
  \\
  &= (a'', \bar\partial_A^*v'')_{L^2(X)}
    \quad\text{(by \eqref{eq:Itoh_1985_p_850_reality_condition_su(E)_2-forms})}
  \\
  &= \overline{(\bar\partial_A^*v'',a'')}_{L^2(X)},
\end{align*}
and thus
\begin{align*}
  (\bar\partial_A^*v'', a'')_{L^2(X)} + (\partial_A^*v', a')_{L^2(X)}
  &=
    (\bar\partial_A^*v'', a'')_{L^2(X)} + \overline{(\bar\partial_A^*v'',a'')}_{L^2(X)}
  \\
  &=
  2\Real(\bar\partial_A^*v'', a'')_{L^2(X)}.
\end{align*}
Substituting the preceding identity gives
\begin{align*}
  (\hat d_{A,\Phi}^{1,*}(v'+v''), (a,\phi))_{L^2(X)}
  &= 2\Real(\bar\partial_A^*v'', a'')_{L^2(X)} - 2\Real(v''\varphi, \tau)_{L^2(X)} - 2\Real(\psi, v''\sigma)_{L^2(X)}
  \\
  &= 2\Real(\bar\partial_A^*v'', a'')_{L^2(X)} - 2\Real(v''\varphi, \tau)_{L^2(X)} - 2\Real(v''\sigma, \psi)_{L^2(X)},
\end{align*}    
for all $(a,\phi)\in\Omega^1(\su(E))\oplus\Omega^0(W_\can^+\otimes E)$, and this completes our verification of \eqref{eq:Itoh_1985_2-7_SO3_monopole_complex_Kaehler_20_02_slE_forms}.

To verify \eqref{eq:Itoh_1985_2-8_SO3_monopole_complex_Kaehler_11omega_slE_forms}, we note that $\xi \in \Omega^0(\su(E))$ and observe that
\begin{align*}
  {}&(\hat d_{A,\Phi}^{1,*}\xi, (a,\phi))_{L^2(X)}
  \\
    &= (\xi, \hat d_{A,\Phi}^1(a,\phi))_{L^2(X)}
      \quad\text{(by \eqref{eq:Itoh_1985_2-6_SO3_monopole_complex_Kaehler_L2_adjoint})}
  \\
    &= (\xi, \Lambda(\bar\partial_Aa' + \partial_Aa''))_{L^2(X)}
       - \left(\xi, i\left(\sigma\otimes\varphi^* + \varphi\otimes\sigma^*\right)_0\right)_{L^2(X)}
  \\
  &\quad + \left(\xi, i\star\left(\tau\otimes\psi^* + \psi\otimes\tau^*\right)_0\right)_{L^2(X)}
  \quad\text{(by \eqref{eq:Itoh_1985_2-6_SO3_monopole_complex_Kaehler})}
  \\
  &= (\xi, [\Lambda,\bar\partial_A]a' + [\Lambda,\partial_A]a'')_{L^2(X)}                                            
  + \left(i\xi, \sigma\otimes\varphi^*\right)_{L^2(X)}
      + \left(i\xi, \varphi\otimes\sigma^*\right)_{L^2(X)}
  \\
  &\quad - \left(i\xi, \star(\tau\otimes\psi^*)\right)_{L^2(X)}
      - \left(i\xi, \star(\psi\otimes\tau^*)\right)_{L^2(X)},
\end{align*}
for all $(a,\phi)\in\Omega^1(\su(E))\oplus\Omega^0(W_\can^+\otimes E)$. Therefore, 
\begin{align*}
  {}&(\hat d_{A,\Phi}^{1,*}\xi, (a,\phi))_{L^2(X)}
  \\
  &= (\xi, -i\partial_A^*a' + i\bar\partial_A^*a'')_{L^2(X)}
   + \left(i\xi\varphi, \sigma\right)_{L^2(X)}
   + \left(i\xi\sigma, \varphi\right)_{L^2(X)}
  \\
  &\quad - \left(i\xi\psi, \tau\right)_{L^2(X)}
      - \left(i\xi\tau, \psi)\right)_{L^2(X)}
    \quad\text{(by \eqref{eq:Kaehler_identity_commutator_Lambda_del-bar_A_and_Lambda_del_A},
    \eqref{eq:Inner_product_sigma_otimes_bar_varphi_star_and_zeta_equals_sigma_and_zeta_varphi}, and
    \eqref{eq:Inner_product_tau_and_zeta_psi_equals_star_tau_otimes_bar_psi_star_and_zeta})}
  \\
  &= i(\partial_A\xi, a')_{L^2(X)} - i(\bar\partial_A\xi,a'')_{L^2(X)}
       + \left(i\xi\varphi, \sigma\right)_{L^2(X)}
   + \left(\sigma, i\xi\varphi\right)_{L^2(X)}
  \\
  &\quad - \left(i\xi\psi, \tau\right)_{L^2(X)}
    - \left(\tau, i\xi\psi)\right)_{L^2(X)}
    \quad\text{(as $(i\xi)^\dagger = i\xi$ for $\xi\in\Omega^0(\su(E))$}
  \\
    &= (i\partial_A\xi, a')_{L^2(X)} - (i\bar\partial_A\xi, a'')_{L^2(X)}
      + 2\Real(i\xi\varphi, \sigma)_{L^2(X)} - 2\Real(i\xi\psi, \tau)_{L^2(X)},
\end{align*}
for all $(a,\phi)\in\Omega^1(\su(E))\oplus\Omega^0(W_\can^+\otimes E)$. Now
\begin{align*}
  (i\partial_A\xi, a')_{L^2(X)} - (i\bar\partial_A\xi, a'')_{L^2(X)}
  &=
  (i\partial_A\xi, a')_{L^2(X)} - ((a'')^\dagger, (i\bar\partial_A\xi)^\dagger)_{L^2(X)}
  \quad\text{(by \eqref{eq:Inner_product_M_dagger_N_equals_inner_product_N_dagger_M})}
  \\
  &=
  (i\partial_A\xi, a')_{L^2(X)} + ((a'')^\dagger, i\partial_A(\xi^\dagger))_{L^2(X)}
  \quad\text{(by Lemma \ref{lem:PointwiseHermitianAdjoint})}
  \\
  &=
    (i\partial_A\xi, a')_{L^2(X)} + (a', i\partial_A\xi)_{L^2(X)}
  \\
  &= 2\Real(i\partial_A\xi, a')_{L^2(X)},
\end{align*}
for all $(a,\phi)\in\Omega^1(\su(E))\oplus\Omega^0(W_\can^+\otimes E)$. Substituting the preceding identity yields
\begin{align*}
  (\hat d_{A,\Phi}^{1,*}\xi, (a,\phi))_{L^2(X)}
  &=
    2\Real(i\partial_A\xi, a')_{L^2(X)} + 2\Real(i\xi\varphi, \sigma)_{L^2(X)} - 2\Real(i\xi\psi, \tau)_{L^2(X)}
  \\
  &=
    -2\Imag(\partial_A\xi, a')_{L^2(X)} - 2\Imag(\xi\varphi, \sigma)_{L^2(X)} + 2\Imag(\xi\psi, \tau)_{L^2(X)},
\end{align*}
for all $(a,\phi)$, and this completes our verification of \eqref{eq:Itoh_1985_2-8_SO3_monopole_complex_Kaehler_11omega_slE_forms}.

To verify \eqref{eq:Itoh_1985_2-8_negative_spinor_SO3_monopole_complex_Kaehler_L2}, we observe that
\begin{align*}
  (\hat d_{A,\Phi}^{1,*}\nu, (a,\phi))_{L^2(X)}
  &= (\nu, \hat d_{A,\Phi}^1(a,\phi))_{L^2(X)}
    \quad\text{(by \eqref{eq:Itoh_1985_2-6_SO3_monopole_complex_Kaehler_L2_adjoint})}
  \\
  &= \left(\nu, \bar\partial_A\sigma + \bar\partial_A^*\tau + a''\varphi - \star(a'\wedge\star\psi)\right)_{L^2(X)} \quad\text{(by \eqref{eq:Itoh_1985_2-6_SO3_monopole_complex_Kaehler})}
  \\
  &= (\bar\partial_A^*\nu, \sigma)_{L^2(X)} + (\bar\partial_A\nu, \tau)_{L^2(X)} + (\nu, a''\varphi)_{L^2(X)}
  \\
  &\quad - (\nu, \star(a'\wedge\star\psi))_{L^2(X)},
\end{align*}
for all $(a,\phi)\in\Omega^1(\su(E))\oplus\Omega^0(W_\can^+\otimes E)$, as desired. This completes the proof of Claim \ref{claim:Itoh_1985_2-7_and_8_SO3_monopole_complex_Kaehler}.
\end{proof}

We now make the

\begin{claim}
\label{claim:L2_inner_product_nu_and_a''_varphi_and_star_(a'_wedge_star_psi)}
Continue the hypotheses and notation in the proof of Proposition \ref{prop:Itoh_1985_proposition_2-3_SO3_monopole_complex_Kaehler}. Then for each $\nu \in \Omega^0(W_\can^-\otimes E) = \Omega^{0,1}(E)$, the following identities hold for all $a'' \in \Omega^{0,1}(\fsl(E))$ with $a'  = -(a'')^\dagger\in \Omega^{1,0}(\fsl(E))$:
\begin{subequations}
\label{eq:L2_inner_product_nu_and_a''_varphi_and_star_(a'_wedge_star_psi)}  
\begin{align}
  \label{eq:L2_inner_product_nu_and_a''_varphi}
  (\nu, a''\varphi)_{L^2(X)}
  &= ((\nu\otimes\varphi^*)_0, a'')_{L^2(X)},
  \\
  \label{eq:L2_inner_product_nu_and_star_a''_wedge_psi}
  (\nu, \star(a'\wedge\star\psi))_{L^2(X)}
  &= -(\star((\star\psi^*)\wedge\nu)_0), a')_{L^2(X)}.
\end{align}
\end{subequations}
\end{claim}

\begin{proof}[Proof of Claim \ref{claim:L2_inner_product_nu_and_a''_varphi_and_star_(a'_wedge_star_psi)}]
To simplify the expression $(\nu, a''\varphi)_{L^2(X)}$ in \eqref{eq:L2_inner_product_nu_and_a''_varphi}, consider elementary tensors $\nu = \beta\otimes s$ and $a'' = \alpha\otimes\zeta$, where $s\in\Omega^0(E)$ and $\zeta \in \Omega^0(\fsl(E))$ and $\alpha,\beta \in \Omega^{0,1}(X)$, and note that by \eqref{eq:L2_inner_product_pq_forms_E} we have
\[
  (\nu, a''\varphi)_{L^2(X)} = \int_X \star\langle\nu, a''\varphi\rangle_{\Lambda^{0,1}(E)}\,d\vol.
\]
We obtain
\begin{align*}
  \langle a''\varphi, \nu\rangle_{\Lambda^{0,1}(E)}
  &= \langle \alpha\otimes\zeta \varphi, \beta\otimes s\rangle_{\Lambda^{0,1}(E)}
  \\
  &= \langle \alpha, \beta\rangle_{\Lambda^{0,1}(X)}\langle \zeta \varphi, s\rangle_E
    \quad\text{(by \eqref{eq:Huybrechts_definition_4-1-11_and_example_4-1-2})}
  \\
  &= \langle \alpha, \beta\rangle_{\Lambda^{0,1}(X)}\langle \zeta, s\otimes \varphi^*\rangle_{\gl(E)}
    \quad\text{(by \eqref{eq:Hermitian_duals_sections_LambdapqE} and \eqref{eq:Inner_product_sigma_otimes_bar_varphi_star_and_zeta_equals_sigma_and_zeta_varphi})}
  \\
  &= \langle \alpha\otimes\zeta, \beta\otimes s\otimes \varphi^*\rangle_{\Lambda^{0,1}(\gl(E))}
    \quad\text{(by \eqref{eq:Huybrechts_definition_4-1-11_and_example_4-1-2})}
  \\
  &= \langle a'', \nu\otimes \varphi^*\rangle_{\Lambda^{0,1}(\gl(E))}
  \\
  &= \langle a'', (\nu\otimes \varphi^*)_0\rangle_{\Lambda^{0,1}(\gl(E))}
  \\
  &= \langle a'', (\nu\otimes \varphi^*)_0\rangle_{\Lambda^{0,1}(\fsl(E))},
\end{align*}  
and so \eqref{eq:L2_inner_product_nu_and_a''_varphi} follows.

To simplify the expression $(\nu, \star(a'\wedge\star\psi))_{L^2(X)}$ in \eqref{eq:L2_inner_product_nu_and_star_a''_wedge_psi}, consider an elementary tensor $\psi = \varpi\otimes t$, where $t\in\Omega^0(E)$ and $\varpi \in \Omega^{0,2}(X)$, and observe that by \eqref{eq:L2_inner_product_pq_forms_E}
\[
  (\nu, \star(a'\wedge\star\psi))_{L^2(X)} = \int_X \star\langle\nu, \star(a'\wedge\star\psi)\rangle_{\Lambda^{0,1}(E)}\,d\vol.
\]
Noting that
\[
  \psi^* = \langle\cdot,\psi\rangle_{\Lambda^{0,2}(E)} = \langle\cdot,\varpi\otimes t\rangle_{\Lambda^{0,2}(E)} = \bar\varpi\otimes \langle\cdot, t\rangle_E = \bar\varpi\otimes t^*
\]
by \eqref{eq:Complex_conjugate_sections_E_and_Lambda02E} and \eqref{eq:Hermitian_duals_sections_LambdapqE}, and using \eqref{eq:Kobayashi_7-6-11} to give
\[
  a' = -(a'')^\dagger = -\overline{a''}^\intercal = -\overline{\alpha\otimes\zeta^\intercal} = -\bar\alpha\otimes\bar\zeta^\intercal = -\bar\alpha\otimes\zeta^\dagger \in \Omega^{1,0}(\fsl(E)),
\]
we obtain
\begin{align*}
  \left\langle \star(a'\wedge\star\psi), \nu\right\rangle_{\Lambda^{0,1}(E)}
  &= -\left\langle \star((\bar\alpha\otimes\zeta^\dagger)\wedge\star(\varpi\otimes t)), \beta\otimes s\right\rangle_{\Lambda^{0,1}(E)}
  \\  
  &= -\left\langle \star(\bar\alpha\wedge\star\varpi)\otimes \zeta^\dagger t, \beta\otimes s\right\rangle_{\Lambda^{0,1}(E)}
  \\  
  &= -\left\langle \star(\bar\alpha\wedge\star\varpi), \beta\right\rangle_{\Lambda^{0,1}(X)}
    \langle \zeta^\dagger t, s\rangle_E
    \quad\text{(by \eqref{eq:Huybrechts_definition_4-1-11_and_example_4-1-2})}
  \\
  &= \left\langle \bar\alpha\wedge\star\varpi, \star\beta\right\rangle_{\Lambda^{1,2}(X)}
    \left\langle \zeta^\dagger t, s\right\rangle_E
    \quad\text{(by \eqref{eq:Huybrechts_proposition_1-2-20_ii} with $k=1$ and $d=4$)}
  \\
  &= \left\langle \bar\alpha\wedge\star\varpi, \star\beta\right\rangle_{\Lambda^{1,2}(X)}
    \left\langle \zeta^\dagger, s\otimes t^*\right\rangle_{\gl(E)}
    \quad\text{(by \eqref{eq:Hermitian_duals_sections_LambdapqE} and \eqref{eq:Inner_product_sigma_otimes_bar_varphi_star_and_zeta_equals_sigma_and_zeta_varphi})}
  \\
  &= \star(\bar\alpha\wedge\star\varpi\wedge \star (\overline{\star\beta}))
    \left\langle \zeta^\dagger, s\otimes t^*\right\rangle_{\gl(E)}
    \quad\text{(by \eqref{eq:Huybrechts_proposition_1-2-20_i} and \eqref{eq:Huybrechts_page_33_inner_product_complex-valued_forms})}
  \\
  &= \star(\bar\alpha\wedge\star\varpi\wedge \star^2 \bar\beta)
    \left\langle \zeta^\dagger, s\otimes t^*\right\rangle_{\gl(E)}
    \quad\text{(by \eqref{eq:Huybrechts_page_33_and_lemma_1-2-24})}
  \\
  &= -\star(\bar\alpha\wedge\star\varpi\wedge \bar\beta)
    \left\langle \zeta^\dagger, s\otimes t^*\right\rangle_{\gl(E)}
    \quad\text{(by \eqref{eq:Huybrechts_proposition_1-2-20_iii} with $k=1$ and $d=4$)}
  \\
  &= \star(\bar\alpha\wedge \star^2(\star\varpi\wedge \bar\beta))
    \left\langle \zeta^\dagger, s\otimes t^*\right\rangle_{\gl(E)}
    \quad\text{(by \eqref{eq:Huybrechts_proposition_1-2-20_iii} with $k=3$ and $d=4$)}
  \\
  &= \left\langle\bar\alpha, \overline{\star(\star\varpi\wedge \bar\beta)} \right\rangle_{\Lambda^{1,0}(X)}
    \left\langle \zeta^\dagger, s\otimes t^*\right\rangle_{\gl(E)}
    \quad\text{(by \eqref{eq:Huybrechts_proposition_1-2-20_i} and \eqref{eq:Huybrechts_page_33_inner_product_complex-valued_forms})}
  \\
  &= \left\langle\bar\alpha, \star(\star\bar\varpi\wedge \beta) \right\rangle_{\Lambda^{1,0}(X)}
    \left\langle \zeta^\dagger, s\otimes t^*\right\rangle_{\gl(E)}
    \quad\text{(by \eqref{eq:Huybrechts_page_33_and_lemma_1-2-24})}
  \\
  &= \left\langle\bar\alpha\otimes \zeta^\dagger,
    \star(\star\bar\varpi\wedge \beta)\otimes s\otimes t^* \right\rangle_{\Lambda^{1,0}(\gl(E))}
    \quad\text{(by \eqref{eq:Huybrechts_definition_4-1-11_and_example_4-1-2})}
  \\
  &= \left\langle\bar\alpha\otimes \zeta^\dagger,
    \star((\star\bar\varpi\otimes t^*)\wedge (\beta\otimes s)) \right\rangle_{\Lambda^{1,0}(\gl(E))}
  \\
  &= -\left\langle a', \star((\star\psi^*)\wedge\nu)_0) \right\rangle_{\Lambda^{1,0}(\fsl(E))}, 
\end{align*}  
and so \eqref{eq:L2_inner_product_nu_and_star_a''_wedge_psi} follows. This completes the proof of Claim \ref{claim:L2_inner_product_nu_and_a''_varphi_and_star_(a'_wedge_star_psi)}.
\end{proof}

Suppose now that $(v,\nu) \in \bH_{A,\Phi}^2$, so that, using \eqref{eq:Donaldson_Kronheimer_lemma_2-1-57_sum_3_terms} to write $v=\frac{1}{2}(v'+v''+\xi\otimes\omega)$,
\begin{equation}
\label{eq:d_APhi^1star_v+v''+Lxi_is_zero}
\hat d_{A,\Phi}^{1,*}(v, \nu) = 0.
\end{equation}
By the reality condition \eqref{eq:Itoh_1985_p_850_reality_condition_su(E)_2-forms} for any $v \in \Omega^+(\su(E))$, we have
\[
v' = -(v'')^\dagger.
\]
The variational equation \eqref{eq:Itoh_1985_2-7_and_8_SO3_monopole_complex_Kaehler} implicitly defines the elliptic system $\hat d_{A,\Phi}^{1,*}(v,\nu)=0$ for $(v,\nu)$ and hence the corresponding elliptic system for $(v'',\nu,\xi)$. We now seek to write the latter system explicitly and hence verify the isomorphisms \eqref{eq:Itoh_1985_proposition_2-3_SO3_monopole_complex_Kaehler_isomorphism} and \eqref{eq:Itoh_1985_proposition_2-3_SO3_monopole_complex_Kaehler_isomorphism_type1}. First, suppose in \eqref{eq:Itoh_1985_2-7_and_8_SO3_monopole_complex_Kaehler} that $a''=0$ (and thus
$a'=-(a'')^\dagger=0$ by \eqref{eq:Kobayashi_7-6-11}) and $\sigma = 0$, while $\tau \in \Omega^{0,2}(E)$ is allowed to be arbitrary. By combining the identity $d_{A,\Phi}^{1,*}(v,\nu)=0$ and the reality condition inherited from \eqref{eq:Itoh_1985_2-6_SO3_monopole_complex_Kaehler_L2_adjoint} with the system of variational identities \eqref{eq:Itoh_1985_2-7_and_8_SO3_monopole_complex_Kaehler} we obtain
\[
  - \Real(v''\varphi, \tau)_{L^2(X)} + \Imag(\xi\psi, \tau)_{L^2(X)} + \Real(\bar\partial_A\nu, \tau)_{L^2(X)} = 0, \quad\text{for all } \tau \in \Omega^{0,2}(E),
\]
which is equivalent to
\begin{equation}
  \label{eq:d_APhi^1*_(v,nu)_tau_component_is_zero_variational}
  - \Real(v''\varphi, \tau)_{L^2(X)} - \Real(i\xi\psi, \tau)_{L^2(X)} + \Real(\bar\partial_A\nu, \tau)_{L^2(X)} = 0, \quad\text{for all } \tau \in \Omega^{0,2}(E).
\end{equation}
By \eqref{eq:Equivalence_real_complex_inner_product_equations}, the variational equation \eqref{eq:d_APhi^1*_(v,nu)_tau_component_is_zero_variational} is equivalent to
\[
  (\bar\partial_A\nu - v''\varphi - i\xi\psi, \tau)_{L^2(X)} = 0, \quad\text{for all } \tau \in \Omega^{0,2}(E).
\]
The preceding variational equation is in turn equivalent to
\begin{equation}
  \label{eq:1}
  \bar\partial_A\nu - v''\varphi - i\xi\psi = 0 \in \Omega^{0,2}(E).
\end{equation}
Second, suppose that $a''=0$ (and thus $a'=-(a'')^\dagger=0$) and $\tau = 0$ in \eqref{eq:Itoh_1985_2-7_and_8_SO3_monopole_complex_Kaehler}, while $\sigma \in \Omega^0(E)$ is arbitrary. By combining the identity $\hat d_{A,\Phi}^{1,*}(v,\nu)=0$ and the reality condition inherited from \eqref{eq:Itoh_1985_2-6_SO3_monopole_complex_Kaehler_L2_adjoint} with the system of variational identities \eqref{eq:Itoh_1985_2-7_and_8_SO3_monopole_complex_Kaehler} we now obtain
\[
   -\Real(v''\sigma, \psi)_{L^2(X)} - \Imag(\xi\varphi, \sigma)_{L^2(X)} + \Real (\bar\partial_A^*\nu,\sigma)_{L^2(X)} = 0, \quad\text{for all } \sigma \in \Omega^0(E),
\]
which is equivalent to
\begin{equation}
  \label{eq:d_APhi^1*_(v,nu)_sigma_component_is_zero_variational}
  \Real(\psi, v''\sigma)_{L^2(X)} - \Real(i\xi\varphi, \sigma)_{L^2(X)} - \Real (\bar\partial_A^*\nu,\sigma)_{L^2(X)} = 0, \quad\text{for all } \sigma \in \Omega^0(E).
\end{equation}
By \eqref{eq:Equivalence_real_complex_inner_product_equations}, the variational equation \eqref{eq:d_APhi^1*_(v,nu)_sigma_component_is_zero_variational} is equivalent to
\[
  (\psi, v''\sigma)_{L^2(X)} - (i\xi\varphi, \sigma)_{L^2(X)} - (\bar\partial_A^*\nu,\sigma)_{L^2(X)} = 0, \quad\text{for all } \sigma \in \Omega^0(E).
\]
By applying the forthcoming identity \eqref{eq:v''sigma_innerproduct_psi_equals_sigma_innerproduct_other} (and reversing the arguments in the inner products on both sides of \eqref{eq:v''sigma_innerproduct_psi_equals_sigma_innerproduct_other})
we see that the preceding variational equation is equivalent to
\[
  \left(\star\left((v'')^\dagger\wedge \star \psi\right), \sigma\right)_{L^2(X)} - (i\xi\varphi, \sigma)_{L^2(X)} - (\bar\partial_A^*\nu,\sigma)_{L^2(X)} = 0, \quad\text{for all } \sigma \in \Omega^0(E).
\]
The preceding variational equation is in turn equivalent to
\begin{equation}
  \label{eq:2}
  \bar\partial_A^*\nu - \star\left((v'')^\dagger\wedge \star \psi\right) + i\xi\varphi
  = 0 \in \Omega^0(E).
\end{equation}
To simplify the expression $(v''\sigma, \psi)_{L^2(X)}$ in the identity \eqref{eq:d_APhi^1*_(v,nu)_sigma_component_is_zero_variational}, consider elementary tensors $v'' = \chi\otimes\eta$ and $\psi = \varpi\otimes s$, where $\eta \in \Omega^0(\fsl(E))$ and $\chi, \varpi \in \Omega^{0,2}(X)$ and $s \in \Omega^0(E)$, and note by \eqref{eq:L2_inner_product_pq_forms_E} that
\[
  (v''\sigma, \psi)_{L^2(X)} = \int_X \star\langle v''\sigma, \psi \rangle_{\Lambda^{0,2}(E)}\,d\vol,
\]
where $\sigma \in \Omega^0(E)$. Because
\[
  (v'')^\intercal = (\chi\otimes\eta)^\intercal = \chi\otimes\eta^\intercal \in \Omega^{2,0}(\fsl(E)),
\]
we then obtain
\begin{align*}
  \langle v''\sigma, \psi\rangle_{\Lambda^{0,2}(E)}
  &= \langle \chi\otimes\eta \sigma, \varpi\otimes s\rangle_{\Lambda^{0,2}(E)}
  \\
  &= \langle \chi, \varpi\rangle_{\Lambda^{0,2}(X)}\langle \eta \sigma, s\rangle_E
    \quad\text{(by \eqref{eq:Huybrechts_definition_4-1-11_and_example_4-1-2})}
  \\
  &= \langle \chi, \varpi\rangle_{\Lambda^{0,2}(X)}\langle \sigma, \eta^\dagger s\rangle_E
  \\
  &= \left\langle \sigma, \overline{\langle \chi, \varpi\rangle}_{\Lambda^{0,2}(X)}\eta^\dagger s\right\rangle_E  
  \\
  &= \left\langle \sigma, \langle \bar\chi, \bar\varpi\rangle_{\Lambda^{0,2}(X)}\eta^\dagger s\right\rangle_E
    \quad\text{(by \eqref{eq:Complex_conjugation_and_inner_product})}
  \\
  &= \left\langle \sigma, \star(\bar\chi\wedge\star\varpi)\, \eta^\dagger s\right\rangle_E
    \quad\text{(by \eqref{eq:Huybrechts_proposition_1-2-20_i} and \eqref{eq:Huybrechts_page_33_inner_product_complex-valued_forms})}
  \\
  &= \left\langle \sigma, \star\left((\bar\chi\otimes \eta^\dagger) \wedge (\star\varpi\otimes s)\right) \right\rangle_E
  \\
  &\qquad \text{(where ``$\wedge$'' combines $\wedge$ on $\Lambda^\bullet(X,\CC)$ and evaluation $\gl(E)\times E \to E$)}
  \\
  &= \left\langle \sigma, \star\left((\chi\otimes \eta)^\dagger \wedge \star(\varpi\otimes s)\right)\right\rangle_E
  \\
  &= \left\langle \sigma, \star\left((v'')^\dagger\wedge \star \psi\right)\right \rangle_E.
\end{align*}
Therefore,
\begin{multline}
  \label{eq:v''sigma_innerproduct_psi_equals_sigma_innerproduct_other}
  \langle v''\sigma, \psi\rangle_{\Lambda^{0,2}(E)}
  =
  \left\langle \sigma, \star\left((v'')^\dagger\wedge \star \psi\right)\right \rangle_E,
  \\
  \text{for all } v'' \in \Omega^{0,2}(\fsl(E)), \sigma \in \Omega^0(E), \text{ and } \psi \in \Omega^{0,2}(E).
\end{multline}
Third, suppose that $a''\in \Omega^{0,1}(\fsl(E))$ (and thus $a' = -(a'')^\dagger \in \Omega^{1,0}(\fsl(E))$) is arbitrary in \eqref{eq:Itoh_1985_2-7_and_8_SO3_monopole_complex_Kaehler} while $\tau = 0$ and $\sigma = 0$. By combining the identity $\hat d_{A,\Phi}^{1,*}(v,\nu)=0$ and the reality condition inherited from \eqref{eq:Itoh_1985_2-6_SO3_monopole_complex_Kaehler_L2_adjoint} with the system of variational identities \eqref{eq:Itoh_1985_2-7_and_8_SO3_monopole_complex_Kaehler} we obtain
\begin{multline*}
  \Real(\bar\partial_A^*v'',a'')_{L^2(X)} - \Imag(\partial_A\xi, a')_{L^2(X)} + \Real(\nu, a''\varphi)_{L^2(X)} - \Real(\nu, \star(a'\wedge\star\psi))_{L^2(X)} = 0,
  \\
  \quad\text{for all } a'' \in \Omega^{0,1}(\fsl(E)),
\end{multline*}  
which is equivalent to
\begin{multline}
  \label{eq:d_APhi^1*_(v,nu)_a''_component_is_zero_variational}
  \Real(\bar\partial_A^*v'',a'')_{L^2(X)} + \Real(i\partial_A\xi, a')_{L^2(X)}
  \\
  + \Real(\nu, a''\varphi)_{L^2(X)} - \Real(\nu, \star(a'\wedge\star\psi))_{L^2(X)} = 0,
  \\
  \quad\text{for all } a'' \in \Omega^{0,1}(\fsl(E)).
\end{multline}
To simplify the term $(i\partial_A\xi, a')_{L^2(X)}$ in \eqref{eq:d_APhi^1*_(v,nu)_a''_component_is_zero_variational}, we observe that
\begin{align*}
  \langle i\partial_A\xi, a'\rangle_{\Lambda^{1,0}(\fsl(E))}
  &=
  \langle i\partial_A\xi, a'\rangle_{\Lambda^1(\fsl(E))}
  \\
  &=
  \langle (a')^\dagger, (i\partial_A\xi)^\dagger\rangle_{\Lambda^1(\fsl(E))}
  \quad\text{(by \eqref{eq:Inner_product_M_dagger_N_equals_inner_product_N_dagger_M})}
  \\
  &=
  \langle (-a''), (-i)\bar\partial_A\xi^\dagger\rangle_{\Lambda^1(\fsl(E))}
  \quad\text{(by \eqref{eq:Kobayashi_7-6-11} and \eqref{eq:Commute_dagger_and_del_glE})}
  \\
  &=
  \langle a'', i\bar\partial_A(-\xi)\rangle_{\Lambda^1(\fsl(E))},
\end{align*}
and thus
\[
  \langle i\partial_A\xi, a'\rangle_{\Lambda^{1,0}(\fsl(E))}
  =
  -\langle a'', i\bar\partial_A\xi\rangle_{\Lambda^{0,1}(\fsl(E))}.
\]
Hence, by substituting the preceding identity into \eqref{eq:d_APhi^1*_(v,nu)_a''_component_is_zero_variational} and using symmetry of the real part of the Hermitian inner product, we obtain the equivalent variational equation,
\begin{multline}
  \label{eq:d_APhi^1*_(v,nu)_a''_component_is_zero_variational_simplified}
  \Real(\bar\partial_A^*v'',a'')_{L^2(X)} - \Real(i\bar\partial_A\xi, a'')_{L^2(X)}
  \\
  + \Real(\nu, a''\varphi)_{L^2(X)} - \Real (\nu, \star(a'\wedge\star\psi))_{L^2(X)} = 0,
  \\
  \text{for all } a'' \in \Omega^{0,1}(\fsl(E)).
\end{multline}
By \eqref{eq:Equivalence_real_complex_inner_product_equations}, the preceding variational equation is equivalent to
\begin{multline*}
  (\bar\partial_A^*v'',a'')_{L^2(X)} - (i\bar\partial_A\xi, a'')_{L^2(X)}
  \\
  + (\nu, a''\varphi)_{L^2(X)} - (\nu, \star(a'\wedge\star\psi))_{L^2(X)} = 0,
  \\
  \text{for all } a'' \in \Omega^{0,1}(\fsl(E)).
\end{multline*}
By substituting the variational identities \eqref{eq:L2_inner_product_nu_and_a''_varphi_and_star_(a'_wedge_star_psi)},
we see that the preceding variational equation is equivalent to
\begin{multline*}
  (\bar\partial_A^*v'',a'')_{L^2(X)} - (i\bar\partial_A\xi, a'')_{L^2(X)}
  \\
  + ((\nu\otimes\varphi^*)_0, a'')_{L^2(X)}
  + ((\star((\star\psi^*)\wedge\nu)_0)^\dagger, a'')_{L^2(X)} = 0,
  \\
  \text{for all } a'' \in \Omega^{0,1}(\fsl(E)).
\end{multline*}
The preceding variational equation is in turn equivalent to
\[
  \bar\partial_A^*v'' - i\bar\partial_A\xi + (\nu\otimes\varphi^*)_0
  + (\star((\star\psi^*)\wedge\nu)_0)^\dagger = 0 \in \Omega^{0,1}(\fsl(E)).
\]
By combining the preceding equation with equations \eqref{eq:1} and \eqref{eq:2}, we obtain the system
\begin{subequations}
  \label{eq:4-1_and_2_and_3}
  \begin{align}
  \label{eq:4-3}
  \bar\partial_A^*v'' - i\bar\partial_A\xi + (\nu\otimes\varphi^*)_0
    + (\star((\star\psi^*)\wedge\nu)_0)^\dagger &= 0 \in \Omega^{0,1}(\fsl(E)),
  \\                                                
  \label{eq:4-1}
  \bar\partial_A\nu - v''\varphi - i\xi\psi &= 0 \in \Omega^{0,2}(E),
  \\
  \label{eq:4-2}
  \bar\partial_A^*\nu - \star\left((v'')^\dagger\wedge \star \psi\right) + i\xi\varphi
  &= 0 \in \Omega^0(E).
\end{align}
\end{subequations}
Therefore, we have shown that if
$(v,\nu)\in\bH_{A,\Phi}^2$ where $\Phi=(\varphi,\psi)$ and
\[
  v = \frac{1}{2}(v' + v'') + \frac{1}{2}\xi\otimes\omega
  \in \Omega^{2,0}(\fsl(E)) \oplus \Omega^{0,2}(\fsl(E))\oplus \Omega^0(\su(E))\otimes\omega,
\]
as in \eqref{eq:Donaldson_Kronheimer_lemma_2-1-57_sum_3_terms}, then $(v'',\nu,\xi)$ satisfy
\eqref{eq:4-1_and_2_and_3}. Conversely, if $(v'',\nu,\xi)$ satisfies \eqref{eq:4-1_and_2_and_3},
then the pair $(v,\nu)$ defined by $(v'',\nu,\xi)$ as in \eqref{eq:Donaldson_Kronheimer_lemma_2-1-57_sum_3_terms} can be seen to satisfy $d_{A,\Phi}^{1,*}(v,\nu)=0$ as follows.  The argument yielding \eqref{eq:1} implies that
$d_{A,\Phi}^{1,*}(v,\nu)$ is orthogonal to elements  of the form
\[
(0,0,\tau)\in \Om^1(\su(E))\oplus\Om^0(E)\oplus \Om^{0,2}(E).
\]
The argument yielding \eqref{eq:2} implies that
$d_{A,\Phi}^{1,*}(v,\nu)$ is orthogonal to elements  of the form
\[
(0,\si,0)\in \Om^1(\su(E))\oplus\Om^0(E)\oplus \Om^{0,2}(E).
\]
Finally, the argument yielding \eqref{eq:4-3} implies that
$d_{A,\Phi}^{1,*}(v,\nu)$ is orthogonal to elements of the form
\[
(a,0,0)\in \Om^1(\su(E))\oplus\Om^0(E)\oplus \Om^{0,2}(E).
\]
Because the preceding three types of elements span $\Om^1(\su(E))\oplus\Om^0(E)\oplus \Om^{0,2}(E)$,
we see that if  $(v'',\nu,\xi)$ satisfies \eqref{eq:4-1_and_2_and_3},
then the pair $(v,\nu)$ defined by $(v'',\nu,\xi)$ as in \eqref{eq:Donaldson_Kronheimer_lemma_2-1-57_sum_3_terms} satisfies $d_{A,\Phi}^{1,*}(v,\nu)=0$.
This completes the proof of the isomorphism \eqref{eq:Itoh_1985_proposition_2-3_SO3_monopole_complex_Kaehler_isomorphism}.

We now verify the isomorphism \eqref{eq:Itoh_1985_proposition_2-3_SO3_monopole_complex_Kaehler_isomorphism_type1} under the additional assumption that $(A,\Phi)$ is type $1$.
Setting $\psi=0$ in \eqref{eq:4-1_and_2_and_3} and combining \eqref{eq:4-1} and \eqref{eq:4-2} yields the system
\begin{subequations}
\label{eq:4-3_type1System}
\begin{align}
  \label{eq:4-3_type1}
  \bar\partial_A^*v'' + (\nu\otimes\varphi^*)_0 - i\bar\partial_A\xi
  &= 0 \in \Omega^{0,1}(\fsl(E)),
  \\
  \label{eq:4-1_and_2_type1}
  \bar\partial_A\nu + \bar\partial_A^*\nu - v''\varphi + i\xi\varphi
  &= 0 \in \Omega^0(E)\oplus \Omega^{0,2}(E).
\end{align}
\end{subequations}
We will now show that
\begin{inparaenum}[\itshape i\upshape)]
\item the terms $i\bar\rd_A\xi$ in \eqref{eq:4-3_type1} and $i\xi\varphi$ in \eqref{eq:4-1_and_2_type1} vanish, so \eqref{eq:4-3_type1System} reduces to the system of equations     \eqref{eq:H_dbar_APhi^01_part_explicit_type1} defining $\bH_{\bar\partial_{A,(\varphi,0)}}^2$, and
\item the section $\xi\in\Omega^0(\su(E))$ belongs to $\bH_{A,\Phi}^0$.
\end{inparaenum}  
Because $\bar\partial_{A,(\varphi,0)}^1\circ\bar\partial_{A,(\varphi,0)}^0 = 0$ by \eqref{eq:d_squared_zero_pre-holomorphic_pair_complex}, we see that
\begin{multline}
  \label{eq:Range_dbarAvarphipsi^0_and_range_dbarAvarphipsi^1star_L2_orthogonal}
  (\bar\partial_{A,(\varphi,0)}^0\zeta, \bar\partial_{A,(\varphi,0)}^{1,*}(v'',\nu))_{L^2(X)}
  \\
  =
  (\bar\partial_{A,(\varphi,0)}^1\circ\bar\partial_{A,(\varphi,0)}^0\zeta, (v'',\nu))_{L^2(X)} = 0, \quad\text{for all } \zeta \in \Omega^0(\fsl(E)),
\end{multline}
while by \eqref{eq:d0StableComplex}, we have
\[
  \bar\partial_{A,(\varphi,0)}^0\zeta
  =
  (\bar\partial_A\zeta, - \zeta\varphi) \in \Omega^{0,1}(\fsl(E)) \oplus \Omega^0(E),
  \quad\text{for all } \zeta \in \Omega^0(\fsl(E)).
\]
By \eqref{eq:H_dbar_APhi^01_part_explicit_type1}, we find that
the component of $\bar\partial_{A,(\varphi,0)}^{1,0}(v'',\nu)$ in $\Om^{0,1}(\fsl(E))\oplus \Om^0(E)$ is given by
\[
  \bar\partial_{A,(\varphi,0)}^{1,*}(v'',\nu) \cap \left(\Omega^{0,1}(\fsl(E)) \oplus \Omega^0(E)\right)
  =
  (\bar\partial_A^*v'' + (\nu\otimes\varphi^*)_0, \bar\partial_A^*\nu).
\]  
By substituting the preceding identities into the $L^2$-orthogonality result \eqref{eq:Range_dbarAvarphipsi^0_and_range_dbarAvarphipsi^1star_L2_orthogonal}, we obtain
\[
  \left( (\bar\partial_A\zeta, - \zeta\varphi),
    (\bar\partial_A^*v'' + (\nu\otimes\varphi^*)_0, \bar\partial_A^*\nu) \right)_{L^2(X)} = 0,
  \quad\text{for all } \zeta \in \Omega^0(\fsl(E)),
\]
that is
\begin{equation}
  \label{eq:Range_dbarAvarphipsi^0_and_range_dbarAvarphipsi^1star_L2_orthogonal_simplified}
  (\bar\partial_A\zeta, \bar\partial_A^*v'' + (\nu\otimes\varphi^*)_0 )_{L^2(X)}
  - (\zeta\varphi,  \bar\partial_A^*\nu)_{L^2(X)}
  = 0, \quad\text{for all } \zeta \in \Omega^0(\fsl(E)).
\end{equation}
Now substituting the identity given by \eqref{eq:4-3_type1},
\[
  \bar\partial_A^*v'' + (\nu\otimes\varphi^*)_0 = \bar\partial_A(i\xi),
\]
into \eqref{eq:Range_dbarAvarphipsi^0_and_range_dbarAvarphipsi^1star_L2_orthogonal_simplified} and choosing $\zeta = i\xi$ yields
\[
  (\bar\partial_A(i\xi), \bar\partial_A(i\xi))_{L^2(X)}
  - (i\xi\varphi,  \bar\partial_A^*\nu)_{L^2(X)}
  = 0,
\]
or equivalently,
\begin{equation}
  \label{eq:Range_dbarAvarphipsi^0_and_range_dbarAvarphipsi^1star_L2_orthogonal_simplified_ixi}
  \|\bar\partial_A\xi\|_{L^2(X)}^2  - (i\xi\varphi,  \bar\partial_A^*\nu)_{L^2(X)}  = 0.
\end{equation}
By taking the imaginary part of the identity \eqref{eq:Range_dbarAvarphipsi^0_and_range_dbarAvarphipsi^1star_L2_orthogonal_simplified_ixi}, we see that
\[
  \Imag(i\xi\varphi,  \bar\partial_A^*\nu)_{L^2(X)} = 0
\]
and hence by \eqref{eq:Equivalence_imaginary_complex_inner_product_equations}, we obtain
\begin{equation}
  \label{eq:xi_varphi_and_d-bar_star_nu_L2-orthogonal}
  (i\xi\varphi,  \bar\partial_A^*\nu)_{L^2(X)} = 0.
\end{equation}
Substituting the $L^2$ orthogonality result \eqref{eq:xi_varphi_and_d-bar_star_nu_L2-orthogonal} into the identity \eqref{eq:Range_dbarAvarphipsi^0_and_range_dbarAvarphipsi^1star_L2_orthogonal_simplified_ixi} gives
\[
  \|\bar\partial_A\xi\|_{L^2(X)}^2 = 0
\]
and thus
\[
  \bar\partial_A\xi = 0 \in \Omega^{0,1}(\fsl(E)).
\]
Thus, \eqref{eq:4-3_type1} reduces to
\begin{equation}
\label{eq:4-3_type1_RemovedTerm}
  \bar\partial_A^*v'' + (\nu\otimes\varphi^*)_0
  = 0 \in \Omega^{0,1}(\fsl(E)).
\end{equation}
In addition, by $ \bar\partial_A\xi = 0$, the identity \eqref{eq:Commute_dagger_and_del_glE}, and the fact that $\xi^\dagger = -\xi$ since $\xi\in\Omega^0(\su(E))$, we also have
\[
  \partial_A\xi = (\bar\partial_A(\xi^\dagger))^\dagger = -(\bar\partial_A\xi)^\dagger
  = 0  \in \Omega^{1,0}(\fsl(E))
\]
and so by \eqref{eq:d_A_sum_components_almost_complex_manifold_integrable},
\begin{equation}
\label{eq:VanishingCovariantOfExtraTerm}
  d_A\xi = \partial_A\xi + \bar\partial_A\xi = 0  \in \Omega^1(\su(E)).
\end{equation}
Because the elements $i\xi\varphi$ and $\bar\partial_A^*\nu$ of $\Omega^0(E)$ are $L^2$-orthogonal by \eqref{eq:xi_varphi_and_d-bar_star_nu_L2-orthogonal} and because the component of \eqref{eq:4-1_and_2_type1} in $\Om^0(E)$ yields
\[
  \bar\partial_A^*\nu + i\xi\varphi = 0 \in \Omega^0(E)
\]
we thus have
\[
  \bar\partial_A^*\nu = 0 \quad\text{and}\quad i\xi\varphi = 0 \in \Omega^0(E).
\]
The equality $i\xi\varphi=0$ implies that \eqref{eq:4-1_and_2_type1} reduces to
\begin{equation}
\label{eq:4-1_and_2_type1_RemovedTerm}
  \bar\partial_A\nu + \bar\partial_A^*\nu - v''\varphi
  = 0 \in \Omega^0(E)\oplus \Omega^{0,2}(E).
\end{equation}
In addition $i\xi\varphi=0$ and \eqref{eq:VanishingCovariantOfExtraTerm} imply that
\[
  d_A\xi - \xi\varphi = 0 \in \Omega^1(\su(E)) \oplus \Omega^0(E).
\]
Combining the preceding equality with \eqref{eq:4-3_type1_RemovedTerm} and \eqref{eq:4-1_and_2_type1_RemovedTerm} yields
\begin{align}
  \label{eq:4-3_type1_corrected_orthogonal_to_H0APhi}
  \bar\partial_A^*v'' + (\nu\otimes\varphi^*)_0 &= 0 \in \Omega^{0,1}(\fsl(E)),
  \\
  \label{eq:4-1_and_2_type1_corrected_orthogonal_to_H0APhi}
  \bar\partial_A\nu + \bar\partial_A^*\nu - v''\varphi  &= 0 \in \Omega^0(E)\oplus \Omega^{0,2}(E),
  \\
  \label{eq:4-1_and_2_type1_corrected_H0APhi}
  d_A\xi - \xi\varphi &= 0.
\end{align}
Hence, when $\psi=0$, we see that $\xi \in \bH_{A,\Phi}^0$ by
\eqref{eq:4-1_and_2_type1_corrected_H0APhi}, the definition of $\bH_{A,\Phi}^0$ in \eqref{eq:H_APhi^0}, and the definition of $d_{A,\Phi}^0$ in \eqref{eq:d_APhi^0}. The identity
\eqref{eq:4-1_and_2_type1_corrected_orthogonal_to_H0APhi} is equal to \eqref{eq:H_dbar_APhi^01_part_explicit_type1_00_and_02} and the identity
\eqref{eq:4-3_type1_corrected_orthogonal_to_H0APhi} is equal to \eqref{eq:H_dbar_APhi^01_part_explicit_type1_01}, so that $(v'',\nu) \in \bH_{\bar\partial_{A,(\varphi,0)}}^2$ by the defining equations \eqref{eq:H_dbar_APhi^01_part_explicit_type1} for $\bH_{\bar\partial_{A,(\varphi,0)}}^2$.

Conversely, suppose
\[
  (v'',\nu) +\xi\otimes\omega \in \bH_{\bar\partial_{A,(\varphi,\psi)}}^2 \oplus \bH_{A,\Phi}^0\otimes\omega,
\]
so that $v = \frac{1}{2}(v'+v''+\xi\otimes\omega) \in \Omega^+(\su(E))$ by \eqref{eq:Donaldson_Kronheimer_lemma_2-1-57_sum_3_terms} when $v' = -(v'')^\dagger$. By reversing the preceding argument, we discover that $(v,\nu) \in \bH_{A,\Phi}^2$. Injectivity of the real linear map \eqref{eq:Real_linear_map_obstruction_spaces_SO(3)_pairs_to_preholomorphic_pairs} is clear, so
the existence of the canonical real linear isomorphisms \eqref{eq:Itoh_1985_proposition_2-3_SO3_monopole_complex_Kaehler_isomorphism}
and \eqref{eq:Itoh_1985_proposition_2-3_SO3_monopole_complex_Kaehler_isomorphism_type1} now follows and this completes the proof of Proposition \ref{prop:Itoh_1985_proposition_2-3_SO3_monopole_complex_Kaehler}.
\end{proof}

\subsection{Auxiliary results required for proof of isomorphism between second-order cohomology groups}
\label{subsec:Auxiliary_results_linear_algebra_second-order_cohomology}
In this subsection, we gather auxiliary results in required by the proof of Proposition \ref{prop:Itoh_1985_proposition_2-3_SO3_monopole_complex_Kaehler}.

\subsubsection{Real transpose operators, complex conjugate operators, and Hermitian adjoint operators on Hilbert spaces}
\label{subsubsec:Complex_conjugate_transpose_adjoint_operator_Hilbert_space}  
For any complex Hilbert spaces $\fH$ and $\fK$ and operator $M\in\Hom(\fH,\fK)$, we recall that the \emph{Hermitian adjoint operator} $M^\dagger\in\Hom(\fK,\fH)$ is defined by the relation (Kadison and Ringrose \cite[Theorem 2.4.2]{KadisonRingrose1})
\begin{equation}
  \label{eq:Hilbert_space_adjoint}
  \langle M^\dagger k, h\rangle_\fH := \langle k, Mh \rangle_\fK, \quad\text{for all } h \in \fH, k \in \fK.
\end{equation}
With respect to almost complex structures $J_\fH$ on $\fH$ and $J_\fK$ on $\fK$ and corresponding complex conjugations $C_\fH(h) = \bar h$ and $C_\fK(k) = \bar k$, we define the \emph{complex conjugate operator} $\bar M\in\Hom(\fK,\fH)$ by the relation
\begin{equation}
  \label{eq:Hilbert_space_complex_conjugate_operator}
  \langle \bar M h, k\rangle_\fK := \langle \bar k, M\bar h \rangle_\fK, \quad\text{for all } h \in \fH, k \in \fK.
\end{equation}
Note that for all $h \in \fH$ and $k \in \fK$,
\begin{align*}
  \langle \bar M\bar h, \bar k \rangle_\fK
  &= \langle k, Mh \rangle_\fK \quad\text{(by \eqref{eq:Hilbert_space_complex_conjugate_operator})}
  \\
  &= \langle \overline{M h}, \bar k \rangle_\fK
  \\
  &= \overline{\langle M h, k \rangle}_\fK.
\end{align*}
We define the \emph{real transpose} $M^\intercal \in \Hom(\fK,\fH)$ by the relation
\begin{equation}
  \label{eq:Hilbert_space_transpose}
\langle M^\intercal k, h\rangle_\fH := \langle M\bar h, \bar k \rangle_\fK, \quad\text{for all } h \in \fH, k \in \fK.
\end{equation}
Note that for all $h \in \fH$ and $k \in \fK$,
\begin{align*}
  \langle \bar M^\intercal k, h\rangle_\fH
  &= \langle \bar M\bar h, \bar k\rangle_\fK \quad\text{(by \eqref{eq:Hilbert_space_transpose})}
  \\
  &= \langle \overline{M h}, \bar k \rangle_\fK \quad\text{(due to \eqref{eq:Hilbert_space_complex_conjugate_operator})}
  \\
  &= \overline{\langle M h, k \rangle}_\fK
  \\
  &= \langle k, Mh\rangle_\fK
  \\
  &= \langle M^\dagger k, h\rangle_\fK \quad\text{(by \eqref{eq:Hilbert_space_adjoint})},
\end{align*}
and thus
\[
  M^\dagger = \bar M^\intercal,
\]
generalizing the usual definition of Hermitian adjoint of a complex finite-dimensional matrix (Lancaster and Tismenetsky \cite[Section 1.5]{Lancaster_Tismenetsky}).

\subsubsection{Alternative definitions of the $L^2$ adjoints of $\partial_A$ and $\bar\partial_A$ on $\Omega^{p,q}(E)$ using a modified Hodge star operator}
\label{subsubsec:Huybrechts_Wells_definitions_L2_adjoints_del_A_and_del-bar_A}
Recall from Huybrechts \cite[Definition 3.1.3, p. 116]{Huybrechts_2005} that if $(X,g,J)$ is a smooth almost Hermitian manifold, then
\begin{equation}
  \label{eq:Huybrechts_definition_3-1-3}
  \partial^* := -\star\circ\bar\partial\circ\star
  \quad\text{and}\quad
  \bar\partial^* := -\star\circ\partial\circ\star.
\end{equation}
According to Huybrechts \cite[Lemma 3.2.3, p. 125]{Huybrechts_2005}, the operators in \eqref{eq:Huybrechts_definition_3-1-3} are $L^2$ adjoints of $\partial$ and $\bar\partial$, respectively, with respect to the Hermitian $L^2$ inner product on $\Omega^\bullet(X,\CC)$ given by
\begin{equation}
  \label{eq:Huybrechts_page_page_33_pointwise_inner_product_pq_forms}
  \langle\alpha,\beta\rangle_{\Lambda^{p,q}(X)}\star 1
  :=
  \alpha\wedge \star\bar\beta \in \Omega^0(X,\CC),
\end{equation}
for all $\alpha,\beta \in \Omega^\bullet(X,\CC)$ and
\begin{equation}
  \label{eq:Huybrechts_definition_3-2-1}
  (\alpha,\beta)_{L^2(X)}
  :=
  \int_X \langle\alpha,\beta\rangle_{\Lambda^{p,q}(X)}\star 1.
\end{equation}
In some references on complex geometry --- and Huybrechts \cite{Huybrechts_2005} and Wells \cite{Wells3} in particular --- when considering differential forms with values in a Hermitian vector bundle $E$ over $(X,g,J)$, it is customary to replace the $\CC$-linear Hodge $\star:\Omega^{p,q}(X) \to \Omega^{n-q,n-p}(X)$ by (see Huybrechts \cite[Definition 4.1.6, p. 168]{Huybrechts_2005}) the $\CC$-anti-linear Hodge star operator
\begin{equation}
  \label{eq:Huybrechts_definition_4-1-6}
  \bar\star_E:\Omega^{p,q}(X,E) \to \Omega^{n-p,n-q}(X,E^*)
\end{equation}
defined on elementary tensors by
\[
  \bar\star_E(\alpha\otimes s) := (\overline{\star\alpha})\otimes \langle\cdot,s\rangle_E = (\star\bar\alpha)\otimes \langle\cdot,s\rangle_E,
\]
for all $\alpha \in \Omega^{p,q}(X)$ and $s \in \Omega^0(E)$. By Huybrechts \cite[Proposition 1.2.20 (iii), p. 32 and Section 4.1, p. 168]{Huybrechts_2005}, one has (using $d=2n$ and $k=p+q$)
\begin{equation}
  \label{eq:Huybrechts_proposition_4-1-61-2-20_iii_and_p_168}
  \bar\star_E\circ\bar\star_E = (-1)^{p+q} \quad\text{on } \Omega^{p,q}(X,E).
\end{equation}
Note that, by using the $\CC$-linear, isometric isomorphism $E\cong E^{**}$ of Hermitian vector bundles, one also has the $\CC$-anti-linear Hodge star operator
\begin{equation}
  \label{eq:Huybrechts_definition_4-1-6_Estar}
  \bar\star_{E^*}:\Omega^{p,q}(X,E^*) \to \Omega^{n-p,n-q}(X,E).
\end{equation}
Following Huybrechts \cite[Definition 4.1.7, p. 169]{Huybrechts_2005}, one defines
\begin{equation}
  \label{eq:Huybrechts_definition_4-1-7}
  \partial_A^* := -\bar\star_{E^*}\circ\,\partial_A\circ\bar\star_E
  \quad\text{and}\quad
  \bar\partial_A^* := -\bar\star_{E^*}\circ\,\bar\partial_A\circ\bar\star_E,
\end{equation}
where $A$ is a unitary connection on $E$.

According to Huybrechts \cite[Lemma 4.1.12, p. 169]{Huybrechts_2005}, the operators in \eqref{eq:Huybrechts_definition_4-1-7} are $L^2$ adjoints of $\partial_A$ and $\bar\partial_A$, respectively, with respect to the Hermitian $L^2$ inner product on $\Omega^\bullet(X,E)$ (see Huybrechts \cite[Section 4.1, p. 168 and Definition 4.1.11, p. 169]{Huybrechts_2005}). Wells provides the same definitions and results for $\bar\star_E$, $\partial_A^*$, and $\bar\partial_A^*$ in \cite[Chapter V, Section 2, pp. 166--168]{Wells3}.

For $\kappa, \mu\in \Omega^{p,q}(E)$, one defines
\begin{equation}
  \label{eq:Huybrechts_page_168_pointwise_inner_product_pq_forms_E}
  \langle\kappa, \mu\rangle_{\Lambda^{p,q}(E)}\star 1
  :=
  \kappa\wedge \bar\star_E\mu \in \Omega^0(X,\CC),
\end{equation}
where ``$\wedge$'' is exterior product in the form part and the evaluation map $(s,t^*) \mapsto t^*(s) = \langle s,t\rangle_E$ in the bundle part, and
\begin{equation}
  \label{eq:Huybrechts_definition_4-1-11}
  (\kappa, \mu)_{L^2(X)}
  :=
  \int_X \langle\kappa, \mu\rangle_{\Lambda^{p,q}(E)}\star 1.
\end{equation}
Note that if $s, t \in \Omega^0(E)$, then $\langle s, t \rangle_E = \tr(s\otimes t^*)$, where $s\otimes t^* \in \Omega^0(\gl(E))$.

For $\kappa\in \Omega^{p,q}(E)$ and $\mu \in \Omega^{p,q+1}(E)$ one has $\bar\star_E\mu \in \Omega^{n-p, n-q-1}(E)$ by \eqref{eq:Huybrechts_definition_4-1-6} and consequently $\bar\partial_A(\bar\star_E\mu) \in \Omega^{n-p, n-q}(E)$, so we have
\begin{align*}
  (\kappa, \bar\partial_A^*\mu)_{L^2(X)}
  &= \int_X \langle\kappa, \bar\partial_A^*\mu\rangle_{\Lambda^{p,q}(E)}\star 1
    \quad\text{(by \eqref{eq:Huybrechts_definition_4-1-11})}
  \\
  &= - \int_X \langle\kappa, \bar\star_E(\bar\partial_A(\bar\star_E\mu)\rangle_{\Lambda^{p,q}(E)}\star 1
    \quad\text{(by \eqref{eq:Huybrechts_definition_4-1-6})})
  \\
  &= -\int_X \kappa \wedge \bar\star_E^2\bar\partial_A(\bar\star_E\mu) \quad\text{(by \eqref{eq:Huybrechts_page_168_pointwise_inner_product_pq_forms_E})}
  \\
  &= (-)^{p+q} \int_X \kappa \wedge \bar\partial_A(\bar\star_E\mu) \quad\text{(by \cite[Proposition 1.2.20 (iii)]{Huybrechts_2005}} 
  \\
  &\qquad\text{with $d=2n$ and $k=2n-p-q$)}
  \\
  &= \int_X (\bar\partial_A\kappa) \wedge \bar\star_E\mu
    \quad\text{(by Leibnitz Rule)}
  \\
  &\qquad\text{extended to $\bar\partial_A$ on $\Omega^{p,q}(E)$)}    
  \\
  &= \int_X \langle\bar\partial_E\kappa, \mu\rangle_{\Lambda^{p,q+1}(E)}\star 1 \quad\text{(by \eqref{eq:Huybrechts_page_168_pointwise_inner_product_pq_forms_E})}
  \\
  &= (\bar\partial_E\kappa, \mu)_{L^2(X)} \quad\text{(by \eqref{eq:Huybrechts_definition_4-1-11})}, 
\end{align*}
where we used the Leibnitz Rule \eqref{eq:Donaldson_Kronheimer_2-1-45_i_pq-forms} in the form
\[
  \bar\partial(\kappa\wedge\bar\star_E\mu)
  =
  (\bar\partial_A\kappa)\wedge \bar\star_E\mu + (-1)^{p+q}\kappa\wedge\partial_A(\bar\star_E\mu),
\]
and Stokes' Theorem in the form
\[
\int_X \bar\partial(\kappa\wedge\bar\star_E\mu) = \int_X d(\kappa\wedge\bar\star_E\mu) = 0.
\]
We do \emph{not} use the definitions of the $L^2$ adjoints of $\partial_A$ and $\bar\partial_A$ on $\Omega^{p,q}(E)$; rather, we employ the definitions provided in Section \ref{subsubsec:Donaldson_Kronheimer_implied_L2_adjoint_dbar_A_and_delA}.

\subsubsection{Definitions of the $L^2$ adjoints of $\partial_A$ and $\bar\partial_A$ on $\Omega^{p,q}(E)$ using the standard Hodge star operator}
\label{subsubsec:Donaldson_Kronheimer_implied_L2_adjoint_dbar_A_and_delA}
Rather than modify the definition of Hodge star $\star:\Omega^{\bullet,\bullet}(X)\to \Omega^{\bullet,\bullet}(X)$ when replacing the product line bundle $X\times\CC$ by an arbitrary Hermitian vector bundle $E$ over $X$ for the purpose of defining the $L^2$ adjoints of $\partial_A$ and $\bar\partial_A$ on $\Omega^{p,q}(E)$, we shall instead use the definitions of $\partial_A^*$ and $\bar\partial_A^*$ implied by Donaldson and Kronheimer \cite[Equations (2.1.12), (2.1.24), p. 37, Equation (2.1.47), and p. 44]{DK} in their definitions of $d_A$, $d_A^*$, $\partial_A$, and $\bar\partial_A$. Thus, we define (for any merely complex vector bundle $E$)
\begin{equation}
  \label{eq:Complex_linear_Hodge_star_pq_forms_complex_vector_bundle}
  \star:\Omega^{p,q}(E) \to \Omega^{n-q,n-p}(E)
\end{equation}
on elementary tensors by
\[
  \star(\alpha\otimes s) := \star\alpha\otimes s \in  \Omega^{n-q,n-p}(E).
\]
For tensors $\kappa, \mu \in \Omega^{p,q}(E)$ and elementary tensors $\alpha\otimes s, \beta\otimes t\in \Omega^{p,q}(E)$ defined by $\alpha,\beta\in\Omega^{p,q}(X)$ and $s,t \in \Omega^0(E)$, we define the pointwise Hermitian inner product on $\Lambda^{p,q}(E) = \Lambda^{p,q}(X)\otimes E$ by
\begin{equation}
  \label{eq:Pointwise_inner_product_pq_forms_E}
  \begin{aligned}
  \langle\kappa, \mu\rangle_{\Lambda^{p,q}(E)}
  &:=
  \langle\kappa, \mu\rangle_{\Lambda^{p,q}(X)\otimes E},
  \\
  \langle \alpha\otimes s, \beta\otimes t \rangle_{\Lambda^{p,q}(X)\otimes E}
  &:=
  \langle\alpha, \beta\rangle_{\Lambda^{p,q}(X)} \langle s, t\rangle_E,
  \\
  \langle\alpha, \beta\rangle_{\Lambda^{p,q}(X)}
  &:= \star(\alpha\wedge \star\bar\beta) \in \Omega^0(X,\CC),
  \end{aligned}
\end{equation}
and we define the $L^2$ Hermitian inner product on $\Omega^{p,q}(E)$ by
\begin{equation}
  \label{eq:L2_inner_product_pq_forms_E}
  (\kappa, \mu)_{L^2(X)}
  :=
  \int_X \star\langle\kappa, \mu\rangle_{\Lambda^{p,q}(E)}\,d\vol,
  \quad\text{for all } \kappa, \mu \in \Omega^{p,q}(E).
\end{equation}
For $\kappa\in \Omega^{p,q}(E)$ and $\mu \in \Omega^{p,q+1}(E)$ one has $\star\mu \in \Omega^{n-q-1,n-p}(E)$ by \eqref{eq:Huybrechts_page_33_and_lemma_1-2-24} and $d_A(\star\mu) \in \Omega^{2n-q-p}(E)$ and $\partial_A(\star\mu) = \pi_{n-q,n-p}d_A(\star\mu) \in \Omega^{n-q,n-p}(E)$ by the definition
\begin{equation}
  \label{eq:del_A}
    \partial_A = \pi_{p+1,q}d_A:\Omega^{p,q}(E) \to \Omega^{p+1,q}(E),
\end{equation}
so that, using (see Warner \cite[Equation (6.1.2)]{Warner})
\begin{equation}
  \label{eq:Warner_6-1-2}
  d_A^* := (-1)^{d(k+1)+1}\star d_A \star \quad\text{on } \Omega^k(X,\CC),
\end{equation}
with $d=2n$ and $k=p+q+1$, we obtain
\begin{multline*}
  \pi_{p,q}d_A^*\mu = -\pi_{p,q}\star d_A(\star\mu) = -\pi_{p,q}\star d_A\pi_{n-q-1,n-p}(\star\mu)
  \\
  = -\pi_{p,q}\star \pi_{n-q,n-p}d_A\pi_{n-q-1,n-p}(\star\mu) = -\pi_{p,q}\star \partial_A\pi_{n-q-1,n-p}(\star\mu)
  = -\star \partial_A(\star\mu)  \in \Omega^{p,q}(E).
\end{multline*}
Recall that by analogy with Warner \cite[Proposition 6.2]{Warner} for $d^*$ on $\Omega^k(X,\RR)$ that
\begin{equation}
  \label{eq:Warner_6-2-1}
  (a,d_A^*b)_{L^2(X)} = (d_Aa,b)_{L^2(X)} \quad\text{for all } a \in \Omega^k(E), b \in \Omega^{k+1}(E).
\end{equation}
Therefore, if we define (see Huybrechts \cite[Definition 3.1.3, p. 116 and Lemma 3.2.3, p. 125]{Huybrechts_2005} for the case $E=X\times\CC$)
\begin{equation}
  \label{eq:dbar_A_star}
  \bar\partial_A^*\mu := \pi_{p,q}d_A^*\mu = -\star \partial_A(\star\mu) \in \Omega^{p,q}(E),
  \quad\text{for all } \mu \in \Omega^{p,q+1}(E).
\end{equation}
and apply \eqref{eq:Warner_6-2-1} we obtain, for all $\kappa \in \Omega^{p,q}(E)$ and $\mu \in \Omega^{p,q+1}(E)$,
\begin{multline*}
  (\kappa, \bar\partial_A^*\mu)_{L^2(X)} = (\kappa, \pi_{p,q}d_A^*\mu)_{L^2(X)} = (\kappa, d_A^*\mu)_{L^2(X)}
  \\
  = (d_A\kappa, \mu)_{L^2(X)} = (\pi_{p,q+1}d_A\kappa, \mu)_{L^2(X)} = (\bar\partial_A\kappa, \mu)_{L^2(X)},
\end{multline*}
and thus $\bar\partial_A^*$ in \eqref{eq:dbar_A_star} is the $L^2$ adjoint of $\bar\partial_A$.

Similarly, for $\kappa\in \Omega^{p,q}(E)$ and $\mu \in \Omega^{p+1,q}(E)$ one has $\star\mu \in \Omega^{n-q,n-p-1}(E)$ by \eqref{eq:Huybrechts_page_33_and_lemma_1-2-24} and $d_A(\star\mu) \in \Omega^{2n-q-p}(E)$ and $\bar\partial_A(\star\mu) = \pi_{n-q,n-p}d_A(\star\mu) \in \Omega^{n-q,n-p}(E)$ by definition of $\bar\partial_A$, and thus we obtain
\begin{multline*}
  \pi_{p,q}d_A^*\mu = -\pi_{p,q}\star d_A(\star\mu) = -\pi_{p,q}\star d_A\pi_{n-q,n-p-1}(\star\mu)
  \\
  = -\pi_{p,q}\star \pi_{n-q,n-p}d_A\pi_{n-q,n-p-1}(\star\mu) = -\pi_{p,q}\star \bar\partial_A\pi_{n-q,n-p-1}(\star\mu)
  = -\star \bar\partial_A(\star\mu)  \in \Omega^{p,q}(E).
\end{multline*}
Therefore, if we define
\begin{equation}
  \label{eq:del_A_star}
  \partial_A^*\mu := \pi_{p,q}d_A^*\mu = -\star \bar\partial_A(\star\mu) \in \Omega^{p,q}(E),
  \quad\text{for all } \mu \in \Omega^{p+1,q}(E),
\end{equation}
we now find that $\partial_A^*$ in \eqref{eq:del_A_star} is the $L^2$ adjoint of $\partial_A$. This completes our discussion of our expressions for the $L^2$ adjoints of $\partial_A$ and $\bar\partial_A$ on $\Omega^{p,q}(E)$ when computed using the \emph{unmodified} Hodge star $\star:\Lambda^\bullet(E)\to \Lambda^\bullet(E)$.

\subsubsection{Auxiliary properties of real and complex Hilbert spaces}
\label{subsubsec:Auxiliary_properties_real_and_complex_Hilbert_spaces}
Let $\fH$ be a complex Hilbert space with real form $\sH$, so $\fH = \sH\otimes_\RR\CC$. We write $h=h_1+ih_2\in\fH$, where $h_1,h_2\in\sH$. Then
\[
  \langle h,h'\rangle_\fH = \left(\langle h_1,h_1'\rangle_\sH + \langle h_2,h_2'\rangle_\sH \right)
  -  i\left(\langle h_1,h_2'\rangle_\sH + \langle h_1',h_2\rangle_\sH \right),
\]
so that
\[
  \Real\,\langle h,h'\rangle_\fH = \langle h_1,h_1'\rangle_\sH + \langle h_2,h_2'\rangle_\sH.
\]
(This matches the convention in \eqref{eq:Kobayashi_7-6-5_abstract_inner_product_space}.) Given $h \in \fH$, suppose that $\Real\,\langle h,h'\rangle_\fH = 0$ for all $h' \in \fH$. By first choosing $h' = h_1' + i0$ and then $h' = 0 + ih_2'$, respectively, we see that
\[
  \langle h_1,h_1'\rangle_\sH = 0 = \langle h_2,h_2'\rangle_\sH, \quad\text{for all } h_1', h_2' \in \sH,
\]
and so we obtain $h_1=h_2=0$. Hence, for any given $h\in\fH$, the following equations are equivalent:
\begin{equation}
  \label{eq:Equivalence_real_complex_inner_product_equations}
  \Real\,\langle h,h'\rangle_\fH = 0 \quad\text{for all } h' \in \fH
  \iff
  \langle h,h'\rangle_\fH = 0 \quad\text{for all } h' \in \fH.
\end{equation}
Alternatively, just observe that
\[
  \Real\,\langle h,ih'\rangle_\fH = \Real\left(-i\langle h,h'\rangle_\fH\right) = \Imag\,\langle h,h'\rangle_\fH
\]
and the conclusion \eqref{eq:Equivalence_real_complex_inner_product_equations} again follows. We shall also use \eqref{eq:Equivalence_real_complex_inner_product_equations} in the equivalent form
\begin{equation}
  \label{eq:Equivalence_imaginary_complex_inner_product_equations}
  \Imag\,\langle h,h'\rangle_\fH = 0 \quad\text{for all } h' \in \fH
  \iff
  \langle h,h'\rangle_\fH = 0 \quad\text{for all } h' \in \fH,
\end{equation}
obtained from the preceding relation between the real and imaginary parts of the complex inner product $\langle h,h'\rangle_\fH$.

For any complex Hilbert spaces $\fH$ and $\fK$ and operator $M\in\Hom(\fH,\fK)$, its \emph{Toeplitz Decomposition} (see Horn and Johnson \cite[p. 7]{Horn_Johnson_matrix_analysis_2013}) is given by
\[
  M = \pi_HM + \pi_SM = M^H + M^S,
\]
where $\pi_HM := \frac{1}{2}(M + M^\dagger)$ and $\pi_SM := \frac{1}{2}(M - M^\dagger) = -i\pi_H(iM)$, since $\pi_H(iM) =  \frac{1}{2}(iM - iM^\dagger) = i\pi_SM$. 

\section[Elliptic complexes for anti-self-dual and holomorphic curvature equations]{Comparison of elliptic complexes for the anti-self-dual and holomorphic curvature equations over complex K\"ahler surfaces}
\label{sec:Elliptic_deformation_complex_moduli_space_ASD_connections_over_complex_Kaehler_surface}
Our goal in this section is to give an exposition of proofs of the isomorphisms \cite[Section 6.4.2, p. 240]{DK}, \cite{Itoh_1985} of cohomology groups for the elliptic deformation complex for the anti-self-dual equation
\begin{equation}
  \label{eq:Elliptic_deformation_complex_ASD_equation}
  \Omega^0(\su(E)) \xrightarrow{d_A} \Omega^1(\su(E)) \xrightarrow{d_A^+} \Omega^+(\su(E)),
\end{equation}
with the corresponding cohomology groups for the elliptic deformation complex for the holomorphic curvature equation
\begin{equation}
  \label{eq:Elliptic_deformation_complex_holomorphic_curvature_equation}
  \Omega^{0,0}(\fsl(E)) \xrightarrow{\bar\partial_A} \Omega^{0,1}(\fsl(E)) \xrightarrow{\bar\partial_A} \Omega^{0,2}(\fsl(E)),
\end{equation}
given a Hermitian vector bundle $E$ over a complex K\"ahler surface $(X,g,J)$. The harmonic representatives of the cohomology groups for the complex \eqref{eq:Elliptic_deformation_complex_ASD_equation} are defined by the real vector spaces,
\begin{subequations}
\label{eq:H_Abullet}
\begin{align}
  \label{eq:H_A^0}
  \bH_A^0 &:= \Ker\left(d_A:\Omega^0(\su(E)) \to \Omega^1(\su(E))\right),
  \\
  \label{eq:H_A^1}
  \bH_A^1 &:= \Ker\left(d_A^+ + d_A^*:\Omega^1(\su(E)) \to \Omega^+(\su(E))\oplus \Omega^0(\su(E))\right),
  \\
  \label{eq:H_A^2}
  \bH_A^2 &:=\Ker \left(d_A^{+,*}:\Omega^+(\su(E)) \to \Omega^1(\su(E))\right),
\end{align}
\end{subequations}
and similarly, the complex vector spaces
\begin{subequations}
\label{eq:H_dbar_A^0bullet}
\begin{align}
\label{eq:H_dbar_A^00}
  \bH_{\bar\partial_A}^0 &:= \Ker\left(\bar\partial_A:\Omega^{0,0}(\fsl(E)) \to \Omega^{0,1}(\fsl(E)) \right),
  \\
\label{eq:H_dbar_A^01}
  \bH_{\bar\partial_A}^1 &:= \Ker\left(\bar\partial_A + \bar\partial_A^*:\Omega^{0,1}(\fsl(E)) \to \Omega^{0,2}(\fsl(E))\oplus\Omega^{0,0}(\fsl(E)) \right),
  \\
\label{eq:H_dbar_A^02}
  \bH_{\bar\partial_A}^2 &:= \Ker\left(\bar\partial_A^*:\Omega^{0,2}(\fsl(E)) \to \Omega^{0,1}(\fsl(E)) \right),
\end{align}
\end{subequations}
define the harmonic representatives of the cohomology groups for the complex \eqref{eq:Elliptic_deformation_complex_holomorphic_curvature_equation}. In this subsection, we deduce the canonical isomorphisms between cohomology groups for the elliptic complexes \eqref{eq:Elliptic_deformation_complex_ASD_equation} and \eqref{eq:Elliptic_deformation_complex_holomorphic_curvature_equation} as corollaries of the corresponding results for pairs proved earlier.

For complex K\"ahler surfaces, the identification in the forthcoming Lemma \ref{lem:Itoh_1985_proposition_2-4} was stated without proof by Donaldson and Kronheimer in \cite[Section 6.4.2, p. 240]{DK} but previously established using methods similar to ours by Itoh in the proof of his \cite[Proposition 2.4]{Itoh_1985} and by Kobayashi in the proof of his \cite[Theorem 7.2.21, p. 227]{Kobayashi_differential_geometry_complex_vector_bundles}. Proposition \ref{prop:Itoh_1985_proposition_2-4_SO3_monopole_complex_Kaehler} yields the

\begin{lem}[Canonical isomorphism between first-order cohomology groups for elliptic complexes over closed, complex K\"ahler surfaces]
\label{lem:Itoh_1985_proposition_2-4}
Let $(E,h)$ be a smooth Hermitian vector bundle over a closed complex K\"ahler surface $(X,g,J)$. If $A$ is a smooth, projectively anti-self-dual unitary connection on $E$ that induces a fixed unitary connection $A_d$ on $\det E$ with $F_{A_d}^{0,2} = 0$, then the isomorphism of real affine spaces,
\begin{equation}
  \label{eq:Isomorphism_between_smooth_unitary_connections_and_01-connections}
  A_0 + \Omega^1(\su(E)) \ni A \mapsto \bar\partial_A \in \Omega^{0,1}(\fsl(E)) + \bar\partial_{A_0},
\end{equation}
induces a canonical isomorphism of finite-dimensional real vector spaces,
\begin{equation}
\label{eq:Itoh_1985_proposition_2-4_isomorphism}
  \bH_A^1 \cong \bH_{\bar\partial_A}^1,
\end{equation}
where $\bH_A^1$ is as in \eqref{eq:H_A^1} and $\bH_{\bar\partial_A}^1$ is as in \eqref{eq:H_dbar_A^01}.
\end{lem}

Recall that $(p,q)$ forms of different bi-type are $L^2$ orthogonal by \eqref{eq:Almost_hermitian_manifold_pq_forms_pointwise_orthogonal_rs_forms_unless_pq_equals_rs}. As a consequence of the \emph{Hodge Decomposition Theorem} (see Huybrechts \cite[Theorem 3.2.8 and Corollary 3.2.9, p. 128]{Huybrechts_2005}) for a closed, K\"ahler manifold $X$ that there is a direct sum decomposition (see Huybrechts \cite[Proposition 3.2.6 (ii), p. 126 and Corollary 3.2.12, p. 129]{Huybrechts_2005} or Wells \cite[Chapter V, Theorem 5.1, p. 197]{Wells3}) 
\begin{equation}
  \label{eq:Hodge_decomposition}
  \bH^r(X,\CC) = \bigoplus_{p+q=r}\bH^{p,q}(X)
\end{equation}
that is independent of the K\"ahler metric on $X$; moreover, $\bH^r(X,\CC) = \bH^r(X,\RR)\otimes_\RR\CC$ and $\overline{\bH^{p,q}(X)} = \bH^{q,p}(X)$. Recall that
\[
  \bH^{p,q}(X) = \{\alpha \in \Omega^{p,q}(X): \bar\partial\alpha = 0 \text{ and } \bar\partial^*\alpha = 0\}
\]
by \cite[Definitions 3.1.5 and 3.2.4]{Huybrechts_2005}.

In particular, $\bH^{0,0}(X) = \bH^0(X,\RR)\otimes_\RR\CC$. Indeed, to see this we may follow Huybrechts \cite[Proof of Proposition 3.2.6 (ii), p. 126]{Huybrechts_2005} or Wells \cite[Chapter V, Proof of Theorem 5.1, pp. 197--198]{Wells3} and observe that $f \in \bH^{0,0}(X)$ if and only if $\bar\square f = 0$ while $\bar\square = \frac{1}{2}\Delta$ by Huybrechts \cite[Proposition 3.1.12, p. 120]{Huybrechts_2005} or Wells \cite[Chapter V, Theorem 4.7, p. 191]{Wells3}, where $\Delta = dd^* + d^*d$ and $\bar\square =  \bar\partial^* \bar\partial +  \bar\partial\bar\partial^*$ by Wells \cite[Chapter V, Section 4, p. 191]{Wells3}. Thus $\Delta f = d^*df = 0$ and so $df = 0$.

Noting that the condition that $X$ has real dimension four in Proposition \ref{prop:Donaldson_1985_proposition_3_and_Itoh_1985_proposition_4-1_SO3_monopole_almost_Hermitian} and Remark \ref{rmk:Isomorphism_nonAbelianH0_preholomH0} can be relaxed (we only used the fact that $X$ has real dimension four in order to describe the canonical \spinc structure over an almost Hermitian manifold), that proposition immediately yields the

\begin{lem}[Canonical isomorphism between zeroth-order cohomology groups for elliptic complexes over closed, almost Hermitian manifolds]
\label{lem:Donaldson_1985_proposition_3_and_Itoh_1985_proposition_4-1_almost_Hermitian}
Let $E$ be a smooth Hermitian vector bundle over a closed, smooth almost Hermitian manifold $(X,g,J)$. If $A$ is a smooth, projectively anti-self-dual unitary connection on $E$ that induces a fixed unitary connection $A_d$ on $\det E$,
then the isomorphism
\eqref{eq:Isomorphism_between_smooth_unitary_connections_and_01-connections} of real affine spaces induces a canonical isomorphism of complex vector spaces,
\begin{equation}
\label{eq:Itoh_1985_proposition_4-1}  
\bH_A^0\otimes_\RR\CC \cong \bH_{\bar\partial_A}^0, 
\end{equation}
where $\bH_A^0$ is as in \eqref{eq:H_A^0} and $\bH_{\bar\partial_A}^0$ is as in \eqref{eq:H_dbar_A^00}.
\end{lem}

For complex K\"ahler surfaces, the identification \eqref{eq:Itoh_1985_proposition_4-1} is stated by Donaldson and Kronheimer in \cite[Section 6.4.2, p. 240]{DK}, implied by Itoh in \cite[Proposition 4.1]{Itoh_1985}, by Kobayashi in \cite[Theorem 7.2.21, p. 227]{Kobayashi_differential_geometry_complex_vector_bundles}, and is implicit in the proof of \cite[Section 4.3.3, Proposition 3.7, p. 324]{FrM} by Friedman and Morgan.

Recall that with respect to a local oriented orthonormal frame $\{e_i\}$ for $TX$ and $\{e_i^*\}$ for $T^*X$ over an open subset $U \subset X$, we may write $\omega = e_1^*\wedge e_2^* + e_3^*\wedge e_4^*$ with volume form $\frac{1}{2}\omega^2 = e_1^*\wedge e_2^* \wedge e_3^*\wedge e_4^* = d\vol$. The pointwise norm $\omega$ is defined through $|\omega|^2\,d\vol = \omega\wedge *\omega$ by Warner \cite[Exercise 2.13 and Equation (6.1.5)]{Warner}. But $*\omega = \omega$ since $\omega$ is self-dual by direct calculation \cite[Proof of Lemma 2.1.57, p. 47]{DK}. Hence, $|\omega|^2\,d\vol = \frac{1}{2}|\omega|^2\omega^2 = \omega^2$ and thus $|\omega| = \sqrt{2}$ on $U$ and so also on $X$. The pointwise orthogonal projection onto the Hermitian line bundle $\Lambda^0(X,\CC)\omega$ is given by $\frac{1}{2}\langle\cdot,\omega\rangle\omega = \frac{1}{2}\Lambda(\cdot)$ since $\omega/\sqrt{2}$ comprises an orthonormal frame and $\Lambda = \langle\cdot,\omega\rangle = L^*$, where $L = (\cdot)\otimes\omega$ (see Section \ref{subsec:Hodge_star_Lefshetz_dual_Lefshetz_operators}). By the reality condition \eqref{eq:Itoh_1985_p_850_reality_condition_su(E)_2-forms}, we have
\[
v' = -(v'')^\dagger.
\]
Proposition \ref{prop:Itoh_1985_proposition_2-3_SO3_monopole_complex_Kaehler} yields the 
  
\begin{lem}[Canonical real linear isomorphism between second-order cohomology groups for elliptic complexes over closed, complex K\"ahler surfaces]
\label{lem:Itoh_1985_proposition_2-3_complex_Kaehler}
(See Itoh \cite[Proposition 2.2]{Itoh_1985} and Donaldson \cite[Proof of Proposition 3]{DonASD}.)
Assume the hypotheses of Lemma \ref{lem:Itoh_1985_proposition_2-4}.
Then the isomorphism \eqref{eq:Isomorphism_between_smooth_unitary_connections_and_01-connections} of real affine spaces induces a canonical isomorphism\footnote{This isomorphism is stated without proof by Donaldson and Kronheimer in \cite[Section 6.4.2, p. 240]{DK} and established by more indirect methods by Itoh in the proof of his \cite[Equation (2.16) and Proposition 2.3]{Itoh_1985} and Kobayashi in the proof of his \cite[Theorem 7.2.21, p. 227]{Kobayashi_differential_geometry_complex_vector_bundles}.}
of finite-dimensional real vector spaces,
\begin{equation}
  \label{eq:Itoh_1985_proposition_2-3_complex_Kaehler}
  \bH_A^2 \cong \bH_{\bar\partial_A}^2\oplus \bH_A^0,
\end{equation}
where $\bH_A^0$ and $\bH_A^2$ are as in \eqref{eq:H_A^0} and \eqref{eq:H_A^2}, respectively, and $\bH_{\bar\partial_A}^2$ is as in \eqref{eq:H_dbar_A^02}.
\end{lem}

\section[Elliptic complexes for Hermitian--Einstein and holomorphic curvature equations]{Comparison of elliptic complexes for the projectively Hermitian--Einstein and holomorphic curvature equations over complex K\"ahler manifolds}
\label{sec:Elliptic_complex_moduli_space_HE_connections_over_complex_Kaehler_manifold}
The results of Section \ref{sec:Elliptic_deformation_complex_moduli_space_ASD_connections_over_complex_Kaehler_surface} generalize to give a comparison between the elliptic complexes for the projectively Hermitian--Einstein and holomorphic curvature equations over complex K\"ahler manifolds and hence isomorphisms between their cohomology groups and harmonic spaces when those manifolds are closed. In particular, we have the following version of Kobayashi \cite[Theorem 7.2.21, p. 227]{Kobayashi_differential_geometry_complex_vector_bundles}, slightly modified to include the condition of a fixed holomorphic structure on the complex determinant line bundle and tailored to our statement of the projectively Hermitian--Einstein equations.

\begin{thm}[Canonical real linear isomorphisms between cohomology groups for elliptic complexes for the projectively Hermitian--Einstein and holomorphic curvature equations over closed, complex K\"ahler manifolds]
\label{thm:Kobayashi_7-2-21}
(Compare Kobayashi \cite[Theorem 7.2.21, p. 227]{Kobayashi_differential_geometry_complex_vector_bundles}.)
Let $(E,h)$ be a smooth Hermitian vector bundle over a closed, complex K\"ahler manifold $(X,g,J)$ with fundamental two-form $\omega = g(\cdot,J\cdot)$ as in \eqref{eq:Fundamental_two-form}. If $A$ is a solution to the forthcoming projectively Hermitian--Einstein equations \eqref{eq:Holomorphic_connection}, \eqref{eq:Einstein_connection} that induces a fixed holomorphic structure on the complex line bundle $\det E$, then the canonical isomorphism
\begin{equation}
  \label{eq:Bijection_unitaryconnections_with_01connections}
  \pi_h^{0,1}:\sA(E,h) \ni A \mapsto \bar\partial_A \in \sA^{0,1}(E)
\end{equation}
of real affine spaces (see Sections \ref{sec:SpinuPairsQuotientSpace} and \ref{sec:Elliptic_deformation_complex_holomorphic_pair_equations} and equation \eqref{eq:Affine_space_01-connections}) induces canonical isomorphisms of finite-dimensional real vector spaces,
\begin{subequations}
\label{eq:Kobayashi_7-2-21_HAbullet_isomorphism_HdbarA0bullet}  
\begin{align}
  \label{eq:Kobayashi_7-2-21_HA0complex_isomorphism_HdbarA00}
  \bH_A^0\otimes_\RR\CC &\cong \bH_{\bar\partial_A}^0,
  \\
  \label{eq:Kobayashi_7-2-21_HA1_isomorphism_HdbarA01}
  \bH_A^1 &\cong \bH_{\bar\partial_A}^1,
  \\
  \label{eq:Kobayashi_7-2-21_HA2_isomorphism_HdbarA02+HdbarA00}
  \bH_A^2 &\cong \bH_{\bar\partial_A}^2\oplus \bH_A^0,
  \\
  \label{eq:Kobayashi_7-2-21_HAk_isomorphism_HdbarA0k}
  \bH_A^k &\cong \bH_{\bar\partial_A}^k, \quad\text{for } k = 3, \ldots,n,               
\end{align}
\end{subequations}
where harmonic spaces $\bH_A^\bullet$ are as in the forthcoming \eqref{eq:HE_equation_bHAbullet} and the harmonic spaces $\bH_{\bar\partial_A}^\bullet$ are as in the forthcoming \eqref{eq:bHdbarA0bullet}.
\end{thm}

\begin{proof}
The isomorphism \eqref{eq:Kobayashi_7-2-21_HA0complex_isomorphism_HdbarA00} is given by Lemma \ref{lem:Donaldson_1985_proposition_3_and_Itoh_1985_proposition_4-1_almost_Hermitian}. The isomorphism \eqref{eq:Kobayashi_7-2-21_HA1_isomorphism_HdbarA01} is implied by Lemma \ref{lem:Itoh_1985_proposition_2-4}, noting that the role of the differential $d_A^+$ in the forthcoming \eqref{eq:Linearization_HE_equation_full}, namely
\begin{align*}
  d_A^+:\Omega^1(\su(E)) \ni a \mapsto
  &\frac{1}{2}\left(\partial_A a' + \bar\partial_Aa''\right) + \Lambda(d_Aa)\otimes\omega
  \\
  \in\ &\Omega^{2,0}(\fsl(E)) \oplus \Omega^{0,2}(\fsl(E)) \oplus \Omega^0(\su(E))\otimes\omega
\end{align*}
can be replaced by that of $\hat d_A^+$ in the forthcoming \eqref{eq:Linearization_HE_equation_codomain_Omega0+Omega02}, namely
\[
  \hat d_A^+:\Omega^1(\su(E)) \ni a \mapsto \frac{1}{2}\bar\partial_Aa'' + \Lambda(d_Aa) \in
  \Omega^{0,2}(\fsl(E)) \oplus \Omega^0(\su(E))
\]
and the fact that the proof of Lemma \ref{lem:Itoh_1985_proposition_2-4} does not require $X$ to have complex dimension two. (Our proof of the more general Proposition \ref{prop:Itoh_1985_proposition_2-4_SO3_monopole_complex_Kaehler} only uses a hypothesis that $X$ has complex dimension two in order to describe the canonical \spinc structure over a complex K\"ahler surface: that hypothesis can be replaced by the assumption in higher dimensions that $\Phi = (\varphi,\psi)$, for $\varphi\in\Omega^0(E)$ and $\psi\in\Omega^{0,2}(E)$ in general and $\Phi=(0,0)$ here.) Similarly, Lemma \ref{lem:Itoh_1985_proposition_2-3_complex_Kaehler} yields the isomorphism \eqref{eq:Kobayashi_7-2-21_HA2_isomorphism_HdbarA02+HdbarA00} because our proof of the more general Proposition \ref{prop:Itoh_1985_proposition_2-3_SO3_monopole_complex_Kaehler} only uses a hypothesis that $X$ has complex dimension two in order to describe the canonical \spinc structure over a complex K\"ahler surface and because, when the complex dimension of $X$ becomes three or more, the definitions of $\bH_A^2$ and $\bH_{\bar\partial_A}^2$ both change by imposing the same condition $\bar\partial_A\beta = 0$ on $\beta \in \Omega^{0,2}(\fsl(E))$. The isomorphism \eqref{eq:Kobayashi_7-2-21_HAk_isomorphism_HdbarA0k} is immediate from the definitions of these harmonic spaces since the higher-degree vector spaces and differentials in the associated elliptic complexes are identical. 
\end{proof}  

\section[Elliptic complexes for projective vortex and holomorphic pair equations]{Comparison of elliptic complexes for the projective vortex and holomorphic pair equations over complex K\"ahler manifolds}
\label{sec:Elliptic_complex_moduli_space_projective vortices_over_complex_Kaehler_manifold}
Lastly, the results of Section \ref{sec:Elliptic_complex_moduli_space_HE_connections_over_complex_Kaehler_manifold} generalize in turn to give a comparison between the elliptic complexes for the projective vortex and holomorphic pair equations over complex K\"ahler manifolds and hence isomorphisms between their cohomology groups and harmonic spaces when those manifolds are closed. In particular, we have the following generalization of Theorem \ref{thm:Kobayashi_7-2-21}.

\begin{thm}[Canonical real linear isomorphisms between cohomology groups for elliptic complexes for the projective vortex and holomorphic pair equations over closed, complex K\"ahler manifolds]
\label{thm:Kobayashi_7-2-21_pairs}
Continue the hypotheses of Theorem \ref{thm:Kobayashi_7-2-21}.  If $(A,\varphi)$ is a solution to the
projective vortex equations \eqref{eq:SO(3)_monopole_equations_almost_Hermitian_alpha} that induces a fixed holomorphic structure on the complex line bundle $\det E$ and $(\bar\partial_A,\varphi)$ is the corresponding solution to the holomorphic pair equations \eqref{eq:Holomorphic_pair}, then the isomorphism $\pi_h^{0,1}:\sA(E,h) \cong \sA^{0,1}(E)$ of real affine spaces in \eqref{eq:Bijection_unitaryconnections_with_01connections} and identity operator on $\Omega^0(E)$ induces canonical isomorphisms of finite-dimensional real vector spaces,
\begin{subequations}
\label{eq:HAvarphi^bullet_isomorphic_HdbarAvarphi^0bullet}  
\begin{align}
  \label{eq:HAvarphi^0complexified_isomorphic_HdbarAvarphi^00}
  \bH_{A,\varphi}^0\otimes_\RR\CC &\cong \bH_{\bar\partial_A,\varphi}^0,
  \\
  \label{eq:HAvarphi^1_isomorphic_HdbarAvarphi^01}
  \bH_{A,\varphi}^1 &\cong \bH_{\bar\partial_A,\varphi}^1,
  \\
  \label{eq:HAvarphi^2_isomorphic_HdbarAvarphi^02+HdbarAvarphi^00}
  \bH_{A,\varphi}^2 &\cong \bH_{\bar\partial_A,\varphi}^2\oplus \bH_{A,\varphi}^0,
  \\
  \label{eq:HAvarphi^k_isomorphic_HdbarAvarphi^0k}
  \bH_{A,\varphi}^k &\cong \bH_{\bar\partial_A,\varphi}^k, \quad\text{for } k = 3, \ldots,n,
\end{align}
\end{subequations}
where harmonic spaces $\bH_{A,\varphi}^\bullet$ are as in \eqref{eq:H_Avarphi^bullet} and the harmonic spaces
$\bH_{\bar\partial_A,\varphi}^\bullet$ are as in \eqref{eq:H_dbar_Avarphi^0bullet}.
\end{thm}

\begin{proof}
Consider the assertion \eqref{eq:HAvarphi^0complexified_isomorphic_HdbarAvarphi^00}. If $X$ has complex dimension two, then Proposition \ref{prop:Donaldson_1985_proposition_3_and_Itoh_1985_proposition_4-1_SO3_monopole_almost_Hermitian} yields the isomorphism,
\begin{equation}
\label{eq:Isom_H0_PreHolomorphic_nonAbelian}
\bH_{\bar\partial_A,(\varphi,0)}^0 \cong \bH_{A,(\varphi,0)}^0\otimes_\RR\CC,
\end{equation}
where $\bH_{\bar\partial_A,(\varphi,0)}^0$ is as in \eqref{eq:H_dbar_APhi^00} and $\bH_{A,(\varphi,0)}^0$ is as in \eqref{eq:H_APhi^0} with $\Phi = (\varphi,0)$. By the definitions of $\bar\partial_{A,(\varphi,0)}^0$ in \eqref{eq:d0StableComplex} and $\bar\partial_{A,\varphi}^0$ in \eqref{eq:d0StablePair}, we see that
\begin{multline*}
  \bar\partial_{A,(\varphi,0)}^0\zeta
  =
  \begin{pmatrix}
  \bar\partial_A\zeta
  \\
  -\zeta(\varphi,0)
  \end{pmatrix}
  =
  \begin{pmatrix}
  \bar\partial_A\zeta
  \\
  -\zeta\varphi
  \end{pmatrix}
  =
  \bar\partial_{A,\varphi}^0\zeta
  \in \Omega^{0,1}(\fsl(E))\oplus\Omega^0(E),
  \\
  \text{for all } \zeta \in \Omega^0(\fsl(E)),
\end{multline*}
and thus
\begin{equation}
\label{eq:Isom_H0_PreHolomorphic_Holomorphic}
\bH_{\bar\partial_A,(\varphi,0)}^0 = \bH_{\bar\partial_A,\varphi}^0,
\end{equation}
where $\bH_{\bar\partial_A,\varphi}^0$ is defined in \eqref{eq:H_dbar_Avarphi^0bullet}.

On the other hand, by the definitions of $d_{A,(\varphi,0)}^0$ in \eqref{eq:d_APhi^0} with $\Phi = (\varphi,0)$ and $d_{A,\varphi}^0$ in \eqref{eq:d0_projective_vortex_elliptic_deformation_complex}, we see that
\begin{multline*}
  d_{A,(\varphi,0)}^0\xi
  =
  \begin{pmatrix}
  d_A\xi
  \\
  -\xi(\varphi,0)
  \end{pmatrix}
  =
  \begin{pmatrix}
  d_A\xi
  \\
  -\xi\varphi
  \end{pmatrix}
  =
  d_{A,\varphi}^0\xi
  \in \Omega^1(\su(E))\oplus\Omega^0(E),
  \\
  \text{for all } \xi \in \Omega^0(\su(E)),
\end{multline*}
and thus
\begin{equation}
\label{eq:Identify_d0_of_ProjVortex_and_nonAbelian}
d_{A,(\varphi,0)}^0\xi
=
d_{A,\varphi}^0\xi, \quad \text{for all } \xi \in \Omega^0(\su(E)).
\end{equation}
The preceding equation yields the equality,
\begin{equation}
\label{eq:H0NonAbelian_and_ProjVortex}
  \bH_{A,(\varphi,0)}^0 = \bH_{A,\varphi}^0,
\end{equation}
where $\bH_{A,\varphi}^0$ is defined in \eqref{eq:H_Avarphi^0}.
Therefore, we see that
\begin{align*}
  \bH_{A,\varphi}^0\otimes_\RR\CC &= \bH_{A,(\varphi,0)}^0\otimes_\RR\CC
                                    \quad\text{(by \eqref{eq:H0NonAbelian_and_ProjVortex})}
  \\
                                  &\cong \bH_{\bar\partial_A,(\varphi,0)}^0
                                    \quad\text{(by \eqref{eq:Isom_H0_PreHolomorphic_nonAbelian})}
  \\
  &= \bH_{\bar\partial_A,\varphi}^0 \quad\text{(by \eqref{eq:Isom_H0_PreHolomorphic_Holomorphic})},
\end{align*}
and this yields \eqref{eq:HAvarphi^0complexified_isomorphic_HdbarAvarphi^00} when $X$ has complex dimension two. To obtain \eqref{eq:HAvarphi^0complexified_isomorphic_HdbarAvarphi^00} when $X$ has complex dimension three or more, we recall that our proof of Proposition \ref{prop:Itoh_1985_proposition_2-4_SO3_monopole_complex_Kaehler} only uses a hypothesis that $X$ has complex dimension two in order to describe the canonical \spinc structure over a complex K\"ahler surface: this hypothesis can be replaced by the assumption in higher dimensions that $\Phi = (\varphi,\psi)$, for $\varphi\in\Omega^0(E)$ and $\psi\in\Omega^{0,2}(E)$ in general and by $\Phi = (\varphi,0)$ here.

Consider the assertion \eqref{eq:HAvarphi^1_isomorphic_HdbarAvarphi^01}. If $X$ has complex dimension two, then Proposition \ref{prop:Itoh_1985_proposition_2-4_SO3_monopole_complex_Kaehler} yields the isomorphism of real vector spaces,
\[
\bH_{\bar\partial_A,(\varphi,0)}^1 \cong \bH_{A,(\varphi,0)}^1,
\]
where $\bH_{\bar\partial_A,(\varphi,0)}^1$ is given by \eqref{eq:H_dbar_APhi^01} and $\bH_{A,(\varphi,0)}^0$ is given by \eqref{eq:H_APhi^1} with $\Phi = (\varphi,0)$. By the definitions of $\bar\partial_{A,(\varphi,0)}^1$ in \eqref{eq:d1StableComplex} and $\bar\partial_{A,\varphi}^1$ in \eqref{eq:d1StablePair}, we see that
\begin{multline*}
  \bar\partial_{A,(\varphi,0)}^1(\alpha,\sigma,0)
  =
  \begin{pmatrix}
  \bar\partial_A\alpha
  \\
  \bar\partial_A\sigma + \alpha\varphi
  \end{pmatrix}
  =
  \bar\partial_{A,\varphi}^1(\alpha,\sigma)
  \in \Omega^{0,2}(\fsl(E))\oplus\Omega^{0,1}(E),
  \\
  \text{for all } (\alpha,\sigma) \in \Omega^{0,1}(\fsl(E))\oplus\Omega^0(E).
\end{multline*}
Here, we have used Lemma \ref{lem:Witten_dichotomy_linearized_type1_SO3_monopole_equation} to allow us to restrict our attention to triples
\[
  (\alpha,\sigma,\tau) \in \Omega^{0,1}(\fsl(E))\oplus\Omega^0(E)\oplus\Omega^{0,2}(E)
\]
with $\tau\equiv 0$ when $(\alpha,\sigma,\tau) \in \Ker\bar\partial_{A,(\varphi,0)}^1$. Moreover, from our proof of \eqref{eq:HAvarphi^0complexified_isomorphic_HdbarAvarphi^00}, we have
\[
  \Ran \bar\partial_{A,(\varphi,0)}^0 = \Ran \bar\partial_{A,\varphi}^0
  \subset \Omega^{0,1}(\fsl(E))\oplus\Omega^0(E),
\]
and therefore, by the definition of $\bH_{\bar\partial_A,(\varphi,0)}^1$ in \eqref{eq:H_dbar_APhi^01}
and that of $\bH_{\bar\partial_A,\varphi}^1$ in \eqref{eq:H_dbar_Avarphi^0bullet},
\[
  \bH_{\bar\partial_A,(\varphi,0)}^1 = \bH_{\bar\partial_A,\varphi}^1.
\]
On the other hand, by the definition of $d_{A,(\varphi,0)}^1$ in \eqref{eq:d1OfSO3MonopoleComplex} with $\Phi = (\varphi,0)$ and its detailed expression in \eqref{eq:Itoh_1985_2-18_SO3_monopole_complex_Kaehler} and the definition of $d_{A,\varphi}^1$ in \eqref{eq:d1_projective_vortex_elliptic_deformation_complex}, we see that
\begin{multline*}
  d_{A,(\varphi,0)}^1(a,(\sigma,0))
  =
  \begin{pmatrix}
  \Lambda(\bar\partial_Aa' + \partial_Aa'') - i(\varphi\otimes\sigma^*+\sigma\otimes\varphi^*)_0
  \\
  \bar\partial_Aa''
  \\
  \bar\partial_A\sigma+a''\varphi
  \end{pmatrix}
  \\
  =
  d_{A,\varphi}^1(a,\sigma)
  \in \Omega^0(\su(E))\oplus\Omega^{0,2}(\fsl(E))\oplus\Omega^{0,1}(\fsl(E)),
  \\
  \text{for all } (a,\sigma) \in \Omega^1(\su(E))\oplus\Omega^0(E).
\end{multline*} 
By combining the preceding relation between $d_{A,(\varphi,0)}^1$ and $d_{A,\varphi}^1$ with the identification of $d_{A,(\varphi,0)}^0$ and $d_{A,\varphi}^0$ given in \eqref{eq:Identify_d0_of_ProjVortex_and_nonAbelian}, we obtain the equality
\begin{equation}
\label{eq:H1NonAbelian_and_ProjVortex}
  \bH_{A,(\varphi,0)}^1 = \bH_{A,\varphi}^1,
\end{equation}
where $\bH_{A,(\varphi,0)}^1$ is defined in \eqref{eq:H_APhi^1}
and $\bH_{A,\varphi}^1$ in \eqref{eq:H_Avarphi^1}.

Therefore, by combining the preceding isomorphisms and equalities of real vector spaces, we see that
\[
  \bH_{A,\varphi}^1 = \bH_{A,(\varphi,0)}^1
  \cong \bH_{\bar\partial_A,(\varphi,0)}^1 = \bH_{\bar\partial_A,\varphi}^1,
\]
and this yields \eqref{eq:HAvarphi^1_isomorphic_HdbarAvarphi^01} when $X$ has complex dimension two. To obtain \eqref{eq:HAvarphi^1_isomorphic_HdbarAvarphi^01} when $X$ has complex dimension three or more, we recall that our hypothesis in Proposition \ref{prop:Itoh_1985_proposition_2-4_SO3_monopole_complex_Kaehler} that $X$ has complex dimension two can be replaced by the assumption that $\Phi=(\varphi,0)$.

Consider the assertion \eqref{eq:HAvarphi^2_isomorphic_HdbarAvarphi^02+HdbarAvarphi^00}. If $X$ has complex dimension two, then Proposition \ref{prop:Itoh_1985_proposition_2-3_SO3_monopole_complex_Kaehler} yields the isomorphism of real vector spaces,
\[
  \bH_{A,(\varphi,0)}^2  \cong \bH_{\bar\partial_A,(\varphi,0)}^2 \oplus \bH_{A,(\varphi,0)}^0,
\]
where $\bH_{\bar\partial_A,(\varphi,0)}^2$ is given by \eqref{eq:H_dbar_APhi^02} and $\bH_{A,(\varphi,0)}^2$ is given by \eqref{eq:H_APhi^2} with $\Phi = (\varphi,0)$. From our proof of \eqref{eq:HAvarphi^0complexified_isomorphic_HdbarAvarphi^00}, we recall that
\[
  \bH_{A,(\varphi,0)}^0 = \bH_{A,\varphi}^0, 
\]
and so we focus on the degree-two harmonic spaces. From our proof of \eqref{eq:HAvarphi^1_isomorphic_HdbarAvarphi^01}, we have the equality
\[
  \Ran \bar\partial_{A,(\varphi,0)}^1 = \Ran \bar\partial_{A,\varphi}^1
  \subset \Omega^{0,2}(\fsl(E))\oplus\Omega^{0,1}(E)
\] 
and therefore,
\[
  \bH_{\bar\partial_A,(\varphi,0)}^2 = \bH_{\bar\partial_A,\varphi}^2.
\]
Moreover, from our proof of \eqref{eq:HAvarphi^1_isomorphic_HdbarAvarphi^01}, we also have the equality
\[
  \Ran d_{A,(\varphi,0)}^1 = \Ran d_{A,\varphi}^1
  \subset
  \Omega^0(\su(E))\oplus\Omega^{0,2}(\fsl(E))\oplus\Omega^{0,1}(E),
\] 
and therefore,
\[
  \bH_{A,(\varphi,0)}^2 = \bH_{A,\varphi}^2. 
\]
Therefore, by combining the preceding isomorphisms and equalities of real vector spaces, we see that
\begin{align*}
  \bH_{A,\varphi}^2 &= \bH_{A,(\varphi,0)}^2
  \\
  &\cong \bH_{\bar\partial_A,(\varphi,0)}^2\oplus \bH_{A,(\varphi,0)}^0
  \\
  &= \bH_{\bar\partial_A,\varphi}^2\oplus \bH_{A,\varphi}^0,
\end{align*}
and this yields \eqref{eq:HAvarphi^2_isomorphic_HdbarAvarphi^02+HdbarAvarphi^00} when $X$ has complex dimension two. The proof of \eqref{eq:HAvarphi^2_isomorphic_HdbarAvarphi^02+HdbarAvarphi^00} when $X$ has complex dimension three or more follows from two observations. First, we recall that our hypothesis in Proposition \ref{prop:Itoh_1985_proposition_2-4_SO3_monopole_complex_Kaehler} that $X$ has complex dimension two can be replaced by the assumption that $\Phi=(\varphi,0)$. Second, when the complex dimension of $X$ becomes three or more, the definitions of $\bH_{A,\varphi}^2$ and $\bH_{\bar\partial_A,\varphi}^2$ both change by imposing the same condition,
\[
  \bar\partial_{A,\varphi}^2(\beta,\nu) = 0 \quad\text{on } (\beta,\nu) \in \Omega^{0,2}(\fsl(E))\oplus\Omega^{0,1}(E).
\]
The isomorphism \eqref{eq:HAvarphi^k_isomorphic_HdbarAvarphi^0k} is immediate from the definitions of these harmonic spaces since the higher-degree vector spaces and differentials in the associated elliptic complexes are identical.
\end{proof}  

Theorem \ref{thm:Kobayashi_7-2-21_pairs} implies the following analogue of Lemmas \ref{lem:H0_Of_NonAbelianMonopoleComplex_Vanishes} and \ref{lem:H0_Of_ProjectiveComplex_Vanishes}.

\begin{lem}[Condition for vanishing of $\bH_{\bar\partial_A,\varphi}^0$]
\label{lem:H0_Of_Holomorphic_Vanishes}
Let $(E,h)$ be a rank-two Hermitian vector bundle over a complex, closed, K\"ahler manifold $X$ and let $p\in(n,\infty)$, where $n$ is the complex dimension of $X$. Let $(A,\varphi)\in\sA(E,h)\times W^{1,p}(E)$ be a solution to the projective vortex equations \eqref{eq:SO(3)_monopole_equations_almost_Hermitian_alpha} and $(\bar\rd_A,\varphi)\in\sA^{0,1}(E)\times W^{1,p}$ be the holomorphic pair obeying \eqref{eq:Holomorphic_pair} that is defined by $(A,\varphi)$. If $\varphi\not\equiv 0$, then
\[
\bH_{\bar\partial_A,\varphi}^0=(0),
\]
where $\bH_{\bar\partial_A,\varphi}^0$ is the harmonic space defined in \eqref{eq:H_dbar_Avarphi^0bullet}.
\end{lem}

\begin{proof}
The conclusion follows from  Lemma \ref{lem:H0_Of_ProjectiveComplex_Vanishes} and the isomorphism
\[
  \bH_{A,\varphi}^0\otimes_\RR\CC \cong \bH_{\bar\partial_A,\varphi}^0
\]
given by Theorem \ref{thm:Kobayashi_7-2-21_pairs}.
\end{proof}

\chapter[K\"ahler structure on moduli space of Hermitian--Einstein connections]{Complex structure and K\"ahler metric on the moduli space of projectively Hermitian--Einstein connections over a complex K\"ahler manifold}
\label{chap:Complex_structure_and_Kaehler_metric_moduli_space_HE_connections_complex_Kaehler_manifold}
In preparation for our development of the corresponding results for the moduli space of non-Abelian monopoles, we recall in this chapter some important results on the existence of a complex structure and K\"ahler metric on the moduli space $M(E,h,\omega)$ of projectively Hermitian--Einstein connections on a Hermitian vector bundle $(E,h)$ over a complex K\"ahler manifold $(X,\omega)$, following Itoh \cite{Itoh_1983, Itoh_1985}, Kim \cite{Kim_1987}, and Kobayashi \cite{Kobayashi_differential_geometry_complex_vector_bundles}. The Hitchin--Kobayashi correspondence provides a real analytic isomorphism from $M(E,h,\omega)$ onto the moduli space $\cM(E,\omega)$ of stable holomorphic structures on $E$, an open subspace of the moduli space $\cM(E)$ of all holomorphic structures on $E$. The existence of a K\"ahler metric on the open subspace $M_\reg^*(E,h,\omega)$ comprising smooth points $[A]$ of $M(E,h,\omega)$ represented by non-split connections $A$ on $E$ is due to Itoh \cite{Itoh_1988} and Kobayashi \cite{Kobayashi_differential_geometry_complex_vector_bundles} in more generality.

In Sections \ref{sec:Moduli_space_simple_bundles} and \ref{sec:Moduli_space_stable_bundles}, we define the moduli space of holomorphic bundles and its subspaces, the moduli space of simple and stable holomorphic bundles with a fixed determinant bundle, following Friedman and Morgan \cite[Chapter IV, Section 4]{FrM} and Kobayashi \cite[Section 7.4]{Kobayashi_differential_geometry_complex_vector_bundles}.  We define projectively Hermitian--Einstein connections and their deformation complex in Section \ref{sec:HE_connections} and introduce their moduli space in Section \ref{sec:Moduli_space_HE_connections}.  After defining the Hitchin--Kobayashi correspondence between the moduli space of projectively Hermitian--Einstein connections and the moduli space of semistable holomorphic bundles in Section \ref{sec:HK_correspondence_between_HE_connections_and_semistable_bundles}, we state a result
that the Hitchin--Kobayashi correspondence defines a real analytic embedding of the moduli space of projectively Hermitian--Einstein connections into the moduli space of simple holomorphic bundles
in Section \ref{sec:Moduli_space_HE_connections_open_subspace_moduli_space_simple_holomorphic_structures}.
We review the Marsden--Weinstein symplectic reduction and construction of a symplectic structure on a quotient in Section \ref{sec:Marsden-Weinstein_symplectic_quotient_construction} and then apply this construction to the moduli space of projectively Hermitian--Einstein connections
in Section \ref{sec:Marsden-Weinstein_reduction_moduli_space_HE_connections_symplectic_quotient}, giving an explicit expression for the resulting symplectic form.

\section{Moduli spaces of simple holomorphic vector bundles}
\label{sec:Moduli_space_simple_bundles}
For a complex manifold $X$ of dimension $n$ and a constant $p\in (n,\infty)$, we update our previous definition \eqref{eq:Affine_space_01-connections} of $\sA^{0,1}(E)$ for a complex vector bundle $E$ over a complex surface $X$ in Section \ref{sec:Elliptic_deformation_complex_for_pre-holomorphic_pair_equations_complex_surface} as the affine space of $(0,1)$-connections\footnote{See Section \ref{sec:Elliptic_deformation_complex_holomorphic_pair_equations}.} $\bar\partial_E$ on $E$ of class $W^{1,p}$
by allowing $E$ to be a complex vector bundle over a complex manifold $X$ of any dimension but requiring that each $(0,1)$-connection induce a fixed smooth, holomorphic structure $\bar\partial_{E_d}$ on the complex line bundle $\det E$, that is,
\begin{equation}
  \label{eq:Holomorphic_structure_fixed_determinant}
  \bar\partial_{\det E} = \bar\partial_{E_d} \text{ on } \Omega^0(\det E).
\end{equation}
Given \eqref{eq:Holomorphic_structure_fixed_determinant}, the condition $F_{\bar\partial_E} = 0  \in \Omega^{0,2}(\gl(E))$ (see equation \eqref{eq:Donaldson_Kronheimer_2-1-50})
is equivalent to $(F_{\bar\partial_E})_0 = 0  \in \Omega^{0,2}(\fsl(E))$ since
$F_{E_d} = \bar\partial_{E_d}\circ\bar\partial_{E_d} \equiv 0 \in \Omega^{0,2}(\det E)$ by assumption that $\bar\partial_{E_d}$ is a holomorphic structure on $\det E$, where for $E$ of complex rank $r$,
\[
  (F_{\bar\partial_E})_0 := F_{\bar\partial_E} - \frac{1}{r}(\tr F_{\bar\partial_E})\,\id_E.
\]
Thus, $\sA^{0,1}(E)$ is modeled on the complex Banach space $W^{1,p}(\Lambda^{0,1}(\fsl(E)))$. (See Friedman and Morgan \cite[Section 4.1.2]{FrM} and \cite[Section 4.1.6]{FrM} for further discussion of the case of fixed holomorphic structures on the determinant line bundle.) 

For the principal $\SL(r,\CC)$-bundle $\SL(E)$ over $X$ defined in \eqref{eq:AutomorphismBundles_SL(E)},
the group $W^{2,p}(\SL(E))$ of $W^{2,p}$ determinant one, complex gauge transformations acts on $\sA^{0,1}(E)$ by
\begin{equation}
\label{eq:SL(E)_Action_on_(0,1)conn}
W^{2,p}(\SL(E))\times \sA^{0,1}(E) \ni (v,\bar\partial_E)
\mapsto
v^{-1}\circ \bar\partial_E\circ v \in \sA^{0,1}(E).
\end{equation}
By adapting the proofs \mutatis due to Freed and Uhlenbeck of \cite[Appendix A, Proposition A.2, p. 160 and Proposition A.3, p. 161]{FU}, one can show that $W^{2,p}(\SL(E))$ is a complex Banach Lie group and that its action \eqref{eq:SL(E)_Action_on_(0,1)conn} on $\sA^{0,1}(E)$ is holomorphic. We write $[\bar\partial_E]$ for the gauge-equivalence class of the $(0,1)$-connection $\bar\partial_E$ under the action \eqref{eq:SL(E)_Action_on_(0,1)conn}
\label{page:EquivClassOf(0,1)Connection}
and let $\sB^{0,1}(E) := \sA^{0,1}(E)/W^{2,p}(\SL(E))$
\label{page:Quotient_space_(0,1)_connections}
denote the corresponding configuration space equipped with quotient topology and which, unlike the quotient space of unitary connections, need not be Hausdorff. We let
\begin{equation}
  \label{eq:Moduli_space_holomorphic_structures}
  \cM(E) := \left\{\bar\partial_E \in \sA^{0,1}(E): F_{\bar\partial_E} = 0\right\}/W^{2,p}(\SL(E))
\end{equation}
denote the moduli space of holomorphic structures on $E$ modulo the group $W^{2,p}(\SL(E))$, equipped with the quotient topology. We have the following analogue for $(0,1)$-connections of Definition \ref{defn:Reducible_split_trivial-stabilizer_unitary_connection} \eqref{item:Split_unitary_connection} for split unitary connections.

\begin{defn}[Split $(0,1)$-connections on complex vector bundles over almost complex manifolds]
\label{defn:Split_01-connection}
Let $E$ be a $W^{2,p}$ complex vector bundle of rank $r\geq 2$ over a connected, almost complex manifold of real dimension $2n$ with $p\in(n,\infty)$. A $W^{1,p}$ $(0,1)$-connection $\bar\partial_E$ on $E$ is \emph{split} if $\bar\partial_E = \bar\partial_{E_1}\oplus\partial_{E_2}$ with respect to a decomposition
\[
  E=E_1\oplus E_2
\]
as a direct sum of proper, $W^{2,p}$ complex subbundles of $E$, where $\bar\partial_{E_i}$ is a $W^{1,p}$ $(0,1)$-connection on $E_i$, for $i=1,2$, and is \emph{non-split} otherwise. If the pair $(E,\bar\partial_E)$ is smooth (respectively, analytic), then each pair $(E_i,\bar\partial_{E_i})$ is smooth (respectively, analytic), for $i=1,2$.
\end{defn}

For consistency with customary notation in complex geometry, we shall occasionally write $\End(E)$ for the vector bundle of endomorphisms of a complex vector bundle $E$ over a manifold $X$ and write $\End_0(E)$ for the vector subbundle of trace free endomorphisms of $E$. In this monograph, we generally denote these endomorphism bundles by $\gl(E)$ and $\fsl(E)$, respectively.

\begin{defn}[Simple $(0,1)$-connections on smooth, complex vector bundles over complex manifolds]
\label{defn:Simple_01-connection}
(Compare Friedman and Morgan \cite[Section 4.1.2, Definition 1.5, p. 283]{FrM} or Kobayashi \cite[Section 5.7, p. 159]{Kobayashi_differential_geometry_complex_vector_bundles}.)
Continue the setting of Definition \ref{defn:Split_01-connection}. An $(0,1)$-connection $\bar\partial_E$ on $E$ is \emph{simple} if one (and thus all) of the following equivalent conditions hold:
\begin{enumerate}
\item\label{item:End0(E)_simple_01-connection} The only $W^{2,p}$ section of $\End_0(E)$ commuting with $\bar\partial_E$ is the zero endomorphism.
  
\item\label{item:End(E)_simple_01-connection} The only $W^{2,p}$ sections of $\End(E)$ commuting with $\bar\partial_E$ are of the form $b\,\id_E$, for $b \in \CC$.
  
\item\label{item:GL(E)_simple_01-connection} The only $W^{2,p}$ sections of $\GL(E)$ commuting with $\bar\partial_E$ are of the form $c\,\id_E$, for $c \in \CC^*$.
  
\item\label{item:SL(E)_simple_01-connection} The only $W^{2,p}$ sections of $\SL(E)$ commuting with $\bar\partial_E$ are of the form $\varrho\,\id_E$, where $\varrho$ is an element of the group $C_r$ of $r$-th roots of unity.
\end{enumerate}
\end{defn}

\begin{rmk}[Non-connected base manifolds]
One could extend Definition \ref{defn:Simple_01-connection} to the setting of non-connected base manifolds by requiring the restriction of the $(0,1)$-connection $\bar\partial_E$ to each component to be simple.
However, we shall not need this level of generality.
\end{rmk}

\begin{rmk}[On induced connections and $(0,1)$-connections on $\End(E)$]
\label{rmk:SimplicityAndKernelOfConn_on_EndE}
We describe an alternate formulation of the condition that $w\in\Omega^0(\End(E))$ commutes with $\bar\partial_A$ in Definition \ref{defn:Simple_01-connection}. Just as a connection $A$ on a vector bundle $E$ induces a connection $A^{\ad}$ on $\End(E)$ as described in Remark \ref{rmk:InducedConn_on_EndE}, an $(0,1)$-connection $\bar\partial_E$ on $E$ defines an $(0,1)$-connection $\bar\partial_{\End(E)}$ on $\End(E)$ by
\begin{equation}
\label{eq:(0,1)Connection_Induced_On_Endomorphism_Bundle}
(\bar\partial_{\End(E)} w)(s) = \bar\partial_E(w(s)) - w(\bar\partial_Es),
\quad\text{for all } w\in\Omega^0(\End(E)) \text{ and } s \in \Omega^0(E).
\end{equation}
Thus, $w$ commutes with  $\bar\partial_E$ if and only if $\bar\partial_{\End(E}w =0$. We will often  abbreviate $\bar\partial_{\End(E)}$ by $\bar\partial_E$ when the meaning is clear from the context.
\end{rmk}

We now verify the claimed equivalence of Items \eqref{item:End0(E)_simple_01-connection}, \eqref{item:End(E)_simple_01-connection}, \eqref{item:GL(E)_simple_01-connection}, and \eqref{item:SL(E)_simple_01-connection} in Definition \ref{defn:Simple_01-connection}.

\begin{lem}[Equivalence of conditions for simple $(0,1)$-connections in Definition \ref{defn:Simple_01-connection}]
\label{lem:EquivalenceOfDefinitionOfSimpleConnections}
Continue the setting of Definition \ref{defn:Split_01-connection}. The conditions for an $(0,1)$-connection $\bar\partial_E$ on $E$ to be simple in Definition \ref{defn:Simple_01-connection} are equivalent.
\end{lem}

\begin{proof}
We shall first prove the equivalence of \eqref{item:End0(E)_simple_01-connection} and \eqref{item:End(E)_simple_01-connection}. Assume \eqref{item:End(E)_simple_01-connection} holds and that $w\in\End_0(E)$ commutes with $\bar\partial_E$. Since $\End_0(E)\subset \End(E)$, then \eqref{item:End(E)_simple_01-connection} implies that $w$ has the form $b\,\id_E$, for some $b \in \CC$, and because $\tr_Ew = rb = 0$, we must have $b=0$ and thus $w=0$. This proves that \eqref{item:End(E)_simple_01-connection} $\implies$ \eqref{item:End0(E)_simple_01-connection}. Conversely, assume \eqref{item:End0(E)_simple_01-connection} holds and that $v\in\End(E)$ commutes with $\bar\partial_E$. Define $w := \pi_0v = v - \frac{1}{r}(\tr_Ev)\,\id_E$ to be the trace free part of $v$, where $\pi_0$ is the projection from $\End(E) = \End_0(E) \oplus \ubarCC\,\id_E$ onto $\End_0(E)$.

We claim that $v$ commutes with $\bar\partial_E$ if and only if $w$ commutes with $\bar\partial_E$ and  $\bar\partial(\tr_E v)=0$. Choose a Hermitian metric $h$ on $E$, let $A$ be the Chern connection on $E$ defined by the pair $(\bar\partial_E,h)$, and let $A^{\ad}$ be the unitary connection on $\End(E)$ induced by $A$ as in \eqref{eq:Connection_Induced_On_Endomorphism_Bundle}. The definition of $d_{A^{\ad}}$ in \eqref{eq:Connection_Induced_On_Endomorphism_Bundle} implies that $d_{A^\ad}\id_E=0$.  Because $A^{\ad}$ is a unitary connection on $\End(E)$, the operator $d_{A^{\ad}}$ preserves both the subbundle $\CC\,\id_E$ and its orthogonal complement $\End_0(E)$ in $\End(E)$.  The $(0,1)$-component of $d_{A^\ad}$ given by the $(0,1)$-connection $\bar\partial_{\End(E)}$ defined in \eqref{eq:(0,1)Connection_Induced_On_Endomorphism_Bundle} therefore admits a direct sum decomposition,
\begin{multline}
\label{eq:DirectSumDecompositionOfEnd(0,1)Connection}
\bar\partial_{\End(E)}
=
\bar\partial_{\End_0(E)}\oplus\bar\partial: W^{2,p}(\End_0(E))\oplus W^{2,p}(X,\CC)
\\
\to W^{1,p}(\La^{0,1}\otimes \End_0(E))\oplus W^{1,p}(\La^{0,1}).
\end{multline}
Thus, $\bar\partial_{\End(E)}v=0$ if and only if $\bar\partial_{\End_0(E)}w=0$ and $\bar\partial \tr_E v = 0$.  By Remark \ref{rmk:SimplicityAndKernelOfConn_on_EndE}, the property that $v$ (respectively, $w$) commutes with $\bar\partial_E$ is equivalent to the property that $\bar\partial_{\End(E)}v=0$ (respectively, $\bar\partial_{\End(E)}w=0$). Therefore, $v$ commutes with $\bar\partial_E$ if and only if $w$ commutes with $\bar\partial_E$ and $\bar\partial(\tr_E v) = 0$ and this proves the claim.

Hence, returning to our proof that \eqref{item:End0(E)_simple_01-connection} $\implies$ \eqref{item:End(E)_simple_01-connection}, we see that $w=0$ by \eqref{item:End0(E)_simple_01-connection} and therefore $v = f\,\id_E$ with $rf := \tr_Ev \in W^{2,p}(X,\CC)$ and $\bar\partial f = 0$ since $v$ commutes with $\bar\partial_E$ by assumption. Since $X$ is connected, we see that $f$ is constant
by Cirici and Wilson \cite[Corollary 1, p. 7]{Cirici_Wilson_2020_almost_hermitian_identities}. (Their result is an extension to almost complex manifolds of the fact that if $f$ is a holomorphic function on closed, connected complex manifold, then $f$ is constant.) This proves that \eqref{item:End0(E)_simple_01-connection} $\implies$ \eqref{item:End(E)_simple_01-connection}.
  
Clearly, \eqref{item:End(E)_simple_01-connection} $\implies$ \eqref{item:GL(E)_simple_01-connection} $\implies$ \eqref{item:SL(E)_simple_01-connection} in Definition \ref{defn:Simple_01-connection}. It suffices to directly prove that \eqref{item:SL(E)_simple_01-connection} $\implies$ \eqref{item:End(E)_simple_01-connection}, but it is easier to separately prove that \eqref{item:SL(E)_simple_01-connection} $\implies$ \eqref{item:GL(E)_simple_01-connection} and then \eqref{item:GL(E)_simple_01-connection} $\implies$ \eqref{item:End(E)_simple_01-connection}.

Assume that \eqref{item:SL(E)_simple_01-connection} holds and let $v$ be a $W^{2,p}$ section of $\End(E)$ commuting with $\bar\partial_E$. Thus $\det v$ is a $W^{2,p}$ section of $\End(\det E)$, that is, a $W^{2,p}$ $\CC$-valued function on $X$. Because $v$ commutes with $\bar\partial_E$, we see that $\det v$ commutes with
the $(0,1)$-connection induced by $\bar\partial_E$ on $\det E$, so $\bar\partial(\det v)=0$ and thus $\det v$ is constant by Cirici and Wilson \cite[Corollary 1, p. 7]{Cirici_Wilson_2020_almost_hermitian_identities}. Therefore, we may choose $\lambda\in\CC$ such that $\lambda^r=\det v$. If $\det v \neq 0$ and thus $\lambda\neq 0$, then $\lambda^{-1}v$ is a $W^{2,p}$ section of $\SL(E)$ which commutes with $\bar\partial_E$, so \eqref{item:SL(E)_simple_01-connection} implies that $\lambda^{-1}v=\varrho\,\id_E$ for some $\varrho\in C_r$. Hence, $v = c\,\id_E$ with $c = \lambda\varrho \in \CC^*$.  Thus, \eqref{item:SL(E)_simple_01-connection} implies \eqref{item:GL(E)_simple_01-connection}.

We now assume \eqref{item:GL(E)_simple_01-connection} and prove that \eqref{item:GL(E)_simple_01-connection} $\implies$ \eqref{item:End(E)_simple_01-connection}.  Let $w$ be  a $W^{2,p}$ section of $\End(E)$ that commutes with $\bar\partial_E$.  Then, $\Exp(tw)$ is a $W^{2,p}$ section of $\GL(E)$, for all $t\in\RR$, which also commutes with  $\bar\partial_E$. (Here $\Exp(\cdot)$ denotes the pointwise exponential map from $\End(E)$ to $\GL(E)$.)  Because we are assuming that \eqref{item:GL(E)_simple_01-connection} holds, for all $t\in\RR$ there is an element $k(t)\in\CC^*$ such that $\Exp(tu)=k(t)\,\id_E$. Since the restriction of $\Exp(tw)=k(t)\,\id_E$ to each fiber is a smooth curve in the general linear group and $\GL(E)$ is a smooth fiber bundle, $k(t)$ depends smoothly on $t$. Therefore
\[
w=
\frac{d}{dt}\Exp(t w)|_{t=0}
=k'(0)\,\id_E,
\]
is diagonal, so \eqref{item:End(E)_simple_01-connection} holds with $b = k'(0)$. Thus, \eqref{item:GL(E)_simple_01-connection} $\implies$ \eqref{item:End(E)_simple_01-connection} and this completes the proof of the lemma.
\end{proof}

We have seen in the proof of Lemma \ref{lem:EquivalenceOfDefinitionOfSimpleConnections} that there is a canonical isomorphism of complex vector bundles, $\End(E) \cong \End_0(E)\oplus \ubarCC\,\id_E$.
If $\sE$ is a holomorphic vector bundle over a complex manifold $X$ then, because the holomorphic map of vector bundles $\tr_\sE: \End(\sE) \to \ubarCC$ admits a holomorphic right inverse, $f\mapsto (1/r)f\,\id_\sE$, the exact sequence of holomorphic vector bundles,
\[
\begin{CD}
0 @>>> \End_0(\sE) @>>> \End(\sE) @> \tr_\sE >> \ubarCC @>>> 0
\end{CD}
\]
splits and there is an isomorphism of holomorphic vector bundles, $\End(\sE) \cong \End_0(\sE)\oplus \ubarCC\,\id_E$.
This holomorphic decomposition of $\End(\sE)$ immediately yields the direct sum decomposition of $\bar\partial_{\End(E)}$ appearing in equation \eqref{eq:DirectSumDecompositionOfEnd(0,1)Connection} without relying on the Chern connection as in the proof of Lemma \ref{lem:EquivalenceOfDefinitionOfSimpleConnections}.
Hence, $v\in\Omega^0(\End(E))$ is holomorphic if and only if $\pi_0v\in \Omega^0(\End_0(E))$ and $\tr_Ev\in\Omega^0(\ubarCC)$ are holomorphic.

The following definition of the stabilizer $\bar\partial_E$ is implicit in Definition \ref{defn:Simple_01-connection}.

\begin{defn}[Stabilizer of $\bar\partial_E$]
\label{defn:StabilizerOf(0,1)Conn}
Continue the setting of Definition \ref{defn:Split_01-connection}. The \emph{stabilizer} $\Stab(\bar\partial_E)$ of the $(0,1)$-connection $\bar\partial_E$ under the action \eqref{eq:SL(E)_Action_on_(0,1)conn} is the subgroup of elements of $W^{2,p}(\SL(E))$ commuting with $\bar\partial_E$.
\end{defn}

The stabilizer has the following Lie group structure.

\begin{lem}[Lie algebra and complex structure of $\Stab(\bar\partial_E)$]
\label{lem:Stab(rdE)_is_Lie_Group}
Continue the setting of Definition \ref{defn:Split_01-connection}, but assume in addition that the manifold is complex. The stabilizer $\Stab(\bar\partial_E)$ of an $(0,1)$-connection $\bar\partial_E$ of class $W^{1,p}$ is a complex Lie subgroup of $W^{2,p}(\SL(E))$ with Lie algebra given by (compare the forthcoming definition \eqref{eq:bHdbarA0bullet})
\begin{equation}
\label{eq:bHdbarA0}
\bH_{\bar\partial_E}^0:=
\Ker \left( \bar\partial_{\End(E)}: W^{2,p}(\fsl(E)) \to W^{1,p}(\La^{0,1}(\fsl(E)) \right),
\end{equation}
where the $(0,1)$-connection $\bar\partial_{\End(E)}$ is defined in \eqref{eq:(0,1)Connection_Induced_On_Endomorphism_Bundle}.
\end{lem}

\begin{proof}
According to Hilgert and Neeb \cite[Lemma 10.1.5 (i), p. 361]{Hilgert_Neeb_structure_geometry_lie_groups}, the stabilizer $\Stab(\bar\partial_E)$ is a closed subgroup of $W^{2,p}(\SL(E))$ (and thus a Lie group) and
furthermore, according to Akhiezer \cite[Section 1.2, p. 8, 6th paragraph]{Akhiezer_lie_group_actions_complex_analysis}, it is necessarily a complex Lie group as in Hilgert and Neeb \cite[Definition 15.1.1 (a), p. 566]{Hilgert_Neeb_structure_geometry_lie_groups}. The argument given in
Lemma \ref{lem:LieGroup_Structure_of_Stab(A)} identifying $\bH_A^0$ with the Lie algebra of $\Stab(A)$ translates directly to identify $\bH_{\bar\partial_E}^0$ with the Lie algebra of $\Stab(\bar\partial_E)$.
\end{proof}

We have the following characterization of $(0,1)$-connections with central stabilizer.

\begin{lem}[Simple $(0,1)$-connections and cohomology groups]
\label{lem:Simple_01-connection}
(Compare Kobayashi \cite[Section 7.3, p. 233, second last line]{Kobayashi_differential_geometry_complex_vector_bundles}.) Continue the notation of Definition \ref{defn:Simple_01-connection}.
The following conditions on an $(0,1)$-connection $\bar\partial_E$ on $E$ are equivalent:
\begin{enumerate}
\item
\label{item:Simple_01-connectionSimple}
The $(0,1)$-connection $\bar\partial_E$ is simple in the sense of Definition \ref{defn:Simple_01-connection}.
\item
\label{item:Simple_01-connectionCentralStabilizer}
The stabilizer of $\bar\partial_E$ is central, so $\Stab(\bar\partial_E) = \{\varrho\,\id_E:\varrho \in C_r\}$, where $C_r$ is the group of $r$-th roots of unity and $E$ has complex rank $r$.
\item
\label{item:Simple_01-connectionVanishingH0}
$\bH_{\bar\partial_E}^0 = (0)$, where the harmonic space $\bH_{\bar\partial_E}^0$ is as in \eqref{eq:bHdbarA0}.
\end{enumerate}
\end{lem}

\begin{proof}
The equivalence of Items \eqref{item:Simple_01-connectionSimple} and \eqref{item:Simple_01-connectionCentralStabilizer} follows immediately from Definition \ref{defn:Simple_01-connection} \eqref{item:SL(E)_simple_01-connection} and Definition \ref{defn:StabilizerOf(0,1)Conn}.

We now prove the equivalence of Items \eqref{item:Simple_01-connectionSimple} and \eqref{item:Simple_01-connectionVanishingH0}. An $(0,1)$-connection $\bar\partial_E$ on $E$ is simple in the sense of Definition \ref{defn:Simple_01-connection} \eqref{item:End0(E)_simple_01-connection} if and only if the only $W^{2,p}$ section $\zeta$ of $\fsl(E)$ that commutes with $\bar\partial_E$ is the zero endomorphism of $E$. The definition of $\bar\partial_{\End(E)}$ in \eqref{eq:(0,1)Connection_Induced_On_Endomorphism_Bundle} identifies
sections of $\fsl(E)$ which commute with $\bar\partial_E$ with the kernel of $\bar\partial_{\End(E)}$, which is $\bH_{\bar\partial_E}^0$ by the definition \eqref{eq:bHdbarA0}.  Hence, an $(0,1)$-connection $\bar\partial_E$ is simple if and only if $\bH_{\bar\partial_E}^0=(0)$, giving the equivalence of Items \eqref{item:Simple_01-connectionSimple} and \eqref{item:Simple_01-connectionVanishingH0}.
\end{proof}

We define the following open subspaces of $\cM(E)$ in \eqref{eq:Moduli_space_holomorphic_structures},
\begin{subequations}
  \label{eq:Moduli_subspaces_holomorphic_structures}
\begin{align}
   \label{eq:Moduli_space_simple_holomorphic_structures}
  \cM^*(E) &:= \left\{[\bar\partial_E] \in \cM(E): (E,\bar\partial_E) \text{ is simple} \right\},
  \\
   \label{eq:Moduli_space_regular_holomorphic_structures}
  \cM_\reg(E) &:= \left\{[\bar\partial_E] \in \cM(E): \bH_{\bar\partial_E}^2 = (0) \right\},
  \\
  \label{eq:Moduli_space_regular_simple_holomorphic_structures}
  \cM_\reg^*(E) &:= \cM^*(E) \cap \cM_\reg(E),
\end{align}
\end{subequations}
of simple or regular holomorphic vector bundles with fixed smooth, holomorphic structure on $\det E$ modulo $W^{2,p}(\SL(E))$, where the harmonic space $\bH_{\bar\partial_E}^2$ is defined by \eqref{eq:bHdbarA0bullet}. The property that $\cM_\reg(E)$ is an open subspace of $\cM(E)$ follows from the Implicit Mapping Theorem, while the fact that $\cM^*(E)$ is an open subspace of $\cM(E)$ follows from the forthcoming Lemma \ref{lem:Openness_subspace_simple_01-connections}. By analogy with our definition \eqref{eq:Moduli_space_simple_holomorphic_structures} of $\cM^*(E)$, we define
\[
  \sA^{0,1;*}(E) := \left\{\bar\partial_E \in \sA^{0,1}(E): \bar\partial_E \text{ is simple} \right\}.
\]  
We have the important

\begin{lem}[Openness of the subspace of simple $(0,1)$-connections]
\label{lem:Openness_subspace_simple_01-connections}
Let $E$ be a smooth Hermitian vector bundle over a closed, complex, Hermitian  manifold $X$ of complex dimension $n$ and let $p\in(n,\infty)$. Then the subspace $\sA^{0,1;*}(E)$ of simple $(0,1)$-connections is an open subspace of the affine space $\sA^{0,1}(E)$ of $(0,1)$-connections.
\end{lem}

\begin{proof}
Let $\bar\partial_E \in \sA^{0,1;*}(E)$, so $\bar\partial_E$ is a $W^{1,p}$ simple holomorphic structure on $E$. According to Lemma \ref{lem:Simple_01-connection}, the condition that $\bar\partial_E$ be simple is equivalent to $\bH_{\bar\partial_E}^0 = (0)$, where $\bH_{\bar\partial_E}^0 \subset W^{2,p}(\fsl(E))$ is the finite-dimensional vector subspace defined in \eqref{eq:bHdbarA0}. It follows from Rudin \cite[Theorem 4.12, p. 99]{Rudin} that
\begin{multline*}
  \Ker\left(\bar\partial_E:W^{2,p}(\fsl(E))\to W^{1,p}(\Lambda^{0,1}(\fsl(E)))\right)
  \\
  =
  \Ran\left(\bar\partial_E^*:W^{3,p}(\Lambda^{0,1}\fsl(E)) \to W^{2,p}(\fsl(E))\right)^\perp,
\end{multline*}
where $\perp$ denotes $L^2$-orthogonal complement. Thus, $\bH_{\bar\partial_E}^0 = (0)$ if and only if $\Ran \bar\partial_E^* = W^{2,p}(\fsl(E))$ and the latter is an open condition in the sense that if $\alpha \in W^{1,p}(\Lambda^{0,1}(\fsl(E)))$ is small, then $\Ran(\bar\partial_E+\alpha)^*= W^{2,p}(\fsl(E))$ and this yields the conclusion.
\end{proof}

\begin{thm}[Local topology of the moduli space of holomorphic vector bundles over a complex manifold near a simple bundle]
\label{thm:Kobayashi_7_3_17}  
(See Kobayashi \cite[Theorem 7.3.17, p. 234]{Kobayashi_differential_geometry_complex_vector_bundles}.)
Let $E$ be a smooth Hermitian vector bundle, with a fixed smooth holomorphic structure on $\det E$, over a closed, complex, Hermitian\footnote{We omit the hypothesis that $X$ is K\"ahler since that only appears to be used by Kobayashi to define harmonic spaces $\bH_{\bar\partial_E}^\bullet$ representing the cohomology groups 
  $H_{\bar\partial_E}^\bullet$ defined by the elliptic complex \eqref{eq:Elliptic_deformation_complex_holomorphic_curvature_equation}.} manifold $X$ of complex dimension $n$ and let $p\in(n,\infty)$. If $\bar\partial_E$ is a holomorphic structure on $E$ such that the holomorphic bundle $(E,\bar\partial_E)$ is simple and
\[
  \cS_{\bar\partial_E}
  := \bar\partial_E + \left\{\alpha \in W^{1,p}\left(\Lambda^{0,1}(\fsl(E)\right): F_{\bar\partial_E+\alpha} = 0
    \text{ and } \bar\partial_E^*\alpha = 0\right\},
\]
as in \cite[Equation (7.3.9), p. 233]{Kobayashi_differential_geometry_complex_vector_bundles}, then the quotient map,
\[
  \pi:\cU_{\bar\partial_E} \ni \bar\partial_E + \alpha \mapsto [\bar\partial_E + \alpha] \in \cM(E),
\]
gives a homeomorphism from an open neighborhood $\cU_{\bar\partial_E}$ of $\bar\partial_E$ in $\cS_{\bar\partial_E}$ onto an open neighborhood of $[\bar\partial_E] = \pi(\bar\partial_E) \in \cM(E)$.
\end{thm}

\begin{thm}[Moduli spaces of simple holomorphic vector bundles over complex manifolds]
\label{thm:Kobayashi_7_3_34_simple}  
(See Kobayashi \cite[Theorem 7.3.23, p. 236, Theorem 7.3.34, p. 239, and Corollary 7.3.32, p. 238]{Kobayashi_differential_geometry_complex_vector_bundles}.)
Continue the hypotheses of Theorem \ref{thm:Kobayashi_7_3_17}. Then the following hold:
\begin{enumerate}
\item\label{item:cM*(E)_is_complex_analytic_space}
The moduli subspace $\cM^*(E)$ in \eqref{eq:Moduli_space_simple_holomorphic_structures} of simple holomorphic vector bundles is a complex analytic space\footnote{See Grauert and Remmert \cite[Section 1.1.5, p. 7]{Grauert_Remmert_coherent_analytic_sheaves}.} (possibly non-Hausdorff and non-reduced\footnote{See Grauert and Remmert \cite[Section 1.1.5, p. 8]{Grauert_Remmert_coherent_analytic_sheaves}.}) with Zariski tangent spaces $\bH_{\bar\partial_E}^1$ at points $[\bar\partial_E] \in \cM^*(E)$ and smooth at points $[\bar\partial_E]$ with $\bH_{\bar\partial_E}^2 = 0$, where the harmonic spaces $\bH_{\bar\partial_E}^\bullet$ are defined by \eqref{eq:bHdbarA0bullet}.
  
\item\label{item:cM*(E)_is_complex_manifold}
The moduli subspace $\cM_\reg^*(E)$ in \eqref{eq:Moduli_space_regular_simple_holomorphic_structures} of regular, simple holomorphic vector bundles is a (possibly non-Hausdorff) complex manifold with tangent spaces $\bH_{\bar\partial_E}^1$ at points $[\bar\partial_E]$.
\end{enumerate}
\end{thm}

\begin{rmk}[On the hypotheses of Theorem \ref{thm:Kobayashi_7_3_34_simple}]
\label{rmk:Theorem_Kobayashi_7_3_34_simple_hypotheses}
Our statement of Theorem \ref{thm:Kobayashi_7_3_34_simple} differs slightly from that implied by \cite{Kobayashi_differential_geometry_complex_vector_bundles} since we restrict to holomorphic structures $\bar\partial_E$ that induce a \emph{fixed} holomorphic structure on $\det E$ as in \eqref{eq:Holomorphic_structure_fixed_determinant} and thus we replace $\GL(E)$ in \cite{Kobayashi_differential_geometry_complex_vector_bundles} by $\SL(E)$ in Theorem \ref{thm:Kobayashi_7_3_34_simple} and the forthcoming Theorem \ref{thm:Kobayashi_7_3_34_stable_hausdorff}. (This restriction has no significant impact on proofs of results in \cite{Kobayashi_differential_geometry_complex_vector_bundles}; see the discussion by Friedman and Morgan \cite[Section 4.1.6, p. 300]{FrM}.)

We note that Theorem \ref{thm:Kobayashi_7_3_17} does not imply that the moduli space $\cM^*(E)$ is Hausdorff, where both $\cM(E)$ and the subspace $\cM^*(E)\subset \cM(E)$ are equipped with the quotient topology: there could be points $[\bar\partial_E]$ and $[\bar\partial_{E'}]$ in $\cM^*(E)$ such that $[\bar\partial_{E'}]$ does not belong to the open neighborhood $\pi(\cU_{\bar\partial_E})$ in $\cM^*(E)$ provided by Theorem \ref{thm:Kobayashi_7_3_17}, but where $[\bar\partial_E]$ and $[\bar\partial_{E'}]$ cannot be separated by open neighborhoods in $\cM^*(E)$ as occurs in the example of the line with two origins in Lee \cite[Exercise 3.16, p. 83]{Lee_john_topological_manifolds}. According to Kobayashi \cite[Proposition 7.3.37, p. 239]{Kobayashi_differential_geometry_complex_vector_bundles}, if both of the points $[\bar\partial_E]$ and $[\bar\partial_{E'}]$ are semistable and at least one of the points is stable, then they can be separated by open neighborhoods in $\cM(E)$. See also the forthcoming Theorem \ref{thm:Kobayashi_7_4_20}.
\end{rmk}

\section{Moduli spaces of stable holomorphic vector bundles}
\label{sec:Moduli_space_stable_bundles}
As customary, for a holomorphic structure $\bar\partial_E$ on a complex vector bundle over a complex manifold, we let $\sE$ denote the holomorphic vector bundle $(E,\bar\partial_E)$ or its sheaf of holomorphic sections.
\label{page:SheafOfHolomorphicSections}

\begin{defn}[Stable and semistable holomorphic vector bundles over complex manifolds]
\label{defn:Stable_holomorphic_structure}
(See Donaldson and Kronheimer \cite[Section 6.1.4, p. 215]{DK}, Friedman and Morgan \cite[Section 4.3.2, Definition 3.3, p. 322]{FrM}, Kobayashi \cite[Section 5.7, p. 154]{Kobayashi_differential_geometry_complex_vector_bundles}, L\"ubke and Teleman \cite[Chapter 1, p. 1]{Lubke_Teleman_2006}, Mumford \cite{Mumford_1963}, \cite[Appendix to Chapter 5, Section C, Definition, p. 224]{Mumford_Fogarty_Kirwan_geometric_invariant_theory}, Takemoto \cite{Takemoto_1972, Takemoto_1973}, or Wells \cite[Appendix, Section 2.1, p. 247]{Wells3}.)
Let $E$ be a smooth complex vector bundle over a closed, complex K\"ahler manifold $(X,\omega)$ of complex dimension $n$. The \emph{normalized degree} or \emph{slope} of $E$ is the real number
\begin{equation}
  \label{eq:Slope}
  \mu(E) := \frac{\deg E}{\rank_\CC E},
\end{equation}
where the $\deg E$ is defined by \eqref{eq:Degree}. A holomorphic vector bundle $\sE = (E,\bar\partial_E)$ is \emph{semistable} if for every holomorphic subbundle $\sF \subset \sE$ with $\rank_\CC\sF > 0$ one has $\mu(\sF) \leq \mu(\sE)$ and $\sE$ is \emph{stable} if for every holomorphic bundle $\sF\subset\sE$ with $0 < \rank_\CC\sF < \rank_\CC\sE$ one has $\mu(\sF) < \mu(\sE)$. We call $\sE$ \emph{strictly semistable} if it is semistable but not stable. The holomorphic vector bundle $\sE$ is \emph{polystable} if
\begin{equation}
  \label{eq:Polystable_sum}
\sE = \sE_1 \oplus \cdots\oplus \sE_q
\end{equation}
as a direct sum of stable holomorphic subbundles with $\mu(\sE_k)=\mu(\sE)$ for $k=1,\ldots,q$, where $1\leq q\leq \rank_\CC\sE$. One says that $\sE$ is
\begin{inparaenum}[\itshape i\upshape)]
\item \emph{unstable} if it is not semistable,
\item \emph{strictly semistable} if it is semistable but not polystable, and
\item \emph{strictly polystable} if it is polystable but not stable.
\end{inparaenum}
\end{defn}

In Donaldson and Kronheimer \cite[Section 6.1.4, p. 215]{DK} and Friedman and Morgan \cite[Section 4.3.2, Definition 3.3, p. 322]{FrM}, the authors replace the condition that $\sF\subset\sE$ be a holomorphic vector subbundle by the condition that $\sF$ be a holomorphic bundle on $X$ that admits a non-zero map $\sF\to\sE$ of holomorphic vector bundles. We recall that if $\sE$ is stable, then $\sE$ is simple (see Friedman and Morgan \cite[Section 4.3.2, Proposition 3.4, p. 323]{FrM} or Kobayashi \cite[Corollary 5.7.14, p. 159]{Kobayashi_differential_geometry_complex_vector_bundles}). When $(X,\omega)$ is a Hodge manifold, then $\omega$ will be the curvature of an ample Hermitian line bundle $(L,h)$ \cite[Section 4.3.2, p. 322]{FrM}.

\begin{thm}[Moduli space of stable holomorphic vector bundles over a complex K\"ahler manifold]
\label{thm:Kobayashi_7_3_34_stable_hausdorff}  
(See Kobayashi \cite[Corollary 5.7.14, p. 159, Proposition 7.3.37, p. 239, and Remark 7.3.38, p. 240]{Kobayashi_differential_geometry_complex_vector_bundles}.)
Continue the hypotheses of Theorem \ref{thm:Kobayashi_7_3_34_simple} and assume further that $X$ is K\"ahler with K\"ahler form $\omega$. Then the moduli space of \emph{stable} holomorphic vector bundles with fixed holomorphic structure on the determinant line bundle,
\begin{equation}
  \label{eq:Moduli_space_omega-stable_holomorphic_structures}
  \cM(E,\omega) := \left\{[\bar\partial_E] \in \cM(E): (E,\bar\partial_E) \text{ is stable} \right\},
\end{equation}
is a Hausdorff topological space and a subset of $\cM^*(E)$.
\end{thm}

\begin{rmk}[Openness in Theorem \ref{thm:Kobayashi_7_3_34_stable_hausdorff}]
\label{rmk:Theorem_Kobayashi_7_3_34_stable_hausdorff_openness}  
Kobayashi notes in \cite[Remark 7.3.38, p. 240]{Kobayashi_differential_geometry_complex_vector_bundles} that there does not yet appear\footnote{We are grateful, however, to Robert Friedman for pointing out that this should still be possible.} to be a direct analytic proof of openness of $\cM(E,\omega)$ in $\cM(E)$, although there is an algebraic proof of openness due to Maruyama \cite{Maruyama_1976} for the moduli space of Gieseker-stable holomorphic bundles over an algebraic variety. L\"ubke and Teleman \cite[Section 4.4, pp. 112--113]{Lubke_Teleman_1995} remark that $\cM^*(E)$ need not be Hausdorff but is complex analytic.
\end{rmk}

Theorem \ref{thm:Kobayashi_7_3_34_stable_hausdorff} is reproved as part of the forthcoming Theorem \ref{thm:Kobayashi_7_4_20} with the aid of the Hitchin--Kobayashi Correspondence Theorem \ref{thm:Hitchin-Kobayashi_correspondence_Hermitian-Einstein_connections_stable_bundles}, which Kobayashi does not use to prove Theorem \ref{thm:Kobayashi_7_3_34_stable_hausdorff}.

\section[Hermitian--Einstein connections over almost Hermitian manifolds]{Projective Hermitian--Einstein connections on Hermitian vector bundles over almost Hermitian manifolds}
\label{sec:HE_connections}
We begin by recalling the well-known

\begin{defn}[Hermitian--Einstein connection on a Hermitian vector bundle over an almost Hermitian manifold]
\label{defn:HE_connection}
(See Donaldson and Kronheimer \cite[Section 6.1.4, p. 215]{DK}, Kobayashi \cite[Section 4.1]{Kobayashi_differential_geometry_complex_vector_bundles}.)
Let $(E,h)$ be a smooth Hermitian vector bundle over an almost Hermitian manifold $(X,g,J)$ with fundamental two-form $\omega = g(\cdot,J\cdot)$ as in \eqref{eq:Fundamental_two-form}. A smooth unitary connection $A$ on $E$ is \emph{Hermitian--Einstein} (or \emph{Hermitian--Yang--Mills}) if
\begin{align}
  \label{eq:Holomorphic_connection}
  F_A^{0,2} &= 0 \in \Omega^{0,2}(\fu(E)),
  \\
  \label{eq:Einstein_connection}
  \Lambda F_A &= -i\lambda\,\id_E \in \Omega^0(\fu(E)).
\end{align}  
for some constant $\lambda \in \RR$.
\end{defn}

Equation \eqref{eq:Einstein_connection} is called the \emph{Einstein condition} (see Kobayashi \cite[Section 4.1, p. 92]{Kobayashi_differential_geometry_complex_vector_bundles}) when $\lambda$ is a constant (as in Definition \ref{defn:HE_connection}) and the \emph{weak Einstein condition} when $\lambda\in C^\infty(X,\RR)$ is a real, smooth function on $X$ (see Kobayashi \cite[Equation (3.1.37), p. 56]{Kobayashi_differential_geometry_complex_vector_bundles}). By Kobayashi \cite[Proposition 4.2.4, p. 96]{Kobayashi_differential_geometry_complex_vector_bundles}
or L\"ubke and Teleman \cite[Lemma 2.1.5 (ii), p. 47]{Lubke_Teleman_1995}, the weak Einstein condition is equivalent to the Einstein condition modulo a conformal change (unique up to a homothety) of the Hermitian metric $h$ on $E$. The term $\lambda$ is called the \emph{Einstein factor}.

Because of \eqref{eq:Holomorphic_connection} and the fact that $A$ is unitary, we also have $F_A^{2,0}=0$. A unitary connection $A$ may be regarded as the Chern connection defined by an almost holomorphic structure $\bar\partial_E$ on $E$ and conversely $\bar\partial_A$ defines an almost holomorphic structure on $E$ (see Donaldson and Kronheimer \cite[Lemma 2.1.54, p. 45]{DK}, L\"ubke and Teleman \cite[Section 7.1]{Lubke_Teleman_1995}, \cite[Proposition 1.1.19, p. 27]{Lubke_Teleman_1995}). Even when the fundamental two-form $\omega \in \Omega^{1,1}(X,\RR)$ is not closed, it is convenient to still use the expression  \eqref{eq:DegreeIntegral} for $\deg E$ (compare De Bartolomeis and Tian \cite[Section 3, p. 244]{DeBartolomeis_Tian_1996}, who assume that $d\omega^{n-1}=0$)
\[
  \deg E = \int_X \frac{i}{2\pi}(\tr_E F_A)\wedge\omega^{n-1},
\]
keeping in mind that when $\omega$ is not closed, $\deg E$ is not independent of the choice of unitary connection induced by $A$ on the Hermitian line bundle $\det E$ or $\omega$, unlike the case where $\omega$ is closed, where $\deg E$ depends only on the cohomology classes $c_1(E)$ and $[\omega]$.
When $\alpha \in \Omega^{1,1}(X)$ and $X$ has complex dimension $n$, we recall that (see Kobayashi \cite[Equation (3.1.18), p. 52]{Kobayashi_differential_geometry_complex_vector_bundles})
\begin{equation}
  \label{eq:Kobayashi_3-1-18}
  \alpha\wedge\omega^{n-1} = \frac{1}{n}(\Lambda\alpha)\omega^n.
\end{equation}
By \eqref{eq:Holomorphic_connection}, we see that $F_A$ and thus also $\tr_E F_A$ have type $(1,1)$ and so by \eqref{eq:Kobayashi_3-1-18}, we have
\[
  (\tr_E F_A)\wedge\omega^{n-1} = \frac{1}{n}(\Lambda (\tr_E F_A))\omega^n.
\]
If $E$ has complex rank $r$, then (see Kobayashi \cite[Proposition 4.2.1, p. 96]{Kobayashi_differential_geometry_complex_vector_bundles} for the case when $(X,\omega)$ is K\"ahler) a solution $A$ to the \emph{Einstein condition} \eqref{eq:Einstein_connection} (see Kobayashi \cite[Equation (3.1.37), p. 56]{Kobayashi_differential_geometry_complex_vector_bundles}) has (constant) Einstein factor $\lambda$ determined by
\begin{align*}
  \deg E
  &=  \frac{i}{2\pi} \int_X (\tr_E F_A)\wedge\omega^{n-1} \quad\text{(by \eqref{eq:DegreeIntegral})}
  \\
  &= \frac{i}{2n\pi} \int_X \Lambda(\tr_E F_A)\omega^n \quad\text{(by \eqref{eq:Kobayashi_3-1-18})}
  \\
  &= \frac{i}{2n\pi} \int_X (-i\lambda r)\omega^n  \quad\text{(by \eqref{eq:Einstein_connection})}
  \\
  &= \frac{r}{2n\pi} \int_X \lambda\omega^n
  \\
  &= \frac{\lambda r(n-1)!}{2\pi} \int_X \frac{\omega^n}{n!}.
\end{align*}
By analogy with the expression for the volume of an almost K\"ahler manifold $(X,g,J)$ when $\omega$ is closed (see Huybrechts \cite[Exercise 3.1.13, p. 123]{Huybrechts_2005}), we define
\begin{equation}
  \label{eq:Pseudovolume_almost_Hermitian_manifold}
  \vol_\omega X := \int_X\frac{\omega^n}{n!}.
\end{equation}
(When $(X,g,J)$ is K\"ahler with corresponding symplectic form $\omega$, then the expression on the right-hand side of \eqref{eq:Pseudovolume_almost_Hermitian_manifold} is equal to the volume of $(X,g)$ defined by the Riemannian metric $g$ --- see Griffiths and Harris \cite[Section 0.2, p. 31, Wirtinger Theorem]{GriffithsHarris} or Huybrechts \cite[Exercise 3.1.8, p. 123]{Huybrechts_2005}.)
Hence, the preceding equation for $\deg E$ becomes
\[
  \deg E = \frac{\lambda r(n-1)!}{2\pi} \vol_\omega X
\]
and so the Einstein factor is given by
\begin{equation}
  \label{eq:Einstein_factor}
  \lambda = \frac{2\pi}{(n-1)!\vol_\omega X}\frac{\deg E}{\rank_\CC E},
\end{equation}
which agrees with the expressions in Kobayashi \cite[Equation (5.8.4), p. 163]{Kobayashi_differential_geometry_complex_vector_bundles}, L\"ubke and Teleman \cite[Lemma 2.1.8, p. 49]{Lubke_Teleman_1995}, \cite[Chapter 1, p. 1]{Lubke_Teleman_2006}, and, when $n=2$, with that in Donaldson and Kronheimer \cite[Section 6.1.4, p. 215]{DK}. For arbitrary $n$ and $r$, the single equation \eqref{eq:Einstein_connection} is equivalent to the pair of equations,
\begin{subequations}
  \label{eq:Einstein_connection_system}
  \begin{align}
    \label{eq:Einstein_connection_trace_part}
    \Lambda \tr F_A &= -ir\lambda \in \Omega^0(X,i\RR),
    \\
    \label{eq:Einstein_connection_tracefree_part}
    \Lambda (F_A)_0 &= 0 \in \Omega^0(\su(E)),
  \end{align}
\end{subequations}  
where the tracefree component of $F_A$ is defined by
\begin{equation}
  \label{eq:Curvature_tracefree_part}
  (F_A)_0 := F_A - \frac{1}{r}(\tr F_A)\id_E \in \Omega^2(\su(E)).
\end{equation}
By analogy with the decomposition of $\Omega^2(X,\CC)$ in \eqref{eq:Donaldson_Kronheimer_lemma_2-1-57} when $X$ has dimension two, one can define (see Kobayashi \cite[Section 7.2, p. 226]{Kobayashi_differential_geometry_complex_vector_bundles}) 
\begin{equation}
  \Omega^+(X,\CC) := \Omega^{2,0}(X) \oplus \Omega^0(X,\CC)\otimes\omega \oplus \Omega^{0,2}(X),
\end{equation}
and (see Kobayashi \cite[Equation (7.4.1), p. 240]{Kobayashi_differential_geometry_complex_vector_bundles})
\begin{equation}
  F_A ^+ := F_A^{2,0} + (\Lambda F_A)\frac{\omega}{2} + F_A^{0,2} \in \Omega^+(\fu(E)).
\end{equation}
We thus observe that the Hermitian--Einstein equations \eqref{eq:Holomorphic_connection}, \eqref{eq:Einstein_connection} are equivalent to
\begin{equation}
  \label{eq:FA_plus_HE_connection}
  F_A^+ = -\frac{i\lambda}{2}\omega\otimes\id_E \in \Omega^+(\fu(E)),
\end{equation}
which is invariant under the action of $\U(E)$. Suppose now that we fix a unitary connection $A_d$ on the Hermitian line bundle $\det E$ and that we restrict our attention to unitary connections $A$ on $E$ that induce the fixed connection $A_d$ on $\det E$, that is,
\begin{equation}
  \label{eq:Unitary_connection_detE_fixed}
  \nabla_{A_{\det E}} = \nabla_{A_d} \quad\text{on } \det E.
\end{equation}
With this restriction, we say that a connection $A$ is \emph{projectively Hermitian--Einstein} if in Definition \ref{defn:HE_connection} the equation \eqref{eq:Holomorphic_connection} is replaced by 
\begin{equation}
  \label{eq:Holomorphic_connection_fixed_unitary_determinant_connection}
  (F_A^{0,2})_0 = 0 \in \Omega^{0,2}(\su(E)).
\end{equation}
Observe that $A_d$ on $\det E$ automatically obeys the Einstein condition \eqref{eq:Einstein_connection} since $A_{\det E} = A_d$ by \eqref{eq:Unitary_connection_detE_fixed}, so
\[
  F_{A_d} = F_{A_{\det E}} = \tr_E F_A
\]
and thus, for $\lambda$ as in \eqref{eq:Einstein_factor},
\[
  \Lambda F_{A_d} = -ir\lambda \in \Omega^0(X,i\RR).
\]  
In the case of an almost Hermitian manifold $(X,g,J)$ of real dimension four, we see  by \eqref{eq:Donaldson_Kronheimer_lemma_2-1-57} that $A$ is a projectively Hermitian--Einstein connection on $E$ if and only if $A$ induces a $g$-anti-self-dual connection on $\su(E)$ and the induced connection on $\det E$ obeys \eqref{eq:Einstein_connection}.

If we assume now that the unitary connection $A_d$ on $\det E$ obeys
\begin{equation}
  \label{eq:Unitary_connection_detE_with_F02_zero}
  F_{A_d}^{0,2} = 0 \in \Omega^{0,2}(X,i\RR),
\end{equation}
then equations \eqref{eq:Holomorphic_connection} and \eqref{eq:Holomorphic_connection_fixed_unitary_determinant_connection} are equivalent since \eqref{eq:Unitary_connection_detE_with_F02_zero} yields
\begin{equation}
\label{eq:(0,2)ComponentOfCurvatureOfUnitaryConnectionTracefree}
  F_A^{0,2} = (F_A^{0,2})_0,
\end{equation}
and equation \eqref{eq:FA_plus_HE_connection} is invariant under the action of $\SU(E)$.

Focussing now on the case of a \emph{complex} Hermitian manifold $(X,g,J)$ and fixed unitary connection $A_d$ on $\det E$ with $F_{A_d}^{0,2} = 0$ as in \eqref{eq:Unitary_connection_detE_with_F02_zero}, the sequence of operators associated with a solution $A$ to the system \eqref{eq:Einstein_connection}, \eqref{eq:Holomorphic_connection_fixed_unitary_determinant_connection} modulo the action of $\SU(E)$ is given by (see Kobayashi \cite[Proposition 4.2.19, p. 226, and Lemma 4.2.20, p. 227]{Kobayashi_differential_geometry_complex_vector_bundles})
\begin{multline}
  \label{eq:HE_connection_fixed_unitary_determinant_connection_elliptic_deformation_complex}
  0 \xrightarrow{} \Omega^0(\su(E)) \xrightarrow{d_A} \Omega^1(\su(E)) \xrightarrow{\hat d_A^+}
  \Omega^0(\su(E))\oplus \Omega^{0,2}(\fsl(E))
  \\
  \xrightarrow{\bar\partial_A\circ\pi_{0,2}} \Omega^{0,3}(\fsl(E)) \xrightarrow{\bar\partial_A}
  \cdots
  \xrightarrow{\bar\partial_A} \Omega^{0,n}(\fsl(E)) \xrightarrow{} 0
\end{multline}
where $\pi_{0,2}$ is the projection from the direct sum $\Omega^0(\su(E))\oplus \Omega^{0,2}(\fsl(E))$ onto the subspace $\Omega^{0,2}(\fsl(E))$. For convenience of notation,  we have written $d_A$ in \eqref{eq:HE_connection_fixed_unitary_determinant_connection_elliptic_deformation_complex} for the covariant derivative induced on $\su(E)$ by the connection $A$ as given in \eqref{eq:Connection_Induced_On_Endomorphism_Bundle}. The harmonic representatives $\bH_A^\bullet$ for the cohomology groups $H_A^\bullet$ of \eqref{eq:HE_connection_fixed_unitary_determinant_connection_elliptic_deformation_complex} defined in Kobayashi \cite[Theorem 7.2.21, p. 227]{Kobayashi_differential_geometry_complex_vector_bundles} by analogy with those of the $g$-anti-self-dual equation when $X$ has real dimension four (see Donaldson and Kronheimer \cite[Equation (4.2.26), p. 138]{DK}). The differential $\hat d_A^+$ is the linearization of the Hermitian--Einstein equations \eqref{eq:Holomorphic_connection} (or equivalently \eqref{eq:Holomorphic_connection_fixed_unitary_determinant_connection} with the condition \eqref{eq:Unitary_connection_detE_with_F02_zero}) and \eqref{eq:Einstein_connection} (or equivalently \eqref{eq:Einstein_connection_tracefree_part} with the condition \eqref{eq:Einstein_connection_trace_part}), so that
\begin{equation}
  \label{eq:Linearization_HE_equation_codomain_Omega0+Omega02}
  \hat d_A^+a := \frac{1}{2}\bar\partial_A a'' + \Lambda(d_Aa) \in \Omega^{0,2}(\fsl(E)) \oplus \Omega^0(\su(E)),
  \quad\text{for all } a \in \Omega^1(\su(E)).
\end{equation}
We note that Kobayashi \cite[Section 7.2, p. 226]{Kobayashi_differential_geometry_complex_vector_bundles} prefers to use the following more traditional generalization of the operator $d_A^+$ in Atiyah, Hitchin and Singer \cite[p. 444]{AHS} or Donaldson and Kronheimer \cite[Equation (1.1.17) and p. 135]{DK} when $X$ has real dimension four:
\begin{multline}
  \label{eq:Linearization_HE_equation_full}
  d_A^+a := \frac{1}{2}\left(\partial_A a' + \bar\partial_A a''\right) + \Lambda(d_Aa)\otimes\omega
  \\
  \in \Omega^{2,0}(\fsl(E))\oplus \Omega^{0,2}(\fsl(E))\oplus\Omega^0(\su(E))\otimes\omega,
  \\
  \text{for all } a \in \Omega^1(\su(E)).
\end{multline}
The fact that the Hermitian--Einstein equations \eqref{eq:Holomorphic_connection}, \eqref{eq:Einstein_connection} are invariant under the action of $C^\infty(\SU(E))$ ensures that
\[
  \hat d_A^+\circ d_A = 0 \quad\text{on } \Omega^0(\su(E))
\]
while the fact that $F_A^{0,2}=0$ ensures that
\begin{align*}
  \bar\partial_A\circ\pi_{0,2}\circ \hat d_A^+  &= 0 \quad\text{on } \Omega^1(\su(E)),
  \\
  \bar\partial_A\circ \bar\partial_A &= 0 \quad\text{on } \Omega^{0,k}(\fsl(E)),
                                       \quad\text{for } k = 2, \ldots, n-2.
\end{align*}
The preceding identities ensure that the sequence \eqref{eq:HE_connection_fixed_unitary_determinant_connection_elliptic_deformation_complex} is a complex and by Kobayashi \cite[Lemma 7.2.20, p. 227]{Kobayashi_differential_geometry_complex_vector_bundles}, this complex is elliptic since the associated symbol sequence is exact.

We define the real vector spaces $\bH_A^\bullet$ to be the harmonic representatives for the cohomology groups $H_A^\bullet$ of the elliptic complex \eqref{eq:HE_connection_fixed_unitary_determinant_connection_elliptic_deformation_complex}:
\begin{subequations}
  \label{eq:HE_equation_bHAbullet}
  \begin{align}
    \label{eq:HE_equation_bHA0}
    \bH_A^0 &:= \Ker d_A\cap \Omega^0(\su(E)),
    \\
    \label{eq:HE_equation_bHA1}
    \bH_A^1 &:= \Ker (\hat d_A^+ + d_A^*)\cap \Omega^1(\su(E)),
    \\
    \label{eq:HE_equation_bHA2}
    \bH_A^2 &:= \Ker (\bar\partial_A\circ\pi_{0,2} + \hat d_A^{+,*})\cap
              \left(\Omega^0(\su(E))\oplus\Omega^{0,2}(\fsl(E))\right),
    \\
    \label{eq:HE_equation_bHAk}
    \bH_A^k &:= \Ker (\bar\partial_A + \bar\partial_A^*)\cap \Omega^{0,k}(\fsl(E)),
              \quad\text{for } k = 3,\ldots,n,
  \end{align}    
\end{subequations}
with the convention that $\Omega^{0,k}(\fsl(E)) = (0)$ if $k>n$.

We define the complex vector spaces $\bH_{\bar\partial_A}^\bullet$ to be the harmonic representatives for the cohomology groups $H_{\bar\partial_A}^\bullet$ of the elliptic complex (see Wells \cite[Chapter IV, Example 5.7, p. 151]{Wells3})
\begin{equation}
  \label{eq:Complex_vector_bundle_dbar_A_elliptic_deformation_fixed_det_holomorphic_structure}
  0 \xrightarrow{} \Omega^0(\fsl(E)) \xrightarrow{\bar\partial_A} \Omega^{0,1}(\fsl(E)) \xrightarrow{} 
  \cdots \xrightarrow{\bar\partial_A} \Omega^{0,n-1}(\fsl(E))
  \xrightarrow{\bar\partial_A} \Omega^{0,n}(\fsl(E)) \xrightarrow{} 0,
\end{equation}
where $\bar\partial_A$ on $E$ induces the fixed holomorphic structure $\bar\partial_{A_{\det}} = \bar\partial_{A_d}$ on $\det E$:
\begin{equation}
  \label{eq:bHdbarA0bullet}
    \bH_{\bar\partial_A}^k := \Ker (\bar\partial_A + \bar\partial_A^*)\cap \Omega^{0,k}(\fsl(E)),
              \quad\text{for } k = 0,\ldots,n,
\end{equation}
with the convention that $\Omega^{0,k}(\fsl(E)) = (0)$ if $k<0$ or $k>n$. For convenience of notation,  we have written $\bar\partial_A$ in \eqref{eq:Complex_vector_bundle_dbar_A_elliptic_deformation_fixed_det_holomorphic_structure} for the covariant derivative induced on $\fsl(E)$ by the $(0,1)$-connection $\bar\partial_A$ as given in \eqref{eq:(0,1)Connection_Induced_On_Endomorphism_Bundle}.

Unlike Kobayashi in \cite[Theorem 7.2.21, p. 227]{Kobayashi_differential_geometry_complex_vector_bundles}, we restrict our attention to the case of fixed unitary connection or holomorphic structure on $\det E$ (see Remark \ref{rmk:Theorem_Kobayashi_7_3_34_simple_hypotheses}). To examine the question of existence of unitary connections $A_d$ on $\det E$ with $F_{A_d}^{0,2} = 0$, we recall the

\begin{prop}[Classification of complex line bundles over a closed manifold]
Let $X$ be a closed, topological manifold.
\begin{enumerate}
\item The group of complex line bundles modulo continuous equivalence is isomorphic to the group $H^1(X,\sC_X^*) \cong H^2(X,\ZZ)$, where $\sC_X^*$ is the sheaf of nowhere-vanishing continuous functions on $X$, that is, continuous functions taking values in the multiplicative group $\CC^* = \CC\less\{0\}$.  
 
\item If $X$ is a complex manifold, then the group of holomorphic line bundles modulo holomorphic equivalence is isomorphic to the group $\Pic(X) = H^1(X,\sO_X^*)$, where $\sO_X^*$ is the sheaf of nowhere-vanishing holomorphic functions on $X$.
\end{enumerate}
\end{prop}

Suppose that $L$ is a holomorphic line bundle over a closed, complex manifold $X$. We recall that the exact sequence of sheaves (see Barth, Hulek, Peters Van de Ven \cite[Section 1.6]{Barth_Hulek_Peters_Van_de_Ven_compact_complex_surfaces}, Ciliberto \cite[Section 1.2]{Ciliberto_classification_complex_algebraic_surfaces} or Griffiths and Harris \cite[p. 139]{GriffithsHarris})
\[
  0 \xrightarrow{\iota} \ZZ_X \xrightarrow{} \sO_X \xrightarrow{\exp} \sO_X^* \xrightarrow{} 0
\]
where $\iota$ is inclusion and $f\mapsto \exp(2\pi f)$ for $f \in \sO_X$, induces a long exact sequence in cohomology and a boundary map $\delta$,
\[
  \cdots \xrightarrow{} H^1(X,\ZZ) \xrightarrow{} H^1(X,\sO_X) \xrightarrow{} H^1(X,\sO_X^*) \xrightarrow{\delta} H^2(X,\ZZ)\xrightarrow{} \cdots
\]
and $c_1(L) = \delta(L) \in H^2(X,\ZZ)$ is the first Chern class of $L$ in the Picard group $\Pic(X) = H^1(X,\sO_X^*)$ (see \cite[Proposition 1.6.1]{Barth_Hulek_Peters_Van_de_Ven_compact_complex_surfaces} or Hirzebruch \cite[Theorem 4.3.1]{Hirzebruch_topological_methods_algebraic_geometry}). If $L$ has a Hermitian metric, then there is a unique unitary connection $A_L$ on $L$ such that $\bar\partial_{A_L} = \bar\partial_L$ (see Kobayashi \cite[Proposition 1.4.9, p. 11]{Kobayashi_differential_geometry_complex_vector_bundles}) and thus $F_{A_L}^{0,2} = \bar\partial_{A_L}\circ \bar\partial_{A_L} = 0$, so $F_{A_L}$ has type $(1,1)$ and \cite[p. 39]{DK}, \cite[p. 141]{GriffithsHarris}
\[
  c_1(L) = \left[\frac{i}{2\pi}F_{A_L}\right] \in H^2(X,\ZZ).
\]
Recall that on a closed, K\"ahler manifold (see Wells \cite[Chapter V, Theorem 5.1, p. 197]{Wells3})
\[
  H^2(X,\CC) = H^{2,0}(X) \oplus H^{1,1}(X) \oplus H^{0,2}(X) 
\]
and $\overline{H^{2,0}(X)} = H^{0,2}(X)$, where for any non-negative integers $p,q$,
\[
  H^{p,q}(X) \cong \Ker\square \cap \Omega^{p,q}(X),
\]
and $\square = \partial^*\partial + \partial\partial^*$ (see Wells \cite[Chapter V, Section 4, pp. 191 and 195]{Wells3}).
The group $\Pic(X)$ fits into an exact sequence
\[
  0 \xrightarrow{} T \xrightarrow{} \Pic(X) \xrightarrow{\delta} H^2(X,\ZZ)\cap H^{1,1}(X) \xrightarrow{} 0
\]  
where $T = \Pic^0(X) = \Ker\delta$ is a complex torus of dimension $h^{1,0}$ and
\[
  \Pic(X)/\Pic^0(X) = H^{1,1}(X)\cap \iota^*H^2(X,\ZZ)
\]
is the \emph{N\'eron--Severi subgroup} \cite[pp. 27 and 143]{Barth_Hulek_Peters_Van_de_Ven_compact_complex_surfaces}, where $\iota^*:H^2(X,\ZZ)\to H^2(X,\CC)$ is induced by the embedding $\iota:\ZZ_X\to\CC_X$ of sheaves.
By the Lefshetz theorem on $(1,1)$ classes \cite[Theorem IV.2.13 143]{Barth_Hulek_Peters_Van_de_Ven_compact_complex_surfaces}, the range of $\delta:\Pic(X)\to H^2(X,\CC)$ is given by the set of cohomology classes which are representable by closed, real differential two-forms of type $(1,1)$. If $X$ is a compact, complex algebraic surface and $\widetilde X$ is the blow-up of $X$ at a point, then $\Pic(\widetilde X) = \Pic(X)\oplus\ZZ[e]$, where $e \subset \widetilde X$ is the exceptional divisor (see \cite[Section 1.7]{Ciliberto_classification_complex_algebraic_surfaces} or Hartshorne \cite[Proposition 5.3.2]{Hartshorne_algebraic_geometry}).

Recall that if $D \subset X$ is a divisor, then (see Huybrechts \cite[Proposition 4.4.13, p. 202]{Huybrechts_2005})
\[
  c_1(\sO_X(D)) = \PD[D] \in H^2(X,\ZZ),
\]
where $\PD[D]$ denotes the Poincar\'e dual of the fundamental class of $D$. Writing $\PD[e] = e^*$, we have
\[
  H^2(\widetilde X,\ZZ) = H^2(X,\ZZ)\oplus \ZZ e^*.
\]
Hence, if $c_1(L) = ke^*$ for $k\in\ZZ$, then $c_1(L)$ is represented by an integer multiple of $\omega_{\CC\PP^2}$, the K\"ahler form on $\CC\PP^2$. (This can also be deduced from Griffiths and Harris \cite[pp. 185--187]{GriffithsHarris}.) More generally, if $X$ is closed, K\"ahler manifold and $D\subset X$ is a complex analytic subvariety of complex codimension $1$, then its Poincar\'e dual $[D] \in H^2(X,\CC)$ has type $(1,1)$ by Griffiths and Harris \cite[pp. 163--164]{GriffithsHarris}. In particular, if $X$ is a complex projective manifold then $e \subset \widetilde X$ is complex analytic subvariety of complex codimension $1$ and so $e^* \in H^{1,1}(X)$.

\section{Moduli spaces of projectively Hermitian--Einstein connections}
\label{sec:Moduli_space_HE_connections}
Given a smooth Hermitian vector bundle $(E,h)$ over a smooth manifold $X$ of dimension $d$ and a constant $p\in (d/2,\infty)$, we recall from Section \ref{sec:Gauge_transformations_stabilizers_unitary_connections} that $\sA(E,h)$ is the affine space of $W^{1,p}$ unitary connections $A$ on $E$ that induce a fixed smooth, unitary connection $A_d$ on the Hermitian line bundle $\det E$, as in \eqref{eq:Unitary_connection_detE_fixed} and so
\[
  A_{\det} = A_d \text{ on } \det E.
\]
When $(X,g,J)$ is an almost Hermitian manifold with fundamental two-form $\omega = g(\cdot,J\cdot)$ as in \eqref{eq:Fundamental_two-form}, we let
\begin{equation}
  \label{eq:Moduli_space_HE_connections}
  M(E,h,\omega) := \left\{A \in \sA(E,h): \eqref{eq:Einstein_connection}, \eqref{eq:Holomorphic_connection_fixed_unitary_determinant_connection} \text{ hold} \right\}/W^{2,p}(\SU(E))
\end{equation}
denote the moduli space of projectively Hermitian--Einstein connections modulo the Banach Lie group $W^{2,p}(\SU(E))$ of determinant-one, unitary automorphisms of $E$, where the action is given in \eqref{eq:W2pSUE_action_on_AEh}. We also define the following open subspaces,
\begin{equation}
  \label{eq:Moduli_space_HE_connections_irreducible}  
  M^*(E,h,\omega) := \left\{[A] \in M(E,h,\omega): A \text{ is non-split} \right\},
\end{equation}
and, when $(X,g,J)$ is complex Hermitian,
\begin{subequations}
  \label{eq:Moduli_space_HE_connections_open_subspaces}
\begin{align}
  \label{eq:Moduli_space_HE_connections_regular}
  M_\reg(E,h,\omega) &:= \left\{[A] \in M(E,h,\omega): \bH_A^2 = (0) \right\},
  \\
  \label{eq:Moduli_space_HE_connections_irreducible_regular}
  M_\reg^*(E,h,\omega) &:= M^*(E,h,\omega) \cap M_\reg(E,h,\omega),
\end{align}
\end{subequations}
where $\bH_A^2$ is defined by the elliptic deformation complex \eqref{eq:HE_connection_fixed_unitary_determinant_connection_elliptic_deformation_complex}. The property that $M_\reg(E,h,\omega)$ is an open subspace of $M(E,h,\omega)$ follows from the Implicit Mapping Theorem, while the fact that $M^*(E,h,\omega)$ is an open subspace of $M(E,h,\omega)$ follows from the forthcoming

\begin{lem}[Openness of the subspace of unitary connections with minimal stabilizer and moduli subspace of non-split projectively Hermitian--Einstein connections]
\label{lem:Openness_moduli_subspace_non-split_projectively_Hermitian-Einstein_connections}
Let $(E,h)$ be a smooth Hermitian vector bundle of rank $r$ over a connected, smooth Riemannian manifold $(X,g)$ of dimension $d$, and $p\in (d/2,\infty)$ be a constant, and $A_d$ be a fixed smooth, unitary connection on the Hermitian line bundle $\det E$. Then the following hold:
\begin{enumerate}
\item\label{item:Openness_subspace_unitary_connections_minimal_stabilizer}
The subspace $\sA^{C_r}(E,h)$ of connections in $\sA(E,h)$ with minimal stabilizer, $C_r$, is open in $\sA(E,h)$, where $C_r$ is the group of $r$-th roots of unity.

\item\label{item:Openness_quotient_subspace_unitary_connections_minimal_stabilizer}
The quotient subspace $\sB^{C_r}(E,h)$ of gauge-equivalence classes of connections in $\sB(E,h) := \sA(E,h)/W^{2,p}(\SU(E))$ with minimal stabilizer, $C_r$, is open in $\sB(E,h)$ equipped with the quotient topology.

\item\label{item:Openness_moduli_subspace_non-split_projectively_Hermitian-Einstein_connections}
The moduli subspace $M^*(E,h,\omega)$ in \eqref{eq:Moduli_space_HE_connections_irreducible}  of gauge-equivalence classes of non-split projectively Hermitian--Einstein connections is open in $M(E,h,\omega)$.
\end{enumerate}
\end{lem}

\begin{rmk}[Openness of the subspace of unitary connections with minimal stabilizer]
\label{rmk:Openness_subspace_unitary_connections_minimal_stabilizer}
Item \eqref{item:Openness_subspace_unitary_connections_minimal_stabilizer} in Lemma \ref{lem:Openness_moduli_subspace_non-split_projectively_Hermitian-Einstein_connections} is often stated without proof in the literature --- see, for example, Donaldson and Kronheimer \cite[Section 4.2.2, p. 132, line $-15$]{DK} or Singer \cite{Singer_1978, Singer_1981}. While it would appear to immediately follow from a more general result described by Rudolph and Schmidt \cite[Corollary 8.3.6, p. 646 and Remark 8.2.9 (2), p. 642]{Rudolph_Schmidt_differential_geometry_mathematical_physics_part2}, their proof seems to us to be incomplete since those authors do not allow for the possibility of gauge transformations $u \in \Stab(A)$ that may be far from the identity. Lawson \cite[Chapter II, Theorem 10.4 (2), p. 33, and pp. 36--37]{Lawson} proves Item \eqref{item:Openness_subspace_unitary_connections_minimal_stabilizer} in the special case where $E$ has rank two.
Kondracki and Rogulski prove more general results of this kind as \cite[Proposition 4.3.1, p. 50, Theorem 4.3.5, p. 55, and Proposition 4.3.7, p. 56]{Kondracki_Rogulski_1986}: see Remark \ref{rmk:Kondracki_Rogulski_openness_results}. As we shall see below, Item \eqref{item:Openness_subspace_unitary_connections_minimal_stabilizer} in Lemma \ref{lem:Openness_moduli_subspace_non-split_projectively_Hermitian-Einstein_connections} may be proved several different ways and we include them since these proofs are of independent interest.
\end{rmk}

\begin{proof}[First proof of Lemma \ref{lem:Openness_moduli_subspace_non-split_projectively_Hermitian-Einstein_connections}
\eqref{item:Openness_subspace_unitary_connections_minimal_stabilizer}]
We can write $\sA(E,h)$ as a disjoint union:
\[
  \sA(E,h) =  \sA^{C_r}(E,h) \sqcup \sA^{>C_r}(E,h),
\]
where $\sA^{>C_r}(E,h) := \{A \in \sA(E,h): \Stab(A) \supsetneq C_r\,\id_E\}$. We claim that $\sA^{>C_r}(E,h)$ is a closed subspace of $\sA(E,h)$. Let $\{A_k\}_{k\in\NN} \subset \sA^{>C_r}(E,h)$ be a sequence of connections that converge in $W^{1,p}(X)$ to a limit $A_\infty\in \sA(E,h)$ as $k\to\infty$. Since $\Stab(A_k) \supsetneq C_r\,\id_E$, there exist elements of $\Stab(A_k)$ that do not belong to $C_r\,\id_E$, for each $k\in\NN$. We now argue that these elements of $\Stab(A_k)\less(C_r\,\id_E)$ can be chosen so that their distance from the center $C_r\,\id_E$ is uniformly bounded from below by a positive constant. By choosing a point $x\in X$ and an $h$-orthonormal basis for the fiber $E_x$, we can identify $\Stab(A_k)$ with the centralizer in $\SU(E_x) \cong \SU(r)$ of the holonomy group of $A_k$ at $x$ (see, for example, Rudolph and Schmidt \cite[Theorem 6.1.5, p. 467]{Rudolph_Schmidt_differential_geometry_mathematical_physics_part2}). For any $u,v \in \SU(r)$, we define
\begin{equation}
\label{eq:Defining_Orbit_Distance}
  \dist_{\SU(r)}(u,v) := \inf_{w\in\SU(r)}\|u-wvw^{-1}\|,
\end{equation}
where $\|\cdot\|$ denotes the operator norm on $\End(E_x) \cong \End(\CC^r)$ induced by the Hermitian inner product $h$ on $E_x \cong \CC^r$. For any closed subset $S \subset \SU(r)$, we define
\[
  \dist_{\SU(r)}(u,S) := \inf_{s\in S}\dist_{\SU(r)}(u,s).
\]
Since $\Stab(A_k) \supsetneq C_r\,\id_E$, there is an element $u_k\in \Stab(A_k)\less (C_r\,\id_E)$ for each $k \in \NN$. Because the orbit of $u_k$ in $\SU(r)$,
\[
  \left\{ v u_k v^{-1}: v \in \SU(r) \right\},
\]
is closed (since it is the image of the compact group $\SU(r)$ under a continuous map to a Hausdorff topological space) and disjoint from the compact set $C_r\,\id_E$, there is a constant $\delta_k = \delta_k(g,h) \in (0,1]$ for each $k \in \NN$ such that
\[
  \dist_{\SU(r)}(u_k,C_r)\ge \delta_k,
\]
where $\|\cdot\|$ denotes the operator norm on $\End(\CC^r)$ induced by the standard Hermitian inner product on $\CC^r$. Now observe that if $\Stab(A_k)$ and $\Stab(A_\ell)$ are conjugate in $\SU(r)$, we can choose $u_k$ and $u_\ell$ to be conjugate in $\SU(r)$.  The definition of $\dist_{\SU(r)}$ in \eqref{eq:Defining_Orbit_Distance} implies that
\[
  \dist_{\SU(r)}(u_k,C_r)=\dist_{\SU(r)}(u_\ell,C_r),
\]
so we can choose $\delta_k=\delta_\ell$. By Lemma \ref{lem:Howe_Subgroups_of_SU(n)}, there are finitely many conjugacy classes of Howe subgroups of $\SU(r)$ and, as in the proof of Lemma \ref{lem:Kronheimer_Lemma2.2}, every stabilizer $\Stab(A_k)$ is conjugate to a Howe subgroup of $\SU(r)$.
Thus, we can choose a constant $\delta = \delta(g,h) \in (0,1]$ and, for each $k \in \NN$, an element $u_k\in\Stab(A_k)$ such that
\begin{equation}
  \label{eq:Dist_uk_Cr_geq_delta}
  \dist_{\SU(r)}(u_k,C_r) \geq \delta, \quad\text{for all } k \in \NN.
\end{equation}
Because $\SU(r)$ is compact, we can select a subsequence $\{k'\}\subset \{k\}$ such that, after relabeling, the sequence $\{u_k\}_{k\in\NN} \subset \SU(r)$ converges to a limit $u_\infty \in \SU(r)$. Hence, when viewed as elements of $W^{2,p}(\SU(E))$, the sequence $\{u_k\}_{k\in\NN} \subset W^{2,p}(\SU(E))$ converges to a limit $u_\infty \in W^{2,p}(\SU(E))$.

Since the action of $W^{2,p}(\SU(E))$ on $\sA(E,h)$ is continuous (in fact, smooth as explained following equation \eqref{eq:W2pSUE_action_on_AEh}), we have
\[
  A_\infty
  = \lim_{k\to\infty}A_k
  = \lim_{k\to\infty}u_k\cdot A_k
  = \left(\lim_{k\to\infty}u_k\right)\cdot\left(\lim_{k\to\infty}A_k\right)
  = u_\infty\cdot A_\infty,
\]
and thus $u_\infty \in \Stab(A_\infty)$. By taking the limit as $k\to\infty$, the inequality \eqref{eq:Dist_uk_Cr_geq_delta} yields
\[
  \dist_{\SU(r)}(u_\infty,C_r) \geq \delta,
\]
and so $u_\infty \notin C_r$. Therefore, $A_\infty \in \sA^{>C_r}(E,h)$ and thus $\sA^{>C_r}(E,h)$ is a closed subspace of $\sA(E,h)$, as claimed. Consequently, $\sA^{C_r}(E,h)$ is an open subspace of $\sA(E,h)$.
\end{proof}

\begin{rmk}[Generalizations of Lemma \ref{lem:Openness_moduli_subspace_non-split_projectively_Hermitian-Einstein_connections} \eqref{item:Openness_subspace_unitary_connections_minimal_stabilizer}]
\label{rmk:Kondracki_Rogulski_openness_results}
Continue the notation of the first proof of Lemma \ref{lem:Openness_moduli_subspace_non-split_projectively_Hermitian-Einstein_connections} \eqref{item:Openness_subspace_unitary_connections_minimal_stabilizer}.
For any Lie subgroup $S \subset \sG_E$, Kondracki and Rogulski define \cite[Section 4.3, p. 49]{Kondracki_Rogulski_1986}
\begin{align*}
  \sA^S(E,h) &:= \left\{A \in \sA(E,h): \Stab(A) = S\right\},
  \\
  \sA^{\geq S}(E,h) &:= \left\{A \in \sA(E,h): S \subset \Stab(A)\right\}.
\end{align*}
(Here, we reverse the notation of Rudolph and Schmidt \cite[Remark 8.2.9 (2), p. 642]{Rudolph_Schmidt_differential_geometry_mathematical_physics_part2}.) Kondracki and Rogulski prove in \cite[Proposition 4.3.1, p. 50 and Theorem 4.3.5, p. 55]{Kondracki_Rogulski_1986} that $\sA^S(E,h)$ is open in $\sA^{\geq S}(E,h)$. They also prove in \cite[Proposition 4.3.7, p. 56]{Kondracki_Rogulski_1986} that $\sA^{(S)}(E,h)$ is open in $\sG_E\cdot\sA^{\geq S}(E,h)$, where $\sA^{(S)}(E,h)$ is the subspace of all connections in $\sA(E,h)$ whose stabilizers are conjugate to $S$.
\end{rmk}  

Before presenting several alternative proofs of Lemma \ref{lem:Openness_moduli_subspace_non-split_projectively_Hermitian-Einstein_connections} \eqref{item:Openness_subspace_unitary_connections_minimal_stabilizer}, we dispose of the

\begin{proof}[Proofs of Lemma \ref{lem:Openness_moduli_subspace_non-split_projectively_Hermitian-Einstein_connections} \eqref{item:Openness_quotient_subspace_unitary_connections_minimal_stabilizer} and \eqref{item:Openness_moduli_subspace_non-split_projectively_Hermitian-Einstein_connections}]
Item \eqref{item:Openness_quotient_subspace_unitary_connections_minimal_stabilizer} follows immediately from the definition of the quotient topology, see Munkres \cite[Section 22, pp. 137-138]{Munkres_topology_second_edition},  on $\sB(E,h)$ (a set $U\subset \sB(E,h)$ is open if and only if $\pi^{-1}(U)\subset\sA(E,h)$ is open where $\pi:\sA(E,h) \to \sB(E,h)$ is the quotient map), the fact that $\pi^{-1}(\sB^{C_r}(E,h))=\sA^{C_r}(E,h)$, and Item \eqref{item:Openness_subspace_unitary_connections_minimal_stabilizer}.

Consider Item \eqref{item:Openness_moduli_subspace_non-split_projectively_Hermitian-Einstein_connections}. Define
\[
  M^{C_r}(E,h,\omega) := M(E,h,\omega)\cap \sB^{C_r}(E,h)
\]
and observe that $M^{C_r}(E,h,\omega)$ is open in $M(E,h,\omega)$ by Item \eqref{item:Openness_quotient_subspace_unitary_connections_minimal_stabilizer}. For each point $[A] \in M(E,h,\omega)$, the $W^{1,p}$ connection $A$ is gauge-equivalent to a smooth connection, so we may assume without loss of generality that $A$ is smooth (see the forthcoming Remark \ref{rmk:Hitchin-Kobayashi_correspondence_Hermitian-Einstein_connections_regularity}). By Corollary \ref{cor:Split_unitary_A_and_Lie_Alg_of_Stab(A)} \eqref{item:A_smooth_and_bHA_non-zero_implies_A_split} and \eqref{item:A_split_implies_bHA_non-zero}, the connection $A$ is non-split if and only if $\bH_A^0 = (0)$. By Corollary \ref{cor:Split_unitary_A_and_Lie_Alg_of_Stab(A)} \eqref{item:HA0_is_Lie_Alg_of_Stab(A)_Irreducible}, the connection $A$ has minimal stabilizer $C_r$ if and only if $\bH_A^0 = (0)$. Hence, we obtain the equality
\[
  M^*(E,h,\omega) = M^{C_r}(E,h,\omega),
\]
and this yields the conclusion.
\end{proof}

We now present several other proofs of Lemma \ref{lem:Openness_moduli_subspace_non-split_projectively_Hermitian-Einstein_connections} \eqref{item:Openness_subspace_unitary_connections_minimal_stabilizer}.

\begin{rmk}[Principal Orbit Type Theorem]
\label{rmk:PrincipalOrbitTypeTheorem}
The following proof is essentially that of the Principal Orbit Type Theorem for the proper action of a Lie group $G$ on a smooth manifold $M$ --- see Meinrenken \cite[Theorem 1.32, p. 15]{Meinrenken_group_actions_manifolds_lecture_notes} or tom Dieck \cite[Theorem 5.14, p. 42]{tomDieck}. The majority of the work required to prove that theorem lies in establishing the existence of a $G$-equivariant tubular neighborhood of the orbit of a point in $M$, as in the argument provided below.
\end{rmk}

\begin{proof}[Second proof of Lemma \ref{lem:Openness_moduli_subspace_non-split_projectively_Hermitian-Einstein_connections} \eqref{item:Openness_subspace_unitary_connections_minimal_stabilizer}]
We shall apply the method of proof of Kondracki and Rogulski \cite[Corollary 3.3.5 (ii), p. 41 and Proposition 4.3.1, p. 50]{Kondracki_Rogulski_1986}.
We shall denote $\sG_E := W^{2,p}(\SU(E))$ and $\tilde\sG_E := W^{2,p}(\SU(E))/(C_r\,\id_E)$ for convenience. For any $A \in \sA(E,h)$, we have $\sG_E\cdot A = \tilde\sG_E\cdot A$. We claim that if $A \in \sA^{C_r}(E,h)$, then the orbit map
\begin{equation}
  \label{eq:Orbit_map_A}
  \iota_A:\tilde\sG_E \ni u \mapsto u\cdot A \in \tilde\sG_E\cdot A \subset \sA(E,h)
\end{equation}
is a smooth embedding and $\tilde\sG_E\cdot A$ is an embedded smooth submanifold of $\sA(E,h)$, where we write $u\cdot A = u^*A$ for convenience for each $u\in \sG_E$. We first show that $\iota_A$ in \eqref{eq:Orbit_map_A} is injective. If $u\cdot A = v\cdot A$ for $u,v\in \sG_E$, then
\[
  A = (uu^{-1})\cdot A = (uu^{-1})^*A = (u^{-1})^*u^*A = (u^{-1})^*v^*A = (vu^{-1})^*A = (vu^{-1})\cdot A,
\]
so $vu^{-1} \in \Stab(A) = C_r\,\id_E$, and thus $u = v \in \tilde\sG_E$, so $\iota_A$ is injective. The map $\iota_A$ in \eqref{eq:Orbit_map_A} is clearly smooth. We now show that $\iota_A$ is an immersion (alternatively, see  Kondracki and Rogulski \cite[Lemma 3.2.5, p. 37]{Kondracki_Rogulski_1986}). The differential of $\iota_A$ at $\id_E \in \sG_E$ is given by (see the discussion around \eqref{eq:d0_projective_vortex_elliptic_deformation_complex})
\[
  (D\iota_A)(\id_E):T_{\id_E}\sG_E = W^{2,p}(\su(E)) \ni \xi \mapsto d_A\xi \in W^{1,p}(T^*X\otimes\su(E)) = T_A\sA(E,h).
\]  
Because $\Stab(A) = C_r\,\id_E$, Lemma \ref{lem:LieGroup_Structure_of_Stab(A)} yields
\[
  \Ker d_A\cap W^{2,p}(\su(E)) = T_{\id_E}\Stab(A) = (0)
\]
and so $\Ker(D\iota_A)(\id_E) = (0)$. For arbitrary $u\in\sG_E$, we note that $\iota_A(u) = \iota_{u\cdot A}(\id_E)$. Thus,
\[
  (D\iota_A)(u) = (D\iota_{u\cdot A})(\id_E) = d_{u\cdot A}.
\]
If $s \in \Stab(A)$, then
\[
  (u^{-1}su)\cdot (u\cdot A) = (u^{-1}su)^*u^*A = (uu^{-1}su)^*A = u^*s^*A = u^*A = u\cdot A,
\]
and so $u^{-1}su \in \Stab(u\cdot A)$. Therefore, $\Stab(u\cdot A) = u^{-1}(\Stab(A))u$ and hence $\Stab(u\cdot A) = C_r\,\id_E$ and $\Ker d_{u\cdot A}\cap W^{2,p}(\su(E)) = T_{\id_E}\Stab(u\cdot A) = (0)$. Consequently, $\Ker(D\iota_A)(u) = (0)$ for all $u\in\sG_E$ and so the map $\iota_A$ in \eqref{eq:Orbit_map_A} is an immersion.
Moreover, the map $\iota_A$ in \eqref{eq:Orbit_map_A} is proper since $\SU(r)$ is compact (see Freed and Uhlenbeck \cite[Proof of Corollary, pp. 50--51]{FU}, Kondracki and Rogulski \cite[Corollary 2.4.11, p. 30]{Kondracki_Rogulski_1986}, or Lawson \cite[Proof of Theorem 10.4 (1), p. 36]{Lawson}). Consequently, the map $\iota_A$ is closed (see Abraham, Marsden, and Ratiu \cite[Section 3.5, p. 201]{AMR} or Lee \cite[Theorem A.57]{Lee_john_smooth_manifolds}) and so is an embedding by Abraham, Marsden, and Ratiu \cite[Theorem 3.5.9, p. 201]{AMR}. Hence, the map $\iota_A$ in \eqref{eq:Orbit_map_A} is a closed, injective smooth immersion and thus a smooth embedding (as in Abraham, Marsden, and Ratiu \cite[Definition 3.5.9, p. 201]{AMR}) of the manifold $\sG_E\cdot A$ into the affine space $\sA(E,h)$.

We define a subbundle $N_{\sG_E\cdot A}$ of the restriction of the tangent bundle $T\sA(E,h)$ to the orbit $\sG_E\cdot A$ of $A$ by
\[
N_{\sG_E\cdot A}
=
\{ (u\cdot A,a)\in\sA(E,h)\times W^{1,p}(T^*X\otimes \su(E): d_{u\cdot A}^*a=0 \,\text{for $u\in \sG_E$}\}.
\]
We see that $N_{\sG_E\cdot A}$ is a normal bundle for the orbit because the fiber of $N_{\sG_E\cdot A}$ over a point $A' = u\cdot A \in \sG_E\cdot A = \sG_E\cdot A'$ is equal to
\[
  \left(T_{A'}\sG_E\cdot A'\right)^\perp = \left(\Ran d_{A'}\right)^\perp = \Ker d_{A'}^*
  \subset W^{1,p}(T^*X\otimes \su(E)) = T_{A'}\sA(E,h),
\]
where ``$\perp$'' denotes $L^2$-orthogonal complement and we use the abbreviations
\begin{align*}
  \Ran d_{A'} &= d_{A'}\left(W^{2,p}(\su(E))\right),
  \\
  \Ker d_{A'}^* &= \Ker \left(d_{A'}^*: W^{1,p}(T^*X\otimes \su(E)) \to L^p(\su(E))\right).
\end{align*}
We argue that $N_{\sG_E\cdot A}$ is a smooth, $\sG_E$-equivariant vector bundle as follows. First, because the map
\[
W^{1,p}(T^*X\otimes \su(E))\ni a\mapsto u^{-1}a u\in W^{1,p}(T^*X\otimes \su(E))
\]
is an isometry when $u\in W^{w,p}(\SU(E))$, we see that $a$ is $L^2$-orthogonal to the image of $d_A$ if and only if $u^{-1}a u$ is $L^2$-orthogonal to the image of $u^{-1}\circ d_A\circ u=d_{u\cdot A}$.  Thus, $a\in \Ker d_A^*$ if and only if $u^{-1}au\in \Ker d_{u\cdot A}^*$. Hence, the map
\[
\tilde\sG_E\times \Ker d_A^* \ni(u,a)\mapsto (u\cdot A,u^{-1}au)\in N_{\sG_E\cdot A}
\]
gives a smooth trivialization of $N_{\sG_E\cdot A}$. Thus, $N_{\sG_E\cdot A}$ is a smooth vector bundle. Moreover, if $\sG_E$ acts on  $N_{\sG_E\cdot A}$ by
\[
\sG_E\times N_{\sG_E\cdot A} \ni (u,(A',a)) \mapsto (u\cdot A',u^{-1}au)\in N_{\sG_E\cdot A},
\]
then the exponential map (defined as in Kondracki and Rogulski \cite[Section 3.3, p. 39]{Kondracki_Rogulski_1986})
\[
  \exp: N_{\sG_E\cdot A} \ni (A',a') \mapsto A'+a' \in \sA(E,h).
\]
is $\sG_E$-equivariant. To see this, observe that
\[
  d_{u^*(A+a)}=u^{-1}\circ d_{A+a}\circ u=u^{-1}\circ (d_A+a)\circ u=u^{-1}\circ d_A\circ u + u^{-1}au=d_{u^*A}+u^{-1}au
\]
Hence, writing $A' = u\cdot A$ and $a' = u^{-1}au$, we obtain
\[
  \exp(A', a') = A'+a' = u\cdot A + u^{-1}au = u\cdot(A+a) = u\cdot\exp(A,a). 
\]
Moreover, $N_{\sG_E\cdot A}$ is clearly a smooth vector subbundle of $T\sA(E,h)\restriction\sG_E\cdot A$, with $\sG_E$-equivariant bundle projection $\pi:N_{\sG_E\cdot A} \to \sG_E\cdot A$ given by the restriction to $\sG_E\cdot A$ of the $\sG_E$-equivariant projection associated to the tangent bundle,
\[
  \pi:T\sA(E,h) \to \sA(E,h).
\]
Observe that $(D\exp)(A,0)$ is the identity operator on $T_A\sA(E,h) = \Ran d_A\oplus\Ker d_A^*$ and so there is an open neighborhood $O_A$ of $(A,0)$ in $N_{\sG_E\cdot A}$ such that $\exp\restriction O_A$ is a diffeomorphism onto an open neighborhood of $A$ in $\sA(E,h)$. Define
\[
  O_{\sG_E\cdot A} := \bigcup_{u\in\sG_E}u\cdot O_A.
\]
Note that $u\cdot O_A = O_{u\cdot A}$, for any $u\in\sG_E$. Hence, $\exp\restriction O_{u\cdot A} = u\cdot\exp\restriction O_A$ is a diffeomorphism onto an open neighborhood of $u\cdot A$ in $\sA(E,h)$, for any $u\in\sG_E$. Therefore, we obtain a
$\sG_E$-invariant
open neighborhood $O_{\sG_E\cdot A}$ of the zero section of $N_{\sG_E\cdot A}$ such that $\sT := \exp\restriction O_{\sG_E\cdot A}$ is a $\sG_E$-equivariant diffeomorphism from $O_{\sG_E\cdot A}$ onto a $\sG_E$-invariant open neighborhood of $\sG_E\cdot A$ in $\sA(E,h)$. (See Kondracki and Rogulski \cite[Lemma 3.3.3, p. 40]{Kondracki_Rogulski_1986} for a similar result.)

We now proceed as in the proof of Kondracki and Rogulski \cite[Corollary 3.3.5 (ii), pp. 41--42]{Kondracki_Rogulski_1986}. Suppose that $A \in \sA(E,h)$ and $A' \in \sG_E\cdot A$. For each $A' \in \sG_E\cdot A$, we define a submanifold,
\begin{equation}
  \label{eq:Kondracki_Rogulski_A'_slice}
  \sS_{A'}
  :=
  \sT\left(T_{A'}\sA(E,h)\cap O_{\sG_E\cdot A}\right)
  =
  \exp\left(T_{A'}\sA(E,h)\cap O_{\sG_E\cdot A}\right)
  \subset
  \sA(E,h),
\end{equation}
so $\sS_{A'}$ is transverse to $\sG_E\cdot A$ at $A'$.

Consider $A_1 \subset O_{\sG_E\cdot A}$, noting that $O_{\sG_E\cdot A}$ is an open neighborhood of $A$ in $\sA(E,h)$, and let $u \in \Stab(A_1)$, so $u\cdot A_1 = A_1$. Thus, $A_1 \in \sS_{A'}$ by definition \eqref{eq:Kondracki_Rogulski_A'_slice} for some $A' \in \sG_E\cdot A$ and so
\[
  A' = \pi\left(\sT^{-1}(A_1)\right) 
  = \pi\left(\sT^{-1}(u\cdot A_1)\right)
  = \pi\left(u\cdot\sT^{-1}(A_1)\right)
  = u\cdot\pi\left(\sT^{-1}(A_1)\right) = u\cdot A'.
\]
Consequently, $u \in \Stab(A') = \Stab(A) = C_r\,\id_E$. Therefore, $\Stab(A_1) = C_r\,\id_E$ since $u$ was arbitrary. Since $A_1\in O_{\sG_E\cdot A}$ was arbitrary, we conclude that $O_{\sG_E\cdot A}\subset\sA^{C_r}(E,h)$ and thus $\sA^{C_r}(E,h)$ is open in $\sA(E,h)$.
\end{proof}

\begin{rmk}[Slice property of $\sS_{A'}$ in \eqref{eq:Kondracki_Rogulski_A'_slice}]
\label{rmk:Slice_property_sSA'}  
More generally, the submanifold $\sS_{A'}$ has the property in Item \eqref{item:Duistermaat_Kolk_2-3-1_GlobalTube} of the Definition \ref{defn:Duistermaat_Kolk_2-3-1} of a slice through a point for the action of a Lie group on a smooth manifold due to Duistermaat and Kolk: If $u \in \sG_E$ and $A_1, u(A_1) \in \sS_{A'}$, then $u \in \Stab(A')$. To see this, observe that the definition of $\sS_{A'}$ as the image under $\sT$ of an open neighborhood of the origin in the fiber $N_{\sG_E\cdot A}|_{A'}$ of $N_{\sG_E\cdot A}$ over $A'$ implies that for any $A''\in \sS_{A'}$ we have $\sT^{-1}(A'')\in N_{\sG_E\cdot A}|_{A'}$. Hence, we see that
\begin{multline*}
  A' = \pi\left(\sT^{-1}(A_1)\right) 
  = \pi\left(\sT^{-1}\left(\sS_{A'}\right)\right)
  = \pi\left(\sT^{-1}(u\cdot A_1)\right)
  \\
  = \pi\left(u\cdot\sT^{-1}(A_1)\right)
  = u\cdot\pi\left(\sT^{-1}(A_1)\right) = u\cdot A',
\end{multline*}
where we have used the fact that $\pi$ and $\sT$ are $\sG_E$-equivariant. Thus, $u \in \Stab(A')$, as claimed.
\end{rmk}

\begin{proof}[Third proof of Lemma \ref{lem:Openness_moduli_subspace_non-split_projectively_Hermitian-Einstein_connections} \eqref{item:Openness_subspace_unitary_connections_minimal_stabilizer}]
Let $A\in \sA^{C_r}(E,h)$, so $A$ is a $W^{1,p}$ unitary connection on $E$ that induces the fixed smooth, unitary connection $A_d$ on the Hermitian line bundle $\det E$ and has $\Stab(A) \cong C_r$. By Corollary \ref{cor:Split_unitary_A_and_Lie_Alg_of_Stab(A)} \eqref{item:HA0_is_Lie_Alg_of_Stab(A)_Irreducible}, the connection $A$ has minimal stabilizer $C_r$ if and only if $\bH_A^0 = (0)$, where $\bH_A^0 \subset W^{2,p}(\su(E))$ is the finite-dimensional vector subspace defined in \eqref{eq:HE_equation_bHA0}. It follows from Rudin \cite[Theorem 4.12, p. 99]{Rudin} that
\begin{multline*}
  \Ker\left(d_A:W^{2,p}(\su(E))\to W^{1,p}(T^*\otimes\su(E))\right)
  \\
  =
  \Ran\left(d_A^*:W^{3,p}(T^*\otimes\su(E)) \to W^{2,p}(\su(E))\right)^\perp,
\end{multline*}
where $\perp$ denotes $L^2$-orthogonal complement. Thus, $\bH_A^0 = (0)$ if and only if $\Ran d_A^* = W^{2,p}(\su(E))$ and the latter is an open condition in the sense that if $a \in W^{1,p}(T^*X\otimes \su(E))$ is small, then $\Ran(d_A+a)^*= W^{2,p}(\su(E))$ and this proves Item \eqref{item:Openness_subspace_unitary_connections_minimal_stabilizer}. 
\end{proof}

\begin{proof}[Fourth proof of Lemma \ref{lem:Openness_moduli_subspace_non-split_projectively_Hermitian-Einstein_connections} \eqref{item:Openness_subspace_unitary_connections_minimal_stabilizer}]
We lightly adapt the argument due to Lawson for \cite[Chapter II, Theorem 10.4 (2), p. 33, and pp. 36--37]{Lawson}. By Corollary \ref{cor:Split_unitary_A_and_Lie_Alg_of_Stab(A)} \eqref{item:HA0_is_Lie_Alg_of_Stab(A)_Irreducible}, the connection $A$ has minimal stabilizer $C_r$ if and only if $\bH_A^0 = (0)$. We claim that $\sA(E,h)\less \sA^{C_r}(E,h) = \{A \in \sA(E,h): \bH_A^0 \neq (0)\}$ is a closed subspace of $\sA(E,h)$. Let $\{A_k\}_{k\in\NN} \subset \sA(E,h)\less \sA^{C_r}(E,h)$ be a sequence that converges in $W^{1,p}$ to a limit $A_\infty \in \sA(E,h)$. By definition of $\{A_k\}_{k\in\NN}$, there is a sequence $\{\xi_k\}_{k\in\NN} \subset W^{2,p}(\su(E))$ such that $d_{A_k}\xi_k = 0$ and $\|\xi_k\|_{L^p(X)} = 1$, for all $k\in\NN$. Observe that
\[
  d_{A_\infty}\xi_k = d_{A_k}\xi_k + [A_\infty - A_k,\xi_k] = [A_\infty - A_k,\xi_k],
\]
for all $k\in\NN$. By the Sobolev Multiplication Theorem (see, for example, Freed and Uhlenbeck \cite[Chapter 6, pp. 95--96]{FU}),
\begin{align*}
  \|d_{A_\infty}\xi_k\|_{W^{1,p}(X)}
  &\leq c_0\|A_\infty - A_k\|_{W^{1,p}(X)} \|\xi_k\|_{W^{2,p}(X)}
  \\
  &\leq c_0\|A_\infty - A_k\|_{W^{1,p}(X)}, \quad\text{for all } k \in \NN,
\end{align*}
for a constant $c_0 = c_0(g,h,p) \in [1,\infty)$. The final inequality above follows from the fact that $\|\xi_k\|_{L^p(X)} = 1$, the identity $d_{A_k}^*d_{A_k}\xi_k = 0$, convergence in $W^{1,p}$ of $A_k$ to $A_\infty$, and the $L^p$ estimate for the Laplace operators $d_{A_k}^*d_{A_k}$ (see Gilbarg and Trudinger \cite[Theorem 9.11, p. 235]{GT}),
\[
  \|\xi_k\|_{W_{A_k}^{2,p}(X)} \leq c_0\left(\|d_{A_k}^*d_{A_k}\xi_k\|_{L^p(X)} + \|\xi_k\|_{L^p(X)}\right)
  = c_0, \quad\text{for all } k \in \NN.
\]  
Thus, $\|d_{A_\infty}\xi_k\|_{W_{A_\infty}^{1,p}(X)} \to 0$ as $k\to\infty$ and we may assume without loss of generality by reindexing the sequence that $\|d_{A_\infty}\xi_k\|_{W_{A_\infty}^{1,p}(X)} \leq 1$, for all $k\in\NN$. Hence,
\[
  \|\xi_k\|_{W_{A_\infty}^{2,p}(X)}
  \leq
  \|d_{A_\infty}\xi_k\|_{W_{A_\infty}^{1,p}(X)} + \|\xi_k\|_{L^p(X)}
  \leq 2, \quad\text{for all } k \in \NN.
\]
The Rellich--Kondrachov Compactness Theorem (see Adams and Fournier \cite[Theorem 6.3, p. 168]{AdamsFournier}) implies that the continuous embedding $W^{2,p}(X) \subset W^{1,p}(X)$ is compact and so, after passing to a subsequence, we obtain that $\xi_k \to \xi_\infty$ in $W^{1,p}$ as $k\to\infty$. Therefore, $\|\xi_\infty\|_{L^p(X)} = 1$ and
\[
  d_{A_\infty}\xi_\infty = \lim_{k\to\infty}d_{A_\infty}\xi_k = 0.
\]
By elliptic regularity for the Laplace operator $d_{A_\infty}^*d_{A_\infty}$, we see that $\xi_\infty\in W^{2,p}(\su(E))$. Consequently, $\bH_{A_\infty}^0 \neq (0)$ and this completes the proof.
\end{proof}

\begin{proof}[Fifth proof of Lemma \ref{lem:Openness_moduli_subspace_non-split_projectively_Hermitian-Einstein_connections} \eqref{item:Openness_subspace_unitary_connections_minimal_stabilizer}]
Let $A \in \sA^{C_r}(E,h)$, so $\Stab(A) \cong C_r$ and $A$ has central stabilizer in the sense of Definition \ref{defn:Reducible_split_trivial-stabilizer_unitary_connection} \eqref{item:central-stabilizer_unitary_connection}. Consider an open neighborhood $\UU_A \subset \sA(E,h)$, a connection $A+a \in \UU_A$, and suppose $u \in W^{2,p}(\SU(E))$ is such that
\begin{equation}
  \label{eq:u_in_SUE_stabilizer_d_A+a}
  u^*(d_A+a)
  = u^{-1}\circ(d_A+a)\circ u
  = d_A+a,
\end{equation}
where $u^*(d_A+a)$ is defined by the action of $u\in W^{2,p}(\SU(E))$ in \eqref{eq:Donaldson-Kronheimer_2-1-7_pullback}. We claim that $u = \varrho\,\id_E$, for some $\varrho \in C_r$. To prove the claim, we shall adapt part of the proof of Freed and Uhlenbeck \cite[Corollary, pp. 50--51]{FU}. The identity \eqref{eq:u_in_SUE_stabilizer_d_A+a} may be rewritten as
\[
  (d_A+a)\circ u - u\circ(d_A+a) = 0,
\]
or equivalently $d_{A+a}u = 0$ (see Remark \ref{rmk:SimplicityAndKernelOfConn_on_EndE}), or equivalently 
\[
  d_Au = -[a, u] = u\circ a - a\circ u.
\]
It will be very convenient to view $W^{2,p}(\SU(E)) \subset W^{2,p}(\gl(E))$, via the fiberwise embeddings $\SU(E_x) \subset \End(E_x) = \gl(E_x)$, where $\SU(E_x) \cong \SU(r)$ and $\gl(E_x) \cong \gl(r,\CC)$, for each $x\in X$. Note that $W^{2,p}(\gl(E)) \cong W^{2,p}(\fsl(E))\oplus W^{2,p}(X,\CC)$ via the isomorphism $\zeta \mapsto \zeta_0 + f\,\id_E$, where $f := \frac{1}{r}\tr_E\zeta$ and $\zeta_0 := \zeta - f\,\id_E$. Moreover, $W^{2,p}(\fsl(E)) \cong W^{2,p}(\su(E)) \oplus W^{2,p}(i\su(E))$ via the isomorphism $\zeta \mapsto \frac{1}{2}(\zeta-\zeta^\dagger) + \frac{1}{2}(\zeta+\zeta^\dagger)$. By assumption, $\Stab(A) \cong C_r$ and so
\[
  \Ker d_A\cap W^{2,p}(\su(E)) = \bH_A^0 = T_{\id_E}\Stab(A) = (0).
\]
Consequently, we obtain
\[
  \Ker d_A\cap W^{2,p}(\fsl(E)) = (0).
\]
Therefore,
\[
  \Ker d_A \cap W^{2,p}(\gl(E)) \cong \Ker d \cap W^{2,p}(X,\CC).
\]  
Consider the $L^2$-orthogonal decomposition,
\[
  W^{2,p}(\gl(E)) = \Ker d_A \cap W^{2,p}(\gl(E)) \oplus \left(\Ker d_A \cap W^{2,p}(\gl(E))\right)^\perp.
\]
We write $u = f\,\id_E + u_0$ with respect to this orthogonal decomposition and obtain
\[
  d_Au_0 = u_0\circ a - a\circ u_0 \quad\text{and}\quad df = 0.
\]
Thus, $f \equiv c \in \CC$, a constant since $X$ is connected by hypothesis. We may schematically write
\[
  d_A^*d_Au_0 = d_A^*(u_0\circ a - a\circ u_0) = u_0\otimes \nabla_A a + a\otimes\nabla_Au_0.
\]
Observe that, for a constant $c_0 = c_0(g,h) \in [1,\infty)$ depending at most on the Riemannian metric $g$ on $X$ and Hermitian metric $h$ on $E$,
\[
  \|d_A^*d_Au_0\|_{L^p(X)} \leq c_0\|u_0\|_{L^\infty(X)}\|\nabla_Aa\|_{L^p(X)}
  + c_0\|\nabla_Au_0\|_{L^{2p}(X)}\|a\|_{L^{2p}(X)}.
\]
By applying the continuous Sobolev embeddings, $W^{2,p}(X) \subset L^\infty(X)$ and $W^{1,p}(X) \subset W^{2p}(X)$ for $d/2 < p < \infty$ (see Adams and Fournier \cite[Theorem 4.12, p. 85]{AdamsFournier}), we obtain
\[
  \|d_A^*d_Au_0\|_{L^p(X)} \leq c_1\|a\|_{W^{1,p}(X)}\|u_0\|_{W^{2,p}(X)},
\]
for a constant $c_1 = c_1(g,h,p) \in [1,\infty)$. Therefore, for constants $c_2,c_3\in [1,\infty)$ depending at most on $g,h,p$ and because $u_0$ is $L^2$-orthogonal to $\Ker d_A \cap W^{2,p}(\gl(E))$, we see that
\[
  \|u_0\|_{W^{2,p}(X)} \leq c_2\|d_A^*d_Au_0\|_{L^p(X)} \leq c_3\|a\|_{W^{1,p}(X)}\|u_0\|_{W^{2,p}(X)}.
\]
By shrinking $\UU_A$ enough to give $\|a\|_{W^{1,p}(X)} \leq 1/(2c_3)$, the preceding inequality yields a contradiction, $\|u_0\|_{W^{2,p}(X)} \leq (1/2)\|u_0\|_{W^{2,p}(X)}$, unless $u_0 \equiv 0$. Therefore, $u_0 \equiv 0$ and $u = c\,\id_E$, for a constant $c \in \CC$. Because $\det u = 1$, we must have $c^r = 1$, so $c$ is an $r$-th root of unity and $u \in C_r$. This completes the proof.
\end{proof}

We next have the

\begin{thm}[Moduli space of non-split, regular projectively Hermitian--Einstein connections is a real analytic manifold]
\label{thm:Kobayashi_7_4_19}
(See Kobayashi \cite[Theorem 7.4.19, p. 243]{Kobayashi_differential_geometry_complex_vector_bundles} and Kim \cite{Kim_1987}.)
Let $E$ be a smooth Hermitian vector bundle over a closed, complex K\"ahler manifold $(X,g,J)$ of dimension $n$ with K\"ahler form $\omega = g(\cdot,J\cdot)$ and let $p\in(n,\infty)$ be a constant. Then the following hold:
\begin{enumerate}
\item\label{eq:M*(E,h,omega)_is_real_analytic_space}
  The moduli space $M^*(E,h,\omega)$ in \eqref{eq:Moduli_space_HE_connections_irreducible} of non-split projectively Hermitian--Einstein connections modulo $W^{2,p}(\SU(E))$ is a (possibly reduced) real analytic space with Zariski tangent spaces $\bH_A^1$ at points $[A]$ and smooth at points at points $[A]$ with $\bH_A^2=(0)$, where the harmonic spaces $\bH_A^\bullet$ are defined by the elliptic deformation complex \eqref{eq:HE_connection_fixed_unitary_determinant_connection_elliptic_deformation_complex}.

\item\label{eq:Mreg*(E,h,omega)_is_real_manifold}
  The moduli space $M_\reg^*(E,h,\omega)$ in \eqref{eq:Moduli_space_HE_connections_irreducible_regular} is a real analytic manifold with tangent spaces $\bH_A^1$ at points $[A] \in M_\reg^*(E,h,\omega)$.
\end{enumerate}
\end{thm}

\begin{rmk}[On the hypothesis in Theorem \ref{thm:Kobayashi_7_4_19} that $(X,g,J)$ is a complex K\"ahler manifold]
\label{rmk:Theorem_Kobayashi_7_4_19_almost_Hermitian}  
It is likely that the conclusions of Theorem \ref{thm:Kobayashi_7_4_19} would continue to hold if $(X,g,J)$ were only assumed to be a complex Hermitian or even almost Hermitian manifold.
\end{rmk}

\begin{rmk}[Version of Theorem \eqref{thm:Kobayashi_7_4_19}  due to Itoh]
\label{rmk:Other_versions_theorem_Kobayashi_7_4_19}  
Itoh \cite[Theorem 2]{Itoh_1985} proves a version of Theorem \ref{thm:Kobayashi_7_4_19} under stronger hypotheses, where he assumes that $X$ has dimension two, $\det E$ is holomorphically trivial, and restricts to the subset of points $[A]$ with $\bH_A^2 = (0)$.
\end{rmk}

\begin{rmk}[Complex analytic version of Theorem \ref{thm:Kobayashi_7_4_19}]
\label{rmk:Complex_analytic_version_theorem_Kobayashi_7_4_19}  
Itoh and Kobayashi state complex analytic versions of Theorem \ref{thm:Kobayashi_7_4_19}, in the sense that they further assert that $M^*(E,h,\omega)$ is a \emph{complex} analytic space and that $M_\reg^*(E,h,\omega)$ is a \emph{complex} manifold whereas we prefer, following Friedman and Morgan in \cite[Section 4.3.4, Theorem 3.9, p. 328]{FrM} and L\"ubke and Teleman in \cite[Theorem 4.4.1, p. 114]{Lubke_Teleman_1995}, to view $M^*(E,h,\omega)$ as a \emph{real} analytic space that is isomorphic (in the sense of real analytic spaces) to the open subset $\cM(E,\omega)$ of the \emph{complex} analytic space $\cM(E)$ given by the moduli subspace of stable, holomorphic vector bundles. Kobayashi considers equation \eqref{eq:Einstein_connection} as defined on the subspace of the complex  $\sA^{0,1}(E)$ of $(0,1)$-connections $\bar\partial_A$ on $E$ that are holomorphic (obeying $F_{\bar\partial_A} = 0$) whereas we (and Friedman and Morgan) consider the pair of equations \eqref{eq:Einstein_connection}, \eqref{eq:Holomorphic_connection_fixed_unitary_determinant_connection} as defined on the real affine space $\sA(E,h)$ of unitary connections $A$ on $E$. 
\end{rmk}

\section[Hitchin--Kobayashi correspondence for Hermitian--Einstein connections]{Hitchin--Kobayashi correspondence between projectively Hermitian--Einstein connections and polystable holomorphic vector bundles}
\label{sec:HK_correspondence_between_HE_connections_and_semistable_bundles}
The statement of the Hitchin--Kobayashi correspondence that we require is given by

\begin{thm}[Hitchin--Kobayashi correspondence between projectively Hermitian--Einstein connections and polystable holomorphic bundles]
\label{thm:Hitchin-Kobayashi_correspondence_Hermitian-Einstein_connections_stable_bundles}
Let $(E,h)$ be a complex, smooth Hermitian vector bundle over a closed complex K\"ahler manifold $(X,\omega)$ and $A_d$ be a smooth unitary connection on the Hermitian line bundle $\det E$ such that $F_{A_d}^{0,2} = 0$.
\begin{enumerate}
\item\label{item:HEconnection_implies_polystable_bundle}
\emph{(Hermitian--Einstein $\implies$ polystable.)} If $A$ is a smooth, projectively Hermitian--Einstein connection on $E$ in the sense of Definition \ref{defn:HE_connection}, then $(E,\bar\partial_A)$ is a polystable holomorphic bundle in the sense of Definition \ref{defn:Stable_holomorphic_structure}. Moreover, the following hold:

\begin{enumerate}
\item\label{item:HEconnection_implies_polystable_bundle_non-split}
If $A$ is non-split in the sense Definition \ref{defn:Reducible_split_trivial-stabilizer_unitary_connection} \eqref{item:Split_unitary_connection}, then $(E,\bar\partial_A)$ is stable.

\item\label{item:HEconnection_implies_polystable_bundle_split}
If $A$ is split in the sense Definition \ref{defn:Reducible_split_trivial-stabilizer_unitary_connection} \eqref{item:Split_unitary_connection}, then $(E,\bar\partial_A)$ is strictly polystable and for each summand $(E_k,\bar\partial_{A_k})$ in the direct sum \eqref{eq:Polystable_sum}, the Chern connection $A_k$ defined by $(h,\bar\partial_{A_k})$ has the same Einstein factor $\lambda$ as that of $A$.
\end{enumerate}

\item\label{item:Polystable_bundle_implies_HEconnection}
\emph{(Polystable $\implies$ Hermitian--Einstein.)} If $(E,\bar\partial_E)$ is a polystable holomorphic bundle in the sense of Definition \ref{defn:Stable_holomorphic_structure} that induces the holomorphic structure $\bar\partial_{A_d}$ on $\det E$, then there is a unique $\SU(E)$-orbit of a projectively Hermitian--Einstein connection $A$ on $E$ such that $A_{\det E} = A_d$ and $(E,\bar\partial_A)$ is isomorphic to $(E,\bar\partial_E)$ by an $\SL(E)$-gauge transformation. Moreover, the following hold:
\begin{enumerate}
\item\label{item:Polystable_bundle_implies_HEconnection_stable}
If $(E,\bar\partial_E)$ is stable, then $A$ is non-split in the sense of Definition \ref{defn:Reducible_split_trivial-stabilizer_unitary_connection} \eqref{item:Split_unitary_connection}.

\item\label{item:Polystable_bundle_implies_HEconnection_strictly_polystable}
If $(E,\bar\partial_E)$ is strictly polystable, then $A$ is split in the sense of Definition \ref{defn:Reducible_split_trivial-stabilizer_unitary_connection} \eqref{item:Split_unitary_connection}.
\end{enumerate}  
\end{enumerate}
\end{thm}

\begin{rmk}[On the hypothesis in Theorem \ref{thm:Hitchin-Kobayashi_correspondence_Hermitian-Einstein_connections_stable_bundles} that $(X,\omega)$ is a complex K\"ahler manifold]
\label{rmk:Hitchin-Kobayashi_correspondence_Hermitian-Einstein_connections_stable_bundles_non-Kaehler}
It is likely that the hypotheses in Theorem \ref{thm:Hitchin-Kobayashi_correspondence_Hermitian-Einstein_connections_stable_bundles} could be relaxed. For example, L\"ubke and Teleman \cite[Theorem 4.1.1, p. 92]{Lubke_Teleman_1995} obtain a version of this result when $(X,g,J)$ is only assumed to be complex Hermitian, while De Bartolomeis and Tian \cite[Theorem 0.1, p. 231]{DeBartolomeis_Tian_1996} prove a weaker version when $(X,g,J)$ is only assumed to be almost Hermitian.
\end{rmk}  

\begin{rmk}[Regularity of projectively Hermitian--Einstein and integrable $(0,1)$-connections in the Hitchin--Kobayashi correspondence]
\label{rmk:Hitchin-Kobayashi_correspondence_Hermitian-Einstein_connections_regularity}  
In Item \eqref{item:HEconnection_implies_polystable_bundle} of Theorem \ref{thm:Hitchin-Kobayashi_correspondence_Hermitian-Einstein_connections_stable_bundles}, if $A$ is a $W^{1,p}$ projectively Hermitian--Einstein connection, then it is immediate that $\bar\partial_A$ is a $W^{1,p}$ integrable $(0,1)$-connection. Conversely, although this is less obvious, if $\bar\partial_E$ is a $W^{1,p}$ integrable $(0,1)$-connection, then $A$ is a $W^{1,p}$ projectively Hermitian--Einstein connection. One approach to proving this conclusion that regularity is preserved by the Hitchin--Kobayashi correspondence is to examine its proof via gradient flow: for example, see Donaldson \cite{DonASD}, Donaldson and Kronheimer \cite{DK}, Jacob \cite{Jacob_2015conm}, and Kobayashi \cite{Kobayashi_differential_geometry_complex_vector_bundles}. See Feehan \cite{Feehan_yang_mills_gradient_flow_v4} for details concerning regularity of gradient flows and their limits. Preservation of regularity can also be understood by examining proofs of Hitchin--Kobayashi correspondences via the continuity method, as in Li and Yau \cite{Li_Yau_1987}, L\"ubke and Teleman \cite{Lubke_Teleman_1995}, and Uhlenbeck and Yau \cite{Uhlenbeck_Yau_1986}. Similar remarks apply if the connections are $W^{k,p}$ for $k\geq 2$ or $C^l$ for $l\geq 1$ or real analytic. We also recall that if $A$ is a $W^{1,p}$ Yang--Mills connection, then there is a $W^{2,p}$ unitary gauge transformation $u$ such that $u^*A$ is smooth by Wehrheim \cite[Theorem 9.4, p. 143]{Wehrheim_2004}. Hence, although the projectively Hermitian--Einstein connection $A$ produced by Item \eqref{item:Polystable_bundle_implies_HEconnection} of Theorem \ref{thm:Hitchin-Kobayashi_correspondence_Hermitian-Einstein_connections_stable_bundles} may only be of class $W^{1,p}$, there is a $W^{2,p}$ determinant one, unitary gauge transformation $u$ such that $u^*A$ is smooth.
\end{rmk}  

\begin{proof}[Proof of Theorem \ref{thm:Hitchin-Kobayashi_correspondence_Hermitian-Einstein_connections_stable_bundles}]
Item \eqref{item:HEconnection_implies_polystable_bundle} is stated without proof by Friedman and Morgan as \cite[Section 4.3.2, Theorem 3.5, p. 323]{FrM} (without a separation of the cases where $A$ is split or non-split) and proved by Donaldson and Kronheimer \cite[Proposition 6.1.11, p. 214]{DK} and by
Kobayashi \cite{Kobayashi_1982} and L\"ubke \cite{Lubke_1983}. In \cite[Proposition 6.1.11, p. 214]{DK}, \cite[Section 4.3.2, Theorem 3.5, p. 323]{FrM}, the authors assume that $E$ has complex rank two, $X$ is a surface, and $(\det E, \bar\partial_{A_d})$ is holomorphically trivial. However, Kobayashi \cite[Theorem 5.8.3, p. 163]{Kobayashi_differential_geometry_complex_vector_bundles} and L\"ubke and Teleman \cite[Theorem 2.3.2 and Remark 2.3.3, p. 55]{Lubke_Teleman_1995} provide statements and expositions of the proof of Item \eqref{item:HEconnection_implies_polystable_bundle} that allow $X$ to have any dimension and $E$ to have any rank. The assertion implicit in Item \eqref{item:HEconnection_implies_polystable_bundle_split} that the slope of $(E_k,\bar\partial_{E_k})$ is equal to the slope of $(E,\bar\partial_E)$ for all $k$ follows from Kobayashi \cite[Theorem 5.8.3, p. 163]{Kobayashi_differential_geometry_complex_vector_bundles} or L\"ubke and Teleman \cite[Theorem 2.3.2 and Remark 2.3.3, p. 55]{Lubke_Teleman_1995}, their conclusions that $E_k$ has the same Einstein factor \eqref{eq:Einstein_factor} as $E$ for all $k$, and the expression \eqref{eq:Slope} for the slope.

Consider Item \eqref{item:Polystable_bundle_implies_HEconnection}. The case of strictly polystable $(E, \bar\partial_E)$ in Item \eqref{item:Polystable_bundle_implies_HEconnection_strictly_polystable} is described by Friedman and Morgan in \cite[Section 4.3.2, paragraph prior to Theorem 3.6, p. 323]{FrM} for the case where $(\det E, \bar\partial_{A_d})$ is holomorphically trivial. The case of stable $(E, \bar\partial_E)$ is stated without proof by Friedman and Morgan as \cite[Section 4.3.2, Theorem 3.6 and Remark (1), p. 323]{FrM} and proved by Donaldson \cite[Theorem 1]{DonASD}, under the additional hypotheses that $X$ is a smooth projective surface having a Hodge metric and corresponding ample line bundle. Item \eqref{item:Polystable_bundle_implies_HEconnection} is proved by Donaldson and Kronheimer \cite[Theorem 6.1.5, p. 210 and Section 6.1.4, p. 215]{DK}), under the additional hypothesis that $(\det E, \bar\partial_{A_d})$ is holomorphically trivial. The case of stable $(E, \bar\partial_E)$ in Item \eqref{item:Polystable_bundle_implies_HEconnection_stable} is also proved by Li and Yau \cite[Theorem 1]{Li_Yau_1987} and Uhlenbeck and Yau \cite[Theorem, p. S262, or Theorem 4.1, p. S278]{Uhlenbeck_Yau_1986}, \cite{Uhlenbeck_Yau_1989}, but they allow $X$ to have any dimension and $E$ to have any rank.

Kobayashi proves Item \eqref{item:Polystable_bundle_implies_HEconnection} in \cite[Theorem 6.10.13, p. 214]{Kobayashi_differential_geometry_complex_vector_bundles} under the additional hypothesis that there exists an ample line bundle $L$ over $X$ such that $\omega$ represents the Chern class $c_1(L)$, but allows $X$ to have any dimension and $E$ to have any rank. L\"ubke and Teleman prove Item \eqref{item:Polystable_bundle_implies_HEconnection_stable} for stable $(E, \bar\partial_E)$ as \cite[Theorem 3.0.1, p. 61]{Lubke_Teleman_1995}, but also allow $X$ to have any dimension and $E$ to have any rank.

More precisely, L\"ubke and Teleman prove the existence of a Hermitian \emph{metric} $h'$ on $E$ that is Hermitian--Einstein in the sense of \cite[Definition 2.1.2, p. 46]{Lubke_Teleman_1995}, that is, the associated $h'$-unitary Chern connection $A'$ defined by the pair $(\bar\partial_E,h')$ is Hermitian--Einstein in the sense of Definition \ref{defn:HE_connection}. Moreover, the Hermitian metric $h'$ is unique up to multiplication by a positive constant. The induced Hermitian metric $h_{\det E}$ on $\det E$ is also Hermitian--Einstein by L\"ubke and Teleman \cite[Lemma 2.1.4 (iv), p. 47]{Lubke_Teleman_1995}. Two Hermitian metrics $h$ and $h'$ are, in general, related by a $\GL(E)$-gauge transformation of $E$ but we claim that here they are actually related by an $\SL(E)$-gauge transformation $v$ of $E$. To see this, observe that $A_d$ is by hypothesis the Chern connection defined by the pair $(\bar\partial_{\det E},h_{\det E})$ and is itself Hermitian--Einstein, so the given Hermitian metric $h_{\det E}$ on $\det E$ is Hermitian--Einstein. By the uniqueness assertion in \cite[Theorem 3.0.1, p. 61]{Lubke_Teleman_1995}, we may suppose (after scaling $h'$ by a positive constant) that $h_{\det E} = h_{\det E}'$ and this proves the claim. Consequently, the $h$-unitary Chern connection $A = v(A')$ defined by the pair $(\bar\partial_E,h)$ is Hermitian--Einstein in the sense of Definition \ref{defn:HE_connection} and projectively Hermitian--Einstein since the Chern connection on $\det E$ (namely, $A_d$) defined by the pair $(\bar\partial_{\det E},h_{\det E})$ is equal to the $h$-unitary connection on $\det E$ induced by the Chern connection on $E$ (namely, $A$) defined by the pair $(\bar\partial_E,h)$.

To prove Item \eqref{item:Polystable_bundle_implies_HEconnection_strictly_polystable} for $(E,\bar\partial_A)$ strictly polystable, recall from Definition \ref{defn:Stable_holomorphic_structure} that
\[
  (E,\bar\partial_E) = (E_1,\bar\partial_{E_1})\oplus\cdots\oplus (E_q,\bar\partial_{E_q}),
\]
for an integer $q$ in the range $2\leq q\leq r$, where the holomorphic subbundles $(E_k,\bar\partial_{E_k})$ are stable for $k=1,\ldots,q$. Therefore, Item \eqref{item:Polystable_bundle_implies_HEconnection_stable} for  stable $(E_k,\bar\partial_{E_k})$ provides Hermitian--Einstein connections $A_k$ on $E_k$ such that $(E_k,\bar\partial_{E_k})$ is isomorphic to $(E_k,\bar\partial_{A_k})$ via an $\SL(E_k)$-gauge transformation $v_k$ and $A_k$ is the Chern connection associated with $(\bar\partial_{A_k},h_{E_k})$ for $k=1,\ldots,q$. Consequently, there is an $\SL(E)$-gauge transformation $v=v_1\oplus\cdots\oplus v_q$ which yields the isomorphism,
\[
  (E,\bar\partial_E)
  \cong
  (E_1,\bar\partial_{A_1})\oplus\cdots\oplus (E_q,\bar\partial_{A_q}).
\]
We define a split Hermitian--Einstein connection
\[
  A := A_1\oplus\cdots\oplus A_q \quad\text{on } E = E_1\oplus\cdots\oplus E_q
\]
and denote $\bar\partial_A := \bar\partial_{A_1}\oplus\cdots\oplus\bar\partial_{A_q}$, so $A$ is the Chern connection associated with $(\bar\partial_A,h)$. Because $\bar\partial_E$ induces $\bar\partial_{A_d}$ on $\det E$ and $\bar\partial_E$ is equivalent to $\bar\partial_A$ via an $\SL(E)$-gauge transformation, then $\bar\partial_A$ also induces $(\bar\partial_A)_{\det E} = \bar\partial_{A_d}$ on $\det E$. The Chern connection $A$ associated with $(\bar\partial_A,h)$ on $E$ induces the Chern connection associated with $((\bar\partial_A)_{\det E},h_{\det E}) = (\bar\partial_{A_d},h_{\det E})$, namely $A_d$.

Because $(E,\bar\partial_E)$ is polystable, then each summand $(E_k,\bar\partial_{E_k})$, and thus also  each $(E_k,\bar\partial_{A_k})$, has slope equal to the slope of $(E,\bar\partial_E)$ by Definition \ref{defn:Stable_holomorphic_structure} for $k=1,\ldots,q$. Hence, the expression \eqref{eq:Einstein_factor} for the Einstein factor and expression \eqref{eq:Slope} for the slope implies that each Chern connection $A_k$ has the same Einstein factor, $\lambda$. But $A = A_1\oplus\cdots\oplus A_q$ and each $A_k$ is Hermitian--Einstein, so we obtain
\[
  F_A^{0,2} = F_{A_1}^{0,2}\oplus\cdots\oplus F_{A_q}^{0,2} = 0
\]
and
\[
  \Lambda F_A
  = \Lambda F_{A_1}\oplus\cdots\oplus \Lambda F_{A_q}
  = -i\lambda\,\id_{E_1} \oplus\cdots\oplus -i\lambda\,\id_{E_q}
  = -i\lambda\,\id_E.
\]
Therefore, $A$ is a Hermitian--Einstein connection on $E$ by Definition \ref{defn:HE_connection} and, since $A$ induces $A_d$ on $\det E$, then $A$ is a projectively Hermitian--Einstein connection on $E$. This completes the proof of Theorem \ref{thm:Hitchin-Kobayashi_correspondence_Hermitian-Einstein_connections_stable_bundles}.
\end{proof}

The following result provides context for the statement of Item \eqref{item:Polystable_bundle_implies_HEconnection} in Theorem \ref{thm:Hitchin-Kobayashi_correspondence_Hermitian-Einstein_connections_stable_bundles}. 

\begin{prop}[Structure of complex rank two, strictly semistable vector bundles over smooth complex projective varieties]
\label{prop:Friedman_4-21}  
(See Friedman \cite[Chapter 4, Proposition 21, p. 98]{FriedmanBundleBook}.)  
Let $X$ be a smooth complex projective variety of dimension $d$ and $H$ be an ample line bundle over $X$. Let $E$ be a complex rank $r$ vector bundle over $X$ with slope $\mu(E) := \frac{1}{r}c_1(E)\cdot H^{d-1}$ (see Friedman \cite[Chapter 4, Definition 1, p. 85]{FriedmanBundleBook}). If $r=2$ and $E$ is strictly semistable (see Friedman \cite[Chapter 4, Definition 4, p. 87]{FriedmanBundleBook} or Definition \ref{defn:Stable_holomorphic_structure}), then exactly one of the following holds:
\begin{enumerate}
\item There is a unique complex line subbundle $L_1\subset E$ with $\mu(L_1)=\mu(E)$. The quotient $E/L_1$ is torsion free and $E$ is canonically defined as an extension,
  \[
    0 \xrightarrow{} L_1 \xrightarrow{} E \xrightarrow{} L_2\otimes\sI_Z
  \]
  where $L_2$ is a complex line bundle over $X$ and $\sI_Z$ is the ideal sheaf of a divisor $Z \subset X$.
  
\item There are exactly two complex line subbundles $L_i\subset E$ for $i=1,2$ with $\mu(L_i)=\mu(E)$ for $i=1,2$ and $E = L_1\oplus L_2$ as a direct sum of complex line subbundles (so $E$ is polystable in the sense of Definition \ref{defn:Stable_holomorphic_structure}).
  
\item $E = L\oplus L$ as a direct sum of complex line subbundles and there are infinitely many complex line subbundles with slope $\mu(E)$, exactly corresponding to the choice of a line in $H^0(E\otimes L^{-1})$. 
\end{enumerate}
More precisely, the following holds: Suppose that $E$ is an arbitrary complex rank-two vector bundle over $X$ that is given as an extension,
\[
    0 \xrightarrow{} L_1 \xrightarrow{} E \xrightarrow{} L_2\otimes\sI_Z
  \]
where $L_i$ are complex line bundles over $X$ for $i=1,2$ and $\sI_Z$ is the ideal sheaf of a divisor $Z \subset X$ and $\mu(L_1)=\mu(E)$. Then $E$ is semistable and either $L_1$ is the unique destabilizing subbundle with torsion free quotient or $Z = \varnothing$ and $E = L_1\oplus L_2$ as a direct sum of complex line subbundles, that is, the extension splits.
\end{prop}

We have the following analogue of Proposition \ref{prop:ReducibleUnitaryFromH0} and Freed and Uhlenbeck \cite[Theorem 3.1, p. 47]{FU} for the case of a unitary connection $A$, adapted here to the case of a holomorphic $(0,1)$-connection $\bar\partial_E$ on a complex vector bundle $E$ over a closed, complex K\"ahler manifold such that the holomorphic bundle $(E,\bar\partial_E)$ is semistable.

\begin{lem}[Stabilizers of polystable holomorphic $(0,1)$-connections]
\label{lem:Split_(0,1)-connection_and_H0}
Let $(E,h)$ be a smooth Hermitian vector bundle over a closed, connected, complex K\"ahler manifold $(X,\omega)$ and $\bar\partial_E$ be a smooth, integrable $(0,1)$-connection on $E$. Let $\Stab(\bar\partial_E)$ denote the stabilizer of $\bar\partial_E$ in $W^{2,p}(\SL(E))$ and assume that $(E,\bar\partial_E)$ is polystable in the sense of Definition \ref{defn:Stable_holomorphic_structure}. Let $\bH_{\bar\partial_E}^0$ be the harmonic space in
\eqref{eq:bHdbarA0bullet} for the elliptic complex \eqref{eq:Complex_vector_bundle_dbar_A_elliptic_deformation_fixed_det_holomorphic_structure} defined by $\bar\partial_E$. The following are equivalent:
\begin{enumerate}
\item\label{item:Split_(0,1)-connection_and_H0_nonzero} $\bH_{\bar\partial_E}^0 \neq (0)$.
  
\item\label{item:Split_(0,1)-connection_and_H0_is_C} $\CC\subseteq \bH_{\bar\partial_E}^0 $.
  
\item\label{item:Split_(0,1)-connection_and_H0_splitting} $E=E_1\oplus E_2$ as a proper\footnote{In the sense that $E_i\subsetneq E$ for $i=1,2$.} direct sum of complex vector bundles and $\bar\partial_E$ is split with respect to this decomposition, so $\bar\partial_E = \bar\partial_{E_1}\oplus \bar\partial_{E_2}$.
  
\item\label{item:Split_(0,1)-connection_and_H0_stabilizer_contains_C} $\Stab(\bar\partial_E)$ contains a Lie subgroup isomorphic to $\CC^*$.
\end{enumerate}
If $E$ has rank two and $\Stab(\bar\partial_E)$ is a proper subgroup\footnote{More precisely, $\Stab(\bar\partial_E)$ is naturally isomorphic to a subgroup $G \subset \SL(r,\CC)$ and we assume that $G\neq \SL(2,\CC)$.} of $\SL(2,\CC)$ then Item \eqref{item:Split_(0,1)-connection_and_H0_nonzero} and the following are equivalent:
\begin{enumerate}
\setcounter{enumi}{4}
\item
\label{item:Split_(0,1)-connection_on_Rank2_and_H0_is_C}
 $\bH_{\bar\partial_E}^0\cong \CC$.
\item
\label{item:Split_(0,1)-connection_on_Rank2_and_H0_stabilizer_is_C} 
$\Stab(\bar\partial_E) \cong \CC^*$.
\end{enumerate}
\end{lem}

\begin{proof}
Let $A$ be the smooth projectively Hermitian--Einstein connection on $(E,h)$ in the sense of Definition \ref{defn:HE_connection} that is provided by the Hitchin--Kobayashi Correspondence Theorem \ref{thm:Hitchin-Kobayashi_correspondence_Hermitian-Einstein_connections_stable_bundles} and Remark \ref{rmk:Hitchin-Kobayashi_correspondence_Hermitian-Einstein_connections_regularity},
so that $(E,\bar\partial_A)$ is isomorphic to $(E,\bar\partial_E)$ by a gauge transformation $v\in W^{2,p}(\SL(E))$ and $\bar\partial_A = v(\bar\partial_E)$. To simplify notation, we shall assume without loss of generality in our proof of Lemma \ref{lem:Split_(0,1)-connection_and_H0} that $v=\id_E$ and thus $\bar\partial_A = \bar\partial_E$. (It is straightforward to relax this assumption.)

Item \eqref{item:Split_(0,1)-connection_and_H0_nonzero} implies Item \eqref{item:Split_(0,1)-connection_and_H0_is_C} because $\bH_{\bar\partial_E}^0$ is a complex vector space.

Assume that Item \eqref{item:Split_(0,1)-connection_and_H0_is_C} holds. If $\bH_{\bar\partial_E}^0 \neq (0)$ and $\bH_A^0$ is the harmonic space in \eqref{eq:HE_equation_bHAbullet} defined by $A$, then $\bH_A^0 \neq (0)$ since $\bH_{\bar\partial_E}^0 \cong \bH_A^0\otimes_\RR\CC$ by Theorem \ref{thm:Kobayashi_7-2-21}.   Consequently, there exists a non-zero section $\xi \in W^{2,p}(\su(E))$ with $d_A\xi = 0$. Hence, Proposition \ref{prop:ReducibleUnitaryFromH0} implies that $A$ is a split unitary connection on $E$, so $A=A_1\oplus A_2$ with respect to a splitting $E=E_1\oplus E_2$ as an orthogonal direct sum of Hermitian vector bundles. Therefore, Theorem \ref{thm:Hitchin-Kobayashi_correspondence_Hermitian-Einstein_connections_stable_bundles} implies that $\bar\partial_E = \bar\partial_{E_1}\oplus \bar\partial_{E_2}$ with respect to the splitting $E=E_1\oplus E_2$ as a direct sum of complex vector bundles and this verifies Item \eqref{item:Split_(0,1)-connection_and_H0_splitting}.

Assume that Item \eqref{item:Split_(0,1)-connection_and_H0_splitting} holds. For $i=1,2$, let $r_i$ be the rank of the complex vector bundle $E_i$ appearing in the decomposition $E=E_1\oplus E_2$. Observe that for any $\lambda \in \CC^*$, we may construct a complex gauge transformation $v \in W^{2,p}(\SL(E))$ by defining
\[
  v := \lambda^{r_2}\,\id_{E_1} \oplus \lambda^{-r_1}\,\id_{E_2},
\]
with respect to the given splitting $E=E_1\oplus E_2$ as a direct sum of complex vector bundles. Because $\bar\partial_E = \bar\partial_{E_1}\oplus \bar\partial_{E_2}$ with respect to this splitting, we obtain that $v(\bar\partial_E) = \bar\partial_E$ and so $v \in \Stab(\bar\partial_E)$. Thus, $\Stab(\bar\partial_E)$ contains a complex Lie subgroup isomorphic to $\CC^*$ and this verifies Item \eqref{item:Split_(0,1)-connection_and_H0_stabilizer_contains_C}.

Assume that Item \eqref{item:Split_(0,1)-connection_and_H0_stabilizer_contains_C} holds.
Because $\bH_{\bar\partial_E}^0$ is the Lie algebra of $\Stab(\bar\partial_E)$ and $\Stab(\bar\partial_E)$ contains a subgroup isomorphic to $\CC^*$, the Lie algebra $\bH_{\bar\partial_E}^0$ must contain a subspace isomorphic to $\CC$ and thus is non-zero.  This verifies Item \eqref{item:Split_(0,1)-connection_and_H0_nonzero} and hence the equivalence of Items
\eqref{item:Split_(0,1)-connection_and_H0_nonzero} through \eqref{item:Split_(0,1)-connection_and_H0_stabilizer_contains_C}.

We now assume that $E$ has rank two and Item \eqref{item:Split_(0,1)-connection_and_H0_nonzero} holds.
Since $\bH_A^0$ is the Lie algebra of the stabilizer $\Stab(A)$ of $A$ in $W^{2,p}(\SU(E))$ by Lemma \ref{lem:LieGroup_Structure_of_Stab(A)}, that stabilizer is necessarily isomorphic to $\{\pm 1\}$, $S^1$, or $\SU(2)$, so we must have $\Stab(A)\cong S^1$ or $\SU(2)$ and $\bH_A^0 \cong \RR$ or $\su(2)$, since $\bH_A^0 \neq (0)$. If $\bH_A^0 \cong \su(2)$, we would have $\bH_{\bar\partial_E}^0 \cong \fsl(2,\CC)$ and thus $\Stab(\bar\partial_E) \cong \SL(2,\CC)$ (since $\bH_{\bar\partial_E}^0$ is the Lie algebra of $\Stab(\bar\partial_E)$ by the proof of Item \eqref{item:Split_(0,1)-connection_and_H0_stabilizer_contains_C}), which contradicts our hypothesis on $\Stab(\bar\partial_E)$. Therefore, $\bH_{\bar\partial_E}^0 \cong \CC$ and this verifies Item \eqref{item:Split_(0,1)-connection_on_Rank2_and_H0_is_C}.

Assume that Item \eqref{item:Split_(0,1)-connection_on_Rank2_and_H0_is_C} holds. Thus, $\bH_{\bar\partial_E}^0 \cong \CC$ and since $\bH_{\bar\partial_E}^0$ is the Lie algebra of $\Stab(\bar\partial_E)$, we obtain $\Stab(\bar\partial_E) \cong \CC^*$ and this verifies Item \eqref{item:Split_(0,1)-connection_on_Rank2_and_H0_stabilizer_is_C}. Lastly, if Item \eqref{item:Split_(0,1)-connection_on_Rank2_and_H0_stabilizer_is_C} holds, then $\bH_{\bar\partial_E}^0\neq(0)$ and thus Item \eqref{item:Split_(0,1)-connection_and_H0_nonzero} holds.

This completes the proof of Lemma \ref{lem:Split_(0,1)-connection_and_H0}.
\end{proof}

\section[Hermitian--Einstein connections and holomorphic vector bundles]{Real analytic embedding of the moduli space of Hermitian--Einstein connections into the moduli space of holomorphic vector bundles}
\label{sec:Moduli_space_HE_connections_open_subspace_moduli_space_simple_holomorphic_structures}
We have the

\begin{thm}[Real analytic embedding of the moduli space of non-split projectively Hermitian--Einstein connections into the moduli space of simple holomorphic vector bundles]
\label{thm:Kobayashi_7_4_20}
(See L\"ubke and Teleman \cite[Theorem 4.4.1, p. 114 and Corollary 4.4.4, p. 118]{Lubke_Teleman_1995} and Miyajima \cite[Corollary, p. 319]{Miyajima_1989}.)  
Let $(E,h)$ be a smooth Hermitian vector bundle over a closed complex K\"ahler manifold $(X,\omega)$ and $A_d$ be a smooth unitary connection on the Hermitian line bundle $\det E$ such that $F_{A_d}^{0,2} = 0$. Then the following hold:
\begin{enumerate}
\item\label{item:Kobayashi_7_4_20_analytic_embedding}
There is an embedding in the sense of real analytic spaces from the moduli space of non-split projectively Hermitian--Einstein connections in \eqref{eq:Moduli_space_HE_connections} onto an open subspace of the moduli space $\cM^*(E)$ of simple holomorphic vector bundles in \eqref{eq:Moduli_space_simple_holomorphic_structures},
\begin{equation}
  \label{eq:Analytic_embedding_HE_moduli_space_irreducibles_onto_stable_bundles}
  M^*(E,h,\omega) \hookrightarrow \cM^*(E),
\end{equation}
given by the moduli subspace $\cM(E,\omega)$ of stable holomorphic vector bundles in \eqref{eq:Moduli_space_omega-stable_holomorphic_structures}.

\item\label{item:Kobayashi_7_4_20_complex_analytic_space}
The moduli space $\cM(E,\omega)$ is a (Hausdorff) complex analytic space.

\item\label{item:Kobayashi_7_4_20_set-theoretic_embedding} There is a set-theoretic injection from the moduli space $M(E,h,\omega)$ of projectively Hermitian--Einstein connections in \eqref{eq:Moduli_space_HE_connections},
\begin{equation}
  \label{eq:Set_theoretic_embedding_HE_moduli_space_onto_polystable_bundles}
  M(E,h,\omega) \hookrightarrow \cM(E),
\end{equation}
whose image is given by the moduli subset $\cM_\ps(E,\omega)$ of polystable holomorphic vector bundles in \eqref{eq:Moduli_space_omega-stable_holomorphic_structures}.
\end{enumerate}
\end{thm}

\begin{proof}
Although proofs due to L\"ubke and Teleman and due to Miyajima are included in the cited references, we shall give alternative proofs here because we shall subsequently rely on its generalization to the case of projective vortices and holomorphic pairs. The Hitchin--Kobayashi Correspondence Theorem \ref{thm:Hitchin-Kobayashi_correspondence_Hermitian-Einstein_connections_stable_bundles} immediately yields the set-theoretic bijections,
\begin{subequations}
  \begin{align}
    \label{eq:Set_theoretic_bijection_HE_moduli_space_and_polystable_bundles}
    M(E,h,\omega) &\leftrightarrow \cM_\ps(E,\omega),
    \\
    \label{eq:Set_theoretic_bijection_HE_moduli_space_non-split_and_stable_bundles}
    M^*(E,h,\omega) &\leftrightarrow \cM(E,\omega),
  \end{align}
\end{subequations}  
Thus, Item \eqref{item:Kobayashi_7_4_20_set-theoretic_embedding} is given by \eqref{eq:Set_theoretic_bijection_HE_moduli_space_and_polystable_bundles}.

We now prove Item \eqref{item:Kobayashi_7_4_20_analytic_embedding}, which asserts that 
\begin{inparaenum}[\itshape i\upshape)]
\item the set-theoretic bijection \eqref{eq:Set_theoretic_bijection_HE_moduli_space_and_polystable_bundles} is an isomorphism of real analytic spaces, and
\item $\cM(E,\omega)$ is an open subspace of the moduli space $\cM^*(E)$.
\end{inparaenum}
In \cite[Corollary 5.7.14, p. 159]{Kobayashi_differential_geometry_complex_vector_bundles}, Kobayashi proves the set-theoretic inclusion $\cM(E,\omega) \subset \cM^*(E)$ and in\footnote{In \cite[Theorem 7.4.20, p. 243]{Kobayashi_differential_geometry_complex_vector_bundles}, a hypothesis that the Hermitian--Einstein connections are irreducible is omitted, although Kobayashi uses this property in his proof \cite[Section 7.4, p. 245, second paragraph]{Kobayashi_differential_geometry_complex_vector_bundles}.} \cite[Theorem 7.4.20, p. 243]{Kobayashi_differential_geometry_complex_vector_bundles}, he proves that $\cM(E,\omega)$, identified with $M^*(E,h,\omega)$ under the image of the bijection \eqref{eq:Set_theoretic_bijection_HE_moduli_space_non-split_and_stable_bundles}, is an open subset of $\cM(E)$ and thus an open subset of $\cM^*(E)$.

Friedman and Morgan \cite[Section 4.3.4, Theorem 3.9, p. 329]{FrM} give a different proof that $\cM(E,\omega)$ is an open subset of $\cM^*(E)$ via their comparison of local Kuranishi models for open neighborhoods of points in the moduli spaces of $g$-anti-self-dual connections and stable holomorphic structures on $E$. We adapt their argument to our setting. Let $[\bar\partial_E] \in \cM(E,\omega)$ and let $[A] \in M^*(E,h,\omega)$ be the corresponding point given by the bijection \eqref{eq:Set_theoretic_bijection_HE_moduli_space_non-split_and_stable_bundles}.

By the analogue of the forthcoming Theorem \ref{thm:Local_Kuranishi_model_for_moduli_space_projective_vortices_StabAvarphi_idE} for a projective vortex with trivial stabilizer, an open neighborhood of $[A]$ in $M^*(E,h,\omega)$ is parametrized by a
\label{page:Kuranishi_model_K(A)}
\emph{Kuranishi model} $K(A)$ comprising the real analytic set $\beps(\bkappa^{-1}(0))$, where $\beps$ is a real analytic embedding from an open neighborhood $N_A \subset \bH_A^1$ of the origin into an open neighborhood of $A$ in the slice $A+\Ker d_A^*\cap W^{1,p}(T^*X\otimes\su(E))$ and $\bkappa:\bH_A^1 \supset N_A \to \bH_A^2$ is a real analytic map. The map $\pi\circ\beps:K(A) \to M(E,h,\omega)$ induced by composition of $\beps:\bkappa^{-1}(0) \to \sA(E,h)$ with the quotient map $\pi:\sA(E,h) \to \sA(E,h)/W^{2,p}(\SU(E))$ (endowed with the quotient topology) is open and so its image is an open neighborhood of $[A]$ in $M(E,h,\omega)$ (and thus its intersection with $M^*(E,h,\omega)$ is also open). Since $\pi\circ\beps$ is also continuous and injective, it is a homeomorphism.

By the forthcoming Theorem \ref{thm:Local_Kuranishi_model_for_simple_point_cH(E)} (or its analogue \ref{thm:Local_Kuranishi_model_for_strongly_simple_point_cP(E)} for a strongly simple holomorphic pair), to each representative $\bar\partial_E$ for a point in $\cM^*(E)$ we may associate a
\label{page:Kuranishi_model_cKdbarE}
\emph{Kuranishi model} $\cK(\bar\partial_E)$ comprising the complex analytic set $\bgamma(\bchi^{-1}(0))$, where $\bgamma$ is a holomorphic embedding from an open neighborhood $N_{\bar\partial_E} \subset \bH_{\bar\partial_E}^1$ of the origin into an open neighborhood of $\bar\partial_E$ in the slice $\bar\partial_E+\Ker \bar\partial_E^*\cap W^{1,p}(\Lambda^{0,1}(\fsl(E))$ and $\bchi:\bH_{\bar\partial_E}^1 \supset N_{\bar\partial_E} \to \bH_{\bar\partial_E}^2$ is a holomorphic map.

The map $\pi\circ\bgamma:\cK(\bar\partial_E) \to \cM(E)$ induced by composition of $\bgamma:\bchi^{-1}(0) \to \sA^{0,1}(E)$ with the quotient map $\pi:\sA^{0,1}(E) \to \sA^{0,1}(E)/W^{2,p}(\SL(E))$ (endowed with the quotient topology) is open and so its image is an open neighborhood of $[\bar\partial_E]$ in $\cM(E)$ (and thus also in $\cM^*(E)$). Unlike in the case of the Kuranishi model $K(A)$ discussed above, it is less clear that the map
$\pi\circ\bgamma$ is injective, although if it is injective, it will necessarily be a homeomorphism (see Friedman and Morgan \cite[Section 4.1.5, p. 300, second paragraph]{FrM}) since it is also clearly continuous. In his proof of \cite[Theoreem 7.3.37, p. 234]{Kobayashi_differential_geometry_complex_vector_bundles}, Kobayashi provides an analytic proof that $\pi\circ\bgamma$ is injective, but because his argument appears not to easily translate from the case of simple holomorphic structures to the case of strongly simple holomorphic pairs, we shall instead rely on the argument due to Friedman and Morgan described below.

By Friedman and Morgan \cite[Section 4.3.4, Theorem 3.8, p. 326]{FrM} or the analogue of the more general Theorem \ref{thm:Friedman_Morgan_4-3-8_projective_vortices} --- our corresponding result for pairs --- the canonical isomorphism \eqref{eq:Bijection_unitaryconnections_with_01connections} of real Banach affine spaces,
\[
  \pi_h^{0,1}:\sA(E,h) \ni A' \mapsto \bar\partial_{A'} \in \sA^{0,1}(E),
\]
induces an isomorphism of real analytic spaces from $K(A)$ onto $\cK(\bar\partial_E)$ and the bijection \eqref{eq:Set_theoretic_bijection_HE_moduli_space_non-split_and_stable_bundles} of quotient spaces. Therefore, the composite map,
\[
  M^*(E,h,\omega) \supset \left.\beps(\bkappa^{-1}(0))\right/W^{2,p}(\SU(E))
\to
\left.\bgamma(\bchi^{-1}(0))\right/W^{2,p}(\SL(E)) \subset \cM(E),
\]
is injective and thus a homeomorphism onto its image. Hence, the image $\cM(E,\omega)$ of $M^*(E,h,\omega)$ is an open subset of $\cM(E)$ (and thus also $\cM^*(E)$). This completes the proof of Item \eqref{item:Kobayashi_7_4_20_analytic_embedding}.

It remains to prove Item \eqref{item:Kobayashi_7_4_20_complex_analytic_space}. Just as in the proof of \cite[Section 4.3.4, Theorem 3.9, pp. 329--330]{FrM} due to Friedman and Morgan, the collection of Kuranishi models $\cK(\bar\partial_E)$ define local complex analytic charts on $\cM(E,\omega)$ and they induce transition maps on overlaps that are isomorphisms of complex analytic spaces (see Friedman and Morgan \cite[Section 4.1.5, Proposition 1.16, p. 296]{FrM}); when $\cM(E,\omega)$ is replaced by $\cM_\reg(E,\omega)$, the fact that the transition maps are holomorphic is also proved by Itoh \cite[Proposition 5.1, p. 859]{Itoh_1985} and Kim \cite[Lemma 2.3, p. 148]{Kim_1987}. Therefore, $\cM(E,\omega)$ is a complex analytic space in the sense of Grauert and Remmert \cite[Section 1.1.5, p. 7]{Grauert_Remmert_coherent_analytic_sheaves}. This proves Item \eqref{item:Kobayashi_7_4_20_complex_analytic_space} and completes the proof of Theorem \ref{thm:Kobayashi_7_4_20}.
\end{proof}  

\begin{rmk}[Openness in Theorem \ref{thm:Kobayashi_7_4_20}]
\label{rmk:Friedman_Morgan_4-3-9_openness}
Itoh proves the openness result in Item \eqref{item:Kobayashi_7_4_20_analytic_embedding} via \cite[Proposition 4.3, p. 857]{Itoh_1985}. While he assumes that $E$ has complex rank two, $X$ has complex dimension two, and $c_1(E)=0$, those assumptions do not appear to be required by his argument.
\end{rmk}

\begin{rmk}[On the hypotheses of Theorems \ref{thm:Kobayashi_7_4_19} and \ref{thm:Kobayashi_7_4_20}]
\label{rmk:Friedman_Morgan_4-3-9_isomorphism}
Friedman and Morgan prove versions of Theorems \ref{thm:Kobayashi_7_4_19} and \ref{thm:Kobayashi_7_4_20} as part of their \cite[Section 4.3.4, Theorem 3.9, p. 329]{FrM}, albeit under stronger hypotheses. Suppose that $X$ is a closed, smooth, complex projective surface, $g$ is a Hodge metric on $X$ with corresponding ample line bundle $L$, and $E$ is an $L$-stable, rank-two complex vector bundle over $X$ with holomorphically trivial determinant $\det E$. Then they prove that the moduli space $\cM(E,L)$ of $L$-stable (in the sense of \cite[Section 4.3.2, p. 322]{FrM}), holomorphic vector bundles $(E,\bar\partial_A)$ modulo the group $W^{2,p}(\SL(E))$ of determinant-one, complex automorphisms of $E$ is a complex analytic space and that the map $M^*(E,h,\omega) \ni [A] \mapsto [\bar\partial_A] \in \cM(E,L)$ is an isomorphism of real analytic spaces. 
\end{rmk}

\begin{rmk}[Complexifications of stabilizers]
\label{rmk:Complexification_stabilizers}  
It is known that $\SL(r,\CC)$ is the universal complexification of $\SU(r)$ (see Hilgert and Neeb \cite[Section 15.2, p. 570]{Hilgert_Neeb_structure_geometry_lie_groups}) and from this one could prove that $W^{2,p}(\SL(r,\CC))$ is the universal complexification of $W^{2,p}(\SU(r))$, essentially as asserted by Atiyah and Bott \cite[Section 8, p. 570, line 12]{Atiyah_Bott_1983}. We refer the reader to Hilgert and Neeb \cite[Definition 15.1.2 and Theorem 15.1.4, p. 566]{Hilgert_Neeb_structure_geometry_lie_groups} for the definition of a universal complexification of a Lie group and its existence and uniqueness. If $A$ is the Chern connection on $(E,h)$ defined by an $(0,1)$-connection $\bar\partial_E$ on $E$, then the stabilizer of $d_A = \partial_A + \bar\partial_A$ in $W^{2,p}(\SU(r))$ is equal to the stabilizer of $\bar\partial_E = \bar\partial_A$ in $W^{2,p}(\SU(r))$. The Banach Lie groups $W^{2,p}(\SU(r))$ and $W^{2,p}(\SL(r,\CC))$ both act (real analytically) on $\sA^{0,1}(E)$ and thus it would be tempting to assume that the stabilizer of $\bar\partial_E$ in $W^{2,p}(\SL(r))$ is equal to the universal complexification of the stabilizer of $\bar\partial_E$ in $W^{2,p}(\SU(r))$, especially since $\bH_A^0\otimes_\RR\CC = \bH_{\bar\partial_A}^0$ by Theorem \ref{thm:Kobayashi_7-2-21} when $A$ is a projectively Hermitian--Einstein connection on $(E,h)$. Unfortunately, elementary examples show that if $G_\CC$ is the universal complexification of a Lie group $G$, and both $G$ and $G_\CC$ act smoothly on a smooth manifold $M$, and $G_{x_0}$ is the stabilizer in $G$ of a point $x_0\in M$, then the universal complexification $(G_{x_0})_\CC$ of $G_{x_0}$ need not coincide with the stabilizer $(G_\CC)_{x_0}$ of $x_0$ in $G_\CC$. (We are grateful to Karl--Hermann Neeb for this observation.) However, while it does not hold in general, the equality $(G_{x_0})_\CC = (G_\CC)_{x_0}$ holds in certain interesting cases: see Heinzer and Huckleberry \cite[Section 3.2, p. 331]{Heinzner_Huckleberry_1999} and Teleman \cite[Proposition 1.3]{Teleman_2004}.
\end{rmk}

\section[Marsden--Weinstein symplectic reduction and symplectic quotients]{Marsden--Weinstein symplectic reduction and construction of symplectic quotients}
\label{sec:Marsden-Weinstein_symplectic_quotient_construction}
Following Kobayashi \cite[Section 7.5]{Kobayashi_differential_geometry_complex_vector_bundles}, we summarize the Marsden--Weinstein symplectic reduction method and the construction of symplectic quotients. For complementary expositions, see Donaldson and Kronheimer \cite[Section 6.5.1, p. 244]{DK}, Marsden and Weinstein \cite{Marsden_Weinstein_1974}, McDuff and Salamon \cite[Proposition 5.4.13, p. 224]{McDuffSalamonSympTop3}, and \cite[Appendix, Section 7.1]{Wells3} by Garc{\'\i}a--Prada. Lerman, Montgomery, and Sjamaar \cite{Lerman_Montgomery_Sjamaar_1993, Lerman_Sjamaar_1991, Lerman_Sjamaar_1996} provide versions of Marsden--Weinstein reduction and the construction of symplectic quotients for singular symplectic spaces. We begin with the

\begin{defn}[(Weak) symplectic form on a Banach manifold]
\label{defn:Symplectic_form_Banach_manifold}
(See Kobayashi \cite[Equation (7.5.1), p. 246]{Kobayashi_differential_geometry_complex_vector_bundles}; compare Abraham, Marsden, and Ratiu \cite[Definition 5.2.12, p. 353]{AMR} for the definition of weak and strong Riemannian metrics on a smooth manifold.)  
Let $V$ be a real, smooth Banach manifold. One calls $\omega \in \Omega^2(V,\RR)$ a \emph{weak symplectic form} on $V$ if it obeys the following conditions:
\begin{enumerate}  
\item
\label{item:Symplectic_form_Banach_manifold_continuous}  
For each $x \in V$, the function $\omega_x : T_xV \times T_xV \to \RR$ is continuous.
\item
\label{item:Symplectic_form_Banach_manifold_non-degenerate}
For each $x \in V$, the bilinear form $\omega_x$ is \emph{weakly non-degenerate}, that is, if $\omega_x(u, v) = 0$ for all $v \in T_xV$, then $u = 0 \in T_xV$.
\item
\label{item:Symplectic_form_Banach_manifold_smooth}
$\omega_x$ is $C^\infty$ in $x \in V$.
\item
\label{item:Symplectic_form_Banach_manifold_closed}
$\omega$ is closed.
\end{enumerate}
One calls $\omega$ a \emph{symplectic form} on $V$ if $\omega_x$ is \emph{(strongly) non-degenerate} for each $x\in V$, that is, the map
\[
  \omega_x:T_xV \ni v \mapsto \omega_x(\cdot,v) \in T_x^*V
\]
is a isomorphism from the Banach space $T_xV$ onto its continuous dual space $T_x^*V$, for each $x\in V$.
\end{defn}

More briefly, $\omega$ is a weak symplectic form in the sense of Definition \ref{defn:Symplectic_form_Banach_manifold} if it is a closed, weakly non-degenerate two-form on $V$. For each $x \in V$, the (linear) map
\[
  T_xV \ni u \mapsto \omega_x(u,\cdot) \in T_x^*V
\]
is continuous by Item \eqref{item:Symplectic_form_Banach_manifold_continuous}  and injective by Item \eqref{item:Symplectic_form_Banach_manifold_non-degenerate} but not necessarily bijective (unless, for example, $V$ is finite-dimensional or a Hilbert space whose inner product is defined by $\omega$ and an almost complex structure).

\begin{defn}[Symplectic maps]
\label{defn:Symplectic_map}
Continue the notation of Definition \ref{defn:Symplectic_form_Banach_manifold}. A smooth map $f:V\to V$ is called \emph{symplectic} if $f^*\omega = \omega$.
\end{defn}

A vector field $a \in \Omega^0(TV)$ is called an \emph{infinitesimal symplectic transformation} if $\sL_a\omega = 0$, where $\sL_a\omega$ is the Lie derivative of $\omega$ with respect to $a$ and defined in the usual way (see Lee \cite[Equation (12.8)]{Lee_john_smooth_manifolds}). Cartan's magic formula gives (see Lee \cite[Theorem 14.35]{Lee_john_smooth_manifolds})
\[
  \sL_a\omega = \iota_ad\omega + d\iota_a\omega.
\]
Since $\omega$ is closed and thus $d\omega=0$, the condition $\sL_a\omega = 0$ thus reduces to $d\iota_a\omega = 0$, that is, the one-form $\iota_a\omega$ is closed, where the interior product is defined as usual by $\iota_a\omega(b) = \omega(a,b)$ for any vector field $b \in \Omega^0(TV)$. (We shall shortly consider the case where this one-form is exact.)

Let $G$ be a real Banach Lie group acting on $V$ as a group of symplectic transformations. Let $\fg$ be the Banach Lie algebra of $G$, and $\fg^*$ its dual Banach space. A \emph{moment map} \label{page:Moment_map} for the action of $G$ on $V$ is a smooth map
\[
  \mu : V \to \fg^*
\]
such that 
\begin{equation}
\label{eq:Kobayashi_7-5-2}
\langle \xi, d\mu(x)v\rangle = \iota_{a_\xi(x)}\omega(v), \quad\text{for all } \xi \in \fg \text{ and }v \in T_xV,
\end{equation}
where
\[
  d\mu(x) : T_xV \to \fg^*
\]
is the derivative of $\mu$ at $x$, and $\langle\cdot,\cdot\rangle:\fg\times\fg^*\to\RR$ denotes the canonical pairing between $\fg$ and $\fg^*$, and $a_\xi \in \Omega^0(TV)$ is the smooth vector field on $V$ defined by $\xi \in \fg$ through the action of $G$ on $V$, namely $\rho:G \times V \to V$, by
\begin{equation}
\label{eq:DefineInfinitesimalLieGroupAction}
  a_\xi(x) := (D_1\rho)(\id_G,x)\xi, \quad\text{for all } x\in V.
\end{equation}
The identity \eqref{eq:Kobayashi_7-5-2} may be written more concisely as
\[
  d\langle \xi, \mu\rangle = \iota_{a_\xi}\omega,
\]
where $\langle\xi,\mu\rangle \in C^\infty(V,\RR)$ is the smooth function obtained from $\xi$ and $\mu$ by the canonical pairing between $\fg$ and its dual $\fg^*$. Given a moment map $\mu$, one defines the
\label{page:co-moment_map}
\emph{co-moment map} $\mu^*:\fg \to C^\infty(V,\RR)$ by (see Donaldson and Kronheimer \cite[Section 6.5, p. 244]{DK})
\[
  \mu^*(\xi) = \langle\xi,\mu\rangle, \quad\text{for all } \xi \in \fg,
\]
in which case the identity \eqref{eq:Kobayashi_7-5-2} may be also written as
\[
  d\mu^*(\xi) = \iota_{a_\xi}\omega.
\]  
A moment map $\mu$ may or may not exist. Its existence means that the closed one-form $\iota_{a_\xi}\omega$ is exact and is equal to $d\langle \xi,\mu\rangle$. One can show that a moment map is unique up to an additive factor: If $\mu'$ is another moment map, then $\mu'-\mu$ is a (constant) element of $\fg^*$.

\begin{hyp}[Marsden--Weinstein symplectic reduction conditions on the moment map]
\label{hyp:Marsden-Weinstein_symplectic_quotient_conditions}
(See Kobayashi \cite[Conditions (a), (b), (c), pp. 246--247]{Kobayashi_differential_geometry_complex_vector_bundles}.)
Let $(V,\omega)$ be as in Definition \ref{defn:Symplectic_form_Banach_manifold} and let $G$ be a real Banach Lie group acting symplectically on $V$.  We shall impose the following conditions on a moment map $\mu : V \to \fg^*$ as in \eqref{eq:Kobayashi_7-5-2}, where $\fg$ is the Lie algebra of $G$ and $\fg^*$ its continuous dual space.
\begin{enumerate}
\item
\label{item:Marsden-Weinstein_symplectic_quotient_conditions_equivariance}  
Assume that $\mu$ is \emph{equivariant} with respect to the \emph{coadjoint action} of $G$ in the sense that
\begin{equation}
\label{eq:Kobayashi_7-5-3}  
\mu(g(x)) = (\Ad(g))^*(\mu(x)), \quad\text{for all } g \in G \text{ and } x \in V.
\end{equation}
Then $G$ leaves $\mu^{-1}(0) \subset V$ invariant and the quotient
\begin{equation}
\label{eq:Kobayashi_7-5-4}  
W = \mu^{-1}(0)/G
\end{equation}
is called the \emph{reduced phase space}. 

\item
\label{item:Marsden-Weinstein_symplectic_quotient_conditions_weakly_regular_value}
Assume that $0 \in\fg^*$ is a \emph{weakly regular value} of $\mu$ in the sense that $\mu^{-1}(0)$ is an embedded submanifold of $V$ and for every $x \in \mu^{-1}(0)$, the inclusion $T_x\mu^{-1}(0) \subseteq \Ker d\mu(x)$ is an equality.

\item
\label{item:Marsden-Weinstein_symplectic_quotient_conditions_slice}
Assume that the action of $G$ on $\mu^{-1}(0)$ is free and that at each point $x \in \mu^{-1}(0)$ there is a \emph{local slice} $S_x \subset \mu^{-1}(0)$ for the action, that is, an open submanifold $S_x \subset \mu^{-1}(0)$ through $x$ which is transversal to the orbit $G\cdot x \subset V$ in the sense that
\[
  T_x\mu^{-1}(0) = T_xS_x + T_x(G\cdot x).
\]  
\end{enumerate}  
\end{hyp}

As Kobayashi notes \cite[Section 7.5, p. 247]{Kobayashi_differential_geometry_complex_vector_bundles}, if $0 \in \fg^*$ is a regular value of $\mu:V \to \fg^*$, so $d\mu(x) : T_xV \to \fg^*$ is surjective for every $x \in \mu^{-1}(0)$, then the Implicit Mapping Theorem guarantees that Condition \eqref{item:Marsden-Weinstein_symplectic_quotient_conditions_weakly_regular_value} holds in Hypothesis \ref{hyp:Marsden-Weinstein_symplectic_quotient_conditions}. Recall that if $\Ad : G \to \Aut(\fg)$ is the adjoint representation of $G$, then the \emph{coadjoint representation} $\Ad^* : G \to \Aut(\fg^*)$ (see Kirillov \cite[Section 4.2.2, Example 4.11, pp. 55--56]{Kirillov_2008}) is defined by $\Ad^*(g) :=\Ad(g^{-1})^*$ for all $g \in G$, that is,
\begin{equation}
  \label{eq:Coadjoint_representation}
  \langle \xi, (\Ad(g))^*\alpha \rangle = \langle \Ad(g^{-1})\xi, \alpha \rangle
\end{equation}
for all $g \in G$, $\xi \in \fg$, and $\alpha \in \fg^*$. We give the following criterion for the regularity of a moment map $\mu$ at point $x\in \mu^{-1}(0)$. (Compare McDuff and Salamon \cite[Lemma 5.2.5, p. 206]{McDuffSalamonSympTop3} and Liberman and Marle \cite[Chapter IV, Corollary 2.3, p. 196]{Libermann_Marle_symplectic_geometry_analytical_mechanics}.)

\begin{lem}[Regularity criterion for moment map]
\label{lem:MomentMapRegularityCriterion}
Let $(V,\omega)$ be as in Definition \ref{defn:Symplectic_form_Banach_manifold} and let $G$ be a real Banach Lie group with Lie algebra $\fg$ acting symplectically on $V$. If $L_x:\fg\to T_x V$ is the linear map defined by $L_x(\xi)=a_\xi(x)$, for each $x\in V$ and where $a_\xi$ is defined in \eqref{eq:DefineInfinitesimalLieGroupAction}, then
\begin{equation}
\label{eq:CharacterizeCokernelOfDerivativeOfMomentMap}
\Ker L_x = \left\{\xi\in\fg: \langle \xi,d\mu(x)v\rangle=0, \text{ for all } v\in T_x V\right\}.
\end{equation}
If $d\mu(x) : T_xV \to \fg^*$ is surjective, then $\Ker L_x = (0)$. Conversely, if the Lie algebra $\fg$ of $G$ is a reflexive Banach space, $\Ran d\mu(x)$ is a closed subspace, and $\Ker L_x = (0)$, then $d\mu(x) : T_xV \to \fg^*$ is surjective.
\end{lem}

\begin{proof}
From the definition of $L_x$, we have $\xi\in\Ker L_x$ if and only if $a_\xi=0$. The weak non-degeneracy of $\om$ (see Item \eqref{item:Symplectic_form_Banach_manifold_non-degenerate} of Definition \ref{defn:Symplectic_form_Banach_manifold}) implies that $\xi\in\Ker L_x$ if and only if
$(\iota_{a_\xi}\om)_x=0$.  The defining property of the moment map \eqref{eq:Kobayashi_7-5-2} then implies that $\xi\in \Ker L_x$ if and only if $\langle \xi,d\mu(x)v\rangle=0$ for all $v\in T_x V$, and this proves the equality \eqref{eq:CharacterizeCokernelOfDerivativeOfMomentMap}.

If $d\mu(x)$ is surjective, then \eqref{eq:CharacterizeCokernelOfDerivativeOfMomentMap} implies that $\xi\in\Ker L_x$ if and only if $\langle \xi,\xi^*\rangle=0$ for all $\xi^*\in\fg^*$ which implies $\xi=0$ by the Hahn--Banach theorem. We now prove the converse statement, assuming that $\fg$ is reflexive and the range of $d\mu(x)$ is closed.  Because the range is closed, there is a non-zero $\xi_0^{**}\in \fg^{**}$ with $\langle \xi^*,\xi_0^{**}\rangle=0$ for all $\xi^*$ in the range of $d\mu(x)$ by Brezis \cite[Corollary 1.8]{Brezis}. Because $\fg$ is reflexive by hypothesis, there is a vector  $\xi_0\in\fg$ that maps to $\xi_0^{**}$ under the isomorphism $\fg\to\fg^{**}$. Because $\langle \xi^*,\xi_0^{**}\rangle=0$ for all $\xi^*$ in the range of $d\mu(x)$, we obtain that $\langle \xi_0,\xi^*\rangle=0$ for all $\xi^*$ in the range of $d\mu(x)$. The characterization \eqref{eq:CharacterizeCokernelOfDerivativeOfMomentMap} implies that $\xi_0\in\Ker L_x$. Because $\xi_0$ maps to the non-zero $\xi_0^{**}$ under the isomorphism $\fg\to\fg^{**}$, $\xi_0$ is non-zero and so $\xi_0\in\Ker L_x$ contradicts the assumption that $\Ker L_x=(0)$.
\end{proof}

\begin{rmk}[Marsden--Weinstein symplectic reduction conditions]
\label{rmk:Marsden-Weinstein_symplectic_quotient_conditions}
More generally, regarding Item \eqref{item:Marsden-Weinstein_symplectic_quotient_conditions_equivariance}, for any $\alpha \in \fg^*$ one can define
\[ 
  W_\alpha = \mu^{-1}(\alpha)/G_\alpha
\]  
where $G_\alpha = \{g \in G: (\Ad(g))^*\alpha = \alpha\}$.

The following diagram illustrates the relationship between $V$, $\mu^{-1}(0)$, and $W$:
\begin{equation}
\label{eq:Kobayashi_7-5-5}  
\begin{CD}
  \mu^{-1}(0) &@>{\jmath}>> V
  \\
  @VV{\pi}V &
  \\
  W = \mu^{-1}(0)/G
\end{CD}  
\end{equation}
where $\jmath$ is the natural injection and $\pi$ is the projection.

Regarding Item \eqref{item:Marsden-Weinstein_symplectic_quotient_conditions_slice}, if we choose $S_x$ sufficiently small, then the projection $\pi:\mu^{-1}(0) \to W$ defines a homeomorphism of $S_x$ onto an open subset $\pi(S_x) \subset W$. This defines a local coordinate chart on $W$ and ensures that $W$ is a manifold, which may or may not be Hausdorff. In order to have a Hausdorff manifold, one must further assume that the action of $G$ on $\mu^{-1}(0)$ is \emph{proper}. In applications, one may need to consider the case where the action of $G$ may not be proper. 
\end{rmk}

We can now state the

\begin{thm}[Marsden--Weinstein reduction and symplectic quotient]
\label{thm:Marsden-Weinstein_symplectic_quotient}
(See Kobayashi \cite[Theorem 7.5.8, p. 247]{Kobayashi_differential_geometry_complex_vector_bundles} or Marsden and Weinstein \cite{Marsden_Weinstein_1974}.)
Let $V$ be a Banach manifold with a symplectic form $\omega_V$ in the sense of Definition \ref{defn:Symplectic_form_Banach_manifold}. Let $G$ be a Banach Lie group acting on $V$. If there is a moment map $\mu: V \to \fg^*$ satisfying Hypothesis \ref{hyp:Marsden-Weinstein_symplectic_quotient_conditions}, then there is a unique symplectic form $\omega_W$ on the reduced phase space $W = \mu^{-1}(0)/G$ such that
\[
  \pi^*\omega_W = \jmath^*\omega_V \quad\text{on } \mu^{-1}(0),
\]
where $\jmath$ and $\pi$ are as in the diagram \eqref{eq:Kobayashi_7-5-5}.
\end{thm}

\begin{rmk}[Construction of weakly non-degenerate two-forms on quotients]
\label{rmk:Marsden-Weinstein_non-degenerate_two-form_quotient}
If $\omega_V$ is merely a weakly non-degenerate two-form on $V$, so it obeys the properties of Definition \ref{defn:Symplectic_form_Banach_manifold} except for the condition \eqref{item:Symplectic_form_Banach_manifold_closed} that it be closed, and we only ask that $\omega_W$ be a weakly non-degenerate two-form on $W$, then Theorem \ref{thm:Marsden-Weinstein_symplectic_quotient} admits a relatively simple proof when $V$ admits a $G$-invariant almost complex structure $J$ that is compatible with $\omega$ and so $V$ admits a $G$-invariant weak Riemannian metric $g_V$ such that $\omega_V = g_V(\cdot,J_V\cdot)$. Since $J_V$ is $G$-invariant, then $J_V$ defines an almost complex structure $J_W$ on the quotient $W = \mu^{-1}(0)/G$, where $\omega_V$ and the smooth map $\mu:V\to\fg^*$ are related by \eqref{eq:Kobayashi_7-5-2}. Similarly, because $g_V$ is $G$-invariant, it defines a weak Riemannian metric $g_W$ on the quotient $W = \mu^{-1}(0)/G$ following the construction by Groisser and Parker \cite[Section 1, pp. 505--512]{GroisserParkerGeometryDefinite}. The weak Riemannian metric $g_W$ and almost complex structure $J_W$ are compatible, and thus define a weakly non-degenerate two-form $\omega_W = g_W(\cdot,J_W\cdot)$ on $W$.
\end{rmk}

\section[Moduli space of Hermitian--Einstein connections as symplectic quotient]{Marsden--Weinstein symplectic reduction and the moduli space of projectively Hermitian--Einstein connections as a symplectic quotient}
\label{sec:Marsden-Weinstein_reduction_moduli_space_HE_connections_symplectic_quotient}
We adapt and extend the discussion and results due to Donaldson and Kronheimer \cite[Section 6.5.3, p. 250]{DK} and Kobayashi \cite[Section 7.6]{Kobayashi_differential_geometry_complex_vector_bundles}. We construct a weak symplectic form on the space of connections and compute the moment map for the action of the group of special unitary gauge transformations in Section \ref{subsec:Symplectic_form_moment_map_affine_space_unitary_connections}.
In Section \ref{subsec:Kuranishi_model_holomorphic_structure}, we construct a Kuranishi model, in the sense of Lemma \ref{lem:EquivariantKuranishiLemma}, for an open neighborhood of a given holomorphic structure in the affine space of all holomorphic structures with fixed determinant. We apply the Marsden--Weinstein Symplectic Reduction Theorem \ref{thm:Marsden-Weinstein_symplectic_quotient} to prove that the moduli space of regular, non-split projective vortices is a complex K\"ahler manifold in Section \ref{subsec:Application_Marsden-Weinstein_theorem_construction_Kaehler_metric_holomorphic_structures}.

\subsection[Symplectic form and moment map on affine space of unitary connections]{Symplectic form and moment map on affine space of unitary connections on a Hermitian vector bundle over an almost Hermitian manifold}
\label{subsec:Symplectic_form_moment_map_affine_space_unitary_connections}
We begin with the

\begin{prop}[Weak symplectic form on the affine space of unitary connections on a Hermitian vector bundle over an almost Hermitian manifold]
\label{prop:Donaldson_Kronheimer_6-5-7_almost_Hermitian}
(See Donaldson and Kronheimer \cite[Equation (6.5.7), p. 250]{DK} and Kobayashi \cite[Equation (7.6.22), p. 253]{Kobayashi_differential_geometry_complex_vector_bundles}.) 
Let $(E,h)$ be a smooth Hermitian vector bundle over a closed, smooth almost Hermitian manifold $(X,g,J)$ with fundamental two-form $\omega = g(\cdot,J\cdot)$ as in \eqref{eq:Fundamental_two-form} and complex dimension $n$. If $p \in (n,\infty)$ is a constant, then the expression\footnote{In \cite[Equation (6.5.7), p. 250]{DK}, Donaldson and Kronheimer employ a different normalization factor in the integral expression on in \eqref{eq:Kobayashi_7-6-22}.}
\begin{equation}
  \label{eq:Kobayashi_7-6-22}
  \bomega(a,b) := n\int_X \tr(a\wedge b) \wedge \omega^{n-1}, \quad\text{for all } a,b \in W^{1,p}(T^*X\otimes\su(E))
\end{equation}
defines a weak symplectic form (Definition \ref{defn:Symplectic_form_Banach_manifold}) on the real Banach space $W^{1,p}(T^*X\otimes\su(E))$. 
\end{prop}

\begin{proof}
When $(X,\omega)$ is symplectic, the conclusion is stated by Donaldson and Kronheimer without proof in \cite[p. 250]{DK} and when $(X,\omega)$ is K\"ahler, the conclusion is proved by Kobayashi\footnote{The result is stated with a typographical error, since $\Phi^{-1}$ should be replaced by $\Phi^{n-1}$, where Kobayashi uses $\Phi$ to denote the K\"ahler form.} \cite[Equation (7.6.22), p. 253]{Kobayashi_differential_geometry_complex_vector_bundles}. We shall adapt Kobayashi's proof from \cite[Section 7.6]{Kobayashi_differential_geometry_complex_vector_bundles}. The almost complex structure $J$, Riemannian metric $g$, and symplectic form $\omega$ on $TX$ obey \cite[Equations (7.6.6) and (7.6.8), p. 251]{Kobayashi_differential_geometry_complex_vector_bundles}
\begin{align}
  \label{eq:Kobayashi_7-6-6}  
  g(Ju,Jv) &= g(u,v), 
  \\
  \label{eq:Kobayashi_7-6-8}
  \omega(u,v) &= g(u,Jv), \quad\text{for all } x \in X \text{ and } u, v \in T_xX.
\end{align}
We extend $J$ from $TX$ to $T_\CC X = TX\otimes_\RR\CC$ by $\CC$-linearity to decompose
\[
  T_\CC X  = T_\CC X^{1,0} \oplus T_\CC X^{0,1}
\]
into the complex $i$ and $-i$ eigen-subbundles of $T_\CC X$, respectively, with an isomorphism of real vector bundles given by
\begin{equation}
\label{eq:Kobayashi_7-6-3_TX}
  TX \ni u \mapsto z = \varphi(u) = \frac{1}{2}(u - iJu) \in T_\CC X^{1,0}.
\end{equation}
The Hermitian metric on the complex vector bundle $T_\CC X^{1,0}$ is related to $g$ and $\omega$ by \cite[Equations (7.6.5) and (7.6.7), p. 251]{Kobayashi_differential_geometry_complex_vector_bundles}
\begin{align}
  \label{eq:Kobayashi_7-6-5}  
  g(u,v) &= h(z,w) + h(w,z) = 2\Real\,h(z,w) = 2\Real\,h(\varphi(u),\varphi(v)), 
  \\
  \label{eq:Kobayashi_7-6-7}  
  \omega(u,v) &= \frac{1}{i}(h(z,w) - h(w,z)) = 2\Imag\,h(z,w) = 2\Imag\,h(\varphi(u),\varphi(v)),
  \\
  \notag
  &\qquad\text{for all } x \in X \text{ and } u, v \in T_xX  \text{ and } z=\varphi(u), w=\varphi(v) \in (T_\CC X^{1,0})_x. 
\end{align}
We also write $\Lambda_\CC^1(X) = T_\CC^*X = T^*X\otimes_\RR\CC = \Lambda^1(X)\otimes_\RR\CC$ and use the almost complex structure $J^*$ on $T_\CC^*X$ induced by the conjugate $\CC$-linear isomorphism $h:T_\CC X \ni z \mapsto h(\cdot,z) \in T_\CC^*X$ of Riemannian vector bundles (so $J^*\theta = Jh^{-1}(\theta)$) to decompose
\[
  T_\CC^*X  = T_\CC^*X^{1,0} \oplus T_\CC^*X^{0,1}
\]
into the complex $-i$ and $i$ eigen-subbundles of $T_\CC^*X$, respectively (following Kobayashi's convention \cite[Equation (7.6.12), p. 252]{Kobayashi_differential_geometry_complex_vector_bundles}) or, as customary,
\[
  \Lambda_\CC^1(X) = \Lambda^{1,0}(X) \oplus \Lambda^{0,1}(X).
\]  
For $p\in (n,\infty)$, define a real Banach space by \cite[Equation (7.6.9), p. 251]{Kobayashi_differential_geometry_complex_vector_bundles}
\begin{equation}
\label{eq:Kobayashi_7-6-9}  
  \sX := W^{1,p}(T^*X\otimes\su(E)).
\end{equation}
Recall from the Sobolev Embedding Theorem (see Adams and Fournier \cite[Theorem 4.12, p. 85]{AdamsFournier}) that $W^{1,p}(X) \subset L^{2p}(X)$ is a continuous embedding\footnote{The conclusion is immediate for $p\geq 2n$ while for $p<2n$ we must require that $p^* = 2np/(2n-p) \geq 2p$ or, equivalently, $p\geq n$.} for all $p \in [n,\infty)$. Given $a \in \sX$, we have the decomposition (see \eqref{eq:Decompose_a_in_Omega1suE_into_10_and_01_components})
\[
  a = \frac{1}{2}\left(a' + a''\right) \quad\text{over } X
\]
where (see Kobayashi \cite[Equation (7.6.10), p. 251]{Kobayashi_differential_geometry_complex_vector_bundles})
\begin{equation}
\label{eq:Kobayashi_7-6-10}
\begin{aligned}
  a' &\in  W^{1,p}(\Lambda^{1,0}(\fsl(E))),
  \\
  a'' &\in  W^{1,p}(\Lambda^{0,1}(\fsl(E))).
\end{aligned}
\end{equation}
As noted in \eqref{eq:Kobayashi_7-6-11}, the condition that $a$ be skew-Hermitian is equivalent to the condition \eqref{eq:Kobayashi_7-6-11}, namely,
\[ 
  a' = -(a'')^\dagger,
\]
where $a^\dagger = \bar{a}^\intercal$ is defined by taking the complex, conjugate transpose of complex matrix-valued representatives with respect to any local frame for $E$ that is orthonormal with respect to the Hermitian metric on $E$. Define an almost complex structure $\bJ$ on $\sX$ by Kobayashi \cite[Equation (7.6.12), p. 252]{Kobayashi_differential_geometry_complex_vector_bundles} (this sign convention is opposite to that of Itoh \cite[Equation (4.1)]{Itoh_1988})
\begin{equation}
\label{eq:Kobayashi_7-6-12}  
  \bJ a' := -ia' \quad\text{and}\quad \bJ a'' := ia''
\end{equation}
and define the complex Banach spaces by \cite[Equation (7.6.13), p. 252]{Kobayashi_differential_geometry_complex_vector_bundles}
\begin{equation}
\label{eq:Kobayashi_7-6-13}
\begin{aligned}
  \sZ &:=  W^{1,p}(\Lambda^{0,1}(\fsl(E))),
  \\
  \bar\sZ &:=  W^{1,p}(\Lambda^{1,0}(\fsl(E))).
\end{aligned}
\end{equation}
Hence, there is an isomorphism of real Banach spaces given by \cite[Equations (7.6.3) and (7.6.14), pp. 251--252]{Kobayashi_differential_geometry_complex_vector_bundles}
\begin{equation}
\label{eq:Kobayashi_7-6-3_appliedin_7-6-14}
\sX \ni a \mapsto \frac{1}{2}a'' = \frac{1}{2}(a-i\bJ a) \in \sZ,
\end{equation}
where we write $a''/2$ rather than $a''$ in \eqref{eq:Kobayashi_7-6-3_appliedin_7-6-14} because of our convention \eqref{eq:Decompose_a_in_Omega1suE_into_10_and_01_components} of writing $a = (a'+a'')/2$.
For $\alpha, \beta \in \sZ$, we define their \emph{local pointwise} Hermitian inner product $h(\alpha,\beta) = \langle\alpha,\beta\rangle$ over $X$ by \cite[Equations (7.6.18) and (7.6.19), p. 252]{Kobayashi_differential_geometry_complex_vector_bundles}
\begin{multline}
\label{eq:Kobayashi_7-6-18_and_19}  
  \langle\alpha,\beta\rangle := \Lambda\left(\frac{1}{i}\tr_E(\alpha\wedge\beta^\dagger)\right)
  \quad\text{and}\quad  \langle\alpha,\beta\rangle\omega^n = \frac{n}{i}\tr_E(\alpha\wedge\beta^\dagger)\wedge\omega^{n-1} \quad\text{over } X,
  \\
  \text{for all } \alpha, \beta \in W^{1,p}(\Lambda^{0,1}(\fsl(E))).
\end{multline}
(By Kobayashi \cite[Equations (7.6.16) and (7.6.18), p. 252]{Kobayashi_differential_geometry_complex_vector_bundles}, one can see that the definition \eqref{eq:Kobayashi_7-6-18_and_19} gives the same Hermitian inner product as that induced on the fibers of $T^*X\otimes\su(E)$ by the Riemannian metric on $T^*X$ and the Hermitian inner product $h$ on $E$.) We define their \emph{global} ($L^2$) Hermitian inner product on $\sZ$ by setting \cite[Equation (7.6.20), p. 252]{Kobayashi_differential_geometry_complex_vector_bundles}
\begin{align*}    
  \bh(\alpha,\beta) &:= \int_X\langle\alpha,\beta\rangle\omega^n
  \\
                    &= \int_X\Lambda\left(\frac{1}{i}\tr_E(\alpha\wedge\beta^\dagger\right)\omega^n
  \\
                    &= \int_X\frac{n}{i}\tr_E(\alpha\wedge\beta^\dagger)\wedge\omega^{n-1},
\end{align*}
that is,
\begin{equation}
\label{eq:Kobayashi_7-6-20}    
\bh(\alpha,\beta) := \int_X\frac{n}{i}\tr_E(\alpha\wedge\beta^\dagger)\wedge\omega^{n-1},
\quad\text{for all } \alpha, \beta \in W^{1,p}(\Lambda^{0,1}(\fsl(E))).
\end{equation}
The corresponding real global ($L^2$) weak inner product (in the sense of Abraham, Marsden, and Ratiu \cite[Definition 5.2.12, p. 352]{AMR}) on the real Banach space $W^{1,p}(X;T^*X\otimes\fu(E))$ is given by \cite[Equation (7.6.21), p. 252]{Kobayashi_differential_geometry_complex_vector_bundles}
\begin{align*}  
  \bg(a,b) &= \bh(a'',b'') + \bh(b'',a'') \quad\text{(by analogue of \eqref{eq:Kobayashi_7-6-5})}
  \\
           &= \int_X\frac{n}{i}\tr_E(a''\wedge (b'')^\dagger + b''\wedge(a'')^\dagger)\wedge\omega^{n-1}
 \quad\text{(by \eqref{eq:Kobayashi_7-6-20})}             
  \\
           &= \int_X\frac{n}{i}\tr_E(-a''\wedge b' - b''\wedge a')\wedge\omega^{n-1}
 \quad\text{(by \eqref{eq:Kobayashi_7-6-11})}             
  \\
  &= \int_X\frac{n}{i}\tr_E(a'\wedge b'' - a''\wedge b')\wedge\omega^{n-1},
\end{align*}
that is,
\begin{equation}
  \label{eq:Kobayashi_7-6-21}
   \bg(a,b) = \int_X\frac{n}{i}\tr_E(a'\wedge b'' - a''\wedge b')\wedge\omega^{n-1}.
\end{equation}  
The corresponding weak symplectic form (in the sense of Definition \ref{defn:Symplectic_form_Banach_manifold})
on the real Banach space $W^{1,p}(T^*X\otimes\su(E))$ is thus given by (compare \cite[Equation (7.6.8), p. 251, and Equation (7.6.22), p. 253]{Kobayashi_differential_geometry_complex_vector_bundles})
\[ 
\bomega(a,b) = \bg(a,\bJ b) = \int_X n\tr_E(a\wedge b)\wedge\omega^{n-1}
\quad\text{(by \eqref{eq:Kobayashi_7-6-21} and analogue of \eqref{eq:Kobayashi_7-6-8})},
\]
which agrees with the expression stated in \eqref{eq:Kobayashi_7-6-22}, since
\begin{align*}
  \bg(a,\bJ b) &= \int_X\frac{n}{i}\tr_E(a'\wedge (\bJ b)'' - a''\wedge (\bJ b)')\wedge\omega^{n-1}
  \\
               &= \int_X\frac{n}{i}\tr_E(a'\wedge ib'' - a''\wedge (-ib'))\wedge\omega^{n-1}
  \\
               &= \int_X n\tr_E(a'\wedge b'' + a''\wedge b')\wedge\omega^{n-1},                 
\end{align*}
noting that $\bJ b = -ib'+ib''$ with $(\bJ b)' = -ib'$ and $(\bJ b)'' = ib''$ (matching the definitions in Kobayashi \cite[Equations (7.6.2), (7.6.12), and (7.6.13), pp. 251--252]{Kobayashi_differential_geometry_complex_vector_bundles}). This completes the proof of Proposition \ref{prop:Donaldson_Kronheimer_6-5-7_almost_Hermitian}.
\end{proof}  

The Riemannian metric $\bg$ on $\sA(E,h)$ given by \eqref{eq:Kobayashi_7-6-21} restricts to define a Riemannian metric $\iota^*\bg = \bg \restriction \bV$ on any smoothly embedded Banach submanifold, $\biota:\bV \hookrightarrow \sA(E,h)$. Similarly, we obtain a closed two-form $\iota^*\bomega = \bomega \restriction \bV$ upon restriction to $\bV$. We shall subsequently choose $\bV$ so that $\bJ$ also restricts from $\sA(E,h)$ to $\bV$ and, in that situation, we have the following identity from \cite[Equation (7.6.9), p. 251]{Kobayashi_differential_geometry_complex_vector_bundles},
\[
  \iota^*\bomega(a,b) = \iota^*\bg(a,(\iota^*\bJ) b), \quad\text{for all } A \in \bV \text{ and } a,b \in T_A\bV,
\]
which we may more simply write as
\begin{equation}
  \label{eq:Kobayashi_7-6-22_bomega_as_fundamental_2-form_for_bg_and_bJ}
  \bomega(a,b) = \bg(a,\bJ b), \quad\text{for all } A \in \bV \text{ and } a,b \in T_A\bV,
\end{equation}
where $\bomega$ is given by \eqref{eq:Kobayashi_7-6-22}. Next, we have the following modification of similar, but not identical results in \cite{DK} and \cite{Kobayashi_differential_geometry_complex_vector_bundles}.

\begin{prop}[Equivariant moment map for the action of the Banach Lie group of determinant-one, unitary gauge transformations on the affine space of unitary connections over an almost K\"ahler manifold]
\label{prop:Donaldson_Kronheimer_6-5-8_almost_Kaehler}
(See Donaldson and Kronheimer \cite[Proposition 6.5.8, p. 251]{DK} for $\SU(E)$ gauge transformations and Kobayashi \cite[Equation (7.6.23), p. 253]{Kobayashi_differential_geometry_complex_vector_bundles} for $\U(E)$ gauge transformations.) Continue the hypotheses of Proposition \ref{prop:Donaldson_Kronheimer_6-5-7_almost_Hermitian}, but assume now that $\omega$ is closed, and thus $(X,g,J)$ is almost K\"ahler, and fix a smooth, unitary connection $A_d$ on $\det E$. Then
\begin{equation}
\label{eq:Moment_map_action_unitary_det_one_gauge_transformations_affine_space_unitary_connections}  
\bmu:\sA(E,h) \to \left(W^{2,p}(\su(E))\right)^*,
\end{equation}
where
\begin{equation}
\label{eq:Moment_map_action_unitary_det_one_gauge_transformations_affine_space_unitary_connections_expression}
  \langle \xi, \bmu(A) \rangle
  :=
  n\int_X \tr_E((\Lambda F_A + i\lambda\,\id_E)\xi)\,\omega^n
\quad\text{for all } \xi \in W^{2,p}(\su(E)),
\end{equation}
is an equivariant moment map in the sense of \eqref{eq:Kobayashi_7-5-3} for the action of the Banach Lie group $W^{2,p}(\SU(E))$ of determinant-one, unitary $W^{2,p}$ automorphisms of $E$ for the (weak) symplectic form $\bomega$ given by \eqref{eq:Kobayashi_7-6-22} on the Banach affine space $\sA(E,h) = A_0 + W^{1,p}(T^*X\otimes \su(E))$, where $A_0$ is a smooth, unitary connection on $E$ that induces $A_d$ on $\det E$ and $\lambda$ is the Einstein constant \eqref{eq:Einstein_factor} appearing in the equation \eqref{eq:Einstein_connection} for a projectively Hermitian--Einstein connection $A$ on $E$.
\end{prop}

Proposition \ref{prop:Donaldson_Kronheimer_6-5-8_almost_Kaehler} for unitary connections is a special case of the forthcoming more general Proposition \ref{prop:Donaldson_Kronheimer_6-5-8_almost_Kaehler_type1_pairs} for unitary pairs and which is proved in Section \ref{sec:Marsden-Weinstein_reduction_moduli_space_SO3_monopoles_symplectic_quotient}.

\begin{rmk}[Role of the Einstein constant]
\label{rmk:Moment_map_Einstein_constant}  
Since the Einstein constant $\lambda$ is determined by the connection $A_{\det} = A_d$ that the unitary connection $A$ induces on $\det E$, which is fixed for all $A\in\sA(E,h)$, one can replace the expression $\Lambda F_A + i\lambda\,\id_E$ in Proposition \ref{prop:Donaldson_Kronheimer_6-5-8_almost_Kaehler} by its trace-free component $\Lambda(F_A)_0$ since the trace component will be identically zero for all $A\in\sA(E,h)$ and thus \eqref{eq:Moment_map_action_unitary_det_one_gauge_transformations_affine_space_unitary_connections_expression} becomes
\begin{equation}
\label{eq:Moment_map_action_unitary_det_one_gauge_transformations_affine_space_unitary_connections_expression_tracefree}
  \langle \xi, \bmu(A) \rangle
  =
  n\int_X \tr_E((\Lambda F_A)_0\xi)\,\omega^n
\quad\text{for all } \xi \in W^{2,p}(\su(E)).
\end{equation}
The preceding simplification also follows from the facts that $\xi \in W^{2,p}(\su(E))$ and $\fu(E)\cong i\underline{\RR}\oplus\su(E)$ as an orthogonal direct sum of Riemannian vector bundles. We shall use the expression \eqref{eq:Moment_map_action_unitary_det_one_gauge_transformations_affine_space_unitary_connections_expression_tracefree} for $\bmu(A)$ in Section \ref{sec:Marsden-Weinstein_reduction_moduli_space_SO3_monopoles_symplectic_quotient}.
\end{rmk}

Note that in the statement of Proposition \ref{prop:Donaldson_Kronheimer_6-5-8_almost_Kaehler} we may write (compare Donaldson and Kronheimer \cite[Equation (6.5.9), p. 251]{DK})
\[
  \tr_E((\Lambda F_A+i\lambda\,\id_E)\xi)\,\omega^n
  =
  n\tr_E\left(\left(F_A+\frac{i\lambda}{n}\,\id_E\otimes\omega\right)\xi\right)\wedge\omega^{n-1}.
\]
We now give the

\begin{proof}[Proof of Proposition \ref{prop:Donaldson_Kronheimer_6-5-8_almost_Kaehler}]
We adapt and add detail to the formal proof of \cite[Proposition 6.5.8, p. 251]{DK} due to Donaldson and Kronheimer and the verification of the moment map property of \cite[Equation (7.6.23), p. 253]{Kobayashi_differential_geometry_complex_vector_bundles} due to Kobayashi. Note that
\[
  W^{2,p}(\su(E)) = T_{\id_E}W^{2,p}(\SU(E))
\]
is the Lie algebra of $W^{2,p}(\SU(E))$ (see Freed and Uhlenbeck \cite[Proposition A.2]{FU}. For $\xi \in W^{2,p}(\su(E))$, the derivative of the real, smooth function $\langle\xi,\bmu\rangle$ on $\sA(E,h)$ at a connection $A\in\sA(E,h)$ in a direction $a \in T_A\sA(E,h) = W^{1,p}(T^*X\otimes\su(E))$ is
\begin{equation}
\label{eq:Derivative_moment_map}
  \langle \xi, d\bmu(A)a\rangle = n\int_X \tr_E(d_Aa)\xi)\wedge\omega^{n-1}.
\end{equation}
Integrating by parts and using the fact that $\omega$ is closed, we can rewrite the integral on the right-hand side of equation \eqref{eq:Derivative_moment_map} as
\begin{equation}
\label{eq:Derivative_moment_map_integration_by_parts}
  \langle \xi, d\bmu(A)a\rangle = -n\int_X \tr_E(a\wedge d_A\xi)\wedge\omega^{n-1}.
\end{equation}
Recall that there is a smooth (right) action $W^{2,p}(\SU(E))\times\sA(E,h) \ni (u,A) \mapsto u^*(A) \in \sA(E,h)$
with (this is the case $\Phi=0$ in equations \eqref{eq:Differential_of_SU(E)_GaugeAction} and \eqref{eq:d_APhi^0})
\[
  \left.\frac{d}{dt}u_t^*(A)\right|_{t=0} = d_A\xi,
\]
for $u_t := \Exp(t\xi) \in W^{2,p}(\SU(E))$ with $t\in (-\eps,\eps)$ and $\xi \in W^{2,p}(\su(E))$. Thus, each $\xi \in W^{2,p}(\su(E))$ defines a smooth vector field $\bX_\xi$ on $\sA(E,h)$ given by
\begin{equation}
\label{eq:Vector_field_affine_space_connections_defined_by_Lie_algebra_group_gauge_transformations}
\bX_\xi(A) := d_A\xi \in T_A\sA(E,h), \quad\text{for all } \xi \in W^{2,p}(\su(E)) \text{ and } A\in\sA(E,h).
\end{equation}  
Hence, by comparing the expression for $\bomega$ in \eqref{eq:Kobayashi_7-6-22} and the right-hand side of \eqref{eq:Derivative_moment_map_integration_by_parts}, we see that
\begin{multline*}
  \langle \xi, d\bmu(A)a\rangle = -\bomega(a,\bX_\xi(A)) = \bomega(\bX_\xi(A),a),
  \\
  \text{for all } \xi \in W^{2,p}(\su(E)), A\in\sA(E,h)\text{ and } a \in T_A\sA(E,h).
\end{multline*}
We have thus verified the moment map condition \eqref{eq:Kobayashi_7-5-2}, namely
\[
  \left\langle\xi, d\bmu(A)a\right\rangle = \bomega(\bX_\xi(A),a),
  \text{for all } \xi \in W^{2,p}(\su(E)), A\in\sA(E,h)\text{ and } a \in T_A\sA(E,h),
\]
that is,
\[
  \left\langle \xi, d\bmu(A) \right\rangle = \iota_{\bX_\xi(A)}\bomega,
  \quad\text{for all } \xi \in W^{2,p}(\su(E)) \text{ and } A\in\sA(E,h).
\]
By adapting Kobayashi's verification of his \cite[Equation (7.6.32), p. 255]{Kobayashi_differential_geometry_complex_vector_bundles} and using
\[
  F_{u(A)} = u^{-1}F_A u \quad\text{and}\quad \Lambda F_{u(A)} = u^{-1}(\Lambda F_A)u,
  \quad\text{for all } u \in W^{2,p}(\SU(E)),
\]
just as in \cite[Equation (7.6.30), p. 254]{Kobayashi_differential_geometry_complex_vector_bundles}, we see that, for any $u \in W^{2,p}(\SU(E))$ and $A \in \sA(E,h)$,
\begin{align*}
  \langle \xi, d\bmu(u(A)) \rangle
  &=
    n\int_X \tr_E((\Lambda F_{u(A)}+i\lambda\,\id_E)\xi)\,\omega^n
    \quad\text{(by
\eqref{eq:Moment_map_action_unitary_det_one_gauge_transformations_affine_space_unitary_connections_expression})}
  \\
  &= n\int_X \tr_E\left(u^{-1}(\Lambda F_A+i\lambda\,\id_E)u\xi\right)\omega^n
  \\
  &= n\int_X \tr_E\left((\Lambda F_A+i\lambda\,\id_E)u\xi u^{-1}\right)\omega^n
  \\
  &= \langle u\xi u^{-1}, d\bmu(u(A)) \rangle \quad\text{(by \eqref{eq:Moment_map_action_unitary_det_one_gauge_transformations_affine_space_unitary_connections})},
    \quad \text{for all } \xi \in W^{2,p}(\su(E)),
\end{align*}
that is, using the definition \eqref{eq:Coadjoint_representation} of the coadjoint representation,
\begin{equation}
  \label{eq:Kobayashi_7-6-32}
  \bmu(u(A)) = ((\Ad(u))^*\bmu(A), \quad \text{for all } u \in W^{2,p}(\SU(E)) \text{ and } A \in \sA(E,h).
\end{equation}
This yields the equivariance with respect to the action of $W^{2,p}(\SU(E))$ required by \eqref{eq:Kobayashi_7-5-3} and completes the proof of Proposition \ref{prop:Donaldson_Kronheimer_6-5-8_almost_Kaehler}.
\end{proof}

By analogy with Donaldson and Kronheimer \cite[p. 251]{DK} for $\lambda = 0$ and Kobayashi \cite[Equation (7.6.33), p. 255]{Kobayashi_differential_geometry_complex_vector_bundles} for $\lambda\in\RR$, we observe that the expression \eqref{eq:Moment_map_action_unitary_det_one_gauge_transformations_affine_space_unitary_connections} yields the identification
\begin{equation}
  \label{eq:Kobayashi_7-6-33}
  \bmu^{-1}(0)\cap\sA(E,h) = \left\{ A \in \sA(E,h): \Lambda F_A = -i\lambda\,\id_E\right\},
\end{equation}
that is, $A \in \bmu^{-1}(0)\cap\sA(E,h)$ if and only if $A$ obeys the Einstein condition \eqref{eq:Einstein_connection}.

Following Kobayashi \cite[Section 7.1, p. 220, line 16]{Kobayashi_differential_geometry_complex_vector_bundles} and by analogy with Donaldson and Kronheimer\footnote{They denote this subset by $\sA^{1,1}(E,h)$ and assume that $E$ has complex rank two, $\det E$ is holomorphically trivial, and $X$ has complex dimension two.} \cite[p. 209]{DK}, we consider the subset of the Banach affine space $\sA(E,h)$ of $W^{1,p}$ unitary connections $A$ on $E$ that induce a fixed smooth, unitary connection $A_d$ on the complex line bundle $\det E$, as in \eqref{eq:Unitary_connection_detE_fixed}, and obey the trace-free component of equation \eqref{eq:Holomorphic_connection}, namely
\begin{equation}
  \label{eq:Subset_unitary_connections_defining_holomorphic_structure_on_E}
  \sH(E,h) := \left\{A \in \sA(E,h): (F_A^{2,0})_0 = 0 \right\},
\end{equation}
so $(F_A^{2,0})_0=0$ too and $(F_A)_0$ has type $(1,1)$ if $A \in \sH(E,h)$.

Hence, by \eqref{eq:Kobayashi_7-6-33} and Definition \ref{defn:HE_connection}, we see that $A$ is a projectively Hermitian--Einstein connection on $E$ if and only if it belongs to the subspace
\begin{equation}
  \label{eq:HE_connections_as_zero_set_moment_map_on_1-1_unitary_connections}
  \bmu^{-1}(0)\cap\sH(E,h) = \left\{ A \in \sH(E,h): \Lambda F_A = -i\lambda\,\id_E\right\},
\end{equation}
where $\sH(E,h)$ is as in \eqref{eq:Subset_unitary_connections_defining_holomorphic_structure_on_E}. (This identification is valid regardless of whether $J$ is integrable or $\omega$ is closed, so we may allow $(X,g,J)$ to be almost Hermitian.)

\subsection[Kuranishi model for a holomorphic structure]{Kuranishi model for an open neighborhood of a holomorphic structure}
\label{subsec:Kuranishi_model_holomorphic_structure}
We consider the subset of the Banach affine space $\sA^{0,1}(E)$ of $(0,1)$-connections $\bar\partial_E$ of class $W^{1,p}$ on $E$ that induce a fixed smooth, holomorphic structure $\bar\partial_{E_d}$ on the complex line bundle $\det E$, as in \eqref{eq:Holomorphic_structure_fixed_determinant}, and obey $F_{\bar\partial_E} = 0$ for $F_{\bar\partial_E} = \bar\partial_E\circ\bar\partial_E$ as in \eqref{eq:Donaldson_Kronheimer_2-1-50}, namely
\begin{equation}
  \label{eq:Subset_01_connections_defining_holomorphic_structure_on_E}
  \cH(E) := \left\{\bar\partial_E \in \sA^{0,1}(E): F_{\bar\partial_E} = 0 \right\},
\end{equation}
by close analogy with Kobayashi \cite[Section 7.1, p. 219, line 1]{Kobayashi_differential_geometry_complex_vector_bundles}. We further define
the following open subspaces of $\cH(E$,
\begin{subequations}
  \label{eq:cH*E_cHregE_cHreg*E}
  \begin{align}
    \label{eq:cH*E}
    \cH^*(E)
    &:= \left\{\bar\partial_E \in \cH(E): \bH_{\bar\partial_E}^0 = (0) \right\},
    \\
    \label{eq:cHregE}
    \cH_\reg(E)
    &:= \left\{\bar\partial_E \in \cH(E): \bH_{\bar\partial_E}^2 = (0) \right\},
    \\
    \cH_\reg^*(E)
    \label{eq:cHreg*E}
    &:= \cH^*(E)\cap \cH_\reg(E),
\end{align}
\end{subequations}
where the harmonic spaces $\bH_{\bar\partial_E}^\bullet$ are defined by \eqref{eq:bHdbarA0bullet}. Recall from Lemma \ref{lem:Simple_01-connection} that $\bH_{\bar\partial_E}^0 = (0)$ if and only if $\bar\partial_E$ is simple. The openness of the subspaces \eqref{eq:cH*E_cHregE_cHreg*E} in $\cH(E)$ follows from the fact that $\cM(E)$ in \eqref{eq:Moduli_space_holomorphic_structures} is the quotient,
\[
  \cM(E) = \left.\cH(E)\right/W^{2,p}(\SL(E)),
\]
and the fact that $\cM^*(E)$ in \eqref{eq:Moduli_space_simple_holomorphic_structures} and $\cM_\reg(E)$ in \eqref{eq:Moduli_space_regular_holomorphic_structures} are proved in Section \ref{sec:Moduli_space_simple_bundles} to be open subspaces of $\cM(E)$. Our main goal in this section is to prove

\begin{thm}[Local Kuranishi models for points in $\cH(E)$]
\label{thm:Local_Kuranishi_model_for_simple_point_cH(E)}
Let $E$ be a Hermitian vector bundle over a complex, Hermitian manifold $X$ and $\bar\partial_{E_d}$ be a fixed smooth, holomorphic structure on the complex line bundle $\det E$. If $p\in(n,\infty)$, where $n$ is the dimension of $X$, and $\bar\partial_E \in \cH(E)$ as in \eqref{eq:Subset_01_connections_defining_holomorphic_structure_on_E} and
\begin{equation}
  \label{eq:dbar_E_slice}
  S_{\bar\partial_E} := \bar\partial_E + \Ker\bar\partial_E^*\cap W^{1,p}(\Lambda^{0,1}(\fsl(E))),
\end{equation}
then there are open neighborhoods $U_{\bar\partial_E} \subset S_{\bar\partial_E}$ of $\bar\partial_E$ and $N_{\bar\partial_E} \subset \bH_{\bar\partial_E}^1$ of the origin such that the following hold:
\begin{enumerate}
\item\label{item:Kuranishi_model_dbarA_cH''(E)}
There are a holomorphic embedding $\bgamma:N_{\bar\partial_E}\to U_{\bar\partial_E}$ and a holomorphic map $\bchi:N_{\bar\partial_E}\to \bH_{\bar\partial_E}^2$ such that $\bgamma(0) = \bar\partial_E$ and $\bchi(0)=0$ and
\begin{equation}
  \label{eq:Kuranishi_model_dbarA_cH''(E)}
  \cH(E) \cap U_{\bar\partial_E} = \bgamma(\bchi^{-1}(0)\cap N_{\bar\partial_E}).
\end{equation}
The Taylor expansion of $\bchi(\tau)$ has $D\bchi(0) = 0$ and second-order term proportional to $\Pi_{\bar\partial_E}(\tau\wedge\tau)$, where $\Pi_{\bar\partial_E}$ is $L^2$-orthogonal projection from $L^p(\Lambda^{0,2}(\fsl(E)))$ onto $\bH_{\bar\partial_E}^2$.

\item\label{item:cM_dbarE^vir_is_bgamma(N_dbarE)}
For $L^2$-orthogonal projection $\Pi_{\Ran\bar\partial_E}$ onto the range of $\bar\partial_E$ in $L^p(\Lambda^{0,2}(\fsl(E)))$ and
\begin{equation}
  \label{eq:Atiyah_Hitchin_Singer_family_page_446_holomorphic_structures}
  \cH_{\bar\partial_E}^\vir(E)
  :=
  \left\{\bar\partial_E + \alpha \in \sA^{0,1}(E):
    \Pi_{\Ran\bar\partial_E} F_{\bar\partial_E + \alpha} = 0\right\}
\end{equation}
then the subset $\cH_{\bar\partial_E}^\vir(E) \cap U_{\bar\partial_E}$ is an embedded complex submanifold of $U_{\bar\partial_E}$ with tangent space $T_{\bar\partial_E}(\cH_{\bar\partial_E}^\vir \cap U_{\bar\partial_E}) = \bH_{\bar\partial_E}^1$ and
\begin{equation}
  \label{eq:cM_dbarE^vir_is_bgamma(N_dbarE)}
  \cH_{\bar\partial_E}^\vir \cap U_{\bar\partial_E} = \bgamma(N_{\bar\partial_E}). 
\end{equation}

\item\label{item:cH(E)_cap_UU_dbarE_biholomorphic_cH(E)_cap_U_dbarE_times_UidE}
If $\bar\partial_E$ is a simple point, then there is an open neighborhood $U_{\id_E} \subset W^{2,p}(\SL(E))/C_r$ of the identity $\id_E$ such that the natural map,
\begin{equation}
  \label{eq:dbar_E_slice_biholomorphic_map}
  U_{\bar\partial_E} \times U_{\id_E} \ni (\bar\partial_E + \alpha, v)
  \mapsto v^*(\bar\partial_E + \alpha) \in \UU_{\bar\partial_E},
\end{equation} 
is a biholomorphic map onto an open neighborhood $\UU_{\bar\partial_E} \subset \sA^{0,1}(E)$ of $\bar\partial_E$ and restricts to a biholomorphic map,
\[
  \left(\cH(E) \cap U_{\bar\partial_E}\right) \times U_{\id_E} \cong \cH(E) \cap \UU_{\bar\partial_E}.
\]
(We recall that $C_r \subset \SL(r,\CC)$ is the subgroup of $r$-th roots of unity when $E$ has rank $r$.)
\end{enumerate}
In particular, $\cH_\reg^*(E)$ in \eqref{eq:cHreg*E} is an embedded complex submanifold of $\sA^{0,1}(E)$.
\end{thm}  

By definition \eqref{eq:Subset_01_connections_defining_holomorphic_structure_on_E} of $\cH(E)$, we have
\[
  \cH(E)\cap S_{\bar\partial_E}
  = \left\{\bar\partial_E + \alpha \in S_{\bar\partial_E}: F_{\bar\partial_E + \alpha} = 0\right\},
\]
and so the definition \eqref{eq:Atiyah_Hitchin_Singer_family_page_446_holomorphic_structures} yields
\begin{equation}
  \cH(E)\cap U_{\bar\partial_E}
  \label{eq:cH''E_cap_slice_equals_zero_locus_F_cap_PhidbarA}
    = \left\{\bar\partial_E + \alpha \in \cH_{\bar\partial_E}^\vir(E)\cap U_{\bar\partial_E}:
    \left(\id-\Pi_{\Ran\bar\partial_E}\right) F_{\bar\partial_E + \alpha} = 0\right\}.
\end{equation}
The maps $\bgamma$ and $\bchi$ in Item \eqref{item:Kuranishi_model_dbarA_cH''(E)} of Theorem \ref{thm:Local_Kuranishi_model_for_simple_point_cH(E)} define a finite-dimensional, local complex analytic model space $\bgamma(\bchi^{-1}(0)\cap N_{\bar\partial_E})$ in the sense of Grauert and Remmert \cite[Section 1.1.2, p. 3]{Grauert_Remmert_coherent_analytic_sheaves}. The equality \eqref{eq:Kuranishi_model_dbarA_cH''(E)} thus gives a characterization of $\cH(E)\cap U_{\bar\partial_E}$ as a finite-dimensional, complex analytic model space.

Friedman and Morgan \cite[Section 4.1.5, p. 296, Remark]{FrM} describe a modification of the Kuranishi model constructed in Theorem \ref{thm:Local_Kuranishi_model_for_simple_point_cH(E)} that allows for points $\bar\partial_E$ that are not simple. Indeed, the open neighborhood $N_{\bar\partial_E}$ can be chosen to be invariant under the action of $\Stab(\bar\partial_E)/\CC^*$ and the map $\bchi$ is equivariant with respect to the action of $\Stab(\bar\partial_E)/\CC^*$.

Before proceeding to the proof of Theorem \ref{thm:Local_Kuranishi_model_for_simple_point_cH(E)}, we make some preliminary comments. Our proof of Theorem \ref{thm:Local_Kuranishi_model_for_simple_point_cH(E)} is motivated by the development of Friedman and Morgan \cite[Section 4.1.5, pp. 294--295, and Section 4.1.6, pp. 300-302]{FrM} and Kobayashi \cite[Section 7.3]{Kobayashi_differential_geometry_complex_vector_bundles} of the Kuranishi model, but modified using an analysis suggested by Taubes \cite{TauIndef}.

According to the Bianchi Identity \eqref{eq:Holomorphic_curvature_Bianchi_identity} for an $(0,1)$-connection $\bar\partial_E\in\sA^{0,1}(E)$, we have $\bar\partial_EF_{\bar\partial_E} = 0$. Thus, we obtain
\begin{align*}
  \bar\partial_EF_{\bar\partial_E+\alpha}
  &= \bar\partial_E\left(F_{\bar\partial_E} + \bar\partial_E\alpha + \alpha\wedge\alpha\right)
  \\
  &= \bar\partial_EF_{\bar\partial_E} + (\bar\partial_E\circ\bar\partial_E)\alpha
    + \bar\partial_E\alpha\wedge\alpha - \alpha\wedge\bar\partial_E\alpha
  \\
  &= F_{\bar\partial_E}\alpha + \bar\partial_E\alpha\wedge\alpha - \alpha\wedge\bar\partial_E\alpha
    \quad\text{(by \eqref{eq:Donaldson_Kronheimer_2-1-50} and \eqref{eq:Holomorphic_curvature_Bianchi_identity})}.
\end{align*}
Therefore, using the fact that wedge product of a complex-valued two-form and one-form commutes and using the anti-symmetry of the Lie bracket on $\fsl(r,\CC)$, where $E$ has complex rank $r$, the preceding equality yields
\begin{equation}
  \label{eq:dbar_E_F_dbar_E+alpha}
  \bar\partial_EF_{\bar\partial_E+\alpha}
  =
  F_{\bar\partial_E}\alpha - 2\alpha\wedge \bar\partial_E\alpha,
    \quad\text{for all } \alpha \in W^{1,p}(\Lambda^{0,1}(\fsl(E))).
\end{equation}
Therefore, if $\bar\partial_E \in \cH(E)$, so that $F_{\bar\partial_E} = 0$ by \eqref{eq:Subset_01_connections_defining_holomorphic_structure_on_E}, we obtain:
\begin{equation}
  \label{eq:F_dbar_E+alpha_Ker_dbar_E}
  F_{\bar\partial_E+\alpha} \in \Ker\bar\partial_E\cap L^p(\Lambda^{0,2}(\fsl(E))),
  \quad \text{for all } \alpha \in \Ker\bar\partial_E\cap W^{1,p}(\Lambda^{0,1}(\fsl(E))).
\end{equation}
Regarding the inclusion \eqref{eq:F_dbar_E+alpha_Ker_dbar_E}, we make the following comments.

\begin{rmk}[On the codomain of the holomorphic curvature map]
\label{rmk:Codomain_holomorphic_curvature_map} 
In the context of this discussion, Friedman and Morgan \cite[Section 4.1.5, p. 294, line 30]{FrM} assert that the ``complex version of the Bianchi identity given after [Definition] 1.5 [see \cite[Section 4.1.2, p. 284, line 14]{FrM} implies that the image of $F$ is actually contained in $\Ker\bar\partial_E$'', where $F$ is the map $\sA^{0,1}\ni \bar\partial_E \mapsto F_{\bar\partial_E} \in L^p(\Lambda^{0,2}(\fsl(E)))$. However, since $\bar\partial_E$ is fixed in the cited application, what would matter is that $\bar\partial_EF_{\bar\partial_E+\alpha} = 0$ for $\alpha \in \Ker\bar\partial_E^*\cap W^{1,p}(\Lambda^{0,1}(\fsl(E)))$. From the derivation of \eqref{eq:F_dbar_E+alpha_Ker_dbar_E}, we see that this inclusion holds when both $F_{\bar\partial_E}=0$ and $\alpha \in \Ker\bar\partial_E\cap W^{1,p}(\Lambda^{0,1}(\fsl(E)))$, but not necessarily for all points in the slice $S_{\bar\partial_E} = \bar\partial_E + \Ker\bar\partial_E^*\cap W^{1,p}(\Lambda^{0,1}(\fsl(E)))$ defined by \eqref{eq:dbar_E_slice}. Of course, if $X$ has complex dimension two, then we have $\bar\partial_EF_{\bar\partial_E+\alpha} = 0$ for all $\alpha \in W^{1,p}(\Lambda^{0,1}(\fsl(E)))$ since $\Lambda^{0,3}(\fsl(E)) = (0)$.
\end{rmk}

We now proceed to the

\begin{proof}[Proof of Theorem \ref{thm:Local_Kuranishi_model_for_simple_point_cH(E)}]
  Hodge theory provides $L^2$-orthogonal decompositions,
\begin{equation}
  \label{eq:Hodge_decomposition_W1p_Lambda01_slE}
\begin{aligned}
  W^{1,p}(\Lambda^{0,1}(\fsl(E)))
  &=
    \bar\partial_E\left(W^{2,p}(\fsl(E))\right)
    \oplus \Ker\bar\partial_E^*\cap W^{1,p}(\Lambda^{0,1}(\fsl(E)))
  \\
  &=
    \bar\partial_E\left(W^{2,p}(\fsl(E))\right)
    \oplus \bH_{\bar\partial_E}^1
    \oplus \bar\partial_E^*\left(W^{2,p}(\Lambda^{0,2}(\fsl(E)))\right)
  \\
  &= \Ker\bar\partial_E\cap W^{1,p}(\Lambda^{0,1}(\fsl(E)))
     \oplus \bar\partial_E^*\left(W^{2,p}(\Lambda^{0,2}(\fsl(E)))\right),  
\end{aligned}
\end{equation}
where $\bH_{\bar\partial_E}^1 = \Ker(\bar\partial_E+\bar\partial_E^*)\cap W^{1,p}(\Lambda^{0,1}(\fsl(E)))$ is defined by \eqref{eq:bHdbarA0bullet}, and
\begin{equation}
  \label{eq:Hodge_decomposition_Lp_Lambda02_slE}
\begin{aligned}
  L^p(\Lambda^{0,2}(\fsl(E)))
  &=
  \bar\partial_E\left(W^{1,p}(\Lambda^{0,1}(\fsl(E)))\right)
  \oplus \Ker\bar\partial_E^*\cap L^p(\Lambda^{0,2}(\fsl(E)))
  \\
  &=
  \bar\partial_E\left(W^{1,p}(\Lambda^{0,1}(\fsl(E)))\right)
  \oplus \bH_{\bar\partial_E}^2
  \oplus \bar\partial_E^*\left(W^{1,p}(\Lambda^{0,3}(\fsl(E)))\right)
  \\
  &= \Ker\bar\partial_E\cap\left(L^p(\Lambda^{0,2}(\fsl(E)))\right)
  \oplus \bar\partial_E^*\left(W^{1,p}(\Lambda^{0,3}(\fsl(E)))\right),
\end{aligned}
\end{equation}
where $\bH_{\bar\partial_E}^2 = \Ker(\bar\partial_E+\bar\partial_E^*)\cap W^{1,p}(\Lambda^{0,2}(\fsl(E)))$ is defined by \eqref{eq:bHdbarA0bullet}. In particular, the Hodge decomposition\eqref{eq:Hodge_decomposition_Lp_Lambda02_slE} gives
\[
  \Ker\bar\partial_E\cap L^p(\Lambda^{0,2}(\fsl(E)))
  =
  \bar\partial_E\left(W^{1,p}(\Lambda^{0,1}(\fsl(E)))\right) \oplus \bH_{\bar\partial_E}^2.
\]
From Friedman and Morgan \cite[Section 4.1.5, p. 294, line 22]{FrM} (see also Freed and Uhlenbeck \cite[Theorem 3.2, p. 49]{FU}), we recall that $S_{\bar\partial_E}$ in \eqref{eq:dbar_E_slice} is a local slice for the action of $W^{2,p}(\SL(E))/C_r$ on $\sA^{0,1}(E)$ at a simple point\footnote{Friedman and Morgan use $W^{2,p}(\SL(E))$ or $W^{2,p}(\GL(E))$ rather than $W^{2,p}(\SL(E))/C_r$ or $W^{2,p}(\GL(E))/\CC^*$, but neither of their choices leads to a bijection in \eqref{eq:dbar_E_slice_biholomorphic_map} because the stabilizer of $\bar\partial_E$ in $W^{2,p}(\SL(E))$ or $W^{2,p}(\GL(E))$ would be isomorphic to $C_r$ or $\CC^*$, respectively.} $\bar\partial_E$.
When $\bar\partial_E$ is simple, the orbit
\[
  W^{2,p}(\SL(E))\cdot \bar\partial_E = \left(W^{2,p}(\SL(E))/C_r\right)\cdot \bar\partial_E
\]
is a smoothly embedded complex submanifold of $\sA^{0,1}(E)$ and $S_{\bar\partial_E}$ is an affine complex submanifold that is transverse to this orbit at the point $\bar\partial_E$ relative to $\sA^{0,1}(E)$ (see \cite[Section 4.1.5, p. 294, line 23]{FrM}).

One way to see that $W^{2,p}(\SL(E))/C_r$ is a complex Lie group is to write
\[
  W^{2,p}(\SL(E))/C_r = W^{2,p}\left(\SL(E)/C_r\right),
\]
and observe that $C_r$ is a discrete, central subgroup of both $\SL(E_x) \cong \SL(r,\CC)$ and $\SU(E_x) \cong \SU(r)$, for all $x\in X$, where the latter isomorphism is defined by the Hermitian metric on $E$. Now $\SL(r,\CC)$ is equal to $\SU(r)_\CC$, the universal complexification of $\SU(r)$ by Hilgert and Neeb \cite[Section 15.2, p. 570]{Hilgert_Neeb_structure_geometry_lie_groups}, and it is well-known that $\SU(r)$ is a connected real Lie group (see, for example, Gupta and Mishra \cite[Theorem 4.5]{Gupta_Mishra_2018}).
According to Hilgert and Neeb \cite[Exercise 15.1.2 (a) (2), p. 570]{Hilgert_Neeb_structure_geometry_lie_groups}, we have
\[
  (\SU(r)/C_r)_\CC = \SU(r)_\CC/C_r = \SL(r,\CC)/C_r
\]
and so $\SL(r,\CC)/C_r$ is a complex Lie group. Thus, $W^{2,p}(\SL(E)/C_r)$ is a complex Lie group.

The Implicit Mapping Theorem (see, for example, Feehan and Maridakis \cite[Theorem F.1, p.  127]{Feehan_Maridakis_Lojasiewicz-Simon_coupled_Yang-Mills}) therefore provides an open neighborhood $U_{\bar\partial_E} \subset S_{\bar\partial_E}$ of $\bar\partial_E$ and an open neighborhood $U_{\id_E} \subset W^{2,p}(\SL(E))$ of the identity $\id_E$ such that the natural map \eqref{eq:dbar_E_slice_biholomorphic_map} is a biholomorphic map onto an open neighborhood $\UU_{\bar\partial_E} \subset \sA^{0,1}(E)$ of the simple point $\bar\partial_E$ (see Friedman and Morgan \cite[Section 4.1.5, p. 294, line 24]{FrM}). This verifies Item \eqref{item:cH(E)_cap_UU_dbarE_biholomorphic_cH(E)_cap_U_dbarE_times_UidE} since the equation $F_{\bar\partial_E}=0$ defining $\cH(E)$ in \eqref{eq:Subset_01_connections_defining_holomorphic_structure_on_E} is invariant under the action of $W^{2,p}(\SL(E))$.

For any (not necessarily simple) point $\bar\partial_E \in \cH(E)$, we have $F_{\bar\partial_E} = \bar\partial_E\circ\bar\partial_E = 0$ by \eqref{eq:Subset_01_connections_defining_holomorphic_structure_on_E}. We restrict the map
\[
  F:\sA^{0,1}(E)\ni\bar\partial_E \mapsto F_{\bar\partial_E} \in L^p(\Lambda^{0,2}(\fsl(E)))
\]
to the slice \eqref{eq:dbar_E_slice} and compose with the $L^2$-orthogonal projection implied by the Hodge decomposition \eqref{eq:Hodge_decomposition_Lp_Lambda02_slE},
\[
  \Pi_{\Ran\bar\partial_E}:
  L^p(\Lambda^{0,2}(\fsl(E)))
  \to
  \bar\partial_E\left(W^{1,p}(\Lambda^{0,1}(\fsl(E)))\right),
\]
to give a holomorphic map,
\begin{equation}
  \label{eq:cF_sA01_to_Ran_dbar_A_cap_LpLambda02_slE}
  \cF:S_{\bar\partial_E} \ni \bar\partial_E+\alpha
  \mapsto
  \Pi_{\Ran\bar\partial_E} F_{\bar\partial_E+\alpha}
  \in
  \bar\partial_E\left(W^{1,p}(\Lambda^{0,1}(\fsl(E)))\right).
\end{equation}
Because
\[
  F_{\bar\partial_E+\alpha}
  = (\bar\partial_E+\alpha)\circ(\bar\partial_E+\alpha)
  = F_{\bar\partial_E} + \bar\partial_E\alpha + \alpha\wedge\alpha,
  \quad\text{for all } \alpha \in W^{1,p}(\Lambda^{0,1}(\fsl(E))),
\]
the differential of $F$ at any point $\bar\partial_E \in \sA^{0,1}(E)$ is given by
\begin{equation}
  \label{eq:Differential_F_W1p_to_Lp}
  (DF)(\bar\partial_E): W^{1,p}(\Lambda^{0,1}(\fsl(E))) \ni \alpha \mapsto \bar\partial_E\alpha
  \in \Ker\bar\partial_E\cap L^p(\Lambda^{0,2}(\fsl(E))).
\end{equation}
The calculation \eqref{eq:Differential_F_W1p_to_Lp} thus gives 
\begin{equation}
  \label{eq:Differential_cF_KerdbarA*_cap_W1p_to_Ran_dbar_A_in_Lp}
  (D\cF)(\bar\partial_E): \Ker\bar\partial_E^*\cap W^{1,p}(\Lambda^{0,1}(\fsl(E)))
  \ni \alpha \mapsto \bar\partial_E\alpha
  \in \bar\partial_E\left(W^{1,p}(\Lambda^{0,1}(\fsl(E)))\right).
\end{equation}
The Hodge decomposition \eqref{eq:Hodge_decomposition_W1p_Lambda01_slE} implies that
\[
  \Ker\bar\partial_E^*\cap W^{1,p}(\Lambda^{0,1}(\fsl(E)))
  =
  \bH_{\bar\partial_E}^1 \oplus \bar\partial_E^*\left(W^{2,p}(\Lambda^{0,2}(\fsl(E)))\right),
\]
and so the partial derivative $(D_2\cF)(\bar\partial_E)$ of $\cF$ in \eqref{eq:cF_sA01_to_Ran_dbar_A_cap_LpLambda02_slE} with respect to the second component of the slice,
\[
  S_{\bar\partial_E}
  =  \bar\partial_E
  +  \bH_{\bar\partial_E}^1 \oplus \bar\partial_E^*\left(W^{2,p}(\Lambda^{0,2}(\fsl(E)))\right),
\]
is an isomorphism, namely
\begin{equation}
  \label{eq:Differential_cF_Ran_dbarA*_in_W1p_to_Ran_dbar_A_in_Lp_isomorphism}
  (D_2\cF)(\bar\partial_E): \bar\partial_E^*\left(W^{2,p}(\Lambda^{0,2}(\fsl(E)))\right)
  \ni \alpha \mapsto \bar\partial_E\alpha
  \in \bar\partial_E\left(W^{1,p}(\Lambda^{0,1}(\fsl(E)))\right).
\end{equation}
Our approach to building a finite-dimensional complex analytic model in the sense of Kuranishi \cite{Kuranishi} that represents all solutions to $F_{\bar\partial_E+\alpha} = 0$ for $\alpha\in\Ker\bar\partial_E^*\cap W^{1,p}(\Lambda^{0,3}(\fsl(E)))$ near the origin is most closely modeled on that of Taubes for the anti-self-dual equation in \cite[Equations (2.7) and (2.8), p. 524 and Theorem 3.2, p. 525]{TauIndef}. Our approach complements but is distinct from that of Kobayashi \cite[Section 7.3]{Kobayashi_differential_geometry_complex_vector_bundles} and more explicit than that of Friedman and Morgan \cite[Section 4.1.5]{FrM}, in particular with regard to the codomain problem mentioned in Remark \ref{rmk:Codomain_holomorphic_curvature_map}. The Hodge decomposition \eqref{eq:Hodge_decomposition_W1p_Lambda01_slE} implies that we may write any $\alpha\in W^{1,p}(\Lambda^{0,1}(\fsl(E)))$ as 
\[
  \alpha = \tau + \bar\partial_E\zeta + \bar\partial_E^*v,
\]
where $\tau \in \bH_{\bar\partial_E}^1$ and $\zeta \in W^{2,p}(\fsl(E))$ and $v \in W^{2,p}(\Lambda^{0,2}(\fsl(E)))$. Since we restrict the map $\bar\partial_E + \alpha \mapsto F_{\bar\partial_E + \alpha}$ to the slice $S_{\bar\partial_E}$, we have $\zeta=0$ in the preceding decomposition. Moreover, we may assume without loss of generality that
\[
  v \in \left(\Ker\bar\partial_E^*\cap W^{2,p}(\Lambda^{0,2}(\fsl(E)))\right)^\perp
  = \bar\partial_E\left(W^{3,p}(\Lambda^{0,1}(\fsl(E)))\right),
\]
where the latter equality is implied by the Hodge decomposition \eqref{eq:Hodge_decomposition_Lp_Lambda02_slE} for $L^p(\Lambda^{0,2}(\fsl(E)))$ and thus we may assume without loss of generality that $v=\bar\partial_E\beta$ for some $\beta \in W^{3,p}(\Lambda^{0,1}(\fsl(E)))$. We compute
\begin{align*}
  F_{\bar\partial_E + \tau + \bar\partial_E^*v}
  &= F_{\bar\partial_E+\tau} + 
  \bar\partial_{E+\tau}\bar\partial_E^*v + \bar\partial_E^*v\wedge\bar\partial_E^*v
  \\
  &= F_{\bar\partial_E+\tau} + (\bar\partial_E+\tau\wedge\cdot)\bar\partial_E^*v + \bar\partial_E^*v\wedge\bar\partial_E^*v
  \\
  &= F_{\bar\partial_E+\tau} + \bar\partial_E\bar\partial_E^*v + \tau\wedge\bar\partial_E^*v + \bar\partial_E^*v\wedge\bar\partial_E^*v.
\end{align*}
Since $\tau \in \bH_{\bar\partial_E}^1 \subset \Ker\bar\partial_E\cap W^{1,p}(\Lambda^{0,1}(\fsl(E)))$, we have by \eqref{eq:F_dbar_E+alpha_Ker_dbar_E} that
\[
  \bar\partial_EF_{\bar\partial_E+\tau} = 0
\]
and so the Hodge decomposition \eqref{eq:Hodge_decomposition_Lp_Lambda02_slE} yields
\[
  F_{\bar\partial_E+\tau} \in \Ker\bar\partial_E\cap L^p(\Lambda^{0,2}(\fsl(E)))
  = \bH_{\bar\partial_E}^2 \oplus \bar\partial_E\left(W^{1,p}(\Lambda^{0,1}(\fsl(E)))\right).
\]
We wish to construct a finite-dimensional, local, complex analytic model space for all solutions in $S_{\bar\partial_E}$ near $\bar\partial_E$ to
\begin{equation}
  \label{eq:Taubes_1984_2-5_fundamental}
  F_{\bar\partial_E + \tau + \bar\partial_E^*v} = 0 \in L^p(\Lambda^{0,2}(\fsl(E))),
\end{equation}
or equivalently,
\begin{equation}
  \label{eq:Taubes_1984_2-5}
  \bar\partial_E\bar\partial_E^*v + \tau\wedge\bar\partial_E^*v + \bar\partial_E^*v\wedge\bar\partial_E^*v
  = -F_{\bar\partial_E+\tau} \in \bH_{\bar\partial_E}^2 \oplus \bar\partial_E\left(W^{1,p}(\Lambda^{0,1}(\fsl(E)))\right).
\end{equation}
Let $G_{\bar\partial_E}$ denote the Green's operator for the Laplace operator, 
\begin{equation}
  \label{eq:Box_dbar_E}
  \Box_{\bar\partial_E} := \bar\partial_E\bar\partial_E^* +  \bar\partial_E^*\bar\partial_E \quad\text{on } \Omega^{0,k}(\fsl(E)), \quad\text{for } k = 0,1,\ldots,n,
\end{equation}
when $X$ has complex dimension $n$, and thus
\[
  G_{\bar\partial_E}\Box_{\bar\partial_E} = \Box_{\bar\partial_E}G_{\bar\partial_E} = \id - \Pi_{\bar\partial_E} = \Pi_{\bar\partial_E}^\perp \quad\text{on } \Omega^{0,k}(\fsl(E)),
\]
where $\Pi_{\bar\partial_E}$ denotes $L^2$-orthogonal projection from $\Omega^{0,k}(\fsl(E))$ onto the harmonic subspace,
\[
  \bH_{\bar\partial_E}^k = \Ker\Box_{\bar\partial_E}\cap \Omega^{0,k}(\fsl(E))
  = \Ker(\bar\partial_E + \bar\partial_{\bar\partial_E}^*)\cap \Omega^{0,k}(\fsl(E)),
\]
as in \eqref{eq:bHdbarA0bullet}. The equation \eqref{eq:Taubes_1984_2-5}, and thus in turn \eqref{eq:Taubes_1984_2-5_fundamental}, may be rewritten as an equivalent pair of equations,
\begin{subequations}
  \label{eq:Taubes_1984_2-7_and_8}
  \begin{align}
    \label{eq:Taubes_1984_2-7}
  \Pi_{\bar\partial_E}^\perp\left(\bar\partial_E\bar\partial_E^*v + \tau\wedge\bar\partial_E^*v + \bar\partial_E^*v\wedge\bar\partial_E^*v + F_{\bar\partial_E+\tau}\right)
  &= 0 \in \bar\partial_E\left(W^{1,p}(\Lambda^{0,1}(\fsl(E)))\right),
    \\
    \label{eq:Taubes_1984_2-8}
  \Pi_{\bar\partial_E}\left(\bar\partial_E\bar\partial_E^*v
  + \tau\wedge\bar\partial_E^*v + \bar\partial_E^*v\wedge\bar\partial_E^*v
  + F_{\bar\partial_E+\tau}\right) &= 0 \in \bH_{\bar\partial_E}^2.
\end{align}
\end{subequations}
Because $\Pi_{\bar\partial_E}^\perp = \id$ on $\bar\partial_E(W^{1,p}(\Lambda^{0,1}(\fsl(E))))$, the quasilinear equation \eqref{eq:Taubes_1984_2-7} is equivalent to
\[
  \bar\partial_E\bar\partial_E^*v
  + \Pi_{\bar\partial_E}^\perp\left(\tau\wedge\bar\partial_E^*v + \bar\partial_E^*v\wedge\bar\partial_E^*v \right)
  =
  -\Pi_{\bar\partial_E}^\perp F_{\bar\partial_E+\tau}
  \in \bar\partial_E\left(W^{1,p}(\Lambda^{0,1}(\fsl(E)))\right).
\]
Since the restriction of the Green's operator,
\[
  G_{\bar\partial_E}:\bar\partial_E(W^{1,p}(\Lambda^{0,1}(\fsl(E))))
  \to \bar\partial_E(W^{3,p}(\Lambda^{0,1}(\fsl(E)))),
\]
is an isomorphism, the quasilinear equation \eqref{eq:Taubes_1984_2-7} is therefore equivalent to
\[
  G_{\bar\partial_E}\bar\partial_E\bar\partial_E^*v + G_{\bar\partial_E}\Pi_{\bar\partial_E}^\perp\left(\tau\wedge\bar\partial_E^*v + \bar\partial_E^*v\wedge\bar\partial_E^*v\right) = -G_{\bar\partial_E}\Pi_{\bar\partial_E}^\perp F_{\bar\partial_E+\tau}.
\]
Because $\Box_{\bar\partial_E}v = \bar\partial_E\bar\partial_E^*v$ for $v\in \bar\partial_E\left(W^{3,p}(\Lambda^{0,1}(\fsl(E)))\right)$ and $G_{\bar\partial_E}\Box_{\bar\partial_E} = \id$ on $\bar\partial_E(W^{3,p}(\Lambda^{0,1}(\fsl(E))))$, the quasilinear equation \eqref{eq:Taubes_1984_2-7} is thus equivalent to
\[
  v + G_{\bar\partial_E}\Pi_{\bar\partial_E}^\perp\left(\tau\wedge\bar\partial_E^*v + \bar\partial_E^*v\wedge\bar\partial_E^*v\right) = -G_{\bar\partial_E}\Pi_{\bar\partial_E}^\perp F_{\bar\partial_E+\tau}.
\]
Hence, for a small enough open neighborhood $N_{\bar\partial_E} \subset \bH_{\bar\partial_E}^1$ of the origin, the Implicit Mapping Theorem (see, for example, Feehan and Maridakis \cite[Theorem F.1, p.  127]{Feehan_Maridakis_Lojasiewicz-Simon_coupled_Yang-Mills}) yields a unique holomorphic map,
\begin{equation}
  \label{eq:bv_holomorphic_structures}
  \bv: N_{\bar\partial_E} \ni \tau \mapsto v(\tau) \in \bar\partial_E(W^{3,p}(\Lambda^{0,1}(\fsl(E)))),
\end{equation}
such that $\bv(0) = 0$ and equation \eqref{eq:Taubes_1984_2-7} holds for all $\tau \in N_{\bar\partial_E}$. We define a holomorphic map
\begin{equation}
  \label{eq:bpsi_holomorphic_structures}
  \bpsi: \bH_{\bar\partial_E}^1 \supset N_{\bar\partial_E}
   \ni \tau \mapsto \bar\partial_E^*v(\tau)
  \in \bar\partial_E^*\left(W^{2,p}(\Lambda^{0,2}(\fsl(E)))\right),
\end{equation}
and observe that $\bpsi(0) = 0$.

We can now complete the verification of Item \eqref{item:Kuranishi_model_dbarA_cH''(E)}. We define holomorphic maps,
\begin{subequations}
  \label{eq:bgamma_and_bchi_holomorphic_structures}
  \begin{align}
  \label{eq:bgamma_holomorphic_structures}
    \bgamma: \bH_{\bar\partial_E}^1 \supset N_{\bar\partial_E} \ni \tau
    \mapsto \bar\partial_E + \tau + \bpsi(\tau) \in S_{\bar\partial_E},
    \\
    \label{eq:bchi_holomorphic_structures}
    \bchi: \bH_{\bar\partial}^1 \supset N_{\bar\partial_E} \ni \tau
    \mapsto \Pi_{\bar\partial_E}F_{\bar\partial_E + \tau + \bpsi(\tau)}
    \in \bH_{\bar\partial_E}^2.
  \end{align}
\end{subequations}
Clearly, $\bar\partial_E \in \bgamma(\bchi^{-1}(0)\cap N_{\bar\partial_E})$ since $F_{\bar\partial_E} = 0$, so $\bchi(0) = F_{\bar\partial_E} = 0$ and $0 \in \bchi^{-1}(0)\cap N_{\bar\partial_E}$ and $\bgamma(0) = \bar\partial_E$. It remains to verify that $D\bchi(0) = 0$ (compare Friedman and Morgan \cite[Section 4.1.5, Proposition 1.20, p. 299]{FrM}). Because $\Pi_{\bar\partial_E}$ is $L^2$-orthogonal projection from $L^p(\Lambda^{0,2}(\fsl(E)))$ and $\bar\partial_E\bar\partial_E^*v \in \bar\partial_E(W^{1,p}(\Lambda^{0,1}(\fsl(E)))$ and $F_{\bar\partial_E} = \tau\wedge\tau$, equation \eqref{eq:Taubes_1984_2-8} and the definitions \eqref{eq:bpsi_holomorphic_structures} and \eqref{eq:bchi_holomorphic_structures} of $\bpsi$ and $\bchi$ yield
\begin{equation}
  \label{eq:bchi_quadratic_form}
  \bchi(\tau) = \Pi_{\bar\partial_E}\left(\tau\wedge\bpsi(\tau) + \bpsi(\tau)\wedge\bpsi(\tau) + \tau\wedge\tau
 \right), \quad\text{for all } \tau \in N_{\bar\partial_E}.
\end{equation}
Since $\bpsi(0) = 0$, the Taylor expansion of $\bchi(\tau)$ has leading order term proportional to $\Pi_{\bar\partial_E}(\tau\wedge\tau)$ and thus $D\bchi(0) = 0$, as claimed.

We claim that, after possibly shrinking the open neighborhood $N_{\bar\partial_E} \subset \bH_{\bar\partial_E}^1$ of the origin and the open neighborhood $U_{\bar\partial_E} \subset S_{\bar\partial_E}$ of $\bar\partial_E$ in the slice \eqref{eq:dbar_E_slice}, the Implicit Mapping Theorem (see, for example, Feehan and Maridakis \cite[Theorem F.1, p.  127]{Feehan_Maridakis_Lojasiewicz-Simon_coupled_Yang-Mills}) yields the equality \eqref{eq:Kuranishi_model_dbarA_cH''(E)}. To see this, suppose first that $\bar\partial_E + \tau + \bpsi(\tau) \in \bgamma(\bchi^{-1}(0)\cap N_{\bar\partial_E})$. The definitions \eqref{eq:bv_holomorphic_structures} and \eqref{eq:bpsi_holomorphic_structures} of $v(\tau)$ and $\bpsi(\tau) = \bar\partial_E^*v(\tau)$ ensure that equation \eqref{eq:Taubes_1984_2-7} holds for all $\tau \in N_{\bar\partial_E}$, while the definition of \eqref{eq:bchi_holomorphic_structures} of $\bchi$ and assumption that $\bchi(\tau) = 0$ ensure that equation \eqref{eq:Taubes_1984_2-8} holds for all $\tau \in N_{\bar\partial_E}$. The pair of equations \eqref{eq:Taubes_1984_2-7_and_8} is equivalent to \eqref{eq:Taubes_1984_2-5} and that in turn is equivalent to \eqref{eq:Taubes_1984_2-5_fundamental}, so that $F_{\bar\partial_E + \tau + \bpsi(\tau)} = 0$ and thus $\bar\partial_E + \tau + \bpsi(\tau) \in \cH(E)\cap U_{\bar\partial_E}$ by definition \eqref{eq:Subset_01_connections_defining_holomorphic_structure_on_E} of $\cH(E)$ and the definition of the slice neighborhood $U_{\bar\partial_E}$ preceding \eqref{eq:dbar_E_slice_biholomorphic_map}. This yields the inclusion
\[
  \bgamma(\bchi^{-1}(0)\cap N_{\bar\partial_E}) \subset \cH(E)\cap U_{\bar\partial_E}.
\]
Conversely, suppose that $\bar\partial_E + \alpha \in \cH(E)\cap U_{\bar\partial_E}$. The definition of the slice neighborhood $U_{\bar\partial_E}$ and the Hodge decomposition \eqref{eq:Hodge_decomposition_W1p_Lambda01_slE} ensure that $\alpha = \tau + \bar\partial_E^*v$, while the definition \eqref{eq:Subset_01_connections_defining_holomorphic_structure_on_E} of $\cH(E)$ implies that $F_{\bar\partial_E + \tau + \bar\partial_E^*v} = 0$. The equivalences of equations just discussed ensure that $\bgamma(\tau) = \bar\partial_E + \tau + \bar\partial_E^*v$ obeys $\bchi(\tau) = 0$ and this yields the inclusion
\[
  \cH(E)\cap U_{\bar\partial_E} \subset \bgamma(\bchi^{-1}(0)\cap N_{\bar\partial_E}),
\]
and hence the equality \eqref{eq:Kuranishi_model_dbarA_cH''(E)}.

We now verify that Item \eqref{item:cM_dbarE^vir_is_bgamma(N_dbarE)} holds. Thus, we claim that
\begin{inparaenum}[\itshape i\upshape)]
\item $\cH_{\bar\partial_E}^\vir(E)\cap U_{\bar\partial_E}$ defined by \eqref{eq:Atiyah_Hitchin_Singer_family_page_446_holomorphic_structures} is a complex embedded submanifold of $U_{\bar\partial_E}$ whose tangent space at $\bar\partial_E$ is equal to $\bH_{\bar\partial_E}^1$, and
\item $\bgamma(N_{\bar\partial_E}) \subset \sA^{0,1}(E)$ is an embedded complex submanifold of $U_{\bar\partial_E}$ whose tangent space at $\bar\partial_E$ is equal to $\bH_{\bar\partial_E}^1$, and
that
\item $\bgamma(N_{\bar\partial_E}) = \cH_{\bar\partial_E}^\vir(E)\cap U_{\bar\partial_E}$, after possibly shrinking the open neighborhoods $N_{\bar\partial_E}$ and $U_{\bar\partial_E}$.
\end{inparaenum}
By definition \eqref{eq:cF_sA01_to_Ran_dbar_A_cap_LpLambda02_slE} of $\cF$, we have
\[
  \cH_{\bar\partial_E}^\vir(E) = \cF^{-1}(0)\cap S_{\bar\partial_E}.
\]
Hence, after possibly shrinking $U_{\bar\partial_E}$ and noting that $\cF$ in \eqref{eq:cF_sA01_to_Ran_dbar_A_cap_LpLambda02_slE} is a submersion by virtue of the isomorphism \eqref{eq:Differential_cF_Ran_dbarA*_in_W1p_to_Ran_dbar_A_in_Lp_isomorphism}, we see that $\cH_{\bar\partial_E}^\vir(E)\cap U_{\bar\partial_E}$ is indeed a complex embedded submanifold of $U_{\bar\partial_E}$ and that its tangent space at $\bar\partial_E$ is
\begin{multline*}
  T_{\bar\partial_E}\left(\cH_{\bar\partial_E}^\vir(E)\cap U_{\bar\partial_E}\right)
  = \Ker D\cF(\bar\partial_E)
  = \Ker\bar\partial_E\cap T_{\bar\partial_E}S_{\bar\partial_E}
  \\
  = \Ker(\bar\partial_E + \bar\partial_E^*)\cap W^{1,p}(\Lambda^{0,1}(\fsl(E)))
  = \bH_{\bar\partial_E}^1.
\end{multline*}
This proves the first claim. On the other hand, because $\bpsi(\tau) \in \bar\partial_E^*(W^{2,p}(\Lambda^{0,2}(\fsl(E))))$ for all $\tau \in N_{\bar\partial_E} \subset \bH_{\bar\partial_E}^1$, the definition \eqref{eq:bgamma_holomorphic_structures} of $\bgamma$ implies that
\[
  \bgamma(N_{\bar\partial_E})
  = \left\{\bar\partial_E + \tau + \bpsi(\tau): \tau \in N_{\bar\partial_E}\right\}
  \subset \sA^{0,1}(E).
\]
is the graph of $\bpsi$ in $\bar\partial_E + \bH_{\bar\partial_E}^1 \oplus \bar\partial_E^*(W^{2,p}(\Lambda^{0,2}(\fsl(E))))$ over $N_{\bar\partial_E}$. Therefore, $\bgamma(N_{\bar\partial_E})$ is a complex embedded submanifold of $U_{\bar\partial_E}$ whose tangent space at $\bar\partial_E$ is
\[
  T_{\bar\partial_E}\bgamma(N_{\bar\partial_E}) = \bH_{\bar\partial_E}^1.
\]
This proves the second claim. Because $\bgamma(N_{\bar\partial_E}) \subset \cH_{\bar\partial_E}^\vir(E)\cap U_{\bar\partial_E}$ we obtain, after possibly shrinking $N_{\bar\partial_E}$ or $U_{\bar\partial_E}$, the inclusion
\[
  \bgamma(N_{\bar\partial_E}) = \cH_{\bar\partial_E}^\vir(E)\cap U_{\bar\partial_E},
\]
which proves the third claim. This completes our verification of Item \eqref{item:cM_dbarE^vir_is_bgamma(N_dbarE)}.

Lastly, we prove the final assertion of Theorem \ref{thm:Local_Kuranishi_model_for_simple_point_cH(E)}. We thus restrict our attention again to points $\bar\partial_E \in \cH(E)$ that are both \emph{regular} in the sense that $\bH_{\bar\partial_E}^2 = (0)$ and \emph{simple} in the sense of Definition \ref{defn:Simple_01-connection}, which is equivalent to $\bH_{\bar\partial_E}^0 = (0)$ by Lemma \ref{lem:Simple_01-connection}. Hence, we restrict our attention to the open subset $\cH_\reg^*(E) \subset \cH(E)$ in \eqref{eq:cHreg*E}. Because $\bH_{\bar\partial_E}^2 = (0)$, the holomorphic maps $\bpsi$ in \eqref{eq:bpsi_holomorphic_structures} and $\bchi$ in \eqref{eq:bchi_holomorphic_structures} are identically zero. It follows that for every point $\bar\partial_E \in \cH_\reg^*(E)$, there is an open neighborhood $\UU_{\bar\partial_E} \subset \sA^{0,1}(E)$ such that $\cH_\reg^*(E)\cap \UU_{\bar\partial_E}$ is biholomorphic to $N_{\bar\partial_E} \times U_{\id_E}$, where $U_{\id_E}$ is the open neighborhood in $W^{2,p}(\SL(E))$ of $\id_E$ appearing in \eqref{eq:dbar_E_slice_biholomorphic_map} and $N_{\bar\partial_E}$ is the open neighborhood in $\bH_{\bar\partial_E}^1$ of the origin appearing in \eqref{eq:bpsi_holomorphic_structures} and \eqref{eq:bgamma_holomorphic_structures}. Hence, $\cH_\reg^*(E)$ is an embedded complex submanifold of $\sA^{0,1}(E)$. This completes the proof of Theorem \ref{thm:Local_Kuranishi_model_for_simple_point_cH(E)}.
\end{proof}

\begin{rmk}[On systems of equations defining local complex analytic model spaces for $\cH(E)$]
\label{rmk:On_systems_equations_defining_local_complex_analytic_model_space_H''(E)}  
Because of the Hodge decomposition for \eqref{eq:Hodge_decomposition_Lp_Lambda02_slE} for $L^p(\Lambda^{0,2}(\fsl(E)))$, the single equation $F_{\bar\partial_E + \alpha} = 0$ for $\bar\partial_E+\alpha \in \UU_{\bar\partial_E}$ is equivalent to a system of three equations,
\begin{align*}
  \Pi_{\Ran\bar\partial_E}F_{\bar\partial_E + \alpha} &= 0 \in \bar\partial_E\left(W^{1,p}(\Lambda^{0,1}(\fsl(E)))\right),
  \\                                                      
  \Pi_{\Ran\bar\partial_E^*}F_{\bar\partial_E + \alpha} &= 0 \in \bar\partial_E^*\left(W^{1,p}(\Lambda^{0,3}(\fsl(E)))\right),
  \\                                                      
  \Pi_{\bar\partial_E}F_{\bar\partial_E + \alpha} &= 0 \in \bH_{\bar\partial_E}^2,                                                       
\end{align*}
and their solution is the approach taken by Kobayashi \cite[Equation (7.3.8), p. 232]{Kobayashi_differential_geometry_complex_vector_bundles}. This equivalence continues to hold, of course, upon restriction to the slice neighborhood $U_{\bar\partial_E}$. However, as we have seen, the preceding system is equivalent to the pair of equations \eqref{eq:Taubes_1984_2-7_and_8} corresponding to the $L^2$-orthogonal decomposition
\[
  \Ker\bar\partial_E\cap L^p(\Lambda^{0,2}(\fsl(E)))
  =
  \bH_{\bar\partial_E}^2 \oplus \bar\partial_E\left(W^{1,p}(\Lambda^{0,1}(\fsl(E)))\right)
\]
implied by \eqref{eq:Hodge_decomposition_Lp_Lambda02_slE}.
\end{rmk}

\subsection[Application of Marsden--Weinstein theorem to construction of K\"ahler metric]{Application of Marsden--Weinstein symplectic reduction theorem to construction of K\"ahler metric on moduli space of holomorphic structures}
\label{subsec:Application_Marsden-Weinstein_theorem_construction_Kaehler_metric_holomorphic_structures}
In preparation for our forthcoming application of the Marsden--Weinstein Symplectic Reduction Theorem \ref{thm:Marsden-Weinstein_symplectic_quotient}, we shall choose
\begin{equation}
  \label{eq:Marsden-Weinstein_V_space_01connections}
  \bV'' := \cH_\reg^*(E)
\end{equation}
in Hypothesis \ref{hyp:Marsden-Weinstein_symplectic_quotient_conditions}, where $\cH_\reg^*(E)$ is as in \eqref{eq:cHreg*E} and defined under the hypotheses of Theorem \ref{thm:Local_Kuranishi_model_for_simple_point_cH(E)}. Our choice of $\bV''$ is different from that of Kobayashi in \cite[Section 7.6, p. 253, line 20]{Kobayashi_differential_geometry_complex_vector_bundles}, as we explain in the following

\begin{rmk}[On the choice of $\bV''$]
Instead of \eqref{eq:Marsden-Weinstein_V_space_01connections}, Kobayashi in \cite[Section 7.6, p. 253, line 20]{Kobayashi_differential_geometry_complex_vector_bundles} chooses
\[
  \bV'' = \left\{v^*\bar\partial_E: \bar\partial_E \in \cH(E) \text{ with } \bH_{\bar\partial_E}^0 = (0) \text{ and } \bH_{\bar\partial_E}^2 = (0) \text{ and } v \in W^{2,p}(\SL(E))\right\}.
\]
There is surely a typographical error in \cite{Kobayashi_differential_geometry_complex_vector_bundles} because the conditions $\bar\partial_E \in \cH(E)$ and $\bH_{\bar\partial_E}^0 = (0)$ and $\bH_{\bar\partial_E}^2 = (0)$ are both invariant under the action of $W^{2,p}(\SL(E))$, so the presence of complex gauge transformations $v$ in the preceding definition of $\bV''$ is redundant. 
\end{rmk}

For the moduli subspace $\cM_\reg^*(E)$ in \eqref{eq:Moduli_space_regular_simple_holomorphic_structures}, one has an identification of sets:
\[
  \cM_\reg^*(E) \subset \bV''/W^{2,p}(\SL(E)).
\]
Suppose now that $(X,g,J)$ is a complex K\"ahler manifold with K\"ahler form $\omega = g(\cdot,J\cdot)$ and that the fixed unitary connection $A_d$ on the Hermitian line bundle $\det E$ has curvature $F_{A_d}$ of type $(1,1)$, so $F_{A_d}^{0,2}=0$. Continuing our analogy with Kobayashi, we follow \cite[Section, p. 253, line 29]{Kobayashi_differential_geometry_complex_vector_bundles} and observe that the K\"ahler metric $\bh$ on $\sA^{0,1}(E)$ given by \eqref{eq:Kobayashi_7-6-20} induces a K\"ahler metric $\bh|_{\bV''} := \bh\restriction\bV''$ on $\bV''$, which we denote simply by $\bh$.

For $\bV''$ as in \eqref{eq:Marsden-Weinstein_V_space_01connections}, we define
\begin{equation}
  \label{eq:Marsden-Weinstein_V_space_unitary_connections}
  \bV := \left(\pi_h^{0,1}\right)^{-1}(\bV'') = \left(\pi_h^{0,1}\right)^{-1}\left(\cH_\reg^*(E)\right),
\end{equation}
so that $\bV$ is an embedded submanifold of the real Banach affine space $\sA(E,h)$ as well via the restriction of the isomorphism \eqref{eq:Bijection_unitaryconnections_with_01connections} of real Banach affine spaces,
\[
  \pi_h^{0,1}:\sA(E,h) \ni A \mapsto \bar\partial_A \in \sA^{0,1}(E),
\]
to a bijection of real analytic subspaces,
\begin{equation}
  \label{eq:Bijection_unitary1-1connections_with_integrable01connections}
  \pi_h^{0,1}:\sH(E,h) \ni A \mapsto \bar\partial_A \in \cH(E),
\end{equation}
where $\sH(E,h)$ is as in \eqref{eq:Subset_unitary_connections_defining_holomorphic_structure_on_E} and $\cH(E)$ is as in \eqref{eq:Subset_01_connections_defining_holomorphic_structure_on_E}. The K\"ahler metric $\bh$ on $\sA^{0,1}(E)$ induces a K\"ahler metric $\bh|_{\bV''} := \bh\restriction\bV''$ on $\bV''$, which we denote simply by $\bh$. We shall write the K\"ahler form $\bomega|_{\bV''} := \bomega\restriction\bV''$ induced by $\bomega$ simply as $\bomega$. Similarly, we may regard $\bV$ as a symplectic manifold and write $\bomega$ for its symplectic form instead of $(\pi_h^{0,1})^*(\bomega|_{\bV''})$.

It is convenient to define
\begin{equation}
  \label{eq:tilde_sG_E}
  \sG_E := W^{2,p}(\SU(E))
  \quad\text{and}\quad
  \tilde\sG_E := W^{2,p}(\SU(E))/C_r.
\end{equation}
We recall that the group $C_r$ of $r$-th roots of unity is equal to the center of $\SU(r)$, where $E$ has complex rank $r$, and acts trivially on $\sA(E,h)$ via the identification $C_r\,\id_E \subset \sG_E$. We note that
\[
  \fg_E := W^{2,p}(\su(E))
\]
is the Lie algebra of both $\tilde\sG_E = W^{2,p}(\SU(E))/C_r$ and $\sG_E = W^{2,p}(\SU(E))$ in \eqref{eq:tilde_sG_E}.

The following lemma is a special case of the forthcoming Lemma \ref{lem:Regular_points_bmu_on_restriction_from_sA(E,h)_times_W1p(E)_to_sP(E,h)}, which applies more generally to projective vortices and which we prove in detail in that context.

\begin{lem}[Projective Hermitian--Einstein connections with vanishing zero-order cohomology groups remain regular points of the moment map upon restriction from $\sA(E,h)$ to $\sH(E,h)$]
\label{lem:Regular_points_bmu_on_restriction_from_sA(E,h)_to_sH(E,h)}
Continue the hypotheses of Proposition \ref{prop:Donaldson_Kronheimer_6-5-7_almost_Hermitian} and assume in addition that the almost complex structure $J$ is integrable, so $(X,g,J)$ is a complex K\"ahler manifold. If $A \in \bmu^{-1}(0)\cap\sH(E,h)$, then the following hold:
\begin{enumerate}
\item $A$ is a regular point of the restriction of the moment map $\bmu$ in \eqref{eq:Moment_map_action_unitary_det_one_gauge_transformations_affine_space_unitary_connections} from $\sA(E,h)$ to the Banach analytic subvariety $\sH(E,h)$,
\begin{equation}
  \label{eq:Moment_map_action_unitary_det_one_gauge_transformations_sH(E,h)}
  \bmu:\sH(E,h) \to \left(W^{2,p}(\su(E))\right)^*,
\end{equation}
in the sense that the differential $d\bmu(A)$ is an epimorphism from the Zariski tangent space $T_A\sH(E,h)$ onto $(W^{2,p}(\su(E)))^*$, if and only if $\bH_A^0 = (0)$, where the harmonic space $\bH_A^0$ is as in \eqref{eq:HE_equation_bHA0}.

\item If $\bH_A^0 = (0)$, then an open neighborhood of $A$ in $\bmu^{-1}(0)\cap\bV$ is an embedded, real analytic submanifold of $\bV$ in \eqref{eq:Marsden-Weinstein_V_space_unitary_connections}.
\end{enumerate}
\end{lem}

We conclude from Lemma \ref{lem:Regular_points_bmu_on_restriction_from_sA(E,h)_to_sH(E,h)} that zero is a regular value of the restriction,
\begin{equation}
  \label{eq:Restriction_of_MomentMap_To_bV_cap_sH*(E,h)}
  \bmu: \bV \to \left(W^{2,p}(\su(E))\right)^*,
\end{equation}
of the moment map $\bmu$ in \eqref{eq:Moment_map_action_unitary_det_one_gauge_transformations_affine_space_unitary_connections}. Indeed, this holds because if $A \in \bmu^{-1}(0)\cap\bV$ for $\bV = (\pi_h^{0,1})^{-1}(\bV'')$ as in \eqref{eq:Marsden-Weinstein_V_space_unitary_connections}, then $A$ is a projectively Hermitian--Einstein connection by the identity \eqref{eq:Kobayashi_7-6-33}. Hence, Theorem \ref{thm:Kobayashi_7-2-21} implies that $\bH_A^0 = (0)$ since $\bH_{\bar\partial_A}^0 = \bH_A^0\otimes_\RR\CC$ and $\bH_{\bar\partial_A}^0 = (0)$ by the definition \eqref{eq:Marsden-Weinstein_V_space_holomorphic_pairs} of $\bV'' = \cH_\reg^*(E)$ and the definition \eqref{eq:cHreg*E} of $\cH_\reg^*(E)$.

Therefore, $0 \in \fg^* := (W^{2,p}(\su(E)))^*$ is a \emph{weakly regular value} of $\bmu$ in \eqref{eq:Restriction_of_MomentMap_To_bV_cap_sH*(E,h)} in the sense of Condition \eqref{item:Marsden-Weinstein_symplectic_quotient_conditions_weakly_regular_value} in Hypothesis \ref{hyp:Marsden-Weinstein_symplectic_quotient_conditions}. Thus, $\bmu^{-1}(0)\cap\bV$ is an embedded submanifold of $\bV$ and for every $A \in \bmu^{-1}(0)\cap\bV$, one has the equality
\[
  T_A\left(\bmu^{-1}(0)\cap \bV \right) = \Ker d\bmu(A).
\]
Thus, Condition \eqref{item:Marsden-Weinstein_symplectic_quotient_conditions_weakly_regular_value} in Hypothesis \ref{hyp:Marsden-Weinstein_symplectic_quotient_conditions} is obeyed with $V$ replaced by $\bV$. Proposition \ref{prop:Donaldson_Kronheimer_6-5-8_almost_Kaehler} verifies that Condition \eqref{item:Marsden-Weinstein_symplectic_quotient_conditions_equivariance} in Hypothesis \ref{hyp:Marsden-Weinstein_symplectic_quotient_conditions} is obeyed with $V$ and $G$ replaced by $\bV$ and $\tilde\sG_E$, respectively.

Each point $A\in \bmu^{-1}(0)\cap \bV$ in the zero set has $\bH_A^0 = (0)$ as explained above. Because $\bH_A^0 = (0)$, the connection $A$ has central stabilizer in $\sG_E$ by Lemma \ref{lem:Kronheimer_Lemma2.2} and thus trivial stabilizer in $\tilde\sG_E$. Hence, the Banach Lie group $\tilde\sG_E$ acts freely on
$\bmu^{-1}(0)\cap \bV$. The Lie group $\tilde\sG_E$ acts properly on $\bmu^{-1}(0)\cap \bV$ since the actions of $W^{2,p}(\SU(E))$ on $\sA(E,h)$ or $\sH(E,h)$ are proper by Kobayashi \cite[Proposition 7.1.14, p. 221]{Kobayashi_differential_geometry_complex_vector_bundles}. (Although Kobayashi considers $\U(E)$ and not $\SU(E)$ gauge transformations in \cite[Proposition 7.1.14, p. 221]{Kobayashi_differential_geometry_complex_vector_bundles}, because $W^{2,p}(\SU(E))$ is a subgroup of $W^{2,p}(\U(E))$, the group acts properly on $\bmu^{-1}(0)\cap \bV$ since Kobayashi shows that this is true for $W^{2,p}(\U(E))$.)

According to Feehan and Maridakis \cite[Corollary 18, p. 19]{Feehan_Maridakis_Lojasiewicz-Simon_coupled_Yang-Mills}, there is a (Coulomb gauge) slice for the action of $\sG_E$ (and hence the action of $\tilde\sG_E$) on $\bmu^{-1}(0)\cap\bV$ in the sense of Condition \eqref{item:Marsden-Weinstein_symplectic_quotient_conditions_slice} with $V$ and $G$ replaced by $\bV$ and $\tilde\sG_E$, respectively. Therefore, all the requirements of Condition \eqref{item:Marsden-Weinstein_symplectic_quotient_conditions_slice} in Hypothesis \ref{hyp:Marsden-Weinstein_symplectic_quotient_conditions} are obeyed. We thus obtain the following modification and enhancement of \cite[Theorem  7.6.36, p. 255]{Kobayashi_differential_geometry_complex_vector_bundles} due to Kobayashi.

\begin{thm}[Moduli space of non-split, regular projectively Hermitian--Einstein connections is a complex K\"ahler manifold]
\label{thm:Kobayashi_7-6-36}
Let $(E,h)$ be a smooth Hermitian vector bundle over a closed, complex K\"ahler manifold $(X,g,J)$ with K\"ahler form $\omega = g(\cdot,J\cdot)$ and complex dimension $n$. Let $p \in (n,\infty)$ be a constant and $A_d$ be a fixed unitary connection on the Hermitian line bundle $\det E$ with curvature obeying $F_{A_d}^{0,2}=0$. Then the moduli space $M_\reg^*(E,h,\omega)$ in \eqref{eq:Moduli_space_HE_connections_irreducible_regular} is a complex K\"ahler manifold with K\"ahler form induced from $\bomega$ in \eqref{eq:Kobayashi_7-6-22}.
\end{thm}

\begin{proof}
According to Theorem \ref{thm:Kobayashi_7_4_20}, the moduli subspace $\cM(E,\omega)\subset \cM^*(E)$ in \eqref{eq:Moduli_space_omega-stable_holomorphic_structures} and \eqref{eq:Moduli_space_simple_holomorphic_structures}, respectively, is an open subspace that is the image of an embedding in the sense of real analytic spaces of the moduli space $M^*(E,h,\omega)$ in \eqref{eq:Moduli_space_HE_connections_irreducible}. By Theorem \ref{thm:Kobayashi_7_3_34_simple} \eqref{item:cM*(E)_is_complex_analytic_space}, the moduli space $\cM^*(E)$ is a complex analytic space in the sense of Grauert and Remmert \cite[Section 1.1.5, p. 7]{Grauert_Remmert_coherent_analytic_sheaves} and so the open subspace $\cM(E,\omega)$ is also a complex analytic space. We claim that the open subset of smooth points,
\begin{equation}
  \label{eq:Moduli_space_omega-stable_holomorphic_structures_regular}
  \cM_\reg(E,\omega)
  :=
  \left\{[\bar\partial_A] \in \cM(E,\omega): \bH_{\bar\partial_A}^2 = (0) \right\},
\end{equation}
is a complex manifold. Indeed, this follows from the facts that $M_\reg^*(E,h,\omega)$ embeds onto an open
subspace\footnote{This openness result is proved by Itoh as \cite[Proposition 4.3, p. 857]{Itoh_1985} in the restricted setting of anti-self-dual connections $A$ on rank-two Hermitian vector bundles $E$ with $c_1(E)=0$ over complex K\"ahler surfaces, but his method of proof is quite general and also applies in the present setting.} of $\cM_\reg^*(E)$ in \eqref{eq:Moduli_space_regular_simple_holomorphic_structures} that is identified with $\cM_\reg(E,\omega)$ by Theorem \ref{thm:Kobayashi_7_4_20} \eqref{item:Kobayashi_7_4_20_analytic_embedding}, and $\cM_\reg^*(E)$ is a complex manifold by Theorem \ref{thm:Kobayashi_7_3_34_simple} \eqref{item:cM*(E)_is_complex_manifold}.

For each point $[A]\in M^*(E,h,\omega)$ in \eqref{eq:Moduli_space_HE_connections_irreducible}, the superscript ``$*$'' indicates that the unitary connection $A$ is non-split, so Corollary \ref{cor:Split_unitary_A_and_Lie_Alg_of_Stab(A)} \eqref{item:A_smooth_and_bHA_non-zero_implies_A_split} implies that $\bH_A^0=(0)$. In addition, for each point $[A]\in M(E,h,\omega)$, Theorem \ref{thm:Kobayashi_7-2-21} ensures, when $\bH_A^0 = (0)$, that $\bH_A^2 = (0)$ if and only if  $\bH_{\bar\partial_A}^2 = (0)$. Theorem \ref{thm:Kobayashi_7_4_19} \eqref{eq:M*(E,h,omega)_is_real_analytic_space} asserts that $M_\reg^*(E,h,\omega)$ in \eqref{eq:Moduli_space_HE_connections_irreducible_regular} is a real analytic manifold. By combining the preceding observations, we see that Theorem \ref{thm:Kobayashi_7_4_20} induces a real analytic embedding of the real analytic manifold $M_\reg^*(E,h,\omega)$ onto the complex manifold $\cM_\reg(E,\omega)$. The moduli space $M_\reg^*(E,h,\omega)$ thus inherits a complex structure from $\cM_\reg(E,\omega)$ via the preceding isomorphism of real analytic manifolds.

We claim that we have the following identification of sets:
\begin{equation}
  \label{eq:Mreg*Ehomega_equals_symplectic_quotient}
  M_\reg^*(E,h,\omega) = \left.\left(\bmu^{-1}(0)\cap \bV\right)\right/W^{2,p}(\SU(E)).
\end{equation}
To verify \eqref{eq:Mreg*Ehomega_equals_symplectic_quotient}, first let $[A]\in M_\reg^*(E,h,\omega)$. From the definition of $M_\reg^*(E,h,\omega)$ in \eqref{eq:Moduli_space_HE_connections_irreducible_regular}, the connection $A$ is non-split and thus $\bH_A^0 = (0)$ by Corollary \ref{cor:Split_unitary_A_and_Lie_Alg_of_Stab(A)} \eqref{item:A_smooth_and_bHA_non-zero_implies_A_split} as noted above, and moreover, $\bH_A^2 = (0)$. Since $A$ is projectively Hermitian--Einstein, then $(F_A^{0,2})_0 = 0$ and thus $F_A^{0,2} = 0$ by hypothesis that $F_{A_d}^{0,2} = 0$. Hence, $F_{\bar\partial_A}=0$, so $\bar\partial_A \in \cH(E)$, where $\cH(E)$ is as in \eqref{eq:Subset_01_connections_defining_holomorphic_structure_on_E}. Since $A$ is projectively Hermitian--Einstein, Theorem \ref{thm:Kobayashi_7-2-21} ensures that $\bH_{\bar\partial_A}^0 = \bH_A^0\otimes_\RR\CC = (0)$ and $\bH_{\bar\partial_A}^2 = (0)$ since $\bH_A^2 = (0)$. Therefore, $\bar\partial_A \in \cH_\reg^*(E)$ by definition \eqref{eq:cHreg*E} of $\cH_\reg^*(E)$ and so $A \in \bV$ by the definition \eqref{eq:Marsden-Weinstein_V_space_unitary_connections} of $\bV = (\pi_h^{0,1})^{-1}\left(\cH_\reg^*(E)\right)$ as an embedded submanifold of $\sA(E,h)$. Furthermore, because $A$ is projectively Hermitian--Einstein, the identity \eqref{eq:HE_connections_as_zero_set_moment_map_on_1-1_unitary_connections} ensures that $A \in \bmu^{-1}(0)\cap\sA(E,h)$ and hence we have the inclusion of sets,
\[
  M_\reg^*(E,h,\omega) \subset \left.\left(\bmu^{-1}(0)\cap \bV\right)\right/W^{2,p}(\SU(E)).
\]
Conversely, suppose that $A \in \bmu^{-1}(0)\cap \bV$. By the definition \eqref{eq:Marsden-Weinstein_V_space_unitary_connections} of $\bV = (\pi_h^{0,1})^{-1}\left(\cH_\reg^*(E)\right)$ and definition \eqref{eq:cHreg*E} of $\cH_\reg^*(E)$ and the identity \eqref{eq:HE_connections_as_zero_set_moment_map_on_1-1_unitary_connections}, we see that $A$ is projectively Hermitian--Einstein. Furthermore, by the definition \eqref{eq:Marsden-Weinstein_V_space_unitary_connections} of $\bV = (\pi_h^{0,1})^{-1}\left(\cH_\reg^*(E)\right)$, and the definition \eqref{eq:cHreg*E} of $\cH_\reg^*(E)$, we have $\bH_{\bar\partial_A}^0 = (0)$ and $\bH_{\bar\partial_A}^2 = (0)$. Hence, Theorem \ref{thm:Kobayashi_7-2-21} ensures that $\bH_A^0 = (0)$, since $\bH_{\bar\partial_A}^0 = (0)$ and $\bH_{\bar\partial_A}^0 = \bH_A^0\otimes_\RR\CC$, and thus also ensures that $\bH_A^2 = (0)$. Because $\bH_A^0 = (0)$, the connection $A$ must be non-split by Corollary \ref{cor:Split_unitary_A_and_Lie_Alg_of_Stab(A)} \eqref{item:A_smooth_and_bHA_non-zero_implies_A_split}. Therefore, $[A] \in M_\reg^*(E,h,\omega)$ by definition \eqref{eq:Moduli_space_HE_connections_irreducible_regular} of $M_\reg^*(E,h,\omega)$ and hence we have the reverse inclusion of sets,
\[
  \left.\left(\bmu^{-1}(0)\cap \bV\right)\right/W^{2,p}(\SU(E)) \subset M_\reg^*(E,h,\omega),
\]
and this proves the claimed equality \eqref{eq:Mreg*Ehomega_equals_symplectic_quotient}. 

The fact that $\bomega$ in \eqref{eq:Kobayashi_7-6-22} or equivalently \eqref{eq:Kobayashi_7-6-22_bomega_as_fundamental_2-form_for_bg_and_bJ} induces a K\"ahler form on the quotient
\[
  \left(\bmu^{-1}(0)\cap \bV\right)/W^{2,p}(\SU(E))
\]
follows from the Marsden--Weinstein Symplectic Reduction Theorem \ref{thm:Marsden-Weinstein_symplectic_quotient} (with $V$, $G$, and $\omega_V$ replaced by $\bV$,  $\tilde\sG_E$, and $\bomega|_\bV := \bomega\restriction\bV$, respectively) and our verification of its hypotheses in the preceding paragraphs. In particular, we obtain a K\"ahler form on $M_\reg^*(E,h,\omega)$ compatible with the induced complex structure and Riemannian metric on $M_\reg^*(E,h,\omega)$.
\end{proof}

\begin{rmk}[Direct proof of the K\"ahler property for the moduli space of non-split, regular, projectively Hermitian--Einstein connections]
\label{rmk:Kaehler_moduli_space_irreducible_regular_projectively_Hermitian-Einstein_connections}
When $(X,g,J)$ is a closed, complex K\"ahler surface and $(E,h)$ is a Hermitian vector bundle over $X$ with trivial determinant line bundle $\det E$, Itoh \cite[Lemma 4.1, p. 21 and Proposition 4.2, p. 22]{Itoh_1988} provides a direct proof that the $L^2$ metric on the corresponding moduli space $M_\reg^*(E,h,\omega)$ of non-split, regular, $g$-anti-self-dual $W^{1,p}$ connections on $E$ modulo $W^{2,p}(\SU(E))$ is K\"ahler,  without appealing to the Marsden--Weinstein Symplectic Reduction Theorem \ref{thm:Marsden-Weinstein_symplectic_quotient}. His method should generalize to prove the K\"ahler property for the $L^2$ metric in Theorem \ref{thm:Kobayashi_7-6-36}, as Kobayashi implies in \cite[Section 7.6, p. 253, line 14]{Kobayashi_differential_geometry_complex_vector_bundles}.
\end{rmk}

\chapter[K\"ahler structure on moduli space of projective vortices]{Complex structure and K\"ahler metric on the moduli space of projective vortices over a complex K\"ahler manifold}
\label{chap:Complex_structure_and_Kaehler_metric_moduli_space_SO3_monopoles_complex_Kaehler_surface}
In this chapter, we prove the existence of a complex structure and K\"ahler metric on the open subset of smooth points of the moduli space $\sM(E,h,\om)$ of projective vortices on a Hermitian rank two vector bundle $(E,h)$ over a closed, complex K\"ahler manifold $(X,\omega)$ defined in the forthcoming \eqref{eq:Moduli_space_projective_vortices}. The Hitchin--Kobayashi correspondence provides a real analytic diffeomorphism from the open subspace of smooth points of $\sM(E,h,\om)$ onto the open subset of smooth points of the moduli space $\fM(E,\omega)$ of stable holomorphic pairs on $E$. We prove the existence of a complex structure on the open subset of smooth points of $\fM(E,\omega)$ by adapting \mutatis the proofs of the corresponding results in Chapter \ref{chap:Complex_structure_and_Kaehler_metric_moduli_space_HE_connections_complex_Kaehler_manifold} due to
Itoh \cite{Itoh_1983, Itoh_1985} for the moduli space $M(E,h,\omega)$ of anti-self-dual connections on a rank-two Hermitian vector bundle $(E,h)$ (with holomorphically trivial determinant) over a closed, complex K\"ahler surface and Kobayashi \cite[Section 7.4]{Kobayashi_differential_geometry_complex_vector_bundles} and Kim \cite{Kim_1987} for the moduli space $M(E,h,\omega)$ of projectively Hermitian--Einstein connections on a Hermitian vector bundle $(E,h)$ over a closed, complex K\"ahler manifold. We prove the existence of a K\"ahler metric on the open subset of smooth points of $\fM(E,\omega)$ by adapting proofs of the corresponding results in Chapter \ref{chap:Complex_structure_and_Kaehler_metric_moduli_space_HE_connections_complex_Kaehler_manifold} due to Itoh \cite{Itoh_1988} for the moduli space of anti-self-dual connections on a rank-two Hermitian bundle $(E,h)$ (with holomorphically trivial determinant) over a complex K\"ahler surface and by Kobayashi \cite[Section 7.6]{Kobayashi_differential_geometry_complex_vector_bundles}) for the moduli space of projectively Hermitian--Einstein connections on a Hermitian vector bundle $(E,h)$ over a closed, complex K\"ahler manifold.
The construction of a symplectic form on the moduli space of stable vector bundles over a smooth algebraic surface is described by Huybrechts and Lehn in \cite[Chapter 10]{Huybrechts_Lehn_geometry_moduli_spaces_sheaves}. We apply the Marsden--Weinstein Symplectic Reduction Theorem \ref{thm:Marsden-Weinstein_symplectic_quotient} to prove that the open subset of smooth points of $\sM(E,h,\om)$ is a K\"ahler manifold and identify a Hamiltonian function for the $S^1$ action on $\sM(E,h,\om)$ given by scalar multiplication on the section, which leads to a proof of Theorem \ref{mainthm:IdentifyCriticalPoints}.

In Sections \ref{sec:Moduli_space_simple_strongly_simple_pairs} and \ref{sec:Moduli_space_stable_pairs}, we define the moduli subspaces of (strongly) simple and stable holomorphic pairs of vector bundles and sections and the properties of semistability and polystability of such pairs. We define the moduli space of projective vortices in Section \ref{sec:Moduli_space_projective_vortices}, giving the Zariski tangent space of this moduli space in Theorem \ref{thm:Kobayashi_7_4_19_projective_vortices}. Section \ref{sec:HK_correspondence_between_SO3_monopoles_and_semistable_pairs} contains our discussion of the Hitchin--Kobayashi correspondence between projective vortices and semistable holomorphic pairs of vector bundles and sections. In Section \ref{sec:Moduli_space_SO3_monopoles_open_subspace_moduli_space_simple_holomorphic_pairs}, we describe how this correspondence gives a real analytic embedding. Finally, in Section \ref{sec:Marsden-Weinstein_reduction_moduli_space_SO3_monopoles_symplectic_quotient}, we build a complex, K\"ahler structure on the subspace of regular points of the moduli space of
projective vortices over a complex, K\"ahler manifold.  We then use this construction to
prove Theorem \ref{thm:Critical_points_Hitchin_Hamiltonian_function_moduli_space_projective_vortices} and hence Theorem \ref{mainthm:IdentifyCriticalPoints}.

\section[Moduli spaces of simple and strongly simple holomorphic pairs]{Moduli spaces of simple and strongly simple holomorphic pairs of vector bundles and sections}
\label{sec:Moduli_space_simple_strongly_simple_pairs}
We continue the assumptions and notation in Section \ref{sec:Moduli_space_simple_bundles} and extend its development from the case of holomorphic structures to holomorphic pairs. We refer to elements $(\bar\partial_E,\varphi) \in \sA^{0,1}(E)\times W^{1,p}(E)$ as a \emph{$(0,1)$-pairs}.
\label{page:(0,1)-pair}
The group $W^{2,p}(\SL(E))$  of complex, determinant-one, $W^{2,p}$ automorphisms of $E$ acts holomorphically on $\sA^{0,1}(E)\times W^{1,p}(E)$ by
\begin{multline}
\label{eq:SL(E)ActionOn(0,1)Pairs}
W^{2,p}(\SL(E))\times
\left( \sA^{0,1}(E)\times W^{1,p}(E)\right)
\ni
\left( v,(\bar\rd_E,\varphi)\right)
\\
\mapsto
v^*(\bar\rd_E,\varphi)
:=
(v^{-1}\circ \bar\rd_E\circ v,v^{-1}\varphi)
\in
\sA^{0,1}(E)\times W^{1,p}(E).
\end{multline}
Let $[\bar\partial_E,\varphi]$ denote the gauge-equivalence class of the $(0,1)$-pair $(\bar\partial_E,\varphi)$
\label{page:EquivClassOf(0,1)Pair} under the action  in \eqref{eq:SL(E)ActionOn(0,1)Pairs}.

By analogy with our definition \eqref{eq:Moduli_space_holomorphic_structures} of the moduli space $\cM(E)$ of holomorphic structures on $E$, let
\begin{equation}
  \label{eq:Moduli_space_holomorphic_pairs}
  \fM(E) := \left\{(\bar\partial_E,\varphi) \in \sA^{0,1}(E)\times W^{1,p}(E): F_{\bar\partial_E} = 0 \text{ and } \bar\partial_E\varphi = 0 \right\}/W^{2,p}(\SL(E))
\end{equation}
denote the moduli space of holomorphic pairs for $E$ with fixed smooth, holomorphic structure on $\det E$ modulo the action of $ W^{2,p}(\SL(E))$ given by \eqref{eq:SL(E)_Action_on_(0,1)conn}, and equipped with the quotient topology.

We define the \emph{stabilizer} of $(\bar\partial_E,\varphi)$ in $W^{2,p}(\SL(E))$ to be
\begin{equation}
  \label{eq:Stabilizer_holomorphic_pair}
  \Stab(\bar\partial_E,\varphi)
  :=
  \left\{v \in W^{2,p}(\SL(E)): 
  v^*(\bar\partial_E,\varphi) = (\bar\partial_E,\varphi) \right\}.
\end{equation}
The stabilizer $\Stab(\bar\partial_E,\varphi)$ has the following Lie group structure, similar to that of $\Stab(\bar\partial_E)$ appearing in Lemma \ref{lem:Stab(rdE)_is_Lie_Group}.

\begin{lem}[Lie algebra and complex structure of $\Stab(\bar\partial_E,\varphi)$]
\label{lem:Stab(rdE,varphi)_is_Lie_Group}
Let $E$ be a smooth complex vector bundle over a connected, complex manifold of dimension $n$ and let $p\in (n,\infty)$. If $(\bar\partial_E,\varphi) \in \sA^{0,1}(E)\times W^{1,p}(E)$, then its stabilizer $\Stab(\bar\partial_E,\varphi)$ is a complex Lie group with Lie algebra given by the harmonic space $\bH_{\bar\partial_E,\varphi}^0$ in \eqref{eq:H_dbar_Avarphi^0bullet}.
\end{lem}

\begin{proof}
The proof is identical to that of Lemma \ref{lem:Stab(rdE)_is_Lie_Group}.
\end{proof}

Keeping in mind the pullback action \eqref{eq:SL(E)ActionOn(0,1)Pairs} of $W^{2,p}(\SL(E))$ on $\sA^{0,1}(E)\times W^{1,p}(E)$, we say that an endomorphism $v\in W^{2,p}(\End(E))$ \emph{commutes} with a $W^{1,p}$ pair $(\bar\partial_E,\varphi)$ if
\begin{equation}
  \label{eq:Pair_commuting_with_endomorphism}
  (\bar\partial_E\circ v,\varphi) = (v\circ\bar\partial_E, v\varphi).
\end{equation}
If $v\in W^{2,p}(\SL(E))$, then $v$ commutes with $(\bar\partial_E,\varphi)$ if and only if $v\in\Stab(\bar\partial_E,\varphi)$. By analogy with Definition \ref{defn:Simple_01-connection}, we make the

\begin{defn}[Simple and strongly simple $(0,1)$-pairs for smooth, complex vector bundles over almost Hermitian manifolds]
\label{defn:Simple_pair}
(Compare Flenner and L\"ubke \cite[Section 1, p. 100, line 2]{Flenner_Lubke_2002}, L\"ubke and Teleman \cite[Definition 3.1.1, p. 23 and Section 6.2.1, p. 64, line 20]{Lubke_Teleman_2006}, Okonek and Teleman \cite[Definition 5.2]{OTQuaternion}, and \cite[Definition 1.2, p. 67]{Suyama_1996}.)
Let $E$ be a smooth complex vector bundle over a connected, smooth almost Hermitian manifold of real dimension $2n$, let $p \in (n,\infty)$, and let $(\bar\partial_E,\varphi) \in \sA^{0,1}(E)\times W^{1,p}(E)$.
\begin{enumerate}
\item\label{item:SL(E)_simple_01-pair}
If $\varphi\equiv 0$, the zero-section pair $(\bar\partial_E,0)$ is \emph{simple} if $\bar\partial_E$ is simple in the sense of Definition \ref{defn:Simple_01-connection}.

\item\label{item:Strongly_simple_pair}
The pair $(\bar\partial_E,\varphi)$ is \emph{strongly simple} if the only $W^{2,p}$ section of $\SL(E)$ that commutes with $(\bar\partial_E,\varphi)$ in the sense of \eqref{eq:Pair_commuting_with_endomorphism} is the identity endomorphism, $\id_E$.
\end{enumerate}
\end{defn}

In Definition \ref{defn:Simple_pair} \eqref{item:Strongly_simple_pair}, the condition for $(\bar\partial_E,\varphi)$ to be strongly simple is equivalent to $\Stab(\bar\partial_E,\varphi) = \{\id_E\}$: we adopt terminology of L\"ubke and Teleman, though other authors call such a pair `simple'. We note the following criterion for strong simplicity of holomorphic pairs that are split in
the sense of the forthcoming Definition \ref{defn:Split_(0,1)-pair}.

\begin{lem}[Conditions for strong simplicity of rank-two, non-zero-section holomorphic pairs]
\label{lem:Rank2_Split_HolomPairsSimple}
Let $E$ be a complex, rank-two vector bundle over a closed, connected, complex K\"ahler manifold $(X,\omega)$. Let $(\bar\partial_E,\varphi)$ be a holomorphic pair on $E$ with $\varphi\not\equiv 0$ and which is split in the sense of the forthcoming Definition \ref{defn:Split_(0,1)-pair}, so $E$ admits a decomposition $E=L_1\oplus L_2$ as a direct sum of complex line bundles and $\bar\rd_E=\bar\rd_{L_1}\oplus\bar\rd_{L_2}$, where $\bar\rd_{L_i}$ is a holomorphic structure on $L_i$ for $i=1,2$, and $\varphi$ is a section of $L_1$.  Let $\sL_i$ denote the sheaves of holomorphic sections of $(L_i,\bar\partial_{L_i})$. Then the following hold:
\begin{enumerate}
\item\label{item:Rank2_Split_HolomPairsSimpleCase}
The pair $(\bar\partial_E,\varphi)$ is strongly simple in the sense of Definition \ref{defn:Simple_pair} if and only if $H^0(X;\sL_1\otimes\sL_2^*)=(0)$.

\item\label{item:Rank2_Split_HolomPairsSimpleCase_DegreeCriterion}
If $\deg L_1<\deg L_2$, then $(\bar\partial_E,\varphi)$ is strongly simple in the sense of Definition \ref{defn:Simple_pair}.
\end{enumerate}
\end{lem}

\begin{proof}
Consider Item \eqref{item:Rank2_Split_HolomPairsSimpleCase}. Suppose that $v\in W^{2,p}(\SL(E))$ commutes with $(\bar\partial_E,\varphi)$. Write
\begin{equation}
\label{eq:Decomposition_Of_Element_of_Stabilizer}
v= \begin{pmatrix} v_{11} & v_{12} \\ v_{21} & v_{22} \end{pmatrix}
\end{equation}
with respect to the decomposition $E=L_1\oplus L_2$, so that $v_{ij}\in W^{2,p}(L_i\otimes L_j^*)$. The condition \eqref{eq:Pair_commuting_with_endomorphism} that $v\varphi=\varphi$ and our hypothesis that $\varphi\not\equiv 0$ imply that $v_{11}=1$ and $v_{21}=0$. Because $v\in W^{2,p}(\SL(E))$, we have $\det v=1$ and so we obtain $v_{22}=1$. Since $v$ commutes with $\bar\partial_E$, we see that $v$ is a holomorphic section of $\End(E)$ and thus $v_{12}$ is a holomorphic section of $L_1\otimes L_2^*$. Hence, the map
\begin{equation}
\label{eq:Stabilizer_Bijection}
H^0(X;\sL_1\otimes\sL_2^*)\ni v_{12}\mapsto
v=\begin{pmatrix} 1 & v_{12} \\ 0 & 1 \end{pmatrix} \in W^{2,p}(\SL(E)),
\end{equation}
gives a group isomorphism between $H^0(X;\sL_1\otimes\sL_2^*)$, where the group operation is addition, and the elements of $W^{2,p}(\SL(E))$ which commute with $(\bar\partial_E,\varphi)$, where the group operation is composition of gauge transformations. This yields Item \eqref{item:Rank2_Split_HolomPairsSimpleCase}.

Consider Item \eqref{item:Rank2_Split_HolomPairsSimpleCase_DegreeCriterion}. The hypothesis $\deg L_1<\deg L_2$ in this case implies that $H^0(X;\sL_1\otimes\sL_2^*)=(0)$ by Lemma \ref{lem:NegDegLineBundlesNoSections}. Therefore, Item \eqref{item:Rank2_Split_HolomPairsSimpleCase_DegreeCriterion} follows from Item \eqref{item:Rank2_Split_HolomPairsSimpleCase}.
\end{proof}

\begin{rmk}[Comparison of Definitions \ref{defn:Simple_01-connection} and \ref{defn:Simple_pair}]
\label{rmk:SimplePairsHave_GL(E)_Stabilizers}
We explain why Definition \ref{defn:Simple_pair} only contains a condition on sections of $\SL(E)$ and does not contain equivalent conditions on sections of $\GL(E)$ and $\End(E)$ as in Definition \ref{defn:Simple_01-connection}.

We say that an $(0,1)$-pair $(\bar\partial_E,\varphi)$ is \emph{$\GL(E)$-strongly simple} if the only $W^{2,p}$ section of $\GL(E)$ which commutes with $(\bar\partial_E,\varphi)$ is the identity isomorphism, $\id_E$. If a pair is $\GL(E)$-strongly simple, then it is strongly simple in the sense of Definition \ref{defn:Simple_pair}.
as in Definition \ref{defn:Simple_pair}.

However, the converse is not true. Let $(\bar\partial_E,\varphi)$ be a holomorphic pair satisfying the hypotheses of Lemma \ref{lem:Rank2_Split_HolomPairsSimple}. Assume that $H^0(X;\sL_1\otimes\sL_2^*)=(0)$, so
$(\bar\partial_E,\varphi)$ is strongly simple by Item \eqref{item:Rank2_Split_HolomPairsSimpleCase} of
Lemma \ref{lem:Rank2_Split_HolomPairsSimple}. If $L_1$ and $L_2$ are the complex line bundles appearing in the hypotheses of Lemma \ref{lem:Rank2_Split_HolomPairsSimple}, then the existence of gauge transformations,
\[
v_\lambda := \id_{L_1}\oplus \lambda\,\id_{L_2} \in W^{2,p}(\GL(E)), \quad\text{for } \lambda\in\CC^*,
\]
which commute with $(\bar\rd_E,\varphi)$ demonstrates that such a pair is not $\GL(E)$-strongly simple. Hence, $\GL(E)$-strongly simple is not equivalent to strongly simple.
\end{rmk}

\begin{rmk}[Strongly simple $(0,1)$-pairs $(\bar\partial_E,\varphi)$ satisfy $\bH_{\bar\partial_E,\varphi}^0 = (0)$]
\label{rmk:SimplePairsAndH0=0}
Suppose that a pair $(\bar\partial_E,\varphi)$ is strongly simple in the sense of Definition
\ref{defn:Simple_pair}, so that $\Stab(\bar\partial_E,\varphi) = \{\id_E\} \subset W^{2,p}(\SL(E))$. By Lemma \ref{lem:Stab(rdE)_is_Lie_Group}, the harmonic space $\bH_{\bar\partial_E,\varphi}^0$ in \eqref{eq:H_dbar_Avarphi^0bullet} is the Lie algebra of $\Stab(\bar\partial_E,\varphi)$ and thus $\bH_{\bar\partial_E,\varphi}^0=(0)$.
\end{rmk}

By analogy with \eqref{eq:Moduli_subspaces_holomorphic_structures}, we define the following subspaces of $\fM(E)$,
\begin{subequations}
  \label{eq:Moduli_subspaces_holomorphic_pairs}
\begin{align}
  \label{eq:Moduli_space_strongly_simple_holomorphic_pairs}
  \fM^{**}(E) &:= \left\{[\bar\partial_E,\varphi] \in \fM(E):
                    (\bar\partial_E,\varphi) \text{ is strongly simple} \right\},
  \\
  \label{eq:Moduli_space_holomorphic_pairs_NonZeroSection}
  \fM^0(E)
  &:=
    \left\{[\bar\partial_E,\varphi]\in \fM(E): \varphi\not\equiv 0\right\},
  \\
   \label{eq:Moduli_space_regular_holomorphic_pairs}
  \fM_\reg(E) &:= \left\{[\bar\partial_E,\varphi] \in \fM(E):
                        \bH_{\bar\partial_E,\varphi}^2 = (0) \right\},
\end{align}
\end{subequations}
where the harmonic space $\bH_{\bar\partial_E,\varphi}^2$ is defined by \eqref{eq:H_dbar_Avarphi^0bullet}. The subspace $\fM^0(E)$ is an open subspace of $\fM(E)$ because it is the complement of the closed subspace of points $[A,\varphi]\in \fM(E)$ with $\varphi\equiv 0$. The subspace $\fM_\reg(E)$ is an open subspace of $\fM(E)$ as a consequence of the Implicit Mapping Theorem. It is unclear whether or not $\fM^{**}(E)$ is an open subspace of $\fM(E)$ in general, but the forthcoming Lemma \ref{lem:Openness_subspace_strongly_simple_01_pairs} provides related openness properties under additional conditions. The method of proof of Theorem \ref{thm:Kobayashi_7_3_17} extends \mutatis to give the

\begin{thm}[Local topology of the moduli space of holomorphic pairs over a complex manifold near a simple pair]
\label{thm:Kobayashi_7_3_17_pairs}  
Let $E$ be a smooth Hermitian vector bundle, with a fixed smooth holomorphic structure on $\det E$, over a closed, Hermitian, complex manifold $X$ of complex dimension $n$ and let $p\in(n,\infty)$. If $(\bar\partial_E,\varphi)$ is a holomorphic pair on $E$ that is strongly simple in the sense of Definition \ref{defn:Simple_pair} and
\begin{multline*}
  \cS_{\bar\partial_E,\varphi}
  := (\bar\partial_E,\varphi) + \left\{(\alpha,\phi) \in W^{1,p}\left(\Lambda^{0,1}(\fsl(E))\oplus E\right): F_{\bar\partial_E+\alpha} = 0 \text{ and }\right.
  \\
  \left.(\bar\partial_E+\alpha)(\varphi + \phi) = 0 \text{ and }
    \bar\partial_{E,\varphi}^{0,*}(\alpha,\phi) = 0\right\},
\end{multline*}
where $\bar\partial_{E,\varphi}^0$ is as in \eqref{eq:d0StablePair} and $\bar\partial_{E,\varphi}^{*,0}$ is its $L^2$ adjoint, then the quotient map,
\[
  \pi:\cU_{\bar\partial_E,\varphi} \ni (\bar\partial_E,\varphi) + (\alpha,\phi) \mapsto [\bar\partial_E + \alpha, \varphi+\phi] \in \fM(E),
\]
gives a topological epimorphism from an open neighborhood $\cU_{\bar\partial_E,\varphi}$ of the pair $(\bar\partial_E,\varphi)$ in $\cS_{\bar\partial_E,\varphi}$ onto an open neighborhood of $[\bar\partial_E,\varphi] = \pi(\bar\partial_E,\varphi) \in \fM(E)$.
\end{thm}

\begin{rmk}[Topology of the moduli space of holomorphic pairs over a complex manifold near points represented by strongly simple pairs or simple pairs]
\label{rmk:Kobayashi_7_3_17_strongly_simple_pair}
Unlike the statement of Theorem \ref{thm:Kobayashi_7_3_17}, we only assert in Theorem \ref{thm:Kobayashi_7_3_17_pairs} that the quotient map $\pi$ gives an epimorphism and not a homeomorphism of topological spaces. This is due to our inability to prove that $\pi$ is injective without additional hypotheses, such as stability of the pair $(\bar\partial_E,\varphi)$. Consequently, we do not prove a direct analogue of Theorem \ref{thm:Kobayashi_7_3_17} for pairs. Instead, we prove an analogue for pairs of the results in Theorem \ref{thm:Kobayashi_7_3_17} under additional hypotheses in the forthcoming Corollary \ref{cor:Lubke_Teleman_6-3-7_and_Kobayashi_7_3_17_pair}.

We would expect an analogue of Theorem \ref{thm:Kobayashi_7_3_17_pairs} to hold when the pair $(\bar\partial_E,0)$ is simple in the sense of Definition \ref{defn:Simple_pair} \eqref{item:SL(E)_simple_01-pair}, with the caveat that such a pair has a stabilizer in $W^{2,p}(\SL(E))$ that is isomorphic to $C_r$ (the group of $r$-th roots of unity) and so is not strongly simple. In particular, we would expect $\fM(E)$ to have an orbifold singularity near such a point. Compare L\"ubke and Teleman \cite[Remark 6.1.6 (2), p. 60]{Lubke_Teleman_2006}.
\end{rmk}

\begin{rmk}[Openness of the subspace of strongly simple holomorphic pairs]
\label{rmk:Openness_subspace_strongly_simple_holomorphic_pairs}  
The question of openness in the affine space $\sA^{0,1}(E)\times W^{1,p}(E)$ of the subspace of strongly simple $(0,1)$-pairs is more involved than that answered by Lemma \ref{lem:Openness_subspace_simple_01-connections}, giving openness in the affine space $\sA^{0,1}(E)$ of the subspace of simple $(0,1)$-connections. The forthcoming Lemma \ref{lem:Openness_subspace_strongly_simple_01_pairs} suffices for our applications in this monograph. Suyama states a more general result in his \cite[Proposition 1.2, p. 67]{Suyama_1996}, giving openness in the space of holomorphic pairs of the subspace of strongly simple holomorphic pairs, though his proof appears to be incomplete: in his \cite[Lemma 2.2, p. 70]{Suyama_1996}, he states that a pair $(\bar\partial_E,\varphi)$ is simple if and only if $\Ker\bar\partial_{E,\varphi}^0 = (0)$, but his proof only gives the ``only if'' part of that assertion and overlooks the possibility that $\Stab(\bar\partial_E,\varphi)$ could be a discrete but non-trivial Lie group. Flenner and L\"ubke provide a related result in \cite[Corollary 3.3, p. 105]{Flenner_Lubke_2002}.
\end{rmk}

\section[Moduli spaces of stable holomorphic pairs]{Moduli spaces of stable holomorphic pairs of vector bundles and sections}
\label{sec:Moduli_space_stable_pairs}
We recall the following definitions due to Bradlow \cite[Definition 2.1.3, p. 177]{Bradlow_1991}. Let $\sE$ denote the sheaf of sections of a holomorphic vector bundle $(E,\bar\partial_E)$ over a complex, K\"ahler manifold. Let $\Coh(\sE)$
\label{page:ProperCoherentSubsheaves}
 be the set of coherent\footnote{We adopt the convention in Bradlow and Garc{\'{\i}}a--Prada \cite[Definition 6.1, p. 580]{BradlowGP} of using coherent rather than reflexive subsheaves.} subsheaves $\sF$ of $\sE$ with $0<\rank\sF<\rank\sE$ and, if $\varphi \in \sE$, let $\Coh_\varphi(\sE) := \{\sF\in\Coh(\sE):\varphi\in\sF\}$
 \label{page:ProperCoherentSubsheavesContainingSection}. Define
\begin{subequations}
\label{eq:Bradlow_mu_M_m}  
\begin{align}
  \label{eq:Bradlow_mu_M}  
  \mu_M(\sE) &:= \max\left(\mu(\sE), \sup_{\sF \in \Coh(\sE)}\mu(\sF)\right),
  \\
  \label{eq:Bradlow_mu_m}  
  \mu_m(\varphi)
             &:= \inf_{\sF \in \Coh_\varphi(\sE)} \frac{\rank\sE\cdot\mu(\sE) - \rank\sF\cdot\mu(\sF)}{\rank\sE - \rank\sF}.
\end{align}
\end{subequations}
Given a short exact sequence of sheaves,
\[
  0 \to \sF \to \sE \to \sE/\sF \to 0,
\]
then $\deg\sE = \deg\sF + \deg(\sE/\sF)$ and $\rank\sE = \rank\sF + \rank(\sE/\sF)$.
Hence, the slope of the quotient is
\[
  \mu(\sE/\sF) = \frac{\deg(\sE/\sF)}{\rank(\sE/\sF)} = \frac{\deg\sE - \deg\sF}{\rank\sE - \rank\sF}
  = \frac{\rank\sE\cdot\mu(\sE) - \rank\sF\cdot\mu(\sF)}{\rank\sE - \rank\sF}.
\]
Hence, one may rewrite the definition \eqref{eq:Bradlow_mu_m} in the simpler form
\begin{equation}
  \label{eq:Bradlow_mu_m_simpler} 
  \mu_m(\varphi) = \inf_{\sF \in \Coh_\varphi(\sE)} \mu(\sE/\sF).
\end{equation}
We have the following generalization of Definition \ref{defn:Stable_holomorphic_structure}.

\begin{defn}[Stable pairs of holomorphic vector bundles and sections over complex K\"ahler manifolds]
\label{defn:StablePair_arbitrary_rank}
(See Bradlow, Daskalopoulos, and Wentworth \cite[Definition 1.1, p. 731]{BradlowDaskWent}, Bradlow \cite[Definition 2.1.5, p. 179]{Bradlow_1991}, Bradlow and Garc{\'{\i}}a--Prada \cite[Definition 6.1, p. 580]{BradlowGP}, and Okonek and Teleman \cite[Section 3, p. 903]{OTVortex}, \cite[Definition 5.3, p. 375]{OTQuaternion}.)
Let $\sE$ denote the sheaf of sections of a holomorphic vector bundle $(E,\bar\partial_E)$ over a closed, complex, K\"ahler manifold. One says that $(\sE,\varphi)$ is \emph{$\varphi$-stable} if
\begin{equation}
  \label{eq:Bradlow_phi_stable}
  \mu_M(\sE) < \mu_m(\varphi).
\end{equation}
One says that $(\sE,\varphi)$ is $\lambda$-\emph{stable} if 
\begin{equation}
\label{eq:Bradlow_lambda_stable}
\mu_M(\sE) < \lambda < \mu_m(\varphi)
\end{equation}
for a real constant $\lambda$. One says that $(\sE,\varphi)$ is $\lambda$-\emph{semistable} if the pair is $\lambda$-stable or
\[
  \mu_M(\sE) = \lambda = \mu_m(\varphi).
\]
One says that $(\sE,\varphi)$ is $\lambda$-\emph{polystable} if $(\sE,\varphi)$ is $\lambda$-\emph{stable} or $\sE = \sE'\oplus\sE''$ as a direct sum of holomorphic subsheaves such that $\varphi$ is a section of $\sE'$ and $(\sE',\varphi)$ is $\lambda$-stable and $\sE''$ is polystable of slope $\lambda$ in the sense of Definition \ref{defn:Stable_holomorphic_structure}. When the value of the constant $\lambda$ is unimportant, we refer to pairs as being \emph{stable}, \emph{polystable}, or \emph{semistable}.
One says that $(\sE,\varphi)$ is
\begin{inparaenum}[\itshape i\upshape)]
\item \emph{unstable} if it is not semistable,
\item \emph{strictly semistable} if it is semistable but not polystable, and
\item \emph{strictly polystable} if it is polystable but not stable.
\end{inparaenum}
\end{defn}

\begin{rmk}[Gieseker--Maruyama and slope stability conditions for pairs]
\label{rmk:Gieseker-stable-pair}
Huybrechts and Lehn \cite[Definition 1.1, p. 69]{Huybrechts_Lehn_1995jag}, \cite[pp. 297--298]{Huybrechts_Lehn_1995ijm} and Y. Lin \cite[Definition 2.1, p. 129]{Lin_2018} employ a version of \emph{stability in the sense of Gieseker} \cite{Gieseker_1977} and \emph{Maruyama} \cite{Maruyama_1978} for vector bundles that is different from the version of \emph{slope stability} (due to Mumford \cite{Mumford_1963}, \cite[Appendix to Chapter 5, Section C, Definition, p. 224]{Mumford_Fogarty_Kirwan_geometric_invariant_theory} and Takemoto \cite{Takemoto_1972} for vector bundles) described in Definition \ref{defn:StablePair_arbitrary_rank}, motivated by their construction of the Gieseker compactification of moduli spaces of pairs over smooth, projective algebraic varieties.
\end{rmk}

\begin{rmk}[Stability for pairs generalizes stability for sheaves]
\label{rmk:Stability_pairs_generalizes_stability_sheaves}
Continue the setting of  Definition \ref{defn:StablePair_arbitrary_rank}. If $\sE$ is stable in the sense of Definition \ref{defn:Stable_holomorphic_structure}, so that $\mu(\sF) < \mu(\sE)$ for all $\sF\in\Coh(\sE)$, then $\mu_M(\sE) = \mu(\sE)$ by \eqref{eq:Bradlow_mu_M}. For all $\sF \in \Coh_\varphi(\sE)$, we have by \eqref{eq:Bradlow_mu_m} that
\begin{multline*}
  \frac{\rank\sE\cdot\mu(\sE) - \rank\sF\cdot\mu(\sF)}{\rank\sE - \rank\sF}
  \\
  =
  \frac{(\rank\sE- \rank\sF)\cdot\mu(\sE) + \rank\sF\cdot(\mu(\sE) - \mu(\sF))}{\rank\sE - \rank\sF}
  >
  \mu(\sE).
\end{multline*}
Hence, $\mu_m(\varphi) > \mu_M(\sE)$ and so $(\sE,\varphi)$ is a stable pair in the sense of Definition \ref{defn:StablePair_arbitrary_rank}.
\end{rmk}

\begin{rmk}[Stability for rank-two, split holomorphic pairs]
\label{rmk:Stability_of_Split_Pairs}
Suppose that $\sE$ has rank two, with $\sE = \sL_1\oplus\sL_2$ as a direct sum of rank one coherent sheaves, so $\deg\sE = \deg\sL_1 + \deg\sL_2$ and
\begin{equation}
  \label{eq:Slope_E_split}
  \mu(\sE) = \frac{1}{2}(\deg\sL_1 + \deg\sL_2),
\end{equation}  
and that $\varphi$ is a section of $L_1$. We consider different cases.

\setcounter{case}{0}
\begin{case}[$\deg\sL_1 \leq \deg\sL_2$ and $\varphi\not\equiv 0$]
\label{case:degL_1Smaller_SectionNonZero}
Every $\sF\in\Coh(\sE)$ is a subsheaf of the locally free sheaf $\sE$ and is thus torsion-free. In addition, every $\sF\in\Coh(\sE)$ has rank one when $\sE$ has rank two and so $\sF$ is automatically stable by Definition \ref{defn:StablePair_arbitrary_rank}. We claim that $\deg\sF\leq\deg\sL_2$ when $\deg\sL_1 \leq \deg\sL_2$. To see this, let $f_i:\sF\to\sL_i$ be the morphism of sheaves defined by the composition $\sF\to\sE\to\sL_i$ for $i=1,2$, where the first map is inclusion and the second map is the natural projection. If $\deg\sF>\deg\sL_i$, then $f_i=0$ by Kobayashi \cite[Proposition 5.7.11 (1), p. 158]{Kobayashi_differential_geometry_complex_vector_bundles}. Because $\sF$ is a subsheaf of $\sE$, at least one of the two maps $f_i$ must be non-zero. By assumption, $\deg\sL_1\leq\deg\sL_2$ and so if $\deg\sF>\deg\sL_2$, then we would have $\deg\sF>\deg\sL_i$ for both $i=1,2$ and so $f_i=0$ for both $i=1,2$, a contradiction. Therefore, we have 
\[
  \sup_{\sF \in \Coh(\sE)}\mu(\sF) = \mu(\sL_2).
\]
Moreover, because $\deg\sL_1 \leq \deg\sL_2$, we have $\mu(\sE) \leq \deg\sL_2 = \mu(\sL_2)$ by \eqref{eq:Slope_E_split} and thus \eqref{eq:Bradlow_mu_M} yields
\[
  \mu_M(\sE) = \max\left(\mu(\sE),\mu(\sL_2)\right) = \mu(\sL_2).
\]
By the proof of the preceding claim, if $\sF\in\Coh_\varphi(\sE)$ with $\deg\sF>\deg\sL_1$, then we would have $f_1=0$.  However, $\varphi$ is a section of $\sF$ since $\sF\in\Coh_\varphi(\sE)$ and the equality $\varphi=f_1\circ\varphi$ would contradict $f_1=0$ when $\varphi\not\equiv 0$. Hence, for all $\sF\in\Coh_\varphi(\sE)$, we have $\deg\sF\le\deg\sL_1$ and thus
\[
  \mu(\sE/\sF) = \deg(\sE/\sF) = \deg(\sE) - \deg(\sF) = \deg\sL_1+\deg\sL_2-\deg\sF \geq \deg\sL_2,
\]
where the third equality follows from \eqref{eq:Slope_E_split}. The preceding inequality, the fact that $\sF=\sL_1$ yields an equality above, and the expression \eqref{eq:Bradlow_mu_m_simpler} for $\mu_m(\varphi)$ gives
\[
  \mu_m(\varphi) = \inf_{\sF \in \Coh_\varphi(\sE)} \mu(\sE/\sF) = \mu(\sL_2).
\]  
Therefore, the preceding equalities yield
\begin{equation}
  \label{eq:Split_pair_muL2_semistable}
  \mu_M(\sE) = \mu(\sL_2) = \mu_m(\varphi),
\end{equation}
and so $(\sE,\varphi)$ is a $\mu(\sL_2)$-semistable but not stable pair in the sense of Definition \ref{defn:StablePair_arbitrary_rank}.

\begin{subcase}[$\deg\sL_1 = \deg\sL_2$ and $\varphi\not\equiv 0$]
By assumption, $\deg\sL_1 = \deg\sL_2$, so \eqref{eq:Bradlow_mu_M} for $\sL_1$ in place of $\sE$ yields
\[
  \mu_M(\sL_1) = \mu(\sL_1) = \mu(\sL_2),
\]
and thus $(\sL_1,\varphi)$ is a $\mu(\sL_2)$-semistable but not $\mu(\sL_2)$-stable pair in the sense of Definition \ref{defn:StablePair_arbitrary_rank}. (The upper $\mu_m(\sL_1,\varphi)$-limit in \eqref{eq:Bradlow_lambda_stable} is vacuous since $\sL_1$ has rank one.) Since $\sL_2$ has slope $\mu(\sL_2)$, we see that $(\sE,\varphi)$ is a $\mu(\sL_2)$-semistable but not polystable pair in the sense of Definition \ref{defn:StablePair_arbitrary_rank}.
\end{subcase}

\begin{subcase}[$\deg\sL_1 < \deg\sL_2$ and $\varphi\not\equiv 0$]
Equation \eqref{eq:Bradlow_mu_M} for $\sL_1$ in place of $\sE$ and the assumption $\deg\sL_1 < \deg\sL_2$ yield
\[
  \mu_M(\sL_1) = \mu(\sL_1) < \mu(\sL_2),
\]
and so $(\sL_1,\varphi)$ is a $\mu(\sL_2)$-stable pair in the sense of Definition \ref{defn:StablePair_arbitrary_rank}. (Again, the upper $\mu_m(\sL_1,\varphi)$-limit in \eqref{eq:Bradlow_lambda_stable} is vacuous since $\sL_1$ has rank one.) Since $\sL_2$ has slope $\mu(\sL_2)$, we see that $(\sE,\varphi)$ is a $\mu(\sL_2)$-polystable but not stable pair in the sense of Definition \ref{defn:StablePair_arbitrary_rank}.
\end{subcase}
\end{case}

\begin{case}[$\deg\sL_1 = \deg\sL_2$ and $\varphi\equiv 0$]
Because $\deg\sL_1 = \deg\sL_2$ and $\varphi\equiv 0$ by assumption, we see that $\sE = \sE'\oplus\sE''$ with $\sE' = 0$ and $\sE'' = \sL_1\oplus\sL_2$ is a polystable sheaf in the sense of Definition \ref{defn:Stable_holomorphic_structure} with slope $\mu(\sL_2)$, and so $(\sE,0)$ is a $\mu(\sL_2)$-polystable pair in the sense of Definition \ref{defn:StablePair_arbitrary_rank}. The equalities \eqref{eq:Split_pair_muL2_semistable} continue to hold for this case, so $(\sE,0)$ is not a stable pair.
\end{case}  

\begin{case}[$\deg\sL_1 > \deg\sL_2$]
By assumption $\deg\sL_1 > \deg\sL_2$ and so the argument from Case \ref{case:degL_1Smaller_SectionNonZero} implies that
\[
  \sup_{\sF \in \Coh(\sE)}\mu(\sF) = \mu(\sL_1).
\]
Because $\deg\sL_1 < \deg\sL_2$, we have $\mu(\sE) < \deg\sL_1 = \mu(\sL_1)$ by \eqref{eq:Slope_E_split}, and thus equation \eqref{eq:Bradlow_mu_M} yields
\[
  \mu_M(\sE) = \max\left(\mu(\sE),\mu(\sL_1)\right) = \mu(\sL_1).
\]
The argument from Case \ref{case:degL_1Smaller_SectionNonZero} also implies that for all $\sF\in\Coh_\varphi(\sE)$ we have $\deg\sF\le\deg\sL_1$ and thus
\[
  \mu(\sE/\sF) \ge\mu(\sL_2),
\]
with equality when $\sF=\sL_1$. The definition \eqref{eq:Bradlow_mu_m_simpler} thus gives
\[
  \mu_m(\varphi) = \mu(\sL_2).
\]
Hence, the preceding equalities and inequalities yield
\[
  \mu_M(\sE) = \mu(\sL_1) > \mu(\sL_2) = \mu_m(\varphi),
\]
and so $(\sE,\varphi)$ is an unstable pair in the sense of Definition \ref{defn:StablePair_arbitrary_rank}.
\end{case}

The preceding simplifications when $\sE$ has complex rank two partly reflect those of Okonek and Teleman \cite[Definition 5.4, p. 375]{OTQuaternion} and  L\"ubke and Teleman \cite[Definition 6.3.2, p. 71]{Lubke_Teleman_2006}.
\end{rmk}

\begin{rmk}[Closure of orbits]
\label{rmk:Closure_of_orbits} 
Note that if $(\sE,\varphi)$ is a strictly $\mu(\sL_2)$-semistable pair as in Remark \ref{rmk:Stability_pairs_generalizes_stability_sheaves}, so $\deg\sL_1 = \deg\sL_2$ and $\varphi\not\equiv 0$, and we act on $(\bar\partial_E,\varphi)$ by a sequence of $\SL(E)$ gauge transformations of the form
\[
  v_t := \begin{pmatrix} t & 0 \\ 0 & t^{-1} \end{pmatrix}, \quad\text{for } t \in \CC^*,
\]
then $v_t(\bar\partial_E,\varphi)$ converges to $(\bar\partial_E,0)$ as $t\to 0$. Hence, the point $[\sE,0]$ is in the closure of the $\SL(E)$-orbit of the semistable point $[\sE,\varphi]$. If we applied the same argument to a strictly $\mu(\sL_2)$-polystable point $[\sE,\varphi]$, so $\deg\sL_1 < \deg\sL_2$ and $\varphi\not\equiv 0$, the point $[\sE,0]$ would be unstable since $\mu_M(\sE)=\mu(\sL_2) > \mu(\sL_1)=\mu_m(0)$.
\end{rmk}

\begin{rmk}[$\mu(\sE)$-polystable pairs]
\label{rmk:mu(sE)_-polystable pairs}
Suppose that $(\sE,\varphi)$ is a $\lambda$-polystable pair in Definition \ref{defn:StablePair_arbitrary_rank}. If $(\sE,\varphi)$ is a $\lambda$-stable, then the definition \eqref{eq:Bradlow_mu_M} and inequalities \eqref{eq:Bradlow_lambda_stable} imply that $\lambda \geq \mu(\sE)$. If $(\sE,\varphi)$ is strictly $\lambda$-polystable, then $\sE = \sE'\oplus\sE''$ as a direct sum of holomorphic subsheaves such that $\varphi$ is a section of $\sE'$ and $(\sE',\varphi)$ is $\lambda$-stable and $\sE''$ is polystable of slope $\lambda$ in the sense of Definition \ref{defn:Stable_holomorphic_structure}, and thus
\[
  \mu_M(\sE') < \lambda < \mu_m(\sE',\varphi) \quad\text{and}\quad \mu(\sE'') = \lambda,
\]
Therefore, $\lambda \geq \mu(\sE')$ by definition \eqref{eq:Bradlow_mu_M} of $\mu_M(\sE')$ and we see that
\begin{align*}
  \mu(\sE) &= \mu(\sE'\oplus\sE'') = \frac{\deg(\sE'\oplus\sE'')}{\rank(\sE'\oplus\sE'')}
  = \frac{\deg\sE' + \deg\sE''}{\rank\sE' + \rank\sE''}
  \\
  &= \frac{\deg\sE'}{\rank\sE'}\frac{\rank\sE'}{\rank\sE' + \rank\sE''}
  + \frac{\deg\sE''}{\rank\sE''}\frac{\rank\sE''}{\rank\sE' + \rank\sE''}
  \\
  &= \mu(\sE')\frac{\rank\sE'}{\rank\sE' + \rank\sE''} + \mu(\sE'')\frac{\rank\sE''}{\rank\sE' + \rank\sE''}
  \\
           &\leq \lambda\frac{\rank\sE'}{\rank\sE' + \rank\sE''} + \lambda\frac{\rank\sE''}{\rank\sE' + \rank\sE''}
             = \lambda.
\end{align*}
Consequently, for any $\lambda$-polystable pair $(\sE,\varphi)$, we conclude that
\begin{equation}
  \label{eq:lambda_geq_mu(sE)_for_lambda_polystable_pair}
  \lambda \geq \mu(\sE).
\end{equation}  
We observe that $\lambda = \mu(\sE)$ in \eqref{eq:lambda_geq_mu(sE)_for_lambda_polystable_pair} if and only if $\varphi \equiv 0$ by the forthcoming identity \eqref{eq:Okonek-Teleman_Bradlow-Garcia-Prada_lambda_mu(sE)_L2norm_varphi_square_identity}.
\end{rmk}  

The result for stable holomorphic bundles analogous to the forthcoming lemma appears in Friedman \cite[Chapter 4, Corollary 8, p. 88]{FriedmanBundleBook}.

\begin{lem}[Stable $\implies$ strongly simple for non-zero-section, holomorphic pairs]
\label{lem:Stable_implies_strongly_simple_for_01-pairs}
Let $E$ be a complex vector bundle over a closed, complex K\"ahler manifold $(X,\om)$. If $(\bar\partial_E,\varphi)$ is a holomorphic pair with $\varphi\not\equiv 0$ that is stable as in Definition \ref{defn:StablePair_arbitrary_rank}, then it is strongly simple as in Definition \ref{defn:Simple_pair}.
\end{lem}

\begin{proof}
We adapt the proof due to Dowker \cite[Proposition 1.3.4. p. 12]{DowkerThesis}, which is based in turn on that of Thaddeus \cite[Lemma 1.6, p. 320]{ThaddeusVerlinde}.  Let $\sE$ denote the sheaf of $\bar\partial_E$-holomorphic sections of $E$. By hypothesis, $(\bar\rd_E,\varphi)$ is stable and so there is a constant $\lambda$ that obeys \eqref{eq:Bradlow_lambda_stable}, namely
\[
\mu_M(\sE) < \lambda < \mu_m(\varphi).
\]
By Definition \ref{defn:Simple_pair}, the pair $(\bar\partial_E,\varphi)$ is strongly simple if and only if the only $W^{2,p}$ section $v$ of $\SL(E)$ that commutes with $(\bar\partial_E,\varphi)$ in the sense of \eqref{eq:Pair_commuting_with_endomorphism} is the identity, $\id_E$. Suppose $v\in W^{2,p}(\SL(E))$ commutes with $(\bar\partial_E,\varphi)$, write $w=v-\id_E$, and observe that \eqref{eq:Pair_commuting_with_endomorphism} yields
\[
  \bar\partial_E\circ w = w\circ\bar\partial_E \quad\text{and}\quad w\varphi = 0.
\]
Consequently, for the $(0,1)$-connection $\bar\partial_{\End(E)}$ on $\End(E)$ defined 
in \eqref{eq:(0,1)Connection_Induced_On_Endomorphism_Bundle},
$\bar\partial_{\End(E)}w = [\bar\partial_E,w] = 0$, so $w$ is a $\bar\partial_{\End(E)}$-holomorphic section of $\End(E)$ and $\varphi\in\Ker w$. We claim that $(\bar\partial_E,\varphi)$ is strongly simple and suppose, to get a contradiction, that $w\not\equiv 0$. Define the subsheaves $\sK \subset \sE$ and $\sJ \subset \sE$ to be the kernel and image, respectively, of $w$. By Grauert and Remmert \cite[Annex, Section 4.2, Consequence 2, p. 237]{Grauert_Remmert_coherent_analytic_sheaves}, the sheaves $\sK$ and $\sJ$ are coherent. Let $r_\sE$, $r_\sK$, and $r_\sJ$ denote the ranks of $\sE$, $\sK$, and $\sJ$ respectively. Because $\varphi$ is a non-zero section of $\sK$ and because $w$ is not the zero endomorphism, we have
\begin{equation}
\label{eq:RankInequality1}
0<r_\sK <r_\sE.
\end{equation}
The inequalities \eqref{eq:RankInequality1} and the fact that $\varphi$ is a section of $\sK$ imply that $\sK\in\Coh_\varphi(\sE)$. The expression for $\mu_m(\varphi)$ in \eqref{eq:Bradlow_mu_m_simpler} thus yields
\begin{equation}
\label{eq:LowerBoundOnMuOfImage}
\mu_m(\varphi)\le \mu(\sE/\sK).
\end{equation}
The short exact sequence (see Grauert and Remmert \cite[Annex, Section 2.2, p. 230]{Grauert_Remmert_coherent_analytic_sheaves}),
\[
0 \to \sK \to \sE \to \sJ \to 0
\]
yields the following equality of ranks,
\begin{equation}
\label{eq:RankEquality}
r_\sE=r_\sK+r_\sJ.
\end{equation}
By combining the preceding equality with \eqref{eq:RankInequality1} we see that $0<r_\sJ< r_\sE$ and therefore, Definition \ref{eq:Bradlow_mu_M} yields
\begin{equation}
\label{eq:UpperBoundOnMuOfImage}
\mu(\sJ)\leq\mu_M(\sE).
\end{equation}
However, the endomorphism $w$ defines an isomorphism $\sE/\sK\to \sJ$ and thus
\begin{equation}
\label{eq:SlopeEquality}
\mu(\sE/\sK)=\mu(\sJ).
\end{equation}
By combining the preceding equality with \eqref{eq:Bradlow_lambda_stable}, \eqref{eq:LowerBoundOnMuOfImage}, and \eqref{eq:UpperBoundOnMuOfImage}, we obtain
\[
\mu(\sJ) \leq \mu_M(\sE) <\lambda < \mu_m(\varphi)\leq \mu(\sE/\sK) = \mu(\sJ),
\]
a contradiction. Hence, $w$ must be the zero endomorphism and so $v=\id_E$ and $(\bar\partial_E,\varphi)$ is strongly simple in the sense of Definition \ref{defn:Simple_pair}. 
\end{proof}

When $E$ has rank two, Lemma \ref{lem:Stable_implies_strongly_simple_for_01-pairs} may be strengthened to include all polystable, non-zero-section pairs.

\begin{lem}[Polystable, non-zero-section pairs are strongly simple on complex rank-two vector bundles]
\label{lem:Polystable_implies_strongly_simple_for_01-pairs_on_Rank2_bundles}
Let $E$ be a rank-two, complex vector bundle over a closed, connected, complex K\"ahler manifold $(X,\omega)$. If $(\bar\partial_E,\varphi)$ is a holomorphic pair with $\varphi\not\equiv 0$ that is polystable as in Definition \ref{defn:StablePair_arbitrary_rank}, then it is strongly simple as in Definition \ref{defn:Simple_pair}.
\end{lem}

\begin{proof}
The assertion follows from Lemma \ref{lem:Stable_implies_strongly_simple_for_01-pairs} when $(\bar\partial_E,\varphi)$ is stable and $E$ has arbitrary rank, so we only need to consider the strictly polystable case, where $(\bar\partial_E,\varphi)$ is polystable but not stable as in Definition \ref{defn:StablePair_arbitrary_rank} and where $E$ has rank two. Thus, we assume that $E$ admits a decomposition $E=L_1\oplus L_2$ as a direct sum of complex line bundles and that $\bar\partial_E=\bar\partial_{L_1}\oplus\bar\partial_{L_2}$, where $\bar\partial_{L_i}$ is a holomorphic structure on $L_i$ for $i=1,2$, and $\varphi$ is a section of $L_1$.

We write $\sE$, $\sL_1$, and $\sL_2$ for the coherent sheaves of holomorphic sections of $E$, $L_1$, and $L_2$, respectively. By Definition \ref{defn:StablePair_arbitrary_rank} and Remark \ref{rmk:Stability_of_Split_Pairs}, the pair $(\sE,\varphi)$ is $\lambda$-polystable with $\lambda=\mu(\sL_2)$ and
\begin{equation}
\label{eq:DegreeL1<DegreeL2}
\deg \sL_1 < \deg \sL_2.
\end{equation}
The conclusion thus follows from Item \eqref{item:Rank2_Split_HolomPairsSimpleCase} of Lemma \ref{lem:Rank2_Split_HolomPairsSimple}
\end{proof}

The stability condition is preserved by the $\SL(E)$ action on holomorphic pairs. By analogy with \eqref{eq:Moduli_space_omega-stable_holomorphic_structures}, we define the following subspaces of the moduli space of holomorphic pairs with fixed holomorphic determinant line bundle in \eqref{eq:Moduli_space_holomorphic_pairs}
\begin{subequations}
\label{eq:Moduli_spaces_stable_and_polystable_holomorphic_pairs}
\begin{align}
  \label{eq:Moduli_space_stable_holomorphic_pairs}
  \fM(E,\omega)
  &:=
    \left\{[\bar\partial_E,\varphi]\in \fM(E): \text{$(\bar\partial_A,\varphi)$ is stable}\right\},
  \\
  \label{eq:Moduli_space_polystable_holomorphic_pairs}
  \fM_\ps(E,\omega)
  &:=
    \left\{[\bar\partial_E,\varphi]\in \fM(E): \text{$(\bar\partial_A,\varphi)$ is polystable}\right\}.
\end{align}
\end{subequations}
We refine the definitions \eqref{eq:Moduli_subspaces_holomorphic_pairs} to give the following subspaces of $\fM(E,\omega)$,
\begin{subequations}
  \label{eq:Moduli_subspaces_stable_holomorphic_pairs}
  \begin{align}
    \label{eq:Moduli_space_non-zero-section_stable_holomorphic_pairs}
    \fM^0(E,\omega) &:= \fM(E,\omega)\cap \fM^0(E),
    \\
    \label{eq:Moduli_space_regular_non-zero-section_stable_holomorphic_pairs}
    \fM_\reg^0(E,\omega) &:= \fM^0(E,\omega)\cap \fM_\reg(E),
\end{align}
\end{subequations}
as well as the following subspaces of $\fM_\ps(E,\omega)$:
\begin{subequations}
  \label{eq:Moduli_subspaces_polystable_holomorphic_pairs}
  \begin{align}
  \label{eq:Moduli_space_non-zero-section_polystable_holomorphic_pairs}
  \fM_\ps^0(E,\omega) &:= \fM_\ps(E,\omega)\cap \fM^0(E),
    \\
  \label{eq:Moduli_space_regular_non-zero-section_polystable_holomorphic_pairs}
  \fM_{\ps,\reg}^0(E,\omega) &:= \fM_\ps^0(E,\omega)\cap \fM_\reg(E).
\end{align}
\end{subequations}
We have the following partial analogue of Lemma \ref{lem:Openness_subspace_simple_01-connections}.

\begin{lem}[Openness of moduli subspaces of strongly simple holomorphic pairs]
\label{lem:Openness_subspace_strongly_simple_01_pairs}
If $E$ is a smooth complex vector bundle over a closed, connected, complex K\"ahler manifold $(X,\omega)$, then the following hold:
\begin{enumerate}
\item\label{item:Openness_subspace_strongly_simple_stable_pairs}
The moduli subspace $\fM^0(E,\omega)$ in \eqref{eq:Moduli_space_non-zero-section_stable_holomorphic_pairs} of non-zero-section, stable holomorphic pairs is an open subspace of the moduli space $\fM(E,\omega)$ in \eqref{eq:Moduli_space_stable_holomorphic_pairs} of stable holomorphic pairs and a subset of the moduli space $\fM^{**}(E)$ of strongly simple holomorphic pairs in \eqref{eq:Moduli_space_strongly_simple_holomorphic_pairs}:
\begin{equation}
  \label{eq:Moduli_space_non-zero-section_stable_holomorphic_pairs_subset_strongly_simple_pairs}
  \fM^0(E,\omega) \subset \fM^{**}(E).
\end{equation}
Moreover, the moduli subspace $\fM_\reg^0(E,\omega)$ in \eqref{eq:Moduli_space_regular_non-zero-section_stable_holomorphic_pairs} is open in $\fM(E,\omega)$.

\item\label{item:Openness_subspace_strongly_simple_polystable_rank-2_pairs} If $E$ has complex rank two, then the
subspace $\fM_\ps^0(E,\omega)$ in \eqref{eq:Moduli_space_non-zero-section_polystable_holomorphic_pairs} of non-zero-section, polystable holomorphic pairs is an open subspace of the moduli space $\fM_\ps(E,\omega)$ in \eqref{eq:Moduli_space_polystable_holomorphic_pairs} of polystable holomorphic pairs and a subset of the moduli space $\fM^{**}(E)$ of strongly simple holomorphic pairs in \eqref{eq:Moduli_space_strongly_simple_holomorphic_pairs}:
\begin{equation}
  \label{eq:Moduli_space_non-zero-section_polystable_rank-2_holomorphic_pairs_subset_strongly_simple_pairs}
  \fM_\ps^0(E,\omega) \subset \fM^{**}(E).
\end{equation}
Moreover, the moduli subspace $\fM_{\ps,\reg}^0(E,\omega)$ in \eqref{eq:Moduli_space_regular_non-zero-section_polystable_holomorphic_pairs} is open in $\fM_\ps(E,\omega)$.
\end{enumerate}
\end{lem}

\begin{proof}
Consider Item \eqref{item:Openness_subspace_strongly_simple_stable_pairs}. Lemma  \ref{lem:Stable_implies_strongly_simple_for_01-pairs} implies that if $(\bar\partial_E,\varphi)$ is a non-zero-section, stable holomorphic pair, then it is strongly simple. This yields the inclusion $\fM^0(E,\omega) \subset \fM^{**}(E)$ in \eqref{eq:Moduli_space_non-zero-section_stable_holomorphic_pairs_subset_strongly_simple_pairs}. The moduli subspaces $\fM^0(E,\omega)$ in \eqref{eq:Moduli_space_non-zero-section_stable_holomorphic_pairs} and $\fM_\reg^0(E,\omega)$ in \eqref{eq:Moduli_space_regular_non-zero-section_stable_holomorphic_pairs} are open in $\fM(E,\omega)$, since the moduli subspaces $\fM^0(E)$ in \eqref{eq:Moduli_space_holomorphic_pairs_NonZeroSection} and $\fM_\reg(E)$ in \eqref{eq:Moduli_space_regular_holomorphic_pairs} are open in $\fM(E)$.

Consider Item \eqref{item:Openness_subspace_strongly_simple_polystable_rank-2_pairs}. Lemma \ref{lem:Polystable_implies_strongly_simple_for_01-pairs_on_Rank2_bundles} implies that if $E$ has complex rank two and $(\bar\partial_E,\varphi)$ is a non-zero-section, polystable holomorphic pair, then it is strongly simple. This yields the inclusion $\fM_\ps^0(E,\omega) \subset \fM^{**}(E)$ given by
\eqref{eq:Moduli_space_non-zero-section_polystable_rank-2_holomorphic_pairs_subset_strongly_simple_pairs}. The subspaces $\fM_\ps^0(E,\omega)$ in \eqref{eq:Moduli_space_non-zero-section_polystable_holomorphic_pairs} and $\fM_{\ps,\reg}^0(E,\omega)$ in \eqref{eq:Moduli_space_regular_non-zero-section_polystable_holomorphic_pairs}
are open in $\fM_\ps(E,\omega)$ by the same argument as used in Item \eqref{item:Openness_subspace_strongly_simple_stable_pairs}.
\end{proof}

Naturally, one would expect an analogue of Theorem \ref{thm:Kobayashi_7_3_34_stable_hausdorff}, asserting that $\fM(E,\omega)$ is a Hausdorff topological space. We shall cite an assertion to this effect (and other properties) in the forthcoming Theorem \ref{thm:Lubke_Teleman_6-3-7}, whose relies on the Hitchin--Kobayashi correspondence that we discuss next (see the forthcoming Theorem \ref{thm:HitchinKobayashiCorrespondenceForPairs}).

\section{Moduli spaces of projective vortices}
\label{sec:Moduli_space_projective_vortices}
For a smooth Hermitian vector bundle $(E,h)$ over a smooth manifold $X$ of dimension $d$ and a constant $p\in (d/2,\infty)$, we consider the affine space of \emph{unitary pairs} $\sA(E,h)\times W^{1,p}(E)$
\label{page:Unitary_pairs} consisting of pairs of $W^{1,p}$ sections $\varphi$ of $E$ and $W^{1,p}$ unitary connections $A$ on $E$ that induce a fixed smooth, unitary connection $A_d$ on the Hermitian line bundle $\det E$, as in \eqref{eq:Unitary_connection_detE_fixed}. The group of special unitary gauge transformations, $W^{2,p}(\SU(E))$, acts on this affine space by
\begin{multline}
\label{eq:DefineSU(E)ActionOnUnitaryPairs}
W^{2,p}(\SU(E))\times \sA(E,h)\times W^{1,p}(E)
\ni
\left( u, (A,\varphi)\right)
\\
\mapsto u^*(A,\varphi) := \left(u^*A,u^{-1}\varphi\right)
\in \sA(E,h)\times W^{1,p}(E),
\end{multline}
where the connection $u^*A$ has covariant derivative $\nabla_{u^*A}=u^{-1}\circ\nabla_A\circ u$
as discussed in Remark \ref{rmk:PushforwardPullbackNotation}. We will consider subspaces of the configuration space of $W^{1,p}$ unitary pairs endowed with the quotient topology,
\begin{equation}
\label{eq:ConfigurationSpaceForProjectiveVortices}
 \sC(E,h):=  \left.\left(\sA(E,h)\times W^{1,p}(E)\right)\right/W^{2,p}(\SU(E)),
\end{equation}
where the action of $W^{2,p}(\SU(E))$ in \eqref{eq:ConfigurationSpaceForProjectiveVortices} is defined in
\eqref{eq:DefineSU(E)ActionOnUnitaryPairs}.
We write elements of $\sC(E,h)$ as $[A,\varphi]$.
\label{page:GaugeEquivClassOfUnitaryPair}

The \emph{stabilizer} of a unitary pair $(A,\varphi)\in \sA(E,h)\times W^{1,p}(E)$ in $W^{2,p}(\SU(E))$ is
\begin{equation}
\label{eq:DefineStabilizerOfUnitaryPair}
\Stab(A,\varphi) := \left\{u\in W^{2,p}(\SU(E)): u^*(A,\varphi)=(A,\varphi)\right\}.
\end{equation}
We note the following Lie group structure of $\Stab(A,\varphi)$.

\begin{lem}[Lie group structure of $\Stab(A,\varphi)$]
\label{lem:LieGroupStructureOfStab(A,varhi)}
Let $(E,h)$ be a smooth Hermitian vector bundle over a smooth manifold $X$ of dimension $d$ and let $p\in (d/2,\infty)$. If $(A,\varphi)\in \sA(E,h)\times W^{1,p}(E)$, then the stabilizer $\Stab(A,\varphi)$ is a Lie group with Lie algebra given by the harmonic space $\bH_{A,\varphi}^0$ defined in \eqref{eq:H_Avarphi^bullet}.
\end{lem}
 
\begin{proof}
The proof is identical to that of Lemma \ref{lem:LieGroupStructureOfStab(A,Phi)}.
\end{proof}
 
We have the following analogue of Definition \ref{defn:Split_trivial_central-stabilizer_spinor_pair} \eqref{item:Split_spinor_pair} for non-Abelian monopoles and Definition \ref{defn:Reducible_split_trivial-stabilizer_unitary_connection} \eqref{item:Split_unitary_connection} for unitary connections.

\begin{defn}[Split, trivial-stabilizer, and central-stabilizer unitary pairs]
\label{defn:Split_trivial_central-stabilizer_unitary_pair}
(Compare Feehan and Leness \cite[Definition 2.2, p. 64]{FL2a} for the case $r=2$.)
Let $(E,h)$ be a smooth Hermitian vector bundle of complex rank $r\geq 2$ over a smooth, connected Riemannian manifold of real dimension four and let $(A,\varphi)\in\sA(E,h)\times W^{1,p}(E)$.
\begin{enumerate}
\item\label{item:Trivial-stabilizer_unitary_pair}
$(A,\varphi)$ has \emph{trivial stabilizer} if the stabilizer group $\Stab(A,\varphi)$ of $(A,\varphi)$ in $W^{2,p}(\SU(E))$ is $\{\id_E\}$ and has \emph{non-trivial stabilizer} otherwise.

\item\label{item:central-stabilizer_unitary_pair}
$(A,\varphi)$ has \emph{central stabilizer} if the stabilizer group $\Stab(A,\varphi)$ of $(A,\varphi)$ in $W^{2,p}(\SU(E))$ is equal to the center $Z(\SU(r)) = C_r$ and has \emph{non-central stabilizer} otherwise.

\item\label{item:Split_unitary_pair}
 $(A,\varphi)$ is \emph{split} if the connection $A$ is split as in Definition \ref{defn:Reducible_split_trivial-stabilizer_unitary_connection} \eqref{item:Split_unitary_connection} and $\varphi\in W^{1,p}( E_1)$. If no such splitting exists, then $(A,\varphi)$ is called \emph{non-split}. We say that $(A,\varphi)$ is \emph{split with respect to the decomposition $E=E_1\oplus E_2$} in \eqref{eq:BasicSplitting} when we wish to specify the bundle splitting.

\item\label{item:Zero-section_unitary_pair}
$(A,\varphi)$ is a \emph{zero-section unitary pair} if $\varphi\equiv 0$.  
\end{enumerate}
\end{defn}

\begin{rmk}[Central stabilizer implies $\bH_{A,\varphi}^0=(0)$ for unitary pairs]
\label{rmk:Central_stabilizer_unitary_pairs_and_H0=0}
If a pair $(A,\varphi)$ has central stabilizer as in Definition \ref{defn:Split_trivial_central-stabilizer_unitary_pair} \eqref{item:central-stabilizer_unitary_pair}, then $\Stab(A,\varphi) = C_r\,\id_E$, a discrete subgroup of $W^{2,p}(\SU(E))$. The Lie algebra of $\Stab(A,\varphi)$ is $\bH_{A,\varphi}^0$ by Lemma \ref{lem:LieGroupStructureOfStab(A,varhi)}. Consequently, if $(A,\varphi)$ has central stabilizer, then $\bH_{A,\varphi}^0=(0)$ since a Lie group is discrete if and only if its Lie algebra is zero.
\end{rmk}

We note the following restriction on $\Stab(A,\varphi)$ when the bundle $E$ has rank two.

\begin{lem}[Stabilizers of non-zero-section unitary pairs on rank-two Hermitian vector bundles]
\label{lem:Stabilizers_non-zero-section_unitary_pairs}
Let $(E,h)$ be a rank-two, Hermitian vector bundle over a connected smooth manifold $X$ of real dimension $d$ and $p\in(d/2,\infty)$ be a constant. If $A$ is a $W^{1,p}$ unitary connection on $E$, and $\varphi\not\equiv 0$ is a $W^{1,p}$ section of $E$, and $u\in W^{2,p}(\SU(E))$ is such that $u(A,\varphi) = (A,\varphi)$ via the action \eqref{eq:DefineSU(E)ActionOnUnitaryPairs}, then $u=\id_E$.
\end{lem}

\begin{proof}
The proof is identical to that of Lemma \ref{lem:NonZeroSection_Spinu_pairs_Have_Trivial_Stabilizer}.
\end{proof}

\begin{rmk}[On the necessity of the assumption that $E$ has rank two in Lemma \ref{lem:Stabilizers_non-zero-section_unitary_pairs}]
\label{lem:NonTrivialDiscreteStabilizer}
We give an example of a unitary pair $(A,\varphi)$ with non-trivial but discrete stabilizer. Let $(E,h)$ be a  rank-three Hermitian vector bundle and let $A\in\sA(E,h)$ admit a splitting $A=A_1\oplus A_2$ with respect to an orthogonal decomposition $E=L_1\oplus E_2$, where $L_1$ is a complex line bundle and $E_2$ a rank-two Hermitian vector bundle.  Let $A_1$ be a unitary connection on $L_1$ and let $A_2$ be a unitary connection with central stabilizer on $E_2$.  Thus,
\[
\Stab(A)
=
\left\{ e^{2i\theta}\,\id_{L_1} \oplus e^{-i\theta}\,\id_{E_2}: e^{i\theta}\in S^1\right\}.
\]
If $\varphi\not\equiv 0$ is a section of $L_1\subset E$, then $\Stab(A,\varphi)\subset\Stab(A)$ is given
by
\[
\Stab(A,\varphi)
=
\{ \id_{L_1} \oplus (\pm\id_{E_2})\}.
\]
Thus, although $\Stab(A,\varphi)$ is non-trivial, it is discrete and so its Lie algebra $\bH_{A,\varphi}^0$ is zero.  Hence, the converse of the implication discussed in Remark \ref{rmk:Central_stabilizer_unitary_pairs_and_H0=0}, that if $(A,\varphi)$ has central stabilizer then $\bH_{A,\varphi}^0=(0)$, does not hold when the rank of $E$ is greater than two.
\end{rmk}

When $(X,g,J)$ is a smooth almost Hermitian manifold with fundamental two-form $\omega = g(\cdot,J\cdot)$ as in \eqref{eq:Fundamental_two-form}, we let
\begin{equation}
  \label{eq:Moduli_space_projective_vortices}
  \sM(E,h,\omega) := \left\{(A,\varphi) \in \sA(E,h)\times W^{1,p}(E): \eqref{eq:SO(3)_monopole_equations_almost_Hermitian_alpha} \text{ hold} \right\}/W^{2,p}(\SU(E))
\end{equation}
denote the \emph{moduli space of projective vortices} modulo the Banach Lie group $W^{2,p}(\SU(E))$ of determinant-one, unitary automorphisms of $E$ of class $W^{2,p}$, noting that the $(1,1)$-determinant condition \eqref{eq:Einstein_connection_trace_part} is obeyed by virtue of our choice that $A$ obeys \eqref{eq:Unitary_connection_detE_fixed} on $\det E$.

Just as it does for the standard circle action \eqref{eq:S1ZAction} on the affine space $\sA(E,h)\times W^{1,p}(W^+\otimes E)$ of \spinu pairs, scalar multiplication on the section defines the \emph{standard circle action on the affine space of unitary pairs}, 
\begin{equation}
\label{eq:S1_Action_On_AffineSpaceForProjectiveVortices}
S^1\times \left( \sA(E,h)\times W^{1,p}(E)\right) \ni \left(e^{i\theta},(A,\varphi)\right)
\mapsto
(A,e^{i\theta}\varphi) \in \sA(E,h)\times W^{1,p}(E).
\end{equation}
Because the circle action \eqref{eq:S1_Action_On_AffineSpaceForProjectiveVortices} commutes\footnote{See  the discussion around tom Dieck \cite[Proposition 3.4, p. 23]{tomDieck} and \cite[Exercise 17, p. 31]{tomDieck} for further information about when group actions descend to a quotient space.} with the action \eqref{eq:DefineSU(E)ActionOnUnitaryPairs} of $W^{2,p}(\SU(E))$, it induces an action
on the configuration space \eqref{eq:ConfigurationSpaceForProjectiveVortices},
\begin{equation}
\label{eq:S1_Action_On_ConfigurationSpaceForProjectiveVortices}
S^1\times \sC(E,h) \ni \left(e^{i\theta},[A,\varphi]\right)
\mapsto
[A,e^{i\theta}\varphi] \in \sC(E,h),
\end{equation}
which we refer to as \emph{standard circle action on the quotient space} of unitary pairs. We note that $\sM(E,h,\omega)$ is closed under the action \eqref{eq:S1_Action_On_ConfigurationSpaceForProjectiveVortices}. The proof of Proposition \ref{prop:FixedPointsOfS1ActionOnSpinuQuotientSpace}, which characterizes the fixed points of the $S^1$ action in Definition \ref{defn:UnitaryZActionOnAffine} on the quotient space $\sC_\ft$ of \spinu pairs, adapts \mutatis to give the

\begin{prop}[Fixed points of the $S^1$ action in $\sC(E,h)$]
\label{prop:FixedPointsOfS1ActionOnUnitaryQuotientSpace}
Let $(E,h)$ be a smooth Hermitian vector bundle over a connected, smooth manifold $X$ of dimension $d$, and $p\in (d/2,\infty)$ be a constant, and $A_d$ be a fixed smooth, unitary connection on the Hermitian line bundle $\det E$.   If $[A,\varphi] \in \sC(E,h)$ obeys
\begin{enumerate}
\item
\label{item:UnitaryPairFixedPointsZeroSection}
$\varphi\equiv 0$, so $(A,0)$ is a zero-section pair as in Definition \ref{defn:Split_trivial_central-stabilizer_unitary_pair} \eqref{item:Zero-section_unitary_pair}, or
\item
\label{item:UnitaryPairFixedPointsSplit}
$(A,\varphi)$ is a split unitary pair as in Definition \ref{defn:Split_trivial_central-stabilizer_unitary_pair} \eqref{item:Split_unitary_pair},
\end{enumerate}
then $[A,\varphi]$ is a fixed point of the $S^1$ action \eqref{eq:S1_Action_On_ConfigurationSpaceForProjectiveVortices} on $\sC(E,h)$. Conversely, if a fixed point of that $S^1$ action on $\sC(E,h)$ is represented by smooth pair $(A,\Phi)$, then it obeys Condition \eqref{item:SO(3)MonopoleFixedPointsZeroSection} or \eqref{item:SO(3)MonopoleFixedPointsReducible}. 
\end{prop}

We define the following open subspaces,
\begin{subequations}
\label{eq:Moduli_space_projective_vortices_non-split_non-zero-section}  
\begin{align}
  \label{eq:Moduli_space_projective_vortices_StabAvarphi_idE}
  \sM^{**}(E,h,\omega) &:= \left\{(A,\varphi) \in \sM(E,h,\omega): \Stab(A,\varphi)=\{\id_E\} \right\},
  \\
  \label{eq:Moduli_space_projective_vortices_non-split}
  \sM^*(E,h,\omega) &:= \left\{[A,\varphi] \in \sM(E,h,\omega): A \text{ is non-split} \right\},
  \\
  \label{eq:Moduli_space_projective_vortices_non-zero_section}
  \sM^0(E,h,\omega) &:= \left\{[A,\varphi] \in \sM(E,h,\omega): \varphi\not\equiv 0 \right\},
  \\
  \label{eq:Moduli_space_projective_vortices_non-split-non-zero-section}
  \sM^{*,0}(E,h,\omega) &:= \sM^*(E,h,\omega) \cap \sM^0(E,h,\omega),
\end{align}
\end{subequations}
and, when $(X,g,J)$ is complex Hermitian,
\begin{subequations}
\label{eq:Moduli_space_projective_vortices_non-split_non-zero-section_regular}  
\begin{align}
  \label{eq:Moduli_space_projective_vortices_regular}
  \sM_\reg(E,h,\omega) &:= \left\{[A,\varphi] \in \sM(E,h,\omega): \bH_{A,\varphi}^2 = (0) \right\},
  \\
  \label{eq:Moduli_space_projective_vortices_StabAvarphi_idE_regular}
  \sM_\reg^{**}(E,h,\omega) &:= \sM^{**}(E,h,\omega) \cap \sM_\reg(E,h,\omega),
  \\
  \label{eq:Moduli_space_projective_vortices_non-split-non-zero-section_regular}
  \sM_\reg^{*,0}(E,h,\omega) &:= \sM^*(E,h,\omega) \cap \sM^0(E,h,\omega) \cap \sM_\reg(E,h,\omega),
  \\                         
  \label{eq:Moduli_space_projective_vortices_non-zero_section_regular}
  \sM_\reg^0(E,h,\omega) &:= \sM^0(E,h,\omega) \cap \sM_\reg(E,h,\omega),
\end{align}
\end{subequations}
where the harmonic space $\bH_{A,\varphi}^2$ in \eqref{eq:H_Avarphi^2} is defined by the elliptic complex  \eqref{eq:Projective_vortex_elliptic_deformation_complex} for the projective vortex equations \eqref{eq:SO(3)_monopole_equations_almost_Hermitian_alpha}. The verification that $\sM^0(E,h,\omega)$ in \eqref{eq:Moduli_space_projective_vortices_non-zero_section} and $\sM_\reg(E,h,\omega)$ in \eqref{eq:Moduli_space_projective_vortices_regular} are open subspaces of $\sM(E,h,\omega)$ is the same as that provided for openness of the subspace $\fM^0(E)$ and $\fM_\reg(E)$ of $\fM(E)$ in \eqref{eq:Moduli_subspaces_holomorphic_pairs}. The proof of Lemma \ref{lem:Openness_moduli_subspace_non-split_projectively_Hermitian-Einstein_connections} adapts \mutatis to yield the

\begin{lem}[Openness of the configuration subspace of unitary pairs and moduli subspace of projectively vortices with trivial stabilizer]
\label{lem:Openness_configuration_subspace_unitary_pairs_trivial_stabilizer}
Let $(E,h)$ be a smooth Hermitian vector bundle over a connected, smooth manifold $X$ of dimension $d$, and $p\in (d/2,\infty)$ be a constant, and $A_d$ be a fixed smooth, unitary connection on the Hermitian line bundle $\det E$. Then the following hold:
\begin{enumerate}
\item\label{item:Openness_subspace_unitary_pairs_minimal_stabilizer}
The subspace $(\sA(E,h)\times W^{1,p}(E))^{**}$ of pairs in $\sA(E,h)\times W^{1,p}(E)$ with minimal stabilizer, $\{\id_E\}$, is open in $\sA(E,h)\times W^{1,p}(E)$.

\item\label{item:Openness_quotient_subspace_unitary_pairs_minimal_stabilizer}
The quotient subspace,
\begin{equation}
\label{eq:ConfigurationSpaceForProjectiveVortices**}
 \sC^{**}(E,h) := \left\{[A,\varphi] \in \sC(E,h): \Stab(A,\varphi) = \{\id_E\}\right\},
\end{equation}
of gauge-equivalence classes of $W^{1,p}$ unitary pairs with minimal stabilizer, $\{\id_E\}$, is open in $\sC(E,h)$ equipped with the quotient topology.

\item\label{item:Openness_moduli_subspace_non-split_projective_vortices}
The moduli subspace,
\[
  \sM^{**}(E,h,\omega) = \sM(E,h,\omega)\cap \sC^{**}(E,h),
\]
in \eqref{eq:Moduli_space_projective_vortices_StabAvarphi_idE} is an open subspace of $\sM(E,h,\omega)$ in \eqref{eq:Moduli_space_projective_vortices}.
\end{enumerate}
\end{lem}

\begin{proof}
The first and second proofs of Lemma \ref{lem:Openness_moduli_subspace_non-split_projectively_Hermitian-Einstein_connections}
\eqref{item:Openness_subspace_unitary_connections_minimal_stabilizer} adapt \mutatis to prove Item \eqref{item:Openness_subspace_unitary_pairs_minimal_stabilizer}, noting that if $(A,\varphi) \in \sA(E,h)\times W^{1,p}(E)$, then
\[
  \Stab(A,\varphi) = \Stab(A)\cap\Stab(r_\varphi),
\]
where (by partial analogy with \eqref{eq:Phi_Omega0suE_to_Omega0V+}) we define
\[
  r_\varphi:W^{2,p}(\SU(E)) \ni u \mapsto u\varphi \in W^{1,p}(E), 
\]
keeping in mind the definition \eqref{eq:DefineStabilizerOfUnitaryPair} of $\Stab(A,\varphi)$ in terms of the action \eqref{eq:DefineSU(E)ActionOnUnitaryPairs} of $W^{2,p}(\SU(E))$ on $\sA(E,h)\times W^{1,p}(E)$.

The proof of Item \eqref{item:Openness_quotient_subspace_unitary_pairs_minimal_stabilizer} follows \mutatis the proof of the corresponding result in Lemma \ref{lem:Openness_moduli_subspace_non-split_projectively_Hermitian-Einstein_connections}, while Item \eqref{item:Openness_moduli_subspace_non-split_projective_vortices} is an immediate consequence of Item \eqref{item:Openness_quotient_subspace_unitary_pairs_minimal_stabilizer}.
\end{proof}

Next, we have the

\begin{lem}[Openness of the moduli subspace of projective vortices with non-split connections]
\label{lem:Openness_moduli_subspace_projective_vortices_non-split_connections}
Let $(E,h)$ be a smooth Hermitian vector bundle over a smooth almost Hermitian manifold $(X,g,J)$ of real dimension $2n$, let $A_d$ be a fixed smooth, unitary connection on the Hermitian line bundle $\det E$, and let $p\in (n,\infty)$ be a constant. Then $\sM^*(E,h,\omega)$ in \eqref{eq:Moduli_space_projective_vortices_non-split} is an open subspace of $\sM(E,h,\omega)$ in \eqref{eq:Moduli_space_projective_vortices} equipped with the quotient topology.
\end{lem}

\begin{proof}
Let $[A,\varphi] \in \sM^*(E,h,\omega)$, so $(A,\varphi)$ is a $W^{1,p}$ projective vortex and $A$ is non-split. 
By Lemma \ref{lem:Regularity_For_ProjVortices}, we may assume without loss of generality that $(A,\varphi)$ is smooth since, if not, it is equivalent via $u\in W^{2,p}(\SU(E))$-gauge transformation to a smooth projective vortex $u^*(A,\varphi)$. According to Items \eqref{item:A_smooth_and_bHA_non-zero_implies_A_split} and \eqref{item:A_split_implies_bHA_non-zero} in Corollary \ref{cor:Split_unitary_A_and_Lie_Alg_of_Stab(A)}, the condition that $A$ be non-split is equivalent to $\bH_A^0 = (0)$, where $\bH_A^0$ is defined as the finite-dimensional kernel of $d_A:W^{2,p}(\su(E))\to W^{1,p}(T^*X\otimes \su(E))$ in \eqref{eq:HE_equation_bHA0}. It follows from Rudin \cite[Theorem 4.12, p. 99]{Rudin} that
\begin{multline*}
  \Ker\left(d_A:W^{2,p}(\su(E))\to W^{1,p}(T^*X\otimes \su(E))\right)
  \\
  =
  \Ran\left(d_A^*:W^{3,p}(T^*X\otimes \su(E))\to W^{2,p}(\su(E))\right)^\perp,
\end{multline*}
where $\perp$ denotes $L^2$-orthogonal complement. Thus, $\bH_A^0 = (0)$ if and only if $\Ran d_A^* = W^{2,p}(\su(E))$ and the latter is an open condition in the sense that if $a \in W^{1,p}(T^*X\otimes \su(E))$ is small, then $\Ran d_{A+a}^* = W^{2,p}(\su(E))$ and this yields the conclusion.
\end{proof}  

We have the following inclusions and equalities of moduli subspaces.

\begin{lem}[Inclusions and equalities of moduli subspaces of projective vortices]
\label{lem:sM*0(E,h,omega)_subset_sM**(E,h,omega)}
Let $(E,h)$ be a smooth Hermitian vector bundle over a connected, smooth, almost Hermitian manifold $(X,g,J)$ of real dimension $2n$ and fundamental two-form $\omega = g(\cdot,J\cdot)$ as in \eqref{eq:Fundamental_two-form}, and $p\in (n,\infty)$ be a constant, and $A_d$ be a fixed smooth, unitary connection on the Hermitian line bundle $\det E$. Then
\begin{equation}
  \label{eq:sC**(E,h)_subset_sC0(E,h)}
  \sC^{**}(E,h) \subset \sC^0(E,h),
\end{equation}
where $\sC^{**}(E,h)$ is as in \eqref{eq:ConfigurationSpaceForProjectiveVortices**} and
\begin{equation}
\label{eq:ConfigurationSpaceForProjectiveVorticesnon-zero-section}
 \sC^0(E,h) := \left\{[A,\varphi] \in \sC(E,h): \varphi \not\equiv 0\right\}.
\end{equation}
Moreover, the following inclusion and equality hold:
\begin{enumerate}
\item\label{item:Subspace_trivial_stabilizer_equals_non-split_non-zero-section_projective vortices}
The moduli subspace $\sM^{*,0}(E,h,\omega)$ in \eqref{eq:Moduli_space_projective_vortices_non-split-non-zero-section} of non-split, non-zero-section, projective vortices is a subset of the moduli subspace $\sM^{**}(E,h,\omega)$ in \eqref{eq:Moduli_space_projective_vortices_StabAvarphi_idE} of projective vortices with trivial stabilizer:
\begin{equation}
  \label{eq:sM*0(E,h,omega)_subset_sM**(E,h,omega)}
  \sM^{*,0}(E,h,\omega) \subset \sM^{**}(E,h,\omega).
\end{equation}

\item\label{item:Subspace_trivial_stabilizer_equals_non-zero-section_rank-2_unitary_pairs}
If $E$ has complex rank two, then the quotient subspace $\sC^0(E,h)$ in \eqref{eq:ConfigurationSpaceForProjectiveVorticesnon-zero-section} of non-zero-section unitary pairs is equal to the moduli subspace $\sC^{**}(E,h)$ in \eqref{eq:ConfigurationSpaceForProjectiveVortices**} of unitary pairs with trivial stabilizer:
\begin{equation}
  \label{eq:sC0(E,h)_equals_sC**(E,h)_E_rank2}
  \sC^0(E,h) = \sC^{**}(E,h).
\end{equation}
\end{enumerate}
\end{lem}

\begin{proof}
For any unitary pair $(A,\varphi)$ with $\Stab(A,\varphi)=\{\id_E\}$ we must have $\varphi\not\equiv 0$ because a zero-section pair $(A,0)$ has stabilizer equal to $\Stab(A)$ in $W^{2,p}(\SU(E))$ and $\Stab(A)$ always contains a subgroup isomorphic to the center $Z(\SU(r)) = C_r$, the group of $r$-th roots of unity, when $E$ has complex rank $r$. This observation yields the inclusion \eqref{eq:sC**(E,h)_subset_sC0(E,h)}.

Consider Item \eqref{item:Subspace_trivial_stabilizer_equals_non-split_non-zero-section_projective vortices}. Let $[A,\varphi]\in \sM^{*,0}(E,h,\omega)$, so that $A$ is not split and $\varphi\not\equiv 0$. We may assume without loss of generality that the pair $(A,\varphi)$ is smooth by Lemma \ref{lem:Regularity_For_ProjVortices}. If $A$ is not split, then Corollary
\ref{cor:Split_unitary_A_and_Lie_Alg_of_Stab(A)} \eqref{item:A_smooth_and_bHA_non-zero_implies_A_split} implies that $\bH_A^0=(0)$. Thus, Corollary \ref{cor:Split_unitary_A_and_Lie_Alg_of_Stab(A)} \eqref{item:HA0_is_Lie_Alg_of_Stab(A)_Irreducible} implies that $A$ has central stabilizer in the sense of
Definition \ref{defn:Reducible_split_trivial-stabilizer_unitary_connection} and so $\Stab(A) = C_r\,\id_E$, where $C_r$ is the group of $r$-th roots of unity.  If $u\in\Stab(A,\varphi)$, then the inclusion
$\Stab(A,\varphi)\subset\Stab(A)$ implies that $u=\varrho\,\id_E$, for some $\varrho\in C_r$. Because $\varphi\not\equiv 0$, the identity $u\varphi = \varphi$ forces $\varrho=1$.  Hence, $\Stab(A,\varphi)=\{\id_E\}$ and so $[A,\varphi]\in \sM^{**}(E,h,\omega)$. Since the point $[A,\varphi]\in \sM^{*,0}(E,h,\omega)$ was arbitrary, this yields the inclusion \eqref{eq:sM*0(E,h,omega)_subset_sM**(E,h,omega)}.

Consider Item \eqref{item:Subspace_trivial_stabilizer_equals_non-zero-section_rank-2_unitary_pairs}. If $[A,\varphi]\in \sC^0(E,h)$, and $E$ has complex rank two, then Lemma \ref{lem:Stabilizers_non-zero-section_unitary_pairs} implies that $(A,\varphi)$ has stabilizer $\Stab(A,\varphi)=\{\id_E\}$ in $W^{2,p}(\SU(E))$, and therefore
\[
  \sC^0(E,h,\omega) \subset \sC^{**}(E,h,\omega).
\]  
Thus, by combining the preceding inclusion with the inclusion \eqref{eq:sC**(E,h)_subset_sC0(E,h)}, we obtain the equality \eqref{eq:sC0(E,h)_equals_sC**(E,h)_E_rank2}. This completes the proof of the lemma.
\end{proof}

We have the following analogue of Donaldson and Kronheimer \cite[Proposition 4.2.9, p. 132]{DK}, Freed and Uhlenbeck \cite[Theorem 3.2, p. 49, and Corollary, p. 50]{FU}, or Lawson \cite[Chapter II, Theorem 10.4, p. 33]{Lawson} for unitary pairs rather than unitary connections, real analytic Banach manifolds rather than smooth Hilbert manifolds, and $X$ of real dimension $d \geq 2$ rather than dimension four. See also Feehan and Maridakis \cite[Theorem 16, p. 18 and Corollary 18, p. 19]{Feehan_Maridakis_Lojasiewicz-Simon_coupled_Yang-Mills}.

\begin{thm}[Topology of the quotient space of unitary pairs with trivial stabilizer over a Riemannian manifold]
\label{thm:DK_prop_5-2-9_FU_corollary_page_50_unitary_pairs}
Let $(E,h)$ be a Hermitian vector bundle over a smooth Riemannian manifold $(X,g)$ and $A_d$ be a fixed smooth, unitary connection on the Hermitian line bundle $\det E$. If $p\in(d/2,\infty)$, where $d$ is the real dimension of $X$, and $(A,\varphi) \in \sA(E,h)\times W^{1,p}(E)$ has trivial stabilizer as in Definition \ref{defn:Split_trivial_central-stabilizer_unitary_pair} (\ref{item:Trivial-stabilizer_unitary_pair}), so $\Stab(A,\varphi) = \{\id_E\}$ and
\begin{equation}
  \label{eq:Unitary_pair_Avarphi_linear_slice}
  S_{A,\varphi}
  := (A,\varphi) + \Ker d_{A,\varphi}^{0,*}\cap W^{1,p}\left(T^*X\otimes\su(E)\oplus E\right),
\end{equation}
where $d_{A,\varphi}^0$ is as in \eqref{eq:d0_projective_vortex_elliptic_deformation_complex} and $d_{A,\varphi}^{*,0}$ is its $L^2$ adjoint, then there is an open neighborhood $U_{\id_E} \subset W^{2,p}(\SU(E))$ of the identity $\id_E$ such that the natural map,
\begin{equation}
  \label{eq:Avarphi_slice_real_analytic_diffeomorphism}
  U_{A,\varphi} \times U_{\id_E} \ni \left(A + a, \varphi + \phi, u\right)
  \mapsto u^*(A + a, \varphi + \phi) \in \UU_{A,\varphi},
\end{equation} 
is a real analytic embedding onto an open neighborhood $\UU_{A,\varphi} \subset \sA(E,h)\times W^{1,p}(E)$ of $(A,\varphi)$ and the quotient map,
\begin{equation}
  \label{eq:Unitary_pair_Avarphi_linear_slice_neighborhood_quotient_map}
  \pi:U_{A,\varphi} \ni (a,\phi) \mapsto [A + a, \varphi+\phi] \in \sC^{**}(E,h),
\end{equation}
gives a homeomorphism from an open neighborhood $U_{A,\varphi}$ of the pair $(A,\varphi)$ in $S_{A,\varphi}$ onto an open neighborhood of $[A,\varphi] = \pi(A,\varphi) \in \sC^{**}(E,h)$. Moreover, $\sC^{**}(E,h)$ is a Hausdorff, real analytic Banach manifold. 
\end{thm}

The proof of Theorem \ref{thm:Local_Kuranishi_model_for_simple_point_cH(E)} yields \mutatis the following analogue for projective vortices.

\begin{thm}[Local Kuranishi models for points in $\sM^{**}(E,h,\omega)$]
\label{thm:Local_Kuranishi_model_for_moduli_space_projective_vortices_StabAvarphi_idE}
Let $(E,h)$ be a Hermitian vector bundle over a complex Hermitian manifold $(X,g,J)$ and $A_d$ be a fixed smooth, unitary connection on the Hermitian line bundle $\det E$. If $p\in(n,\infty)$, where $n$ is the complex dimension of $X$, and $(A,\varphi)$ is a solution to the projective vortex equations \eqref{eq:SO(3)_monopole_equations_almost_Hermitian_alpha} with trivial stabilizer as in Definition \ref{defn:Split_trivial_central-stabilizer_unitary_pair} (\ref{item:Trivial-stabilizer_unitary_pair}), so $\Stab(A,\varphi) = \{\id_E\}$, and $S_{A,\varphi}$ is the affine slice through $(A,\varphi)$ defined in \eqref{eq:Unitary_pair_Avarphi_linear_slice},
\[
  S_{A,\varphi} = (A,\varphi) + \Ker d_{A,\varphi}^{0,*}\cap W^{1,p}\left(T^*X\otimes\su(E)\oplus E\right),
\]
and
\begin{multline}
  \label{eq:Projective_vortex_map}
  \sS:\sA(E,h)\times W^{1,p}(E) \ni (A,\varphi)
  \\
  \mapsto
  \left(\Lambda (F_A)_0 - \frac{i}{2}(\varphi\otimes\varphi^*)_0, (F_A^{0,2})_0, \bar\partial_A\varphi\right)
  \in L^p\left( \su(E) \oplus \Lambda^{0,2}(\fsl(E)) \oplus \Lambda^{0,1}(E) \right)
\end{multline}
is the real analytic \emph{projective vortex map} defined by the projective vortex equations \eqref{eq:SO(3)_monopole_equations_almost_Hermitian_alpha}, then there are open neighborhoods $U_{A,\varphi} \subset S_{A,\varphi}$ of $(A,\varphi)$ and $N_{A,\varphi} \subset \bH_{A,\varphi}^1$ of the origin such that the following hold:
\begin{enumerate}
\item\label{item:Kuranishi_model_Avarphi_moduli_space_projective_vortices_StabAvarphi_idE}
There are a real analytic embedding,
\begin{equation}
  \label{eq:Kuranishi_embedding_map_projective_vortex}
  \beps:\bH_{A,\varphi}^1 \supset N_{A,\varphi} \to U_{A,\varphi} \subset S_{A,\varphi}
\end{equation}
such that $\beps(0,0) = (A,\varphi)$ and a real analytic map
\begin{equation}
  \label{eq:Kuranishi_obstruction_map_projective_vortex}
  \bkappa:\bH_{A,\varphi}^1 \supset N_{A,\varphi} \to \bH_{A,\varphi}^2,
\end{equation}
such that $\bkappa(0,0)=0$ and $D\bkappa(0,0) = 0$ and
\begin{equation}
  \label{eq:Kuranishi_model_Avarphi_moduli_space_projective_vortices_StabAvarphi_idE}
  \sS^{-1}(0) \cap U_{A,\varphi} = \beps(\bkappa^{-1}(0)\cap N_{A,\varphi}).
\end{equation}

\item\label{item:Kuranishi_model_Avarphi_moduli_space_projective_vortices^vir_is_bepsN_Avarphi)}
For $L^2$-orthogonal projection $\Pi_{\Ran d_{A,\varphi}^1}$ onto the range of $d_{A,\varphi}^1$ in
\[
  L^p\left(\su(E)\oplus\Lambda^{0,2}(\fsl(E))\oplus\Lambda^{0,1}(E)\right)
\]
and
\begin{equation}
  \label{eq:Atiyah_Hitchin_Singer_family_page_446_projective_vortices}
  \sS^{-1}(0)_{A,\varphi}^\vir
  :=
  \left\{(A,\varphi) + (a,\phi) \in \sA(E,h)\times W^{1,p}(E): \Pi_{\Ran d_{A,\varphi}^1}\sS(A+a,\varphi+\phi) = 0\right\}
\end{equation}
the subset $\sS^{-1}(0)_{A,\varphi}^\vir \cap U_{A,\varphi}$ is an embedded real analytic submanifold of $U_{A,\varphi}$ with tangent space $T_{A,\varphi}(\sS^{-1}(0)_{A,\varphi}^\vir \cap U_{A,\varphi}) = \bH_{A,\varphi}^1$ and
\begin{equation}
  \label{item:Kuranishi_model_Avarphi_moduli_space_projective_vortices^vir_is_beps(N_Avarphi)}
  \sS^{-1}(0)_{A,\varphi}^\vir \cap U_{A,\varphi} = \beps(N_{A,\varphi}). 
\end{equation}

\item\label{item:Avarphi_local_virtual_moduli_space_projective_vortices}
The quotient map $\pi$ in \eqref{eq:Unitary_pair_Avarphi_linear_slice_neighborhood_quotient_map} yields a homeomorphism from $\sS^{-1}(0)_{A,\varphi}^\vir \cap U_{A,\varphi}$ onto the \emph{local real virtual moduli space of projective vortices},
\begin{equation}
  \label{eq:sM_Avarphi^vir_real}
  \sM_{A,\varphi}^{\vir,\RR}(E,h,\omega)
  :=
  \pi\left(\sS^{-1}(0)_{A,\varphi}^\vir\cap U_{A,\varphi}\right) \subset \sC(E,h),
\end{equation}
and gives $\sM_{A,\varphi}^{\vir,\RR}(E,h,\omega)$ the structure of a real analytic manifold, with Zariski tangent space $\bH_{A,\varphi}^1$ at $[A,\varphi]$ and containing $\sM(E,h,\omega) \cap \pi(U_{A,\varphi})$ as a real analytic subspace.
\end{enumerate}
\end{thm}

Given Theorems \ref{thm:DK_prop_5-2-9_FU_corollary_page_50_unitary_pairs} and \ref{thm:Local_Kuranishi_model_for_moduli_space_projective_vortices_StabAvarphi_idE}, one can quickly show that the method of proof of Theorem \ref{thm:Kobayashi_7_4_19} and the equality \eqref{eq:sC0(E,h)_equals_sC**(E,h)_E_rank2} when $E$ has rank two yield


\begin{thm}[Moduli space of regular projective vortices with trivial stabilizer is a real analytic manifold]
\label{thm:Kobayashi_7_4_19_projective_vortices}
Let $(E,h)$ be a smooth Hermitian vector bundle over a closed, complex Hermitian manifold $(X,g,J)$
of complex dimension $n$ with fundamental two-form $\omega = g(\cdot,J\cdot)$ as in \eqref{eq:Fundamental_two-form}
and let $p\in(n,\infty)$ be a constant. Then the following hold:
\begin{enumerate}
\item\label{item:Kobayashi_7_4_19_projective_vortices_StabAvarphi_idE}
The moduli subspace $\sM^{**}(E,h,\omega)$ in \eqref{eq:Moduli_space_projective_vortices_StabAvarphi_idE} is a (possibly reduced) real analytic space with Zariski tangent spaces $\bH_{A,\varphi}^1$ at points $[A,\varphi]$ and smooth at points at points $[A,\varphi]$ with $\bH_{A,\varphi}^2 = 0$, where the harmonic spaces $\bH_{A,\varphi}^\bullet$ are defined in \eqref{eq:H_Avarphi^bullet} for the elliptic complex \eqref{eq:Projective_vortex_elliptic_deformation_complex} for the projective vortex equations \eqref{eq:SO(3)_monopole_equations_almost_Hermitian_alpha}.

\item\label{item:Kobayashi_7_4_19_projective_vortices_StabAvarphi_idE_regular} The moduli subspace
  $\sM_\reg^{**}(E,h,\omega)$ in \eqref{eq:Moduli_space_projective_vortices_StabAvarphi_idE_regular}is a real analytic manifold with tangent spaces $\bH_{A,\varphi}^1$ at points $[A,\varphi] \in \sM_\reg^{**}(E,h,\omega)$.

\item\label{item:Kobayashi_7_4_19_projective_vortices_rank-2_non-zero-section}
  If $E$ has rank two, then Items \eqref{item:Kobayashi_7_4_19_projective_vortices_StabAvarphi_idE} and \eqref{item:Kobayashi_7_4_19_projective_vortices_StabAvarphi_idE_regular} hold with $\sM^{**}(E,h,\omega)$ and $\sM_\reg^{**}(E,h,\omega)$ replaced by $\sM^0E,h,\omega)$ in \eqref{eq:Moduli_space_projective_vortices_non-zero_section} and $\sM_\reg^0(E,h,\omega)$ in \eqref{eq:Moduli_space_projective_vortices_non-zero_section_regular}, respectively.
\end{enumerate}
\end{thm}

\begin{rmk}[Uhlenbeck compactifications of moduli spaces of projective vortices]
\label{rmk:Uhlenbeck_compactification_moduli_space_vortices}
When $X$ has complex dimension two in Theorem \ref{thm:Kobayashi_7_4_19_projective_vortices}, then existence of an Uhlenbeck compactification for the moduli space $\sM(E,h,\omega)$ of projective vortices should follow by methods similar to those used to prove the existence of such a compactification for the moduli space of non-Abelian monopoles in \cite[Thorem 1.1, p. 270]{FL1}. When $X$ has complex dimension three or more, Tian and Yang \cite{Tian_Yang_2002} have proved the existence of an Uhlenbeck-type compactification for the moduli space of vortices (see \cite[Equations (1.1) and (1.2)]{Tian_Yang_2002}), building on earlier results due to Tian \cite{TianGTCalGeom} for the moduli space of projectively Hermitian--Einstein connections.
\end{rmk}

\begin{rmk}[Topology of the moduli space of projective vortices over a complex K\"ahler manifold near a zero-section point]
\label{rmk:Kobayashi_7_3_17_at_irreducible_zero-section_points}
For the reasons noted in Remark \ref{rmk:Kobayashi_7_3_17_strongly_simple_pair} for the moduli space $\fM(E)$ of holomorphic pairs in \eqref{eq:Moduli_space_holomorphic_pairs}, we would not expect an analogue of Item \eqref{item:Kobayashi_7_4_19_projective_vortices_StabAvarphi_idE} in Theorem \ref{thm:Kobayashi_7_4_19_projective_vortices} to hold for the moduli space $\sM^*(E,h,\omega)$ of non-split projective vortices in \eqref{eq:Moduli_space_projective_vortices_non-split}. Indeed, a non-split zero-section pair $(A,0)$ has a stabilizer in $W^{2,p}(\SU(E))$ that is isomorphic, assuming $X$ is connected, to $C_r$ (the group of $r$-th roots of unity) and so $\sM^*(E,h,\omega)$ would have an orbifold singularity at the zero-section point $[A,0]$.
\end{rmk}

\section[Hitchin--Kobayashi correspondence between vortices and holomorphic pairs]{Hitchin--Kobayashi correspondence between projective vortices and semistable holomorphic pairs of vector bundles and sections}
\label{sec:HK_correspondence_between_SO3_monopoles_and_semistable_pairs}
We have the following analogue for $(0,1)$-pairs of Definition \ref{defn:Split_trivial_central-stabilizer_unitary_pair} \eqref{item:Split_unitary_pair} for unitary pairs and of Definition \ref{defn:Split_01-connection} for split $(0,1)$-connections.

\begin{defn}[Split $(0,1)$-pairs]
\label{defn:Split_(0,1)-pair}
Let $E$ be a $W^{2,p}$ complex vector bundle of rank $r\geq 2$ over a closed, almost complex manifold of real dimension $2n$ with $p\in(n,\infty)$. A $W^{1,p}$ $(0,1)$-connection $\bar\partial_A$ on $E$ is \emph{split} if $\bar\partial_A = \bar\partial_{A_1}\oplus\bar\partial_{A_2}$ with respect to a decomposition
\begin{equation}
  \label{eq:E_equals_L_1_oplus_L_2}
  E=E_1\oplus E_2
\end{equation}
as a direct sum of proper, $W^{2,p}$ complex subbundles of $E$, where $\bar\partial_{A_i}$ is a $W^{1,p}$ $(0,1)$-connection on $E_i$, for $i=1,2$, and $\varphi$ is a section of $E_1$. The pair is \emph{non-split} otherwise. If the pair $(E,\bar\partial_A)$ is smooth (respectively, analytic), then the pair $(E_i,\bar\partial_{A_i})$ is smooth (respectively, analytic) for $i=1,2$.
\end{defn}

We have the following analogue of Theorem \ref{thm:Hitchin-Kobayashi_correspondence_Hermitian-Einstein_connections_stable_bundles}. The result is essentially due to Bradlow \cite{Bradlow_1991}, but there are differences since Bradlow considers the vortex equations and we consider the projective vortex equations.

\begin{thm}[Hitchin--Kobayashi correspondence between projective vortices and polystable holomorphic pairs]
\label{thm:HitchinKobayashiCorrespondenceForPairs}
(See Bradlow and Garc{\'{\i}}a-Prada \cite[Theorem 6.4, p. 581 and Theorem 6.9, p. 583]{BradlowGP} and Okonek and Teleman in \cite[Theorem 6.2, p. 377]{OTQuaternion}, \cite[Theorem 3.3, p. 903]{OTVortex}.)
Let $(E,h)$ be a smooth Hermitian vector bundle over a closed, complex K\"ahler manifold $(X,\om)$ and $A_d$ be a fixed smooth, unitary connection on the Hermitian line bundle $\det E$ with $F_{A_d}^{0,2}=0$. Then the following hold:
\begin{enumerate}
\item \emph{(Projective vortex $\implies$ polystable holomorphic pair.)}
\label{item:Projective_vortex_implies_polystable_holomorphic pair}
If $(A,\varphi)$ is a solution to the projective vortex equations \eqref{eq:SO(3)_monopole_equations_almost_Hermitian_alpha} and $A$ obeys the determinant condition \eqref{eq:Unitary_connection_detE_fixed}, so $A$ induces $A_d$ on $\det E$, then $(\bar\partial_A,\varphi)$ is an $(0,1)$-pair such that $\bar\partial_A$ obeys the determinant condition \eqref{eq:Holomorphic_structure_fixed_determinant}, so $\bar\partial_A$ induces $\bar\partial_{A_d}$ on $\det E$, and $(\bar\partial_A,\varphi)$ is polystable in the sense of Definition \ref{defn:StablePair_arbitrary_rank}. Moreover, the following hold:
\begin{enumerate}
\item\label{item:Projective_vortex_implies_polystable_holomorphic_pair_non-split}
If $(A,\varphi)$ is non-split as in Definition \ref{defn:Split_trivial_central-stabilizer_unitary_pair} (\ref{item:Split_unitary_pair}), then $(\bar\partial_A,\varphi)$ is $\varphi$-stable.

\item\label{item:Projective_vortex_implies_polystable_holomorphic_pair_split}
If $(A,\varphi)$ is split as in Definition \ref{defn:Split_trivial_central-stabilizer_unitary_pair} (\ref{item:Split_unitary_pair}) and $(A_1,\varphi)$ is non-split,
then $(\bar\partial_A,\varphi)$ is split as in Definition \ref{defn:Split_(0,1)-pair} and $(\bar\partial_{A_1},\varphi)$ is $\varphi$-stable while $(E_2,\bar\partial_{A_2})$ is polystable as in Definition \ref{defn:Stable_holomorphic_structure}.
\end{enumerate}

\item \emph{(Polystable holomorphic pair $\implies$ projective vortex.)}
\label{item:HitchinKobayashiCorrespondenceForPairs_holomorphic_pair}
Let $(\bar\partial_E,\varphi)$ be a holomorphic pair obeying the determinant condition \eqref{eq:Holomorphic_structure_fixed_determinant}, so $\bar\partial_A$ induces $\bar\partial_{A_d}$ on $\det E$. If $(\bar\partial_E,\varphi)$ is polystable in the sense of Definition \ref{defn:StablePair_arbitrary_rank}, then there is a unique $\SU(E)$-orbit of a solution $(A',\varphi')$ to the projective vortex equations \eqref{eq:SO(3)_monopole_equations_almost_Hermitian_alpha}  such that $A$ induces $A_d$ on $\det E$ and $(\bar\partial_{A'},\varphi')$ is isomorphic to $(\bar\partial_E,\varphi)$ by an $\SL(E)$-gauge transformation.
Moreover, the following hold:
\begin{enumerate}
\item\label{item:HitchinKobayashiCorrespondenceForPairs_holomorphic_pair_non-split}
  If $(\bar\partial_E,\varphi)$ is non-split as in Definition \ref{defn:Split_(0,1)-pair}, then $(A',\varphi')$ is non-split as in Definition \ref{defn:Split_trivial_central-stabilizer_unitary_pair} (\ref{item:Split_unitary_pair}).

\item\label{item:HitchinKobayashiCorrespondenceForPairs_holomorphic_pair_split}
  If $(\bar\partial_E,\varphi)$ is split as in Definition \ref{defn:Split_(0,1)-pair} and $(\bar\partial_{E_1},\varphi)$ is non-split, then $(A_1',\varphi')$ is non-split as in Definition \ref{defn:Split_trivial_central-stabilizer_unitary_pair} (\ref{item:Split_unitary_pair}) while $A_2$ may be non-split or split as in Definition \ref{defn:Reducible_split_trivial-stabilizer_unitary_connection}.
\end{enumerate}
\end{enumerate}
\end{thm}

\begin{rmk}[Regularity of projective vortices and integrable $(0,1)$-pairs in the Hitchin--Kobayashi correspondence]
\label{rmk:HitchinKobayashiCorrespondenceForPairs_regularity}
Our comments in Remark \ref{rmk:Hitchin-Kobayashi_correspondence_Hermitian-Einstein_connections_regularity} apply almost verbatim to the conclusions of Theorem \ref{thm:HitchinKobayashiCorrespondenceForPairs} as well. In Item \eqref{item:Projective_vortex_implies_polystable_holomorphic pair}  of Theorem \ref{thm:HitchinKobayashiCorrespondenceForPairs}, if $(A,\varphi)$ is a $W^{1,p}$ projective vortex, then it is immediate that $(\bar\partial_A,\varphi)$ is a $W^{1,p}$ integrable $(0,1)$-pair. Conversely, although this is less obvious, if $(\bar\partial_E,\varphi)$ is a $W^{1,p}$ integrable $(0,1)$-pair, then $(A,\varphi')$ is a $W^{1,p}$ projective vortex. One approach to proving this conclusion that regularity is preserved by the Hitchin--Kobayashi correspondence is to examine its proof via gradient flow: for example, see Bradlow \cite{Bradlow_1991} and Simpson \cite{Simpson_1988}. Similar remarks apply if the pairs are $W^{k,p}$ for $k\geq 2$ or $C^l$ for $l\geq 1$ or real analytic. Moreover, although the projective vortex $(A,\varphi')$ produced by Item \eqref{item:HitchinKobayashiCorrespondenceForPairs_holomorphic_pair} of Theorem \ref{thm:HitchinKobayashiCorrespondenceForPairs} may only be of class $W^{1,p}$, Lemma \ref{lem:Regularity_For_ProjVortices} yields a $W^{2,p}$ determinant one, unitary gauge transformation $u$ such that $u^*(A,\varphi')$ is smooth.
\end{rmk}

\begin{rmk}[On the hypothesis in Theorem \ref{thm:HitchinKobayashiCorrespondenceForPairs} that $(X,\omega)$ is a complex K\"ahler manifold]
\label{rmk:HitchinKobayashiCorrespondenceForPairs_non-Kaehler}
It is likely that the hypotheses in Theorem \ref{thm:HitchinKobayashiCorrespondenceForPairs} could be relaxed, just as we noted in Remark \ref{rmk:Hitchin-Kobayashi_correspondence_Hermitian-Einstein_connections_stable_bundles_non-Kaehler} for the Hitchin--Kobayashi correspondence between Hermitian--Einstein connections and polystable holomorphic structures. For example, L\"ubke and Teleman \cite[Theorem 6.2.1, p. 65]{Lubke_Teleman_2006} obtain a version of this result when $(X,g,J)$ is only assumed to be complex Hermitian, while Jacob \cite[Theorem 1, p. 119]{Jacob_2015conm} proves a version for Higgs pairs when $(X,g,J)$ is only assumed to be complex Hermitian equipped with a Gauduchon metric.
\end{rmk}

Dowker \cite{DowkerThesis} provides a proof of much of Theorem \ref{thm:HitchinKobayashiCorrespondenceForPairs} when $X$ has dimension two and $E$ has complex rank two, but his definitions of stability and semistability (see \cite[Definitions 1.3.2 and 1.3.5]{DowkerThesis}) appear to be weaker than those of Bradlow and so we avoid relying on his work. We shall instead apply results due to Okonek and Teleman \cite{OTQuaternion, OTVortex} and ideas explained to us by Richard Wentworth \cite{Wentworth_2022private}. We begin by recalling the following important results due to Bradlow, though we translate them to the setting of pairs $(A,\varphi)$ of unitary connections and sections from the setting of pairs $(H,\varphi)$ of Hermitian metrics and sections. The concepts of stability in Definition \ref{defn:StablePair_arbitrary_rank} reflect the conclusions of the forthcoming

\begin{thm}
\label{thm:Bradlow_2-1-6}  
(See Bradlow \cite[Theorem 2.1.6, p. 179]{Bradlow_1991}.)  
Let $(E,h)$ be a smooth Hermitian vector bundle over a closed, complex K\"ahler manifold $(X,\omega)$. If $(A,\varphi)$ is a smooth solution to the \emph{vortex equations} for some positive parameter $\tau$, namely
\begin{subequations}
  \label{eq:Vortex_equations}
  \begin{align}
    \label{eq:Vortex_equation_LambdaFA}
    \Lambda F_A &= \frac{i}{2}\varphi\otimes\varphi^* - \frac{i}{2}\tau\,\id_E,
    \\
    \label{eq:Vortex_equation_FA02}
    F_A^{0,2} &= 0,
    \\
    \label{eq:Vortex_equation_dbarA_varphi}
    \bar\partial_A\varphi &= 0,
  \end{align}
\end{subequations}
then exactly one of the following two cases hold:
\begin{enumerate}
\item\label{item:Bradlow_2-1-6_a} $(\bar\partial_A,\varphi)$ is $\varphi$-stable as in Definition \ref{defn:StablePair_arbitrary_rank} and
  \begin{equation}
    \label{eq:Bradlow_tau_stable}
    \mu_M(E) < \frac{\tau\vol(X)}{4\pi} < \mu_m(\varphi).
  \end{equation}
    
\item\label{item:Bradlow_2-1-6_b} $(\bar\partial_A,\varphi)$ is split as in Definition \ref{defn:Split_(0,1)-pair}, so there is a splitting $(E,\bar\partial_A) = (E_1,\bar\partial_{A_1}) \oplus (E_2,\bar\partial_{A_2})$ as a direct sum of holomorphic vector bundles and $\varphi$ is a section of $E_1$.  
\end{enumerate}
In Case \eqref{item:Bradlow_2-1-6_b}, the pair $(\bar\partial_{A_1},\varphi)$ is $\varphi$-stable as in Definition \ref{defn:StablePair_arbitrary_rank} and satisfies \eqref{eq:Bradlow_tau_stable}, while the holomorphic vector bundle $(E_2,\bar\partial_{A_2})$ is polystable as in Definition \ref{defn:Stable_holomorphic_structure} with slope $\tau\vol(X)/4\pi$.
\end{thm}

\begin{thm}
\label{thm:Bradlow_3-1-1}
(See Bradlow \cite[Theorem 3.1.1, p. 184]{Bradlow_1991}.)
Let $(E,h)$ be a smooth Hermitian vector bundle over a closed, complex K\"ahler manifold $(X,\omega)$. If $(\bar\partial_E,\varphi)$ is a smooth, holomorphic pair that is $\varphi$-stable in the sense of Definition \ref{defn:StablePair_arbitrary_rank} and $\tau$ is a real constant obeying $\mu_M < \tau\vol(X)/4\pi < \mu_m(\varphi)$, then there exists a smooth Hermitian metric $H$ on $E$ and smooth $H$-unitary connection $A$ on $E$ such that $\bar\partial_A = \bar\partial_E$ and a $\U(E)$-gauge transformation $u$ such that $A'$ is $h$-unitary and the pair $(A',\varphi') = u^*(A,\varphi)$ solves the vortex equations \eqref{eq:Vortex_equations}.
\end{thm}

A solution to the projective vortex equations \eqref{eq:SO(3)_monopole_equations_almost_Hermitian_alpha} does not yield a solution to the vortex equations \eqref{eq:Vortex_equations} for some positive parameter $\tau$, but rather a weaker version of those equations. Let $(E,h)$ be a smooth Hermitian vector bundle over a closed, complex K\"ahler manifold $(X,\omega)$ and $A$ be a smooth unitary connection on $E$ and $\varphi$ be a smooth section of $E$. According to Okonek and Teleman \cite[Proof of Theorem 6.2, p. 179]{OTQuaternion}, \cite[Section 3 p. 902]{OTVortex}, a pair $(A,\varphi)$ is a smooth solution to the \emph{weak vortex equations} for some smooth function $t:X\to\RR$ if
\begin{subequations}
  \label{eq:Weak_vortex_equations}
  \begin{align}
    \label{eq:Weak_vortex_equation_LambdaFA}
    \Lambda F_A &= \frac{i}{2}\varphi\otimes\varphi^* - \frac{i}{2}t\,\id_E,
    \\
    \label{eq:Weak_vortex_equation_FA02}
    F_A^{0,2} &= 0,
    \\
    \label{eq:Weak_vortex_equation_dbarA_varphi}
    \bar\partial_A\varphi &= 0.
  \end{align}
\end{subequations}
Clearly, the system \eqref{eq:Weak_vortex_equations} is identical to \eqref{eq:Vortex_equations} except that in \eqref{eq:Weak_vortex_equation_LambdaFA} the function $f$ need not be a positive constant as assumed in \eqref{eq:Vortex_equation_LambdaFA}.

\begin{proof}[Proof of Theorem \ref{thm:HitchinKobayashiCorrespondenceForPairs}]
Consider Item \eqref{item:Projective_vortex_implies_polystable_holomorphic pair}, so $(A,\varphi)$ is a solution to the projective vortex equations \eqref{eq:SO(3)_monopole_equations_almost_Hermitian_alpha}. Following Okonek and Teleman \cite[Proof of Theorem 6.2, p. 377]{OTQuaternion} and Wentworth \cite{Wentworth_2022private}, we claim that $(A,\varphi)$ is therefore a solution to the weak vortex equations \eqref{eq:Weak_vortex_equations}. Indeed, by \eqref{eq:SO(3)_monopole_equations_(1,1)_curvature_alpha} there is a smooth function $t:X\to\RR$ such that
\[
  \Lambda F_A = \frac{i}{2}\varphi\otimes\varphi^* - \frac{i}{2}t\,\id_E,
\]
since we may write
\[
  (F_A)_0 = F_A - \frac{1}{r}(\tr_EF_A)\,\id_E
  \quad\text{and}\quad
  (\varphi\otimes\varphi^*)_0 = \varphi\otimes\varphi^* - \frac{1}{r}|\varphi|_E^2\,\id_E,
\]
noting that $\tr_EF_A = F_{A_d}$ and, with respect to a local orthonormal frame $\{e_1,\ldots,e_r\}$ for $E$,
\[
  \tr_E(\varphi\otimes\varphi^*) = \sum_{j=1}^r \langle \varphi\otimes\varphi^*(e_j),e_j\rangle_E = \sum_{j=1}^r\langle \varphi,e_j\rangle \langle e_j,\varphi\rangle = |\varphi|_E^2.
\]
Substituting these identities into \eqref{eq:SO(3)_monopole_equations_(1,1)_curvature_alpha} gives
\[
  \Lambda F_A - \frac{1}{r}\Lambda(\tr_EF_A)\,\id_E
  =
  \frac{i}{2}\varphi\otimes\varphi^* - \frac{i}{2r}|\varphi|_E^2\,\id_E,
\]
and thus
\[
  \Lambda F_A = \frac{i}{2}\varphi\otimes\varphi^*
  - \frac{i}{2r}\left(2i\Lambda(\tr_EF_A)\,\id_E + |\varphi|_E^2\right)\,\id_E.
\]
Therefore, we choose
\[
  t := \frac{1}{r}\left(2i\Lambda\tr_EF_A + |\varphi|_E^2\right).
\]
By \eqref{eq:SO(3)_monopole_equations_(0,2)_curvature_zero_alpha}, we have $(F_A^{0,2})_0 = 0$ and because $A$ induces $A_d$ on $\det E$ with $F_{A_d}^{0,2} = 0$ by assumption, we have $F_A^{0,2} = 0$ and so $A$ solves \eqref{eq:Weak_vortex_equation_FA02}. Equations \eqref{eq:SO(3)_monopole_equations_Dirac_almost_Hermitian_alpha} and \eqref{eq:Weak_vortex_equation_dbarA_varphi} are identical and so $(A,\varphi)$ solves \eqref{eq:Weak_vortex_equation_dbarA_varphi}. Thus, $(A,\varphi)$ solves the weak vortex equations  \eqref{eq:Weak_vortex_equations}, as claimed. 

According to Okonek and Teleman \cite[Proof of Theorem 6.2, p. 377]{OTQuaternion} and Bradlow and Garc{\'{\i}}a-Prada \cite[Theorem 6.4, p. 581]{BradlowGP}, the $(0,1)$-pair $(\bar\partial_A,\varphi)$ is $\lambda$-polystable for the parameter given by the average value of the function $t$, namely
\begin{equation}
  \label{eq:Okonek-Teleman_Bradlow-Garcia-Prada_lambda_mu(sE)_L2norm_varphi_square_identity}
  \lambda = \frac{(n-1)!}{4\pi}\int_Xt\,d\vol_\omega = \mu(\sE) + \frac{(n-1)!}{4\pi r}\|\varphi\|_{L^2(X)}^2,
\end{equation}
noting the definition \eqref{eq:Slope} of $\mu(\sE)$ and expression
\eqref{eq:DegreeIntegral} for $\deg\sE$:
\[
  \mu(\sE) = \frac{1}{r}\int_X\frac{i}{2\pi}\frac{1}{n}(\Lambda\tr_EF_A)\omega^n
  = \frac{(n-1)!}{4\pi}\int_X \frac{1}{r}(2i\Lambda\tr_EF_A)\,d\vol_\omega.
\]  
This assertion follows because, as the preceding authors note, replacing the constant $\tau$ by the function $t$ does not impact the proof of Theorem \ref{thm:Bradlow_2-1-6}. The main conclusion of Item \eqref{item:Projective_vortex_implies_polystable_holomorphic pair} thus follows from Theorem \ref{thm:Bradlow_2-1-6}.

Consider Item \eqref{item:Projective_vortex_implies_polystable_holomorphic_pair_non-split}. Because the pair $(A,\varphi)$ is not split, Case \eqref{item:Bradlow_2-1-6_b} of Theorem \ref{thm:Bradlow_2-1-6} cannot hold. Therefore, Case \eqref{item:Bradlow_2-1-6_a} of Theorem \ref{thm:Bradlow_2-1-6} applies and the pair $(\bar\partial_A,\varphi)$ must be $\varphi$-stable and this verifies Item \eqref{item:Projective_vortex_implies_polystable_holomorphic_pair_non-split}. Consider Item \eqref{item:Projective_vortex_implies_polystable_holomorphic_pair_split}. Because the pair $(A,\varphi)$ is split, we must have by Definition \ref{defn:Split_trivial_central-stabilizer_unitary_pair} \eqref{item:Split_unitary_pair} that $A=A_1\oplus A_2$ with respect to a decomposition $E=E_1\oplus E_2$ as an orthogonal direct sum of Hermitian vector bundles and $\varphi$ is a section of $E_1$. By assumption, the pair $(A_1,\varphi)$ is non-split and because $(A_1,\varphi)$ is also a solution to the projective vortex equations \eqref{eq:SO(3)_monopole_equations_almost_Hermitian_alpha}, Item \eqref{item:Projective_vortex_implies_polystable_holomorphic_pair_non-split} implies that the pair $(\bar\partial_{A_1},\varphi)$ must be $\varphi$-stable. The unitary connection $A_2$ is projectively Hermitian--Einstein and so the pair $(E_2,\bar\partial_{A_2})$ must be polystable by Item \eqref{item:HEconnection_implies_polystable_bundle} of Theorem \ref{thm:Hitchin-Kobayashi_correspondence_Hermitian-Einstein_connections_stable_bundles}. This verifies Item \eqref{item:Projective_vortex_implies_polystable_holomorphic_pair_split}.

Consider Item \eqref{item:HitchinKobayashiCorrespondenceForPairs_holomorphic_pair}, so $(\bar\partial_E,\varphi)$ is a polystable, holomorphic pair. According to Okonek and Teleman \cite[Theorem 6.2, p. 377]{OTQuaternion} (and which they verify by modifying Bradlow's proof of Theorem \ref{thm:Bradlow_3-1-1}), there is a Hermitian metric $H$ on $E$ that induces the same Hermitian metric $h_{\det E}$ on the complex line bundle $\det E$ as that of the given Hermitian metric $h$ and such that the $H$-unitary Chern connection $A$ defined by $(\bar\partial_E,H)$ and section $\varphi$ solve the projective vortex equations \eqref{eq:SO(3)_monopole_equations_almost_Hermitian_alpha}. By assumption, $A_d$ is the $h_{\det E}$-unitary Chern connection on $\det E$ defined by $\bar\partial_{\det E}$ and so $A$ induces $A_d$ on $\det E$.
If $(\bar\partial_E,\varphi)$ is a stable pair, then $H$ and thus $A$ are unique by Okonek and Teleman \cite[Theorem 6.2, p. 377]{OTQuaternion}. If $(\bar\partial_E,\varphi)$ is strictly polystable, then $\varphi \equiv 0$ and $(E,\bar\partial_E)$ is strictly polystable and Theorem \ref{thm:Hitchin-Kobayashi_correspondence_Hermitian-Einstein_connections_stable_bundles} 
\eqref{item:Polystable_bundle_implies_HEconnection} implies that the $\SU(E)$-orbit of $(A',0)$ is uniquely determined by $(E,\bar\partial_E)$. There is an $\SU(E)$-gauge transformation $u$ such that $A' = u^*A$ is an $h$-unitary connection on $E$ inducing $A_d$ on $\det E$ and the pair $(A',\varphi') = u^*(A,\varphi)$ solves the projective vortex equations \eqref{eq:SO(3)_monopole_equations_almost_Hermitian_alpha}. This proves the main conclusion of Item \eqref{item:HitchinKobayashiCorrespondenceForPairs_holomorphic_pair}. The dichotomy represented by Items \eqref{item:HitchinKobayashiCorrespondenceForPairs_holomorphic_pair_non-split} and \eqref{item:HitchinKobayashiCorrespondenceForPairs_holomorphic_pair_split} follows from the dichotomy provided by Items \eqref{item:Projective_vortex_implies_polystable_holomorphic_pair_non-split} and \eqref{item:Projective_vortex_implies_polystable_holomorphic_pair_split}.
\end{proof}

We have the following refinement of Theorem \ref{thm:HitchinKobayashiCorrespondenceForPairs} when $E$ has rank two and the pairs are split.

\begin{lem}[Hitchin--Kobayashi correspondence for split projective vortices and polystable holomorphic pairs]
\label{lem:HKIdentificationOfSplitPairs}
Let $(E,h)$ be a rank-two, smooth Hermitian vector bundle over a closed, connected complex K\"ahler manifold $(X,\om)$. Let $A_d$ be a fixed smooth, unitary connection on the Hermitian line bundle $\det E$ with $F_{A_d}^{0,2}=0$. Then the following hold:
\begin{enumerate}
\item
\label{item:HKIdentificationOfSplitPairs1}
If $(A,\varphi)$ is a solution to the projective vortex equations
\eqref{eq:SO(3)_monopole_equations_almost_Hermitian_alpha} which is split as in Definition \ref{defn:Split_trivial_central-stabilizer_unitary_pair} (\ref{item:Split_unitary_pair}) with respect to the decomposition $E=L_1\oplus L_2$ as an orthogonal direct sum of Hermitian line bundles, then the pair $(\bar\partial_A,\varphi)$ is split as in Definition \ref{defn:Split_(0,1)-pair} with respect to the decomposition $E = L_1\oplus L_2$ as a direct sum of complex line bundles.

\item
\label{item:HKIdentificationOfSplitPairs2}
If $(\bar\partial_E,\varphi)$ is a holomorphic pair that is polystable as in Definition \ref{defn:StablePair_arbitrary_rank} and split with respect to a decomposition
$E=L_1\oplus L_2$ as a direct sum of complex vector bundles, then the $\SL(E)$-orbit of $(\bar\partial_E,\varphi)$ contains the $\SU(E)$-orbit of a holomorphic pair $(\bar\partial_A,\varphi')$, where $(A,\varphi')$ is a solution to the projective vortex equations that is split with respect to the decomposition $E=L_1\oplus L_2$ as an orthogonal direct sum of Hermitian line bundles.
\end{enumerate}
\end{lem}

\begin{proof}
Consider Item \eqref{item:HKIdentificationOfSplitPairs1}. By Definition \ref{defn:Split_trivial_central-stabilizer_unitary_pair} \eqref{item:Split_unitary_pair}, we
can write $A=A_1\oplus A_2$, where $A_i$ is a unitary connection on $L_i$ for $i=1,2$ and $\varphi$ is a section of $L_1$. Thus, $\bar\rd_A=\bar\rd_{A_1}\oplus \bar\rd_{A_2}$ so $(\bar\rd_A,\varphi)$ is split in the sense of Definition \ref{defn:Split_(0,1)-pair}.  This verifies Item \eqref{item:HKIdentificationOfSplitPairs1}.

Consider Item \eqref{item:HKIdentificationOfSplitPairs2}. Let $A''$ be the Chern connection defined by
$(\bar\partial_E,h)$, so $\bar\partial_{A''}=\bar\partial_E$. By Item \eqref{item:HitchinKobayashiCorrespondenceForPairs_holomorphic_pair} of Theorem \ref{thm:HitchinKobayashiCorrespondenceForPairs}, there is an $\SL(E)$-gauge transformation $v$
such that $v(A'',\varphi)=(A,\varphi')$ is a solution to the projective vortex equations. By the proof of Dowker \cite[Proposition 2.2.2, p. 17]{DowkerThesis}, the gauge transformation $v$ has the form
\[
v= e^{\la/2}\, \id_{L_1}\oplus e^{-\la/2}\, \id_{L_2},
\]
where $\la$ is a real-valued function on $X$. Because $\bar\partial_E = \bar\partial_{L_1}\oplus \bar\partial_{L_2}$
is split with respect to the decomposition $E=L_1\oplus L_2$, the uniqueness of the Chern connection (see Kobayashi \cite[Proposition 1.4.9, p. 11]{Kobayashi_differential_geometry_complex_vector_bundles}) implies that $A''=A_1\oplus A_2$, where $A_i$ is the Chern connection defined by $(\bar\partial_{L_i}, h_{L_i})$, for $i=1,2$. Because $A''$ is split with respect to the decomposition $E=L_1\oplus L_2$ as an orthogonal direct sum of Hermitian line bundles, the preceding expression for $v$ implies that $v^*A''$ is also split with respect to the decomposition $E=L_1\oplus L_2$ as an orthogonal direct sum of Hermitian line bundles and $u^{-1}(\varphi)$ is a section of $L_1$. Thus, $(A,\varphi')=v(A'',\varphi)$ is the desired split solution to the projective vortex equations. This verifies Item \eqref{item:HKIdentificationOfSplitPairs2} and completes the proof of the lemma.
\end{proof}

\section[Real analytic embeddings of moduli spaces of projective vortices]{Real analytic embeddings of moduli spaces of projective vortices and non-Abelian monopoles into moduli spaces of pairs of  holomorphic pairs}
\label{sec:Moduli_space_SO3_monopoles_open_subspace_moduli_space_simple_holomorphic_pairs}
We have the following analogue of Theorem \ref{thm:Kobayashi_7_4_20} and generalization of L\"ubke and Teleman \cite[Theorem 6.3.7, p. 72]{Lubke_Teleman_2006} to allow $E$ to have arbitrary complex rank.

\begin{thm}[Real analytic embedding of the moduli space of projective vortices into the moduli space of holomorphic pairs]
\label{thm:Lubke_Teleman_6-3-7}
(See L\"ubke and Teleman \cite[Theorem 6.3.7, p. 72]{Lubke_Teleman_2006}.)  
Let $(E,h)$ be a smooth Hermitian vector bundle over a closed complex K\"ahler manifold $(X,\omega)$ and $A_d$ be a smooth unitary connection on the Hermitian line bundle $\det E$ such that $F_{A_d}^{0,2} = 0$. If $X$ has complex dimension $n$ and $p\in(n,\infty)$ is a constant, then the bijection
\begin{equation}
\label{eq:Map_W1p_unitary_pairs_to_01-pairs}
  \sA(E,h)\times W^{1,p}(E) \ni (A,\varphi) \mapsto (\bar\partial_A,\varphi) \in \sA^{0,1}(E)\times W^{1,p}(E)
\end{equation}
between $W^{1,p}$ pairs induces the following properties:
\begin{enumerate}
\item\label{item:Lubke_Teleman_6-3-7_analytic_embedding_rank-2_non-zero_section}
If $E$ has rank two, then there is an embedding in the sense of real analytic spaces
from the moduli subspace $\sM^0(E,h,\omega)$ of non-zero-section projective vortices in \eqref{eq:Moduli_space_projective_vortices_non-zero_section} onto an open subspace of the moduli space $\fM(E)$ of holomorphic pairs in \eqref{eq:Moduli_space_holomorphic_pairs},
\begin{equation}
  \label{eq:Lubke_Teleman_6-3-7_analytic_embedding_rank-2_non-zero_section}
  \sM^0(E,h,\omega) \hookrightarrow \fM(E),
\end{equation}
given by the moduli subspace $\fM_\ps^0(E,\omega)$ of non-zero-section, polystable holomorphic pairs in \eqref{eq:Moduli_space_non-zero-section_polystable_holomorphic_pairs} and which is contained in $\fM^{**}(E)$ in \eqref{eq:Moduli_space_strongly_simple_holomorphic_pairs}.
  
\item\label{item:Lubke_Teleman_6-3-7_analytic_embedding_non-split-nonzero-section}
There is an embedding in the sense of real analytic spaces from the moduli subspace $\sM^{*,0}(E,h,\omega)$ of non-split, non-zero-section projective vortices in
\eqref{eq:Moduli_space_projective_vortices_non-split-non-zero-section} onto an open subspace of the moduli space
$\fM(E)$ of holomorphic pairs in \eqref{eq:Moduli_space_holomorphic_pairs},
\begin{equation}
  \label{eq:Lubke_Teleman_6-3-7_analytic_embedding_non-split-nonzero-section}
  \sM^{*,0}(E,h,\omega) \hookrightarrow \fM(E),
\end{equation}
given by the moduli subspace $\fM^0(E,\omega)$ of non-zero-section, stable holomorphic pairs in \eqref{eq:Moduli_space_non-zero-section_stable_holomorphic_pairs} and which is contained in $\fM^{**}(E)$ in \eqref{eq:Moduli_space_strongly_simple_holomorphic_pairs}.

\item\label{item:Lubke_Teleman_6-3-7_set-theoretic_embedding}
There are set-theoretic bijections from
\begin{enumerate}
\item The moduli space $\sM^*(E,h,\omega)$ of non-split projective vortices in \eqref{eq:Moduli_space_projective_vortices_non-split} onto a subset of the moduli space \eqref{eq:Moduli_space_holomorphic_pairs} of holomorphic pairs,
\begin{equation}
  \label{eq:Lubke_Teleman_6-3-7_set-theoretic_embedding_non-split_stable}
  \sM^*(E,h,\omega) \hookrightarrow \fM(E),
\end{equation}
given by the moduli subset $\fM(E,\omega)$ of stable holomorphic pairs in \eqref{eq:Moduli_space_stable_holomorphic_pairs}.

\item The moduli space $\sM(E,h,\omega)$ of projective vortices in \eqref{eq:Moduli_space_projective_vortices} onto a subset of the moduli space \eqref{eq:Moduli_space_holomorphic_pairs} of holomorphic pairs,
\begin{equation}
  \label{eq:Lubke_Teleman_6-3-7_set-theoretic_embedding}
  \sM(E,h,\omega) \hookrightarrow \fM(E),
\end{equation}
given by the moduli subset $\fM_\ps(E,\omega)$ of polystable holomorphic pairs in \eqref{eq:Moduli_space_polystable_holomorphic_pairs}.
\end{enumerate}
\end{enumerate}
\end{thm}

\begin{rmk}
\label{rmk:Lubke_Teleman_6-3-7_embedding_topological_spaces}
L\"ubke and Teleman \cite[Theorem 6.3.7, p. 72]{Lubke_Teleman_2006} further assert that the set-theoretic injection \eqref{eq:Lubke_Teleman_6-3-7_set-theoretic_embedding} actually yields a homeomorphism from $\sM(E,h,\omega)$ onto $\fM_\ps(E,\omega)$.
\end{rmk}

Theorem \ref{thm:Lubke_Teleman_6-3-7} is proved in Section \ref{sec:Proofs_theorem_Lubke_Teleman_6-3-7_and_corollary}. Teleman states a version of Theorem \ref{thm:Lubke_Teleman_6-3-7} in \cite[Theorem 2.2.6]{TelemanNonabelian} and remarks that it can be proved by the same methods as in Donaldson and Kronheimer \cite{DK}, L\"ubke and Teleman \cite{Lubke_Teleman_1995}, and Okonek and Teleman \cite{OTVortex}. A result similar to Theorem \ref{thm:Lubke_Teleman_6-3-7} is stated by Okonek and Teleman in \cite[Theorem 4.1]{OTVortex} and its proof is outlined by them in \cite[pp. 905--906]{OTVortex}. The references \cite{DK} and \cite{Lubke_Teleman_1995} are concerned with the relationship between moduli spaces of anti-self-dual connections over complex K\"ahler surfaces \cite{DK} and, more generally, projectively Hermitian--Einstein connections over complex Hermitian manifolds and moduli spaces of stable holomorphic vector bundles.

Theorem \ref{thm:Lubke_Teleman_6-3-7} yields the following partial analogue for moduli spaces of holomorphic pairs of Theorem \ref{thm:Kobayashi_7_3_17} for holomorphic vector bundles and which we had previewed in Remark \ref{rmk:Kobayashi_7_3_17_strongly_simple_pair}.

\begin{cor}[Moduli spaces of stable or polystable holomorphic pairs over complex manifolds]
\label{cor:Lubke_Teleman_6-3-7_and_Kobayashi_7_3_17_pair}
Continue the hypotheses of Theorem \ref{thm:Lubke_Teleman_6-3-7}. Then the following hold:
\begin{enumerate}
\item\label{item:fM_ps^0(E,omega)_rank-2_is_complex_analytic_space}
If $E$ has complex rank two, then the moduli subspace $\fM_\ps^0(E,\omega)$ in \eqref{eq:Moduli_space_regular_non-zero-section_polystable_holomorphic_pairs} of non-zero section, polystable holomorphic vector pairs is a (Hausdorff) complex analytic space with Zariski tangent spaces $\bH_{\bar\partial_E,\varphi}^1$ at points $[\bar\partial_E,\varphi] \in \fM_\ps^0(E,\omega)$ and smooth at points $[\bar\partial_E,\varphi]$ with $\bH_{\bar\partial_E,\varphi}^2 = 0$, where the harmonic spaces $\bH_{\bar\partial_E,\varphi}^\bullet$ are defined by \eqref{eq:H_dbar_Avarphi^0bullet}. In particular, the moduli subspace $\fM_{\ps,\reg}^0(E,\omega)$ of non-zero section, regular, polystable holomorphic pairs is a complex manifold with tangent spaces $\bH_{\bar\partial_E,\varphi}^1$ at points $[\bar\partial_E,\varphi]$ and the embedding \eqref{eq:Lubke_Teleman_6-3-7_analytic_embedding_rank-2_non-zero_section} restricts to an isomorphism of real analytic manifolds,
\begin{equation}
  \label{eq:Lubke_Teleman_6-3-7_analytic_isomorphism_rank-2_non-zero_section_regular}
  \sM_\reg^0(E,h,\omega) \cong \fM_{\ps,\reg}^0(E,\omega) \subset \fM^{**}(E). 
\end{equation}

\item\label{item:fM^0(E,omega)_is_complex_analytic_space}
  The moduli subspace $\fM^0(E,\omega)$ in \eqref{eq:Moduli_space_non-zero-section_stable_holomorphic_pairs} of non-zero section, stable holomorphic vector pairs is a (Hausdorff) complex analytic space with Zariski tangent spaces $\bH_{\bar\partial_E,\varphi}^1$ at points $[\bar\partial_E,\varphi] \in \fM^0(E,\omega)$ and smooth at points $[\bar\partial_E,\varphi]$ with $\bH_{\bar\partial_E,\varphi}^2 = 0$, where the harmonic spaces $\bH_{\bar\partial_E,\varphi}^\bullet$ are defined by \eqref{eq:H_dbar_Avarphi^0bullet}. In particular, the moduli subspace $\fM_\reg^0(E,\omega)$ in \eqref{eq:Moduli_space_regular_non-zero-section_stable_holomorphic_pairs} of non-zero section, regular, stable holomorphic pairs is a complex manifold with tangent spaces $\bH_{\bar\partial_E,\varphi}^1$ at points $[\bar\partial_E,\varphi]$ and the embedding \eqref{eq:Lubke_Teleman_6-3-7_analytic_embedding_non-split-nonzero-section} restricts to an isomorphism of real analytic manifolds,
\begin{equation}
  \label{eq:Lubke_Teleman_6-3-7_analytic_isomorphism_non-split-nonzero-section_regular}
   \sM_\reg^{*,0}(E,h,\omega) \cong \fM_\reg^0(E,\omega) \subset \fM^{**}(E).
\end{equation}
\end{enumerate}
\end{cor}

Corollary \ref{cor:Lubke_Teleman_6-3-7_and_Kobayashi_7_3_17_pair} is proved in Section \ref{sec:Proofs_theorem_Lubke_Teleman_6-3-7_and_corollary}. We now specialize Theorem \ref{thm:Lubke_Teleman_6-3-7} to the case of complex surfaces.

\begin{thm}[Real analytic embedding of the moduli space of non-Abelian monopoles into the moduli space of holomorphic pairs]
\label{thm:Lubke_Teleman_6-3-10}
(See L\"ubke and Teleman \cite[Theorem 6.3.10, p. 74]{Lubke_Teleman_2006}.)
Let $(E,h)$ be a rank-two, smooth Hermitian vector bundle over a closed complex K\"ahler surface $(X,g,J)$ with K\"ahler form $\omega = g(\cdot,J\cdot)$ and $A_d$ be a smooth unitary connection on the Hermitian line bundle $\det E$ such that $F_{A_d}^{0,2} = 0$. Define the rank-two, smooth Hermitian vector bundle $F$ over $X$ by $F=K_X\otimes E^*$ where $K_X$ is the canonical line bundle over $X$ defined in \eqref{eq:DefineCanonicalLineBundle}. For $\ft = (\rho_\can, W_\can, E)$, the moduli space $\sM_\ft$ of non-Abelian monopoles is the union of two closed subspaces,
\[
  \sM_\ft = \sM_{\ft,1}\cup\sM_{\ft,2},
\]
comprising moduli spaces of non-Abelian monopoles of type $1$ and $2$, respectively.
\begin{enumerate}
\item
\label{item:Lubke_Teleman_6-3-10_RealAnalyticEmbedding}
The map \eqref{eq:Map_W1p_unitary_pairs_to_01-pairs} of $W^{1,p}$ pairs induces embeddings in the sense of real analytic spaces from the moduli subspaces $\sM_{\ft,1}^0$ and $\sM_{\ft,2}^0$ of non-zero-section non-Abelian monopoles of type $1$ and type $2$, respectively, onto open subspaces of the moduli spaces $\fM(E)$ and $\fM(F)$, respectively, of holomorphic pairs in \eqref{eq:Moduli_space_holomorphic_pairs},
\begin{equation}
  \label{eq:Lubke_Teleman_6-3-10_analytic_embedding}
  \sM_{\ft,1}^0 \hookrightarrow \fM(E),
  \quad
  \sM_{\ft,2}^0 \hookrightarrow \fM(F),
\end{equation}
given by the moduli subspaces $\fM_\ps^0(E,\omega)$ and $\fM_\ps^0(F,\omega)$, respectively, of non-zero-section, polystable holomorphic pairs, and which are contained by \eqref{eq:Moduli_space_non-zero-section_polystable_rank-2_holomorphic_pairs_subset_strongly_simple_pairs} in the moduli subspaces $\fM^{**}(E)$ and $\fM^{**}(F)$, respectively, of strongly simple holomorphic pairs.

\item
\label{item:Lubke_Teleman_6-3-10_SetTheoreticEmbedding}
The map \eqref{eq:Map_W1p_unitary_pairs_to_01-pairs} of $W^{1,p}$ pairs induces set-theoretic injections from $\sM_{\ft,1}$ and $\sM_{\ft,2}$ onto subsets of the moduli spaces \eqref{eq:Moduli_space_holomorphic_pairs} of holomorphic pairs,
\begin{equation}
  \label{eq:Lubke_Teleman_6-3-10_set-theoretic_embedding_}
  \sM_{\ft,1} \hookrightarrow \fM(E),
  \quad
  \sM_{\ft,2} \hookrightarrow \fM(F),
\end{equation}
given by the moduli subsets $\fM_\ps(E,\omega)$ and $\fM_\ps(F,\om)$, respectively, of polystable holomorphic pairs as in \eqref{eq:Moduli_space_polystable_holomorphic_pairs}.
\end{enumerate}
\end{thm}

\begin{proof}
The identification $W_{\can} = \Lambda^{0,0}(X)\oplus \Lambda^{0,2}(X)$ given by \eqref{eq:Canonical_spinc_bundles} and continuous embedding of Banach spaces,
\[
  W^{1,p}(E) \ni \varphi \mapsto (\varphi,0) \in W^{1,p}(W_{\can}\otimes E) = W^{1,p}\left(E \oplus \Lambda^{0,2}(E)\right)
\]
defines a real analytic embedding of Banach manifolds,
\[
  j_{\ft,1}:\sC^0(E,h)=\left(\sA(E,h)\times W^{1,p}(E)\less\{0\}\right)/W^{p,2}(\SU(E)) \to \sC_\ft^0,
\]
where $\sC^0(E,h)$ is the open subspace \eqref{eq:ConfigurationSpaceForProjectiveVorticesnon-zero-section} of gauge-equivalence classes of non-zero section pairs in $\sC(E,h)$. By Lemma \ref{lem:Okonek_Teleman_1995_3-1} and Remark \ref{rmk:Projective_vortices_type_1_monopoles}, 
the pullback  by $j_{\ft,1}$ of the map \eqref{eq:PerturbedSO3MonopoleEquation_map} defined by the non-Abelian monopole equations \eqref{eq:PerturbedSO3MonopoleEquations} is equal to the map defining the projective vortex equations \eqref{eq:SO(3)_monopole_equations_almost_Hermitian_alpha}. Hence, the bijection
$\sM^0(E,h,\om) \leftrightarrow \sM_{\ft,1}^0$ defined by $j_{\ft,1}$ is a real analytic isomorphism. The results stated for $\sM_{\ft,1}^0$ and $\sM_{\ft,1}$ now follow from Items \eqref{item:Lubke_Teleman_6-3-7_analytic_embedding_rank-2_non-zero_section} and \eqref{item:Lubke_Teleman_6-3-7_set-theoretic_embedding} in Theorem \ref{thm:Lubke_Teleman_6-3-7}, respectively.

By Remarks \ref{rmk:Projective_vortices_type_2_monopoles} and \ref{rmk:Decoupling_Dirac_equation_and_equivalence_type_I_and_II_solutions}, the subspace $\sM_{\ft,2}^0$ is identified with $\sM^0(F,h,\om)$ when $F= K_X\otimes E^*$. Hence, the conclusions for $\sM_{\ft,2}^0$ and $\sM_{\ft,2}$ follow from the same argument as those for $\sM_{\ft,1}^0$ and $\sM_{\ft,1}$.
\end{proof}

Let $E$ be a complex vector bundle over a smooth, almost complex manifold $X$. We recall that $\sA^{0,1}(E)$ denotes the affine space of $(0,1)$-connections of class $W^{1,p}$ that induce a fixed, integrable $(0,1)$-connection $\bar\rd_{A_d}$ on the complex line bundle $\det E$ as in \eqref{eq:Holomorphic_structure_fixed_determinant}. We define
\begin{equation}
\label{eq:(0,1)PairQuotientSpaces}
\sC^{0,1}(E):=\left( \sA^{0,1}(E)\times W^{1,p}(E)\right)/W^{2,p}(\SL(E)),
\end{equation}
to be the corresponding configuration space equipped with quotient topology, where the $\SL(E)$-action in \eqref{eq:(0,1)PairQuotientSpaces} is defined in \eqref{eq:SL(E)ActionOn(0,1)Pairs}, and note that this quotient space need not be Hausdorff, unlike the quotient space of unitary pairs.

\begin{defn}[Standard $\CC^*$ action on affine and quotient spaces of $(0,1)$-pairs]
\label{defn:CStarActionOn(0,1)Pairs}
The \emph{standard $\CC^*$ action on the affine space of $(0,1)$-pairs} is defined by
\begin{equation}
  \label{eq:CZActionOnAffine}
  \CC^*\times \sA^{0,1}(E)\times \Omega^0(E) \ni
  \left( \la,(\bar\rd_A,\varphi)\right)\mapsto (\bar\rd_A,\la\varphi)
  \in \sA^{0,1}(E)\times \Omega^0(E).
\end{equation}
The action \eqref{eq:CZActionOnAffine} commutes with the action of $\SL(E)$ gauge transformations given in \eqref{eq:SL(E)ActionOn(0,1)Pairs} and thus induces the \emph{standard  action on the quotient space $\sC^{0,1}(E)$} given by
\begin{equation}
\label{eq:CActionsOn(0,1)PairQuotients}
\CC^*  \times \sC^{0,1}(E) \ni \left(\la, [\bar\rd_A,\varphi]\right)\mapsto [\bar\rd_A,\la\varphi] \in \sC^{0,1}(E).
\end{equation}
\end{defn}

The $S^1$ action on $\sC(E,h)$ in \eqref{eq:S1_Action_On_ConfigurationSpaceForProjectiveVortices} could be extended to a $\CC^*$ action like that in \eqref{eq:CActionsOn(0,1)PairQuotients} but we will not need that extension. Furthermore, the bijection \eqref{eq:Map_W1p_unitary_pairs_to_01-pairs} is equivariant with respect to the $\SU(E)$ action on the domain and the $\SU(E)$ action on the range (as a subgroup of $\SL(E)$) as well as with respect to the $\CC^*$-action given by scalar multiplication on the section. Hence, the bijection \eqref{eq:Map_W1p_unitary_pairs_to_01-pairs} defines a $\CC^*$-equivariant map $\sC(E,h)\to \sC^{0,1}(E)$.

We now identify the fixed points of the $\CC^*$ action \eqref{eq:CActionsOn(0,1)PairQuotients} on $\sC^{0,1}(E)$. 

\begin{lem}[Fixed points of $\CC^*$ action $\sC^{0,1}(E)$]
\label{lem:FixedPointsOfC*ActionOnQuotientSpace01Pairs}
Let $E$ be a complex vector bundle over a connected, almost complex manifold and let $[\bar\rd_A,\varphi]\in \sC^{0,1}(E)$. If $\varphi\equiv 0$ or $(\bar\rd_A,\varphi)$ is a split pair in the sense of Definition \ref{defn:Split_(0,1)-pair}, then $[\bar\rd_A,\varphi]$ is a fixed point of the $\CC^*$ action \eqref{eq:CActionsOn(0,1)PairQuotients}. 
\end{lem}

\begin{proof}
If $\varphi\equiv 0$ or $(\bar\rd_A,\varphi)$ is a split $(0,1)$-pair, then $[\bar\rd_A,\varphi]$ is a fixed point of the $\CC^*$ action by the same argument given in the proof of the analogous assertion in Proposition \ref{prop:FixedPointsOfS1ActionOnSpinuQuotientSpace}.
\end{proof}

We have the following complete characterization of the fixed points of the $\CC^*$ action on the moduli space of polystable holomorphic pairs when the rank of the bundle is two.

\begin{lem}[Fixed points of $\CC^*$ action on $\fM_\ps(E,\omega)$]
\label{lem:FixedPointsOfC*ActionOnPolystablePairs}
Let $(E,h)$ be a rank-two, smooth Hermitian vector bundle over a closed complex K\"ahler manifold $(X,\omega)$.
The moduli space $\fM_\ps(E,\omega)$ of polystable holomorphic pairs in \eqref{eq:Moduli_space_polystable_holomorphic_pairs} is closed under the $\CC^*$ action  \eqref{eq:CActionsOn(0,1)PairQuotients} and the fixed points of this action on  $\fM_\ps(E,\omega)$ are represented by smooth pairs $(\bar\partial_E,\varphi)$ such that
\begin{enumerate}
\item
\label{item:PolystableFixedPointsZeroSection}
$\varphi\equiv 0$, so $(\bar\partial_E,\varphi)$ is a zero-section pair, or
\item
\label{item:PolystableFixedPointsSplit}
$\varphi\not\equiv 0$ and $(\bar\partial_E,\varphi)$ is split as in Definition \ref{defn:Split_(0,1)-pair} with respect to a decomposition $E=L_1\oplus L_2$, where $\deg L_1 < \deg L_2$.
\end{enumerate}
\end{lem}

\begin{proof}
The condition defining polystability in Definition \ref{defn:StablePair_arbitrary_rank} is clearly invariant under scalar multiplication by $\CC^*$ on the section, so the moduli space $\fM_\ps(E,\omega)$ of polystable holomorphic pairs in \eqref{eq:Moduli_space_polystable_holomorphic_pairs} is closed under the $\CC^*$ action. Lemma \ref{lem:FixedPointsOfC*ActionOnQuotientSpace01Pairs} implies that the gauge-equivalence classes of zero-section and split pairs in $\fM_\ps(E,\omega)$ are fixed points of the $\CC^*$ action \eqref{eq:CActionsOn(0,1)PairQuotients}.

To see that any fixed point of the action on $\fM_\ps(E,\omega)$ must be the gauge-equivalence class of a zero-section or split pair, we first observe that if $[\bar\partial_E,\varphi]\in \sC^{0,1}(E)$ is a fixed point of the $\CC^*$ action \eqref{eq:CActionsOn(0,1)PairQuotients}, then it is a fixed point of the $S^1$ action on
$\sC^{0,1}(E)$ defined by the $\CC^*$ action \eqref{eq:CActionsOn(0,1)PairQuotients} by considering $S^1$ as a subgroup of $\CC^*$.  Next, the bijection \eqref{eq:Map_W1p_unitary_pairs_to_01-pairs} is $S^1$-equivariant with respect to the actions defined by scalar multiplication on the section and so the set-theoretic embedding  in \eqref{eq:Lubke_Teleman_6-3-7_set-theoretic_embedding},
\[
\sM(E,h,\omega) \hookrightarrow \fM(E),
\]
is $S^1$-equivariant with respect to the $S^1$ action \eqref{eq:S1_Action_On_ConfigurationSpaceForProjectiveVortices} on the domain and the $S^1$ action on the codomain \eqref{eq:CActionsOn(0,1)PairQuotients} implied by considering $S^1$ as a subgroup of $\CC^*$. Thus, if a point in $\fM_\ps(E,\omega)$ is a fixed point of the $\CC^*$ action \eqref{eq:CActionsOn(0,1)PairQuotients}, then it is the image under the map \eqref{eq:Lubke_Teleman_6-3-7_set-theoretic_embedding} of a fixed point of the $S^1$ action
\eqref{eq:S1_Action_On_ConfigurationSpaceForProjectiveVortices} on $\sM(E,h,\omega)$.  By Proposition \ref{prop:FixedPointsOfS1ActionOnUnitaryQuotientSpace}, the fixed points of the $S^1$ action \eqref{eq:S1_Action_On_ConfigurationSpaceForProjectiveVortices} on $\sM(E,h,\omega)$ are the gauge-equivalence classes of zero-section and split projective vortices. The embedding \eqref{eq:Lubke_Teleman_6-3-7_set-theoretic_embedding} sends
\begin{inparaenum}[\itshape i\upshape)]
\item gauge-equivalence classes of zero-section projective vortices to gauge-equivalence classes of zero-section holomorphic pairs and
\item gauge-equivalence classes of split projective vortices to gauge-equivalence classes of split holomorphic pairs.
\end{inparaenum}
Thus, fixed points of the $\CC^*$ action \eqref{eq:CActionsOn(0,1)PairQuotients} on $\fM_\ps(E,\omega)$ are gauge-equivalence classes of zero-section holomorphic pairs or of non-zero-section split holomorphic pairs. The conclusion for non-zero-section, split polystable holomorphic pairs that $\deg L_1 < \deg L_2$ follows from Remark \ref{rmk:Stability_of_Split_Pairs}.
\end{proof}

\begin{rmk}[Fixed points of the circle action on the configuration space of unitary pairs and the circle-invariant moduli subspace of projective vortices]
\label{rmk:Fixed_points_circle_action_configuration_space_and_moduli_subspace_non-Abelian_monopoles}
Fixed points of the $S^1$ action on the moduli space of projective vortices $\sM(E,\om,h)$ are given by the intersection of $\sM(E,\om,h)$ with the fixed points of the $S^1$ action on $\sC(E,h)$ identified in
Proposition \ref{prop:FixedPointsOfS1ActionOnUnitaryQuotientSpace}.
This contrasts with the behavior of critical points of a smooth function: $[A,\varphi]\in\sM(E,\om,h)$ can be a critical point of the restriction of a smooth function $f$ on $\sC(E,h)$ to the subspace $\sM(E,\om,h)$ without being a critical point of $f$ on $\sC(E,h)$.  As discussed in Chapter \ref{chap:Hessian_restriction_smooth_function_submanifold_Euclidean_space}, finding critical points of the restriction of a function to a subspace can be considerably more difficult than finding critical points of the function in the ambient space.  The bijection between fixed points of the $S^1$ action and critical points of a Hamiltonian function given in Theorem \ref{mainthm:Frankel_almost_Hermitian} allows us to
find critical points using  the simpler criterion of finding fixed points.
\end{rmk}

\section[Moduli space of projective vortices as a symplectic quotient]{Marsden--Weinstein symplectic reduction and the moduli space of projective vortices as a symplectic quotient}
\label{sec:Marsden-Weinstein_reduction_moduli_space_SO3_monopoles_symplectic_quotient}
We now construct a symplectic structure on $\sM^0(E,h,\omega)$, the moduli space of non-zero-section projective vortices
on a Hermitian vector bundle $(E,h)$ over a complex, K\"ahler manifold $(X,g,J)$ with K\"ahler two-form $\omega = g(\cdot,J\cdot)$. We define a weak symplectic structure on the affine space of unitary pairs and compute a moment map for the action of the group of special unitary gauge transformations on this space in
Section \ref{subsec:Symplectic_form_moment_map_affine_space_unitary_pairs}. In Section \ref{subsec:Properties_moment_map_affine_space_unitary_pairs}, we characterize the regular points of this moment map.
We construct a Kuranishi model for an open neighborhood of a holomorphic pair in the moduli space of holomorphic pairs in Section \ref{subsec:Kuranishi_model_holomorphic_pair}. We apply the Marsden--Weinstein Symplectic Reduction Theorem \ref{thm:Marsden-Weinstein_symplectic_quotient} to the moduli space of projective vortices
in Sections \ref{subsec:Marsden-Weinstein_reduction_moduli_space_SO3_monopoles_symplectic_quotient} and \ref{subsec:Marsden-Weinstein_reduction_moduli_space_SO3_monopoles_symplectic_quotient_singular_points}, 
considering separately and respectively the cases where $[A,\varphi]$ is a regular or singular point of the moduli space $\sM(E,h,\omega)$. In Section \ref{subsec:Characterization_critical_points_Hamiltonian_function_circle_action_on_projective_vortices}, we discuss the characterization of critical points of the Hamiltonian function for the circle action on the moduli space of projective vortices.

\subsection[Symplectic form and moment map on affine space of unitary pairs]{Symplectic form and moment map on affine space of unitary pairs on a Hermitian vector bundle over an almost Hermitian manifold}
\label{subsec:Symplectic_form_moment_map_affine_space_unitary_pairs}
In this subsection, we describe how the open subset of smooth points of the moduli space of projective vortices over a complex, K\"ahler manifold inherits a K\"ahler metric. A similar result was proved by L\"ubke and Teleman as \cite[Theorem 6.2.8, p. 68]{Lubke_Teleman_2006} for the moduli space of type $1$ non-Abelian monopoles over a complex, K\"ahler surface.

Let $(E,h)$ be a  Hermitian vector bundle over an almost complex manifold $(X,J)$. We recall from \eqref{eq:Kobayashi_7-6-12} that the almost complex structure $J$ on $TX$ induces an almost complex structure on the real vector space $\Omega^1(X) = \Omega^0(T^*X)$ and thus the real vector space $\Omega^1(\su(E)) = \Omega^0(T^*X\otimes\su(E))$ while the vector space $\Omega^0(E)$ has an almost complex structure induced by multiplication by $\CC$ on the fibers of $E$. Therefore, the vector space $\Omega^1(\su(E))\oplus \Omega^0(E)$ has an almost complex structure defined by
\begin{equation}
  \label{eq:Almost_complex_structure_affine_space_pairs_unitary_connections_and_sections}
  \bJ(a,\phi) := (-ia' + ia'', i\phi),
\end{equation}
for all $a = \frac{1}{2}(a'+a'') \in \Omega^1(\su(E))$ with
$a'' \in \Omega^{0,1}(\fsl(E))$ and $a' = -(a'')^\dagger \in \Omega^{1,0}(\fsl(E))$,
following our convention \eqref{eq:Decompose_a_in_Omega1suE_into_10_and_01_components}, and $\phi \in \Omega^0(E)$. Similarly, we may define a (weak or $L^2$) Riemannian metric $\bg$ that is compatible with $\bJ$ on these vector spaces by extending its definition on $\Omega^1(\su(E))$ to give (following the convention of Kobayashi \cite[Equation (7.6.5), p. 261]{Kobayashi_differential_geometry_complex_vector_bundles} relating Riemannian and Hermitian inner products),
\begin{equation}
  \label{eq:L2_metric_affine_space_pairs_unitary_connections_and_sections}
  \bg((a_1,\phi_1), (a_2,\phi_2)) := (a_1,a_2)_{L^2(X)} + \Real(\phi_1,\phi_2)_{L^2(X)},
\end{equation}
for all $a_j \in \Omega^1(\su(E))$ and $\phi_j \in \Omega^0(E)$, for $j=1,2$. We then have the following analogue of Proposition \ref{prop:Donaldson_Kronheimer_6-5-7_almost_Hermitian}.

\begin{prop}[Weak symplectic form on the affine space of unitary pairs on a Hermitian vector bundle over an almost Hermitian manifold]
\label{prop:Donaldson_Kronheimer_6-5-7_almost_Hermitian_pairs}
Let $(E,h)$ be a smooth Hermitian vector bundle over a closed, smooth almost Hermitian manifold $(X,g,J)$ of real dimension $n$. If $p \in (n,\infty)$ is a constant, then the expression
\begin{multline}
  \label{eq:Kobayashi_7-6-22_pairs}
  \bomega((a_1,\phi_1),(a_2,\phi_2)) := \bg((a_1,\phi_1), \bJ(a_2,\phi_2)),
  \\
  \text{for all } a_1,a_2 \in W^{1,p}(T^*X\otimes\su(E))
  \text{ and } \phi_1,\phi_2 \in W^{1,p}(E),
\end{multline}
defines a weak symplectic form in the sense of Definition \ref{defn:Symplectic_form_Banach_manifold} on the real Banach space $W^{1,p}(T^*X\otimes\su(E)) \oplus W^{1,p}(E)$.
\end{prop}

We observe that the Banach Lie group $W^{2,p}(\SU(E))$ acts isometrically on the Banach affine space $\sA(E,h)\times W^{1,p}(E)$ with respect to the Riemannian metric $\bg$ in \eqref{eq:L2_metric_affine_space_pairs_unitary_connections_and_sections}
as well as symplectically in the sense of Definition \ref{defn:Symplectic_map}:
\begin{align*}
  \bg(u(a_1,\phi_1), u(a_2,\phi_2)) &= \bg((a_1,\phi_1), (a_2,\phi_2)),
  \\
  \bomega(u(a_1,\phi_1),u(a_2,\phi_2)) &= \bomega((a_1,\phi_1),(a_2,\phi_2)),
  \\
                                    &\qquad \text{for all } u \in W^{2,p}(\SU(E)), a_1,a_2 \in W^{1,p}(T^*X\otimes\su(E)),
                                      \\
  &\qquad \text{ and } \phi_1,\phi_2 \in W^{1,p}(E).
\end{align*}
Note that the action of $W^{2,p}(\SU(E))$ in the preceding equalities is the action on the tangent space of $\sA(E,h)\times W^{1,p}(E)$ given by the derivative of the action \eqref{eq:DefineSU(E)ActionOnUnitaryPairs} on the affine space.  This action also preserves the almost complex structure $\bJ$ in \eqref{eq:Almost_complex_structure_affine_space_pairs_unitary_connections_and_sections}.

There is a K\"ahler metric $\bh$ on $\sA^{0,1}(E)\times W^{1,p}(E)$ given by the following generalization of the definition \eqref{eq:Kobayashi_7-6-20} of a K\"ahler metric on $\sA^{0,1}(E)$,
\begin{multline}
\label{eq:Kobayashi_7-6-20_pairs}    
\bh((\alpha_1,\phi_1),(\alpha_2,\phi_2))
:= \int_X\frac{n}{i}\tr_E(\alpha_1\wedge\alpha_2^\dagger)\wedge\omega^{n-1}
+ \int_X \langle\phi_1,\phi_2\rangle_h\,\omega^n,
\\
\text{for all } (\alpha_i, \phi_i) \in W^{1,p}\left(\Lambda^{0,1}(\fsl(E))\oplus E\right),
\text{ for } i=1,2.
\end{multline}
We have the following analogue of Proposition \ref{prop:Donaldson_Kronheimer_6-5-8_almost_Kaehler}.

\begin{prop}[Equivariant moment map for action of Banach Lie group of determinant-one, unitary gauge transformations on the affine space of unitary pairs over an almost K\"ahler manifold]
\label{prop:Donaldson_Kronheimer_6-5-8_almost_Kaehler_type1_pairs}
(Compare Bradlow and Garc{\'\i}a--Prada \cite[Equation (6.1), p. 580, and last two paragraphs of Section 6.1]{BradlowGP}, Okonek and Teleman \cite[Proposition 3.2]{OTVortex} and \cite[Discussions following Propositions 3.1.4 and 3.21]{OTSurvey}, and Teleman \cite[Section 1.2.1 and Remark 2.2.3]{TelemanNonabelian}.)  
Continue the hypotheses of Proposition \ref{prop:Donaldson_Kronheimer_6-5-7_almost_Hermitian_pairs}, but assume now that $\omega$ is closed, and thus $(X,g,J)$ is almost K\"ahler, and fix a smooth, unitary connection $A_d$ on $\det E$. Then
\begin{equation}
\label{eq:Moment_map_action_unitary_det_one_gauge_transformations_affine_space_pairs}  
\bmu:\sA(E,h)\times W^{1,p}(E) \to \left(W^{2,p}(\su(E))\right)^*,
\end{equation}
where
\begin{equation}
\label{eq:Moment_map_action_unitary_det_one_gauge_transformations_affine_space_pairs_expression}
  \langle \xi, \bmu(A,\varphi) \rangle
  := \left(\xi, \Lambda (F_A)_0 - \frac{i}{2}(\varphi\otimes\varphi^*)_0\right)_{L^2(X)},
\quad\text{for all } \xi \in W^{2,p}(\su(E)),
\end{equation}
is an equivariant moment map in the sense of \eqref{eq:Kobayashi_7-5-3} for the action of the Banach Lie group $W^{2,p}(\SU(E))$ of determinant-one, unitary $W^{2,p}$ automorphisms of $E$ for the (weak) symplectic form $\bomega$ given by \eqref{eq:Kobayashi_7-6-22_pairs} on the Banach affine space $\sA(E,h) \times W^{1,p}(E)$, where $\sA(E,h) = A_0 + W^{1,p}(T^*X\otimes \su(E))$, and $A_0$ is a smooth, unitary connection on $E$ that induces $A_d$ on $\det E$.
\end{prop}

\begin{rmk}[Moment map for the affine space of pairs and the non-Abelian monopole equations]
\label{rmk:Moment_map_SO3_monopoles}  
Okonek and Teleman \cite[Section 6.2]{OTMasterCoupling} imply that there is a moment map for the affine space of pairs $\sA(E,h) \times W^{1,p}(W^+\otimes E)$ with symplectic form defined by the Riemannian metric and compatible almost complex structure on this space, but we shall restrict our attention to the simpler case of type $1$ pairs $(A,\varphi)$ provided by Proposition \ref{prop:Donaldson_Kronheimer_6-5-8_almost_Kaehler_type1_pairs} when $X$ has real dimension four.
\end{rmk}

\begin{proof}[Proof of Proposition \ref{prop:Donaldson_Kronheimer_6-5-8_almost_Kaehler_type1_pairs}]
For convenience, let us write $\bmu(A,\varphi) := \bmu_1(A)+\bmu_2(\varphi)$, where we set 
\begin{align*}  
  \bmu_1(A) &:= \Lambda (F_A)_0, \quad\text{for all } A \in \sA(E,h),
  \\
  \bmu_2(\varphi) &:= - \frac{i}{2}(\varphi\otimes\varphi^*)_0, \quad\text{for all } \varphi \in W^{1,p}(E).
\end{align*}
Let $\bomega_1$ denote the weak symplectic form on $\sA(E,h)$ given by \eqref{eq:Kobayashi_7-6-22} and let $\bX^1$ be the smooth section of $\Hom(L(\sG_E),T\sA(E,h)) \cong T\sA(E,h)\otimes L(\sG_E)^*$ generated by the smooth action of the Banach Lie group $\sG_E = W^{2,p}(\SU(E))$ on the Banach affine space $\sA(E,h)$, so by \eqref{eq:Differential_of_SU(E)_GaugeAction} we have
\[
  \bX_\xi^1(A) = d_A\xi \in T_A\sA(E,h), \quad\text{for all } A \in \sA(E,h) \text{ and } \xi \in W^{2,p}(\su(E)),
\]
where $\sG_E$ has Lie algebra given by
\[
  L(\sG_E) := T_{\id_E}\sG_E = T_{\id}W^{2,p}(\SU(E)) = W^{2,p}(\su(E)).
\]
Proposition \ref{prop:Donaldson_Kronheimer_6-5-8_almost_Kaehler} and the identity \eqref{eq:Moment_map_action_unitary_det_one_gauge_transformations_affine_space_unitary_connections_expression_tracefree} in Remark \ref{rmk:Moment_map_Einstein_constant} imply that
\begin{multline}
\label{eq:Moment_map_action_unitary_det_one_gauge_transformations_affine_space_pairs_expression_1}
  \langle\xi,d\bmu_1(A)a\rangle = \iota_{\bX_\xi^1(A)}\bomega_1(a),
  \\
  \text{for all } \xi \in W^{2,p}(\su(E)) \text{ and } a \in T_A\sA(E,h) = W^{1,p}(T^*X\otimes\su(E)).
\end{multline}
Indeed, to see this, observe that the definition of $\bmu_1$ gives
\[
  \langle \xi,d\bmu_1(A)a\rangle = (\xi,\Lambda d_Aa)_{L^2(X)},
  \quad\text{for all } \xi \in W^{2,p}(\su(E)) \text{ and }
  a \in W^{1,p}(T^*X\otimes\su(E)).
\]
On the other hand, according to \eqref{eq:Kobayashi_7-6-22}, we have
\begin{multline*}
  \bomega_1(\bX_\xi^1(A),a) = n\int_X \tr_E(d_A\xi\wedge a)\wedge\omega^{n-1}
  = -n\int_X \tr_E(\xi d_Aa)\wedge\omega^{n-1}
  \\
  = -\int_X \tr_E(\xi\Lambda d_Aa)\wedge\omega^n
  = -(\xi,\Lambda d_Aa)_{L^2(X)} = \langle\xi,d\bmu_1(A)a\rangle,
  \\
  \text{for all } \xi \in W^{2,p}(\su(E)) \text{ and }
  a \in W^{1,p}(T^*X\otimes\su(E)).
\end{multline*}
To obtain the second equality above, we used integration by parts and the fact that $d\omega=0$:
\[
  d\left(\tr_E(\xi a)\wedge\omega^{n-2}\right)
  =
  \tr_E(d_A\xi\wedge a)\wedge\omega^{n-2} + \tr_E(\xi d_Aa)\wedge\omega^{n-2}.
\]
Therefore, $\bmu_1$ and $\bomega_1$ obey the key defining relation \eqref{eq:Kobayashi_7-5-2} for a moment map for the action of $W^{2,p}(\SU(E))$ on $\sA(E,h)$ and this verifies the claimed identity \eqref{eq:Moment_map_action_unitary_det_one_gauge_transformations_affine_space_pairs_expression_1}.

Next, we claim that
\begin{equation}
\label{eq:Moment_map_action_unitary_det_one_gauge_transformations_affine_space_pairs_expression_2}
  \langle\xi,d\bmu_2(\varphi)\sigma\rangle = \iota_{\bX_\xi^2(\varphi)}\bomega_2(\sigma),
  \quad\text{for all } \xi \in W^{2,p}(\su(E)) \text{ and } \sigma \in W^{1,p}(E),
\end{equation}
where, by the convention of Kobayashi \cite[Equation (7.6.8), p. 251]{Kobayashi_differential_geometry_complex_vector_bundles},
\[
  \bomega_2(\sigma_1,\sigma_2) := \Real(\sigma_1,i\sigma_2)_{L^2(X)} = \frac{1}{n!}\int_X\Real\langle\sigma_1,i\sigma_2\rangle_E\, \omega^n,
  \quad\text{for all } \sigma_1,\sigma_2 \in W^{1,p}(E),
\]
and $\bX^2$ is the section of $\Hom(L(\sG_E),TW^{1,p}(E)) \cong TW^{1,p}(E)\otimes L(\sG_E)^*$ generated by the action of $W^{2,p}(\SU(E))$ on $W^{1,p}(E)$, so \eqref{eq:Differential_of_SU(E)_GaugeAction} gives
\[
  \bX_\xi^2(\varphi) = -\xi\varphi \in T_\varphi W^{1,p}(E) = W^{1,p}(E),
  \quad\text{for all } \text{for all } \varphi \in W^{1,p}(E) \text{ and } \xi \in W^{2,p}(\su(E)).
\]
Therefore,
\[
  \bomega_2(\bX_\xi^2(\varphi),\sigma) = -\Real(\xi\varphi,i\sigma)_{L^2(X)},
  \quad\text{for all } \xi \in W^{2,p}(\su(E)) \text{ and } \sigma \in W^{1,p}(E).
\]
The definition of $\bmu_2$ and the facts that $\xi \in W^{2,p}(\su(E))$ with $\xi^\dagger = -\xi$ and $\fu(E) \cong i\underline{\RR}\oplus\su(E)$ as an orthogonal direct sum of Riemannian vector bundles give
\begin{align*}
  \langle \xi,d\bmu_2(\varphi)\sigma\rangle
  &= -\frac{1}{2}(\xi,i\sigma\otimes\varphi^* + i\varphi\otimes\sigma^*)_{L^2(X)}
  \\
  &= -\frac{1}{2}(\xi\varphi,i\sigma)_{L^2(X)} - \frac{1}{2}(\xi\sigma,i\varphi)_{L^2(X)}
  \\
  &= -\frac{1}{2}(\xi\varphi,i\sigma)_{L^2(X)} + \frac{1}{2}(\sigma,\xi i\varphi)_{L^2(X)}
  \\
  &= -\frac{1}{2}(\xi\varphi,i\sigma)_{L^2(X)} - \frac{1}{2}(i\sigma,\xi\varphi)_{L^2(X)}
  \\
  &= -\Real(\xi\varphi,i\sigma)_{L^2(X)},
  \quad\text{for all } \xi \in W^{2,p}(\su(E)) \text{ and } \sigma \in W^{1,p}(E).
\end{align*}
Therefore, $\bmu_2$ and $\bomega_2$ obey the key defining relation \eqref{eq:Kobayashi_7-5-2} for a moment map for the action of $W^{2,p}(\SU(E))$ on $W^{1,p}(E)$ and this verifies the claimed identity \eqref{eq:Moment_map_action_unitary_det_one_gauge_transformations_affine_space_pairs_expression_2}.

By combining the two identities \eqref{eq:Moment_map_action_unitary_det_one_gauge_transformations_affine_space_pairs_expression_1} and \eqref{eq:Moment_map_action_unitary_det_one_gauge_transformations_affine_space_pairs_expression_2}, we see that for any pair $(A,\varphi) \in \sA(E,h)\times W^{1,p}(E)$ we have
\begin{multline*}
  \langle \xi,d\bmu(A,\varphi)(a,\sigma)\rangle
  =
  \langle \xi,d\bmu_1(A)a\rangle + \langle \xi,d\bmu_1(\varphi)\sigma\rangle
  \\
  =
  \bomega_1(\bX_\xi^1(A),a) + \bomega_2(\bX_\xi^2(\varphi),\sigma)
  =
  \bomega(\bX_\xi(A,\varphi),(a,\sigma)),
  \\
  \text{for all } \xi \in W^{2,p}(\su(E)) \text{ and } (a,\sigma) \in W^{1,p}(\su(E))\oplus W^{1,p}(E).
\end{multline*}
Here, $\bX := \bX^1+\bX^2$ is the section of
\[
  \Hom\left(L(\sG_E),T\left(\sA(E,h)\times W^{1,p}(E)\right)\right)
  \cong \left(T\sA(E,h)\times TW^{1,p}(E)\right)\otimes L(\sG_E)^*
\]
generated by the action of $W^{2,p}(\SU(E))$ on the Banach affine space $\sA(E,h)\times W^{1,p}(E)$, so
\begin{align*}
  \bX_\xi(A,\varphi)
  &= \bX_\xi^1(A) + \bX_\xi^2(\varphi)
  \\
  &= (d_A\xi,-\xi\varphi) \in T_A\sA(E,h)\oplus T_\varphi W^{1,p}(E) = W^{1,p}(\su(E))\oplus W^{1,p}(E),
  \\
  &\qquad\text{for all } \varphi \in (A,\varphi) \in \sA(E,h)\times W^{1,p}(E) \text{ and } \xi \in W^{2,p}(\su(E)).
\end{align*}
Therefore, $\bmu$ and $\bomega$ obey the key defining relation \eqref{eq:Kobayashi_7-5-2} for a moment map for the action of $W^{2,p}(\SU(E))$ on $\sA(E,h)\times W^{1,p}(E)$ and this verifies the claimed identity \eqref{eq:Moment_map_action_unitary_det_one_gauge_transformations_affine_space_pairs_expression}.

Equivariance in the sense of \eqref{eq:Kobayashi_7-5-3} of $\bmu_1$ with respect to the action of $W^{2,p}(\SU(E))$ required by \eqref{eq:Kobayashi_7-5-3} follows from Proposition \ref{prop:Donaldson_Kronheimer_6-5-8_almost_Kaehler} and equivariance of $\bmu_2$ and thus $\bmu$ with respect to the action of $W^{2,p}(\SU(E))$ follows from the proof of Proposition \ref{prop:Donaldson_Kronheimer_6-5-8_almost_Kaehler}.
\end{proof}

It is convenient to define the following analogue of the set in $\sH(E,h)$ in \eqref{eq:Subset_unitary_connections_defining_holomorphic_structure_on_E}, but now for pairs,
\begin{equation}
  \label{eq:Subset_pre-holomorphic_pairs}
  \sP(E,h)
  =
  \left\{(A,\varphi) \in \sA(E,h)\times W^{1,p}(E): (F_A^{0,2})_0 = 0 \text{ and } \bar\partial_A\varphi = 0 \right\}.
\end{equation}
The expression \eqref{eq:Moment_map_action_unitary_det_one_gauge_transformations_affine_space_pairs_expression} yields the identification
\begin{equation}
  \label{eq:Kobayashi_7-6-33_pairs}
  \bmu^{-1}(0)\cap\sP(E,h)
  =
  \left\{ (A,\varphi) \in \sP(E,h): \Lambda (F_A)_0 = \frac{i}{2}(\varphi\otimes\varphi^*)_0 \right\},
\end{equation}
that is, $(A,\varphi) \in \bmu^{-1}(0)\cap\sP(E,h)$ if and only if $(A,\varphi)\in\sP(E,h)$ obeys the projective vortex equation \eqref{eq:SO(3)_monopole_equations_(1,1)_curvature_alpha}; the condition $A_{\det E} = A_d$ is assured by the assumption that all $A\in\sA(E,h)$ obey the determinant condition \eqref{eq:Unitary_connection_detE_fixed}.

Hence, by \eqref{eq:SO(3)_monopole_equations_almost_Hermitian_alpha} we see that $(A,\varphi)$ is a projective vortex (or non-Abelian monopole of type $1$ when $X$ has complex dimension two) on $E$ if and only if it belongs to the subspace $\bmu^{-1}(0)\cap\sP(E,h)$ in \eqref{eq:Kobayashi_7-6-33_pairs}. (Again, this identification is valid regardless of whether $J$ is integrable or $\omega$ is closed, so we may allow $(X,g,J)$ to be almost Hermitian.)

\subsection[Properties of the moment map on the space of unitary pairs]{Properties of the moment map on the Banach affine space of unitary pairs}
\label{subsec:Properties_moment_map_affine_space_unitary_pairs}
By a minor modification of the short argument due to Kobayashi following \cite[Equation (7.6.34), p. 255]{Kobayashi_differential_geometry_complex_vector_bundles}
(see the discussion around \eqref{eq:Restriction_of_MomentMap_To_bV_cap_sH*(E,h)} and Lemma \ref{lem:MomentMapRegularityCriterion}), we obtain the

\begin{lem}[Regular points of the moment map on the affine space of pairs]
\label{lem:Regular_value_moment_map_affine_space_pairs}  
Continue the hypotheses of Proposition \ref{prop:Donaldson_Kronheimer_6-5-8_almost_Kaehler_type1_pairs} and assume in addition that the almost complex structure $J$ is integrable, so $(X,g,J)$ is a complex K\"ahler manifold. If the $\Omega^0(\su(E))$-component of the harmonic space $\bH_{A,\varphi}^2$ in \eqref{eq:H_Avarphi^bullet} vanishes, then $(A,\varphi) \in \bmu^{-1}(0)\cap(\sA(E,h)\times W^{1,p}(E))$ is a regular point of the moment map $\bmu$ in \eqref{eq:Moment_map_action_unitary_det_one_gauge_transformations_affine_space_pairs} and an open neighborhood of $(A,\varphi)$ in $\bmu^{-1}(0)\cap(\sA(E,h)\times W^{1,p}(E))$ is an embedded, complex submanifold of the affine space $\sA(E,h)\times W^{1,p}(E)$.
\end{lem}

\begin{proof}
From the definition of the moment map $\bmu$, we have, for any pair $(A,\varphi) \in \sA(E,h)\times W^{1,p}(E)$,
\[
  \langle \xi, \bmu(A,\varphi) \rangle
  = \left(\xi, \Lambda (F_A)_0 - \frac{i}{2}(\varphi\otimes\varphi^*)_0\right)_{L^2(X)},
\quad\text{for all } \xi \in W^{2,p}(\su(E)).
\]
It is convenient to define the co-moment map implied by the definition of $\bmu$ in \eqref{eq:Moment_map_action_unitary_det_one_gauge_transformations_affine_space_pairs}, 
\begin{equation}
  \label{eq:Co-moment_map_action_unitary_det_one_gauge_transformations_affine_space_pairs}
  \bmu^*:\sA(E,h)\times W^{1,p}(E) \ni (A,\varphi) \mapsto \Lambda (F_A)_0 - \frac{i}{2}(\varphi\otimes\varphi^*)_0 \in L^p(\su(E)),
\end{equation}
where we note that, for $p' = p/(p-1) \in (1,n/(n-1))$, there is a continuous embedding of Sobolev spaces,
\[
  L^p(\su(E)) \subset W^{-2,p'}(\su(E)),
\]
since $p\in (n,\infty)$ by hypothesis and $2n$ is the real dimension of $X$. (Recall that the Sobolev Embedding Theorem (see Adams and Fournier \cite[Theorem 4.12, p. 85]{AdamsFournier}) yields a continuous embedding $W^{2,p}(X) \subset C^0(X)$ of Sobolev spaces and, of course, there is a continuous embedding $C^0(X) \subset L^{p'}$ for $p' = p/(p-1) \in (1,n/(n-1))$. If $T:X_1\to X_2$ is a bounded linear map of Banach spaces, then the adjoint $T^*:X_2^*\to X_1^*$ is a bounded linear map of Banach spaces by Rudin \cite[Theorem 4.10, p. 98]{Rudin}, where $X_i^*$ is the normed dual of $X_i$ for $i=1,2$, and if $\Ran T$ is dense in $X_2$, then $T^*$ is injective by \cite[Corollary 4.12 (b), p. 99]{Rudin}. Because the ranges of the preceding Sobolev embeddings are dense by \cite[Theorem 3.17, p. 68]{AdamsFournier}, duality thus yields a continuous embedding of Sobolev spaces, $L^p(X) \subset W^{-2,p'}(X)$, as claimed.) Hence, we obtain
\begin{multline*}
  d\bmu^*(A,\varphi)(a,\sigma)
  =
  \Lambda d_Aa - \frac{i}{2}(\varphi\otimes\sigma^* + \sigma\otimes\varphi^*)_0,
  \\
  \text{for all } (a,\sigma) \in W^{1,p}(T^*X\otimes\su(E)\oplus E).
\end{multline*}
On the other hand, from the definition \eqref{eq:d1_projective_vortex_elliptic_deformation_complex} of the differential $d_{A,\varphi}^1$, its $\Omega^0(\su(E))$-component $d_{A,\varphi}^{1;0}$ is given by
\begin{multline*}
d_{A,\varphi}^{1;0}(a,\sigma)
=
\Lambda(\bar\partial_Aa' + \partial_Aa'') - i(\varphi\otimes\sigma^*+\sigma\otimes\varphi^*)_0 \in L^p(\su(E)),
\\
\text{for all } (a,\sigma) \in W^{1,p}(T^*X\otimes\su(E)\oplus E),
\end{multline*}
where we write $a=\frac{1}{2}(a'+a'')$, with $a''\in  W^{1,p}(\Lambda^{0,1}(\fsl(E)))$ and $a'=-(a'')^\dagger \in W^{1,p}(\Lambda^{1,0}(\fsl(E)))$. Since the almost complex structure $J$ is integrable, we have from \eqref{eq:d_A_sum_components_almost_complex_manifold_integrable} that
\[
  d_Aa
  =
  \frac{1}{2}(\partial_A + \bar\partial_A)(a'+a'')
  =
  \frac{1}{2}(\partial_A a' + \bar\partial_Aa' + \partial_Aa'' + \bar\partial_Aa''),
\]
and thus, by definition of the dual Lefshetz operator $\Lambda$ in Section \ref{subsec:Hodge_star_Lefshetz_dual_Lefshetz_operators},
\[
  \Lambda d_Aa = \frac{1}{2}\Lambda(\bar\partial_Aa' + \partial_Aa''),
\]
since $\partial_A a' \in \Omega^{2,0}(\fsl(E))$ and $\bar\partial_A a'' \in \Omega^{0,2}(\fsl(E))$. Therefore, by combining the preceding identities, we see that
\begin{multline}
  \label{eq:dbmu*_Avarphi_asigma}
  d\bmu^*(A,\varphi)(a,\sigma)
  =
  \frac{1}{2}\Lambda(\bar\partial_Aa' + \partial_Aa'') - \frac{i}{2}(\varphi\otimes\sigma^* + \sigma\otimes\varphi^*)_0,
  \\
  \text{for all } (a,\sigma) \in W^{1,p}(T^*X\otimes\su(E)\oplus E).
\end{multline}
Consequently, we obtain
\begin{equation}
  \label{eq:dbmu_Avarphi_equals_dAvarphi1;0}
  d\bmu^*(A,\varphi)(a,\sigma) = \frac{1}{2}d_{A,\varphi}^{1;0}(a,\sigma),
  \quad\text{for all } (a,\sigma) \in W^{1,p}(T^*X\otimes\su(E)\oplus E).
\end{equation}
Now assume that $(A,\varphi) \in \bmu^{-1}(0)\cap(\sA(E,h)\times W^{1,p}(E))$. By hypothesis, the $\Omega^0(\su(E))$-component of the harmonic space $\bH_{A,\varphi}^2$ in \eqref{eq:H_Avarphi^bullet} vanishes and so the following operator is surjective:
\[
  d_{A,\varphi}^{1;0}:W^{1,p}(T^*X\otimes\su(E)\oplus E) \to L^p(\su(E)).
\]
By the preceding identity, the differential
\[
  d\bmu^*(A,\varphi):W^{1,p}(T^*X\otimes\su(E)\oplus E) \to L^p(\su(E))
\]
is also surjective and the differential
\[
  d\bmu(A,\varphi):W^{1,p}(T^*X\otimes\su(E)\oplus E) \to \left(W^{2,p}(\su(E))\right)^*
\]
is surjective (see Lemma \ref{lem:MomentMapRegularityCriterion}), noting that $L^p(\su(E)) \subset W^{2,p}(\su(E))$ and so $W^{2,p}(\su(E))^* \subset (L^p(\su(E)))^*$. Since the pair $(A,\varphi) \in \bmu^{-1}(0)\cap(\sA(E,h)\times W^{1,p}(E))$ was arbitrary, the conclusion follows.
\end{proof}

The proofs of Proposition \ref{prop:Donaldson_Kronheimer_6-5-8_almost_Kaehler_type1_pairs} and Lemma \ref{lem:Regular_value_moment_map_affine_space_pairs} yield the

\begin{lem}[Projective vortices with vanishing zero-order cohomology groups are regular points of the moment map]
\label{lem:Pairs_trivial_stabilizer_group_regular_points_moment_map}
Continue the hypotheses of Lemma \ref{lem:Regular_value_moment_map_affine_space_pairs}. Then a point $(A,\varphi) \in \bmu^{-1}(0)\cap(\sA(E,h)\times W^{1,p}(E))$ is a regular point of the moment map $\bmu$ in \eqref{eq:Moment_map_action_unitary_det_one_gauge_transformations_affine_space_pairs} if and only if $\bH_{A,\varphi}^0 = (0)$, where the harmonic space $\bH_{A,\varphi}^0$ is as in \eqref{eq:H_APhi^0}, and, in that situation, an open neighborhood of $(A,\varphi)$ in $\bmu^{-1}(0)\cap(\sA(E,h)\times W^{1,p}(E))$ is an embedded, real analytic submanifold of the Banach affine space $\sA(E,h)\times W^{1,p}(E)$.  
\end{lem}

\begin{proof}
We recall from the proof of Proposition \ref{prop:Donaldson_Kronheimer_6-5-8_almost_Kaehler_type1_pairs} that
\begin{multline*}
  \langle \xi, d\bmu(A,\varphi)(a,\sigma) \rangle = - (\xi,\Lambda d_Aa)_{L^2(X)} - \Real(\xi\varphi,i\sigma)_{L^2(X)},
  \\
  \text{for all } \xi \in W^{2,p}(\su(E)) \text{ and } (a,\sigma) \in W^{1,p}(T^*X\otimes \su(E) \oplus E).
\end{multline*}
By applying integration by parts as in the proof of Proposition \ref{prop:Donaldson_Kronheimer_6-5-8_almost_Kaehler_type1_pairs}, we see that
\begin{multline*}
  (\xi,\Lambda d_Aa)_{L^2(X)}
  = \int_X \tr_E(\xi\Lambda d_Aa)\wedge\omega^n
  = n\int_X \tr_E(\xi d_Aa)\wedge\omega^{n-1}
  \\
  = -n\int_X \tr_E(d_A\xi \wedge a)\wedge\omega^{n-1}
  = -(d_A\xi,a)_{L^2(X)}.
\end{multline*}
Consequently, we obtain
\begin{multline*}
  \langle \xi, d\bmu(A,\varphi)(a,\sigma) \rangle = (d_A\xi,a)_{L^2(X)} - \Real(\xi\varphi,i\sigma)_{L^2(X)},
  \\
  \text{for all } \xi \in W^{2,p}(\su(E)) \text{ and } (a,\sigma) \in W^{1,p}(T^*X\otimes \su(E) \oplus E).
\end{multline*}
Suppose now that $\xi \in W^{2,p}(\su(E))$ obeys $\langle \xi, d\bmu(A,\varphi)(a,\sigma) \rangle = 0$, for all $(a,\sigma) \in W^{1,p}(T^*X\otimes \su(E) \oplus E)$. The preceding identity implies that this condition on $\xi$ is equivalent to
\[
  d_{A,\varphi}^0\xi = (d_A\xi,-\xi\varphi) = 0 \in W^{1,p}(T^*X\otimes \su(E) \oplus E),
\]
that is, $\xi \in \bH_{A,\varphi}^0$. Hence, $(A,\varphi) \in \bmu^{-1}(0)\cap(\sA(E,h)\times W^{1,p}(E))$ is a regular point of $\bmu$ if and only if $\bH_{A,\varphi}^0 = (0)$.
\end{proof}

\begin{rmk}[Alternative proof that projective vortices with vanishing zero-order cohomology groups are regular points of the moment map]
\label{rmk:Proof_projective_vortices_vanishing_zero-order_cohomology_regular_points_moment_map} 
Although we shall gave a separate, direct proof above the conclusion of Lemma \ref{lem:Pairs_trivial_stabilizer_group_regular_points_moment_map} for a projective vortex $(A,\varphi)$ can also be inferred from Proposition \ref{prop:Itoh_1985_proposition_2-3_SO3_monopole_complex_Kaehler} and Remark \ref{rmk:Relation_H_APhi^2_and_H_dbar_APhi^2}, when $X$ has complex dimension two, and from Theorem \ref{thm:Kobayashi_7-2-21_pairs}, when $X$ has arbitrary complex dimension. Those results establish the canonical isomorphism $\bH_{A,\varphi}^2 \cong \bH_{A,\varphi}^0\oplus\bH_{\bar\partial_A,\varphi}^2$, where $\bH_{A,\varphi}^0 = \Ker d_{A,\varphi}^0\cap\Omega^0(\su(E))$ and $\bH_{\bar\partial_A,\varphi}^2 \subset\Omega^{0,2}(\fsl(E))\oplus\Omega^{0,1}(E)$. When $X$ has dimension two, the harmonic space $\bH_{A,\varphi}^2$ is the $L^2$-orthogonal complement in
\[
  \Omega^0(\su(E))\oplus\Omega^{0,2}(\fsl(E))\oplus\Omega^{0,1}(E)
\]
of the range of the linearization $d\sS(A,\varphi) = d_{A,\varphi}^1$ in the definition \eqref{eq:Projective_vortex_map} of the real analytic, $W^{2,p}(\SU(E))$-equivariant projective vortex map, namely
\begin{multline*}
  \sS(A,\varphi)
  =
  \left(\Lambda (F_A)_0 - \frac{i}{2}(\varphi\otimes\varphi^*)_0, (F_A^{0,2})_0, \bar\partial_A\varphi\right),
  \\
  \text{for all } (A,\varphi) \in \sA(E,h)\times W^{1,p}(E).
\end{multline*}
In higher dimensions, $\bH_{A,\varphi}^2$ is the $L^2$-orthogonal complement in
\[
  \Ker d_{A,\varphi}^2\cap(\Omega^0(\su(E))\oplus\Omega^{0,2}(\fsl(E))\oplus\Omega^{0,1}(E))
\]
of the range of $d_{A,\varphi}^1$. By the definitions \eqref{eq:Projective_vortex_map} of $\sS$ and  \eqref{eq:Co-moment_map_action_unitary_det_one_gauge_transformations_affine_space_pairs} of $\bmu^*$
we see that
\[
  \bmu^*(A,\varphi) = \sS_0(A,\varphi), \quad\text{for all } (A,\varphi) \in \sA(E,h)\times\Omega^0(E),
\]
where we write $\sS = (\sS_0,\sS_{0,2},\sS_{0,1})$, with the components of $\sS$ being defined by the projections of $\sS$ onto the summands $\Omega^0(\su(E))$, $\Omega^{0,2}(\fsl(E))$, and $\Omega^{0,1}(E)$, respectively, of the codomain of $\sS$. Hence,
\[
  d\bmu^*(A,\varphi) = d\sS_0(A,\varphi), \quad\text{for all } (A,\varphi) \in \sA(E,h)\times\Omega^0(E).
\]
The differential $d\bmu(A,\varphi)$ is surjective if and only if $d\bmu^*(A,\varphi)$ is surjective. The preceding identity implies that the latter condition is equivalent to $d\sS_0(A,\varphi)$ being surjective and that is in turn is equivalent to $\bH_{A,\varphi}^0 = (0)$, since $\bH_{A,\varphi}^0$ is the $\Omega^0(\su(E))$-component of $\bH_{A,\varphi}^2$ by definition of $\sS$.
\end{rmk}

\subsection[Kuranishi model for a holomorphic pair]{Kuranishi model for an open neighborhood of a holomorphic pair}
\label{subsec:Kuranishi_model_holomorphic_pair}
We consider the subset of the Banach affine space $\sA^{0,1}(E)\times W^{1,p}(E)$ of $(0,1)$-pairs $(\bar\partial_E,\varphi)$ of class $W^{1,p}$ on $E$ that induce a fixed smooth, holomorphic structure $\bar\partial_{E_d}$ on the complex line bundle $\det E$, as in \eqref{eq:Holomorphic_structure_fixed_determinant}, and obey $F_{\bar\partial_E} = 0$ and $\bar\partial_E\varphi = 0$, for $F_{\bar\partial_E} = \bar\partial_E\circ\bar\partial_E$ as in \eqref{eq:Donaldson_Kronheimer_2-1-50}, and thus define the following analogue of the set $\cH(E)$ in \eqref{eq:Subset_01_connections_defining_holomorphic_structure_on_E}, but now for pairs:
\begin{equation}
  \label{eq:Subset_01-pairs_defining_holomorphic_pair_on_E}
  \cP(E) := \left\{(\bar\partial_E,\varphi) \in \sA^{0,1}(E)\times W^{1,p}(E): F_{\bar\partial_E} = 0 \text{ and } \bar\partial_E\varphi = 0 \right\}.
\end{equation}
According to the Bianchi Identity \eqref{eq:Holomorphic_pair_Bianchi_identity} for $(0,1)$-pairs, we have
\[
  \bar\partial_{E,\varphi}^2(F_{\bar\partial_E},\bar\partial_E\varphi) = (0,0)
  \in \Omega^{0,3}(\fsl(E)) \times \Omega^{0,2}(E),
\]
where the differential $\bar\partial_{E,\varphi}^2$ is as in \eqref{eq:dkStablePair}. Recall from \eqref{eq:d1StablePair} that
\[
  \bar\partial_{E,\varphi}^1(\alpha,\phi) = (\bar\partial_E\alpha, \bar\partial_E\phi + \alpha\varphi)
\]
and that if $(\alpha,\phi) \in \bH_{\bar\partial_E,\varphi}^1$, then $\bar\partial_{E,\varphi}^1(\alpha,\phi) = 0$. As a consequence of the identity \eqref{eq:(0,1)Connection_Induced_On_Endomorphism_Bundle}, we have
\[
  \bar\partial_E(\alpha\phi) = (\bar\partial_{\End(E)}\alpha)\phi - \alpha\wedge\bar\partial_E\phi,
\]
and which we abbreviate by writing $\bar\partial_E(\alpha\phi) = (\bar\partial_E\alpha)\phi - \alpha\wedge\bar\partial_E\phi$. We proceed by analogy with the calculation leading to \eqref{eq:dbar_E_F_dbar_E+alpha}:
\begin{align*}
  &\bar\partial_{E,\varphi}^2\left(F_{\bar\partial_E+\alpha},(\bar\partial_E+\alpha)(\varphi+\phi)\right)
  \\
  &\quad = \bar\partial_{E,\varphi}^2\left(F_{\bar\partial_E} + \bar\partial_E\alpha + \alpha\wedge\alpha, \bar\partial_E\varphi + \bar\partial_E\phi + \alpha(\varphi+\phi)\right)
  \\
  &\quad = \bar\partial_{E,\varphi}^2\left(\bar\partial_E\alpha + \alpha\wedge\alpha, \bar\partial_E\phi + \alpha(\varphi+\phi)\right) \quad\text{(by \eqref{eq:Holomorphic_pair_Bianchi_identity})}
  \\
  &\quad = \left(\bar\partial_E\bar\partial_E\alpha + \bar\partial_E(\alpha\wedge\alpha), \bar\partial_E\bar\partial_E\phi + \bar\partial_E(\alpha(\varphi+\phi))
    - (\bar\partial_E\alpha + \alpha\wedge\alpha)\varphi \right)
  \\
  &\quad\qquad\text{(by \eqref{eq:dkStablePair} with $k=2$)}
  \\
  &\quad = \left(F_{\bar\partial_E}\alpha - 2\alpha\wedge\bar\partial_E\alpha, F_{\bar\partial_E}\phi + \bar\partial_E(\alpha(\varphi+\phi)) - (\bar\partial_E\alpha + \alpha\wedge\alpha)\varphi \right)
  \\
  &\quad = \left(F_{\bar\partial_E}\alpha - 2\alpha\wedge\bar\partial_E\alpha, F_{\bar\partial_E}\phi + (\bar\partial_E\alpha)(\varphi+\phi) - \alpha\wedge\bar\partial_E(\varphi+\phi) - (\bar\partial_E\alpha)\varphi
    - \alpha\wedge\alpha\varphi\right)
  \\
  &\quad = \left(F_{\bar\partial_E}\alpha - 2\alpha\wedge\bar\partial_E\alpha, F_{\bar\partial_E}\phi + (\bar\partial_E\alpha)\phi - \alpha\wedge\bar\partial_E\varphi - \alpha\wedge\bar\partial_E\phi - \alpha\wedge\alpha\varphi\right)
  \\
  &\quad = \left(F_{\bar\partial_E}\alpha, F_{\bar\partial_E}\phi + (\bar\partial_E\alpha)\phi
    - \alpha\wedge\bar\partial_E\varphi\right)
    - \left(2\alpha\wedge\bar\partial_E\alpha, \alpha\wedge\bar\partial_E\phi
    + \alpha\wedge\alpha\varphi\right)    
  \\
  &\quad = \left(F_{\bar\partial_E}\alpha, F_{\bar\partial_E}\phi + (\bar\partial_E\alpha)\phi
    - \alpha\wedge\bar\partial_E\varphi\right)
    - \alpha\wedge\left(\bar\partial_E\alpha, \bar\partial_E\phi + \alpha\varphi\right)
    - \left(\alpha\wedge\bar\partial_E\alpha, 0\right).  
\end{align*}
Therefore, by substituting the expression \eqref{eq:d1StablePair} for $\bar\partial_{E,\varphi}^1(\alpha,\phi)$, we obtain
\begin{multline}
  \label{eq:dbar_Evarphi2_F_dbar_E+alpha_dbar_E+alpha(varphi+phi)}
  \bar\partial_{E,\varphi}^2\left(F_{\bar\partial_E+\alpha},(\bar\partial_E+\alpha)(\varphi+\phi)\right)
  \\
  = \left(F_{\bar\partial_E}\alpha, F_{\bar\partial_E}\phi + (\bar\partial_E\alpha)\phi
    - \alpha\wedge\bar\partial_E\varphi\right)
    - \alpha\wedge\bar\partial_{E,\varphi}^1(\alpha,\phi)
    - \left(\alpha\wedge\bar\partial_E\alpha, 0\right). 
\end{multline}
Hence, for $(\bar\partial_E,\varphi) \in \cP(E)$, so $F_{\bar\partial_E} = 0$ and $\bar\partial_E\varphi = 0$ by definition \eqref{eq:Subset_01-pairs_defining_holomorphic_pair_on_E} of $\cP(E)$, and for $(\alpha,\phi)$ obeying $\bar\partial_{E,\varphi}^1(\alpha,\phi) = 0$, which also gives $\bar\partial_E\alpha = 0$, we obtain the following analogue of \eqref{eq:F_dbar_E+alpha_Ker_dbar_E}:
\begin{multline}
  \label{eq:F_dbar_E+alpha_and_dbar_E_varphi_in_Ker_dbar_Evarphi}
  \left(F_{\bar\partial_E+\alpha}, \bar\partial_{E+\alpha}(\varphi+\phi)\right)
  \in \Ker\bar\partial_{E,\varphi}^1\cap L^p\left(\Lambda^{0,2}(\fsl(E)) \oplus \Lambda^{0,1}(E) \right),
  \\
  \text{for all } (\alpha,\phi) \in \Ker\bar\partial_{E,\varphi}^1\cap W^{1,p}(\Lambda^{0,1}(\fsl(E))\oplus E).
\end{multline}
Given the important inclusion \eqref{eq:F_dbar_E+alpha_and_dbar_E_varphi_in_Ker_dbar_Evarphi}, the proof of Theorem \ref{thm:Local_Kuranishi_model_for_simple_point_cH(E)} yields \mutatis the following analogue for holomorphic pairs; for an alternative approach to the construction of the Kuranishi model provided by Theorem \ref{thm:Local_Kuranishi_model_for_strongly_simple_point_cP(E)}, we refer the reader to Section \ref{subsec:Kuranishi_model_moduli_space_holomorphic_pairs}.

\begin{thm}[Local Kuranishi models for points in $\cP(E)$]
\label{thm:Local_Kuranishi_model_for_strongly_simple_point_cP(E)}
Let $E$ be Hermitian vector bundle over a complex, Hermitian manifold $X$ and $\bar\partial_{E_d}$ be a fixed smooth, holomorphic structure on the complex line bundle $\det E$. If $p\in(n,\infty)$, where $n$ is the complex dimension of $X$, and $(\bar\partial_E,\varphi) \in \cP(E)$ as in \eqref{eq:Subset_01-pairs_defining_holomorphic_pair_on_E} and
\begin{equation}
  \label{eq:dbar_Evarphi_slice}
  S_{\bar\partial_E,\varphi} := (\bar\partial_E,\varphi) + \Ker\bar\partial_{E,\varphi}^{0,*}\cap W^{1,p}\left(\Lambda^{0,1}(\fsl(E))\oplus E\right),
\end{equation}
then there are open neighborhoods $U_{\bar\partial_E,\varphi} \subset S_{\bar\partial_E,\varphi}$ of $(\bar\partial_E,\varphi)$ and $N_{\bar\partial_E,\varphi} \subset \bH_{\bar\partial_E,\varphi}^1$ of the origin such that the following hold:
\begin{enumerate}
\item\label{item:Kuranishi_model_dbarEvarphi_cP(E)}
There are a holomorphic embedding,
\begin{equation}
  \label{eq:Kuranishi_embedding_map_holomorphic_pairs_trivial_stabilizer}
  \bgamma:\bH_{\bar\partial_E,\varphi}^1 \supset N_{\bar\partial_E,\varphi}
  \to U_{\bar\partial_E,\varphi} \subset S_{\bar\partial_E,\varphi},
\end{equation}
such that $\bgamma(0,0) = (\bar\partial_E,\varphi)$ and a holomorphic map,
\begin{equation}
  \label{eq:Kuranishi_obstruction_map_holomorphic_pairs_trivial_stabilizer}
  \bchi:\bH_{\bar\partial_E,\varphi}^1 \supset N_{\bar\partial_E,\varphi}
  \to \bH_{\bar\partial_E,\varphi}^2,
\end{equation}
such that $\bchi(0,0)=0$ and
\begin{equation}
  \label{eq:Kuranishi_model_dbarEvarphi_cP(E)}
  \cP(E) \cap U_{\bar\partial_E,\varphi} = \bgamma(\bchi^{-1}(0)\cap N_{\bar\partial_E,\varphi}).
\end{equation}
The Taylor expansion of $\bchi(\tau,\sigma)$ has $D\bchi(0,0) = 0$ and second-order term comparable\footnote{In the sense that the coefficients of the quadratic form depend at most on the pair $(\bar\partial_E,\varphi)$ and the geometry of the complex Hermitian manifold $X$.} to $\Pi_{\bar\partial_E,\varphi}(\tau\wedge\tau, \tau\sigma)$, where $\Pi_{\bar\partial_E,\varphi}$ is $L^2$-orthogonal projection from $L^p(\Lambda^{0,2}(\fsl(E))\oplus\Lambda^{0,1}(E))$ onto $\bH_{\bar\partial_E,\varphi}^2$.

\item\label{item:cP_dbarEvarphi^vir_is_bgamma(N_dbarEvarphi)}
For $L^2$-orthogonal projection $\Pi_{\Ran\bar\partial_{E,\varphi}^1}$ onto the range of $\bar\partial_{E,\varphi}^1$ in $L^p(\Lambda^{0,2}(\fsl(E))\oplus\Lambda^{0,1}(E))$ and
\begin{multline}
  \label{eq:Atiyah_Hitchin_Singer_family_page_446_holomorphic_pairs}
  \cP_{\bar\partial_E,\varphi}^\vir(E)
  :=
  \left\{(\bar\partial_E,\varphi) + (\alpha,\phi) \in \sA^{0,1}(E) \times W^{1,p}(E): \right.
  \\
  \left.\Pi_{\Ran\bar\partial_{E,\varphi}^1}\left(F_{\bar\partial_E + \alpha}, (\bar\partial_E + \alpha)(\varphi + \phi)\right) = 0\right\},
\end{multline}
the subset $\cP_{\bar\partial_E,\varphi}^\vir(E) \cap U_{\bar\partial_E,\varphi}$ is an embedded complex submanifold of $U_{\bar\partial_E,\varphi}$ with tangent space $T_{\bar\partial_E,\varphi}(\cP_{\bar\partial_E,\varphi}^\vir \cap U_{\bar\partial_E,\varphi}) = \bH_{\bar\partial_E,\varphi}^1$ and
\begin{equation}
  \label{eq:cP_dbarEvarphi^vir_is_bgamma(N_dbarEvarphi)}
  \cP_{\bar\partial_E,\varphi}^\vir \cap U_{\bar\partial_E,\varphi} = \bgamma(N_{\bar\partial_E,\varphi}). 
\end{equation}

\item\label{item:cP(E)_cap_UU_dbarEvarphi_biholomorphic_cP(E)_cap_U_dbarEvarphi_times_UidE}
If $(\bar\partial_E,\varphi)$ is a strongly simple point as in Definition \ref{defn:Simple_pair}, so $\Stab(\bar\partial_E,\varphi) = \{\id_E\}$, then there is an open neighborhood $U_{\id_E} \subset W^{2,p}(\SL(E))$ of the identity $\id_E$ such that the natural map,
\begin{equation}
  \label{eq:dbar_Evarphi_slice_biholomorphic_map}
  U_{\bar\partial_E,\varphi} \times U_{\id_E} \ni \left(\bar\partial_E + \alpha, \varphi + \phi, v\right)
  \mapsto v^*(\bar\partial_E + \alpha, \varphi + \phi) \in \UU_{\bar\partial_E,\varphi},
\end{equation} 
is a biholomorphic map onto an open neighborhood $\UU_{\bar\partial_E,\varphi} \subset \sA^{0,1}(E)\times W^{1,p}(E)$ of $(\bar\partial_E,\varphi)$ and restricts to a biholomorphic map,
\[
  \left(\cP(E) \cap U_{\bar\partial_E,\varphi}\right) \times U_{\id_E} \cong \cP(E) \cap \UU_{\bar\partial_E,\varphi}.
\]  
\end{enumerate}
\end{thm}

We define the following subspaces,
\begin{subequations}
  \label{eq:Subspaces_strongly_simple_or_non-zero-section_or_regular_holomorphic_pairs}
  \begin{align}
    \label{eq:cP**(E)}
    \cP^{**}(E)
    &:=
      \left\{(\bar\partial_E,\varphi) \in \cP(E):
      (\bar\partial_E,\varphi) \text{ is strongly simple}\right\},
    \\
    \label{eq:cP0(E)}
    \cP^0(E)
    &:=
      \left\{(\bar\partial_E,\varphi) \in \cP(E): \varphi \not\equiv 0\right\},
    \\
    \label{eq:cPreg(E)}  
    \cP_\reg(E)
    &:=
    \left\{(\bar\partial_E,\varphi) \in \cP(E): 
    \bH_{\bar\partial_E,\varphi}^2 = (0) \right\}.  
  \end{align}
\end{subequations}
where the harmonic space $\bH_{\bar\partial_E,\varphi}^2$ is defined in \eqref{eq:H_dbar_Avarphi^0bullet}. The set $\cP_\reg(E)$ in \eqref{eq:cPreg(E)} is an open subspace of $\cP(E)$ as a consequence of the Implicit Mapping Theorem and the subspace $\cP^0(E)$ in \eqref{eq:cP0(E)} is an open subspace of $\cP(E)$ because it is the complement of the closed subspace of zero section holomorphic pairs. Assume now that $(X,\omega)$ is a complex K\"ahler manifold and define  
\begin{subequations}
  \label{eq:Subspaces_stable_polystable_holomorphic_pairs}
  \begin{align}
    \label{eq:cP(E,omega)}
    \cP(E,\omega)
    &:=
      \left\{(\bar\partial_E,\varphi) \in \cP(E):
      (\bar\partial_E,\varphi) \text{ is stable}\right\},
    \\
    \label{eq:cPps(E,omega)}
    \cP_\ps(E,\omega)
    &:=
      \left\{(\bar\partial_E,\varphi) \in \cP(E):
      (\bar\partial_E,\varphi) \text{ is polystable}\right\}.
\end{align}
\end{subequations}
We now set
\begin{subequations}
  \label{eq:Subspaces_strongly_simple_or_non-zero-section_stable_or_polystable_holomorphic_pairs}
  \begin{align}
    \label{eq:cP0reg(E,omega)}
    \cP_\reg^0(E,\omega) &:= \cP^0(E) \cap \cP(E,\omega) \cap \cP_\reg(E),
    \\
    \label{eq:cP0psreg(E,omega)}
    \cP_{\ps,\reg}^0(E,\omega) &:= \cP^0(E) \cap \cP_\ps(E,\omega) \cap \cP_\reg(E).
  \end{align}
\end{subequations}
According to Theorem \ref{thm:Lubke_Teleman_6-3-7} \eqref{item:Lubke_Teleman_6-3-7_analytic_embedding_non-split-nonzero-section}, the moduli subspace $\fM^0(E,\omega)$ in \eqref{eq:Moduli_space_stable_holomorphic_pairs} is open in $\fM(E)$ and because
\[
  \fM^0(E,\omega) = \left.\cP^0(E,\omega)\right/W^{2,p}(\SL(E)),
\]
we see that $\cP^0(E,\omega)$ is an open subspace of $\cP(E)$. By Lemma \ref{lem:Openness_subspace_strongly_simple_01_pairs} \eqref{item:Openness_subspace_strongly_simple_stable_pairs}, we have $\fM^0(E,\omega) \subset \fM^{**}(E)$ and so $\cP^0(E,\omega) \subset \cP^{**}(E)$. Therefore, the subspace $\cP_\reg^0(E,\omega)$ is also open in $\cP(E)$ and contained in $\cP^{**}(E)$. As a consequence of Theorem \ref{thm:Local_Kuranishi_model_for_strongly_simple_point_cP(E)} \eqref{item:cP(E)_cap_UU_dbarEvarphi_biholomorphic_cP(E)_cap_U_dbarEvarphi_times_UidE}, the subspace
$\cP_\reg^0(E,\omega)$ is an embedded complex submanifold of $\sA^{0,1}(E)\times W^{1,p}(E)$.

If $E$ has complex rank two, Theorem \ref{thm:Lubke_Teleman_6-3-7} \eqref{item:Lubke_Teleman_6-3-7_analytic_embedding_rank-2_non-zero_section} implies that the moduli subspace $\fM_\ps^0(E,\omega)$ in \eqref{eq:Moduli_space_regular_non-zero-section_polystable_holomorphic_pairs} is open in $\fM(E)$ and because
\[
  \fM_\ps^0(E,\omega) = \left.\cP_\ps^0(E,\omega)\right/W^{2,p}(\SL(E)),
\]
we see that $\cP_\ps^0(E,\omega)$ is an open subspace of $\cP(E)$. According to Lemma \ref{lem:Openness_subspace_strongly_simple_01_pairs} \eqref{item:Openness_subspace_strongly_simple_polystable_rank-2_pairs}, we have $\fM_\ps^0(E,\omega) \subset \fM^{**}(E)$ and so $\cP_\ps^0(E,\omega) \subset \cP^{**}(E)$. Therefore, the subspace $\cP_{\ps,\reg}(E,\omega)$ is also open in $\cP(E)$ and contained in $\cP^{**}(E)$. As a consequence of Theorem \ref{thm:Local_Kuranishi_model_for_strongly_simple_point_cP(E)} \eqref{item:cP(E)_cap_UU_dbarEvarphi_biholomorphic_cP(E)_cap_U_dbarEvarphi_times_UidE}, the subspace
$\cP_{\ps,\reg}^0(E,\omega)$ is an embedded complex submanifold of $\sA^{0,1}(E)\times W^{1,p}(E)$.

\subsection[Symplectic structure on the moduli space of regular projective vortices]{Marsden--Weinstein reduction and symplectic structure for the moduli space of regular projective vortices}
\label{subsec:Marsden-Weinstein_reduction_moduli_space_SO3_monopoles_symplectic_quotient}
Suppose now that $(E,h)$ is a Hermitian vector bundle over a complex K\"ahler manifold $(X,g,J)$ of dimension $n$ with K\"ahler form $\omega = g(\cdot,J\cdot)$ and that the fixed unitary connection $A_d$ on the Hermitian line bundle $\det E$ has curvature $F_{A_d}$ of type $(1,1)$, so $F_{A_d}^{0,2}=0$.

In preparation for our forthcoming application of the Marsden--Weinstein Symplectic Reduction Theorem \ref{thm:Marsden-Weinstein_symplectic_quotient}, we shall choose (by analogy with our previous choice of $\bV''$ in \eqref{eq:Marsden-Weinstein_V_space_01connections})
\begin{equation}
  \label{eq:Marsden-Weinstein_V_space_holomorphic_pairs}
  \bV''
  :=
  \begin{cases}
    \cP_\reg^0(E,\omega), &\text{if } \rank_\CC E \geq 2,
    \\
    \cP_{\ps,\reg}^0(E,\omega), &\text{if } \rank_\CC E = 2,
  \end{cases}  
\end{equation}
in Hypothesis \ref{hyp:Marsden-Weinstein_symplectic_quotient_conditions}, where $\cP_\reg^0(E,\omega)$ is as in \eqref{eq:cP0reg(E,omega)} and $\cP_{\ps,\reg}^0(E,\omega)$ is as in \eqref{eq:cP0psreg(E,omega)}. By our discussion in Section \ref{subsec:Kuranishi_model_holomorphic_pair} of the properties of the subspaces $\cP_\reg^0(E,\omega)$ and $\cP_{\ps,\reg}^0(E,\omega)$, we note that $\bV''$  is an embedded complex submanifold of $\sA^{0,1}(E)\times W^{1,p}(E)$ and is contained in $\cP^{**}(E)$ in \eqref{eq:cP**(E)}:
\begin{subequations}
  \begin{gather}
    \label{eq:bV''_embedded_complex_submanifold_affine_space}
    \bV'' \hookrightarrow \sA^{0,1}(E)\times W^{1,p}(E),
    \\
    \label{eq:bV''_in_cP**(E)}
    \bV'' \subset \cP^{**}(E).
  \end{gather}
\end{subequations}
For $\bV''$ as in \eqref{eq:Marsden-Weinstein_V_space_holomorphic_pairs}, we define
\begin{equation}
  \label{eq:Marsden-Weinstein_V_space_unitary_pairs}
  \bV := \left(\pi_h^{0,1}\right)^{-1}(\bV''),
\end{equation}
so that $\bV$ is an embedded submanifold of the real Banach affine space $\sA(E,h)\times W^{1,p}(E)$ as well via the restriction of the canonical isomorphism of real Banach affine spaces (extending the canonical isomorphism $\pi_h^{0,1}:\sA(E,h) \to \sA^{0,1}(E)$ of real Banach affine spaces in \eqref{eq:Bijection_unitaryconnections_with_01connections}),
\begin{equation}
  \label{eq:Bijection_unitarypairs_with_01pairs}
  \pi_h^{0,1}:\sA(E,h) \times W^{1,p}(E) \ni (A,\varphi) \mapsto (\bar\partial_A,\varphi) \in \sA^{0,1}(E) \times W^{1,p}(E),
\end{equation}
to a bijection of real analytic subspaces,
\begin{equation}
  \label{eq:Bijection_unitary1-1pairs_with_integrable01pairs}
  \pi_h^{0,1}:\sP(E,h) \ni (A,\varphi) \mapsto (\bar\partial_A,\varphi) \in \cP(E),
\end{equation}
where $\sP(E,h)$ is as in \eqref{eq:Subset_pre-holomorphic_pairs} and $\cP(E)$ is as in \eqref{eq:Subset_01-pairs_defining_holomorphic_pair_on_E}.

The K\"ahler metric $\bh$ on $\sA^{0,1}(E)\times W^{1,p}(E)$ in \eqref{eq:Kobayashi_7-6-20_pairs} is related to the almost complex structure $\bJ$ in \eqref{eq:Almost_complex_structure_affine_space_pairs_unitary_connections_and_sections}, (weak or $L^2$) Riemannian metric $\bg$ in \eqref{eq:L2_metric_affine_space_pairs_unitary_connections_and_sections}, and symplectic form $\bomega$ in \eqref{eq:Kobayashi_7-6-22_pairs} on $\sA(E,h)\times W^{1,p}(E)$ in the usual way. The K\"ahler metric $\bh$ on $\sA^{0,1}(E)\times W^{1,p}(E)$ induces a K\"ahler metric $\bh|_{\bV''} := \bh\restriction\bV''$ on the embedded complex submanifold $\bV''$, which we denote simply by $\bh$. We shall write the K\"ahler form $\bomega|_{\bV''} := \bomega\restriction\bV''$ induced by $\bomega$ simply as $\bomega$. Similarly, we may regard $\bV$ as a symplectic manifold and write $\bomega$ for its symplectic form instead of $(\pi_h^{0,1})^*(\bomega|_{\bV''})$. 

Since the moment map $\bmu$ on $\bV$ is the restriction of the moment map $\bmu$ in \eqref{eq:Moment_map_action_unitary_det_one_gauge_transformations_affine_space_pairs} on the Banach affine space $\sA(E,h)\times W^{1,p}(E)$, we shall need the following refinement of Lemma \ref{lem:Regular_value_moment_map_affine_space_pairs}.

\begin{lem}[Regular points of $\bmu$ remain regular upon restriction from $\sA(E,h)\times W^{1,p}(E)$ to $\sP(E,h)$]
\label{lem:Regular_points_bmu_remain_regular_on_restriction_from_sA(E,h)_times_W1p(E)_to_sP(E,h)}
Continue the hypotheses of Lemma \ref{lem:Regular_value_moment_map_affine_space_pairs}. If $(A,\varphi) \in \sP(E,h)$ as in \eqref{eq:Subset_pre-holomorphic_pairs} is a regular point of the moment map $\bmu$ in \eqref{eq:Moment_map_action_unitary_det_one_gauge_transformations_affine_space_pairs},
\[
  \bmu:\sA(E,h)\times W^{1,p}(E) \to \left(W^{2,p}(\su(E))\right)^*,
\]  
then $(A,\varphi)$ is also a regular point of the restriction of the moment map $\bmu$ to the Banach analytic subvariety $\sP(E,h)$,
\begin{equation}
  \label{eq:Moment_map_action_unitary_det_one_gauge_transformations_sP(E,h)}
  \bmu:\sP(E,h) \to \left(W^{2,p}(\su(E))\right)^*,
\end{equation}
in the sense that the differential $d\bmu(A,\varphi)$ is an epimorphism from the Zariski tangent space $T_{A,\varphi}\sP(E,h)$ onto $(W^{2,p}(\su(E)))^*$.
\end{lem}

\begin{proof}
From the definition \eqref{eq:d1_projective_vortex_elliptic_deformation_complex} of the differential $d_{A,\varphi}^1$, we may write
\[
  d_{A,\varphi}^1 = d_{A,\varphi}^{1;0} + d_{A,\varphi}^{1;1}:
  W^{1,p}\left(T^*X\otimes\su(E)\oplus E\right)
  \to
  L^p(\su(E)) \oplus L^p\left(\Lambda^{0,2}(\fsl(E)) \oplus \Lambda^{0,1}(E)\right),
\]
where $d_{A,\varphi}^{1;0}$ and $d_{A,\varphi}^{1;1}$ correspond to the $\Omega^0(\su(E))$ and $\Omega^{0,2}(\fsl(E))\oplus \Omega^{0,1}(E)$ components of $d_{A,\varphi}^1$, respectively:
\begin{subequations}
  \label{eq:dAvarphi1;0_and_1;1}
  \begin{align}
    \label{eq:dAvarphi1;0}
    d_{A,\varphi}^{1;0}: W^{1,p}\left(T^*X\otimes\su(E)\oplus E\right)
    &\to L^p(\su(E)),
    \\
     \label{eq:dAvarphi1;1}
    d_{A,\varphi}^{1;1}: W^{1,p}\left(T^*X\otimes\su(E)\oplus E\right)
    &\to L^p\left(\Lambda^{0,2}(\fsl(E)) \oplus \Lambda^{0,1}(E)\right).
  \end{align}
\end{subequations}
By definition \eqref{eq:Subset_pre-holomorphic_pairs} of $\sP(E,h)$, we have
\begin{multline*}
  \sP(E,h)
  =
  \left\{(A,\varphi) \in \sA(E,h)\times W^{1,p}(E): \right.
    \\
    \left. F_{\bar\partial_A} = 0 \in 
    L^p(\Lambda^{0,2}(\fsl(E))) \text{ and }
    \bar\partial_A\varphi = 0 \in L^p(\Lambda^{0,1}(E)) \right\},
\end{multline*}
and so its Zariski tangent space at $(A,\varphi)$ is 
\[
  T_{A,\varphi}\sP(E,h)
  =
  \Ker d_{A,\varphi}^{1;1}\cap W^{1,p}\left(T^*X\otimes\su(E)\oplus E\right).
\]
We claim that if the operator $d_{A,\varphi}^{1;0}$ in \eqref{eq:dAvarphi1;0} is surjective, then its restriction to the closed subspace $T_{A,\varphi}\sP(E,h)$,
\[
  d_{A,\varphi}^{1;0}:\Ker d_{A,\varphi}^{1;1}\cap W^{1,p}\left(T^*X\otimes\su(E)\oplus E\right) \to L^p(\su(E))
\]
is also surjective. To see this, suppose more generally that $T:\bE \to \bF$ is a bounded, linear operator between Banach spaces and that $\bF = \bF_0\oplus \bF_1$ as a direct sum of closed linear subspaces. Define $\bE_i := T^{-1}(\bF_i)$ for $i=0,1$, so that $\bE_0 + \bE_1 \subset \bE$. Write $T_i := \pi_i\circ T$, with $\pi_i:\bF\to \bF_i$ denoting the projection operator for $i=0,1$. If $x_0 \in \bE_0$, then $Tx_0 \in \bF_0$ (by definition of $\bE_0$) while $T_1x_0 = \pi_1Tx_0 = 0$ (since $Tx_0\in \bF_0$), so that $x_0 \in \Ker T_1$. We conclude that $\bE_0 \subset \Ker T_1$. In particular, if $T_0:\bE \to \bF_0$ is surjective, then $T_0:\Ker T_1 \to \bF_0$ is also surjective. The claim now follows by applying this observation from linear algebra to the operator $d_{A,\varphi}^1 = d_{A,\varphi}^{1;0} + d_{A,\varphi}^{1;1}$.

We conclude that $d_{A,\varphi}^{1;0}$ is surjective upon restriction to $T_{A,\varphi}\sP(E,h)$. The identity \eqref{eq:dbmu_Avarphi_equals_dAvarphi1;0} therefore implies that the restriction of the differential,
\[
  d\bmu^*(A,\varphi):T_{A,\varphi}\sP(E,h) \to L^p(\su(E))
\]
is surjective. Hence, the restriction of the differential,
\[
  d\bmu(A,\varphi):T_{A,\varphi}\sP(E,h) \to \left(W^{2,p}(\su(E))\right)^*
\]
is also surjective and this completes the proof of the lemma.
\end{proof}

Similarly, because the moment map $\bmu$ on $\bV$ is the restriction of the moment map $\bmu$ on the Banach affine space $\sA(E,h)\times W^{1,p}(E)$, we shall need the following refinement of Lemma \ref{lem:Pairs_trivial_stabilizer_group_regular_points_moment_map} and whose proof is an immediate consequence of the proof of that lemma.

\begin{lem}[Projective vortices with vanishing zero-order cohomology groups remain regular points of the moment map upon restriction from $\sA(E,h)\times W^{1,p}(E)$ to $\sP(E,h)$]
\label{lem:Regular_points_bmu_on_restriction_from_sA(E,h)_times_W1p(E)_to_sP(E,h)}
Continue the hypotheses of Lemma \ref{lem:Regular_value_moment_map_affine_space_pairs}. If $(A,\varphi) \in \bmu^{-1}(0)\cap\sP(E,h)$, then the following hold:
\begin{enumerate}
\item $(A,\varphi)$ is a regular point of the restricted moment map $\bmu$ in \eqref{eq:Moment_map_action_unitary_det_one_gauge_transformations_sP(E,h)} if and only if $\bH_{A,\varphi}^0 = (0)$, where the harmonic space $\bH_{A,\varphi}^0$ is as in \eqref{eq:H_APhi^0}.

\item If $\bH_{A,\varphi}^0 = (0)$, then an open neighborhood of $(A,\varphi)$ in $\bmu^{-1}(0)\cap\bV$ is an embedded, real analytic submanifold of $\bV$ in \eqref{eq:Marsden-Weinstein_V_space_unitary_pairs}.
\end{enumerate}
\end{lem}

We conclude from Lemma \ref{lem:Regular_points_bmu_on_restriction_from_sA(E,h)_times_W1p(E)_to_sP(E,h)} that zero is a regular value of the restriction,
\[
  \bmu: \bV \to \left(W^{2,p}(\su(E))\right)^*,
\]
of the moment map $\bmu$ in \eqref{eq:Moment_map_action_unitary_det_one_gauge_transformations_affine_space_pairs}. Indeed, this holds because if $(A,\varphi) \in \bmu^{-1}(0)\cap\bV$ for $\bV = (\pi_h^{0,1})^{-1}(\bV'')$ as in \eqref{eq:Marsden-Weinstein_V_space_unitary_pairs}, then $(A,\varphi)$ is a projective vortex by the identity \eqref{eq:Kobayashi_7-6-33_pairs}. Hence, Theorem \ref{thm:Kobayashi_7-2-21_pairs} implies that $\bH_{A,\varphi}^0 = (0)$ since $\bH_{\bar\partial_A,\varphi}^0 = \bH_{A,\varphi}^0\otimes_\RR\CC$ and $\bH_{\bar\partial_A,\varphi}^0 = (0)$ by definition \eqref{eq:Marsden-Weinstein_V_space_holomorphic_pairs} of $\bV''$, the fact that $\bH_{\bar\partial_A,\varphi}^0$ is the Lie algebra of $\Stab(\bar\partial_A,\varphi)$ (see Lemma \ref{lem:Stab(rdE,varphi)_is_Lie_Group}), and the fact that $\Stab(\bar\partial_A,\varphi) = \{\id_E\}$ by the inclusion $\bV'' \subset \cP^{**}(E)$ in \eqref{eq:bV''_in_cP**(E)} and definition \eqref{eq:cP**(E)} of $\cP^{**}(E)$.

Therefore, $0 \in \fg^* := (W^{2,p}(\su(E)))^*$ is a weakly regular value of $\bmu$ on $\bV$ in the sense of Condition \eqref{item:Marsden-Weinstein_symplectic_quotient_conditions_weakly_regular_value} in Hypothesis \ref{hyp:Marsden-Weinstein_symplectic_quotient_conditions} for the Marsden--Weinstein Symplectic Reduction Theorem \ref{thm:Marsden-Weinstein_symplectic_quotient}. Thus, $\bmu^{-1}(0)\cap\bV$ is an embedded submanifold of $\bV$ and for every $(A,\varphi) \in \bmu^{-1}(0)\cap\bV$, one has the equality
\[
  T_{A,\varphi}\left(\bmu^{-1}(0)\cap \bV \right) = \Ker d\bmu(A,\varphi).
\]
Thus, Condition \eqref{item:Marsden-Weinstein_symplectic_quotient_conditions_weakly_regular_value} in Hypothesis \ref{hyp:Marsden-Weinstein_symplectic_quotient_conditions} is obeyed with $\fg$ as above and $V$ and $\mu$ replaced by $\bV$ and $\bmu$, respectively.

The group $\sG_E = W^{2,p}(\SU(E))$ acts freely on the open subspace $\sP^{**}(E,h)\subset \sP(E,h)$ given by
\begin{equation}
  \label{eq:sP**(E,h)} 
  \sP^{**}(E,h) := \left\{(A,\varphi) \in \sP(E,h): \Stab(A,\varphi) = \{\id_E\} \right\},
\end{equation}
by analogy with our definition \eqref{eq:Moduli_space_projective_vortices_StabAvarphi_idE} of the open subspace
$\sM^{**}(E,h,\omega) \subset \sM(E,h,\omega)$. Here, openness of $\sP^{**}(E,h)$ in $\sP(E,h)$ is a consequence of the fact that $\sC^{**}(E,h)$ is an open subspace of $\sC(E,h)$ with the quotient topology by Lemma \ref{lem:Openness_configuration_subspace_unitary_pairs_trivial_stabilizer}, so
\[
  \left.\sP^{**}(E,h)\right/W^{2,p}(\SU(E)) = \left(\sP(E,h)/W^{2,p}(\SU(E))\right)\cap \sC^{**}(E,h)
\]
is an open subspace of $\sP(E,h)/W^{2,p}(\SU(E))$ with the quotient topology.

By partial analogy with our definition \eqref{eq:Moduli_space_projective_vortices_non-split-non-zero-section} of $\sM^{*,0}(E,h,\omega)$, we define the following open subspace of $\sP(E,h)$,
\begin{equation}
  \label{eq:sPnatural0(E,h)} 
  \sP^{\natural,0}(E,h) := \left\{(A,\varphi) \in \sP(E,h): \bH_A^0 = (0) \text{ and } \varphi\not\equiv 0 \right\}.
\end{equation}
Here, openness of $\sP^{\natural,0}(E,h)$ in $\sP(E,h)$ follows from the facts that
\begin{inparaenum}[\itshape i\upshape)]
\item $\sP^0(E,h)$ is an open subspace $\sP(E,h)$, where
\[
  \sP^0(E,h) := \left\{(A,\varphi) \in \sP(E,h): \varphi\not\equiv 0 \right\},
\]
and \item $\bH_A^0 = (0)$ is an open condition on the affine space $\sA(E,h)$ of connections $A$ by the proof of Lemma \ref{lem:Openness_moduli_subspace_projective_vortices_non-split_connections},
\item $\sP^\natural(E,h)$ is an open subspace $\sP(E,h)$, and \item $\sP^{\natural,0}(E,h)  = \sP^\natural(E,h) \cap \sP^0(E,h)$.
\end{inparaenum}

Note that if $A$ is a $W^{1,p}$ unitary connection obeying $\bH_A^0 = (0)$, then $A$ is non-split by Corollary \ref{cor:Split_unitary_A_and_Lie_Alg_of_Stab(A)} \eqref{item:A_split_implies_bHA_non-zero}. Consequently, we obtain the inclusion
\begin{equation}
  \label{eq:sPnatural0(E,h,omega)_subset_sP*0(E,h,omega)}
  \sP^{\natural,0}(E,h) \subset \sP^{*,0}(E,h),
\end{equation} 
where
\begin{equation}
  \label{eq:sP*0(E,h)} 
  \sP^{*,0}(E,h) := \left\{(A,\varphi) \in \sP(E,h): A \text{ is non-split and } \varphi\not\equiv 0 \right\}.
\end{equation}
Conversely, Corollary \ref{cor:Split_unitary_A_and_Lie_Alg_of_Stab(A)} \eqref{item:A_smooth_and_bHA_non-zero_implies_A_split} asserts that, when $A$ is smooth, if $A$ is non-split then $\bH_A^0 = (0)$. Since every $W^{1,p}$ pair $(A',\varphi')$ representing a point in $\sM(E,h,\omega)$ is equivalent via a gauge transformation in $W^{2,p}(\SU(E))$ to a smooth pair $(A,\varphi)$, the preceding inclusion restricts to an equality of moduli subspaces,
\[
  \sM^{\natural,0}(E,h,\omega) = \sM^{*,0}(E,h,\omega),
\]  
where $\sM^{*,0}(E,h,\omega)$ is as in \eqref{eq:Moduli_space_projective_vortices_non-split-non-zero-section} and
\begin{equation}
  \label{eq:sMnatural0(E,h)} 
  \sM^{\natural,0}(E,h,\omega) := \left\{[A,\varphi] \in \sM(E,h,\omega): \bH_A^0 = (0) \text{ and } \varphi\not\equiv 0 \right\}.
\end{equation}
One sees that $\sM^{\natural,0}(E,h,\omega)$ is an open subspace of $\sM(E,h,\omega)$ by the fact that $\bH_A^0 = (0)$ is an open condition from the proof of Lemma \ref{lem:Openness_moduli_subspace_projective_vortices_non-split_connections} and the fact that $\varphi\not\equiv 0$ is an open condition.

The proof of the inclusion \eqref{eq:sM*0(E,h,omega)_subset_sM**(E,h,omega)} also verifies that the following inclusion holds:
\[
  \sP^{*,0}(E,h) \subset \sP^{**}(E,h).
\]
Therefore, by combining the preceding inclusion with the inclusion \eqref{eq:sPnatural0(E,h,omega)_subset_sP*0(E,h,omega)} we obtain
\begin{equation}
  \label{eq:sPnatural0(E,h,omega)_subset_sP**(E,h,omega)} 
  \sP^{\natural,0}(E,h) \subset \sP^{**}(E,h),
\end{equation}
Hence, $\sG_E$ acts freely on the open subspace,
\begin{equation}
  \label{eq:bVnatural0}
  \bV^{\natural,0} := \bV \cap \sP^{\natural,0}(E,h).
\end{equation}
of the manifold $\bV$ in \eqref{eq:Marsden-Weinstein_V_space_holomorphic_pairs}, a requirement of Condition \eqref{item:Marsden-Weinstein_symplectic_quotient_conditions_slice} in Hypothesis \ref{hyp:Marsden-Weinstein_symplectic_quotient_conditions} with $G$, $V$, and $\mu$ replaced by $\sG_E$, $\bV^{\natural,0}$, and $\bmu$, respectively. 

Furthermore, the Banach Lie group $\sG_E$ acts properly on $\bmu^{-1}(0)\cap\bV^{\natural,0}$ since the actions of $W^{2,p}(\SU(E))$ on $\sA(E,h)\times W^{1,p}(E)$ or $\sP(E,h)$ are proper by modifying Kobayashi's proof of \cite[Proposition 7.1.14, p. 221]{Kobayashi_differential_geometry_complex_vector_bundles}. (As discussed in Section \ref{sec:Marsden-Weinstein_reduction_moduli_space_HE_connections_symplectic_quotient}, Kobayashi uses the group $W^{2,p}(\U(E))$ rather than $W^{2,p}(\SU(E))$ but this poses no additional difficulties because $W^{2,p}(\SU(E))$ is a subgroup of $W^{2,p}(\U(E))$. We obtain properness of the action of $\sG_E$ on unitary pairs from the properness of the action on $\sA(E,h)$: if $u_n^*(A_n,\varphi_n)$ and $(A_n,\varphi_n)$ converge in $W^{1,p}$, then $u_nA_n$ and $A_n$ converge in $W^{1,p}$ and so a subsequence of $u_n$ converges in $W^{2,p}$.) According to Feehan and Maridakis \cite[Corollary 18, p. 19, and remarks for pairs, p. 20]{Feehan_Maridakis_Lojasiewicz-Simon_coupled_Yang-Mills}, there is a (Coulomb gauge) slice for the action of $\sG_E$ on $\bmu^{-1}(0)\cap\bV^{\natural,0}$ in the sense of Condition \eqref{item:Marsden-Weinstein_symplectic_quotient_conditions_slice}.
Therefore, all the requirements of Condition \eqref{item:Marsden-Weinstein_symplectic_quotient_conditions_slice} in Hypothesis \ref{hyp:Marsden-Weinstein_symplectic_quotient_conditions} are obeyed with $G$, $V$, and $\mu$ replaced by $\sG_E$, $\bV^{\natural,0}$, and $\bmu$, respectively.

The analogue of Proposition \ref{prop:Donaldson_Kronheimer_6-5-8_almost_Kaehler} for $\sA(E,h)\times W^{1,p}(E)$ in place of $\sA(E,h)$ verifies that Condition \eqref{item:Marsden-Weinstein_symplectic_quotient_conditions_equivariance} in Hypothesis \ref{hyp:Marsden-Weinstein_symplectic_quotient_conditions} is obeyed with $G$, $V$, and $\mu$ replaced by $\sG_E$, $\bV^{\natural,0}$, and $\bmu$, respectively.

\begin{rmk}[Open subspaces of the affine space of unitary pairs when $E$ has complex rank two]
\label{rmk:Simplification_E_rank_two_subspaces_affine_space_unitary_pairs}
By analogy with our customary definitions, we set     
\begin{equation}
  \label{eq:sP0}
  \sP^0(E,h) := \left\{(A,\varphi) \in \sP(E,h): \varphi\not\equiv 0 \right\}.
\end{equation}
If $\varphi\equiv 0$, then $\Stab(A,0) \supset C_r\,\id_E$, where $r$ is the rank of $E$ and $C_r$ is the group of $r$-th roots of unity. Conversely, if $\Stab(A,\varphi) = \{\id_E\}$, we must have $\varphi\not\equiv 0$ and thus we have the inclusion
\[
  \sP^{**}(E,h) \subset \sP^0(E,h).
\]
If we restrict to the case where $E$ has complex rank two and $\varphi\not\equiv 0$, then $\Stab(A,\varphi) = \{\id_E\}$ by Lemma \ref{lem:Stabilizers_non-zero-section_unitary_pairs} and so we obtain the reverse inclusion
\[
  \sP^0(E,h) \subset \sP^{**}(E,h).
\]
Hence, when $E$ has complex rank two, we obtain the equality
\[
  \sP^0(E,h) = \sP^{**}(E,h).
\]
Compare our derivation of the equality \eqref{eq:sC0(E,h)_equals_sC**(E,h)_E_rank2}.
\end{rmk}

We now prove the following analogue of Theorem \ref{thm:Kobayashi_7-6-36}.

\begin{thm}[Moduli space of regular projective vortices as a complex K\"ahler manifold]
\label{thm:Kobayashi_7-6-36_pairs}
Let $(E,h)$ be a smooth Hermitian vector bundle over a closed, complex K\"ahler manifold $(X,g,J)$ with K\"ahler form $\omega = g(\cdot,J\cdot)$ and complex dimension $n$. Let $p \in (n,\infty)$ be a constant and $A_d$ be a fixed unitary connection on the Hermitian line bundle $\det E$ with curvature obeying $F_{A_d}^{0,2}=0$. Then the following hold:
\begin{enumerate}
\item\label{item:sMreg*0(E.h,omega)_is_complex_Kaehler}
  The moduli subspace $\sM_\reg^{*,0}(E,h,\omega)$ in \eqref{eq:Moduli_space_projective_vortices_non-split-non-zero-section_regular} is a complex K\"ahler manifold with K\"ahler form induced from $\bomega$ in \eqref{eq:Kobayashi_7-6-22_pairs}.

\item\label{item:sMreg0(E.h,omega)_is_complex_Kaehler_if_E_rank_2}
  If $E$ has rank two, then the moduli subspace $\sM_\reg^0(E,h,\omega)$ in \eqref{eq:Moduli_space_projective_vortices_non-zero_section_regular} is a complex K\"ahler manifold with K\"ahler form induced from $\bomega$ in \eqref{eq:Kobayashi_7-6-22_pairs}.
\end{enumerate}
\end{thm}

\begin{proof}
We adapt the proof of Theorem \ref{thm:Kobayashi_7-6-36}.

Consider Item \eqref{item:sMreg*0(E.h,omega)_is_complex_Kaehler}. According to Corollary \ref{cor:Lubke_Teleman_6-3-7_and_Kobayashi_7_3_17_pair} \eqref{item:fM^0(E,omega)_is_complex_analytic_space}, the moduli subspace $\fM_\reg^0(E,\omega)$ in \eqref{eq:Moduli_space_regular_non-zero-section_stable_holomorphic_pairs} of non-zero section, regular, stable holomorphic pairs is a complex manifold and there is an isomorphism \eqref{eq:Lubke_Teleman_6-3-7_analytic_isomorphism_non-split-nonzero-section_regular} of real analytic manifolds,
\[
  \sM_\reg^{*,0}(E,h,\omega) \cong \fM_\reg^0(E,\omega).
\]
Hence, the moduli space $\sM_\reg^{**}(E,h,\omega)$ inherits a complex structure from $\fM_\reg^0(E,\omega)$ via the preceding isomorphism of real analytic manifolds.

We claim that we have the following identification of sets:
\begin{equation}
  \label{eq:sMreg*0Ehomega_equals_symplectic_quotient}
  \sM_\reg^{*,0}(E,h,\omega) = \left.(\bmu^{-1}(0)\cap \bV^{\natural,0})\right/W^{2,p}(\SU(E)),
\end{equation}
where $\bV^{\natural,0}$ is defined in \eqref{eq:bVnatural0}.

To verify \eqref{eq:sMreg*0Ehomega_equals_symplectic_quotient}, first let $[A,\varphi]\in \sM_\reg^{*,0}(E,h,\omega)$. From the definition of $\sM_\reg^{*,0}(E,h,\omega)$ in \eqref{eq:Moduli_space_projective_vortices_non-split-non-zero-section_regular}, we have $\bH_{A,\varphi}^2 = (0)$. Since $(A,\varphi) \in \sP(E,h)$ in \eqref{eq:Subset_pre-holomorphic_pairs}, then $(\bar\partial_A,\varphi) \in \cP(E)$ in \eqref{eq:Subset_01-pairs_defining_holomorphic_pair_on_E}. Moreover, since $(A,\varphi)$ is a projective vortex and $\bH_{A,\varphi}^2 = (0)$, Theorem \ref{thm:Kobayashi_7-2-21_pairs} implies that $\bH_{\bar\partial_A,\varphi}^2 = (0)$. The isomorphism
\[
  \sM_\reg^{*,0}(E,h,\omega) \ni [A,\varphi] \leftrightarrow [\bar\partial_A,\varphi] \in \fM_\reg^0(E,\omega)
\]
given by \eqref{eq:Lubke_Teleman_6-3-7_analytic_isomorphism_non-split-nonzero-section_regular} implies that $(\bar\partial_A,\varphi) \in \bV''$ by the definition \eqref{eq:Marsden-Weinstein_V_space_holomorphic_pairs} of $\bV'' = \cP_\reg^0(E,\omega)$, and so $(A,\varphi) \in \bV$ by the definition $\bV = (\pi_h^{0,1})^{-1}(\bV'') \subset \sA(E,h)\times W^{1,p}(E)$ in \eqref{eq:Marsden-Weinstein_V_space_unitary_pairs}. Because $A$ is non-split and $\varphi\not\equiv 0$ by definition \eqref{eq:Moduli_space_projective_vortices_non-split-non-zero-section} of $\sM^{*,0}(E,h,\omega)$, then $(A,\varphi) \in \bV^{\natural,0} = \bV\cap\sP^{\natural,0}(E,h)$ by the definition \eqref{eq:bVnatural0} of $\bV^{\natural,0}$, the definition \eqref{eq:sPnatural0(E,h)} of $\sP^{\natural,0}(E,h)$, the inclusion $\sP^{\natural,0}(E,h) \subset \sP^{*,0}(E,h)$ in \eqref{eq:sPnatural0(E,h,omega)_subset_sP*0(E,h,omega)}, and the definition \eqref{eq:sP*0(E,h)} of $\sP^{*,0}(E,h)$. Furthermore, because $(A,\varphi)$ is a projective vortex, the identity \eqref{eq:Kobayashi_7-6-33_pairs} ensures that $(A,\varphi) \in \bmu^{-1}(0)\cap (\sA(E,h)\times W^{1,p}(E))$ and hence we have the inclusion of sets,
\[
  \sM_\reg^{*,0}(E,h,\omega) \subset \left.(\bmu^{-1}(0)\cap \bV^{\natural,0})\right/W^{2,p}(\SU(E)).
\]
Conversely, suppose that $(A,\varphi) \in \bmu^{-1}(0)\cap \bV^{\natural,0}$. By the definition $\bV'' = \cP_\reg^0(E,\omega)$ in \eqref{eq:Marsden-Weinstein_V_space_holomorphic_pairs} and the definition $\bV = (\pi_h^{0,1})^{-1}(\bV'') \subset \sA(E,h)\times W^{1,p}(E)$ in \eqref{eq:Marsden-Weinstein_V_space_unitary_pairs}, and the identity \eqref{eq:Kobayashi_7-6-33_pairs}, we see that $(A,\varphi)$ is a projective vortex with $\bH_{\bar\partial_A,\varphi}^2 = (0)$. Furthermore, by the definition $\bV^{\natural,0} = \bV\cap\sP^{\natural,0}(E,h)$ in \eqref{eq:bVnatural0}, the definition \eqref{eq:sPnatural0(E,h)} of $\sP^{\natural,0}(E,h)$, the inclusion $\sP^{\natural,0}(E,h) \subset \sP^{*,0}(E,h)$ in \eqref{eq:sPnatural0(E,h,omega)_subset_sP*0(E,h,omega)}, and the definition \eqref{eq:sP*0(E,h)} of $\sP^{*,0}(E,h)$, we must have that $A$ is non-split and $\varphi\not\equiv 0$. Hence, we obtain $\Stab(A,\varphi) = \{\id_E\}$ by the inclusion \eqref{eq:sM*0(E,h,omega)_subset_sM**(E,h,omega)} and definition \eqref{eq:Moduli_space_projective_vortices_StabAvarphi_idE} of $\sM^{**}(E,h,\omega)$. Thus, $\bH_{A,\varphi}^0 = (0)$, since $\bH_{A,\varphi}^0$ is the Lie algebra of $\Stab(A,\varphi)$
by Lemma \ref{lem:LieGroupStructureOfStab(A,varhi)}, and because $\bH_{\bar\partial_A,\varphi}^2 = (0)$, Theorem \ref{thm:Kobayashi_7-2-21_pairs} ensures that $\bH_{A,\varphi}^2 = (0)$. Therefore, $[A,\varphi] \in \sM_\reg^{*,0}(E,h,\omega)$ by definition \eqref{eq:Moduli_space_projective_vortices_non-split-non-zero-section_regular} of $\sM_\reg^{*,0}(E,h,\omega)$ and hence we have the reverse inclusion of sets,
\[
  \left.(\bmu^{-1}(0)\cap \bV^{\natural,0})\right/W^{2,p}(\SU(E)) \subset \sM_\reg^{*,0}(E,h,\omega),
\]
and this proves the claimed equality \eqref{eq:sMreg*0Ehomega_equals_symplectic_quotient}. 

The fact that $\bomega$ in \eqref{eq:Kobayashi_7-6-22_pairs} induces a K\"ahler form on the quotient
\[
  \left(\bmu^{-1}(0)\cap \bV^{\natural,0}\right)/W^{2,p}(\SU(E))
\]
follows from the Marsden--Weinstein Symplectic Reduction Theorem \ref{thm:Marsden-Weinstein_symplectic_quotient} (with $V$, $G$, and $\omega_V$ replaced by $\bV^{\natural,0}$,  $\sG_E$, and $\bomega|_{\bV^{\natural,0}} := \bomega\restriction\bV^{\natural,0}$, respectively) and our verification of its hypotheses in the paragraphs preceding the statement of the theorem. In particular, we obtain a K\"ahler form on $\sM_\reg^{**}(E,h,\omega)$ compatible with the induced complex structure and Riemannian metric on $\sM_\reg^{**}(E,h,\omega)$.

Consider Item \eqref{item:sMreg0(E.h,omega)_is_complex_Kaehler_if_E_rank_2}, where we assume that $E$ has rank two. According to Corollary \ref{cor:Lubke_Teleman_6-3-7_and_Kobayashi_7_3_17_pair} \eqref{item:fM_ps^0(E,omega)_rank-2_is_complex_analytic_space}, the moduli subspace $\fM_{\ps,\reg}^0(E,\omega)$ in \eqref{eq:Moduli_space_regular_non-zero-section_stable_holomorphic_pairs} of non-zero section, regular, polystable holomorphic pairs is a complex manifold and there is an isomorphism \eqref{eq:Lubke_Teleman_6-3-7_analytic_isomorphism_rank-2_non-zero_section_regular} of real analytic manifolds,
\[
  \sM_\reg^0(E,h,\omega) \cong \fM_{\ps,\reg}^0(E,\omega).
\]
Hence, the moduli space $\sM_\reg^0(E,h,\omega)$ inherits a complex structure from $\fM_{\ps,\reg}^0(E,\omega)$ via the preceding isomorphism of real analytic manifolds.

The proof of the identification of sets,
\begin{equation}
  \label{eq:sMreg0Ehomega_equals_symplectic_quotient_Erank2}
  \sM_\reg^0(E,h,\omega) = \left.(\bmu^{-1}(0)\cap \bV^0)\right/W^{2,p}(\SU(E)),
\end{equation}
where 
\begin{equation}
  \label{eq:Define_bV0}
  \bV^0 := \bV \cap \sP^0(E,h)
\end{equation}
and $\sP^0(E,h)$ is defined by \eqref{eq:sP0}, is similar though not identical to the proof of \eqref{eq:sMreg*0Ehomega_equals_symplectic_quotient}. Let $[A,\varphi]\in \sM_\reg^0(E,h,\omega)$. From the definition of $\sM_\reg^0(E,h,\omega)$ in \eqref{eq:Moduli_space_projective_vortices_non-zero_section_regular}, we have $\bH_{A,\varphi}^2 = (0)$. Since $(A,\varphi) \in \sP(E,h)$ in \eqref{eq:Subset_pre-holomorphic_pairs}, then $(\bar\partial_A,\varphi) \in \cP(E)$ in \eqref{eq:Subset_01-pairs_defining_holomorphic_pair_on_E}. Moreover, since $(A,\varphi)$ is a projective vortex and $\bH_{A,\varphi}^2 = (0)$, Theorem \ref{thm:Kobayashi_7-2-21_pairs} implies that $\bH_{\bar\partial_A,\varphi}^2 = (0)$. The isomorphism
\[
  \sM_\reg^0(E,h,\omega) \ni [A,\varphi] \leftrightarrow [\bar\partial_A,\varphi] \in \fM_{\ps,\reg}^0(E,\omega)
\]
given by \eqref{eq:Lubke_Teleman_6-3-7_analytic_isomorphism_rank-2_non-zero_section_regular} implies that $(\bar\partial_A,\varphi) \in \bV''$ by the definition \eqref{eq:Marsden-Weinstein_V_space_holomorphic_pairs} of $\bV'' = \cP_{\ps,\reg}^0(E,\omega)$, and so $(A,\varphi) \in \bV$ by the definition $\bV = (\pi_h^{0,1})^{-1}(\bV'') \subset \sA(E,h)\times W^{1,p}(E)$ in \eqref{eq:Marsden-Weinstein_V_space_unitary_pairs}. As $\varphi\not\equiv 0$ by definition \eqref{eq:Moduli_space_projective_vortices_non-zero_section} of $\sM^0(E,h,\omega)$, then
$(A,\varphi) \in \bV^0$ by the definition $\bV^0 = \bV\cap\sP^0(E,h)$ 
in \eqref{eq:Define_bV0} and the definition of $\sP^0(E,h)$ in \eqref{eq:sP0}. Furthermore, because $(A,\varphi)$ is a projective vortex, the identity \eqref{eq:Kobayashi_7-6-33_pairs} ensures that $(A,\varphi) \in \bmu^{-1}(0)\cap (\sA(E,h)\times W^{1,p}(E))$ and hence we have the inclusion of sets,
\[
  \sM_\reg^0(E,h,\omega) \subset \left.(\bmu^{-1}(0)\cap \bV^0)\right/W^{2,p}(\SU(E)).
\]
Conversely, suppose that $(A,\varphi) \in \bmu^{-1}(0)\cap \bV^0$. By the definition $\bV'' = \cP_{\ps,\reg}^0(E,\omega)$ in \eqref{eq:Marsden-Weinstein_V_space_holomorphic_pairs} and the definition $\bV = (\pi_h^{0,1})^{-1}(\bV'') \subset \sA(E,h)\times W^{1,p}(E)$ in \eqref{eq:Marsden-Weinstein_V_space_unitary_pairs}, and the identity \eqref{eq:Kobayashi_7-6-33_pairs}, we see that $(A,\varphi)$ is a projective vortex with $\bH_{\bar\partial_A,\varphi}^2 = (0)$. Furthermore, by the definition $\bV^0 = \bV\cap\sP^0(E,h)$ in \eqref{eq:Define_bV0} and the definition of $\sP^0(E,h)$ in \eqref{eq:sP0}, we must have that $\varphi\not\equiv 0$. Hence, we obtain $\Stab(A,\varphi) = \{\id_E\}$ by the equality \eqref{eq:sC0(E,h)_equals_sC**(E,h)_E_rank2} and definition \eqref{eq:Moduli_space_projective_vortices_StabAvarphi_idE} of $\sM^{**}(E,h,\omega)$. Thus, $\bH_{A,\varphi}^0 = (0)$, since $\bH_{A,\varphi}^0$ is the Lie algebra of $\Stab(A,\varphi)$ by Lemma \ref{lem:LieGroupStructureOfStab(A,varhi)}, and because $\bH_{\bar\partial_A,\varphi}^2 = (0)$, Theorem \ref{thm:Kobayashi_7-2-21_pairs} ensures that $\bH_{A,\varphi}^2 = (0)$. Therefore, $[A,\varphi] \in \sM_\reg^0(E,h,\omega)$ by definition \eqref{eq:Moduli_space_projective_vortices_non-zero_section_regular} of $\sM_\reg^0(E,h,\omega)$ and hence we have the reverse inclusion of sets,
\[
  \left.(\bmu^{-1}(0)\cap \bV^0)\right/W^{2,p}(\SU(E)) \subset \sM_\reg^0(E,h,\omega),
\]
and this proves the claimed equality \eqref{eq:sMreg0Ehomega_equals_symplectic_quotient_Erank2}. 

The remainder of the proof of Item \eqref{item:sMreg0(E.h,omega)_is_complex_Kaehler_if_E_rank_2} is identical to that of Item \eqref{item:sMreg*0(E.h,omega)_is_complex_Kaehler}. This completes the proof of Theorem \ref{thm:Kobayashi_7-6-36_pairs}.
\end{proof}

\begin{rmk}[Direct proof of the K\"ahler property for the moduli space of non-zero-section, regular projective vortices]
\label{rmk:Kaehler_moduli_space_non-zero-section_regular_projective_vortices}
The approach due to Itoh \cite[Section 4]{Itoh_1988} that we discussed in Remark \ref{rmk:Kaehler_moduli_space_irreducible_regular_projectively_Hermitian-Einstein_connections} to prove the K\"ahler property for the moduli space of non-split, regular, projectively Hermitian--Einstein connections --- without appealing to the Marsden--Weinstein Symplectic Reduction Theorem \ref{thm:Marsden-Weinstein_symplectic_quotient} --- should generalize to prove Theorem \ref{thm:Kobayashi_7-6-36_pairs}.
\end{rmk}

\subsection[Symplectic  structure for moduli space of projective vortices near singular points]{Marsden--Weinstein reduction and circle-invariant symplectic  structure for the moduli space of projective vortices near singular points}
\label{subsec:Marsden-Weinstein_reduction_moduli_space_SO3_monopoles_symplectic_quotient_singular_points}
We shall need an analogue of Theorem \ref{thm:Kobayashi_7-6-36_pairs} that applies to a circle-invariant open neighborhood of a point $[A,\varphi] \in  \sM^{**}(E,h,\omega)$ in \eqref{eq:Moduli_space_projective_vortices_StabAvarphi_idE}, where $\bH_{A,\varphi}^2$ in \eqref{eq:H_Avarphi^2} may be non-zero. In our applications in Section \ref{subsec:Characterization_critical_points_Hamiltonian_function_circle_action_on_projective_vortices}, we shall need a circle-invariant symplectic structure on a finite-dimensional manifold $\sM_{A,\varphi}^{\vir,\CC}(E,h,\omega) \subset \sC^{**}(E,h)$, as in the forthcoming definition \eqref{eq:sM_Avarphi^vir_circle_invariant}, that contains a circle-invariant open neighborhood of $[A,\varphi]$ in $\sM^{**}(E,h,\omega)$ through a circle-equivariant embedding of topological manifolds and an open neighborhood of $[A,\varphi]$ in $\sM_\reg^{**}(E,h,\omega)$ through a circle-equivariant embedding of symplectic manifolds. For this purpose, we shall use a construction that generalizes one due to Taubes \cite[Definition 8.1, p. 370]{Taubes_1987} for Yang--Mills connections over closed, four-dimensional, smooth Riemannian manifolds, an intrinsic version of his construction in \cite[Equations (2.7) and (2.8), p. 524]{TauIndef}.

Recall that $(\bar\partial_E,\varphi)$ is a holomorphic pair on $E$ if and only if it belongs to the zero locus of the map $\fS$ in \eqref{eq:Holomorphic_pair_map}, namely
\[
  \fS:\sA^{0,1}(E)\times W^{1,p}(E)
  \ni (\bar\partial_E,\varphi) \mapsto (F_{\bar\partial_E},\bar\partial_E\varphi) \in
  L^p\left(\Lambda^{0,2}(\fsl(E))\oplus \Lambda^{0,1}(E)\right).
\]
The preceding map is $\CC^*$-equivariant with respect to the standard action \eqref{eq:CZActionOnAffine}
induced by complex scalar multiplication on $W^{1,p}(E)$ in the domain and complex scalar multiplication on $L^p(\Lambda^{0,1}(E))$ in the codomain.

By analogy with the definition \eqref{eq:Box_dbar_E} of $\Box_{\bar\partial_E}$, consider the following Laplace
operators defined by the elliptic complex \eqref{eq:Holomorphic_pair_elliptic_complex} for a smooth holomorphic pair $(\bar\partial_E,\varphi)$ and integers $k=0,\ldots,n$, where $X$ has complex dimension $n$:
\begin{equation}
  \label{eq:Box_dbar_Evarphi}
  \Box_{\bar\partial_E,\varphi}^k
  :=
  \bar\partial_{E,\varphi}^{k,*}\bar\partial_{E,\varphi}^k + \bar\partial_{E,\varphi}^{k-1}\bar\partial_{E,\varphi}^{k-1,*}
  \quad\text{on } \Omega^{0,k}(\fsl(E)) \oplus \Omega^{0,k-1}(E).
\end{equation}
Although the operators $\Box_{\bar\partial_E,\varphi}^k$ have Sobolev and not smooth coefficients, they still have discrete spectra of eigenvalues that are contained in $[0,\infty)$. (See Feehan and Maridakis \cite[Proposition 2.2.3, p. 23]{Feehan_Maridakis_Lojasiewicz-Simon_coupled_Yang-Mills} for the statement and proof of a similar result for a Laplace operator defined by a unitary pair --- compactness of the bundle structure group is not required in the proof.) Let $\delta = \delta(\bar\partial_E,\varphi)$ be the minimum of the least positive eigenvalues of $\Box_{\bar\partial_E,\varphi}^k$ for $k=0,\ldots,n$. Define 
\[
  \Pi_{\delta;E,\varphi}^k:
  \Omega^{0,k}(\fsl(E)) \oplus \Omega^{0,k-1}(E)
  \to
  \Omega^{0,k}(\fsl(E)) \oplus \Omega^{0,k-1}(E)
\]
to be the finite-rank $L^2$-orthogonal spectral projection onto the span of the eigenvectors of $\Box_{\bar\partial_E,\varphi}^k$ with eigenvalue strictly less than $\delta(\bar\partial_E,\varphi)/2$, for $k=0,\ldots,n$.

Observe that if $(w,\nu) \in \Omega^{0,k}(\fsl(E)) \oplus \Omega^{0,k-1}(E)$ is an eigenvector of $\Box_{\bar\partial_E,\varphi}^k$ with eigenvalue $\mu$, then $\bar\partial_{E,\varphi}^k(w,\nu)$ is an eigenvector of $\Box_{\bar\partial_E,\varphi}^{k+1}$ with the same eigenvalue for $k=0,\ldots,n$:
\begin{align*}
  \Box_{\bar\partial_E,\varphi}^{k+1}\bar\partial_{E,\varphi}^k(w,\nu)
  &=
  \left(\bar\partial_{E,\varphi}^{k+1,*}\bar\partial_{E,\varphi}^{k+1}
    + \bar\partial_{E,\varphi}^k\bar\partial_{E,\varphi}^{k,*}\right)\bar\partial_{E,\varphi}^k(w,\nu)
  \\
  &=
  \bar\partial_{E,\varphi}^k\bar\partial_{E,\varphi}^{k,*}\bar\partial_{E,\varphi}^k(w,\nu)
  \\
  &=
  \bar\partial_{E,\varphi}^k\left(\bar\partial_{E,\varphi}^{k,*}\bar\partial_{E,\varphi}^k(w,\nu)
    + \bar\partial_{E,\varphi}^{k-1}\bar\partial_{E,\varphi}^{k-1,*}(w,\nu)\right)
  \\
  &=
  \bar\partial_{E,\varphi}^k\Box_{\bar\partial_E,\varphi}^k(w,\nu)
  \\
  &=
  \mu\bar\partial_{E,\varphi}^k(w,\nu).
\end{align*}
(We use the fact that $\bar\partial_{E,\varphi}^{k+1}\circ \bar\partial_{E,\varphi}^k = 0$ for $k=0,\ldots,n$ to obtain the previous identity.) Hence, the differentials $\bar\partial_{E,\varphi}^k$ and spectral projections $\Pi_{\delta;E,\varphi}^k$ commute in the following sense:
\[
  \bar\partial_{E,\varphi}^k\Pi_{\delta;E,\varphi}^k = \Pi_{\delta;E,\varphi}^{k+1}\bar\partial_{E,\varphi}^k,
  \quad\text{for } k = 0,1,\ldots,n.
\]  
Recall the $L^2$-orthogonal decomposition implied by Rudin \cite[Theorem 4.12, p. 99]{Rudin},
\[
  \Omega^{0,k}(\fsl(E)) \oplus \Omega^{0,k-1}(E)
  =
  \Ker\bar\partial_{E,\varphi}^k \oplus \Ran\bar\partial_{E,\varphi}^{k,*},
\]
and the equality
\[
   \Ker\bar\partial_{E,\varphi}^k = \Ran\bar\partial_{E,\varphi}^{k-1} \oplus \bH_{\bar\partial_E,\varphi}^k,
\]
implied by the definition \eqref{eq:H_dbar_Avarphi^0bullet} of the harmonic spaces $\bH_{\bar\partial_E,\varphi}^k$,
\[
  \bH_{\bar\partial_E,\varphi}^k
  = \Ker\left(\bar\partial_{E,\varphi}^k + \bar\partial_{E,\varphi}^{k-1,*} \right)
  = \Ker\Box_{\bar\partial_E,\varphi}^k.
\]
For $p \in (n,\infty)$ and a small enough open neighborhood of a given holomorphic pair $(\bar\partial_E,\varphi)$,
\[
  \UU_{\bar\partial_E,\varphi} \subset \sA^{0,1}(E) \times W^{1,p}(E),
\]
standard results on eigenvalue perturbation theory for linear operators (see Kato \cite{Kato}) yield the following isomorphisms for $k = 0,\ldots,n$ implied by $L^2$-orthogonal projections:
\[
  \Ran \Pi_{\delta;E+\alpha,\varphi+\phi}^k
  \cong
  \bH_{\bar\partial_E,\varphi}^k,
  \quad\text{for all } (\bar\partial_E+\alpha,\varphi+\phi) \in \UU_{\bar\partial_E,\varphi}.
\]
Since the preceding isomorphisms are equivariant with respect to the action \eqref{eq:SL(E)ActionOn(0,1)Pairs} of $W^{2,p}(\SU(E))$ (viewed as a Lie subgroup of $W^{2,p}(\SL(E))$) on $\sA^{0,1}(E) \times W^{1,p}(E)$ and the induced actions on the domains and codomains of the preceding isomorphisms, we may assume that $\UU_{\bar\partial_E,\varphi}$ is $W^{2,p}(\SU(E))$-invariant by replacing it with its orbit under the action of $W^{2,p}(\SU(E))$. We write
\[
  \Pi_{\delta;E+\alpha,\varphi+\phi}^{k,\perp}
  :=
  \id - \Pi_{\delta;E+\alpha,\varphi+\phi}^k,
  \quad\text{for all } (\bar\partial_E+\alpha,\varphi+\phi) \in \UU_{\bar\partial_E,\varphi},
\]
where $\id$ is the identity operator on $\Omega^{0,k}(\fsl(E)) \oplus \Omega^{0,k-1}(E)$. The sequence of operators $\bar\partial_{E+\alpha,\varphi+\phi}^k$ for $k=0,\ldots,n$ is only guaranteed to form a complex when $(\alpha,\phi) = (0,0)$, since $(\bar\partial_E,\varphi)$ is assumed to be a holomorphic pair, in which case the sequence of operators $\bar\partial_{E,\varphi}^k$ does form a complex for $k=0,\ldots,n$. 

It will be useful to consider the rolled-up complex obtained via the operators,
\[
  \bar D_{E+\alpha,\varphi+\phi}^k
  = \bar\partial_{E+\alpha,\varphi+\phi}^k + \bar\partial_{E+\alpha,\varphi+\phi}^{k-1,*}:
  \Omega^k \to \Omega^{k+1} \oplus \Omega^{k-1},
  \quad\text{for } k=0,\ldots,n,
\]
where we abbreviate $\Omega^k := \Omega^{0,k}(\fsl(E)) \oplus \Omega^{0,k-1}(E)$, and rolled-up operator,
\[
  \bar D_{E+\alpha,\varphi+\phi} = \bigoplus_{k \text{ odd}} \bar D_{E+\alpha,\varphi+\phi}^k :
  \bigoplus_{k \text{ odd}}\Omega^k \to \bigoplus_{l \text{ even}}\Omega^l.
\]
The operator $\bar D_{E,\varphi}$ has finite-dimensional kernel and cokernel given by
\begin{equation}
  \label{eq:Ker_and_coker_barD_Evarphi}
  \Ker \bar D_{E,\varphi} = \bigoplus_{k \text{ odd}} \bH_{\bar\partial_E,\varphi}^k
  \quad\text{and}\quad
  \Coker \bar D_{E,\varphi} = \bigoplus_{l \text{ even}} \bH_{\bar\partial_E,\varphi}^l,
\end{equation}
where the cokernel is defined to be the $L^2$-orthogonal complement of the range. 

We compose the operators $\bar D_{E+\alpha,\varphi+\phi}^k$ with $L^2$-orthogonal projections to define
\begin{multline*}
  \bar D_{\delta;E+\alpha,\varphi+\phi}^k
  :=
  \left(\Pi_{\delta;E+\alpha,\varphi+\phi}^{k+1,\perp}\bar\partial_{E+\alpha,\varphi+\phi}^k
    + \Pi_{\delta;E+\alpha,\varphi+\phi}^{k-1,\perp}\bar\partial_{E+\alpha,\varphi+\phi}^{k-1,*}\right)
  \Pi_{\delta;E+\alpha,\varphi+\phi}^{k,\perp}
  \\
  \Omega^k \to \Omega^{k+1} \oplus \Omega^{k-1},
  \quad\text{for } k=0,\ldots,n,
\end{multline*}
and thus define
\[
  \bar D_{\delta;E+\alpha,\varphi+\phi} = \bigoplus_{k \text{ odd}} \bar D_{\delta;E+\alpha,\varphi+\phi}^k :
  \bigoplus_{k \text{ odd}}\Omega^k \to \bigoplus_{l \text{ even}}\Omega^l.
\]
If $(\alpha,\phi) = (0,0)$, we obtain $\bar D_{\delta;E,\varphi}^k = \bar D_{E,\varphi}^k = \bar\partial_{E,\varphi}^k + \bar\partial_{E,\varphi}^{k-1,*}$. If $(\alpha,\phi) \neq (0,0)$ but $(\bar\partial_E+\alpha,\varphi+\phi) \in \UU_{\bar\partial_E,\varphi}$, the rolled-up operator is Fredholm with kernel and cokernel given by
\begin{align*}
  \Ker \bar D_{\delta;E+\alpha,\varphi+\phi}
  &= \bigoplus_{k \text{ odd}} \Ran \Pi_{\delta;E+\alpha,\varphi+\phi}^k
  \cong \bigoplus_{k \text{ odd}} \bH_{\bar\partial_E,\varphi}^k,
  \\
  \Coker \bar D_{\delta;E+\alpha,\varphi+\phi}
  &= \bigoplus_{l \text{ even}} \Ran \Pi_{\delta;E+\alpha,\varphi+\phi}^l
                                            \cong \bigoplus_{l \text{ even}} \bH_{\bar\partial_E,\varphi}^l.
\end{align*}
We obtain holomorphic finite-rank vector bundles over $\UU_{\bar\partial_E,\varphi}$ whose fibers $\Ran\Pi_\delta^k(\bar\partial_E+\alpha,\varphi+\phi)$ over a pair $(\bar\partial_E+\alpha,\varphi+\phi)$ are isomorphic to $\bH_{\bar\partial_E,\varphi}^k$, and these vector bundles are biholomorphic to the product vector bundles, $\UU_{\bar\partial_E,\varphi} \times \bH_{\bar\partial_E,\varphi}^k$ for $k = 0,\ldots,n$. We shall apply the preceding observations in the transversality argument below.

By replacing $\UU_{\bar\partial_E,\varphi}$ with its $S^1$ orbit, we may assume that $\UU_{\bar\partial_E,\varphi}$ is also invariant under the action of $S^1$ induced by the standard $\CC^*$ action \eqref{eq:CZActionOnAffine} on $\sA^{0,1}(E)\times W^{1,p}(E)$, recalling that this $S^1$ action commutes with the action of $W^{2,p}(\SU(E))$. For convenience in the argument below, we shall also write
\[
  \Pi_{\delta;E+\alpha,\varphi+\phi}^2
  =
  \Pi_\delta^2\left(\bar\partial_E+\alpha,\varphi+\phi\right),
  \quad\text{for all } (\bar\partial_E+\alpha,\varphi+\phi) \in \UU_{\bar\partial_E,\varphi}.
\]
In partial analogy with our definition of $\widehat\fS$ in the forthcoming \eqref{eq:Holomorphic_pair_map_Banach_spaces_orthogonal_projections}, we consider the map
\begin{multline}
  \label{eq:Holomorphic_pair_map_intrinsic_virtual}
  \fS_\delta:\sA^{0,1}(E)\times W^{1,p}(E) \supset
  \UU_{\bar\partial_E,\varphi} 
  \ni (\bar\partial_E+\alpha,\varphi+\sigma)
  \\
  \mapsto
  \Pi_\delta^{2,\perp}\left(\bar\partial_E+\alpha,\varphi+\phi\right)
  \fS\left(\bar\partial_E+\alpha,\varphi+\phi\right)
  \\
  \in L^p\left(\Lambda^{0,2}(\fsl(E))\oplus \Lambda^{0,1}(E)\right).
\end{multline}
We observe that $\fS_\delta$ is circle-equivariant with respect to the action of $S^1$ on the domain and the action of $S^1$ on $L^p(\Lambda^{0,1}(E))$ in the codomain induced by complex scalar multiplication on $E$:
\begin{multline*}
  \Pi_\delta^{2,\perp}\left(\bar\partial_E+\alpha,e^{i\theta}(\varphi+\phi)\right)
  \fS\left(\bar\partial_E+\alpha,e^{i\theta}(\varphi+\phi)\right)
  \\
  = e^{i\theta}\cdot \Pi_\delta^{2,\perp}\left(\bar\partial_E+\alpha,\varphi+\phi\right)
  \fS\left(\bar\partial_E+\alpha,\varphi+\phi\right),
  \\
  \text{for all } e^{i\theta} \in S^1 \text{ and }
  (\bar\partial_E+\alpha,\varphi+\phi) \in \UU_{\bar\partial_E,\varphi}.
\end{multline*}
Similarly, $\fS_\delta$ is $W^{2,p}(\SU(E))$-equivariant with respect to the action \eqref{eq:SL(E)ActionOn(0,1)Pairs} on the domain $\UU_{\bar\partial_E,\varphi} \subset \sA^{0,1}(E) \times W^{1,p}(E)$ and the action on $\Omega^2$ induced by the usual actions on $\fsl(E)$ and $E$:
\begin{multline*}
  \Pi_\delta^{2,\perp}\left(u\cdot\left(\bar\partial_E+\alpha,\varphi+\phi\right)\right)
  \fS\left(u\cdot\left(\bar\partial_E+\alpha,\varphi+\phi\right)\right)
  \\
  = u\cdot\left(\Pi_\delta^{2,\perp}\left(\bar\partial_E+\alpha,\varphi+\phi\right)
    \fS\left(\bar\partial_E+\alpha,\varphi+\phi\right)\right),
  \\
  \text{for all } u \in W^{2,p}(\SU(E)) \text{ and }
  (\bar\partial_E+\alpha,\varphi+\phi) \in \UU_{\bar\partial_E,\varphi}.
\end{multline*}
Since $\fS(\bar\partial_E,\varphi) = (0,0)$ by assumption that $(\bar\partial_E,\varphi)$ is a holomorphic pair, we see that $\fS_\delta$ has differential
\[
  d\fS_\delta(\bar\partial_E,\varphi):
  L^p\left(\fsl(\Lambda^{0,1}(E))\oplus E\right) \to L^p\left(\Lambda^{0,2}(\fsl(E))\oplus \Lambda^{0,1}(E)\right)
\]
given by
\[
  d\fS_\delta(\bar\partial_E,\varphi)
  =
  \Pi_\delta^{2,\perp}(\bar\partial_E,\varphi)d\fS(\bar\partial_E,\varphi)
  =
  \Pi_{E,\varphi}^{2,\perp}\bar\partial_{E,\varphi}^1
  =
  \bar\partial_{E,\varphi}^1\Pi_{E,\varphi}^{1,\perp},
\]
and therefore, since $\bar\partial_{E,\varphi}^1 = 0$ on $\Ran\Pi_{E,\varphi}^1 = \bH_{\bar\partial_E,\varphi}^1$,
\begin{equation}
  \label{eq:dfS_delta}
  d\fS_\delta(\bar\partial_E,\varphi) = \bar\partial_{E,\varphi}^1.
\end{equation}
Given a holomorphic pair $(\bar\partial_E,\varphi)$, we consider
\begin{multline*}
  \fS_\delta^{-1}(0) = \left\{(\bar\partial_E+\alpha,\varphi+\phi) \in \UU_{\bar\partial_E,\varphi}:
   \fS_\delta(\bar\partial_E+\alpha,\varphi+\phi) = (0,0) \right.
  \\
  \in \left. L^p\left(\Lambda^{0,2}(\fsl(E))\oplus \Lambda^{0,1}(E)\right)\right\}.
\end{multline*}
We claim that $\fS_\delta^{-1}(0)$ is an embedded complex submanifold of $\UU_{\bar\partial_E,\varphi}$. To see this, we begin by observing that 
\begin{multline*}
  \fS_\delta^{-1}(0) = \left\{(\bar\partial_E+\alpha,\varphi+\phi) \in \UU_{\bar\partial_E,\varphi}:
   \fS(\bar\partial_E+\alpha,\varphi+\phi) \in \Ran\Pi_\delta^2\left(\bar\partial_E+\alpha,\varphi+\phi\right) \right.
  \\
  \subset \left. L^p\left(\Lambda^{0,2}(\fsl(E))\oplus \Lambda^{0,1}(E)\right)\right\},
\end{multline*}
and thus
\[
  \fS_\delta^{-1}(0) = \fS^{-1}(\Xi_\delta^2),
\]
where $\Xi_\delta^2$ is the finite-rank holomorphic vector bundle over $\UU_{\bar\partial_E,\varphi}$ with fibers $\Ran\Pi_{\delta;E+\alpha,\varphi+\phi}^2$ over each pair $(\bar\partial_E+\alpha,\varphi+\phi) \in \UU_{\bar\partial_E,\varphi}$, and which is biholomorphic to the finite-rank product vector bundle, $\UU_{\bar\partial_E,\varphi} \times \bH_{\bar\partial_E,\varphi}^2$.

We augment the map $\fS$ by a holomorphic family of linear operators parametrized by $\UU_{\bar\partial_E,\varphi}$ to obtain a holomorphic map,
\[
  \widetilde\fS = \fS\times \bigoplus_{\text{odd } k > 1} D_\delta:
  \UU_{\bar\partial_E,\varphi}\times \bigoplus_{\text{odd } k > 1}W^{1,p}(\Lambda^k)
  \to
  \bigoplus_{\text{even } l > 0}L^p(\Lambda^l),
\]
where we abbreviate $\Lambda^k := \Lambda^{0,k}(\fsl(E))\oplus \Lambda^{0,k-1}(E)$, and define
\[
  D_\delta(\bar\partial_E+\alpha,\varphi+\phi)
  :=
  \bar D_{\delta;E+\alpha,\varphi+\phi}.
\]
We have the following more explicit expression for values of $\widetilde\fS$ and which helps illustrate the reason for the constraint $k>1$ in the sum over odd integers $k$ in the domain:
\begin{align*}
  \widetilde\fS\left(\left(\bar\partial_E+\alpha,\varphi+\phi\right),(\bw,\bnu)\right)
  &=
  \fS\left(\bar\partial_E+\alpha,\varphi+\phi\right)
  +
  \bar D_{\delta;E+\alpha,\varphi+\phi}(\bw,\bnu),
  \\
  &\qquad\text{for all }
  \left(\bar\partial_E+\alpha,\varphi+\phi\right) \in \UU_{\bar\partial_E,\varphi}
  \\
  &\qquad\qquad\text{and}\quad
  (\bw,\bnu) \in \bigoplus_{\text{odd } k > 1}W^{1,p}(\Lambda^k).
\end{align*}
By applying \eqref{eq:dfS_delta}, we see that $\widetilde\fS$ has differential at $(\bar\partial_E,\varphi)$ given by
\[
  d\widetilde\fS\left(\left(\bar\partial_E,\varphi\right),(0,0)\right)
  =
  \bar\partial_{E,\varphi}^1 + \bar D_{E,\varphi}:
  \bigoplus_{\text{odd } k}W^{1,p}(\Lambda^k) \to \bigoplus_{\text{even } l > 0} L^p(\Lambda^l).
\]
For the same reasoning underlying \eqref{eq:Ker_and_coker_barD_Evarphi}, we see that the preceding operator has cokernel
\[
  \Coker d\widetilde\fS\left(\left(\bar\partial_E,\varphi\right),(0,0)\right)
  =
  \bigoplus_{\text{even } l > 0} \bH_{\bar\partial_E,\varphi}^l.
\]
Hence, the augmented map $\widetilde\fS$ is transverse at $(\bar\partial_E,\varphi)$ to the finite-rank holomorphic vector subbundle of the product vector bundle,
\[
  \Xi_\delta = \bigoplus_{\text{even } l > 0}\Xi_\delta^l
  \subset
  \UU_{\bar\partial_E,\varphi} \times \bigoplus_{\text{even } l > 0}L^p(\Lambda^l),
\]
with fibers over each pair $(\bar\partial_E+\alpha,\varphi+\phi) \in \UU_{\bar\partial_E,\varphi}$ given by
\[
  \bigoplus_{\text{even } l > 0}\Ran\Pi_{\delta;E+\alpha,\varphi+\phi}^l.
\]
The vector subbundle $\Xi_\delta^l$ has fiber over a pair $(\bar\partial_E+\alpha,\varphi+\phi) \in \UU_{\bar\partial_E,\varphi}$ given by $\Ran\Pi_{\delta;E+\alpha,\varphi+\phi}^l$ for each even integer $l>0$. Moreover, $\Xi_\delta$ is isomorphic to the product vector bundle,
\[
  \UU_{\bar\partial_E,\varphi} \times \bigoplus_{\text{even } l>0}\bH_{\bar\partial_E,\varphi}^l.
\]
The preceding local triviality of $\Xi_\delta$ over $\UU_{\bar\partial_E,\varphi}$ (and consequently its manifold structure) allows us to appeal to standard results on transversality for complex analytic maps of Banach manifolds with respect to embedded complex submanifolds of the codomain. (See, for example, Lee \cite[Theorem 6.30 (a), p. 144]{Lee_john_smooth_manifolds} in the case of smooth maps of finite-dimensional manifolds and Lang \cite[Chapter II, Proposition 2.4, p. 27]{Lang_introduction_differential_topology} in the case of smooth or analytic maps of Banach manifolds.) The preimage $\widetilde\fS^{-1}(\Xi_\delta)$ is thus an embedded complex analytic submanifold that is biholomorphic to the product,
\[
  \fS^{-1}(\Xi_\delta^2) \times \bigoplus_{\text{odd } k > 1}\bH_{\bar\partial_E,\varphi}^k.
\]
In particular, $\fS_\delta^{-1}(0) = \fS^{-1}(\Xi_\delta^2)$ is an embedded complex analytic submanifold of $\UU_{\bar\partial_E,\varphi}$, as claimed. We thus choose
\[
  \bV^{\vir,\prime\prime} := \fS_\delta^{-1}(0),
\]
noting that this is an embedded complex analytic submanifold of $\UU_{\bar\partial_E,\varphi}$.

\begin{rmk}[Alternative argument for transversality of sections of Banach vector bundles]
\label{rmk:Alternative_argument_transversality_sections_Banach_vector_bundles} 
In the preceding argument, we relied on the easily-verified fact that $\Xi_\delta^l$ (for even $l>0$) and $\Xi_\delta$ are holomorphic vector bundles over $\UU_{\bar\partial_E,\varphi}$ that are isomorphic to product vector bundles. Instead of proving transversality of the augmented map $\widetilde\fS$ with respect to the finite-rank, holomorphic vector bundle $\Xi_\delta$, we could instead have tried to prove transversality of $\fS$ with respect to a complex analytic, infinite-rank vector bundle $\Xi_\delta^2\oplus \Ran\bar\partial_{\delta;\bullet}^{2,*}$ over $\UU_{\bar\partial_E,\varphi}$ with fibers
\[
  \Ran\Pi_{\delta;E+\alpha,\varphi+\phi}^2 \oplus \Ran\bar\partial_{\delta;E+\alpha,\varphi+\phi}^{2,*}.
\]
This approach involves proving that $\Xi_\delta^2\oplus \Ran\bar\partial_{\delta;\bullet}^{2,*}$ is isomorphic to the vector subbundle of
\[
  \UU_{\bar\partial_E,\varphi} \times L^p\left(\Lambda^{0,2}(\fsl(E)) \oplus \Lambda^{0,1}(E)\right)
\]
given by the product vector bundle,
\[
 \UU_{\bar\partial_E,\varphi} \times \bH_{\bar\partial_E,\varphi}^2 \oplus \Ran\bar\partial_{E,\varphi}^{2,*}.
\]
To accomplish this, one would verify that there are complex analytic, infinite-rank vector bundles,
\[
  \Ran\bar\partial_{\delta;\bullet}^{k-1},
  \quad \Ker\bar\partial_{\delta;\bullet}^k,
  \quad\text{and}\quad \Ran\bar\partial_{\delta;\bullet}^{k,*},
\]  
with fibers over pairs $(\bar\partial_E+\alpha,\varphi+\phi) \in \UU_{\bar\partial_E,\varphi}$ given by
\[
  \Ran\bar\partial_{\delta;E+\alpha,\varphi+\phi}^{k-1},
  \quad
  \Ker\bar\partial_{\delta;E+\alpha,\varphi+\phi}^k,
  \quad\text{and}\quad
  \Ran\bar\partial_{\delta;E+\alpha,\varphi+\phi}^{k,*},
\]
that are isomorphic to the product vector subbundles of
\[
  \UU_{\bar\partial_E,\varphi} \times L^p\left(\Lambda^{0,k}(\fsl(E)) \oplus \Lambda^{0,k-1}(E)\right)
\]
given by
\[
  \UU_{\bar\partial_E,\varphi} \times \Ran\bar\partial_{E,\varphi}^{k-1},
  \quad
  \UU_{\bar\partial_E,\varphi} \times \Ker\bar\partial_{E,\varphi}^k
  \quad\text{and}\quad
  \quad \UU_{\bar\partial_E,\varphi} \times \Ran\bar\partial_{E,\varphi}^{k,*}.
\]
For this purpose, we use the fact that each differential $\bar\partial_{E,\varphi}^k$ is an isomorphism from $(\Ker\bar\partial_{E,\varphi}^k)^\perp$ onto the closed subspace $\Ran\bar\partial_{E,\varphi}^k$, for $k = 0,1,2$.
\end{rmk}  

We now assume that $(A,\varphi)$ is a projective vortex with trivial stabilizer and that $(\bar\partial_E,\varphi) = \pi_h^{0,1}(A,\varphi)$ and observe that
\begin{equation}
\label{eq:VirtualSpaceForProjectiveVortices}
  \bV^\vir := \left(\pi_h^{0,1}\right)^{-1}\left(\bV^{\vir,\prime\prime}\right)
\end{equation}
is an embedded real analytic submanifold of the open neighborhood $\UU_{A,\varphi}$ of $(A,\varphi)$ in $\sA(E,h)\times W^{1,p}(E)$ given by
\[
  \UU_{A,\varphi}
  :=
  \left(\pi_h^{0,1}\right)^{-1}\left(\UU_{\bar\partial_A,\varphi}\right).
\]
The map $\pi_h^{0,1}$ in \eqref{eq:Bijection_unitarypairs_with_01pairs} is an $S^1$-equivariant isomorphism of real Banach affine spaces with respect to the standard $S^1$ action \eqref{eq:S1_Action_On_AffineSpaceForProjectiveVortices} on the domain and the $S^1$ action implied by the standard $\CC^*$ action \eqref{eq:CZActionOnAffine} on the codomain; moreover, these $S^1$ actions commute with the actions of $W^{2,p}(\SU(E))$ on the domain and codomain. We define
\begin{equation}
  \label{eq:sM_Avarphi^vir_circle_invariant}
  \sM_\symp^\vir(E,h,\omega) := \left.\left(\bmu^{-1}(0)\cap\bV^\vir\right)\right/W^{2,p}(\SU(E))
\end{equation}
to be the \emph{local symplectic virtual moduli space of projective vortices}. The forthcoming Theorem \ref{thm:Kobayashi_7-6-36_pairs_local_virtual_moduli_space} provides a circle-invariant symplectic structure on $\sM_\symp^\vir(E,h,\omega)$. (We shall see in Section \ref{sec:Complex_Kaehler_structure_moduli_space_projective_vortices_near_singular_points} how to construct a local complex virtual moduli space of projective vortices, $\sM_{A,\varphi}^{\vir,\CC}(E,h,\omega)$, via the comparison of Kuranishi models in Section \ref{subsec:Friedman-Morgan_4-3-4} yields a complex K\"ahler structure on $\sM_{A,\varphi}^{\vir,\CC}(E,h,\omega)$, but its circle invariance is less transparent because of the need to choose a representative $(A,\varphi)$ for a point $[A,\varphi]$ in the quotient space.)

\begin{thm}[Local virtual moduli space of non-zero-section projective vortices as a circle-invariant symplectic manifold]
\label{thm:Kobayashi_7-6-36_pairs_local_virtual_moduli_space}
Continue the hypotheses of Theorem \ref{thm:Kobayashi_7-6-36_pairs}. If $[A,\varphi]\in \sM^{**}(E,h,\omega)$ as in \eqref{eq:Moduli_space_projective_vortices_StabAvarphi_idE}, then there is an open neighborhood $\UU_{A,\varphi}$ of the pair $(A,\varphi)$ in $\sA(E,h)\times W^{1,p}(E)$ that is $S^1$-invariant with respect to the circle action \eqref{eq:S1_Action_On_AffineSpaceForProjectiveVortices} and such that the following hold.
\begin{enumerate}
\item\label{item:Local_virtual_moduli_space_properties}
The local virtual moduli space $\sM_\symp^\vir(E,h,\omega)$ in \eqref{eq:sM_Avarphi^vir_circle_invariant} is an embedded real analytic submanifold of $\sC^{**}(E,h)$ as in \eqref{eq:ConfigurationSpaceForProjectiveVortices**}, is an $S^1$-invariant symplectic manifold with respect to the circle action \eqref{eq:S1_Action_On_ConfigurationSpaceForProjectiveVortices} on $\sC(E,h)$ with symplectic form induced from $\bomega$ in \eqref{eq:Kobayashi_7-6-22_pairs}, and has tangent space at $[A,\varphi]$ represented by $\bH_{A,\varphi}^1$ in \eqref{eq:H_Avarphi^1}.

\item\label{item:Local_virtual_moduli_space_identities}
  The following equality holds,
\[
  \sM(E,h,\omega)\cap\pi(\UU_{A,\varphi}) = \sM(E,h,\omega)\cap \sM_\symp^\vir(E,h,\omega),
\]
and $\sM(E,h,\omega)\cap \pi(\UU_{A,\varphi})$ is an $S^1$-equivariantly embedded topological submanifold of $\sM_\symp^\vir(E,h,\omega)$, where $\pi:\sA(E,h)\times W^{1,p}(E) \to \sC(E,h)$ is the quotient map associated to the definition \eqref{eq:ConfigurationSpaceForProjectiveVortices} of the quotient space $\sC(E,h)$. Moreover, the preceding equality restricts to equalities,
\begin{align*}
  \sM^{**}(E,h,\omega)\cap \pi(\UU_{A,\varphi})
  &=
    \sM^{**}(E,h,\omega)\cap \sM_\symp^\vir(E,h,\omega),
  \\
   \sM_\reg^{**}(E,h,\omega)\cap \pi(\UU_{A,\varphi})
  &=
  \sM_\reg^{**}(E,h,\omega)\cap \sM_\symp^\vir(E,h,\omega),
\end{align*}  
and $\sM_\reg^{**}(E,h,\omega)\cap \pi(\UU_{A,\varphi})$ is an $S^1$-equivariantly embedded symplectic submanifold of $\sM_\symp^\vir(E,h,\omega)$, where $\sM^{**}(E,h,\omega)$ is as in \eqref{eq:Moduli_space_projective_vortices_StabAvarphi_idE} and $\sM_\reg^{**}(E,h,\omega)$ is as in \eqref{eq:Moduli_space_projective_vortices_StabAvarphi_idE_regular}.
\end{enumerate}
\end{thm}

\begin{proof}
We follow the construction and notation preceding the statement of the theorem. Consider Item \eqref{item:Local_virtual_moduli_space_properties}. By hypothesis, $\Stab(A,\varphi) = \{\id_E\}$ and so Lemma \ref{lem:Openness_configuration_subspace_unitary_pairs_trivial_stabilizer} \eqref{item:Openness_subspace_unitary_pairs_minimal_stabilizer} implies that $\Stab(A+a,\varphi+\varphi) = \{\id_E\}$ for all pairs $(A+a,\varphi+\varphi)$ in a small enough open neighborhood $\UU_{A,\varphi}$ of $(A,\varphi)$. Hence, $\pi(\UU_{A,\varphi}) \subset \sC^{**}(E,h)$ and $W^{2,p}(\SU(E))$ acts freely on the real analytic embedded submanifold $\bV^\vir$ defined in \eqref{eq:VirtualSpaceForProjectiveVortices}.

If $(A,\varphi) \in \bmu^{-1}(0)\cap\sP(E,h)$ as in \eqref{eq:Subset_pre-holomorphic_pairs} is a regular point of the moment map $\bmu$ in \eqref{eq:Moment_map_action_unitary_det_one_gauge_transformations_affine_space_pairs},
\[
  \bmu:\sA(E,h)\times W^{1,p}(E) \to \left(W^{2,p}(\su(E))\right)^*,
\]  
then the proof of Lemma \ref{lem:Regular_points_bmu_remain_regular_on_restriction_from_sA(E,h)_times_W1p(E)_to_sP(E,h)} implies that $(A,\varphi)$ is also a regular point of the restriction of the moment map $\bmu$ to the embedded real analytic submanifold $\bV^\vir$,
\[
  \bmu:\bV^\vir \to \left(W^{2,p}(\su(E))\right)^*.
\]
Therefore, after possibly shrinking $\UU_{A,\varphi}$ and then replacing it by its orbit under the action of $S^1\times W^{2,p}(\SU(E))$ defined by \eqref{eq:DefineSU(E)ActionOnUnitaryPairs} and \eqref{eq:S1_Action_On_AffineSpaceForProjectiveVortices},
we obtain that $0$ is a regular value of $\bmu$ on $\bV^\vir$ and the preceding moment map satisfies Hypothesis \ref{hyp:Marsden-Weinstein_symplectic_quotient_conditions} and $\sM_\symp^\vir(E,h,\omega)$ is an embedded real analytic submanifold of the real analytic Banach manifold $\sC^{**}(E,h)$. (Because the actions of gauge transformations $u \in W^{2,p}(\SU(E))$ and $e^{i\theta} \in S^1$ on $\sA(E,h)\times W^{1,p}(E)$ commute, we have
\[
  \Stab(A,e^{i\theta}\varphi) = \Stab(A,\varphi) = \{\id_E\},
\]
and so the Banach manifold $\sC^{**}(E,h)$ is indeed a circle-invariant subspace of $\sC(E,h)$.) In particular, the Marsden--Weinstein Symplectic Reduction Theorem \ref{thm:Marsden-Weinstein_symplectic_quotient} implies that $\sM_\symp^\vir(E,h,\omega)$ inherits a symplectic structure $\bomega$ from the weak symplectic structure $\bomega$ on $\sA(E,h)\times W^{1,p}(E)$ defined in \eqref{eq:Kobayashi_7-6-22_pairs}.

By the definitions \eqref{eq:dAvarphi1;1} of $d_{A,\varphi}^{1;1}$ and \eqref{eq:d1StablePair} of $\bar\partial_{A,\varphi}^1$, we see that 
\begin{equation}
  \label{eq:dAvarphi1;1_equals_dbar_Avarphi1_pi_h01}
  d_{A,\varphi}^{1;1} = \bar\partial_{A,\varphi}^1\circ d\pi_h^{0,1}:
  \Omega^1(\su(E))\oplus\Omega^0(E) \to \Omega^{0,2}(\fsl(E))\oplus\Omega^{0,1}(E),
\end{equation}
where the definition \eqref{eq:Bijection_unitarypairs_with_01pairs} of $\pi_h^{0,1}$ yields the following isomorphism of Fr\'echet spaces:
\[
  d\pi_h^{0,1}:\Omega^1(\su(E))\oplus\Omega^0(E) \to \Omega^{0,1}(\fsl(E))\oplus\Omega^0(E).
\]  
By definition \eqref{eq:sM_Avarphi^vir_circle_invariant} of $\sM_\symp^\vir(E,h,\omega)$ and the existence of a local affine, Coulomb-gauge slice $S_{A,\varphi} = (A,\varphi) + \Ker d_{A,\varphi}^{0,*}$ through $(A,\varphi)$ for the action of $W^{2,p}(\su(E))$ on $\sA(E,h)\times W^{1,p}(E)$, we obtain that
\begin{align*}
  T_{A,\varphi}\sM_\symp^\vir(E,h,\omega)
  &= \Ker d\bmu(A,\varphi)\cap \Ker\left(d\fS_\delta(\bar\partial_A,\varphi)\circ d\pi_h^{0,1}\right)
    \cap \Ker d_{A,\varphi}^{0,*}
  \\
  &=
    \Ker d\bmu^*(A,\varphi)\cap \Ker\left(\bar\partial_{A,\varphi}^1\circ d\pi_h^{0,1}\right)
     \cap \Ker d_{A,\varphi}^{0,*}  \quad\text{(by \eqref{eq:dfS_delta})}
  \\
  &=
    \Ker d_{A,\varphi}^{1;0} \cap \Ker d_{A,\varphi}^{1;1}
    \cap \Ker d_{A,\varphi}^{0,*} \quad\text{(by \eqref{eq:dbmu_Avarphi_equals_dAvarphi1;0}
    and \eqref{eq:dAvarphi1;1_equals_dbar_Avarphi1_pi_h01})}
  \\
  &=
    \Ker d_{A,\varphi}^1 \cap \Ker d_{A,\varphi}^{0,*} \quad\text{(by \eqref{eq:dAvarphi1;0_and_1;1})}
  \\
  &= \bH_{A,\varphi}^1 \quad\text{(by \eqref{eq:H_Avarphi^bullet})},
\end{align*}  
as claimed. The preceding arguments also imply that $\sM_\symp^\vir(E,h,\omega)$ is an $S^1$-invariant subspace of $\sC(E,h)$ and this completes the proof of Item \eqref{item:Local_virtual_moduli_space_properties}. The statements about set-theoretic inclusions and circle-equivariant embeddings in Item \eqref{item:Local_virtual_moduli_space_identities} are clear.
\end{proof}

\begin{rmk}[K\"ahler spaces and Theorem \ref{thm:Kobayashi_7-6-36_pairs_local_virtual_moduli_space}]
\label{rmk:Kaehler_spaces}  
We may view the construction in Theorem \ref{thm:Kobayashi_7-6-36_pairs_local_virtual_moduli_space} as a replacement of the concept of a \emph{K\"ahler space}, originally due to Grauert \cite{Grauert_1962} and later modified by Moishezon \cite{Moishezon_1975} --- a class of complex analytic spaces equipped with suitably defined K\"ahler metrics. We refer to
%
%
Bingener \cite{Bingener_1983a, Bingener_1983b}, Fischer \cite{Fischer_complex_analytic_geometry}, Fujiki \cite{Fujiki_1983}, Grauert and Remmert \cite{Grauert_Remmert_coherent_analytic_sheaves}, H\"ormander \cite{Hormander_introduction_complex_analysis_several_variables}, Koll\'ar \cite{Kollar_lectures_resolution_singularities}, Treger \cite{Treger_2016arxiv}, and Varouchas \cite{Varouchas_1989} for related results.
\end{rmk}

\subsection{Characterization of critical points of the Hamiltonian function for the circle action on the moduli space of projective vortices}
\label{subsec:Characterization_critical_points_Hamiltonian_function_circle_action_on_projective_vortices}
We begin with the

\begin{defn}[Critical point of a smooth function on the moduli space of projective vortices]
\label{defn:Critical_point_Hitchin_Hamiltonian_function_moduli_space_projective_vortices}
Let $(E,h)$ be a smooth Hermitian vector bundle over a closed, almost Hermitian manifold $(X,g,J)$ with fundamental two-form $\omega = g(\cdot,J\cdot)$ as in \eqref{eq:Fundamental_two-form} and real dimension $2n$. Let $p \in (n,\infty)$ be a constant and $A_d$ be a fixed unitary connection on the Hermitian line bundle $\det E$. If $[A,\varphi]$ is a point in the moduli space $\sM^{**}(E,h,\omega)$ in \eqref{eq:Moduli_space_projective_vortices_StabAvarphi_idE}, then $[A,\varphi]$ is a \emph{critical point} of the restriction to $\sM^{**}(E,h,\omega)$ of the smooth function, 
  \begin{equation}
    \label{eq:Hamiltonian_function}
    f: \sC^{**}(E,h) \ni [A,\varphi]
  \mapsto \frac{1}{2}\|\varphi\|_{L^2(X)}^2 \in \RR,
\end{equation}
if one has
\begin{equation}
  \label{eq:Critical_point_Hitchin_Hamiltonian_function_moduli_projective_vortices}
  T_{[A,\varphi]}\sM^{**}(E,h,\omega) \subset \Ker df[A,\varphi]
\end{equation}
where $T_{[A,\varphi]}\sM^{**}(E,h,\omega) \cong \bH_{A,\varphi}^1$ is the \emph{Zariski tangent space} to $\sM^{**}(E,h,\omega)$ at $[A,\varphi]$.
\end{defn}

To see that $f$ in \eqref{eq:Hamiltonian_function} is a \emph{Hamiltonian function} for the circle action \eqref{eq:S1_Action_On_ConfigurationSpaceForProjectiveVortices} in the sense of \eqref{eq:MomentMap}, observe that the vector field $\bX$ on $\sA(E,h)\times W^{1,p}(E)$ generating the action of $e^{i\theta}\in S^1$ induced by scalar multiplication on $W^{1,p}(E)$ is given by
\begin{multline*}
  \bX_{A,\varphi}
  = \left.\frac{d}{d\theta}e^{i\theta}\cdot(A,\varphi)\right|_{\theta=0}
  = \left.\frac{d}{d\theta}(A,e^{i\theta}\varphi)\right|_{\theta=0}
  \\
  = (0,i\varphi)
  = \bJ_{A,\varphi}(0,\varphi)
  \in W^{1,p}(T^*X\otimes\su(E) \oplus E).
\end{multline*}
For convenience, we record the equality
\begin{equation}
\label{eq:S1GeneratorExpression}
 \bX_{A,\varphi}
 =
 (0,i\varphi), \quad\text{for all } (A,\varphi) \in \sA(E,h)\times W^{1,p}(E).
\end{equation}
On the other hand, for $(a,\phi) \in W^{1,p}(T^*X\otimes\su(E) \oplus E)$, we have
\begin{multline*}
  df(A,\varphi)(a,\phi)
  = \left.\frac{d}{dt}f(A+ta,\varphi+t\phi)\right|_{t=0}
  = \frac{1}{2}\left.\frac{d}{dt}\|\varphi+t\phi\|_{L^2(X)}^2\right|_{t=0}
  \\
  = \frac{1}{2}\left((\varphi,\phi)_{L^2(X)} + (\phi,\varphi)_{L^2(X)}\right)
  = \Real(\phi,\varphi)_{L^2(X)}.
\end{multline*}
Using the expressions for $\bX_{A,\varphi}$ and $\bomega$ in \eqref{eq:Kobayashi_7-6-22_pairs} at the point $(A,\varphi)$ in terms of $\bJ$ in \eqref{eq:Almost_complex_structure_affine_space_pairs_unitary_connections_and_sections} and $\bg$ in \eqref{eq:L2_metric_affine_space_pairs_unitary_connections_and_sections} and the compatibility of $\bJ$ and $\bg$, we see that
\begin{align*}
  \iota_{\bX_{A,\varphi}}\bomega_{A,\varphi}(a,\phi)
  &= \bomega(\bX_{A,\varphi},(a,\phi))
  \\
  &= \bomega(\bJ(0,\varphi),(a,\phi))
  \\
  &= \bg(\bJ(0,\varphi),\bJ(a,\phi))
  \\
  &= \bg((0,\varphi),(a,\phi))
  \\
  &= (0,a)_{L^2(X)} + \Real(\varphi,\phi)_{L^2(X)},
\end{align*}
and therefore we obtain
\begin{multline}
  \label{eq:iota_bX_bomega}
  \iota_{\bX_{A,\varphi}}\bomega_{A,\varphi}(a,\phi)
  = \Real(\phi,\varphi)_{L^2(X)},
  \\
  \text{for all } (A,\varphi) \in \sA(E,h)\times W^{1,p}(E)
  \text{ and } (a,\phi) \in W^{1,p}(T^*X\otimes\su(E)\oplus E).
\end{multline}
Thus, since $(a,\phi) \in W^{1,p}(T^*X\otimes\su(E) \oplus E)$ was arbitrary, we see that
\begin{equation}
  \label{eq:Moment_map_affine_space_unitary_pairs}
  df(A,\varphi) = \iota_{\bX_{A,\varphi}}\bomega_{A,\varphi}
  = \left(\iota_\bX\bomega\right)_{A,\varphi},
  \quad\text{for all } (A,\varphi) \in \sA(E,h)\times W^{1,p}(E).
\end{equation}
Consequently, $f$ in \eqref{eq:Hamiltonian_function} is a Hamiltonian function in the sense of \eqref{eq:MomentMap} for the natural circle action on $\sA(E,h)\times W^{1,p}(E)$.

The function $f$ in \eqref{eq:Hamiltonian_function} is well-defined on the quotient $\sC(E,h)$ in \eqref{eq:ConfigurationSpaceForProjectiveVortices}. According to Theorems \ref{thm:Kobayashi_7-6-36_pairs} and \ref{thm:Kobayashi_7-6-36_pairs_local_virtual_moduli_space}, when we restrict to the case where $(X,g,J)$ is a complex, K\"ahler manifold with K\"ahler form $\omega=g(\cdot,J\cdot)$ and the curvature of $A_d$ on $\det E$ obeys $F_{A_d}^{0,2}=0$, the moduli spaces $\sM_\reg^{**}(E,h,\omega)$ in \eqref{eq:Moduli_space_projective_vortices_StabAvarphi_idE_regular} and $\sM_\symp^\vir(E,h,\omega)$ in \eqref{eq:sM_Avarphi^vir_circle_invariant} are complex, K\"ahler manifolds with K\"ahler form $\bomega$ determined by \eqref{eq:Kobayashi_7-6-22_pairs}. We claim that the relation \eqref{eq:Moment_map_affine_space_unitary_pairs} descends to the quotient to give
\begin{equation}
  \label{eq:Moment_map_moduli_space_projective_vortices}
  df[A,\varphi] = \iota_{\bX_{[A,\varphi]}}\bomega_{[A,\varphi]}
   = \left(\iota_\bX\bomega\right)_{[A,\varphi]},
\end{equation}
for all points $[A,\Phi]$ in the quotient spaces $\sM_\reg^{**}(E,h,\omega)$ or $\sM_\symp^\vir(E,h,\omega)$, where $\bX$ is the vector field generating the circle action \eqref{eq:S1_Action_On_ConfigurationSpaceForProjectiveVortices} on the open subspace given by the Banach manifold
\[
  \sC^{**}(E,h) := \left\{[A,\varphi] \in \sC(E,h): \Stab(A,\varphi) = \{\id_E\} \right\} \subset \sC(E,h).
\]
Indeed, to prove this claim, we observe that the vector field $\bX$ on $\sA(E,h)\times W^{1,p}(E)$ defined by \eqref{eq:S1GeneratorExpression} is invariant under the action of $u\in W^{2,p}(\SU(E))$ since 
\begin{multline*}
  \bX_{u(A,\varphi)} = (0,iu(\varphi)) = u_*(0,i\varphi) \in W^{1,p}(T^*X\otimes\su(E)\oplus E),
  \\
  \text{for all } (A,\varphi) \in \sA(E,h)\times W^{1,p}(E).
\end{multline*}
We have already noted that, by definition, $f:\sA(E,h)\times W^{1,p}(E) \to \RR$ is invariant under the action of $W^{2,p}(\SU(E))$, and so the one-form $df$ on $\sA(E,h)\times W^{1,p}(E)$ is invariant too under the action of $u\in W^{2,p}(\SU(E))$:
\begin{multline*}
  df(u(A,\varphi))(u_*(a,\phi)) = df(A,\varphi)(a,\phi),
  \\
  \text{for all } (A,\varphi) \in \sA(E,h)\times W^{1,p}(E)
  \text{ and } (a,\phi) \in W^{1,p}(T^*X\otimes\su(E)\oplus E).
\end{multline*}
Next, we observe that the expression \eqref{eq:iota_bX_bomega} implies that the one-form $\iota_\bX\bomega$ on $\sA(E,h)\times W^{1,p}(E)$ is also invariant under the action of $u\in W^{2,p}(\SU(E))$: 
\begin{align*}
  \left(\iota_\bX\bomega\right)_{u(A,\varphi)}(u_*(a,\phi))
  &= \iota_{\bX_{u(A,\varphi)}}\bomega_{u(A,\varphi)}(u_*(a,\phi))
  \\
  &= \Real(u(\phi),u(\varphi))_{L^2(X)}
  = \Real(\phi,\varphi)_{L^2(X)}
  = \left(\iota_\bX\bomega\right)_{A,\varphi}(a,\phi),
  \\
  &\qquad \text{for all } (A,\varphi) \in \sA(E,h)\times W^{1,p}(E)
  \\
  &\qquad \text{and } (a,\phi) \in W^{1,p}(T^*X\otimes\su(E)\oplus E).
\end{align*}
By combining the preceding observations and noting that $\bomega$ in \eqref{eq:Moment_map_moduli_space_projective_vortices} is the well-defined two-form on the quotient spaces
$\sM_\reg^{**}(E,h,\omega)$ and $\sM_\symp^\vir(E,h,\omega)$ provided by Theorems \ref{thm:Kobayashi_7-6-36_pairs} and \ref{thm:Kobayashi_7-6-36_pairs_local_virtual_moduli_space}, we see that the identity \eqref{eq:Moment_map_moduli_space_projective_vortices} holds, as claimed.  Therefore, $f$ in \eqref{eq:Hamiltonian_function} is a Hamiltonian function for the induced circle action on
$\sM_\reg^{**}(E,h,\omega)$ or $\sM_\symp^\vir(E,h,\omega)$.

\begin{thm}[Characterization of critical points of the Hamiltonian function on the moduli space of projective vortices]
\label{thm:Critical_points_Hitchin_Hamiltonian_function_moduli_space_projective_vortices}
Let $(E,h)$ be a smooth Hermitian vector bundle over a closed, complex K\"ahler manifold $(X,g,J)$ with K\"ahler form $\omega = g(\cdot,J\cdot)$ and complex dimension $n$. Let $p \in (n,\infty)$ be a constant and $A_d$ be a fixed unitary connection on the Hermitian line bundle $\det E$ with curvature obeying $F_{A_d}^{0,2}=0$. Then a point $[A,\varphi]$ is a critical point (in the sense of Definition \ref{defn:Critical_point_Hitchin_Hamiltonian_function_moduli_space_projective_vortices}) of the restriction of the Hamiltonian function \eqref{eq:Hamiltonian_function} to $\sM^{**}(E,h,\omega)$ in \eqref{eq:Moduli_space_projective_vortices_StabAvarphi_idE},
\[
  f: \sM^{**}(E,h,\omega) \ni [A,\varphi] \mapsto \frac{1}{2}\|\varphi\|_{L^2(X)}^2 \in \RR,
\]
if and only if $[A,\varphi]$ is a fixed point of the $S^1$ action on $\sM^{**}(E,h,\omega)$.
\end{thm}

\begin{proof}
Let $[A,\varphi]$ be a fixed point of the $S^1$ action on $\sM^{**}(E,h,\omega)$, so it is also a fixed point of the $S^1$ action on $\sM_\symp^\vir(E,h,\omega)$ defined by the action \eqref{eq:S1_Action_On_ConfigurationSpaceForProjectiveVortices}. By Theorem \ref{thm:Kobayashi_7-6-36_pairs_local_virtual_moduli_space}, the local virtual moduli space $\sM_\symp^\vir(E,h,\omega)$ is an almost K\"ahler manifold and so Frankel's Theorem \ref{mainthm:Frankel_almost_Hermitian} implies that $[A,\varphi]$ is thus a critical point of $f:\sM_\symp^\vir(E,h,\omega)\to\RR$ in the usual sense and hence also a critical point of $f:\sM^{**}(E,h,\omega)\to\RR$ in the sense of Definition \ref{defn:Critical_point_Hitchin_Hamiltonian_function_moduli_space_projective_vortices} because, by construction, the tangent space to $\sM_\symp^\vir(E,h,\omega)$ at $[A,\varphi]$ is equal to the Zariski tangent space to
$\sM^{**}(E,h,\omega)$ at $[A,\varphi]$.

Conversely, let $[A,\varphi]$ be a critical point of $f:\sM^{**}(E,h,\omega)\to\RR$ in the sense of Definition \ref{defn:Critical_point_Hitchin_Hamiltonian_function_moduli_space_projective_vortices}. Consequently, $[A,\varphi]$ is a critical point of $f:\sM_\symp^\vir(E,h,\omega)\to\RR$ in the usual sense and Frankel's Theorem \ref{mainthm:Frankel_almost_Hermitian} implies that $[A,\varphi]$ is a fixed point of the $S^1$ action on $\sM_\symp^\vir(E,h,\omega)$ and thus also a fixed point of the $S^1$ action on $\sM^{**}(E,h,\omega)$.
\end{proof}

\begin{proof}[Proof of Theorem \ref{mainthm:IdentifyCriticalPoints}]
Let $(E,h)$ be a smooth, rank-two Hermitian vector bundle over $X$ with the property that the \spinu structure $\ft$ is given by $\ft=(\rho,W_\can\otimes E)$, where $(\rho,W_\can)$ is the canonical \spinc structure of Definition \ref{defn:Canonical_spinc_bundles}. Let $[A,\Phi]\in \sM_\ft^0$ be a critical point of the smooth function $f$ in \eqref{eq:Hitchin_function} with $\Phi\not\equiv 0$. By Lemma \ref{lem:NonZeroSection_Spinu_pairs_Have_Trivial_Stabilizer}
(the analogue for \spinu pairs of Lemma \ref{lem:Stabilizers_non-zero-section_unitary_pairs}) and the hypotheses that $E$ has complex rank two and $\Phi\not\equiv 0$, we have $\Stab(A,\Phi)=\{\id_E\}$. Trivial stabilizers are required by the restriction in Theorem
\ref{thm:Critical_points_Hitchin_Hamiltonian_function_moduli_space_projective_vortices} of $f$ to $\sM^{**}(E,h,\omega)$. Because $E$ has complex rank two, we have $\sM^{**}(E,h,\omega) = \sM^0(E,h,\omega)$ by the equality \eqref{eq:sC0(E,h)_equals_sC**(E,h)_E_rank2} implied by Lemma
\ref{lem:sM*0(E,h,omega)_subset_sM**(E,h,omega)} \eqref{item:Subspace_trivial_stabilizer_equals_non-zero-section_rank-2_unitary_pairs}.

The identifications of the moduli space $\sM^0(E,h,\omega)$ of non-zero-section projective vortices on $E$ with the moduli space of type $1$ non-zero-section non-Abelian monopoles $\sM_{\ft,1}^0$ and the moduli space $\sM^0(F,h,\omega)$ of non-zero-section projective vortices on $F:=E^*\otimes K_X$ with the moduli space of type $2$ non-zero-section non-Abelian monopoles $\sM_{\ft,2}^0$ give real analytic isomorphisms
\[
j_{\ft,1}: \sM^0(E,h,\omega) \to \sM_{\ft,1}^0
\quad\text{and}\quad
j_{\ft,2}: \sM^0(F,h,\omega) \to \sM_{\ft,2}^0.
\]
Both $j_{\ft,1}$ and $j_{\ft,2}$ are $S^1$-equivariant with respect to the actions defined by scalar multiplication on sections and both maps preserve the $L^2$ norms of the sections.

We shall temporarily write the Hamiltonian function defined in \eqref{eq:Hitchin_function} as $f_\ft:\sM_\ft^0\to \RR$ and write the Hamiltonian functions implied by \eqref{eq:Hamiltonian_function} as $f_E:\sM^0(E,h,\omega)\to\RR$ and $f_F:\sM^0(F,h,\omega)\to \RR$. Because $j_{\ft,1}$ and $j_{\ft,2}$ preserve the $L^2$-norms of the sections, we obtain $f_E=f_\ft\circ j_{\ft,1}$ and $f_F=f_\ft\circ j_{\ft,2}$. Thus, if $[A,\Phi]=j_{\ft,1}([A,\varphi_1])$ and $[A,\Phi]$ is a critical point of $f_\ft$, then $[A,\varphi_1]$ is a critical point of $f_E$. Theorem \ref{thm:Critical_points_Hitchin_Hamiltonian_function_moduli_space_projective_vortices} implies that $[A,\varphi_1]$ is a fixed point of the $S^1$ action. By the $S^1$-equivariance of $j_{\ft,1}$, we see that $[A,\Phi]$ is a fixed point of the $S^1$ action on $\sM_\ft^0$. Theorem \ref{mainthm:IdentifyCriticalPoints} follows from the characterization of fixed points of the $S^1$ action in Proposition \ref{prop:FixedPointsOfS1ActionOnSpinuQuotientSpace} and Feehan and Leness \cite[Lemma 3.12, p. 95]{FL2a}. A similar argument proves Theorem \ref{mainthm:IdentifyCriticalPoints} when $[A,\Phi]=j_{\ft,2}([A,\varphi_2])$.
\end{proof}

\chapter{Real analytic isomorphism of moduli spaces}
\label{chap:Analytic_isomorphism_moduli_spaces}
Our primary goal in this chapter is to prove Theorem \ref{thm:Lubke_Teleman_6-3-7}, and hence Theorem \ref{thm:Lubke_Teleman_6-3-10} as a corollary of Theorems \ref{thm:Lubke_Teleman_6-3-7} and \ref{thm:HitchinKobayashiCorrespondenceForPairs} and Lemmas \ref{lem:SO3_monopole_equations_almost_Kaehler_manifold} and \ref{lem:Okonek_Teleman_1995_3-1}. While Theorem \ref{thm:Lubke_Teleman_6-3-7} is stated by L\"ubke and Teleman as \cite[Theorem 6.3.7, p. 72]{Lubke_Teleman_2006}, it is deduced as a corollary of a universal Hitchin--Kobayashi correspondence developed by them in \cite{Lubke_Teleman_2006}. However, we believe that it is valuable to provide a direct, self-contained proof and we include one in this chapter modeled on that of \cite[Section 4.3.4, Theorem 3.9, p. 329]{FrM} due to Friedman and Morgan, a special case of Theorem \ref{thm:Kobayashi_7_4_20}. Friedman and Morgan attribute some of their arguments to Douady \cite{Douady_1966_sem_bourbaki, Douady_1966aifg}, who is also a source for Miyajima \cite{Miyajima_1989}.

In Section \ref{sec:Real_complex_analytic_functions_Banach_spaces}, we review the construction of complex analytic spaces associated to Fredholm maps of complex Banach spaces. In Section \ref{sec:Slice_action_complex_gauge_transformations_Kuranishi_model_moduli_space_holomorphic_pairs}, we present a slice theorem for the action of the group of complex gauge transformations on the affine space of $(0,1)$-pairs and construct Kuranishi models for open neighborhoods of points in the moduli space of holomorphic pairs.  We prove the existence of a local moduli functor for holomorphic pairs in Section \ref{subsec:Friedman-Morgan_4-2} and discuss the global moduli functor in Section \ref{subsec:Friedman-Morgan_4-4-1}. In Section \ref{subsec:Friedman-Morgan_4-3-4} we prove the existence of a real analytic isomorphism between Kuranishi models for open neighborhoods of points in the moduli spaces of projective vortices and holomorphic pairs, respectively. We use this comparison of Kuranishi models to construct a virtual moduli space of projective vortices as a complex K\"ahler manifold in Section \ref{sec:Complex_Kaehler_structure_moduli_space_projective_vortices_near_singular_points}.
We conclude in Section \ref{sec:Proofs_theorem_Lubke_Teleman_6-3-7_and_corollary} by proving Theorem \ref{thm:Lubke_Teleman_6-3-7}, which gives real analytic embeddings of the moduli space of non-zero-section projective vortices into the moduli space of holomorphic pairs, and proving Corollary \ref{cor:Lubke_Teleman_6-3-7_and_Kobayashi_7_3_17_pair}.

\section[Complex analytic space associated to a Fredholm map]{Real and complex analytic functions on Banach spaces and the complex analytic space associated to a Fredholm map}
\label{sec:Real_complex_analytic_functions_Banach_spaces}
While much of the background material discussed by Friedman and Morgan in \cite[Sections 4.1.3 and 4.1.4]{FrM} is by now standard, some of the concepts that they discuss are less well-known and so we shall quote their expositions here and refer to \cite[Section 4.1.3 and 4.1.4]{FrM} for more standard material and, also following those authors, refer to the reader to Douady \cite{Douady_1966aifg} for additional details.

Let $H_1$ and $H_2$ be Banach spaces over $\KK=\RR$ or $\CC$. We refer to \cite[Section 4.1.3]{FrM} for the definition of the Banach space of bounded homogeneous polynomials of degree $n$ from $H_1$ to $H_2$ and its identification with the Banach space of all symmetric bounded multilinear maps from $H_1$ to $H_2$, along with the definition of analytic maps in terms of convergent power series, called \emph{holomorphic maps} when $\KK=\CC$.

\begin{lem}
\label{lem:Friedman_Morgan_4-1-8}  
(See Friedman and Morgan \cite[Section 4.1.3, Lemma 1.8, p. 286]{FrM}.)  
Let $B \subset H_1$ be a ball of radius $R$ about the origin, and let $f: B \to H_2$ be a bounded continuous function such that, for all linear maps $\varphi: \CC \to H_1$ and all elements $\alpha \in H_2^*$, where $H_2^*$ is the continuous dual space of bounded linear functions on $H_2$, the map $\alpha\circ f\circ\varphi$ is a holomorphic function from $\varphi^{-1}(B) \subset \CC$ to $\CC$. Then $f$ is analytic and its radius of convergence is at least $R$. In particular, if $f$ is complex differentiable at each point of $B$, then it is
analytic.
\end{lem}

Lemma \ref{lem:Friedman_Morgan_4-1-8} is used to prove the

\begin{thm}[Inverse and implicit mapping theorems for complex analytic maps of complex Banach spaces]
\label{thm:Friedman_Morgan_4-1-9}    
(See Friedman and Morgan \cite[Section 4.1.3, Theorem 1.9, p. 287]{FrM}.)  
Let $H$ and $H'$ be Banach spaces over $\KK$, and let $U$ be an open subset of $H$.
\begin{enumerate}
\item Suppose that $F: U \to H'$ is an analytic map and that the differential of $F$ is an isomorphism at a point $z \in U$. Then $F(U)$ contains an open set $V$ about the point $w = F(z)$, and there exists an analytic map $G:V\to U$ such that $G(w) = z$ and such that $G\circ F = \id_H$ and $F\circ G = \id_{H'}$.

\item Suppose that $F: U \to H'$ is an analytic map and that $H$ is a topological direct sum $H = H_1 \oplus H_2$. Suppose that $z = z_1 + z_2 \in U$, where $z_i \in H_i$ for $i=1,2$, that $F(z) = 0$, and that the differential of $F$ restricted to $H_2$ is an isomorphism from $H_2$ to $H'$. Then there are open neighborhoods $U_1 \subset U\cap H_1$ of $z_1$ and $U_2 \subset U\cap H_2$ of $z_2$ and an analytic map $\varphi: U_1 \to U_2$, such that the graph of $\varphi$ is contained in $U$ and such that $(U_1 \times U_2) \cap F^{-1}(0)$ is the graph of $\varphi$.
\end{enumerate}
\end{thm}

We also refer the reader to Berger \cite[Definition 2.3.1]{Berger_1977}, Deimling \cite[Definition 15.1]{Deimling_1985}, or Zeidler \cite[Definition 8.8]{Zeidler_nfaa_v1} for definitions of analytic maps of Banach spaces and to Whittlesey \cite{Whittlesey_1965} for a development of their properties. For other expositions of the proof of Theorem \ref{thm:Friedman_Morgan_4-1-9}, we note that statements and proofs of the Inverse Mapping Theorem for analytic maps of Banach spaces are provided by Berger \cite[Corollary 3.3.2]{Berger_1977} (complex), Deimling \cite[Theorem 4.15.3]{Deimling_1985} (real or complex), and Zeidler \cite[Corollary 4.37]{Zeidler_nfaa_v1} (real or complex). The corresponding Implicit Mapping Theorems for analytic maps can be proved in the standard way as corollaries, for example \cite[Theorem 2.5.7, p. 121]{AMR} and \cite[Theorem 4.H]{Zeidler_nfaa_v1}.

We now turn to a summary of some less well-known concepts discussed in \cite[Section 4.1.3]{FrM} that are required in order to understand and extend Kuranishi's Theorem \cite{Kuranishi} on the deformation of complex structures on a closed, complex manifold to the case of complex structures on a complex vector bundle over a closed, complex manifold. For this purpose, we shall need the concept of a \emph{complex analytic space} and, for the purpose of the proof of Theorem \ref{thm:Lubke_Teleman_6-3-7}, that of a \emph{real analytic space}.

There are several references for \emph{complex} analytic spaces, \label{page:Complex_analytic_space} including Abhyankar \cite{Abhyankar_local_analytic_geometry}, Aroca, Hironaka, and Vicente \cite{Aroca_Hironaka_Vicente_complex_analytic_desingularization} (based on the earlier three-volume series by Aroca, Hironaka, and Vicente \cite{Hironaka_infinitely_near_singular_points, Aroca_Hironaka_Vicente_desingularization_theorems, Aroca_Hironaka_Vicente_theory_maximal_contact}), Bierstone and Milman \cite{Bierstone_Milman_1997}, \cite[Section 2]{Bierstone_Milman_1989} (and references therein), Chirka \cite{Chirka_1989}, Fischer \cite{Fischer_complex_analytic_geometry}, Grauert and Remmert \cite{Grauert_Remmert_coherent_analytic_sheaves}, Griffiths and Harris \cite{GriffithsHarris}, Gunning and Rossi \cite{Gunning_Rossi_analytic_functions_several_complex_variables}, and Narasimhan \cite{Narasimhan_introduction_theory_analytic_spaces}.

For the analogous development of \emph{real} analytic spaces\label{page:Real_analytic_space},  references are more limited and not widely available, although they include Hironaka \cite{Hironaka_intro_real-analytic_sets_maps, Hironaka_1973} and Guaraldo, Macr\`\i, and Tancredi \cite{Guaraldo_Macri_Tancredi_topics_real_analytic_spaces}, together with Chriestenson \cite{ChriestensonThesis} and references therein. A few references such as \cite{Guaraldo_Macri_Tancredi_topics_real_analytic_spaces} consider analytic spaces over real or complex analytic spaces simultaneously whenever possible while Onishchik \cite{Onishchik} allows $\KK$ to be a complete non-discretely normed field of characteristic zero, essentially the view taken as well by Bierstone and Milman \cite{Bierstone_Milman_1997}.

Let $T$ be a (finite-dimensional) complex analytic space or real analytic ringed space. (See Grauert and Remmert \cite[Section 1.1.5, p. 7]{Grauert_Remmert_coherent_analytic_sheaves} for the definition of a complex analytic space.) We now define morphisms from $T$ into a Banach space. As in Friedman and Morgan \cite[Section 4.1.3]{FrM}, which we follow very closely, we only consider $\KK=\CC$ and note that the modifications required for $\KK=\RR$ are straightforward.

We first define a local model for a holomorphic map from $T$ to $H$ where $H$ is a complex Banach space.
Locally, a finite-dimensional complex analytic space $T$ is a complex analytic subspace\footnote{These are usually just called \emph{complex spaces} --- see Friedman and Morgan \cite[p. 11]{FrM} or Grauert and Remmert \cite[Section 1.1.5, p. 7]{Grauert_Remmert_coherent_analytic_sheaves} --- but we shall need to consider both complex analytic and real analytic spaces in our applications, thus $\KK$-analytic spaces where $\KK=\RR$ or $\CC$.}  of an open set $B \subset \CC^N$ for some integer $N$. 
Given such an embedding $T \subset B$, a \emph{local model} for a holomorphic map $f$ is the restriction of a holomorphic map $g$ from an open neighborhood of $T$ in $B$ to $H$. Two maps $g_1: B \to H$ and $g_2: B \to H$ are \emph{equal} on $T$ if $(g_1-g_2)^*\sO_H \subset \sI_T$, where $\sO_H$ is the sheaf of germs of holomorphic functions on $H$ (see Grauert and Remmert \cite[Section 1.1.2, p. 3]{Grauert_Remmert_coherent_analytic_sheaves}) and $\sI_T$ is the ideal sheaf of $T$ in $B$. (Thus, for each $t\in T$ there are an open neighborhood $U \subset \CC^N$ of $t$ (with $U\subset B$) and holomorphic functions $f_1,\ldots, f_k \in \sO(U)$ such that (see Grauert and Remmert \cite[pp. 4, 7, and 224]{Grauert_Remmert_coherent_analytic_sheaves})
\[
  (\sI_T)_U = \sI_T\restriction U = \sO_U f_1 + \cdots + \sO_U f_k,
\]
and therefore, locally,
\[
  T\cap U = \{z\in U: f_1(z) = \cdots = f_k(z) = 0\},
\]
while $T = \supp(\sO_B/\sI_T)$, globally.) Hence, there is a well-defined morphism of ringed spaces from $T$ to $H$. Thus, given two embeddings $T \subset B_i$, for $i = 1,2$, if $f: T \to H$ is locally the restriction of a holomorphic map $g_1: B_1 \to H$, then it is also locally the restriction of a holomorphic map $g_2: B_2 \to H$, and moreover there is a holomorphic map $\alpha:B_1 \to B_2$, inducing the given embedding on $T$, and such that $g_1$ and $g_2\circ\alpha$ are equal on $T$. (This statement can also be used as a definition of when two maps $g_i: B_i \to H$ for $i=1,2$ are equal on $T$.) Thus, two maps are equal on $T$ if and only if they induce the same morphism of ringed spaces from $T$ to $H$. Using this definition of equality, one can thus define a global holomorphic map from $T$ to $H$; locally it should be of the form described above.

We see that composition of holomorphic maps is well-defined as follows. Suppose that $g_1: B \to H$ and $g_2: B \to H$ are two holomorphic maps such that $g_1$ and $g_2$ are equal on $T$. Let $H'$ be a complex Banach space and  $F: H \to H'$  a holomorphic map defined on a neighborhood of $\Imag g_1 \cup \Imag g_2$. Then $F \circ g_1$ is equal to $F \circ g_2$ on $T$. Thus if $f$ denotes the restriction of $g_i$ to $T$, then $F \circ f$ is well-defined.

\begin{defn}[Containment in the sense of complex analytic spaces]
\label{defn:Friedman_Morgan_page_288}
(See Friedman and Morgan \cite[Section 4.1.3, p. 288]{FrM}.)
Suppose that $f: T \to H$ is a morphism from the finite-dimensional complex analytic space $T$ to the complex Banach space $H$, and that $F: U \to H'$ is a holomorphic map from a neighborhood $U$ of $\Imag f$ to another complex Banach space $H'$. One says that $F^{-1}(0)$ \emph{contains} $\Imag f$ \emph{in the sense of complex analytic spaces} if the following holds: Locally, there is an embedding of $T$ into an open subset $B$ of $\CC^N$ and a holomorphic map $g: B \to H$ such that the restriction of $g$ to $T$ is $f$, and we require that, for every $h$ in the maximal ideal of $0$ on $H'$, that is, for every holomorphic function $h$ defined on an open neighborhood of $0$ and vanishing at $0$, the composition $h \circ F \circ g$ lies in $\sI_T$. Equivalently, the map $F \circ f$ is equal to zero on $T$.
\end{defn}

The preceding definition does not depend on the choice of an embedding of $T$ into $B$ or of the extension of $f$ to a map $g$ on $B$. There are various possible definitions for the vanishing of $F$ on $\Imag f$, but when $T$ is finite-dimensional, they all agree:

\begin{lem}
\label{lem:Friedman_Morgan_4-1-10}  
(See Friedman and Morgan \cite[Section 4.1.3, Lemma 1.10, p. 288]{FrM}.)  
Let $g: B \to H$ be an analytic map from an open subset $B$ of $\CC^N$ to the Banach space $H$. Let $F: H \to H'$ be an analytic map from an open neighborhood of $\Imag g$ into a Banach space $H'$. Then the following three ideal sheaves in $\sO_B$ are equal:
\begin{enumerate}
\item The ideal sheaf generated by pullbacks of the form $h \circ F \circ g$, where $h$ is a bounded, linear function on $H'$.

\item The ideal sheaf generated by pullbacks of the form $h \circ F \circ g$, where $h$ is a holomorphic function defined in a neighborhood of $0 \in H'$ and vanishing at $0$.

\item The ideal sheaf generated by all functions of the form $\langle \alpha, F\circ g\rangle$, where $\alpha$ is a holomorphic function from an open set in $B$ to $(H')^*$, the dual space of bounded linear functions on $H'$, and $\langle \cdot, \cdot \rangle$ denotes the natural pairing between $H'$ and $(H')^*$.
\end{enumerate}
\end{lem}

The following lemma is used frequently by Friedman and Morgan in \cite[Sections 4.2 and 4.3]{FrM}:

\begin{lem}
\label{lem:Friedman_Morgan_4-1-11}  
(See Friedman and Morgan \cite[Section 4.1.3, Lemma 1.11, p. 289]{FrM}.)  
Let $f: T \to H$ be a morphism from the complex analytic space $T$ to the complex Banach space $H$. Suppose that $H = H_1 \oplus H_2$, where the $H_i$ are closed subspaces for $i=1,2$. Finally suppose that the image of $f$ is contained in $\pi_2^{-1}(0)$ in the sense of complex analytic spaces, where $\pi_i: H \to H_i$ is the projection for $i=1,2$. Then $f = \iota_1 \circ \pi_1 \circ f$, where $\iota_i: H_i \hookrightarrow H$ is the inclusion for $i=1,2$. In other words, the morphism $f$ factors through the inclusion of $H_1$ in $H$.
\end{lem}

\begin{defn}[Germ of a set at a point]
\label{defn:Germ_set}
(See Huybrechts \cite[Definition 1.1.21, p. 18]{Huybrechts_2005}.)
Suppose that $A_1, A_2 \subset \CC^N$ are subsets. One says that $A_1$ and $A_2$ have the \emph{same germ at the origin} if there exists an open neighborhood $U$ of the origin such that $A_1 \cap U = A_2 \cap U$.
The \emph{germ of a set $A$ in $\CC^N$ at the origin}, often denoted $(A,0)$, is the equivalence class of $A$ under the preceding relation. 
\end{defn}

If $A_1, A_2$ in Definition \ref{defn:Germ_set} are complex analytic sets (see Grauert and Remmert \cite[Section 4.1.1]{Grauert_Remmert_coherent_analytic_sheaves}), they define the same analytic germ at the origin if the topological subsets $A_1, A_2$ define the same at the origin as in Definition \ref{defn:Germ_set}. The
\label{page:Germ_complex_analytic_set}
\emph{germ of a complex analytic set $A$ at the origin in $\CC^N$} is the equivalence class of $A$ under the relation in Definition \ref{defn:Germ_set}. One may define the
\label{page:Germ_complex_analytic_space}
\emph{germ of a complex analytic space at a point} in exactly the same way as in Definition \ref{defn:Germ_set} for the germ of an analytic set in $\CC^N$ at the origin. Lastly, we recall the construction of the complex analytic space defined by a Fredholm map. 

\begin{thm}
\label{thm:Friedman_Morgan_4-1-13}    
(See Friedman and Morgan \cite[Section 4.1.4, Theorem 1.13, p. 290]{FrM}.)  
Let $H_1$ and $H_2$ be two complex Banach spaces, $U$ an open neighborhood of the origin in $H_1$, and $F:U \to H_2$ a holomorphic map. Suppose that $F(0) = 0$ and that the differential $dF(0)$ of $F$ at the origin is Fredholm.
\begin{enumerate}
\item There is a natural structure of the germ of a finite-dimensional complex analytic space $T$ on the germ of $F^{-1}(0)$ at the origin.
  
\item With the complex structure given in the first item, $T$ is contained in $F^{-1}(0)$ in the sense of complex analytic spaces.

\item The germ $T$ has the following universal property: If $S$ is the germ of a finite-dimensional complex analytic space and there is a holomorphic map $f: S \to H_1$ which is contained in $F^{-1}(0)$ in the sense of complex spaces, then there is a unique map $\tilde f: S \to T$ such that the map $f = \iota \circ \tilde f$, where $\iota:T\to H_1$ is the inclusion map.
\end{enumerate}
\end{thm}

In \cite[Section 4.1.4, Theorem 1.14, p. 292]{FrM}, Friedman and Morgan provide a generalization of Theorem \ref{thm:Friedman_Morgan_4-1-13} from the case of Fredholm maps of Banach spaces to the case of Fredholm sections of Banach vector bundles. We shall need this generalization for our proof of Theorem \ref{thm:Friedman_Morgan_4-3-8_projective_vortices} and so we quote their result here.

Following Friedman and Morgan \cite[Section 4.1.4, p. 292]{FrM}, we remark that a holomorphic map $F: H_1 \supset U \to H_2$ (as in Theorem \ref{thm:Friedman_Morgan_4-1-13}) may be viewed as a section of the product Banach space bundle $H_2 \times U \to U$. More generally, one can consider an arbitrary Banach vector bundle over a neighborhood $U$ of the origin in $H_1$. Since we are only concerned with germs near the origin, we may assume that this vector bundle is the trivial bundle $\pi: H_2 \times U \to U$. Let $s$ be an  analytic section of $H_2 \times U\to U$ that is Fredholm at the origin. That is, if we write $s = (F,\id_U)$, then $s$ is Fredholm at the origin if and only if $dF(0)$ is Fredholm. Consequently, $\Ker ds(0)$ is finite dimensional and there is a finite-dimensional vector subspace $V_0$ of the fiber $\pi^{-1}(0) = \Ker d\pi(0) = H_2$ of the bundle over $0 \in U \subset H_1$ such that the projection of $\Ran ds(0)$ to $H_2$ is closed and the linear span of this range with $V_0$ is equal to $H_2$.  Suppose now that we are given a finite-dimensional analytic vector subbundle $V \subset H_2\times U$ with the property that the linear span of the fiber of $V$ over $0 \in U$ and $\Ran ds(0) \subset H_2$ is equal to $H_2$, that is, $V$ is transverse to $s$ over $0\in U$. After possibly shrinking $U$, standard methods imply that $s^{-1}(V) = B$ is a finite-dimensional analytic submanifold of $U$. By the construction of $B$, the restriction of $s$ to $B$ defines a section of $V\restriction B$, which we shall continue to denote by $s$. Thus, $s$ is a section of the finite-dimensional analytic vector bundle $V\restriction B$ over $B$, and so there is a natural structure of a complex analytic space $T$ on $s^{-1}(0)$.

\begin{thm}
\label{thm:Friedman_Morgan_4-1-14}    
(See Friedman and Morgan \cite[Section 4.1.4, Theorem 1.14, p. 292]{FrM}.)  
Assume the notation of the preceding paragraph.
\begin{enumerate}
\item The structure of the germ of a finite-dimensional complex space $T$ on the germ of $s^{-1}(0)$ at $0$ is independent of the choice of the finite-dimensional subbundle $V$.
  
\item With the complex structure given in the preceding item, $T$ is contained in $s^{-1}(0)$ in the sense of complex analytic spaces.
  
\item The germ $T$ has the following universal property: If $S$ is the germ of a finite-dimensional complex analytic space and there is a holomorphic map $f: S \to U$ which is contained in $s^{-1}(0)$ in the sense of complex analytic spaces, then there is a unique map $\tilde f: S \to T$ such that $f = \iota \circ \tilde f$, where $\iota:T \to H_1$ is the inclusion map.
\end{enumerate}
\end{thm}

\section[Local Kuranishi model for moduli space of holomorphic pairs]{Slice for the action of the group of complex gauge transformations and local Kuranishi model for the moduli space of holomorphic pairs}
\label{sec:Slice_action_complex_gauge_transformations_Kuranishi_model_moduli_space_holomorphic_pairs}
In this section, we adapt the development in \cite[Sections 4.1.5 and 4.1.6]{FrM} due to Friedman and Morgan of the local Kuranishi model for the moduli space of holomorphic structures on complex vector bundles to the case of holomorphic pairs. We assume the setting of Section \ref{sec:Elliptic_deformation_complex_holomorphic_pair_equations}, so $E$ is a complex vector bundle over a complex manifold $X$ of complex dimension $n$, and furthermore we require $X$ to be closed.

We review the relevant groups of gauge transformations in Section \ref{subsec:Group_complex_gauge_transformations} and discuss the notion of a slice and of a local slice for the action of a Lie group on a smooth manifold in Section \ref{subsec:Slice_action_Lie_group_smooth_manifold}.  We construct a local slice for the action of the group of complex, determinant-one gauge transformations on the affine space of $(0,1)$-pairs in
Section \ref{subsec:Coulomb_gauge_slice_condition_for_01-pairs}.  In Section \ref{subsec:Kuranishi_model_moduli_space_holomorphic_pairs}, we construct a Kuranishi model associated to a point in the moduli space of holomorphic pairs and prove that it is independent of the choice of slice.

\subsection{Groups of complex gauge transformations}
\label{subsec:Group_complex_gauge_transformations}
Let $E$ be a smooth complex vector bundle over a closed, complex manifold $X$. Recall that complex gauge transformations are defined as sections of the bundles defined in \eqref{eq:AutomorphismBundles},
\begin{align*}
  \GL(E)& := \{u\in \End(E): u(x)\in\GL(E_x), \text{ for all } x\in X\},
\\
  \SL(E)& := \{u\in \End(E): u(x)\in\SL(E_x), \text{ for all } x\in X\},
\end{align*}
both smooth principal fiber bundles over $X$ with structure groups $\GL(r,\CC)$ and $\SL(r,\CC)$, respectively, when $E$ has complex rank $r$. Recall that $\GL(r,\CC)$ is a \emph{complex analytic group} (see Lee \cite[Sections 1.2 and 1.3]{Lee_structure_complex_lie_groups}) and $\SL(r,\CC)$ is a connected, \emph{complex Lie subgroup} of $\GL(r,\CC)$ (see Lee \cite[Sections 1.2 and 1.3]{Lee_structure_complex_lie_groups}). The exponential map $\exp:\fg \to G$ of a complex analytic group is complex analytic (see Lee \cite[Theorem 1.15]{Lee_structure_complex_lie_groups}) and surjective when $G=\GL(r,\CC)$ (see Hall \cite[Exercises 2.9 and 2.10]{Hall_lie_groups_algebras_representations}).

By analogy with Freed and Uhlenbeck \cite[Appendix A]{FU}, we consider the groups of smooth and $W^{2,p}$ sections of $\GL(E)$ and $\SL(E)$, denoted by $C^\infty(\GL(E))$ and $W^{2,p}(\GL(E))$, respectively, and similarly for $\SL(E)$, where $p\in(n,\infty)$ if $X$ has complex dimension $n$. The complex analytic exponential map $\exp:\gl(r,\CC)\to\GL(r,\CC)$ induces maps
\[
  \Exp:\Omega^0(\gl(E)) \to C^\infty(\GL(E))
  \quad\text{and}\quad
  \Exp:W^{2,p}(\gl(E)) \to W^{2,p}(\GL(E)),
\]
and similarly for $\exp:\fsl(r,\CC)\to\SL(r,\CC)$. Here, $\Omega^0(\gl(E)) = C^\infty(\gl(E))$, a complex vector space of $C^\infty$ sections of the smooth complex vector bundle $\gl(E)$ over $X$, while $W^{2,p}(\gl(E))$ is a complex Banach space of $W^{2,p}$ sections of the smooth complex vector bundle $\gl(E)$ over $X$. The proof due to Freed and Uhlenbeck of their \cite[Proposition A.2]{FU} extends to show that the exponential map $\Exp$ is complex analytic on an open neighborhood of the origin in $W^{2,p}(\gl(E))$ and that $W^{2,p}(\GL(E))$ is a (Banach) \emph{complex Lie group} in the sense of Lee \cite[Section 1.2]{Lee_structure_complex_lie_groups}, with complex Lie algebra $W^{2,p}(\gl(E))$. Similarly, $W^{2,p}(\SL(E))$ is a (Banach) complex Lie group, with complex Lie algebra $W^{2,p}(\fsl(E))$.

The groups $C^\infty(\GL(E))$ and $C^\infty(\SL(E))$, along with their Sobolev completions, $W^{2,p}(\GL(E))$ and $W^{2,p}(\SL(E))$, are often called \emph{complex gauge groups} or \emph{groups of complex gauge transformations} and denoted by $\sG_E^\CC$ \cite[Section 4.1.2, Definition 1.7, p. 285]{FrM}. Before proceeding to discuss the distinction between slice theorems for the actions of $W^{2,p}(\SU(E))$ on $\sA(E,h)$ versus $W^{2,p}(\SL(E))$ on $\sA^{0,1}(E)$, and similarly for pairs, we recall general facts concerning the existence of slices for the action of Lie groups on smooth manifolds.

\subsection{Slice for the action of a Lie group on a smooth manifold}
\label{subsec:Slice_action_Lie_group_smooth_manifold}
Duistermaat and Kolk \cite[Chapter 2]{DuistermaatLieGroups} provide an exceptionally clear discussion of the action of Lie groups on smooth manifolds and we recall some essential relevant points here. As in \cite[Chapter 2]{DuistermaatLieGroups}, a $C^k$ map of manifolds allows for finite integers $k$, thus $1\leq k<\infty$, or $k=\infty$ (smooth) or $k=\omega$ (analytic).

\begin{lem}
\label{lem:Duistermaat_Kolk_2-1-1}  
(See Duistermaat and Kolk \cite[Lemma 2.1.1]{DuistermaatLieGroups}.)  
Let $\rho: G \times M \to M$ be a $C^k$ action (for $k \geq 1$) of a Lie group $G$ on a smooth manifold $M$. For a point $x_0 \in M$, let $S \subset M$ be an embedded $C^k$ submanifold through $x_0$ such that
\begin{equation}
  \label{eq:Duistermaat_Kolk_2-1-1}
  T_{x_0}M = \Ran D_1\rho(\id,x_0) \oplus T_{x_0}S,
\end{equation}
where $\id\in G$ is the identity element and $D_1\rho(\id,x_0):\fg \to T_{x_0}M$ denotes the partial derivative of $\rho$ in the $G$ direction, and
$\fg = T_{\id}G$ is the Lie algebra of $G$. Let $C \subset G$ be an embedded $C^k$ submanifold through $\id_G$ such that
\begin{equation}
  \label{eq:Duistermaat_Kolk_2-1-2}
  \fg = \fg_{x_0} \oplus T_{\id}C,
\end{equation}
where $\fg_{x_0} \subset \fg$ is the Lie algebra of the stabilizer subgroup $G_{x_0}\subset G$ of the point $x_0$. Then there are an open neighborhood $C_0 \subset C$ of $\id_G$, an open neighborhood $S_0 \subset S$ of $x_0$, and an open neighborhood $M_0 \subset M$ of $x_0$ such that
\[
  \rho_0 := \rho\restriction (C_0\times S_0) \to M
\]
is a $C^k$ diffeomorphism from $C_0\times S_0$ onto $M_0$.
\end{lem}

Recall that $G_{x_0} := \{g\in G: \rho(g,x_0)=x_0\} = \Stab(x_0)$ is the \emph{stabilizer} subgroup of $G$ for the point $x_0$ and is a Lie subgroup of $G$ \cite[p. 94]{DuistermaatLieGroups}.

\begin{defn}[Slice for the action of a Lie group on a smooth manifold]
\label{defn:Duistermaat_Kolk_2-3-1}    
(See Duistermaat and Kolk \cite[Definition 2.3.1]{DuistermaatLieGroups}.)  
Let $\rho: G \times M \to M$ be a $C^k$ action (for $k \geq 1$) of a Lie group $G$ on a smooth manifold $M$. A $C^k$ \emph{slice} at $x_0 \in M$ for the action $\rho$ is an embedded $C^k$ submanifold $S \subset M$ through $x_0$ such that
\begin{enumerate}
\item
\label{item:Duistermaat_Kolk_2-3-1_SplitTangentSpace}
$T_{x_0}M = \Ran D_1\rho(\id,x_0) \oplus T_{x_0}S$ and $T_xM = \Ran D_1\rho(\id,x) + T_xS$ for all $x\in S$.
\item
\label{item:Duistermaat_Kolk_2-3-1_StabInvariant}
$S$ is $G_{x_0}$-invariant.
\item
\label{item:Duistermaat_Kolk_2-3-1_GlobalTube}
If $x\in S$ and $g \in G$ and $\rho(g,x) \in S$, then $g \in G_{x_0}$.
\end{enumerate}
\end{defn}

It follows that the inclusion map, $S \to M$, induces a bijective mapping, and thus a homeomorphism,
\[
  S/G_{x_0} \ni G_{x_0}\cdot x \mapsto G\cdot x \in M/G
\]
from the quotient space $S/G_{x_0}$ onto an open neighborhood of $G\cdot x_0$ in the quotient space $M/G$. The action of $G_{x_0}$ on $S$ necessarily has $x_0$ as a fixed point.

\begin{rmk}[Local slice for the action of a Lie group on a smooth manifold]
\label{rmk:Duistermaat_Kolk_2-3-1_local_slice}  
(Compare Haydys \cite[Section 6.6, Definition 178]{Haydys_introduction_to_gauge_theory} or Morgan \cite[Remark 4.5.6, p. 62]{MorganSWNotes}.)
We call $S$ in Definition \ref{defn:Duistermaat_Kolk_2-3-1} a \emph{local slice} if it obeys Conditions \eqref{item:Duistermaat_Kolk_2-3-1_SplitTangentSpace} and \eqref{item:Duistermaat_Kolk_2-3-1_StabInvariant}, but not necessarily Condition \eqref{item:Duistermaat_Kolk_2-3-1_GlobalTube}.   
\end{rmk}  

The quotient topology on $M/G$ is the one for which $V$ is an open subset of $M/G$ if and only if the $G$-invariant subset $\pi^{-1}(V)$ is open in $M$, where $\pi:M\to M/G$ is the canonical projection. In general, the quotient topology need not be Hausdorff. However, the quotient topology on $M/G$ is Hausdorff if and only if the \emph{orbit relation} (or \emph{graph}) $\{(x, y) \in M \times M: y \in G\cdot x\}$ is a closed subset of $M \times M$ (see Duistermaat and Kolk \cite[Lemma 1.11.2]{DuistermaatLieGroups}). An action $\rho:G\times M\to M$ is said to be \emph{proper} if the induced map
\[
  G\times M \ni (g,x) \mapsto (g\cdot x, x) \in M\times M
\]
is proper and, according to Duistermaat and Kolk \cite[Lemma 1.11.3]{DuistermaatLieGroups}, a proper and continuous action $\rho$ yields a Hausdorff quotient space $M/G$.

\begin{defn}[Proper action of a Lie group at a point]
\label{defn:Duistermaat_Kolk_2-3-2}    
(See Duistermaat and Kolk \cite[Definition 2.3.2]{DuistermaatLieGroups}.)
Continue the notation of Definition \ref{defn:Duistermaat_Kolk_2-3-1}. The action $\rho$ is said to be \emph{proper at $x_0$} if for every sequence $\{x_j\} \subset M$ and $\{g_j\} \subset G$ such that $\lim_{j\to\infty} x_j = x_0$ and $\lim_{j\to\infty} g_j\cdot x_j = x_0$, there is a subsequence $\{j(k)\} \subset \{j\}$ such that $\{g_{j(k)}\}$ converges in $G$ as $k\to\infty$.
\end{defn}

A Lie group $G$ acts properly on a manifold $M$ if and only if the action is proper at each point $x \in M$ and the topology of the quotient space $M/G$ is Hausdorff (see Duistermaat and Kolk \cite[pp. 103--104 and Lemma 1.11.3]{DuistermaatLieGroups}).

\begin{thm}[Existence of a slice for the action of a Lie group on a smooth manifold]
\label{thm:Duistermaat_Kolk_2-3-3}    
(See Duistermaat and Kolk \cite[Theorem 2.3.3]{DuistermaatLieGroups}.)
Let $\rho$ be a $C^k$ action ($k \geq 1$) of a Lie group $G$ on a smooth manifold $M$, and suppose that the action is proper at $x_0 \in M$. Then there exists a $C^k$ slice $S$ at $x_0$ for the action $\rho$.
\end{thm}

\subsection{Coulomb gauge for the action of the group of complex gauge transformations on the affine space of (0,1) pairs}
\label{subsec:Coulomb_gauge_slice_condition_for_01-pairs}
In \cite[Theorems 16 and 17, p. 18, and Corollary 18, p. 19]{Feehan_Maridakis_Lojasiewicz-Simon_coupled_Yang-Mills}, when $(E,h)$ is a Hermitian vector bundle over a closed, smooth Riemannian manifold $X$ and $\det (E,h)$ has a fixed unitary connection, Feehan and Maridakis proved a Coulomb gauge slice theorem for the action of $W^{2,p}(\SU(E))$ on the affine space of pairs, $\sA(E,h)\times W^{1,p}(E)$, where
$\sA(E,h) = A_0 + W^{1,p}(T^*X\otimes \su(E))$, and $p\in (d/2,\infty)$ is a constant, $X$ has real dimension $d\geq 2$, and $A_0$ is a fixed, smooth unitary connection on $E$.

We shall need an analogue of this result for the action of $W^{2,p}(\SL(E))$ on the affine space of $(0,1)$-pairs, $\sA^{0,1}(E)\times W^{1,p}(E)$. Recall from the definition in \eqref{eq:H_dbar_Avarphi^0bullet} of the harmonic representative space of the elliptic complex of holomorphic pairs \eqref{eq:Holomorphic_pair_elliptic_complex} that
\[
  \bH_{\bar\partial_E,\varphi}^0
  =
  \Ker\left(\bar\partial_{E,\varphi}^0:W^{2,p}(\fsl(E))\to W^{1,p}(\Lambda^{0,1}(\fsl(E)))\right)
\]
is a complex linear subspace and so also is its $L^2$-orthogonal complement.  By Lemma \ref{lem:Stab(rdE,varphi)_is_Lie_Group}, the complex vector space $\bH_{\bar\partial_E,\varphi}^0$ is the Lie algebra of the stabilizer subgroup $\Stab(\bar\partial_E,\varphi)$ of $(\bar\partial_E,\varphi)$ in \eqref{eq:Stabilizer_holomorphic_pair}. It will be convenient to define
\begin{multline}
  \label{eq:Subgroup_W2pSLE_Exponential_orthogonal_complement_kernel}
  \bC_{\bar\partial_E,\varphi}(\delta) := \left\{\Exp\zeta: \|\zeta\|_{W^{2,p}(X)} < \delta \text{ and } \right.
\\
\left. \zeta \in \left(\Ker\left(\bar\partial_{E,\varphi}^0:W^{2,p}(\fsl(E))\to W^{1,p}(\Lambda^{0,1}(\fsl(E)))\right) \right)^\perp\right\} \subset W^{2,p}(\SL(E)),
\end{multline}
for any $(0,1)$-pair $(\bar\partial_E,\varphi)$ on $E$ of class $W^{1,p}$ and constant $\delta \in (0,\infty]$. For small enough $\delta\in(0,1]$, the subset $\bC_{\bar\partial_E,\varphi}(\delta) \subset W^{2,p}(\SL(E))$ is an embedded complex submanifold since $\Exp$ is a complex analytic embedding on a small enough open neighborhood of the origin in $W^{2,p}(\fsl(E))$ while, as noted above, so also is its $L^2$-orthogonal complement. Thus, we have the following analogue of the direct sum decomposition \eqref{eq:Duistermaat_Kolk_2-1-2},
\[
  W^{2,p}(\fsl(E)) = \bH_{\bar\partial_E,\varphi}^0 \oplus T_{\id_E}\bC_{\bar\partial_E,\varphi}(\delta),
\]
noting that $W^{2,p}(\fsl(E)) = T_{\id_E}W^{2,p}(\SL(E))$ is the Lie algebra of the Banach Lie group $W^{2,p}(\SL(E))$.

We can now state the desired analogue of Lemma \ref{lem:Duistermaat_Kolk_2-1-1}  and the well-known Coulomb gauge theorem due to Freed and Uhlenbeck \cite[Theorem 3.2]{FU} for the affine space of connections on $E$ modulo $W^{2,\ell}(\SU(E))$, where $(X,g)$ was a Riemannian four-manifold, $\ell\geq 3$ was an integer, $E$ had complex rank two, and the reference connection $A$ had stabilizer $\Stab(A) \cong \ZZ/2\ZZ$. The result below extends the simpler version in Theorem \ref{thm:Local_Kuranishi_model_for_strongly_simple_point_cP(E)} \eqref{item:cP(E)_cap_UU_dbarEvarphi_biholomorphic_cP(E)_cap_U_dbarEvarphi_times_UidE}, where we assumed that $\Stab(\bar\partial_E,\varphi) = \{\id_E\}$.

\begin{thm}[Existence of $W^{2,p}$ Coulomb gauge transformations for $(0,1)$-pairs of class $W^{1,p}$]
\label{thm:Existence_of_complex_gauge_transformation_to_Coulomb_01_pairs}
Let $(E,h)$ be a Hermitian vector bundle over a closed, complex, Hermitian manifold $(X,g,J)$. Let $\bar\partial_{E_d}$ be a fixed, smooth $(0,1)$-connection on the Hermitian line bundle $\det E$ and $p\in(n,\infty)$ be a constant, where $X$ has complex dimension $n$. If $(\bar\partial_E,\varphi)$ is an $(0,1)$-pair on $E$ of class $W^{1,p}$ that induces $\bar\partial_{E_d}$ on $\det E$, then there are a $\Stab(\bar\partial_E,\varphi)$-invariant open neighborhood of the pair $(\bar\partial_E,\varphi)$,
\[
  \UU_{\bar\partial_E,\varphi}
  \subset
  \sA^{0,1}(E) \times W^{1,p}(E), 
\]
a constant $\delta = \delta(\bar\partial_E,\varphi,g,h,p) \in (0,1]$, and a complex analytic map,
\[
  \UU_{\bar\partial_E,\varphi}
  \ni (\bar\partial_E+\alpha,\varphi+\sigma) \mapsto
  v_{\bar\partial_E+\alpha,\varphi+\sigma} \in \bC_{\bar\partial_E,\varphi}(\delta),
\]
such that, where we abbreviate $v = v_{\bar\partial_E+\alpha,\varphi+\sigma}$,
\begin{equation}  
\label{eq:Coulomb_gauge_slice_condition_for_01-pairs}
v\cdot\left(\bar\partial_E+\alpha,\varphi + \sigma\right)
\\
\in S_{\bar\partial_E,\varphi},
\end{equation}
and $v \in W^{2,p}(\SL(E))$ acts on $\sA^{0,1}(E)\times W^{1,p}(E)$ via \eqref{eq:SL(E)ActionOn(0,1)Pairs}, and we recall from \eqref{eq:dbar_Evarphi_slice} that 
\[
  S_{\bar\partial_E,\varphi} = (\bar\partial_E,\varphi) + \Ker\left(\bar\partial_{E,\varphi}^{0,*}:W^{1,p}(\Lambda^{0,1}(\fsl(E))\oplus E)
\to W^{1,p}(\fsl(E))\right),
\]
the operator $\bar\partial_{E,\varphi}^0$ is as in \eqref{eq:d0StablePair}, and the map
\begin{multline*}
  \UU_{\bar\partial_E,\varphi}
  \ni
  \left(\bar\partial_E+\alpha,\varphi + \sigma\right)
  \mapsto
  \left(v_{\bar\partial_E+\alpha,\varphi+\sigma}\cdot\left(\bar\partial_E+\alpha,\varphi + \sigma\right),
    v_{\bar\partial_E+\alpha,\varphi+\sigma} \right)
  \\
  \in
  U_{\bar\partial_E,\varphi} \times U_{\id_E}
  \subset
  S_{\bar\partial_E,\varphi} \times  \bC_{\bar\partial_E,\varphi}(\delta)
\end{multline*}
is a complex analytic embedding onto the product of an open neighborhood $U_{\bar\partial_E,\varphi}$ of the pair $(\bar\partial_E,\varphi)$ in $S_{\bar\partial_E,\varphi}$ and an open neighborhood $U_{\id_E}$ of the identity $\id_E$ in $\bC_{\bar\partial_E,\varphi}(\delta)$. 
\end{thm}

\begin{rmk}[Dependence of open neighborhood on auxiliary data used to define Sobolev space norms]
\label{rmk:Dependence_constants_on_auxiliary_data_defining_Sobolev_space_norms}  
In addition to their dependence on $\bar\partial_E,\varphi,g,h,p$, the open neighborhood in Theorem \ref{thm:Existence_of_complex_gauge_transformation_to_Coulomb_01_pairs} may also depend, for example, on a choice of smooth connection on $E$ used to define covariant derivatives of sections of $E$. If $\bar\partial_E$ is smooth, then we may choose $A$ to be the smooth Chern connection defined by $\bar\partial_E$ and the Hermitian metric $h$ on $E$. If $\bar\partial_E$ is not smooth, we may choose a nearby $(0,1)$-connection that is smooth and use that and $h$ to define a smooth Chern connection for the purpose of defining Sobolev space norms.
\end{rmk}

\begin{rmk}[Applications of Theorem \ref{thm:Existence_of_complex_gauge_transformation_to_Coulomb_01_pairs}]
In practice, we shall mainly apply Theorem \ref{thm:Existence_of_complex_gauge_transformation_to_Coulomb_01_pairs} to describe a neighborhood of a pair $(\bar\partial_E,\varphi)$ with $\Stab(\bar\partial_E,\varphi) = \{\id_E\}$ and we can replace $\bC_{\bar\partial_E,\varphi}(\delta)$ by $W^{2,p}(\SL(E))$ in the statement.
\end{rmk}  

Theorem \ref{thm:Existence_of_complex_gauge_transformation_to_Coulomb_01_pairs} is proved by using the Implicit Mapping Theorem for complex analytic maps of Banach spaces (for example, see Feehan and Maridakis \cite[Theorem F.1, p. 127]{Feehan_Maridakis_Lojasiewicz-Simon_coupled_Yang-Mills} for the statement and proof of a quantitative version) and, consequently, the gauge transformations $u$ depend complex analytically on the pairs $(\bar\partial_E+\alpha,\varphi + \sigma)$. Recall from \eqref{eq:(0,1)PairQuotientSpaces} that $\sC^{0,1}(E) = (\sA^{0,1}(E)\times W^{1,p}(E))/W^{2,p}(\SL(E))$ is the quotient space of $(0,1)$-pairs of class $W^{1,p}$.

While we might expect an analogue of Feehan and Maridakis \cite[Corollary 18, p. 19]{Feehan_Maridakis_Lojasiewicz-Simon_coupled_Yang-Mills} to hold, namely that $\sC^{0,1;**}(E)$ is a complex analytic Banach manifold, the standard proof that quotient spaces are Hausdorff (see Freed and Uhlenbeck \cite[Corollary, p. 50]{FU}) fails here since the structure group $\SL(r,\CC)$ is \emph{non-compact},  (where $r$ is the complex rank of $E$), unlike the case of $\SU(2)$ considered in \cite{FU} (their argument would extend without significant change from $\SU(2)$ to any compact Lie group $G$) and so the \emph{graph} of the action,
\[
  \Gamma := \left\{\left((\bar\partial_E,\varphi),v\cdot(\bar\partial_E,\varphi)\right):
  (\bar\partial_E,\varphi) \in \sA^{0,1}(E)\times W^{1,p}(E) \text{ and } v \in W^{2,p}(\SL(E) \right\},
\]
need not be closed in $\sA^{0,1}(E)\times W^{1,p}(E)$. Theorem \ref{thm:Existence_of_complex_gauge_transformation_to_Coulomb_01_pairs} does not give a (global) slice for the action of $W^{2,p}(\SL(E))$ on $\sA^{0,1}(E)\times W^{1,p}(E)$ because the condition in Item \eqref{item:Duistermaat_Kolk_2-3-1_GlobalTube} of the Definition \ref{defn:Duistermaat_Kolk_2-3-1} of a slice might fail since the group element in that condition need not be close to the identity. Therefore, following Friedman and Morgan \cite[Sections 4.1.5 and 4.1.6]{FrM}, Itoh \cite[Sections 3 and 4]{Itoh_1985}, and Kobayashi \cite[Section 7.3]{Kobayashi_differential_geometry_complex_vector_bundles}, we shall focus instead on the local Kuranishi models defined by the Coulomb gauge slice condition \eqref{eq:Coulomb_gauge_slice_condition_for_01-pairs} and the holomorphic pair equations \eqref{eq:Holomorphic_pair}.

\subsection{Local Kuranishi model for the moduli space of holomorphic pairs}
\label{subsec:Kuranishi_model_moduli_space_holomorphic_pairs}
We continue the hypotheses of Theorem \ref{thm:Existence_of_complex_gauge_transformation_to_Coulomb_01_pairs}. As
in Friedman and Morgan \cite[Chapter IV]{FrM}, we shall suppress explicit notation for Sobolev $L^p$  completions with $p\in (n,\infty)$, unless required for clarity, in order to avoid notational clutter and instead write
\[
  \sA^{0,1}(E)\times \Omega^0(E)
  =
  \left(\bar\partial_E+\Omega^{0,1}(\fsl(E))\right) \times \Omega^0(E)
  \quad\text{and}\quad C^\infty(\SL(E))
\]
for the affine space of $(0,1)$-pairs (inducing a fixed $(0,1)$-connection on the complex line bundle $\det E$) and the group of complex, determinant-one automorphisms of $E$, where $\bar\partial_E$ is smooth $(0,1)$-connection on $E$. Recall that the left-hand side of the holomorphic pair equations \eqref{eq:Holomorphic_pair} define the complex-analytic holomorphic pair map \eqref{eq:Holomorphic_pair_map},
\begin{multline*}
  \fS:\sA^{0,1}(E)\times \Omega^0(E) \ni (\bar\partial_E+\alpha,\varphi+\sigma)
  \\
  \mapsto
  \left( F_{\bar\partial_E+\alpha},(\bar\partial_E+\alpha)(\varphi+\sigma) \right)
  \in
  \Omega^{0,2}(\fsl(E)) \oplus \Omega^{0,1}(E).
\end{multline*}
Thanks to the Bianchi Identity \eqref{eq:Holomorphic_pair_Bianchi_identity}, the map \eqref{eq:Holomorphic_pair_map} may also be viewed as a morphism of sheaves,
\begin{multline}
  \label{eq:Holomorphic_pair_and_Coulomb_gauge_map_sheaves}
  \fS:\sA^{0,1}(E)\times \Omega^0(E)
  \ni (\bar\partial_E+\alpha,\varphi+\sigma)
  \\
  \mapsto
  \left( F_{\bar\partial_E+\alpha},(\bar\partial_E+\alpha)(\varphi+\sigma) \right),
  \\
  \in \Ker\bar\partial_{E+\alpha,\varphi+\sigma}^2\cap\left(\Omega^{0,2}(\fsl(E)) \oplus \Omega^{0,1}(E)\right),
\end{multline}
where $\bar\partial_{E+\alpha,\varphi+\sigma}^2$ is defined in \eqref{eq:dkStablePair}. Assume now that $(\bar\partial_E,\varphi)$ is a \emph{holomorphic pair}, so $(\bar\partial_E,\varphi) \in \fS^{-1}(0)$. The linearization of the map \eqref{eq:Holomorphic_pair_map} at $(\bar\partial_E,\varphi)$ is given by
\[
  d\fS(\bar\partial_E,\varphi):\Omega^{0,2}(\fsl(E)) \oplus \Omega^0(E) \ni (\alpha,\sigma)
  \mapsto
  \bar\partial_{E,\varphi}^1(\alpha,\sigma)
  \in
  \Omega^{0,2}(\fsl(E)) \oplus \Omega^{0,1}(E),
\]
where the operator $\bar\partial_{E,\varphi}^1$ is as in \eqref{eq:d1StablePair}. Because the operator $\bar\partial_{E,\varphi}^1$ is a differential in an elliptic complex \eqref{eq:Holomorphic_pair_elliptic_complex} when $(\bar\partial_E,\varphi)$ is a holomorphic pair, then
\[
  \Ran\bar\partial_{E,\varphi}^1 \subset \Ker \bar\partial_{E,\varphi}^2.
\]
Hence, we obtain a linear operator,
\[
  d\fS(\bar\partial_E,\varphi)
  =
  \bar\partial_{E,\varphi}^1:
  \Omega^{0,1}(\fsl(E)) \oplus \Omega^0(E) 
  \\
  \to \Ker \bar\partial_{E,\varphi}^2\cap\left(\Omega^{0,2}(\fsl(E)) \oplus \Omega^{0,1}(E)\right),
\]
that we shall now verify to be Fredholm upon restriction to the tangent space to the Coulomb-gauge, affine slice $S_{\bar\partial_E,\varphi}$ in \eqref{eq:dbar_Evarphi_slice} through the pair $(\bar\partial_E,\varphi)$:
\begin{equation}
  \label{eq:Tangent_space_dbar_Evarphi_slice}
  T_{\bar\partial_E,\varphi}S_{\bar\partial_E,\varphi}
  =
  \bar\partial_{E,\varphi}^{0,*}\cap \left(\Omega^{0,1}(\fsl(E)) \oplus \Omega^0(E)\right).
\end{equation}
The Hodge decomposition for elliptic complexes (see Gilkey \cite[Theorem 1.5.2]{Gilkey2}) yields $L^2$-orthogonal decompositions,
\begin{subequations}
\label{eq:Hodge_decomposition_holomorphic_pair_complex}
\begin{align}
  \label{eq:Hodge_decomposition_holomorphic_pair_complex_Omega^0}
  \Omega^0(\fsl(E))
  &=
  \bH_{\bar\partial_E,\varphi}^0 \oplus \Ran\bar\partial_{E,\varphi}^{0,*},
  \\
  \label{eq:Hodge_decomposition_holomorphic_pair_complex_Omega^01_0}
  \Omega^{0,1}(\fsl(E)) \oplus \Omega^0(E)
  &=
    \bH_{\bar\partial_E,\varphi}^1 \oplus \Ran\bar\partial_{E,\varphi}^0 \oplus \Ran\bar\partial_{E,\varphi}^{1,*},
  \\
  \label{eq:Hodge_decomposition_holomorphic_pair_complex_Omega^02_01}
  \Omega^{0,2}(\fsl(E)) \oplus \Omega^{0,1}(E)
  &=
    \bH_{\bar\partial_E,\varphi}^2 \oplus \Ran\bar\partial_{E,\varphi}^1 \oplus \Ran\bar\partial_{E,\varphi}^{2,*},
  \\
  \label{eq:Hodge_decomposition_holomorphic_pair_complex_Omega^03_02}
  \Omega^{0,2}(E)
  &=
    \bH_{\bar\partial_E,\varphi}^3 \oplus \Ran\bar\partial_{E,\varphi}^2,
\end{align}
\end{subequations}
where the harmonic spaces are as in \eqref{eq:H_dbar_Avarphi^0bullet} (for $k=0,1,2,3$),
\begin{align*}
  \bH_{\bar\partial_E,\varphi}^0
  &=
    \Ker\bar\partial_{E,\varphi}^0,
  \\
  \bH_{\bar\partial_E,\varphi}^1
  &=
    \Ker\bar\partial_{E,\varphi}^1 \cap \Ker\bar\partial_{E,\varphi}^{0,*},
  \\
  \bH_{\bar\partial_E,\varphi}^2
  &=
    \Ker\bar\partial_{E,\varphi}^2 \cap \Ker\bar\partial_{E,\varphi}^{1,*},
  \\
  \bH_{\bar\partial_E,\varphi}^3
  &=
    \Ker\bar\partial_{E,\varphi}^{2,*}.    
\end{align*}
Clearly, the expression for $d\fS(\bar\partial_E,\varphi)$ yields
\[
  \Ker d\fS(\bar\partial_E,\varphi) \cap \Ker \bar\partial_{E,\varphi}^{0,*}
  =
  \Ker\left(\bar\partial_{E,\varphi}^1 + \bar\partial_{E,\varphi}^{0,*}\right)
  =
  \bH_{\bar\partial_E,\varphi}^1.
\]
Moreover,
\[
  \Ker\left(\bar\partial_{E,\varphi}^2 + \bar\partial_{E,\varphi}^{1,*}\right)
  \cap\left(\Omega^{0,2}(\fsl(E)) \oplus \Omega^{0,1}(E) \right)
    =
  \bH_{\bar\partial_E,\varphi}^2,
\]
and so the expression for $d\fS(\bar\partial_E,\varphi)$ also yields
\begin{align*}
  \left(\Ran d\fS(\bar\partial_E,\varphi)\right)^\perp
  &=
    \left(\Ran\left(\bar\partial_{E,\varphi}^1 +\bar\partial_{E,\varphi}^{0,*}\right)\right)^\perp
    \cap\left(\Ker \bar\partial_{E,\varphi}^2 \oplus \Ran\bar\partial_{E,\varphi}^{0,*}\right)
  \\
  &=
    \Ran\left(\bar\partial_{E,\varphi}^1\right)^\perp \cap\Ker \bar\partial_{E,\varphi}^2
   \oplus
    \left(\Ran\bar\partial_{E,\varphi}^{0,*}\right)^\perp
    \cap\Ran\bar\partial_{E,\varphi}^{0,*}
  \\
  &=
    \Ker\bar\partial_{E,\varphi}^{1,*} \cap\Ker \bar\partial_{E,\varphi}^2
  \\
  &=
  \bH_{\bar\partial_E,\varphi}^2.
\end{align*}
Hence, the restriction of the operator $d\fS(\bar\partial_E,\varphi)$ to the tangent space \eqref{eq:Tangent_space_dbar_Evarphi_slice} to the slice $S_{\bar\partial_E,\varphi}$ has finite-dimensional kernel and cokernel and so is Fredholm (see Abramovich and Aliprantis \cite[Definition 4.37 and Lemma 4.38, p. 156]{Abramovich_Aliprantis_2002}).

The Kuranishi Method \cite{Kuranishi} (via the Implicit Mapping Theorem \ref{thm:Friedman_Morgan_4-1-9} for complex analytic maps of complex Banach spaces and the forthcoming Corollary \ref{cor:Friedman_Morgan_1999_1-2_restricted_to_Coulomb_gauge_slice_neighborhood} or by extending the construction in Theorem \ref{thm:Local_Kuranishi_model_for_strongly_simple_point_cP(E)}) yields a $\Stab(\bar\partial_E,\varphi)$-invariant open neighborhood $\UU_{\bar\partial_E,\varphi}$ of $(\bar\partial_E,\varphi)$ in $\sA^{0,1}(E)\times\Omega^0(E)$ and a corresponding $\Stab(\bar\partial_E,\varphi)$-invariant open neighborhood $N_{\bar\partial_E,\varphi}$ of the origin in $\bH_{\bar\partial_E,\varphi}^1$, so
\[
  (\bar\partial_E,\varphi) + N_{\bar\partial_E,\varphi}
  =
  \left( (\bar\partial_E,\varphi) + \bH_{\bar\partial_E,\varphi}^1 \right) \cap \UU_{\bar\partial_E,\varphi},
\]
and a $\Stab(\bar\partial_E,\varphi)$-equivariant, complex analytic embedding (see Lemma \ref{lem:EquivariantKuranishiLemma})
\begin{equation}
  \label{eq:Kuranishi_embedding_map_holomorphic_pairs}
  \bgamma: \bH_{\bar\partial_E,\varphi}^1\supset N_{\bar\partial_E,\varphi} \to \sA^{0,1}(E)\times\Omega^0(E)
\end{equation}
with $\bgamma(0,0) = (\bar\partial_E,\varphi)$ and a $\Stab(\bar\partial_E,\varphi)$-equivariant, complex analytic map
\begin{equation}
  \label{eq:Kuranishi_obstruction_map_holomorphic_pairs}
  \bchi: \bH_{\bar\partial_E,\varphi}^1\supset N_{\bar\partial_E,\varphi} \to \bH_{\bar\partial_E,\varphi}^2,
\end{equation}
such that
\[
  \fS(\bgamma(\alpha,\sigma)) = 0
  \iff
  \bchi(\alpha,\sigma) = 0,
  \quad\text{for all } (\alpha,\sigma) \in N_{\bar\partial_E,\varphi}^1.
\]
The definitions of $\bgamma$ and $\bchi$ in \eqref{eq:Kuranishi_embedding_map_holomorphic_pairs} and \eqref{eq:Kuranishi_obstruction_map_holomorphic_pairs}, respectively, slightly extend those in \eqref{eq:Kuranishi_embedding_map_holomorphic_pairs_trivial_stabilizer} and \eqref{eq:Kuranishi_obstruction_map_holomorphic_pairs_trivial_stabilizer}, respectively, that were already provided by the construction of Theorem \ref{thm:Local_Kuranishi_model_for_strongly_simple_point_cP(E)} in the sense that here we allow non-trivial stabilizer groups, $\Stab(\bar\partial_E,\varphi)$. Below, we restate Lemma \ref{lem:Friedman_Morgan_1999_1-2} in a way that is more useful for our present application.

\begin{cor}[Local equivalence of holomorphic pair map with composition of holomorphic pair map and orthogonal projection]
\label{cor:Friedman_Morgan_1999_1-2_restricted_to_Coulomb_gauge_slice_neighborhood}
Continue the hypotheses of Theorem \ref{thm:Existence_of_complex_gauge_transformation_to_Coulomb_01_pairs} and assume further that $(\bar\partial_E,\varphi)$ is a smooth holomorphic pair in the sense of \eqref{eq:Holomorphic_pair}. Then after possibly shrinking the open neighborhood $\UU_{\bar\partial_E,\varphi}$ of the pair $(\bar\partial_E,\varphi)$ in $\sA^{0,1}(E)\times\Omega^0(E)$ provided by Theorem \ref{thm:Existence_of_complex_gauge_transformation_to_Coulomb_01_pairs}, the zero locus of the map $\fS$ in \eqref{eq:Holomorphic_pair_map} is equal to the zero locus of the map
\begin{multline}
  \label{eq:Holomorphic_pair_map_Banach_spaces_orthogonal_projections}
  \widehat\fS:\sA^{0,1}(E)\times \Omega^0(E) \supset
  \UU_{\bar\partial_E,\varphi} 
  \ni (\bar\partial_E+\alpha,\varphi+\sigma)
  \\
  \mapsto
  \Pi_{\bar\partial_E,\varphi}^{0,2}\left(F_{\bar\partial_E + \alpha},(\bar\partial_E+\alpha)(\varphi+\sigma)\right),
  \\
  \in \Ker\bar\partial_{E,\varphi}^2\cap\left(\Omega^{0,2}(\fsl(E)) \oplus \Omega^{0,1}(E)\right),
\end{multline}
where
\begin{equation}
  \label{eq:L^2_orthogonal_projection_Ker_dbar_Evarphi^2}
  \Pi_{\bar\partial_E,\varphi}^{0,2}: \Omega^{0,2}(\fsl(E)) \oplus \Omega^{0,1}(E)
  \to
  \Ker\bar\partial_{E,\varphi}^2\cap\left(\Omega^{0,2}(\fsl(E)) \oplus \Omega^{0,1}(E)\right)
\end{equation}
is $L^2$-orthogonal projection.
\end{cor}

The advantage of the map $\widehat\fS$ in \eqref{eq:Holomorphic_pair_map_Banach_spaces_orthogonal_projections} over the equivalent map $\fS$ in \eqref{eq:Holomorphic_pair_map} is that it defines a Fredholm map of fixed Banach spaces, given a choice of holomorphic pair $(\bar\partial_E,\varphi)$, rather than a Fredholm map of sheaves, where the Fredholm property is obtained upon restriction to the tangent space at $(\bar\partial_E,\varphi)$ to the slice $S_{\bar\partial_E,\varphi}$. We have the following analogue for holomorphic pairs of Friedman and Morgan \cite[Section 4.1.5, Definition 1.15, p. 295]{FrM} for holomorphic bundles.

\begin{defn}[Kuranishi model space for a smooth holomorphic pair]
\label{defn:Friedman_Morgan_4-1-15_holomorphic_pairs}
Continue the hypotheses of Theorem \ref{thm:Existence_of_complex_gauge_transformation_to_Coulomb_01_pairs} and assume further that $(\bar\partial_E,\varphi)$ is a smooth holomorphic pair in the sense of \eqref{eq:Holomorphic_pair}. The holomorphic embedding $\bgamma$ in \eqref{eq:Kuranishi_embedding_map_holomorphic_pairs} is the \emph{Kuranishi embedding map} and the holomorphic map $\bchi$ in \eqref{eq:Kuranishi_obstruction_map_holomorphic_pairs} is the \emph{Kuranishi obstruction map} for $(\bar\partial_E,\varphi)$, and the germ of the complex analytic set $\bchi^{-1}(0)$ at the origin in $\bH_{\bar\partial_E,\varphi}^1$, with its natural ringed space structure, is the \emph{Kuranishi model} for $(\bar\partial_E,\varphi)$. We denote this germ of a complex analytic space by $\fK(\bar\partial_E,\varphi)$.
\end{defn}

Note that Definition \ref{defn:Friedman_Morgan_4-1-15_holomorphic_pairs} yields finite-dimensional, complex analytic subsets,
\[
  \bchi^{-1}(0) \subset \bH_{\bar\partial_E,\varphi}^1
  \quad\text{and}\quad
  \bgamma(\bchi^{-1}(0)) \subset \fS^{-1}(0) \subset \sA^{0,1}(E)\times\Omega^0(E).
\]
The Kuranishi model is simplest when $\Stab(\bar\partial_E,\varphi) = \{\id_E\}$ (compare Friedman and Morgan \cite[Section 4.1.5, Definition 1.15, p. 295]{FrM}), though our Definition \ref{defn:Friedman_Morgan_4-1-15_holomorphic_pairs} does not make this assumption. When $\Stab(\bar\partial_E,\varphi)$ is non-trivial, the domain of the map $\bchi$ in \eqref{eq:Kuranishi_obstruction_map_holomorphic_pairs} may be chosen to be $\Stab(\bar\partial_E,\varphi)$-invariant and $\bchi$ to be $\Stab(\bar\partial_E,\varphi)$-equivariant; see Friedman and Morgan \cite[Section 4.1.5, Remark, p. 296]{FrM}.

We have the following analogue of Definition \ref{defn:Friedman_Morgan_4-1-15_holomorphic_pairs} in the context of projective vortices, though we restrict to the case of trivial stabilizer.

\begin{defn}[Kuranishi model space for a smooth projective vortex]
\label{defn:Friedman_Morgan_4-1-15_projective_vortex}
Let $(E,h)$ be a Hermitian vector bundle over a complex, Hermitian manifold $(X,g,J)$ with fundamental two-form $\omega = g(\cdot,J\cdot)$ as in \eqref{eq:Fundamental_two-form} and let $A_d$ be a fixed smooth, unitary connection on the Hermitian line bundle $\det E$. Let $(A,\varphi)$ is a solution to the projective vortex equations \eqref{eq:SO(3)_monopole_equations_almost_Hermitian_alpha} with trivial stabilizer as in Definition \ref{defn:Split_trivial_central-stabilizer_unitary_pair} (\ref{item:Trivial-stabilizer_unitary_pair}), so $\Stab(A,\varphi) = \{\id_E\}$, and such that $A$ induces $A_d$ on $\det E$. The real analytic embedding $\beps$ in \eqref{eq:Kuranishi_embedding_map_projective_vortex} and the real analytic map $\bkappa$ in \eqref{eq:Kuranishi_obstruction_map_projective_vortex} provided by Theorem \ref{thm:Local_Kuranishi_model_for_moduli_space_projective_vortices_StabAvarphi_idE} are the \emph{Kuranishi embedding map} and \emph{Kuranishi obstruction map} for $(A,\varphi)$, respectively. The germ of the real analytic set $\bkappa^{-1}(0)$ at the origin in $\bH_{E,\varphi}^1$, with its natural ringed space structure, is the \emph{Kuranishi model} for $(A,\varphi)$. We denote this germ of a real analytic space by $\sK(A,\varphi)$.
\end{defn}

Our Definition \ref{defn:Friedman_Morgan_4-1-15_holomorphic_pairs} potentially relies on choices that should not be intrinsic to a complex analytic space, namely the data used to define the local slice for the $\SL(E)$-action, such as Hermitian metrics on the complex bundle $E$ and complex manifold $X$. The proof of \cite[Section 4.1.5, Proposition 1.16, p. 296]{FrM} due to Friedman and Morgan adapts to give the

\begin{prop}[Independence of Kuranishi models from choices of local slices for the action of the automorphism group]
\label{prop:Friedman_Morgan_4-1-16_holomorphic_pairs}
Continue the hypotheses of Theorem \ref{thm:Existence_of_complex_gauge_transformation_to_Coulomb_01_pairs} and assume further that $(\bar\partial_E,\varphi)$ is a smooth holomorphic pair in the sense of \eqref{eq:Holomorphic_pair}. Then the Kuranishi models for $(\bar\partial_E,\varphi)$ defined by two different local slices for the action of $W^{2,p}(\SL(E))$ are isomorphic.
\end{prop}

\begin{proof}
We include the proof since we shall need to appeal to similar arguments later in this chapter. If $S_{\bar\partial_E,\varphi}' \subset \sA^{0,1}(E)\times W^{1,p}(E)$ is another local slice (in the sense of Remark \ref{rmk:Duistermaat_Kolk_2-3-1_local_slice})
through $(\bar\partial_E,\varphi)$ for the action of $W^{2,p}(\SL(E))$, then there are an open neighborhood $U_{\bar\partial_E,\varphi}' \subset S_{\bar\partial_E,\varphi}'$ of $(\bar\partial_E,\varphi)$, an open neighborhood $U_{\bar\partial_E,\varphi} = \UU_{\bar\partial_E,\varphi} \cap S_{\bar\partial_E,\varphi}$ of $(\bar\partial_E,\varphi)$ in the affine, Coulomb-gauge, local slice $S_{\bar\partial_E,\varphi}$ given by Theorem \ref{thm:Existence_of_complex_gauge_transformation_to_Coulomb_01_pairs}, and a biholomorphic map:
\[
  \bv: U_{\bar\partial_E,\varphi} \ni (\bar\partial_E+\alpha,\varphi+\sigma)
  \mapsto
  v\cdot(\bar\partial_E+\alpha,\varphi+\sigma) \in U_{\bar\partial_E,\varphi}',
\]
where $v \in W^{2,p}(\SL(E))$ acts on $\sA^{0,1}(E)\times W^{1,p}(E)$ via \eqref{eq:SL(E)ActionOn(0,1)Pairs} and we abbreviate
\[
  v = v_{\bar\partial_E+\alpha,\varphi+\sigma} = \bv(\bar\partial_E+\alpha,\varphi+\sigma).
\]  
Recall that the map $\fS$ is $W^{2,p}(\SL(E))$-equivariant, so
\begin{multline*}
  \fS\left(v\cdot(\bar\partial_E+\alpha,\varphi+\sigma)\right)
  =
  v\cdot\fS\left(\bar\partial_E+\alpha,\varphi+\sigma\right)
  \in
  L^p\left(\Lambda^{0,2}(\fsl(E)) \oplus \Lambda^{0,1}(E)\right),
  \\
  \text{for all } v \in W^{2,p}(\SL(E)).
\end{multline*}
The following map,
\begin{multline*}
  L^p\left(\Lambda^{0,2}(\fsl(E)) \oplus \Lambda^{0,1}(E)\right) \times U_{\bar\partial_E,\varphi}
  \ni
  \left((\beta,\nu),(\bar\partial_E+\alpha,\varphi+\sigma)\right)
  \\
  \mapsto \left(v_{\bar\partial_E+\alpha,\varphi+\sigma}^{-1}\cdot(\beta,\nu),(\bar\partial_E+\alpha,\varphi+\sigma)\right)
  \in
  L^p\left(\Lambda^{0,2}(\fsl(E)) \oplus \Lambda^{0,1}(E)\right) \times U_{\bar\partial_E,\varphi},
\end{multline*}
is a bundle automorphism of the product Banach vector bundle, which we shall also denote by $\bv^{-1}$.  The construction of the Kuranishi model via the affine, Coulomb-gauge, local slice $S_{\bar\partial_E,\varphi}$ corresponds to taking the inverse image of the product subbundle,
\[
  \bH_{\bar\partial_E,\varphi}^2 \times U_{\bar\partial_E,\varphi}
  \subset
  L^p\left(\Lambda^{0,2}(\fsl(E)) \oplus \Lambda^{0,1}(E)\right) \times U_{\bar\partial_E,\varphi},
\]
under the section associated to $\fS$. The Kuranishi model associated to the local slice $S_{\bar\partial_E,\varphi}'$ corresponds to taking the inverse image of the finite-dimensional subbundle which is the image of the product bundle under $\bv^{-1}$. By Theorem \ref{thm:Friedman_Morgan_4-1-14}, the associated Kuranishi models are isomorphic.
\end{proof}

Next, we have the following analogue for holomorphic pairs of Friedman and Morgan \cite[Section 4.1.5, Proposition 1.17, p. 297]{FrM} for holomorphic bundles.

\begin{prop}
\label{prop:Friedman_Morgan_4-1-17_holomorphic_pairs}
Let $(T, t_0)$ be a germ of a complex analytic space and $f: T \to \sA^{0,1}(E)\times W^{1,p}(E)$ be a holomorphic map with $f(t_0) = (\bar\partial_E,\varphi)$ and such that $f(T)$ is contained in $\fS^{-1}(0)$ in the sense of complex analytic spaces. Then the projection
\[
  \UU_{\bar\partial_E,\varphi} = U_{\bar\partial_E,\varphi}\times U_{\id_E}
  \to U_{\bar\partial_E,\varphi}
\]
defines a morphism of germs of complex analytic spaces from $T$ to
\[
  \fK(\bar\partial_E,\varphi) \subseteq N_{\bar\partial_E,\varphi} \subset U_{\bar\partial_E,\varphi},
\]
where $\fK(\bar\partial_E,\varphi)$ is as in the Kuranishi model of Definition \ref{defn:Friedman_Morgan_4-1-15_holomorphic_pairs} and
$N_{\bar\partial_E,\varphi} \subset \bH_{\bar\partial_E,\varphi}^1$ is an open neighborhood of the origin and
the preceding open neighborhoods are as in Theorem \ref{thm:Existence_of_complex_gauge_transformation_to_Coulomb_01_pairs}:
\begin{align*}
  \UU_{\bar\partial_E,\varphi} &\subset \sA^{0,1}(E) \times W^{1,p}(E),
  \\
  U_{\bar\partial_E,\varphi} &= S_{\bar\partial_E,\varphi} \cap \UU_{\bar\partial_E,\varphi},
  \\
  U_{\id_E} &\subset W^{2,p}(\SL(E)).
\end{align*}
More precisely, there is a map of complex analytic spaces,
\[
  \bv:T \ni t \mapsto v_t = \bv(t) \in U_{\id_E} \subset W^{2,p}(\SL(E)),
\]
so that the composite map
\[
  \bv\cdot f:T \ni t \mapsto v_t\cdot f(t) \in \fK(\bar\partial_E,\varphi) \subseteq N_{\bar\partial_E,\varphi}
\]
gives the desired morphism, where $v_t \in W^{2,p}(\SL(E))$ acts on $\sA^{0,1}(E)\times W^{1,p}(E)$ via \eqref{eq:SL(E)ActionOn(0,1)Pairs}.
\end{prop}

\begin{proof}
We first show that an open neighborhood of $(\bar\partial_E,\varphi)$ in $\fS^{-1}(0)$ is contained in an embedded complex submanifold of $\UU_{\bar\partial_E,\varphi}$ that is biholomorphic to an open neighborhood of $(0,\id_E)$ in
\[
  \bH_{\bar\partial_E,\varphi}^1 \times U_{\id_E},
\]
such that the biholomorphic map restricts to the inclusion
\[
  \fK(\bar\partial_E,\varphi) \subset \bH_{\bar\partial_E,\varphi}^1.
\]
We first observe that the following map is holomorphic:
\begin{multline*}
  S_{\partial_E,\varphi} \times  U_{\id_E}
  =
  \left((\bar\partial_E,\varphi) + \bH_{\bar\partial_E,\varphi}^1 \oplus \Ran\bar\partial_{E,\varphi}^{2,*}\right)
  \times  U_{\id_E}
  \ni \left(\left(\bar\partial_E+\alpha,\varphi+\sigma\right), v\right)
  \\
  \mapsto
  \fS\left(v\cdot\left(\bar\partial_E+\alpha,\varphi+\sigma\right)\right)
  \in
  L^p\left(\Lambda^{0,2}(\fsl(E))\oplus \Lambda^{0,1}(E)\right).
\end{multline*}
Note from the Hodge decompositions \eqref{eq:Hodge_decomposition_holomorphic_pair_complex_Omega^01_0} and \eqref{eq:Hodge_decomposition_holomorphic_pair_complex_Omega^02_01} that there are $L^2$-orthogonal decompositions:
\begin{align*}
  W^{1,p}\left(\Lambda^{0,1}(\fsl(E)) \oplus E\right)
  &=
    \bH_{\bar\partial_E,\varphi}^1 \oplus \Ran\bar\partial_{E,\varphi}^0 \oplus \Ran\bar\partial_{E,\varphi}^{1,*},
  \\
  L^p\left(\Lambda^{0,2}(\fsl(E)) \oplus \Lambda^{0,1}(E) \right)
  &=
    \bH_{\bar\partial_E,\varphi}^2 \oplus \Ran\bar\partial_{E,\varphi}^1 \oplus \Ran\bar\partial_{E,\varphi}^{2,*}.
\end{align*}
By an application of the Implicit Mapping Theorem similar to that used in proving the existence of the Kuranishi obstruction map $\bchi$ in \eqref{eq:Kuranishi_obstruction_map_holomorphic_pairs}, there is a holomorphic map
\[
  \tilde\bgamma_\perp: \bH_{\bar\partial_E,\varphi}^1 \times U_{\id_E}
  \supset \widetilde U_{\bar\partial_E,\varphi}
  \to
  \Ran\bar\partial_{E,\varphi}^{1,*}
\]
from an open neighborhood $\widetilde U_{\bar\partial_E,\varphi}$ of $(0,\id_E)$ such that
\[
  \fS\left(v\cdot\left(\left(\bar\partial_E,\varphi\right) + (\alpha,\sigma) + \tilde\bgamma_\perp(\alpha,\sigma)\right)\right)
  \in
  \bH_{\bar\partial_E,\varphi}^2,
  \quad\text{for all } ((\alpha,\sigma),v) \in \widetilde U_{\bar\partial_E,\varphi},
\]
and the graph of $\tilde\bgamma_\perp$ contains all pairs $((\alpha',\sigma'),v)$ in an open neighborhood of $(0,\id_E)$ for which
\[
  \fS\left(v\cdot\left(\bar\partial_E+\alpha',\varphi+\sigma'\right)\right) = (0,0).
\]
One can show
that $\tilde\bgamma_\perp$ restricts to $\bgamma_\perp$, when we write the Kuranishi embedding map $\bgamma$ in \eqref{eq:Kuranishi_embedding_map_holomorphic_pairs} as
\[
  \bgamma = \bgamma_\parallel + \bgamma_\perp: \bH_{\bar\partial_E,\varphi}^1\supset N_{\bar\partial_E,\varphi}
  \to
  (\bar\partial_E,\varphi) + \bH_{\bar\partial_E,\varphi}^1 \oplus \Ran\bar\partial_{E,\varphi}^{1,*}
  \subset \sA^{0,1}(E)\times W^{1,p}(E),
\]
with
\[
   \bgamma_\parallel: \bH_{\bar\partial_E,\varphi}^1\supset N_{\bar\partial_E,\varphi}
  \to
  (\bar\partial_E,\varphi) + \bH_{\bar\partial_E,\varphi}^1
  \subset \sA^{0,1}(E)\times W^{1,p}(E).
\]
Define the composition
\[
  \tilde\bchi: \widetilde U_{\bar\partial_E,\varphi} \ni ((\alpha,\sigma),v)
  \mapsto
  \fS\left(v\cdot\left( \left(\bar\partial_E,\varphi\right) + (\alpha,\sigma) + \tilde\bgamma_\perp(\alpha,\sigma) \right)\right)
  \in \bH_{\bar\partial_E,\varphi}^2.
\]
In particular, $\tilde\bchi$ takes values in the finite-dimensional complex vector space $\bH_{\bar\partial_E,\varphi}^2$. It follows from Lemma \ref{lem:Friedman_Morgan_4-1-11} that the hypothesis that the image of $f:T \to \sA^{0,1}(E)\times\Omega^0(E)$ is contained in $\fS^{-1}(0)$ in the sense of complex analytic spaces is equivalent to saying that the image of $f$ is contained in the graph of $\tilde\bgamma_\perp$ and that, moreover, if we replace $f$ by the corresponding map $\tilde f:T\to \widetilde U_{\bar\partial_E,\varphi}$, then $\tilde\bchi\circ\tilde f = 0$ on $T$. We also have the following facts about $\tilde\bchi$:
\begin{itemize}
\item For each $v \in U_{\id_E}$, the gauge equivariance of $\fS$ implies that there is an inclusion of complex analytic spaces,
\[
  \fK\left(\bar\partial_E,\varphi\right) \times \{v\} \subset \tilde\bchi^{-1}(0).
\]
\item The restriction of $\tilde\bchi$ to the slice $N_{\bar\partial_E,\varphi}\times\{\id_E\}$ defines a holomorphic map from a finite-dimensional open neighborhood $N_{\bar\partial_E,\varphi}$ of the origin in $\bH_{\bar\partial_E,\varphi}^1$ into $\bH_{\bar\partial_E,\varphi}^2$ that is equal to the Kuranishi obstruction map $\bchi$ in \eqref{eq:Kuranishi_obstruction_map_holomorphic_pairs}.
\end{itemize}
The conclusion of the proposition is then a consequence of the forthcoming Corollary \ref{cor:Friedman_Morgan_4-1-19} to Lemma \ref{lem:Friedman_Morgan_4-1-18}.
\end{proof}

The proof of \cite[Section 4.1.5, Proposition 1.17, p. 297]{FrM} due to Friedman and Morgan, and hence its adaptation to prove Proposition \ref{prop:Friedman_Morgan_4-1-17_holomorphic_pairs} relies on two results about complex analytic spaces whose statements and proofs (included in \cite[Section 4.1.5]{FrM}) they attribute to Douady \cite{Douady_1966_sem_bourbaki}. We include them here for completeness.

\begin{lem}
\label{lem:Friedman_Morgan_4-1-18}  
(See Friedman and Morgan \cite[Section 4.1.5, Lemma 1.18, p. 297]{FrM} and Douady \cite{Douady_1966_sem_bourbaki}.)
Let $\fK$ be the germ of a complex analytic space at $0$ inside $\CC^n$, defined in a neighborhood of the origin by the holomorphic functions $g_1, \ldots, g_s$. Let $U'$ be an open subset of a complex Banach space. Suppose that $\tilde\bchi$ is a holomorphic map from a neighborhood $\widetilde U$ of the origin in $\CC^n\times U'$ to a Banach space $\sX$ with $\tilde\bchi(0) = 0$ such that the following hold:
\begin{enumerate}
\item There exist holomorphic functions $h_1,\ldots, h_s$ defined in a neighborhood of the origin in $\sX$ such that, for all $i$,
\[
  h_i \circ \tilde\bchi(z, 0) = g_i(z).
\]

\item For all $v \in U'$ sufficiently close to $0$, the complex analytic subspace $\fK \times \{v\}$ of $\widetilde U$ is contained in $\tilde\bchi^{-1}(0)$ in the sense of complex analytic spaces. 
\end{enumerate}
Then there exist a neighborhood $N$ of $0$ in $\CC^n$ and a neighborhood $N'$ of $0$ in $U'$ and holomorphic functions $\alpha_{ij}(u,v)$ on $N \times N'$ such that
\[
    \sum_j \alpha_{ij}(u,v)(h_j\circ\tilde\bchi)(u, v) = g_i(u).
\]
\end{lem}

\begin{cor}
\label{cor:Friedman_Morgan_4-1-19}
(See Friedman and Morgan \cite[Section 4.1.5, Corollary 1.19, p. 298]{FrM} and Douady \cite{Douady_1966_sem_bourbaki}.)
Suppose that $\fK$ and $\tilde\bchi$ satisfy the hypotheses of Lemma \ref{lem:Friedman_Morgan_4-1-18}. Let $f: T\to \widetilde U$ be a holomorphic map from the germ of a complex space $(T, t_0)$ to $\widetilde U$ with $f(t_0) = 0$, and suppose that $\tilde\bchi^{-1}(0)$ contains $\Imag f$ in the sense of complex analytic spaces. Then the image of the map $T \to \CC^n$ given by composing $f$ with the projection to the first factor is contained in $\fK$ in the sense of complex analytic spaces.
\end{cor}

The first condition of Lemma \ref{lem:Friedman_Morgan_4-1-18} is equivalent to saying that $\tilde\bchi^{-1}(0)\cap(\CC^n\times\{0\})$ is contained in $\fK$ in the sense of complex analytic spaces, and hence, assuming the second condition, is equal to $\fK$. We now consider general properties of the map $\bchi$. The proof of the following proposition requires only a  minor modification of the proof of \cite[Section 4.1.5, Proposition 1.20, p. 299]{FrM}.

\begin{prop}
\label{prop:Friedman_Morgan_4-1-20}    
(See Friedman and Morgan \cite[Section 4.1.5, Proposition 1.20, p. 299]{FrM} for the case of holomorphic bundles.)  
Let $\bchi$ be the Kuranishi obstruction map \eqref{eq:Kuranishi_obstruction_map_holomorphic_pairs}. Then the following hold:
\begin{enumerate}
\item $\bchi(0) = 0$.
\item The differential of $\bchi$ at the origin is identically zero.
\end{enumerate}
\end{prop}

\section{Local moduli functor for holomorphic pairs}
\label{subsec:Friedman-Morgan_4-2}
Throughout this section, $E$ denotes a complex vector bundle over a closed, complex manifold $X$. Our development of the local moduli functor for holomorphic pairs $(\bar\partial_E,\varphi)$ on $E$ and proof of its representability is based on that of Friedman and Morgan in \cite[Section 4.2]{FrM} for holomorphic structures $\bar\partial_E$ on $E$.

We define the local moduli functor and the concepts of representability and universal object in Section \ref{subsec:Friedman-Morgan_4-2-1}. We prove that the local moduli functor is representable at regular points in Section \ref{subsec:Friedman-Morgan_4-2-2}. The proof that the local moduli functor is representable in general, without restriction to regular points, appears in Sections \ref{subsec:Friedman-Morgan_4-2-3} and \ref{subsec:Friedman-Morgan_4-2-4}. 

\subsection{Definition of the local moduli functor for holomorphic pairs}
\label{subsec:Friedman-Morgan_4-2-1}
We adapt the discussion due to Friedman and Morgan in \cite[Section 4.2.1]{FrM} from the case of holomorphic bundles to holomorphic pairs and, in particular, refer the reader to Friedman and Morgan \cite[Section 4.2.1]{FrM} for the definition of \emph{(global) moduli functors} \cite[Section 4.2.1, Definition 2.1, p. 303]{FrM} and \emph{fine moduli spaces} \cite[Section 4.2.1, p. 304]{FrM} for holomorphic bundles and a discussion of what it means for global moduli functors to be representable and why one cannot expect that property for the moduli space of holomorphic bundles, let alone holomorphic pairs. For a discussion of moduli functors and representability from an algebro-geometric perspective (semistable sheaves over projective schemes over an algebraically closed field), we refer to Huybrechts and Lehn \cite[Chapter 4]{Huybrechts_Lehn_geometry_moduli_spaces_sheaves}.

Although many of the definitions and results of this section hold for general categories, we will state them for the category of complex analytic germs or complex analytic spaces for the sake of simplicity. The objects in the category of complex analytic spaces are complex analytic spaces (see Grauert and Remmert \cite[Section 1.1.5, p. 7]{Grauert_Remmert_coherent_analytic_sheaves}), and the morphisms in this category are holomorphic maps (see \cite[Section 1.1.4, p. 6]{Grauert_Remmert_coherent_analytic_sheaves}). Similarly, the objects in the category of germs of complex analytic spaces are the germs of complex analytic spaces and the morphisms are the germs of holomorphic maps. We will often define functors on the category of germs of complex analytic spaces by defining the value of the functor on a complex analytic space representing the germ, trusting this does not cause confusion.

\begin{defn}[The contravariant morphism functor from the category of complex analytic spaces to the category of sets]
\label{defn:Friedman_Morgan_4-2-3}
(See Friedman and Morgan in \cite[Section 4.2.1, Definition 2.3, p. 304]{FrM}.)
Given a complex analytic space $Z$, the \emph{contravariant morphism functor $h_Z$ from the category of complex analytic spaces to the category of sets} is defined as follows:
\begin{enumerate}
\item If $T$ is a complex analytic space, then $h_Z(T)$ is the set of morphisms from $T$ to $Z$.
\item For a morphism $f: T \to T'$, the corresponding map $h_Z(f): h_Z(T') \to h_Z(T)$ is defined by the pullback $f^*$.
\end{enumerate}
If $Z$ is a complex analytic germ, the functor $h_Z$ is defined in the same way.
\end{defn}

The following definitions are given for general categories in, for example, Lang \cite[Chapter I, Section 7, p. 30]{LangAlgebra}, Rotman \cite[Section 1.2, p. 26]{Rotman_Intro_to_HomologicalAlg_2009}, Fantechi et al. \cite[Definition 2.1, p. 14 and Definition 2.2, p. 15]{Fantechi_Goettsche_Illusie_Kleiman_Nitsure_Vistoli_FGA}, or Kashiwara and Schapira \cite[Definition 1.1.6, p. 25]{Kashiwara_Schapira_sheaves_manifolds}.

\begin{defn}[Representable functor and universal object]
\label{defn:RepresentableFunctor_UniversalObject}
A contravariant functor $F$ from the category of complex analytic spaces to the category of sets is \emph{representable} if it is isomorphic to the contravariant morphism functor $h_Z$ in Definition \ref{defn:Friedman_Morgan_4-2-3} for some complex analytic space $Z$, which in that case is unique up to a unique isomorphism.

Let $F$ be a contravariant functor as above.  A \emph{universal object for $F$} is a pair $(Z,\sE)$, where  $Z$ is a complex analytic space and $\sE$ is an element of $F(Z)$ with the property that the natural transformation $h_Z \to F$ given, for any complex analytic space $Y$, by
\[
h_Z(Y)\ni f\mapsto F(f)(\sE) \in F(Y)
\]
is a bijection of sets.
\end{defn}

\begin{exmp}[Property of being a universal object]
\label{exmp:UniversalObj}
Let $F$ be the contravariant functor from the category of complex analytic spaces to the category of sets obtained by defining $F(Y)$ to be the set of holomorphic vector bundles over the complex analytic space $Y$. If $Z$ is a complex analytic space and $\sE\in F(Z)$, that is, $\sE\to Z$ is a holomorphic vector bundle, then the pair $(Z,\sE)$ is a universal object for $F$ if for every complex analytic space $Y$, the natural transformation $h_Z\to F$ given by
\[
h_Z(Y) \ni f\mapsto f^*\sE\in F(Y),
\]
is a bijection of sets. Thus this natural transformation is an isomorphism of functors and we see that if $F$ has a universal object, then $F$ is a representable functor.
\end{exmp}

By generalizing the discussion in Example \ref{exmp:UniversalObj}, we see that if there is a universal object for a functor then the functor is representable. The opposite implication, that if a functor is representable then it has a universal object for general categories is given, for example, in Fantechi et al. \cite[Proposition 2.3, p. 15]{Fantechi_Goettsche_Illusie_Kleiman_Nitsure_Vistoli_FGA}.

Our goal in this section is proving that the following functor is representable.

\begin{defn}[Local moduli functor]
\label{defn:Friedman_Morgan_4-2-4}  
(Compare Friedman and Morgan in \cite[Section 4.2.1, Definition 2.4, p. 304]{FrM}.) Let $E$ be a smooth complex vector bundle over a closed, complex manifold $X$ with $\det E \cong X\times\CC$ (as smooth complex line bundles). Assume that $(\bar\partial_{E_0},\varphi_0)$ is a holomorphic pair on $E$ in the sense of \eqref{eq:Holomorphic_pair} and that $\bar\partial_{E_0}$ induces the usual $\bar\partial$-operator on $\det E$, so $\det E \cong \sO_X$ (thus, $\det E$ is holomorphically trivial), and that the pair $(\bar\partial_{E_0},\varphi_0)$ is strongly simple in the sense of Definition \ref{defn:Simple_pair}. The \emph{local moduli functor} $\bM_{\bar\partial_{E_0},\varphi_0}$ \label{page:Local_moduli_functor} at $(\bar\partial_{E_0},\varphi_0)$, from germs of complex analytic spaces to sets is defined as follows. For a germ $(T, t_0)$ of a complex analytic space, consider triples consisting of
\begin{enumerate}
\item A holomorphic vector bundle $\sE$ over $T \times X$ with $\det\sE \cong \sO_{T\times X}$ (trivial determinant).
\item An isomorphism of holomorphic vector bundles,
  \[
  \Upsilon: \sE\restriction (\{t_0\}\times X) \cong (E,\bar\partial_{E_0}).
  \]
\item A holomorphic section $\Phi \in \Omega^0(T\times X,\sE)$ such that
  \[
    \Upsilon\circ\Phi\restriction (\{t_0\}\times X) = \varphi_0 \in \Omega^0(E).
  \] 
\end{enumerate}
We say that two such triples, $(\sE_1,\Phi_1,\Upsilon_1)$ and $(\sE_2,\Phi_2,\Upsilon_2)$, are \emph{equivalent} if there is an isomorphism $\Xi:\sE_1 \to \sE_2$ such that
\[
  \Upsilon_2 \circ \Xi = \Upsilon_1 \quad\text{and}\quad \Xi \circ \Phi_1 = \Phi_2.
\]
We define $\bM_{\bar\partial_{E_0},\varphi_0}(T,t_0)$ to be the set of equivalence classes of these pairs. (As in \cite[Section 4.2.1, Definition 2.4, p. 304]{FrM}, we do not have to allow for twisting by the pullback of a line bundle on $T$ since all such line bundles are holomorphically trivial.) If $f: (T, t_0) \to (T', t_0')$ is an analytic map, then the map $\bM_{\bar\partial_{E_0},\varphi_0}(f)$ is defined by the pullback $(f\times\id_X)^*$.
\end{defn}

We prove that $\bM_{\bar\partial_{E_0},\varphi_0}$ is representable in Theorem \ref{thm:Friedman_Morgan_4-2-5_holomorphic_pairs} by constructing a universal object for $\bM_{\bar\partial_{E_0},\varphi_0}$, as in Example \ref{exmp:UniversalObj}.

\begin{rmk}[Determinant line bundles that are holomorphically trivial versus holomorphically fixed]
\label{rmk:Determinant_line_bundles_holomorphically_trivial_versus_fixed}
In our applications, we require $\det E$ to be a complex line bundle with a fixed holomorphic structure. At beginning of \cite[Section 4.1.6]{FrM}, Friedman and Morgan also restrict their attention to $\bar\partial_E$-operators on $E$ for which the induced operator on the complex line bundle $\det E$ is fixed. They then further assume that $\det E$ is trivial and that the induced operator on $\det E$ is the usual $\bar\partial$-operator (that is, $\det E$ is holomorphically trivial), remarking that the general case is handled by minor modifications. To preserve consistency with Friedman and Morgan \cite[Section 4.2.1, especially p. 304]{FrM} and for the purpose of our proof of the forthcoming Theorem \ref{thm:Friedman_Morgan_4-2-5_holomorphic_pairs}, we shall thus also assume that $\det E \cong X\times\CC$ as complex line bundles and that $\bar\partial_E$-operators on $E$ are required to induce the usual $\bar\partial$-operator on $X\times\CC$, so $(\det E,\bar\partial_{\det E})$ is holomorphically trivial. We shall leave the proof of the generalization from the case of holomorphically trivial to holomorphically fixed determinant line bundles to the reader.
\end{rmk}  

\begin{rmk}[Equivalence versus gauge transformations]
\label{rmk:LocalModuliEquivalence}
The notion of equivalence in Definition \ref{defn:Friedman_Morgan_4-2-4} is weaker than gauge equivalence by arbitrary smooth, determinant-one automorphisms of $E$. Because the local moduli functor is defined on germs, we can assume that the space $T$ is arbitrarily small. Thus, continuity of the isomorphism $\Xi$ implies that the restriction of $\Xi$ to $\{t\}\times X$ is a gauge transformation close to the fixed gauge transformation $\Upsilon_2^{-1}\circ\Upsilon_1$.
\end{rmk}

We can now state an analogue for holomorphic pairs of Kuranishi's Theorem \cite{Kuranishi}, proved by Friedman and Morgan for holomorphic bundles \cite[Section 4.2.1, Theorem 2.5, p. 305]{FrM}, on the existence of a local moduli space for the moduli functor; although it does not have a role in the proof of our main results in Chapter \ref{chap:Introduction}, we include it and its proof for completeness.

\begin{thm}[Kuranishi's Theorem on the existence of a local moduli space for the moduli functor for holomorphic pairs]
\label{thm:Friedman_Morgan_4-2-5_holomorphic_pairs}
Assume the hypotheses of Definition \ref{defn:Friedman_Morgan_4-2-4} and let $\fK(\bar\partial_{E_0},\varphi_0)$ be the Kuranishi model for $(\bar\partial_{E_0},\varphi_0)$. Then there exists a holomorphic vector bundle $\sE$ over $\fK(\bar\partial_{E_0},\varphi_0)\times X$ and a holomorphic section $\Phi\in\Omega^0(\fK(\bar\partial_{E_0},\varphi_0)\times X; \sE)$, such that $\det\sE$ is holomorphically trivial, with the property that the local moduli functor $\bM_{\bar\partial_{E_0},\varphi_0}$ is represented (in the sense of Remark \ref{rmk:Representability_local_moduli_functor}) by the germ $(\fK(\bar\partial_{E_0},\varphi_0),(\bar\partial_{E_0},\varphi_0))$ of the Kuranishi model at $(\bar\partial_{E_0},\varphi_0)$ together with $(\sE,\Phi)$.
\end{thm}

The property of the local moduli functor in Theorem \ref{thm:Friedman_Morgan_4-2-5_holomorphic_pairs} being representable is explained by Friedman and Morgan in more generality in \cite[Section 4.2.1, p. 304]{FrM} for moduli spaces of holomorphic vector bundles and moduli functors and we adapt their discussion in the following

\begin{rmk}[Representability of the local moduli functor]
\label{rmk:Representability_local_moduli_functor}
Assume the hypotheses of Definition \ref{defn:Friedman_Morgan_4-2-4}. Then the local moduli functor $\bM_{\bar\partial_{E_0},\varphi_0}$ is representable (in the sense of Definition \ref{defn:Friedman_Morgan_4-2-3}) if there exists a complex analytic space $\fK$ and a holomorphic vector bundle $\sE$ over $\fK\times X$ and holomorphic section $\Phi\in\Omega^0(\fK\times X;\sE)$, whose restriction to each slice $\{t\} \times X$ for $t\in\fK$ is a strongly simple holomorphic pair, with the following property: Given any family of strongly simple\footnote{Friedman and Morgan write `stable' rather than `simple' in the corresponding sentence \cite[Section 4.2.1, p. 304, line 22]{FrM} in their context of holomorphic bundles, but that appears to be a typographical error.} holomorphic pairs $(\sF,\Psi)$ over $T\times X$ parametrized by a complex analytic space $T$ with\footnote{We employ the abbreviations $\sF(t_0)=\sF|_{\{t_0\}\times X}$ and $\Psi(t_0)=\Psi|_{\{t_0\}\times X}$ here.}
$\sF(t_0) = (E,\bar\partial_{E_0})$ and $\Psi(t_0) = \varphi_0$ for some point $t_0\in T$ then, after possibly replacing $T$ by an open neighborhood of $t_0$, there are a unique morphism $f:T \to \fK$ of complex analytic spaces and a holomorphic line bundle $H$ over $T$ with holomorphic section $\eta \in \Omega^0(T;H)$ such that $\sF = (f\times\id_X)^*\sE\otimes\pr_\fK^*H$, where $\pr_\fK:\fK\times X\to\fK$ is projection onto the first factor, and $\Psi = (f\times\id_X)^*\Phi\otimes \pr_\fK^*\eta$.
\end{rmk}

We note that Theorem \ref{thm:Friedman_Morgan_4-2-5_holomorphic_pairs} only gives representability of the \emph{local} moduli functor $\bM_{\bar\partial_{E_0},\varphi_0}$.  As described in Friedman and Morgan \cite[Section 4.2.1, p. 304, paragraph following Definition 2.3]{FrM} for the case of simple holomorphic structures, if the \emph{global} moduli functor $\bM_E$ were representable, then there would be a universal object $(\fM,(\sE,\Phi))$ representing $\bM_E$ in the sense of Definition \ref{defn:RepresentableFunctor_UniversalObject}.  According to Friedman and Morgan \cite[Chapter IV, Theorem 4.3, p. 335]{FrM}, representability of their local moduli functor implies the moduli space of stable holomorphic bundles is a \emph{coarse} moduli space (see Definition \ref{defn:FrM_IV.4.2}).

\subsection{Proof of Kuranishi's Theorem in the smooth case}
\label{subsec:Friedman-Morgan_4-2-2}
In \cite[Section 4.2.2]{FrM}, Friedman and Morgan first prove \cite[Section 4.2.1, Theorem 2.5, p. 305]{FrM} in the special case where the Kuranishi model for a holomorphic bundle $(E, \bar\partial_E)$ is smooth at the point $[\bar\partial_E]$. We shall prove the following analogue for holomorphic pairs of their result.

\begin{thm}[Kuranishi's Theorem on the existence of a local moduli space for the moduli functor for holomorphic pairs near a smooth point]
\label{thm:Friedman_Morgan_4-2-5_holomorphic_pairs_smooth_point}
Assume the hypotheses of Theorem \ref{thm:Friedman_Morgan_4-2-5_holomorphic_pairs} and, in addition, that the reduction $\fK(\bar\partial_{E_0},\varphi_0)_\red$ of the Kuranishi model $\fK(\bar\partial_{E_0},\varphi_0)$ is smooth at the point $(\bar\partial_{E_0},\varphi_0)$. Then the conclusions of Theorem \ref{thm:Friedman_Morgan_4-2-5_holomorphic_pairs} hold. 
\end{thm}

We postpone the proof of Theorem \ref{thm:Friedman_Morgan_4-2-5_holomorphic_pairs_smooth_point}. Recall that a complex analytic space $T$ is \emph{smooth} at a point $t\in T$ if there is an open neighborhood $U$ of $t$ in $T$ such that the restriction of $T$ to that open neighborhood, $(U,\sO_T\restriction U)$, is isomorphic (in the sense of complex analytic spaces) to $(D,\sO_D)$, where $D \subset \CC^n$ is an open neighborhood of the origin (see Grauert and Remmert \cite[Section 1.1.4 and p. 8]{Grauert_Remmert_coherent_analytic_sheaves}). More explicitly, because $T$ is a complex analytic space, there are open neighborhoods $U \subset T$ of $t$ and $D\subset\CC^N$ of the origin and finitely many holomorphic functions $f_1,\ldots,f_k \in \sO(D)$ such that $(U,\sO_T\restriction U)$ is isomorphic to the \emph{complex model space} $(Y,\sO_Y)$, where $Y := \supp(\sO_D/\sI)$ and $\sO_Y := \sO_D/\sI$ with $\sI = \sO_D f_1 + \cdots + \sO_D f_k$ (see Grauert and Remmert \cite[Sections 1.1.2, 1.1.4, and 1.1.5]{Grauert_Remmert_coherent_analytic_sheaves}); then $T$ is \emph{smooth} at $t$ if and only if
\begin{equation}
  \label{eq:Smoothness_criterion_analytic_space}
  \rank_t(f_1,\ldots,f_k) : = \rank_\CC\left(\frac{\partial f_i}{\partial z_j}(t)\right)= k,
\end{equation}
where $z_1,\ldots,z_N$ are complex coordinates on $\CC^N$, in which case an open neighborhood of $t$ in $T$ is a complex manifold of dimension $n=N-k$ (see Griffiths and Harris \cite[p. 20]{GriffithsHarris}).

One says that a complex analytic space $T$ is \emph{reduced} at a point $t\in T$ if the stalk $\sO_{T,t}$ of its structure sheaf $\sO_T$ is a reduced ring (so it does not contain any non-zero nilpotent germs $s_t \in \sO_{T,t}$ such that $s_t^m=0$ for some integer $m\geq 2$) and \emph{non-reduced} at $t$ otherwise; $T$ is a \emph{reduced complex analytic space} if $T$ is reduced at all its points and \emph{non-reduced} otherwise (see Grauert and Remmert \cite[Section 1.1.5, p. 8]{Grauert_Remmert_coherent_analytic_sheaves}).

If $T$ is a complex analytic space with structure sheaf $\sO_T$, then $\sN_T \subset \sO_T$ denotes the \emph{nilradical} of $\sO_T$ and thus, by definition, the stalk of $\sN_{T,t} \subset \sO_{T,t}$ at each point $t\in T$ is the ideal of all nilpotent germs in $\sO_{T,t}$ (see Grauert and Remmert \cite[Section 4.2.5]{Grauert_Remmert_coherent_analytic_sheaves}). Thus $T$ is reduced if and only if $\sN_T$ is the zero ideal  (see Grauert and Remmert \cite[Section 4.3.3]{Grauert_Remmert_coherent_analytic_sheaves}).

If $T$ is a complex analytic space, one denotes
\[
  \sO_{T,\red} := \sO_T/\sN_T \quad\text{and}\quad T_\red := (T,\sO_{T,\red}), 
\]
and calls $T_\red$ the \emph{reduction} of $T$; this is a closed, reduced complex analytic subspace of $T$ and the embedding map $\red: T_\red \to T$ (called the \emph{reduction map}) comprises the identity map $T \to T$ of topological spaces and the residue epimorphism $\sO_T \to \sO_T/\sN_T$ (see Grauert and Remmert \cite[Section 4.3.2]{Grauert_Remmert_coherent_analytic_sheaves}). Every holomorphic map $f: S \to T$ between complex analytic spaces canonically determines a holomorphic map $f_\red: S_\red \to T_\red$ of the reductions and the diagram
\[
  \begin{CD}
  S_\red &@>{f_\red}>> T_\red
  \\
  @VV{\red}V &@VV{\red}V
  \\
  S &@>{f}>> T
\end{CD} 
\]  
is commutative (see Grauert and Remmert \cite[Section 4.3.2]{Grauert_Remmert_coherent_analytic_sheaves}). Note that $T$ is reduced if and only if $T = T_\red$ (see Grauert and Remmert \cite[Section 4.3.3]{Grauert_Remmert_coherent_analytic_sheaves}).

We begin with the construction of the holomorphic structure on the
universal bundle over $\fK(\bar\partial_{E_0},\varphi_0)\times X$. We have the following analogue for holomorphic pairs of \cite[Section 4.2.2, Proposition 2.6, p. 306]{FrM} for holomorphic bundles due to Friedman and Morgan.

\begin{prop}
\label{prop:Friedman_Morgan_4-2-6_holomorphic_pairs}
Continue the hypotheses of Theorem \ref{thm:Friedman_Morgan_4-2-5_holomorphic_pairs}. Suppose that $\fK(\bar\partial_{E_0},\varphi_0)$ is smooth, or more generally that its reduction $\fK(\bar\partial_{E_0},\varphi_0)_\red$ is smooth. Then there is a holomorphic structure $\sE$ on the smooth complex vector bundle $\pr_X^*E$ over $\fK(\bar\partial_{E_0},\varphi_0)_\red\times X$, where $\pr_X:\fK(\bar\partial_{E_0},\varphi_0)_\red\times X \to X$ denotes the projection onto the second factor, and a holomorphic section $\Phi \in \Omega^0(\fK(\bar\partial_{E_0},\varphi_0)_\red\times X;\sE)$ such that, for each $t \in \fK(\bar\partial_{E_0},\varphi_0)_\red$, the following hold:
\begin{enumerate}
\item The holomorphic structure on $E$ defined by restriction of $\pr_X^*E$ to the slice $\{t\} \times X$ is isomorphic to the holomorphic structure on $E$ defined by the point $t$.
  
\item The preceding isomorphism identifies the restriction of $\Phi$ to the slice $\{t\} \times X$ with the holomorphic section of $E$ defined by $t$.  
\end{enumerate}
If $(\sF,\Psi)$ is another holomorphic vector bundle with holomorphic section over $\fK(\bar\partial_{E_0},\varphi_0)_\red\times X$ having these properties, then $\sF$ is isomorphic to $\sE$  as a holomorphic bundle and this isomorphism identifies $\Psi$ with $\Phi$.
\end{prop}

\begin{proof}
  The assertions in the proposition that do not involve the section $\Phi$ are provided by \cite[Section 4.2.2, Proposition 2.6, p. 306]{FrM} and so it suffices to prove the additional assertions regarding the section $\Phi$. For clarity, we adopt the notation of the proof of \cite[Section 4.2.2, Proposition 2.6, p. 306]{FrM} as far as practical. As in \cite{FrM}, we first show the existence of the required holomorphic structure $\sE$ and section $\Phi$ locally over open sets of the form $\fK(\bar\partial_{E_0},\varphi_0)_\red\times N$, where $N \subset X$ is an open subset such that $E\restriction N \cong N \times \CC^r$ (isomorphic as smooth complex vector bundles). We let $t_1, \ldots, t_m$ be local holomorphic coordinates for $\fK(\bar\partial_{E_0},\varphi_0)_\red$ near $(\bar\partial_{E_0},\varphi_0)$. The inclusion
\[
  \fK(\bar\partial_{E_0},\varphi_0)_\red
  \subset
  \Omega^{0,1}(\fsl(E))\times\Omega^0(E) 
\]
and the identification
\[
  \sA^{0,1}(E) \times \Omega^0(E) = (\bar\partial_{E_0},\varphi_0) + \Omega^{0,1}(\fsl(E))\times\Omega^0(E)
\]
now define a family of $(0,1)$-pairs $(\bar\partial_E(t),\varphi(t))$ on $E$, depending holomorphically on $t \in \fK(\bar\partial_{E_0},\varphi_0)_\red$, with $\bar\partial_E(t)^2 = 0$ and $\bar\partial_E(t)\varphi(t) = 0$.

Let $\{s_i\}_{i=1}^r$ be a basis of $C^\infty$ sections for $E\restriction N$ defined by the preceding trivialization, and let $\{\pr_X^*s_i\}_{i=1}^r$ be the corresponding local basis of $C^\infty$ sections for $\pr_X^*E\restriction (\fK(\bar\partial_{E_0},\varphi_0)_\red\times N)$. Let $\sigma$ be a $C^\infty$ section of $\pr_X^*E$ over $\fK(\bar\partial_{E_0},\varphi_0)_\red\times N$ and write $\sigma = \sum_{i=1}^r f_i\pr_X^*s_i$, where $f_i\in C^\infty(\fK(\bar\partial_{E_0},\varphi_0)_\red\times N;\CC)$ for $i=1,\ldots,r$. We set
\[
  \bar\partial_\sE\sigma := \bar\partial_E(t)\sigma + \sum_{i=1}^r \bar\partial_tf_i \otimes \pr_X^*s_i, 
\]
where the operator $\bar\partial_t$ means that we compute $\bar\partial$ with respect to the $t$-variables only while $\bar\partial_E(t)\sigma$ is a section of the smooth vector bundle $\pr_X^*(\Lambda^{0,1}(E))$ whose coefficients with respect to the local basis $\{\pr_X^*s_i\}_{i=1}^r$ depend holomorphically on $t$ (since $\bar\partial_E(t)$ depends holomorphically of $t$).

The operator $\bar\partial_\sE$ can be shown to be independent of the choice of local basis $\{s_i\}_{i=1}^r$,
and thus defines an operator $\bar\partial_\sE$ on $\pr_X^*E$ over $\fK(\bar\partial_{E_0},\varphi_0)_\red\times X$. Moreover, $\bar\partial_\sE$ is a $(0,1)$-connection in the sense of \eqref{eq:Donaldson_Kronheimer_2-1-45_i} since the operator is clearly linear over $\CC$ and if $h$ is a smooth, complex-valued function on $\fK(\bar\partial_{E_0},\varphi_0)_\red\times X$ and $\sigma$ is a smooth section of $\pr_X^*E$ as above, then
\begin{align*}
  \bar\partial_\sE(h\sigma)
  &=
    \bar\partial_E(t)(h\sigma) + \sum_{i=1}^r \bar\partial_t(hf_i) \otimes \pr_X^*s_i
  \\
  &= \bar\partial_x h\otimes \sigma + h\bar\partial_E(t)\sigma
    + \sum_{i=1}^r \left( (\bar\partial_th)f_i + h\bar\partial_tf_i\right) \otimes \pr_X^*s_i
  \\
  &= \bar\partial_x h\otimes \sigma + \bar\partial_th\otimes\left(\sum_{i=1}^r f_i \otimes \pr_X^*s_i\right)
    + h\bar\partial_E(t)\sigma + h\sum_{i=1}^r \bar\partial_tf_i \otimes \pr_X^*s_i
  \\
  &= \bar\partial h\otimes \sigma + h\bar\partial_\sE\sigma,
\end{align*}
where the operator $\bar\partial_x$ means that we compute $\bar\partial$ with respect to the $x$-variables only, where $x\in X$.  We now verify that $\bar\partial_\sE^2=0$. We first observe that
\begin{align*}
  \bar\partial_\sE\sigma
  &=
    \bar\partial_E(t)\left(\sum_{i=1}^r f_i\pr_X^*s_i\right)
    + \sum_{i=1}^r \bar\partial_tf_i \otimes \pr_X^*s_i
  \\
  &= \sum_{i=1}^r \bar\partial_x f_i\otimes\pr_X^*s_i + f_i\bar\partial_E(t)\pr_X^*s_i
    + \bar\partial_tf_i \otimes \pr_X^*s_i
  \\
  &= \sum_{i=1}^r \left(\bar\partial_x f_i + \bar\partial_tf_i\right)\otimes\pr_X^*s_i
    + f_i\otimes \pr_X^*\bar\partial_E(t)s_i,
 \\
  &= \sum_{i=1}^r \bar\partial f_i\otimes\pr_X^*s_i + f_i\otimes \pr_X^*\bar\partial_E(t)s_i,
    \quad\text{for all } \sigma \in \Omega^0(\pr_X^*E).
\end{align*}
Consequently,
\begin{align*}
  \bar\partial_\sE^2\sigma
  &=
    \sum_{i=1}^r \bar\partial_\sE\left(\bar\partial f_i\otimes\pr_X^*s_i + f_i\otimes \pr_X^*\bar\partial_E(t)s_i\right)
  \\
  &= \sum_{i=1}^r \bar\partial^2 f_i\otimes\pr_X^*s_i - \bar\partial f_i\wedge \bar\partial_E(t)\pr_X^*s_i
    + \bar\partial f_i\wedge \pr_X^*\bar\partial_E(t)s_i
    + f_i\otimes \bar\partial_E(t)\pr_X^*\bar\partial_E(t)s_i
  \\
  &\qquad\text{(by Leibnitz Rule \eqref{eq:Donaldson_Kronheimer_2-1-45_i} for an $(0,1)$-connection)}
  \\
  &= \sum_{i=1}^r \bar\partial^2 f_i\otimes\pr_X^*s_i - \bar\partial f_i\wedge \pr_X^*\bar\partial_E(t)s_i
    + \bar\partial f_i\otimes \pr_X^*\bar\partial_E(t)s_i
    + f_i\otimes \pr_X^*\bar\partial_E(t)^2s_i
  \\
  &= 0, \quad\text{for all } \sigma \in \Omega^0(\pr_X^*E).
\end{align*}
Hence, the $(0,1)$-connection $\bar\partial_\sE$ is integrable and so there is a unique holomorphic structure on the smooth vector bundle $\pr_X^*E$ for which $\bar\partial_\sE$ is the $\bar\partial$-operator (see, for example, Friedman and Morgan \cite[Section 4.1.2, Lemma 1.6, p. 284]{FrM}.) We now define
\begin{equation}
  \label{eq:Definition_holomorphic_section_universal_holomorphic_bundle_smooth_point}
  \Phi(t,x) := \varphi(t)(x) \in E_x, \quad\text{for all } (t,x) \in \fK(\bar\partial_{E_0},\varphi_0)_\red\times X.
\end{equation}
Because $\Phi(t,x)$ is separately $C^\infty$ in $t$ and $x$, then $\Phi$ is a $C^\infty$ section of $\pr_X^*E$ over $\fK(\bar\partial_{E_0},\varphi_0)_\red\times X$ and we may write $\Phi = \sum_i h_i\pr_X^*s_i$. Over $\fK(\bar\partial_{E_0},\varphi_0)_\red\times N$, we have
\[
  \bar\partial_\sE\Phi = \bar\partial_E(t)\varphi(t) + \sum_{i=1}^r \bar\partial_th_i \otimes \pr_X^*s_i = 0,
\]
where the final equality follows from the fact that $\bar\partial_E(t)\varphi(t)=0$ and the fact that $\varphi(t)$ and hence the functions $h_i$ depend holomorphically on $t$, so that $\bar\partial_th_i = 0$ for all $i$. The remaining assertions involving $\Phi$ follow immediately.
\end{proof}

Given a holomorphic vector bundle $\sF$ over a product $T\times X$, where $T$ is a complex manifold, we would like to see that, at least locally on $T$, the bundle $\sF$ is induced by a holomorphic map from $T$ to $\fK(\bar\partial_{E_0},\varphi_0)$. The following lemma shows that $\sF$ is locally trivial in a $C^\infty$ sense.

\begin{lem}
\label{lem:Friedman_Morgan_4-2-7}
(See Friedman and Morgan \cite[Section 4.2.2, Lemma 2.7, p. 306]{FrM}.)
Let $X$ be a complex manifold, $(E,\bar\partial_E)$ be a holomorphic vector bundle over $X$, and $\sF$ be a holomorphic vector bundle over $T\times X$, where $T$ is a complex manifold and $t_0 \in T$ is a point. We assume that we are given an isomorphism $\iota:\sF\restriction(\{t_0\}\times X) \cong (E,\bar\partial_E)$ of holomorphic vector bundles.
Then, possibly after replacing $T$ by an open neighborhood of $t_0$, there is a $C^\infty$ isomorphism $\rho:\sF \to \pr_X^*E$ of smooth complex vector bundles with the following properties:
\begin{enumerate}
\item The restriction of $\rho$ to $\sF\restriction(\{t_0\}\times X)$ is equal to $\iota$.
  
\item For each $t\in T$, the restriction of $\rho$ to $\sF\restriction(\{t\}\times X)$ is holomorphic.
\end{enumerate}
If in addition we are given trivializations of $\det\sF$ and $\det E$ which agree under $\iota$, then we can choose $\rho$ to identify the given trivialization of $\det\sF$ with the trivialization of $\det\pr_X^*E$ induced by the trivialization of $\det E$.
\end{lem}

We can now adapt the argument due to Friedman and Morgan \cite[Section 4.2.2, pp. 307--308]{FrM} from the case of holomorphic bundles to holomorphic pairs and proceed to the

\begin{proof}[Proof of Theorem \ref{thm:Friedman_Morgan_4-2-5_holomorphic_pairs_smooth_point}]
Let $(\sE,\Phi)$ be the pair comprising the holomorphic bundle and section over $\fK(\bar\partial_{E_0},\varphi_0)_\red\times X$ provided by Proposition \ref{prop:Friedman_Morgan_4-2-6_holomorphic_pairs}. We need to show that the local moduli functor $\bM_{\bar\partial_{E_0},\varphi_0}$ is represented, in the sense of Remark \ref{rmk:Representability_local_moduli_functor}, by the germ
\[
  \left(\fK(\bar\partial_{E_0},\varphi_0),(\bar\partial_{E_0},\varphi_0)\right)
\]
of the Kuranishi model at $(\bar\partial_{E_0},\varphi_0)$ together with $(\sE,\Phi)$.
  
Therefore, suppose that (as in Remark \ref{rmk:Representability_local_moduli_functor}) we are given a complex manifold $T$, a point $t_0\in T$, and a holomorphic vector bundle $\sF$ and section $\Psi$ over $T\times X$. Fix an isomorphism $\iota:\sF\restriction(\{t_0\}\times X) \cong (E,\bar\partial_{E_0})$ of holomorphic vector bundles. Using the $C^\infty$ isomorphism $\rho$ provided by Lemma \ref{lem:Friedman_Morgan_4-2-7}, we can view the $\bar\partial$-operator on $\sF$ as an $(0,1)$-connection $D$ on $\pr_X^*E$ and the section $\Psi$ of $\sF$ as a section of $\pr_X^*E$. There is also the pullback $\pr_X^*\bar\partial_{E_0}$ on $\pr_X^*E$ of the $(0,1)$-connection $\bar\partial_{E_0}$ on $E$ and the pullback $\pr_X^*\varphi_0$ to $\pr_X^*E$ of the section $\varphi_0$ of $E$. We then obtain the difference\footnote{In \cite[Section 4.2.2, p. 307]{FrM}, there is a minor typographical error, where ``$S$'' should be replaced by ``$M$'' throughout.}
\[
  (D-\pr_X^*\bar\partial_{E_0}, \Phi - \pr_X^*\varphi_0)
  \in
  \Omega^{0,1}(T\times X; \pr_X^*\fsl(E)) \oplus \Omega^0(T\times X; \pr_X^*E).
\]
We can project this pair onto  $\Om^0(T\times X;\pr_X^*(\Lambda^{0,1}(\fsl(E)))\oplus \pr_X^*E)$ using the direct sum decomposition of vector spaces,
\begin{multline*}
  \Omega^{0,1}(T\times X; \pr_X^*\fsl(E)) \oplus \Omega^0(T\times X; \pr_X^*E)
  \\
  = 
  \Om^0(T\times X;(\pr_T^*\Lambda^{0,1}\oplus \pr_X^*\Lambda^{0,1})\otimes \pr_X^*\fsl(E))
  \oplus
  \Om^0(T\times X;\pr_X^*E)
  \\
  =
  \Om^0(T\times X;\pr_X^*(\Lambda^{0,1}(\fsl(E)))\oplus \pr_X^*E)
\end{multline*}
As pointed out in \cite[Section 4.2.2, p. 307]{FrM}, this is equivalent to restricting the $(0,1)$-form component of the pair to tangent vectors pointing in the $X$ direction\footnote{Our notation differs from that in \cite{FrM} where they project the $(0,1)$-form component to $\pr_X^*(\Omega^{0,1}(X; \fsl(E))$, elements of which would be constant with respect to $t\in T$.}. Thus, one sees (as in \cite[Section 4.2.2, p. 307]{FrM}) that this projection of
\[
  (D-\pr_X^*\bar\partial_{E_0}, \Phi - \pr_X^*\varphi_0)
\]
is a $C^\infty$ section $s$ of $\Om^0(T\times X;\pr_X^*(\Lambda^{0,1}(\fsl(E)))\oplus \pr_X^*E)$, which thus defines a map
\[
  f: T\to \Omega^{0,1}(X;\fsl(E))\oplus \Omega^0(X;E)
\]
by the formula
\[
  f(t) := s(t,\cdot), \quad\text{for all } t \in T.
\]
Using the identification of $\Omega^{0,1}(X;\fsl(E))$ with the affine space $\sA^{0,1}(E)$ given by taking $\bar\partial_{E_0}$ as the origin in $\sA^{0,1}(E)$, we can also view $f$ as taking values in $\sA^{0,1}(E)\times \Omega^0(X;E)$.

The image of $f$ is contained in the subspace $\fS^{-1}(0)$ corresponding to holomorphic $(0,1)$-pairs. Friedman and Morgan note that a straightforward argument shows that $f$ is holomorphic in their setting of holomorphic bundles and so equally straightforward in our setting of holomorphic pairs. Thus we have defined a holomorphic map
\[
  f: T \to \sA^{0,1}(E)\times\Omega^0(E), \quad\text{with } f(t_0) = (\bar\partial_{E_0}, \varphi_0).
\]
Consequently, possibly after shrinking $T$, there is a holomorphic map from $T$ to $\fK(\bar\partial_{E_0}, \varphi_0)$, which we shall continue to denote by $f$. Since $T$ is reduced, the image of $f$ lies in the closed subspace $\fK(\bar\partial_{E_0}, \varphi_0)_\red$ of $\fK(\bar\partial_{E_0}, \varphi_0)$.

The proof that $\sF \cong (f\times\id_X)^*\sE$ follows without change from that in \cite[Section 4.2.2, pp. 307--308]{FrM} and proof of the equality $\Psi = (f\times\id_X)^*\Phi$ presents no new difficulties. The proof of uniqueness of the map $f$ follows exactly as in \cite[Section 4.2.2, pp. 307--308]{FrM}.
\end{proof}

\subsection{Proof of Kuranishi's Theorem in the general  case: Construction of the universal holomorphic bundle and section}
\label{subsec:Friedman-Morgan_4-2-3}
In this section, we begin the proof of Theorem \ref{thm:Friedman_Morgan_4-2-5_holomorphic_pairs}, extending the argument used to prove the simpler Theorem \ref{thm:Friedman_Morgan_4-2-5_holomorphic_pairs_smooth_point}, where the complex analytic space $\fK(\bar\partial_{E_0}, \varphi_0)_\red$ was assumed to be smooth at $[\bar\partial_{E_0}, \varphi_0]$. Again, we shall adapt the proof of Friedman and Morgan (see \cite[Section 4.2.3]{FrM})  from the case of holomorphic bundles to holomorphic pairs, noting that they in turn base their argument on that of Douady \cite{Douady_1966_sem_bourbaki}. In this section, we shall construct a universal holomorphic vector bundle over $\fK(\bar\partial_{E_0}, \varphi_0)\times X$ and prove

\begin{thm}[Existence of a universal holomorphic vector bundle and holomorphic section]
\label{thm:Friedman_Morgan_4-2-8_holomorphic_pairs}
Assume the hypotheses of Theorem \ref{thm:Friedman_Morgan_4-2-5_holomorphic_pairs}. Then there exists a holomorphic vector bundle structure $\sE$ on the vector bundle $\pr_X^*E$ over $\fK(\bar\partial_{E_0}, \varphi_0)\times X$ and a holomorphic section $\Phi \in \Omega^0(\fK(\bar\partial_{E_0}, \varphi_0)\times X;\sE)$, such that, for every $t \in \fK(\bar\partial_{E_0}, \varphi_0)$, the following hold:
\begin{enumerate}
\item The holomorphic structure on $E$ defined by restriction of $\pr_X^*E$ to the slice $\{t\} \times X$ is isomorphic to the holomorphic structure on $E$ defined by the point
\[
  t \in \fK(\bar\partial_{E_0}, \varphi_0) \subset \fS^{-1}(0) \subset \sA^{0,1}(E)\times \Omega^0(E).
\]
\item The preceding isomorphism identifies the restriction of $\Phi$ to the slice $\{t\} \times X$ with the holomorphic section of $E$ defined by $t$.  
\end{enumerate}
\end{thm}

Our proof of Theorem \ref{thm:Friedman_Morgan_4-2-8_holomorphic_pairs}, which we give further below, is modeled on that of Friedman and Morgan for their \cite[Section 4.2.3, Theorem 2.8, p. 308]{FrM}. Rather than reproduce their proof\footnote{Given in \cite[Section 4.2.3, pp. 308--314]{FrM}, with many expository details.} in full, we shall instead focus on the construction of the holomorphic section $\Phi$ and adapt their arguments in \cite[Section 4.2.3, pp. 313--314]{FrM}. Recall that if $T$ is a complex analytic space with structure sheaf $\sO_T$, then sections of $\sO_T$ are customarily called \emph{holomorphic functions} (see Grauert and Remmert \cite[p. 9]{Grauert_Remmert_coherent_analytic_sheaves} or Gunning and Rossi \cite[p. 150]{Gunning_Rossi_analytic_functions_several_complex_variables}); a morphism of complex analytic spaces is called a \emph{holomorphic map} and an isomorphism is called a \emph{biholomorphic map} (see \cite[p. 7]{Grauert_Remmert_coherent_analytic_sheaves} or \cite[p. 150]{Gunning_Rossi_analytic_functions_several_complex_variables}).

Suppose that $T$ is a complex model space, as in the discussion around the smoothness criterion \eqref{eq:Smoothness_criterion_analytic_space}, so there are a domain $D \subset \CC^N$ and finitely generated ideal $\sI \subset \sO_D$ such that $\sO_T = \sO_D/\sI$ and $T \subset D$ is a closed subset. If $s$ is a section of $\sO_T$, then (see Grauert and Remmert \cite[Section 1.1.6, p. 9]{Grauert_Remmert_coherent_analytic_sheaves}) there is a holomorphic function $f \in \sO_D(D)$ such that $s = \check f\restriction T$, where $\check f$ denotes the image of $f$ in $\sO_D(D)/\sI(D)$. (In other words, $s$ is holomorphic at a point of the closed subset $T\subset D$ if and only if it is the restriction of a holomorphic function on an open neighborhood of that point in $\CC^N$.) The section $\Phi$ of the holomorphic vector bundle in the statements of Theorems \ref{thm:Friedman_Morgan_4-2-5_holomorphic_pairs} or \ref{thm:Friedman_Morgan_4-2-8_holomorphic_pairs} is understood to be \emph{holomorphic} in this sense.

\begin{proof}[Proof of Theorem \ref{thm:Friedman_Morgan_4-2-8_holomorphic_pairs}]
The holomorphic vector bundle structure $\sE$ on $\pr_X^*E$ is provided by \cite[Section 4.2.3, Theorem 2.8, p. 308]{FrM}. Indeed, as in \cite[Section 4.2.3, p. 313]{FrM}, we may assume without loss of generality that $E$ and $X$ are equipped with real analytic, Hermitian metrics and that the pair $(\bar\partial_{E_0}, \varphi_0)$ is real analytic. Consequently, the pairs
\[
    (\alpha,\sigma) \in \bH_{\bar\partial_{E_0}, \varphi_0}^1 \subset
    \Omega^{0,1}(\fsl(E))\oplus \Omega^0(E)
\]
are real analytic by a well-known result of Morrey \cite{Morrey_1958_part1, Morrey_1958_part2} since they are solutions to a system of nonlinear elliptic partial differential equations with real analytic coefficients. Theorem \ref{thm:Friedman_Morgan_4-2-5_holomorphic_pairs_smooth_point} now applies with $\fK(\bar\partial_{E_0}, \varphi_0)$ replaced by an open neighborhood of the origin (corresponding to the point $(\bar\partial_{E_0}, \varphi_0)$) in the harmonic space $\bH_{\bar\partial_{E_0}, \varphi_0}^1$.

The conclusions of Theorem \ref{thm:Friedman_Morgan_4-2-8_holomorphic_pairs} regarding $\sE$ are obtained in \cite[Section 4.2.3, pp. 313--314]{FrM} essentially by restriction to the complex analytic subspace
\[
  \fK(\bar\partial_{E_0}, \varphi_0) \subset \bH_{\bar\partial_{E_0}, \varphi_0}^1.
\]
When $\fK(\bar\partial_{E_0}, \varphi_0)$ is smooth at the point $[\bar\partial_{E_0}, \varphi_0]$, we defined the section $\Phi$ in \eqref{eq:Definition_holomorphic_section_universal_holomorphic_bundle_smooth_point} and established that it is holomorphic with the required properties in our proof of Theorem \ref{thm:Friedman_Morgan_4-2-5_holomorphic_pairs_smooth_point}. When $\fK(\bar\partial_{E_0}, \varphi_0)$ is only assumed to be a complex analytic space, we obtain a holomorphic section $\Phi$ by restriction to $\fK(\bar\partial_{E_0}, \varphi_0)$ of the holomorphic section on an open neighborhood of the origin in the affine harmonic space $\bH_{\bar\partial_{E_0}, \varphi_0}^1$.
\end{proof}


\subsection{Proof of Kuranishi's Theorem in the general  case: Conclusion}
\label{subsec:Friedman-Morgan_4-2-4}
In this section, we conclude the proof of Theorem \ref{thm:Friedman_Morgan_4-2-5_holomorphic_pairs}, again adapting the exposition due to Friedman and Morgan \cite[Section 4.2.4]{FrM} from the case of holomorphic bundles to that of holomorphic pairs. Their argument is lengthy (occupying \cite[Section 4.2.3, pp. 314--318]{FrM}), so we shall just describe the additional details needed to construct the universal holomorphic section.

\begin{proof}[Conclusion of the proof of Theorem \ref{thm:Friedman_Morgan_4-2-5_holomorphic_pairs}]
We assume that we are given a holomorphic pair $(\bar\partial_{E_0},\varphi_0)$ on a complex vector bundle $E$ over a closed, complex manifold $X$; a holomorphic trivialization of $(E,\bar\partial_{E_0})$; and a complex analytic space $T$ and point $t_0 \in T$. Let $\sF \to T\times X$ be a holomorphic vector bundle and $\Psi \in \Omega^0(T\times X;\sF)$ be a holomorphic section such that the restriction of the pair $(\sF,\Psi)$ to every slice $\{t\} \times T$ is a simple pair, together with a holomorphic isomorphism $\psi: \sF\restriction(\{t_0\}\times X) \to (E,\bar\partial_{E_0})$ that identifies $\Psi(t_0,\cdot)$ with $\varphi_0$. We shall also suppose that $\det\sF$ is holomorphically trivial and fix a trivialization of $\det\sF$ that is compatible with the trivialization of $\det\sF$ induced by $\psi$. As in \cite[Section 4.2.4]{FrM}, there are three remaining steps to the proof of Theorem \ref{thm:Friedman_Morgan_4-2-5_holomorphic_pairs}:
\begin{enumerate}
\item Show that the holomorphic pair $(\sF,\Psi)$ defines a holomorphic map $\tilde f$ from $T$ to $\sA^{0,1}(E)\times\Omega^0(E)$ whose image $\Imag\tilde f$ is contained in $\fS^{-1}(0)$ in the sense of complex analytic spaces. Then, by Proposition \ref{prop:Friedman_Morgan_4-1-17_holomorphic_pairs}, possibly after shrinking $T$, there is a holomorphic map $u:T \to W^{2,p}(\SL(E)))$ such that the map $f = u^*\tilde f$ has image $\Imag f$ contained in $\fS^{-1}(0)$ in the sense of complex analytic spaces.
  
\item Show that $\sF \cong (f \times \id_X)^*\sE$ as holomorphic vector bundles over an open neighborhood of $t_0$ and that this isomorphism identifies $\Psi$ with $\Phi$.
  
\item Show that, if $f_1$ and $f_2$ are two holomorphic maps from $T$ to $\fK(\bar\partial_{E_0},\varphi_0)$ such that there is an automorphism of holomorphic vector bundles $(f_1\times\id_X)^*\sE \cong (f_2\times\id_X)^*\sE$ that identifies $\Phi$ with itself then, possibly after shrinking $T$, one has $f_1=f_2$.
\end{enumerate}

\begin{step}[$\sF$ defines a holomorphic map $\tilde f:T\to \sA^{0,1}(E)\times\Omega^0(E)$ such that $\Imag f \subset \fS^{-1}(0)$ in the sense of complex analytic spaces]
In the case of a holomorphic vector bundle $\sF$, rather than holomorphic pair $(\sF,\Psi)$ as here, Friedman and Morgan complete this step in \cite[Section 4.2.4, pp. 315--317]{FrM}. No significant change in their argument is required by consideration of the holomorphic section $\Psi$, with their application of \cite[Section 4.1.5, Proposition 1.17, p. 297]{FrM} is replaced by an application of our Proposition \ref{prop:Friedman_Morgan_4-1-17_holomorphic_pairs}.
\end{step}

\begin{step}[$(\sF,\Psi) \cong (f \times \id_X)^*(\sE,\Phi)$ as holomorphic pairs over an open neighborhood of $t_0$ in $T$]
In the case of a holomorphic vector bundle $\sF$, rather than a holomorphic pair $(\sF,\Psi)$ as here, Friedman and Morgan complete this step in \cite[Section 4.2.4, pp. 317--318]{FrM}. Again, no significant change in their argument is required when passing from consideration of $\sF$ alone to that of the pair $(\sF,\Psi)$. 
\end{step}

\begin{step}[For holomorphic maps $f_1, f_2:T \to\fK(\bar\partial_{E_0},\varphi_0)$ that induce an isomorphism of holomorphic pairs, $(f_1\times\id_X)^*(\sE,\Phi) \cong (f_2\times\id_X)^*(\sE,\Phi)$, show that $f_1=f_2$]
Friedman and Morgan complete this short step in \cite[Section 4.2.4, p. 318]{FrM} and no change is required when passing from consideration of $\sF$ alone to that of the pair $(\sF,\Psi)$. 
\end{step}

This completes the proof of Theorem \ref{thm:Friedman_Morgan_4-2-5_holomorphic_pairs}.
\end{proof}

\section{Global moduli functor for holomorphic pairs}
\label{subsec:Friedman-Morgan_4-4-1}
We compare the local moduli functor in Definition \ref{defn:Friedman_Morgan_4-2-4} with the following adaptation of \cite[Section 4.4.1, Definition 4.1, p. 334]{FrM} due to Friedman and Morgan for the global moduli functor for holomorphic vector bundles.

\begin{defn}[Global moduli functor for non-zero-section, stable holomorphic pairs]
\label{defn:FrM_Definition_IV.4.1_pairs}
Let $E$ be a complex vector bundle over closed, complex K\"ahler manifold $(X,\omega)$ with $\det E \cong X\times\CC$ (as smooth complex vector bundles). The \emph{(global) moduli functor} $\bM_E^0$ \label{page:Global_moduli_functor} for non-zero-section, stable holomorphic pairs from the category of complex analytic spaces to the category of sets is defined as follows. For any complex analytic space $T$, consider triples consisting of
\begin{enumerate}
\item A holomorphic vector bundle $\sE$ over $T \times X$ with $\det\sE \cong \sO_{T\times X}$ (trivial determinant).
  
\item An isomorphism of holomorphic vector bundles,
\[
  \Upsilon: \sE\restriction (\{t\}\times X) \cong (E,\bar\partial_{E_t}), 
\]
where $\bar\partial_{E_t}$ is a holomorphic structure on $E$.

\item A holomorphic section $\Phi \in \Omega^0(T\times X,\sE)$ such that
\[
  \Upsilon\circ\Phi\restriction (\{t\}\times X) = \varphi_t \in \Omega^0(E),
\]
where $\varphi_t \not\equiv 0$ and $\varphi_t$ is holomorphic with respect to $\bar\partial_{E_t}$ and the pair $(\bar\partial_{E_t},\varphi_t)$ is stable in the sense of Definition \ref{defn:StablePair_arbitrary_rank}.
\end{enumerate}
We define $\bM_E^0(T)$ to be the set of equivalence classes of these pairs, where equivalence is in the sense of Definition \ref{defn:Friedman_Morgan_4-2-4}.
\end{defn}

We adapt \cite[Section 4.4.1, Definition 4.2, p. 335]{FrM} due to Friedman and Morgan for holomorphic vector bundles to give the

\begin{defn}[Coarse moduli space of non-zero-section, stable holomorphic pairs]
\label{defn:FrM_IV.4.2}
(Compare Friedman and Morgan \cite[Section 4.4.1, Definition 4.2, p. 335]{FrM}.)    
A \emph{coarse moduli space} for the moduli functor $\bM_E^0$ in Definition \ref{defn:FrM_Definition_IV.4.1_pairs} is
a complex analytic space $\fM$ with the following properties:
\begin{enumerate}
\item
There is a morphism of functors $f:\bM_E^0\to h_\fM$, where $h_\fM$ is the functor given in Definition \ref{defn:Friedman_Morgan_4-2-3}.
\item
If $T$ is any complex analytic space and $f':\bM_E^0\to h_N$ is a morphism of functors, then there is a unique morphism
$\phi:\bM_E^0\to N$ which induces a natural morphism $\phi_*:h_\fM\to h_N$ such that $f'=\phi_*\circ f$.
\item
The points of $\fM$ are in one-to-one correspondence with the set of isomorphism classes of non-zero-section, stable holomorphic pairs $(\bar\partial_E,\varphi)$ on $E$.
\end{enumerate}
\end{defn}

Because the proof of the corresponding result for stable holomorphic bundles in Friedman and Morgan \cite[Section 4.4.1, Theorem 4.3, p. 335]{FrM} is primarily in the setting of category theory, it translates to give the following analogue for stable holomorphic pairs; it does not have a role in the proof of our main results in Chapter \ref{chap:Introduction}, but is included for completeness.

\begin{thm}[Coarse moduli space of non-zero-section, stable holomorphic pairs]
\label{thm:FrM_IV.4.3}  
Let $E$ be a smooth, complex vector bundle with trivial determinant over a closed, complex K\"ahler surface $(X,\omega)$ and let $\bar\partial_{E_d}$ be a fixed holomorphic structure on the complex line bundle $\det E$. Then the moduli space $\fM^0(E,\omega)$ of non-zero-section, stable holomorphic pairs is a coarse moduli space for the global moduli functor $\bM_E^0$.
\end{thm}

\begin{proof}
The primary deviation from the proof of \cite[Section 4.4.1, Theorem 4.3, p. 335]{FrM} is that we apply our Kuranishi Theorem \ref{thm:Friedman_Morgan_4-2-5_holomorphic_pairs} for holomorphic pairs rather than \cite[Section 4.2.1, Theorem 2.5, p. 305]{FrM}. We restrict to non-zero-section, stable holomorphic pairs in order to make use of the fact that such pairs are strongly simple by Lemma \ref{lem:Stable_implies_strongly_simple_for_01-pairs}, as required by the hypotheses of Theorem \ref{thm:Friedman_Morgan_4-2-5_holomorphic_pairs}.
\end{proof}

See Remark \ref{rmk:Determinant_line_bundles_holomorphically_trivial_versus_fixed} for a discussion of modifications needed to pass from the case of determinant line bundles that are holomorphically trivial, as assumed in this section, to holomorphically fixed.

\section[Comparison of Kuranishi models]{Comparison of Kuranishi models for projective vortices and holomorphic pairs}
\label{subsec:Friedman-Morgan_4-3-4}
In this section we extend the arguments due to Friedman and Morgan in \cite[Section 4.3.4]{FrM}, which compare the analytic structures on the moduli space of anti-self-dual connections $A$ on a Hermitian vector bundle $E$ over a smooth, closed, complex projective surface $S$ and the moduli space of holomorphic structures $\bar\partial_A$ on $E$, to a comparison of the analytic structures on the moduli space of projective vortices $(A,\varphi)$ on $E$ over a closed, K\"ahler manifold $X$ and the moduli space of holomorphic pairs $(\bar\partial_A,\varphi)$ on $E$. As Friedman and Morgan observe at the beginning of \cite[Section 4.3.4]{FrM}, when $[A]$ is a smooth point of its moduli space with trivial stabilizer in $C^\infty(\SU(E))$, the argument is straightforward and the same is true for a point $[A,\varphi]$ in the moduli space of projective vortices. Indeed, from our comparison of the
elliptic deformation complexes provided by Theorem \ref{thm:Kobayashi_7-2-21_pairs}
if the harmonic spaces $\bH_{A,\varphi}^0$ and $\bH_{A,\varphi}^2$ in \eqref{eq:H_Avarphi^bullet} are each zero, then the harmonic spaces $\bH_{\bar\partial_A,\varphi}^0$ and $\bH_{\bar\partial_A,\varphi}^2$ in \eqref{eq:H_dbar_Avarphi^0bullet} vanish also and, moreover, the map $(A,\varphi) \mapsto (\bar\partial_A,\varphi)$ induces an isomorphism $\bH_{A,\varphi}^1 \cong \bH_{\bar\partial_A,\varphi}^1$ of real vector spaces, defined by \eqref{eq:H_Avarphi^bullet} and \eqref{eq:H_dbar_Avarphi^0bullet} respectively. The Inverse Mapping Theorem for real analytic maps (see Section \ref{sec:Real_complex_analytic_functions_Banach_spaces})
now yields a real analytic diffeomorphism from an open neighborhood of $[A,\varphi]$ in its moduli space onto an open neighborhood of $[\bar\partial_A,\varphi]$ in its moduli space. Thus, by analogy with \cite[Section 4.3.4]{FrM}, we shall focus on the case where $[A,\varphi]$ need not be a smooth point
in its moduli space.

We first prove the following analogue of \cite[Section 4.3.4, Theorem 3.8, p. 326]{FrM}, attributed by Friedman and Morgan to Donaldson \cite[Section 3, pp. 154--157]{DonHCobord}.

\begin{thm}[Real analytic isomorphism from the Kuranishi model associated to a projective vortex onto the Kuranishi model associated to a holomorphic pair]
\label{thm:Friedman_Morgan_4-3-8_projective_vortices}
Assume the hypotheses of Theorem \ref{thm:Lubke_Teleman_6-3-7} and let $(A,\varphi)$ be a smooth projective vortex on $E$ with stabilizer $\{\id_E\}$ in $W^{2,p}(\SU(E))$. Let $(\bar\partial_A,\varphi)$ be the corresponding smooth holomorphic pair on $E$ and assume that it is strongly simple, so has stabilizer $\{\id_E\}$ in $W^{2,p}(\SL(E))$. Then there are open neighborhoods of the origin $N_{\bar\partial_A,\varphi} \subset \bH_{\bar\partial_A,\varphi}^1$ and $N_{A,\varphi}\subset \bH_{A,\varphi}^1$ such that the following assertions hold. Let $\beps:N_{A,\varphi} \to \sA(E,h) \times \Omega^0(E)$ and $\bkappa:N_{A,\varphi} \to \bH_{A,\varphi}^2$ be the real analytic Kuranishi embedding and obstruction maps provided by Theorem \ref{thm:Local_Kuranishi_model_for_moduli_space_projective_vortices_StabAvarphi_idE} associated to the affine local slice \eqref{eq:Unitary_pair_Avarphi_linear_slice}, namely
\[
  S_{A,\varphi} = (A,\varphi) + \Ker d_{A,\varphi}^{0,*} \subset \sA(E,h) \times \Omega^0(E),
\]
and the projective vortex equations \eqref{eq:SO(3)_monopole_equations_almost_Hermitian_alpha}, and denote the germ of this real analytic space by $\sK(A,\varphi)$. Similarly, let $\bgamma:N_{\bar\partial_A,\varphi} \to \sA^{0,1}(E) \times \Omega^0(E)$ and $\bchi:N_{\bar\partial_A,\varphi} \to \bH_{\bar\partial_A,\varphi}^2$ be the complex analytic Kuranishi embedding and obstruction maps given in Definition \ref{defn:Friedman_Morgan_4-1-15_holomorphic_pairs} and associated to the affine local slice \eqref{eq:dbar_Evarphi_slice} through the pair $(\bar\partial_A,\varphi)$, namely
\[
  S_{\bar\partial_A,\varphi}
  =
  (\bar\partial_A,\varphi) + \Ker \bar\partial_{A,\varphi}^{0,*}\cap\left(\Omega^{0,1}(\fsl(E))\oplus\Omega^0(E)\right)
   \subset \sA^{0,1}(E) \times \Omega^0(E),
\]
and the holomorphic pair equations \eqref{eq:Holomorphic_pair}. Then the
isomorphism \eqref{eq:Bijection_unitarypairs_with_01pairs} of real affine spaces,
\[
  \pi_h^{0,1}:\sA(E,h)\times\Omega^0(E) \ni (A,\varphi) \mapsto (\bar\partial_A,\varphi)  \in \sA^{0,1}(E)\times\Omega^0(E),
\]
induces an isomorphism of real analytic spaces from the Kuranishi model $\sK(A,\varphi)$ onto the Kuranishi model $\fK(\bar\partial_A,\varphi)$ in Definition \ref{defn:Friedman_Morgan_4-1-15_holomorphic_pairs} and thus an isomorphism of real analytic sets,
\begin{equation}
  \label{eq:Real_analytic_isomorphism_beps(preimage_bkappa_0)_onto_bgamma(preimage_bchi_0)}
  \beps\left(\bkappa^{-1}(0)\right) \cong \bgamma\left(\bchi^{-1}(0)\right). 
\end{equation}
\end{thm}

\begin{proof}
The argument is essentially that of Friedman and Morgan but there are sufficient differences in details to justify adapting their exposition in full. The Kuranishi model in Definition \ref{defn:Friedman_Morgan_4-1-15_holomorphic_pairs} provides a local, complex analytic model space for the holomorphic pair $(\bar\partial_A,\varphi)$ corresponding to the open neighborhood $N_{\bar\partial_A,\varphi}$ of the origin in $\bH_{\bar\partial_A,\varphi}^1$ via the holomorphic embedding $\bgamma:N_{\bar\partial_A,\varphi}\to \sA^{0,1}(E)\times\Omega^0(E)$ and complex analytic obstruction map $\bchi:\bH_{\bar\partial_A,\varphi}^1\supset N_{\bar\partial_A,\varphi}\to\bH_{\bar\partial_A,\varphi}^2$ defined by applying the Implicit Mapping Theorem and Kuranishi Method near $(\bar\partial_A,\varphi)$ to the restriction of the complex analytic map $\fS$ in \eqref{eq:Holomorphic_pair_and_Coulomb_gauge_map_sheaves} or equivalently $\widehat\fS$ in \eqref{eq:Holomorphic_pair_map_Banach_spaces_orthogonal_projections} to the affine local slice $S_{\bar\partial_A,\varphi}$.

The corresponding Kuranishi model for the projective vortex $(A,\varphi)$ provided by Theorem \ref{thm:Local_Kuranishi_model_for_moduli_space_projective_vortices_StabAvarphi_idE} is a local, real analytic model space for an open neighborhood of the point $[A,\varphi]$ in $\sM^{**}(E,h,\omega)$ corresponding to the open neighborhood $N_{A,\varphi}$ of the origin in $\bH_{A,\varphi}^1$ via the real analytic embedding $\beps:\bH_{A,\varphi}^1\supset N_{A,\varphi}\to \sA(E,h)\times\Omega^0(E)$ and real analytic obstruction map $\bkappa:\bH_{A,\varphi}^1\supset N_{A,\varphi}\to\bH_{A,\varphi}^2$. This model is defined by applying the Implicit Mapping Theorem and Kuranishi Method near $(A,\varphi)$ to the real analytic map $\sS$ in \eqref{eq:Projective_vortex_map} defined by the restriction of the left-hand side of the projective vortex equations \eqref{eq:SO(3)_monopole_equations_almost_Hermitian_alpha} to the affine local slice $S_{A,\varphi}$, which we write here as
\begin{multline}
  \label{eq:Projective_vortex_and_Coulomb_gauge_map_sheaves}
  \sS:\sA(E,h)\times \Omega^0(E) \supset S_{A,\varphi} \ni (A+a,\varphi+\sigma)
  \\
  \mapsto
  \left( (\Lambda F_{A+a})_0 - \frac{i}{2}((\varphi+\sigma)\otimes(\varphi+\sigma)^*)_0,
    (F_{A+a}^{0,2})_0, \bar\partial_{A+a}(\varphi+\sigma) \right)
  \\
  \in \Omega^0(\su(E))\oplus \Omega^{0,2}(\fsl(E)) \oplus \Omega^{0,1}(E).
\end{multline}
Equivalently, and courtesy of the forthcoming Lemma \ref{lem:Friedman_Morgan_1999_1-2_projective_vortices} --- which is the analogue for projective vortices of Corollary \ref{cor:Friedman_Morgan_1999_1-2_restricted_to_Coulomb_gauge_slice_neighborhood} --- one may apply the Implicit Mapping Theorem and Kuranishi Method near $(A,\varphi)$ to the real analytic map
\begin{multline}
  \label{eq:Projective_vortex_and_Coulomb_gauge_map_Banach_spaces_orthogonal_projections}
  \hat\sS:\sA(E,h)\times \Omega^0(E ) \supset S_{A,\varphi} \ni (A+a,\varphi+\sigma)
  \\
  \mapsto
  \left( (\Lambda F_{A+a})_0 - \frac{i}{2}((\varphi+\sigma)\otimes(\varphi+\sigma)^*)_0,
    \Pi_{A,\varphi}^2\left((F_{A+a}^{0,2})_0, \bar\partial_{A+a}(\varphi+\sigma)\right) \right)
  \\
  \in \Omega^0(\su(E))\oplus \Ker d_{A,\varphi}^2\cap\left(\Omega^{0,2}(\fsl(E)) \oplus \Omega^{0,1}(E)\right),
\end{multline}
where the differential $d_{A,\varphi}^2$ is defined in \eqref{eq:d2_projective_vortex} and
\[
  \Pi_{A,\varphi}^2: \Omega^{0,2}(\fsl(E)) \oplus \Omega^{0,1}(E)
  \to
  \Ker d_{A,\varphi}^2\cap\left(\Omega^{0,2}(\fsl(E)) \oplus \Omega^{0,1}(E)\right)
\]
is $L^2$-orthogonal projection.  The equality between $\bar\rd_{A,\varphi}^2$ (defined in \eqref{eq:dStablePair}) and the restriction of
$d_{A,\varphi}^2$ to $\Omega^{0,2}(\fsl(E)) \times \Omega^{0,1}(E)$ in the definition \eqref{eq:d2_projective_vortex} implies that
\begin{equation}
\label{eq:EqualityBetweenProjVotexAndHolPairProjections}
\Pi_{A,\varphi}^2 = \Pi_{\bar\rd_A,\varphi}^{0,2},
\end{equation}
where $\Pi_{\bar\partial_A,\varphi}^{0,2}$ is the $L^2$-orthogonal projection defined in Lemma \ref{lem:Friedman_Morgan_1999_1-2} onto the closed subspace $\Ker\bar\partial_{A,\varphi}^2$:
\[
  \Pi_{\bar\partial_A,\varphi}^{0,2}: \Omega^{0,2}(\fsl(E)) \oplus \Omega^{0,1}(E)
  \to
  \Ker\bar\partial_{A,\varphi}^2 \subset\Omega^{0,2}(\fsl(E)) \oplus \Omega^{0,1}(E).
\]
As in Friedman and Morgan \cite[Proof of Theorem 3.8, Section 4.3.4, pp. 326--328]{FrM}, to facilitate the comparison of the Kuranishi models $(\beps,\bkappa)$ and $(\bgamma,\bchi)$, we refine our previous definition of the affine slice $S_{A,\varphi}$ in \eqref{eq:Unitary_pair_Avarphi_linear_slice} by replacing it with the real analytic subvariety,
\begin{multline}
  \label{eq:S_Avarphi_prime_nonlinear_slice}
  S_{A,\varphi}' := \{A,\varphi\}
  + \left\{(a,\sigma) \in \Omega^1(\su(E))\oplus\Omega^0(E): d_{A,\varphi}^{0,*}(a,\sigma) = 0 \right.
  \\
  \left. \text{and } \pi_0\sS(A+a,\varphi+\sigma) = 0 \right\},
\end{multline}
defined by the nonlinear map $\sS$ in \eqref{eq:Projective_vortex_and_Coulomb_gauge_map_sheaves} and the $L^2$-orthogonal projection \eqref{eq:pi_0}, namely
\[
  \pi_0:\Omega^0(\su(E))\oplus \Omega^{0,2}(\fsl(E)) \oplus \Omega^{0,1}(E)  \to \Omega^0(\su(E)).
\]
As before, we denote $\pi_0^\perp := \id - \pi_0$, where here ``$\id$'' means the identity operator on the vector space $\Omega^0(\su(E))\oplus \Omega^{0,2}(\fsl(E)) \oplus \Omega^{0,1}(E)$. Thus, we can write $\sS = \sS_0 + \sS_1$ (see also Remark \ref{rmk:Proof_projective_vortices_vanishing_zero-order_cohomology_regular_points_moment_map}), where
\begin{equation}
\label{eq:DefineComponentsOfProjectiveVortexEquation}
\sS_0 := \pi_0\sS\quad\text{and}\quad \sS_1 := \sS - \sS_0.
\end{equation}
Hence, we see from the definitions of $\sS$ in \eqref{eq:DefineComponentsOfProjectiveVortexEquation} and $\sS$ in \eqref{eq:Projective_vortex_and_Coulomb_gauge_map_sheaves} and $\pi_0$ in \eqref{eq:pi_0} and the co-moment map $\bmu^*:\sA(E,h)\times\Omega^0(E) \to \Omega^0(\su(E))$ in  \eqref{eq:Co-moment_map_action_unitary_det_one_gauge_transformations_affine_space_pairs} that
\[
  \pi_0\sS(A+a,\varphi+\sigma)
  =
  (\Lambda F_{A+a})_0 - \frac{i}{2}((\varphi+\sigma)\otimes(\varphi+\sigma)^*)_0
  =
  \bmu^*(A+a,\varphi+\sigma),
\]
where $\bmu^*$ is associated to the moment map $\bmu:\sA(E,h)\times\Omega^0(E) \to \Omega^0(\su(E))^*$ in \eqref{eq:Moment_map_action_unitary_det_one_gauge_transformations_affine_space_pairs}. 

\begin{claim}
\label{claim:S_Avarphi_prime_embedded_real_analytic_submanifold_near_Avarphi}
There is an open neighborhood $U_{A,\varphi}'$ of $(A,\varphi)$ in the real analytic subvariety $S_{A,\varphi}'$ in \eqref{eq:S_Avarphi_prime_nonlinear_slice} that is an embedded, real analytic submanifold of an open neighborhood $U_{A,\varphi}$ of $(A,\varphi)$ in the affine, Coulomb-gauge, local slice $S_{A,\varphi}$ in \eqref{eq:Unitary_pair_Avarphi_linear_slice}. Moreover,
\begin{equation}
  \label{eq:Zero_locus_sS1_on_UAvarphi_prime_equals_zero_locus_sS1_on_UAvarphi}
  \iota_1 \circ \sS_1 = \sS \circ \iota_{A,\varphi} \quad\text{on } U_{A,\varphi}',
\end{equation}
where $\iota_{A,\varphi}:U_{A,\varphi}' \hookrightarrow U_{A,\varphi}$ is the embedding, and the following diagram commutes,
\[
  \begin{CD}
  U_{A,\varphi}' &@>{\sS_1}>> \Omega^{0,2}(\fsl(E))\oplus\Omega^{0,1}(E)
  \\
  @VV{\iota_{A,\varphi}}V &@VV{\iota_1}V
  \\
  U_{A,\varphi} &@>{\sS}>> \Omega^0(\su(E)) \oplus\Omega^{0,2}(\fsl(E))\oplus\Omega^{0,1}(E)
  \end{CD}
\]
where $\iota_1$ is the canonical inclusion, $(w,\nu) \mapsto (0,w,\nu)$, for all $(w,\nu) \in \Omega^{0,2}(\fsl(E))\oplus\Omega^{0,1}(E)$.
\end{claim}

\begin{proof}
We have $(A,\varphi) \in \bmu^{-1}(0)$, because $(A,\varphi)$ is a projective vortex by hypothesis, which is equivalent to $\sS(A,\varphi) = 0$ for $\sS$ as in \eqref{eq:Projective_vortex_map}, and the facts that $\sS_0 = \bmu^*$ and $\bmu(A,\varphi) = 0 \iff \bmu^*(A,\varphi) = 0$. Lemma \ref{lem:Pairs_trivial_stabilizer_group_regular_points_moment_map} (see also Remark \ref{rmk:Proof_projective_vortices_vanishing_zero-order_cohomology_regular_points_moment_map}) implies that $(A,\varphi)$ is a regular point of the moment map $\bmu$ and thus a regular point of co-moment map $\bmu^* = \sS_0:\sA(E,h)\times\Omega^0(E) \to \Omega^0(\su(E))$ if and only if $\bH_{A,\varphi}^0 = (0)$.

Because $d_{A,\varphi}^1 = d\sS(A,\varphi)$ by \eqref{eq:d1_projective_vortex_elliptic_deformation_complex} and the definition of $\sS$ in \eqref{eq:Projective_vortex_map} and as $d_{A,\varphi}^1\circ d_{A,\varphi}^0 = 0$ by \eqref{eq:d_squared_zero_projective_vortex_elliptic_deformation_complex}, we see that
\[
  d\sS(A,\varphi) = 0 \quad\text{on } \Ran d_{A,\varphi}^0.
\]
Therefore, since $\sS_0 = \pi_0\sS$, where $\pi_0$ is the projection \eqref{eq:pi_0}, we obtain
\[
  d\sS_0(A,\varphi) = 0 \quad\text{on } \Ran d_{A,\varphi}^0.
\]
Due to the $L^2$-orthogonal decomposition,
\[
  T_{A,\varphi}\left(\sA(E,h)\times\Omega^0(E)\right)
  =
  \Omega^1(\su(E))\oplus \Omega^0(E)
  =
  \Ran d_{A,\varphi}^0 \oplus \Ker d_{A,\varphi}^{0,*},
\]
the operator $d\sS_0(A,\varphi):\Omega^1(\su(E))\oplus \Omega^0(E) \to \Omega^0(\su(E))$ is surjective if and only if its restriction to the tangent space $T_{A,\varphi}S_{A,\varphi}$ of the affine slice $S_{A,\varphi}$ is surjective:
\[
  d\sS_0(A,\varphi):\Ker d_{A,\varphi}^{0,*}\cap\left(\Omega^1(\su(E))\oplus \Omega^0(E)\right) \to \Omega^0(\su(E)).
\]
Consequently, $(A,\varphi)$ is a regular point of the real analytic map, $\sS_0:S_{A,\varphi} \to \Omega^0(\su(E))$, if and only if $\bH_{A,\varphi}^0 = (0)$, which holds by our hypothesis that $(A,\varphi)$ is strongly simple and the fact that $\bH_{A,\varphi}^0$ is the Lie algebra of $\Stab(A,\varphi)$ by Lemma \ref{lem:LieGroupStructureOfStab(A,varhi)}. Therefore, by the Implicit Mapping Theorem for real analytic maps of Banach spaces (see Section \ref{sec:Real_complex_analytic_functions_Banach_spaces}), there is an open neighborhood of $(A,\varphi)$ in the nonlinear slice $S_{A,\varphi}'$ that is an embedded, real analytic submanifold of an open neighborhood of $(A,\varphi)$ in the affine slice $S_{A,\varphi}$, as asserted.

By definition of $S_{A,\varphi}'$ in \eqref{eq:S_Avarphi_prime_nonlinear_slice}, we have $\sS_0 = 0$ on $S_{A,\varphi}'$ and thus $\iota_1\circ\sS_1 = \sS\circ\iota_{A,\varphi}$ on $S_{A,\varphi}'$ and this yields the identity \eqref{eq:Zero_locus_sS1_on_UAvarphi_prime_equals_zero_locus_sS1_on_UAvarphi} and commutative diagram. This completes the proof of Claim \ref{claim:S_Avarphi_prime_embedded_real_analytic_submanifold_near_Avarphi}.
\end{proof}

As a consequence of Claim \ref{claim:S_Avarphi_prime_embedded_real_analytic_submanifold_near_Avarphi} --- especially the fact that $\sS_0 = 0$ on $S_{A,\varphi}'$ --- the Kuranishi model associated to the pair $(A,\varphi)$ obtained by applying the Implicit Mapping Theorem to the map $\sS$ in \eqref{eq:Projective_vortex_and_Coulomb_gauge_map_sheaves}, namely
\[
  \sS: S_{A,\varphi} \to \Omega^0(\su(E))\oplus \Omega^{0,2}(\fsl(E)) \oplus \Omega^{0,1}(E),
\]  
is identical to the Kuranishi model obtained by applying the Implicit Mapping Theorem to the restriction of $\sS$ in \eqref{eq:Projective_vortex_map} to the open neighborhood $U_{A,\varphi}'$ of $(A,\varphi)$ in the nonlinear local slice $S_{A,\varphi}'$ in \eqref{eq:S_Avarphi_prime_nonlinear_slice} to the map
\[
  \sS_1: S_{A,\varphi}' \to \Omega^{0,2}(\fsl(E)) \oplus \Omega^{0,1}(E).
\]
This means, in particular, that both applications of the Implicit Mapping Theorem yield the same maps \eqref{eq:Kuranishi_embedding_map_projective_vortex} and \eqref{eq:Kuranishi_obstruction_map_projective_vortex}, respectively,
\[
  \beps:\bH_{A,\varphi}^1 \supset N_{A,\varphi} \to \sA(E,h)\times\Omega^0(E)
  \quad\text{and}\quad
  \bkappa:\bH_{A,\varphi}^1 \supset N_{A,\varphi} \to \bH_{A,\varphi}^2,
\]
and thus define the same Kuranishi model, $\sK(A,\varphi)$, as in Definition  \ref{defn:Friedman_Morgan_4-1-15_projective_vortex}.

Because we have assumed, in the hypotheses of Theorem  \ref{thm:Lubke_Teleman_6-3-7}, that $F_{A_d}^{0,2}=0$, where $A_d$ is the fixed unitary connection induced on $\det E$ by all $A'\in\sA(E,h)$, we have by \eqref{eq:(0,2)ComponentOfCurvatureOfUnitaryConnectionTracefree} that
\[
  F_{A'}^{0,2} = (F_{A'}^{0,2})_0, \quad\text{for all } A'\in \sA(E,h).
\]
We observe that the equation for a solution $(A+a,\varphi+\sigma)$ to the $\Omega^{0,2}(\fsl(E)) \oplus \Omega^{0,1}(E)$-component of the projective vortex equations \eqref{eq:Projective_vortex_and_Coulomb_gauge_map_sheaves} near a given solution $(A,\varphi)$ (see the definition \eqref{eq:DefineComponentsOfProjectiveVortexEquation} of $\sS_1$ and Remark \ref{rmk:Proof_projective_vortices_vanishing_zero-order_cohomology_regular_points_moment_map}),
\begin{multline*}
  \sS_1(A+a,\varphi+\sigma)
  = (\sS_{0,2},\sS_{0,1})(A+a,\varphi+\sigma)
  \\
  = (0,0) \in \Omega^{0,2}(\fsl(E)) \oplus \Omega^{0,1}(E),
  \quad\text{for all } (a,\sigma) \in S_{A,\varphi},
\end{multline*}
is equal to the equation for a solution (with $\pi_h^{0,1}$ as in \eqref{eq:Bijection_unitarypairs_with_01pairs}),
\[
  (\bar\partial_A+a'',\varphi+\sigma) = \pi_h^{0,1}(A+a,\varphi+\sigma),
\]
to the holomorphic pair equations \eqref{eq:Holomorphic_pair} near the corresponding holomorphic pair $(\bar\partial_A,\varphi)$, namely
\[
  \fS(\bar\partial_A+a'',\varphi+\sigma) = (0,0) \in \Omega^{0,2}(\fsl(E)) \oplus \Omega^{0,1}(E).
\]
In writing the preceding equation, it is useful to recall the definition \eqref{eq:Holomorphic_pair_map} of the holomorphic pair map $\fS$:
\begin{multline*}
  \fS:\sA^{0,1}(E)\times \Omega^0(E) \ni (\bar\partial_A+\alpha,\varphi+\sigma)
  \\
  \mapsto
  \left(F_{\bar\partial_A+\alpha},(\bar\partial_A+\alpha)(\varphi+\sigma)\right)
  \in
  \Omega^{0,2}(\fsl(E)) \oplus \Omega^{0,1}(E).
\end{multline*}
We now compare the maps $\fS$ and $\sS_1$ restricted from $S_{A,\varphi}$ to the nonlinear local slice $S_{A,\varphi}'$ in \eqref{eq:S_Avarphi_prime_nonlinear_slice}: for all $(a,\sigma) \in S_{A,\varphi}'$, we have
\begin{align*}
  \sS_1(A+a,\varphi+\sigma)
  &= \left( (F_{A+a}^{0,2})_0, \bar\partial_{A+a}(\varphi+\sigma) \right)
  \\
  &= \left( F_{A+a}^{0,2}, \bar\partial_{A+a}(\varphi+\sigma) \right)
  \quad\text{(by \eqref{eq:(0,2)ComponentOfCurvatureOfUnitaryConnectionTracefree})}
  \\
  &= \left( F_{\bar\partial_{A+a}}, \bar\partial_{A+a}(\varphi+\sigma) \right)
    \quad\text{(by \eqref{eq:Holomorphic_curvature})}
  \\
  &= \left( F_{\bar\partial_A+a''}, (\bar\partial_A+a'')(\varphi+\sigma) \right)    
  \\
  &= \fS(\bar\partial_A + a'', \varphi+\sigma)
    \quad\text{(by definition \eqref{eq:Holomorphic_pair_map} of $\fS$)}
  \\
  &= \fS\left(\pi_h^{0,1}(A + a, \varphi+\sigma)\right)
    \quad\text{(by definition \eqref{eq:Bijection_unitarypairs_with_01pairs} of $\pi_h^{0,1}$)}.   
\end{align*}
In summary, we have shown that
\begin{equation}
  \label{eq:sS1_equals_fS1_on_slice_SAvarphi_prime}
  \sS_1 = \fS\circ\pi_h^{0,1} \quad\text{on } U_{A,\varphi}'.
\end{equation}
It remains to take into account the fact that the equation $\fS(\bar\partial_A+\alpha,\varphi+\sigma') = 0$ is solved near a given solution $(\bar\partial_A,\varphi)$ for $(\alpha,\sigma')$ in the affine local slice
$S_{\bar\partial_A,\varphi}$, rather than $(a'',\sigma)$ in the nonlinear local slice $\pi_h^{0,1}S_{A,\varphi}'$ described in the following

\begin{claim}
\label{claim:pih01S_Avarphi_prime_slice_near_dbarAvarphi}  
The subset $\pi_h^{0,1}S_{A,\varphi}'$ is a real analytic subvariety of $\sA^{0,1}(E)\times\Omega^0(E)$ and there is an open neighborhood of $(\bar\partial_A,\varphi)$ in $\pi_h^{0,1}S_{A,\varphi}'$ that is an embedded real analytic submanifold of an open neighborhood of $(\bar\partial_A,\varphi)$ in $\sA^{0,1}(E)\times\Omega^0(E)$ and a local slice for the action of $C^\infty(\SL(E))$ on $\sA^{0,1}(E)\times\Omega^0(E)$.
\end{claim}  

\begin{proof}
Because $\pi_h^{0,1}$ in \eqref{eq:Bijection_unitarypairs_with_01pairs} is an isomorphism of real Banach affine spaces, the conclusions --- aside from the assertion that $\pi_h^{0,1}S_{A,\varphi}'$ is a local slice --- follow from Claim \ref{claim:S_Avarphi_prime_embedded_real_analytic_submanifold_near_Avarphi}. For the slice property, we will show that the tangent space to $\pi_h^{0,1}S_{A,\varphi}'$ at $(\bar\partial_A,\varphi)$ equals $T_{\bar\partial_A,\varphi}S_{\bar\partial_A,\varphi}$, that is,
\begin{equation}
  \label{eq:Tangent_space_pih01_S_Avarphi_prime}
  T_{\bar\partial_A,\varphi}\pi_h^{0,1}S_{A,\varphi}'
  =
  \pi_h^{0,1}T_{A,\varphi}S_{A,\varphi}'
  =
  T_{\bar\partial_A,\varphi}S_{\bar\partial_A,\varphi}.
\end{equation}
The first equality in \eqref{eq:Tangent_space_pih01_S_Avarphi_prime} follows from the facts that $\pi_h^{0,1}$ is linear and $(\bar\partial_A,\varphi) = \pi_h^{0,1}(A,\varphi)$. Observe that
\begin{multline*}
  T_{\bar\partial_A,\varphi}S_{\bar\partial_A,\varphi}
  =
  T_{\bar\partial_A,\varphi}\left((\bar\partial_A,\varphi) + \Ker\bar\partial_{A,\varphi}^{0,*}\cap
    \left(\Omega^{0,1}(\fsl(E))\oplus \Omega^0(E)\right) \right)
  \\
  =
  \Ker\bar\partial_{A,\varphi}^{0,*}\cap \left(\Omega^{0,1}(\fsl(E))\oplus \Omega^0(E)\right),  
\end{multline*}
and so, more explicitly,
\[
  T_{\bar\partial_A,\varphi}S_{\bar\partial_A,\varphi}
  =
  \left\{(\alpha,\sigma) \in \Omega^{0,1}(\fsl(E))\oplus \Omega^0(E):
    \bar\partial_{A,\varphi}^{0,*}(\alpha,\sigma) = 0\right\}.
\]
Similarly,
\begin{multline*}
  T_{A,\varphi}S_{A,\varphi}'
  =
  T_{A,\varphi}\left( \sS_0^{-1}(0) \cap\left((A,\varphi) + \Ker d_{A,\varphi}^{0,*}\cap
      \left(\Omega^1(\su(E))\oplus \Omega^0(E)\right)\right)\right)
  \\
  =
  \Ker d\sS_0(A,\varphi) \cap \Ker d_{A,\varphi}^{0,*}\cap\left(\Omega^1(\su(E))\oplus \Omega^0(E)\right).
\end{multline*}
and so, more explicitly,
\[
  T_{A,\varphi}S_{A,\varphi}'
  =
  \left\{(a,\sigma) \in \Omega^1(\su(E))\oplus \Omega^0(E):
    d_{A,\varphi}^{0,*}(a,\sigma) = 0 \text{ and } d\sS_0(A,\varphi)(a,\sigma) = 0\right\}.
\]
From equation \eqref{eq:dbmu*_Avarphi_asigma} and the fact that $\sS_0 = \bmu^*$, we have
\begin{multline*}
  d\sS_0(A,\varphi)(a,\sigma)
  =
  \frac{1}{2}\Lambda(\bar\partial_Aa' + \partial_Aa'') - \frac{i}{2}(\varphi\otimes\sigma^* + \sigma\otimes\varphi^*)_0,
  \\
  \text{for all } (a,\sigma) \in W^{1,p}(T^*X\otimes\su(E)\oplus E).
\end{multline*}
Therefore, $d\sS_0(A,\varphi)(a,\sigma) = 0$ if and only if
\[
  \Lambda(\bar\partial_Aa' + \partial_Aa'') - i(\varphi\otimes\sigma^* + \sigma\otimes\varphi^*)_0  = 0.
\]
By the K\"ahler Identities \eqref{eq:Kaehler_identity_commutator_Lambda_del-bar_A_and_Lambda_del_A}, namely,
\[
  [\Lambda,\bar\partial_A] = -i\partial_A^* \quad\text{and}\quad [\Lambda,\partial_A] = i\bar\partial_A^*,
\]  
the preceding identity is equivalent to
\[
  \bar\partial_A\Lambda a' -i\partial_A^*a' + \partial_A\Lambda a'' + i\bar\partial_A^*a''
  - i(\varphi\otimes\sigma^* + \sigma\otimes\varphi^*)_0  = 0.
\]
But $a' \in \Omega^{1,0}(\fsl(E))$ and $a'' \in \Omega^{0,1}(\fsl(E))$, thus $\Lambda a' = 0$ and $\Lambda a'' = 0$ and so the preceding identity is equivalent to
\[
  -i\partial_A^*a' + i\bar\partial_A^*a'' - i(\varphi\otimes\sigma^* + \sigma\otimes\varphi^*)_0  = 0,
\]
that is, after dividing by the factor $i$,
\[
  -\partial_A^*a' + \bar\partial_A^*a'' - (\varphi\otimes\sigma^* + \sigma\otimes\varphi^*)_0  = 0.
\] 
Equation \eqref{eq:Real_Coulomb_slice_SO3_monopole_complex_Kaehler_del_A_and_del-bar_A_adjoints_tensor_products} in Lemma \ref{lem:Real_Coulomb_slice_SO3_monopole_complex_Kaehler_del_A_and_del-bar_A_adjoints_tensor_products} provides a version of the Coulomb-gauge slice condition that is equivalent to the usual identity $d_{A,\varphi}^{0,*}(a,\sigma) = 0$:
\[
  \partial_A^*a' + \bar\partial_A^*a'' - \left(\sigma\otimes\varphi^* - \varphi\otimes\sigma^*\right)_0 = 0.
\]
(Although our analysis in Section \ref{sec:Isomorphisms_between_first-order_cohomology_groups} concerned non-Abelian monopoles over almost Hermitian four-manifolds rather than projective vortices over almost Hermitian manifolds of arbitrary dimension, that has no impact on the application here.) By adding the preceding two identities and dividing by $2$, we obtain
\[
  \bar\partial_A^*a'' - (\sigma\otimes\varphi^*)_0 = 0.
\]
The expression \eqref{eq:R_varphi_star_sigma_is_sigma_tensor_varphi_star_tracefree} for $R_\varphi^*$ ensures that the preceding identity is equivalent to
\[
  \bar\partial_A^*a'' - R_\varphi^*\sigma = 0,
\]
where $R_\varphi^*:\Omega^0(E) \to \Omega^0(\fsl(E))$ is defined by the adjoint of $R_\varphi$ in \eqref{eq:Right_multiplication_of_section_slE_by_section_E}. But the derivation of equation \eqref{eq:H_dbar_APhi^01_explicit_0} indicates (after substituting $\tau=0$) that the preceding identity is equivalent to
\[
  \bar\partial_{A,\varphi}^{0,*}(a'',\varphi) = 0.
\]
We conclude that the second equality in \eqref{eq:Tangent_space_pih01_S_Avarphi_prime} holds and this completes the proof of Claim \ref{claim:pih01S_Avarphi_prime_slice_near_dbarAvarphi}.
\end{proof}

Since $\pi_h^{0,1}S_{A,\varphi}'$ and $S_{\bar\partial_A,\varphi}$ are both local slices for the action of
$C^\infty(\SL(E))$, Theorem \ref{thm:Existence_of_complex_gauge_transformation_to_Coulomb_01_pairs} implies that, after possibly shrinking the open neighborhood $U_{A,\varphi}'$ of $(A,\varphi)$ in $S_{A,\varphi}'$ given by Claim \ref{claim:S_Avarphi_prime_embedded_real_analytic_submanifold_near_Avarphi}, there is a real analytic map,
\begin{equation}
  \label{eq:bv_UAvarphi_prime_to_SL(E)}
  \bv: U_{A,\varphi}' \ni (A+a,\varphi+\sigma)
  \mapsto  v_{A+a,\varphi+\sigma} = \bv(A+a,\varphi+\sigma) \in C^\infty(\SL(E)),
\end{equation}
with $\bv(A,\varphi) = \id_E$ and $D\bv(A,\varphi) = 0$ such that the resulting map,
\begin{multline}
  \label{eq:u_pi01_SAvarphi_prime_to_sA01_Omega0E}
 \Theta:U_{A,\varphi}' \ni (A+a,\varphi+\sigma) \mapsto 
  \left(\bv\cdot\pi_h^{0,1}\right)(A+a,\varphi+\sigma)
\\
= v\cdot\left(\bar\partial_{A+a},\varphi+\sigma\right)
  \in U_{\bar\partial_A,\varphi} \subset \sA^{0,1}(E)\times\Omega^0(E),
\end{multline}
gives a real analytic isomorphism from $U_{A,\varphi}'$ onto an open neighborhood $U_{\bar\partial_A,\varphi}$ of $(\bar\partial_A,\varphi)$ in $S_{\bar\partial_A,\varphi}$ with differential (compare \eqref{eq:Tangent_space_pih01_S_Avarphi_prime})
\[
  D\Theta(A,\varphi) = \pi_h^{0,1}:T_{A,\varphi}S_{A,\varphi}' \to T_{\bar\partial_A,\varphi}S_{\bar\partial_A,\varphi}.
\]
In \eqref{eq:u_pi01_SAvarphi_prime_to_sA01_Omega0E}, we abbreviate $v = v_{A+a,\varphi+\sigma}$ and recall that gauge transformations $v \in C^\infty(\SL(E))$ act on pairs $(\bar\partial_{A+a},\varphi+\sigma)$ as in \eqref{eq:SL(E)ActionOn(0,1)Pairs}.

Without loss of generality, we can suppose that the open neighborhood $N_{\bar\partial_A,\varphi} \subset \bH_{\bar\partial_A,\varphi}^1$ of the origin on which the holomorphic embedding and obstruction maps $(\bgamma,\bchi)$ are defined corresponds under $\pi_h^{0,1}$ to the open neighborhood $N_{A,\varphi} \subset \bH_{A,\varphi}^1$ on which the real analytic embedding and obstruction maps $(\beps,\bkappa)$ are defined. Thus, we have two Kuranishi obstruction maps --- a holomorphic map $\bchi:N_{\bar\partial_A,\varphi}\to\bH_{\bar\partial_A,\varphi}^2$ and a real analytic map from $\bkappa:N_{A,\varphi} \to \bH_{A,\varphi}^2$ --- and, by analogy with Friedman and Morgan \cite[Section 4.3.4, p. 328]{FrM}, a diagram (which need not be commutative because of the gauge transformations, $v$) that summarizes the relationships between these real analytic maps:
\[
  \begin{CD}
  \bH_{A,\varphi}^1 \supset N_{A,\varphi} &@>{\beps}>> U_{A,\varphi}' &@>{\sS_1}>> \Omega^{0,2}(\fsl(E))\oplus\Omega^{0,1}(E)
  \\
  @VV{\pi_h^{0,1}}V &@VV{\Theta}V &@|
  \\
  \bH_{\bar\partial_A,\varphi}^1 \supset N_{\bar\partial_A,\varphi} &@>{\bgamma}>> U_{\bar\partial_A,\varphi} &@>{\fS}>> \Omega^{0,2}(\fsl(E))\oplus\Omega^{0,1}(E)
  \end{CD}
\]
Here, $\bgamma$ denotes the embedding of the open neighborhood $N_{\bar\partial_A,\varphi}$ of the origin in the vector space $\bH_{\bar\partial_A,\varphi}^1$ into the affine slice $S_{\bar\partial_A,\varphi}$ provided by the Kuranishi model $\fK(\bar\partial_A,\varphi)$ in Definition \ref{defn:Friedman_Morgan_4-1-15_holomorphic_pairs}. Similarly, $\beps$ denotes the embedding of the open neighborhood $N_{A,\varphi}$ of the origin in the vector space $\bH_{A,\varphi}^1$ into the real analytic subvariety $S_{A,\varphi}'$ of the affine slice $S_{A,\varphi}$ provided by the Kuranishi model $\sK(A,\varphi)$ in Definition \ref{defn:Friedman_Morgan_4-1-15_projective_vortex}. Indeed, we can compare $\sS_1$ and $\fS\circ\Theta$ as follows, noting that $a''=\pi_h^{0,1}a$. For all $(A+a,\varphi+\sigma) \in U_{A,\varphi}'$, we have
\begin{align*}
  (\fS\circ\Theta)(A+a,\varphi+\sigma)
  &= \fS\left(\left(\bv\cdot\pi_h^{0,1}\right)(A+a,\varphi+\sigma)\right)
    \quad\text{(by definition \eqref{eq:u_pi01_SAvarphi_prime_to_sA01_Omega0E} of $\Theta$)}
  \\
  &= \fS(v\cdot\pi_h^{0,1}(A+a,\varphi+\sigma))
    \quad\text{(by definition \eqref{eq:bv_UAvarphi_prime_to_SL(E)} of $\bv$)}
  \\
  &= v\cdot\fS\left(\pi_h^{0,1}(A+a,\varphi+\sigma)\right)
  \\
  &\qquad\text{(by $C^\infty(\SL(E))$-equivariance of $\fS$ in \eqref{eq:Holomorphic_pair_map})}
  \\
  &= v\cdot\sS_1(A+a,\varphi+\sigma)
    \quad\text{(by \eqref{eq:sS1_equals_fS1_on_slice_SAvarphi_prime})}
  \\
  &= (\bv\cdot\sS_1)(A+a,\varphi+\sigma)
  \quad\text{(by definition \eqref{eq:bv_UAvarphi_prime_to_SL(E)} of $\bv$)},  
\end{align*}
where we abbreviate $v = v_{A+a,\varphi+\sigma}$. In the preceding calculation, we write $v\cdot\sS_1$ and $v\cdot\fS$ for the action of $C^\infty(\SL(E))$ on $(w,\nu)\in\Om^{0,2}(\fsl(E))\oplus\Om^{0,1}(E)$ given by $v(w,\nu)=(v^{-1}w v,v^{-1}\nu)$. In summary, we have shown that
\begin{equation}
\label{eq:fS_Composed_with_bh_Equals_sS1_and_GaugeChange}
\fS\circ\Theta
=
\bv\cdot\sS_1 \quad\text{on } U_{A,\varphi}'.
\end{equation}
The equalities \eqref{eq:EqualityBetweenProjVotexAndHolPairProjections} and \eqref{eq:fS_Composed_with_bh_Equals_sS1_and_GaugeChange} imply that the following identity also holds:
\begin{equation}
\label{eq:fS_Composed_with_bh_Equals_sS1_and_GaugeChange_projections}
\Pi_{A,\varphi}^2\sS_1
=
\Pi_{\bar\rd_A,\varphi}^{0,2}\bv^{-1}\cdot(\fS\circ \Theta) \quad\text{on } S_{A,\varphi}'.
\end{equation}
(Recall that the projection $\Pi_{A,\varphi}^2$ is defined for unitary pairs $(A,\varphi)$, as is the map $\sS_1$, while the projection $\Pi_{\bar\rd_A,\varphi}^{0,2}$ is defined for $(0,1)$-pairs $(\bar\partial_A,\varphi)$, as is the map $\fS$; they are compared in the identity \eqref{eq:EqualityBetweenProjVotexAndHolPairProjections}.)
It follows, as in the proof of Proposition \ref{prop:Friedman_Morgan_4-1-16_holomorphic_pairs}, that the Kuranishi model for a projective vortex and the Kuranishi model for a holomorphic pair can both be viewed as the zero locus of a Fredholm section of a real analytic vector bundle over $S_{A,\varphi}'$ with respect to two different trivializations. Theorem \ref{thm:Friedman_Morgan_4-1-14} implies that the real analytic structure defined on this zero locus is independent of the choice of trivialization. Hence, the two Kuranishi models are isomorphic as real analytic model spaces. This completes the proof of Theorem \ref{thm:Friedman_Morgan_4-3-8_projective_vortices}.
\end{proof}

Our proof of Theorem \ref{thm:Friedman_Morgan_4-3-8_projective_vortices} relied on the following analogue for projective vortices of Lemma \ref{lem:Friedman_Morgan_1999_1-2} and Corollary \ref{cor:Friedman_Morgan_1999_1-2_restricted_to_Coulomb_gauge_slice_neighborhood} for the holomorphic pair and holomorphic pair and Coulomb gauge equations, respectively.

\begin{lem}[Local equivalence of projective vortex map with composition of projective vortex map and orthogonal projection]
\label{lem:Friedman_Morgan_1999_1-2_projective_vortices}
Continue the hypotheses and notation of Theorem \ref{thm:Friedman_Morgan_4-3-8_projective_vortices}, but relax the hypothesis that $(A,\varphi)$ has trivial stabilizer in $W^{2,p}(\SU(E))$. Then, after possibly shrinking the open neighborhood $U_{A,\varphi}$ of $(A,\varphi)$ in the affine local slice $S_{A,\varphi}$, one has
\[
  \sS^{-1}(0)\cap U_{A,\varphi} = \hat\sS^{-1}(0)\cap U_{A,\varphi},
\]  
where $\sS$ is as in \eqref{eq:Projective_vortex_and_Coulomb_gauge_map_sheaves} and $\hat\sS$ is as in \eqref{eq:Projective_vortex_and_Coulomb_gauge_map_Banach_spaces_orthogonal_projections}.
\end{lem}

\begin{proof}
Clearly, if a pair $(A+a,\varphi+\sigma)$ is in the zero locus of $\sS$, then it is also contained in the zero locus of $\hat\sS$. The verification that a pair $(A+a,\varphi+\sigma)$ in the zero locus of $\hat\sS$ is also contained in the zero locus of $\sS$, possibly after shrinking $S_{A,\varphi}$ to an open neighborhood of the origin, follows \mutatis that of Lemma \ref{lem:Friedman_Morgan_1999_1-2}.
\end{proof}

\section[Complex K\"ahler virtual moduli space of projective vortices]{Complex K\"ahler virtual moduli space of projective vortices near singular points}
\label{sec:Complex_Kaehler_structure_moduli_space_projective_vortices_near_singular_points}
In this section, we describe how the comparison of Kuranishi models in Theorem \ref{thm:Friedman_Morgan_4-3-8_projective_vortices} naturally leads to the construction in Section \ref{subsec:Complex_Kaehler_local_virtual_moduli_space_projective_vortices} of $\sM_{A,\varphi}^{\vir,\CC}(E,h,\omega)$ in the forthcoming definition \eqref{eq:sM_Avarphi^vir_complex} as a local virtual moduli space of projective vortices as a complex K\"ahler manifold. Theorem \ref{thm:Kobayashi_7-6-36_pairs_local_virtual_moduli_space} provided a local virtual moduli space of projective vortices $\sM_\symp^\vir(E,h,\omega)$ in \eqref{eq:sM_Avarphi^vir_circle_invariant} equipped with an $S^1$-invariant symplectic form, but its proof does not conclude that $\sM_\symp^\vir(E,h,\omega)$ is a complex, K\"ahler manifold. While our description in this section does not conclude that $\sM_{A,\varphi}^{\vir,\CC}(E,h,\omega)$ is an $S^1$-invariant complex K\"ahler manifold in general, it does so when $[A,\varphi]\in\sM^0(E,h,\omega)$ is a fixed point of the standard $S^1$ action, as we explain in Section \ref{subsec:Circle-invariant_complex_Kaehler_local_virtual_moduli_space_projective_vortices_near_fixed_points}.

\subsection{Complex K\"ahler local virtual moduli space of projective vortices}
\label{subsec:Complex_Kaehler_local_virtual_moduli_space_projective_vortices}
Continue the hypotheses of Theorem \ref{thm:Friedman_Morgan_4-3-8_projective_vortices}, so $(A,\varphi)$ is a smooth representative of a point $[A,\varphi] \in \sM^{**}(E,h,\omega)$. Lemma \ref{lem:sM*0(E,h,omega)_subset_sM**(E,h,omega)} implies that $\sC^{**}(E,h) \subset \sC^0(E,h)$ and so $\varphi\not\equiv 0$. Theorem \ref{thm:Friedman_Morgan_4-3-8_projective_vortices} asserts that the canonical isomorphism \eqref{eq:Bijection_unitarypairs_with_01pairs} of real Banach affine spaces,
\[
  \pi_h^{0,1}:\sA(E,h)\times W^{1,p}(E) \to \sA^{0,1}(E)\times W^{1,p}(E),
\]
induces the isomorphism \eqref{eq:Real_analytic_isomorphism_beps(preimage_bkappa_0)_onto_bgamma(preimage_bchi_0)} of real analytic sets,
\[
  \beps\left(\kappa^{-1}(0)\right) \to \bgamma\left(\bchi^{-1}(0)\right),
\]
induced by an isomorphism of real analytic spaces from the Kuranishi model $\sK(A,\varphi)$ in Definition \ref{defn:Friedman_Morgan_4-1-15_projective_vortex} and the Kuranishi model $\fK(\bar\partial_A,\varphi)$ in Definition \ref{defn:Friedman_Morgan_4-1-15_holomorphic_pairs}. The complex analytic model space $\fK(\bar\partial_A,\varphi)$ provides an embedded complex submanifold,
\[
  \bgamma\left(N_{\bar\partial_A,\varphi}\right)
  \hookrightarrow \sA^{0,1}(E)\times W^{1,p}(E),
\]
and consequently the (weak) K\"ahler two-form on the complex Banach affine space $\sA^{0,1}(E)\times W^{1,p}(E)$ given by $\bomega$ in \eqref{eq:Kobayashi_7-6-22_pairs} restricts to a K\"ahler two-form on $\bgamma(N_{\bar\partial_A,\varphi})$ and which we continue to denote by $\bomega$. By construction of the complex analytic space $\fK(\bar\partial_A,\varphi)$, the Zariski tangent space $\bH_{\bar\partial_A,\varphi}^1$ to $\bgamma(\bchi^{-1}(0))$ at the origin is equal to the true tangent space to $\bgamma(N_{\bar\partial_A,\varphi})$ at the pair $(\bar\partial_A,\varphi)$:
\[
  T_{\bar\partial_A,\varphi}\bgamma\left(N_{\bar\partial_A,\varphi}\right) = \bH_{\bar\partial_A,\varphi}^1.
\]
The isomorphism \eqref{eq:Real_analytic_isomorphism_beps(preimage_bkappa_0)_onto_bgamma(preimage_bchi_0)} of real analytic sets is induced by the real analytic diffeomorphism \eqref{eq:u_pi01_SAvarphi_prime_to_sA01_Omega0E},
\[
  \Theta: U_{A,\varphi}' \to U_{\bar\partial_A,\varphi},
\]
from the embedded real analytic submanifold $U_{A,\varphi}'$ provided by Claim \ref{claim:S_Avarphi_prime_embedded_real_analytic_submanifold_near_Avarphi} --- namely, an open neighborhood of $(A,\varphi)$ in the real analytic subvariety $S_{A,\varphi}'$ in \eqref{eq:S_Avarphi_prime_nonlinear_slice} of the real Banach affine space $\sA(E,h)\times W^{1,p}(E)$ --- onto an open neighborhood $U_{\bar\partial_A,\varphi}$ of $(\bar\partial_A,\varphi)$ in the complex Banach affine subspace $S_{\bar\partial_A,\varphi}$ in \eqref{eq:dbar_Evarphi_slice} of the complex Banach affine space $\sA^{0,1}(E)\times W^{1,p}(E)$. We thus obtain an embedded real analytic submanifold,
\[
  \Theta^{-1}\left(\bgamma\left(N_{\bar\partial_A,\varphi}\right)\right)
  \hookrightarrow \sA(E,h)\times W^{1,p}(E).
\]
Moreover, the complex structure $\bJ$ and K\"ahler two-form $\bomega$ on $\bgamma(N_{\bar\partial_A,\varphi})$ pull back under $\Theta$ to a complex structure and K\"ahler two-form on $\Theta^{-1}(\bgamma(N_{\bar\partial_A,\varphi}))$, and which we continue to denote by $\bJ$ and $\bomega$, respectively. The tangent space to $\Theta^{-1}(\bgamma(N_{\bar\partial_A,\varphi}))$ at $(A,\varphi)$ obeys
\[
  T_{A,\varphi}\Theta^{-1}\left(\bgamma\left(N_{\bar\partial_A,\varphi}\right)\right)
  \cong
  \bH_{\bar\partial_A,\varphi}^1,
\]
where the isomorphism is defined by the differential $D\Theta(\bar\partial_A,\varphi)$. Theorem \ref{thm:Kobayashi_7-2-21_pairs} implies that
\[
  \bH_{\bar\partial_A,\varphi}^1 \cong \bH_{A,\varphi}^1.
\]
As an aside, we note that $\bH_{A,\varphi}^1$ is isomorphic to the Zariski tangent space to $\sM(E,h,\omega)$ at the point $[A,\varphi]$ and the true tangent spaces at $[A,\varphi]$ to the local real virtual moduli space of projective vortices $\sM_{A,\varphi}^{\vir,\RR}(E,h,\omega)$ in \eqref{eq:sM_Avarphi^vir_real} and the local symplectic virtual moduli space of projective vortices $\sM_\symp^\vir(E,h,\omega)$ in \eqref{eq:sM_Avarphi^vir_circle_invariant}.
  
By the definition \eqref{eq:S_Avarphi_prime_nonlinear_slice} of the nonlinear slice $S_{A,\varphi}'$ through $(A,\varphi)$ and the definition \eqref{eq:Unitary_pair_Avarphi_linear_slice} of the affine slice $S_{A,\varphi}$ through $(A,\varphi)$, we see that
\[
  S_{A,\varphi}' \subset S_{A,\varphi}.
\]
By hypothesis in Theorem \ref{thm:Friedman_Morgan_4-3-8_projective_vortices} that $[A,\varphi]\in \sM^{**}(E,h,\omega)$ as in \eqref{eq:Moduli_space_projective_vortices_StabAvarphi_idE}, we have $\Stab(A,\varphi) = \{\id_E\}$. The openness of $\sC^{**}(E,h)$ in $\sC(E,h)$ provided by Lemma \ref{lem:Openness_configuration_subspace_unitary_pairs_trivial_stabilizer} \eqref{item:Openness_quotient_subspace_unitary_pairs_minimal_stabilizer} ensures,
after possibly shrinking $N_{\bar\partial_A,\varphi}$, that the image of $\Theta^{-1}(\bgamma(N_{\bar\partial_A,\varphi}))$ under the quotient map $\pi:\sA(E,h)\times W^{1,p}(E) \to \sC(E,h)$ associated to the definition \eqref{eq:ConfigurationSpaceForProjectiveVortices} of the quotient space $\sC(E,h)$ is contained in $\sC^{**}(E,h)$. Hence, Theorem \ref{thm:DK_prop_5-2-9_FU_corollary_page_50_unitary_pairs} implies, after possibly further shrinking $N_{\bar\partial_A,\varphi}$, that the quotient map $\pi$ restricts to an embedding of real analytic manifolds,
\[
  \pi:\Theta^{-1}\left(\bgamma\left(N_{\bar\partial_A,\varphi}\right)\right)
  \hookrightarrow \sC^{**}(E,h).
\]
We define the following finite-dimensional, embedded smooth real submanifold,
\begin{equation}
  \label{eq:sM_Avarphi^vir_complex}
  \sM_{[A,\varphi]}^{\vir,\CC}(E,h,\omega)
  := \pi\left(\Theta^{-1}\left(\bgamma\left(N_{\bar\partial_A,\varphi}\right)\right)\right)
  \subset \sC^{**}(E,h).
\end{equation}
The manifold $\sM_{[A,\varphi]}^{\vir,\CC}(E,h,\omega)$ thus inherits (abusing notation as above) a complex structure $\bJ$ and K\"ahler metric $\bomega$ from the complex, K\"ahler manifold $\bgamma(N_{\bar\partial_A,\varphi})$, through the preceding real analytic embeddings $\pi$ and $\Theta$, and its tangent space at $[A,\varphi]$ obeys
\[
  T_{[A,\varphi]}\sM_{[A,\varphi]}^{\vir,\CC}(E,h,\omega)
  \cong
  \bH_{\bar\partial_A,\varphi}^1,
\]
where the isomorphism is defined by the composition of the differentials $D\pi(A,\varphi)$ and $D\Theta(\bar\partial_A,\varphi)$. We have proved the following refinement of Theorem \ref{thm:Kobayashi_7-6-36_pairs_local_virtual_moduli_space}:

\begin{thm}[Local virtual moduli space of non-zero-section projective vortices as a complex K\"ahler manifold]
\label{thm:Kobayashi_7-6-36_pairs_local_virtual_moduli_space_Kaehler}
Continue the hypotheses and notation of Theorem \ref{thm:Friedman_Morgan_4-3-8_projective_vortices}. Then the local virtual moduli space $\sM_{[A,\varphi]}^{\vir,\CC}(E,h,\omega)$ in \eqref{eq:sM_Avarphi^vir_complex} is an embedded real analytic submanifold of $\sC^{**}(E,h)$ as in \eqref{eq:ConfigurationSpaceForProjectiveVortices**}, is a complex K\"ahler manifold with K\"ahler two-form induced from $\bomega$ in \eqref{eq:Kobayashi_7-6-22_pairs}, has tangent space at $[A,\varphi]$ represented by $\bH_{A,\varphi}^1$ in \eqref{eq:H_Avarphi^1}, the following equality holds,
\[
  \sM(E,h,\omega)\cap\pi(\beps(N_{A,\varphi})) = \sM(E,h,\omega)\cap \sM_{[A,\varphi]}^{\vir,\CC}(E,h,\omega),
\]
and $\sM(E,h,\omega)\cap \pi(\beps(N_{A,\varphi}))$ is an embedded topological submanifold of $\sM_{[A,\varphi]}^{\vir,\CC}(E,h,\omega)$, where $\pi:\sA(E,h)\times W^{1,p}(E) \to \sC(E,h)$ is the quotient map associated to the definition \eqref{eq:ConfigurationSpaceForProjectiveVortices} of the quotient space $\sC(E,h)$.
\end{thm}

\subsection{Circle-invariant complex K\"ahler local virtual moduli space of projective vortices near fixed points of the circle action}
\label{subsec:Circle-invariant_complex_Kaehler_local_virtual_moduli_space_projective_vortices_near_fixed_points}
We now suppose that $[A,\varphi] \in \sM^0(E,h,\omega)$ is a fixed-point of the standard circle action \eqref{eq:S1_Action_On_ConfigurationSpaceForProjectiveVortices} on $\sC(E,h)$. Because $\varphi\not\equiv 0$, Proposition \ref{prop:FixedPointsOfS1ActionOnUnitaryQuotientSpace} implies that the pair $(A,\varphi)$ is split as in Definition \ref{defn:Split_trivial_central-stabilizer_unitary_pair} \eqref{item:Split_unitary_pair}, so $A=A_1\oplus A_2$ with respect to a decomposition $E=E_1\oplus E_2$ as an orthogonal direct sum of proper Hermitian vector subbundles and $\varphi \in \Omega^0(E_1)$. According to the forthcoming identity \eqref{eq:UnitaryS1ActionsRelation_rank-r}, the standard circle action $\rho_Z$ in \eqref{eq:S1_Action_On_AffineSpaceForProjectiveVortices} is gauge-equivalent via the family $\rho_{\SU}$ of gauge transformations in the forthcoming \eqref{eq:DefineUnitaryS1sActionsAtReducible_rank-r} to the circle action $\rho_2$ in the forthcoming  \eqref{eq:DefineUnitaryS12ActionsAtReducible_rank-r} defined by scalar multiplication on the Hermitian vector subbundle $E_2$:
\[
  \rho_2(e^{ir\theta})=\rho_{\SU}(e^{-i\theta})\rho_Z(e^{ir_2\theta}), \quad\text{for all } e^{i\theta} \in S^1,
\]
where $E$ has complex rank $r$ and $E_i$ has complex rank $r_i$ for $i=1,2$. By the forthcoming Lemma \ref{lem:EqualityOfUnitaryS1ActionsOnQuotientSpace} and Remark \ref{rmk:ExtensionOfEqualityOfUnitaryS1ActionsOnQuotientSpace}, we conclude that $(A,\varphi)$ is a fixed point of the following analogue of the circle action $\rho_2^\sA$ in \eqref{eq:UnitaryS12ActionOnNonAbelianPairsAffineSpace} on the affine space of unitary rather than \spinu pairs:
\begin{multline*}
\rho_2^\sA: S^1\times \sA(E,h)\times W^{1,p}(E) \ni \left(e^{i\theta},(A,\varphi)\right)
\\
\mapsto \left(\rho_2(e^{-i\theta})^* A,\rho_2(e^{i\theta})\varphi\right) \in \sA(E,h)\times W^{1,p}(E).
\end{multline*}
Similarly, the standard $\CC^*$ action $\rho_Z^\CC$ in the forthcoming definition \eqref{eq:DefineS1ZActionsAtReducible} is gauge-equivalent via the family $\rho_{\SL}$ of gauge transformations in \eqref{eq:DefineS1sActionsAtReducible} to the $\CC^*$ action $\rho_2^\CC$ in \eqref{eq:DefineS12ActionsAtReducible} defined by scalar multiplication on the vector subbundle $E_2$, that is, the forthcoming identity \eqref{eq:ComparingS1ActionsAtReducible} holds:
\[
  \rho_2^\CC(\lambda^r)=\rho_{\SL}(\lambda)^{-1}\rho_Z^\CC(\lambda^{r_2}), \quad\text{for all } \lambda\in\CC^*.
\]
The canonical isomorphism $\pi_h^{0,1}$ in \eqref{eq:Bijection_unitarypairs_with_01pairs} of Banach affine spaces is $S^1$-equivariant with respect to the $S^1$ action $\rho_2^\sA$ on $\sA(E,h)\times W^{1,p}(E)$ and the circle action implied by the $\CC^*$ action $\rho_2^{\sA,\CC}$ in the forthcoming definition \eqref{eq:S12ActionOnAffine} on the affine space of $(0,1)$-pairs:
\begin{multline*}
\rho_2^{\sA,\CC}:\CC^*\times \sA^{0,1}(E)\times W^{1,p}(E) \ni \left(\lambda,(\bar\partial_E,\varphi)\right)
\\
\mapsto
\left( \rho_2^\CC(\lambda)\circ\bar\partial_E\circ \rho_2^\CC(\lambda)^{-1}, \rho_2^\CC(\lambda)\varphi\right)
\in \sA^{0,1}(E)\times W^{1,p}(E).
\end{multline*}
In particular, the forthcoming Lemma \ref{lem:S1_Equivariance_of_CohomologyIsoms} and Remark \ref{rmk:Generalizing_S1_Equivariance_of_CohomologyIsoms} imply that the isomorphisms $\bH_{A,\varphi}^k \cong \bH_{\bar\partial_A,\varphi}^k$ provided by Theorem \ref{thm:Kobayashi_7-2-21_pairs} are $S^1$-equivariant with respect to the $S^1$ actions implied by $\rho_2^\sA$ and $\rho_2^{\sA,\CC}$ on the domain and codomain, respectively. The proof of Lemma \ref{lem:EquivariantKuranishiLemma} implies that the maps $(\beps,\bkappa)$ defining the Kuranishi model $\sK(A,\varphi)$ and the maps $(\bgamma,\bchi)$ defining the Kuranishi model $\fK(\bar\partial_A,\varphi)$ can be chosen to be $S^1$-equivariant: see the proof of Lemma \ref{lem:S1EquivariantKuranishiModelForReducibleIn_non_AbelianMonopoleModuli} for a very detailed argument in the case of type $1$ non-Abelian monopoles over four-dimensional manifolds. It follows easily that the isomorphism of real analytic spaces from $\sK(A,\varphi)$ onto $\fK(\bar\partial_A,\varphi)$ provided by Theorem \ref{thm:Friedman_Morgan_4-3-8_projective_vortices} is $S^1$-equivariant with respect to the preceding $S^1$ actions. Consequently, the local virtual moduli space $\sM_{[A,\varphi]}^{\vir,\CC}(E,h,\omega)$ in \eqref{eq:sM_Avarphi^vir_complex} is an $S^1$-invariant, complex K\"ahler manifold and thus we have proved the following refinement of Theorem \ref{thm:Kobayashi_7-6-36_pairs_local_virtual_moduli_space_Kaehler}:

\begin{thm}[Local virtual moduli space of non-zero-section projective vortices as an $S^1$-invariant, complex K\"ahler manifold]
\label{thm:Kobayashi_7-6-36_pairs_local_virtual_moduli_space_Kaehler_S1-invariant}
Continue the hypotheses and notation of Theorem \ref{thm:Friedman_Morgan_4-3-8_projective_vortices} and assume further that $[A,\varphi]$ is a fixed-point of the standard circle action \eqref{eq:S1_Action_On_AffineSpaceForProjectiveVortices} on $\sC(E,h)$. Then in addition to the conclusions of Theorem \ref{thm:Kobayashi_7-6-36_pairs_local_virtual_moduli_space_Kaehler}, we obtain that the complex K\"ahler manifold $\sM_{[A,\varphi]}^{\vir,\CC}(E,h,\omega)$ is $S^1$-invariant with respect to the standard circle action \eqref{eq:S1_Action_On_ConfigurationSpaceForProjectiveVortices} on $\sC(E,h)$ and the topological embedding of $\sM(E,h,\omega)\cap \pi(\beps(N_{A,\varphi}))$ into $\sM_{[A,\varphi]}^{\vir,\CC}(E,h,\omega)$ is $S^1$-equivariant with respect to the action \eqref{eq:S1_Action_On_ConfigurationSpaceForProjectiveVortices}.
\end{thm}

\section[Embedding of moduli spaces of projective vortices into holomorphic pairs]{Proof of real analytic embedding of the moduli space of projective vortices into the moduli space of holomorphic pairs}
\label{sec:Proofs_theorem_Lubke_Teleman_6-3-7_and_corollary}
In this section, we prove Theorem \ref{thm:Lubke_Teleman_6-3-7} --- an analogue for projective vortices and holomorphic pairs on Hermitian vector bundles over complex K\"ahler manifolds of \cite[Section 4.3.9, Theorem 3.9, p. 328]{FrM}, due to Friedman and Morgan, for anti-self-dual connections and holomorphic structures on Hermitian vector bundles over complex K\"ahler surfaces --- and Corollary \ref{cor:Lubke_Teleman_6-3-7_and_Kobayashi_7_3_17_pair}.

\begin{proof}[Proof of Theorem \ref{thm:Lubke_Teleman_6-3-7}]
We adapt our proof of Theorem \ref{thm:Kobayashi_7_4_20} for Hermitian--Einstein connections and holomorphic structures on Hermitian vector bundles over complex K\"ahler manifolds, which is based in turn on the proof of \cite[Section 4.3.9, Theorem 3.9, p. 328]{FrM} due to Friedman and Morgan.

Consider Item \eqref{item:Lubke_Teleman_6-3-7_set-theoretic_embedding}. The Hitchin--Kobayashi Correspondence Theorem \ref{thm:HitchinKobayashiCorrespondenceForPairs} for projective vortices immediately yields the following set-theoretic bijections,
\begin{subequations}
  \label{eq:Set_theoretic_bijections_projective_vortices_and_polystable_pairs}
  \begin{align}
    \label{eq:Set_theoretic_bijection_non-zero-section_projective_vortices_and_polystable_pairs}
    \sM^0(E,h,\omega) &\leftrightarrow \fM_\ps^0(E,\omega),
    \\                    
    \label{eq:Set_theoretic_bijection_projective_vortices_and_polystable_pairs}
    \sM(E,h,\omega) &\leftrightarrow \fM_\ps(E,\omega),                      
    \\
    \label{eq:Set_theoretic_bijection_non-split_non-zero-section_projective_vortices_and_stable_pairs}
    \sM^{*,0}(E,h,\omega) &\leftrightarrow \fM^0(E,\omega),
    \\
    \label{eq:Set_theoretic_bijection_non-split_projective_vortices_and_stable_pairs}
    \sM^*(E,h,\omega) &\leftrightarrow \fM(E,\omega).
  \end{align}
\end{subequations}
Thus, Item \eqref{item:Lubke_Teleman_6-3-7_set-theoretic_embedding} follows from \eqref{eq:Set_theoretic_bijection_projective_vortices_and_polystable_pairs} and \eqref{eq:Set_theoretic_bijection_non-split_projective_vortices_and_stable_pairs}.

Consider Item \eqref{item:Lubke_Teleman_6-3-7_analytic_embedding_non-split-nonzero-section}, which asserts that 
\begin{inparaenum}[\itshape i\upshape)]
\item the set-theoretic bijection \eqref{eq:Set_theoretic_bijection_non-split_non-zero-section_projective_vortices_and_stable_pairs} is an isomorphism of real analytic spaces, 
\item $\fM^0(E,\omega)$ is an open subspace of the moduli space $\fM(E)$, and
\item $\fM^0(E,\omega) \subset \fM^{**}(E)$.  
\end{inparaenum}  
Let $[\bar\partial_E,\varphi] \in \fM^0(E,\omega)$ and let $[A,\varphi] \in \sM^{*,0}(E,h,\omega)$ be the corresponding point given by the bijection \eqref{eq:Set_theoretic_bijection_non-split_non-zero-section_projective_vortices_and_stable_pairs}.

By Theorem \ref{thm:Local_Kuranishi_model_for_moduli_space_projective_vortices_StabAvarphi_idE} for a projective vortex with trivial stabilizer, an open neighborhood of $[A,\varphi]$ in $\sM^{*,0}(E,h,\omega)$ is parametrized by a Kuranishi model $\sK(A,\varphi)$ in Definition \ref{defn:Friedman_Morgan_4-1-15_projective_vortex}, comprising the real analytic set $\beps(\bkappa^{-1}(0))$, where $\beps$ is a real analytic embedding from an open neighborhood $N_{A,\varphi} \subset \bH_{A,\varphi}^1$ of the origin into an open neighborhood of $(A,\varphi)$ in the slice \eqref{eq:Unitary_pair_Avarphi_linear_slice}, namely
\[
  S_{A,\varphi} = (A,\varphi)+\Ker d_{A,\varphi}^{0,*}\cap W^{1,p}(T^*X\otimes\su(E)\oplus E),
\]
and $\bkappa:\bH_{A,\varphi}^1 \supset N_{A,\varphi} \to \bH_{A,\varphi}^2$ is a real analytic map. The map $\pi\circ\beps:\sK(A,\varphi) \to \sM(E,h,\omega)$ induced by composition of $\beps:\bkappa^{-1}(0) \to \sA(E,h)\times W^{1,p}(E)$ with the quotient map
\[
  \pi:\sA(E,h)\times W^{1,p}(E) \to \sC(E,h)
\]
is open (where $\sC(E,h)$ is endowed with the quotient topology) and so its image is an open neighborhood of $[A,\varphi]$ in $\sM(E,h,\omega)$ (and thus its intersection with $\sM^{*,0}(E,h,\omega)$ is also open). Since $\pi\circ\beps$ is also continuous and injective, it is a homeomorphism.

By Theorem \ref{thm:Local_Kuranishi_model_for_strongly_simple_point_cP(E)} for a strongly simple holomorphic pair, to each representative $(\bar\partial_E,\varphi)$ for a point in $\fM^{**}(E)$ we may associate a Kuranishi model $\fK(\bar\partial_E,\varphi)$ comprising the complex analytic set $\bgamma(\bchi^{-1}(0))$, where $\bgamma$ is a holomorphic embedding from an open neighborhood $N_{\bar\partial_E,\varphi} \subset \bH_{\bar\partial_E,\varphi}^1$ of the origin into an open neighborhood of $(\bar\partial_E,\varphi)$ in the slice \eqref{eq:dbar_Evarphi_slice}, namely
\[
  S_{\bar\partial_E,\varphi}
  =
  (\bar\partial_E,\varphi)
  + \Ker \bar\partial_{E,\varphi}^{0,*}\cap W^{1,p}\left(\Lambda^{0,1}(\fsl(E))\oplus E\right),
\]
and $\bchi:\bH_{\bar\partial_E,\varphi}^1 \supset N_{\bar\partial_E,\varphi} \to \bH_{\bar\partial_E,\varphi}^2$ is a holomorphic map.

The map $\pi\circ\bgamma:\fK(\bar\partial_E,\varphi) \to \fM(E)$ induced by composition of $\bgamma:\bchi^{-1}(0) \to \sA^{0,1}(E)\times W^{1,p}(E)$ with the quotient map
\[
  \pi:\sA^{0,1}(E)\times W^{1,p}(E) \to \sC^{0,1}(E)
\]
is open (where $\sC(E,h)$ is endowed with the quotient topology) and so its image is an open neighborhood of $[\bar\partial_E,\varphi]$ in $\fM(E)$; however $\pi\circ\bgamma$ need not be injective, although if it is injective, it will necessarily be a homeomorphism since it is also clearly continuous.

By Theorem \ref{thm:Friedman_Morgan_4-3-8_projective_vortices}, the canonical isomorphism of real Banach affine spaces,
\[
  \pi_h^{0,1}:\sA(E,h)\times W^{1,p}(E) \ni (A',\varphi')
  \mapsto (\bar\partial_{A'},\varphi) \in \sA^{0,1}(E)\times W^{1,p}(E),
\]
induces an isomorphism of real analytic spaces from $\sK(A,\varphi)$ onto $\fK(\bar\partial_E,\varphi)$ and the bijection \eqref{eq:Set_theoretic_bijection_non-split_non-zero-section_projective_vortices_and_stable_pairs} of quotient spaces. Therefore, the composite map,
\[
  \sM^{*,0}(E,h,\omega) \supset \left.\beps(\bkappa^{-1}(0))\right/W^{2,p}(\SU(E))
\to
\left.\bgamma(\bchi^{-1}(0))\right/W^{2,p}(\SL(E)) \subset \fM(E),
\]
is injective and thus a homeomorphism onto its image. Hence, the image $\fM^0(E,\omega)$ of $M^{*,0}(E,h,\omega)$ is an open subset of $\fM(E)$.

Lastly, the inclusion $\fM^0(E,\omega) \subset \fM^{**}(E)$ is provided by Lemma \ref{lem:Openness_subspace_strongly_simple_01_pairs} \eqref{item:Openness_subspace_strongly_simple_stable_pairs}.
This completes the proof of Item \eqref{item:Lubke_Teleman_6-3-7_analytic_embedding_non-split-nonzero-section}.

Consider Item \eqref{item:Lubke_Teleman_6-3-7_analytic_embedding_rank-2_non-zero_section}. The argument here is virtually identical to the proof of Item \eqref{item:Lubke_Teleman_6-3-7_analytic_embedding_non-split-nonzero-section}, except that we replace the role of the bijection \eqref{eq:Set_theoretic_bijection_non-split_non-zero-section_projective_vortices_and_stable_pairs} by that of \eqref{eq:Set_theoretic_bijection_non-zero-section_projective_vortices_and_polystable_pairs} and replace the role of Lemma \ref{lem:Openness_subspace_strongly_simple_01_pairs} \eqref{item:Openness_subspace_strongly_simple_stable_pairs} by that of Lemma \ref{lem:Openness_subspace_strongly_simple_01_pairs} \eqref{item:Openness_subspace_strongly_simple_polystable_rank-2_pairs}. This completes the proof of the theorem.
\end{proof}

We conclude with the

\begin{proof}[Proof of Corollary \ref{cor:Lubke_Teleman_6-3-7_and_Kobayashi_7_3_17_pair}]
We first consider Item \eqref{item:fM^0(E,omega)_is_complex_analytic_space} and then indicate the changes needed to prove Item \eqref{item:fM_ps^0(E,omega)_rank-2_is_complex_analytic_space}.

By the same argument as in the proof of Item  \eqref{item:Kobayashi_7_4_20_complex_analytic_space} in Theorem \ref{thm:Kobayashi_7_4_20}, the collection of Kuranishi models $\fK(\bar\partial_E,\varphi)$ in Definition \ref{defn:Friedman_Morgan_4-1-15_holomorphic_pairs} define local complex analytic charts on $\fM^0(E,\omega) \subset \fM^{**}(E)$ and they induce transition maps on overlaps that are isomorphisms of complex analytic spaces. Therefore, $\fM^0(E,\omega)$ is a complex analytic space in the sense of Grauert and Remmert \cite[Section 1.1.5, p. 7]{Grauert_Remmert_coherent_analytic_sheaves}. 
  
The open subspace of smooth points, namely $\fM_\reg^0(E,\omega)$ in \eqref{eq:Moduli_space_regular_non-zero-section_stable_holomorphic_pairs}, is a complex manifold as a consequence of the construction of the Kuranishi models $\fK(\bar\partial_E,\varphi)$ in Definition \ref{defn:Friedman_Morgan_4-1-15_holomorphic_pairs}.

According to Item \eqref{item:Lubke_Teleman_6-3-7_analytic_embedding_non-split-nonzero-section} in Theorem \ref{thm:Lubke_Teleman_6-3-7}, there is an isomorphism in the sense of real analytic spaces from the moduli subspace $\sM^{*,0}(E,h,\omega)$ of non-split, non-zero-section projective vortices in
\eqref{eq:Moduli_space_projective_vortices_non-split-non-zero-section} onto the moduli subspace $\fM^0(E,\omega)$ of non-zero-section, stable holomorphic pairs in \eqref{eq:Moduli_space_non-zero-section_stable_holomorphic_pairs}:
\[
  \sM^{*,0}(E,h,\omega) \cong \fM^0(E,\omega).
\]
We claim that the preceding isomorphism restricts to an isomorphism \eqref{eq:Lubke_Teleman_6-3-7_analytic_isomorphism_non-split-nonzero-section_regular} in the sense of real analytic spaces of the open subspaces of regular points, namely:
\[
  \sM_\reg^{*,0}(E,h,\omega) \cong \fM_\reg^0(E,\omega).
\]
Let $[A,\varphi] \in \sM_\reg^{*,0}(E,h,\omega)$. From the definition of $\sM_\reg^{*,0}(E,h,\omega)$ in \eqref{eq:Moduli_space_projective_vortices_non-split-non-zero-section_regular}, we have $\bH_{A,\varphi}^2 = (0)$. Since $(A,\varphi)$ is a projective vortex and $\bH_{A,\varphi}^2 = (0)$, Theorem \ref{thm:Kobayashi_7-2-21_pairs} implies that $\bH_{\bar\partial_A,\varphi}^2 = (0)$ and thus $[\bar\partial_A,\varphi] \in \fM_\reg^0(E,\omega)$ by definition \eqref{eq:Moduli_space_regular_non-zero-section_stable_holomorphic_pairs}. Conversely, suppose $[\bar\partial_A,\varphi] \in \fM_\reg^0(E,\omega)$ and let $[A,\varphi]$ be the corresponding point in $\sM^{*,0}(E,h,\omega)$. We have $\Stab(A,\varphi) = \{\id_E\}$ via the inclusion $\sM^{*,0}(E,h,\omega) \subset \sM^{**}(E,h,\omega)$ provided by \eqref{eq:sM*0(E,h,omega)_subset_sM**(E,h,omega)} in Lemma \ref{lem:sM*0(E,h,omega)_subset_sM**(E,h,omega)} and by the definition \eqref{eq:Moduli_space_projective_vortices_StabAvarphi_idE} of $\sM^{**}(E,h,\omega)$. Thus, $\bH_{A,\varphi}^0 = (0)$, since $\bH_{A,\varphi}^0$ is the Lie algebra of $\Stab(A,\varphi)$
by Lemma \ref{lem:LieGroupStructureOfStab(A,varhi)}. Because $\bH_{\bar\partial_A,\varphi}^2 = (0)$, Theorem \ref{thm:Kobayashi_7-2-21_pairs} ensures that $\bH_{A,\varphi}^2 = (0)$. Therefore, $[A,\varphi] \in \sM_\reg^{*,0}(E,h,\omega)$ by definition \eqref{eq:Moduli_space_projective_vortices_non-split-non-zero-section_regular} of $\sM_\reg^{*,0}(E,h,\omega)$ and this verifies the isomorphism \eqref{eq:Lubke_Teleman_6-3-7_analytic_isomorphism_non-split-nonzero-section_regular}.

By adapting the argument that $\fM^0(E,\omega) \subset \fM^{**}(E)$ is a complex analytic space, we see that the moduli subspace $\sM^{*,0}(E,h,\omega) \subset \sM^{**}(E,h,\omega)$ is a real analytic space since the collection of Kuranishi models $\sK(A,\varphi)$ in Definition \ref{defn:Friedman_Morgan_4-1-15_projective_vortex} define local real analytic charts on $\sM^{*,0}(E,h,\omega)$ and they induce transition maps on overlaps that are isomorphisms of real analytic spaces. Therefore, $\sM^{*,0}(E,h,\omega)$ is a real analytic space.

The open subspace of smooth points, namely $\sM_\reg^{*,0}(E,h,\omega)$ in \eqref{eq:Moduli_space_projective_vortices_non-zero_section_regular}, is a complex manifold as a consequence of the construction of the Kuranishi models $\sK(A,\varphi)$ in Definition \ref{defn:Friedman_Morgan_4-1-15_projective_vortex}. This proves Item \eqref{item:fM^0(E,omega)_is_complex_analytic_space}.

Consider Item \eqref{item:fM_ps^0(E,omega)_rank-2_is_complex_analytic_space}, so $E$ has complex rank two. The proof of Item \eqref{item:fM^0(E,omega)_is_complex_analytic_space} now shows that $\fM_\ps^0(E,\omega) = \fM_\ps^{**}(E,\omega)$ is a complex analytic space and that $\fM_{\ps,\reg}^0(E,\omega) = \fM_{\ps,\reg}^{**}(E,\omega)$ is a complex manifold. To verify the isomorphism \eqref{eq:Lubke_Teleman_6-3-7_analytic_isomorphism_rank-2_non-zero_section_regular}, we replace the
role in the proof of the isomorphism \eqref{eq:Lubke_Teleman_6-3-7_analytic_isomorphism_non-split-nonzero-section_regular} played by the inclusion $\sM^{*,0}(E,h,\omega) \subset \sM^{**}(E,h,\omega)$ in \eqref{eq:sM*0(E,h,omega)_subset_sM**(E,h,omega)} by the equality $\sM^0(E,h,\omega) = \sM^{**}(E,h,\omega)$ from \eqref{eq:sC0(E,h)_equals_sC**(E,h)_E_rank2} in Lemma \ref{lem:sM*0(E,h,omega)_subset_sM**(E,h,omega)}. The proofs of the remaining parts of Item \eqref{item:fM_ps^0(E,omega)_rank-2_is_complex_analytic_space} are the same as for the corresponding parts of Item \eqref{item:fM^0(E,omega)_is_complex_analytic_space}. This completes the proof of the corollary.
\end{proof}

\chapter{Calculation of virtual Morse--Bott index via Hirzebruch--Riemann--Roch theorem}
\label{chap:VirtualMorseIndexComputation}
In this chapter, we prove Theorem \ref{mainthm:MorseIndexAtReduciblesOnKahler} and Corollary \ref{maincor:MorseIndexAtReduciblesOnKahlerWithSO3MonopoleCharacteristicClasses}
by computing the virtual Morse--Bott signature \eqref{eq:Virtual_Morse-Bott_signature_moduli_space_non-abelian_monopoles} of the Hamiltonian function \eqref{eq:Hitchin_function} on the moduli space $\sM_\ft$ of non-Abelian monopoles over a complex K\"ahler surface at a point $[A,\Phi]\in\sM_\ft$ represented by a pair $(A,\Phi)$ which is split in the sense of Definition \ref{defn:Split_trivial_central-stabilizer_spinor_pair}. We will compute the virtual Morse--Bott index \eqref{eq:Virtual_Morse-Bott_index_moduli_space_non-abelian_monopoles}, namely
\[
\dim_\RR\bH_{A,\Phi}^{-,1}-\dim_\RR\bH_{A,\Phi}^{-,2},
\]
where, for $k=0,1,2$, the real vector space $\bH_{A,\Phi}^{-,k}$ is the subspace of the harmonic subspace $\bH_{A,\Phi}^k$ of the elliptic deformation complex \eqref{eq:SO3MonopoleDefComplex} for non-Abelian monopoles on which the $S^1$ action has negative weight. We will construct an $S^1$-equivariant Kuranishi model for the neighborhood of a split pair $[A,\Phi]$ in $\sM_\ft$ and describe  $S^1$ actions on $\bH_{A,\Phi}^k$ in Section \ref{sec:UnitaryS1EquivariantKuranishiModelForReducibleNonAbelianMonopole}.
However, as described in Item \eqref{item:Isomorphic_S1_Representations_Real} of Lemma \ref{lem:Isomorphic_S1_Representations} and in Proposition \ref{prop:Direct_Sum_Decomposition_of_Real_S1_Representations}, the sign of the weight of a real representation is not well-defined.  Hence, the negative-weight subspaces $\bH_{A,\Phi}^{-,k}$ are not defined by the $S^1$ action alone. By Lemma \ref{lem:S1_Invariant_J_Is_MultiplicationBy_i}, we see that such weights are well-defined if there is an $S^1$-invariant almost complex structure, in the sense of \eqref{eq:Circle_invariant_(1,1)-tensor}, on the representation space.
In this monograph, the almost complex structure on $\bH_{A,\Phi}^k$ is defined, when $X$ is a K\"ahler surface, by the isomorphism between the elliptic complexes of non-Abelian monopoles and projective vortices given in Sections \ref{sec:Isomorphisms_between_first-order_cohomology_groups}, \ref{sec:Isomorphisms_between_zeroth-order_cohomology_groups}, and \ref{sec:Isomorphisms_between_second-order_cohomology_groups}.
By Corollary \ref{cor:Comparison_elliptic_complexes_cohomology_groups_pre-holomorphic_and_holomorphic_pairs}, the elliptic complex for projective vortices is isomorphic to the elliptic complex for holomorphic pairs, where the almost complex structure is obvious. We will thus work with the elliptic complex \eqref{eq:Holomorphic_pair_elliptic_complex} for holomorphic pairs. We compute the $S^1$ action induced on this elliptic complex by a gauge-transformed modification of the action on $(0,1)$-pairs given by scalar multiplication on the section in Sections \ref{sec:CircleActionsOnAffineSpace} and \ref{subsec:S1EquivComplex}. In Section \ref{sec:EquivarianceOfIsomorphisms}, we verify the $S^1$-equivariance of the isomorphisms proven in Chapter \ref{chap:Elliptic_deformation_complex_moduli_space_SO(3)_monopoles_over_almost_Hermitian_four-manifold}
and Corollaries
\ref{cor:Comparison_elliptic_complexes_cohomology_groups_pre-holomorphic_and_holomorphic_pairs}
and
\ref{cor:Vanishing_third-order_cohomology_group_holomorphic_pair_elliptic_complex_Kaehler_surface}
between the harmonic spaces $\bH_{A,\Phi}^k$ for the elliptic deformation complex for non-Abelian monopoles and the harmonic spaces $\bH_{\bar\rd_A,\varphi}^k$ for the elliptic complex for holomorphic pairs, for $k=0,1,2$. We give a decomposition of the elliptic complex \eqref{eq:Holomorphic_pair_elliptic_complex} into a direct sum of elliptic subcomplexes on which the modified $S^1$ action defined in Section \ref{sec:CircleActionsOnAffineSpace} acts with zero, positive, and negative weight. We identify the negative-weight subcomplex of  the elliptic complex for holomorphic pairs in Section \ref{sec:NormalTangentialSplittingOfDefComplex}.  If we write $\bH_{\bar\rd_A,\varphi}^{-,k}$ for the harmonic spaces of this negative-weight subcomplex, then the virtual Morse--Bott index will be equal to
\begin{equation}
\label{eq:Holomorphic_VirtualIndex_IntroductionExpression}
\dim_\RR \bH_{\bar\rd_A,\varphi}^{-,1}
-\dim_\RR \bH_{\bar\rd_A,\varphi}^{-,2}.
\end{equation}
We apply the Hirzebruch--Riemann--Roch Theorem to compute the Euler characteristic of the negative-weight subcomplex in Section \ref{sec:IndOfSubcomplex}. The expression for the virtual Morse--Bott index in
Theorem \ref{mainthm:MorseIndexAtReduciblesOnKahler} follows from this computation, the observation that $\bH_{\bar\rd_A,\varphi}^{-,k}=0$ for $k\neq 1,2$, and a discussion of how the $S^1$-equivariant Kuranishi model for a neighborhood of $[A,\Phi]$ in $\sM_\ft$ gives a local real analytic model space (by analogy with the definition in Grauert and Remmert \cite[Section 1.1.2, pp. 3--4]{Grauert_Remmert_analytic_local_algebras} for a complex model space). Similar computations involving the positive-weight and zero-weight subcomplexes of the elliptic complex for holomorphic pairs give the expressions for the co-index and nullity of the Hamiltonian function, completing the proof of Theorem \ref{mainthm:MorseIndexAtReduciblesOnKahler}.
Corollary \ref{maincor:MorseIndexAtReduciblesOnKahlerWithSO3MonopoleCharacteristicClasses} follows from Theorem \ref{mainthm:MorseIndexAtReduciblesOnKahler} and a computation of the characteristic classes of the \spinu and \spinc structures in the statement of the corollary. In Section \ref{sec:CheckingIndexComputation}, we verify that the sum of the indices of the zero-weight, positive-weight, and negative-weight subcomplexes is equal to the index of the total complex. In Section \ref{sec:RestrictionOfHolomorphicPairsMapToWeightSubspaces}, we compute the restriction of the holomorphic-pair map to the weight spaces defined in Section \ref{sec:NormalTangentialSplittingOfDefComplex}.

\begin{rmk}[Hypotheses in Chapter \ref{chap:VirtualMorseIndexComputation} on the dimension of $X$, the rank of $E$, the choice of spin${}^u$ rather than unitary pairs, and the choice of non-Abelian monopoles rather than projective vortices]
\label{rmk:DimensionRankAssumptions}
Some of the results in this chapter, in particular those of Sections \ref{sec:UnitaryS1EquivariantKuranishiModelForReducibleNonAbelianMonopole} and \ref{sec:EquivarianceOfIsomorphisms}, are developed in terms of spin${}^u$ pairs and type $1$ non-Abelian monopoles on rank-two \spinu structures over manifolds of real dimension four for simplicity of exposition and for future applications. However, many of the results for such \spinu pairs or type $1$ non-Abelian monopoles hold more generally for unitary pairs or projective vortices, respectively, on Hermitian vector bundles $E$ of higher rank over manifolds of arbitrary real even dimension. We give more precise statements about such generalizations for the results of Section \ref{sec:UnitaryS1EquivariantKuranishiModelForReducibleNonAbelianMonopole} in Remarks \ref{rmk:ExtensionOfEqualityOfUnitaryS1ActionsOnQuotientSpace},
\ref{rmk:GeneralizationOf_S1_Equivariance_Of_nonAbelianMonopole_Def_Complex}, and
\ref{rmk:GeneralizingKuranishiModels}.
 We discuss the generalization of Lemma \ref{lem:S1_Equivariance_of_CohomologyIsoms} in Remark \ref{rmk:Generalizing_S1_Equivariance_of_CohomologyIsoms}.

The calculations in Sections \ref{sec:CircleActionsOnAffineSpace} and \ref{subsec:S1EquivComplex} are carried out for complex vector bundles $E$ of arbitrary rank greater than or equal to two and complex base manifolds $X$ of any dimension $n$.

The calculations in Section \ref{sec:NormalTangentialSplittingOfDefComplex}
also allow complex base manifolds $X$ of any dimension $n$, but assume that $E$ has complex rank two. Generalizations to the case of vector bundles of rank greater than two are possible, but the statements become more complicated.

In Section \ref{sec:IndOfSubcomplex}, we assume that $(\bar\partial_A,\varphi)$ is a holomorphic pair on a vector bundle over a complex K\"ahler \emph{surface} because the Todd class appearing in Theorem \ref{thm:HRR} has a simple expression in this case (see Remark \ref{rmk:DifficultyInRRComputationInHigherDimension}). However, the computations of the Euler characteristics in Lemmas \ref{lem:IndexOfRedStablePairsDefComplexNormalMinus} and \ref{lem:IndexOfRedStablePairsDefComplexNormalPlus} using the Hirzebruch--Riemann-Roch index formula in Theorem \ref{thm:HRR} would extend to the case of holomorphic pairs over higher-dimensional complex K\"ahler manifolds.  Extensions of the results of Section \ref{sec:IndOfSubcomplex} to higher-rank vector bundles would encounter complications similar to those that would occur in Section \ref{sec:NormalTangentialSplittingOfDefComplex} for higher-rank vector bundles, as well as the difficulties of computing Chern characters of higher-rank vector bundles.
\end{rmk}

\section[Circle action on the Kuranishi model]{Circle action on the Kuranishi model for an open neighborhood of a point represented by a split non-Abelian monopole}
\label{sec:UnitaryS1EquivariantKuranishiModelForReducibleNonAbelianMonopole}
Let $(A,\Phi)$ be a split non-Abelian monopole on the \spinu structure $\ft=(\rho,W\otimes E)$.
In this section, we construct an $S^1$-equivariant Kuranishi model for an open neighborhood of $[A,\Phi]$ in $\sM_\ft$. A Kuranishi model for points in the moduli space of projective vortices appears in Theorem  \ref{thm:Friedman_Morgan_4-3-8_projective_vortices} and one for points in the moduli space of holomorphic pairs in  Definition \ref{defn:Friedman_Morgan_4-1-15_holomorphic_pairs}.

The equivariant Kuranishi model in Lemma \ref{lem:EquivariantKuranishiLemma} gives an $S^1$-equivariant, real analytic isomorphism between an open neighborhood of a point in the zero locus of a Fredholm map between Banach spaces and an open neighborhood of the zero locus of a map between finite-dimensional vector spaces. To apply this result to the point $[A,\Phi]$ in the real analytic manifold $\sC_\ft^0$ rather than to a subspace of a Banach space, we first describe an open neighborhood of $[A,\Phi]$ in $\sM_\ft$ in a chart for an open neighborhood of $[A,\Phi]$ in $\sC_\ft^0$. We obtain such a chart from the inverse of the composition of the projection map,
\begin{equation}
\label{eq:QuotientMapOfPairs}
\pi:\sA(E,h)\times\Omega^0(W^+\otimes E) \to \sC_\ft
\end{equation}
with the embedding
\begin{equation}
\label{eq:LocalSliceBall}
\bB_{A,\Phi}\ni (a,\phi)\mapsto (A+a,\Phi+\phi)\in
\sA(E,h)\times\Om^0(W^+\otimes E),
\end{equation}
where $\bB_{A,\Phi}\subset\Ker d_{A,\Phi}^{0,*}$ is the ball defined in \eqref{eq:SliceBall}. (As in Section \ref{subsec:Kuranishi_model_moduli_space_holomorphic_pairs}, we continue our convention of suppressing explicit notation for Sobolev $L^p$ completions unless required for clarity, with $p\in (n,\infty)$ when $X$ has complex dimension $n$, in order to avoid notational clutter.) An open neighborhood of $[A,\Phi]$ in $\sM_\ft$ is then real analytically isomorphic to an open neighborhood of the origin in the zero locus of the Fredholm map $\fS_{A,\Phi}$ given by the composition of the embedding \eqref{eq:LocalSliceBall} with the map $\fS$ defined in \eqref{eq:PerturbedSO3MonopoleEquation_map} by the non-Abelian monopole equations.

To apply the equivariant Kuranishi model in Lemma \ref{lem:EquivariantKuranishiLemma} to the map $\fS_{A,\Phi}$, we must determine the $S^1$ action on $\bB_{A,\Phi}$ which makes the chart described above equivariant. While the projection map \eqref{eq:QuotientMapOfPairs} is $S^1$-equivariant with respect to the standard action \eqref{eq:S1ZAction} on its domain and the standard action \eqref{eq:S1ZActionOnQuotientSpace} on $\sC_\ft$, the embedding \eqref{eq:LocalSliceBall} is not.  Indeed, unless $\Phi\equiv 0$, the point $(A,\Phi)$ is not a fixed point of the $S^1$ action \eqref{eq:S1ZAction}, so the image of the embedding \eqref{eq:LocalSliceBall}  is not closed under this $S^1$ action. Thus we cannot apply Lemma \ref{lem:EquivariantKuranishiLemma} to the map $\fS_{A,\Phi}$ with the standard $S^1$ action on the domain of $\fS_{A,\Phi}$.

We remedy this problem in Section \ref{subsec:UnitaryS1ActionOnAffineSpace} by introducing an apparently new $S^1$ action (see \eqref{eq:UnitaryS12ActionOnNonAbelianPairsAffineSpace}) on the affine space of spin${}^u$s. This $S^1$ action differs from the standard $S^1$ action by gauge transformations and thus induces the same $S^1$ action, up to a positive multiplicity, as the standard $S^1$ action \eqref{eq:S1ZActionOnQuotientSpace}
on the quotient space. In addition, $(A,\Phi)$ is a fixed point of this modified $S^1$ action. (Compare the work of Hitchin \cite[Section 7, p. 95]{Hitchin_1987} and of Gothen \cite[Section 2.3.2]{GothenThesis}.)

We show in Section \ref{subsec:UnitaryS1ActionOnS1DeformationComplex} that this modification of the action \eqref{eq:S1ZAction} defines an $S^1$-equivariant structure in the sense of Definition \ref{defn:GEquivariantStructure} on the elliptic deformation complex for the non-Abelian monopole equations.
This will imply that the image of the embedding \eqref{eq:LocalSliceBall} is closed under the modified action \eqref{eq:UnitaryS12ActionOnNonAbelianPairsAffineSpace}.  That fact and the $S^1$-equivariance of the map $\fS_{A,\Phi}$ allow us to construct the desired $S^1$-equivariant Kuranishi model in Section \ref{subsec:EquivKuranishiModelOfReducibleNonAbelianMonopole} using these modified $S^1$ actions.

\subsection{Modified circle action on the affine space of spin${}^u$ pairs}
\label{subsec:UnitaryS1ActionOnAffineSpace}
Let $\ft=(\rho,W\otimes E)$ be a \spinu structure over a closed, oriented, smooth four-dimensional manifold $X$. We begin by introducing a modification of the $S^1$ action \eqref{eq:S1ZAction}  on pairs on $\ft$. We assume that $E$ splits as a direct sum of Hermitian line bundles, $E=L_1\oplus L_2$. Recall that $\rho_2:S^1\to C^\infty(\U(E))$ is the homomorphism defined in \eqref{eq:DefineUnitaryS12ActionsAtReducible} using the splitting $E=L_1\oplus L_2$,
\[
\rho_2(e^{i\theta})=\id_{L_1}\oplus e^{i\theta}\, \id_{L_2}, \quad\text{for all } e^{i\theta} \in S^1.
\]
We define an $S^1$ action on the affine space of pairs by
\begin{multline}
\label{eq:UnitaryS12ActionOnNonAbelianPairsAffineSpace}
\rho_2^\sA: S^1\times \sA(E,h)\times\Om^0(W^+\otimes E) \ni \left(e^{i\theta},(A,\Phi)\right)
\\
\mapsto \left(\rho_2(e^{-i\theta})^* A,\rho_2(e^{i\theta})\Phi\right) \in \sA(E,h)\times\Om^0(W^+\otimes E),
\end{multline}
where, by the pullback convention described in Remark \ref{rmk:PushforwardPullbackNotation},
$\nabla_{\rho_2(e^{-i\theta})^* A} = \rho_2(e^{i\theta})\circ\nabla_A\circ \rho_2(e^{-i\theta})$. We define $\rho_2^\sA$ as in \eqref{eq:UnitaryS12ActionOnNonAbelianPairsAffineSpace} rather than by transposing $e^{i\theta}$ and $e^{-i\theta}$ so that, as we show in the following lemma, the $S^1$ action \eqref{eq:UnitaryS12ActionOnNonAbelianPairsAffineSpace}
defines an action on $\sC_\ft$ which is the same, up to positive multiplicity, as the standard action given in \eqref{eq:S1ZActionOnQuotientSpace}.

\begin{lem}[Gauge-equivalence of modified and standard circle actions on the affine space of \spinu pairs up to positive multiplicity]
\label{lem:EqualityOfUnitaryS1ActionsOnQuotientSpace}
Let $\ft=(\rho,W\otimes E)$ be a spin${}^u$ structure over an oriented, smooth four-dimensional Riemannian manifold $(X,g)$. If $E$ splits as a direct sum of complex line bundles, $E=L_1\oplus L_2$ then, for all $(A,\Phi)\in\sA(E,h)\times\Om^0(W^+\otimes E)$ and $e^{i\theta}\in S^1$,
\begin{equation}
\label{eq:L2andZActionsEqualOnQuotient}
\left[ \rho_2^\sA(e^{2i\theta})(A,\Phi)\right]
=
\left[ A,e^{i\theta}\Phi\right].
\end{equation}
\end{lem}

\begin{rmk}[Generalizations of Lemma \ref{lem:EqualityOfUnitaryS1ActionsOnQuotientSpace} to unitary pairs on vector bundles of arbitrary rank over manifolds of arbitrary dimension]
\label{rmk:ExtensionOfEqualityOfUnitaryS1ActionsOnQuotientSpace}
The proof of Lemma \ref{lem:EqualityOfUnitaryS1ActionsOnQuotientSpace} adapts \mutatis to prove the analogous result for unitary pairs $(A,\varphi)\in\sA(E,h)\times \Omega^0(E)$ over a smooth manifold of arbitrary dimension when $E$ splits as a direct sum of Hermitian line bundles.  The result also extends to the case when $E$ has higher rank and admits an orthogonal splitting $E=E_1\oplus E_2$, where $E_i$ is a proper Hermitian vector subbundle of rank $r_i$ for $i=1,2$.  When $E$ has rank $r\geq 2$, the relations \eqref{eq:UnitaryS1ActionsRelation} and \eqref{eq:L2andZActionsEqualOnQuotient} are replaced by the relations,
\begin{align}
  \label{eq:UnitaryS1ActionsRelation_rank-r}
  \rho_2(e^{ir\theta})
  &=
    \rho_{\SU}(e^{-i\theta})\rho_Z(e^{ir_2\theta}),
  \\
  \label{eq:L2andZActionsEqualOnQuotient_rank-r}
  [\rho_2^{\sA}(e^{ir\theta})(A,\varphi)]
  &=
    [A, e^{ir_2\theta}\varphi], \quad\text{for all } e^{i\theta}\in S^1,
\end{align}
and the maps $\rho_{\SU}:S^1\to C^\infty(\SU(E))$ and $\rho_2:S^1\to C^\infty(\U(E))$ in \eqref{eq:DefineUnitaryS1ActionsAtReducible} are replaced by
\begin{subequations}
\label{eq:DefineUnitaryS1ActionsAtReducible_rank-r}
\begin{align}
\label{eq:DefineUnitaryS1sActionsAtReducible_rank-r}
\rho_{\SU}(e^{i\theta})&:=e^{ir_2\theta}\,\id_{E_1}\oplus e^{-ir_1\theta}\,\id_{E_2},
\\
\label{eq:DefineUnitaryS12ActionsAtReducible_rank-r}
\rho_2(e^{i\theta})&:=\id_{E_1}\oplus e^{i\theta}\,\id_{E_2}, \quad\text{for all } e^{i\theta} \in S^1.
\end{align}
\end{subequations}
Compare \eqref{eq:UnitaryS1ActionsRelation_rank-r} with the identity \eqref{eq:ComparingS1ActionsAtReducible} in the forthcoming Lemma \ref{lem:RelateS1Actions} for $\CC^*$ actions.
\end{rmk}

\begin{proof}[Proof of Lemma \ref{lem:EqualityOfUnitaryS1ActionsOnQuotientSpace}]
Recall from \eqref{eq:UnitaryS1ActionsRelation} that the three homomorphisms $\rho_2$, $\rho_Z$, and $\rho_{\SU}$ from $S^1$ to $C^\infty(\U(E))$ defined in \eqref{eq:DefineUnitaryS1ActionsAtReducible} are related by
\[
\rho_2(e^{2i\theta})
=
\rho_{\SU}(e^{-i\theta})\rho_Z(e^{i\theta})
\quad\text{for all } e^{i\theta}\in S^1.
\]
For any $e^{i\theta}\in S^1$ and any $A\in \sA(E,h)$,
\begin{equation}
\label{eq:S1ZInStabilizer}
\rho_Z(e^{i\theta})^*A=A,
\end{equation}
because the image of $\rho_Z$ lies in the stabilizer of every connection. We obtain
\begin{align*}
\left[ \rho_2^\sA(e^{2i\theta})(A,\Phi)\right]
&=
\left[ \rho_2(e^{-2i\theta})^* A,\rho_2(e^{2i\theta})\Phi\right]
\quad\text{(by the definition \eqref{eq:UnitaryS12ActionOnNonAbelianPairsAffineSpace})}
\\
&=
\left[ \left(\rho_{\SU}(e^{i\theta})\rho_Z(e^{-i\theta})\right)^* A,\rho_{\SU}(e^{-i\theta})\rho_Z(e^{i\theta})\Phi\right]
\quad\text{(by \eqref{eq:UnitaryS1ActionsRelation})}
\\
&=\left[ \rho_{\SU}(e^{i\theta})^* A,\rho_{\SU}(e^{-i\theta})e^{i\theta}\Phi\right]
\quad\text{(by \eqref{eq:S1ZInStabilizer})}
\\
&=\left[ \rho_{\SU}(e^{i\theta})^*(A,e^{i\theta}\Phi\right]
\quad\text{(by \eqref{eq:PullPushActionsOnUnitaryPairs})}
\\
&=
\left[ A,e^{i\theta}\Phi\right]
\quad\text{(because $\rho_{\SU}(e^{-i\theta})\in C^\infty(\SU(E))$).}
\end{align*}
This proves \eqref{eq:L2andZActionsEqualOnQuotient} and hence the lemma.
\end{proof}

\begin{rmk}[On the definition of the modified circle action]
\label{rmk:CommentOnPositiveMultiplicity}
While the actions induced on $\sC_\ft$ by $\rho_2^\sA$ and by the standard action \eqref{eq:S1ZAction} differ by a factor of two as shown in \eqref{eq:L2andZActionsEqualOnQuotient}, to compute the virtual Morse--Bott index \eqref{eq:Virtual_Morse-Bott_index_moduli_space_non-abelian_monopoles} we only need to know the dimensions of the subspaces of $\bH_{A,\Phi}^k$ on which the $S^1$ action has negative weight.  Because it is only the \emph{sign} of the weight and not the magnitude which is important, the factor of two does not affect our computation of the virtual Morse--Bott index. However, the fact that the weights differ by a \emph{positive} factor is essential and this determines our use of the definition of $\rho_2^\sA$ given in \eqref{eq:UnitaryS12ActionOnNonAbelianPairsAffineSpace} rather than a definition with the roles of $e^{i\theta}$ and $e^{-i\theta}$ reversed.
\end{rmk}

While a split spin${}^u$ pair is a fixed point of the standard $S^1$ action \eqref{eq:S1ZActionOnQuotientSpace} on the quotient space $\sC_\ft$, a split spin${}^u$ pair is not a fixed point of the standard $S^1$ action  \eqref{eq:S1ZAction} on the affine space unless $\Phi\equiv 0$.  As we show in the following lemma, split \spinu pairs are fixed points of the $S^1$ action \eqref{eq:UnitaryS12ActionOnNonAbelianPairsAffineSpace} on the affine space of \spinu pairs.

\begin{lem}[Circle action on affine space of \spinu pairs]
\label{lem:S12ActionOnNonAbelianPairsAffineSpace}
Continue the hypotheses of Lemma \ref{lem:EqualityOfUnitaryS1ActionsOnQuotientSpace}. If $(A,\Phi)$ is a spin${}^u$ pair that is split in the sense of Definition \ref{defn:Split_trivial_central-stabilizer_spinor_pair} with respect to the orthogonal splitting $E = L_1 \oplus L_2$ as a direct sum of Hermitian line bundles and $\Phi=(\Phi_1,0)$ where $\Phi_1\in\Om^0(W^+\otimes L_1)$, then $(A,\Phi)$ is a fixed point of the $S^1$ action \eqref{eq:UnitaryS12ActionOnNonAbelianPairsAffineSpace}. If
\begin{equation}
\label{eq:UnitaryS12ActionOnTangentSpaceOfFixedPoint}
(D_2\rho_2^\sA):S^1\times \Om^1(\su(E))\times\Om^0(W^+\otimes E) \to \Om^1(\su(E))\times\Om^0(W^+\otimes E)
\end{equation}
denotes the linear action on the tangent space at $(A,\Phi)$ to the affine space of spin${}^u$ pairs induced\footnote{As described in Section \ref{sec:CircleActionsFixedPtsHessians}.} by the derivative of $\rho_2^\sA$ in the directions tangent to the affine space in the sense of \eqref{eq:Circle_action_tangent_bundle}, then
\begin{equation}
\label{eq:ExplicitUnitaryS12ActionOnTangentSpaceOfFixedPoint}
(D_2\rho_2^\sA)(e^{i\theta})\left(a,\phi\right)
=
\left( \rho_2(e^{i\theta})a \rho_2(e^{-i\theta}),\rho_2(e^{i\theta})\phi\right),
\end{equation}
for all $e^{i\theta}\in S^1$, and $a\in\Om^1(\su(E))$, and $\phi\in\Om^0(W^+\otimes E)$.
\end{lem}

\begin{rmk}[Generalization of Lemma \ref{lem:S12ActionOnNonAbelianPairsAffineSpace}]
\label{rmk:GeneralizationOfS12ActionOnNonAbelianPairsAffineSpace}
Lemma \ref{lem:S12ActionOnNonAbelianPairsAffineSpace} also holds for split unitary pairs on complex bundles $E$ of arbitrary rank over manifolds of arbitrary dimension.
\end{rmk}

\begin{proof}[Proof of Lemma \ref{lem:S12ActionOnNonAbelianPairsAffineSpace}]
By the definition of a split pair, we have $A=A_1\oplus A_2$, where $A_i$ is a unitary connection on the
Hermitian line bundle $L_i$ for $i=1,2$. By the definition of $\rho_2$ in \eqref{eq:DefineUnitaryS12ActionsAtReducible}, the gauge transformation  $\rho_2(e^{i\theta})$ preserves the splitting $L_1\oplus L_2$ and acts as the identity on $L_1$ and multiplication by $e^{i\theta}$ on $L_2$.  Hence, for all $e^{i\theta}\in S^1$, we have $\rho_2(e^{-i\theta})^*(A_1\oplus A_2)=A_1\oplus A_2$. Moreover, $\rho_2(e^{i\theta})(\Phi_1,0)=(\Phi_1,0)$. Thus, $(A,\Phi)$ is a fixed point of the action \eqref{eq:UnitaryS12ActionOnNonAbelianPairsAffineSpace}.

We claim that if $(A',\Phi') = (A,\Phi)+(a,\phi) \in\sA(E,h)\times\Om^0(W^+\otimes E)$, then
\begin{equation}
\label{eq:ExpansionOfUnitaryS12Action}
\rho_2^\sA(e^{i\theta})(A',\Phi')
=
(A,\Phi)+\left( \rho_2(e^{i\theta})a \rho_2(e^{-i\theta}),\rho_2(e^{i\theta})\phi\right).
\end{equation}
To see this, we compute
\begin{align*}
\rho_2^\sA(e^{i\theta})(A',\Phi')
&=
\rho_2^\sA(e^{i\theta})(A+a,\Phi+\phi)
\\
&=
\left( \rho_2(e^{-i\theta})^*(A+a),\rho_2(e^{i\theta})(\Phi+\phi)\right)
  \\
  &=
    \left(\rho_2(e^{-i\theta})^*A + \rho_2(e^{i\theta})a\rho_2(e^{-i\theta}),                                                                          \rho_2(e^{i\theta})\Phi+\rho_2(e^{i\theta})\phi)\right)
    \quad\text{(by \eqref{eq:UnitaryS12ActionOnNonAbelianPairsAffineSpace})}
\\
&=
\left(A + \rho_2(e^{i\theta})a \rho_2(e^{-i\theta}),\Phi+\rho_2(e^{i\theta})\phi)\right),
\end{align*}
where the last equality follows because $(A,\Phi)$ is a fixed point of the action \eqref{eq:UnitaryS12ActionOnNonAbelianPairsAffineSpace}.
The preceding identities prove the claim \eqref{eq:ExpansionOfUnitaryS12Action}. The expression for the derivative of $\rho_2^\sA$ in \eqref{eq:ExplicitUnitaryS12ActionOnTangentSpaceOfFixedPoint} follows from \eqref{eq:ExpansionOfUnitaryS12Action}.
\end{proof}

\subsection{Circle-equivariant structure on the elliptic deformation complex for a non-Abelian monopole}
\label{subsec:UnitaryS1ActionOnS1DeformationComplex}
We now prove that the elliptic deformation complex \eqref{eq:SO3MonopoleDefComplex} for a non-Abelian monopole $(A,\Phi)$
has an $S^1$-equivariant structure in the following sense.

\begin{defn}[Group-equivariant structures on complexes]
  \label{defn:GEquivariantStructure}
Let $\KK=\RR$ or $\CC$, and $G$ be a group, and $(\sF_\bullet,d_\bullet)$ be a complex,
\[
\dots
 \xrightarrow{d_{k-1}} \sF_k
 \xrightarrow{d_{k}} \sF_{k+1}
 \xrightarrow{d_{k+1}} \sF_{k+2}
 \xrightarrow{d_{k+2}} \dots,
\]
where each $\sF_k$ is a $\KK$-vector space and the differential $d_k \in \Hom_\KK(\sF_k,\sF_{k+1})$ obeys $d_{k+1}\circ d_k = 0$ for $k\in\ZZ$. A \emph{$G$-equivariant structure}  on $(\sF_\bullet,d_\bullet)$ is collection of $\KK$-linear group representations $\rho_k:G\to \GL(\sF_k,\KK)$, for $k\in\ZZ$, such that
\begin{equation}
\label{eq:EquivariantDiffDefinition}
\rho_{k+1}(g)\circ d_k=d_k\circ\rho_k(g), \quad\text{for all $g\in G$ and $k\in\ZZ$}.
\end{equation}
A \emph{real $G$-equivariant structure} is a \emph{$G$-equivariant structure} with $\KK=\RR$ and a
\emph{complex $G$-equivariant structure} is a \emph{$G$-equivariant structure} with $\KK=\CC$.
\end{defn}

\begin{rmk}[$G$-equivariant structures and $G$-modules]
By the equivalence between $\CC$-linear representations of a group $G$ and modules over the group ring $\CC G$ (see Rotman \cite[Proposition 8.37]{RotmanAlgebra_2002}), a complex has a complex $G$-equivariant structure if and only if it is a chain complex of $\CC G$-modules in the sense of Weibel \cite[Definition 1.1.1]{Weibel_introduction_homological_algebra}.
\end{rmk}

We recall that the elliptic deformation complex \eqref{eq:SO3MonopoleDefComplex} for a non-Abelian monopole $(A,\Phi)$ is given by
\[
\begin{CD}
0
@>>>
\Omega^0(\su(E))
@> d_{A,\Phi}^0>>
\begin{matrix}
\Omega^1(\su(E))
\\
\oplus
\\
\Omega^0(W^+\oplus E)
\end{matrix}
@> d_{A,\Phi}^1 >>
\begin{matrix}
\Omega^+(\su(E))
\\
\oplus
\\
\Omega^0(W^-\otimes E)
\end{matrix}
@>>>
0
\end{CD}
\]
We recall also that $\rho_2:S^1\to C^\infty(\U(E))$ is the homomorphism defined in \eqref{eq:DefineUnitaryS12ActionsAtReducible} using the splitting $E=L_1\oplus L_2$ as a direct sum of Hermitian line bundles.  We define linear $S^1$ actions on the three vector spaces appearing in the complex \eqref{eq:SO3MonopoleDefComplex} by
\begin{subequations}
\label{eq:S1ActionOnSO(3)MonopolesDefComplexSpaces}
\begin{equation}
\label{eq:S1ActionOnSO(3)MonopolesDefComplexSpacesZeroTerm}
S^1\times \Omega^0(\su(E)) \ni (e^{i\theta},\xi) \mapsto \rho_2(e^{i\theta})\xi \rho_2(e^{i\theta})^{-1} \in \Om^0(\su(E)),
\end{equation}
\begin{multline}
\label{eq:S1ActionOnSO(3)MonopolesDefComplexSpacesFirstTerm}
S^1\times\left( \Om^1(\su(E))\oplus\Om^0(W^+\otimes E)\right)
\ni (e^{i\theta},(a,\phi))
\\
\mapsto \left(\rho_2(e^{i\theta})a \rho_2(e^{i\theta})^{-1},\rho_2(e^{i\theta})\phi \right)
\in \Om^1(\su(E))\oplus\Om^0(W^+\otimes E),
\end{multline}
\begin{multline}
\label{eq:S1ActionOnSO(3)MonopolesDefComplexSpacesSecondTerm}
S^1\times\left( \Om^+(\su(E))\oplus\Om^0(W^-\otimes E)\right)
\ni (e^{i\theta},(v,\nu))
\\
\mapsto
\left(\rho_2(e^{i\theta})v \rho_2(e^{i\theta})^{-1},\rho_2(e^{i\theta})\nu \right)
\in \Om^+(\su(E))\oplus\Om^0(W^-\otimes E).
\end{multline}
\end{subequations}
Observe that the $S^1$ action \eqref{eq:S1ActionOnSO(3)MonopolesDefComplexSpacesFirstTerm}
is equal to the $S^1$ action given in \eqref{eq:ExplicitUnitaryS12ActionOnTangentSpaceOfFixedPoint}
and the action \eqref{eq:S1ActionOnSO(3)MonopolesDefComplexSpacesSecondTerm} is equal to the $S^1$ action defined in \eqref{eq:S1L2ActionOnnonAbelianMonopoleMapCodomain}.

\begin{lem}[Circle-equivariant structure on the elliptic deformation complex for non-Abelian monopoles]
\label{lem:S1_Equivariance_Of_nonAbelianMonopole_Def_Complex}
Let $\fs=(\rho,W\otimes L_1)$ be a spin${}^c$ structure and $\ft=(\rho,W\otimes E)$ be a spin${}^u$ structure over a four-dimensional, oriented, smooth Riemannian manifold $(X,g)$. If $(A,\Phi)\in\sA(E,h)\times \Omega^0(W^+\otimes E)$ is a non-Abelian monopole that is split in the sense of Definition \ref{defn:Split_trivial_central-stabilizer_spinor_pair} with respect to a decomposition $E=L_1\oplus L_2$ as an orthogonal direct sum of Hermitian line bundles and $\Phi=(\Phi_1,0)$ with $\Phi_1\in\Om^0(W^+\otimes L_1)$, then the group actions \eqref{eq:S1ActionOnSO(3)MonopolesDefComplexSpaces} define an $S^1$-equivariant structure in the sense of Definition \ref{defn:GEquivariantStructure} on the elliptic deformation complex \eqref{eq:SO3MonopoleDefComplex} for a non-Abelian monopole.
\end{lem}

\begin{rmk}[Generalizations of Lemma \ref{lem:S1_Equivariance_Of_nonAbelianMonopole_Def_Complex} and Corollaries  \ref{cor:S1_Action_on_HarmonicSections_Of_nonAbelianMonopole_Def_Complex} and
 \ref{cor:UnitaryS12ActionOnSlice}]
\label{rmk:GeneralizationOf_S1_Equivariance_Of_nonAbelianMonopole_Def_Complex}
The proof of Lemma \ref{lem:S1_Equivariance_Of_nonAbelianMonopole_Def_Complex} adapts \mutatis to prove that the elliptic complex \eqref{eq:Projective_vortex_elliptic_deformation_complex} for a split projective vortex on a vector bundle of higher rank over a complex manifold of arbitrary dimension admits an $S^1$-equivariant structure. One shows this by proving that the differentials $d_{A,\varphi}^k$ for $k=0,\dots,n$ (defined in
\eqref{eq:d0_projective_vortex_elliptic_deformation_complex},
\eqref{eq:d1_projective_vortex_elliptic_deformation_complex},
\eqref{eq:d2_projective_vortex}, and \eqref{eq:dk_projective_vortex}) in the elliptic complex for a projective vortex
satisfy the forthcoming gauge-equivariance condition \eqref{eq:GaugeEquivOfnonAbelianMonopoleDifferentials}.  Corollaries \ref{cor:S1_Action_on_HarmonicSections_Of_nonAbelianMonopole_Def_Complex} and
 \ref{cor:UnitaryS12ActionOnSlice} admit the same generalization.
\end{rmk}

\begin{proof}[Proof of Lemma \ref{lem:S1_Equivariance_Of_nonAbelianMonopole_Def_Complex}]
For any gauge transformation $u\in C^\infty(\U(E))$, the differentials $d_{A,\Phi}^k$ in the non-Abelian monopole elliptic deformation complex \eqref{eq:SO3MonopoleDefComplex} satisfy
\begin{equation}
\label{eq:GaugeEquivOfnonAbelianMonopoleDifferentials}
d_{u^*(A,\Phi)}^k
=
u^{-1}\circ d_{A,\Phi}^k\circ u, \quad\text{for } k = 0,1,
\end{equation}
where $u^*(A,\Phi)=(u^*A,u^{-1}\Phi)$ as defined in \eqref{eq:GaugeActionOnSpinuPairs} with the convention described following \eqref{eq:GaugeActionOnSpinuPairs} that $d_{u^*A}=u^{-1}\circ d_A\circ u$. One can see that \eqref{eq:GaugeEquivOfnonAbelianMonopoleDifferentials} holds by direct computation from the definitions of the differentials in \eqref{eq:d_APhi^0} and \eqref{eq:d1OfSO3MonopoleComplex}. Observe that $C^\infty(\U(E))$ acts on $\Om^k(\su(E))$ by
\[
C^\infty(\U(E)) \times \in\Om^k(\su(E)) \ni (u,\omega)
  \mapsto \ad(u)\om= u\om u^{-1}
  \in \Om^k(\su(E))
\]
and therefore
\begin{equation}
\label{eq:GaugeTransformationAndCovariantDerivative}
(u^{-1}\circ d_A\circ u)\om
=
\ad(u^{-1})\left( d_A (\ad(u)\om)\right).
\end{equation}
We further note that $C^\infty(\U(E))$ acts on linear spaces of pairs by
\begin{multline}
  \label{eq:GaugeTransformationActionOnPairs}
  C^\infty(\U(E)) \times \Om^k(\su(E))\oplus\Om^0(W^\pm\otimes E) \ni (u,(\om,\phi))
  \\
  \mapsto u(\om,\phi) := (\ad(u)\om,u\phi)
  \in \Om^k(\su(E))\oplus\Om^0(W^\pm\otimes E).
\end{multline}
Observe that the $S^1$ actions \eqref{eq:S1ActionOnSO(3)MonopolesDefComplexSpaces} are induced by
the action \eqref{eq:GaugeTransformationActionOnPairs} with $u=\rho_2(e^{i\theta})$. Hence, for all $\xi\in\Om^0(\su(E))$ and $u\in C^\infty(\U(E))$ we have
\begin{align*}
d_{u^*(A,\Phi)}^0 \xi &=
                        (d_{u^*A}\xi, -\xi u^{-1}\Phi)
                        \\
                        &\qquad\text{(by the definition of $d_{A,\Phi}^0$ in \eqref{eq:d_APhi^0} and $u^*(A,\Phi)=(u^*A,u^{-1}\Phi)$)}
\\
&=
     \left(( u^{-1}\circ d_A\circ u)\xi, -u^{-1}(\ad(u)\xi)\Phi\right)
     \quad\text{(by \eqref{eq:GaugeEquivOfnonAbelianMonopoleDifferentials})}
\\
                     &=
\left(\ad(u^{-1})\left( d_A(\ad(u)\xi)\right),-u^{-1}((\ad(u)\xi)\Phi)\right)
\quad\text{(by \eqref{eq:GaugeTransformationAndCovariantDerivative})}
\\
&=
u^{-1}\left(d_A(\ad(u)\xi),-(\ad(u)\xi)\Phi\right)
\quad\text{(by \eqref{eq:GaugeTransformationActionOnPairs})}
\\
&=
u^{-1} \left(d_{A,\Phi}^0(\ad(u)\xi)\right)
\quad\text{(by \eqref{eq:d_APhi^0})}
\\
&=
(u^{-1}\circ d_{A,\Phi}^0\circ u )\xi,
\end{align*}
which proves that \eqref{eq:GaugeEquivOfnonAbelianMonopoleDifferentials} holds for $k=0$.
We give an alternative proof of \eqref{eq:GaugeEquivOfnonAbelianMonopoleDifferentials} for $k=0$ in
Remark \ref{rmk:AlternativeProofOfEquivarianceofd0} by using the fact that $d_{A,\Phi}^0$ is the derivative of the action of the group of gauge transformations at the identity.

A similar computation proves that \eqref{eq:GaugeEquivOfnonAbelianMonopoleDifferentials} holds for $k=1$.  Indeed, for all $a\in\Om^1(\su(E))$, $\phi\in \Om^0(W^+\otimes E)$, and $u\in C^\infty(\U(E))$ we have
\begin{align*}
d_{u^*(A,\Phi)}^1 (a,\phi)&=
\begin{pmatrix}
d_{u^*A}^+a -(u^{-1}\Phi\otimes \phi^*+\phi\otimes (u^{-1}\Phi)^*)_{00}
\\
D_{u^*A}\phi+\rho(a) u^{-1}\Phi
\end{pmatrix}
\\
& \quad\quad\text{(by definition of $d_{A,\Phi}^1$ in  \eqref{eq:d1OfSO3MonopoleComplex} and $u^*(A,\Phi)=(u^*A,u^{-1}\Phi)$)}
\\
&=
\begin{pmatrix}
\ad(u^{-1})\left( d_A^+ (\ad(u)a)\right) - u^{-1}\left( \Phi\otimes (u\phi)^*+u\phi\otimes \Phi^* \right)_{00}u
\\
u^{-1} D_A (u\phi) +u^{-1}\rho(\ad(u)a)\Phi
\end{pmatrix}
  \\
&\quad
\quad\text{(by \eqref{eq:GaugeTransformationAndCovariantDerivative} and $(u\phi)^* = u^{-1}\phi^*$, with $\phi^* = \langle\cdot,\phi\rangle_{W^+\otimes E})$ by \eqref{eq:SO(3)_monopole_equations})}
\\
&=
\begin{pmatrix}
\ad(u^{-1})\left( d_A^+ (\ad(u)a) -(\Phi\otimes (u\phi)^*+u\phi\otimes \Phi^*)_{00}\right)
\\
u^{-1}\left( D_A(u\phi)+ \rho(\ad(u)a)\Phi\right)
\end{pmatrix}
\\
&=
(u^{-1}\circ d_{A,\Phi}^1\circ u)(a,\phi)
\quad\text{(by \eqref{eq:GaugeTransformationActionOnPairs})},
\end{align*}
which proves that \eqref{eq:GaugeEquivOfnonAbelianMonopoleDifferentials} holds for $k=1$. Alternately, one  can observe that $d_{A,\Phi}^1$ is the derivative of the map $\fS$ defined in
\eqref{eq:PerturbedSO3MonopoleEquation_map} by the non-Abelian monopole equations and use the equivariance of
$\fS$ described in  Lemma \ref{lem:S1ActionsForReduciblePairGerm},
noting that the derivative of the action \eqref{eq:S1L2ActionOnnonAbelianMonopoleMapCodomain} on the domain of
$\fS$ is given by \eqref{eq:S1ActionOnSO(3)MonopolesDefComplexSpacesFirstTerm} (as shown in Lemma \ref{lem:S12ActionOnNonAbelianPairsAffineSpace}) and the action \eqref{eq:S1L2ActionOnnonAbelianMonopoleMapCodomain} on the codomain of $\fS$ equals the action \eqref{eq:S1ActionOnSO(3)MonopolesDefComplexSpacesSecondTerm}.

Because $(A,\Phi)$ is a fixed point of the action $\rho_2^\sA$ by Lemma \ref{lem:S12ActionOnNonAbelianPairsAffineSpace}, the definition of the action $\rho_2^\sA$ in \eqref{eq:UnitaryS12ActionOnNonAbelianPairsAffineSpace} implies that
\[
(A,\Phi)
=
\rho_2^\sA(e^{i\theta})(A,\Phi)
=
(\rho_2(e^{i\theta})_*A,\rho_2(e^{i\theta})\Phi),\quad\text{for all $e^{i\theta}\in S^1$}.
\]
Therefore, the circle equivariance of $d_{A,\Phi}^k$ given in \eqref{eq:GaugeEquivOfnonAbelianMonopoleDifferentials} implies that
\[
d_{A,\Phi}^k
=
\rho_2 (e^{-i\theta})\circ d_{A,\Phi}^k\circ\rho_2(e^{i\theta}),
\quad\text{for all } e^{i\theta}\in S^1.
\]
The preceding equality and the observation that the $S^1$ actions in \eqref{eq:S1ActionOnSO(3)MonopolesDefComplexSpaces} are induced by
the action \eqref{eq:GaugeTransformationActionOnPairs} with $u=\rho_2(e^{i\theta})$ imply that the differentials $d_{A,\Phi}^k$ satisfy the equalities \eqref{eq:EquivariantDiffDefinition} required to make the complex $S^1$-equivariant with respect to the $S^1$ actions \eqref{eq:S1ActionOnSO(3)MonopolesDefComplexSpaces}.
\end{proof}

We now present an alternative proof of \eqref{eq:GaugeEquivOfnonAbelianMonopoleDifferentials} for $k=0$.

\begin{rmk}[Alternative proof of \eqref{eq:GaugeEquivOfnonAbelianMonopoleDifferentials} for $k=0$]
\label{rmk:AlternativeProofOfEquivarianceofd0}
By the same argument which led to \eqref{eq:ExpansionOfUnitaryS12Action}, we see that for
any spin${}^u$ pair $(A',\Phi')=(A,\Phi)+(a,\phi)$, where $(a,\phi)\in\Om^1(\su(E))\times\Om^0(W^+\otimes E)$,
and any $u\in\Om^0(\U(E))$, we have
\[
u^*(A',\Phi')=u^*(A,\Phi)+ (u^{-1}au,u^{-1}\phi).
\]
Hence, the differential of the right (pull-back) action of $u$,
\[
R_u:\sA(E,h)\times\Om^0(W^+\otimes E)\ni (A',\Phi')\mapsto 
u^*(A',\Phi')\in \sA(E,h)\times\Om^0(W^+\otimes E),
\]
at $(A,\Phi)$ is given by
\begin{equation}
\label{eq:DifferentialOfPullbackMap}
DR_u(A,\Phi)(a,\phi)=(u^{-1}au,u^{-1}\phi),
\quad\text{for all $(a,\phi)\in\Om^1(\su(E))\times\Om^0(W^+\otimes E)$.}
\end{equation}
Thus, for any spin${}^u$ pair $(A',\Phi')$, combining \eqref{eq:Differential_of_SU(E)_GaugeAction} and \eqref{eq:d_APhi^0} yields
\begin{equation}
\label{eq:d0IsDerivativeOfGaugeGroupAction}
d_{A',\Phi'}^0\xi
=
\left.\frac{d}{dt}\left(\Exp(t\xi)\right)^*(A',\Phi')\right|_{t=0}, 
\quad\text{for all } \xi\in\Om^0(\su(E)).
\end{equation}
Thus,
\begin{align*}
d_{u^*(A,\Phi)}^0\xi
  &=
    \left.\frac{d}{dt}\left(\Exp(t\xi)\right)^*u^*(A',\Phi')\right|_{t=0}
\quad\text{(by \eqref{eq:d0IsDerivativeOfGaugeGroupAction})}
\\
  &=
    \left.\frac{d}{dt}\left(u\Exp(t\xi)\right)^*(A',\Phi')\right|_{t=0}
\\
  &=
    \left.\frac{d}{dt}\left(\Exp(t\ad(u)\xi)u\right)^*(A',\Phi')\right|_{t=0}
\\
&=
\left.\frac{d}{dt}u^*\left(\Exp(t\ad(u)\xi)\right)^*(A',\Phi')\right|_{t=0}
\\
  &=
    DR_u(A,\Phi)\left.\frac{d}{dt}\left(\Exp(t\ad(u)\xi)\right)^*(A',\Phi')\right|_{t=0}
\quad\text{(by the Chain Rule)}
\\
&=
u^{-1} d_{A,\Phi}^0(\ad(u)\xi)
\quad\text{(by \eqref{eq:DifferentialOfPullbackMap} and \eqref{eq:d0IsDerivativeOfGaugeGroupAction})}
\\
&=
(u^{-1}\circ d_{A,\Phi}^0\circ u )\xi.
\end{align*}
This completes the alternative proof of  \eqref{eq:GaugeEquivOfnonAbelianMonopoleDifferentials} for $k=0$.
\end{rmk}

Lemma \ref{lem:S1_Equivariance_Of_nonAbelianMonopole_Def_Complex}, the definition of the harmonic spaces $\bH_{A,\Phi}^k$ in \eqref{eq:H_APhi^bullet}, and the fact that the actions  \eqref{eq:S1ActionOnSO(3)MonopolesDefComplexSpaces} are isometries with respect to the $L^2$ inner products yield the following

\begin{cor}[Circle actions on harmonic spaces for the elliptic deformation complex of a split non-Abelian monopole]
\label{cor:S1_Action_on_HarmonicSections_Of_nonAbelianMonopole_Def_Complex}
Continue the assumptions of Lemma \ref{lem:S1_Equivariance_Of_nonAbelianMonopole_Def_Complex} and assume that $X$ is closed. Then the circle actions \eqref{eq:S1ActionOnSO(3)MonopolesDefComplexSpaces} induce circle actions on the harmonic spaces $\bH_{A,\Phi}^k$ in \eqref{eq:H_APhi^bullet}, for $k=0,1,2$.
\end{cor}

We now show that the slice space \eqref{eq:LocalSliceBall} is closed under the action  \eqref{eq:UnitaryS12ActionOnNonAbelianPairsAffineSpace}.

\begin{cor}[Modified $S^1$ action on $\bB_{A,\Phi}$]
\label{cor:UnitaryS12ActionOnSlice}
Continue the assumptions of Lemma \ref{lem:S1_Equivariance_Of_nonAbelianMonopole_Def_Complex} assume that $X$ is closed.
Then the ball $\bB_{A,\Phi}\subset \Ker d_{A,\Phi}^{0,*}$ defined in \eqref{eq:SliceBall} is closed under the action
\eqref{eq:S1ActionOnSO(3)MonopolesDefComplexSpacesFirstTerm} and the embedding,
\begin{equation}
\label{eq:SliceBallEmbeddingInQuotientSpace}
\bB_{A,\Phi}\ni (a,\phi)\mapsto (A+a,\Phi+\phi)\in\sA(E,h)\times \Om^0(W^+\otimes E)
\end{equation}
is equivariant with respect to the action \eqref{eq:S1ActionOnSO(3)MonopolesDefComplexSpacesFirstTerm}
on its domain and the action \eqref{eq:UnitaryS12ActionOnNonAbelianPairsAffineSpace} on its codomain.
\end{cor}

\begin{proof}
Recall from its definition in \eqref{eq:SliceBall} that $\bB_{A,\Phi}$ is a ball with respect to the $W^{1,p}_A$ norm around the origin in $\Ker d_{A,\Phi}^{0,*}$.
Lemma \ref{lem:S1_Equivariance_Of_nonAbelianMonopole_Def_Complex} and the fact that the  circle actions \eqref{eq:S1ActionOnSO(3)MonopolesDefComplexSpaces} are isometries imply that $d_{A,\Phi}^{0,*}$ is equivariant with respect to these actions. This $S^1$-equivariance of $d_{A,\Phi}^{0,*}$ implies that $\Ker d_{A,\Phi}^{0,*}$ is closed under the action \eqref{eq:S1ActionOnSO(3)MonopolesDefComplexSpacesFirstTerm}. Because the $S^1$ action \eqref{eq:S1ActionOnSO(3)MonopolesDefComplexSpacesFirstTerm} preserves the $W^{1,p}_A$ norm, the ball $\bB_{A,\Phi}$ in $\Ker d_{A,\Phi}^{0,*}$ is closed under this action. The asserted equivariance of the embedding \eqref{eq:SliceBallEmbeddingInQuotientSpace} then follows from the
expression in \eqref{eq:ExpansionOfUnitaryS12Action} for the action \eqref{eq:UnitaryS12ActionOnNonAbelianPairsAffineSpace}.
\end{proof}

\subsection{Application to the construction of circle-equivariant Kuranishi models}
\label{subsec:EquivKuranishiModelOfReducibleNonAbelianMonopole}
We construct a Kuranishi model for an open neighborhood of a point $[A,\Phi]$ in $\sM_\ft$ represented by a split pair $(A,\Phi)$, focusing on its circle-equivariance properties. The constructions of circle-equivariant Kuranishi models associated to projective vortices and holomorphic pairs were outlined in Section \ref{subsec:Circle-invariant_complex_Kaehler_local_virtual_moduli_space_projective_vortices_near_fixed_points}.

In order to apply our abstract result for equivariant Kuranishi models (Lemma \ref{lem:EquivariantKuranishiLemma}), we first show that the restriction of $d_{A,\Phi}^1$ to $\Ker d_{A,\Phi}^{0,*}$ defines splittings of the domain and codomain of this restriction of $d_{A,\Phi}^1$ that are closed under the $S^1$ actions, as in
the hypotheses of Lemma \ref{lem:EquivariantKuranishiLemma}.

\begin{lem}[Circle invariance of the splittings of the domain and codomain of the differential $d_{A,\Phi}^1$]
\label{lem:FredholmMapAndGEquivariantSplitting}
Continue the assumptions of Lemma \ref{lem:S1_Equivariance_Of_nonAbelianMonopole_Def_Complex} and assume that $X$ is closed. Recall that $d_{A,\Phi}^{0,*}$ is the $L^2$-adjoint of the operator $d_{A,\Phi}^0$ defined in \eqref{eq:d_APhi^0}. Then the subspace,
\[
E_1:=\Ker d_{A,\Phi}^{0,*}\subset W^{1,p}(T^*X\otimes \su(E))\oplus W^{1,p}(W^+\otimes E),
\]
is closed under the $S^1$ action \eqref{eq:S1ActionOnSO(3)MonopolesDefComplexSpacesFirstTerm}. If we write $L$ for the restriction of the operator $d_{A,\Phi}^1$ defined in \eqref{eq:d1OfSO3MonopoleComplex} to $E_1$ and denote
\[
E_2:= L^p(\La^+\otimes\su(E))\oplus L^p(W^-\otimes E),
\]
and\footnote{This involves no loss of generality because, according to Feehan and Leness \cite[Proposition 3.7, p. 323]{FL1}, any $W^{2,2}$ non-Abelian monopole $(A,\Phi)$ is equivalent via a $W^{2,3}$ unitary gauge transformation to a $C^\infty$ non-Abelian monopole and the argument for a $W^{1,p}$ non-Abelian monopole is very similar.} the pair $(A,\Phi)$ is $C^\infty$, then the Banach spaces $E_1$ and $E_2$ admit splittings
\begin{align*}
E_1= \Ker L \oplus F_1,
\\
E_2=\Ran L \oplus C_2,
\end{align*}
where $F_1$ is the $L^2$-orthogonal complement of $\Ker L$ in $E_1$ and $C_2$ is the $L^2$-orthogonal complement of $\Ran L$ in $E_2$. The subspaces $\Ker d_{A,\Phi}^1$ and $F_1$ are closed under the $S^1$ action \eqref{eq:S1ActionOnSO(3)MonopolesDefComplexSpacesFirstTerm} on $E_1$ and the subspaces $\Ran d_{A,\Phi}^1$ and $C_2$ are closed under the $S^1$ action \eqref{eq:S1ActionOnSO(3)MonopolesDefComplexSpacesSecondTerm}  on $E_2$.
\end{lem}

Before proving Lemma \ref{lem:FredholmMapAndGEquivariantSplitting}, we recall the following elementary technical result.

\begin{lem}[Group invariance of the orthogonal complement of an invariant subspace]
\label{lem:HilbertSpaceClosedUnderAction}
(See Br\"ocker and tom Dieck \cite[Chapter II, Proof of Proposition 1.9, p. 68]{BrockertomDieck} or Knapp \cite[Chapter IV, Equation (4.5), p. 239]{Knapp_1986}.)  
Let $\rho:G\to \Or(H)$ be an orthogonal representation of a group $G$ on a Hilbert space $H$. If $V \subset H$ is a linear subspace that is $\rho(G)$-invariant, then the orthogonal complement $V^\perp$ of $V$ in $H$ is also $\rho(G)$-invariant.
\end{lem}

\begin{proof}
By the hypothesis that $\rho$ is an orthogonal representation (and a group homomorphism), $\rho(g)^* = \rho(g)^{-1} = \rho(g^{-1})$, for all $g \in G$, where the relation $\langle\rho(g)^*v,w\rangle_H = \langle v,\rho(g)w\rangle_H$, for all $v,w\in V$ and $g\in G$, defines the adjoint $\rho(g)^*\in\Or(H)$ of $\rho(g)$ with respect to the inner product $\langle\cdot,\cdot\rangle_H$ on $H$. If $w\in V^\perp$, then for all $v\in V$ and $g\in G$ we have $\rho(g^{-1})v \in V$ since $V$ is $\rho(G)$-invariant, and thus
\[
  \langle v,\rho(g)w\rangle_H
  = \left\langle \rho(g)^*v,w\right\rangle_H
  = \left\langle\rho(g^{-1})v,w\right\rangle_H = 0.
\]
Hence, $\rho(g)w \in V^\perp$, for all $w\in V^\perp$ and $g\in G$, as asserted.
\end{proof}

We now give the

\begin{proof}[Proof of Lemma \ref{lem:FredholmMapAndGEquivariantSplitting}]
Because $d_{A,\Phi}^0$ is $S^1$-equivariant with respect to the $S^1$ actions
\eqref{eq:S1ActionOnSO(3)MonopolesDefComplexSpacesZeroTerm} on its domain
and \eqref{eq:S1ActionOnSO(3)MonopolesDefComplexSpacesFirstTerm} on its codomain by Lemma \ref{lem:S1_Equivariance_Of_nonAbelianMonopole_Def_Complex} and because these $S^1$ actions are orthogonal with respect to the $L^2$ inner product, the operator $d_{A,\Phi}^{0,*}$ is also $S^1$-equivariant.  Hence, the kernel of $d_{A,\Phi}^{0,*}$ is invariant under the $S^1$ action \eqref{eq:S1ActionOnSO(3)MonopolesDefComplexSpacesFirstTerm}.

The equivariance of $d_{A,\Phi}^1$ with respect to the $S^1$ action \eqref{eq:S1ActionOnSO(3)MonopolesDefComplexSpacesFirstTerm} on its domain and the $S^1$ action
\eqref{eq:S1ActionOnSO(3)MonopolesDefComplexSpacesSecondTerm} on its codomain (shown in Lemma \ref{lem:S1_Equivariance_Of_nonAbelianMonopole_Def_Complex}) imply that $\Ker L$ and $\Ran L$ are invariant under the appropriate $S^1$ actions.

By hypothesis, the pair $(A,\Phi)$ is $C^\infty$ and the first-order differential operator
\begin{multline*}
  d_{A,\Phi}^1 + d_{A,\Phi}^{0,*}:
  W^{1,p}(T^*X\otimes \su(E))\oplus W^{1,p}(W^+\otimes E)
  \\
  \to
  L^p(\su(E)) \oplus L^p(\La^+\otimes \su(E))\oplus L^p(W^-\otimes E)
\end{multline*}
is elliptic with $C^\infty$ coefficients and Fredholm (see Feehan \cite{Feehan_yang_mills_gradient_flow_v4}, Gilkey \cite[Lemma 1.4.5]{Gilkey2}, or Hormander \cite[Theorem 19.2.1]{Hormander_v3}). Consequently, we have the following inclusions of finite-dimensional subspaces,
\[
  \Ker\left(d_{A,\Phi}^1 + d_{A,\Phi}^{0,*}\right)
  \subset
   C^\infty(T^*X\otimes \su(E))\oplus C^\infty(W^+\otimes E)
\]
and
\begin{multline*}
  \left(\Ran\left(d_{A,\Phi}^1 + d_{A,\Phi}^{0,*}\right)\right)^\perp
  =
  \Ker\left(d_{A,\Phi}^1 + d_{A,\Phi}^{0,*}\right)^*
  =
  \Ker\left(d_{A,\Phi}^{1,*} + d_{A,\Phi}^0\right)
  \\
  \subset
  C^\infty(\su(E)) \oplus C^\infty(\La^+\otimes \su(E))\oplus C^\infty(W^-\otimes E).
\end{multline*}
Moreover, there are Banach space splittings that are $L^2$-orthogonal,
\[
  W^{1,p}(T^*X\otimes \su(E))\oplus W^{1,p}(W^+\otimes E)
  =
  \Ker\left(d_{A,\Phi}^1 + d_{A,\Phi}^{0,*}\right)
  \oplus
  \left(\Ker\left(d_{A,\Phi}^1 + d_{A,\Phi}^{0,*}\right)\right)^\perp
\]
and, noting that $d_{A,\Phi}^1 + d_{A,\Phi}^{0,*}$ has closed range,
\[
  L^p(\su(E)) \oplus L^p(\La^+\otimes \su(E))\oplus L^p(W^-\otimes E)
  =
  \Ran\left(d_{A,\Phi}^1 + d_{A,\Phi}^{0,*}\right)
  \oplus
  \left(\Ran\left(d_{A,\Phi}^1 + d_{A,\Phi}^{0,*}\right)\right)^\perp.
\]
Since $L = (d_{A,\Phi}^1 + d_{A,\Phi}^{0,*})\restriction \Ker d_{A,\Phi}^{0,*}$, then $L$ also has closed range and we obtain Banach space splittings that are $L^2$-orthogonal:
\[
  E_1 = \Ker L \oplus (\Ker L)^\perp
  \quad\text{and}\quad
  E_2 = \Ran L \oplus (\Ran L)^\perp.
\]
In particular, $E_1 = \Ker L \oplus F_1$ and $E_2 = \Ran L \oplus F_2$ are Banach space splittings that are $L^2$-orthogonal, where $F_1:= (\Ker L)^\perp\cap E_1$ and $C_2 := (\Ran L)^\perp\cap E_2$.

Since $\Ker L$ and $\Ran L$ are invariant under the $S^1$ actions \eqref{eq:S1ActionOnSO(3)MonopolesDefComplexSpacesFirstTerm} and \eqref{eq:S1ActionOnSO(3)MonopolesDefComplexSpacesSecondTerm}, respectively, Lemma \ref{lem:HilbertSpaceClosedUnderAction} implies that $F_1$ and $C_2$ are invariant under the $S^1$ actions \eqref{eq:S1ActionOnSO(3)MonopolesDefComplexSpacesFirstTerm} and \eqref{eq:S1ActionOnSO(3)MonopolesDefComplexSpacesSecondTerm}, respectively. 
\end{proof}

We have previously constructed Kuranishi models associated to projective vortices in Definition \ref{defn:Friedman_Morgan_4-1-15_projective_vortex} and associated to holomorphic pairs in Definition \ref{defn:Friedman_Morgan_4-1-15_holomorphic_pairs} and outlined their $S^1$-equivariance properties in Section \ref{subsec:Circle-invariant_complex_Kaehler_local_virtual_moduli_space_projective_vortices_near_fixed_points}. We now construct an $S^1$-equivariant Kuranishi model associated to a non-Abelian monopole.

\begin{lem}[Circle-equivariance of the local Kuranishi model for an open neighborhood of a point in the moduli space of non-Abelian monopoles represented by a split pair]
\label{lem:S1EquivariantKuranishiModelForReducibleIn_non_AbelianMonopoleModuli}
Continue the assumptions of Lemma \ref{lem:S1_Equivariance_Of_nonAbelianMonopole_Def_Complex} and assume that $X$ is closed and $\Phi\not\equiv 0$, so that $\Stab(A,\Phi) = \{\id_E\}$ by Lemma \ref{lem:NonZeroSection_Spinu_pairs_Have_Trivial_Stabilizer}. Then, with respect to the $S^1$ actions on $\bH_{A,\Phi}^1$ and $\bH_{A,\Phi}^2$ given by Corollary \ref{cor:S1_Action_on_HarmonicSections_Of_nonAbelianMonopole_Def_Complex}, there are a circle-invariant open neighborhood $\sU\subset\bH_{A,\Phi}^1$ of the origin and an $S^1$-equivariant analytic map,
\[
\bkappa:\bH_{A,\Phi}^1 \supset \sU \to \bH_{A,\Phi}^2,
\]
and an analytic embedding $\beps:\bH_{A,\Phi}^1 \supset \sU\to\sC_\ft$ that is $S^1$-equivariant with respect to the $S^1$ action on $\bH_{A,\Phi}^1$ given by Corollary \ref{cor:S1_Action_on_HarmonicSections_Of_nonAbelianMonopole_Def_Complex}
and the $S^1$ action \eqref{eq:S12ActionsOnSpinuQuotient} on $\sC_\ft$ such that $\beps$ gives an $S^1$-equivariant homeomorphism between $\bkappa^{-1}(0)$ and an open neighborhood of $[A,\Phi]$ in $\sM_\ft$.
\end{lem}

\begin{rmk}[Generalizations of Lemma \ref{lem:S1EquivariantKuranishiModelForReducibleIn_non_AbelianMonopoleModuli} to projective vortices on vector bundles of arbitrary rank over manifolds of arbitrary dimension]
\label{rmk:GeneralizingKuranishiModels}
The Kuranishi model for a projective vortex on a vector bundle of arbitrary rank over a complex Hermitian manifold of arbitrary dimension is provided by Theorem \ref{thm:Local_Kuranishi_model_for_moduli_space_projective_vortices_StabAvarphi_idE}. If the projective vortex $(A,\varphi)$ is split with respect to an orthogonal decomposition $E=E_1\oplus E_2$ as a direct sum of proper Hermitian vector subbundles and $\Stab(A,\varphi)=\{\id_E\}$, then the ingredients in the proof of Lemma \ref{lem:S1EquivariantKuranishiModelForReducibleIn_non_AbelianMonopoleModuli} may be used to show that the Kuranishi model in Theorem \ref{thm:Local_Kuranishi_model_for_moduli_space_projective_vortices_StabAvarphi_idE} is $S^1$-equivariant with respect to the $S^1$ actions on the harmonic spaces $\bH_{A,\varphi}^1$ and $\bH_{A,\varphi}^2$ discussed in Remark \ref{rmk:GeneralizationOf_S1_Equivariance_Of_nonAbelianMonopole_Def_Complex}: see Section \ref{subsec:Circle-invariant_complex_Kaehler_local_virtual_moduli_space_projective_vortices_near_fixed_points} for an outline of the proof of $S^1$-equivariance.
\end{rmk}

\begin{rmk}[Holomorphic embedding and obstruction maps in Lemma \ref{lem:S1EquivariantKuranishiModelForReducibleIn_non_AbelianMonopoleModuli}]
\label{rmk:S1EquivariantKuranishiModelForReducibleIn_non_AbelianMonopoleModuli_holomorphic_maps}  
By the identification of type $1$ non-Abelian monopoles over complex K\"ahler surfaces with projective vortices (see Remark \ref{rmk:Projective_vortices_type_1_monopoles}), the pair $(\beps,\bkappa)$ of real analytic maps in Lemma \ref{lem:S1EquivariantKuranishiModelForReducibleIn_non_AbelianMonopoleModuli} may be assumed to be holomorphic via the isomorphism of (germs of) real analytic spaces provided by Theorem \ref{thm:Friedman_Morgan_4-3-8_projective_vortices} between the Kuranishi model $\sK(A,\varphi)$ for a projective vortex and the Kuranishi model $\fK(\bar\partial_A,\varphi)$ for the corresponding holomorphic pair.
\end{rmk}

\begin{proof}[Proof of Lemma \ref{lem:S1EquivariantKuranishiModelForReducibleIn_non_AbelianMonopoleModuli}]
Because $\Phi\not\equiv 0$ and $E$ has rank two, Lemma \ref{lem:NonZeroSection_Spinu_pairs_Have_Trivial_Stabilizer} implies that $\Stab(A,\Phi)=\{\id_E\}$. Let
\[
\bB_{A,\Phi}\subset \Ker d_{A,\Phi}^{0,*}
\]
be the ball defined in \eqref{eq:SliceBall} with the property that the map 
\begin{equation}
\label{eq:SliceBallEmbeddingInQuotient}
\bB_{A,\Phi}\ni (a,\phi)\mapsto [A+a,\Phi+\phi]\in\sC_\ft^0
\end{equation}
defines an embedding $\bB_{A,\Phi}\to \sC_\ft^0$.  By Feehan and Maridakis  \cite[Corollary 18, p. 19]{Feehan_Maridakis_Lojasiewicz-Simon_coupled_Yang-Mills}, the inverse map to the embedding \eqref{eq:SliceBallEmbeddingInQuotient} provides a chart for an analytic\footnote{The cited result only applies to the quotient space of connections but, as noted following \cite[Corollary 18, p. 19]{Feehan_Maridakis_Lojasiewicz-Simon_coupled_Yang-Mills}, the proof given there translates easily to the quotient space of pairs.} Banach manifold structure on $\sC_\ft^0$. Thus, if we write $\fS_{A,\Phi}$ for the map
\begin{equation}
\label{eq:PU(2)MonopoleEquationOnSlice}
\bB_{A,\Phi} \ni (a,\phi)
\mapsto
\fS_{A,\Phi}(a,\phi)
:=
\fS(A+a,\Phi+\varphi) \in L^p(\La^+\otimes\su(E))\oplus L^p(W^-\otimes E),
\end{equation}
where $\fS$ is defined in \eqref{eq:PerturbedSO3MonopoleEquation_map}, then, after possibly shrinking $\bB_{A,\Phi}$,  the embedding \eqref{eq:SliceBallEmbeddingInQuotient} gives an analytic embedding of $\bB_{A,\Phi}\cap\fS_{A,\Phi}^{-1}(0)$ onto an open neighborhood of $[A,\Phi]$ in $\sM_\ft$.

By Corollary \ref{cor:UnitaryS12ActionOnSlice}, the embedding \eqref{eq:SliceBallEmbeddingInQuotient} is $S^1$-equivariant with respect to the action \eqref{eq:S1ActionOnSO(3)MonopolesDefComplexSpacesFirstTerm} on $\bB_{A,\Phi}$ and the $S^1$ action  \eqref{eq:S12ActionsOnSpinuQuotient} on $\sC_\ft$. Hence, the analytic embedding of $\bB_{A,\Phi}\cap\fS_{A,\Phi}^{-1}(0)$ onto an open neighborhood of $[A,\Phi]$ in $\sM_\ft$ is $S^1$-equivariant. Finally, we note that $\fS_{A,\Phi}$ is $S^1$-equivariant with respect to the $S^1$ action 
\eqref{eq:S1ActionOnSO(3)MonopolesDefComplexSpacesFirstTerm}
on its domain and the $S^1$ action \eqref{eq:S1ActionOnSO(3)MonopolesDefComplexSpacesSecondTerm} on its codomain.  By Corollary \ref{cor:S1_Action_on_HarmonicSections_Of_nonAbelianMonopole_Def_Complex}, these are the  $S^1$ actions on $\bH_{A,\Phi}^1$ and $\bH_{A,\Phi}^2$ described in the statement of that corollary.
Lemma \ref{lem:FredholmMapAndGEquivariantSplitting} implies that we may apply Lemma \ref{lem:EquivariantKuranishiLemma} to the map $\fS_{A,\Phi}$ which gives
the existence of the maps $\bkappa$ and $\beps$ with the properties in the statement of the lemma.
\end{proof}

\section{$\CC^*$ actions on the affine space of (0,1)-pairs}
\label{sec:CircleActionsOnAffineSpace}
As noted in the Introduction to this chapter, to compute the signs of weights of the $S^1$ actions on the harmonic spaces $\bH_{A,\Phi}^k$ in the Kuranishi model provided by Lemma \ref{lem:S1EquivariantKuranishiModelForReducibleIn_non_AbelianMonopoleModuli}, we need these actions to be \emph{complex} representations.  We will accomplish this by showing that there is an $S^1$-equivariant isomorphism between the harmonic spaces for the non-Abelian monopole elliptic deformation complex and
the harmonic spaces for the elliptic complex of holomorphic pairs and using the obvious almost complex structure on the latter spaces. In this section, we begin this process by identifying $\CC^*$ actions (and thus by restriction to $S^1\subset\CC^*$, we identify $S^1$ actions) on the affine space of $(0,1)$-pairs which play the same role in computing actions on the harmonic spaces of the elliptic complex of holomorphic pairs vas the action \eqref{eq:UnitaryS12ActionOnNonAbelianPairsAffineSpace} did in computing actions on the harmonic spaces of the elliptic deformation complex of non-Abelian monopoles.

Let $E$ be a rank-$r$, smooth complex vector bundle over a complex manifold and let $(\bar\rd_E,\varphi)$ be a holomorphic pair on $E$ which is split in the sense of Definition \ref{defn:Split_(0,1)-pair} with respect to a decomposition $E=E_1\oplus E_2$ as in \eqref{eq:E_equals_L_1_oplus_L_2} as a direct sum of proper complex vector subbundles, $E_1$ and $E_2$ of complex ranks $r_1$ and $r_2$, respectively, and $\varphi\in\Omega^0(E_1)$. Lemma \ref{lem:FixedPointsOfC*ActionOnQuotientSpace01Pairs} implies that the point $[\bar\rd_E,\varphi]\in \sC^{0,1}(E)$ is a fixed point of the standard $\CC^*$ action \eqref{eq:CActionsOn(0,1)PairQuotients} on $\sC^{0,1}(E)$.  If $\varphi\not\equiv 0$, then the pair $(\bar\rd_E,\varphi)$ is not a fixed point of the action \eqref{eq:CZActionOnAffine} on the affine space $\sA^{0,1}(E)\times \Omega^0(E)$. Thus, for the same reasons motivating our definition of the modified $S^1$ action \eqref{eq:UnitaryS12ActionOnNonAbelianPairsAffineSpace} on the affine space of \spinu pairs,
we now define a modified $\CC^*$ action for which $(\bar\rd_E,\varphi)$ will be a fixed point.

\subsection{Modified $\CC^*$ action on the affine space of (0,1)-pairs}
\label{sec:NewS1Action}
Assume that $E$ is split as in \eqref{eq:E_equals_L_1_oplus_L_2} and consider the three subgroups of $C^\infty(\GL(E))$ given by the images of the homomorphisms $\CC^*\to C^\infty(\GL(E))$ defined by
\begin{subequations}
\label{eq:DefineS1ActionsAtReducible}
\begin{align}
\label{eq:DefineS1ZActionsAtReducible}
\rho_Z^\CC(\la)&:=\la\,\id_{E_1}\oplus \la\,\id_{E_2},
  \\
\label{eq:DefineS1sActionsAtReducible}
 \rho_{\SL}(\la)&:=\la^{r_2}\,\id_{E_1}\oplus \la^{-r_1}\,\id_{E_2},
\\
\label{eq:DefineS12ActionsAtReducible}
\rho_2^\CC(\la)&:=\id_{E_1}\oplus \la\,\id_{E_2}, \quad\text{for all } \lambda \in \CC^*,
\end{align}
\end{subequations}
where $E_i$ has complex rank $r_i$ for $i=1,2$. The restrictions of the homomorphisms \eqref{eq:DefineS1ActionsAtReducible} to $S^1\subset\CC^*$ are equal to the homomorphisms defined previously in \eqref{eq:DefineUnitaryS1ActionsAtReducible} (when $E$ has rank two) and in \eqref{eq:DefineUnitaryS1ActionsAtReducible_rank-r} (when $E$ has arbitrary rank).

We consider the modified $\CC^*$ action on the affine space of $(0,1)$-pairs defined by
\begin{multline}
\label{eq:S12ActionOnAffine}
\rho_2^{\sA,\CC}:\CC^*\times \sA^{0,1}(E)\times \Omega^0(E)
\ni \left(\la,(\bar\rd_E,\varphi)\right)
\\
\mapsto
\left( \rho_2^\CC(\la)\circ\bar\rd_E\circ \rho_2^\CC(\la)^{-1}, \rho_2^\CC(\la)\varphi\right)
\in
\sA^{0,1}(E)\times \Omega^0(E).
\end{multline}
The action \eqref{eq:S12ActionOnAffine} is analogous to the action $\rho_2^{\sA}$ on spin${}^u$ pairs \eqref{eq:UnitaryS12ActionOnNonAbelianPairsAffineSpace}.

\begin{rmk}[On the definition of the modified $\CC^*$ action]
\label{rmk:Choice_Of_rhoAffine_Action}
The $\CC^*$ action $\rho_2^{\sA,\CC}$ is defined as in \eqref{eq:S12ActionOnAffine} to ensure that the map between type $1$ \spinu pairs and $(0,1)$-pairs, $(A,\Phi)=(A,(\varphi,0))\mapsto (\bar\rd_A,\varphi)$, is $S^1$-equivariant with respect to the $S^1$ action on \spinu pairs given by $\rho_2^\sA$ and the $S^1$ action on $(0,1)$-pairs defined by the composition of $\rho_2^{\sA,\CC}$ with the inclusion $S^1\subset\CC^*$. We use this $S^1$-equivariance in Lemma \ref{lem:S1_Equivariance_of_CohomologyIsoms} to prove that the isomorphisms between the harmonic spaces for the elliptic deformation complex for non-Abelian monopoles and the harmonic spaces for the elliptic complex for holomorphic pairs are $S^1$-equivariant, as required for the reasons discussed in the introduction to this chapter. The same comments apply to the map between unitary pairs and $(0,1)$-pairs, $(A,\varphi)\mapsto (\bar\rd_A,\varphi)$, and the isomorphisms described in the forthcoming Remark \ref{rmk:Generalizing_S1_Equivariance_of_CohomologyIsoms} between the harmonic spaces for the elliptic complex for projective vortices and the harmonic spaces for the elliptic complex for holomorphic pairs.
\end{rmk}

In Lemma \ref{lem:EqualityOfUnitaryS1ActionsOnQuotientSpace}, we showed that the actions induced by $\rho_2^{\sA}$ and $\rho_Z$ on the quotient $\sC_\ft$ were the same, up to positive multiplicity.  In the following, we note that the actions $\rho_2^{\sA,\CC}$ and $\rho_Z^\CC$ induce the same actions on $\sC^{0,1}(E)$ up to positive multiplicity.

\begin{lem}[Gauge-equivalence of modified and standard $\CC^*$ actions on the affine space of $(0,1)$-pairs up to positive multiplicity]
\label{lem:RelateS1Actions}
Let $E$ be a rank-$r$, smooth, complex vector bundle over a complex manifold. If $E$ is split as in \eqref{eq:E_equals_L_1_oplus_L_2}, so $E=E_1\oplus E_2$ as a direct sum of proper complex vector subbundles $E_i$ of rank $r_i$ for $i=1,2$, then the quotient map,
\begin{equation}
\label{eq:QuotientMapOf(0,1)Pairs}
\pi_\sC:\sA^{0,1}(E)\times\Om^0(E) \ni (\bar\rd_E,\varphi) \mapsto [\bar\rd_E,\varphi] \in \sC^{0,1}(E),
\end{equation}
is $\CC^*$-equivariant with respect to the $\CC^*$ action $\rho_2^{\sA,\CC}$ with multiplicity $r$ on the affine space and the standard action \eqref{eq:CActionsOn(0,1)PairQuotients}
with multiplicity $r_2$ on the quotient space $\sC^{0,1}(E)$, that is,
\begin{equation}
\label{eq:C2andCZActionsEqualOnQuotient}
\pi_\sC\left( \rho_2^{\sA,\CC}(\la^r,(\bar\rd_E,\varphi)\right)=[\bar\rd_E,\la^{r_2}\varphi],
\quad\text{for all } \lambda \in \CC^*.
\end{equation}
\end{lem}

\begin{proof}
The three homomorphisms \eqref{eq:DefineS1ActionsAtReducible} are related by the identity
\begin{equation}
\label{eq:ComparingS1ActionsAtReducible}
\rho_2^\CC(\la^r)=\rho_{\SL}(\la)^{-1}\rho_Z^\CC(\la^{r_2}), \quad\text{for all } \la\in\CC^*.
\end{equation}
Because the image of $\rho_{\SL}$ lies in $C^\infty(\SL(E))$, Equation \eqref{eq:ComparingS1ActionsAtReducible} implies that for any pair $(\bar\rd_E,\varphi)$ we have
\begin{align*}
&\pi_\sC\left( \rho_2^{\sA,\CC}\left(\la^r,(\bar\rd_E,\varphi)\right)\right)
  \\
  &\quad=
[ \rho_2^{\sA,\CC}\left(\la^r,(\bar\rd_E,\varphi)\right)]
\quad\text{(by definition of $\pi$ in \eqref{eq:QuotientMapOf(0,1)Pairs})}
\\
&\quad=
[\rho_2^\CC(\la^r)\circ\bar\rd_E\circ\rho_2^\CC(\la)^{-r},\rho_2^\CC(\la^r)\varphi]
\quad\text{(by the definition of $\rho_2^{\sA,\CC}$ in \eqref{eq:S12ActionOnAffine})}
\\
&\quad=
[\rho_{\SL}(\la)^{-1}\rho_Z^\CC(\la^{r_2})\circ \bar\rd_E\circ\rho_Z^\CC(\la^{r_2})^{-1}\rho_{\SL}(\la),
\rho_{\SL}(\la)^{-1}\rho_Z^\CC(\la^{r_2})\varphi]
     \quad\text{(by \eqref{eq:ComparingS1ActionsAtReducible})}
  \\
  &\quad=
[\rho_{\SL}(\la)^{-1}\circ \bar\rd_E\circ\rho_{\SL}(\la), \rho_{\SL}(\la)^{-1}\rho_Z^\CC(\la^{r_2})\varphi]
\\
&\qquad\text{(because $\rho_Z^\CC(\la)\circ\bar\rd_E\circ\rho_Z^\CC(\la)^{-1}=\bar\rd_E$)}
  \\
&\quad=
 [\rho_{\SL}(\la)^* ( \bar\rd_E,\rho_Z^\CC(\la^{r_2})\varphi)]
\quad\text{(by definition of action in \eqref{eq:SL(E)ActionOn(0,1)Pairs})}
 \\
&\quad=
[ \bar\rd_E,\rho_Z^\CC(\la^{r_2})\varphi]
\quad\text{(because $\rho_{\SL}(\lambda)\in C^\infty(\SL(E))$)}
\\
&\quad =
[\bar\rd_E,\la^{r_2}\varphi]\quad\text{(by \eqref{eq:DefineS1ZActionsAtReducible})},
\end{align*}
which yields Equation \eqref{eq:C2andCZActionsEqualOnQuotient}.
\end{proof}

\subsection{$\CC^*$ actions on tangent spaces to fixed points in the  affine space of (0,1)-pairs}
As in \eqref{eq:FixedPointS1Action}, the action \eqref{eq:S12ActionOnAffine} of  $\rho_2^{\sA,\CC}$ on the affine space $\sA^{0,1}(E)\times \Omega^0(E)$ defines an action on the tangent space of $\sA^{0,1}(E)\times \Omega^0(E)$ at fixed points of the action.  Because $\sA^{0,1}(E)\times \Omega^0(E)$ is an affine space, the tangent space
at a pair $(\bar\rd_E,\varphi)$ is identified with the underlying vector space $\sE_1$ appearing in the deformation complex \eqref{eq:Holomorphic_pair_elliptic_complex} (see also \eqref{eq:AbbreviatedHolomorphicPairDefComplex} and\eqref{eq:DefineSpacesOfAbbreviatedHolomorphicPairDefComplex}),
\begin{equation}
\label{eq:AffineTangent}
T_{\bar\partial_E,\varphi} \left( \sA^{0,1}(E)\times \Omega^0(E)\right)
=
\Omega^{0,1}(\fsl(E))\oplus \Omega^0(E).
\end{equation}
We now give an explicit expression for the $\CC^*$ action induced by $\rho_2^{\sA,\CC}$ on the tangent space of a fixed point and give an expression for $\rho_2^{\sA,\CC}$  in terms of this derivative. The analogous result for the action $\rho_2^\sA$ was given in Lemma \ref{lem:S12ActionOnNonAbelianPairsAffineSpace}.

\begin{lem}[Induced $\CC^*$ action on the tangent space to the affine space of $0,1)$ pairs at a fixed point]
\label{lem:ActionOnTangentofAffineSpace}
Let $E$ be a rank-$r$, smooth, complex vector bundle over a complex manifold. If $(\bar\partial_E,\varphi)$ is a split pair on $E$ with respect to a decomposition $E=E_1\oplus E_2$ in the sense of Definition \ref{defn:Split_(0,1)-pair}, so $\varphi\in \Omega^0(E_1)$, then $(\bar\partial_E,\varphi)$ is a fixed point of the $\CC^*$ action $\rho_2^{\sA,\CC}$  on $\sA^{0,1}(E)\times \Omega^0(E)$ given in \eqref{eq:S12ActionOnAffine}. The action $D_2\rho_2^{\sA,\CC}$  induced by the derivative of $\rho_2^{\sA,\CC}$ on the tangent space at $(\bar\rd_E,\varphi)$ \eqref{eq:AffineTangent} (in the sense of \eqref{eq:Circle_action_tangent_bundle}) is
\begin{multline}
\label{eq:S1L2ActionOnTangent}
(D_2\rho_2^{\sA,\CC})\left( \la,(\alpha,\sigma)\right)
=
\left( \rho_2^\CC(\la) \alpha\rho_2^\CC(\la)^{-1},\rho_2^\CC(\la)\si \right),
\\
\text{for all } \lambda \in \CC^* \text{ and } (\alpha,\sigma) \in \Omega^{0,1}(\fsl(E))\oplus \Omega^0(E).
\end{multline}
\end{lem}

\begin{proof}
The split pair $(\bar\rd_E,\varphi)=(\bar\rd_{E_1}\oplus\bar\rd_{E_2},\varphi)$ is a fixed point for the action \eqref{eq:S12ActionOnAffine} because the constant gauge transformations $\rho_2^\CC(\la)$ satisfy
\begin{equation}
\label{eq:StabilizerOfSplit01Connection}
\rho_2^\CC(\la)\circ (\bar\rd_{E_1}\oplus\bar\rd_{E_2})\circ\rho_2^\CC(\la)^{-1}
=
\bar\rd_{E_1}\oplus\bar\rd_{E_2},
\quad\text{for all } \lambda \in \CC^*,
\end{equation}
and because the component of $\varphi$ in $E_2$ is identically zero, so
\begin{equation}
\label{eq:SL(E)StabilizerOfSection}
\rho_2^\CC(\la)\varphi=\varphi,\quad\text{for all } \lambda \in \CC^*.
\end{equation}
 Using the equality,
\begin{equation}
\label{eq:(0,1)PairsUnderlyingVectorSpace}
\sA^{0,1}(E)\times \Omega^0(E)=(\bar\rd_E,\varphi)+\left( \Om^{0,1}(\fsl(E))\times \Omega^0(E)\right),
\end{equation}
we claim that the $\CC^*$ action \eqref{eq:S12ActionOnAffine} satisfies
\begin{multline}
\label{eq:CC*ActionOnAffine}
\rho_2^{\sA,\CC}\left(\la,\left( (\bar\rd_E,\varphi)+(\alpha,\sigma)\right)\right)
=
(\bar\rd_E,\varphi)
+
\left( \rho_2^\CC(\la) \alpha\rho_2^\CC(\la)^{-1},\rho_2^\CC(\la)\si\right)
,
\\
\quad\text{for all } \lambda \in \CC^* \text{ and } (\alpha,\sigma) \in \Omega^{0,1}(\fsl(E))\oplus \Omega^0(E).
\end{multline}
We prove that \eqref{eq:CC*ActionOnAffine} holds by using \eqref{eq:(0,1)PairsUnderlyingVectorSpace} to write $(\bar\rd_{E'},\varphi') = (\bar\rd_E,\varphi)+(\alpha,\sigma)$ and observing that
\begin{align*}
&\rho_2^{\sA,\CC}\left( \la, (\bar\rd_{E'},\varphi')\right)
  \\
  &\quad=
\left( \rho_2^\CC(\la)\circ\bar\partial_{E'}\circ \rho_2^\CC(\la)^{-1},\rho_2^\CC(\la)\varphi'\right)
\quad\text{(by definition of $\rho_2^{\sA,\CC}$ in \eqref{eq:S12ActionOnAffine})}
\\
&\quad=
\left( \rho_2^\CC(\la)\circ\left(\bar\rd_E+\alpha \right)\circ \rho_2^\CC(\la)^{-1},\rho_2^\CC(\la)(\varphi+\sigma)\right)
\\
&\quad=
\left( \rho_2^\CC(\la)\circ \bar\rd_E\circ \rho_2^\CC(\la)^{-1},\rho_2^\CC(\la)\varphi\right)
+
\left( \rho_2^\CC(\la) \alpha\rho_2^\CC(\la)^{-1},\rho_2^\CC(\la)\si\right)
  \\
&\quad=
(\bar\rd_E,\varphi)
+
\left( \rho_2^\CC(\la) \alpha\rho_2^\CC(\la)^{-1},\rho_2^\CC(\la)\si\right)
\quad\text{(by \eqref{eq:StabilizerOfSplit01Connection} and \eqref{eq:SL(E)StabilizerOfSection}).}
\end{align*}
The expression for the derivative \eqref{eq:S12ActionOnAffine} follows immediately from
Equation \eqref{eq:CC*ActionOnAffine}.
\end{proof}

\section{$\CC^*$-equivariance of the elliptic complex for holomorphic pairs}
\label{subsec:S1EquivComplex}
In this section, we compute the $\CC^*$ action on the elliptic complex for holomorphic pairs which will appear in the desired $S^1$-equivariant isomorphism between the harmonic spaces for the elliptic complex for non-Abelian monopoles and the harmonic spaces for the elliptic complex for holomorphic pairs.

We continue to assume that $E$ is a rank-$r$, smooth, complex vector bundle over a complex manifold, of complex dimension $n$. Let $(\bar\partial_E,\varphi)$ be a  holomorphic pair on $E$ which is split with respect to a decomposition $E=E_1\oplus E_2$ in the sense of Definition \ref{defn:Split_(0,1)-pair}. In particular, we assume  that $\varphi$ is a section of $E_1$. We now describe a $\CC^*$-equivariant structure (in the sense of Definition \ref{defn:GEquivariantStructure}) on the elliptic complex \eqref{eq:Holomorphic_pair_elliptic_complex} for a holomorphic pair at $(\bar\partial_E,\varphi)$.
We will write this elliptic complex as
\begin{equation}
\label{eq:AbbreviatedHolomorphicPairDefComplex}
\begin{CD}
0
@>>>
\sE_0
@>\bar\rd_{E,\varphi}^0 >>
\sE_1
@>\bar\rd_{E,\varphi}^1 >>
\sE_2
@>\bar\rd_{E,\varphi}^2 >>
\dots
@>\bar\rd_{E,\varphi}^n >>
\sE_{n+1}
@>>>
0
\end{CD}
\end{equation}
where
\begin{equation}
\label{eq:DefineSpacesOfAbbreviatedHolomorphicPairDefComplex}
\sE_k :=\Om^{0,k}(\fsl(E))\oplus \Om^{0,k-1}(E),\quad\text{for } k=0,\ldots,n+1,
\end{equation}
and the operators $\bar\rd_{E,\varphi}^k$ are as in \eqref{eq:dStablePair}, with the convention that $\Om^{0,k}(\fsl(E)) = (0)$ and $\Om^{0,k}(E) = (0)$ when $k<0$ or $k>n$. We define $\CC^*$ actions 
\begin{equation}
\label{eq:CActionOnHolomorphicPairsDeformationComplex}
\rho_{2,k}:\CC^*\times\sE_k\to\sE_k
\end{equation}
by setting
\begin{multline}
\label{eq:CActionOnHolomorphicPairsDeformationComplexCases}
\rho_{2,k}\left(\la, \left( w,\nu\right)\right)
:= \left( \rho_2^\CC(\la) w \rho_2^\CC(\la)^{-1},\rho_2^\CC(\la)\nu\right),
\\
\text{for all } w\in\Om^{0,k}(\fsl(E)) \text{ and } \nu\in\Om^{0,k-1}(E)
 \text{ and } \lambda \in \CC \text{ and } k \in \ZZ,
\end{multline}
where $\rho_2^\CC(\lambda)$ is as in \eqref{eq:DefineS12ActionsAtReducible}. We note that the action $\rho_{2,1}$ is equal to the action $D_2\rho_2^{\sA,\CC}$ appearing in \eqref{eq:S1L2ActionOnTangent}. To make expressions such as the forthcoming \eqref{eq:CEquivarianceOfDefComplexOperators} easier to read, we will abuse notation and write $\rho_{2,k}(\la)$ for the $\CC$-linear map defined by $\rho_{2,k}(\la,\cdot):\sE_k\to\sE_k$. Just as  the $S^1$-equivariance of the elliptic deformation complex for non-Abelian monopoles was established in Lemma \ref{lem:S1_Equivariance_Of_nonAbelianMonopole_Def_Complex}, the $\CC^*$-equivariance in the sense of Definition \ref{defn:GEquivariantStructure} of the complex \eqref{eq:AbbreviatedHolomorphicPairDefComplex} with respect to the actions \eqref{eq:CActionOnHolomorphicPairsDeformationComplex} is established in the following

\begin{lem}[$\CC^*$-equivariant structure on the elliptic complex for holomorphic pairs]
\label{lem:S1EquivariantDeformationComplex}
Let $E$ be a smooth, rank-$r$, complex vector bundle over a complex manifold, of complex dimension $n$, and $(\bar\partial_E,\varphi)$ be a holomorphic pair on  $E$ which is split with respect to a decomposition $E=E_1\oplus E_2$ in the sense of Definition \ref{defn:Split_(0,1)-pair}, so that $\varphi\in\Omega^0(E_1)$.  Then the differentials in the elliptic complex \eqref{eq:AbbreviatedHolomorphicPairDefComplex} are equivariant with respect to the $\CC^*$ actions \eqref{eq:CActionOnHolomorphicPairsDeformationComplex} in the sense that
\begin{equation}
\label{eq:CEquivarianceOfDefComplexOperators}
\bar\partial_{E,\varphi}^k\circ \rho_{2,k}(\la)
=
\rho_{2,k+1}(\la)\circ \bar\partial_{E,\varphi}^k,\quad\text{for $k=0,\dots,n$}.
\end{equation}
The elliptic complex \eqref{eq:AbbreviatedHolomorphicPairDefComplex} then has a $\CC^*$-equivariant
structure in the sense of Definition \ref{defn:GEquivariantStructure} with respect to the actions \eqref{eq:CActionOnHolomorphicPairsDeformationComplex}. In addition, if the base manifold is closed then the $L^2$ adjoints of the differentials in the elliptic complex \eqref{eq:AbbreviatedHolomorphicPairDefComplex} are also equivariant with respect to these $\CC^*$ actions:
\begin{equation}
\label{eq:CEquivarianceOfDefComplexAdjOperatorsd}
\bar\partial_{E,\varphi}^{k,*}\circ \rho_{2,k+1}(\la)
=
\rho_{2,k}(\la)\circ \bar\partial_{E,\varphi}^{k,*},\quad\text{for $k=0,\dots,n$}.
\end{equation}
\end{lem}

\begin{proof}
Let $k\in\ZZ$ and $\lambda\in\CC^*$. Because $\rho_2^\CC(\la)$ acts trivially on the split pair $(\bar\rd_E,\varphi)$, we have
\begin{subequations}
\label{eq:Rho2StabilizesPair}
\begin{align}
\label{eq:Rho2CommutesWithPartialA}
\bar\rd_E\circ \rho_2^\CC(\la)&=\rho_2^\CC(\la)\circ\bar\rd_E,
\\
\label{eq:Rho2StabilizesSection}
\rho_2^\CC(\la)\varphi&=\varphi.
\end{align}
\end{subequations}
Equation \eqref{eq:CEquivarianceOfDefComplexOperators} then follows from  the preceding equality and the definition of $\bar\rd_{E,\varphi}^k$ in \eqref{eq:dStablePair}:
\begin{align*}
\bar\rd_{E,\varphi}^k\rho_{2,k}(\la) \begin{pmatrix} w \\ \nu \end{pmatrix}
&=
\bar\rd_{E,\varphi}^k \begin{pmatrix} \rho_2^\CC(\la)w \rho_2^\CC(\la)^{-1} \\ \rho_2^\CC(\la)\nu \end{pmatrix}
\quad\text{(by \eqref{eq:CActionOnHolomorphicPairsDeformationComplexCases})}
\\
&=
\begin{pmatrix}
  \bar\partial_E \left(\rho_2^\CC(\la)w \rho_2^\CC(\la)^{-1}\right)
  \\
  \bar\partial_E\left(\rho_2^\CC(\la)\nu\right) + (-1)^{k-1}\rho_2^\CC(\la)w \rho_2^\CC(\la)^{-1}\varphi
\end{pmatrix}
\quad\text{(by \eqref{eq:dkStablePair})}
\\
&=
\begin{pmatrix}
   \rho_2^\CC(\la)\left(\bar\partial_Ew\right) \rho_2^\CC(\la)^{-1}
  \\
  \rho_2^\CC(\la)\bar\partial_E\nu + (-1)^{k-1}\rho_2^\CC(\la)w\varphi
\end{pmatrix}
  \quad\text{(by \eqref{eq:Rho2CommutesWithPartialA} and \eqref{eq:Rho2StabilizesSection})}
  \\
&= \rho_{2,k+1}(\la)\begin{pmatrix}
  \bar\partial_Ew
  \\
  \bar\partial_E\nu + (-1)^{k-1}w\varphi
\end{pmatrix}
\quad\text{(by \eqref{eq:CActionOnHolomorphicPairsDeformationComplexCases})}  
\\
&=
\rho_{2,k+1}(\la) \bar\rd_{E,\varphi}^k \begin{pmatrix} w \\ \nu \end{pmatrix}
\quad\text{(by \eqref{eq:dkStablePair}),}
\end{align*}
as asserted in \eqref{eq:CEquivarianceOfDefComplexOperators}.

Because the $\CC^*$ action \eqref{eq:CActionOnHolomorphicPairsDeformationComplexCases} on $\sE_k$ is given by scalar multiplication on the complex line bundle $L_2$, the $L^2$ adjoint of $\rho_{2,k}(\la) \in \End(\sE_k)$ is given by $\rho_{2,k}(\bar\la)\in \End(\sE_k)$, we obtain
\begin{equation}
\label{eq:AdjointOfRho2}
\left( \rho_{2,k}(\la) \begin{pmatrix} w \\ \nu \end{pmatrix},
\begin{pmatrix} w' \\ \nu'\end{pmatrix} \right)_{L^2(X)}
=
\left(  \begin{pmatrix} w \\ \nu \end{pmatrix},
\rho_{2,k}(\bar\la)\begin{pmatrix} w' \\ \nu'\end{pmatrix} \right)_{L^2(X)}.
\end{equation}
When $X$ is closed, the identity \eqref{eq:AdjointOfRho2} and the $\CC^*$-equivariance \eqref{eq:CEquivarianceOfDefComplexOperators} of the operators $\bar\rd^k_{E,\varphi}$ yield
\begin{align*}
{}&\left( \bar\rd_{E,\varphi}^{k,*} \left(\rho_{2,k+1}(\la) \begin{pmatrix} w \\ \nu \end{pmatrix}\right), \begin{pmatrix} w' \\ \nu'\end{pmatrix} \right)_{L^2(X)}
\\
{}&\quad=
\left(   \begin{pmatrix} w \\ \nu \end{pmatrix},
\rho_{2,k+1}(\bar\la)\bar\rd_{E,\varphi}^k\begin{pmatrix} w' \\ \nu'\end{pmatrix} \right)_{L^2(X)}
\quad\text{(by definition of $L^2$ adjoint and \eqref{eq:AdjointOfRho2})}
\\
{}&\quad=
\left(   \begin{pmatrix} w \\ \nu \end{pmatrix},
\bar\rd_{E,\varphi}^k\rho_{2,k}(\bar\la)\begin{pmatrix} w' \\ \nu'\end{pmatrix} \right)_{L^2(X)}
\quad\text{(by \eqref{eq:CEquivarianceOfDefComplexOperators})}
\\
{}&\quad=
\left( \rho_{2,k}(\la) \bar\rd_{E,\varphi}^{k,*} \begin{pmatrix} w \\ \nu \end{pmatrix},
\begin{pmatrix} w' \\ \nu'\end{pmatrix} \right)_{L^2(X)}
\quad\text{(by definition of $L^2$ adjoint and \eqref{eq:AdjointOfRho2}).}
\end{align*}
The preceding equalities give \eqref{eq:CEquivarianceOfDefComplexAdjOperatorsd}.
\end{proof}

Lemma \ref{lem:S1EquivariantDeformationComplex} implies that the $\CC^*$ action on $\sE_k$ defines an action on the harmonic subspace $\bH_{\bar\rd_E,\varphi}^k \subset \sE_k$, defined in \eqref{eq:H_dbar_Avarphi^0bullet} when the base complex manifold is closed. The following is analogous to Corollary \ref{cor:S1_Action_on_HarmonicSections_Of_nonAbelianMonopole_Def_Complex}.

\begin{cor}[$\CC^*$ actions on harmonic spaces for the elliptic complex of a split holomorphic pair]
\label{cor:CActionOnHarmonicSpacesOfHolomPairsComplex}
Continue the hypotheses and notation of Lemma \ref{lem:S1EquivariantDeformationComplex} and assume, in addition, that the base complex manifold is closed. Then the $\CC^*$ action \eqref{eq:CActionOnHolomorphicPairsDeformationComplexCases} preserves the harmonic subspaces $\bH_{\bar\rd_E,\varphi}^k \subset \sE_k$ defined by \eqref{eq:H_dbar_Avarphi^0bullet}, for $k=0,\dots,n+1$.
\end{cor}

\section{Circle equivariance of the isomorphisms between harmonic spaces}
\label{sec:EquivarianceOfIsomorphisms}
We now show that the isomorphisms between the harmonic spaces for the elliptic deformation complex for non-Abelian monopoles and those for the elliptic complex for holomorphic pairs are equivariant with respect to the $S^1$ actions on these spaces defined in the preceding sections.

If $(A,\Phi)$ is a type $1$ non-Abelian monopole with $\Phi\not\equiv 0$, then $\bH_{A,\Phi}^0=(0)$ by Lemma \ref{lem:H0_Of_NonAbelianMonopoleComplex_Vanishes}. 
Because $(A,\Phi)$ is type 1, we can write $(A,\Phi)=(A,\varphi)$ where $\varphi\in\Om^0(E)$
and $(A,\varphi)$ is a projective vortex.
Then equations \eqref{eq:Itoh_1985_proposition_2-4_SO3_monopole_complex_Kaehler_isomorphism_simplified} in Remark \ref{rmk:Difference_widehatH_dbar_APhi^01_and_H_dbar_APhi^01} and
\eqref{eq:H2_NonAbelianMonopole_HolomorphicPair} in Remark \ref{rmk:Relation_H_APhi^2_and_H_dbar_APhi^2} yield isomorphisms of real vector spaces,
\begin{equation}
\label{eq:Isomorphisms_of_non_AbelianHarmonics_with_HolomorphPairHarmonics}
\bH_{A,\Phi}^k \cong \bH_{\bar\partial_A,\varphi}^k\quad\text{for } k=0,1,2,
\end{equation}
where $\bH_{A,\Phi}^k$ in \eqref{eq:H_APhi^bullet} is a harmonic space for the elliptic deformation complex for a non-Abelian monopole defined and $\bH_{\bar\partial_A,\varphi}^k$ in \eqref{eq:H_dbar_Avarphi^0bullet} denotes a harmonic space for the elliptic complex for holomorphic pairs.

We will show that these isomorphisms are $S^1$-equivariant with respect to the action defined in this section on the non-Abelian monopole deformation complex \eqref{eq:SO3MonopoleDefComplex} and the action \eqref{eq:CActionOnHolomorphicPairsDeformationComplex} on the elliptic complex \eqref{eq:H_dbar_Avarphi^0bullet} for holomorphic pairs.  We note that the $\CC^*$ action defined on
$\bH_{\bar\rd_A,\varphi}^k$ by Corollary \ref{cor:CActionOnHarmonicSpacesOfHolomPairsComplex}
defines  an $S^1$ action by restriction to the subgroup $S^1$ of $\CC^*$.

\begin{lem}[Equivariance of the isomorphisms between the harmonic spaces for elliptic complexes for non-Abelian monopoles and holomorphic pairs]
\label{lem:S1_Equivariance_of_CohomologyIsoms}
Continue the hypotheses of Lemma \ref{lem:S1_Equivariance_Of_nonAbelianMonopole_Def_Complex} and assume in addition that $X$ is a closed, K\"ahler surface, let $(A,\Phi)$ be a type $1$ non-Abelian monopole, and let $(\bar\rd_A,\varphi)$ be the holomorphic pair defined by $(A,\Phi)$ as in \eqref{eq:del_bar_A} and \eqref{eq:Phi_coupled_spinor_is_pair_sections_E_oplus_02E}. If $\Phi\not\equiv 0$, then
$\bH_{A,\Phi}^0=(0)$ and the isomorphisms
\[
\bH_{A,\Phi}^k \cong \bH_{\bar\rd_A,\varphi}^k, \quad\text{for } k=1,2,
\]
defined by \eqref{eq:Isomorphisms_of_non_AbelianHarmonics_with_HolomorphPairHarmonics} for $k=1,2$ are $S^1$-equivariant with respect to the $S^1$ actions on $\bH_{A,\Phi}^k$ given by Corollary \ref{cor:S1_Action_on_HarmonicSections_Of_nonAbelianMonopole_Def_Complex} and on
$\bH_{\bar\rd_A,\varphi}^k$ given by \eqref{eq:CActionOnHolomorphicPairsDeformationComplex} as described in Corollary \ref{cor:CActionOnHarmonicSpacesOfHolomPairsComplex}.
\end{lem}

\begin{rmk}[Generalization of Lemma \ref{lem:S1_Equivariance_of_CohomologyIsoms} to projective vortices on Hermitian vector bundles of higher rank over K\"ahler manifolds of arbitrary dimension]
\label{rmk:Generalizing_S1_Equivariance_of_CohomologyIsoms}
Lemma \ref{lem:S1_Equivariance_of_CohomologyIsoms} easily extends to the case where the type $1$ non-Abelian monopole $(A,\Phi)$ is replaced by a projective vortex $(A,\varphi)$ on a Hermitian vector bundle of arbitrary rank over a K\"ahler manifold of arbitrary dimension $n$. If the projective vortex is split in the sense of Definition \ref{defn:Split_trivial_central-stabilizer_unitary_pair} \eqref{item:Split_unitary_pair} with respect to an orthogonal decomposition $E=E_1\oplus E_2$ of proper Hermitian vector subbundles, $E_1$ and $E_2$, then the isomorphisms of harmonic spaces provided by Theorem \ref{thm:Kobayashi_7-2-21_pairs},
\[
  \bH_{A,\varphi}^k \cong \bH_{\bar\rd_A,\varphi}^k, \quad\text{for } 0,1,\dots,n,
\]
are $S^1$-equivariant with respect to the $S^1$ action on $\bH_{A,\varphi}^k$ 
defined by analogy with the action \eqref{eq:S1ActionOnSO(3)MonopolesDefComplexSpaces},
as discussed in Remark \ref{rmk:GeneralizationOf_S1_Equivariance_Of_nonAbelianMonopole_Def_Complex},
and the $S^1$ action on $\bH_{\bar\partial_A,\varphi}^k$ given by 
\eqref{eq:CActionOnHolomorphicPairsDeformationComplex} as described in Corollary \ref{cor:CActionOnHarmonicSpacesOfHolomPairsComplex}.
\end{rmk}

\begin{proof}[Proof of Lemma \ref{lem:S1_Equivariance_of_CohomologyIsoms}]
For $k=1$, the isomorphism \eqref{eq:Isomorphisms_of_non_AbelianHarmonics_with_HolomorphPairHarmonics} is given by
the map \eqref{eq:Real_linear_map_tangent_spaces_SO(3)_pairs_to_preholomorphic_pairs}
and the identification of $\bH_{\bar\partial_A,\varphi}^1$ with the harmonic space of the pre-holomorphic pair $\bH_{\bar\partial_A,(\varphi,0)}^1$ in
\eqref{eq:Cohomologies_complexes_pre-holomorphic_pair_type1_and_holomorphic_pair_n_is_2_simplified}.
The  inverse of the isomorphism \eqref{eq:Itoh_1985_proposition_2-4_SO3_monopole_complex_Kaehler_isomorphism} is induced by the real linear map
\begin{multline}
  \label{eq:H1IsomorphismMap}
  \sE_1=\Om^{0,1}(\fsl(E))\oplus\Om^0(E) \ni
  (\alpha,\sigma)
  \mapsto
  I_1(\alpha,\sigma)
:=
\left(\frac{1}{2}(\alpha-\alpha^\dagger),(\si,0)\right)
\\
\in \Om^1(\su(E))\oplus\Om^0(E)\oplus\Om^{0,2}(E) = \Omega^1(\su(E))\oplus\Omega^0(W^+\otimes E),
\end{multline}
where the factor of $1/2$ appears because of our convention in \eqref{eq:Decompose_a_in_Omega1suE_into_10_and_01_components} that $a=(1/2)(a'+a'')$. The restriction of $I_1$ to $\bH_{\bar\rd_A,\varphi}$ gives the isomorphism
\eqref{eq:Isomorphisms_of_non_AbelianHarmonics_with_HolomorphPairHarmonics} for $k=1$.

The $\CC^*$ actions on $\bH_{\bar\rd_A,\varphi}^k$ given in Corollary \ref{cor:CActionOnHarmonicSpacesOfHolomPairsComplex} are defined by the actions $\rho_{2,k}$ in \eqref{eq:CActionOnHolomorphicPairsDeformationComplexCases} on the complex vector spaces $\sE_k$ in \eqref{eq:DefineSpacesOfAbbreviatedHolomorphicPairDefComplex}).  The 
restriction of homomorphism $\rho_2^\CC$ (defined in \eqref{eq:DefineS12ActionsAtReducible}) to $S^1\subset\CC^*$ appearing in the definition of $\rho_{2,k}$ is equal to the homomorphism $\rho_2$ in \eqref{eq:DefineUnitaryS12ActionsAtReducible},
\begin{equation}
\label{eq:RestrictionOfS12ActionToS1}
\rho_2^\CC(e^{i\theta})=\rho_2(e^{i\theta}), \quad\text{for all } e^{i\theta}\in S^1.
\end{equation}
Thus, for $e^{i\theta}\in S^1$ and $(\alpha,\sigma)\in\sE_1$, we obtain
\begin{align*}
I_1\left(\rho_{2,1}(e^{i\theta})(\alpha,\sigma) \right)
&=
I_1\left( \left( \rho_2(e^{i\theta})\alpha\rho_2(e^{i\theta})^{-1},\rho_2(e^{i\theta})\si\right)\right)
\quad\text{(by definitions \eqref{eq:CActionOnHolomorphicPairsDeformationComplexCases} and \eqref{eq:RestrictionOfS12ActionToS1})}
\\
&=
\left(
    \frac{1}{2}
    \left(
    \rho_2(e^{i\theta})\alpha\rho_2(e^{i\theta})^{-1} - \left( \rho_2(e^{i\theta})\alpha\rho_2(e^{i\theta})^{-1}\right)^\dagger
    \right),
    (\rho_2(e^{i\theta})\si,0)
     \right)
  \\
  &\qquad\text{(by definition \eqref{eq:H1IsomorphismMap})}
\\
&=
\left(
     \rho_2(e^{i\theta})\frac{1}{2}\left(\alpha-\alpha^\dagger\right)\rho_2(e^{i\theta})^{-1},
     \rho_2(e^{i\theta})(\si,0)
     \right)
     \\
&\qquad\text{(because $\rho_2(e^{i\theta})^\dagger=\rho_2(e^{i\theta})^{-1}$, for all $e^{i\theta}\in S^1$).}
\end{align*}
This last expression above coincides with the expression \eqref{eq:S1ActionOnSO(3)MonopolesDefComplexSpacesFirstTerm} defining the $S^1$ action on $I_1(\alpha,\sigma)$. Hence, for $k=1$, the isomorphism
\eqref{eq:Isomorphisms_of_non_AbelianHarmonics_with_HolomorphPairHarmonics} has the claimed $S^1$-equivariance.

We now prove the statement for $k=2$. The inverse of the isomorphism
 \eqref{eq:Itoh_1985_proposition_2-3_SO3_monopole_complex_Kaehler_isomorphism_type1} is induced by the real linear map
\begin{multline}
\label{eq:H2IsomorphismMap}
\sE_2=\Om^{0,2}(\fsl(E))\oplus\Om^{0,1}(E) \ni (w,\nu)
\mapsto
I_2(w,\nu) := \left(\frac{1}{2}(w-w^\dagger),\nu)\right)
\\
\in
\Omega^+(\su(E))\oplus\Omega^{0,1}(E) = \Omega^+(\su(E))\oplus\Omega^0(W^-\otimes E),
\end{multline}
where the factor of $1/2$ is due to our convention from equation \eqref{eq:Donaldson_Kronheimer_lemma_2-1-57_sum_3_terms}. As in the case for $k=1$, \eqref{eq:Cohomologies_complexes_pre-holomorphic_pair_type1_and_holomorphic_pair_n_is_2_simplified} implies that
\[
\bH_{\bar\rd_A,(\varphi,0)}^2
=
\bH_{\bar\rd_A,\varphi}^2.
\]
Hence, the restriction of $I_2$ to $\bH_{\bar\rd_A,\varphi}^2$ gives the isomorphism
\eqref{eq:Isomorphisms_of_non_AbelianHarmonics_with_HolomorphPairHarmonics} for $k=2$. For $e^{i\theta}\in S^1$ and $(w,\nu)\in \sE_2$, we compute that
\begin{align*}
I_2\left(\rho_{2,2}(e^{i\theta})(w,\nu)\right)
 &=
I_2\left( \rho_2(e^{i\theta})w\rho_2(e^{i\theta})^{-1},\rho_2(e^{i\theta})\nu\right)
\quad\text{(by definition \eqref{eq:CActionOnHolomorphicPairsDeformationComplexCases} and \eqref{eq:RestrictionOfS12ActionToS1})}
\\
 &=
\left( \frac{1}{2}\left( \rho_2(e^{i\theta})w \rho_2(e^{i\theta})^{-1} -
    \left(\rho_2(e^{i\theta}) w \rho_2(e^{i\theta})^{-1}\right)^\dagger\right),
    \rho_2(e^{i\theta})\nu
   \right)
  \\
  &\qquad\text{(by definition \eqref{eq:H2IsomorphismMap})}
\\
&=
\left( \rho_2(e^{i\theta})\frac{1}{2}\left(w - w^\dagger\right) \rho_2(e^{i\theta})^{-1},
    \rho_2(e^{i\theta})\nu
     \right)
     \\
&\qquad\text{(because $\rho_2(e^{i\theta})^\dagger=\rho_2(e^{i\theta})^{-1}$ for $e^{i\theta}\in S^1$).}
\end{align*}
The final expression above coincides with the expression \eqref{eq:S1ActionOnSO(3)MonopolesDefComplexSpacesSecondTerm} defining the $S^1$ action on $I_2(w,\nu)$. Hence, for $k=2$,  the isomorphism
\eqref{eq:Isomorphisms_of_non_AbelianHarmonics_with_HolomorphPairHarmonics}
is $S^1$-equivariant and this completes the proof of the lemma.
\end{proof}

\section[Weight decomposition of elliptic complex into weight subcomplexes]{Weight decomposition of the elliptic  complex for a split holomorphic pair into weight subcomplexes}
\label{sec:NormalTangentialSplittingOfDefComplex}
We now assume that $E$ has complex rank two and describe the decomposition of the elliptic complex \eqref{eq:AbbreviatedHolomorphicPairDefComplex} for a split holomorphic pair into a direct sum of subcomplexes on which the $\CC^*$ action \eqref{eq:CActionOnHolomorphicPairsDeformationComplex} has zero weight, positive weight, and negative weight. A splitting $E=L_1\oplus L_2$ of a complex rank-two vector bundle over a complex manifold into a direct sum of two complex line bundles, $L_1$ and $L_2$, induces a decomposition
\begin{equation}
\label{eq:slDecomp}
\fsl(E)
\cong
\ubarCC\oplus (L_1\otimes L_2^*)\oplus (L_2\otimes L_1^*),
\end{equation}
where $\ubarCC=X\times\CC$ is the product complex line bundle. Just as in the definition of the weight of an $S^1$ action following equation \eqref{eq:FixedPointTangentSpaceDecomp}, we define the \emph{weight} of a $\CC^*$ action $\rho:\CC^*\times V\to V$ on a complex vector space $V$ to be an integer $m$ if $\rho(\lambda)(v)=\lambda^m v$, for all $\lambda\in\CC^*$ and $v\in V$. 

The forthcoming Lemma \ref{lem:WeightSpaceDecomp} serves to motivate the constructions in this section. We apply Lemma \ref{lem:WeightSpaceDecomp} and its method of proof in Sections \ref{subsec:Zero-weight_subcomplex}, \ref{subsec:Positive-weight_subcomplex}, and \ref{subsec:The_negative-weight_subcomplex} to define the subcomplexes of the holomorphic pair elliptic complex \eqref{eq:AbbreviatedHolomorphicPairDefComplex} on which the $\CC^*$ action \eqref{eq:CActionOnHolomorphicPairsDeformationComplex} has zero, positive, and negative weight, respectively.  In Section \ref{subsec:Decomposition_elliptic_deformation_complex_holomorphic_pairs}, we prove that the elliptic complex \eqref{eq:AbbreviatedHolomorphicPairDefComplex} decomposes into a direct sum of these three subcomplexes.

\begin{lem}
\label{lem:WeightSpaceDecomp}
Let $E$ be a complex rank-two vector bundle over a complex manifold which admits a decomposition $E=L_1\oplus L_2$ as a direct sum of complex line bundles. Then the isomorphism
\begin{equation}
\label{eq:DeformationComplexBundleDecomposition}
\left(\ubarCC\oplus L_1 \right)
\oplus
\left(L_2\otimes L_1^*\oplus L_2 \right)
\oplus
\left(L_1\otimes L_2^* \right)
\cong
\fsl(E)\oplus E
\end{equation}
is $\CC^*$-equivariant with respect to the $\CC^*$ action on $\fsl(E)\oplus E$ given by
\[
  \left( \la, (\zeta,\sigma)\right)\mapsto \left( \rho_2^\CC(\la) \zeta \rho_2^\CC(\la)^{-1},\rho_2^\CC(\la)\sigma\right),
  \quad\text{for all $\lambda\in\CC^*$, $\zeta\in\fsl(E)$, and $\sigma\in E$,}
\]
and the $\CC^*$ action on $\left(\ubarCC\oplus L_1 \right)\oplus \left(L_2\otimes L_1^*\oplus L_2 \right)
\oplus\left(L_1\otimes L_2^* \right)$ given by
\[
\left( 
\la, (\zeta_\CC,\si_1), (\alpha_{21},\si_2), \alpha_{12}
\right)
\mapsto
\left(
 (\zeta_\CC,\si_1), (\la\alpha_{21},\la\si_2), \la^{-1}\alpha_{12}
\right), \quad\text{for all } \lambda\in\CC^*,
\]
and for all $(\zeta_\CC,\si_1)\in \ubarCC\oplus L_1$, and $(\alpha_{21},\si_2)\in L_2\otimes L_1^*\oplus L_2$, and $\alpha_{12}\in L_1\otimes L_2^*$, where $\rho_2^\CC:\CC^*\to\GL(E)$ is the homomorphism \eqref{eq:DefineS12ActionsAtReducible}.
\end{lem}

\begin{proof}
We write the isomorphism \eqref{eq:DeformationComplexBundleDecomposition} as
\[
\left( (\zeta_\CC,\si_1), (\alpha_{21},\si_2), \alpha_{12}\right)
\mapsto
(\zeta,\sigma)
:=
\left( 
    \begin{pmatrix}
    \zeta_\CC & \alpha_{21} \\ \alpha_{12} & -\zeta_\CC
    \end{pmatrix},
    \begin{pmatrix}
    \si_1 \\ \si_2
    \end{pmatrix}
\right),
\]
where on the right-hand side we have written the element of $\fsl(E)$ as a matrix  and the element of $E$ as a column vector, both with respect with respect to the splitting $E=L_1\oplus L_2$. Under the isomorphism \eqref{eq:DeformationComplexBundleDecomposition}, we see that
\[
\left(
 (\zeta_\CC,\si_1), (\la\alpha_{21},\la\si_2), \la^{-1}\alpha_{12}
\right)
\mapsto
\left( \begin{pmatrix}
\zeta_\CC & \la\alpha_{21} \\ \la^{-1}\alpha_{12} & -\zeta_\CC
\end{pmatrix},
\begin{pmatrix}
\si_1 \\ \la\si_2
\end{pmatrix}
\right), \quad\text{for all } \lambda \in \CC^*,
\]
and we observe that
\begin{align*}
\left( \begin{pmatrix}
\zeta_\CC & \la\alpha_{21} \\ \la^{-1}\alpha_{12} & -\zeta_\CC
\end{pmatrix},
\begin{pmatrix}
\si_1 \\ \la\si_2
\end{pmatrix}
\right)
&=
\left( 
    \begin{pmatrix} 1 & 0 \\ 0 & \la \end{pmatrix}
    \begin{pmatrix} \zeta_\CC & \alpha_{21} \\ \alpha_{12} & -\zeta_\CC
    \end{pmatrix}
    \begin{pmatrix} 1 & 0 \\ 0 & \la^{-1} \end{pmatrix},
    \begin{pmatrix} 1 & 0 \\ 0 & \la \end{pmatrix}
    \begin{pmatrix}
    \si_1 \\ \si_2
    \end{pmatrix}
\right)
\\
&=
\left(
    \rho_2^\CC(\la)
    \begin{pmatrix} \zeta_\CC & \alpha_{21} \\ \alpha_{12} & -\zeta_\CC
    \end{pmatrix}
    \rho_2^\CC(\la)^{-1},
    \rho_2^\CC(\la)
    \begin{pmatrix}
    \si_1 \\ \si_2
    \end{pmatrix}
\right), \quad\text{for all } \lambda \in \CC^*,
\end{align*}
as asserted.
\end{proof}

We use the decomposition \eqref{eq:DeformationComplexBundleDecomposition} to define subcomplexes of the elliptic complex \eqref{eq:AbbreviatedHolomorphicPairDefComplex} using decompositions of the vector spaces $\sE_k$ appearing in \eqref{eq:DefineSpacesOfAbbreviatedHolomorphicPairDefComplex} into subspaces on which $\CC^*$ acts with zero, positive, and negative weights:
\begin{equation}
\label{eq:DefineSpacesOfAbbreviatedHolomorphicPairDefComplexWeightDecomp}
\sE_k=\sE_k^0\oplus \sE_k^+\oplus \sE_k^-, \quad\text{for } k = 0, \ldots, n+1.
\end{equation}
The $\CC^*$-equivariance \eqref{eq:CEquivarianceOfDefComplexOperators} of the operators $\bar\rd_{E,\varphi}^k$ established in Lemma \ref{lem:S1EquivariantDeformationComplex} then implies that they are given by the direct sum of their restrictions to the subspaces in \eqref{eq:DefineSpacesOfAbbreviatedHolomorphicPairDefComplexWeightDecomp}
and so the elliptic complex \eqref{eq:AbbreviatedHolomorphicPairDefComplex} is a direct sum of the elliptic complex formed by these vector spaces.

\subsection{Zero-weight subcomplex}
\label{subsec:Zero-weight_subcomplex}
We continue to assume that the holomorphic pair $(\bar\rd_A,\varphi)$ is split with respect to a decomposition $E=L_1\oplus L_2$ over a complex manifold of complex dimension $n$.  We write $\bar\rd_E=\bar\rd_{L_1}\oplus \bar\rd_{L_2}$, where $\bar\rd_{L_i}$ is a holomorphic structure on the complex line bundle $L_i$, and recall that $\varphi\in\Om^0(L_1)$. We define the \emph{tangential} or \emph{zero-weight} subcomplex of the elliptic complex \eqref{eq:AbbreviatedHolomorphicPairDefComplex} for holomorphic pairs by
\begin{equation}
\label{eq:ZeroWeightPairDefComplex}
\begin{CD}
0
@>>>
\sE_0^0
@>\bar\rd_{E,\varphi}^{0,0} >>
\sE_1^0
@>\bar\rd_{E,\varphi}^{0,1} >>
\sE_2^0
@>\bar\rd_{E,\varphi}^{0,2} >>
\dots
@>\bar\rd_{E,\varphi}^{0,n} >>
\sE_{n+1}^0
@>>>
0
\end{CD}
\end{equation}
where the complex vector spaces $\sE_k^0$ are defined by
\begin{equation}
\label{eq:SpacesOfWeightZeroSubcomplex}
\sE_k^0 :=\Om^{0,k}(\ubarCC)\oplus \Om^{0,k-1}(L_1), \quad\text{for $k\in\ZZ$},
\end{equation}
with the convention that $\Om^{0,k}(\underline{\CC}) = (0)$ and $\Om^{0,k}(L_1) = (0)$ when $k<0$ or $k>n$,
and the operators $\bar\rd_{E,\varphi}^{0,k}$ are defined by
\begin{multline}
\label{eq:dZeroWeightComplex}
\bar\partial_{E,\varphi}^{0,k}(\alpha_\CC,\sigma_1)
:=
\begin{pmatrix}
\bar\partial\alpha_\CC
\\
\bar\partial_{L_1}\sigma_1 + (-1)^{k-1}\alpha_\CC\varphi
\end{pmatrix},
\\
\text{for $k=0,1,\dots,n$ and all } (\alpha_\CC,\sigma_1) \in \Omega^{0,k}(\ubarCC)\oplus \Omega^{0,k-1}(L_1).
\end{multline}
The fact that the operators defined in \eqref{eq:dZeroWeightComplex} satisfy $\bar\rd_{E,\varphi}^{0,k}\circ\bar\rd_{E,\varphi}^{0,k-1}=0$ and thus ensure that \eqref{eq:ZeroWeightPairDefComplex} forms a complex will be shown in the forthcoming Corollary \ref{cor:ZeroAndPositiveAreComplexes}. If we further choose smooth Hermitian metrics on the complex vector bundle $E$ and on the complex base manifold and assume that the manifold is closed, we may define the harmonic spaces of the zero-weight subcomplex by setting
\begin{equation}
  \label{eq:HarmonicSpacesForZeroWeight}
  \bH_{\bar\partial_E,\varphi}^{0,k}:=\Ker\left( \bar\rd_{E,\varphi}^{0,k}+ \bar\rd_{E,\varphi}^{0,k-1,*}\right),
\quad\text{for $k=0,1,\dots,n+1$,}
\end{equation}
where $\bar\rd_{E,\varphi}^{0,k,*}$ is the $L^2$-adjoint of $\bar\rd_{E,\varphi}^{0,k}$.

\subsection{Positive-weight subcomplex}
\label{subsec:Positive-weight_subcomplex}
We define the \emph{positive-weight} subcomplex of the elliptic complex \eqref{eq:AbbreviatedHolomorphicPairDefComplex} for holomorphic pairs by
\begin{equation}
  \label{eq:PositiveWeightPairDefComplex}
\begin{CD}
0
@>>>
\sE_0^+
@>\bar\rd_{E,\varphi}^{+,0} >>
\sE_1^+
@>\bar\rd_{E,\varphi}^{+,1} >>
\sE_2^+
@>\bar\rd_{E,\varphi}^{+,2} >>
\dots
@>\bar\rd_{E,\varphi}^{+,n} >>
\sE_{n+1}^+
@>>>
0
\end{CD}
\end{equation}
where the complex vector spaces $\sE_k^+$ are defined by
\begin{equation}
\label{eq:SpacesOfWeight+Subcomplex}
\sE_k^+:=\Om^{0,k}(L_2\otimes L_1^*)\oplus \Om^{0,k-1}(L_2), \quad\text{for $k\in\ZZ$,}
\end{equation}
with the convention that $\Om^{0,k}(L_2\otimes L_1^*) = (0)$ and $\Om^{0,k}(L_2) = (0)$ when $k<0$ or $k>n$. The operators $\bar\rd_{E,\varphi}^{+,k}$ are defined by
\begin{multline}
\label{eq:dPlusWeightComplex}
\bar\partial_{E,\varphi}^{+,k}(\alpha_{21},\sigma_2)
:=
\begin{pmatrix}
\bar\partial\alpha_\CC
\\
\bar\partial_{L_2\otimes L_1^*}\sigma_2 + (-1)^{k-1}\alpha_{21}\varphi
\end{pmatrix},
\\
\text{for $k=0,1,\dots,n$ and all }
(\alpha_{21},\sigma_2) \in \Omega^{0,k}(L_2\otimes L_1^*)\oplus \Omega^{0,k-1}(L_2),
\end{multline}
where $\bar\partial_{L_2\otimes L_1^*}$ denotes the holomorphic structure on the complex line bundle $L_1\otimes L_2^*$ induced by the holomorphic structure $\bar\partial_{L_2}$ on the complex line bundle $L_2$ and the holomorphic structure $\bar\partial_{L_1^*}$ on the complex dual line bundle $L_1^*$ induced by $\bar\partial_{L_1}$ on $L_1$ (see Napier and Ramachandran \cite[Proposition 3.1.10]{Napier_Ramachandran_introduction_riemann_surfaces}).

As with the operators in the zero-weight subcomplex, we will show that the definition of the operators \eqref{eq:dPlusWeightComplex} ensures that \eqref{eq:PositiveWeightPairDefComplex} forms a complex in the forthcoming Corollary \ref{cor:ZeroAndPositiveAreComplexes}. If we further choose smooth Hermitian metrics on the complex vector bundle $E$ and on the base complex manifold and assume that the manifold is closed, we may define the harmonic spaces for the positive-weight subcomplex by 
\begin{equation}
  \label{eq:HarmonicSpacesForPositiveWeight}
\bH_{\bar\partial_E,\varphi}^{+,k} :=\Ker\left( \bar\rd_{E,\varphi}^{+,k}+ \bar\rd_{E,\varphi}^{+,k-1,*}\right),
\quad\text{for $k=0,1,\dots,n+1$,}
\end{equation}
where $\bar\rd_{E,\varphi}^{+,k,*}$ is the $L^2$-adjoint of $\bar\rd_{E,\varphi}^{+,k}$. 

\subsection{Negative-weight subcomplex}
\label{subsec:The_negative-weight_subcomplex}
We define the \emph{negative-weight} subcomplex of the elliptic complex \eqref{eq:AbbreviatedHolomorphicPairDefComplex} for holomorphic pairs by
\begin{equation}
  \label{eq:NegativeWeightPairDefComplex}
\begin{CD}
0
@>>>
\sE_0^-
@>\bar\rd_{E,\varphi}^{-,0} >>
\sE_1^+
@>\bar\rd_{E,\varphi}^{-,1} >>
\sE_2^+
@>\bar\rd_{E,\varphi}^{-.2} >>
\dots
@>\bar\rd_{E,\varphi}^{-,n} >>
\sE_{n+1}^-
@>>>
0
\end{CD}
\end{equation}
where the complex vector spaces $\sE_k^-$ are defined by
\begin{equation}
\label{eq:SpacesOfWeight-Subcomplex}
\sE_k^-:=\Om^{0,k}(L_1\otimes L_2^*), \quad\text{for $k\in\ZZ$,}
\end{equation}
with the convention that $\Om^{0,k}(L_1\otimes L_2^*) = (0)$ when $k<0$ or $k>n$. We define the operators $\bar\rd_{E,\varphi}^{-,k}$ by
\begin{align}
  \label{eq:dSNegWeightPair}
\bar\partial_{E,\varphi}^{-,k}
  &:=
  \bar\rd_{L_1\otimes L_2^*}, \quad\text{for $k=0,\dots,n$},
\end{align}
where $\bar\partial_{L_1\otimes L_2^*}$ denotes the holomorphic structure on the complex line bundle $L_2\otimes L_1^*$ induced by the holomorphic structure $\bar\partial_{L_1}$ on the complex line bundle $L_1$ and the holomorphic structure $\bar\partial_{L_2^*}$ on the complex dual line bundle $L_2^*$ induced by $\bar\partial_{L_2}$ on $L_2$. Because $\bar\rd_{L_1\otimes L_2^*}$ is a holomorphic structure, the definition of the operators \eqref{eq:dSNegWeightPair} ensures that \eqref{eq:NegativeWeightPairDefComplex} is indeed a complex. If we further choose smooth Hermitian metrics on the complex vector bundle $E$ and on the base complex manifold and assume that the manifold is closed, we may define the harmonic spaces of the negative-weight subcomplex by
\begin{equation}
  \label{eq:HarmonicSpacesForNegativeWeight}
\bH_{\bar\partial_E,\varphi}^{-,k}:=\Ker\left( \bar\rd_{E,\varphi}^{-,k}+ \bar\rd_{E,\varphi}^{-,k-1,*}\right),
\quad\text{for $k=0,1,\dots,n+1$,}
\end{equation}
where $\bar\rd_{E,\varphi}^{-,k,*}$ is the $L^2$-adjoint of $\bar\rd_{E,\varphi}^{-,k}$.

\subsection{Decomposition of the elliptic complex for holomorphic pairs}
\label{subsec:Decomposition_elliptic_deformation_complex_holomorphic_pairs}
We first describe a decomposition of the complex vector spaces in the elliptic complex \eqref{eq:AbbreviatedHolomorphicPairDefComplex} into zero, positive, and negative-weight subbundles for the $\CC^*$ action.

\begin{lem}[Weight decomposition of the vector spaces in the elliptic complex for a split holomorphic pair]
\label{lem:TheWeightSpaceDecomposition}
Continue the hypotheses and notation of Lemma \ref{lem:S1EquivariantDeformationComplex}. Then the complex vector spaces $\sE_k$ appearing in the elliptic complex \eqref{eq:AbbreviatedHolomorphicPairDefComplex} for holomorphic pairs admit direct-sum decompositions
\begin{equation}
\label{eq:DefineSpacesOfAbbreviatedHolomorphicPairDefComplexWeightDecomp1}
\sE_k=\sE_k^0\oplus \sE_k^+\oplus \sE_k^-, \quad\text{for } k = 0,1,\ldots,n+1,
\end{equation}
where $\sE_k^0$ is defined in \eqref{eq:SpacesOfWeightZeroSubcomplex}, and $\sE_k^+$ in \eqref{eq:SpacesOfWeight+Subcomplex}, and $\sE_k^-$ in \eqref{eq:SpacesOfWeight-Subcomplex}.
The decomposition \eqref{eq:DefineSpacesOfAbbreviatedHolomorphicPairDefComplexWeightDecomp1} has the following properties:
\begin{enumerate}
\item
\label{item:TheWeightSpaceDecompositionOrthogonality}
The decomposition \eqref{eq:DefineSpacesOfAbbreviatedHolomorphicPairDefComplexWeightDecomp1} is $L^2$-orthogonal.
\item
\label{item:TheWeightSpaceDecomposition}
The summands in the decomposition \eqref{eq:DefineSpacesOfAbbreviatedHolomorphicPairDefComplexWeightDecomp1}
are weight spaces for the $\CC^*$ action on $\sE_k$ defined in \eqref{eq:CActionOnHolomorphicPairsDeformationComplex}, with the action having weight zero on $\sE_k^0$, weight one on $\sE_k^+$, and weight negative one on $\sE_k^-$.
\end{enumerate}
\end{lem}

\begin{proof}
The decomposition  \eqref{eq:DefineSpacesOfAbbreviatedHolomorphicPairDefComplexWeightDecomp1} follows from
the splitting $E=L_1\oplus L_2$, the splitting of $\fsl(E)$ in \eqref{eq:slDecomp}, and the definitions of the vector spaces $\sE_k^0$ and $\sE_k^+$ and $\sE_k^-$ in \eqref{eq:SpacesOfWeightZeroSubcomplex} and \eqref{eq:SpacesOfWeight+Subcomplex} and \eqref{eq:SpacesOfWeight-Subcomplex}, respectively. The subspaces in \eqref{eq:DefineSpacesOfAbbreviatedHolomorphicPairDefComplexWeightDecomp1} are $L^2$-orthogonal because the direct-sum decomposition of $\fsl(E)\oplus E$ in \eqref{eq:DeformationComplexBundleDecomposition} is pointwise orthogonal. Finally, the characterization of the subspaces as weight spaces follows from Lemma \ref{lem:WeightSpaceDecomp} and the definitions \eqref{eq:SpacesOfWeightZeroSubcomplex} and \eqref{eq:SpacesOfWeight+Subcomplex} and \eqref{eq:SpacesOfWeight-Subcomplex}.
\end{proof}

We now describe the corresponding decomposition for the operators $\bar\rd_{E,\varphi}^k$.

\begin{lem}[Weight decomposition of the differentials in the elliptic complex for a split holomorphic pair]
\label{lem:WeightSpaceDecompositionOfHolomorphicPairsDeformationComplex}
Continue the hypotheses and notation of Lemma \ref{lem:S1EquivariantDeformationComplex}.  Then the operators $\bar\rd_{E,\varphi}^k$ in the elliptic complex \eqref{eq:AbbreviatedHolomorphicPairDefComplex} for holomorphic pairs admit direct-sum decompositions
\begin{equation}
\label{eq:DirectSumOfOperators}
\bar\rd_{E,\varphi}^k=
\bar\rd_{E,\varphi}^{0,k}
\oplus
\bar\rd_{E,\varphi}^{+,k}
\oplus
\bar\rd_{E,\varphi}^{-,k},
\quad\text{for } k = 0,1,\ldots,n,
\end{equation}
with respect to the decomposition \eqref{eq:DefineSpacesOfAbbreviatedHolomorphicPairDefComplexWeightDecomp1}
of the complex vector spaces $\sE_k$.
\end{lem}

\begin{proof}
By comparing the definitions \eqref{eq:dZeroWeightComplex} and \eqref{eq:dPlusWeightComplex} and \eqref{eq:dSNegWeightPair} of the operators $\bar\rd_{E,\varphi}^{0,k}$ and $\bar\rd_{E,\varphi}^{+,k}$ and $\bar\rd_{E,\varphi}^{-,k}$ with the definition \eqref{eq:dkStablePair} of $\bar\rd_{E,\varphi}^k$, respectively, we see that $\bar\rd_{E,\varphi}^{0,k}$ and $\bar\rd_{E,\varphi}^{+,k}$ and $\bar\rd_{E,\varphi}^{-,k}$ are equal to
the restrictions of $\bar\rd_{E,\varphi}^k$ to $\sE_k^0$ and $\sE_k^+$ and $\sE_k^-$, respectively. The $\CC^*$-equivariance \eqref{eq:CEquivarianceOfDefComplexOperators} of $\bar\rd_{E,\varphi}^k$ established in Lemma \ref{lem:S1EquivariantDeformationComplex} and the characterization of $\sE_k^0$ and $\sE_k^+$ and $\sE_k^-$ as weight spaces with distinct weights in Item \eqref{item:TheWeightSpaceDecomposition} of Lemma \ref{lem:TheWeightSpaceDecomposition} imply that
\[
\bar\rd_{E,\varphi}^k\sE_k^0 \subset \sE_{k+1}^0,\quad
\bar\rd_{E,\varphi}^k\sE_k^+ \subset \sE_{k+1}^+,\ \text{and}\quad
\bar\rd_{E,\varphi}^k\sE_k^- \subset \sE_{k+1}^-.
\]
Thus, each operator $\bar\rd_{E,\varphi}^k$  has the asserted direct-sum decomposition.
\end{proof}

By combining the direct-sum decomposition of $\bar\rd_{E,\varphi}^k$ in \eqref{eq:DirectSumOfOperators} with the equality $\bar\rd_{E,\varphi}^{k+1}\circ\bar\rd_{E,\varphi}^k=0$, we obtain the

\begin{cor}[Weight decomposition of the elliptic complex of a split holomorphic pair into weight subcomplexes]
\label{cor:ZeroAndPositiveAreComplexes}
Continue the hypotheses and notation of Lemma \ref{lem:S1EquivariantDeformationComplex}. Then
\[
\bar\rd_{E,\varphi}^{0,k+1}\circ \bar\rd_{E,\varphi}^{0,k}=0,\quad
\bar\rd_{E,\varphi}^{+,k+1}\circ \bar\rd_{E,\varphi}^{+,k}=0,\quad \text{and}\quad
\bar\rd_{E,\varphi}^{-,k+1}\circ \bar\rd_{E,\varphi}^{-,k}=0,
\quad\text{for } k=0,\dots,n-1.
\]
Thus, the elliptic complex \eqref{eq:AbbreviatedHolomorphicPairDefComplex} admits a decomposition into a direct sum of the subcomplexes
\[
(\sE_k,\bar\rd_{E,\varphi}^k)
=
(\sE_k^0,\bar\rd_{E,\varphi}^{0,k})
\oplus
(\sE_k^+,\bar\rd_{E,\varphi}^{+,k})
\oplus
(\sE_k^-,\bar\rd_{E,\varphi}^{-,k}),
\]
where $(\sE_k^0,\bar\rd_{E,\varphi}^{0,k})$ is defined in \eqref{eq:ZeroWeightPairDefComplex},
$(\sE_k^+,\bar\rd_{E,\varphi}^{+,k})$ is defined in \eqref{eq:PositiveWeightPairDefComplex},
and
$(\sE_k^-,\bar\rd_{E,\varphi}^{-,k})$ is defined in \eqref{eq:NegativeWeightPairDefComplex}.
\end{cor}

Together, Lemmas \ref{lem:TheWeightSpaceDecomposition} and \ref{lem:WeightSpaceDecompositionOfHolomorphicPairsDeformationComplex} and Corollary \ref{cor:ZeroAndPositiveAreComplexes} imply that the elliptic complex \eqref{eq:AbbreviatedHolomorphicPairDefComplex} for holomorphic pairs admits a  decomposition as a complex into a direct sum of the zero-weight, positive-weight, and negative-weight subcomplexes.  
By \eqref{eq:EllipticCohomologyDecomposition}, this implies that the associated harmonic spaces admit a similar decomposition (see Rotman \cite[Exercise 10.30]{RotmanAlgebra_2002}).

\begin{cor}[Weight decomposition of the harmonic spaces for the elliptic complex of a split holomorphic pair]
\label{cor:WeightOfS1ActionOnCohomology}
Continue the hypotheses and notation of Lemma \ref{lem:S1EquivariantDeformationComplex} and assume further that the complex vector bundle $E$ and the complex base manifold are endowed with smooth Hermitian metrics and that the manifold is closed. Then the harmonic spaces \eqref{eq:H_dbar_Avarphi^0bullet} for the elliptic complex \eqref{eq:AbbreviatedHolomorphicPairDefComplex} for holomorphic pairs admit direct-sum decompositions
\begin{equation}
  \label{eq:WeightSpaceDecompOfDeformationCohomology}
\bH^k_{\bar\rd_E,\varphi}
=
\bH^{0,k}_{\bar\rd_E,\varphi}\oplus \bH^{+,k}_{\bar\rd_E,\varphi} \oplus \bH^{-,k}_{\bar\rd_E,\varphi},
\quad\text{for $k=0,\dots,n+1$},
\end{equation}
where $\bH^{0,k}_{\bar\rd_E,\varphi}$ is as in \eqref{eq:HarmonicSpacesForZeroWeight}, and
$\bH^{+,k}_{\bar\rd_E,\varphi}$ is as in \eqref{eq:HarmonicSpacesForPositiveWeight}, and
$\bH^{-,k}_{\bar\rd_E,\varphi}$ is as in \eqref{eq:HarmonicSpacesForNegativeWeight}.
The complex vector spaces $\bH^{0,k}_{\bar\rd_E,\varphi}$ and $\bH^{+,k}_{\bar\rd_E,\varphi}$ and $\bH^{-,k}_{\bar\rd_E,\varphi}$ are, respectively, the zero, positive, and negative-weight subspaces for the harmonic spaces of the elliptic complex \eqref{eq:AbbreviatedHolomorphicPairDefComplex} with respect to the $\CC^*$ action on $\bH^k_{\bar\rd_E,\varphi}$ defined in Corollary \ref{cor:CActionOnHarmonicSpacesOfHolomPairsComplex}.
\end{cor}

\section{Euler characteristics of the subcomplexes}
\label{sec:IndOfSubcomplex}
We compute the Euler characteristics of the positive, negative, and zero-weight complexes \eqref{eq:PositiveWeightPairDefComplex}, \eqref{eq:NegativeWeightPairDefComplex}, and \eqref{eq:ZeroWeightPairDefComplex}, respectively, by applying the Hirzebruch--Riemann--Roch
Theorem (see Friedman \cite[Chapter 1, Theorem 9, p. 15]{FriedmanBundleBook} or Wells \cite[Chapter IV, Theorem 5.8, p. 153]{Wells3}).  Recall from Wells \cite[Chapter IV, Example 5.5, p. 150]{Wells3} that given a complex vector bundle $F$ with holomorphic structure $\bar\partial_F$ over a complex manifold $X$ of complex dimension $n$, the Euler characteristic of the elliptic complex associated with $(F,\bar\partial_F)$ is defined to be (see \eqref{eq:Index_elliptic_complex})\[
  \Euler_\CC(\bar\partial_F) :=
\sum_{k=0}^n (-1)^k \dim_\CC \bH_{\bar\partial_F}^k,
\]
where
\[
\bH_{\bar\partial_F}^k
:= \Ker\left( \bar\partial_F+\bar\partial_F^*:\Om^{0,k}(X;F) \to \Om^{0,k+1}(X;F)\oplus \Om^{0,k-1}(X;F)\right).
\]
We will use the following form of the Hirzebruch--Riemann--Roch Theorem. Recall from Wells \cite[Chapter IV, Section 5, p. 152]{Wells3} or Lawson and Michelson \cite[Chapter III, Equation (11.22), p. 234 and Example 11.11, p. 230]{LM} that the Chern character $\ch(\sE)$ and the Todd class $\operatorname{Todd}(X)$ are formal power series in the characteristic
classes of $\sE$ and $TX$, respectively.

\begin{thm}[Hirzebruch--Riemann--Roch Theorem]
\label{thm:HRR}
(See Hirzebruch \cite[Theorem 21.1.1]{Hirzebruch_topological_methods_algebraic_geometry}.)
Let $(F,\bar\partial_F)$ be a holomorphic vector bundle over a closed, complex manifold $X$. Then
\begin{equation}
\label{eq:GeneralHRR}
\Euler_\CC(\bar\partial_F)
=
\left\langle \ch(F)\smile\operatorname{Todd}(X),[X]\right\rangle,
\end{equation}
where $\ch(F)$ is the Chern character of $F$ and $\operatorname{Todd}(X)$ is the Todd class of $X$. If $(L,\bar\partial_L)$ is a holomorphic line bundle over a closed, complex surface $X$, then
\begin{equation}
\label{eq:SimpleHRR}
\Euler_\CC(\bar\partial_L)
=
\chi_h(X)+\frac{1}{2}c_1(L)\cdot c_1(X) +\frac{1}{2}c_1(L)^2,
\end{equation}
where $\chi_h(X)$ is the holomorphic Euler characteristic \eqref{eq:DefineTopCharNumbersOf4Manifold} of $X$.
\end{thm}

\begin{rmk}[Specialization from complex manifolds to complex surfaces]
\label{rmk:DifficultyInRRComputationInHigherDimension}
We shall now specialize to the case where $X$ has complex dimension $n=2$ since the computations involving the Todd class for higher-dimensional complex manifolds appearing in the formula \eqref{eq:GeneralHRR} become more involved, even if still tractable.
\end{rmk}

\begin{lem}[Euler characteristic of the negative-weight subcomplex]
\label{lem:IndexOfRedStablePairsDefComplexNormalMinus}
Continue the hypotheses and notation of Lemma \ref{lem:S1EquivariantDeformationComplex} and assume further that the complex vector bundle $E$ and the complex base manifold $X$ are endowed with smooth Hermitian metrics and that $X$ is closed with complex dimension two. Then the Euler characteristic of the negative-weight subcomplex \eqref{eq:NegativeWeightPairDefComplex} is given by
\begin{align*}
  \Euler_\CC\left(\bar\rd_{E,\varphi}^{-,\bullet}\right)
  &:=
    \sum_{k=0}^2 (-1)^k\dim_\CC\bH^{-,k}_{\bar\rd_E,\varphi}
  \\
  &\,= \Euler_\CC\left( \bar\partial_{L_1\otimes L_2^*}\right)
  \\
  &\,= \chi_h(X)
    +
    \frac{1}{2}\left( c_1(L_1)-c_1(L_2)\right)\cdot c_1(X)
    +
    \frac{1}{2}\left( c_1(L_1)-c_1(L_2)\right)^2.
\end{align*}
\end{lem}

\begin{proof}
Because the operators $\bar\rd_{E,\varphi}^{-,k}$ in the negative-weight complex \eqref{eq:NegativeWeightPairDefComplex} are given by a holomorphic structure $\bar\partial_{L_1\otimes L_2^*}$ on the complex line bundle $L_1\otimes L_2^*$ according to \eqref{eq:dSNegWeightPair}, we see that
\[
  \Euler_\CC\left(\bar\rd_{E,\varphi}^{-,\bullet}\right)
  =
  \Euler_\CC\left( \bar\partial_{L_1\otimes L_2^*}\right)
\]
and hence the formula for the Euler characteristic follows immediately from equation \eqref{eq:SimpleHRR}.
\end{proof}

\begin{lem}[Euler characteristic of the positive-weight subcomplex]
\label{lem:IndexOfRedStablePairsDefComplexNormalPlus}
Continue the hypotheses and notation of Lemma \ref{lem:IndexOfRedStablePairsDefComplexNormalMinus}. Then the Euler characteristic of the positive weight subcomplex \eqref{eq:PositiveWeightPairDefComplex} is
\begin{align*}
  \Euler_\CC\left(\bar\rd_{E,\varphi}^{+,\bullet}\right)
  &:=
    \sum_{i=0}^2 (-1)^k\dim_\CC\bH^{+,k}_{\bar\rd_E,\varphi}
  \\
  &\,=
    \Euler_\CC(\bar\rd_{L_2\otimes L_1^*})-\Euler_\CC(\bar\rd_{L_2})
  \\
  &\,= \frac{1}{2}\left( c_1(L_2)-c_1(L_1)\right)\cdot c_1(X)
    +
    \frac{1}{2}\left( c_1(L_2)-c_1(L_1)\right)^2
  \\
  &\quad
    -\frac{1}{2}c_1(L_2)\cdot c_1(X) -\frac{1}{2}c_1(L_2)^2.
\end{align*}
\end{lem}

\begin{proof}
We first claim that changing the differentials in the complex \eqref{eq:PositiveWeightPairDefComplex} through the homotopy $\partial_{E,t\varphi}^{+,k}$ for $t\in[0,1]$ does not change the Euler characteristic of the complex. The pairs $(\bar\rd_E,t\varphi)$ with $t\in[0,1]$ are still split and holomorphic, so the operators $\partial_{E,t\varphi}^{+,k}$ still define a complex.  According to Gilkey \cite[Equation (1.5.18)]{Gilkey2}, the index of the `rolled-up' elliptic operator \cite[Equation (1.5.1)]{Gilkey2},
\[
Q_{E,t\varphi} :\bigoplus_{k \text{ even}}\sE_k^+ \to \bigoplus_{k\text{ odd}}\sE_k^+,
\]
where
\[
  Q_{E,t\varphi} := \bigoplus_{k\text{ even}} \left(\bar\rd_{E,t\varphi}^{+,k} + \bar\rd_{E,t\varphi}^{+,k-1,*}\right),
\]  
is equal to the Euler characteristic of the complex $(\sE_\bullet, \partial_{E,t\varphi}^{+,\bullet})$. Because the leading symbols of the differential operators $\partial_{E,t\varphi}^{+,k}$ are independent of $t$, the index of the operator $Q_{E,t\varphi}$ is constant with respect to $t$ by Gilkey \cite[Lemma 1.4.5]{Gilkey2}. Hence, the Euler characteristic of the complex \eqref{eq:PositiveWeightPairDefComplex} is equal to that of the complex $(\sE_\bullet, \partial_{E,0}^{+,\bullet})$.

When $\varphi=0$, the complex \eqref{eq:PositiveWeightPairDefComplex} splits into a direct sum of the subcomplex
\begin{equation}
\label{eq:PositiveWeightPairDefComplexComponent1}
\begin{CD}
0 @>>>
\Om^0(L_2\otimes L_1^*)
@> \bar\partial_{L_2\otimes L_1^*}>>
\Om^{0,1}(L_2\otimes L_1^*)
@> \bar\rd_{L_2\otimes L_1^*} >>
\Om^{0,2}(L_2\otimes L_1^*)
@>>>
0
\end{CD}
\end{equation}
and the subcomplex
\begin{equation}
\label{eq:PositiveWeightPairDefComplexComponent2}
\begin{CD}
0 @>>>
\Om^0(L_2)
@> \bar\partial_{L_2}>>
\Om^{0,1}(L_2)
@> \bar\partial_{L_2} >>
\Om^{0,2}(L_2)
@>>>
0
\end{CD}
\end{equation}
By the Hirzebruch--Riemann--Roch formula \eqref{eq:SimpleHRR}, the Euler characteristic of the complex
\eqref{eq:PositiveWeightPairDefComplexComponent1} is
\[
  \Euler_\CC\left(\bar\rd_{L_2\otimes L_1^*}\right)
  =
  \frac{1}{2}\left( c_1(L_2)-c_1(L_1)\right)^2
  +
  \frac{1}{2}\left( c_1(L_2)-c_1(L_1)\right)\cdot c_1(X)
  +
  \chi_h(X),
\]
and the Euler characteristic of the complex \eqref{eq:PositiveWeightPairDefComplexComponent2} is 
\[
  \Euler_\CC\left(\bar\rd_{L_2}\right)
  =
  \frac{1}{2} c_1(L_2)^2 +\frac{1}{2}c_1(L_2)\cdot c_1(X) + \chi_h(X).
\]
Because the terms of the complex \eqref{eq:PositiveWeightPairDefComplexComponent2} appear in the complex \eqref{eq:PositiveWeightPairDefComplex} with the opposite parity, the Euler characteristic of \eqref{eq:PositiveWeightPairDefComplex} is equal to the difference between the Euler characteristics of the complexes \eqref{eq:PositiveWeightPairDefComplexComponent1} and \eqref{eq:PositiveWeightPairDefComplexComponent2}, and this gives the desired formula.
\end{proof}

\begin{lem}[Euler characteristic of the zero-weight subcomplex]
\label{lem:IndexOfRedStablePairsDefComplexNormalZero}
Continue the hypotheses and notation of Lemma \ref{lem:IndexOfRedStablePairsDefComplexNormalMinus}. Then the Euler characteristic of the zero-weight subcomplex \eqref{eq:ZeroWeightPairDefComplex} is
\begin{align*}
  \Euler_\CC\left(\bar\rd_{E,\varphi}^{0,\bullet}\right)
  &:=
    \sum_{i=0}^2 (-1)^k\dim_\CC\bH^{0,k}_{\bar\rd_E,\varphi}
  \\
  &\,=
    \Euler_\CC(\bar\rd)-\Euler_\CC(\bar\rd_{L_1})
  \\
  &\,=  -\frac{1}{2} c_1(L_1)^2 -\frac{1}{2}c_1(L_1)\cdot c_1(X).
\end{align*}
\end{lem}

\begin{proof}
Applying the argument used in the proof of Lemma \ref{lem:IndexOfRedStablePairsDefComplexNormalPlus} to the zero-weight complex \eqref{eq:ZeroWeightPairDefComplex} gives
\[
\Euler_\CC\left(\bar\rd_{E,\varphi}^{0,\bullet}\right)
=
\Euler_\CC\left(\bar\rd_{E,0}^{0,\bullet}\right).
\]
When $\varphi=0$, the complex \eqref{eq:ZeroWeightPairDefComplex} splits into a direct sum of the subcomplex
\begin{equation}
\label{eq:ZeroWeightPairDefComplexComponent1}
\begin{CD}
0 @>>>
\Om^0(\ubarCC)
@> \bar\partial >>
\Om^{0,1}(\ubarCC)
@> \bar\rd >>
\Om^{0,2}(\ubarCC)
@>>>
0
\end{CD}
\end{equation}
and the subcomplex
\begin{equation}
\label{eq:ZeroWeightPairDefComplexComponent2}
\begin{CD}
0 @>>>
\Om^0(L_1)
@> \bar\partial_{L_1}>>
\Om^{0,1}(L_1)
@> \bar\partial_{L_1} >>
\Om^{0,2}(L_1)
@>>>
0
\end{CD}
\end{equation}
By the Hirzebruch--Riemann--Roch index formula \eqref{eq:SimpleHRR}, the Euler characteristic of the complex
\eqref{eq:ZeroWeightPairDefComplexComponent1} is
\[
  \Euler_\CC\left(\bar\rd\right)
  =
  \chi_h(X),
\]
and the Euler characteristic of the complex \eqref{eq:ZeroWeightPairDefComplexComponent2} is
\[
  \Euler_\CC\left(\bar\rd_{L_1}\right)
  =
  \frac{1}{2} c_1(L_1)^2 +\frac{1}{2}c_1(L_1)\cdot c_1(X) + \chi_h(X).
\]
Because the terms of the complex \eqref{eq:ZeroWeightPairDefComplexComponent2} appear in the complex \eqref{eq:ZeroWeightPairDefComplex} with the opposite parity, the Euler characteristic of \eqref{eq:ZeroWeightPairDefComplex} is equal to the difference of the Euler characteristics of the complexes \eqref{eq:ZeroWeightPairDefComplexComponent1} and \eqref{eq:ZeroWeightPairDefComplexComponent2}, and this gives the desired formula.
\end{proof}

We can now complete the

\begin{proof}[Proof of Theorem \ref{mainthm:MorseIndexAtReduciblesOnKahler}]
A  non-Abelian monopole $(A,\Phi)$ on a  K\"ahler surface defines a pre-holomorphic pair $(A,\Phi)$ with $\Phi=(\varphi,\psi)\in\Om^0(E)\oplus \Om^{0,2}(E)$ by Lemma \ref{lem:SO3_monopole_equations_almost_Kaehler_manifold}.  Because $(A,\Phi)$ is type $1$ by hypothesis, we have $\Phi=(\varphi,0)$ and Remark \ref{rmk:Holomorphic_pair} implies that $(\rd_A,\varphi)$ defines a holomorphic pair (and the gauge-equivalence class $[\bar\rd_A,\varphi]$ is the image of $[A,\Phi]$ under the embedding provided by Theorem \ref{thm:Lubke_Teleman_6-3-10}). More explicitly, a type $1$ non-Abelian monopole $(A,\Phi)$ defines a solution $(A,\varphi)$ of the projective vortex equations by Remark \ref{rmk:Projective_vortices_type_1_monopoles} and thus $(\bar\rd_A,\varphi)$ is a holomorphic pair by Theorem \ref{thm:HitchinKobayashiCorrespondenceForPairs}. Because $(A,\varphi)$ is a split solution to the projective vortex equations with respect to the splitting $E=L_1\oplus L_2$, Item \eqref{item:HKIdentificationOfSplitPairs1} of Lemma \ref{lem:HKIdentificationOfSplitPairs} implies that $(\bar\rd_A,\varphi)$ is a split pair with respect to the decomposition $E=L_1\oplus L_2$.

We now show that the virtual Morse--Bott index \eqref{eq:Virtual_Morse-Bott_index_moduli_space_non-abelian_monopoles} for the Hamiltonian function $f$ in \eqref{eq:Hitchin_function} can be computed in terms of the Euler characteristic of the negative-weight complex defined in \eqref{eq:NegativeWeightPairDefComplex}. By Lemma \ref{lem:S1EquivariantKuranishiModelForReducibleIn_non_AbelianMonopoleModuli} and Remark \ref{rmk:S1EquivariantKuranishiModelForReducibleIn_non_AbelianMonopoleModuli_holomorphic_maps}, the Kuranishi model for an open neighborhood of $[A,\Phi]$ in $\sM_\ft^0$,
\[
\bkappa:\bH_{A,\Phi}^1 \supset\sU\to\bH_{A,\Phi}^2 \quad\text{and}\quad \beps:\bH_{A,\Phi}^1 \supset\sU\to\sC_\ft,
\]
is defined by $S^1$-equivariant holomorphic maps, where $S^1$ acts on $\bH_{A,\Phi}^1$ and $\bH_{A,\Phi}^2$ by the actions given in Corollary \ref{cor:S1_Action_on_HarmonicSections_Of_nonAbelianMonopole_Def_Complex} and on $\sC_\ft$ by the action \eqref{eq:S12ActionsOnSpinuPreQuotient}. Hence, the $S^1$ actions on $\bH_{A,\Phi}^1$ and $\bH_{A,\Phi}^2$ determine the virtual Morse--Bott index  as described in Section \ref{sec:Virtual_Morse-Bott_index_Hamiltonian_function_circle_action_complex_analytic_space}.
By Lemma \ref{lem:S1_Equivariance_of_CohomologyIsoms}, the isomorphisms
\eqref{eq:Isomorphisms_of_non_AbelianHarmonics_with_HolomorphPairHarmonics} of real vector spaces,
\[
I_k: \bH_{\bar\rd_A,\varphi}^k
\cong
\bH_{A,\Phi}^k \quad\text{for $k=1,2$},
\]
are $S^1$-equivariant with respect to the $S^1$ action on $\bH_{\bar\rd_A,\varphi}^k$ given in Corollary \ref{cor:CActionOnHarmonicSpacesOfHolomPairsComplex} and the $S^1$ action on $\bH_{A,\Phi}^k$ given in Corollary \ref{cor:S1_Action_on_HarmonicSections_Of_nonAbelianMonopole_Def_Complex}. In addition, these isomorphisms and the almost complex structures on $\bH_{\bar\rd_A,\varphi}^k$ (the fact that $\bH_{\bar\rd_A,\varphi}^k$ is a complex vector space is noted prior to its definition in \eqref{eq:H_dbar_Avarphi^0bullet}) define almost complex structures on $\bH_{A,\Phi}^k$, for $k=1,2$. Hence, the negative-weight subspace of $\bH_{A,\Phi}^k$ for the $S^1$ actions given in
Corollary \ref{cor:S1_Action_on_HarmonicSections_Of_nonAbelianMonopole_Def_Complex} with respect to the almost complex structures given by the isomorphisms with $\bH_{\bar\rd_A,\varphi}^k$ are given by
\begin{equation}
\label{eq:DefineNegativeWeightHarmonicSubspaceForNonAbelianComplex}
\bH_{A,\Phi}^{-,k}:=I_k\left( \bH_{\bar\rd_A,\varphi}^{-,k}\right), \quad\text{for } k=1,2,
\end{equation}
where $\bH_{\bar\rd_A,\varphi}^{-,k}$ is defined in \eqref{eq:HarmonicSpacesForNegativeWeight}. (Recall that we needed the $S^1$ action to be invariant with respect to an almost complex structure in order to resolve
the ambiguity of the signs of the weights of a real representation described in Item \eqref{item:Isomorphic_S1_Representations_Real} of Lemma \ref{lem:Isomorphic_S1_Representations}
and in Proposition \ref{prop:Direct_Sum_Decomposition_of_Real_S1_Representations}.) The preceding isomorphisms and the definition
\eqref{eq:Virtual_Morse-Bott_index_moduli_space_non-abelian_monopoles} of the virtual Morse--Bott index yield
\begin{equation}
\label{eq:MBIndex1}
\la_{[A,\Phi]}^-(f)= 2\left(\dim_\CC \bH_{\bar\rd_A,\varphi}^{-,1}-\dim_\CC \bH_{\bar\rd_A,\varphi}^{-,2}\right).
\end{equation}
Because $\Phi\not\equiv 0$, Lemma \ref{lem:H0_Of_NonAbelianMonopoleComplex_Vanishes} implies that $\bH_{A,\Phi}^0$ is trivial and so, by \eqref{eq:Isomorphisms_of_non_AbelianHarmonics_with_HolomorphPairHarmonics}, the space
$\bH_{\bar\rd_A,\varphi}^0$ also vanishes. Therefore, any subspace of $\bH_{\bar\rd_A,\varphi}^0$ must vanish, so that
\[
\bH_{\bar\rd_A,\varphi}^{-,0}=0.
\]
Similarly, Corollary \ref{cor:Vanishing_third-order_cohomology_group_holomorphic_pair_elliptic_complex_Kaehler_surface} implies that $\bH_{\bar\rd_A,\varphi}^3=0$, and so
\[
\bH_{\bar\rd_A,\varphi}^{-,3}=0.
\]
Therefore, substituting the preceding vanishing results into the definition of the Euler characteristic yields
\[
  \Euler_\CC\left(\bar\rd_{A,\varphi}^\bullet\right)
  =
  \sum_{k=0}^3(-1)^k\dim_\CC \bH_{\bar\rd_A,\varphi}^{-,k}
  =
  -\dim_\CC \bH_{\bar\rd_A,\varphi}^{-,1} + \dim_\CC \bH_{\bar\rd_A,\varphi}^{-,2}.
\]  
Combining the preceding equality with \eqref{eq:MBIndex1} gives
\begin{equation}
\label{eq:MBIndexEqualsIndexOfDefComplex}
  \la_{[A,\Phi]}^-(f) = -2\Euler_\CC\left(\bar\rd_{A,\varphi}^{-,\bullet}\right).
\end{equation}
By applying Lemma \ref{lem:IndexOfRedStablePairsDefComplexNormalMinus} to compute the Euler characteristic on the right-hand side of \eqref{eq:MBIndexEqualsIndexOfDefComplex}, we obtain
\[
  \lambda_{[A,\Phi]}^-(f)
  =
  -2\chi_h(X) -\left( c_1(L_1) -c_1(L_2)\right)\cdot c_1(X) -\left( c_1(L_1)-c_1(L_2)\right)^2.
\]
This verifies \eqref{eq:MorseIndexAtReduciblesOnKahler}.

The argument giving \eqref{eq:MBIndex1} also implies
\begin{subequations}
\label{eq:CoIndexAndNullityEqualTruncatedEulerChar}
\begin{align}
\label{eq:CoIndexAndNullityEquaTruncatedlEulerCharPositiveWeight}
\la_{[A,\Phi]}^+(f)& = 2\left(\dim_\CC \bH_{\bar\rd_A,\varphi}^{+,1}-\dim_\CC \bH_{\bar\rd_A,\varphi}^{+,2}\right),
\\
\label{eq:CoIndexAndNullityEqualTruncatedEulerCharZeroWeight}
\la_{[A,\Phi]}^0(f)& = 2\left(\dim_\CC \bH_{\bar\rd_A,\varphi}^{0,1}-\dim_\CC \bH_{\bar\rd_A,\varphi}^{0,2}\right).
\end{align}
\end{subequations}
As before, $\bH_{\bar\rd_A,\varphi}^0=0$ implies that $\bH_{\bar\rd_A,\varphi}^{+,0}=0$ and
$\bH_{\bar\rd_A,\varphi}^{0,0}=0$.  Similarly, $\bH_{\bar\rd_A,\varphi}^3=0$ implies that $\bH_{\bar\rd_A,\varphi}^{+,3}=0$ and
$\bH_{\bar\rd_A,\varphi}^{0,3}=0$.  Then equations \eqref{eq:CoIndexAndNullityEqualTruncatedEulerChar} imply that
\begin{subequations}
\label{eq:CoIndexAndNullityEqualEulerChar}
\begin{align}
\label{eq:CoIndexAndNullityEqualEulerCharPositiveWeight}
\la_{[A,\Phi]}^+(f)& = -2\Euler_\CC\left(\bar\rd_{A,\varphi}^{+,\bullet}\right),
\\
\label{eq:CoIndexAndNullityEqualEulerCharZeroWeight}
\la_{[A,\Phi]}^0(f)& = -2\Euler_\CC\left(\bar\rd_{A,\varphi}^{0,\bullet}\right).
\end{align}
\end{subequations}
Combining \eqref{eq:CoIndexAndNullityEqualEulerCharPositiveWeight} with Lemma \ref{lem:IndexOfRedStablePairsDefComplexNormalPlus} yields \eqref{eq:MorseCoIndexAtReduciblesOnKahler}.
Combining \eqref{eq:CoIndexAndNullityEqualEulerCharZeroWeight} with Lemma \ref{lem:IndexOfRedStablePairsDefComplexNormalZero} yields \eqref{eq:MorseNullIndexAtReduciblesOnKahler}.
This completes the proof of Theorem \ref{mainthm:MorseIndexAtReduciblesOnKahler}.
\end{proof}

By expressing the first Chern classes of $L_1$ and $L_2$ in terms of the characteristic classes of \spinc and \spinu structures, Theorem \ref{mainthm:MorseIndexAtReduciblesOnKahler} yields the

\begin{proof}[Proof of Corollary \ref{maincor:MorseIndexAtReduciblesOnKahlerWithSO3MonopoleCharacteristicClasses}]
Recall that $\fs_{\can}=(\rho_{\can},W_{\can})$ denotes the canonical \spinc structure over $X$ described in
Definition \ref{defn:Canonical_spinc_bundles}.  By Feehan and Leness \cite[Lemma 2.3, p. 64]{FL2a}, there is a complex, rank-two vector bundle $E$ over $X$ such that we can write the \spinu structure in the statement of the corollary as $\ft=(\rho_{\can},W_{\can},E)$.  By Kronheimer and Mrowka \cite[Proposition 1.1.1, p. 3, and discussion on p. 5]{KMBook}, there is a complex line bundle $L_1$ over $X$ such that we can write the \spinc structure in the statement of the corollary as $\fs=(\rho_{\can},W_{\can}\otimes L_1)$. By the hypotheses of Theorem \ref{mainthm:MorseIndexAtReduciblesOnKahler}, the bundle $E$ admits a decomposition as an orthogonal direct sum of Hermitian line bundles $E=L_1\oplus L_2$, where $L_2\cong L_1^*\otimes\det E$.
We will compute  $c_1(L_1)$ and $c_1(L_2)$ in terms of $c_1(\fs)$ and $c_1(\ft)$, the characteristic classes defined in \eqref{eq:DefineChernClassOfSpinc} and \eqref{eq:SpinUCharacteristics}, respectively.  The conclusion of Corollary \ref{maincor:MorseIndexAtReduciblesOnKahlerWithSO3MonopoleCharacteristicClasses} will then follow by applying these computations to the expression for virtual Morse--Bott index given in \eqref{eq:MorseIndexAtReduciblesOnKahler}.

The equalities $c_1(W_{\can}^+)=-c_1(K_X)$ and $c_1(K_X)=-c_1(X)$ (where $c_1(X)$ denotes the first Chern class of the holomorphic tangent bundle $T^{1,0}(X)$) from Morgan \cite[Section 7.1, p. 110]{MorganSWNotes} and Kotschick \cite[Fact 2.1]{KotschickSW} and $c_1(\ft)=c_1(W_{\can}^+)+c_1(E)$ by its definition in \eqref{eq:SpinUCharacteristics}
(where we note that the definition $V^+=W_{\can}^+\otimes E$ implies that
$\det V^+ \cong (\det W_{\can}^+)^{\otimes 2} \otimes (\det E)^{\otimes 2}$ and so $c_1(V^+)=2c_1(W_{\can}^+)+2c_1(E)$) and $c_1(E)=c_1(L_1)+c_1(L_2)$ yield the following identities:
\begin{subequations}
\label{eq:SpincForStableSplitting}
\begin{align}
\label{eq:SpincForStableSplitting1}
c_1(\fs)&=c_1(X)+2c_1(L_1),
\\
\label{eq:c1SpinuForStableSplitting}
c_1(\ft)&=c_1(X)+c_1(L_1)+c_1(L_2).
\end{align}
\end{subequations}
Subtracting \eqref{eq:c1SpinuForStableSplitting} from \eqref{eq:SpincForStableSplitting1} yields
\begin{equation}
\label{eq:SpincForStableSplitting2}
c_1(L_1)-c_1(L_2)=c_1(\fs)-c_1(\ft).
\end{equation}
We then compute
\begin{align*}
\la_{[A,\Phi]}^-(f)
&=
-2\chi_h(X) -\left( c_1(L_1) -c_1(L_2)\right)\cdot c_1(X) -\left( c_1(L_1)-c_1(L_2)\right)^2
\quad\text{(by Theorem \ref{mainthm:MorseIndexAtReduciblesOnKahler})}
\\
&=
-2\chi_h(X) -\left( c_1(\fs)-c_1(\ft)\right)\cdot c_1(X) -\left( c_1(\fs)-c_1(\ft)\right)^2
\quad\text{(by \eqref{eq:SpincForStableSplitting2}).}
\end{align*}
This completes the proof of equation \eqref{eq:MorseIndexAtReduciblesOnKahlerSpinNotationType1} and hence the proof of  Corollary \ref{maincor:MorseIndexAtReduciblesOnKahlerWithSO3MonopoleCharacteristicClasses}.
\end{proof}

\section{Checking the index computations}
\label{sec:CheckingIndexComputation}
As a test of the accuracy of the computations of Section \ref{sec:IndOfSubcomplex}
(in particular  Lemmas \ref{lem:IndexOfRedStablePairsDefComplexNormalMinus},
\ref{lem:IndexOfRedStablePairsDefComplexNormalPlus}, and \ref{lem:IndexOfRedStablePairsDefComplexNormalZero}), we now verify that the following equality holds:
\begin{equation}
  \label{eq:IndicesSumToNonAbelianMonopoleDimension}
\expdim\sM_\ft
=
\la_{[A,\Phi]}^0
+
\la_{[A,\Phi]}^+
+
\la_{[A,\Phi]}^-.
\end{equation}
The direct-sum decomposition \eqref{eq:WeightSpaceDecompOfDeformationCohomology} in Corollary \ref{cor:WeightOfS1ActionOnCohomology} of the elliptic complex for holomorphic pairs implies that
\begin{equation}
  \label{eq:EulerOfDecomposition}
\Euler_\CC(\bar\rd_{\bar\rd_A,\phi}^\bullet)
=
\Euler_\CC(\bar\rd_{\bar\rd_A,\phi}^{0,\bullet})
+
\Euler_\CC(\bar\rd_{\bar\rd_A,\phi}^{+,\bullet})
+
\Euler_\CC(\bar\rd_{\bar\rd_A,\phi}^{-,\bullet}).
\end{equation}
Combining \eqref{eq:EulerOfDecomposition} with the isomorphism \eqref{eq:Isomorphisms_of_non_AbelianHarmonics_with_HolomorphPairHarmonics}
and Proposition \ref{eq:IdentifyExpDimOfnonAbelianMonopoleModuliSpaceWithEulerCharacteristic} implies that
\eqref{eq:IndicesSumToNonAbelianMonopoleDimension} will hold if our computations of the indexes appearing in
\eqref{eq:IndicesSumToNonAbelianMonopoleDimension} are correct.
Throughout this computation, we shall assume that $(A,\Phi)$ is a split non-Abelian monopole on a \spinu structure $\ft=(\rho,W_{\can}\otimes E)$, where $E$ splits as a direct sum of Hermitian line bundles $E=L_1\oplus L_2$.

The expressions for the indices in \eqref{eq:MorseIndexAtReduciblesOnKahler}, \eqref{eq:MorseCoIndexAtReduciblesOnKahler}, and \eqref{eq:MorseNullIndexAtReduciblesOnKahler}
derived in Lemmas \ref{lem:IndexOfRedStablePairsDefComplexNormalMinus},
\ref{lem:IndexOfRedStablePairsDefComplexNormalPlus},
and \ref{lem:IndexOfRedStablePairsDefComplexNormalZero}, respectively, imply that
\begin{align*}
&\la_{[A,\Phi]}^0
+
\la_{[A,\Phi]}^+
+
\la_{[A,\Phi]}^-
\\
&\quad=
-2\chi_h(X) -\left(c_1(L_1)-c_1(L_2)\right)\cdot c_1(X) - \left( c_1(L_1)-c_1(L_2)\right)^2
\\
&\qquad
-\left( c_1(L_2)-c_1(L_1)\right)\cdot c_1(X) - \left( c_1(L_1)-c_1(L_2)\right)^2 +c_1(L_2)\cdot c_1(X) +c_1(L_2)^2
\\
&\qquad
+c_1(L_1)^2+c_1(L_1)\cdot c_1(X),
\end{align*}
and thus, after simplifying the last expression on the right-hand side above,
\begin{equation}
\label{eq:SumOfIndexCoIndexNullIndex}
\begin{aligned}
&\la_{[A,\Phi]}^0
+
\la_{[A,\Phi]}^+
+
\la_{[A,\Phi]}^-
\\
&\quad=
-2\chi_h(X) -2\left( c_1(L_1)-c_1(L_2)\right)^2+ c_1(L_1)^2+c_1(L_2)^2+\left( c_1(L_1)+c_1(L_2)\right)\cdot c_1(X).
\end{aligned}
\end{equation}
We now compute the left-hand-side of \eqref{eq:IndicesSumToNonAbelianMonopoleDimension}.  By \eqref{eq:Transv},
\begin{align*}
\expdim\sM_\ft
&=
d_a(\ft) + 2n_a(\ft)
\\
&=
-2p_1(\ft)-6\chi_h(X)
+\frac{1}{2}p_1(\ft) + \frac{1}{2}c_1(\ft)^2 -\frac{1}{2}\left( c_1(X)^2 -8\chi_h(X)\right),
\end{align*}
and thus, after simplifying the last expression on the right-hand side above,
\begin{equation}
\label{eq:ExpDim1}
\expdim\sM_\ft
=
-\frac{3}{2}p_1(\ft)+ \frac{1}{2}c_1(\ft)^2 -\frac{1}{2} c_1(X)^2-2\chi_h(X).
\end{equation}
Because $E=L_1\oplus L_2$, so $c_1(E)=c_1(L_1)+c_1(L_2)$ and $c_2(E)=c_1(L_1)c_1(L_2)$, we obtain
\begin{align*}
p_1(\ft)&= p_1(\su(E))\quad\text{(by \eqref{eq:SpinUCharacteristics} and \eqref{eq:SpinAssociatedBundles})}
\\
&= c_1(E)^2-4c_2(E)\quad\text{(by Donaldson and Kronheimer \cite[Equation (2.1.39), p. 42]{DK})}
\\
&=\left(c_1(L_1)+c_1(L_2)\right)^2-4c_1(L_1)c_1(L_2)
\\
&=\left(c_1(L_1)-c_1(L_2)\right)^2.
\end{align*}
Substituting $p_1(\ft)=(c_1(L_1)-c_1(L_2)^2$ and
the equality
$c_1(\ft)= c_1(X)+c_1(L_1)+c_1(L_2)$ from \eqref{eq:c1SpinuForStableSplitting} into \eqref{eq:ExpDim1} gives
\begin{align*}
\expdim\sM_\ft
&=
-\frac{3}{2}\left( c_1(L_1)-c_1(L_2)\right)^2
+
\frac{1}{2}c_1(X)^2+\left( c_1(L_1)+c_1(L_2)\right)\cdot c_1(X)
\\
&\qquad
+\frac{1}{2}\left(c_1(L_1)+c_1(L_2)\right)^2
-\frac{1}{2} c_1(X)^2-2\chi_h(X),
\end{align*}
and thus, after simplifying the last expression on the right-hand side above,
\begin{equation}
\label{eq:DimensionForNonAbelianModuliInCanonicalSplitting}
\begin{aligned}
\expdim\sM_\ft
&=
-2\chi_h(X) -\frac{3}{2}\left( c_1(L_1)-c_1(L_2)\right)^2+\frac{1}{2}\left(c_1(L_1)+c_1(L_2)\right)^2
\\
&\qquad+\left( c_1(L_1)+c_1(L_2)\right)\cdot c_1(X).
\end{aligned}
\end{equation}
Combining \eqref{eq:SumOfIndexCoIndexNullIndex} and \eqref{eq:DimensionForNonAbelianModuliInCanonicalSplitting} yields
\begin{align*}
&\expdim\sM_\ft - \left(\la_{[A,\Phi]}^0
+
\la_{[A,\Phi]}^+
+
\la_{[A,\Phi]}^- \right)
\\
&\quad =
\frac{1}{2}\left( c_1(L_1)-c_1(L_2)\right)^2
+
\frac{1}{2}\left( c_1(L_1)+c_1(L_2)\right)^2
-c_1(L_1)^2-c_1(L_2)^2
\\
&\quad=0.
\end{align*}
Hence, the computations of Lemmas \ref{lem:IndexOfRedStablePairsDefComplexNormalMinus},
\ref{lem:IndexOfRedStablePairsDefComplexNormalPlus},
and \ref{lem:IndexOfRedStablePairsDefComplexNormalZero} do satisfy \eqref{eq:IndicesSumToNonAbelianMonopoleDimension}, as expected.

\section{Restriction of the holomorphic-pair map to weight subspaces}
\label{sec:RestrictionOfHolomorphicPairsMapToWeightSubspaces}
As in the preceding sections, let $E$ be a complex rank-two vector bundle over a complex manifold and let $(\bar\rd_E,\varphi)$ be an $(0,1)$-pair on $E$ which is split with respect to a decomposition $E=L_1\oplus L_2$ as a direct sum of complex line bundles, $L_1$ and $L_2$, in the sense of Definition \ref{defn:Split_(0,1)-pair}, so $\varphi\in\Omega^0(L_1)$.
We write
\[
  (\rd_E,\varphi)
  =
  \left(\bar\rd_{L_1}\oplus\bar\rd_{L_2},\varphi\oplus 0\right),
\]
where $\bar\rd_{L_j}$ is an $(0,1)$-connection on $L_j$ for $j=1,2$. In Lemma \ref{lem:WeightSpaceDecompositionOfHolomorphicPairsDeformationComplex}, we showed that the linearization $\bar\rd_{A,\phi}^1$ of the holomorphic-pair map \eqref{eq:Holomorphic_pair_map} admitted a direct-sum decomposition with respect to the weight subspaces.  We now show how the image of the restriction of the holomorphic-pair map \eqref{eq:Holomorphic_pair_map} to each of these weight subspaces is contained in the corresponding subspace of $\sE_2$. (We explain the significance of this observation in Feehan \cite{Feehan_analytic_spaces}.)

Recall from
\eqref{eq:DefineSpacesOfAbbreviatedHolomorphicPairDefComplex} that we defined
\[
\sE_1=\Om^{0,1}(\fsl(E))\oplus \Om^0(E)
\quad\text{and}\quad
\sE_2=\Om^{0,2}(\fsl(E))\oplus \Om^{0,1}(E).
\]
By \eqref{eq:(0,1)PairsUnderlyingVectorSpace}, we have
\[
\sA^{0,1}(E)\times \Omega^0(E)=(\bar\rd_E,\varphi)+\left( \Om^{0,1}(E)\times \Omega^0(E)\right),
\]
and so we can consider the holomorphic-pair map $\fS$ in \eqref{eq:Holomorphic_pair_map},
\[
  (\rd_{E'},\varphi')\mapsto \fS(\rd_{E'},\varphi') = \left( F_{E'},\bar\rd_{E'}\varphi'\right),
\]
after identifying $(\bar\rd_E,\varphi)$ with the origin, as a map
\[
\fS_0:\sE_1\to\sE_2.
\]
More explicitly, we have
\begin{equation}
\label{eq:HolomorphicPairMapBetweenVectorSpaces}
\fS_0(\alpha,\sigma)
=
\left(
\bar\rd_E\alpha +\alpha\wedge\alpha ,
\bar\rd_E\si+\alpha\varphi+\alpha\si
\right),
\quad\text{for all } (\alpha,\sigma) \in \Om^{0,1}(\fsl(E))\oplus \Om^0(E).
\end{equation}
Consequently, the maps $\fS_0$ and $\fS$ are related by (see Kobayashi \cite[Equation (7.2.1), p. 223]{Kobayashi})
\begin{align*}
  \fS(\bar\partial_E+\alpha, \varphi+\sigma)
  &=
\left(
F_{\bar\partial_E+\alpha},(\bar\partial_E+\alpha)\left(\varphi+\si\right)
    \right)
    \\
&=
\left(
F_{\bar\partial_E},\bar\rd_E\varphi
\right)
+
\left(
\bar\rd_E\alpha +\alpha\wedge\alpha ,
\bar\rd_E\si+\alpha\varphi+\alpha\si
\right)
\\
&=\left(
F_{\bar\partial_E},\bar\rd_E\varphi
\right)
+
     \fS_0(\alpha,\sigma),
     \quad\text{for all } (\alpha,\sigma) \in \Om^{0,1}(\fsl(E))\oplus \Om^0(E).
\end{align*}
The weight-subspace decompositions of $\sE_1$ and $\sE_2$ in \eqref{eq:DefineSpacesOfAbbreviatedHolomorphicPairDefComplexWeightDecomp1} give decompositions of the domain and range of $\fS_0$ as
\[
\sE_k=\sE_k^0\oplus \sE_k^+ \oplus \sE_k^-,\quad\text{for } k=1,2,
\]
where, as in
\eqref{eq:SpacesOfWeightZeroSubcomplex},
\eqref{eq:SpacesOfWeight+Subcomplex}, and
\eqref{eq:SpacesOfWeight-Subcomplex}, we have
\begin{align*}
\sE_1^0&=\Om^{0,1}(\ubarCC)\oplus \Om^0(L_1),
\\
\sE_1^+&= \Om^{0,1}(L_2\otimes L_1^*)\oplus \Om^0(L_2),
\\
\sE_1^-&= \Om^{0,1}(L_1\otimes L_2^*),
\\
\sE_2^0&=\Om^{0,2}(\ubarCC)\oplus \Om^{0,1}(L_1),
\\
\sE_2^+&=\Om^{0,2}(L_2\otimes L_1^*)\oplus \Om^{0,1}(L_2),
\\
\sE_2^-&=\Om^{0,2}(L_1\otimes L_2^*).
\end{align*}
We claim that $\fS_0$ maps each of the subspaces $\sE_1^0$, $\sE_1^+$, and $\sE_1^-$ to the corresponding subspace of $\sE_2$.

\begin{lem}[Restriction of the holomorphic-pair map to weight subspaces]
\label{lem:WeightSpaceDecompositionOfHolomorphicPairsMap}
Let $E$ be a complex rank-two vector bundle over a complex manifold and let $(\bar\rd_E,\varphi)$ be
an $(0,1)$-pair
on $E$ which is split with respect to a decomposition $E=L_1\oplus L_2$ as a direct sum of complex line bundles, $L_1$ and $L_2$, in the sense of Definition \ref{defn:Split_(0,1)-pair}, so $\varphi\in\Omega^0(L_1)$. Let $\sE_k^0$, $\sE_k^+$, and $\sE_k^-$ be the weight subspaces of $\sE_k$, defined in \eqref{eq:DefineSpacesOfAbbreviatedHolomorphicPairDefComplexWeightDecomp1} for $k=1,2$.
Let $\fS_0:\sE_1\to\sE_2$ be the holomorphic pair map \eqref{eq:HolomorphicPairMapBetweenVectorSpaces} induced from the definition \eqref{eq:Holomorphic_pair_map} of the holomorphic pair map $\fS$ and choice of $(\bar\rd_E,\varphi)$ as the origin of the affine space domain. Then the following inclusions hold:
\begin{equation}
\label{eq:WeightSpaceDecompOfHolomorphicPairsMap}
\fS_0(\sE_1^0) \subset \sE_2^0,
\quad
\fS_0(\sE_1^+) \subset \sE_2^+,
\quad
\fS_0(\sE_1^-) \subset \sE_2^-.
\end{equation}
\end{lem}

\begin{proof}
Recalling that $\bar\rd_{E,\varphi}^1$ is the linearization \eqref{eq:d1StablePair} of the holomorphic pair map $\fS$ (and thus also $\fS_0$) at $(\rd_E,\phi)$, we may write
\begin{align*}
\fS_0(\alpha,\sigma)
&=
\left(
\bar\rd_E\alpha +\alpha\wedge\alpha ,
\bar\rd_E\si+\alpha\varphi+\alpha\si
\right)
\\
&=
\bar\rd_{E,\varphi}^1\left(\alpha,\si\right)
+
\left(\alpha\wedge\alpha,\alpha\si\right),
\quad\text{for all } (\alpha,\sigma) \in \Om^{0,1}(\fsl(E))\oplus \Om^0(E).     
\end{align*}
We denote the quadratic part of the preceding expression for $\fS_0(\alpha,\sigma)$ by
\[
  q(\alpha,\sigma)
  :=
\left(
\alpha\wedge\alpha ,
\alpha\si
\right),
\]
and thus we have
\begin{equation}
\label{eq:HolomorphicPairsMapDecompsitionIntoLinearAndQuadratic}
\fS_0(\alpha,\sigma)
=
\bar\rd_{E,\varphi}^1(\alpha,\sigma)
+
q(\alpha,\sigma),
\quad\text{for all } (\alpha,\sigma) \in \Om^{0,1}(\fsl(E))\oplus \Om^0(E).  
\end{equation}
By Lemma \ref{lem:WeightSpaceDecompositionOfHolomorphicPairsDeformationComplex},
the linear term $\bar\rd_{E,\varphi}^1$ in $\fS_0$ respects the weight-subspace decompositions of $\sE_1$ and $\sE_2$:
\begin{equation}
\label{eq:HolomorphicPairsMapLinearTermInclusions}
\bar\rd_{E,\varphi}^1\left( \sE_1^0\right)
\subset
\sE_2^0,\quad
\bar\rd_{E,\varphi}^1\left( \sE_1^+\right)
\subset
\sE_2^+,\quad
\bar\rd_{E,\varphi}^1\left( \sE_1^-\right)
\subset
\sE_2^-.
\end{equation}
Hence, the inclusions \eqref{eq:HolomorphicPairsMapDecompsitionIntoLinearAndQuadratic} for the linear term imply that the inclusions \eqref{eq:WeightSpaceDecompOfHolomorphicPairsMap} will follow from the following inclusions for the quadratic terms:
\begin{equation}
\label{eq:HolomorphicPairsMapQuadraticTermInclusions}
q(\sE_1^0) \subset \sE_2^0,
\quad
q(\sE_1^+) \subset \sE_2^+,
\quad
q(\sE_1^-) \subset \sE_2^-.
\end{equation}
We now verify the inclusions \eqref{eq:HolomorphicPairsMapQuadraticTermInclusions}. We may write
elements $(\alpha,\sigma) \in\Om^{0,1}(\fsl(E)) \oplus \si\in\Om^0(E)$ as
\[
\alpha=
\begin{pmatrix}
\zeta_\CC & \alpha_{12}
\\
\alpha_{21} & -\zeta_\CC
\end{pmatrix}
\in\Om^{0,1}(\fsl(E))
\quad\text{and}\quad
\si=
\begin{pmatrix}
\si_1 \\ \si_2
\end{pmatrix}
\in\Om^0(E),
\]
where
$\zeta_\CC\in\Om^{0,1}(\ubarCC)$, and $\alpha_{12}\in\Om^{0,1}(L_1\otimes L_2^*)$, and $\alpha_{21}\in\Om^{0,1}(L_2\otimes L_1^*)$, and $\si_i\in\Om^0(L_i)$. Thus,
\begin{align*}
\sE_1^0
&=
\left\{
\left(\begin{pmatrix}
\zeta_\CC & 0
\\
0 & -\zeta_\CC
\end{pmatrix},
\begin{pmatrix}
\si_1 \\ 0
\end{pmatrix}
\right)\in\sE_1
\right\},
\\
\sE_1^+&=
\left\{
\left(\begin{pmatrix}
0 & 0
\\
\alpha_{21} & 0
\end{pmatrix},
\begin{pmatrix}
0 \\ \si_2
\end{pmatrix}
\right)\in\sE_1
\right\},
\\
\sE_1^-&=
\left\{
\left(\begin{pmatrix}
0 & \alpha_{12}
\\
0 & 0
\end{pmatrix},
\begin{pmatrix}
0 \\ 0
\end{pmatrix}
\right)\in\sE_1
\right\}.
\end{align*}
Similarly, we may write elements $(\eta,\psi) \in\Om^{0,2}(\fsl(E)) \oplus \Om^{0,1}(E)$ as
\[
\eta=\begin{pmatrix}
\eta_\CC & \eta_{12}
\\
\eta_{21} & -\eta_\CC
\end{pmatrix}
\quad\text{and}\quad
\psi
=
\begin{pmatrix}
\psi_1 \\ \psi_2
\end{pmatrix},
\]
where $\eta_\CC\in\Om^{0,2}(\ubarCC)$, and $\eta_{12}\in\Om^{0,2}(L_1\otimes L_2^*)$, and $\eta_{21}\in\Om^{0,2}(L_2\otimes L_1^*)$, and $\psi_1\in\Om^{0,1}(L_1)$,
and $\psi_2\in\Om^{0,1}(L_2)$. Thus,
\begin{align*}
\sE_2^0
&=
\left\{
\left(\begin{pmatrix}
\eta_\CC & 0
\\
0 & -\eta_\CC
\end{pmatrix},
\begin{pmatrix}
\psi_1 \\ 0
\end{pmatrix}
\right)\in\sE_2
\right\},
\\
\sE_2^+&=
\left\{
\left(\begin{pmatrix}
0 & 0
\\
\eta_{21} & 0
\end{pmatrix},
\begin{pmatrix}
0 \\ \psi_2
\end{pmatrix}
\right)\in\sE_2
\right\},
\\
\sE_2^-&=
\left\{
\left(\begin{pmatrix}
0 & \eta_{12}
\\
0 & 0
\end{pmatrix},
\begin{pmatrix}
0 \\ 0
\end{pmatrix}
\right)\in\sE_2
\right\}.
\end{align*}
Therefore, we may compute the quadratic term $q$ as
\begin{align*}
q(\alpha,\sigma)&=
\left(
\alpha\wedge\alpha ,
\alpha\si
\right)
\\
&=
\left(
\begin{pmatrix}
\zeta_\CC & \alpha_{12}
\\
\alpha_{21} & -\zeta_\CC
\end{pmatrix}
\wedge
\begin{pmatrix}
\zeta_\CC & \alpha_{12}
\\
\alpha_{21} & -\zeta_\CC
\end{pmatrix},
\begin{pmatrix}
\zeta_\CC & \alpha_{12}
\\
\alpha_{21} & -\zeta_\CC
\end{pmatrix}
\begin{pmatrix}
\si_1
\\
\si_2
\end{pmatrix}
\right)
\\
&=
\left(
\begin{pmatrix}
\alpha_{12}\wedge\alpha_{21} & 2\zeta_\CC\wedge\alpha_{12}
\\
2\alpha_{21}\wedge\zeta_\CC & \alpha_{21}\wedge\alpha_{12}
\end{pmatrix},
\begin{pmatrix}
\zeta_\CC\si_1 +\alpha_{12}\si_2
\\
\alpha_{21}\si_1 -\zeta_\CC\si_2
\end{pmatrix}
\right),
\quad\text{for all } (\alpha,\sigma) \in \Om^{0,1}(\fsl(E))\oplus \Om^0(E).    
\end{align*}
For $(\zeta_\CC,\si_1)\in \sE_1^0$, we compute
\[
q(\zeta_\CC,\si_1)
=
\left(
\begin{pmatrix}
0 & 0 \\ 0 & 0
\end{pmatrix},
\begin{pmatrix}
\zeta_\CC\si_1
\\
0
\end{pmatrix}
\right)
\in
\sE_2^0.
\]
For $(\alpha_{21},\si_2)\in \sE_1^+$, we compute
\[
q(\alpha_{21},\si_2)
=
\left(
\begin{pmatrix}
0 & 0 \\ 0 & 0
\end{pmatrix},
\begin{pmatrix}
0
\\
0
\end{pmatrix}
\right)
\in \sE_2^+.
\]
Finally, for $(\alpha_{12},0)\in \sE_1^-$, we compute
\[
q(\alpha_{12},0)
=
\left(
\begin{pmatrix}
0 & 0 \\ 0 & 0
\end{pmatrix},
\begin{pmatrix}
0
\\
0
\end{pmatrix}
\right)
\in \sE_2^-.
\]
This completes the verification of the inclusions \eqref{eq:HolomorphicPairsMapQuadraticTermInclusions} and hence the proof of the lemma.
\end{proof}

\begin{rmk}[Non-triviality of the quadratic term on $\sE_1$]
As we can see from the proof of Lemma \ref{lem:WeightSpaceDecompositionOfHolomorphicPairsMap}, we do not have
\[
q(\sE_1^+\oplus \sE_1^-)\subset \sE_2^+\oplus \sE_2^-,
\]
as the diagonal terms $\alpha_{12}\wedge\alpha_{21}$  and $\alpha_{21}\wedge\alpha_{12}$ and the term $\alpha_{12}\si_2$ appearing in the expressions for $q$ indicate.
\end{rmk}

\chapter[Bubbling and virtual Morse--Bott indices for compactified moduli spaces]{Bubbling and virtual Morse--Bott indices for compactified moduli spaces of non-Abelian monopoles over four-manifolds of Seiberg--Witten simple type}
\label{chap:Bubbling}
The remaining steps in our program to prove the Bogomolov--Miyaoka--Yau inequality for four-manifolds of Seiberg--Witten simple type (Conjecture \ref{conj:BMY_Seiberg-Witten}) are briefly summarized below:
\begin{enumerate}
\item\label{item:Step_Gieseker} Use the paradigm described in Section \ref{sec:Morse_theory_existence_anti-self-dual_connections} to prove the Bogomolov--Miyaoka--Yau inequality for compact, complex surfaces of general type (originally due to Miyaoka \cite{Miyaoka_1977} and Yau \cite{YauPNAS, Yau}) via the results of the present monograph, our companion monograph \cite{Feehan_analytic_spaces} on Morse--Bott theory for analytic spaces, and our ongoing work on a generalization of the Gieseker compactification \cite{Gieseker_1977} for the moduli space of stable holomorphic bundles to stable holomorphic pairs over smooth, complex projective surfaces. 

\item\label{item:Step_SW_simple_type_top_level} Extend the main results of our current monograph by replacing the hypothesis that the four-manifold is a compact, complex K\"ahler surface by the hypothesis that it is a standard four-manifold of Seiberg--Witten simple type with a non-vanishing Seiberg--Witten invariant.

\item\label{item:Step_SW_simple_type_all_levels} Extend the main results of Step \ref{item:Step_Gieseker} using an analytic version of the Gieseker compactification to replace a hypothesis that the four-manifold is a smooth, complex projective surface by the hypothesis that it is a standard four-manifold of Seiberg--Witten simple type with at least one non-vanishing Seiberg--Witten invariant.  
\end{enumerate}

Seiberg--Witten points in lower levels of the Uhlenbeck compactification $\bar\sM_\ft$ pose additional difficulties. Our results in Section \ref{sec:Morse_theory_moduli_space_PU2_monopoles_Kaehler} do not apply to them and so we outline approaches in this chapter that should allow extensions that do yield analogues of the results in Section \ref{sec:Morse_theory_moduli_space_PU2_monopoles_Kaehler} for Seiberg--Witten points in all levels of $\bar\sM_\ft$, first in the case where $X$ is a smooth, complex projective surface and then more generally when $X$ is a standard four-manifold of Seiberg--Witten simple type with a non-zero Seiberg--Witten invariant as in Conjecture \ref{conj:BMY_Seiberg-Witten}.

\section{Computation of virtual Morse--Bott indices for Seiberg--Witten points at the boundaries of compactifications of the moduli space of non-Abelian monopoles}
\label{sec:Computation_virtual_Morse_indices_lower-level_Seiberg-Witten_points}
Let $\ft$ be a \spinu structure on a standard four-manifold $X$, where $\ft$ and $X$ satisfy the hypotheses of Theorem \ref{mainthm:ExistenceOfSpinuForFlow}.  By Theorem \ref{mainthm:ExistenceOfSpinuForFlow}, $\sM_\ft^{*,0}$ is non-empty and we can start the Hitchin gradient flow $[A(t),\Phi(t)]$ from an initial point in $\sM_\ft^{*,0}$. There are then three possible outcomes:

First, the flow may converge to a possibly ideal zero-section point $[A_\infty,0,\bx] \in \bar\sM_\ft$ and we will obtain a projectively anti-self-dual  connection $A_\infty$ on a bundle $E_\ell$, with $c_1(E_\ell)=c_1(E)$ and $c_2(E_\ell)=c_2(E)-\ell$ for some $l\geq 0$, that obeys \eqref{eq:p1_lower_bound} since $p_1(\su(E_\ell))=p_1(\su(E))+4\ell\geq p_1(\su(E))$. (This will also prove that $\bar M_\kappa^w(g)$ is non-empty.)

Second, the flow may converge to a possibly ideal Seiberg--Witten point $[A_\infty,\Phi_\infty,\bx]$ in the Uhlenbeck (or \emph{analytical}) compactification $\bar\sM_\ft^A = \bar\sM_\ft$, where $[A_\infty,\Phi_\infty]$ is a Seiberg--Witten point in the ambient non-Abelian monopole moduli space $\sM_{\ft(\ell)}$ for some integer $\ell\geq 0$. If $\sM_{\ft(\ell)}$ has positive expected dimension, then by our discussion in Section \ref{subsec:Feasibility_SO(3)-monopole_cobordism_method}, the moduli space $\sM_{\ft(\ell)}^{*,0}$ will be non-empty and we can restart the downward Hitchin gradient flow from an initial point in $\sM_{\ft(\ell)}^{*,0}$ since $[A_\infty,\Phi_\infty]$ will again have positive virtual
Morse index by our discussions in Section \ref{sec:Morse_theory_moduli_space_PU2_monopoles_Kaehler}.

Third, $\sM_{\ft(\ell)}$ may have negative expected dimension and so $M_\fs$ is an \emph{isolated} subset of $\sM_{\ft(\ell)}$ since $\sM_{\ft(\ell)}^{*,0}$ is empty by \cite{FeehanGenericMetric, TelemanGenericMetric}
and $M_\fs$ is disjoint from the  moduli subspace $\bar M_{\kappa-\ell}^w(g) \subset \bar\sM_{\ft(\ell)}$ of anti-self-dual connections. A careful study of this most difficult case requires extensions from anti-self-dual connections to non-Abelian monopoles of our bubble-tree analysis
\cite{FeehanGeometry, Feehan_yang_mills_gradient_flow_v4} and
gluing theory \cite{FL9}, and corresponding local Kuranishi model for an open neighborhood of
$[A_\infty,\Phi_\infty,\bx]$ in the analytical compactification $\bar\sM_\ft^A$  of the moduli space of non-Abelian monopoles; this should allow us to compute the virtual Morse index of $[A_\infty,\Phi_\infty,\bx]$ and examine when an open neighborhood of the unstable subvariety in $\bar\sM_\ft^A$ through
$[A_\infty,\Phi_\infty,\bx]$ meets $\sM_\ft^{*,0}$. Increasing $p_1(\su(E))$ above the lower bound in \eqref{eq:p1_lower_bound} increases the tendency of isolated Seiberg--Witten moduli spaces in this third case to occur.

\section{Gieseker compactification of the moduli space of stable pairs over a complex K\"ahler surface}
\label{sec:GiesekerCompactificationOfStablePairs}
Because of the subtlety of the analysis required to describe unstable manifolds on neighborhoods of lower level Seiberg--Witten moduli subspaces in the Uhlenbeck compactification of the moduli space of non-Abelian monopoles, as described in Section \ref{sec:Computation_virtual_Morse_indices_lower-level_Seiberg-Witten_points}, it is helpful to consider a further specialization of our program from the setting of Conjecture \ref{conj:BMY_Seiberg-Witten}, where $X$ is allowed to be a standard four-manifold of Seiberg--Witten simple type with a non-zero Seiberg--Witten invariant, to the case where $X$ is a smooth, complex projective surface. Of course, if $X$ is further assumed to be a compact, smooth complex surface of general type, then the Bogomolov--Miyaoka--Yau inequality \eqref{eq:BMY} holds by Miyaoka \cite[Theorem 4, p. 234]{Miyaoka_1977} and Yau \cite{YauPNAS, Yau}, so this specialization provides way to partially validate our strategy to prove Conjecture \ref{conj:BMY_Seiberg-Witten} in full generality.

When $X$ is a smooth, complex projective surface, we may replace the role of the Uhlenbeck compactification by that of a compactification of the open moduli subspace of slope stable holomorphic pairs $(\sE,\varphi)$ in the complex, projective moduli space of pairs of coherent sheaves and sections that are semistable in the sense of Gieseker \cite{Gieseker_1977} and Maruyama \cite{Maruyama_1977, Maruyama_1978}. This framework should allow us to analyze the unstable manifolds on open neighborhoods of points that are strictly semistable in the sense of Gieseker and Maruyama using an approach that is broadly similar to the one used in this monograph for open neighborhoods of points that are slope polystable. This is work in progress with Wentworth \cite{Feehan_Leness_Wentworth_virtual_morse_theory_stable_pairs_bmy_kaehler}.

In this setting, a pair, or rather a triple, $(\sE,\varphi,\eps)$ comprises a coherent sheaf $\sE$ with fixed Hilbert polynomial, a homomorphism $\eps:\det\sE\to\sL$ (called an \emph{orientation}) from the determinant of $\sE$ to a fixed rank one sheaf $\sL$, and a homomorphism $\varphi:\sO_X \to \sE$. Such pairs are required to satisfy semistability conditions inspired by those of Gieseker and Maruyama for coherent sheaves. Our construction in \cite{Feehan_Leness_Wentworth_virtual_morse_theory_stable_pairs_bmy_kaehler} is similar to that of Okonek, Schmitt, and Teleman \cite{OTMasterPairs}, except that rather than consider \emph{framings} $\phi:\sE \to \sO_X$ in the terminology of Huybrechts and Lehn \cite{Huybrechts_Lehn_1995jag, Huybrechts_Lehn_1995ijm}, we consider \emph{co-framings} $\varphi:\sO_X \to \sE$ in the terminology of Flenner and L\"ubke \cite{Flenner_Lubke_2002}. Our construction is also similar to that of Lin \cite{LinThesis, Lin_2018} and Wandel \cite{Wandel_2015}, who consider co-framings rather than framings, but who only construct the Gieseker moduli space of semistable, \emph{unoriented} pairs $(\sE, \varphi)$. Thus, our approach combines features of that of Okonek, Schmitt, and Teleman \cite{OTMasterPairs} and of Lin \cite{LinThesis, Lin_2018} and Wandel \cite{Wandel_2015}, while completing many details omitted by Dowker \cite{DowkerThesis} in his Ph.D. thesis.

\section{Extension to the case of standard four-manifolds of Seiberg--Witten simple type}
As noted in \cite{KMThom}, smooth four-manifolds of Seiberg--Witten simple type are \emph{almost complex} and so $(X,g,J)$ is \emph{almost Hermitian}, with almost complex structure $J$ having a Nijenhuis tensor $N$ that may be non-zero and fundamental $2$-form $\omega=g(\cdot,J\cdot)$ that is non-degenerate but need not be closed. In future work, we hope to generalize our calculations in Section \ref{sec:Morse_theory_moduli_space_PU2_monopoles_Kaehler} for virtual Morse--Bott indices on moduli spaces $\sM_\ft$ of non-Abelian monopoles over complex K\"ahler surfaces to four-manifolds of Seiberg--Witten simple type with a strategy inspired by
\begin{inparaenum}[\itshape i\upshape)]
\item use of extrinsic perturbations of the Seiberg--Witten monopole equations to extend Witten's calculation of Seiberg--Witten invariants for complex K\"ahler surfaces to \emph{symplectic four-manifolds} by Taubes \cite{TauSymp, TauSympMore} and to complex \emph{non-K\"ahler} surfaces by Biquard \cite{Biquard_1998},
\item use of approximations to extend results for complex K\"ahler manifolds to \emph{symplectic manifolds} by Donaldson \cite{DonSympAlmostCx}, and
\item existence of almost complex structures with $L^p$-small Nijenhuis tensors, following Evans \cite{Evans_2012} for the case $p=2$ and Fernandez, Shin, and Wilson \cite{Fernandez_Shin_Wilson_2021} for the conjectural case $p=\infty$.
\end{inparaenum}

\section[Analytic compactification of the moduli space of non-Abelian monopoles]{Analytic compactification of the moduli space of non-Abelian monopoles over a smooth four-manifold}
As discussed in Section \ref{sec:Computation_virtual_Morse_indices_lower-level_Seiberg-Witten_points},
if a lower-level limit point is not an isolated point of $\sM_{\ft(\ell)}$ and is not a local minimum, one can
restart the gradient flow and continue it within $\sM_{\ft(\ell)}$.  However, if a lower-level limit point is isolated or a local minimum, this procedure is not available and a closer analysis of the gradient flow in an open neighborhood of such potential limit points is required.

To describe, at least approximately, the gradient flow in the neighborhood of a lower-level point, we may use the gluing maps of Feehan and Leness \cite{FL3,FL9}. For further details of this construction, see Feehan and Leness \cite[Section 3.1]{FL3}, \cite[Chapter 6]{FL5}, or (for an analogous construction for the moduli space of anti-self-dual connections) Friedman and Morgan \cite[Section 3.4]{FrM}. As described in \cite{FLLevelOne,FL5}, these maps consist of a smoothly-stratified topological space $\sU$ which we refer to as a \emph{virtual neighborhood}, a vector bundle $\sV\to\sU$, and an obstruction section, $\bchi:\sU\to\sV$. The gluing map $\bga$ gives a smoothly-stratified, $S^1$-equivariant embedding of $\sU$ into $\bar\sC_\ft$ such that the restriction of $\bga$ to the zero-locus $\bchi^{-1}(0)$ gives a smoothly-stratified homeomorphism of $\bchi^{-1}(0)$ onto an open neighborhood of a lower-level point in $\bar\sM_\ft$.

Let $([A_0,\Phi_0],\bx)\in\sM_{\ft(\ell)}\times\Sym^\ell(X)$ be a point with $\ell \geq 1$. We now describe an open neighborhood of $([A_0,\Phi_0],\bx)$ in the virtual neighborhood $\sU$. Let $\sO_{A_0,\Phi_0}$ be an open neighborhood of $[A_0,\Phi_0]$ in $\sM_{\ft(\ell)}$. Assume that $\bx\in\Sym^\ell(X)$ is given by the distinct points $\{x_1,\dots,x_r\}$ where $x_i$ has multiplicity $\ka_i$. Then an open neighborhood of $([A_0,\Phi_0],\bx)$ in the virtual neighborhood $\sU$ is given by
\begin{equation}
\label{eq:GluingNgh}
\sU_{A_0,\Phi_0,\bx}
:=
\left.\left(
\sO_{A_0,\Phi_0}
\times
\prod_{i=1}^r \sO_i \times \barM^{s,\diamond}_{\ka}\times (0,\eps_i)
\right)\right/S^1,
\end{equation}
where $\sO_i\subset T_{x_i}X$ is an open neighborhood of the origin, and $\barM^{s,\diamond}_{\ka}$ is the Uhlenbeck compactification of the mass- and scale-centered, framed, charge $\ka_i$ instantons on $S^4$, and the intervals $[0,\eps_i)$ represent the `scales' of the instantons being glued in at the points $x_i$. As a scale approaches zero, the pair approaches a lower level. The $S^1$ quotient in \eqref{eq:GluingNgh} is given by $S^1$ acting diagonally by the action of the stabilizer of $(A_0,\Phi_0)$ on $\sO_{A_0,\Phi_0}$ and on the frames in each factor of $\barM^{s,\diamond}_{\ka}$.

The gradient of the Hamiltonian function is given by $-\bJ\xi$, where $\bJ$ is an almost complex structure and $\xi$ is the vector field generated by the $S^1$ action on $\sC_\ft$. The $S^1$-equivariance of the gluing map gives an exact description of $\xi$ on the space $\sU_{A_0,\Phi_0,\bx}$.  We can only give an approximate description of the action of the almost complex structure $\bJ$ on $\sU_{A_0,\Phi_0,\bx}$ although careful analysis may yield an asymptotic description. This description leads us to conjecture that the gradient of the Hamiltonian function on the virtual neighborhood \eqref{eq:GluingNgh} is approximately orthogonal to the directions given by the scale parameters in $[0,\eps_i)$. The gradient flow on the virtual neighborhood $\sU_{A_0,\Phi_0,\bx}$ would therefore be parallel to the lower strata in the sense that the scale parameters would be unchanged under the flow.  

The analysis of the flow on $\bar\sM_\ft\cap \sU_{A_0,\Phi_0,\bx}$, in particular a description of the
stable and unstable manifolds, would require an asymptotic (as the scales go to zero) analysis of the behavior of $\bJ X$ that takes into account the factors of $\sO_i \times \barM^{s,\diamond}_{\ka}\times (0,\eps_i)$ in \eqref{eq:GluingNgh}. In Feehan and Leness \cite{FL9}, we construct gluing maps whose virtual neighborhoods are manifolds with corners; the smooth structure so provided should prove a useful tool in this analysis.

\backmatter

%
%

\bibliography{/Users/pfeehan/Dropbox/LATEX/Bibinputs/master,/Users/pfeehan/Dropbox/LATEX/Bibinputs/mfpde}
\bibliographystyle{amsplain-nodash}

%
%

\chapter*{Notation Index}

\begin{tabular}{ r p{5in} }
  $(A,\Phi)$
  & spin${}^u$ pair;
  page \pageref{page:SpinuPairs}
  \\
  $[A,\Phi]$
  & gauge equivalence class of spin${}^u$ pair;
  page \pageref{page:GaugeEquivClassOfSpinuPair}
  \\
  $(A,\varphi)$
  & unitary pair;
  page \pageref{page:Unitary_pairs}
  \\
  $[A,\varphi]$
  & gauge equivalence class of unitary pair;
  page \pageref{page:GaugeEquivClassOfUnitaryPair}
  \\
  $\sA(E,h)$
  & Affine space of unitary connections on a Hermitian vector bundle $(E,h)$ over a smooth manifold,
  possibly inducing a fixed unitary connection on the Hermitian line bundle $\det(E,h)$;
  Section \ref{sec:SpinuPairsQuotientSpace} prior to Equation \eqref{eq:SpinUConfiguration}, page
   \pageref{page:UnitaryConnectionsOverManifoldOfDim_d}
  \\
  $\sA^{0,1}(E)$
  & Affine space of $(0,1)$-connections on a complex vector bundle $E$ over a complex manifold, 
   possibly inducing a fixed $(0,1)$-connection on the complex line bundle $\det E$;
    Equation \eqref{eq:Affine_space_01-connections}
  \\
  $\sB(E,h)$
  &Quotient space of unitary connections on a Hermitian vector bundle $E$ over a smooth manifold, inducing a fixed connection on the Hermitian line bundle $\det E$; page \pageref{page:Quotient_space_unitary_connections}
  \\
  $\sB^{0,1}(E)$
  &Quotient space of $(0,1)$-connections on a complex vector bundle $E$ over a complex manifold, 
   possibly inducing a fixed $(0,1)$-connection on the complex line bundle $\det E$; page \pageref{page:Quotient_space_(0,1)_connections}
  \\
  $B(x),\widehat B(X)$
  & Basic classes and expanded basic classes;
  Equations \eqref{eq:SetOfBasicClasses} and \eqref{eq:NonEmptyReducibleMonopoles}
  \\
  $C^\infty(E)$
  & Vector space of smooth sections of a smooth vector bundle $E$ over a smooth manifold
  \\
  $\Coh(\sE)$
  & Set of coherent subsheaves $\sF$ of $\sE$ with $0<\rank\sF<\rank\sE$;
  page \pageref{page:ProperCoherentSubsheaves}
  \\
  $\Coh_\varphi(\sE)$
  & Set of coherent subsheaves $\sF$ of $\sE$ with $0<\rank\sF<\rank\sE$ that contain a section $\varphi$ of $\sE$;
  page \pageref{page:ProperCoherentSubsheavesContainingSection}
  \\
  $\Crit f$
  & Critical set of $f$;
  Definition \ref{maindefn:Morse-Bott_function}
  \\
  $\sC(E,h)$
  & Quotient space of unitary pairs on a Hermitian vector bundle $E$ over a smooth manifold, inducing a fixed connection on the Hermitian line bundle $\det E$; Equation \eqref{eq:ConfigurationSpaceForProjectiveVortices}
  \\
  $\sC^{**}(E,h)$
  & Subspace of $\sC(E,h)$ given by unitary pairs with trivial stabilizer; Equation \eqref{eq:ConfigurationSpaceForProjectiveVortices**}
  \\
  $\sC^0(E,h)$
  & Subspace of $\sC(E,h)$ given by non-zero-section unitary pairs; Equation \eqref{eq:ConfigurationSpaceForProjectiveVorticesnon-zero-section}
  \\
  $\sC^{0,1}(E)$
  & Quotient space of $(0,1)$-pairs on a complex vector bundle $E$ over a complex manifold, inducing a fixed $(0,1)$-connection on the complex line bundle $\det E$;
  Equation \eqref{eq:(0,1)PairQuotientSpaces}
  \\
  $\sC_\fs$
  & Quotient space of pairs for a spin${}^c$ structure $\fs$;
  Equation \eqref{eq:SpincConfig}
  \\
  $\sC_\fs^0$
  & Subspace of $\sC_\fs$ given by non-zero-section \spinc pairs;
  Equation \eqref{eq:SpincConfigNonZero}
  \\
  $\sC_\ft$
  & Quotient space of pairs for a spin${}^u$ structure $\ft$;
    Equation \eqref{eq:SpinUConfiguration}
  \\
  $\sC_\ft^{**}$
  & Subspace of $\sC_\ft$ given by \spinu pairs with trivial stabilizer;
    Equation \eqref{eq:eq:SpinuQuotientSpaceSubspacese_trivial_stabilizer}
  \\
  $\sC_\ft^*$
  & Subspace of $\sC_\ft$ given by non-split \spinu pairs;
    Equation \eqref{eq:eq:SpinuQuotientSpaceSubspacese_irreducible}
  \\
  $\sC^0_\ft$
  & Subspace of $\sC_\ft$ given by non-zero-section \spinu pairs;
    Equation \eqref{eq:eq:SpinuQuotientSpaceSubspaces_non-zero-section}
  \\
  $\sC_\ft^{*,0}$
  & Subspace of $\sC_\ft$ given by non-split, non-zero-section \spinu pairs;
  Equation \eqref{eq:eq:SpinuQuotientSpaceSubspaces_irredicible_non-zero-section} 
  \\
  $c_i(E)$
  & $i$-th Chern class of a complex vector bundle $E$;
  \\
  $C_n$
  & Cyclic group of order $n$;
  page \pageref{page:CyclicGroup}
  \\
  $c_i(X)$
  & Chern classes of the holomorphic tangent bundle $T^{1,0}(X)$ of an almost complex manifold $X$
  \\
  $c_1(\fs)$
  & First Chern class of  spin${}^c$ structure $\fs$;
  Equation \eqref{eq:DefineChernClassOfSpinc}
  \\
  $c_1(\ft)$
  & First Chern class of  spin${}^u$ structure $\ft$;
  Equation \eqref{eq:SpinUCharacteristics}
  \\
  $c_1(X)^2$
  & First Chern class squared of $X$ (homotopy invariant);
  Equation \eqref{eq:DefineTopCharNumbersOf4Manifold}
  \\
  $\chi_h(X)$
  & Holomorphic Euler characteristic of $X$ (homotopy invariant);
  Equation \eqref{eq:DefineTopCharNumbersOf4Manifold}
  \\
  $d_a(\ft)$
  & Expected dimension of moduli space of anti-self-dual connections;
  Equation \eqref{eq:Transv}
  \\
  $d(\fs)$
  & Expected dimension of moduli space of Seiberg--Witten monopoles;
  Equation \eqref{eq:DimSW}
  \\
  $D_A$, $D_B$
  & Dirac operators of spin${}^u$ and spin${}^c$ connections $A$ and $B$;
  pages \pageref{Dirac_operator} and \pageref{page:SpincDiracOperator}
  \\
  $\bar\partial_E$
  & $(0,1)$-connection;
  Equations \eqref{eq:Donaldson_Kronheimer_2-1-48}
  \\
  $[\bar\partial_E]$
  & Gauge equivalence class of $(0,1)$-connection;
  page \pageref{page:EquivClassOf(0,1)Connection}
  \\
  $(\bar\partial_E,\varphi)$
  & $(0,1)$-pair;
  Definition \ref{defn:Simple_pair}
  \\
  $[\bar\partial_E,\varphi]$
  & Gauge equivalence class of $(0,1)$-pair;
  page \pageref{page:EquivClassOf(0,1)Pair}
 \end{tabular}

\begin{tabular}{ r p{5in} }
  $d_{A,\Phi}^k$
  & Differential operators in elliptic complex of non-Abelian monopoles;
  Equations \eqref{eq:d_APhi^0} and \eqref{eq:d1OfSO3MonopoleComplex}
  \\
  $d_{A,\varphi}^k$
  & Differential operators in elliptic complex of projective vortices;
  Equations \eqref{eq:d0_projective_vortex_elliptic_deformation_complex},
  \eqref{eq:d1_projective_vortex_elliptic_deformation_complex},
  \eqref{eq:d2_projective_vortex}, and \eqref{eq:dk_projective_vortex}
\\
  $\bar\partial_{A,(\varphi,\psi)}^k$
  & Differential operators in elliptic complex of pre-holomorphic pairs;
  Equations \eqref{eq:d0StableComplex} and \eqref{eq:d1StableComplex}
  \\
  $\bar\partial_{E,\varphi}^k$
  & Differential operators in elliptic complex of holomorphic pairs;
  Equation \eqref{eq:dStablePair}
  \\
  $\hat d_A^+$
  & Linearization of Hermitian--Einstein equation;
  Equation \eqref{eq:Linearization_HE_equation_codomain_Omega0+Omega02}
  \\
  $\deg E$
  & Degree of vector bundle $E$ with respect to K\"ahler form $\om$;
  Equation \eqref{eq:DegreeIntegral}
  \\
  $\sE$
  & Holomorphic vector bundle or its sheaf of holomorphic sections;
  page \pageref{page:SheafOfHolomorphicSections}
  \\
  $e(X)$
  & Euler number of $X$
  \\
  $\Euler_\KK(d_\bullet)$
  & Index of an elliptic complex;
  Equation \eqref{eq:Index_elliptic_complex}
  \\
  $\fg_\ft$
  & Vector bundle of tracefree, skew-Hermitian endomorphisms associated to spin${}^u$ structure $\ft$;
  Equation \eqref{eq:SpinAssociatedBundles}
  \\
  $\bg(a_1,a_2)$
  & Weak ($L^2$) Riemannian metric on affine space of unitary connections $\sA(E,h)$;
  Equation \eqref{eq:Kobayashi_7-6-21}
  \\
  $\bg((a_1,\phi_1), (a_2,\phi_2))$
  & Weak ($L^2$) Riemannian metric on affine space of unitary pairs $\sA(E,h)\times W^{1,p}(E)$;
  Equation \eqref{eq:L2_metric_affine_space_pairs_unitary_connections_and_sections}
  \\
  $\GL(E)$
  & Principle bundle of frames for a vector bundle $E$;
  Equation \eqref{eq:AutomorphismBundles}
  \\
  $\grad_g f$
  & Gradient vector field of function $f$ with respect to Riemannian metric $g$;
  Equation \eqref{eq:DefineGradient}
  \\
  $\bh(a_1,a_2)$
  & Weak ($L^2$) Hermitian metric on affine space of unitary connections $\sA(E,h)$;
    Equation \eqref{eq:Kobayashi_7-6-20}
  \\
  $\bh((\alpha_1,\phi_1),(\alpha_2,\phi_2))$
  & Weak ($L^2$) Hermitian metric on affine space of unitary pairs $\sA(E,h)\times W^{1,p}(E)$;
    Equation \eqref{eq:Kobayashi_7-6-20_pairs}    
  \\
  $\Hess_g f(p)$
  & Hessian operator of function $f$ with respect to Riemannian metric $g$;
  Equation \eqref{eq:Hessian_operator}
  \\
  $\hess f(p)$
  & Hessian bilinear form of $f$;
    Equation \eqref{eq:Hessian_bilinear_form}
  \\
   $\bH_A^k$
  & Harmonic space for the elliptic complex of anti-self-dual connections on a four-manifold
  or projectively Hermitian--Einstein  connections on a K\"ahler manifold
  Equations \eqref{eq:H_Abullet} and \eqref{eq:HE_equation_bHAbullet}
  \\
  $\bH_{\bar\partial_E}^k$
  & Harmonic space for the elliptic complex of holomorphic structures;
  Equations \eqref{eq:H_dbar_A^0bullet} for $n=2$ and \eqref{eq:bHdbarA0bullet} for arbitrary $n$
  \\
  $\bH_{A,\Phi}^k$
  & Harmonic space of the elliptic complex of non-Abelian monopoles;
  Equation \eqref{eq:H_APhi^bullet}
    \\
  $\bH_{A,\phi}^k$
  & Harmonic space of the elliptic complex of projective vortices;
  Equation \eqref{eq:H_Avarphi^bullet}
  \\
  $\bH_{\bar\partial_A,(\varphi,\psi)}^k$
  & Harmonic space of the elliptic complex of pre-holomorphic pairs;
  Equation \eqref{eq:H_dbar_APhi^0bullet}
    \\
  $\bH_{\bar\partial_E,\varphi}^k$
  & Harmonic space of the elliptic complex of holomorphic pairs;
  Equation \eqref{eq:H_dbar_Avarphi^0bullet}
  \\
  $\bH_{\bar\partial_A,\varphi}^{0,k}$
  & Harmonic space of the zero-weight subcomplex of the elliptic complex of holomorphic pairs;
  Equation \eqref{eq:HarmonicSpacesForZeroWeight}
  \\
  $\bH_{\bar\partial_A,\varphi}^{+,k}$
  & Harmonic space of the positive-weight subcomplex of the elliptic complex of holomorphic pairs;
  Equation \eqref{eq:HarmonicSpacesForPositiveWeight}
  \\
  $\bH_{\bar\partial_A,\varphi}^{-,k}$
  & Harmonic space of the negative-weight subcomplex of the elliptic complex of holomorphic pairs;
  Equation \eqref{eq:HarmonicSpacesForNegativeWeight}
  \\
  $\iota_{\fs,\ft}$
  & Embedding of quotient spaces $\sC_\fs$ into $\sC_\ft$;
  Equation \eqref{eq:DefnOfIotaOnQuotient}
  \\
  $\tilde\iota_{\fs,\ft}$
  & Embedding of affine space $\sA(L,h_L)\times W^{1,p}(W^+\otimes L_1)$ into $\sA(E,h)\times W^{1,p}(W^+\otimes E)$;
  Equation \eqref{eq:DefnOfIota}
  \\
  $\sI\!\!\sM_{\ft}$
  & Moduli space of ideal non-Abelian monopoles;
  Equation \eqref{eq:idealmonopoles}
  \\
  $\bJ(a)$
  & Almost complex structure on affine space of unitary connections $\sA(E,h)$; Equation \eqref{eq:Almost_complex_structure_affine_space_pairs_unitary_connections_and_sections}
  \\
  $\bJ(a,\phi)$
  & Almost complex structure on affine space of unitary pairs $\sA(E,h)\times W^{1,p}(E)$;
    Equation \eqref{eq:Almost_complex_structure_affine_space_pairs_unitary_connections_and_sections}
\end{tabular}

\begin{tabular}{ r p{5in} }
  $K(A)$
  & Kuranishi model associated to a Hermitian--Einstein connection, $A$;
    page \pageref{page:Kuranishi_model_K(A)}
  \\
  $\sK(A,\varphi)$
  & Kuranishi model associated to a projective vortex, $(A,\varphi)$;
    page \pageref{defn:Friedman_Morgan_4-1-15_projective_vortex}
  \\
  $\cK(\bar\partial_E)$
  & Kuranishi model associated to a holomorphic structure, $\bar\partial_E$;
    page \pageref{page:Kuranishi_model_cKdbarE}
  \\
  $\fK(\bar\partial_E,\varphi)$
  & Kuranishi model associated to a holomorphic pair, $(\bar\partial_E,\varphi)$;
    page \pageref{defn:Friedman_Morgan_4-1-15_holomorphic_pairs}
  \\
  $K_X$
  & Canonical line bundle of an almost complex manifold;
  Equation \eqref{eq:DefineCanonicalLineBundle}
  \\
  $\kappa(E)$
  & Topological charge of the vector bundle $E$;
  page \pageref{page:kaChargeOfBundle}
  \\
  $\lambda_p^-(f),\lambda_p^+(f),\lambda_p^0(f)$
  & Virtual Morse--Bott index, co-index, and nullity;
  Equation \eqref{eq:Virtual_Morse-Bott_signature}
  \\
  $\lambda^-(\ft,\fs)$
  & Formal Morse--Bott index associated to spin${}^u$ structure $\ft$ and spin${}^c$ structure $\fs$;
  Equation \eqref{eq:FormalMorseIndexIntroThm}
  \\
  $\ell(\ft,\fs)$
  & Level of Seiberg--Witten moduli space $M_\fs$ in moduli space of ideal non-Abelian monopoles $\sI\!\!\sM_{\ft}$;
  Equation \eqref{eq:ReducibleLevel}
  \\
  $L^p(E)$
  & Banach space of $L^p$ sections of a Hermitian vector bundle $E$;
  Equation \eqref{eq:Aubin_definition_2-3}
  \\
  $M_\fs$
  & Moduli space of Seiberg--Witten monopoles for a \spinc structure $\fs$ over a Riemannian four-manifold;
  Equation \eqref{eq:Moduli_space_Seiberg-Witten_monopoles}
  \\
  $M_\kappa^w(X,g)$
  & Moduli space of projective $g$-anti-self-dual connections over a Riemannian four-manifold $(X,g)$;
  Equation \eqref{eq:ASDModuliSpace}
  \\
  $\sM_\ft$
  & Moduli space of non-Abelian monopoles for a \spinu structure $\ft$ over a Riemannian four-manifold;
  Equation \eqref{eq:Moduli_space_SO(3)_monopoles}
  \\
  $\sM_\ft^*$, $\sM_\ft^0$, $\sM_\ft^{*,0}$
  & Subspaces of moduli space of non-Abelian monopoles for a \spinu structure $\ft$ over a Riemannian four-manifold; Equation \eqref{eq:PU2MonopoleSubspaces}
  \\
  $\sM_{\ft,\reg}$
  & Subspace of regular points of $\sM_\ft$;
  page \pageref{page:RegularPointsOfNonAbelianMonopoleModuli}  
  \\
  $\sM_{\ft,1}$, $\sM_{\ft,2}$
  & Moduli spaces of type $1$ and type $2$ non-Abelian monopoles for a \spinu structure $\ft$ over an almost Hermitian four-manifold; Remarks \ref{rmk:Projective_vortices_type_1_monopoles} and \ref{rmk:Projective_vortices_type_2_monopoles}
  \\
  $\bar\sM_{\ft}$
  & Uhlenbeck compactification of moduli space of non-Abelian monopoles $\sM_\ft$; Prior to equation \eqref{eq:idealmonopoles}
  \\
  $\cM(E)$ & Moduli space of holomorphic structures on a complex vector bundle $E$ over a complex manifold; Equation \eqref{eq:Moduli_space_holomorphic_structures}
  \\
  $\cM(E,\omega)$
  & Moduli space of stable holomorphic structures on a complex vector bundle $E$ over a complex, K\"ahler manifold $(X,\omega)$; Equation \eqref{eq:Moduli_space_omega-stable_holomorphic_structures}
  \\
 $\cM_\ps(E,\omega)$
 & Moduli space of polystable holomorphic structures on a complex vector bundle $E$ over a complex, K\"ahler manifold $(X,\omega)$; Theorem \ref{thm:Kobayashi_7_4_20}
  \\
  $M(E,h,\omega)$
  & Moduli space of projectively Hermitian--Einstein connections on a Hermitian vector bundle $(E,h)$ over an almost Hermitian manifold $(X,\omega)$; Equation \eqref{eq:Moduli_space_HE_connections}
  \\
  $\fM(E)$
  & Moduli space of holomorphic pairs on a complex vector bundle $E$ over a complex manifold;
  Equation \eqref{eq:Moduli_space_holomorphic_pairs}
  \\
  $\fM^{**}(E)$
  & Moduli space of strongly simple holomorphic pairs on a complex vector bundle $E$ over a complex manifold;
  Equation \eqref{eq:Moduli_space_strongly_simple_holomorphic_pairs}
  \\
  $\fM^0(E)$
  & Moduli space of non-zero-section holomorphic pairs on a complex vector bundle $E$ over a complex manifold;
  Equation \eqref{eq:Moduli_space_holomorphic_pairs_NonZeroSection}
  \\
  $\fM_\reg(E)$
  & Moduli space of regular holomorphic pairs on a complex vector bundle $E$ over a complex manifold;
    Equation \eqref{eq:Moduli_space_regular_holomorphic_pairs}
  \\
  $\fM(E,\omega)$
  & Moduli space of stable holomorphic pairs on a complex vector bundle $E$ over a complex K\"ahler manifold $(X,\omega)$; Equation \eqref{eq:Moduli_space_stable_holomorphic_pairs}
  \\
  $\fM_\ps(E,\omega)$
  & Moduli space of polystable holomorphic pairs on a complex vector bundle $E$ over a complex K\"ahler manifold  $(X,\omega)$; Equation \eqref{eq:Moduli_space_polystable_holomorphic_pairs}
  \\
  $\sM(E,h,\omega)$
  & Moduli space of projective vortices on a Hermitian vector bundle $(E,h)$ over an almost Hermitian manifold $(X,\omega)$; Equation \eqref{eq:Moduli_space_projective_vortices}
  \\
  $\sM^*(E,h,\omega)$
  & Moduli space of non-split projective vortices on a Hermitian vector bundle $(E,h)$ over an almost Hermitian manifold $(X,\omega)$; Equation \eqref{eq:Moduli_space_projective_vortices_non-split}
  \\
  $\sM^0(E,h,\omega)$
  & Moduli space of non-zero-section projective vortices on a Hermitian vector bundle $(E,h)$ over an almost Hermitian manifold $(X,\omega)$; Equation \eqref{eq:Moduli_space_projective_vortices_non-zero_section}
  \end{tabular}

\begin{tabular}{ r p{5in} }
  $\sM^{*,0}(E,h,\omega)$
  & Moduli space of non-split, non-zero-section projective vortices on a Hermitian vector bundle $(E,h)$ over an almost Hermitian manifold $(X,\omega)$; Equation \eqref{eq:Moduli_space_projective_vortices_non-split-non-zero-section}
  \\
  $\sM^{**}(E,h,\omega)$
  & Moduli space of projective vortices with trivial stabilizer on a Hermitian vector bundle $(E,h)$ over an almost Hermitian manifold $(X,\omega)$; Equation \eqref{eq:Moduli_space_projective_vortices_StabAvarphi_idE}
  \\
  $\sM_\reg(E,h,\omega)$
  & Moduli space of regular projective vortices on a Hermitian vector bundle $(E,h)$ over a complex Hermitian manifold $(X,\omega)$; Equation \eqref{eq:Moduli_space_projective_vortices_regular}
  \\
  $\sM_\reg^0(E,h,\omega)$
  & Moduli space of regular, non-zero-section projective vortices on a Hermitian vector bundle $(E,h)$ over a complex Hermitian manifold $(X,\omega)$; Equation \eqref{eq:Moduli_space_projective_vortices_non-zero_section_regular}
  \\
  $\sM_\reg^{**}(E,h,\omega)$
  & Moduli space of regular, trivial stabilizer projective vortices on a Hermitian vector bundle $(E,h)$ over a complex Hermitian manifold $(X,\omega)$; Equation \eqref{eq:Moduli_space_projective_vortices_StabAvarphi_idE_regular}
  \\
  $\bM_E^0$
  & Global moduli functor; Definition \ref{defn:FrM_Definition_IV.4.1_pairs}
  \\
  $\bM_{\bar\partial_E,\varphi}$
  & Local moduli functor; Definition \ref{defn:Friedman_Morgan_4-2-4}
  \\
  $\mu(E)$
  & Slope of vector bundle $E$;
  Equation \eqref{eq:Slope}
  \\
  $\mu_M(\sE)$;
  & Bounded slope of a sheaf $\sE$;
  Equation \eqref{eq:Bradlow_mu_M}
  \\
  $\mu_\varphi(\sE)$
  & Slope of a pair $(\sE,\varphi)$;
  Equation \eqref{eq:Bradlow_mu_m}
  \\
  $\bmu(A)$
  & Moment map for the action of $W^{2,p}(\SU(E))$ on the affine space of unitary connections $\sA(E,h)$; Equation \eqref{eq:Moment_map_action_unitary_det_one_gauge_transformations_affine_space_unitary_connections}
   \\
  $\bmu(A,\varphi)$
  & Moment map for the action of $W^{2,p}(\SU(E))$ on the affine space of unitary pairs $\sA(E,h\times W^{1,p}(E))$; Equation \eqref{eq:Moment_map_action_unitary_det_one_gauge_transformations_affine_space_pairs}  
  \\
  $\bmu^*(A,\varphi)$
  & Co-moment map for the action of $W^{2,p}(\SU(E))$ on the affine space of unitary pairs $\sA(E,h\times W^{1,p}(E))$; Equation \eqref{eq:Co-moment_map_action_unitary_det_one_gauge_transformations_affine_space_pairs}
  \\
  $n_a(\ft)$
  & Expected codimension of moduli space of anti-self-dual connections in $\sM_\ft$;
  Equation \eqref{eq:Transv}
  \\
  $\Lambda^p(E)$, $\Lambda^{p,q}(E)$
  & Vector spaces of smooth $p$-forms or $(p,q)$-forms taking values in
    a smooth vector bundle $E$ over a smooth or almost complex manifold;
    Equation \eqref{eq:VectorBundleValued(pq)forms}
  \\
  $\Omega^p(E)$, $\Omega^{p,q}(E)$
  & Fr\'echet spaces of smooth sections of $\Lambda^p(E)$ or $\Lambda^{p,q}(E)$;
    Page \pageref{eq:VectorBundleValued(pq)forms}
  \\
  $\bomega(a_1,a_2)$
  & Weak ($L^2$) symplectic form on affine space of unitary connections $\sA(E,h)$;
  Equation \eqref{eq:Kobayashi_7-6-22}
  \\
  $\bomega((a_1,\phi_1),(a_2,\phi_2))$
  & Weak ($L^2$) symplectic form on affine space of unitary pairs $\sA(E,h)\times W^{1,p}(E)$;
    Equation \eqref{eq:Kobayashi_7-6-22_pairs}
  \\
  $\pi_h^{0,1}$
  & Canonical isomorphisms of Banach affine space of unitary connections or pairs; Equations \eqref{eq:Bijection_unitaryconnections_with_01connections} or \eqref{eq:Bijection_unitarypairs_with_01pairs}
  \\
  $p_1(V)$
  & First Pontrjagin class of a real vector bundle $V$;
  \\
  $p_1(\ft)$
  & First Pontrjagin class of  spin${}^u$ structure $\ft$;
  Equation \eqref{eq:SpinUCharacteristics}
  \\
  $Q_X$
  & Intersection form of four-manifold $X$;
  page \pageref{page:IntersectionForm}
  \\
  $(\rho,W)$
  & Spin${}^c$ structure;
  page \pageref{page:Spinc_structure}
  \\
  $(\rho,W\otimes E)$
  & Spin${}^u$ structure;
  page \pageref{page:Spinu_structure}
  \\
  $\rho_2$, $\rho_Z$, $\rho_{\SU}$
  & $S^1$ actions associated to an orthogonal decomposition $E=L_1\oplus L_2$ of a Hermitian vector bundle;
  Equation \eqref{eq:DefineUnitaryS1ZActionsAtReducible}
  \\
  $\rho_2^\sA$
  & $S^1$ action on affine space of spin${}^u$ pairs $\sA(E,h)\times W^{1,p}(E,h)$;
  Equation \eqref{eq:UnitaryS12ActionOnNonAbelianPairsAffineSpace}
  \\
  $\rho_2^\CC$, $\rho_{\SL}$, $\rho_Z^\CC$
  & $\CC^*$ actions defined by a direct-sum decomposition $E=E_1\oplus E_2$ of a complex vector bundle;
    Equation \eqref{eq:DefineS1ActionsAtReducible}
  \\
  $\rho_2^{\sA,\CC}$
  & $\CC^*$ action defined by a direct-sum decomposition $E=E_1\oplus E_2$ of a complex vector bundle;
    Equation \eqref{eq:S12ActionOnAffine}
  \\
  $\sS$
  & Non-Abelian monopole or projective vortex map; Equation \eqref{eq:PerturbedSO3MonopoleEquation_map} or \eqref{eq:Projective_vortex_map}
  \\
  $\fS$
  & Holomorphic pair map; Equation \eqref{eq:Holomorphic_pair_map}
  \\
  $S_{A,\varphi}$
  & Coulomb-gauge local affine slice through a unitary pair $(A,\varphi)$; Equation \eqref{eq:Unitary_pair_Avarphi_linear_slice}
  \\
  $S_{\bar\partial_E,\varphi}$
  & Coulomb-gauge local affine slice through an $(0,1)$-pair $(\bar\partial_E,\varphi)$; Equation \eqref{eq:dbar_Evarphi_slice}
  \\
  $\SW_X(\fs)$
  & Seiberg--Witten invariant associated to a spin${}^c$ structure $\fs$;
  page \pageref{page:SW_Invariant}
  \\
  $\|s\|_{W_\nabla^{k,p}(M)}$
  & Sobolev norm of a section $s$ of a vector bundle;
  Equation \eqref{eq:Aubin_definition_2-3}
  \\
  $\SL(E)$
  & Bundle of special linear automorphisms of a complex vector bundle $E$;
  Equation \eqref{eq:AutomorphismBundles}
  \end{tabular}

\begin{tabular}{ r p{5in} }
  $\Stab(A)$
  & Stabilizer in $W^{2,p}(\SU(E))$ of a unitary connection $A$;
  Equation \eqref{eq:DefineStabilizerOfConnection_in_SU(E)}
  \\
  $\Stab(A,\Phi)$
  & Stabilizer in $W^{2,p}(\SU(E))$ of a spin${}^u$ pair $(A,\Phi)$;
  Equation \eqref{eq:Define_Stabilizer_of_Spinu_Pair_in_SU(E)}
  \\
  $\Stab(A,\varphi)$
  & Stabilizer in $W^{2,p}(\SU(E))$ of unitary pair $(A,\varphi)$;
  Equation \eqref{eq:DefineStabilizerOfUnitaryPair}
  \\
  $\Stab(\bar\partial_E)$
  & Stabilizer in $W^{2,p}(\SL(E))$ of an $(0,1)$-connection $\bar\partial_E$:
  Definition \ref{defn:StabilizerOf(0,1)Conn}
  \\
  $\Stab(\bar\partial_E,\varphi)$
  & Stabilizer in $W^{2,p}(\SL(E))$ of an $(0,1)$-pair $(\bar\partial_E,\varphi)$;
    Equation \eqref{eq:Stabilizer_holomorphic_pair}
  \\
  $\SU(E)$
  & Bundle of special unitary automorphisms of a Hermitian vector bundle $E$;
   Equation \eqref{eq:AutomorphismBundles_SU(E)}
   \\
  $\fs$
  & Spin${}^c$ structure;
  page \pageref{page:Spinc(X)}
  \\
  $\fs_{\can}$
  & Canonical spin${}^c$ structure of an almost complex manifold;
  Definition \ref{defn:Canonical_spinc_bundles}
  \\
  $\Spinc(X)$
  & Spin${}^c$ structures on $X$;
  page \pageref{page:Spinc(X)}
  \\
  $\su(E)$
  & Bundle of skew-Hermitian, trace-free endomorphisms of a Hermitian vector bundle $E$;
  page \pageref{page:su(V)}
  \\
  $\ft$
  & Spin${}^u$ structure;
  page \pageref{page:Spinu_structure}
    \\
  $\U(E)$
  & Bundle of unitary automorphisms of a Hermitian vector bundle $E$;
   Equation \eqref{eq:AutomorphismBundles_U(E)}
  \\
  $\fu(E)$
  & Bundle of skew-Hermitian endomorphisms of a Hermitian vector bundle $E$;
  page \pageref{page:u(V)}
  \\
  $w_2(\ft)$
  & Second Stiefel--Whitney class of  spin${}^u$ structure $\ft$;
  Equation \eqref{eq:SpinUCharacteristics}
  \\
  $\widetilde X$
  & Smooth blow-up of a four-manifold $X$;
  page \pageref{page:BlowUp}
  \\
  $W^{k,p}(E)$
  & Space of Sobolev sections of vector bundle $E$;
  Equation \eqref{eq:Aubin_definition_2-3}
\end{tabular}

%
%

\begin{theindex}
  \item $(0,1)$-connection, \pageref{eq:Donaldson_Kronheimer_2-1-48}
  \item $(0,1)$-pair, \pageref{page:(0,1)-pair}
  \item Adjoint differential operator, $L^2$ or formal \pageref{page:Adjoint_Differential_Operator}  
  \item Almost complex manifold, \pageref{page:Almost_complex_manifold}
  \item Almost complex structure, \pageref{page:Almost_complex_structure}
  \item Almost Hermitian manifold, \pageref{page:Almost_Hermitian_manifold}
  \item Almost K\"ahler manifold, \pageref{page:Almost_Kaehler_manifold}
  \item Almost symplectic manifold, \pageref{page:almost_symplectic_manifold}
  \item Analytic space
    \subitem complex analytic space, \pageref{page:Complex_analytic_space}
    \subitem real analytic space,\pageref{page:Real_analytic_space}
  \item Basic class, \pageref{page:Basic_class}
  \item Bogomolov--Miyaoka--Yau inequality, \pageref{eq:BMY}
  \item Bogomolov--Miyaoka--Yau conjecture for four-manifolds of Seiberg--Witten simple type,
    \pageref{conj:BMY_Seiberg-Witten}
  \item Canonical line bundle, \pageref{eq:DefineCanonicalLineBundle}
  \item Canonical spin${}^c$ structure, \pageref{defn:Canonical_spinc_bundles}
  \item Characteristic element in $H^2(X;\ZZ)$ of a four-manifold $X$, \pageref{page:characteristic}  
  \item Chern connection, \pageref{rmk:Chern_connection}
  \item Circle-invariant tensor, \pageref{page:circle_invariant}
  \item Clifford modules, \pageref{sec:SpincuStr}
  \item Clifford multiplication, \pageref{eq:CliffordMapDefn}
  \item Coadjoint action, \pageref{item:Marsden-Weinstein_symplectic_quotient_conditions_equivariance}
  \item Co-moment map, \pageref{page:co-moment_map} 
  \item Compatible triples of almost complex structures, metrics, and non-degenerate two-forms, \pageref{page:Compatible_triple}
  \item Configuration space,
    \subitem see Quotient space
  \item Critical point on complex analytic space, \pageref{page:Critical_point}
  \item $\om$-degree,
    \subitem see $\om$-slope
  \item Elliptic complex, \pageref{page:Elliptic_complex},
    \subitem Cohomology of an elliptic complex, \pageref{eq:Cohomology_groups_elliptic_complex}
    \subitem Euler characteristic of an elliptic complex, \pageref{eq:Index_elliptic_complex}
    \subitem 
    Harmonic spaces of an elliptic complex, \pageref{eq:Harmonic_spaces_elliptic_complex}
    \subitem Holomorphic pairs complex, \pageref{eq:Holomorphic_pair_elliptic_complex}
        \subsubitem Negative-weight subcomplex, \pageref{eq:NegativeWeightPairDefComplex}
        \subsubitem Positive-weight subcomplex, \pageref{eq:PositiveWeightPairDefComplex}
        \subsubitem Tangential subcomplex, \pageref{eq:ZeroWeightPairDefComplex}
        \subsubitem Zero-weight subcomplex --- see Tangential subcomplex
    \subitem Non-Abelian monopole complex, \pageref{eq:SO3MonopoleDefComplex}
    \subitem Pre-holomorphic pair complex, \pageref{eq:Pre-holomorphic_pair_elliptic_complex}
    \subitem Projective vortex complex, \pageref{eq:Projective_vortex_elliptic_deformation_complex}
  \item Equivariant Kuranishi model, \pageref{lem:EquivariantKuranishiLemma}
    \subitem Equivariant embedding and obstruction maps, \pageref{page:Embedding_obstruction_Kuranishi}  
  \item Equivariant complex, \pageref{defn:GEquivariantStructure}
  \item Frankel's Theorem, \pageref{mainthm:Frankel_almost_Hermitian}
  \item Fredholm complex, \pageref{sec:Fredholm_complexes_Euler_characteristics}
  \item Gauge transformations,
    \subitem on affine space of $(0,1)$-connections, \pageref{eq:SL(E)_Action_on_(0,1)conn}
    \subitem on affine space of $(0,1)$-pairs, \pageref{eq:SL(E)ActionOn(0,1)Pairs}
    \subitem on affine space of \spinc pairs, \pageref{eq:SWGaugeGroupAction}
    \subitem on affine space of \spinu pairs, \pageref{eq:GaugeActionOnSpinuPairs}
    \subitem on affine space of unitary connections, \pageref{eq:W2pSUE_action_on_AEh}
    \subitem on affine space of unitary pairs, \pageref{eq:DefineSU(E)ActionOnUnitaryPairs}
  \item $G$-equivariant structure on a complex,
    \subitem see Equivariant complex
  \item Germ at a point,
    \subitem Complex analytic set \pageref{page:Germ_complex_analytic_set}
    \subitem Complex analytic space \pageref{page:Germ_complex_analytic_space}
    \subitem Set \pageref{defn:Germ_set}
  \item Hamiltonian circle action, \pageref{mainthm:Frankel_almost_Hermitian}
  \item Hamiltonian vector field, \pageref{mainthm:Frankel_almost_Hermitian}
  \item Hermitian Clifford module,
    \subitem see Clifford module
  \item Hermitian--Einstein connection, \pageref{defn:HE_connection}
    \subitem projectively Hermitian--Einstein connection, \pageref{eq:Holomorphic_connection_fixed_unitary_determinant_connection}
  \item Hessian form, \pageref{eq:Hessian_bilinear_form}
  \item Hessian operator, \pageref{eq:Hessian_operator}
  \item Hirzebruch--Riemann--Roch Theorem, \pageref{thm:HRR}
  \item Hitchin Hamiltonian function, \pageref{page:Hitchin_Function}
  \item Hitchin--Kobayashi correspondence,
    \subitem Hermitian--Einstein connections and stable holomorphic bundles, \pageref{sec:HK_correspondence_between_HE_connections_and_semistable_bundles}
    \subitem Projective vortices and stable holomorphic pairs, \pageref{sec:HK_correspondence_between_SO3_monopoles_and_semistable_pairs}
  \item Hodge star, \pageref{page:Hodge_star}
  \item Holomorphic pair, \pageref{eq:Holomorphic_pair}
  \item Holomorphic structure, \pageref{eq:dbarE_squared_is_zero}
  \item Kuranishi model,
    \subitem for a Hermitian--Einstein connection, \pageref{page:Kuranishi_model_K(A)}
    \subitem for a projective vortex, \pageref{defn:Friedman_Morgan_4-1-15_projective_vortex}
    \subitem for a holomorphic structure, \pageref{page:Kuranishi_model_cKdbarE}
    \subitem for a holomorphic pair, \pageref{defn:Friedman_Morgan_4-1-15_holomorphic_pairs}  
  \item Kuranishi's Theorem on existence of local moduli functor for holomorphic pairs, \pageref{thm:Friedman_Morgan_4-2-5_holomorphic_pairs}
  \item Lefshetz operator, \pageref{page:Lefschetz_operator}
  \item Local slice,
    \subitem see Slice
  \item Marsden--Weinstein Symplectic Reduction Theorem, \pageref{thm:Marsden-Weinstein_symplectic_quotient}
  \item Moduli functor,
    \subitem Global, \pageref{page:Global_moduli_functor}
    \subitem Local, \pageref{page:Local_moduli_functor}
  \item Moduli space, 
    \subitem Anti-self-dual connections, \pageref{eq:ASDModuliSpace}
    \subitem Holomorphic pairs, \pageref{eq:Moduli_space_holomorphic_pairs}
    \subitem Holomorphic structures, \pageref{eq:Moduli_space_holomorphic_structures}
    \subitem Non-Abelian monopoles, \pageref{eq:Moduli_space_SO(3)_monopoles}
    \subitem Projectively Hermitian--Einstein connections, \pageref{eq:Moduli_space_HE_connections}
    \subitem Projective vortices, \pageref{eq:Moduli_space_projective_vortices}
    \subitem Seiberg--Witten monopoles, \pageref{eq:Moduli_space_Seiberg-Witten_monopoles}
    \subitem Stable holomorphic bundles, \pageref{eq:Moduli_space_omega-stable_holomorphic_structures}
    \subitem Stable holomorphic pairs, \pageref{eq:Moduli_space_stable_holomorphic_pairs}
  \item Moment map, \pageref{page:Moment_map}
  \item Morse--Bott index, co-index, and nullity, \pageref{page:Morse--Bott_index}
  \item Morse--Bott function, \pageref{page:Morse-Bott_function}
  \item Nijenhuis tensor, \pageref{subsec:Nijenhuis_tensor_components_exterior_derivative}
  \item Non-Abelian monopole equations,
    \subitem Perturbed, \pageref{eq:PerturbedSO3MonopoleEquations}
    \subitem Unperturbed, \pageref{eq:SO(3)_monopole_equations}
    \subitem Type $1$, \pageref{page:Type1_and_Type2_NonAbelianMonopole}
    \subitem Type $2$, \pageref{page:Type1_and_Type2_NonAbelianMonopole}
  \item Partial connection,
    \subitem see $(0,1)$-connection
  \item Polystable,
    \subitem Holomorphic bundle, \pageref{defn:Stable_holomorphic_structure}
    \subitem Holomorphic pair, \pageref{defn:StablePair_arbitrary_rank}  
  \item Pre-holomorphic pair equations, \pageref{rmk:Pre-holomorphic_pairs}
  \item Projective vortex, \pageref{eq:SO(3)_monopole_equations_almost_Hermitian_alpha}
  \item Quotient space,
    \subitem of $(0,1)$-connections, \pageref{page:Quotient_space_(0,1)_connections}
    \subitem of $(0,1)$-pairs, \pageref{eq:(0,1)PairQuotientSpaces}
    \subitem of \spinc pairs, \pageref{eq:SpincConfig}
    \subitem of \spinu pairs, \pageref{eq:SpinUConfiguration}
    \subitem of unitary connections, \pageref{page:Quotient_space_unitary_connections}
    \subitem of unitary pairs, \pageref{eq:ConfigurationSpaceForProjectiveVortices}  
  \item Rank of spin${}^u$ structure, \pageref{page:Rank_Of_Spinu}
  \item Split,
    \subitem $(0,1)$-pair, \pageref{defn:Split_(0,1)-pair}
    \subitem Non-Abelian monopole, \pageref{defn:Split_trivial_central-stabilizer_spinor_pair}
    \subitem Unitary connection, \pageref{defn:Reducible_split_trivial-stabilizer_unitary_connection}
    \subitem Unitary pair, \pageref{defn:Split_trivial_central-stabilizer_unitary_pair}
  \item Representable functor, \pageref{defn:Friedman_Morgan_4-2-3}
  \item Seiberg--Witten monopole, \pageref{eq:SeibergWitten}
  \item Seiberg--Witten simple type, \pageref{page:Seiberg_Witten_simple_type}
  \item Semiconnection,
    \subitem see $(0,1)$-connection
  \item Semistable,
    \subitem Holomorphic bundle, \pageref{defn:Stable_holomorphic_structure}
    \subitem Holomorphic pair, \pageref{defn:StablePair_arbitrary_rank}  
  \item Simple,
    \subitem $(0,1)$-connection, \pageref{defn:Simple_01-connection}
    \subitem $(0,1)$-pair, \pageref{defn:Simple_pair}
  \item Slice,
    \subitem for action of Lie group on smooth manifold, \pageref{defn:Duistermaat_Kolk_2-3-1}
    \subitem for action of gauge transformations on unitary pairs, \pageref{eq:Unitary_pair_Avarphi_linear_slice}
    \subitem for action of gauge transformations on $(0,1)$-pairs, \pageref{eq:dbar_Evarphi_slice}
    \subitem local slice for action of Lie group on smooth manifold, \pageref{defn:Duistermaat_Kolk_2-3-1}
  \item $\om$-slope, \pageref{eq:Slope}
  \item Sobolev norm for section of a Riemannian vector bundle over a Riemannian manifold, \pageref{eq:Aubin_definition_2-3}  
  \item Spin connection, \pageref{eq:SpinConnection}
  \item Spin${}^c$ structure, \pageref{page:Spinc_structure}
  \item Spin${}^u$ structure, \pageref{page:Spinu_structure}
  \item Stabilizer,
    \subitem of an $(0,1)$-connection, \pageref{defn:StabilizerOf(0,1)Conn}
    \subitem of an $(0,1)$-pair, \pageref{eq:Stabilizer_holomorphic_pair}
    \subitem of a \spinu pair, \pageref{eq:Define_Stabilizer_of_Spinu_Pair_in_SU(E)}
    \subitem of a unitary connection, \pageref{eq:DefineStabilizerOfConnection_in_SU(E)}
    \subitem of a unitary pair, \pageref{eq:DefineStabilizerOfUnitaryPair}
  \item Stable,
    \subitem Holomorphic bundle, \pageref{defn:Stable_holomorphic_structure}
    \subitem Holomorphic pair, \pageref{defn:StablePair_arbitrary_rank}
  \item Standard action,
    \subitem of $\CC^*$ on affine or quotient space of $(0,1)$-pairs, \pageref{defn:CStarActionOn(0,1)Pairs}
    \subitem of $S^1$ on affine or quotient space of \spinu pairs, \pageref{defn:UnitaryZActionOnAffine}
    \subitem of $S^1$ on affine or quotient space of unitary pairs, \pageref{eq:S1_Action_On_AffineSpaceForProjectiveVortices}
  \item Standard four-manifold, \pageref{defn:Standard}
  \item Strongly simple holomorphic pair, \pageref{item:Strongly_simple_pair}  
  \item Symplectic manifold, \pageref{page:symplectic_manifold}
  \item Uhlenbeck convergence and Uhlenbeck topology, \pageref{defn:UhlenbeckConvergence}
  \item Unitary pair, \pageref{page:Unitary_pairs}
  \item Universal object, \pageref{defn:RepresentableFunctor_UniversalObject}
  \item Unstable,
    \subitem Holomorphic bundle, \pageref{defn:Stable_holomorphic_structure}
    \subitem Holomorphic pair, \pageref{defn:StablePair_arbitrary_rank}    
  \item Virtual Morse--Bott index, \pageref{sec:Virtual_Morse-Bott_index_Hamiltonian_function_circle_action_complex_analytic_space}
  \item Weak symplectic form, \pageref{defn:Symplectic_form_Banach_manifold}  
  \item Witten's dichotomy, \pageref{page:Wittens_Dichotomy}
  \item Weight of a circle action on a smooth manifold, \pageref{page:weight_circle_action}
  \item Zariski tangent space, \pageref{Zariski_Tangent_Space}
  \item Zero-section \spinu pair, \pageref{defn:ZeroSectionSO(3)Pair}
\end{theindex}

\end{document}